\input amstex
\documentstyle{amsppt}


\pagewidth{14 true cm}

\pageheight{24 true cm}

\NoRunningHeads \NoBlackBoxes
\define\a{\alpha}

\define\g{\gamma}

\define\BC{\Bbb C}

\define\M{\Bbb M}
\define\N{\Bbb N}
\define\C{\Bbb C}
\define\CO{\Cal O}
\define\ol{\overline}

\define\la{\lambda}

\define\f{\frac}

\define\th{\theta}

\define\om{\omega}

\define\de{\delta}

\define\ve{\varepsilon}
\define\bs{\bigskip}
\define\ms{\medskip}

\define\si{\sigma}
\define\vv{\vskip 3 true mm}

\topmatter
\title 
Irreducibility algorithm for the Weierstrass polynomials of two
complex variables and the Puiseux expansions: Part[A], Part[B], Part[C]
\endtitle

\subjclass Primary 32S15, 14E15  \endsubjclass


\author Chunghyuk Kang
\bigskip\centerline{To Youngsang Yu(wife), Sehgyung(daughter) and Sehwon(son)}
\endauthor


\address
Department of Mathematics,
Seoul National University, Seoul 151--742, Korea
\endaddress
\address
E-mail address: chkang\@math.snu.ac.kr
\endaddress

\endtopmatter

\document

\baselineskip = 0.4 true cm plus 3pt minus 3pt


 \vfill \pagebreak

 $$\align
&\text{\bf \centerline{$\underline{\text{\bf{Preface}}}$}}\\
&\quad \qquad \qquad \qquad \qquad \\
\endalign$$

It is very fundamental to study irreducible plane curve
singularities in algebraic geometry. The contents of this book
consist of three parts, called Part[A], Part[B] and Part[C] with
good appendix. Our aim is to prove by Part[B] and Part[C] that a
complete irreducibility algorithm for the Weierstrass polynomial of
two complex variables and the corresponding standard Puiseux expansions 
in Part[A] can be explicitly and rigorously computable in an elementary way, as
follows. For brevity, Weierstrass polynomials may be written by
W-polys throughout this book.    \ms

\noindent$\underline{\text{\rm {\bf By The 1st Algorithm of
Part[A],} we can find explicit algorithm for computing a one-}}$

\noindent$\underline{\text{\rm to-one function from Family(1)(the
family of some irreducible W-polys of two complex }}$

\noindent$\underline{\text{\rm  variables of the recursive type)
onto Family(2)(the family of the standard Puiseux expansions),}}$

\noindent$\underline{\text{\rm  as far as the multiplicity sequences of
irreducible plane curve singularities are concerned.}}$

For convenience, each element of Family(1) is called the standard Puiseux
W-poly of two complex variables of the recursive type throughout this book, 
because  {\bf in Part[B]} each element of Family(1) among the irreducible W-polys 
will be shown to satisfy the same kind of property as the standard Puiseux expansion 
among the Puiseux expansion does. As good consequence, {\bf Family(1) gives a new 
classification on the family of the irreducible Weierstrass polynomials of
two complex variables} because it was well-known that Family(2) has already given one classification on the family of the Puiseux expansions. 
The proof for {\bf The 1st Algorithm} will be completely done {\bf in Part[B]}. \ms

\noindent$\underline{\text{\rm {\bf By The 2nd Algorithm of
Part[A]}, we can find a complete and explicit algorithm for}}$

\noindent$\underline{\text{\rm computing irreducibility
criterion of all the W-polys of two complex variables. }}$ To
compute irreducibility criterion for germs of analytic functions
of two complex variables, without loss of generality, {\bf in
Part[A]} it suffices to find {\bf The 2nd Algorithm} for computing
completely irreducible W-polys from all the W-polys of two complex
variables by using the Weierstrass preparation theorem and the
Weierstrass division theorem. 

{\bf By The 2nd Algorithm again}, as soon as any given W-poly of two
complex variables $f\in \BC\{y,z\}$ is found to be
irreducible in $\BC\{y,z\}$, without using The 1st Algorithm we can find
$g\in \text{\rm Family(1)}$ such that $f$ and $g$ have a homeomorphic resolution.
The proof for {\bf The 2nd Algorithm} will be completely done {\bf in
Part[C]}.\ms

\noindent$\underline{\text{\rm {\bf By The 3rd Algorithm of Part[A]}
again, we can find explicit algorithm for computing }}$

\noindent$\underline{\text{\rm the standard Puiseux expansion which
has the same multiplicity sequence as the zero set of  }}$

\noindent$\underline{\text{\rm any given irreducible W-poly of two
complex variables does.}}$ Equivalently, as soon as any given W-poly
$f(y,z)$ of two complex variables is found to be irreducible in
$\BC\{y,z\}$ by {\bf The 2nd Algorithm}, {\bf in Part[A]} we can
find {\bf The 3rd Algorithm} for computing the standard Puiseux
expansion, denoted by $C(t)$, such that $C(t)$ and the zero set
$f(y,z)=0$ at $0\in \BC^2$ have the same multiplicity sequence,
using {\bf The 1st Algorithm}. The proof for {\bf The 3rd Algorithm}
will be completely done {\bf in Part[C]}. \ms

{\bf In Appendix,} as very good applications of Part[A], we can find
an elementary explicit

\noindent$\underline{\text{\rm algorithm for computing a one-to-one
function from Family(2) onto Family(3)(the family of}}$

\noindent$\underline{\text{\rm all the multiplicity sequences
defining irreducible plane curve singularities) and so on. }}$
\bs\bs\bs\bs\bs\bs\bs\bs

$\underline{\text{\rm Key words and phrases}}.$

the standard Puisuex expansions, the divisors under the standard
resolutions, the multiplicity sequences, the standard irreducible
W-polys of two complex variables of the recursive type,  the
Euclidean multiplicity sequence for two positive integers,
quasisingularity, the Weierstrass preparation theorem, and the
Weierstrass division theorem, the division algorithm for the
W-polys.

 \vfill \pagebreak

\centerline{$\underline{\text{\bf{TABLE OF CONTENTS}}}$}

$$\align
\quad \text{\bf \centerline{Introduction  \quad \qquad}} \\
\noindent \text{\bf(For brevity, Weierstrass polynomials may be
written by W-polys.)}\qquad \qquad \qquad  \\
\endalign$$
$$\align
\quad &\text{\bf \centerline{Part[A](Part[A1], Part[A2])}}\\
&\text{\bf How to compute irreducible
W-polys of two complex variables}\\
&\text{\bf  and their standard Puiseux expansions
with no need of proofs}\\
\endalign$$ \ms

\roster {\bf Part[A1] In preparation for the representation of
irreducibility algorithms {\indent}for the W-polys of two complex
variables } \ms

\item"{\S 1.0.}" The definition of the new terminology for three
explicit algorithms for all the irreducible W-polys of two complex
variables and some remarks \dotfill 20  \vv

\item"\S {1.1}." The definition for {\bf the standard irreducible W-polys
of two complex variables of the recursive type}(the standard Puiseux
W-polys of two complex variables of the recursive type) \dotfill
13\vv

\item"\S {1.2}." The new terminology and notations in preparation for
studying three families with equivalence relations \dotfill 14 \vv

\item"\S {1.3}."  What does an equivalence relation of any two elements in
{\bf Family(1)(the family of the standard Puiseux W-polys of two
complex variables of the recursive type)} mean? \dotfill 16 \vv
\endroster \ms

{\bf Part[A2] The rigorous representation of explicit irreducibility algorithms 
for the W-polys of two complex
variables with examples and without proofs  and related
topics in the Puiseux expansions }

\roster
\item"\S 1.4."  {\bf The 1st Algorithm} for computing
a one-to-one function between {\bf Family(1) and Family(2)(the
family of the standard Puiseux expansions)} with its examples
\dotfill 17 \vv

\item"\S {1.5}."  {\bf The division algorithm} for W-polys in
preparation for the computations of The 2nd Algorithm and The 3rd
Algorithm \dotfill 21 \vv

\item"\S {1.6}."  {\bf Two fundamental lemmas} for the representation of the
local defining equations of irreducible plane curve singularities
\dotfill 38 \vv

\item"\S {1.7}." {\bf Irreducibility criterion of W-polys of two
complex variables}(A generalized representation of irreducible
{W}-polys of two complex variables defining irreducible plane curve
singularities) \dotfill 28 \vv

\item"\S {1.8}." {\bf The 2nd Algorithm} for computing completely irreducible
W-polys from all the W-polys of two complex variables. \dotfill
33\vv

\item"\S {1.9}." {\bf The 3rd Algorithm} for computing the corresponding standard
Puiseux expansion from any irreducible W-poly of two complex
variables with respect to the multiplicity sequences \dotfill 47 \vv

\item "\S {1.10}." {\bf Examples for The 2nd Algorithm and
The 3rd Algorithm} \dotfill 49 \vv
\endroster \bs

$$\align
\text{\bf \centerline{Part[B](Part[B1],\dots, Part[B5])}}\\
 \text{\bf Explicit algorithm for computing the correspondence
between the irreducible} \\
\text{\bf W-polys of two complex variables and the Puiseux
expansions with proofs} \quad \quad
\endalign$$

{\bf Part[B1] Foundations  \vv }

\roster
\item"\S 2." New definitions for quasisingularity and a
generalization of one coordinate patch covering of the local
coordinates used in the standard resolution process of irreducible
plane curve singularities \dotfill 65 \vv

\item"\S3." The representation for the local defining equations
of irreducible plane curve singularities which have either the same
multiplicity sequence or the homeomorphic resolution under the
standard resolution as the curve singularity
\text{(\{$z^n+y^k=0$\})} does and its generalizations \dotfill 63\vv

\item"\S 4." The proofs of Theorem 3.6 and Theorem 3.7 in $\S {3}$
\dotfill 71 \vv
\endroster\ms

{\bf Part[B2] How to construct the Puiseux convergent power
series of the recursive type in $\BC\{y,z\}$(irreducible
W-polys of two complex variables of recursive types in
$\BC\{y,z\}$) \vv }

\roster
\item"\S 5." The construction of special types of polynomials
of two complex variables defining all irreducible plane curve
singularities at the origin in $\BC^2$ \dotfill 89 \vv

\item"\S 6." The proofs of Theorem 5.0 with five sublemmas and
corollaries in $\S {5}$ \dotfill 84 \vv

\item"\S 7." For any two Puiseux convergent power series $g_r$ and
$\phi_\rho$ of the recursive types of two complex variables, how to
compute the necessary and sufficient condition for \text{$\phi_\rho
\buildrel \text{{\rm divisor}} \over \sim g_r$} under the standard
resolutions and their classifications \dotfill 97 \vv
\endroster\bs

\indent{\bf {Part[B3] Review on the definition of the Puiseux
pairs and the multiplicity {\indent}sequences \vv}}

\roster
\item"\S 8." The definition of the multiplicity and Puiseux
exponents(the Puiseux pairs) and the review of theorems about an
equivalence of the multiplicity and Puiseux exponents and the
multiplicity sequence for irreducible curves  \dotfill  112 \vv

\item"\S 9." New definition of the join of subsequences of a finite
sequence and its application to the representation of the
multiplicity sequences for irreducible plane curve singularities
\dotfill 120 \vv
\endroster\bs

\indent{\bf Part[B4] In preparation for the proof of The 1st
Algorithm  \vv}

\roster
\item"\S 10." To find the necessary and sufficient condition for any
two Puiseux convergent power series of recursive types of two
complex variables to have the same multiplicity sequence and their
classifications \dotfill 123 \vv
\endroster \bs

{\bf Part[B5] The 1st Algorithm with proofs  \vv}

\roster
\item"\S 11." {\bf The 1st Algorithm} for computing
 a one-to-one function from Family(1) onto Family(2) with proofs
\dotfill 134\vv
\endroster \bs

$$\align
\quad &\text{\bf \centerline{Part[C](Part[C1],\dots, Part[C4])}}\\
&\text{\bf A complete and explicit irreducibility algorithms for the
W-polys of two}\\
&\text{\bf complex variables with proofs and related topics in the
Puiseux expansions}
\endalign$$

{\bf Part[C1] The new terminology and notations in preparation
for finding irreducibility criterion of W-polys in
$\BC\{y,z\}$  \vv}

\roster
\item "\S 12."  The generalized standard irreducible W-polys of the
recursive r-type with Theorem 12.0 \dotfill 156 \vv

\item"\S 13." The proofs of Theorem 12.0 with five sublemmas and
corollaries in $\S {12}$ \dotfill 151\vv

\item"\S 14." How to compute the divisor of the total transform of
irreducible plane curve singularities defined by the quasi-Puiseux
convergent power series of the recursive type of two complex
variables under the standard resolution \dotfill 175 \vv \endroster
\bs

{\bf Part[C2] The division algorithm for the W-polys}

\roster
\item"\S 15."  {\bf The division algorithm for the W-polys} \dotfill 220 \vv
\endroster \bs

{\bf Part[C3] Irreducibility criterion of W-polys of two complex
variables  \vv}

\roster
\item"\S {16}." In preparation for irreducibility criterion of W-polys of two
complex variables \dotfill 224 \vv

\item"\S {17}." {\bf Irreducibility criterion of W-polys of two
complex variables(A generalized representation of irreducible
{W}-polys of two complex variables)} \dotfill 229 \vv

\item"\S 18." {The proofs of Proposition 16.7,
Proposition 16.8 and Theorem 16.6 in $\S {17}$} \dotfill 236\vv
\endroster \bs

{\bf Part[C4] The 2nd Algorithm for computing irreducible W-polys
from all the W-polys of two complex variables and The 3rd
algorithm for computing the corresponding standard Puiseux
expansion from any irreducible W-poly of two complex
variables \vv}

\roster
\item"\S {19}." {\bf The 2nd Algorithm} for computing completely irreducible
W-polys from all the W-polys of two complex variables with proofs
\dotfill 258\vv

\item"\S {20}." {\bf The 3rd Algorithm} for computing the corresponding standard
Puiseux expansion from any irreducible W-poly of two complex
variables with respect to the multiplicity sequences with proofs
\dotfill 258 \vv
\endroster

$$\align
& \centerline{\bf Appendix } \\
& \qquad \qquad\qquad \qquad\text{\bf How to apply the
algorithms  of Part[A] to finding the} \\
&  \qquad \qquad\qquad \qquad
\text{\bf algorithms in Appendix A, Appendix B and Appendix C.} \\
\endalign$$

\roster \qquad \qquad \item"{\bf Appendix A:}" The
$\alpha$-algorithm for finding a one-to-one function from Family(2)
onto Family(3) \dotfill 260 \vv \endroster

\roster \qquad \qquad \item"{\bf Appendix B:}"  The
$\beta$-algorithm for finding a one-to-one function from Family(2)
onto Family(4) \dotfill 274 \vv \endroster

\roster\qquad \qquad \item"{\bf Appendix C:}"  The
$\gamma$-algorithm for finding a one-to-one function from Family(2)
onto Family(5) \dotfill 279 \vv
\endroster \ms

\vfill  \ms

\pagebreak

$$\align
\indent &\text{\bf \centerline{INTRODUCTION(Aug.01)}}\\
&\qquad \qquad \qquad \qquad\qquad \qquad\\
\endalign$$

It is very fundamental to study what problems can be computed in
irreducible plane curve singularities in Algebraic geometry? For
brevity, Weierstrass polynomials may be written by W-polys
throughout this book.

In general, irreducible plane curve singularities are given by
either complex parametrizations or germs of analytic functions of two
complex variables.

Whenever all the irreducible plane curve singularities are given by complex 
parametrizations, it has been already known by
Theorem B(Theorem 8.5([Br],[Bu],[Z1])) and Theorem A(Theorem 8.10([K2])) that 
the irreducible plane curve singularities can be classified
by the standard Puisuex expansions as far as the multiplicity sequences of
irreducible plane curve singularities are concerned.

To find an explicit algorithm for finding the correspondence between the Puiseux
expansions and the W-polys of two complex variables in this book later,
we will use Theorem 8.10([K2]), 
instead of the well-known theorem(Theorem 8.5), because Theorem B is not
applicable to compute such an algorithm:

{\bf Theorem A:} Whenever any two irreducible parametrization have
the same Puiseux pairs(i.e., the same standard Puisuex
expansions) by a nonsingular change of the parametrization, then
they have the same multiplicity sequences, and conversely.

{\bf Theorem B:} As far as arbitrary Puiseux expansion of
irreducible plane curve singularities is concerned, any two
irreducible plane curve singularities have the same topological
types if and only if they have the same Puiseux pairs.

$(*)$ For example, it has been well-known by Theorem
8.3(Enriques-Chisini) that an algorithm for finding a one-to-one
correspondence between Family(2)(the family of the standard Puiseux
expansions) and Family(3)(the family of all the multiplicity
sequences of irreducible plane curves with isolated singularity
under the standard resolution) is computable. But, even if it is
small, {\bf as an application of the Algorithms in this book,} we
can get a complete, explicit and elementary algorithm for computing
a one-to-one function from Family(2) onto Family(3)(See Theorem A.4
and Theorem A.5 of Appendix A).

On the other hand, it has been not known at all that the
classification of the family of germs of analytic functions at
$0\in\C^2$, which are irreducible in $\BC\{y,z\}={}_2\CO$, the ring
of convergent power series at $0\in\BC^2$, with isolated singularity
at the origin, can be computed as far as the multiplicity sequences
of irreducible plane curve singularities are concerned. Also, it has
been not known yet how any equivalence of the irreducible germs of
analytic functions in $\BC\{y,z\}$ with isolated singularity at the
origin and the Puiseux expansions can be computed.

In fact, when I stayed at Purdue University in 1988, S. S. Abhyankar
gave me his long manuscript[Ab3], called "Irreducibility criterion
for germs of analytic functions of two complex variables", which was
his answer to T. C. Kuo's question.  It was very difficult for me to
find any algorithm for computing the irreducibility criterion from
his manuscript. But, for a long time, the theme of the above
manuscript has been made me very interesting.

{\bf To compute irreducibility criterion for germs of analytic
functions of two complex variables, by the Weierstrass preparation
theorem it suffices to find explicit algorithms for computing
irreducibility criterion of all the W-polys of two complex
variables.}\ms

{\bf The following are very interesting problems which can be solvable by three algorithms, 
called The 1st Algorithm, The 2nd Algorithm and The 3rd Algorithm of Part[A], which 
will be found in this book, for the first time.} \ms
$\underline{\text{\bf {Problem(A)}}}$  Whenever any given W-poly
$f(y,z)$ of two complex variables is found to be irreducible in
$\BC\{y,z\}$, then we can compute the standard Puiseux
expansion, denoted by $C(t)$, such that the curve $C(t)$ and the zero set
$f(y,z)=0$ have the same multiplicity sequence? \ms

$\underline{\text{\bf {Problem(B)}}}$  Let the
parametrization for arbitrary irreducible curve $C(t)$ be given by
$$\align 
\text{$y(t)=t^n$ and
$z(t)=c_1t^{k_1}+c_2t^{k_2}+\cdots=c_1t^{k_1}(1+H(t))$,} \tag {$*$}
\endalign$$ 
where $1<n$, $1<k_1<k_2<\cdots, $ and the $c_i$ are nonzero
complex numbers, and either $n\le k_1$ or $n\ge k_1$ and $H(t)$ is just the substitution.
Then, we can compute a W-poly
$f(y,z)$ of two complex variables such that  the curve $C(t)$ and the zero set
$f(y,z)=0$ have the same multiplicity sequence? \ms

The contents of this book consist of three parts, called Part[A],
Part[B] and Part[C] with good appendix. Our aim is to prove by Part[B] and Part[C] that 
a complete irreducibility algorithm for the W-poly of
two complex variables and the corresponding standard Puiseux expansions 
in Part[A] can be explicitly and rigorously computable in an elementary way. \ms

In more detail, Part[A] is divided by Part[A1] and Part[A2];
 Part[B] is divided by Part[B1],\dots, Part[B5];
 Part[C] is divided by Part[C1],\dots, Part[C4]. \ms

Now, in order to solve the above problems completely, computationally and
rigorously, we are going to prove by Part[B] and Part[C]
that the definition for the new teminology appearing in \S 1.1 of Part[A1] 
must be well-defined, and so The 1st Algorithm in \S 1.4 of Part[A2], The 2nd 
Algorithm in \S 1.8 of Part[A2] and The 3rd Algorithm in \S 1.9 of Part[A2]
can be completely and rigorously computable in an elementary, as follows: \ms

\noindent$\underline{\text{\rm {\bf By the definition for the new teminology in \S 1.1 of Part[A1],} {\rm  the standard Puiseux   }}}$
  
\noindent$\underline{\text{\rm {\rm W-poly of   two complex variables of the
 recursive r-type in {\S 1.1} is the new irreducible W-poly   } }}$

\noindent$\underline{\text{\rm  of two  complex variables of ${}_2\CO$ being well-defined by Definition 1.1,  and   Family(1) is called  }}$

\noindent$\underline{\text{\rm the 1st family, consisting of the standard Puiseux W-polys of 
two complex variables of }}$ 

\noindent$\underline{\text{\rm  recursive r-type.}}$

Since we can prove by Definition [B4] of Definition 5.0.0 and Theorem 5.1 that each element 
of Family(1) is an irreducible W-poly in $z$ with coefficients in $\BC\{y\}$, it will be redefined by the standard 
irreducible(Puiseux) W-poly of two complex variables of the recursive r-type, and 
for notation Family(1) is called the 1st family of all the standard irreducible
W-polys of two complex variables of the recursive r-type by Definition 1.1 and Theorem 5.1  throughout this book. \ms

\noindent$\underline{\text{\rm {\bf By the 1st Algorithm in
\S 1.4 of Part[A2],} we can find explicit algorithm for computing }}$

\noindent$\underline{\text{\rm a one-to-one function from Family(1) 
onto Family(2)(the family of the standard Puiseux  }}$

\noindent$\underline{\text{\rm expansions), as far as the multiplicity sequences of
irreducible plane curve singularities are  }}$

\noindent$\underline{\text{\rm  concerned.}}$

The proof for {\bf The 1st Algorithm} will be completely done {\bf in
Part[B5]}. \ms

\noindent$\underline{\text{\rm {\bf By the 2nd Algorithm in
\S 1.8 of Part[A2],} we can find a complete and explicit }}$

\noindent$\underline{\text{\rm algorithm for computing the irreducibility
criterion of all the W-polys of two complex variables. }}$ 

To compute the irreducibility criterion for germs of analytic functions
of two complex variables, without loss of generality, {\bf in
Part[A]} it suffices to find {\bf The 2nd Algorithm} for computing
completely irreducible W-polys from all the W-polys of two complex
variables by using the Weierstrass preparation theorem and the
Weierstrass division theorem. 

By The 2nd Algorithm again, as soon as any given W-poly of two
complex variables $f\in \BC\{y,z\}$ is found to be
irreducible in $\BC\{y,z\}$, {\bf without using The 1st Algorithm,} we can compute
$g\in \text{\rm Family(1)}$ such that $f$ and $g$ have a homeomorphic resolution.

The proof for {\bf The 2nd Algorithm} will be completely done {\bf in
Part[C]}.\ms

\noindent$\underline{\text{\rm {\bf By The 3rd Algorithm in
\S 1.9 of Part[A2],} we can compute explicit algorithm   }}$

\noindent$\underline{\text{\rm  for computing the standard Puiseux expansion which
has the same multiplicity sequence as    }}$

\noindent$\underline{\text{\rm the zero set of any given irreducible W-poly of two
complex variables does.}}$ 

First, as soon as any given W-poly
$f(y,z)$ of two complex variables is found to be irreducible in
$\BC\{y,z\}$ by {\bf The 2nd Algorithm again} {\bf in Part[A],} 
we can find $g\in \text{\rm Family(1)}$ such that $f$ and $g$ have a 
homeomorphic resolution. After then, {\bf by The 1st Algorithm} we can compute 
the standard Puiseux expansion, denoted by $C(t)$, such that $C(t)$ and the zero set
$g(y,z)=0$ at $0\in \BC^2$ have the same multiplicity sequence.

The proof for {\bf The 3rd Algorithm} will be completely done {\bf in Part[C]}. \ms

Before we finish writing the above three algorithms without proofs in Part[A2] of 
Part[A] completely and explicitly, in advance we will give an introduction to 
three algorihms by [I], [II] and [III] without proofs in {\bf Introduction}, 
in preparation for writing three algorithms without proofs in Part[A2] of Part[A].
After then, use the above algorithms and without proofs we will 
find the method how to compute irreducibility criterion of all the W-polys 
$f(y,z)\in \BC\{y,z\}$ of two complex variables, as follows: Note that the proofs 
of these algorithms will be completely finished in Part[B] and in Part[C].

Moreover, we will show how to classify all the irreducible W-polys of two complex
variables with respect to the standard Puiseux expansions, using
{\bf three explicit algorithms in Part[A]}. \bs

\noindent{\bf [I] The 1st Algorithm(The algorithm for finding a
one-to-one function between Family(1) and Family(2)) in Part[A]} \ms

Note that The 1st Algorithm consists of  the first half of The 1st
Algorithm(Theorem 1.4) and the second half of The 1st Algorithm(Theorem 1.6) in Part [A2]. \ms 

\noindent$\underline{\text{\bf [I-1] The first half of The 1st
Algorithm(Theorem 1.4) and Example 1.4.1 }}$

\proclaim{Theorem 1.4(Part[A])(The first half of The 1st Algorithm: an
algorithm for finding a one-to-one function $\phi$ from Family(1)
into Family(2))}

$\underline{\text{\bf {Assumptions}}}$ Let $g_r\in \BC\{y,z\}$ be
the standard Puiseux W-poly of the recursive r-type in $z$ in {\rm
Family(1)} satisfying six conditions with the same notations as in
{\rm Definition 1.1}. \ms

$\underline{\text{\bf Conclusions}}$

\noindent$\underline{\text{\rm {\bf } {\rm By explicit algorithm in
(1.4.1),} we can compute the standard Puiseux expansion $C_r(t)$}}$

\noindent$\underline{\text{\rm for the curve $C$ such that the zero
set $(V(g_r))$ and $C_r(t)$ have the same multiplicity sequence.}}$
\ms

\noindent$\underline{\text{\bf(Algorithm 1.4.1 for Theorem 1.4)}}$
$$\align
 \text{$C_r(t):=$} & \left\{\eqalign{y&=t^n \cr
z&=t^{\alpha_1}+t^{\alpha_2}+\cdots +t^{\alpha_r}, \cr} \right. \tag
1.4.1
 \\
 \text{such that}  \quad (1.4.1.1)& \quad
 n =n_1n_2\cdots n_r \quad \text{and} \quad
 \alpha_1 =\beta_{1,1}n_2\cdots n_r,  \\
\quad (1.4.1.2)& \quad \alpha_j =\alpha_{j-1}+
\widehat{\Delta}_jn_{j+1}n_{j+2}\cdots  n_r,
 \endalign$$
\quad where $\widehat{\Delta}_j
=\Delta_j(\beta_{j,k})^j_{k=1}-n_jn_{j-1}\Delta_{j-1}(\beta_{j-1,k})^{j-1}_{k=1}>0$
for $2\le j\le r$ and $\Delta_1(t)=t$. \quad {$\square$}
\endproclaim \ms

\noindent{\bf Example 1.4.1 for Theorem 1.4:}

Let the polynomial $g_3$ in $\BC[y,z]$ be given as follows:
$$
 \text{$g_1=z^3+y^4$, $g_2=g_1^{5}+y^{18}z^{2}$,
$g_3=g_2^{3}+y^{61}z^{1}$.} \tag 1.4.2
$$

\roster Then, it is clear by either Definition 1.1 or Definition 1.2
that either $g_3$ is the standard Puiseux W-poly of the recursive
3rd type or $g_3\in \text{\rm Family[1]}$. Now, it is easy to
compute by (1.4.1) in Algorithm 1.4.1 for Theorem 1.4 that the
standard Puiseux expansion for $C_3(t)$ such that the zero set
$(V(g_3))$ and $C_3(t)$ have the same multiplicity sequence, is
given by
$$
\text{$C_3(t):=$}  \left\{\eqalign{ y= &t^{45} \cr z=
&t^{60}+t^{66}+t^{71}. \cr } \right. \tag 1.4.3
$$
\endroster

\noindent$\underline{\text{\bf [I-2] The second half of The 1st
Algorithm(Theorem 1.6) and Example 1.6.3 }}$ \ms

To use the the algorithm in Theorem 1.6, we need to prove Sublemma
1.5, whose proof is trivial, because there is a well-known theorem
for the Euclidean algorithm that for any two positive integers $A$
and $B$, we find two integers $\gamma$ and $\delta$ such that
$\gcd(A,B)=\gamma A+\delta B$, noting that the proof of Sublemma 1.5
just follows from a well-known theorem. \ms

\definition{\bf Sublemma 1.5.(Corollary 7.6) for Theorem 1.6}

$\underline{\text{\bf{Assumptions}}}$ Let $A\ge 2$ and $B\ge 2$ be
integers with $\gcd(A,B)=1$. Let $p$ be an integer such that $p>nAB$
for some integer $n\ge 2$.

$\underline{\text{\bf {Conclusions}}}$ We can compute a unique pair
of two integers $s_1$ and $t_1$ such that $p=s_1A+t_1 B$ with $0\le
s_1<B$ and $t_1>A$. \quad $\square$
\enddefinition \ms

\proclaim{Theorem 1.6(Part[A])(Theorem 11.4:Algorithm for finding a unique
element of Family(1) corresponding to any standard Puiseux
expansion of Family(2))}

$\underline{\text{\bf {Assumptions}}}$ Let $C_r(t)$ be the standard Puiseux
expansion of the $r$-type, defined by
$$\align
(1.6.1) \qquad \text{$C_r(t):=$} \left\{\eqalign{ y=&t^n, \cr
z=&t^{\alpha_1}+t^{\alpha_2}+\cdots +t^{\alpha_r}, \cr} \right.
\\
\text{where} \quad
 2\le n <\alpha_1<\alpha_2<\cdots <\alpha_r  & \quad  \text{and} \\
 n >d_1>d_2>\cdots
 >d_r=1  \quad   \text{with} & \quad
\text{$d_i=\gcd(n,\alpha_1,\dots,\alpha_i)$, $1\le i\le r$.}\qquad
\qquad
\endalign$$

$\underline{\text{\bf Conclusions}}$ To compute the standard Puiseux
W-poly $g_r$ of the recursive $r$-type such that the zero set
$V(g_r)$ and $C_r(t)$ have the same multiplicity sequence, is to
find {\bf explicit algorithm(Algorithm 1.6.2)}, using a finite
number $\frac{r(r-1)}{2}$ of {\rm Sublemma 1.5(Corollary 7.6)}, as
soon as the standard Puiseux W-poly $g_r$ of the recursive $r$-type
satisfies the same kind of properties and notations as in {\rm
Definition 1.1}, and also as in {\rm(1.4.1.1)} and {\rm(1.4.1.2)} of
{\rm Algorithm 1.4.1} for {\rm Theorem 1.4}, for notation.
\endproclaim \ms

\noindent{\bf Example 1.6.3 for Theorem 1.6:}

Let the parametrization $C(t)$ for the curve $C$ be given as
follows:
$$
 \text{$y= t^{100}$
and $z= t^{250}+t^{375}+t^{410}+t^{417}$(See page 517 of [Bri-Kn])}
\tag 1.6.2
$$

It is easy to compute that $C(t)$ is the standard Puiseux expansion.
Then, it can be proved by following explicit algorithm{\bf
(Algorithm 1.6.2)} for Theorem 1.6 that we can compute the unique
standard irreducible W-poly $f(y,z)$ in z with coefficients in
$\BC\{y\}$ such that the zero set $f(y,z)=0$ and the parametrization
$C(t)$  have the same multiplicity sequence at the origin, where the
construction of $f(y,z)=g_4(y,z)$ is as follows:

$f=g_4=g_3^{5}+y^{300}zg_1g_2$ where $g_1=g_1(y,z)=z^2+y^5$,
$g_2=g_2(y,z)=g_1^{2}+y^{10}z$ and $g_3=g_3(y,z)=g_2^{5}+y^{58}g_1$
and $g_4=g_4(y,z)$ and $f=f(y,z)$. \bs \bs 

\noindent {\bf Irreduciblity plane curve singularities in terms of

\noindent{ Puiseux expansions and irreducible Weierstrass polynomials }

\vfill \pagebreak

\noindent{\bf [II]*** The 2nd Algorithm(A complete and explicit
algorithm for computing irreducibility criterion of the
W-polys of two complex variables) in Part[A]} \ms

Observe the following:

\noindent{\bf(a)} {\bf The statement of Theorem 0.3(Theorem1.13: Irreducibility 
criterion of the W-polys of two complex variables)} is defined by the necessary 
and sufficient condition for the W-polys of two complex variables to be irreducible
in $\BC\{y,z\}$. \ms

\noindent{\bf(b)} In order to write Theorem 0.3 throughout this book, 
it is most surprising and helpful that it suffices to find the method how to apply 
{\bf Theorem 0.1(The Weierstrass division algorithm for the W-polys) with Theorem 0.2} 
to solve the desird problem, using Theorem 12.0(The generalized standard irreducible 
W-polys of the recursive r-type) and Theorem 1.15(The 2nd algorithm), 
which is the generalization of Lemma 1.12.1(The fundamental algorithm for finding 
irreducibility criterion of any W-poly in $f\in \BC\{y\}[z]$ with $f\in$ the type$[2]$
under the standard resolution and its generalizations). Therefore, such an algorithm 
must be explicitly and rigorously computable in an elementary way. \ms

\noindent{\bf(c)} {\bf The 2nd Algorithm consists of Theorem 1.15 in \S 1.8 of Part[A2]} 
because Theorem 1.15 with Theorem 12.0 gives the method how to compute irreducibility 
criterion of W-polys of two complex variables. \ms   

\noindent{\bf(d)} {\bf As an application of The 2nd Algorithm}, we may classify all the irreducible W-polys with respect to Family(1),  as far as the homeomorphic resolutions 
are concerned. \ms

\noindent$\underline{\text{\rm  
In order to find irreducibility criterion for germs of analytic functions of two complex}}$

\noindent$\underline{\text{\rm variables, by Definition 1.1 in \S 1.1 of Part[A1] and 
Theorem 1.7 in \S 1.5 of Part[A2],  }}$

\noindent$\underline{\text{\rm  it suffices to find The 2nd Algorithm for 
irreducibility criterion of the W-polys of two }}$

\noindent$\underline{\text{\rm complex variables defined by equations in (0.1.1), 
using the  W-division theorem rather  }}$

\noindent$\underline{\text{\rm than a nonsingular change of coordinates of the
point in $\C^2$. }}$ \ms

Throughout this book, in preparation for finding irreducibility criterion of any
W-polys of two complex variables, first of all, it is needed to prove 
the following theorem:

\proclaim{Theorem 0.1(Theorem 1.8: The Weierstrass division algorithm for the
W-polys)}

$\underline{\text{\bf {Assumptions}}}$ \quad Let $f\in \BC\{y\}[z]$
be an arbitrary $W$-poly of degree $n\ge 2$ in $z$. Without loss of
generality, we may assume that $f$ satisfies the following form:
$$\align
 f=z^n+\sum^{n-2}_{i=0} a_iy^{\a_i}z^i \quad \text{with $a_{n-1}$ identically zero}, 
\tag 0.1.1 \endalign$$
where for $0\le i\le n-2$, each $a_i=a_i(y)$ is a unit in
${}_2\CO_0$  if exists, and the $\a_i$ are
positive integers. Assume that 
$f(y,z)$ may not irreducible in $\BC\{y,z\}$. 

In addition, assume that we have the following:

\noindent {\rm(0.1.2)} \qquad \quad   \quad $2\le n\le\a_0$.  \bs

$\underline{\text{\bf {Conclusions}}}$

$\underline{\text{For any sequence of integers, $\{n_k\ge 2: k=1,2,\dots,{\ell}\}$  
with  $n=\Pi^{\ell}_{k=1}n_k$,}}$ we can compute a unique sequence
of $W$-polys in $z$, $\{f_k:k=1,2,\dots,{\ell}\}$
such that $f_k\in \BC\{y\}[z]$ is a $W$-poly in $z$ with coefficients in $\C\{y\}$,
satisfying the fact, denoted by
$\underline{\text{\bf \rm Fact[I]}}$:
Let $f_{-1}=y$, $f_0=z$. \ms

$\underline{\text{\rm Fact[I]}}$

For each $k=1,2,\dots,{\ell}$ and $\ell\le r$, $f_k\in \BC\{y\}[z]$ can be
uniquely written in the form
$$\align
& f_k=f^{n_k}_{k-1}+\sum^{n_k-2}_{i=0} R_{k,i}f^i_{k-1}\in 
\C\{f_{-1},f_0,\dots,f_{k-2}\}[f_{k-1}]  \tag 0.1.3\\
& \text{with \quad $R_{k,{n_k}-1}$ identically zero} \quad \text{and} \quad
\text{with \quad $f=f_{\ell}$,} \\
\endalign$$
satisfying the following:
\roster
\item "(1)" For each $k$,
$f_k\in \C\{y\}[z]$ is a $W$-poly of degree
$\Pi^k_{t=1}n_t$ in $z$ with coefficients in $\C\{y\}$.

\item "{}(1a)" Let $k$ and $i$ be fixed with $1\le k\le
{\ell}$ and $0\le i\le n_k-2$, and if exists, then $R_{1,i}=R_{1,i}(y)$ is a
nonunit in $\C\{y\}$ and for each $k\ge 2$, $R_{k,i}=R_{k,i}(y,z)\in
\C\{y\}[z]$ is a polynomial of degree $<\Pi^{k-1}_{t=1}n_t$ in $z$
with $R_{k,i}(0,z)=0$.\ms

\item "(2)" For each fixed $k=1,2,\dots,{\ell}$, $f_k=f_k(f_{-1},f_0,\dots,
f_{k-1})\in \C\{f_{-1},f_0,\dots,f_{k-2}\}[f_{k-1}]$ of {\rm(0.1.3)} is a W-poly in $f_{k-1}$
with coefficients in $\C\{f_{-1},f_0,\dots,f_{k-2}\}$,
considering $f_{-1},f_0,\dots,f_{k-1}$ as independent complex $({k}+1)$-variables
at $0\in \C^{{k}+1}$ where $f_{-1}=y$ and $f_0=z$ if necessary,
with the following property {\rm (2a)}:

\item"{}(2a)"  Let $k$ and $i$ be fixed with $1\le k\le {\ell}$
and $0\le i\le n_k-2$, and for any nonzero monomial
$\Pi^k_{t=1}f^{\de_t}_{t-2}$ in $R_{k,i}=R_{k,i}(f_{-1},f_0,\dots,
f_{k-2})\in \C\{f_{-1},f_0,\dots,f_{k-2}\}$, $\de_1>0$ and
$0\le \de_t<n_{t-1}$ for $t=1,2,\dots, k$. 
\endroster 
\endproclaim \ms

\definition{Definition 0.2.0(Definition 2.4 and Definition 2.5)} 
Suppose that $f\in\BC\{y,z\}$ has an isolated singular point at $0\in \BC^2$. 
It is said that \text{$f\in$the type $[j]$} or belongs to the type $[j]$ under the
standard resolution if $f$ satisfies the following properties (a)
and (b) after m iterations of blow-ups which is needed only to get
the standard resolution of the singularity of the curve defined by
$f$: Sometimes, we write \text{$V(f)\in$the type $[j]$}, instead of
\text{$f\in$the type $[j]$}.

(a) There are exactly j exceptional curves of the first kind with
$j\le m$, each of which has three distinct intersections with other
exceptional curves and the proper transform.

(b) Each of the remaining (m-j) exceptional curves rather than the
above j exceptional curves in (a) has at most two distinct
intersection points with other exceptional curves and the proper
transform.
\enddefinition
\ms

\definition{Remark 0.2.0.1}

(i) If $f$ has a nonsingular point at the origin, we say that $f\in$
the type $[0]$.

(ii) If $f\in$ the type $[j]$ for an integer $j\ge 1$, then $f$
may not be irreducible by the following example:

 Let $f=(z+y)(z+2y)(z+3y)+y^6$. Then $f\in$ the type $[1]$, but
$f$ is not irreducible in $\BC\{y,z\}$.
\enddefinition \ms

\proclaim{Theorem 0.2(How to apply the division algorithm of the W-polys of two complex 
variables to finding irreducibility criterion of all the W-polys of two complex variables)}

$\underline{\text{\bf {Assumptions}}}$ \quad Let $f\in \BC\{y\}[z]$
be an arbitrary $W$-poly of degree $n\ge 2$ in $z$. Without loss of
generality, we may assume that $f$ satisfies the following form:
$$\align
 f=z^n+\sum^{n-2}_{i=0} a_iy^{\a_i}z^i \quad \text{with $a_{n-1}$ identically zero}, 
\tag 0.2.1 \endalign$$
where for $0\le i\le n-2$, each $a_i=a_i(y)$ is a unit in
${}_2\CO_0$ if exists, and the $\a_i$ are positive integers. 

In addition, assume that we have the following:

\noindent {\rm(0.2.2)} \qquad \quad   \quad $2\le n\le\a_0$.  \bs

$\underline{\text{\bf {Conclusions}}}$\ {\bf The necessary and
sufficient condition for $f(y,z)$ to be irreducible in $\BC\{y,z\}$,
is given as follows: For convenience, we use Definition 0.2.0.}

For some positive integer ${\ell}$ with $1\le {\ell}\le r$, there exists uniquely 
a sequence of integers, 
$\{n_k\ge 2: k=1,2,\dots,{\ell}\}$ with $n=\Pi^{\ell}_{k=1}n_k$ and a sequence 
of $W$-polys in $z$, $\{f_k:k=1,2,\dots,{\ell}\}$
such that each $f_k\in \BC\{y\}[z]$ is a $W$-poly of degree
$\Pi^k_{t=1}n_t$ in $z$ with coefficients in $\C\{y\}$,
satisfying two facts, denoted by
$\underline{\text{\bf Fact[I]}}$ and $\underline{\text{\bf Fact[II]}}$, respectively:
Let $f_{-1}=y$, $f_0=z$. \ms

\noindent$\underline{\text{\bf Fact[I]}}$ {\bf By Theorem 0.1, $f_k$ in this theorem 
satisfies the same results as $f_k$ of (0.1.3) in Theorm 0.1 does.} 
Then, given a sequence of integers, $\{n_k\ge 2: k=1,2,\dots,{\ell}\}$ as above, 
each $f_k\in \BC\{y\}[z]$ can be uniquely written in the form
$$\align
\text{\rm (0.2.3)} \qquad\qquad f_k=f^{n_k}_{k-1}+\sum^{n_k-2}_{i=0} R_{k,i}f^i_{k-1}\in 
\C\{f_{-1},f_0,\dots,f_{k-2}\}[f_{k-1}]
\quad \text{with $f=f_{\ell}$} \qquad \qquad\\
\endalign$$
satisfying the same properties and notations as in Theorem 0.1, as follows:
\roster
\item "(1)" For each $k$,
$f_k\in \C\{y\}[z]$ is a $W$-poly of degree
$\Pi^k_{t=1}n_t$ in $z$ with coefficients in $\C\{y\}$.

\item "{}(1a)" Let $k$ and $i$ be fixed with $1\le k\le
{\ell}$ and $0\le i\le n_k-2$, and if exists, then $R_{1,i}=R_{1,i}(y)$ is a
nonunit in $\C\{y\}$ and for each $k\ge 2$, $R_{k,i}=R_{k,i}(y,z)\in
\C\{y\}[z]$ is a polynomial of degree $<\Pi^{k-1}_{t=1}n_t$ in $z$
with $R_{k,i}(0,z)=0$.

\item "(2)" For each fixed $k=1,2,\dots,{\ell}$, $f_k=f_k(f_{-1},f_0,\dots,
f_{k-1})\in \C\{f_{-1},f_0,\dots,f_{k-2}\}[f_{k-1}]$ of {\rm(0.1.3)} is a W-poly in $f_{k-1}$
with coefficients in $\C\{f_{-1},f_0,\dots,f_{k-2}\}$,
considering $f_{-1},f_0,\dots,f_{k-1}$ as independent complex $({k}+1)$-variables
at $0\in \C^{{k}+1}$ where $f_{-1}=y$ and $f_0=z$ if necessary,
with the property {\rm (2a)}:

\item"{}(2a)"  Let $k$ and $i$ be fixed with $1\le k\le {\ell}$
and $0\le i\le n_k-2$, and for any nonzero monomial
$\Pi^k_{t=1}f^{\de_t}_{t-2}$ in $R_{k,i}=R_{k,i}(f_{-1},f_0,\dots,
f_{k-2})\in \C\{f_{-1},f_0,\dots,f_{k-2}\}$, $\de_1>0$ and
$0\le \de_t<n_{t-1}$ for $t=1,2,\dots, k$. 
\endroster \ms

\noindent$\underline{\text{\bf Fact[II]}}$ 
$$\align
[A]  \quad  \quad &\text{$f$ is irreducible in $\BC\{y,z\}$
with $f\in$the type $[{\ell}]$ under the standard resolution}  \qquad \quad\\
 {\iff} \quad &\text{$f_k$ is irreducible in
$\C\{y,z\}$ with $f_k\in$the type $[k]$ under the standard resolution} \\
& \text{for each $k=1,2,\dots,{\ell}$ with $\ell\le r$}. \\
\endalign$$
\endproclaim \ms

\definition{Remark 0.2.1} For the proof of Theorem 0.2, in order to prove that
the necessary and sufficient condition for $f(y,z)$ of Theorem 0.2 to be irreducible
in $\BC\{y,z\}$ is true, there is nothing to prove by Theorem 0.1 that Fact[I] in 
Conclusions  of this theorem is true.
\enddefinition \ms

\proclaim{Theorem 0.3(Irreducibility criterion of
W-polys of two complex variables(a generalized representation of irreducible
{W}-polys of two complex variables))}

$\underline{\text{\bf {Assumptions}}}$ \quad Let $f\in \BC\{y\}[z]$
be an arbitrary $W$-poly of degree $n\ge 2$ in $z$. Without loss of
generality, we may assume that $f$ satisfies the following form:
$$ f=z^n+\sum^{n-2}_{i=0} a_iy^{\a_i}z^i, \tag 0.3.1 $$
where for $0\le i\le n-2$, each $a_i=a_i(y)$ is a unit in
${}_2\CO_0$ if exists, and the $\a_i$ are
positive integers. Note that $a_{n-1}$ is identically zero.  Write
$n=d_2n_1$ and $\a_0 =d_2\a_{1,0,1}$ with $d_2=\gcd(n,\a_0)$. Write
$n=\Pi^{r}_{k=1}q_k$ with positive integers $q_k\ge 2$ for all $k$
where $\Pi^{r}_{k=1}q_k$ is a factorization of prime numbers $q_k$.

In addition, assume that we have the following:

\noindent {\rm(0.3.2)} \qquad \quad   \quad $2\le n\le\a_0$.  \bs

$\underline{\text{\bf {Conclusions}}}$\ {\bf The necessary and
sufficient condition for $f(y,z)$ to be irreducible in $\BC\{y,z\}$,
is given as follows:}

For some positive integer ${\ell}$ with $1\le {\ell}\le r$, there exists uniquely 
a sequence of integers, 
$\{n_k\ge 2: k=1,2,\dots,{\ell}\}$ with $n=\Pi^{\ell}_{k=1}n_k$ and a sequence 
of $W$-polys in $z$, $\{f_k:k=1,2,\dots,{\ell}\}$
such that each $f_k\in \BC\{y\}[z]$ is a $W$-poly of degree
$\Pi^k_{t=1}n_t$ in $z$ with coefficients in $\C\{y\}$,
satisfying two facts, denoted by
$\underline{\text{\bf Fact[I]}}$ and $\underline{\text{\bf Fact[II]}}$, respectively:
Let $f_{-1}=y$, $f_0=z$. \ms

For the proof of Theorem 0.3, note that there is nothing for the proof of Fact[I] of this theorem.

\noindent$\underline{\text{\bf Fact[I]}}$ {\bf By Theorem 0.1, $f_k$ in this theorem 
satisfies the same results as $f_k$ of (0.1.3) in Theorm 0.1 does.} 
So, given a sequence of integers, $\{n_k\ge 2: k=1,2,\dots,{\ell}\}$ as above, 
each $f_k\in \BC\{y\}[z]$ can be uniquely written in the form
$$\align
\text{\rm (0.3.3)} \qquad\qquad f_k=f^{n_k}_{k-1}+\sum^{n_k-2}_{i=0} R_{k,i}f^i_{k-1}\in 
\C\{f_{-1},f_0,\dots,f_{k-2}\}[f_{k-1}]
\quad \text{with $f=f_{\ell}$} \qquad \qquad\\
\endalign$$
satisfying the same properties and notations as $f_k$ of {\rm(0.2.3)} in 
{\rm Theorm 0.2} does.  
\ms

\noindent$\underline{\text{\bf Fact[II]}}$ By {\rm[A]} of {\rm Fact[II]} in {\rm Theorem 0.2}, 
{\rm Fact[II]} of this theorem is given by {\rm [B]}: 
$$\align
\noindent [B] \quad \quad  &\text{$f$ is irreducible in $\BC\{y,z\}$}  
\qquad \qquad\qquad \qquad \qquad \qquad\qquad \qquad\\
\text{\rm{$\iff$}}\quad &\text{$f$ satisfies a finite number {\rm (${\ell}$)} pairs of inequalities
in {\rm Eq(1)}, \dots, {\bf \rm(${\ell}$)}, which}\\ 
&\text{can be shown by the {\rm (${\ell}$)} steps,  {\rm Step(1) of Fact[II]}, \dots, {\rm Step(${\ell}$) of Fact[II]}, as follows:} \\
\endalign$$

\noindent$\underline{\text{\rm Step(1) of Fact[II]}}$ 

${\underline{\text{$f_1$ is irreducible in
$\C\{y,z\}$ with $f_1\in$the type $[1]$ under the standard resolution}}}$
\ms
\noindent$\iff$ $f_1={f_0}^{n_1}+\sum^{n_1-2}_{i=0} R_{1,i}{f_0}^{i}$
of {\rm(0.3.3)} satisfies a pair of inequalities in \text{\rm Eq(1)} of \text{\rm Step(1) of Fact[II]},  
denoted by \text{\rm Eq(1.1)} and \text{\rm Eq(1.2)} of \text{\rm Eq(1)}, as follows: \ms

For $0\le i\le n_1-2$, each $R_{1,i}=b_iy^{\a_{1,i,1}}$
with a unit $b_i\in\C \{y\}$ and a positive integer $\a_{1,i,1}$ if
exists. Denote $A_{1,i}$ by $b_i(0)$ for convenience of notations.

Define a function $\th_1:\N_0\to\N_0$ by $\th_1(t)=t$ where $\N_0$
is the set of nonnegative integers and $\N^k_0$ is k-dimensional cartesian 
product of $\N_0$. \ms

\noindent$\underline{\text{\rm Eq(1) of Step(1)}}$ For all $i=0,1,\dots,{n_1}-2$, a pair of inequalities in {\rm  Eq(1)} can be given
by 
$$\align
\text{\rm Eq(1.1) of Eq(1)} \qquad \qquad &{\th_1(\a_{1,i,1})}>({n_1}-i),  \\
\text{\rm Eq(1.2) of Eq(1)} \qquad \qquad  &\gcd(n_1,\a_{1,0,1})=1 \quad \text{and} \\
&\f{\th_1(\a_{1,i,1})}{n_1-i}=\f{\a_{1,i,1}}{n_1-i}\ge
\f{\a_{1,0,1}}{n_1}=\f{\th_1(\a_{1,0,1})}{n_1}. \qquad \qquad \qquad  \\
\endalign$$

\noindent$\underline{\text{\rm Step(2) of Fact[II]}}$ 

${\underline{\text{$f_2$ is irreducible in
$\C\{y,z\}$ with $f_2\in$the type $[2]$ under the standard resolution}}}$
\ms
\noindent$\iff$ Assuming that $f_2={f_1}^{n_1}+\sum^{n_1-2}_{i=0} R_{2,i}{f_{1}}^i$  
of {\rm(0.3.3)} has the first pair of inequalities in \text{\rm Eq(1)}, 
which was already given by {\rm Step(1)} of {\rm Fact[II]}, $f_2$ must satisfy the second pair of 
inequalities in {\rm {\rm Eq(2)}} of {\rm Step(2)} of {\rm Fact[II]}. \ms

Then, it remains 
to find the second pair of inequalities in \text{\rm Eq(2)} of \text{\rm Step(2)} of {\rm Fact[II]}, 
denoted by \text{\rm Eq(2.1)} and \text{\rm Eq(2.2)} of \text{\rm Eq(2)}, as follows: \ms

For given integers $n_1,\a_{1,0,1}$ and a function $\th_1$ in {\rm
(1)}, define a function $\th_2:\N^2_0\to \N_0$
by
$$
\th_2(t_1,t_2)=t_2\th_1(\a_{1,0,1})+n_1\th_1(t_1)=t_2\a_{1,0,1}+n_1t_1
\quad \text{for each $(t_1,t_2)\in \N^2_0$}. \tag 0.3.4 $$

Since $f_1$ is irreducible in $\C\{y,z\}$ with $f_1 \in$the type $[1]$ under the
standard resolution, 
it is clear that for any two nonzero monomials $y^{\a_1}z^{\a_2}$ and
$y^{\de_1}z^{\de_2}$ in $R_{2i}$ with $i$ fixed,
$$\align
\text{\rm(0.3.4-1)} \qquad \qquad &
\th_2(\a_1,\a_2)=\th_2(\de_1,\de_2)\ \text{if and only if}\
\a_1=\de_1\ \text{and}\ \a_2=\de_2. \qquad \qquad\qquad \qquad \\
 &\text{So, there exists a unique nonzero monomial
$A_{2,i}y^{\a_{2,i,1}}z^{\a_{2,i,2}}$ in $R_{2,i}$} \\
&\text{with a nonzero constant $A_{2,i}$ such that
$\th_2(\a_{2,i,1},
\a_{2,i,2})=\text{$\min$}\{\th_2(\de_1,\de_2)\}$} \\
&\text{for any nonzero monomial $y^{\de_1}z^{\de_2}$ in $R_{2,i}$
with $i$ fixed.}
\endalign$$

\noindent$\underline{\text{\rm Eq(2) of Step(2)}}$
For all $i=0,1,\dots, n_2-2$, a pair of inequalities in {\rm {\rm Eq(2)}} can be given
by 
$$\align
\text{\rm Eq(2.1) of Eq(2)} \qquad \qquad\qquad    &{\th_2(\a_{2,i,k})^2_{k=1}}>({n_2}-i) n_{1}{\a_{1,0,1}},  \\
\text{\rm Eq(2.2) of Eq(2)}   \qquad \qquad\qquad 
&\gcd(n_2,\th_2(\a_{2,0,k})^2_{k=1})=1 \quad \text{and} \qquad \qquad \qquad \qquad\\
&\f{\th_2(\a_{2,i,1}, \a_{2,i,2})}{n_2-i}\ge
\f{\th_2(\a_{2,0,1},\a_{2,0,2})}{n_2}.\qquad \qquad \qquad \qquad
\endalign$$ \ms

\quad By induction assumption on the integer $(m-1)\le {\ell-1}$ and by {\rm Step(2), Step(3), \dots, Step(m-1)},
there exists a sequence $\{f_2,f_3,\dots,f_{m-1}\}$, each element $f_k${\rm ($2\le k\le m-1$)}
of which has the same kind of {\rm(k)} pairs of inequalities 
in {\rm Eq(1)}, {\rm Eq(2)}, \dots, {\rm Eq(k)}, which was already given by {\rm Step(k)} of {\rm Fact[II]}, as $f_2$ has two pairs of inequalities 
in {\rm Eq(1)} and {\rm Eq(2)} of {\rm Step(2)} of {\rm Fact[II]}. \ms

\noindent$\underline{\text{\bf Step(m) of Fact[II]}}$ 

${\underline{\text{$f_m$ is irreducible in
$\C\{y,z\}$ with $f_m\in$the type $[m]$ under the standard resolution}}}$

\noindent$\iff$ Assuming that 
$f_m={f_{m-1}}^{n_{m-1}}+\sum^{n_{m-1}-2}_{i=0} R_{m,i}{f_{m-1}}^i\in \C\{f_{-1},f_0,\dots,f_{m-2}\}[f_{m-1}]$  has the same kind of {\rm(m-1)} pairs of inequalities 
in {\rm Eq(1)}, {\rm Eq(2)}, \dots, {\rm Eq(m-1)}, which was already given by {\rm Step(m-1)} of {\rm Fact[II]}, $f_m$ must satisfy the $m$-th pair of 
inequalities in {\rm Eq(m)} of {\rm Step(m)} of {\rm Fact[II]}. \ms
Then, it remains 
to find the $m$-th pair of inequalities in \text{\rm Eq(m)} of \text{\rm Step(m)}, 
denoted by \text{\rm Eq(m.1)} and \text{\rm Eq(m.2)}, as follows: 

We may assume by induction assumption 
that $f_{m-1}={f_{m-2}}^{n_{m-2}}+\sum^{n_{m-2}-2}_{i=0} R_{m-1,i}{f_{m-2}}^i$ satisfies a finite number {\rm (m-1)} of a pair of 
inequalities in {\rm Eq(1)}, {\rm Eq(2)}, \dots, {\rm Eq(m-1)}. 

Inductively, define $\th_m:\N^m_0\to \N_0$ by
$$
(0.3.5) \quad
\th_m(t_k)^m_{k=1}=t_m\th_{m-1}(\a_{m-1,0,k})^{m-1}_{k=1}
+n_{m-1}\th_{m-1}(t_k)^{m-1}_{k=1} \quad \text{for each
$(t_k)^m_{k=1}\in \N^k_0$},
$$
where recall by induction assumption that for a fixed $i$,
$A_{m-1,i}\Pi^{m-1}_{k=1}f^{\a_{m-1,i,k}}_{k-2}$ is a unique nonzero
monomial in $R_{m-1,i}$ with a constant $A_{m-1,i}$ such that
$$
\th_{m-1}(\a_{m-1,i,k})^{m-1}_{k=1}=\text{$\min$}\{\th_{m-1}(\de_k)^{m-1}_{k=1}\},
\tag 0.3.5-0
$$
for any nonzero monomial $\Pi^{m-1}_{k=1}f^{\de_k}_{k-2}$ in
$R_{m-1,i}$.

Since $f_{m-1}$ is irreducible in $\C\{y,z\}$ with $f_{m-1} \in$the type $[m-1]$ 
under the standard resolution, it is clear that
for any two nonzero monomials $\Pi^m_{k=1}f^{\a_k}_{k-2}$ and
$\Pi^m_{k=1}f^{\de_k}_{k-2}$ in $R_{m,i}$ with $i$ fixed,
$$\align
\text{\rm(0.3.5-1)} \qquad \qquad
&\text{$\th_m(\a_k)^m_{k=1}=\th_m(\de_k)^m_{k=1}$
if and only if $\a_k=\de_k$ for $k=1,2,\dots, m.$} \qquad \qquad\qquad \qquad\\
&\text{So, there exists a unique nonzero-monomial
$A_{m,i}\Pi^m_{k=1}f^{\a_{m,i,k}}_{k-2}$ in $R_{m,i}$} \\
&\text{with a constant $A_{m,i}$ such that
$\th_m(\a_{m,i,k})^m_{k=1}=\text{$\min$}\{\th_m(\de_k)^m_{k=1}\}$}\\
&\text{for any nonzero monomial $\Pi^m_{k=1}f^{\de_k}_{k-2}$ in
$R_{m,i}$.}\\
\endalign$$

\noindent$\underline{\text{\rm Eq(m) of Step(m)}}$ For all \text{$i=0,1,\dots, {n_m}-2$},  
a pair of inequalities in {\rm Eq(m)} can be given by 
$$\align
\text{\rm Eq(m.1) of Eq(m)} \qquad \qquad &{\th_m(\a_{m,i,k})^m_{k=1}}>({n_m}-i)
 n_{m-1}\th_{m-1}(\a_{m-1,0,k})^{m-1}_{k=1},\\
\text{\rm Eq(m.2) of Eq(m)} \qquad \qquad  
&\gcd(n_m,\th_m(\a_{m,0,k})^m_{k=1})=1 \quad \text{and} \qquad\qquad  \\
&\f{\th_m(\a_{m,i,k})^m_{k=1}}{{n_m}-i}\ge
\f{\th_m(\a_{m,0,k})^m_{k=1}}{n_m}. \\  
\endalign$$
\endproclaim \ms

\definition{Remark 0.3.1}  {\rm(1)} It is clear that the proof of Fact[I] in Theorem 0.3
is true. As far as the proof of Fact[II] in Theorem 0.3 is concerned, note that
the sufficiency of the condition of Fact[II] can be proved by Theorem 12.0, 
and the necessity of the condition of Fact[II] can be proved by Theorem 1.15. \ms

{\rm(2)} Moreover, the conclusion in Theorem 1.15 is explicit
algorithm for finding completely irreducible W-polys from all the
W-polys of two complex variables. \ms

{\rm(3)} The converse of Theorem 16.5 can be represented by Theorem 16.6.
Moreover, we can compute irreducibility criterion of W-polys
defining plane curve singularities at the origin in $\C^2$ in the
process of the proof of Theorem 16.6 together with Proposition 16.7
and Proposition 16.8 completely and rigorously, using the Euclidean
algorithm and Theorem 15.4(The Division Algorithm for the W-polys).
$\square$
\enddefinition \ms

\definition{Remark 0.3.2}  Consider the sequence $S=\{f_k: 1\le k\le{\ell}\}$
with $f_{\ell}=f$ where $f_k=f_k(y,z,\dots,f_{k-1})\in
\C\{y,z,f_1,\dots, f_{k-1}\}$, which have the same properties and
notations as we have seen in {\rm(0.3.3)} of the conclusion of
Theorem $0.3$. If $f\in \BC\{y\}[z]$ is irreducible in $\BC\{y,z\}$,
then $f=f_{\ell}(y,z,f_1,\dots,f_{\ell-1})\in
\C\{y,z,f_1,\dots,f_{\ell-2}\}[f_{\ell-1}]$ is an irreducible
$W$-poly of degree $n_{\ell}$ in $f_{\ell-1}$ with coefficients in
$\BC\{y,z,f_1,\dots,f_{\ell-2}\}$ and with multiplicity $n_{\ell}$
at the origin in $\C^{\ell}$. $\square$ \ms
\enddefinition \ms

\proclaim{Corollary 0.4} $\underline{\text{\bf {Assumptions}}}$
\quad Under the same assumption and conclusion as in Theorem $0.3$,
note that $f_k$ is irreducible in $\BC\{y,z\}$ with isolated
singularity at $(0,0)$ in $\BC^2$ for $k\ge 1$. In particular, for
each $k=1,2,\dots,\ell$, let $V(H_{k})=\{(y,z):H_{k}(y,z)=0\}$ be an
analytic variety at $(0,0)$ in $\BC^2$, each of which is defined as
follows:
$$\align
(0.4.3.1) \qquad \qquad \text{\rm(i)} \qquad
H_1&=z^{n_1}+y^{\alpha_{1,0,1}} \quad \text{with $n_1\ge 2$ and
$\alpha_{1,0,1}\ge 2$}. \qquad \qquad \qquad \qquad \\
\text{\rm(ii)} \qquad H_2&=H^{n_2}_1+y^{\alpha_{2,0,1}}z^{\alpha_{2,0,2}}.\\
\qquad  &\ldots\ldots \\
\text{\rm({$\ell$})} \qquad
H_{\ell}&=H^{n_{\ell}}_{\ell-1}+y^{\alpha_{\ell,0,1}}z^{\alpha_{\ell,0,2}}H^{\alpha_{\ell,0,3}}_1\cdots
H^{\alpha_{\ell,0,\ell}}_{\ell-2}. \\
\endalign$$

$\underline{\text{\bf {Conclusions}}}$ Then, we have the following:

\quad ${f_k} \buildrel \text{{\rm multiseq }}
\over \sim H_{k}$ for each $k=1,2,\dots,\ell$.
\endproclaim \ms

\noindent{\bf Remark 0.4.1. } 
By Theorem $5.0$, we have the following: \ms

$H_{j+1}$ is irreducible in
$\BC\{y,z\}$ with $H_{j+1}\in \text{\rm the type[j+1]}$ under the standard resolution  \ms
$\iff$ the following conditions hold: \ms

\noindent(1) $\gcd(n_1,\alpha_{1,0,1})=1$. \ms

\noindent(2)
$\gcd(n_2,\theta_2(\alpha_{2,0,k})^{2}_{k=1})=1$ with $\dfrac{\th_2(\a_{2,0,1},\a_{2,0,2})}{n_2}>n_1\a_{1,0,1}$.\ms

\qquad \qquad \qquad $\ldots\ldots$

\noindent(j+1) $\gcd(n_{j+1},\theta_{j+1}(\alpha_{{j+1},0,k})^{j+1}_{k=1})=1$ with $\dfrac{\th_{j+1}(\a_{{j+1},0,k})^{j+1}_{k=1}}{n_{j+1}}>n_j\theta_j(\alpha_{j,0,k})^{j}_{k=1}$ . 
$\square$ \ms

\noindent{\bf Remark 0.4.2.} 
By Theorem $12.0$, we have the following: \ms 

$f_{j+1}$ is irreducible in
$\BC\{y,z\}$ with $f_{j+1}\in \text{\rm the type[j+1]}$ under the standard resolution \ms
$\iff$ the following conditions hold:

\noindent(1)  $\gcd(n_1,\alpha_{1,0,1})=1$,
$\dfrac{\th_1(\a_{1,i,1})}{n_1-i}=\dfrac{\a_{1,i,1}}{n_1-i}\ge
\dfrac{\a_{1,0,1}}{n_1}=\dfrac{\th_1(\a_{1,0,1})}{n_1}$ for $0\le
i\le n_1-2$.

\noindent(2) $\gcd(n_2,\theta_2(\alpha_{2,0,1},\alpha_{2,0,2}))=1$,
$\dfrac{\th_2(\a_{2,i,1}, \a_{2,i,2})}{n_2-i}\ge
\dfrac{\th_2(\a_{2,0,1},\a_{2,0,2})}{n_2}>n_1\a_{1,0,1}$ for $0\le i\le n_2-2$.
\ms

\qquad \qquad \qquad $\ldots\ldots$

\noindent(j+1)
$\gcd(n_{j+1},\theta_{j+1}(\alpha_{{j+1},0,k})^{j+1}_{k=1})=1$,
$\dfrac{\th_{j+1}(\a_{{j+1},i,k})^{j+1}_{k=1}}{{n_{j+1}}-i}\ge
\dfrac{\th_{j+1}(\a_{{j+1},0,k})^{j+1}_{k=1}}{n_{j+1}}
>n_j\theta_j(\alpha_{j,0,k})^{j}_{k=1}$ for $0\le
i\le n_{j}-2$. $\square$
\ms

Therefore, to find such irreducibility criterion of all the W-polys of two complex
variables, it suffices to solve the following problems in order:\ms

$\underline{\text{\bf {Problem(I)}}}$ Firstly, it is needed to constuct the new terminology
of Definition 1.0.2 by Definition 5.0.0, called the standard Puiseux W-poly
$g_r(y,z)\in\BC\{y\}[z]$ of the recursive type of two complex variables from the family of all
the W-polys of two complex variables.
After then, it will be proved by Theorem 5.0.1 that the standard Puiseux W-poly $g_r(y,z)\in\BC\{y\}[z]$ of the recursive type is an irreducible W-poly in $z$ with coefficients in $\BC\{y\}$. \ms

$\underline{\text{\bf {Problem(II)}}}$ Secondly, it will be proved by Theorem 7.4 that an equivalence of any two standard Puiseux W-polys of the recursive type of of two complex variables implies homeomorphic resolutions under the same standard resolution, and conversely. \ms

$\underline{\text{\bf {Problem(III)}}}$ Thirdly, it remains to show
that the necessary and sufficient condition for any W-poly of two complex variables to be irreducible in $\BC\{y,z\}$ can be uniquely represented by a generalization of the standard Puiseux
W-poly $g_r(y,z)\in\BC\{y\}[z]$ of the recursive type, using the Weierstass division theorem
and an equivalence of homeomorphic resolutions. \ms

In preparation for finding a solution of the above three problems, in Part [A],
we are going to write some definitions, and lemmas and theorems without proof,
relative to irreducible W-polys of $n\ge 2$ complex variables at $0\in\BC^n$.

\newpage

{\bf As an example for The 2nd Algorithm in Part[A],} let $f(y,z)$
be a W-poly of two complex variables, which is given by the
following:
$$\align
\text{\rm $(*1)$} \quad
f(y,z)=&z^{16}+4y^{3}z^{14}+\{4y^5+6y^6\}z^{12}+\{12y^8+4y^9\}z^{10}+
\{6y^{10}+12y^{11}+y^{12}\}z^{8} \quad \\
&+\{12y^{13}+4y^{14}+y^{17}\}z^{6}+\{4y^{15}+6y^{16}+y^{20}\}z^{4}
+\{4y^{18}+y^{22}\}z^{2}\\
&+y^{24}z+\{y^{20}+y^{29}\}.
\endalign$$

\noindent{\bf Example 1.6.3 for Theorem 0.3:}

As it has been seen in Assumptions of Theorem 0.3,
we may begin to assume without loss of generality that $f(y,z)$
is an arbitrary W-poly of two complex variables, 
$$ f(y,z)=z^n+\sum^{n-2}_{i=0} a_iy^{\a_i}z^i, \tag 0.3.1 $$
where for $0\le i\le n-2$, each $a_i=a_i(y)$ is a unit in
${}_2\CO_0$ if exists, and the $\a_i$ are
positive integers.

\noindent$\underline{\text{\rm Step(1) of Fact[II]}}$ 
\ms

\noindent$\underline{\text{\rm Step(2) of Fact[II]}}$ 

\ms

{\bf{Example(Example [I] for Lemma 11 and Lemma 12)(2012-5300.page 30)}

(1) Let $g=g(y,z)$, $h=h(y,z)$, $f=f(y,z)$ and $\phi=\phi(y,z)$ be
in ${\BC}\{y,z\}$, which can be written as follows:
$$g=z^4+2y^5.$$
$$h=z^4+3y^4z+2y^5.$$
$$f=(z^4+2y^5)^3+4y^{10}z(z^4+2y^5)+3y^{15}z=g^3+(4y^{10}z)g+3y^{15}z.$$
$$\phi=(z^4+3y^4z+2y^5)^3+4y^{10}z(z^4+3y^4z+2y^5)+3y^{15}z=h^3+(4y^{10}z)h+3y^{15}z.$$

\text{\bf The 1st Problem:} The aim is to compute that two zero sets
$\{g=0\}$ and $\{h=0\}$ have the same kind of the singularity at
$0$.

\text{\bf The 2nd Problem:} Let $F=F(y,z)=g^3+y^{15}z$. The aim is
to compute that two zero sets $\{f=0\}$ and $\{\phi=0\}$ have the
same kind of the singularity at $0$ as $\{F=0\}$ does at $0$. \bs

{\bf{Example(Example [II]] for Lemma 11 and Lemma 12).}  By $g$ of
Example [I], $h=h(y,z)$, $f=f(y,z)$ and $\phi=\phi(y,z)$ be in
${\BC}\{y,z\}$, which can be rewritten as follows:
$$\align
h&=g+3y^4z+2y^7 =g+\sum_{\alpha,\beta\ge
0}c_{\alpha,\beta}y^{\alpha}z^{\beta},  \quad {with}\quad
4\alpha+5\beta>20  \\
f&=(z^4+2y^5)^3+4y^{10}z(z^4+2y^5)+3y^{15}z \quad \text{with \quad $g=(z^4+2y^5)$ } \\
 & =g^3+R_1g+R_0  \quad \text{with\
   $R_1=4y^{10}z$ and $R_0=3y^{15}z$ {in} ${\BC}\{y,z\}$} \\
     & =g^3+\sum_{\gamma,\delta\ge 0}b_{\gamma,\delta}y^{\gamma}z^{\delta},  \quad {with}\quad
    4\gamma+5\delta>60 \\
          &=z^{12}+6y^5z^8+4y^{10}z^5++12y^{10}z^4+ 8y^{15}z+3y^{15}z +8y^{15}\\
    & \qquad   \text{where \quad $g=(z^4+2y^5)$. }
 \endalign$$

{\bf{Example(Example [III]] for Lemma 11 and Lemma 12).}  By
$g=g(y,z)$ of Example [I], $h=h(y,z)$, $f=f(y,z)$ and
$\phi=\phi(y,z)$ be in ${\BC}\{y,z\}$, which can be rewritten as
follows:
$$\align
\phi&=(z^4+3y^4z+2y^5)^3+4y^{10}z(z^4+3y^4z+2y^5)+3y^{15}z+y^{17} \\
 & =h^3+R_1h+R_0  \quad {with}\quad
    R_i=R_i\{y,z\}\in {\BC}\{y,z\} \\
      &=g^3+9y^4zg^2+(27y^8z^2+4y^{10}z)g+27y^{12}z^3+12y^{14}z^2+3y^{15}z+y^{17}\\
      & =g^3+T_2g^2+T_1g+R_0  \quad {with}\quad
    T_i=T_i\{y,z\}\in {\BC}\{y,z\} \\
     & =g^3+\sum_{\gamma,\delta\ge 0}a_{\gamma,\delta}y^{\gamma}z^{\delta},  \quad {with}\quad
    4\gamma+5\delta>60 \\
       &=z^{12}+9y^4z^9+6y^5z^8+27y^8z^6+(6y^{9}+4y^{10})z^5+12y^{10}z^4\\
    &\quad +27y^{12}z^3+(54y^{13}+12y^{14})z^2+(36y^{14}+11y^{15})z+8y^{15}+y^{17} \\
   \endalign$$
\bs\bs

\noindent{\bf [III] The 3rd Algorithm(A complete and explicit
algorithm for computing the corresponding standard Puiseux expansion
from any irreducible W-poly of two complex variables) in Part[A]}

Note that {\bf The 3rd Algorithm consists of Theorem 1.16 together
with Theorem 1.4.} {\bf As an example for The 3rd Algorithm in
Part[A],} let $f(y,z)$ be a W-poly of two complex variables, which
is given by the above $(*1)$. By the same method just as above, it
can be computed that $f(y,z)$ of $(*1)$ is irreducible in
$\BC\{y,z\}$.

To find the standard Puiseux expansion $C(t)$ such that that $C(t)$
and the zero set $f(y,z)=0$ at $0\in \BC^2$ have the same
multiplicity sequence, we use {\bf The 1st Algorithm}.

In preparation, $(f_1,f)$ is a generalized representation of $f$ in
the sense of Theorem $16.5$ because it was just known that $f(y,z)$
of $(*1)$ can be computed as follows:
$$\align
f&=f^4_{1}+y^{17}z^2f_{1}+y^{24}z+y^{29}
=f^{d_2}_{1}+\sum^{d_2-2}_{i=0}R^{(3)}_{2,i}f_{1} \quad \text{with $f_1=h_{1,3}$,} \tag $*2$ \\
 f_1&=z^4+y^3z^2+y^5=z^{n_1}+\sum^{n_1-2}_{i=0}R^{(3)}_{1,i}z^i, \\
\endalign$$
where $f_1=h_{1,3}$ and $R^{(3)}_{2,3}=0$, $R^{(3)}_{2,2}=0$,
$R^{(3)}_{2,1}=y^{17}z^2$, $R^{(3)}_{2,0}=y^{24}z+y^{29}$ and
$R^{(3)}_{1,3}=0$, $R^{(3)}_{1,2}=y^3$, $R^{(3)}_{1,1}=0$,
$R^{(3)}_{1,0}=y^{5}$, as we have seen in {\rm(Eq.5.2.1)} of
{\rm(Eq.5) in $\S 1.10$.

Next, by  Theorem 1.16 of $\S1.9$, as it is seen in Example 1.10.2
of $\S.1.10$, we can get that $f(y,z)=0$ and $\phi(y,z)=0$ have the
same multiplicity sequence at $0\in \BC^2$ where
$\phi(y,z)={{\phi}_1}^4+y^{24}z$ where
${\phi}_1={\phi}_1(y,z)=(z^4+y^5)$. Since $\phi(y,z)$ is the
standard irreducible W-poly of two complex variables of the
recursive 2-type, then it is easy to compute by the 1st
algorithm(Theorem 1.4) that the standard Puiseux expansion $C(t)$,
which defined by $y(t)=t^{16}$ and $z(t)=t^{20}+t^{41}$, and the
zero set $\phi(y,z)=0$ have the same multiplicity sequence at
$0\in\BC^2$.  So, $C(t)$ and the zero set $f(y,z)=0$ have the same
multiplicity sequence at $0\in\BC^2$. \bs

{\bf By Appendix(as an application of the Algorithms in this
paper)}, for the above example $f\in \BC\{y,z\}$ as it was seen in
an example for The 2nd Algorithm in Part[A], we can show how to
compute the following problems (A), (B), (C),(D) with corresponding
solutions: \ms

\noindent{\bf (A):} $\underline{\text{\bf The problem is to compute
$g_2 \in$ Family[1]}}$ such that $g_2$ and $f$ have the same
multiplicity sequence.  As a unique solution, it was found by The
3rd Algorithm that $g_2=\phi_2$. \ms

\noindent{\bf (B):} $\underline{\text{\bf The problem is to compute
$C(t)\in$ Family[2]}}$ such that $C(t)$ and $f$ have the same
multiplicity sequence.  By The 3rd Algorithm in Part[A], the
standard Puiseux expansion $C(t)$ can be found by $y(t)=t^{16}$ and
$z(t)=t^{20}+t^{41}$. \ms

\noindent{\bf (C):} $\underline{\text{\bf The problem is to compute
$(g_2\circ\tau_{\xi})_{divisor}\in$ Family(4)}}$ such that
$(g_2\circ\tau_{\xi})_{divisor}=V^{(\xi)}(g_2)+\sum^{\xi}_{i=1}e_iE_i$
where $\tau=\tau_{\xi}=\pi_1\circ\pi_2\circ\cdots
\circ\pi_{\xi}:M^{(\xi)}\to \BC^2$ is the composition of a finite
number $\xi$ of successive blow-ups $\pi_i$ at the origin in
$\BC^2$, which is needed only to get the standard resolution of the
singularity of $V(g_2)$ where
$\tau^{-1}_{\xi}(0,0)=E=\cup^{\xi}_{i=1}E_{i}$, called the
exceptional curves of the first kind, is the decomposition into
irreducible components.

For any $g_r \in \text{\rm Family(1)}$, it can be proved by Theorem
14.0 of $\S 14$ of Part[C] that $(g_r\circ\tau_{\xi})_{divisor}$ can
be rewritten by
$$\align
\text{\rm(C.1)}  \qquad  (g_r\circ\tau_{\xi})_{divisor}
=V^{(\xi)}(g_r)+\sum^{\lambda_1}_{i=1}e_iE_i
+\sum^{\lambda_2}_{i=\lambda_1+1}e_iE_i+\dots+
\sum^{\lambda_{r}}_{i=\lambda_{r-1}+1}e_iE_i, \qquad \qquad \qquad
\qquad
\endalign$$
with $\lambda_r=\xi$ where each $e_i$ is the multiplicity of
$g_r\circ\tau_{\xi}$ along $E_i$ for $1\le i\le \xi=\lambda_r$ and
$V^{(\xi)}(g_r)$ is the proper transform of $V(g_r)$ under
$\tau_{\xi}$ and write $\lambda_0=0$ for notation, satisfying two
properties:

Write $L=V^{(\xi)}(g_r)\bigcup ({\bigcup}^{\xi}_{i=1} E_i)$, and for
any $A\subset M^{(\xi)}$ $\overline{A}$ is the closure of $A$ in
$M^{(\xi)}$.

\noindent{\rm(i)} For each $i=1,\dots,\xi$, $E_i\bigcap
\overline{(L-E_i)}$ has at most three distinct points under
$\tau_{\xi}$ in $L$.

\noindent{\rm(ii)}{\rm(iia)} There is a strictly increasing finite
sequence $\{\lambda_i: 1\le i\le r\}$ such that
$E_{\lambda_i}\bigcap \overline{(L-E_{\lambda_i})}$ has exactly
three distinct points under $\tau_{\lambda_r}$ in $L$ for each
$i=1,2,\dots,r$.

{\rm(iib)} For any $j\not\in \{\lambda_i: 1\le i\le r\}$ where $1\le
j\le {\xi}$, $E_j\bigcap \overline{(L-E_j)}$ has at most two
distinct points under $\tau_{\lambda_r}$ in $L$.

Since $\{e_i:i=1,2,\dots,\xi=\lambda_r\}$ is a strictly increasing
sequence and $\{\lambda_i: 1\le i\le r\}$ is a unique sequence where
${\lambda_i}<{\lambda_j}$ for all $i<j$, using the same properties
and notations as in (C.1) it is clear by Definition 9.1 and
Definition 9.2 that
$\{(g_r\circ\tau_{\lambda_r})_{divisor}\}_{seq.}\}$ is well-defined
by
$$\align
\text{\rm(C.2)} \quad
\text{$\{(g_r\circ\tau_{\lambda_r})_{divisor}\}_{seq.}
\equiv\{e_i:i=1,2,\dots,\lambda_r\}\equiv\text{\rm
Join}\{B_1,B_2,\dots,B_r\}$, as sequence,} \qquad
\endalign$$
where for $j=1,2,\dots,r$, the j-th subsequence $B_j$ is written
respectively as follows:
$$\align
\text{\rm(C.3)} \qquad \qquad & B_1=\{e_i: i=1,2,\dots,\lambda_1 \}
\quad
\text{with $1<\lambda_1<\lambda_2<\cdots<\lambda_r$,}   \\
& B_j=\{e_{\lambda_{j-1}+i}:i=1,2,\dots,(\lambda_j-\lambda_{j-1})\}
\quad \text{for $j=2,3,\dots,r$.} \qquad \qquad \qquad \qquad
\endalign$$

So, $\text{\bf {Family(4)}}$ of\text{\rm(C.1)} can be identified
with $\text{\bf Family(4)}_{seq.}$ of \text{\rm(C.4)}
$$\align
\text{\rm(C.4)} \quad \quad &\text{$\underline{\text{\rm{$\text{\bf
Family(4)}_{seq.}$}}=\text{\{$\{(g_r)_{{\text{\rm
divisor}}}\}_{seq.}:g_r\in \text{\rm
{Family(1)}}$ }}$} \\
& \text{$\underline{\text{where $\tau_\xi:M\to \BC^2$ is the
standard resolution of the singularity of {$V(g_r)$}\}}}$.} \qquad
\qquad
\endalign$$

For example, if $f=g_2$, the problem is to compute
$(g_2\circ\tau_{\xi})_{divisor}=V^{(\xi)}(g_2)+\sum^{\lambda_1}_{i=1}e_iE_i
+\sum^{\lambda_2}_{i={\lambda_1}+1}e_iE_i$ by (C.1). So,
$\{(g_2\circ\tau_{\lambda_2})_{divisor}\}_{seq.}
=\{e_i:i=1,2,\dots,\lambda_2\}=\text{\rm Join}\{B_1,B_2\}$, where
$B_1=\{e_i: 1\le i\le \lambda_1=5\}=\{16, 20, 40, 60, 80\}$ with
$e_5=80$, and $B_2=\{e_i: {{\lambda_1}+1}=4\le i\le
{\lambda_2=14}\}=\{84,88,92,96,100,101,202,303,404\}$ with
$e_{14}=404$. \ms

\noindent{\bf (D):} $\underline{\text{\bf The problem is to compute
$\text{\rm Multiseq{$(g_r(y,z)$)}}\in$ Family[3]}}$ for any $g_r \in
\text{\rm Family(1)}$ as we have seen in (C),  the problem is to
compute $\text{\rm Multiseq{$(g_r(y,z)$)}}\in \text{\rm Family[3]}$
where $\text{\rm Multiseq{$(g_r(y,z)$}}$ is defined by the
multiplicity sequence of the zero set \text{\rm V($g_r(y,z)$)} with
isolated singularity under the standard resolution.

For notation, we may assume that $g_r\in \text{\rm Family(1)}$
satisfies the same properties and notations as in (C.1). For any
given $g_r \in \text{\rm Family(1)}$, note that $\{\lambda_j:
1<\lambda_1<\cdots<\lambda_r, 1\le j\le r \}$ is a uniquely
determined finite strictly increasing sequence in (C.1). By the same
way as we have used in (C.4), $\text{\rm {Multiseq(V($g_r$))}}$ in
\text{\rm Family(3)} is well-defined by
$$\align
\text{\rm(D.1})\qquad \qquad  \text{\rm
Multiseq(V($g_r$))}\equiv\{c_i:i=1,2,\dots,\lambda_r\}
\equiv\text{\rm Join}\{P_1,P_2,\dots,P_r\}, \quad\text{as sequence}
\qquad
\endalign$$
for $j=1,2,\dots,r$, each subsequence $P_j$ of which can be uniquely
written as follows:
$$\align
\text{\rm (D.2)}  \qquad \qquad  & P_1=\{c_i: i=1,2,\dots,\lambda_1
\}
 \quad \text{and}  \qquad \qquad\qquad \qquad \\
& P_j=\{c_{\lambda_{j-1}+i}:i=1,2,\dots,(\lambda_j-\lambda_{j-1})\}
\quad \text{for $j=2,3,\dots,r$. \quad {$\square$}} \qquad\qquad
\qquad
\endalign$$
 \ms

By Theorem 1.4 and Theorem 1.6,  for any given standard Puiseux
W-poly $g_r$ of the recursive $ r$-type we can find an algorithm how
to compute the standard Puiseux expansion \text{\rm $C_r(t)$} such
that the zero set $V(g_r)$ and $C_r(t)$ have the same multiplicity
sequence. By Theorem 10.1 or Appendix, it is easy to find an
algorithm for computing $\text{\rm Multiseq{$(C_r(t)$)}}\in
\text{\rm Family[3]}$ where $\text{\rm Multiseq{$(C_r(t)$}}$ is
defined by the multiplicity sequence of the zero set of the standard
Puiseux expansion \text{\rm $C_r(t)$} with isolated singularity
under the standard resolution. Then, it is clear by Definition 9.1,
Definition 9.2 and Theorem 10.1 that $\text{\rm
Multiseq{($g_2$(y,z))}}=\{[16,20],[4,21]\}$.  \ms

\noindent{\bf(E):} $\underline{\text{\bf the problem is to compute
$(g_r\circ\tau_{\lambda_r})_{\text{\rm singular part of
divisor}}\in$ Family(5)}}$ for any $g_r \in \text{\rm Family(1)}$,
called the new family defined by this paper only, such that
$(g_r\circ\tau_{\lambda_r})_{\text{\rm singular part of divisor}}$
can be written as follows:
$$\align
(g_r\circ\tau_{\xi})_{\text{\rm singular part of the divisor}} =
\sum^{r-1}_{i=0}
\{e_{\lambda_i+1}E_{\lambda_i+1}+e_{\lambda_i+s_i}E_{\lambda_i+s_i}\}
\tag {E.1}
\endalign$$
where $\lambda_0=0$, satisfying the following properties: For
convenience of notation, let
${\Omega}^{(1)}=\bigcup^{\lambda_1}_{i=1}{E_{i}}$,
${\Omega}^{(2)}=\bigcup^{\lambda_2}_{i={{\lambda_1}+1}}{E_{i}}$,\dots,
and ${\Omega}^{(r)}=\bigcup^{\lambda_r}_{i={{\lambda_{r-1}}+1}}
{E_{i}}$.

$\underline{\text{\rm Property(1)}}$ \quad For any $E_t\subset
{\Omega}^{(w+1)}$ with $w+1\le r$, $E_t\bigcap
\overline{{\Omega}^{(w+1)}-E_{t}}$ \quad have at most two distinct
points in ${\Omega}^{(w+1)}$. \ms

$\underline{\text{\rm Property(2)}}$ Let $w$ be with $w+1\le r$, and
${\Omega}^{(w+1)}=\bigcup^{\lambda_{w+1}}_{i=\lambda_{w}+1}E_i$.

{\rm(i)} There are two distinct exceptional curves of the first kind
in ${\Omega}^{(w+1)}$, denoted by $E_{{\lambda_w}+1}$ and
$E_{{\lambda_w}+s_w}$ with $1<s_w\le\lambda_{w+1}-{\lambda_w}$, each
of which satisfies the following property: Note that
${\Omega}^{(w+1)}=\bigcup^{\lambda_{w+1}}_{i=\lambda_{w}+1}E_i$.

{\rm(ia)} $E_{\lambda_w+1}\bigcap
\overline{{\Omega}^{(w+1)}-E_{\lambda_w+1}}$ {\quad} has one and
only one intersection point in ${\Omega}^{(w+1)}$.

{\rm(ib)} $E_{{\lambda_w}+s_w}\bigcap
\overline{{\Omega}^{(w+1)}-E_{{\lambda_w}+s_w}}$ {\quad} has one and
only one intersection point in ${\Omega}^{(w+1)}$.

{\rm(ii)} $E_{{\lambda_w}+j}\bigcap
\overline{{\Omega}^{(w+1)}-E_{{\lambda_w}+j}}$ {\quad} has two
distinct intersection points in ${\Omega}^{(w+1)}$ for any $j$ where
$1<j\le\lambda_{w+1}-{\lambda_w}$ and $j\neq {s_w}$. \ms

$\underline{\text{\bf In Appendix $C$,}}$ it will be proved that we
can find an explicit algorithm for finding a one-to-one function
from Family(2) onto {{Family(5)}}(the family of the singular parts
of the divisors defined by the total transforms of irreducible plane
curves with isolated singularity under the standard resolution),
without mentioning anything about Family[4], because
$(g_r\circ\tau_{\lambda_r})_{\text{\rm singular part of divisor}}$
can be identified with a strictly increasing sequence
$\{e_1,e_{s_1};e_{{\lambda_1}+1},e_{{\lambda_1}+s_2};\dots;e_{{\lambda_{r-1}}+1},e_{{\lambda_{r-1}}+s_r}
\}$, consisting of $2r$ elements in
$\{(g_r\circ\tau_{\lambda_r})_{divisor}\}_{seq.}$. For any $g_2 \in
\text{\rm Family[1]}$, we can compute
$(g_2\circ\tau_{\xi})_{\text{\rm singular part of the divisor}}$ as
follows:
$$\align
\text{\rm (E.2)} \quad (g_2\circ\tau_{\xi})_{\text{\rm singular part
of the divisor}}
= V^{(\xi)}(g_2)+16E_1+20E_{s_1}+84E_{{\lambda_1}+1}+101E_{{\lambda_1}+s_2}\\
\endalign$$
by Appendix $C$ where $s_1=2$ and $s_2=6$. Note that $\lambda_1=5$
and $\lambda_2=14$.

Conversely, if $\{(g_2\circ\tau_{\lambda_r})_{divisor}\}_{seq.}$ is
given by  a sequence $Q=\{16,20,84,101\}$, then we can compute by
Appendix $C$ that the standard Puiseux expansion defined by
$y=t^{16}$ and $z=t^{20}+t^{41}$ and $Q$ are equivalent as far as
the multiplicity sequences are concerned.

In preparation for finding a solution of the above three problems, in this
section we are going to write some definitions, and lemmas and theorems without proof,
relative to irreducible W-polys of $n\ge 2$ complex variables at $0\in\BC^n$.

\vfill
\pagebreak

$$\align
\quad &\text{\bf \centerline{Part[A](Part[A1], Part[A2])}}\\
&\text{\bf How to use a complete and explicit irreducibility
algorithm for the W-polys of }\\
&\text{\bf two complex variables without proofs and
related topics in the Puiseux expansions}\\
\endalign$$ \ms

{\bf Part[A1] In preparation for the representation of
irreducibility algorithms for the W-polys of two complex variables}
\ms

{\bf {\S1.0.} The definition of the new terminology for three
explicit algorithms for all the irreducible W-polys of two complex
variables and some remarks} \ms

To solve all the problems in {\bf Introduction} completely, the aim
in $\S 1.0$ is to represent explicit algorithms, consisting of three
algorithms, completely, computationally and rigorously in an
elementary way without proofs, as follows: \ms

\noindent$\underline{\text{\rm By Theorem 1.4 and Theorem 1.6 of
Part[A], we can find The 1st Algorithm for computing  }}$

\noindent$\underline{\text{\rm a one-to-one map from the family of
some irreducible W-polys of two complex variables of }}$

\noindent$\underline{\text{\rm the recursive type(Family(1)) onto
the family of the standard Puiseux expansions(Family(2)).}}$

In $\S1.4$ the aim is how to use The 1st Algorithm for
computing a one-to-one map between Family(1) and Family(2) by
Theorem $1.4$ and Theorem 1.6 of $\S1.4$. Note by Definition 1.1 of $\S1.1$ and
Definition 1.2 of $\S1.2$ that each element of Family(1) is called
the standard Puiseux W-poly of the recursive type throughout this
paper. \ms

\noindent$\underline{\text{\rm By Theorem 1.13 and Theorem 1.15 of
Part[A], we can find The 2nd Algorithm for  }}$

\noindent$\underline{\text{\rm computing irreducibility criterion of all the
W-polys of two complex variables, by using  }}$

\noindent$\underline{\text{\rm Theorem 1.8(The Weierstrass division algorithm 
for the W-polys) of $\S1.5$.}}$

To compute irreducibility criterion for germs of analytic functions
of two complex variables, without loss of generality, it is clear
that it suffices to find The 2nd Algorithm for computing completely
irreducible W-polys from all the W-polys of two complex variables 
by the Weierstrass preparation
theorem and the Weierstrass division theorem. The proof of The 2nd
Algorithm can be given by Theorem 16.6 together with Proposition
16.7 and Proposition 16.8 in Part[C]. \ms

\noindent$\underline{\text{\rm By Theorem 1.16 of Part[A], we can
find The 3rd Algorithm for computing the standard }}$

\noindent$\underline{\text{\rm Puiseux expansion which has the same
multiplicity sequence as the zero set of any given   }}$

\noindent$\underline{\text{\rm irreducible W-poly in $\BC\{y,z\}$ does at
$0\in\BC^2$.}}$

In $\S1.8$, as soon as any W-poly $f(y,z)$ of two complex variables
is found to be irreducible in $\BC\{y,z\}$ by The 2nd Algorithm, in
$\S1.9$ the aim is how to use The 3rd Algorithm for computing the
standard Puiseux W-poly of the recursive type(or the standard
Puiseux expansion) directly which has the same multiplicity sequence
at $0\in\BC^2$ as the analytic variety $\{f(y,z)=0\}$ does at
$0\in\BC^2$, by The 1st Algorithm. Thus, the solution of The 3rd
Algorithm can be given by Theorem 1.16 of $\S1.9$. \ms

\definition{Remark 1.0} {\bf The complete proof of three explicit algorithms will
be done later in the other sections of this book.}

{\rm(a)} In order to succeed in the computation of The 1st
Algorithm, it is very interesting and important by Definition $1.1$
that we can define the new terminology, the standard
irreducible(Puiseux) W-polys of two complex variables of the
recursive type, which will be shown to be well-defined by Theorem
1.3 of $\S1.3$. \ms

{\rm(b)} The proofs of Theorem $1.4$ and Theorem 1.6 will be done by
Theorem 11.2 and Theorem 11.4 of $\S 11$, respectively. The proofs
of Theorem 1.15 and Theorem 1.16 will be done by Theorem 16.6 and
Theorem 19.1, respectively. \ms

{\rm(c)} $\underline{\text{\rm To find The 1st Algorithm with
proof,}}$ we use Theorem A([K2]), instead of the well-known
theorem(Theorem B), because Theorem B is not applicable to compute
The 1st algorithm in this book:

{\bf Theorem A:} Whenever any two irreducible parametrization have
the same Puiseux pairs(that is, the same standard Puisuex
expansions) by a nonsingular change of the parametrization, then
they have the same multiplicity sequences, and conversely. $\square$\ms

{\bf Theorem B:} As far as arbitrary Puiseux expansion of
irreducible plane curve singularities is concerned, any two
irreducible plane curve singularities have the same topological
types if and only if they have the same Puiseux pairs.  $\square$ \ms

Note that Theorem B was proved by K. Brauner[Br], W.Burau[Bu] and
O.Zariski[Z1] and that the proof of Theorem A was done by Theorem
8.8([K2]) and Theorem 8.10([K2])), using $\sigma$-process only,
without using Theorem B. \ms

{\rm(d)} $\underline{\text{\rm To find three algorithms with
proofs}}$, we have not studied the theory of Newton polygon. Then we
did not mention and use the definition of the Newton polygon at all
throughout this book. In preparation for the representation of the
2nd and the 3rd algorithms in $\S1.6$, it is very important to say
that Theorem 1.8(The Division Algorithm for the W-polys) of $\S1.5$
can have an important role of the 2nd and the 3rd algorithms in
$\S1.6$, which will be shown by Theorem 15.4 of $\S 15$ and Theorem
16.6 of $\S 16$. \quad $\square$
\enddefinition \ms

{\bf \S1.1. The definition for {\bf the standard irreducible W-polys
of two complex variables of the recursive type}(the standard Puiseux
W-polys of two complex variables of the recursive type)}\ms

Let $\BC\{z_1,z_2,\dots,z_n\}$ be the ring of either convergent
power series or germs of analytic functions at $0\in\BC^n$, and
${}_n\CO={}_n\CO_0$ denote the ring of germs of holomorphic
functions at $0\in\BC^n$. As the ring, $\BC\{z_1,z_2,\dots,z_n\}$
can be identified with ${}_n\CO$. Also, let $\BC[z_1,z_2,\dots,z_n]$
be the ring of polynomials of n complex variables.

\noindent{\bf(i)} $f\in
\BC\{z_1,z_2,\dots,z_{n-1}\}[z_n]={}_{n-1}\CO[z_n]$ is called 
$\underline{\text{\bf a W-poly of degree $\nu$ in $z_n$}}$ with coefficients in
$\BC\{z_1,z_2,\dots,z_{n-1}\}={}_{n-1}\CO$ if
$f(z_1,z_2,\dots,z_n)={z_n}^{\nu}+a_1{z_n}^{\nu-1}+\cdots+a_{\nu-1}{z_n}+a_{\nu}$
where the $a_j=a_j(z_1,z_2,\dots,z_{n-1})$ are holomorphic functions
at $0\in\BC^{n-1}$ and $a_j(0)=0$ for $1\le j\le \nu$. Note that
${}_{n-1}\CO[z_n]\subseteq{{}_n\CO}$.

\noindent{\bf(ii)} A nonunit $f\in {}_n\CO$ is said to be 
$\underline{\text{\bf reducible}}$ in ${}_n\CO$ if it can be written as 
a product $f=g_1g_2$ where $g_1$ and $g_2$ are not unit in ${}_n\CO$, 
and a nonunit $f\in {}_n\CO$ that is not reducible in ${}_n\CO$ is called  
$\underline{\text{\bf irreducible}}$ in ${}_n\CO$.
 \ms
 
\noindent$\underline{\text{\bf {Lemma 1.0.0(The well-known lemma).}}}$ 
$\underline{\text{\bf Assumptions}}$ Let $f$ in ${}_n\CO$ be the W-poly in $z_n$ with coefficients in ${}_{n-1}\CO$.

$\underline{\text{\bf Conclusions}}$ Then, $f$ is reducible over ${}_{n-1}\CO[z_n]$ if and only if it is reducible over ${}_n\CO$.
If the W-poly $h$ is reducible, its factors are also W-polys, up to units of ${}_n\CO$.
$\square$ \ms

\proclaim{Lemma 1.0.1}  $\underline{\text{\bf Assumptions}}$ Let
$f(y,z)$ be $\in \BC\{y,z\}$ with isolated singularity at $(0,0)$,
defined by the following:
$$\align
\text{$f(y,z)=\sum_{\alpha,\beta\ge
0}c_{\alpha,\beta}y^{\alpha}z^{\beta}$} \quad \text{with
$\alpha+\beta\ge 2$,} \qquad \qquad \tag 1.0.1.1
\endalign$$

$\underline{\text{\bf Conclusions}}$ If $f(y,z)$ is irreducible in
$\BC\{y,z\}$, we have the following:

{\rm(1)} Then, $f=f(y,z)$ can be written in the form:
$$\align
\noindent(1.0.1.2) \qquad \qquad
\text{$f=b_nz^n+b_{n-1}y^{\beta_{n-1}}z^{n-1} + \dots+b_0y^{\beta_0}$
\quad with $n\ge 2$ and $k=\beta_0\ge 2$,} \qquad \qquad
\endalign$$
where each $b_i=b_i(y,z)$ is a unit in $\BC\{y,z\}$ if exists, and
each $\beta_i$ is a positive integer for $i=1,2,\dots,n$, and
also $b_0$ and $b_n$ must be units in $\BC\{y,z\}$. \ms

{\rm(2)} By {\rm the W-preparation theorem}, $f$ of
{\rm(1.0.1.2)} can be rewritten by the following:
$$\align
f&=ug, \tag 1.0.1.3\\
g&=z^n+c_{n-1}y^{\alpha_{n-1}}z^{n-1} +
\dots+c_0y^{\alpha_0}   \text{\quad with $\alpha_0=\beta_0$},
\endalign$$
where $u=u(y,z)$ is a unit in $\BC\{y,z\}$ and $g=g(y,z)$ is a
W-poly in $z$ such
that each $c_i=c_i(y)$ is a unit in $\BC\{y\}$ if exists, and the
$\alpha_i$ are positive integers for $i=1,2,\dots,n$. $\square$
\endproclaim \ms

In order to succeed in the computation of our three algorithm, in
this section it is very interesting and important to prove that we
can construct the new terminology, some irreducible W-polys
of two complex variables of the recursive type,  called
{\rm ``the standard irreducible W-polys of two complex variables of the
recursive type"}
throughout this book, which will be shown to be
equivalent to the standard Puiseux expansion with Explicit Algorithm
of $\S11$, as far as the multiplicity sequences of irreducible plane
curve singularities are concerned. \ms

\definition{Definition 1.0.2(The new terminology)} Let $N_0$ be the set of
nonnegative integers and $N^k_0$ be its $k$-dimensional copy. Let
$r$ be an arbitrary positive integer. 

{\bf } By [B4] of Definition 5.0.0, $g_r\in\BC\{y,z\}$ is 
called $\underline{\text{\bf the standard Puiseux polynomial of }}$ 

\noindent$\underline{\text{\bf two complex variables of the recursive r-type in $z$}}$
if there are sequences $\{X_k:k=1,2,\dots,r\}$ with
$X_k\subset N_0$, $\{g_k:k=1,2,\dots,r\}$ with $g_k\in\BC\{y,z\}$
and $\{\text{${{\Delta}_k}:N^k_0\to N_0$ is an }$
$\text{integer-valued function for $k=1,2,\dots,r$}\}$
satisfying the following $\underline{\text{\rm six conditions}}$: \ms

\noindent$\underline{\text{\rm Six conditions}}$ are denoted by
\text{\rm The 1st ${\text{\rm{Cond}}}^{\text{{\rm(0)}}}$}, $\dots$,
\text{\rm The 6-th ${\text{\rm{Cond}}}^{\text{{\rm(0)}}}$} for
notation. \ms

\noindent $\underline{\text{\rm The 1st
${\text{\rm{Cond}}}^{\text{{\rm(0)}}}$}}$ Let
$\{X_j:j=1,2,\dots,r\}$ with $X_j\subset N_0$ be defined as follows:

\roster
\item"(1) (1a)"
 $X_1=\{n_1,\beta_{1,1}\}$ with $n_1\ge 2$ and $\beta_{1,1}\ge 1$.
 \item"$\quad$ (1b)" $X_j=\{n_j,\beta_{j,1},\beta_{j,2},\dots,\beta_{j,j}\}$
 with $n_j\ge 2$ \quad where $j=2,\dots,r$.
 \endroster
If $j\ge 2$, assume that at least one of
$\beta_{j,1},\beta_{j,2},\dots,\beta_{j,j}$ is nonzero. \ms

\noindent $\underline{\text{\rm The 2nd
${\text{\rm{Cond}}}^{\text{{\rm(0)}}}$}}$ For each $j=1,2,\dots,r$,
let $g_j=g_j(y,z)$ be in $\BC\{y,z\}$, each of which is defined by
the following way:
 \roster
 \item"(2) (2a)" $g_1=z^{n_1}+y^{\beta_{1,1}}$.
  \item"$\quad$ (2b)"
 $g_j=g^{n_j}_{j-1}+y^{\beta_{j,1}}z^{\beta_{j,2}}g^{\beta_{j,3}}_1\cdots
 g^{\beta_{j,j}}_{j-2}$ \quad with $g_{-1}=y$ and $g_0=z$, where $j=2,\dots,r$.
\endroster
\ms

\noindent $\underline{\text{\rm The 3rd
${\text{\rm{Cond}}}^{\text{{\rm(0)}}}$}}$ Let $\{\Delta_k: N^k_0\to
N_0: k=1,2,\dots,r\}$ be a sequence such that each $\Delta_k$ is an
integer-valued function defined by the following:

\roster
\item"(3) (3a)" $\Delta_1(t)=t$ for each $t\in N_0$.
\item"$\quad$ (3b)"
$\Delta_j(t_k)^j_{k=1}=t_j\Delta_{j-1}(\beta_{j-1,k})^{j-1}_{k=1}
+n_{j-1}\Delta_{j-1}(t_k)^{j-1}_{k=1}$ for each $(t_k)^j_{k=1}\in
N^j_0$

\noindent where $j=2,\dots,r$.
\endroster \ms

\noindent $\underline{\text{\rm The 4-th
${\text{\rm{Cond}}}^{\text{{\rm(0)}}}$}}$ The following inequalities
hold: Note that $2\le j\le r$. \roster
\item"(4) (4a)" $\Delta_1(\beta_{1,1})=\beta_{1,1}>0$ with $n_1\ge 2$.
\item"$\quad$ (4b)" $\Delta_j(\beta_{j,k})^j_{k=1}
>n_jn_{j-1}\Delta_{j-1}(\beta_{j-1,k})^{j-1}_{k=1}$
where $j=2,\dots,r$.
\endroster \ms

\noindent $\underline{\text{\rm The 5-th
${\text{\rm{Cond}}}^{\text{{\rm(0)}}}$}}$ \quad The following
inequalities hold:

\roster \item"(5)(5a)" $\gcd(n_j,\Delta_j(\beta_{j,k})^j_{k=1})=1$
for $1\le j\le r$.
\endroster

\noindent $\underline{\text{\rm The 6-th
${\text{\rm{Cond}}}^{\text{{\rm(0)}}}$}}$ \quad The following
inequalities hold: Note that $2\le j\le r$.

\noindent\rm{(6)(6a)}  $2\le n_1<\beta_{1,1}$. \ms
\rm{(6b)} $n_{j}\ge 2$, $\beta_{j,1}>0$, and $0\le
\beta_{j,k}<n_{k-1}$ for $2\le j\le r$ and $2\le k\le j$. 
$\square$ \ms
\enddefinition \ms

\proclaim{Lemma 1.0.3} $\underline{\text{\bf {Assumptions}}}$ 
Let $g_r$ be the standard Puiseux polynomial of two complex 
variables of the recursive r-type in $z$
as in {\rm Definition 1.0.2}. \ms

$\underline{\text{\bf {Conclusions}}}$

{\rm(i)} It is clear that $g_r(y,z)$ in {\rm Definition 1.0.2} is a W-poly
of degree $n_1n_2\cdots n_r$ in $z$ with coefficients in $\BC\{y\}$,
without using \text{\rm The 5-th
${\text{\rm{Cond}}}^{\text{{\rm(0)}}}$} in {\rm Definition 1.0.2}.

{\rm(ii)} It can be proved by {\rm Theorem 5.0} that $g_r\in \BC\{y,z\}$
of {\rm Definition 1.0.2} is irreducible in $\BC\{y,z\}$, without using
\text{\rm The 6-th
${\text{\rm{Cond}}}^{\text{{\rm(0)}}}$} in {\rm Definition 1.0.2.  $\square$
\endproclaim \ms

\definition{Definition 1.1 (by Definition 1.0.2 and Lemma 1.0.3)} Let $g_r\in\BC\{y,z\}$ be 
the standard Puiseux polynomial of the recursive r-type in $z$
as in Definition 1.0.2, satisfying the same properties and
notations as in Definition 1.0.2. Hereafter, $g_r\in\BC\{y,z\}$ is called   
$\underline{\text{\bf the standard irreducible(Puiseux) W-poly of 
two complex variables of the }}$ 

\noindent$\underline{\text{\bf  recursive r-type  
}}$ instead of using this terminology in Definition 1.0.2 throughout this book.
\enddefinition \ms

\definition{Remark 1.1.1}
In preparation for computing The 1st Algorithm, the aim is how to construct 
the standard Puiseux W-polys from the irreducible W-polys as we have seen in the
construction of the standard Puiseux expansions from the Puiseux
expansions. \quad  $\square$
\enddefinition \ms

Assuming that the proof of Theorem 5.0 is done, the contents of Definition 1.1 
can be rewritten as follows.

\definition{Definition 1.1} Let $N_0$ be the set of
nonnegative integers and $N^k_0$ be its $k$-dimensional copy. Let
$r$ be an arbitrary positive integer. \ms

\noindent{\bf } $g_r\in\BC\{y,z\}$ is called $\underline{\text{\rm
the standard irreducible(Puiseux) W-poly of two complex variables }}$ 

\noindent$\underline{\text{\rm of the recursive r-type in $z$}}$
if there are sequences $\{X_k:k=1,2,\dots,r\}$ with
$X_k\subset N_0$, $\{g_k:k=1,2,\dots,r\}$ with $g_k\in\BC\{y,z\}$
and $\{\text{${{\Delta}_k}:N^k_0\to N_0$ is an integer-valued
function for}$ $\text{$k=1,2,\dots,r$}\}$ satisfying the following
$\underline{\text{\rm six conditions}}$: \ms

\noindent$\underline{\text{\bf Six conditions}}$ are denoted by
\text{\rm The 1st ${\text{\rm{Cond}}}^{\text{{\rm(0)}}}$}, $\dots$,
\text{\rm The 6-th ${\text{\rm{Cond}}}^{\text{{\rm(0)}}}$} for
notation. \ms

\noindent $\underline{\text{\rm The 1st
${\text{\rm{Cond}}}^{\text{{\rm(0)}}}$}}$ Let
$\{X_j:j=1,2,\dots,r\}$ with $X_j\subset N_0$ be defined as follows:

\roster
\item"(1) (1a)"
 $X_1=\{n_1,\beta_{1,1}\}$ with $n_1\ge 2$ and $\beta_{1,1}\ge 1$.
 \item"$\quad$ (1b)" $X_j=\{n_j,\beta_{j,1},\beta_{j,2},\dots,\beta_{j,j}\}$
 with $n_j\ge 2$ \quad where $j=2,\dots,r$.
 \endroster
If $j\ge 2$, assume that at least one of
$\beta_{j,1},\beta_{j,2},\dots,\beta_{j,j}$ is nonzero. \ms

\noindent $\underline{\text{\rm The 2nd
${\text{\rm{Cond}}}^{\text{{\rm(0)}}}$}}$ For each $j=1,2,\dots,r$,
let $g_j=g_j(y,z)$ be in $\BC\{y,z\}$, each of which is defined by
the following way:
 \roster
 \item"(2) (2a)" $g_1=z^{n_1}+y^{\beta_{1,1}}$.
  \item"$\quad$ (2b)"
 $g_j=g^{n_j}_{j-1}+y^{\beta_{j,1}}z^{\beta_{j,2}}g^{\beta_{j,3}}_1\cdots
 g^{\beta_{j,j}}_{j-2}$ \quad with $g_{-1}=y$ and $g_0=z$, where $j=2,\dots,r$.
\endroster
\ms

\noindent $\underline{\text{\rm The 3rd
${\text{\rm{Cond}}}^{\text{{\rm(0)}}}$}}$ Let $\{\Delta_k: N^k_0\to
N_0: k=1,2,\dots,r\}$ be a sequence such that each $\Delta_k$ is an
integer-valued function defined by the following:

\roster
\item"(3) (3a)" $\Delta_1(t)=t$ for each $t\in N_0$.
\item"$\quad$ (3b)"
$\Delta_j(t_k)^j_{k=1}=t_j\Delta_{j-1}(\beta_{j-1,k})^{j-1}_{k=1}
+n_{j-1}\Delta_{j-1}(t_k)^{j-1}_{k=1}$ for each $(t_k)^j_{k=1}\in
N^j_0$

\noindent where $j=2,\dots,r$.
\endroster \ms

\noindent $\underline{\text{\rm The 4-th
${\text{\rm{Cond}}}^{\text{{\rm(0)}}}$}}$ The following inequalities
hold: Note that $2\le j\le r$. \roster
\item"(4) (4a)" $\Delta_1(\beta_{1,1})=\beta_{1,1}>0$ with $n_1\ge 2$.
\item"$\quad$ (4b)" $\Delta_j(\beta_{j,k})^j_{k=1}
>n_jn_{j-1}\Delta_{j-1}(\beta_{j-1,k})^{j-1}_{k=1}$
where $j=2,\dots,r$.
\endroster \ms

\noindent $\underline{\text{\rm The 5-th
${\text{\rm{Cond}}}^{\text{{\rm(0)}}}$}}$ \quad The following
inequalities hold:

\roster \item"(5)(5a)" $\gcd(n_j,\Delta_j(\beta_{j,k})^j_{k=1})=1$
for $1\le j\le r$.
\endroster

\noindent $\underline{\text{\rm The 6-th
${\text{\rm{Cond}}}^{\text{{\rm(0)}}}$}}$ \quad The following
inequalities hold: Note that $2\le j\le r$.

\noindent\rm{(6)(6a)}  $2\le n_1<\beta_{1,1}$.

\rm{(6b)} $n_{j}\ge 2$, $\beta_{j,1}>0$, and $0\le
\beta_{j,k}<n_{k-1}$ for $2\le j\le r$ and $2\le k\le j$. \quad
$\square$
\enddefinition \ms

\definition{Remark 1.1.1}

{\rm(i)} It is clear that $g_r(y,z)$ of Definition 1.1 is a W-poly
of degree $n_1n_2\cdots n_r$ in $z$ with coefficients in $\BC\{y\}$,
without using \text{\rm The 5-th
${\text{\rm{Cond}}}^{\text{{\rm(0)}}}$}.

{\rm(ii)} It can be proved by Theorem $5.0$ that $g_r\in \BC\{y,z\}$
of Definition 1.1 is irreducible in $\BC\{y,z\}$, without using
\text{\rm The 6-th ${\text{\rm{Cond}}}^{\text{{\rm(0)}}}$}.

{\rm(iii)} In preparation for computing
The 1st Algorithm, the aim is how to find the difference between the construction 
of the standard Puiseux W-polys from the irreducible W-polys and the
construction of the standard Puiseux expansions from the Puiseux
expansions. \quad  $\square$
\enddefinition \ms

{\bf \S1.2. The new terminology and notations in preparation for
studying three families with equivalence relations}

In preparation for a good success of the main aim, we define three
Families, that is, Family(1), Family(2), Family(3) with equivalence
relations, respectively, as follows: Throughout this book, the
family and the set have the same meaning.

\definition{Definition 1.2}
{\bf[I](1)} For brevity, let $\underline{\text{\bf{Family(0)}}}$ be
the $0$-th family, consisting of all the convergent power series $f
\in \BC\{y,z\}$ such that $f$ is irreducible in $\BC\{y,z\}$ with
isolated singularity at $0\in \BC^2$. \ms 
\noindent{\bf(2)} For convenience of notation, we use the following:

For any two finite sequences $A=\{a_i\in \BC:1\le i\le n\}$ and
$B=\{b_i \in \BC:1\le i\le m\}$, it is said that $A$ and $B$ are
equivalent if and only if either $a_i=b_i$ for $i=1,2,\dots,n=m$ or
${A}\equiv{B}$ as sequence. If $A$ and $B$ are not equivalent, we
write ${A}\not \equiv{B}$ as sequence.\ms

\noindent{\bf[II](1)} $\underline{\text{\bf{Family(1)}}}$ is called
the first family, consisting of all the standard Puiseux W-polys $f
\in \BC\{y,z\}$ of the recursive type with isolated singularity at
$0\in \BC^2$, denoted by \ms

\noindent$\text{\rm(1.2.1)}  \quad
\underline{\text{\rm{Family(1)}}=\{\text{\rm $f$ is arbitrary
standard Puiseux W-poly of the recursive r-type:}}$

\quad\quad $\underline{\text{$f\in \text{\rm Family(0) and $r$ are
arbitrary positive integers}$}}\}$. \ms

Now, in preparation for finding an equivalence relation on Family(1)
by (1c) later, it remains to define arbitrary two elements $g_r$ and
$\phi_{\rho}$ in Family(1) by (1a) and (1b), as follows:\ms

\noindent\text{(\bf 1a)} $g_r\in \BC\{y,z\}$ is called the standard
Puiseux W-poly in $\BC[y,z]$ of the recursive r-type if $g_r\in
\BC\{y,z\}$ satisfies the same properties and notations as we have
seen in Definition 1.1. \ms

\noindent\text{(\bf 1b)} By the same method as we have seen in
either Definition 1.1 or (1a), another element $\phi_{\rho}\in
\text{\rm Family(1)}$ is called the standard Puiseux W-poly
$\phi_{\rho}\in \BC\{y,z\}$ of the recursive {$\rho$}-type at
$(0,0)\in \BC^2$, if there are sequences $\{W_k:k=1,2,\dots,\rho\}$
with $W_k\subset N_0$, $\{\phi_k:k=1,2,\dots,\rho\}$ with $\phi_k\in
\BC\{y,z\}$ and $\{\text{$\omega_k:N^k_0\to N_0$ is an
integer-valued function for}$ $\text{$k=1,2,\dots,\rho$}\}$
satisfying six conditions: \ms

\noindent$\underline{\text{\rm Six conditions for $\phi_\rho$}}$ are
denoted by \text{\rm The 1st
${\text{\rm{Cond}}}^{\text{{\rm(0)}}}$}, $\dots$, \text{\rm The 6-th
${\text{\rm{Cond}}}^{\text{{\rm(0)}}}$}. \ms

\noindent $\underline{\text{\rm The 1st
${\text{\rm{Cond}}}^{\text{{\rm(0)}}}$}}$ Let
$\{W_j:j=1,2,\dots,\rho\}$ with $W_j\subset N_0$ be defined as
follows:

\roster
\item"(1) (1.1)"$W_1=\{\ell_1,\delta_{1,1}\}$ with $\ell_1\ge 2$ and
 $\delta_{1,1}\ge 1$.

\item"(1.2)"
$W_{j}=\{\ell_{j},\delta_{j,1},\delta_{j,2},\dots,\delta_{j,j}\}$
with $\ell_{j}\ge 2$ \quad where $j=2,\dots,\rho$.
\endroster
If $j\ge 2$, then assume that at least one of
$\delta_{j,1},\delta_{j,2},\dots,\delta_{j,j}$ is nonzero. \ms

\noindent $\underline{\text{\rm The 2nd
${\text{\rm{Cond}}}^{\text{{\rm(0)}}}$}}$ For each $j=1,2,\dots,r$,
let $\phi_j=\phi_j(y,z)$ be in $\BC\{y,z\}$, each of which is
defined by the following way:

\roster
\item"(2) (2.1)"$\phi_1=z^{\ell_1}+y^{\delta_{1,1}}$.

\item" (2.2)"
$\phi_{j}=\phi_{j-1}^{\ell_{j}}+y^{\delta_{j,1}}z^{\delta_{j,2}}
      \phi_1^{\delta_{j,3}}\cdots \phi_{j-2}^{\delta_{j,j}}$.
\quad where $j=2,\dots,{\rho}$.
\endroster
\ms

\noindent $\underline{\text{\rm The 3rd
${\text{\rm{Cond}}}^{\text{{\rm(0)}}}$}}$ Let $\{\omega_k: N^k_0\to
N_0: k=1,2,\dots,\rho\}$ be a sequence such that each $\omega_k$ is
an integer-valued function defined by the following:

\roster
\item"(3) (3.1)" $\omega_1(t)=t$ for each $t\in N_0$.

\item"(3.2)"
$\omega_j(t_k)^{j}_{k=1}=t_j\omega_{j-1}(\delta_{j-1,k})^{j-1}_{k=1}
+\ell_{j-1}\omega_{j-1}(t_k)^{j-1}_{k=1}$ for each
$(t_k)^{j}_{k=1}\in N^{j}_0$

\noindent where $j=2,\dots,{\rho}$.
\endroster \ms

\noindent $\underline{\text{\rm The 4-th
${\text{\rm{Cond}}}^{\text{{\rm(0)}}}$}}$ The following inequalities
hold: Note that $2\le j\le \rho$.
\roster
\item"(4) (4.1)"
$\omega_j(\delta_{1,1})=\delta_{1,1}>0$ with
\item" (4.1)"
$\omega_j(\delta_{j,k})^j_{k=1}>\ell_j
\ell_{j-1}\omega_{j-1}(\delta_{j-1,k})^{j-1}_{k=1}$  for $2\le j\le
\rho$.
\endroster \ms

\noindent $\underline{\text{\rm The 5-th
${\text{\rm{Cond}}}^{\text{{\rm(0)}}}$}}$ The following inequalities
hold:

\roster
\item"(5) (5.1)"
$\gcd(\ell_{j},\omega_{j}(\delta_{j,k})^{\rho}_{k=1})=1$ for $1\le
j\le \rho$.   \endroster \ms

\noindent$\underline{\text{\rm The 6-th
${\text{\rm{Cond}}}^{\text{{\rm(0)}}}$}}$ \quad The following
inequalities hold: Note that $2\le j\le \rho$.

\roster

\item"(6)(6.1)"  $2\le {\ell}_1<\delta_{1,1}$.

\item"(6.2)"  ${\ell}_{j}\ge 2$, $\delta_{j,1}>0$, and $0\le
\delta_{j,k}<{\ell}_{k-1}$ for $2\le j\le \rho$ and $2\le k\le j$.
\endroster \ms

\noindent{\bf(1c)} For any standard Puiseux W-poly $g_r\in
\BC\{y,z\}$ of the recursive r-type in (1a) and any standard Puiseux
W-poly $\phi_{\rho}\in \BC\{y,z\}$ of the recursive ${\rho}$-type in
(1b), it is said that $g_r$ and $\phi_{\rho}$ are equivalent,
denoted by $g_r\equiv\phi_\rho$ in Family(1), if the following are
satisfied:
$$\align
(1.2.2) \qquad \qquad &     \text{$n_j=\ell_j$  \quad\qquad for
each $j=1,2,\dots,r=\rho$, \quad  and}     \\
 & \text{$\beta_{j,k}=\delta_{j,k}$
 \qquad  for each $j=1,2,\dots,r$ \quad {and} \quad for all $k=1,2,\dots,j.$}
\qquad \qquad\\
 \endalign$$ \ms

\noindent{\bf(2)} $\underline{\text{\bf{Family(2)}}}$ is the 2nd
family, consisting of all the irreducible curves with the standard
Puisuex expansions, denoted by \ms

\noindent(1.2.3)\quad \text{Family(2)=\{$C_r(t)$:$C_r(t)$ is the
standard Puiseux expansion of the $r$-type for any $r\in \N$\}}. \bs

In more detail, we define the standard Puiseux expansion $C_r(t)$ of
the $r$-type by (1.2.4), and also the standard Puiseux expansion
$C_s(t)$ of the $s$-type by (1.2.5), respectively. After then, we
will define an equivalence relation for any two standard Puisuex
expansions $C_r(t)$ of the $r$-type and $C_s(t)$ of the $s$-type by
(1.2.6) of (2c).

Note by Definition 8.1 that the parametrization $C(t)$ for arbitrary
irreducible plane curve $C$ can be defined by $y(t)=t^n$ and
$z(t)=c_1t^{k_1}+c_2t^{k_2}+\cdots=c_1t^{k_1}(1+H(t))$, where $1<n$,
$1<k_1<k_2<\cdots, $ and the $c_i$ are nonzero complex numbers and
$H(t)$ is just the substitution.

\noindent{\bf (2a)} The standard Puiseux expansion $C_r(t)$ of the
$r$-type for the curve $C$ is as follows:
$$\align
(1.2.4) \qquad \text{$C_r(t):=$} \left\{\eqalign{ y=&t^n, \cr
z=&t^{\alpha_1}+t^{\alpha_2}+\cdots +t^{\alpha_r}, \cr} \right. \\
\text{where} \quad
 2\le n <\alpha_1<\alpha_2<\cdots <\alpha_r  & \quad  \text{and} \\
 n >d_1>d_2>\cdots
 >d_r=1  \quad   \text{with} & \quad
\text{$d_i=\gcd(n,\alpha_1,\dots,\alpha_i)$, $1\le i\le r$.} \qquad
\qquad
\endalign$$

\noindent{\bf (2b)} The standard Puiseux expansion $C_s(t)$ of the
$s$-type for the curve $C'$ is as follows:
$$\align
(1.2.5) \qquad \text{$C_s(t):=$} \left\{\eqalign{ y=&t^m, \cr
z=&t^{\beta_1}+t^{\beta_2}+\cdots +t^{\beta_s}, \cr} \right. \\
\text{where} \quad
 2\le m <\beta_1<\beta_2<\cdots <\beta_s & \quad \text{and} \\
 m >\overline{d}_1>\overline{d}_2>\cdots
 >\overline{d}_s=1  \quad \text{with} & \quad
\text{${\overline{d}}_i=\gcd(m,\beta_1,\dots,\beta_i)$, $1\le i\le
s$.} \qquad \qquad
\endalign$$

\noindent{\bf (2c)} Whenever the standard Puiseux expansions
$C_r(t)$ of the $r$-type in (1.2.4) and $C_s(t)$ of the $s$-type in
(1.2.5) are chosen arbitrary, then it is said that $C_r(t)$ and
$C_s(t)$ are equivalent if the following conditions are satisfied:
$$\align
\text{$n=m$ \quad and \quad $\alpha_i=\beta_i$ {\quad} for {\quad}
$i=1,2,\dots,r=s$.} \tag 1.2.6
\endalign$$ \ms

\noindent{\bf(3)} $\underline{\text{\bf{Family(3)}}}$ is the 3rd
family, consisting of all the multiplicity sequences of irreducible
plane curves with isolated singularity under the standard
resolution, denoted by \ms

\noindent(1.2.7) \qquad
\text{\rm{Family(3)}}=\text{\{Multiseq(V(f)):$f\in \text{\rm
{Family(0)}}$ and f is irreducible in in ${}_2\CO$\}} \ms

{\noindent}where for any $f\in \text{\rm {Family(0)}}$, we define
$\text{\rm Multiseq(V(f))}$ by the multiplicity sequence of $f$, and
next an equivalence relation for any two multiplicity sequences in
Family(3), as follows: \ms

\noindent{\bf (3a)} For any $f\in \text{\rm {Family(0)}}$, to define
the multiplicity sequence of $f$, we may assume for notation that
$\tau_{\xi}:M^{(\xi)}\to \BC^2$ with $\tau=\tau_{\xi}$ is the
standard resolution of the singularity of $V(f)$ as we have used for
$f\in \text{\rm {Family(0)}}$, which is the the composition of a
finite number $\xi$ of successive blow-ups $\pi_i$ at the origin in
$\BC^2$. Let $c_0$ be the multiplicity of this curve germ $f$ at
this point. If we blow up once, then we again find at most one
singularity. Let $c_1$ be the multiplicity of the curve of the germ
blown up once, $c_2$ be the multiplicity of the curve of the germ
blown up twice, and continue to the standard resolution. The
sequence ends with a sequences of ones. The sequences of these
multiplicities, $\{c_0,c_1,\dots,c_{{\xi}-1}\}$, where the last one
is not is not counted, is then the multiplicity sequence.

Then, $\text{\rm Multiseq(V(f))}$ is written as follows:
$$\align
\text{\rm Multiseq(V(f))}=\{c_i:i=0,1,\dots,{\xi}-1\}. \tag 1.2.8
\endalign$$

\noindent{\bf(3b)} For any $f$ and $g$ in Family(0), an equivalence
relation for any two multiplicity sequences \text{\rm
Multiseq(V(f))} and \text{\rm Multiseq(V(g))} in Family(3) is
defined as follows:
$$\align
& \text{$f$ and $g$ in Family(0) have the same multiplicity sequence} \tag 1.2.9\\
 \text{$\iff$} \qquad & either \quad  \text{$f \buildrel \text{{\rm
multiseq}} \over \sim g$} \quad or \quad \text{$V(f) \buildrel
\text{{\rm multiseq}} \over \sim V(g)$} \quad \text{at
$0\in \BC^2$}\\
 \text{$\iff$} \qquad & \text{\text{\rm
Multiseq(V(f))}{$\equiv$}\text{\rm Multiseq(V(g))} as sequence.}
\qquad  \square \qquad \qquad \qquad
\endalign$$
\enddefinition \ms

{\bf \S1.3. What does an equivalence relation of any two elements in
Family(1)(the family of the standard Puiseux W-polys in $\BC\{y,z\}$
of the recursive type) mean?} \ms

In order to succeed in the computation of The 1st Algorithm for
finding a one-to-one correspondence between Family(1) and Family(2),
the aim is to study the equivalent class of the standard Puiseux
W-polys in Family(1) with respect to the multiplicity sequences.

\proclaim{Theorem 1.3} $\underline{\text{\bf {Assumptions}}}$ Let
$r$ and $\rho$ be arbitrary positive integers. By the same way as in
Definition 1.1, let $g_r$ be the standard irreducible W-poly in $z$
of the recursive r-type, satisfying the same properties and
notations as in Definition 1.1. Also, let $\phi_{\rho}$ be the
standard irreducible W-poly in $z$ of the recursive $\rho$-type,
satisfying the same properties and notations as in Definition 1.2.

\noindent $\underline{\text{\bf Conclusions}}$ \quad Then, we have
the following:
$$\align
(1.3.1) \qquad \qquad  & \text{$g_r$ and $\phi_{\rho}$ have the same
multiplicity sequence.}  \\
  \iff \quad &
 \text{$n_j=\ell_j$ and $\beta_{j,k}=\delta_{j,k}$ for all
$j=1,2,\dots,r=\rho$ and all
$k=1,2,\dots,j$.} \qquad \qquad \\
   \quad & \text{That is, $g_r$ and $\phi_{\rho}$ are equivalent in the sense of
Definition 1.2.}\\
\endalign$$

Moreover, it can be easily proved by Theorem $7.3$ that the
following holds:
$$\align
& \text{$g_r$ and $\phi_{\rho}$ have the same multiplicity
sequence} \tag 1.3.2 \\
 \iff \quad &
\text{$g_r \buildrel \text{{\rm divisor}} \over \sim \phi_\rho$
under the standard resolutions. \quad $\square$}
\endalign$$
\endproclaim

\definition{Remark 1.3.1} {\rm(i)} Theorem 1.3 can be
proved by Theorem 7.3 and Theorem 10.2 where the new terminology of
(1.3.2) is defined by Definition 2.4 and Definition 2.6.

{\rm(ii)} Without assuming that both $2\le n_1<\beta_{1,1}$ and
$2\le \ell_1<\delta_{1,1}$, it can be easily proved that the
conclusion of Theorem 1.3 may not be true by the following example:
$$\align
& g_1=z^3+y^8 \qquad \text{and} \tag 1.3.1.1\\
& \phi_2=\phi^3_1+y^2z^3 \quad \text{with \quad $\phi_1=z+y^2$,}
\endalign$$
because $g_1$ and $\phi_{2}$ have the same multiplicity sequence,
and also they have the same divisor under two standard resolutions,
but the condition in (1.3.1) does not hold. \quad $\square$
\enddefinition \ms

\noindent{\bf Remark 1.3.2.} It will be proved by Theorem $1.3$ and
Theorem A(Theorem 8.10([K2])) that we can compute one-to-one
function $F$ from Family(1) into Family(2). It will be proved later
that such a function \text{F:Family(1)$\mapsto$ Family(2)} must be
onto.  \quad $\square$  \bs

\vfill \pagebreak

{\bf Part[A2] The rigorous representation of explicit
irreducibility algorithms for the Weierstrass polynomials of two
complex variables with examples and without proofs and related
topics in the Puiseux expansions} \bs

{\bf \S1.4. The 1st Algorithm for computing a one-to-one function
between Family(1) and Family(2)(the family of the standard Puiseux
expansions) with its examples} \ms

In this section, the first half of The 1st Algorithm can be given by
Algorithm 1.4.1 for Theorem 1.4 with Example 1.4.1, and also the
second half of the 1st Algorithm can be given by Algorithm 1.6.2 for
Theorem 1.6 with Example 1.6.3. \ms

{\bf \S1.4.A. The first half of The 1st Algorithm(Theorem 1.4)}

\proclaim{Theorem 1.4(Theorem 11.2:Algorithm for finding a
one-to-one function from Family(1) into Family(2))}

$\underline{\text{\bf {Assumptions}}}$ Let $g_r\in \BC\{y,z\}$ be
the standard Puiseux W-poly of the recursive r-type in $z$ in {\rm
Family(1)} satisfying six conditions with the same notations as in
{\rm Definition 1.1}. \ms

$\underline{\text{\bf Conclusions}}$

\noindent$\underline{\text{\rm {\bf [I]} {\rm By explicit algorithm
in (1.4.1),} we can compute the standard Puiseux expansion}}$

\noindent$\underline{\text{\rm $C(g_r:t)$ in Family(2) such that
\text{$\text{\rm{Multiseq}}(V(g_r))\equiv \text{\rm
Multiseq}(C(g_r:t))$  \text{\rm as sequence:}}}}$ \ms

\noindent$\underline{\text{\rm(Algorithm 1.4.1 for Theorem 1.4)}}$
$$\align
 \text{$C(g_r:t):=$} & \left\{\eqalign{y&=t^n \cr
z&=t^{\alpha_1}+t^{\alpha_2}+\cdots +t^{\alpha_r}, \cr} \right. \tag
1.4.1
 \\
 \text{such that}  \qquad
 n &=n_1n_2\cdots n_r \quad \text{and} \quad
 \alpha_1 =\beta_{1,1}n_2\cdots n_r, \\
 \alpha_j &=\alpha_{j-1}+
\widehat{\Delta}_jn_{j+1}n_{j+2}\cdots
 n_r, \\
 \endalign$$
where $\widehat{\Delta}_j
=\Delta_j(\beta_{j,k})^j_{k=1}-n_jn_{j-1}\Delta_{j-1}(\beta_{j-1,k})^{j-1}_{k=1}>0$
for $2\le j\le r$ and $\Delta_1(t)=t$. \ms

\noindent{\bf [II]} Let \text{\rm $\Psi$
:Family(1)$\rightarrow$Family(2)} be a function defined by
$\Psi(g_r)=C(g_r:t)$ for any $g_r$ in {\rm Family(1)}. Then, $\Psi$
is a one-to-one function from {\rm {Family(1)}} into {\rm
Family(2)}. \quad $\square$
\endproclaim \ms

\noindent{\bf Remark 1.4.0.} {\rm(a)} Note that the parametrization
in (1.4.1) satisfies the following:

\roster

\item"({}a.1)"  $n<\alpha_1<\alpha_2<\alpha_3<\cdots <\alpha_r$.

\item"(a.2)"  $n> d_1>d_2>\cdots >d_r=1$ with $d_i=
\gcd(n,\alpha_1,\dots,\alpha_i)$, where $d_i=n_{i+1}n_{i+2}\cdots
n_r$ for $i=1,2,\dots,r-1$. So, $\alpha_j-\alpha_{j-1}=
\widehat{\Delta}_jd_j.$
\item"(b)"  $\Psi$, being given by a one-to-one function from
{\rm {Family(1)}} into {\rm Family(2)} by Algorithm 1.4.1 for
Theorem 1.4, can be proved to onto map which is computable by
Algorithm 1.6.2 for Theorem 1.6.  \quad $\square$
\endroster\ms

\noindent{\bf Example 1.4.1 for Theorem 1.4:}

Let the polynomial $g_3$ in $\BC[y,z]$ be given as follows:
$$
 \text{$g_1=z^3+y^4$, $g_2=g_1^{5}+y^{18}z^{2}$,
$g_3=g_2^{3}+y^{61}z^{1}$.} \tag 1.4.2
$$
\roster

Then, $g_3$ is the standard Puiseux W-poly of the recursive 3rd type
because of the following computations {\rm(i)}, {\rm(ii)},
{\rm(iii)} and {\rm(iv)}:

\item"(i)" By Definition $1.1$ and (1.4.2), $n_1=3$, $n_2=5$, $n_3=3$,
$\Delta_1(\beta_{1,1})=\beta_{1,1}=4$,
$\Delta_2(\beta_{2,1},\beta_{2,2})=62$, and
$\Delta_3(\beta_{3,1},\beta_{3,2},\beta_{3,3})=\beta_{3,3}\Delta_2(\beta_{2,1},\beta_{2,2})
+n_2\Delta_2(\beta_{3,1},\beta_{3,2})=5\cdot 187=935$.

\item"(ii)"  $n_1=3<\beta_{1,1}=4$,
$\Delta_2(\beta_{2,1},\beta_{2,2})=62>n_2n_1\Delta_1(\beta_{1,1})=60$,
and $\Delta_3(\beta_{3,1},\beta_{3,2},\beta_{3,3})=935
>n_3n_2\Delta_2(\beta_{2,1},\beta_{2,2})=930$.

\item"(iii)"$\gcd(n_1,\beta_{1,1})=\gcd(3,4)=1$,
$\gcd(n_2,\Delta_2(\beta_{2,k})^2_{k=1})=\gcd(5,62)=1$ and

{\noindent}$\gcd(n_3,\Delta_3(\beta_{3,k})^3_{k=1})=\gcd(3,935)=1$.

\item"(iv)" $n_{j}\ge 2$, $\beta_{j,1}>0$ for all $j=1,2,3$. Also,
$0\le \beta_{j,k}<n_{k-1}$ for $2\le j\le 3$ and $2\le k\le j$. \ms

{\noindent}Note by (i), (ii), (iii) and (iv) and by Theorem $5.0$
that $g_3$ is irreducible in $\BC\{y,z\}$.
\endroster\ms

Now, it is easy to compute by (1.4.1) in Algorithm 1.4.1 for Theorem
1.4 that the standard Puiseux expansion for $C_3(t)$ such that
\text{$V(g_r)\equiv C_3(t)$ \text{\rm (multi. seq.)}} is given by
$$
\text{$C_3(t):=$}  \left\{\eqalign{ y= &t^{45} \cr z=
&t^{60}+t^{66}+t^{71}. \cr } \right. \tag 1.4.3
$$
because of the following computations {\rm(a)} and {\rm(b)}:
\roster
\item"(a)" $n=n_1n_2n_3=45$, $\alpha_1=\beta_{1,1}n_2n_3=60$, and
$\alpha_2-\alpha_1=\widehat{\Delta}_2n_3=2\cdot 3=6$ implies that
$\alpha_2=66$ because
$\widehat{\Delta}_2=\Delta_2(\beta_{2,1},\beta_{2,2})-n_2n_1\beta_{1,1}=2$.

\item"(b)" $\alpha_3-\alpha_2=\widehat{\Delta}_3=5$ implies that
$\alpha_3=71 $ because
$\widehat{\Delta}_3=\Delta_3(\beta_{3,1},\beta_{3,2},\beta_{3,3})
-n_3n_2\Delta_2(\beta_{2,1},\beta_{2,2})=935-3\cdot 5\cdot 62=5$.
\quad $\square$
\endroster \ms

{\bf \S1.4.B. The second half of The 1st Algorithm(Theorem 1.6)}

\definition{\bf Sublemma 1.5.(Corollary 7.6) for Theorem 1.6}

$\underline{\text{\bf{Assumptions}}}$ Let $A\ge 2$ and $B\ge 2$ be
integers with $\gcd(A,B)=1$. Let $p$ be an integer such that $p>nAB$
for some integer $n\ge 2$.

$\underline{\text{\bf {Conclusions}}}$ We can compute a unique pair
of two integers $s_1$ and $t_1$ such that $p=s_1A+t_1 B$ with $0\le
s_1<B$ and $t_1>A$. \quad $\square$
\enddefinition \ms

\proclaim{Theorem 1.6(Theorem 11.4:Algorithm for finding the unique
element of Family(1) corresponding to any given standard Puiseux
expansion of Family(2))}

$\underline{\text{\bf {Assumptions}}}$ Let the standard Puiseux
expansion of the $r$-type $C_r(t)$ be given by
$$\align
(1.6.1) \qquad \text{$C_r(t):=$} \left\{\eqalign{ y=&t^n, \cr
z=&t^{\alpha_1}+t^{\alpha_2}+\cdots +t^{\alpha_r}, \cr} \right.
\\
\text{where} \quad
 2\le n <\alpha_1<\alpha_2<\cdots <\alpha_r  & \quad  \text{and} \\
 n >d_1>d_2>\cdots
 >d_r=1  \quad   \text{with} & \quad
\text{$d_i=\gcd(n,\alpha_1,\dots,\alpha_i)$, $1\le i\le r$.}\qquad
\qquad
\endalign$$

$\underline{\text{\bf Conclusions}}$ To compute the standard Puiseux
W-poly $g_r$ of the recursive $r$-type with \text{$V(g_r) \buildrel
\text{{\rm multiseq}} \over \sim C_r(t)$} \text{at $0\in \BC^2$} is
to find {\bf explicit algorithm(Algorithm 1.6.2)}, using a finite
number $\frac{r(r-1)}{2}$ of Sublemma 1.5(Corollary 7.6), as soon as
the standard Puiseux W-poly $g_r$ of the recursive $r$-type
satisfies the same kind of properties and notations as in Definition
$1.1$, for notation. \ms

\noindent{\bf (Algorithm 1.6.2 for Theorem 1.6)} To compute an
algorithm for finding one and only one standard Puiseux W-poly
$g_r\in \BC\{y\}[z]$ of the recursive $r$-type such that
\text{$V(g_r) \buildrel \text{{\rm multiseq}} \over \sim C_r(t)$} at
$0\in \BC^2$, we may assume that the above $g_r$ satisfies the same
properties and notations as $g_r$ of Definition $1.1$ does.

To find such an algorithm, using Step(1) and Step(2), it suffices to
compute the family of sets $\{X_j:1\le j\le r\}$, satisfying
following properties:
\roster
\item"(1)(1a)" $X_1=\{n_1,\beta_{1,1}\}$ with $2\le
n_1<\beta_{1,1}$. \ms

\item"{}(1b)" $X_j=\{n_j,\beta_{j,1},\beta_{j,2},\dots,\beta_{j,j}\}$ with
$n_j\ge 2$ for $j=2,\dots,r$, satisfying the six conditions in
Definition $1.1$ and the following equations in {\rm(1.6.2)}:
$$\align
\text{\rm(1.6.2)(i)} \quad &\text{$n=n_1d_1$ and
$\alpha_1=\beta_{1,1}d_1$ with
$d_1=\gcd(n,\alpha_1)$} \\
\text{\rm{\quad(ii)}} \quad &\text{$d_{j-1}=n_jd_{j}$ and
$\alpha_j-\alpha_{j-1}= \widehat{\Delta}_jd_j$ with
$d_j=\gcd(d_{j-1},\alpha_j-\alpha_{j-1})$
for $2\le j\le r$}, \\
&\text{where $\widehat{\Delta}_j
=\Delta_j(\beta_{j,k})^j_{k=1}-n_jn_{j-1}\Delta_{j-1}(\beta_{j-1,k})^{j-1}_{k=1}$.}
\endalign$$
\endroster
Note that all $n_i\ge 2$, and so $n>d_1>\cdots>d_r=1$. It may be
assumed by Remark 1.4.0 that the equations in {\rm(1.6.2)} of
Theorem 1.6 and the equations in {\rm(1.4.1)} are the same. \ms

\noindent$\underline{\text{\bf {Step(1) for Algorithm 1.6.2.}}}$ We
can easily compute a finite set $\{(n_j,\widehat{\Delta}_j)\in N^2:
j=1,2,\dots,r\}$ of unique pairs, each of which satisfies the
following: \roster
\item"(1.1)" Let $d_1 =\gcd(n,\alpha_1)$. It is easy to
compute a unique pair $(n_1, \widehat{\Delta}_1)\in N^2$ such that
$n=n_1d_1$ and $\alpha_1=\widehat{\Delta}_1d_1$ with $\gcd(n_1,
\widehat{\Delta}_1)=1$. Note that $\widehat{\Delta}_1>n_1\ge 2$
because $n>d_1$, and write $\widehat{\Delta}_1=\beta_{1,1}$ for
notation.\ms

\item"(1.2)" Let $d_j =\gcd(d_{j-1},\alpha_j-\alpha_{j-1})$
for $2\le j\le r$. It is easy to compute a unique pair $(n_j,
\widehat{\Delta}_j) \in N^2$ such that $d_{j-1}=n_jd_j$ and
$\alpha_j-\alpha_{j-1}=\widehat{\Delta}_jd_j$ with
$\gcd(d_j,\widehat{\Delta}_j)=1$. Note that
$d_j=\gcd(n,\alpha_1,\alpha_2,\dots,\alpha_j)$. \endroster

\noindent $\underline{\text{\bf {Step(2) for Algorithm 1.6.2.}}}$
Let $\{(n_j,\widehat{\Delta}_j)\in N^2:j=1,2,\dots,r\}$ be already
given by {\rm Step(1)} where $\gcd(n_j,\widehat{\Delta}_j)=1$ for
$1\le j\le r$.

\roster
\item"(2.1)" With $\{n_1, \widehat{\Delta}_1\}$ in $(1.1)$ of
Step (1), let $\Delta_1:N_0\to N_0$ be a function defined by
$$\ \Delta_1(t)=t. $$

We can compute a solution $\beta_{1,1}=\widehat{\Delta}_1$  such
that $\Delta_1(\beta_{1,1})=\beta_{1,1}>n_1\ge 2$. \ms

\item"(2.2)" With $(n_2, \widehat{\Delta}_2)$ in $(1.2)$ of
Step (1), let $\Delta_2:N^2_0\to N_0$ be a function defined by
$$\split
\Delta_2(t_1,t_2) &=t_2\beta_{1,1}+n_1t_1 \quad with \quad
\text{$\beta_{1,1}=\widehat{\Delta}_1$}.
\endsplit$$

Given $p_2=n_2n_1\beta_{1,1}+\widehat{\Delta}_2$, we can compute a
unique pair $(a_2,b_2)$ in $N^2$ such that
$a_2\beta_{1,1}+b_2n_1=p_2$ with $b_2>\beta_{1,1}$ and $0\le
a_2<n_1${\rm(by Sublemma 1.5, because $p_2>2n_1\beta_{1,1}$ and
$\gcd(n_1,\beta_{1,1})=1$)}.

Write $\beta_{2,1}=a_2$ and $\beta_{2,2}=b_2$ with $0\le a_2<n_1$
and $b_2>\beta_{1,1}$, and then
$p_2=\Delta_2(\beta_{2,1},\beta_{2,2})$. \ms

\item"(2.3)" With $\{n_j, \widehat{\Delta}_j\}$ in $(1.2)$ of
Step (1), suppose we have proved that the following are true:

For each $j=2,3,\dots,{\ell}$ with $\ell<r$, assume that $d_j>1$,
and then use the induction assumption on the positive integer {j}.
Given $p_j=
\widehat{\Delta}_j+n_jn_{j-1}\Delta_{j-1}(\beta_{j-1,k})^{j-1}_{k=1}$
for $j=2,3,\dots,\ell$, we may assume by a finite number (j-1) of
use of Sublemma 1.5 that we can compute a unique sequence
$X_j=\{\beta_{j,k}:k=1,2,\dots,j\}\subset N_0$ such that
$\Delta_{j}(\beta_{j,k})^{j}_{k=1}=p_j$ satisfying the six
conditions in Definition 1.1. \ms

\item"(2.3.0)" Then, a computational algorithm for
$\{\beta_{\ell+1,k}:k=1,2,\dots,\ell+1\}\subset N_0$ such that
$\Delta_{\ell+1}(\beta_{{\ell+1},k})^{{\ell+1}}_{k=1}=p_{\ell+1}$
with the six conditions in Definition 1.1 can be easily represented
by {\rm(i)}, {\rm(ii)} and {\rm(iii)}.

We use Sublemma 1.5 $\ell$-times, as follows:

With $\{n_{\ell+1},\widehat{\Delta}_{\ell+1}\}\in N^2$ in $(1,2)$ of
Step (1), let $\Delta_{\ell+1}:N^{\ell+1}_0\to N_0$ be a function
defined by
$$
\align \Delta_{\ell+1}(t_k)^{\ell+1}_{k=1}&
=t_{\ell+1}\Delta_{\ell}(\beta_{{\ell},k})^{{\ell}}_{k=1}
+n_{\ell}\Delta_{\ell}(t_k)^{\ell}_{k=1}.
\endalign
$$

{\bf(i)} Given
$p_{{\ell}+1}=n_{{\ell}+1}n_{{\ell}}\Delta_{{\ell}}(\beta_{{\ell},k})^{{\ell}}_{k=1}
+\widehat{\Delta}_{{\ell}+1}$, we can compute a unique pair
$(a_{{\ell}+1},b_{{\ell}+1})$ in $N^2$ such that
$p_{{\ell}+1}=a_{{\ell}+1}\Delta_{{\ell}}(\beta_{{\ell},k})^{\ell}_{k=1}
+b_{{\ell}+1}n_{\ell}$ with
$b_{{\ell}+1}>\Delta_{{\ell}}(\beta_{{\ell},k})^{\ell}_{k=1}$ and
$0\le a_{{\ell}+1}=\beta_{{\ell}+1,{\ell}+1}<n_{\ell}$(by Sublemma
1.5 once, because $p_{\ell+1}
>2n_{\ell}\Delta_{\ell}(\beta_{\ell,k})^{\ell}_{k=1}$
and $\gcd(n_{\ell},\Delta_{\ell}(\beta_{\ell,k})^{\ell}_{k=1})=1$).

{\bf(ii)} Since
$b_{{\ell}+1}>\Delta_{{\ell}}(\beta_{{\ell},k})^{\ell}_{k=1}\ge
2n_{\ell}n_{\ell-1}\Delta_{\ell-1}(\beta_{\ell-1,k})^{\ell-1}_{k=1}$,
then by induction on the positive integer $\ell<s$, and  by a finite
number $(\ell-1)$ of use of Sublemma 1.5, we can compute a unique
sequence $(\beta_{\ell+1,k})^{\ell}_{k=1}\subset N^{\ell}_0$ such
that
$$\Delta_{\ell}(\beta_{{\ell}+1,k})^{\ell}_{k=1}=b_{{\ell}+1}
\quad \text{with $\beta_{{\ell}+1,1}>0$ and $0\le
\beta_{{\ell}+1,k}<n_{k-1}$ for $2\le k\le {\ell}$}.
$$

{\bf (iii)} By {\rm(i)} and {\rm(ii)}, by a finite number {$\ell$}
of use of Sublemma 1.5, we can compute a unique sequence
$\{\beta_{{\ell+1},k}:k=1,2,\dots,{\ell+1}\}\subset N_0$ such that
$p_{{\ell}+1}=a_{{\ell}+1}\Delta_{{\ell}}(\beta_{{\ell},k})^{\ell}_{k=1}
+b_{{\ell}+1}n_{\ell}
=\beta_{{\ell}+1,{\ell}+1}\Delta_{{\ell}}(\beta_{{\ell},k})^{\ell}_{k=1}
+n_{\ell}\Delta_{\ell}(\beta_{{\ell}+1,k})^{\ell}_{k=1}=
\Delta_{\ell+1}(\beta_{{\ell+1},k})^{\ell+1}_{k=1}$ satisfying the
six conditions in Definition 1.1. \quad $\square$
\endroster
\endproclaim \ms

{\bf Example 1.6.3 for Theorem 1.6(See page 517 of [Bri-Kn]).} Let
the parametrization $C_4(t)$ for the Puiseux expansion be given by
an example in page 517 of [Bri-Kn.].
$$
\text{$C_4(t):=$}  \left\{\eqalign{ y= &t^{100} \cr z=
&t^{250}+t^{375}+t^{410}+t^{417}. \cr } \right. \tag 1.6.3
$$

This is the standard Puiseux expansion because of the  computations
{\rm(i)} and {\rm(ii)}:

{\rm(i)} $n<\alpha_1<\cdots<\alpha_4$ where $n=100$, $\alpha_1=250$,
$\alpha_2=375$, $\alpha_3=410$ and $\alpha_4=417$.

{\rm(ii)} $n=100>d_1=50>d_2=25>d_3=5>d_4=1$ where
$d_1=\gcd(n,\alpha_1)$, $d_2=\gcd(d_1,\alpha_2-\alpha_1)$,
$d_3=\gcd(d_2,\alpha_3-\alpha_2)$ and
$d_4=\gcd(d_3,\alpha_4-\alpha_3)$.

Now, the problem is how to find a one and only one $g_4\in \text{\rm
Family(1)}$ such that \text{$V(g_4) \buildrel \text{{\rm multiseq}}
\over \sim C_4(t)$} at the origin in $\BC^2$. For the solution of
the above problem, by a finite number
$\frac{r(r-1)}{2}=\frac{4\cdot3}{2}$ of use of Sublemma 1.5, it
suffices to follow {\bf (Algorithm 1.6.2 for Theorem 1.6)}.

After the following computation is done, the above polynomial $g_4$
of the recursive $4$-th type is as follows:
$$
(1.6.4) \qquad \qquad \text{$g_1=z^2+y^5$, $g_2=g_1^{2}+y^{10}z$,
$g_3=g_2^{5}+y^{58}g_1$ and $g_4=g_3^{5}+y^{300}zg_1g_2$.} \qquad
\qquad
$$
\ms

\noindent$\underline{\text{\bf {Step(1) for Algorithm 1.6.2.}}}$ We
can compute a finite set $\{(n_j,\widehat{\Delta}_j)\in N^2:
j=1,2,3,4\}$ of pairs, each of which satisfies the following: Recall
that $d_j=\gcd(d_{j-1},\alpha_j-\alpha_{j-1})$ for $1\le j\le 4$
where $d_{0}=n$ and $\alpha_{0}=0$.

{\rm(1)} Let $n=d_1n_1$ and $\alpha_1=d_1\widehat{\Delta}_1$. Then,
$d_1=50$, $n_1=2$ and $\gamma_{11}=\widehat{\Delta}_1=5$.

{\rm(2)} Let $d_1=d_2n_2$ and
$\alpha_2-\alpha_1=d_2\widehat{\Delta}_2$. Then, $d_2=25$, $n_2=2$
and $\widehat{\Delta}_2=5$.

{\rm(3)} Let $d_2=d_3n_3$ and
$\alpha_3-\alpha_2=d_3\widehat{\Delta}_3$. Then, $d_3=5$, $n_3=5$
and $\widehat{\Delta}_3=7$.

{\rm(4)} Let $d_3=d_4n_4$ and
$\alpha_4-\alpha_3=d_4\widehat{\Delta}_4=7$. \ms

\noindent$\underline{\text{\bf {Step(2) for Algorithm 1.6.2.}}}$ Let
$\{(d_j,\widehat{\Delta}_j)\in N^2: j=1,2,3,4\}$ of pairs be given
by {\rm Step (1)} where $\gcd(n_j,\widehat{\Delta}_j)=1$ for $1\le
j\le 4$.

By the Euclidean algorithm in Sublemma 1.5(Corollary $7.6$), for
each $j=2,3,4$, we compute a finite unique sequence
$\{\beta_{j,1},\beta_{j,2},\dots,\beta_{j,j}\}$ with six conditions
in Definition 1.1 for the standard Puiseux W-poly in Family(1), as
follows: \roster
\item"(i)" To compute $g_1$, $n_1=2$ and
$\beta_{1,1}=5$. So, $g_1=z^{2}+y^{5}$. \ms

\item"(ii)" To compute $g_2$, we use Sublemma 1.5 once. Note that
$\Delta_2(\beta_{2,1},\beta_{2,2})=\beta_{1,1}\beta_{2,2}+n_1\beta_{2,1}
=5\beta_{2,2}+2\beta_{2,1}$ where
$\Delta_2(\beta_{2,1},\beta_{2,2})=p_2=n_2n_1\beta_{1,1}+\widehat{\Delta}_2
=2\cdot2\cdot5+5=25$.

Since $\Delta_2(\beta_{2,1},\beta_{2,2})>n_2n_1\beta_{1,1}$, we can
compute a unique solution $\{\beta_{2,1},\beta_{2,2}\}\subset N_0$
such that $\Delta_2(\beta_{2,1},\beta_{2,2})=25$ with
$\beta_{2,2}<n_1=2$. By Sublemma 1.5, $\beta_{2,1}=10$ and
$\beta_{2,2}=1$. So, $g_2=g_1^{2}+y^{10}z^1$. \ms

\item"(iii)" To compute $g_3$, we use Sublemma 1.5 twice. Note that
$\Delta_3(\beta_{3,k})^3_{k=1}
=\beta_{3,3}\Delta_2(\beta_{2,1},\beta_{2,2})+n_2\Delta_2(\beta_{3,1},\beta_{3,2})
=25\beta_{3,3}+2\Delta_2(\beta_{3,1},\beta_{3,2})>n_3n_2\Delta_2(\beta_{2,1},\beta_{2,2})$
where $\Delta_3(\beta_{3,k})^3_{k=1}=p_3=
n_3n_2\Delta_2(\beta_{2,,1},\beta_{2,2})+\widehat{\Delta}_3=5\cdot
2\cdot 25+7=257$. \ms

\item"(iii-a)" To compute $\beta_{3,3}$, since
$\gcd(n_2,\Delta_2(\beta_{2,1},\beta_{2,2}))=1$, by Sublemma 1.5 we
can compute a unique pair $(a_3,b_3)=(1,116)$ such that
$25a_3+2b_3=257$ where $0\le 1=a_3=\beta_{3,3}<n_2$ and
$b_3=116>\Delta_2(\beta_{2,1},\beta_{2,2})$. Then, $\beta_{3,3}=1$
and $\Delta_2(\beta_{3,1},\beta_{3,2})=b_3=116$. \ms

\item"(iii-b)" Since
$116=\Delta_2(\beta_{3,1},\beta_{3,2})=2\beta_{3,1}+5\beta_{3,2}
>\Delta_2(\beta_{2,1},\beta_{2,2})$ by (iii-a),
then by Sublemma 1.5 we can compute a unique solution
$\{\beta_{3,1},\beta_{3,2}\}\subset N_0$ such that
$\Delta_2(\beta_{3,1},\beta_{3,2})=116$ with $\beta_{3,2}=0<n_1=2$
and $\beta_{3,1}=58$.

So, by (iii-a) and (iii-b), $\beta_{3,1}=58$, $\beta_{3,2}=0<n_1$
and $\beta_{3,3}=1<n_2$. So, $g_3=g_2^{5}+y^{58}g_1^{1}$. \ms

\item"(iv)" To compute $g_4$, we use Sublemma 1.5 three times. It is clear
that $\Delta_4(\beta_{4,k})^4_{k=1}
=\beta_{4,4}\Delta_3(\beta_{3,k})^3_{k=1}
+n_3\Delta_3(\beta_{4,k})^3_{k=1}= 257\beta_{4,4}
+5\Delta_3(\beta_{4,1},\beta_{4,2},\beta_{4,3})
>n_4n_3\Delta_3(\beta_{3,k})^3_{k=1}$
where $\Delta_4(\beta_{4,k})^4_{k=1}=p_4=
n_4n_3\Delta_3(\beta_{3,k})^3_{k=1}+\widehat{\Delta}_3=5\cdot 5\cdot
257+7$. \ms

\item"(iv-a)" To compute $\beta_{4,4}$, since
$\gcd(n_3,\Delta_3(\beta_{3,k})^3_{k=1})=1$, by Sublemma 1.5 once we
can compute a unique pair $(a_4,b_4)=(1,1235)$ such that $257a_4
+5b_4=5\cdot 5\cdot 257+7$ where $0\le 1=a_4=\beta_{4,4}<n_3$ and
$b_4=1235>\Delta_3(\beta_{3,k})^3_{k=1}$. \ms

\item"(iv-b)" To compute  a unique solution $\{\beta_{4,k};k=1,2,3\}$
 such that $b_4=\Delta_3(\beta_{4,k})^3_{k=1}>\Delta_3(\beta_{3,k})^3_{k=1}$
with  $\beta_{4,1}>0$, $\beta_{4,2}<n_1$ and $\beta_{4,3}<n_2$, then
using Sublemma 1.5 twice and the same method as we have used in
(iii),  we can compute $\beta_{4,1}=300$, $\beta_{4,2}=1<n_1$ and
$\beta_{4,3}=1<n_2$.

So, by (iv-a) and (iv-b), $\beta_{4,1}=300$, $\beta_{4,2}=1$,
$\beta_{4,3}=1<n_2$ and $\beta_{4,4}=1$.
\endroster \ms

Thus, by (i), (ii), (iii) and (iv), the above polynomial
$f(y,z)=g_4$ of the recursive type is as follows:

$f=g_4=g_3^{5}+y^{300}z^1g_1^{1}g_2^1$ where $g_1=z^2+y^5$,
$g_2=g_1^{2}+y^{10}z^1$ and $g_3=g_2^{5}+y^{58}g_1^{1}$.

Therefore,  $f$ is the unique standard Puiseux W-poly of the
recursive type because of \text{\rm The 6-th
${\text{\rm{Cond}}}^{\text{{\rm(0)}}}$}, and so
$f$ is the W-poly of the recursive type in $z$.
 \quad $\square$ \bs

{\bf \S1.5. The Weierstrass division algorithm for the
W-polys in preparation for finding The 2nd Algorithm and The 3rd Algorithm} \ms

As in Definition 15.0 with Notation 15.0.1, the Weierstrass
preparation theorem and the Weierstrass division theorem can be
written by The WPT and The WDT respectively, for brevity of
notation. In order to succeed in the computations of the 2nd and the
3rd algorithms in $\S1.6$, in this section it is very important to
say without any other proofs that The Weierstrass division algorithm for the
W-polys(Theorem 1.8(Theorem 15.4) with two sublemmas can have an
important role of the 2nd and the 3rd algorithms in $\S1.6$, which
will be shown in $\S 16$, later.

\proclaim{Theorem 1.7(Theorem 15.2: The WDT for the W-polys)}

$\underline{\text{\bf{Assumptions}}}$ Let $h\in {}_{n-1}\CO[z_n]$ be
a $W$-poly of degree $\nu
>0$ in $z_n$. Let $f\in {}_{n-1}\CO[z_n]$ be a $W$-poly of degree $\mu\ge \nu$
in $z_n$, and $\ell$ be a positive
integer with $\ell\nu \le \mu <(\ell+1)\nu$. \ms

$\underline{\text{\bf Conclusions}}$

{\rm(1)} Then, $f$ can be written uniquely in the form
$$
f=\sum^{\ell}_{i=0} r_ih^i \quad \text{with \quad $h^{0}=1$}, \tag
1.7.1
$$
where if $\mu \ge \ell\nu$ then for $i=0,1,\dots,{\ell}-1$, each
$r_i\in {}_{n-1}\CO[z_n]$ is a polynomial of degree $<\nu$ in $z_n$
with $r_i(0,\dots, 0, z_n)$ identically zero, and if $\mu =\ell\nu$
then $r_{\ell}$ is equal to one and if $\mu
>\ell\nu$ then $r_{\ell}\in {}_{n-1}\CO[z_n]$ is a $W$-poly
of degree $\mu -\ell\nu <\nu$ in $z_n$. \ms

{\rm(2)} In addition, suppose $h\in {}_{n-1}\CO[z_n]$ has a
multiplicity $\nu >0$ at $0\in \BC^n$ and $f\in {}_{n-1}\CO[z_n]$
has a multiplicity $\mu \ge \nu$ at $0\in \BC^n$. Then, the above
representation $f=\sum^{\ell}_{i=0}r_ih^i$ of {\rm (1.7.1)}
satisfies the property such that for $i=0,1,\dots,{\ell}-1$, each
$r_i$ has a multiplicity $\ge \mu -i\nu$ at $0\in \BC^n$ and such
that if $\mu = \ell\nu$ then $r_\ell$ has a zero multiplicity at
$0\in \BC^n$. $\square$
\endproclaim \ms

\proclaim{Corollary 1.7.1(Theorem 15.2: The WDT for the W-polys)}

$\underline{\text{\bf{Assumptions}}}$ Let $f=z^n+\sum^{n-1}_{i=0}
a_iy^{\a_i}z^i$ be a $W$-poly of degree $n\ge 2$ in $z$ where for
$0\le i\le n-1$, each $a_i=a_i(y)$ is a unit in $\C\{y,z\}$ and 
the $\a_i$ a positive integer, if exists. In particular,
let $h=z+\f{a_{n-1}y^{\a_{n-1}}}n$ be a $W$-poly of degree $1$ in $z$ where 
$a_{n-1}=a_{n-1}(y)$ is a unit in $\C\{y,z\}$ and the $\a_{n-1}$ is 
a positive integer. \ms 

$\underline{\text{\bf Conclusions}}$ Instead of a nonsingular change of coordintes,
using an equation in (1.7.1) where $n={\mu}={\ell}{\nu}=1{d}$ with ${\ell}=1$,
the above $f$ can be rewritten uniquely in the form 
$$\align
(1.7.2) \qquad & f=h^d+\sum^{d-2}_{i=0}r_ih^i\quad {and} \quad r_{d-1}=0 \quad 
with \quad h^{0}=1 \ and \ r=r_0, \\
&  \text{$r_i\in \C\{y\}[z]$ is a polynomial of degree $<{\ell}$ in $z$
for $0\le i\le d-2$, } \qquad \qquad \qquad \qquad
\endalign$$
where for $0\le i\le d-2$, each
$r_i=r_i(y,z)\in \C\{y\}[z]$ is a polynomial of degree $<\ell$ in $z$
with $r_i(0,z)$ identically zero, if exists.
\endproclaim \ms

\proclaim{Theorem 1.8(Theorem 15.4:The Weierstrass division algorithm for the
W-polys)}

$\underline{\text{\bf {Assumptions}}}$ \quad Let $f\in \BC\{y\}[z]$
be an arbitrary $W$-poly of degree $n\ge 2$ in $z$. Without loss of
generality, we may assume by Theorem 1.7 that $f$ satisfies the following form:
$$\align
 f=z^n+\sum^{n-2}_{i=0} a_iy^{\a_i}z^i, \tag 1.8.1 
\endalign$$
where for $0\le i\le n-2$, each $a_i=a_i(y)$ is a unit in ${}_2\CO_0$, if exists, 
and the $\a_i$ are positive integers. Assume that 
$f(y,z)$ may not be irreducible in $\BC\{y,z\}$, and note that $a_{n-1}$ is identically
zero. Write $n=\Pi^{\ell}_{k=1}n_k$ with
positive integers $n_k\ge 2$ for all $k$ where the $n_k$ may not be
the factorization of prime numbers. 

In addition, assume that we have the following:

\noindent {\rm(1.8.2)} \qquad \quad   \quad $2\le n\le\a_0$.  \bs

$\underline{\text{\bf {Conclusions}}}$

$\underline{\text{Given a sequence of positive integers, $\{n_k\ge 2: k=1,2,\dots,{\ell}\}$,}}$  we can compute a unique sequence
of $W$-polys in $z$, $\{f_k:k=1,2,\dots,{\ell}\}$
such that $f_k\in \BC\{y\}[z]$ is a $W$-poly in $z$ with coefficients in $\C\{y\}$,
satisfying the fact, denoted by
$\underline{\text{\bf \rm Fact[I]}}$:
Let $f_{-1}=y$, $f_0=z$. \ms

$\underline{\text{\rm Fact[I]}}$

For each $k=1,2,\dots,{\ell}$ and $\ell\le r$, $f_k\in \BC\{y\}[z]$ can be
uniquely written in the form
$$\align
\text{\rm (1.8.3)} \qquad\qquad f_k=f^{n_k}_{k-1}+\sum^{n_k-2}_{i=0} R_{k,i}f^i_{k-1}\in 
\C\{f_{-1},f_0,\dots,f_{k-2}\}[f_{k-1}]
\quad \text{with $f=f_{\ell}$} \qquad \qquad\\
\endalign$$
satisfying the following:
\roster
\item "(1)" For each $k$,
$f_k\in \C\{y\}[z]$ is a $W$-poly of degree
$\Pi^k_{t=1}n_t$ in $z$ with coefficients in $\C\{y\}$.

\item "{}(1a)" Let $k$ and $i$ be fixed with $1\le k\le
{\ell}$ and $0\le i\le n_k-2$, and if exists, then $R_{1,i}=R_{1,i}(y)$ is a
nonunit in $\C\{y\}$ and for each $k\ge 2$, $R_{k,i}=R_{k,i}(y,z)\in
\C\{y\}[z]$ is a polynomial of degree $<\Pi^{k-1}_{t=1}n_t$ in $z$
with $R_{k,i}(0,z)=0$.\ms

\item "(2)" For each fixed $k=1,2,\dots,{\ell}$, $f_k=f_k(f_{-1},f_0,\dots,
f_{k-1})\in \C\{f_{-1},f_0,\dots,f_{k-2}\}[f_{k-1}]$ of {\rm(1.8.3)} is a W-poly in $f_{k-1}$
with coefficients in $\C\{f_{-1},f_0,\dots,f_{k-2}\}$,
considering $f_{-1},f_0,\dots,f_{k-1}$ as independent complex $({k}+1)$-variables
at $0\in \C^{{k}+1}$ where $f_{-1}=y$ and $f_0=z$ if necessary,
with the following property {\rm (2a)}:

\item"{}(2a)"  Let $k$ and $i$ be fixed with $1\le k\le {\ell}$
and $0\le i\le n_k-2$, and for any nonzero monomial
$\Pi^k_{t=1}f^{\de_t}_{t-2}$ in $R_{k,i}=R_{k,i}(f_{-1},f_0,\dots,
f_{k-2})\in \C\{f_{-1},f_0,\dots,f_{k-2}\}$, $\de_1>0$ and
$0\le \de_t<n_{t-1}$ for $t=1,2,\dots, k$. 
\endroster 
\endproclaim \ms

\definition{Remark 1.8.0 for the proof of Theorem 1.8}
If ${\ell}=1$, there is nothing to prove for $f=f_1$. 
If ${\ell}=2$ and $f=f_2$ in Assumptions and Conclusions of Theorem 1.8, for brevity 
this statement will be called Theorem 1.8.1 in the process of induction proof 
for Theorem 1.8. For this case, it suffice to prove Theorem 1.8.1 with Sublemma 1.9 
and Sublemma 1.10. The proof of Theorem 1.8.1  will be done in the beginning 
of the proof of Theorem 15.4 in \S 15 without difficulty. The remaining proof of Theorem 1.8 can be finished in the proof of Theorem 15.4 in \S 15. 
\enddefinition \ms

\proclaim{Theorem 1.8.1(Theorem 15.4: Weierstrass division algorithm for the
W-polys)}

$\underline{\text{\bf Assumptions}}$ Let $f=z^n+\sum^{n-2}_{i=0}
a_iy^{\a_i}z^i$ be a $W$-poly of degree $n\ge 2$ in $z$ where for
$0\le i\le n-2$, each $a_i=a_i(y)$ is a unit in ${}_2\CO_0$ if
exists and the $\a_i$ are positive integers. Assume that $f$ may not
be irreducible in ${}_2\CO_0$, and note that $a_{n-1}$ is
identically zero for convenience. Write $n=\Pi^{\ell}_{k=1}n_k$ with
positive integers $n_k\ge 2$ for all $k$ where the $n_k$ may not be
the factorization of prime numbers. \ms

$\underline{\text{\bf Conclusions}}$ We can compute a unique
$W$-poly in $z$, $f_1\in \C \{y\}[z]$, satisfying the following
notations and properties: Let $f_{-1}=y$ and $f_0=z$.
$$\align
f_1&=f^{n_1}_{0}+\sum^{n_1-2}_{i=0} R_{1,i}f^i_{0}
\quad \text{and} \tag 1.8.1.1 \\
f&=f^{d_{2}}_1 +\sum^{d_{2}-2}_{i=0} S_{2,i}f^i_1,
\endalign$$
such that \quad{\rm(i)} \quad $n=d_{2}n_1$ with $n=d_1$,

\qquad \quad \quad {\rm(ii)} \quad $f_1=f_1(y,z)\in \C\{y\}[z]$ is a
$W$-poly of degree $n_1$ in $z$,

\qquad \quad \quad {\rm(iii)} \quad $f\in \C\{y\}[z,f_1]\subseteq
\C\{y,z\}[f_1]$ is a $W$-poly of degree $d_{2}$ in $f_1$, \ms

{\noindent}considering $f_{-1}=y,f_0=z,f_1$ as independent complex
$(3)$-variables at the origin in $\C^{3}$, with two properties {\rm
(1)} and {\rm (2)}: \ms

\roster
\item"(1)(1a)" Let $i$ be fixed with $0\le i\le n_1-2$. If exists, then
$R_{1,i}=R_{1,i}(y)$ is a nonunit in $\C\{y\}$.

\item"(1b)" Let $i$ be fixed with $0\le i\le d_{2}-2$. Then
$S_{2,i}=S_{2,i}(y,z)\in \C\{y\}[z]$ is a polynomial of degree
$<n_1$ in z and $S_{2,i}(0,z)=0$. \ms

\item"(2)(2a)" Let $i$ be fixed with $0\le i\le n_1-2$.
For any nonzero monomial $y^{\de_1}$ in $R_{1,i}=R_{1,i}(y)\in
\C\{y\}$, $\de_1>0$.

\item"(2b)" Let $i$ be fixed with $0\le i\le d_{2}-2$. For
any nonzero monomial $y^{\de_1}z^{\delta_2}$ in
$S_{2,i}=S_{2,i}(y,z)\in \C\{y\}[z]$, $\delta_1>0$ and
$\delta_2<n_{1}$. \quad $\square$
\endroster
\endproclaim \ms

\noindent{\bf Remark 1.8.1.1} {\rm(a)} It is most interesting and important 
in this book that explicit algorithm for finding a construction of (1.8.1.1) 
in the conclusion of Theorem 1.8.1 can be completely computed from an
equation in (1.9.1) of Sublemma 1.9 and an equation in (1.10.1) of
Sublemma 1.10, using an equation in (1.7.1) of Theorem 1.7.

{\rm(b)} Note that Theorem 1.8(Theorem 15.4) is a
generalization of Theorem 1.8.1 It will be proved by Sublemma
15.5.$\alpha$ and Sublemma 15.5.$\beta$ of $\S15$ that Theorem 15.4 is
true. So, it can be easily proved by Sublemma $1.9$(Sublemma
15.4.$\alpha$) and Sublemma $1.10$(Sublemma 15.5.$\beta$) that Theorem 1.8.1
is true.
$\square$ \ms

\proclaim{Sublemma 1.9 for Theorem 1.8.1(Sublemma 15.5.$\alpha$ for
Theorem 15.4)}

$\underline{\text{\bf Assumptions}}$ Suppose that the same
properties and notations as in the assumption of Theorem 1.8.1 hold.
\ms

$\underline{\text{\bf Conclusions}}$ We show that $h_{1,1}$ and $f$
can be constructed as follows:
$$
\cases
h_{1,1} &=f^{n_{1}}_0 +\sum^{n_{1}-2}_{i=0} R^{(1)}_{1,i}f^i_0, \\
f &=h^{d_{2}}_{1,1}+\sum^{d_{2}-1}_{i=0} T_{2,i}h^i_{1,1},
\endcases \tag 1.9.1
$$
where $h_{1,1}\in \C\{y\}[z]$ is a $W$-poly of degree $n_1$ in z and
$n=n_1d_2$, satisfying the following facts, {\rm Fact(A), Fact(B),
Fact(C), Fact(D) and Fact(E)}.

\roster

\item"$\underline{\text{\rm Fact(A)}}$" For each
$i=0,1,\dots,n_{1}-2$, $R^{(1)}_{1,i}=R_{1,i}^{(1)}(y)\in \C\{y\}$
with $R^{(1)}_{1,i}(0)=0$, if exists.

\item"$\underline{\text{\rm Fact(B)}}$" For each
$i=0,1,\dots,n_{1}-2$, and for any nonzero monomial $y^{\de_1}$ in
$R^{(1)}_{1,i}\in \C\{y\}$, $\de_1>0$.

\item"$\underline{\text{\rm Fact(C)}}$" For each $i=0,1,\dots,
d_{2}-1$, $T_{2,i}=T_{2,i}(y,z)=\sum a_{p,q}y^pz^q$ with a nonzero
constant $a_{p,q}$ such that $p>0$ and $q<n_1$ and that
$T_{2,i}(0,z)=0$.

Moreover, considering $y,z,f_1=h_{1,1}$ as independent complex
$(3)$-variables at the origin in $\C^{3}$, then, $T_{2,i}\in
\C\{y\}[z]\subseteq \C\{y,z\}$ satisfies two facts {\rm Fact(D) and
Fact(E)}.

\item"$\underline{\text{\rm Fact(D)}}$"  For each $i=0,1,\dots, d_{2}-1$,
and for any nonzero monomial $\Pi^{2}_{t=1}f^{\g_t}_{t-2}$ in
$T_{2,i}$, $\g_1>0$ and $\g_2<n_1$.

\item"$\underline{\text{\rm Fact(E)}}$" In particular, if $i=d_{2}-1$ for
$T_{2,i}$ of {\rm Fact(D)}, then $\g_{2}\le n_{1}-2$. $\square$
\endroster
\endproclaim \ms

\proclaim{Sublemma 1.10 for Theorem 1.8.1(Sublemma 15.5.$\beta$ for Theorem
15.4)}

$\underline{\text{\bf Assumptions}}$ Suppose that the same
properties and notations as in the assumption of Theorem 1.8.1 hold.
By the same way as we have seen in {\rm (1.9.1)} of the conclusion
of Sublemma $1.9$, we may assume that $(h_{1,1},f)$ can be
constructed as follows:
$$
\cases
h_{1,1} &=f^{n_{1}}_0 +\sum^{n_{1}-2}_{i=0} R^{(1)}_{1,i}f^i_0, \\
f &=h^{d_{2}}_{1,1}+\sum^{d_{2}-1}_{i=0} T_{2,i}h^i_{1,1},
\endcases \tag 1.10.1
$$
satisfying the facts, denoted by {\rm Fact(A), Fact(B), Fact(C),
Fact(D) and Fact(E)}. For brevity of notation, let $h_1=h_{1,1}$,
$R^{(1)}_i =R^{(1)}_{1,i}$ for $0\le i\le n_{1}-2$ and
$T^{(1)}_i=T^{(1)}_{2,i}=T_{2,i}$ for $0\le i\le d_{2}-1$,
respectively. \ms

$\underline{\text{\bf Conclusions}}$ Then, $(f_{1},f)$ for $f$ in
the conclusion of Theorem 1.8.1 can be constructed as follows:

{\bf Case[I]:} If $T^{(1)}_{d_{2}-1}$ in $(h_1,f)$ is zero, let
$f_{1}=h_1$, $R_{1,i}=R^{(1)}_i$ for $0\le i\le n_{1}-2$ and
$S_{2,i}=T^{(1)}_i$ for $0\le i\le d_{2}-2$, respectively. Then, the
construction of $(f_{1},f)$ has been already finished. \ms

{\bf Case[II]:} If $T^{(1)}_{d_{2}-1}$ is not zero, for finding such
a construction of $(f_{1},f)$, it suffices to follow two steps, {\rm
Step(1)} and {\rm Step(2)}.

$\underline{\text{\bf Step(1) for Case[II]}}$ Then, there is a
sequence of pairs, $H=\{(h_p,f):p=1,2,\dots\}$, each pair of which
can be constructed with five properties, called {\rm Property(1)},
{\rm Property(2)}, {\rm Property(3)}, {\rm Property(4)} and {\rm
Property(5)}, as follows:
$$
(1.10.2)(1.10.2.1)\quad \quad  \cases h_{1} &=f^{n_{1}}_0
+\sum^{n_{1}-2}_{i=0}R^{(1)}_if^i_0
\text{\quad with $R^{(1)}_i=R^{(1)}_{1,i}$ in $(1.10.1)$,}\\
f &=h^{d_{2}}_{1}+\sum^{d_{2}-1}_{i=0} T^{(1)}_ih^i_{1} \text{\quad
with $T^{(1)}_i=T^{(1)}_{2,i}$ in $(1.10.1)$,} \qquad \qquad
\endcases $$
$$\cases h_{2} &=h_1 +\f 1{d_{2}}T^{(1)}_{d_{2}-1}=f^{n_{1}}_0
+\sum^{n_{1}-2}_{i=0}R^{(2)}_if^i_0, \qquad \qquad \quad\\
f &=h^{d_{2}}_{2}+\sum^{d_{2}-1}_{i=0} T^{(2)}_ih^i_{2},
\endcases \tag 1.10.2.2
$$
$$\cases h_{3} &=h_2 +\f 1{d_{2}}T^{(2)}_{d_{2}-1}=f^{n_{1}}_0
+\sum^{n_{1}-2}_{i=0}R^{(3)}_if^i_0, \qquad \qquad \quad\\
f &=h^{d_{2}}_{3}+\sum^{d_{2}-1}_{i=0} T^{(3)}_ih^i_{3},
\endcases \tag 1.10.2.3
$$
$$ {\dots,}\qquad \qquad \qquad \qquad \qquad \qquad
\qquad \qquad \quad  $$

satisfying the following properties and notations: \ms

$\underline{\text{\rm Property(1)}}$ Let $p$ be fixed with $p\ge 1$.
For each $i=0,1,\dots, d_{2}-2$, $R^{(p+1)}_i=R^{(p+1)}_i(y)\in
\C\{y\}$ with $R^{(p+1)}_{i}(0)=0$, if exists. \ms

$\underline{\text{\rm Property(2)}}$ Let $p$ be fixed with $p\ge 1$.
For each $i=0,1,\dots, d_{2}-1$, $T^{(p+1)}_i=T^{(p+1)}_i(y,z)=\sum
a^{(p+1)}_{\alpha,\beta}y^{\alpha}z^{\beta}$ with a nonzero constant
$a^{(p+1)}_{\alpha,\beta}$ such that $\alpha>0$ and $0\le\beta<n_1$
and that $T^{(p+1)}_{i}(0,z)=0$.\ms

Consider $f_{-1}=y,f_0=z$ as independent complex $(2)$-variables at
the origin in $\C^{2}$.

$\underline{\text{\rm Property(3)}}$ Let $p$ and $i$ be fixed with
$p\ge 1$ and $0\le i\le n_{1}-2$. For any nonzero monomial
$f^{\de_1}_{-1}$ in $R^{(p+1)}_{i}=R^{(p+1)}_{i}(y)\in \C\{y\}$,
$\de_1>0$.

$\underline{\text{\rm Property(4)}}$ Let $p$ and $i$ be fixed with
$p\ge 1$ and $0\le i\le d_{2}-1$. For any nonzero monomial
$f^{\de_1}_{-1}f^{\de_2}_{0}$ in
$T^{(p+1)}_i=T_i^{(p+1)}(f_{-1},f_{0})\in \C\{f_{-1},f_0\}$,
$\de_1>0$ and $\de_2<n_{1}$.

$\underline{\text{\rm Property(5)}}$ In particular, if $i=d_{2}-1$
for $T^{(p+1)}_i$ of {\rm Property(4)}, then $\de_{2}\le n_{1}-2$.
\ms

$\underline{\text{\bf Step(2) for Case[II]}}$ By {\rm Step(1)} for
{\rm Case[II]}, there is a pair $(h_{\nu+1},f)\in H$ which satisfies
the following property:

$\underline{\text{\rm Property(6)}}$ There is an integer $\nu\le
\f{n_{1}+1}2$ such that $T^{(p)}_{d_{2}-1}\ne 0$ for
$p=1,2,\dots,\nu$ and
$T^{(\nu+1)}_{d_{2}-1}=T^{(\nu+2)}_{d_{2}-1}=\cdots =0$. That is,
$(h_{\nu},f)\not= (f_{1},f)$ and $(h_{\nu+1},f)= (f_{1},f)$ for an
integer $\nu \le \f{n_{1}+1}2$. $\square$
\endproclaim

{\bf Remark 1.10.1.} {\rm(a)} It is clear by Sublemma $1.9$ that
{\rm Property(1)}, {\rm Property(2)}, {\rm Property(3)}, {\rm
Property(4)} and {\rm Property(5)} for \text{\rm($h_1$,f)} in
{\rm(1.10.2.1)} are equivalent to {\rm Fact(A)}, {\rm Fact(C)}, {\rm
Fact(B)}, {\rm Fact(D)} and {\rm Fact(E)} for
\text{\rm($h_{1,1}$,f)} in {\rm(1.9.1)}, respectively.

{\rm(b)} It is clear by Sublemma $15.5.\alpha$ that $(h_1,f)$ of
Sublemma 1.9 was already constructed with five properties. $\square$
\bs

\vfill \pagebreak

{\bf \S1.6. Two fundamental lemmas for the representation of the
local defining equations of irreducible plane curve singularities}
\ms

In order to succeed in the computations of the 2nd and the 3rd
algorithms in $\S1.7$, $\S1.8$ and $\S1.9$, it is very important to
say without any other proofs that we can write two fundamental
lemmas for the representation of the local defining equations of
irreducible plane curve singularities, that is, Lemma 1.11 and Lemma
1.12 in this section.

\proclaim{Lemma 1.11(The fundamental lemma for the representation of
the local defining equations of irreducible plane curve
singularities)}

$\underline{\text{\bf Assumptions}}$ Let
$f=f(y,z)=b_nz^n+b_0y^{\beta_0}+\sum^{n-1}_{i=1}b_iy^{\beta_i}z^i$
be in $\BC\{y,z\}$ where for $0\le i\le n$, each $b_i=b_i(y,z)$ is a
unit in $\BC\{y,z\}$ if exists and the $\beta_i$ are positive
integers. Let $m$ be the multiplicity of $f$ at $0\in \BC^2$ with
$n\ge 2$ and ${\beta_0}\ge 2$. Let $d=\gcd(n,\beta_0)$, and write
$n=n_1d$ and $\beta_0=\beta_{1,0,1}d$ with
$gcd(n_1,\beta_{1,0,1})=1$. Note that $d$ may be equal to $n$.

In particular, if $b_i(y,z)=b_i(y)$ for all $i$  and $b_n=1$ then
$f(y,z)$ is called a $W$-poly in $z$. \ms

$\underline{\text{\bf Conclusions}}$

{\bf Fact[I]:}   If $f$ is irreducible in ${}_2\CO_0$ then $b_n$ and
$b_0$ are units in $\BC\{y,z\}$, and $m=n$ or $\beta_0$. So, $f$
must satisfy the following necessary condition:
$$
\f{\beta_i}{n-i}\ge \f{\beta_0}{n} \quad \text{for} \quad 0\le i\le
n-1. \tag 1.11.1
$$

For $f$ of $(1.11.1)$, it suffices to consider two cases:

{\rm Case(A)} $\gcd(n,\beta_0)=1$, and {\rm Case(B)}
$\gcd(n,\beta_0)>1$. \ms

{\rm Case(A)} Let $\gcd(n,\beta_0)=1$. Then, $f$ is irreducible in
${}_2\CO_0$ if and only if $(1.11.1)$ holds. In this case, $f\in$the
type$[1]$ in the sense of Definition 2.5. \ms

{\rm Case(B)} Let $\gcd(n,\beta_0)>1$. If $f$ is irreducible in
${}_2\CO_0$, then $f\in$the type$[\ell]$ with ${\ell}\ge 1$ in the
sense of Definition 2.5. \ms

{\bf Fact[II]:} If $f$ is irreducible in ${}_2\CO_0$, f can be
represented as follows:
$$\align
 g_1&=z^{n_1}+ \xi y^{k_1} \quad \text{with} \quad
k_1=\beta_{1,0,1},  \tag 1.11.2\\
f &=A\cdot{g_1}^d +\sum_{\alpha,\beta\ge 0}
c_{\alpha,\beta}y^{\alpha}z^{\beta} \quad \text{with} \quad
n_1\alpha+k_1\beta>n_1k_1d, \\
\endalign$$
where the $c_{\alpha,\beta}$ are nonzero complex numbers for some
nonnegative integers $\alpha$ and $\beta$ such that
$n_1\alpha+\beta_{1,0,1}\beta>n_1\beta_{1,0,1}d$, satisfying the
following properties : \roster
\item "(i)" $A$ and $\xi$ are the unique nonzero complex
numbers such that $A=b_n(0,0)\ne 0$, $dA\xi=b_{n-n_1}(0,0)\ne 0$,
and \quad ${{d}\choose{i}}A\xi^i=b_{n-in_1}(0,0)$ for $1\le i\le d$.

\item "(ii)" $\frac{\beta_i}{n-i}\ge \frac{\beta_{0}}{n}=\frac{\beta_{1,0,1}}{n_1}$ for
$0\le i\le n-1$. \quad $\square$
\endroster
\endproclaim \ms

The proofs of Lemma 1.11 and Lemma 1.11.1 will be done by Theorem $3.2$, Theorem
$3.4$ and Theorem 3.6. \ms

\noindent{\bf Remark 1.11.0.} {\rm(a)} Let $d=\gcd(n,\beta_{0})=1$.
If $f$ is irreducible in ${}_2\CO_0$, $f\buildrel \text{{\rm resol
}} \over \sim z^n+y^{\beta_0}$ at $0\in\BC^2$ in the sense of
Definition 2.4.

{\rm(b)} If $f$ is defined by $f=(z^2+y^3)^2+y^2z^4$ satisfying an
equation in {\rm(1.11.2)}, note that $f$ is not irreducible in
$\BC\{y,z\}$, satisfying an equation in {\rm(1.11.1)}. \quad
$\square$ \ms

\noindent{\bf Lemma 1.11.1} \quad $\underline{\text{\bf
Assumptions}}$ Let $V(f)=\{(y,z): f(y,z)=0\}$ be an analytic variety
at $(0,0)$ in $\BC^2$ with isolated singularity at the origin, which
is written in the form,
$$\align
 g_1&=z^{n_1}+ \xi y^{k_1} \quad \text{with} \quad
k_1=\beta_{1,0,1},  \tag 1.11.3\\
f &=A\cdot{g_1}^d +\sum_{\alpha,\beta\ge 0}
c_{\alpha,\beta}y^{\alpha}z^{\beta} \quad \text{with} \quad
n_1\alpha+k_1\beta>n_1k_1d, \\
\endalign$$
where $A$ and $\xi$ are nonzero complex numbers, and the
$c_{\alpha,\beta}$ are nonzero complex numbers for some nonnegative
integers $\alpha$ and $\beta$ such that $n_1\alpha+k_1\beta>n_1k_1d$
and $2\le n_1<k_1$ and $\gcd(n_1,k_1)=1$. Note that $f$ may not be
irreducible in $\BC\{y,z\}$. \ms

$\underline{\text{\bf Conclusions}}$ We may assume that $g_1$ of
(1.11.3) can be identified with either $g_1$ of (3.6.1) in Theorem
3.6 or $g_1$ of (5.4.0) in Sublemma 5.4. Following the same
properties and notations as in either Theorem 3.6 or Sublemma 5.4,
let $\tau_m=\pi_1\circ\pi_2\circ\cdots\circ\pi_m:M^{(m)}\to\BC^2$ be
the compositions of a finite number $m$ of successive blow-ups
$\pi_i$ which is needed to get the standard resolution of the
singular point of $V(g_1)$. Since $V(g_1)$ has the singular point at
the origin, along $v=0$ $\tau_m:M^{(m)}\to\BC^2$ as a composition of
analytic mappings and $(f\circ\tau_m)_{total}$ can be rewritten in
the following form: Note that $2\le j\le r$.
$$\align
(1.11.4) \qquad \qquad  \tau_m(v,u)&=(y,z)=(v^{n_1}u^a,v^{k_1}u^b),  \\
(g_1\circ\tau_m)_{proper}&=(1+{\xi}u), \\
(f\circ\tau_m)_{total}&=(f\circ\tau_m)(v,u)
=v^{e_{m}}u^{\rho_{m}}(f\circ\tau_m)_{proper}\quad
\text{with  $g_j=f$}, \qquad \qquad \\
(f\circ\tau_m)_{proper}&=A(1+{\xi}u)^{d}+\sum_{\alpha,\beta\ge
0}c_{\alpha,\beta}v^{n_1\alpha+k_1\beta-n_1k_1d}
u^{a\alpha+b\beta-bn_1d}, \qquad \qquad \qquad \qquad\\
\endalign
$$
where \roster \item "(i)" $a$ and $b$ are some nonnegative integers
such that $ak_1-bn_1=1$,

\item "(ii)" $e_{m}=n_1k_1d$ and $\rho_{m}=bn_1d$
and $\rho_{\alpha,\beta}=a\alpha+b\beta-bn_1d\ge 0$,

\item "(iii)" $E_m=\{v=0\}$ is defined by the $m-th$ exceptional
curve of the first kind,

\item "(iv)" $V^{(m)}(g_1)\cap(\cup^m_{i=1}E_i)=V^{(m)}(g_1)\cap
E_m=\{(v,1+{\xi} u)=(0,0)\}$.
\endroster
because note by Theorem 3.6 or Sublemma 5.4 that we can use the same
$\tau_m$ for the composition of the first finite number $m$ of
successive blow-ups in preparation for finding the standard
resolution of the singular point $(0,0)$ of $V(f)$ if exists, as a
reduced variety.

Moreover, it is very interesting that $(f\circ\tau_m)_{proper}$, as
an element in $\BC\{v,1+{\xi} u\}$ with $(v,1+{\xi} u)=(0,0)$ and
$y=v$, has the same properties and notations at $(y,z)=(0,0)$ as
$f(y,z)$ does at $(y,z)=(0,0)$ in Lemma 1.11 in the sense of
Definition 2.6.  $\square$ \ms

\proclaim{Corollary 1.11.2(for Lemma 1.11)} 

$\underline{\text{\bf Assumptions}}$ Assume that $f=b_nz^n+b_0y^{\beta_0}+\sum^{n-1}_{i=1}b_iy^{\beta_i}z^i$ satisfies 
the same assumption as $f$ does in {\rm Lemma 1.11} and {\rm Lemma 1.11.1}. 
Let $f$ be irreducible in ${}_2\CO_0$, and let $\gcd(n,\beta_0)>1$.\ms
$\underline{\text{\bf Conclusions}}$

$\underline{\text{\bf Conclusions}}$

{\bf Fact(1):} $\underline{\text{\rm If $b_{n-1}$ is a unit in $\C\{y,z\}$}}$, 
then it can be found by the above lemma that either $f\in$the type$[1]$ 
or $f\in$ the type$[2]$ in the sense of Definition 2.5, which follows 
from two examples:

\noindent{\rm(a)} If $f=z^2+(1+y)y^{10}+2y^5z=(z+y^5)^2+y^{11}$, then $f\in$
the type$[1]$ in the sense of Definition 2.5. 

\noindent{\rm(b)} If $f=z^4+y^{10}+y^4z^3+2y^{5}z^2=(z^2+y^5)^2+y^4z^3$, then
$f\in$ the type$[2]$ in the sense of Definition 2.5. \ms

{\bf Fact(2):} $\underline{\text{\rm If $b_{n-1}$ is identically zero 
in $\C\{y,z\}$ and $2\le n\le\beta_0$ in addition}}$, 
then it will be found by {\rm Lemma 1.12} that $f\in$the type$[{\ell}]$ 
with $\ell\ge 2$ in the sense of Definition 2.5. 
$\square$ 
\endproclaim \ms

\proclaim{Lemma 1.12(The fundamental algorithm for finding 
irreducibility criterion of any W-poly in $f\in \BC\{y\}[z]$ with
$f\in$ the type$[1]$ in the sense of Definition 2.5 and its generalizations)}

$\underline{\text{\bf Assumptions}}$ Let
$f=f(y,z)=a_nz^n+a_{n-2}y^{\a_{n-2}}z^{n-2}+\cdots
+a_1y^{\a_1}z+a_0y^{\a_0}$ be in $\BC\{y,z\}$ where for $0\le i\le
n-2$, each $a_i=a_i(y,z)$ is a unit in ${}_2\CO_0$ if exists and the
$\a_i$ are positive integers. $\underline{\text{\rm Note that $a_{n-1}$ is identically
zero}}$. Write $n=d_2n_1$ and $\a_0=d_2\a_{1,0,1}$ with
$d_2=\gcd(n,\a_0)$. In addition, assume that we have the following:

\noindent {\rm(1.12.0)} \qquad \quad   \quad $2\le n\le\a_0$.  \bs

$\underline{\text{\bf Conclusions}}$

{\bf Fact:}  If $f$ is irreducible in ${}_2\CO_0$, $f$ must satisfy
the necessary condition:
$$
\f{\alpha_i}{n-i}\ge \f{\alpha_0}{n} \quad \text{for} \quad 0\le
i\le n-2. \tag 1.12.1
$$

Also, $d_2=\gcd(n,\alpha_0)<n$ because $f$ is irreducible in
${}_2\CO_0$ and $a_{n-1}$ is identically zero.

For $f$ of $(1.12.1)$, it suffices to consider two cases:

{\rm Case(A)} $\gcd(n,\alpha_0)=1$, and {\rm Case(B)}
$1<\gcd(n,\alpha_0)<n$. \ms

{\bf Case(A)} The necessary and sufficient condition for $f(y,z)$ to
be irreducible in ${}_2\CO_0$ with $f\in$the $type[1]$ in the sense of
Definition 2.5 is as follows:
$$
\text{$\gcd(n,\a_{0})=1$ \quad {and} \quad $\f{\alpha_i}{n-i}>
\f{\alpha_0}{n}$ \quad \text{for} \quad $0\le i\le n-2$.} \tag
1.12.2
$$
In this case, $f\buildrel \text{{\rm multiseq }} \over \sim
z^n+{\xi}y^{\alpha_0}=g_1$ with $\xi=a_0(0)$. Equivalently, f can be
rewritten as follows:
$$\align
(1.12.2^*) \qquad   f =z^{n}+ \xi y^{\alpha_0}
+\sum_{\alpha,\beta\ge 0} c_{\alpha,\beta}y^{\alpha}z^{\beta} \quad
\text{with} \quad n\alpha+\alpha_0\beta>n\alpha_0 \quad
\text{and} \quad \gcd(n,\a_{0})=1. \qquad  \\
 \endalign$$

Moreover, $V(f) \buildrel \text{{\rm resol}} \over \sim V(g_1)$ at
the origin in $\BC^2$ in the sense of Definition 2.4.  \ms

{\bf Case(B)} Let $1<\gcd(n,\alpha_0)<n$ with $a_{n-1}$ zero. If $f$
is irreducible in ${}_2\CO_0$, then $f\in$the type$[\ell]$ with
${\ell}\ge 2$ in the sense of Definition 2.5. $\square$ \endproclaim
\ms

\proclaim{Lemma 1.12.1(The fundamental algorithm for finding 
irreducibility criterion of any W-poly in $f\in \BC\{y\}[z]$ with
$f\in$ the type$[2]$ in the sense of Definition 2.5 and its
generalizations)}

$\underline{\text{\bf Assumptions}}$ Let
$f=f(y,z)=z^n+a_{n-2}y^{\a_{n-2}}z^{n-2}+\cdots
+a_1y^{\a_1}z+a_0y^{\a_0}$ be a W-poly of degree $n$ in z where for
$0\le i\le n-2$, each $a_i=a_i(y)$ is a unit in $\BC\{y\}$ if exists
and the $\a_i$ are positive integers. Note that $a_{n-1}$ is
identically zero. Write $n=d_2n_1$ and $\a_0=d_2\a_{1,0,1}$ with
$d_2=\gcd(n,\a_0)$. By the same way as we have seen in {\rm(1.12.0)}
of Lemma 1.12, assume that we have the same additioal inequality:

\noindent {\rm(1.12.0)} \qquad \quad   \quad $2\le n\le\a_0$.  \bs

$\underline{\text{\bf Conclusions}}$

{\bf Fact[I]:} If $f$ is irreducible in ${}_2\CO_0$, f can be
represented as follows:
$$\align
 g_1 &=z^{n_1}+ \xi_1 y^{k_1} \quad \text{with} \quad
k_1=\alpha_{1,0,1}=\beta_{1,0,1}, \tag 1.12.3 \\
f &=g_1^{d_2} +\sum_{\alpha,\beta\ge 0}
c_{\alpha,\beta}y^{\alpha}z^{\beta} \quad
\text{with \quad $\alpha>0$ \ {and} \ $\beta\le n-2$},  \\
 \endalign$$
where $\xi_1=\f{1}{d_2}a_{n-n_1}(0)$, $\alpha_{1,0,1}=\a_{n-n_1}$
and the $c_{\alpha,\beta}$ are nonzero complex numbers for some
nonnegative integers $\alpha$ and $\beta$ such that
$n_1\alpha+\alpha_{1,0,1}\beta>n_1\alpha_{1,0,1}d_2$, with the
following:

\roster
\item "(i)" $\xi_1=\f{1}{d_2}a_{n-n_1}(0)$ is the unique nonzero complex
number such that ${\xi_1}^{d_2}=a_n(0)$ and \quad
${{d_2}\choose{i}}{\xi_1}^i=a_{n-in_1}(0)$ for $1\le i\le d_2$.

\item "(ii)" $\frac{\alpha_{n-in_1}}{in_1}= \frac{\alpha_{1,0,1}}{n_1}$ for
$1\le i\le {{d_2}}$.

\item "(iii)" Either $ \f{\a_j}{n-j}> \frac{\alpha_{1,0,1}}{n_1}
=\frac{\alpha_{0}}{n}$ or
$n_1\a_j+\alpha_{1,0,1}j>n_1\alpha_{1,0,1}d_2$ for any $j\not=
in_1$, if exists.
\endroster

{\bf Fact[II]:} If $f$ is irreducible in ${}_2\CO_0$, $f$ of
{\rm(1.12.3)} can be rewritten in the form
$$\align
f=g^{d_2}_1 +\sum^{d_2-1}_{i=0} T_{2,i}g^i_1, \tag 1.12.4
\endalign$$
satisfying the following:

$\underline{\text{\rm(a)}}$ For $i=0,1,\dots,d_{2}-1$ and for any
nonzero monomial $y^{\de_1}z^{\de_2}$ in $T_{2,i}$, $\de_1>0\quad
\text{and}\quad \de_2<n_1$. Without assuming irreducibility of $f$
in $\BC\{y,z\}$, note by Sublemma 1.9 of Theorem 1.8 that a
coefficient $T_{2,d_{2}-1}$ of ${g_{1}}^{d_2-1}$ may not be
identically zero.

$\underline{\text{\rm(b)}}$ Let $i$ be fixed with $0\le i\le
d_{2}-1$. For brevity of notation, let $\th_1: \N_0\to \N_0$ and
$\th_2: \N^2_0\to \N_0$ be integer-valued functions, each of which
is defined respectively, as follows:
$$\align
(1.12.5) \qquad & \text{$\th_1(t)=t$ for each $t\in \N_0$,} \\
& \text{$\th_2(t_1,t_2)=t_2\th_1(k_1) + n_1\th_1(t_1) =
t_2k_1+n_1t_1$ for each $(t_1,t_2)\in \N^2_0$.} \qquad \qquad \qquad \\\\
\endalign$$

By {\rm(1.11.2)} of {\rm Lemma 1.11}, for any nonzero monomial
$y^{\de_1}z^{\delta_2}$ in $T_{2,i}=T_{2,i}(y,z)\in \C\{y\}[z]$,
$\th_2(\delta_1,\delta_2)>n_1k_1(d_2-i)$.

$\underline{\text{\rm(c)}}$ Following the same properties and
notations as we have used in {\rm (1.11.4)} of {\rm Lemma 1.11.1},
let $\tau_m=\pi_1\circ\pi_2\circ\cdots\circ\pi_m:M^{(m)}\to\BC^2$ be
the compositions of a finite number $m$ of successive blow-ups
$\pi_i$ which is needed to get the standard resolution of the
singular point of $V(g_1)$. With the representation in
{\rm(1.11.4)}, $(f\circ\tau_m)_{total}$ with
$(f\circ\tau_m)_{proper}$ can be rewritten as follows:
$$\align
(1.12.6) \qquad  (f\circ\tau_m)_{total}&=(f\circ\tau_m)(v,u)
=v^{e_{m}}u^{\rho_{m}}(f\circ\tau_m)_{proper},\\
(f\circ\tau_m)_{proper}&=(1+{\xi}u)^{d_2}+\sum^{d_2-1}_{i=0}
\ve_{i}v^{\th_2(\beta_{2,i,1},\beta_{2,i,2})+n_1k_1i-n_1k_1d_2}
(1+{\xi}u)^{i}, \qquad \qquad \qquad \qquad\\
\endalign$$
satisfying the property {\rm(v)} with the properties {\rm(i)},
{\rm(ii)}, {\rm(iii)}, {\rm(iv)} in {\rm(1.11.4)}, as follows:

\roster
\item "(v)" If $T_{2,i}$ is nonzero for $i=0,1,\dots,d_2-1$, then there is a
unique nonzero monomial $C_{2,i}y^{\beta_{2,i,1}}z^{\beta_{2,i,2}}$
in $T_{2,i}$ with a constant $C_{2,i}=\ve_{i}(0,0)$ such that
$\th_2(\beta_{2,i,1},\beta_{2,i,2})
=\min\{\th_2(\gamma_1,\gamma_2)\}$ for any nonzero monomial
$y^{\gamma_1}z^{\gamma_2}$ in $T_{2,i}$ where $\ve_{i}=\ve_{i}(v,
1+{\xi}u)$ is a unit in $\BC\{v, 1+{\xi}u\}$ for $i=0,1,\dots,d_2-1$
if exists.
\endroster  \ms

$\underline{\text{\rm(d)}}$ For all $i= 0,1,\dots,d_2-1$, the
following holds:
$$\align
&\gcd(d_2,\th_2(\beta_{2,0,1},\beta_{2,0,2}))\ge 1 \quad \text{and}  \tag 1.12.7\\
& \dfrac{\th_2(\beta_{2,i,1},\beta_{2,i,2})}{d_2-i}\ge
\dfrac{\th_2(\beta_{2,0,1},\beta_{2,0,2})}{d_2}>n_1k_1.
\endalign$$

Then, either $\gcd(d_{2},\th_2(\beta_{2,0,1},\beta_{2,0,2}))=1$ or
$1<\gcd(d_{2},\th_2(\beta_{2,0,1},\beta_{2,0,2}))\le d_{2}$. \ms

$\underline{\text{\rm(1d-1)}}$ Suppose
$\gcd(d_{2},\th_2(\beta_{2,0,1},\beta_{2,0,2}))=1$. Then $f$ is
irreducible in ${}_2\CO_0$ with $f\in $ the type $[2]$ in the sense
of Definition 2.5 if and only if the inequality in {\rm(1.12.6)}
holds and $g_1$ is irreducible in ${}_2\CO_0$ with $g_{1} \in$ the
type $[1]$ in the sense of Definition 2.5.

$\underline{\text{\rm(1d-2)}}$ Suppose
$1<\gcd(d_{2},\th_2(\beta_{2,0,1},\beta_{2,0,2}))\le d_{2}$ in
{\rm(1.12.6)}. To find an irreducible criterion of any W-poly in
$f\in \BC\{y\}[z]$ with $f\in$ the type $[2]$ in the sense of
Definition 2.5, it remains to solve two subcases respectively:

$\underline{\text{\rm Subcase(i) of (1d-2)}}$ Let
$\gcd(d_{2},\th_2(\beta_{2,0,1},\beta_{2,0,2}))=d_{2}$ in
{\rm(1.12.6)}. Then, $f$ is either irreducible or not in
${}_2\CO_0$. If $f$ is irreducible in ${}_2\CO_0$ then $f \in$ the
type $[\ell]$ with $\ell\ge 2$ in the sense of Definition 2.5.

In this case, we can find a necessary and sufficient condition for
$f$ to be irreducible in ${}_2\CO_0$ with $f\in $ the type $[2]$ in
the sense of Definition 2.5, following the method as in {\rm (b) of
Remark 1.12.1.1 for Lemma 1.12.1}.\ms

$\underline{\text{\rm Subcase(ii) of (1d-2)}}$ Let
$1<\th_2(\beta_{2,0,1},\beta_{2,0,2}))< d_{2}$ in {\rm(1.12.6)}.
Then, $f$ is either irreducible or not in ${}_2\CO_0$. If $f$ is
irreducible in ${}_2\CO_0$ then $f \in$ the type $[\ell]$ with
$\ell\ge 3$ in the sense of Definition 2.5. $\square$ \ms
\endproclaim \ms

\definition{Remark 1.12.1.1 for Lemma 1.12.1} {\rm(a)} Assuming that a
coefficient $T_{2,d_{2}-1}$ of ${g_{1}}^{d_2-1}$ is zero, then it is
clear that $f$ is irreducible in ${}_2\CO_0$ with $f\in $ the type
$[2]$ in the sense of Definition 2.5 if and only if
$\gcd(d_{2},\th_2(\beta_{2,0,1},\beta_{2,0,2}))=1$ and (1.12.7)
holds.

{\rm(b)} Assuming that a coefficient $T_{2,d_{2}-1}$ of
${g_{1}}^{d_2-1}$ is nonzero, by Theorem 1.8 and Lemma 1.12.1 we
will find a necessary and sufficient condition for $f$ to be
irreducible in ${}_2\CO_0$ with $f\in $ the type $[2]$ in the sense
of Definition 2.5. First, using the same kind of properties and
notations as in (1.10.2) in Sublemma 1.10 of Theorem 1.8 and Lemma
1.12.1, it suffices to consider the following: Note that
$h_{1,1}=g_1$ with $T^{(1)}_{2,i}=T_{2,i}$ as we have seen in
Sublemma 1.10.
$$\cases
h_{1,{\nu_1}+1}
&=z^{n_1}+\sum^{n_1-2}_{i=0}R^{({\nu_1}+1)}_{1,i}z^i\
{\text{with}\ h_{1,{\nu_1}+1}\buildrel \text{{\rm multiseq }} \over \sim g_1} \\
f&=h^{d_2}_{1,{\nu_1}+1}+\sum^{d_2-2}_{i=0}T^{({\nu_1}+1)}_{2,i}h^i_{1,{\nu_1}+1}\
{\text{with}\ T^{({\nu_1}+1)}_{2,d_2-1}=0,}
\endcases \tag 1.12.8
$$
for an integer $\nu_{1}\ge 0$ where $R^{({\nu_1}+1)}_{1,i}\in
\C\{y\}$ for $0\le i\le n_1-2$, and $T^{({\nu_1}+1)}_{2,i}\in
\C\{y\}[z]$ for $0\le i\le d_2-2$, satisfying the following:

{\rm(i)} For $0\le i\le n_1-2$, each
$R^{({\nu_1}+1)}_{1,i}=b_iy^{\beta^{({\nu_1}+1)}_{1,i,1}}$ with a
unit $b_i\in\C \{y\}$ and a positive integer
$\beta^{({\nu_1}+1)}_{1,i,1}$ if exists. Denote $A_{1,i}$ by
$b_i(0)$ for convenience of notations.

{\rm(ii)} For any nonzero monomial $y^{\g_1}z^{\g_2}$ in
$T^{({\nu_1}+1)}_{2,i}$,
$$
\g_1>0 \quad \text{and} \quad \g_2<n_1. \tag 1.12.9
$$

{\rm(iii)} Using the same kind of properties and notations as in
Lemma 1.12.1, if $T^{({\nu_1}+1)}_{2,i}\ne 0$, let
$C^{({\nu_1}+1)}_{2,i}\Pi^2_{k=1}f^{\beta^{({\nu_1}+1)}_{2,i,k}}_{k-2}$
be a unique nonzero monomial with a constant $C^{({\nu_1}+1)}_{2,i}$
in $T^{({\nu_1}+1)}_{2,i}$ such that
$\th_2(\beta^{({\nu_1}+1)}_{2,i,k})^2_{k=1}=\min\{\th_2(\g_k)^2_{k=1}\}$
for any nonzero monomial $\Pi^2_{k=1}f^{\g_k}_{k-2}$ in
$T^{(p)}_{2,i}$ where $f_{-1}=y$ and $f_0=z$.

Then, $f$ is irreducible in ${}_2\CO_0$ with $f\in $ the type $[2]$
in the sense of Definition 2.5 if and only if the following
inequalities in {\rm(1.12.10)} hold:
$$\align
(1.12.10) \qquad  &\gcd(n_1,\beta^{({\nu_1}+1)}_{1,0,1})=1 \quad
\text{and} \quad \f{\th_1(\beta^{({\nu_1}+1)}_{1,i,1})}{n_1-i}\ge
\f{\th_1(\beta^{({\nu_1}+1)}_{1,0,1}}{n_1} \qquad \text{for $0\le
i\le n_1-2$,} \qquad \qquad \\
&\gcd(d_2,\th_2(\beta^{({\nu_1}+1)}_{2,0,k})^2_{k=1}) =1 \qquad \text{and} \\
&\f{\th_2(\beta^{({\nu_1}+1)}_{2,i,k})^2_{k=1}}{d_2-i} \ge
\f{\th_2(\beta^{({\nu_1}+1)}_{2,0,k})^2_{k=1}}{d_2}>n_1\alpha_{1,0,1}
\quad \text{for $0\le i\le d_2-2$.}
\endalign$$

Also, if $f$ satisfies (1.12.10), $h_{1,{\nu_1}+1}\buildrel
\text{{\rm multiseq}} \over \sim g_1$ and $f\buildrel \text{{\rm
multiseq}} \over \sim
{g_1}^{d_2}+y^{\beta^{({\nu_1}+1)}_{2,0,1}}z^{\beta^{({\nu_1}+1)}_{2,0,2}}$.
$\square$ \ms
\enddefinition

\proclaim{Lemma 1.12.{$\alpha$}} $\underline{\text{\bf
Assumptions}}$ Let $f=f(y,z)=z^n+a_{n-2}y^{\a_{n-2}}z^{n-2}+\cdots
+a_1y^{\a_1}z+a_0y^{\a_0}$ be a W-poly of degree $n$ in z where for
$0\le i\le n-2$, each $a_i=a_i(y)$ is a unit in $\BC\{y\}$ if exists
and the $\a_i$ are positive integers. Note that $a_{n-1}$ is
identically zero. Write $n=d_2n_1$ and $\a_0=d_2\a_{1,0,1}$ with
$d_2=\gcd(n,\a_0)$. In addition, assume that we have the following:

\noindent {\rm(1.12.{$\alpha$}.0)} \qquad \quad   \quad $n>\a_0\ge
2$. \ms

$\underline{\text{\bf Conclusions}}$

{\bf Case(A)} Let $\gcd(n,\alpha_0)=1$ with $a_{n-1}$ zero. The
necessary and sufficient condition for $f(y,z)$ to be irreducible in
${}_2\CO_0$ with $f\in$the type$[1]$ in the sense of Definition 2.5
is as follows:
$$
\text{\quad $\f{\alpha_i}{n-i}\ge \f{\alpha_0}{n}$ \quad \text{for}
\quad $0\le i\le n-2$.} \tag 1.12.11
$$

In this case, $f\buildrel \text{{\rm multiseq }} \over \sim
z^n+{\xi}y^{\alpha_0}=g_1$. Moreover, $V(f) \buildrel \text{{\rm
resol}} \over \sim V(g_1)$ at $0\in\BC^2$ in the sense of Definition
2.4. \ms

{\bf Case(B)} Let $1<\gcd(n,\alpha_0)<\alpha_0$ with $a_{n-1}$ zero.
If $f$ is irreducible in ${}_2\CO_0$, then $f\in$the type$[\ell]$
with ${\ell}\ge 2$ in the sense of Definition 2.5. \ms

{\bf Case(C)} Let $\gcd(n,\alpha_0)=\alpha_0$ with $a_{n-1}$ zero.
If $f$ is irreducible in ${}_2\CO_0$, then $f\in$the type$[\ell]$
with ${\ell}\ge 1$ in the sense of Definition 2.5. $\square$
\endproclaim \ms

\noindent{\bf Remark 1.12.{$\alpha$}.1.} Let $f(y,z)$ of Lemma
1.12.{$\alpha$} be defined by $(z^2+y)^4+y^4z$. Then,
$\gcd(n,\alpha_0)=\alpha_0=4>1$, $f(y,z)$ is irreducible in
${}_2\CO_0$, and $f\in$the type $[1]$ in the sense of Definition
2.5. \quad $\square$ \bs

{\bf{Example(Example [I] for Lemma 11 and Lemma 12)(2012-5300.page 30)}

(1) Let $g=g(y,z)$, $h=h(y,z)$, $f=f(y,z)$ and $\phi=\phi(y,z)$ be
in ${\BC}\{y,z\}$, which can be written as follows:
$$g=z^4+2y^5.$$
$$h=z^4+3y^4z+2y^5.$$
$$f=(z^4+2y^5)^3+4y^{10}z(z^4+2y^5)+3y^{15}z=g^3+(4y^{10}z)g+3y^{15}z.$$
$$\phi=(z^4+3y^4z+2y^5)^3+4y^{10}z(z^4+3y^4z+2y^5)+3y^{15}z=h^3+(4y^{10}z)h+3y^{15}z.$$

\text{\bf The 1st Problem:} The aim is to compute that two zero sets
$\{g=0\}$ and $\{h=0\}$ have the same kind of the singularity at
$0$.

\text{\bf The 2nd Problem:} Let $F=F(y,z)=g^3+y^{15}z$. The aim is
to compute that two zero sets $\{f=0\}$ and $\{\phi=0\}$ have the
same kind of the singularity at $0$ as $\{F=0\}$ does at $0$. \bs

{\bf{Example(Example [II]] for Lemma 11 and Lemma 12).}  By $g$ of
Example [I], $h=h(y,z)$, $f=f(y,z)$ and $\phi=\phi(y,z)$ be in
${\BC}\{y,z\}$, which can be rewritten as follows:
$$\align
h&=g+3y^4z+2y^7 =g+\sum_{\alpha,\beta\ge
0}c_{\alpha,\beta}y^{\alpha}z^{\beta},  \quad {with}\quad
4\alpha+5\beta>20  \\
f&=(z^4+2y^5)^3+4y^{10}z(z^4+2y^5)+3y^{15}z \quad \text{with \quad $g=(z^4+2y^5)$ } \\
 & =g^3+R_1g+R_0  \quad \text{with\
   $R_1=4y^{10}z$ and $R_0=3y^{15}z$ {in} ${\BC}\{y,z\}$} \\
     & =g^3+\sum_{\gamma,\delta\ge 0}b_{\gamma,\delta}y^{\gamma}z^{\delta},  \quad {with}\quad
    4\gamma+5\delta>60 \\
          &=z^{12}+6y^5z^8+4y^{10}z^5++12y^{10}z^4+ 8y^{15}z+3y^{15}z +8y^{15}\\
    & \qquad   \text{where \quad $g=(z^4+2y^5)$. }
 \endalign$$

{\bf{Example(Example [III]] for Lemma 11 and Lemma 12).}  By
$g=g(y,z)$ of Example [I], $h=h(y,z)$, $f=f(y,z)$ and
$\phi=\phi(y,z)$ be in ${\BC}\{y,z\}$, which can be rewritten as
follows:
$$\align
\phi&=(z^4+3y^4z+2y^5)^3+4y^{10}z(z^4+3y^4z+2y^5)+3y^{15}z+y^{17} \\
 & =h^3+R_1h+R_0  \quad {with}\quad
    R_i=R_i\{y,z\}\in {\BC}\{y,z\} \\
      &=g^3+9y^4zg^2+(27y^8z^2+4y^{10}z)g+27y^{12}z^3+12y^{14}z^2+3y^{15}z+y^{17}\\
      & =g^3+T_2g^2+T_1g+R_0  \quad {with}\quad
    T_i=T_i\{y,z\}\in {\BC}\{y,z\} \\
     & =g^3+\sum_{\gamma,\delta\ge 0}a_{\gamma,\delta}y^{\gamma}z^{\delta},  \quad {with}\quad
    4\gamma+5\delta>60 \\
       &=z^{12}+9y^4z^9+6y^5z^8+27y^8z^6+(6y^{9}+4y^{10})z^5+12y^{10}z^4\\
    &\quad +27y^{12}z^3+(54y^{13}+12y^{14})z^2+(36y^{14}+11y^{15})z+8y^{15}+y^{17} \\
   \endalign$$
\bs\bs

{\bf $\S$1.7. Irreducibility criterion of W-polys of two complex
variables(A generalized representation of irreducible {W}-polys of
two complex variables)} \ms

In this section, in order to find irreducibility criterion for germs
of analytic functions of two complex variables, without loss of
generality, it suffices to find the necessary and sufficient
condition for $f(y,z)$ of all the W-polys of two complex variables
to be irreducible in $\BC\{y,z\}$ with $f\in {\text{\rm
type}}[{\ell}]$ in the sense of Definition 2.5 in terms of Theorem
1.13, using Theorem 15.2(The WDT for the W-polys) and Theorem
15.4(The Division Algorithm for the W-polys) and Theorem 12.0. In
$\S$1.8 and $\S$1.9, as an application of this theorem, it will be
found without proof that we can write The 2nd Algorithm and The 3rd
Algorithm. \ms

\proclaim{Theorem 1.13(How to find a generalized representation of
irreducible {W}-polys of two complex variables(Irreducibility
criterion of W-polys of two complex variables))}

$\underline{\text{\bf {Assumptions}}}$ \quad Let $f\in \BC\{y\}[z]$
be an arbitrary $W$-poly of degree $n\ge 2$ in $z$. Without loss of
generality, we may assume that $f$ satisfies the following form:
$$ f=z^n+\sum^{n-2}_{i=0} a_iy^{\a_i}z^i, \tag Eq.1 $$
where for $0\le i\le n-2$, each $a_i=a_i(y)$ is a unit in
${}_2\CO_0$ for $0\le i\le n-2$,  if exists, and the $\a_i$ are
positive integers. Note that $a_{n-1}$ is identically zero.  Write
$n=d_2n_1$ and $\a_0 =d_2\a_{1,0,1}$ with $d_2=\gcd(n,\a_0)$. Write
$n=\Pi^{r}_{k=1}n_k$ with positive integers $n_k\ge 2$ for all $k$
where the $n_k$ may not be the factorization of prime numbers.

In addition, assume that we have the following:

\noindent {\rm(1.13.0)} \qquad \quad   \quad $2\le n\le\a_0$.  \bs

$\underline{\text{\bf {Conclusions}}}$\ {\bf The necessary and
sufficient condition for $f(y,z)$ to be irreducible in $\BC\{y,z\}$
with $f\in {\text{\bf type}}[{\ell}]$ in the sense of Definition
2.5(Theorem 12.0) is as follows:}

By {\rm {Theorem 15.4(The Division Algorithm for $W$-polys)}} for
each $k=1,2,\dots,{\ell}$, $f_k$ an $f$ can be written in the form
$$ \cases
f_k &=f^{n_k}_{k-1}+\sum^{{n_k}-2}_{i=0} R_{k,i}f^i_{k-1}  \\
f_{\ell}& =f^{n_{\ell}}_{\ell-1} +\sum^{n_{\ell}-2}_{i=0}
R_{\ell,i}f^i_{\ell-1}  \quad \text{with $f=f_{\ell}$}
\endcases \tag 1.13.1
$$
where, considering $y,z,f_1,\dots, f_j$ as independent complex
$(j+2)$-variables at the origin in $\C^{j+2}$ with $f_{-1}=y$ and
$f_{0}=z$,

{\rm(i)} $n=\Pi^{\ell}_{k=1}n_k$ with $n_k\ge 2$ for $1\le k\le
{\ell}$;

{\rm(ii)} for each fixed $k$ and for each $i$ with $0\le i\le
n_k-2$, $R_{k,i}\in \C\{y,z,f_1,\dots,f_{k-2}\}$;

{\rm(iii)} for each $k=1,2,\dots,{\ell-1}$, $f_k=f_k(y,z,f_1,\dots,
f_{k-1})\in \C\{y,z,f_1,\dots, f_{k-2}\}[f_{k-1}]$;

{\rm(iv)} $f=f(y,z,f_1,\dots, f_{\ell-1})\in \C\{y,z,f_1,\dots,
f_{\ell-2}\}[f_{\ell-1}]$ with $f=f_{\ell}$;

{\noindent}satisfying a finite number of conditions, each of which
is represented respectively, as follows: \ms

$\underline{ \text{\bf (1) Condition[A] for $f_1(y,z)\in $the
type[1] in the sense of Definition 2.5:}}$

$R_{1,i}\in \C\{y\}$ satisfies the properties {\rm(1a)}, {\rm(1b)}
and {\rm(1c)} for each $i=0,1,\dots,n_1-2,$ and then
$f_1=f_1(y,z)\in \C\{y,z\}$ satisfies the properties {\rm(1d)}. \ms

{\rm(1a)} For $0\le i\le n_1-2$, each $R_{1,i}=b_iy^{\a_{1,i,1}}$
with a unit $b_i\in\C \{y\}$ and a positive integer $\a_{1,i,1}$ if
exists. Denote $A_{1,i}$ by $b_i(0)$ for convenience of notations.

{\rm(1b)} Define a function $\th_1:\N_0\to\N_0$ by $\th_1(t)=t$
where $\N_0$ is the set of nonnegative integers.

{\rm(1c)} ${\th_1(\a_{1,i,1})}>({n_1}-i)$ for all $i=0,1,\dots,
n_1-2$, where $n_1$ is the multiplicity of $f_1$ at $0\in \BC^2$
with $n_1\ge 2$.

{\rm(1d)} For all $i=0,1,\dots, n_1-2$,
$$\align
&\gcd(n_1,\a_{1,0,1})=1 \quad \text{and} \tag 1.13.2\\
&\f{\th_1(\a_{1,i,1})}{n_1-i}=\f{\a_{1,i,1}}{n_1-i}\ge
\f{\a_{1,0,1}}{n_1}=\f{\th_1(\a_{1,0,1})}{n_1}.
\endalign$$

$\underline{ \text{\bf (2) Condition[A] for $f_2(y,z)\in$the type[2]
in the sense of Definition 2.5:}}$

$R_{2,i}\in \C\{y\}[z]$ satisfies the properties {\rm(2a)},
{\rm(2b)} and {\rm(2c)} for each $i=0,1,\dots,n_2-2,$ and then
$f_2=f_2(y,z,f_1)\in \C\{y,z\}[f_1]$ satisfies the properties
{\rm(2d)}. \ms

{\rm(2a)} For any nonzero monomial $y^{\de_1}z^{\de_2}$ in
$R_{2,i}$, $\de_1>0$ and $\de_2<n_1$.

{\rm(2b)} Let $\N^2_0$ be two-dimensional cartesian product of
$\N_0$. For given integers $n_1,\a_{1,0,1}$ and a function $\th_1$
in {\rm Cond[A] of the 1st type}, define a function $\th_2:\N^2_0\to
\N_0$ by
$$
\th_2(t_1,t_2)=t_2\th_1(\a_{1,0,1})+n_1\th_1(t_1)=t_2\a_{1,0,1}+n_1t_1
\quad \text{for each $(t_1,t_2)\in \N^2_0$}. \tag 1.13.3 $$

Then, for any two nonzero monomials $y^{\a_1}z^{\a_2}$ and
$y^{\de_1}z^{\de_2}$ in $R_{2i}$ with $i$ fixed,
$$\align
\text{\rm(1.13.3-1)} \qquad \qquad &
\th_2(\a_1,\a_2)=\th_2(\de_1,\de_2)\ \text{if and only if}\
\a_1=\de_1\ \text{and}\ \a_2=\de_2. \qquad \qquad\qquad \qquad \\
 &\text{So, there exists a unique nonzero monomial
$A_{2,i}y^{\a_{2,i,1}}z^{\a_{2,i,2}}$ in $R_{2,i}$} \\
&\text{with a nonzero constant $A_{2,i}$ such that
$\th_2(\a_{2,i,1},
\a_{2,i,2})=\text{$\min$}\{\th_2(\de_1,\de_2)\}$} \\
&\text{for any nonzero monomial $y^{\de_1}z^{\de_2}$ in $R_{2,i}$
with $i$ fixed.}
\endalign$$

{\rm(2c)} For all $i=0,1,\dots, n_2-2$,
$$\align
{\th_2(\a_{2,i,k})^2_{k=1}}>({n_2}-i) n_{1}{\a_{1,0,1}}. \tag 1.13.4
\endalign$$

{\rm(2d)} For all $i=0,1,\dots, n_2-2$,
$$\align
&\gcd(n_2,\th_2(\a_{2,0,k})^2_{k=1})=1 \quad \text{and} \tag 1.13.5\\
&\f{\th_2(\a_{2,i,1}, \a_{2,i,2})}{n_2-i}\ge
\f{\th_2(\a_{2,0,1},\a_{2,0,2})}{n_2}.
\endalign$$

$\underline{ \text{\bf (3) Condition[A] for $f_m(y,z)\in$the type[m]
in the sense of Definition 2.5:}}$

For each fixed $m$ with $3\le m\le {\ell-1}$, $R_{m,i}\in
\C\{y,z,f_1,\dots,f_{m-2}\}$ satisfies the properties {\rm(3a)},
{\rm(3b)} and {\rm(3c)} for each $i=0,1,\dots,n_m-2$, and then
$f_m=f_m(y,z,f_1,\dots,f_{m-1})\in
\C\{y,z,f_1,\dots,f_{m-2}\}[f_{m-1}]$ satisfies the properties
{\rm(3d)}. \ms

{\rm(3a)} For any nonzero monomial $\Pi^m_{k=1}f^{\de_k}_{k-2}$ in
$R_{m,i}$ with $f_{-1}=y$ and $f_0=z$, $\de_1>0$ and $\de_k<n_{k-1}$
for $k=2,3,\dots,m$. \ms

{\rm(3b)} By induction assumption on the integer $(m-1)\le
{\ell-1}$, there exists a sequence $\{f_3,f_4,\dots,f_{m-1}\}$, each
of which satisfies the same kind of properties and notations as we
have seen in $\underline{ \text{\rm Condition[A] for
$f_2(y,z)\in$the type[2] in the sense of Definition 2.5}}$. Then
inductively, define $\th_m:\N^m_0\to \N_0$ where $\N^m_0$ is its
$m$-dimensional cartesian product by
$$
(1.13.6) \quad
\th_m(t_k)^m_{k=1}=t_m\th_{m-1}(\a_{m-1,0,k})^{m-1}_{k=1}
+n_{m-1}\th_{m-1}(t_k)^{m-1}_{k=1} \quad \text{for each
$(t_k)^m_{k=1}\in \N^k_0$},
$$
where recall by induction assumption that for a fixed $i$,
$A_{m-1,i}\Pi^{m-1}_{k=1}f^{\a_{m-1,i,k}}_{k-2}$ is a unique nonzero
monomial in $R_{m-1,i}$ with a constant $A_{m-1,i}$ such that
$$
\th_{m-1}(\a_{m-1,i,k})^{m-1}_{k=1}=\text{$\min$}\{\th_{m-1}(\de_k)^{m-1}_{k=1}\},
\tag 1.13.7
$$
for any nonzero monomial $\Pi^{m-1}_{k=1}f^{\de_k}_{k-2}$ in
$R_{m-1,i}$.

Then, for any two nonzero monomials $\Pi^m_{k=1}f^{\a_k}_{k-2}$ and
$\Pi^m_{k=1}f^{\de_k}_{k-2}$ in $R_{m,i}$ with $i$ fixed,
$$\align
\text{\rm(1.13.7-1)} \qquad \qquad
&\text{$\th_m(\a_k)^m_{k=1}=\th_m(\de_k)^m_{k=1}$
if and only if $\a_k=\de_k$ for $k=1,2,\dots, m.$} \qquad \qquad\qquad \qquad\\
&\text{So, there exists a unique nonzero-monomial
$A_{m,i}\Pi^m_{k=1}f^{\a_{m,i,k}}_{k-2}$ in $R_{m,i}$} \\
&\text{with a constant $A_{m,i}$ such that
$\th_m(\a_{m,i,k})^m_{k=1}=\text{$\min$}\{\th_m(\de_k)^m_{k=1}\}$}\\
&\text{for any nonzero monomial $\Pi^m_{k=1}f^{\de_k}_{k-2}$ in
$R_{m,i}$.}\\
\endalign$$

{\rm(3c)} For all $i=0,1,\dots, n_m-2$,
$$\align
{\th_m(\a_{m,i,k})^m_{k=1}}>({n_m}-i)
n_{m-1}\th_{m-1}(\a_{m-1,0,k})^{m-1}_{k=1}.\tag 1.13.8
\endalign$$

{\rm(3d)} For all $i=0,1,\dots,n_m-2$,
$$\align
&\gcd(n_m,\th_m(\a_{m,0,k})^m_{k=1})=1 \quad \text{and} \tag 1.13.9\\
&\f{\th_m(\a_{m,i,k})^m_{k=1}}{{n_m}-i}\ge
\f{\th_m(\a_{m,0,k})^m_{k=1}}{n_m}. \quad {\square}
\endalign$$
\endproclaim

\noindent{\bf Remark 1.13.1.} {\rm(1)} There is nothing to prove for
Theorem 1.13, because the sufficiency of the condition in Theorem
1.13 can be proved by Theorem 16.5, and the necessity of the
condition in Theorem 1.13 can be proved by Theorem 16.6, too.

{\rm(2)} The converse of Theorem 16.5 can be represented by Theorem
16.6. Moreover, we can compute irreducibility criterion of W-polys
defining plane curve singularities at the origin in $\C^2$ in the
process of the proof of Theorem 16 together with Proposition 16.7
and Proposition 16.8 completely and rigorously, using the Euclidean
algorithm and Theorem 15.4(The Division Algorithm for the W-polys).
$\square$ \ms

\definition{Remark 1.13.2} Consider the sequence $S=\{f_k: 1\le k\le{\ell}\}$
with $f_{\ell}=f$ where $f_k=f_k(y,z,\dots,f_{k-1})\in
\C\{y,z,f_1,\dots, f_{k-1}\}$, which have the same properties and
notations as we have seen in {\rm(1.13.1)} of the conclusion of
Theorem $1.13$. If $f\in \BC\{y\}[z]$ is irreducible in
$\BC\{y,z\}$, then $f=f_{\ell}(y,z,f_1,\dots,f_{\ell-1})\in
\C\{y,z,f_1,\dots,f_{\ell-2}\}[f_{\ell-1}]$ is an irreducible
$W$-poly of degree $n_{\ell}$ in $f_{\ell-1}$ with coefficients in
$\BC\{y,z,f_1,\dots,f_{\ell-2}\}$ and with multiplicity $n_{\ell}$
at the origin in $\C^{\ell}$. $\square$ \ms
\enddefinition

\proclaim{Corollary 1.13.3} $\underline{\text{\bf {Assumptions}}}$
\quad Under the same assumption and conclusion as in Theorem $1.13$,
note that $f_k$ is irreducible in $\BC\{y,z\}$ with isolated
singularity at $(0,0)$ in $\BC^2$ for $k\ge 1$. In particular, for
each $k=1,2,\dots,\ell$, let $V(H_{k})=\{(y,z):H_{k}(y,z)=0\}$ be an
analytic variety at $(0,0)$ in $\BC^2$, each of which is defined as
follows:
$$\align
(1.13.3.1) \qquad \qquad \text{\rm(i)} \qquad
H_1&=z^{n_1}+y^{\alpha_{1,0,1}} \quad \text{with $n_1\ge 2$ and
$\alpha_{1,0,1}\ge 2$}. \qquad \qquad \qquad \qquad \\
\text{\rm(ii)} \qquad H_2&=H^{n_2}_1+y^{\alpha_{2,0,1}}z^{\alpha_{2,0,2}}.\\
\qquad  &\ldots\ldots \\
\text{\rm({$\ell$})} \qquad
H_{\ell}&=H^{n_{\ell}}_{\ell-1}+y^{\alpha_{\ell,0,1}}z^{\alpha_{\ell,0,2}}H^{\alpha_{\ell,0,3}}_1\cdots
H^{\alpha_{\ell,0,\ell}}_{\ell-2}. \\
\endalign$$

$\underline{\text{\bf {Conclusions}}}$

{\bf Fact[I]:} Then, ${f_k} \buildrel \text{{\rm multiseq }} \over
\sim H_{k}$ for each $k=1,2,\dots,\ell$.

{\bf Fact[II]:} Then, ${f_k} \buildrel \text{{\rm resol }} \over
\sim H_{k}$ for each $k=1,2,\dots,\ell$. \quad $\square$
\endproclaim \ms

\definition{Remark for Corollary 1.13.3}
{\bf(I)} Note by Theorem $5.0$ that $H_{j+1}$ is irreducible in
$\BC\{y,z\}$ with $H_{j+1}\in \text{\rm the type[j+1]}$ in the sense
of Definition of $2.5$ $\iff$ $\gcd(n_1,\alpha_{1,0,1})=1$,
$\gcd(n_2,\theta_2(\alpha_{2,0,k})^{2}_{k=1})=1$,\dots,
$\gcd(n_{j+1},\theta_{j+1}(\alpha_{{j+1},0,k})^{j+1}_{k=1})=1$. \ms

\noindent{\bf(II)} Note that $f_{j+1}$ is irreducible in
$\BC\{y,z\}$ with $f_{j+1}\in \text{\rm the type[j+1]}$ in the sense
of Definition of $2.5$ $\iff$ the following holds:

(1)\quad  $\gcd(n_1,\alpha_{1,0,1})=1$ and
$\dfrac{\th_1(\a_{1,i,1})}{n_1-i}=\dfrac{\a_{1,i,1}}{n_1-i}\ge
\dfrac{\a_{1,0,1}}{n_1}=\dfrac{\th_1(\a_{1,0,1})}{n_1}$ for $0\le
i\le n_1-2$.

(2) \quad $\gcd(n_2,\theta_2(\alpha_{2,0,1},\alpha_{2,0,2}))=1$ and
$\dfrac{\th_2(\a_{2,i,1}, \a_{2,i,2})}{n_2-i}\ge
\dfrac{\th_2(\a_{2,0,1},\a_{2,0,2})}{n_2}$ for $0\le i\le n_2-2$.
\ms

\qquad \qquad \qquad $\ldots\ldots$

(j+1) \quad
$\gcd(n_{j+1},\theta_{j+1}(\alpha_{{j+1},0,k})^{j+1}_{k=1})=1$ and
$\dfrac{\th_{j+1}(\a_{{j+1},i,k})^{j+1}_{k=1}}{{n_{j+1}}-i}\ge
\dfrac{\th_{j+1}(\a_{{j+1},0,k})^{j+1}_{k=1}}{n_{j+1}}$ for $0\le
i\le n_{j+1}-2$. $\square$
\enddefinition \bs

\proclaim{Theorem 1.14(A generalized representation of irreducible
{W}-polys of two complex variables(Irreducibility criterion of
W-polys of two complex variables))}

$\underline{\text{\bf {Assumptions}}}$ \quad Let $f\in \BC\{y\}[z]$
be an arbitrary $W$-poly of degree $n\ge 2$ in $z$. Without loss of
generality, we may assume that $f$ satisfies the following form:
$$ f=z^n+\sum^{n-2}_{i=0} a_iy^{\a_i}z^i, \tag Eq.1 $$
where for $0\le i\le n-2$, each $a_i=a_i(y)$ is a unit in
${}_2\CO_0$ for $0\le i\le n-2$,  if exists, and the $\a_i$ are
positive integers. Note that $a_{n-1}$ is identically zero.  Write
$n=d_2n_1$ and $\a_0 =d_2\a_{1,0,1}$ with $d_2=\gcd(n,\a_0)$. Write
$n=\Pi^{r}_{k=1}n_k$ with positive integers $n_k\ge 2$ for all $k$
where the $n_k$ may not be the factorization of prime numbers.

In preparation for finding the irreducibility criterion for $f(y,z)$
of all the W-polys of two complex variables, it suffices to consider
three cases, respectively:
$$\align
&\text{{\bf Case(${\alpha}$)} \quad  $2\le n<\a_0$,}  \tag 1.14.0\\
 &\text{{\bf Case(${\beta}$)} \quad  $2\le \a_0<n$ with $\gcd(n,\a_0)<\a_0$,}  \\
&\text{{\bf Case(${\gamma}$)} \quad $n=p\a_0$ for an integer $p>0$.} \\
\endalign$$

$\underline{\text{\bf {Conclusions}}}$ \quad {\bf(I)} {\bf Let
Case(${\alpha}$) hold.} If $f(y,z)$ is irreducible in $\BC\{y,z\}$,
note that $2\le n<\a_0$ if and only if $2\le n<\a_0$ with
$\gcd(n,\a_0)<n$.

{\bf Then, the necessary and sufficient condition for $f(y,z)$ to be
irreducible in $\BC\{y,z\}$ with $f\in {\text{\bf type}}[{\ell}]$ in
the sense of Definition 2.5(Theorem 12.0) was already given by the
condition in the conclusion in Theorem 1.13.} \ms

{\bf(II)} {\bf Let Case(${\beta}$) hold.} {\bf Then, the necessary
and sufficient condition for $f(y,z)$ to be irreducible in
$\BC\{y,z\}$ with $f\in {\text{\bf type}}[{\ell}]$ in the sense of
Definition 2.5 is the same as the given condition  by the same
method as we have used in the conclusion in Theorem 1.13.} \ms

{\bf(III)} {\bf Let Case(${\gamma}$) hold.} {\bf Then, the necessary
and sufficient condition for $f(y,z)$ to be irreducible in
$\BC\{y,z\}$ with $f\in {\text{\bf type}}[{\ell-1}]$ in the sense of
Definition 2.5 is as follows:}

By {\rm {Theorem 15.4(The Division Algorithm for $W$-polys)}} for
each $k=1,2,\dots,{\ell-1}$, $f_k$ and $f$ can be written in the
form
$$ \cases
f_k &=f^{n_k}_{k-1}+\sum^{{n_k}-2}_{i=0} R_{k,i}f^i_{k-1}  \\
f& =f^{n_{\ell}}_{\ell-1} +\sum^{n_{\ell}-2}_{i=0}
R_{\ell,i}f^i_{\ell-1}
\endcases \tag 1.14.1
$$
where, considering $y,z,f_1,\dots, f_j$ as independent complex
$(j+2)$-variables at the origin in $\C^{j+2}$ with $f_{-1}=y$ and
$f_{0}=z$,

{\rm(i)} $n=\Pi^{\ell}_{k=1}n_k$ with $n_k\ge 2$ for $1\le k\le
{\ell}$;

{\rm(ii)} for each fixed $k$ and for each $i$ with $0\le i\le
n_k-2$, $R_{k,i}\in \C\{y,z,f_1,\dots,f_{k-2}\}$;

{\rm(iii)} for each $k=1,2,\dots,{\ell-1}$, $f_k=f_k(y,z,f_1,\dots,
f_{k-1})\in \C\{y,z,f_1,\dots, f_{k-2}\}[f_{k-1}]$;

{\rm(iv)} $f=f(y,z,f_1,\dots, f_{\ell-1})\in \C\{y,z,f_1,\dots,
f_{\ell-2}\}[f_{\ell-1}]$ with $f=f_{\ell}$;

{\noindent}satisfying a finite number of conditions, each of which
is represented respectively, as follows: \ms

$\underline{ \text{\bf (1) Condition[A] for $f_1(y,z)\in $the
type[0] in the sense of Definition 2.5:}}$

$R_{1,i}\in \C\{y\}$ satisfies the properties {\rm(1a)}, {\rm(1b)}
and {\rm(1c)} for each $i=0,1,\dots,n_1-2,$ and then
$f_1=f_1(y,z)\in \C\{y,z\}$ satisfies the properties {\rm(1d)}. \ms

{\rm(1a)} For $0\le i\le n_1-2$, each $R_{1,i}=b_iy^{\a_{1,i,1}}$
with a unit $b_i\in\C \{y\}$ and a positive integer $\a_{1,i,1}$ if
exists. Denote $A_{1,i}$ by $b_i(0)$ for convenience of notations.

{\rm(1b)} Define a function $\th_1:\N_0\to\N_0$ by $\th_1(t)=t$
where $\N_0$ is the set of nonnegative integers.

{\rm(1c)} ${\th_1(\a_{1,i,1})}>({\a_{1,0,1}}-i)$ for all
$i=0,1,\dots, n_1-2$, where $\a_{1,0,1}=1$ is the multiplicity of
$f_1$ at $0\in \BC^2$ with $n_1\ge 2$.

{\rm(1d)} For all $i=0,1,\dots, n_1-2$,
$$\align
&\gcd(n_1,\a_{1,0,1})=1 \quad \text{and} \tag 1.14.2\\
&\f{\th_1(\a_{1,i,1})}{n_1-i}=\f{\a_{1,i,1}}{n_1-i}\ge
\f{\a_{1,0,1}}{n_1}=\f{\th_1(\a_{1,0,1})}{n_1}.
\endalign$$

$\underline{ \text{\bf (2) Condition[A] for $f_2(y,z)\in$the type[1]
in the sense of Definition 2.5:}}$

$R_{2,i}\in \C\{y\}[z]$ satisfies the properties {\rm(2a)},
{\rm(2b)} and {\rm(2c)} for each $i=0,1,\dots,n_2-2,$ and then
$f_2=f_2(y,z,f_1)\in \C\{y,z\}[f_1]$ satisfies the properties
{\rm(2d)}. \ms

{\rm(2a)} For any nonzero monomial $y^{\de_1}z^{\de_2}$ in
$R_{2,i}$, $\de_1>0$ and $\de_2<n_1$.

{\rm(2b)} Let $\N^2_0$ be two-dimensional cartesian product of
$\N_0$. For given integers $n_1,\a_{1,0,1}$ and a function $\th_1$
in {\rm Cond[A] of the 1st type}, define a function $\th_2:\N^2_0\to
\N_0$ by
$$
\th_2(t_1,t_2)=t_2\th_1(\a_{1,0,1})+n_1\th_1(t_1)=t_2\a_{1,0,1}+n_1t_1
\quad \text{for each $(t_1,t_2)\in \N^2_0$}. \tag 1.14.3 $$

Then, for any two nonzero monomials $y^{\a_1}z^{\a_2}$ and
$y^{\de_1}z^{\de_2}$ in $R_{2i}$ with $i$ fixed,
$$\align
\text{\rm(1.14.3-1)} \qquad \qquad &
\th_2(\a_1,\a_2)=\th_2(\de_1,\de_2)\ \text{if and only if}\
\a_1=\de_1\ \text{and}\ \a_2=\de_2. \qquad \qquad\qquad \qquad \\
 &\text{So, there exists a unique nonzero monomial
$A_{2,i}y^{\a_{2,i,1}}z^{\a_{2,i,2}}$ in $R_{2,i}$} \\
&\text{with a nonzero constant $A_{2,i}$ such that
$\th_2(\a_{2,i,1},
\a_{2,i,2})=\text{$\min$}\{\th_2(\de_1,\de_2)\}$} \\
&\text{for any nonzero monomial $y^{\de_1}z^{\de_2}$ in $R_{2,i}$
with $i$ fixed.}
\endalign$$

{\rm(2c)} For all $i=0,1,\dots, n_2-2$,
$$\align
{\th_2(\a_{2,i,k})^2_{k=1}}>({n_2}-i) n_{1}{\a_{1,0,1}}. \tag 1.14.4
\endalign$$

{\rm(2d)} For all $i=0,1,\dots, n_2-2$,
$$\align
&\gcd(n_2,\th_2(\a_{2,0,k})^2_{k=1})=1 \quad \text{and} \tag 1.14.5\\
&\f{\th_2(\a_{2,i,1}, \a_{2,i,2})}{n_2-i}\ge
\f{\th_2(\a_{2,0,1},\a_{2,0,2})}{n_2}.
\endalign$$

$\underline{ \text{\bf (3) Condition[A] for $f_m(y,z)\in$the
type[m-1] in the sense of Definition 2.5:}}$

For each fixed $m$ with $3\le m\le {\ell-1}$, $R_{m,i}\in
\C\{y,z,f_1,\dots,f_{m-2}\}$ satisfies the properties {\rm(3a)},
{\rm(3b)} and {\rm(3c)} for each $i=0,1,\dots,n_m-2$, and then
$f_m=f_m(y,z,f_1,\dots,f_{m-1})\in
\C\{y,z,f_1,\dots,f_{m-2}\}[f_{m-1}]$ satisfies the properties
{\rm(3d)}. \ms

{\rm(3a)} For any nonzero monomial $\Pi^m_{k=1}f^{\de_k}_{k-2}$ in
$R_{m,i}$ with $f_{-1}=y$ and $f_0=z$, $\de_1>0$ and $\de_k<n_{k-1}$
for $k=2,3,\dots,m$. \ms

{\rm(3b)} By induction assumption on the integer $(m-1)\le
{\ell-1}$, there exists a sequence $\{f_3,f_4,\dots,f_{m-1}\}$, each
of which satisfies the same kind of properties and notations as we
have seen in $\underline{ \text{\rm Condition[A] for
$f_2(y,z)\in$the type[1] in the sense of Definition 2.5}}$. Then
inductively, define $\th_m:\N^m_0\to \N_0$ where $\N^m_0$ is its
$m$-dimensional cartesian product by
$$
(1.14.6) \quad
\th_m(t_k)^m_{k=1}=t_m\th_{m-1}(\a_{m-1,0,k})^{m-1}_{k=1}
+n_{m-1}\th_{m-1}(t_k)^{m-1}_{k=1} \quad \text{for each
$(t_k)^m_{k=1}\in \N^k_0$},
$$
where recall by induction assumption that for a fixed $i$,
$A_{m-1,i}\Pi^{m-1}_{k=1}f^{\a_{m-1,i,k}}_{k-2}$ is a unique nonzero
monomial in $R_{m-1,i}$ with a constant $A_{m-1,i}$ such that
$$
\th_{m-1}(\a_{m-1,i,k})^{m-1}_{k=1}=\text{$\min$}\{\th_{m-1}(\de_k)^{m-1}_{k=1}\},
\tag 1.14.7
$$
for any nonzero monomial $\Pi^{m-1}_{k=1}f^{\de_k}_{k-2}$ in
$R_{m-1,i}$.

Then, for any two nonzero monomials $\Pi^m_{k=1}f^{\a_k}_{k-2}$ and
$\Pi^m_{k=1}f^{\de_k}_{k-2}$ in $R_{m,i}$ with $i$ fixed,
$$\align
\text{\rm(1.14.7-1)} \qquad \qquad
&\text{$\th_m(\a_k)^m_{k=1}=\th_m(\de_k)^m_{k=1}$
if and only if $\a_k=\de_k$ for $k=1,2,\dots, m.$} \qquad \qquad\qquad \qquad\\
&\text{So, there exists a unique nonzero-monomial
$A_{m,i}\Pi^m_{k=1}f^{\a_{m,i,k}}_{k-2}$ in $R_{m,i}$} \\
&\text{with a constant $A_{m,i}$ such that
$\th_m(\a_{m,i,k})^m_{k=1}=\text{$\min$}\{\th_m(\de_k)^m_{k=1}\}$}\\
&\text{for any nonzero monomial $\Pi^m_{k=1}f^{\de_k}_{k-2}$ in
$R_{m,i}$.}\\
\endalign$$

{\rm(3c)} For all $i=0,1,\dots, n_m-2$,
$$\align
{\th_m(\a_{m,i,k})^m_{k=1}}>({n_m}-i)
n_{m-1}\th_{m-1}(\a_{m-1,0,k})^{m-1}_{k=1}.\tag 1.14.8
\endalign$$

{\rm(3d)} For all $i=0,1,\dots,n_m-2$,
$$\align
&\gcd(n_m,\th_m(\a_{m,0,k})^m_{k=1})=1 \quad \text{and} \tag 1.14.9\\
&\f{\th_m(\a_{m,i,k})^m_{k=1}}{{n_m}-i}\ge
\f{\th_m(\a_{m,0,k})^m_{k=1}}{n_m}. \quad {\square}
\endalign$$
\endproclaim \ms

\proclaim{Corollary 1.14.1} $\underline{\text{\bf {Assumptions}}}$
\quad Suppose that the same assumption as in Theorem $1.14$ holds.
In addition, assume that we have the following:

\noindent {\rm(1.14.0)} \qquad \quad   \quad $n=p\a_0$ for an
integer $p>0$.  \ms

In particular, for each $k=1,2,\dots,\ell$, let
$V(H_{k})=\{(y,z):H_{k}(y,z)=0\}$ be an analytic variety at $(0,0)$
in $\BC^2$, each of which is defined as follows:
$$\align
(1.14.1.1) \qquad \qquad \text{\rm(i)} \qquad
H_1&=z^{n_1}+y^{\alpha_{1,0,1}} \quad \text{with $n_1\ge 2$ and
$\alpha_{1,0,1}\ge 1$}. \qquad \qquad \qquad \qquad \\
\text{\rm(ii)} \qquad H_2&=H^{n_2}_1+y^{\alpha_{2,0,1}}z^{\alpha_{2,0,2}}.\\
\qquad  &\ldots\ldots \\
\text{\rm({$\ell$})} \qquad
H_{\ell}&=H^{n_{\ell}}_{\ell-1}+y^{\alpha_{\ell,0,1}}z^{\alpha_{\ell,0,2}}H^{\alpha_{\ell,0,3}}_1\cdots
H^{\alpha_{\ell,0,\ell}}_{\ell-2}. \\
\endalign$$

$\underline{\text{\bf {Conclusions}}}$

{\bf Fact[I]:} Then, ${f_k} \buildrel \text{{\rm multiseq }} \over
\sim H_{k}$ for each $k=1,2,\dots,\ell$.

{\bf Fact[II]:} Then, ${f_k} \buildrel \text{{\rm resol }} \over
\sim H_{k}$ for each $k=1,2,\dots,\ell$. \quad $\square$
\endproclaim \ms

{\bf $\S$1.8. The 2nd Algorithm for computing completely irreducible
W-polys from all the W-polys of two complex variables.} \ms

Noting that the statements of Theorem 1.15 and Corollary 1.15.1 are
different, observe by either Theorem 1.15 or Corollary 1.15.1 that
The 2nd Algorithm can be solved as follows:

$\underline{\text{\bf Assumptions}}$ Let $f\in \BC\{y\}[z]$ be an
arbitrary $W$-poly of degree $n\ge 2$ in $z$, satisfying the same
properties and notations as in the assumptions of Theorem
1.15(Corollary 1.15.1).

$\underline{\text{\bf Conclusions}}$

{\bf(i)} The aim of Theorem 1.15 is  to prove the following:

We can compute explicitly when $f$ is irreducible in $\BC\{y,z\}$.
If $f$ is irreducible in $\BC\{y,z\}$, we can find $H_{\ell} \in$
Family(1) such that $f\buildrel \text{{\rm multiseq}} \over \sim
H_{\ell}$ for any irreducible W-poly $f\in \BC\{y\}[z]$ $\iff$ by
the induction method on the positive integer $\ell$ we must follow
the computations over all the k-steps, $k=1,2,\dots,{\ell}$, in
Theorem 1.15. Note by Definition 1.2 that $H_{\ell}\in \BC\{y\}[z]$
is called the standard irreducible(Puiseux) W-poly of the recursive
${\ell}$-type in $z$. \ms

{\bf(ii)} The aim of Corollary 1.15.1 is to prove the following:

We can compute explicitly when $f$ is irreducible in $\BC\{y,z\}$.
If $f$ is irreducible in $\BC\{y,z\}$, we can find a generalized
representation of any irreducible W-poly $f\in \BC\{y\}[z]$ with
$f\buildrel \text{{\rm multiseq}} \over \sim H_{\ell}$ in the sense
of Theorem 1.13 $\iff$ by the induction method on the positive
integer $\ell$ we must follow the computations over all the k-steps,
$k=1,2,\dots,{\ell}$, in Corollary 1.15.1. \ms

\proclaim{Theorem 1.15(The 2nd Algorithm:Irreducibility algorithm for 
the Weierstrass polynomials of two complex variables)}

$\underline{\text{\bf Assumptions}}$ \quad Let $f\in \BC\{y\}[z]$ be
an arbitrary $W$-poly of degree $n\ge 2$ in $z$. Without loss of
generality, we may assume that $f$ satisfies the following form:
$$ f=z^n+\sum^{n-2}_{i=0} a_iy^{\a_i}z^i, \tag Eq.1 $$
where for $0\le i\le n-2$, each $a_i=a_i(y)$ is a unit in
${}_2\CO_0$ for $0\le i\le n-2$,  if exists, and the $\a_i$ are
positive integers. Note that $a_{n-1}$ is identically zero.  Write
$n=d_2n_1$ and $\a_0 =d_2\a_{1,0,1}$ with $d_2=\gcd(n,\a_0)$. Write
$n=\Pi^{r}_{k=1}n_k$ with positive integers $n_k\ge 2$ for all $k$
where the $n_k$ may not be the factorization of prime numbers.

In addition, assume that we have the following:
$$  2\le n\le \alpha_0. \tag 1.15.0 $$

$\underline{\text{\bf Conclusions}}$ \quad If $f$ is irreducible in
${}_2\CO_0$ with isolated singularity at $0\in \BC^2$, then $f\in $
the type $[\ell]$ for some $\ell\le r$ in the sense of Definition
2.5. By the induction method on the positive integer $r$, the aim is
to compute an elementary algorithm for finding irreducible W-polys
from all the W-polys in $\BC\{y\}[z]$, using $q$ iterations of the
following steps with $q\le \ell$: {\bf Observe that the statement on
the 3rd step may be omitted if necessary, to simplify the statements
for this theorem by the induction method}. \ms

$\underline{\text{\bf The 1st step: To find the irreducibility
algorithm for $f\in $ the type [1]}}$

$\underline{\text{\bf in the sense of  Definition 2.5.}}$ \ms

With equations in {\rm(Eq.2)}, the aim in this step is how to
compute the necessary and sufficient condition for $f\in $ the type
$[1]$ in the sense of Definition 2.5.

If $f$ is irreducible in ${}_2\CO_0$, then $f$ of {\rm(Eq.1)} must
satisfy the following necessary condition:
$$\align
 \f{\a_i}{n-i}\ge \f{\a_0}n \quad
\text{and} \quad{1\le\gcd(n,\a_0)<n} \quad \text{for \quad $0\le
i\le n-2$.} \tag Eq.2
\endalign$$

If $f$ satisfies {\rm(Eq.2)}, it suffices to consider the following
two cases for the 1st step:

$\underline{\text{\rm  Case(A)}}$  $\gcd(n,\a_0)=1$ and
$\underline{\text{\rm  Case(B)}}$  $1<\gcd(n,\a_0)<n$. \ms

$\underline{\text{\bf Case(A) of The 1st step:}}$ Let
$\gcd(n,\a_0)=1$. Then, $f$ is irreducible in ${}_2\CO_0$ $\iff$
{\rm(Eq.2)} holds. In this case, $f\in $ the type $[1]$ in the sense
of Definition 2.5. \ms

\text{\rm Remark for Case(A) of The 1st step}. The equation in
{\rm(Eq.1)} itself is a generalized representation of $f$ in the
sense of Theorem $1.13$, since $a_{n-1}=0$. \ms

$\underline{\text{\bf Case(B) of The 1st step:}}$ Let
$1<\gcd(n,\a_0)<n$. Then, take the 2nd step. $\square$ \bs

\noindent$\underline{\text{\bf Remark for The 1st step: To find $H_1
\in$ Family(1) such that $f\buildrel \text{{\rm multiseq}} \over
\sim H_1$ for }}$

\noindent$\underline{\text{\bf any $f\in $ the type [1] in the sense
of  Definition 2.5.}}$  \quad In Case(A) of The 1st step, $f\in $
the type $[1]$ in the sense of Definition 2.5 if and only if
$f\buildrel \text{{\rm multiseq}} \over \sim H_1$ where
$H_1=z^{n}+y^{\alpha_{0}}$ with $\gcd(n,\a_0)=1$. $\square$ \bs

$\underline{\text{\bf The 2nd step: To find the irreducibility
algorithm for $f\in $ the type [2]}}$

$\underline{\text{\bf in the sense of  Definition 2.5.}}$ \ms

With the equations, which are defined by {\rm(Eq.3)}, {\rm(Eq.4)}
and {\rm(Eq.5)} in the 2nd step later, the aim in this step is how
to compute the necessary and sufficient condition for $f\in $ the
type $[2]$ in the sense of Definition 2.5.

For the irreducibility algorithm for the 2nd step, it suffices to
consider Case(B) for the 1st step only. Let $1<\gcd(n,\a_0)<n$.

Then, note that $f\in $ the type [$\ell$] for some $\ell\ge 2$ in
the sense of Definition 2.5.

Let $d_2=\gcd(n,\a_0)$, and then write $n=n_1d_2$ and
$\a_0=\alpha_{1,0,1}d_2$. To solve the above problem, if $f$ is
irreducible in ${}_2\CO_0$, it is easy to compute by {\rm(1.12.3)}
of Lemma 1.12.1 that $f$ of {\rm(Eq.1)} must satisfy the following
necessary condition:
$$\align
\text{\rm(Eq.3)}\quad \text{\rm(a)} \quad &g_1=z^{n_1}+\xi_1y^{\alpha_{1,0,1}},  \\
 \text{\rm(b)}\quad &f =g_1^{d_2} +\sum_{\alpha, \beta\ge 0}
c_{\alpha,\beta}y^{\alpha}z^{\beta} \quad \text{with
$n_1\alpha+\alpha_{1,0,1}\beta>n_1\alpha_{1,0,1}d_2$, $\alpha>0$,
$\beta\le n-2$,} \qquad \quad
\endalign$$
where $\xi_1=\f{1}{d_2}a_{n-n_1}(0)$, $\alpha_{1,0,1}=\a_{n-n_1}$
and the $c_{\alpha,\beta}$ are nonzero complex numbers for some
nonnegative integers $\alpha>0$ and $\beta\le n-2$ such that
$n_1\alpha+\alpha_{1,0,1}\beta>n_1\alpha_{1,0,1}d_2$.

Then, $g_1$ must satisfy the following necessary condition:
$$
1<d_2<n\quad \text{and}\quad \xi_1 =\f {1}{d_{2}} a_{n-n_1}(0)\ne 0.
\tag Eq.3.1
$$
Note that $g_1$ is irreducible in ${}_2\CO_0$, since
$\gcd(n_1,\alpha_{1,0,1})=1$.

Apply the WDT(Theorem 1.7) to $f(y,z)$ with a divisor $g_1(y,z)$.
Whether or not $f$ is irreducible in ${}_2\CO_0$, by either
{\rm(1.12.4)} of Lemma 1.12.1 or {\rm(1.9.1)} of Sublemma 1.9,
$(g_1,f)$ can be written in the form
$$\cases
g_1 &=z^{n_1}+\xi_1 y^{\si_1}\ \qquad \qquad \text{{with} \quad
$\si_1=\a_{1,0,1}$ \ {and} \
$\xi_1 =\f {1}{d_{2}} a_{n-n_1}(0)$,} \\
f &=g^{d_{2}}_1+\sum^{d_{2}-1}_{i=1}T^{(1)}_{2,i}g^i_1,
\endcases \tag Eq.4)(Eq.4.1
$$
where $T^{(1)}_{2,i}\in \C\{y\}[z]$ for $0\le i\le d_2-1$,
satisfying the following {\rm(i)} and {\rm(ii)}:

{\rm(i)} For any nonzero monomial $y^{\g_1}z^{\g_2}$ in
$T^{(1)}_{2,i}$,
$$
\g_1>0\quad \text{and}\quad \g_2<n_1. \tag Eq.4.2
$$

{\rm(ii)} Define a function $\th_2 :\N^2_0\to \N_0$ by
$\th_2(t_k)^2_{k=1}=t_2\alpha_{1,0,1}+n_1t_1$.  Then, $\th_2$ is
one-to-one.

If $T^{(1)}_{2,i}\ne 0$, let
$C^{(1)}_{2,i}\Pi^2_{k=1}f^{\beta^{(1)}_{2,i,k}}_{k-2}$ be a unique
nonzero monomial with a constant $C^{(1)}_{2,i}$ in $T^{(1)}_{2,i}$
such that
$\th_2(\beta^{(1)}_{2,i,k})^2_{k=1}=\min\{\th_2(\g_k)^2_{k=1}\}$ for
any nonzero monomial $\Pi^2_{k=1}f^{\g_k}_{k-2}$ in $T^{(1)}_{2,i}$.
Note by Lemma 1.12.1 that the construction of $T^{(1)}_{2,i}$ is
trivial where $f_{-1}=y$ and $f_0=z$. \ms \ms

If $f$ is irreducible in ${}_2\CO_0$, then $f$ must satisfy the
following necessary condition:
$$
\f{\th_2(\beta^{(1)}_{2,i,k})^2_{k=1}}{d_2-i}\ge
\f{\th_2(\beta^{(1)}_{2,0,k})^2_{k=1}}{d_2}>n_1\alpha_{1,0,1} \quad
\text{for \quad $0\le i\le d_2-1$.} \tag Eq.4.3
$$

If $f$ satisfies {\rm(Eq.4.3)}, in order to find an irreducibility
algorithm for $f$, it suffices to consider the following two cases
for the 2nd step:

$\underline{\text{\rm Case(A)}}$
$\gcd(d_2,\th_2(\beta^{(1)}_{2,0,k})^2_{k=1})=1$ and
$\underline{\text{\rm Case(B)}}$
$\gcd(d_2,\th_2(\beta^{(1)}_{2,0,k})^2_{k=1})>1$. \ms

$\underline{\text{\bf Case(A) of the 2nd step:}}$ Let
$\gcd(d_2,\th_2(\beta^{(1)}_{2,0,k})^2_{k=1})=1$. Then, $f$ is
irreducible in ${}_2\CO_0$ $\iff$ {\rm(Eq.4.3)} holds.

In this case, $f\in $ the type $[2]$ in the sense of Definition 2.5.
So $f\buildrel \text{{\rm multiseq}} \over \sim
H^{d_2}_1+\Pi^2_{k=1}H^{\si_{2,k}}_{k-2}$ where $H_{-1}=y$, $H_0=z$,
$H_1=z^{n_1}+y^{\alpha_{1,0,1}}$ and $\si_{2,k}=\beta^{(1)}_{2,0,k}$
for $1\le k\le 2$. \ms

\text{\rm Remark for Case(A) of the 2nd step}. Note that $(g_1,f)$
in {\rm(Eq.4.1)} may not be a generalized representation of $f$ in
the sense of Theorem $1.13$, since $T^{(1)}_{2,d_2-1}$ may not be
zero. To find a generalized representation of $f$ in the sense of
Theorem $1.13$, until we get the same result in the conclusion of
Sublemma 1.10(Sublemma 15.5), it suffices to apply either Sublemma
1.10(Sublemma 15.5) or the same kind of a fine sequence of pairs in
{\rm(Eq.5)} to $(g_1,f)$. \ms

$\underline{\text{\bf Case(B) of the 2nd step:}}$ Let
$\gcd(d_2,\th_2(\beta^{(1)}_{2,0,k})^2_{k=1})>1$. To find a
necessary condition for $f$ to be irreducible in ${}_2\CO_0$, we may
assume that $f$ satisfies {\rm(Eq.2)}, {\rm(Eq.3)} and {\rm(Eq.4)}.

Then, it suffices to consider the following two subcases of Case(B)
for the 2nd step:

$\underline{\text{\rm Subcase(B1)}}$ $T^{(1)}_{2,d_2-1}\ne 0$, and
$\underline{\text{\rm Subcase(B2)}}$ $T^{(1)}_{2,d_2-1}=0$. \ms

$\underline{\text{\bf Subcase(B1) of the 2nd step:}}$ Assume that
$T^{(1)}_{2,d_2-1}\ne 0$. To solve the subcase, first of all, we
must eliminate $T^{(1)}_{2,d_2-1}$ whether or not $f$ is irreducible
in ${}_2\CO_0$. To do it, by Sublemma 1.10 for Theorem 1.8.1
(Sublemma 15.5.$\beta$ of Sublemma 15.5 of the 1st statement for Theorem 15.4), 
we can compute a unique finite sequence of
pairs $\{(h_{1,p},f):1\le p\le {\nu_1}+1\}$ in {\rm(Eq.5)} for a
unique integer ${\nu_1}\le \f{n_1+1}2$, each pair of which can be
written in the form
$$\cases
h_{1,1} &=z^{n_1}+\xi_1y^{\alpha_{1,0,1}} =z^{n_1}+R^{(1)}_{1,0}\
\text{with}\
h_{1,1}=g_1\\
f &=h^{d_2}_{1,1} +\sum^{d_2-1}_{i=0}T^{(1)}_{2,i}h^i_{1,1},
\endcases \tag {Eq.5)(Eq.5.1)(Eq.5.1.1}
$$
and {\quad} for $1\le p\le {\nu_1}-1$,
$$\cases
h_{1,p+1} &=h_{1,p}+\f
1{d_2}T^{(p)}_{2,d_2-1}=z^{n_1}+\sum^{n_1-2}_{i=0}R^{(p+1)}_{1,i}z^i\\
f &=h^{d_2}_{1,p+1}+\sum^{d_2-1}_{i=0}T^{(p+1)}_{2,i}h^i_{1,p+1}\
\text{with}\ T^{(p+1)}_{2,d_2-1}\not=0,
\endcases \tag Eq.5.1.2
$$
and
$$\cases
h_{1,{\nu_1}+1} &=h_{1,{\nu_1}} +\f
1{d_2}T^{({\nu_1})}_{2,d_2-1}=z^{n_1}+\sum^{n_1-2}_{i=0}R^{({\nu_1}+1)}_{1,i}z^i\\
f
&=h^{d_2}_{1,{\nu_1}+1}+\sum^{d_2-2}_{i=0}T^{({\nu_1}+1)}_{2,i}h^i_{1,{\nu_1}+1}\
\text{with}\ T^{({\nu_1}+1)}_{2,d_2-1}=0,
\endcases \tag Eq.5.1.3
$$
where $R^{(p)}_{1,i}\in \C\{y\}$ for $1\le p\le {\nu_1}+1$ and $0\le
i\le n_1-2$; and $T^{(p)}_{2,,i}\in \C\{y\}[z]$ for $1\le p\le
{\nu_1}+1$ and $0\le i\le d_2-1$, satisfying the following
properties:

{\rm(i)} For any nonzero monomial $y^\de$ in $R^{(p)}_{1,i},\ \de
>0$.

{\rm(ii)} For any nonzero monomial $y^{\g_1}z^{\g_2}$ in
$T^{(p)}_{2,i}$,
$$
\g_1>0 \quad \text{and} \quad \g_2<n_1. \tag Eq.5.2
$$

{\rm(iii)} If $T^{(p)}_{2,i}\ne 0$, let
$C^{(p)}_{2,i}\Pi^2_{k=1}f^{\beta^{(p)}_{2,i,k}}_{k-2}$ be a unique
nonzero monomial with a constant $C^{(p)}_{2,i}$ in $T^{(p)}_{2,i}$
such that
$\th_2(\beta^{(p)}_{2,i,k})^2_{k=1}=\min\{\th_2(\g_k)^2_{k=1}\}$ for
any nonzero monomial $\Pi^2_{k=1}f^{\g_k}_{k-2}$ in $T^{(p)}_{2,i}$
where $f_{-1}=y$ and $f_0=z$.

\noindent{\bf {Remark.}} It is easy to compute
$C^{(p)}_{2,i}\Pi^2_{k=1}f^{\beta^{(p)}_{2,i,k}}_{k-2}$ by an
elementary way because $T^{(p)}_{2,i}\in \C\{y\}[z]$ is a polynomial
of finite degree $<n_1$ in $z$  and $\th_2$ is one to one. \ms

If $f$ is irreducible in ${}_2\CO_0$, then $f$ satisfy the following
necessary condition:
$$\align
\text{\rm(Eq.5.3)(Eq.5.3.1)} \quad
&\f{\th_2(\beta^{(p)}_{2,i,k})^2_{k=1}}{d_2-i} \ge
\f{\th_2(\beta^{(p)}_{2,0,k})^2_{k=1}}{d_2}>n_1\alpha_{1,0,1} \quad
\text{for $1\le p\le {\nu_1}$, $0\le i\le {d_2}-1$}.  \qquad \qquad \\
\text{\rm(Eq.5.3.2)}\quad &
\f{\th_2(\beta^{({\nu_1}+1)}_{2,i,k})^2_{k=1}}{d_2-i} \ge
\f{\th_2(\beta^{({\nu_1}+1)}_{2,0,k})^2_{k=1}}{d_2}>n_1\alpha_{1,0,1}
\quad \text{for $0\le i\le d_2-2$.}
\endalign$$

\noindent{\bf {Remark.}} If $f$ satisfies {\rm(Eq.5.3)},
$h_{1,p}\buildrel \text{{\rm multiseq}} \over \sim h_{1,1}=g_1$ for
all $p\ge 1$. \ms

Assuming that $f$ satisfies {\rm(Eq.5.3)}, to find an irreducibility
algorithm for Subcase(B1), it suffices to consider the following two
subcases, Subcase(B1-a) and Subcase(B1-b):

$\underline{\text{\rm Subcase(B1-a)}}$ Let
$\gcd(d_2,\th_2(\beta^{(1)}_{2,0,k})^2_{k=1})>1$ and
$\gcd(d_2,\th_2(\beta^{({\nu_1}+1)}_{2,0,k})^2_{k=1})=1$.

$\underline{\text{\rm Subcase(B1-b)}}$ Let
$\gcd(d_2,\th_2(\beta^{({\nu_1}+1)}_{2,0,k})^2_{k=1})>1$. \ms

$\underline{\text{\bf Subcase(B1-a) of the 2nd step}}$

{\rm(a)} Let $\gcd(d_2,\th_2(\beta^{(1)}_{2,0,k})^2_{k=1})>1$ and
$\gcd(d_2,\th_2(\beta^{({\nu_1}+1)}_{2,0,k})^2_{k=1})=1$. Then, $f$
is irreducible in ${}_2\CO_0$ if and only if {\rm(Eq.5.3.2)} holds.
Thus, $f\in $ the type $[2]$ in the sense of Definition 2.5.

{\rm(b)} Assume that
$\gcd(d_2,\th_2(\beta^{({p}+1)}_{2,0,k})^2_{k=1})=1$ for some $p\le
{\nu_1}+1$. Then, $f$ is irreducible in ${}_2\CO_0$ if and only if
an inequality in {\rm(Eq.5.3.1)} holds without mentioning any
inequality in {\rm(Eq.5.3.2)}. \ms

\text{\rm Remark for Subcase(B1-a) of the 2nd step}.
$(f_1,f)=(h_{1,{\nu_1}+1},f)$ is a generalized representation of $f$
in the sense of Theorem $1.13$ since $T^{({\nu_1}+1)}_{2,d_2-1}=0$.
\ms

$\underline{\text{\bf Subcase(B1-b) of the 2nd step}}$ Let
$\gcd(d_2,\th_2(\beta^{({\nu_1}+1)}_{2,0,k})^2_{k=1})>1$, noting
that $T^{({\nu_1}+1)}_{2,d_2-1}=0$. To find a necessary and
sufficient condition for $f$ to be irreducible in ${}_2\CO_0$, take
the next step. \ms

$\underline{\text{\bf Subcase(B2) of the 2nd step}}$ It was assumed
by this subcase that $T^{(1)}_{2,d_2-1}=0$, noting that
$\gcd(d_2,\th_2(\beta^{(1)}_{2,0,k})^2_{k=1})>1$. To find a
necessary and sufficient condition for $f$ to be irreducible in
${}_2\CO_0$, take the next step. $\square$ \ms

\noindent$\underline{\text{\bf Remark for The 2nd step: To find $H_2
\in$ Family(1) such that $f\buildrel \text{{\rm multiseq}} \over
\sim H_2$ for }}$

\noindent$\underline{\text{\bf any $f\in $ the type [2] in the sense
of Definition 2.5.}}$ \quad In this case, $f\in $ the type $[2]$ in
the sense of Definition 2.5, and so $f\buildrel \text{{\rm
multiseq}} \over \sim H_2=H^{d_2}_1+\Pi^2_{k=1}H^{\si_{2,k}}_{k-2}$
where $H_{-1}=y$, $H_0=z$, $H_1=z^{n_1}+y^{\alpha_{1,0,1}}$ and
$\si_{2,k}=\beta^{({\nu_1}+1)}_{2,0,k}$ for $1\le k\le 2$ and for
some positive integer ${\nu_1}+1$. $\square$ \bs

$\underline{\text{\bf The 3rd step: To find the irreducibility
algorithm for $f\in $ the type [3]}}$

$\underline{\text{\bf in the sense of  Definition 2.5.}}$ \ms

With three equations {\rm(Eq.6)}, {\rm(Eq.7)} and {\rm(Eq.8)} in the
3rd step, the aim in this step is to find the necessary and
sufficient condition for $f\in $ the type $[3]$ in the sense of
Definition 2.5.

For the proof of this step, recall the defining equation
$(h_{1,{\nu_1}+1},f)$ of {\rm(Eq.5)} as we have seen in Subcase(B1)
of the 2nd step:
$$\cases
h_{1,{\nu_1}+1} &=z^{n_1}+\sum^{n_1-2}_{i=0}R^{({\nu_1}+1)}_{1,i}z^i\\
f
&=h^{d_2}_{1,{\nu_1}+1}+\sum^{d_2-2}_{i=0}T^{({\nu_1}+1)}_{2,i}h^i_{1,{\nu_1}+1}\
\text{with}\ T^{({\nu_1}+1)}_{2,d_2-1}=0.
\endcases \tag Eq.5
$$

By either Subcase(B1-b) or Subcase(B2) of the 2nd step, it may be
assumed that

$\gcd(d_2,\th_2(\beta^{({\nu_1}+1)}_{2,0,k})^2_{k=1})>1$ and
$T^{({\nu_1}+1)}_{2,d_2-1}=0$ for some positive integer ${\nu_1}+1$.
\ms

\noindent{\bf Remark.} If $T^{(1)}_{2,d_2-1}=0$ and
$\gcd(d_2,\th_2(\beta^{(1)}_{2,0,k})^2_{k=1})>1$ from $(h_{1,1},f)$
of {\rm(Eq.5.1)} , then $(h_{1,1},f)=(g_1,f)$ can be viewed as
$(h_{1,{\nu_1}+1},f)$ of {\rm(Eq.5)} with ${\nu_1}=0$. \ms

For brevity of notation, $h_{1,{\nu_1}+1},T^{({\nu_1}+1)}_{2,i}$,
$C^{({\nu_1}+1)}_{2,i}$ and $\beta^{({\nu_1}+1)}_{2,i,k}$ can be
replaced by $f_1,T_{2,i},C_{2,i}$ and $\beta_{2,i,k}$ for $0\le i\le
d_2-2$ and $1\le k\le 2$, respectively because $g_1=h_{1,1}\buildrel
\text{{\rm multiseq}} \over \sim h_{1,p}$ for all $p\ge 1$ and
$T^{({\nu_1}+1)}_{2,d_2-1}=0$. Recall that
$f_1=z^{n_1}+\sum^{n_1-2}_{i=0}R_{1,i}z^i$ from
$(h_{1,{\nu_1}+1},f)$ where $R^{({\nu_1}+1)}_{1,i}$ is defined to be
$R_{1,i}$ for each $i$.

Let $d_3=\gcd(d_2,\th_2(\beta^{(1)}_{2,0,k})^2_{k=1})$, and then
write $d_2=n_2d_3$. To solve the above problem, we need to construct
$g_2$ as follows:
$$
g_2=f^{n_2}_1 +\xi_2y^{\si_{2,1}}z^{\si_{2,2}}\in \C\{y\}[z,f_1],
\tag Eq.6
$$
where $\xi_2=\f 1{d_3}C_{2,d_2-n_2}$ and
$\si_{2,k}=\beta_{2,d_2-n_{2,k}}$ for $1\le k\le 2$.

Then, $g_2$ must satisfy the following necessary condition:
$$
1<d_3<d_2,\ \xi_2\ne 0,\
\th_2(\si_{2,k})^2_{k=1}>n_2n_1\beta_{1,0,1}, \tag Eq.6.1
$$
and $\gcd(n_2, \th_2(\si_{2,k})^2_{k=1})=1$.

If $g_2$ of {\rm(Eq.6)} satisfies inequalities in {\rm(Eq.6.1)},
then $g_2(y,z)$ is irreducible in ${}_2\CO_0$ with $g_2\in $ the
type $[2]$ in the sense of Definition 2.5 because $g_2(y,z,f_1)$ can
be viewed as an element of ${}_2\CO_0$.

Apply the WDT to $f(y,z)$ with a divisor $g_2(y,z)$. Then, $(g_2,f)$
can be written in the form
$$\cases
g_2 &=f^{n_2}_1 +\xi_2y^{\si_{2,1}}z^{\si_{2,2}}\\
f &=g^{d_3}_2 +\sum^{d_3-1}_{i=0}T^{(1)}_{3,i}g^i_2,
\endcases \tag Eq.7)(Eq.7.1
$$
where $T^{(1)}_{3,i}\in \C\{y\}[z,f_1]$ for $0\le i\le d_3-1$,
satisfying the following property:

{\rm(i)} For any nonzero monomial $\Pi^3_{k=1}f^{\g_k}_{k-2}$ in
$T^{(1)}_{3,i}$,
$$
\g_1>0,\ \g_2<n_1\quad \text{and}\quad \g_3<n_2. \tag Eq.7.2
$$

{\rm(ii)} Define $\th_3 :\N^{(3)}_0\to \N_0$ by
$\th_3(t_k)^3_{k=1}=t_3\th_2(\si_{2,k})^2_{k=1}+n_2\th_2(t_k)^2_{k=1}$,
and then $\th_3$ is one-to-one.

If $T^{(1)}_{3,i}\ne 0$, let
$C^{(1)}_{3,i}\Pi^3_{k=1}f^{\beta^{(1)}_{3,i,k}}_{k-2}$ be a unique
nonzero monomial with a constant $C^{(1)}_{3,i}$ in $T^{(1)}_{3,i}$
such that
$\th_3(\beta^{(1)}_{3,i,k})^3_{k=1}=\min\{\th_3(\g_k)^3_{k=1}\}$ for
any nonzero monomial $\Pi^3_{k=1}f^{\g_k}_{k-2}$ in $T^{(1)}_{3,i}$.
\ms

\noindent{\bf Remark.} The above construction of $T^{(1)}_{3,i}$ is
as follows: Note that $g_2(y,z)$ is a polynomial in $z$ of degree
$\Pi^2_{k=1}n_k$, and for $0\le i\le d_3-1$, $T^{(1)}_{3,i}\in
\C\{y\}[z]$ is a polynomial in $z$ of degree $<\Pi^2_{k=1}n_k$, and
$f_1(y,z)\in \C\{y\}[z]$ is a polynomial in $z$ of degree $n_1$.
Again, apply the WDT to each $T^{(1)}_{3,i}$ with a divisor
$f_1(y,z)$. Then $T^{(1)}_{3,i}$ can be written just as above.
$\square$ \ms

Since $f$ is irreducible in ${}_2\CO_0$, $f$ must satisfy the
following necessary condition:
$$
\f{\th_3(\beta^{(1)}_{3,i,k})^3_{k=1}}{d_3-i}\ge
\f{\th_3(\beta^{(1)}_{3,0,k})^3_{k=1}}{d_3}>n_2\th_2(\sigma_{2,k})^2_{k=1}.
\text{\quad for $0\le i\le d_3-1$.} \tag Eq.7.3
$$

Recall that $\si_{2,k}=\beta_{2,d_2-n_2,k}$ for $1\le k\le 2$ from
the construction of $g_2$.

Since $f$ satisfies {\rm(Eq.7.3)}, then to find an irreducibility
algorithm for $f$, it suffices to consider the following two
subcases Case(A) and Case(B):

$\underline{\text{\rm Case(A)}}$
$\gcd(d_3,\th_3(\beta^{(1)}_{3,0,k})^3_{k=1})=1$, and
$\underline{\text{\rm Case(B) }}$
$\gcd(d_3,\th_3(\beta^{(1)}_{3,0,k})^3_{k=1})>1$. \ms

$\underline{\text{\bf Case(A) of the 3rd step}}$ Let
$\gcd(d_3,\th_3(\beta^{(1)}_{3,0,k})^3_{k=1})=1$. Then, $f$ is
irreducible in ${}_2\CO_0$ if and only if {\rm(Eq.7.2)} holds. In
this case, if $f$ is irreducible in ${}_2\CO_0$, then $f\in $ the
type $[3]$ in the sense of Definition 2.5 and also $f\buildrel
\text{{\rm multiseq}} \over \sim
H^{d_3}_2+\Pi^3_{k=1}H^{\si_{3,k}}_{k-2}$ where $H_{-1}=y$, $H_0=z$,
$H_1=z^{n_1}+y^{\beta_{1,0,1}}$ and
$H_2=H^{n_2}_1+y^{\si_{2,1}}z^{\si_{2,2}}$ and
$\si_{3,k}=\beta^{(1)}_{3,0,k}$ for $1\le k\le 3$. \ms

\text{\rm Remark for Case(A) of the 3rd step}. Let
$f_1=h_{1,{\nu_1}+1}$ in {\rm(Eq.5)}, and $(g_2,f)$ in {\rm(Eq.7)}.
Then, $(f_1,g_2,f)$ may not be a generalized representation of $f$
in the sense of Theorem $16.5$, since $T^{(1)}_{3,d_3-1}$ may not be
zero. To find a generalized representation of $f$ in the sense of
Theorem $16.5$, until we get the same kind of result in the
conclusion of Sublemma 15.5, it suffices to apply either Sublemma
15.5 or the same kind of a fine sequence of pairs in {\rm(Eq.8)} to
$(g_2,f)$. \ms

$\underline{\text{\bf Case(B) of the 3rd step:}}$ Let
$\gcd(d_3,\th_3(\beta^{(1)}_{3,0,k})^2_{k=1})>1$. To find a
necessary condition for $f$ to be irreducible in ${}_2\CO_0$, we may
assume that $f$ satisfies {\rm(Eq.7)} with
{\rm(Eq.2)},{\rm(Eq.3)},\dots,{\rm(Eq.6)}.

Then, it suffices to consider the following two subcases:

{\rm Subcase(B1)} $T^{(1)}_{3,d_3-1}\ne 0$, and {\rm Subcase(B2)}
$T^{(1)}_{3,d_3-1}=0$. \ms

$\underline{\text{\bf Subcase(B1) of the 3rd step:}}$ Let
$T^{(1)}_{3,d_3-1}\ne 0$. To solve the case, first of all, we must
eliminate $T^{(1)}_{3,d_3-1}$ whether or not $f$ is irreducible in
${}_2\CO_0$. To do it, by Theorem 1.8(Theorem 15.4), we can compute
a unique finite sequence of pairs, $\{(h_{2,p},f):1\le p\le
{\nu_2}+1\}$ in {\rm(Eq.8)} for a unique integer ${\nu_2} \le
\f{n_2+1}2$, each pair of which can be written in the form
$$\cases
h_{2,1} &=f^{n_2}_1 +\xi_2y^{\si_{2,1}}z^{\si_{2,2}}
=f^{n_2}_1+R^{(1)}_{2,0}\ \text{with $h_{2,1}=g_2$} \\
f &=h^{d_3}_{2,1} +\sum^{d_3-1}_{i=0}T^{(1)}_{3,i}h^i_{2,1},
\endcases \tag Eq.8)(Eq.8.1)(Eq.8.1.1
$$
and for $1\le p\le{\nu_2}-1$,
$$\cases
h_{2,p+1} &=h_{2p}+\f
1{d_3}T^{(p)}_{3,d_3-1}=f^{n_2}_1+\sum^{n_2-2}_{i=0}R^{(p+1)}_{2i}f^i_1\\
f &=h^{d_3}_{2,p+1}+\sum^{d_3-1}_{i=0}T^{(p+1)}_{3,i}h^i_{2,p+1}\
\text{with}\ T^{(p+1)}_{3,d_3-1}\not=0,
\endcases \tag Eq.8.1.2
$$
and
$$\cases
h_{2,{\nu_2}+1} &=h_{2,{\nu_2}} +\f
1{d_3}T^{({\nu_2})}_{3,d_3-1}=f^{n_2}_1
+\sum^{n_2-2}_{i=0}R^{({\nu_2}+1)}_{2,i}f^i_1\\
f
&=h^{d_3}_{2,{\nu_2}+1}+\sum^{d_3-2}_{i=0}T^{({\nu_2}+1)}_{3,i}h^i_{2,{\nu_2}+1}\
\text{with}\ T^{({\nu_2}+1)}_{3,d_3-1}=0,
\endcases \tag Eq.8.1.3
$$
where $R^{(p)}_{2,i}\in \C\{y\}[z]$ for $1\le p\le {\nu_2}+1$ and
$0\le i\le n_2-2$; and $T^{(p)}_{3,i}\in \C\{y\}[z,f_1]$ for $1\le
p\le {\nu_2}+1$ and $0\le i\le d_3-1$, satisfying the following
properties: \roster \item"(i)" For any nonzero monomial
$y^{\de_1}z^{\de_2}$ in $R^{(p)}_{2,i}$, $\de_1>0$ and $\de_2<n_1$.
\item"(ii)" For any nonzero monomial $\Pi^3_{k=1}f^{\g_k}_{k-2}$ in
$T^{(p)}_{3,i}$,
$$
\g_1>0,\ \g_2<n_1\quad \text{and}\quad \g_3<n_2. \tag Eq.8.2
$$
\item"(iii)" If $T^{(p)}_{3,i}\ne 0$, let
$C^{(p)}_{3,i}\Pi^3_{k=1}f^{\beta^{(p)}_{3,i,k}}_{k-2}$ be a unique
nonzero monomial with a constant $C^{(p)}_{3,i}$ in $T^{(p)}_{3,i}$
such that
$\th_3(\beta^{(p)}_{3,i,k})^3_{k=1}=\text{min}\{\th_3(\g_k)^3_1\}$
for any nonzero monomial $\Pi^3_{k=1}f^{\g_k}_{k-2}$ in
$T^{(p)}_{3,i}$ where $f_{-1}=y$ and $f_0=z$.
\endroster

\noindent{\bf Remark.} It is possible to compute
$C^{(p)}_{3,i}\Pi^3_{k=1}f^{\beta^{(p)}_{3,i,k}}_{k-2}$ by an
elementary way because $T^{(p)}_{3,i}\in \C\{y\}[z,f_1]$ in $z$ and
$f_1$, recalling that for any nonzero monomial
$\Pi^3_{k=1}f^{\g_k}_{k-2}$ in $T^{(p)}_{3,i}$, $\g_2<n_1$ and
$\g_3<n_2$ and $\th_3$ is one-to-one. \ms

Now, Since $f$ is irreducible in ${}_2\CO_0$, then $f$ must satisfy
the following necessary condition:
$$\align
\text{\rm(Eq.8.3)(Eq.8.3.1)} \quad
\f{\th_3(\beta^{(p)}_{3,i,k})^3_{k=1}}{d_3-i}& \ge
\f{\th_3(\beta^{(p)}_{3,0,k})^3_{k=1}}{d_3}>n_2\th_2(\si_{2,k})^2_{k=1}
\text{\quad for $1\le p\le {\nu_2}$,  $0\le i\le d_3-1$.}  \qquad \qquad \\
\text{\rm(Eq.8.3.2)}
\f{\th_3(\beta^{({\nu_2}+1)}_{3,i,k})^3_{k=1}}{d_3-i}& \ge
\f{\th_3(\beta^{({\nu_2}+1)}_{3,0,k})^3_{k=1}}{d_3}>n_2\th_2(\si_{2,k})^2_{k=1}
\text{\quad for $0\le i\le d_3-2$.}  \qquad \\
\endalign$$

\noindent{\bf Remark.} Since $f$ satisfies {\rm(Eq.8.3)}, then
$h_{2,p}\buildrel \text{{\rm multiseq}} \over \sim h_{2,1}=g_2$ for
all $p\ge 1$. \ms

Since $f$ satisfies {\rm(Eq.8.3)}, to find an irreducibility
algorithm for Subcase(B1), then there exist the following two
possibilities, Subcase(B1-a) and Subcase(B1-b):

$\underline{\text{\rm Subcase(B1-a)}}$ Let
$\gcd(d_3,\th_3(\beta^{(1)}_{3,0,k})^3_{k=1})>1$ and
$\gcd(d_3,\th_3(\beta^{({\nu_2}+1)}_{3,0,k})^3_{k=1})=1$.

$\underline{\text{\rm Subcase(B1-b)}}$ Let
$\gcd(d_3,\th_3(\beta^{({\nu_2}+1)}_{3,0,k})^3_{k=1})>1$. \ms

$\underline{\text{\bf Subcase(B1-a) of the 3rd step}}$

{\rm(a)} Let $\gcd(d_3,\th_3(\beta^{(1)}_{3,0,k})^3_{k=1})>1$ and
$\gcd(d_3,\th_3(\beta^{({\nu_2}+1)}_{3,0,k})^3_{k=1})=1$. Then, $f$
is irreducible in ${}_2\CO_0$ if and only if {\rm(Eq.8.3.2)} holds.
Thus, $f\in $ the type $[3]$ in the sense of Definition 2.5.

{\rm(b)} Assume that
$\gcd(d_3,\th_3(\beta^{({p}+1)}_{3,0,k})^3_{k=1})=1$ for some $p\le
{\nu_2}+1$. Then, $f$ is irreducible in ${}_2\CO_0$ if and only if
an inequality in {\rm(Eq.8.3.1)} holds without mentioning any
inequality in {\rm(Eq.8.3.2)}. \ms

\text{\rm Remark for Subcase(B1-a) of the 3rd step}.
$(f_1,f_2,f)=(f_1,h_{2,{\nu_2}+1},f)$ is a generalized
representation of $f$ in the sense of Theorem $1.13$, since
$T^{({\nu_2}+1)}_{2,d_2-1}=0$. \ms

$\underline{\text{\bf Subcase(B1-b) of the 3rd step}}$ Let
$\gcd(d_3,\th_3(\beta^{({\nu_2}+1)}_{3,0,k})^3_{k=1})>1$, noting
that $T^{({\nu_2}+1)}_{3,d_3-1}=0$. To find a necessary and
sufficient condition for $f$ to be irreducible in ${}_2\CO_0$, take
the next step. \ms

$\underline{\text{\bf Subcase(B2) of the 3rd step}}$ It was assumed
by this subcase that $T^{(1)}_{3,d_3-1}=0$, noting that
$\gcd(d_3,\th_3(\beta^{(1)}_{3,0,k})^3_{k=1})>1$. To find a
necessary and sufficient condition for $f$ to be irreducible in
${}_2\CO_0$, take the next step.  $\square$ \ms

\noindent$\underline{\text{\bf Remark for The 3rd step: To find $H_3
\in$ Family(1) such that $f\buildrel \text{{\rm multiseq}} \over
\sim H_3$ for }}$

\noindent$\underline{\text{\bf any $f\in $ the type [3] in the sense
of  Definition 2.5.}}$ \quad Let
$\gcd(d_3,\th_3(\beta^{({\nu_2}+1)}_{3,0,k})^3_{k=1})=1$. Then, $f$
is irreducible in ${}_2\CO_0$ if and only if {\rm(Eq.7.2)} holds. In
this case, if $f$ is irreducible in ${}_2\CO_0$, then $f\in $ the
type $[3]$ in the sense of Definition 2.5 and also $f\buildrel
\text{{\rm multiseq}} \over \sim
H^{d_3}_2+\Pi^3_{k=1}H^{\si_{3,k}}_{k-2}$ where $H_{-1}=y$, $H_0=z$,
$H_1=z^{n_1}+y^{\beta_{1,0,1}}$ and
$H_2=H^{n_2}_1+y^{\si_{2,1}}z^{\si_{2,2}}$ and
$\si_{3,k}=\beta^{(\nu+1)}_{3,0,k}$ for $1\le k\le 3$ and for some
positive integer ${\nu_2}+1$.  $\square$ \bs

The general case will be proved by induction. Let $f\in \BC\{y\}[z]$
be an arbitrary $W$-poly of degree $n\ge 2$ in $z$. Suppose we have
shown for any integer $q\ge 3$ that for each \text{\rm(j)-th} step,
$1\le j\le q-1$, the irreducibility algorithm for $f\in $ the type
$[j]$ has been found, by using the same kind of properties and
notations as we have seen in the proof of finding the irreducibility
algorithm for the \text{\rm(j)-th} step, $1\le j\le q-1$. Now, we
may assume by induction on $q$ that the irreducibility algorithm for
the \text{\rm(q-1)-th} step can be solvable for $q\ge 3$, and so if
$f$ is irreducible in ${}_2\CO_0$, we may assume without loss of
generality that if $f\in $ the type $[k]$ for an integer $k\ge q$ in
the sense of Definition 2.5 then
$\gcd(d_{q-1},\th_{q-1}(\beta^{({\nu_{q-2}}+1)}_{q-1,0,k})^{q-1}_{k=1})>1$
and $T^{({\nu_{q-2}}+1)}_{q-1,d_{{q-1}-1}}=0$ for some positive
integer ${\nu_{q-2}}+1$. Then, it remains to show that the
generalized irreducibility algorithm for $f\in $ the type $[q]$ can
be written as follows. \ms

$\underline{\text{\bf The q-th step: To find the irreducibility
algorithm for $f\in $ the type [q]}}$

$\underline{\text{\bf in the sense of  Definition 2.5.}}$ \ms

With three equations, which is defined by {\rm(Eq.(3q-3))},
{\rm(Eq.(3q-2))} and {\rm(Eq.(3q-1))} in the q-th step later, the
aim in this step is how to compute the necessary and sufficient
condition for $f\in $ the type $[q]$ in the sense of Definition 2.5.

Suppose by induction on the positive integer (q-1) that $f$ is
irreducible in ${}_2\CO_0$ and let $f$ satisfy a finite number
(3q-4) of conditions, which have been represented by {\rm(Eq.1)},
{\rm(Eq.2)}, {\rm(Eq.3)}, $\dots$, {\rm(Eq.{(3q-6)})},
{\rm(Eq.{(3q-5)})},{\rm(Eq.{(3q-4)})}. \ms

For the proof of this step, recall the defining equation
$(h_{{q-2},{\nu_{q-2}}+1},f)$ of {\rm(Eq.{(3q-4)})} as we have
already seen in Subcase(B1) of the {(q-1)}-th step:
$$\cases
h_{{q-2},{\nu_{q-2}}+1} &=f^{n_{q-2}}_{q-3}+\sum^{n_{q-2}-2}_{i=0}R^{({\nu_{q-2}}+1)}_{q-2,i}f^i_{q-3}\\
f
&=h^{d_{q-1}}_{{q-2},{\nu_{q-2}}+1}+\sum^{d_{q-1}-2}_{i=0}T^{({\nu_{q-2}}+1)}_{q-1,i}h^i_{{q-2},{\nu_{q-2}}+1}\
\text{with}\ T^{({\nu_{q-2}}+1)}_{q-1,d_{q-1}-1}=0.
\endcases \tag Eq.(3q-4)
$$

By either Subcase(B1-b) or Subcase(B2) of the {(q-1)}-th step, it
may be assumed that
$\gcd(d_{q-1},\th_{q-1}(\beta^{({\nu_{q-2}}+1)}_{{q-1},0,k})^{q-1}_{k=1})>1$
and $T^{({\nu_{q-2}}+1)}_{{(q-1)},d_{(q-1)}-1}=0$ for some positive
integer ${\nu_{q-2}}+1$.

\noindent{\bf Remark.} If $T^{(1)}_{{q-1},d_{q-1}-1}=0$ and
$\gcd(d_{q-1},\th_{q-1}(\beta^{(1)}_{{q-1},0,k})^{q-1}_{k=1})>1$ for
Subcase(B2), then $(h_{{q-2},{1}},f)=(g_{q-2},f)$ can be viewed as
$(h_{q-2,\nu+1},f)$ with ${\nu_{q-2}}=0$. \ms

For brevity of notation,
$h_{q-2,{\nu_{q-2}}+1},T^{({\nu_{q-2}}+1)}_{q-1,i}$,
$C^{({\nu_{q-2}}+1)}_{q-1,i}$ and
$\beta^{({\nu_{q-2}}+1)}_{q-1,i,k}$ can be replaced by
$f_{q-2},T_{q-1,i},C_{q-1,i}$ and $\beta_{q-1,i,k}$ for $0\le i\le
d_{q-1}-2$ and $1\le k\le q-1$, respectively because
$g_{q-2}=h_{{q-2},1}\buildrel \text{{\rm multiseq}} \over \sim
h_{{q-2},p}$ for all $p\ge 1$ and
$T^{({\nu_{q-2}}+1)}_{q-1,d_{q-1}-1}=0$. Recall that
$f_{q-2}=f^{n_{q-2}}_{q-3}+\sum^{n_{q-2}-2}_{i=0}R^{({\nu_{q-2}}+1)}_{q-2,i}f^i_{q-3}$
from $(h_{{q-2},{\nu_{q-2}}+1},f)$ where
$R^{({\nu_{q-2}}+1)}_{{q-2},i}$ is defined to be $R_{{q-2},i}$ for
each $i$.

Let $d_{q}=\gcd(d_{q-1},\th_{q-1}(\beta_{{q-1},0,k})^{q-1}_{k=1})$
and then write $d_{q-1}=n_{q-1}d_{q}$. To solve the above problem,
we need to construct $g_{q-1}$ as follows:
$$
g_{q-1}=f^{n_{q-1}}_{q-2}
+\xi_{q-1}\Pi^{q-1}_{k=1}f^{\si_{{q-1},k}}_{k-2}\in
\C\{y\}[z,f_1,\dots,f_{q-2}], \tag Eq.(3q-3)
$$
where $\xi_{q-1}=\f 1{d_{q}}C_{{q-1},d_{q-1}-n_{q-1}}$ and
$\si_{{q-1},k}=\beta_{{q-1},d_{q-1}-n_{q-1},k}$ for $1\le k\le
{q-1}$.

Then, $g_{q-1}$ must satisfy the following necessary condition:
$$
1<d_{q}<d_{q-1},\ \xi_{q-1}\ne 0,\
\th_{q-1}(\si_{{q-1},k})^{q-1}_{k=1}>n_{q-1}n_{q-2}
\th_{q-2}(\alpha_{{q-2},0,k})^{q-2}_{k=1}, \tag Eq.(3q-3).1
$$
and $\gcd(n_{q-1}, \th_{q-1}(\si_{{q-1},k})^{q-1}_{k=1})=1$.

If $g_{q-1}$ of {\rm(Eq.(3q-3))} satisfies inequalities in
{\rm(Eq.(3q-3).1)}, then $g_{q-1}(y,z)$ is irreducible in
${}_{q-1}\CO_0$ with $g_{q-1}\in $ the type$[q-1]$ because
$g_{q-1}(y,z,f_1,\dots,f_{q-2})$ can be viewed as an element of
${}_{q}\CO_0$.

Apply the WDT to $f(y,z)$ with a divisor $g_{q-1}(y,z)$. Then,
$(g_{q-1},f)$ can be written in the form
$$\cases
g_{q-1} &=f^{n_{q-1}}_{q-2}
+\xi_{q-1}\Pi^q_{k=1}f^{\si_{{q-1},k}}_{k-2}\\
f &=g^{d_{q}}_{q-1} +\sum^{d_{q}-1}_{i=0}T^{(1)}_{q,i}g^i_{q-1},
\endcases \tag Eq.(3q-2))(Eq.(3q-2).1
$$
where $T^{(1)}_{q,i}\in \C\{y\}[z,f_1,\dots,f_{q-2}]$ for $0\le i\le
d_q-1$, satisfying the following property:

{\rm(i)} For any nonzero monomial $\Pi^q_{k=1}f^{\g_k}_{k-2}$ in
$T^{(1)}_{q,i}$,
$$
\g_1>0 \quad \text{and} \quad \g_k<n_{k-1}  \quad \text{for} \quad
2\le k\le q. \tag Eq.(3q-2).2
$$

{\rm(ii)} Define $\th_q :\N^{(q)}_0\to \N_0$ by
$\th_q(t_k)^q_{k=1}=t_q\th_{q-1}(\si_{{q-1},k})^{q-1}_{k=1}
+n_{q-1}\th_{q-1}(t_k)^{q-1}_{k=1}$, and then $\th_q$ is one-to-one.

If $T^{(1)}_{q,i}\ne 0$, let
$C^{(1)}_{q,i}\Pi^q_{k=1}f^{\beta^{(1)}_{q,i,k}}_{k-2}$ be a unique
nonzero monomial with a constant $C^{(1)}_{q,i}$ in $T^{(1)}_{q,i}$
such that
$\th_q(\beta^{(1)}_{q,i,k})^q_{k=1}=\min\{\th_q(\g_k)^q_{k=1}\}$ for
any nonzero monomial $\Pi^q_{k=1}f^{\g_k}_{k-2}$ in $T^{(1)}_{q,i}$.
\ms

\noindent{\bf {Remark.}} The above construction of $T^{(1)}_{q,i}$
is as follows: Note that $g_{q-1}(y,z)$ is a polynomial of degree
$\Pi^{q-1}_{k=1}n_k$ in $z$, and for $0\le i\le d_q-1$,
$T^{(1)}_{q,i}\in \C\{y\}[z]$ is a polynomial of degree
$<\Pi^{q-1}_{k=1}n_k$ in $z$, and $f_{q-2}(y,z)\in \C\{y\}[z]$ is a
polynomial of degree $\Pi^{q-2}_{k=1}n_k$ in $z$. Again, apply the
WDT to each $T^{(1)}_{q,i}$ with a divisor $f_{q-2}(y,z)$. Then
$T^{(1)}_{q,i}$ can be written just as above. \ms

Since $f$ is irreducible in ${}_2\CO_0$, $f$ must satisfy the
following necessary condition:
$$
\f{\th_q(\beta^{(1)}_{q,i,k})^q_{k=1}}{d_q-i}\ge
\f{\th_q(\beta^{(1)}_{q,0,k})^q_{k=1}}{d_q}>n_{q-1}\th_{q-1}(\si_{{q-1},k})^{q-1}_{k=1}
\text{\quad for $0\le i\le d_q-1$.} \tag Eq.(3q-2).3
$$

Recall that $\si_{{q-1},k}=\beta_{{q-1},d_{q-1}-n_{q-1},k}$ for
$1\le k\le {q-1}$ from the construction of $g_{q-1}$. \ms

If $f$ satisfies {\rm(Eq.(3q-2).3)}, then to find an irreducibility
algorithm for $f$, it suffices to consider the following two cases:

$\underline{\text{\rm Case(A)}}$
$\gcd(d_q,\th_q(\beta^{(1)}_{q,0,k})^q_{k=1})=1$, and
$\underline{\text{\rm Case(B)}}$
$\gcd(d_q,\th_q(\beta^{(1)}_{q,0,k})^q_{k=1})>1$. \ms

$\underline{\text{\bf Case(A) of the q-th step}}$ Let
$\gcd(d_q,\th_q(\beta^{(1)}_{q,0,k})^q_{k=1})=1$. Then, $f$ is
irreducible in ${}_2\CO_0$ if and only if {\rm(Eq.(3q-2).3)} holds.
In this case, if $f$ is irreducible in ${}_2\CO_0$, then $f\in $ the
type [q] and also $f\buildrel \text{{\rm multiseq}} \over \sim
H^{d_q}_{q-1}+\Pi^q_{k=1}H^{\si_{q,k}}_{k-2}$ where $H_{-1}=y$,
$H_0=z$, $H_1=z^{n_1}+y^{\beta_{1,0,1}}$,
$H_2=H^{n_2}_1+y^{\si_{2,1}}z^{\si_{2,2}}$,
$H_j=H^{n_j}_{j-1}+\Pi^j_{k=1}H^{\si_{j,k}}_{k-2}$ for $3\le j\le
q-2$,
$H_{q-1}=H^{d_{q-1}}_{q-2}+\Pi^{q-1}_{k=1}H^{\si_{{q-1},k}}_{k-2}$,
noting that for each fixed $j=1,2,\dots,q$,
$\si_{j,k}=\beta^{(1)}_{j,0,k}$ for $1\le k\le j$. \ms

\text{\rm Remark for Case(A) of the q-th step}. Let $(g_{q-1},f)$ be
in {\rm(Eq.(3q-2))}. By Theorem $1.13$,
$(f_1,\dots,f_{q-2},g_{q-1},f)$ may not be a generalized
representation of $f$, since $T^{({\nu_{q-1}}+1)}_{q-1,d_{q-1}-1}$
may not be zero. To find a generalized representation of $f$ in the
sense of Theorem $1.13$, until we get the same kind of result in the
conclusion of Sublemma 15.5, it suffices to apply either Sublemma
15.5 or the same kind of a fine sequence of pairs in
{\rm(Eq.(3q-1))} to $(g_{q-1},f)$ \ms

$\underline{\text{\bf Case(B) of the q-th step}}$ Let
$\gcd(d_q,\th_q(\beta^{(1)}_{q,0,k})^q_{k=1})>1$. To find a
necessary condition for $f$ to be irreducible in ${}_2\CO_0$, then
we may assume that $f$ satisfies {\rm(Eq.2)}, {\rm(Eq.3)},\dots,
{\rm(Eq.3q-2)}.

Then, it suffices to consider the following two subcases:

$\underline{\text{\rm Subcase(B1)}}$ $T^{(1)}_{q,d_q-1}\ne 0$, and
$\underline{\text{\rm Subcase(B2)}}$ $T^{(1)}_{q,d_q-1}=0$. \ms

$\underline{\text{\bf Subcase(B1) of the q-th step}}$ Let
$T^{(1)}_{q,d_q-1}\ne 0$. To solve the case, first of all, we must
eliminate $T^{(1)}_{q,d_q-1}$ whether or not $f$ is irreducible in
${}_2\CO_0$. To do it, by Theorem 1.8(Theorem 15.4), we can compute
a unique finite sequence of pairs, $\{(h_{q-1,p},f):1\le p\le
\nu+1\}$ in {\rm(Eq.(3q-1))} for a unique integer ${\nu_{q-1}} \le
\f{n_{q-1}+1}2$, each pair of which can be written in the form
$$\cases
h_{{q-1,}1} &=f^{n_{q-1}}_{q-2}
+\xi_{q-1}\Pi^{q-1}_{k=1}f^{\si_{q-1,k}}_{k-2}
=f^{n_{q-1}}_{q-2}+R^{(1)}_{q-1,0}\ \text{with $h_{q-1,1}=g_{q-1}$}\\
f &=h^{d_q}_{q-1,1} +\sum^{d_q-1}_{i=0}T^{(1)}_{q,i}h^i_{q-1,1},
\endcases \tag Eq.(3q-1))(Eq.(3q-1).1)(Eq.(3q-1).1.1
$$
and for $1\le p\le{\nu_{q-1}}-1$,
$$\cases
h_{q-1,p+1} &=h_{q-1,p}+\f
1{d_q}T^{(p)}_{q,d_q-1}=f^{n_{q-1}}_{q-2}+\sum^{n_{q-1}-2}_{i=0}R^{(p+1)}_{q-1,i}f^i_{q-2}\\
f
&=h^{d_q}_{q-1,p+1}+\sum^{d_q-1}_{i=0}T^{(p+1)}_{q,i}h^i_{q-1,p+1}\
\text{with}\ T^{(p+1)}_{q,d_q-1}\not=0,
\endcases \tag Eq.(3q-1).1.2
$$
and
$$\cases
h_{q-1,{\nu_{q-1}}+1} &=h_{q-1,{\nu_{q-1}}} +\f
1{d_q}T^{({\nu_{q-1}})}_{q,d_q-1}=f^{n_{q-1}}_{q-2}
+\sum^{n_{q-1}-2}_{i=0}R^{({\nu_{q-1}}+1)}_{q-1,i}f^i_{q-2}\\
f
&=h^{d_q}_{q-1,{\nu_{q-1}}+1}+\sum^{d_q-2}_{i=0}T^{({\nu_{q-1}}+1)}_{q,i}h^i_{q-1,{\nu_{q-1}}+1}\
\text{with}\ T^{({\nu_{q-1}}+1)}_{q,d_q-1}=0,
\endcases \tag Eq.(3q-1).1.3
$$

where $\si_{q-1,k}=\beta_{q-1,d_{q-1}-n_{q-1},k}$ for $1\le k\le
q-1$, $R^{(p)}_{q-1,i}\in \C\{y\}[z,f_1,\dots,f_{q-3}]$ for $1\le
p\le {\nu_{q-1}}+1$ and $0\le i\le n_{q-1}-2$; and
$T^{(p)}_{{\nu_{q-1}},i}\in \C\{y\}[z,f_1,\dots,f_{q-2}]$ for $1\le
p\le \nu+1$ and $0\le i\le d_q-1$, satisfying the following
properties:

\roster

\item"(i)" For any nonzero monomial $\Pi^{q-1}_{k=1}f^{\de_k}_{k-2}$ in
$R^{(p)}_{q-1,i}$,

\qquad \qquad \qquad $\de_1>0$ \quad and \quad $\de_k<n_{k-1}$ \quad
for $2\le k\le q-1$.

\item"(ii)" For any nonzero monomial $\Pi^q_{k=1}f^{\g_k}_{k-2}$ in
$T^{(p)}_{q,i}$,
$$
\text{$\g_1>0$ \quad and \quad $\g_q<n_{q-1}$ \quad for $2\le k\le
q$.} \tag Eq.(3q-1).2
$$

\item"(iii)" Define a function $\th_q:\N^q_0\to \N_0$ by
$\th_q(t_k)^q_{k=1}=t_q\th_{q-1}(\si_{q-1,k})^{q-1}_{k=1}
+n_{q-1}\th_{q-1}(t_k)^{q-1}_{k=1}$, and then $\th_{q}$ is
one-to-one. If $T^{(p)}_{q,i}\ne 0$, let
$C^{(p)}_{q,i}\Pi^q_{k=1}f^{\beta^{(p)}_{q,i,k}}_{k-2}$ be a unique
nonzero monomial with a constant $C^{(p)}_{q,i}$ in $T^{(p)}_{q,i}$
such that
$\th_q(\beta^{(p)}_{q,i,k})^q_{k=1}=\min\{\th_q(\g_k)^q_{k=1}\}$ for
any nonzero monomial $\Pi^q_{k=1}f^{\g_k}_{k-2}$ in $T^{(p)}_{q,i}$
where $f_{-1}=y$ and $f_0=z$.
\endroster

\noindent{\bf{Remark.}} It is possible to compute
$C^{(p)}_{q,i}\Pi^q_{k=1}f^{\beta^{(p)}_{q,i,k}}_{k-2}$ by an
elementary way because $T^{(p)}_{q,i}\in
\C\{y\}[z,f_1,f_2,\dots,f_{q-2}]$ in $z$,$f_1$,\dots,$f_{q-2}$,
recalling that for any nonzero monomial
$\Pi^{q}_{k=1}f^{\g_k}_{k-2}$ in $T^{(p)}_{q,i}$, $\g_1>0$ and
$\g_k<n_{k-1}$ for $2\le k\le q$, $\th_q$ is one-to-one. \ms

Since $f$ is irreducible in ${}_2\CO_0$, then $f$ must satisfy the
following necessary condition:
$$\align
&\f{\th_q(\beta^{(p)}_{q,i,k})^q_{k=1}}{d_q-i} \ge
\f{\th_q(\beta^{(p)}_{q,0,k})^q_{k=1}}{d_q}>n_{q-1}\th_{q-1}(\si_{q-1,k})^{q-1}_{k=1}
\text{\ for $1\le p\le {\nu_{q-1}}$, $0\le i\le d_q-1$,} \tag
{(Eq.(3q-1).3)(Eq.(3q-1).3.1)}\\
&\f{\th_q(\beta^{({\nu_{q-1}}+1)}_{q,i,k})^q_{k=1}}{d_q-i} \ge
\f{\th_q(\beta^{({\nu_{q-1}}+1)}_{q,0,k})^q_{k=1}}{d_q}>n_{q-1}\th_{q-1}(\si_{q-1,k})^{q-1}_{k=1}
\text{\quad for $0\le i\le d_q-2$.}  \tag{(Eq.(3q-1).3.2)} \\
\endalign$$

\noindent{\bf{Remark.}} Since $f$ satisfies {\rm(Eq.(3q-1).3)},
$h_{q-1,p}\buildrel \text{{\rm multiseq}} \over \sim
h_{q-1,1}=g_{q-1}$ for all $p\ge 1$. \ms

Since $f$ satisfies {\rm(Eq.(3q-1).3)}, to find an irreducibility
algorithm for Subcase(B1), then there exist the following two
possibilities, Subcase(B1-a) and Subcase(B1-b):

$\underline{\text{\rm Subcase(B1-a)}}$ Let
$\gcd(d_q,\th_q(\beta^{(1)}_{q,0,k})^q_{k=1})>1$ and
$\gcd(d_q,\th_q(\beta^{({\nu_{q-1}}+1)}_{q,0,k})^q_{k=1})=1$.

$\underline{\text{\rm Subcase(B1-b)}}$ Let
$\gcd(d_q,\th_q(\beta^{({\nu_{q-1}}+1)}_{q,0,k})^q_{k=1})>1$. \ms

$\underline{\text{\bf Subcase(B1-a) of the q-th step}}$

{\rm(a)} Let $\gcd(d_q,\th_q(\beta^{(1)}_{q,0,k})^q_{k=1})>1$ and
$\gcd(d_q,\th_q(\beta^{({\nu_{q-1}}+1)}_{q,0,k})^q_{k=1})=1$. Then,
$f$ is irreducible in ${}_2\CO_0$ if and only if {\rm
(Eq.(3q-1).3.2)} holds. Thus, $f\in $ the type $[q]$ in the sense of
Definition 2.5.

{\rm(b)} Assume that
$\gcd(d_q,\th_q(\beta^{(p)}_{q,0,k})^q_{k=1})=1$ for some $p\le
{\nu_{q-1}}+1$. Then, $f$ is irreducible in ${}_2\CO_0$ if and only
if an inequality in {\rm(Eq.(3q-1).3.1)} holds without mentioning
any inequality in {\rm(Eq.(3q-1).3.2)}. \ms

\text{\rm Remark for Subcase(B1-a) of the q-th step}.
$(f_1,\dots,f_{q-2},f_{q-1},f)$ with $f_{q-1}=h_{q-1,{\nu_{q-1}}+1}$
is a generalized representation of irreducible {W}-poly of the
recursive {q}-type for $f$ in the sense of Theorem $1.13$, since
$T^{({\nu_{q-1}}+1)}_{q-1,d_{q-1}-1}=0$. \ms

$\underline{\text{\bf Subcase(B1-b) of the q-th step}}$ Let
$\gcd(d_q,\th_q(\beta^{({\nu_{q-1}}+1)}_{q,0,k})^q_{k=1})>1$, noting
that $T^{({\nu_{q-1}}+1)}_{q,d_q-1}=0$. To find a necessary and
sufficient condition for $f$ to be irreducible in ${}_2\CO_0$, take
the (q+1)-th step. \ms

$\underline{\text{\bf Subcase(B2) of the q-th step}}$ It was assumed
by this subcase that $T^{(1)}_{q,d_q-1}=0$, noting that
$\gcd(d_q,\th_q(\beta^{(1)}_{q,0,k})^q_{k=1})>1$. To find a
necessary and sufficient condition for $f$ to be irreducible in
${}_2\CO_0$, take the (q+1)-th step. \quad $\square$

\noindent$\underline{\text{\bf Remark for the q-th step: To find
$H_q \in$ Family(1) such that $f\buildrel \text{{\rm multiseq}}
\over \sim H_q$ for}}$

\noindent$\underline{\text{\bf any $f\in $ the type [q] in the sense
of  Definition 2.5.}}$  \quad Let
$\gcd(d_q,\th_q(\beta^{({\nu_{q-1}}+1)}_{q,0,k})^q_{k=1})=1$. Then,
$f$ is irreducible in ${}_2\CO_0$ if and only if {\rm(Eq.(3q-1).2)}
holds. In this case, if $f$ is irreducible in ${}_2\CO_0$, then
$f\in $ the type $[q]$ and also $f\buildrel \text{{\rm multiseq}}
\over \sim H^{d_q}_{q-1}+\Pi^q_{k=1}H^{\si_{q,k}}_{k-2}$ where
$H_{-1}=y$, $H_0=z$, $H_1=z^{n_1}+y^{\beta_{1,0,1}}$, and for $3\le
j\le q$,
$H_{j-1}=H^{d_{j-1}}_{j-2}+\Pi^{j-1}_{k=1}H^{\si_{{j-1},k}}_{k-2}\in
{\text{\rm Family[1]}}$, noting that for each fixed
$j=1,2,\dots,q-1$ and for $1\le k\le j$,
$\si_{j,k}=\beta_{j,d_j-n_j,k}$, and for $1\le k\le q-1$
$\si_{q,k}=\beta^{({\nu_{q-1}}+1)}_{q,0,k}$  with some positive
integer ${\nu_{q-1}+1}$, because of the equations {\rm(Eq.(3j-3)},
{\rm(Eq.(3j-2)} and {\rm(Eq.(3j-1)} with $2\le j\le q$ in the
conclusion of Theorem 1.16. Note that
$d_{j+1}=\gcd(d_j,\theta_j(\sigma_{j,k})^j_{k=1})$ with
$d_j=n_jd_{j+1}$ for $j\ge 2$. \quad $\square$
\endproclaim \ms

\proclaim{Corollary 1.15.1 for Theorem 1.15}

$\underline{\text{\bf Assumptions}}$ Suppose that the same
properties and notations as in the assumptions of Theorem 1.15 hold.
\ms

$\underline{\text{\bf Conclusions}}$ \quad If $f$ is irreducible in
${}_2\CO_0$ with isolated singularity at $0\in \BC^2$, then $f\in $
the type $[\ell]$ for some $\ell\le r$ in the sense of Definition
2.5. By the induction method on the positive integer $r$, the aim is
to compute an elementary algorithm for finding irreducible W-polys
from all the W-polys in $\BC\{y\}[z]$, using $q$ iterations of the
following steps with $q\le \ell$: {\bf Observe that the statement on
the 3rd step may be omitted if necessary, to simplify the statements
for this theorem by the induction method}.

Following the same properties and notations as in the conclusions of
Theorem 1.15, to find $f\in $ the type $[\ell]$ for some $\ell\le r$
in the sense of Definition 2.5, there is nothing to prove by Theorem
1.15 that the following are true: \ms

$\underline{\text{\bf The 1st step: To find the irreducibility
algorithm for $f\in $ the type [1]}}$

$\underline{\text{\bf in the sense of  Definition 2.5.}}$ \ms

If $f$ satisfies {\rm(Eq.2)} of Theorem 1.15, it suffices to
consider the following case for the coefficient $a_{n-1}$ of
$z^{n-1}$.

$\underline{\text{\bf Case(I)}}$  $a_{n-1}= 0$.  There are two
subcases only.

$\underline{\text{\rm Subcase(I-1)}}$ $\gcd(n,\alpha_{0})=1$ and
$\underline{\text{\rm Subcase(I-2)}}$ $1<\gcd(n,\alpha_{0})<n$. \ms

$\underline{\text{\bf Subcase(I-1)}}$ Let $\gcd(n,\alpha_{0})=1$.
Then, $f$ is irreducible in ${}_2\CO_0$ with $f\in$ the type $[1]$
$\iff$ the inequality in {\rm(Eq.2)} of Theorem 1.15 holds. \ms

$\underline{\text{\bf Subcase(I-2)}}$ Let $1<\gcd(n,\alpha_{0})<n$.
If $f$ is irreducible in ${}_2\CO_0$, note that $f\in $ the type
$[{\ell}]$ for some $\ell\ge 2$ in the sense of Definition 2.5. To
find a necessary and sufficient condition for $f$ to be irreducible
in ${}_2\CO_0$, take the next step. \ms

$\underline{\text{\bf The 2nd step: To find the irreducibility
algorithm for $f\in $ the type [2]}}$

$\underline{\text{\bf in the sense of  Definition 2.5.}}$ \ms

With three equations, which is defined by {\rm(Eq.3)}, {\rm(Eq.4)}
and {\rm(Eq.5)} in the 2nd step of Theorem 1.15 later, the aim in
this step is how to compute the necessary and sufficient condition
for $f\in $ the type $[2]$ in the sense of Definition 2.5, {\bf
without need of computing any inequality in (Eq.5.3.1)} of Theorem
1.15.

For the irreducibility algorithm for the 2nd step, it suffices to
consider the following subcase only for the 1st step:

$\underline{\text{\rm Subcase(I-2)}}$ Let $1<\gcd(n,\a_0)<n$ and
$a_{n-1}=0$. \ms

If $f$ is irreducible in ${}_2\CO_0$, note that $f\in $ the type
$[\ell]$ for some $\ell\ge 2$ in the sense of Definition 2.5. To
find the irreducibility algorithm for $f\in$ the type $[2]$, using
the inequalities in {\rm(Eq.3)} with {\rm(Eq.3.1)} and {\rm(Eq.4)}
with {\rm(Eq.4.1)} and {\rm(Eq.4.2)} of Theorem 1.15 it suffices to
consider two cases:

$\underline{\text{\rm Case(I)}}$ $T^{(1)}_{2,d_2-1}= 0$, and
$\underline{\text{\rm Case(II)}}$ $T^{(1)}_{2,d_2-1}\ne 0$. \ms

$\underline{\text{\bf Case(I)}}$ Let $T^{(1)}_{2,d_2-1}= 0$. There
are two subcases only.

$\underline{\text{\rm Subcase(I-1)}}$
$\gcd(d_2,\th_2(\beta^{(1)}_{2,0,k})^2_{k=1})=1$ and
$\underline{\text{\rm Subcase(I-2)}}$
$\gcd(d_2,\th_2(\beta^{(1)}_{2,0,k})^2_{k=1})>1$. \ms

$\underline{\text{\bf Subcase(I-1)}}$ Let
$\gcd(d_2,\th_2(\beta^{(1)}_{2,0,k})^2_{k=1})=1$. Then, $f$ is
irreducible in ${}_2\CO_0$ with $f\in$ the type $[2]$ $\iff$ the
inequality in {\rm Eq(4.3)} of Theorem 1.15 holds. \ms

$\underline{\text{\bf Subcase(I-2)}}$ Let
$\gcd(d_2,\th_2(\beta^{(1)}_{2,0,k})^2_{k=1})>1$. If $f$ is
irreducible in ${}_2\CO_0$, note that $f\in $ the type $[\ell]$ for
some $\ell\ge 3$ in the sense of Definition 2.5. To find a necessary
and sufficient condition for $f$ to be irreducible in ${}_2\CO_0$ in
this case, take the next step. To compute the inequality in {\rm
(Eq.4.3)} may not be needed, if necessary. \ms

$\underline{\text{\bf Case(II)}}$ Let $T^{(1)}_{2,d_2-1}\ne 0$. To
solve the subcase, first of all, we must eliminate
$T^{(1)}_{2,d_2-1}$ whether or not $f$ is irreducible in
${}_2\CO_0$. To do it, by Sublemma 1.10 for Theorem 1.8(Sublemma
15.5 for Theorem 15.4), we can compute a unique finite sequence of
pairs $\{(h_{1,p},f):1\le p\le {\nu_1}+1\}$ in {\rm(Eq.5)} of
Theorem 1.15 for a unique integer ${\nu_1}\le \f{n_1+1}2$, which
satisfies the following:

{\rm(i)} $T^{(p)}_{2,d_2-1}\not=0$ for $p=1,2,\dots,{\nu_1}$ and
$T^{(\nu_1+1)}_{2,d_2-1}=0$.

{\rm(ii)} If $f$ satisfies {\rm(Eq.5.1)} and {\rm(Eq.5.2)} of
Theorem 1.15, without computing any inequality in {\rm(Eq.5.3.1)} of
Theorem 1.15, to find an irreducibility algorithm for $f$, it
suffices to consider the following two subcases:

$\underline{\text{\rm Subcase(II-1)}}$
$\gcd(d_2,\th_2(\beta^{(\nu_1+1)}_{2,0,k})^2_{k=1})=1$ and
$\underline{\text{\rm Subcase(II-2)}}$
$\gcd(d_2,\th_2(\beta^{(\nu_1+1)}_{2,0,k})^2_{k=1})>1$. \ms

$\underline{\text{\bf Subcase(II-1)}}$ Let
$\gcd(d_2,\th_2(\beta^{(\nu_1+1)}_{2,0,k})^2_{k=1})=1$. Then, $f$ is
irreducible in ${}_2\CO_0$ with $f\in$ the type $[2]$ $\iff$
{\rm(Eq.5.3.2)} of Theorem 1.15 holds. \ms

$\underline{\text{\bf Subcase(II-2)}}$ Let
$\gcd(d_2,\th_2(\beta^{(\nu_1+1)}_{2,0,k})^2_{k=1})>1$. If $f$ is
irreducible in ${}_2\CO_0$, note that $f\in $ the type $[\ell]$ for
some $\ell\ge 3$ in the sense of Definition 2.5. To find a necessary
and sufficient condition for $f$ to be irreducible in ${}_2\CO_0$ in
this case, take the next step. To compute the inequality in {\rm
Eq(5.3.1)} of Theorem 1.15 may not be needed, if necessary.\ms

$\underline{\text{\bf The 3rd step: To find the irreducibility
algorithm for $f\in $ the type [3]}}$

$\underline{\text{\bf in the sense of  Definition 2.5.}}$ \ms

With three equations, which is defined by {\rm(Eq.6)}, {\rm(Eq.7)}
and {\rm(Eq.8)} of Theorem 1.15 in the 3rd step later, the aim in
this step is how to compute the necessary and sufficient condition
for $f\in $ the type $[3]$ in the sense of Definition 2.5, {\bf
without need of computing any inequalty in (Eq.8.3.1)} of Theorem
1.15.

For the irreducibility algorithm for the 3rd step, it suffices to
consider the following subcase only for the 2nd step:

$\underline{\text{\rm Subcase(II-2)}}$ Let
$\gcd(d_2,\th_2(\beta^{(\nu_1+1)}_{2,0,k})^2_{k=1})>1$ and
$T^{({\nu_1}+1)}_{2,d_2-1}=0$ for some integer $\nu_1\ge 0$.\ms

If $f$ is irreducible in ${}_2\CO_0$ for \text{\rm Subcase(II-2)},
note that $f\in $ the type $[\ell]$ for some $\ell\ge 3$ in the
sense of Definition 2.5. To find the irreducibility algorithm for
$f\in$ the type $[3]$, using the inequalities in {\rm(Eq.6)} with
{\rm(Eq.6.1)} and {\rm(Eq.7)} with {\rm(Eq.7.1)} and {\rm(Eq.7.2)}
of Theorem 1.15 it suffices to consider two new cases:

$\underline{\text{\rm Case(I)}}$ $T^{(1)}_{3,d_3-1}= 0$, and
$\underline{\text{\rm Case(II)}}$ $T^{(1)}_{3,d_3-1}\ne 0$. \ms

$\underline{\text{\bf Case(I)}}$ Let $T^{(1)}_{3,d_3-1}= 0$. There
are two subcases only.

$\underline{\text{\rm Subcase(I-1)}}$
$\gcd(d_3,\th_3(\beta^{(1)}_{3,0,k})^3_{k=1})=1$ and
$\underline{\text{\rm Subcase(I-2)}}$
$\gcd(d_3,\th_3\beta^{(1)}_{3,0,k})^3_{k=1})>1$. \ms

$\underline{\text{\bf Subcase(I-1)}}$ Let
$\gcd(d_3,\th_3(\beta^{(1)}_{3,0,k})^3_{k=1})=1$. Then, $f$ is
irreducible in ${}_2\CO_0$ with $f\in$ the type $[3]$ $\iff$ {\rm
(Eq.7.3)} of Theorem 1.15 holods. \ms

$\underline{\text{\bf Subcase(I-2)}}$ Let
$\gcd(d_3,\th_3(\beta^{(1)}_{3,0,k})^3_{k=1})>1$. If $f$ is
irreducible in ${}_2\CO_0$, note that $f\in $ the type $[\ell]$ for
some $\ell\ge {4}$ in the sense of Definition 2.5. To find a
necessary and sufficient condition for $f$ to be irreducible in
${}_2\CO_0$ in this case, take the 4-th step. To compute the
inequality in {\rm (Eq.7.3)} of Theorem 1.15 may not be needed, if
necessary.\ms

$\underline{\text{\bf Case(II)}}$ Let $T^{(1)}_{3,d_3-1}\ne 0$. To
solve the subcase, first of all, we must eliminate
$T^{(1)}_{3,d_3-1}$ whether or not $f$ is irreducible in
${}_2\CO_0$. To do it, by Sublemma 1.10 for Theorem 1.8(Sublemma
15.5 for Theorem 15.4), we can compute a unique finite sequence of
pairs $\{(h_{2,p},f):1\le p\le {\nu_2}+1\}$ in {\rm(Eq.7)} of
Theorem 1.15 for a unique integer ${\nu_2}\le \f{n_3+1}2$, which
satisfies the following:

{\rm(i)} $T^{(p)}_{3,d_3-1}\not=0$ for $p=1,2,\dots,{\nu_2}$ and
$T^{(\nu_2+1)}_{2,d_3-1}=0$.

{\rm(ii)} If $f$ satisfies {\rm(Eq.8.1)} and {\rm(Eq.8.2)} of
Theorem 1.15, without computing any inequaity in {\rm(Eq.8.3.1)} of
Theorem 1.15, to find an irreducibility algorithm for $f$, it
suffices to consider the following two subcases:

$\underline{\text{\rm Subcase(II-1)}}$
$\gcd(d_3,\th_3(\beta^{(\nu_2+1)}_{3,0,k})^3_{k=1})=1$ and
$\underline{\text{\rm Subcase(II-2)}}$
$\gcd(d_3,\th_3(\beta^{(\nu_2+1)}_{3,0,k})^3_{k=1})>1$. \ms

$\underline{\text{\bf Subcase(II-1)}}$ Let
$\gcd(d_3,\th_3(\beta^{(\nu_2+1)}_{3,0,k})^3_{k=1})=1$. Then, $f$ is
irreducible in ${}_2\CO_0$ with $f\in$ the type $[3]$ $\iff$
{\rm(Eq.8.3.2)} of Theorem 1.15 holds. \ms

$\underline{\text{\bf Subcase(II-2)}}$ Let
$\gcd(d_3,\th_3(\beta^{(\nu_2+1)}_{3,0,k})^3_{k=1})>1$. If $f$ is
irreducible in ${}_2\CO_0$ for \text{\rm Subcase(II-2)}, note that
$f\in $ the type $[\ell]$ for some $\ell\ge 4$ in the sense of
Definition 2.5. To find a necessary and sufficient condition for $f$
to be irreducible in ${}_2\CO_0$ in this case, take the next step.
To compute the inequality in {\rm Eq(8.3.1)} of Theorem 1.15 may not
be needed, if necessary. \ms

The general case will be proved by induction. Let $f\in \BC\{y\}[z]$
be an arbitrary $W$-poly of degree $n\ge 2$ in $z$. Suppose we have
shown for any integer $q\ge 3$ that for each \text{\rm(j)-th} step,
$1\le j\le q-1$, the irreducibility algorithm for $f\in $ the type
$[j]$ has been found, by using the same kind of properties and
notations as we have seen in the proof of finding the irreducibility
algorithm for the \text{\rm(j)-th} step, $1\le j\le q-1$. Now, we
may assume by induction on $q$ that the irreducibility algorithm for
the \text{\rm(q-1)-th} step can be solvable for $q\ge 3$, and so if
$f$ is irreducible in ${}_2\CO_0$, we may assume without loss of
generality that if $f\in $ the type $[k]$ for an integer $k\ge q$ in
the sense of Definition 2.5 then
$\gcd(d_{q-1},\th_{q-1}(\beta^{({\nu_{q-2}}+1)}_{q-1,0,k})^{q-1}_{k=1})>1$
and $T^{({\nu_{q-2}}+1)}_{q-1,d_{{q-1}-1}}=0$ for some positive
integer ${\nu_{q-2}}+1$. Then, it remains to show that the
generalized irreducibility algorithm for $f\in $ the type $[q]$ can
be written as follows. \ms

$\underline{\text{\bf The q-th step: To find the irreducibility
algorithm for $f\in $ the type [q]}}$

$\underline{\text{\bf in the sense of  Definition 2.5.}}$ \ms

With three equations, which is defined by {\rm(Eq.(3q-3))},
{\rm(Eq.(3q-2))} and {\rm(Eq.(3q-1))} in the q-th step of Theorem
1.15 later, the aim in this step is how to compute the necessary and
sufficient condition for $f\in $ the type $[q]$ in the sense of
Definition 2.5, {\bf without need of computing any inequalty in
(Eq.(3q-1).3.1)} of Theorem 1.15.

Suppose by induction on the positive integer (q-1) that $f$ is
irreducible in ${}_2\CO_0$ and let $f$ satisfy a finite number
(3q-4) of conditions, which have been represented by {\rm(Eq.1),
(Eq.2), (Eq.3), $\dots$, (Eq.{(3q-6)}), (Eq.{(3q-5)}),
(Eq.{(3q-4)})} of Theorem 1.15. \ms

For the proof of this step, recall the defining equation
$(h_{{q-2},{\nu_{q-2}}+1},f)$ of {\rm(Eq.{(3q-4)})} of Theorem 1.15
as we have already seen in Subcase(B1) of the {(q-1)}-th step of
Theorem 1.15:
$$\cases
h_{{q-2},{\nu_{q-2}}+1} &=f^{n_{q-2}}_{q-3}
+\sum^{n_{q-2}-2}_{i=0}R^{({\nu_{q-2}}+1)}_{q-2,i}f^i_{q-3}\\
f &=h^{d_{q-1}}_{{q-2},{\nu_{q-2}}+1}
+\sum^{d_{q-1}-2}_{i=0}T^{({\nu_{q-2}}+1)}_{q-1,i}h^i_{{q-2},{\nu_{q-2}}+1}
\endcases \tag Eq.(3q-4)
$$
of Theorem 1.15, noting that
$T^{({\nu_{q-2}}+1)}_{q-1,d_{q-1}-1}=0$. \ms

By either \text{\rm Subcase(I-2)} or \text{\rm Subcase(II-2)} of the
{(q-1)}-th step, it may be assumed that
$\gcd(d_{q-1},\th_{q-1}(\beta^{({\nu_{q-2}}+1)}_{{q-1},0,k})^{q-1}_{k=1})>1$
and $T^{({\nu_{q-2}}+1)}_{{(q-1)},d_{(q-1)}-1}=0$ for some positive
integer ${\nu_{q-2}}+1$.

\noindent{\bf Remark.} If $T^{(1)}_{{q-1},d_{q-1}-1}=0$ and
$\gcd(d_{q-1},\th_{q-1}(\beta^{(1)}_{{q-1},0,k})^{q-1}_{k=1})>1$ for
Subcase(I-2), then $(h_{{q-2},{1}},f)=(g_{q-2},f)$ can be viewed as
$(h_{q-2,\nu_{q-2}+1},f)$ with ${\nu_{q-2}}=0$. $\square$ \ms

For the irreducibility algorithm for the q-th step, it suffices to
consider the following subcase only for the (q-1)-th step:

$\underline{\text{\rm Subcase(II-2)}}$ Let
$\gcd(d_{q-1},\th_{q-1}(\beta^{(\nu_{q-2}+1)}_{{q-1},0,k})^{q-1}_{k=1})>1$
for some integer $\nu_{q-2}\ge 0$.\ms

If $f$ is irreducible in ${}_2\CO_0$ for \text{\rm Subcase(II-2)},
note that $f\in $ the type $[\ell]$ for some $\ell\ge 3$ in the
sense of Definition 2.5. To find the irreducibility algorithm for
$f\in$ the type $[q]$ in the sense of Definition 2.5, using the
inequalities in {\rm(Eq.3q-3)} with {\rm(Eq.(3q-3).1)} and
{\rm(Eq.3q-2)} with {\rm(Eq.(3q-2).1)} and {\rm(Eq.(3q-2).2)} of
Theorem 1.15, it suffices to consider two new cases:

$\underline{\text{\rm Case(I)}}$ $T^{(1)}_{q,d_q-1}= 0$, and
$\underline{\text{\rm Case(II)}}$ $T^{(1)}_{q,d_q-1}\ne 0$. \ms

$\underline{\text{\bf Case(I)}}$ Let $T^{(1)}_{q,d_q-1}= 0$. There
are two subcases only.

$\underline{\text{\rm Subcase(I-1)}}$
$\gcd(d_q,\th_q(\beta^{(1)}_{q,0,k})^q_{k=1})=1$ and
$\underline{\text{\rm Subcase(I-2)}}$
$\gcd(d_q,\th_q(\beta^{(1)}_{q,0,k})^q_{k=1})>1$. \ms

$\underline{\text{\bf Subcase(I-1)}}$ Let
$\gcd(d_q,\th_q(\beta^{(1)}_{q,0,k})^q_{k=1})=1$. Then, $f$ is
irreducible in ${}_2\CO_0$ with $f\in$ the type $[q]$ in the sense
of Definition 2.5 $\iff$ the inequality in {\rm(Eq.(3q-2).3)} of
Theorem 1.15 holds. \ms

$\underline{\text{\bf Subcase(I-2)}}$ Let
$\gcd(d_q,\th_q(\beta^{(1)}_{q,0,k})^q_{k=1})>1$. If $f$ is
irreducible in ${}_2\CO_0$, note that $f\in $ the type $[\ell]$ for
some $\ell\ge {q+1}$ in the sense of Definition 2.5. To find a
necessary and sufficient condition for $f$ to be irreducible in
${}_2\CO_0$ in this case, take the (q+1)-th step. To compute the
inequality in {\rm Eq.(3q-2).3)} of Theorem 1.15 may not be needed,
if necessary. \ms

$\underline{\text{\bf Case(II)}}$ Let $T^{(1)}_{q,{d_q}-1}\ne 0$. To
solve the subcase, first of all, we must eliminate
$T^{(1)}_{q,{d_q}-1}$ whether or not $f$ is irreducible in
${}_2\CO_0$. To do it, by Sublemma 1.10 for Theorem 1.8(Sublemma
15.5 for Theorem 15.4, we can compute a unique finite sequence of
pairs $\{(h_{q-1,p},f):1\le p\le {\nu_{q-1}+1}\}$ in {\rm(Eq.3q-2)}
of Theorem 1.15 for a unique integer $\nu_{q-1}\le \f{{n_q}+1}{2}$,
which satisfies the following:

{\rm(i)} $T^{(p)}_{q,d_q-1}\not=0$ for $p=1,2,\dots,{\nu_{q-1}}$ and
$T^{(\nu_{q-1}+1)}_{2,d_q-1}=0$.

{\rm(ii)} If $f$ satisfies {\rm(Eq.(3q-1).1)} and {\rm(Eq.(3q-1).2)}
of Theorem 1.15, without computing any inequaity in
{\rm(Eq.(3q-1).3.1)} of Theorem 1.15, to find an irreducibility
algorithm for $f$, it suffices to consider the following two
subcases:

\noindent$\underline{\text{\rm Subcase(II-1)}}$
$\gcd(d_q,\th_q(\beta^{(\nu_{q-1}+1)}_{{q-1},0,k})^{q-1}_{k=1})=1$;
$\underline{\text{\rm Subcase(II-2)}}$
$\gcd(d_q,\th_q(\beta^{(\nu_{q-1}+1)}_{{q-1},0,k})^{q-1}_{k=1})>1$.
\ms

$\underline{\text{\bf Subcase(II-1)}}$ Let
$\gcd(d_q,\th_q(\beta^{(\nu_{q-1}+1)}_{{q-1},0,k})^{q-1}_{k=1})=1$.
Then, $f$ is irreducible in ${}_2\CO_0$ with $f\in$ the type $[q]$
in the sense of Definition 2.5 $\iff$ \text{\rm(Eq.(3q-1).3.2)} of
Theorem 1.15 holds. \ms

$\underline{\text{\bf Subcase(II-2)}}$ Let
$\gcd(d_q,\th_q(\beta^{(\nu_{q-1}+1)}_{{q-1},0,k})^{q-1}_{k=1})>1$.
If $f$ is irreducible in ${}_2\CO_0$, note that $f\in $ the type
$[\ell]$ for some $\ell\ge q+1$ in the sense of Definition 2.5. To
find a necessary and sufficient condition for $f$ to be irreducible
in ${}_2\CO_0$ in this case, take the (q+1)-th step.
\text{\rm(Eq.(3q-1).3.1)} of Theorem 1.15 may not be needed if
necessay. $\square$
\endproclaim \ms

\definition{Remark 1.15.2}
Note that the proof of Theorem 1.13 can be done by Theorem 16.5 and
Theorem 16.6. Using the same method as we have used in the process
of the proof of Theorem 16.6 together with Proposition 16.7 and
Proposition 16.8, there is nothing to prove for Theorem 1.15 with
Corollary 1.15.1(The 2nd Algorithm). $\square$ \enddefinition \bs

{\bf $\S$ 1.9. The 3rd Algorithm for computing the corresponding
standard Puiseux expansion from any irreducible W-poly of two
complex variables with respect to the multiplicity sequences} \ms

In preparation for finding the representation of the desired
algorithm, we prefer to use the properties and notations in
Definition 1.2, Theorem 1.4 and Theorem 1.6. By Definition 1.2,
$\underline{\text{\rm{Family(1)}}}$ is  the first family, denoted by
\text{\rm{Family(1)}}= \text{\{$f\in$ Family(0):$f$ is arbitrary
standard}

\noindent \text{Puiseux polynomial of the recursive r-type and $r$
are any positive integers}\} and $\underline{\text{\rm{Family(2)}}}$
is the 2nd family, denoted by \text{\rm{Family(2)}}= \text{\{
$C_r(t)$:$C_r(t)$ is the standard Puiseux expansion of the }

\noindent \text{$r$-type for any $r \in $ N\}}.

Then, we use the following notations in {\rm(i)} and {\rm(ii)}.

\noindent{\rm(i)} For any $f$ and $g$ in Family(0), an equivalence
relation for any two multiplicity sequences \text{\rm
Multiseq(V(f))} and \text{\rm Multiseq(V(g))} in Family(3) is
defined as follows:
$$\align
& \text{\text{\rm Multiseq(V(f))} and \text{\rm Multiseq(V(g))}
are equivalent} \tag $*$\\
 \text{$\iff$} \qquad & \text{$f \buildrel \text{{\rm
multiseq}} \over \sim g$} \quad \text{at
the origin in $\BC^2$} \\
\endalign$$

\noindent{\rm(ii)} For any $f$ in Family(0) and any $C_r(t)$ in
Family(2), an equivalence relation for any two multiplicity
sequences \text{\rm Multiseq(V(f))} and \text{\rm
Multiseq($C_r(t)$)} in Family(3) is defined as follows:
$$\align
& \text{\text{\rm Multiseq(V(f))} and \text{\rm Multiseq($C_r(t)$)}
are equivalent} \tag $**$\\
 \text{$\iff$} \qquad & \text{$f \buildrel \text{{\rm
multiseq}} \over \sim C_r$} \quad \text{at
the origin in $\BC^2$ \qquad $\square$ }\\
\endalign$$  \ms

\proclaim{Theorem 1.16(The 3rd Algorithm: Explicit algorithm for
finding the corresponding standard Puiseux expansion from any given
irreducible W-poly of two complex variables)}

$\underline{\text{\bf Assumptions}}$ Suppose that the same
properties and notations as in the assumptions of Theorem 1.15 hold.
In addition, let $f$ of {\rm(Eq.1)} of Theorem 1.15 be irreducible
in ${}_2\CO_0$ with isolated singularity at $0\in \BC^2$. Write
$n=\Pi^{r}_{k=1}n_k$ with positive integers $n_k\ge 2$ for all $k$
where the $n_k$ may not be the factorization of prime numbers. Then
$f\in $ the type $[\ell]$ for some $1\le {\ell}\le r$ in the sense
of Definition 2.5.\bs

$\underline{\text{\bf Conclusions}}$ Following the same properties
and notations as in the assumptions and the conclusions of Theorem
1.15, to find explicit algorithm for computing the corresponding
standard Puiseux expansion from any irreducible W-poly $f \in
\BC\{y\}[z]$, it suffices to solve The 1st Problem by The 1st
Conclusion and The 2nd Problem by the proof of The 2nd Conclusion
respectively, because Conclusions of this theorem can be divided by
The 1st Conclusion and The 2nd Conclusion, as follows:

$\underline{\text{\rm The 1st Problem}}$ \quad The problem is to
find explicit algorithm for computing the standard Puiseux
polynomial of the recursive r-type $H_q\in {\text{\rm Family[q]}}$
for any $f\in $ the type $\text{\rm[q]}$ in the sense of Definition
2.5 such that $f \buildrel {\text{\rm multiseq}} \over \sim H_q$  at
$0\in\BC^2$ for each $q=1,2,\dots,r$, as an application of Theorem
1.15. \bs

$\underline{\text{\rm The 2nd Problem}}$ \quad The problem is to
find Explicit algorithm for computing standard Puisuex expansion
$C_q\in {\text{\rm Family[2]}}$ for any standard Puiseux polynomial
of the recursive r-type $H_q\in {\text{\rm Family[1]}}$  such that
$H_q \buildrel {\text{\rm multiseq}} \over \sim C_q$  at
$0\in\BC^2$, as an application of Theorem 1.4.  \bs

$\underline{\text{\bf [I] The 1st Conclusion(A solution of The 1st
problem)}}$

Let $f\in $ the type \text{\rm[$\ell$]} for some $\ell\le r$ in the
sense of Definition 2.5, and follow the same properties and
notations as in the assumptions and the conclusions of Theorem 1.15
for a solution of The 1st problem. By the induction method on the
set of the positive integers, the aim is to solve the 1st problem,
using $q$ iterations of the following steps with $q\le \ell$. \ms

$\underline{\text{\bf The 1st step: To find $H_1\in {\text{\rm
Family[1]}}$ for any $f\in $ the type [1] in the sense of  }}$

$\underline{\text{\bf Definition 2.5 such that $f \buildrel
{\text{\rm multiseq}} \over \sim H_1$  at $0\in\BC^2$}}$ \ms

For the aim in this step, by Remark for The 1st step in the
conclusion of Theorem 1.15 we can find $H_1\in {\text{\rm
Family[1]}}$ for any $f\in $ the type $[1]$ in the sense of
Definition 2.5 such that $f\buildrel \text{{\rm multiseq}} \over
\sim H_1$ where $H_1=z^{n}+y^{\alpha_{0}}\in Family[1]$ with
$\gcd(n,\a_0)=1$. $\square$ \bs

$\underline{\text{\bf The 2nd step: To find $H_2\in {\text{\rm
Family[1]}}$ for any $f\in $ the type [2] in the sense of }}$

$\underline{\text{\bf Definition 2.5 such that $f \buildrel
{\text{\rm multiseq}} \over \sim H_2$  at $0\in\BC^2$}}$ \ms

For the aim in this step, by Remark for The 2nd step in the
conclusion of Theorem 1.15, we can find $H_2\in {\text{\rm
Family[1]}}$ for any $f\in $ the type $[2]$ in the sense of
Definition 2.5 such that $f\buildrel \text{{\rm multiseq}} \over
\sim H_2=H^{d_2}_1+\Pi^2_{k=1}H^{\si_{2,k}}_{k-2}$ where $H_{-1}=y$,
$H_0=z$, $H_1=z^{n_1}+y^{\alpha_{1,0,1}}\in {\text{\rm Family[1]}}$
and $\si_{2,k}=\beta^{({\nu_1}+1)}_{2,0,k}$ for $1\le k\le 2$ and
for some integer ${\nu_1}+1>0$, because of the equations,
{\rm(Eq.3)}, {\rm(Eq.4)} and {\rm(Eq.5)} in the conclusion of
Theorem 1.15. $\square$  \ms

$\underline{\text{\bf The q-th step: To find $H_q\in {\text{\rm
Family[1]}}$ for any $f\in $ the type [q] in the sense of }}$

$\underline{\text{\bf Definition 2.5 such that $f \buildrel
{\text{\rm multiseq}} \over \sim H_q$  at $0\in\BC^2$}}$ \ms

For the aim in this step, by Remark for The q-th step in the
conclusion of Theorem 1.15, we can find $H_q\in {\text{\rm
Family[1]}}$ for any $f\in $ the type $[q]$ in the sense of
Definition 2.5 such that $f\buildrel \text{{\rm multiseq}} \over
\sim H_q=H^{d_q}_{q-1}+\Pi^q_{k=1}H^{\si_{q,k}}_{k-2}$ where
$H_{-1}=y$, $H_0=z$, $H_1=z^{n_1}+y^{\alpha_{1,0,1}}$, and for $3\le
j\le q$,
$H_{j-1}=H^{d_{j-1}}_{j-2}+\Pi^{j-1}_{k=1}H^{\si_{{j-1},k}}_{k-2}\in
{\text{\rm Family[1]}}$, noting that for each fixed
$j=1,2,\dots,q-1$ and for $1\le k\le j$,
$\si_{j,k}=\beta_{j,d_j-n_j,k}$, and for $1\le k\le q-1$
$\si_{q,k}=\beta^{({\nu_{q-1}}+1)}_{q,0,k}$  with some positive
integer ${\nu_{q-1}+1}$, because of the equations, {\rm(Eq.(3j-3)},
{\rm(Eq.(3j-2)} and {\rm(Eq.(3j-1)} with $2\le j\le q$ in the
conclusion of Theorem 1.15. Note that
$d_{j+1}=\gcd(d_j,\theta_j(\sigma_{j,k})^j_{k=1})$ with
$d_j=n_jd_{j+1}$ for $j\ge 2$. $\square$ \ms

$\underline{\text{\bf [II] The 2nd Conclusion(A solution of The 2nd
problem)}}$ \bs

In preparation for finding a solution of The 2nd problem, follow the
same properties and notations as in {\rm The 1st Conclusion} in
order to avoid the complexity of the terminology and notations of
the statements between Theorem 1.4 and Theorem 1.16. For any $H_q\in
{\text{\rm Family[1]}}$ of The 1st Conclusion,  it suffices to show
by explicit algorithm in {\rm(1.4.1)} of Theorem 1.4 that we can
rewrite explicit algorithm for finding the standard Puiseux
expansion, denoted by the curve $C(H_q:t)$, which satisfies
\text{\rm $f \buildrel {\text{\rm multiseq}} \over \sim H_q$} at
$0\in\BC^2$ explicitly and rigorously, as in Sublemma 1.17. \ms

\noindent{\bf{Sublemma 1.17(Theorem 1.4:Algorithm for finding a
one-to-one function from Family(1) into Family(2))}}

$\underline{\text{\bf {Assumptions}}}$ Following the same properties
and notations as in the assumption and The 1st Conclusion of Theorem
1.16, recall the following:

{\rm(i)} We may assume by Lemma 1.8 that $f\in \BC\{y\}[z]$ is an
arbitrary $W$-poly of degree $n\ge 2$ in $z$, satisfying the
following form:
$$ f=z^n+\sum^{n-2}_{i=0} a_iy^{\a_i}z^i, \tag Eq.1 $$

{\rm(ii)} In the q-th step in The 1st Conclusion of this
theorem{\rm(A solution of The 1st problem)}, we can compute $H_q\in
{\text{\rm Family[1]}}$ for any $f\in $ the type $[q]$ in the sense
of Definition 2.5 such that \text{\rm $f \buildrel {\text{\rm
multiseq}} \over \sim H_q$} at $0\in\BC^2$, directly. \ms

$\underline{\text{\bf Conclusions}}$

\noindent$\underline{\text{\rm {\bf [I]} {\rm By Algorithm 1.4.1 for
Theorem 1.4,}  we can compute the standard Puiseux expansion}}$

\noindent$\underline{\text{\rm for the curve $C(H_q:t)$ such that
\text{$\text{\rm{Multiseq}}(V(H_q))\equiv \text{\rm
Multiseq}(C(H_q:t))$ \text{\rm as sequence:}}}}$ \ms

\noindent$\underline{\text{\rm(Algorithm 1.4.1 for Theorem 1.4)}}$
$$\align
\text{\rm(Eq.2)} \quad \quad \text{$C(H_q:t):=$} &
\left\{\eqalign{y&=t^n \cr z&=t^{\gamma_1}+t^{\gamma_2}+\cdots
+t^{\gamma_q}, \cr} \right.
 \\
 \text{such that}  \quad
 d_2  &=\gcd(n,\gamma_1) \quad\text{with}\quad n=n_1d_2 \quad\text{and}\quad
 \gamma_1=\alpha_{1,0,1}d_2, \quad \text{and}\\
\text{for $1\le j\le r$,}\quad
d_{j+1}&=\gcd(d_j,\gamma_j-\gamma_{j-1}) \quad\text{with}\quad
d_j=n_jd_{j+1} \quad\text{and}\quad
\gamma_j-\gamma_{j-1}=\widehat{\theta}_jd_j
  \endalign$$
where $d_r=1$ and $\widehat{\theta}_j
=\theta_j(\beta_{j,0,k})^j_{k=1}-n_jn_{j-1}\theta_{j-1}(\beta_{j-1,0,k})^{j-1}_{k=1}$
is a positive integer for $2\le j\le r$ and $\theta_1(t)=t$. Note
that $d_j=n_{j}d_{j+1}$ for $1\le j\le r-1$ . \bs

\noindent$\underline{\text{\bf Remark for the q-th step: To find
$H_q \in$ Family(1) such that $f\buildrel \text{{\rm multiseq}}
\over \sim H_q$ for}}$

\noindent$\underline{\text{\bf any $f\in $ the type [q] in the sense
of Definition 2.5.}}$ \quad Let
$\gcd(d_q,\th_q(\beta^{({\nu_{q-1}}+1)}_{q,0,k})^q_{k=1})=1$. Then,
$f$ is irreducible in ${}_2\CO_0$ if and only if {\rm(Eq.(3q-1).2)}
holds. In this case, if $f$ is irreducible in ${}_2\CO_0$, then
$f\in $ the type [q] in the sense of Definition 2.5 and also
$f\buildrel \text{{\rm multiseq}} \over \sim
H^{d_q}_{q-1}+\Pi^q_{k=1}H^{\si_{q,k}}_{k-2}$ where $H_{-1}=y$,
$H_0=z$, $H_1=z^{n_1}+y^{\beta_{1,0,1}}$,
$H_2=H^{n_2}_1+y^{\si_{2,1}}z^{\si_{2,2}}$,
$H_j=H^{n_j}_{j-1}+\Pi^j_{k=1}H^{\si_{j,k}}_{k-2}$ for $3\le j\le
q-2$,
$H_{q-1}=H^{d_{q-1}}_{q-2}+\Pi^{q-1}_{k=1}H^{\si_{{q-1},k}}_{k-2}$,
noting that for each fixed $j=1,2,\dots,q-2$,
$\si_{j,k}=\beta^{(1)}_{j,0,k}$ for $1\le k\le j$ and
$\si_{q-1,k}=\beta^{({\nu_{q-1}}+1)}_{q-1,0,k}$ for for $1\le k\le
q-1$. $\square$

Assume that $f\in $ the type {\rm[q]} with $q\ge 3$ in the sense of
Definition 2.5.

{\rm(i)} Define a function $\th_2 :\N^2_0\to \N_0$ by
$\th_2(t_k)^2_{k=1}=t_2\th_1(\si_{1,1})+n_1t_1$.
$\sigma_{1,1}=\alpha_{1,0,1}$. \ms

{\rm(ii)} By induction on positive integers j=2,3,\dots,q-1, define
$\th_{j} :\N^{(j)}_0\to \N_0$ by
$\th_{j}(t_k)^{j}_{k=1}=t_{j}\th_{j-1}(\si_{j-1,k})^{j-1}_{k=1}+n_{j-1}\th_{j-1}(t_k)^{j-1}_{k=1}$,
and then $\th_{j}$ is one-to-one. By {\rm(Eq.3j-3)},
$\si_{j-1,k}=\beta_{j-1,{d_{j-1}-n_{j-1}},k}$ for $k=1,2$. \ms

{\rm(q)} Define $\th_q :\N^{(q)}_0\to \N_0$ by
$\th_{q}(t_k)^{q}_{k=1}=t_{q}\th_{q-1}(\si_{q-1,k})^{q-1}_{k=1}+n_{q-1}\th_{q-1}(t_k)^{q-1}_{k=1}$,
and then $\th_{q-1}$ is one-to-one. By {\rm(Eq.3q-3)},
$\si_{q-1,k}=\beta^{(\nu+1)}_{q-1,0,k}$ for $k=1,2$. $\square$
\endproclaim \bs

{\bf \S 1.10. Examples for The 2nd Algorithm(Theorem 1.15(Corollary
1.15.1)) and The 3rd Algorithm(Theorem 1.16)} \ms

{\bf As Examples for The 2nd Algorithm  and The 3rd Algorithm,}  let
$f(y,z)$ be a W-poly of two complex variables, which is given by the
following:
$$\align
\text{\rm (Eq.1)} \quad
f(y,z)=&z^{16}+4y^{3}z^{14}+\{4y^5+6y^6\}z^{12}+\{12y^8+4y^9\}z^{10}+
\{6y^{10}+12y^{11}+y^{12}\}z^{8} \quad \\
&+\{12y^{13}+4y^{14}+y^{17}\}z^{6}+\{4y^{15}+6y^{16}+y^{20}\}z^{4}
+\{4y^{18}+y^{22}\}z^{2}\\
&+y^{24}z+\{y^{20}+y^{29}\}.
\endalign$$

By the same method and notations as in Theorem 1.15(Corollary
1.15.1) and Theorem 1.16, we are going to find The 2nd Algorithm and
The 3rd Algorithm for $f(y,z)$ of (Eq.1). \ms

\noindent$\underline{\text{\bf Example 1.10.1 for The 2nd Algorithm:
To find irreducibility of $f$ of (Eq.1) in ${}_2\CO_0$}}$ \ms

We can find The 2nd algorithm in both Theorem 1.15 and Corollary
1.15.1 to compute irreducibility of $f$ of (Eq.1) in ${}_2\CO_0$. To
find a solution,  we use the 1st method defined by The 2nd Algorithm
in Theorem 1.15,  and the 2nd method defined by The 2nd Algorithm in
Corollary 1.15.1, respectively. \bs

\noindent{\bf [The 1st Method] We use The 2nd algorithm in Theorem
1.15, to find a solution.} \ms

\noindent$\underline{\text{\bf The 1st step: To find the
irreducibility algorithm for $f\in $ the type [${\ell}$] with
{${\ell}\ge 1$}}}$

\noindent$\underline{\text{\bf in the sense of  Definition 2.5.}}$
\ms

If $f$ is irreducible in ${}_2\CO_0$, then $f$ must satisfy the
following necessary condition:
$$\align
 \f{\a_i}{n-i}\ge \f{\a_0}n \quad
\text{and} \quad{1\le\gcd(n,\a_0)<n} \quad \text{for \quad $0\le
i\le n-2$.} \tag Eq.2
\endalign$$
Note that $1<\gcd(n,\a_0)=4<n$ where $n=16$ and $\a_0=20$.

\noindent{\bf {Remark.}} It is clear that $\f{\a_i}{n-i}\ge
\f{\a_0}n=\f{5}{4}$ for $0\le i\le 14$, if exists. Since
$1<\gcd(n,\a_0)<n$, $f\not\in $ the type $[1]$ in the sense of
Definition 2.5, and so either $f\in $ the type $[\ell]$ for $\ell\ge
2$ in the sense of Definition 2.5  or $f$ is not irreducible in
$\BC\{y,z\}$. So, it suffices to follow the second step. \ms

\noindent$\underline{\text{\bf The 2nd step: To find the
irreducibility algorithm for $f\in $ the type [${\ell}$] with
{${\ell}\ge 2$}}}$

\noindent$\underline{\text{\bf in the sense of  Definition 2.5.}}$
\ms

Let $d_2=\gcd(n,\a_0)$, and then write $n=n_1d_2$ and
$\a_0=\beta_{1,0,1}d_2$. Then, $n_1=4$ and $\beta_{1,0,1}=5$ with
$d_2=4$. If $f$ is irreducible in ${}_2\CO_0$, it suffices to show
that $f$ with $g_1$ in (Eq.1) can be represented as follows:
$$\align
\text{(Eq.3)}\quad\text{(a)} \quad
g_1&=z^{n_1}+\xi_1y^{\beta_{1,0,1}} =z^{4}+y^{5} \quad \text{with
$\xi_1=1$}, \qquad \qquad \\
\text{(b)}\quad f &=(z^{n_1}+ \xi_1 y^{\beta_{1,0,1}})^{d_2}
+\sum_{\alpha, \beta\ge 0} c_{\alpha,\beta}y^{\alpha}z^{\beta} \quad
\text{with} \quad n_1\alpha+\beta_{1,0,1}\beta>n_1\beta_{1,0,1}d_2,
\quad
\endalign$$
\noindent{\bf {Remark.}} It is trivial that $(g_1, f)$ satisfies the
inequalities in (Eq.3).

Whether the W-poly $f(y,z)$  is irreducible in $\BC\{y,z\}$ or not,
apply the WDT(Theorem 1.7) to $f(y,z)$ with a divisor $g_1(y,z)$. By
The Division Alorithm for the W-polys(Theorem 1.8), we can show that
$(g_1,f)$ can be written in the form
$$\cases
g_1 &=z^{4}+y^{5} \quad \text{with
$h_{1,1}=g_1$}, \\
f &=g^{d_{2}}_1+\sum^{d_{2}-1}_{i=1}T^{(1)}_{2,i}g^i_1,\
\text{with}\ T^{(1)}_{2,d_2-1}\not=0,
\endcases \tag {Eq.4)(Eq.4.1}
$$
with the following property:
$$\align
\text{\rm(Eq.4.1)(Eq.4.1.1)} \quad
&T^{(1)}_{2,3}=\{4y^{3}z^2+6y^6\},\quad
T^{(1)}_{2,2}=\{4y^{9}z^2-6y^{11}+y^{12}\},\\ \quad
&T^{(1)}_{2,1}=\{-4y^{14}z^2+y^{17}z^2-2y^{17}+y^{20}\},\quad
T^{(1)}_{2,0}=\{y^{24}z+y^{22}-y^{25}+y^{29}\}.\\
\endalign$$

If $f$ is irreducible in ${}_2\CO_0$, it can be proved that
$(g_1,f)$ satisfies the necessary condition:
$$\align
\text{\rm(Eq.4.1)(Eq.4.1.2)}\qquad \qquad
&\f{\th_2(\beta^{(1)}_{2,i,k})^2_{k=1}}{d_2-i}\ge
\f{\th_2(\beta^{(1)}_{2,0,k})^2_{k=1}}{d_2}>n_1\beta_{1,0,1} \quad
\text{for \quad $0\le i\le d_2-1$} \qquad \qquad\\
\text{(Eq.4.1.3)}\qquad \qquad&\gcd(d_2,\th_2(\beta^{(1)}_{2,0,k})^2_{k=1})=4\ge {1}. \\
\endalign$$

\noindent{\bf {Remark.}} It is easy to compute that $(g_1,f)$
satisfies the inequalities in (Eq.4.1.2) and (Eq.4.1.3) because
$\f{\th_2(\beta^{(1)}_{2,i,k})^2_{k=1}}{d_2-i}=22$ for $0\le i\le 3$
and $\gcd(d_2,\th_2(\beta^{(1)}_{2,0,k})^2_{k=1})=\gcd(4,88)=4>1$.

Since $\gcd(d_2,\th_2(\beta^{(1)}_{2,0,k})^2_{k=1})>1$ and the
coefficient $T^{(1)}_{2,d_2-1}=4y^3z^2+6y^6\ne 0$, to find a
necessary and sufficient condition for $f$ to be irreducible in
${}_2\CO_0$, first of all, we can use the same notations and
properties as in $\underline{\text{\rm Case(B) of The 2nd step in
the conclusions of Theorem 1.15}}$. Then, it suffices to consider
the following:

$\underline{\text{\rm Subcase(B1) of Case(B) of the 2nd step:}}$
Since $\gcd(d_2,\th_2(\beta^{(1)}_{2,0,k})^2_{k=1})=4>1$ and
$T^{(1)}_{2,d_2-1}=4y^3z^2+6y^6\ne 0$, to find whether the above
W-poly $f(y,z)$ is irreducible in $\BC\{y,z\}$ or not, by the
WDT(Theorem 1.7) only $(g_1,f)$ can be rewritten by $(h_{1,2},f)$
$$\cases
h_{1,2} &=h_{1,1}+\f
1{d_2}T^{(1)}_{2,d_2-1}=z^{n_1}+\sum^{n_1-2}_{i=0}R^{(2)}_{1,i}z^i\\
f &=h^{d_2}_{1,2}+\sum^{d_2-1}_{i=0}T^{(2)}_{2,i}h^i_{1,2}\
\text{with}\ T^{(2)}_{2,d_2-1}\not=0,
\endcases \tag Eq.5)(Eq.5.1
$$
with the following property:
$$\align
\noindent \text{\rm(Eq.5.1)(Eq.5.1.1)}
  \quad h_{1,2}&=g_1+\f{1}{4}T^{(1)}_{2,d_2-1}
=z^4+y^3z^2+\{y^5+\f{3}{2}y^6\}=z^{n_1}+\sum^{n_1-2}_{i=0}R^{(2)}_{1,i}z^i, \qquad \qquad \qquad\\
\quad \text{and} \quad  T^{(2)}_{2,3}=-6y^6,&
T^{(2)}_{2,2}=\f{27}{2}y^{12},
T^{(2)}_{2,1}=-\f{27}{2}y^{18}+y^{17}z^2, T^{(2)}_{2,0}=y^{24}z+\f{81}{16}y^{24}-y^{25}+y^{29}. \\
\endalign$$

If $f$ is irreducible in ${}_2\CO_0$, it can be proved that
$(h_{1,2},f)$ satisfies the necessary condition. In particular, if
$\gcd(d_2,\th_2(\beta^{(2)}_{2,0,k})^2_{k=1})=1$ in (Eq.5.1.3), $f$
must be irreducible in ${}_2\CO_0$:
$$\align
\text{\rm(Eq.5.1.2)}\qquad \qquad
&\f{\th_2(\beta^{(2)}_{2,i,k})^2_{k=1}}{d_2-i} \ge
\f{\th_2(\beta^{(2)}_{2,0,k})^2_{k=1}}{d_2}>n_1\beta_{1,0,1} \quad
\text{for \quad $0\le i\le d_2-1$} \qquad \qquad\\
\text{(Eq.5.1.3)}\qquad \qquad &\gcd(d_2,\th_2(\beta^{(2)}_{2,0,k})^2_{k=1})> 1. \\
\endalign$$

\noindent{\bf {Remark.}} It is clear that $h_{1,2}\buildrel
\text{{\rm multiseq}} \over \sim h_{1,1}=g_1$ and note that
$(h_{1,2},f)$ satisfies the inequalities in (Eq.5.1.2) and
(Eq.5.1.3) because $\th_2(\beta^{(2)}_{2,3,k})^2_{k=1})=24$,
$\th_2(\beta^{(2)}_{2,2,k})^2_{k=1})=48$,
$\th_2(\beta^{(2)}_{2,1,k})^2_{k=1})=72$, and
$\th_2(\beta^{(2)}_{2,0,k})^2_{k=1})=96$. So, we need the next step
to find that $f$ is irreducible in ${}_2\CO_0$. \ms

$\underline{\text{\rm Subcase(B2) of Case(B) of the 2nd step:}}$
Since $T^{(2)}_{2,d_2-1}=-6y^6\ne 0$ and
$\gcd(d_2,\th_2(\beta^{(2)}_{2,0,k})^2_{k=1})> 1$, to find whether
the above W-poly $f(y,z)$ is irreducible in $\BC\{y,z\}$ or not, if
necessay, for the computation of this step, the proof of the
inequalities in \text{\rm(Eq.5.1.2)} and \text{\rm(Eq.5.1.3)} may
not be needed for \text{\rm(Eq.5.1.1)} of {\bf Subcase(B1) of
Case(B) of the 2nd step}.

To find whether the above W-poly $f(y,z)$ is irreducible in
$\BC\{y,z\}$ or not, by the WDT(Theorem 1.7) only $(h_{1,2},f)$ can
be rewritten by $(h_{1,3},f)$
$$\cases
h_{1,3} &=h_{1,2}+\f
1{d_2}T^{(2)}_{2,d_2-1}=z^{n_1}+\sum^{n_1-2}_{i=0}R^{(3)}_{1,i}z^i\\
f &=h^{d_2}_{1,3}+\sum^{d_2-1}_{i=0}T^{(3)}_{2,i}h^i_{1,2}\
\text{with}\ T^{(2)}_{2,d_2-1}=0,
\endcases \tag Eq.5)(Eq.5.2
$$
with the following property:
$$\align
\text{\rm(Eq.5.2.1)}\qquad \qquad
f&=h^4_{1,3}+y^{17}z^2h^2_{1,3}+y^{24}z+y^{29}
=h^{d_2}_{1,3}+\sum^{d_2-2}_{i=0}T^{(3)}_{2,i}h^i_{1,3}, \qquad \qquad \\
\quad h_{1,3}&=h_{1,2}+\f{1}{4}T^{(2)}_{3,d_2-1}
=h_{1,2}+\f{1}{4}\{-6y^6\}=z^4+y^5+y^3z^2=z^{n_1}+\sum^{n_1-2}_{i=0}R^{(3)}_{1,i}z^i, \qquad \qquad\\
\text{where} \quad & T^{(3)}_{2,3}=0, T^{(3)}_{2,2}=0,
T^{(3)}_{2,1}=y^{17}z^2, T^{(3)}_{2,0}=y^{24}z+y^{29}. \\
\endalign$$

Now, if $f$ is irreducible in ${}_2\CO_0$, it can be proved that
$(h_{1,2},f)$ satisfies the necessary condition:
$$\align
\text{(\rm Eq.5.2.2)} \qquad \qquad
&\f{\th_2(\beta^{(3)}_{2,i,k})^2_{k=1}}{d_2-i}\ge
\f{\th_2(\beta^{(3)}_{2,0,k})^2_{k=1}}{d_2}>n_1\beta_{1,0,1} \quad
\text{for $0\le i\le d_2-2$,} \qquad \qquad \qquad \\
\text{(Eq.5.2.3)}\qquad \qquad &\gcd(d_2,\th_2(\beta^{(3)}_{2,0,k})^2_{k=1})=1. \\
\endalign$$

\noindent{\bf {Remark.}} {\rm(i)} It is clear that $h_{1,3}\buildrel
\text{{\rm multiseq}} \over \sim h_{1,1}=g_1$ and note that
$(h_{1,3},f)$ satisfies the inequalities in (Eq.5.2) because
$\th_2(\beta^{(3)}_{2,1,k})^2_{k=1}=98$, and
$\th_2(\beta^{(3)}_{2,0,k})^2_{k=1}=101$. So, $f$ is irreducible in
${}_2\CO_0$.

{\rm(ii)} It can be computed that $T^{(2)}_{2,d_2-1}=-6y^6\ne 0$ and
$\gcd(d_2,\th_2(\beta^{(3)}_{2,0,k})^2_{k=1})=1$. So
$(f_1,f)=(h_{1,\nu+1},f)=(h_{1,3},f)$ is a generalized
representation of $f$ in the sense of Theorem $1.15$. \quad
$\square$ \bs

\noindent{\bf [The 2nd Method] We use The 2nd algorithm in Corollary
1.15.1, to find a solution.} \ms

To compute irreducibility of $f$ of (Eq.1) in ${}_2\CO_0$ by The 2nd
Algorithm in Corollary 1.15.1, it suffices to consider the
following:

{\rm(a)} Since $T^{(1)}_{2,3}=\{4y^{3}z^2+6y^6\}\ne 0$, in this case
to compute the inequality in {\rm (Eq.4.1.2)} and {\rm(Eq.4.1.2)}
may not be needed, if necessary, in order to find a solution,\ms

{\rm(b)} Since $T^{(2)}_{2,d_2-1}=-6y^6\ne 0$ with $d_2=4$, in this
case to compute the inequality in {\rm(Eq.5.1.2)} and
{\rm(Eq.5.1.3)} may not be needed if necessary, in order to find a
solution. \ms

{\rm(c)} Since $T^{(3)}_{2,d_2-1}=0$ with $d_2=4$, in this case to
compute the inequality in {\rm(Eq.5.2.2)} and {\rm(Eq.5.2.3)} must
be needed if necessary, in order to find a solution. \ms

Thus, without provig the inequality in {\rm (Eq.4.1.2)} and
{\rm(Eq.4.1.2)} of {\rm(a)} and the inequality in {\rm(Eq.5.1.2)}
and {\rm(Eq.5.1.3)} of {\rm(b)} in Theorem 15.5, following the same
methods and notations as in {\rm(c)} we can find the same solution
as in Theorem 15.5. \quad $\square$ \bs

\noindent$\underline{\text{\bf Example 1.10.2 for The 3rd Algorithm:
To find the standard Puiseux expansion}}$

\noindent$\underline{\text{\bf C(f) for an irreducible W-poly $f$ of
(Eq.1) such that $f \buildrel {\text{\rm multiseq}} \over \sim C(f)$
at $0\in\BC^2$}}$ \ms

By Remark for the 2nd step in the conclusion of Theorem 1.15, we can
find $H_2\in {\text{\rm Family[1]}}$ for any $f\in $ the type [2] in
the sense of  Definition 2.5 such that $f\buildrel \text{{\rm
multiseq}} \over \sim H_2=H^{d_2}_1+\Pi^2_{k=1}H^{\si_{2,k}}_{k-2}$
where $H_{-1}=y$, $H_0=z$, $H_1=z^{4}+y^{5}\in {\text{\rm
Family[1]}}$ and $\si_{2,k}=\beta^{({\nu_1}+1)}_{2,0,k}$ for $1\le
k\le 2$ and for some integer ${\nu_1}+1>0$, because of the
equations,(Eq.3), (Eq.4) and (Eq.5) in the conclusion of Theorem
1.15, noting that $d_2=4$, $\si_{2,1}=24$ and $\si_{2,2}=1$. Let
$C(H_q:t)$ be the the standard Puiseux expansion defined by
$$\align
\text{(\rm Eq.6)} \quad \quad \text{$C(H_q:t):=$} &
\left\{\eqalign{y&=t^n \cr z&=t^{\gamma_1}+t^{\gamma_2}+\cdots
+t^{\gamma_q}, \cr} \right.
 \\
 \text{such that}  \quad
 d_2  &=\gcd(n,\gamma_1) \quad\text{with}\quad n=n_1d_2=12 \quad\text{and}\quad
 \gamma_1=\alpha_{1,0,1}d_2=20, \quad \text{and}\\
\quad d_3&=\gcd(d_2,\gamma_2-\gamma_{1})=1 \quad\text{with}\quad
d_2=n_2d_{3} \quad\text{and}\quad
\gamma_2-\gamma_{1}=\widehat{\theta}_2d_3
  \endalign$$
where $d_3=1$ and $\widehat{\theta}_j
=\theta_j(\beta_{j,0,k})^j_{k=1}-n_jn_{j-1}\theta_{j-1}(\beta_{j-1,0,k})^{j-1}_{k=1}$
is a positive integer for $2\le j\le r$ and $\theta_1(t)=t$. Note
that $d_j=n_{j}d_{j+1}$ for $1\le j\le r-1$ . \ms

By Theorem 1.4, it remains to compute $\gamma_2-\gamma_{1}$, which
follows from

$\widehat{\theta}_2d_3=\theta_2(\beta_{2,0,k})^2_{k=1}-n_2n_1\theta_{1}(\beta_{1,0,1})
=\theta_j(\sigma_{2,1},\sigma_{2,2})-n_2n_1\alpha_{1,0,1}=\theta_2(24,1)-5\cdot4\cdot4=21$.

Thus, $C(H_q:t)$ can be computed by $y=t^{16}$ and
$z=t^{20}+t^{41}$. $\square$  \bs

\vfill \pagebreak

$$\align
\quad & \text{\bf \centerline{Part[B](Part[B1],\dots, Part[B5])}}\\
& \text{\bf Explicit algorithm for computing the correspondence
between the irreducible}\\
&\text{\bf W-polys of two complex variables and the Puiseux
expansions with proofs}
\endalign$$ \ms

{\bf Part[B1] Foundations} \ms

{\bf \S 2. New definitions for quasisingularity and a generalization
of one coordinate patch covering of the local coordinates used in
the standard resolution process of irreducible plane curve
singularities} \ms

\definition{Definition 2.1}
Let $\BC\{y,z\}$ be the ring of convergent power series or analytic
functions at $(y,z)=(0,0)$. Let $V(F)=\{(y,z):F=F(y,z)=0\}$ be an
analytic variety at $(y,z)=(0,0)$ where $F$ is in $\BC\{y,z\}$.
Assume that $V(F)$ has an isolated singular point at the origin in
$\BC^2$ as reduced variety. Note that $F$ may have multiple factors
as analytic function at the origin. Let $\pi_1:M\to \BC^2$ be the
blow-up of $\BC^2$ at $(y,z)=(0,0)$, which is the singular point
$(0,0)$ of $V(F)$. Let $U_1=(v_1,u_1)$ and $U_2=(v'_1, u'_1)$ be
coordinate patches for $M$ with $\pi_1(v_1,u_1)=(y,z)=(v_1,v_1u_1)$
and $\pi_1(v'_1,u'_1)=(y,z)=(v'_1u'_1,v'_1)$ where $u'_1=
\frac{1}{u_1}$ and $v'_1=v_1u_1$. For brevity of notations, the
above $F(y,z)$ is square-free in $\BC\{y,z\}$ with isolated
singularity at the origin. Let $e$ be the multiplicity of $V(F)$ at
$(0,0)$, with $e\ge 2$.

Then ${\pi_1}^{-1}(V(F))$, the total transform of $V(F)$ under
$\pi_1$, is given locally by
$(F\circ\pi_1)_{total}=F(v_1,v_1u_1)=v_1^eF_1(v_1,u_1)$ along
$v_1=0$ where $F_1(v_1,u_1)$ is in $\BC\{v_1,u_1\}$ and
$(F\circ\pi_1)_{total}=F(v'_1u'_1,v'_1)={v'_1}^e F_2(v'_1,u'_1)$
along $v'_1=0$ where $F_2(v'_1,u'_1)$ is in $\BC\{v'_1,u'_1\}$. We
call $V^{(1)}(F)$ the proper transform of $V(F)$ under $\pi_1$ where
$V^{(1)}(F)=\{(v_1,u_1):F_1(v_1,u_1)=0\}\cup
\{(v'_1,u'_1):F_2(v'_1,u'_1)=0\}$. We say that
$E_1=\{(v_1,u_1):v_1=0\}\cup\{(v'_1,u'_1):v'_1=0\}$ is an
exceptional curve of the first kind. In this case, each of
$U_1=\{(v_1,u_1)\}$ and $U_2=\{(v'_1,u'_1)\}$ is called one
coordinate patch of the given local coordinates for blow-up $\pi_1$,
respectively. Note that if $F(y,z)$ is irreducible in $\BC\{y,z\}$,
then just one coordinate patch is needed for the study of
$V^{(1)}(F)$.

After $m$ iterations of blow-ups, let
$\tau_m=\pi_1\circ\pi_2\circ\cdots \circ\pi_m:M^{(m)}\to \BC^2$. Let
$V^{(m)}(F)$ be the proper transform of $V(F)$ under $\tau_m$. Let
$E^{(m)}=\tau^{-1}_m(0,0)$. Then $E^{(m)}$ is, by definition, an
exceptional set of the first kind. Let $E^{(m)}= \cup^m_{i=1} E_i$
be the decomposition into irreducible components. Each $E_i$ is
called an exceptional curve of the first kind. Let
$(F\circ\tau_m)_{divisor}$ be the divisor of $F\circ\tau_m$ defined
by $(F\circ\tau_m)_{divisor}=V^{(m)}(F)+\sum^m_{i=1}e_iE_i$ where
each $e_i$ is the multiplicity of $F\circ\tau_m$ along $E_i$ for
$1\le i\le m$.
\enddefinition

Then, we have the following well-known theorem.

\proclaim{Theorem 2.2} Let $V(F)=\{(y,z): F(y,z)=0\}$ be an analytic
variety at $(0,0)$ where $F(y,z)$ is in $\BC\{y,z\}$ with isolated
singularity at the origin. There exists an analytic manifold $M$ by
using the composition of a finite number $m$ of successive blow-ups,
$\tau_m:M\to \BC^2$, such that if $R$ is the set of regular points
on $V(F)$ then $\tau_m:\overline{\tau^{-1}_m(R)}\to V(F)$ is a
resolution of a singular point $(0,0)$ of $V(F)$, where
$\overline{\tau^{-1}_m(R)}$ is the closure of $\tau^{-1}_m(R)$ in
$M$.
\endproclaim

\definition{Remark 2.2.1} {\rm(i)} If $V$ is an analytic set, a
resolution of the singularities $f$ consists of a complex manifold
and a proper analytic map $\pi:M\rightarrow V$ such that $\pi$ is
biholomorphic on the inverse image of $R$, the regular points of
$V$, and such that $\pi^{-1}(R)$ is dense in $M$.

{\rm(ii)} Blow-ups are canonical. Namely, let $\phi:\BC^2\to \BC^2$
be a biholomorphic map, and let ${\pi'}:M'\to \BC^2$ be a blow-up at
$\phi(0,0)$. Then, there is a unique induced biholomorphic map
$\phi:M\to {M'}$ such that $\phi\circ\pi_1={\pi'}\circ\phi'$.
\enddefinition

\proclaim{Corollary 2.3} Under the same assumption of Theorem
$2.2$, after additional blow-ups any two components of
$V^{(m)}(F)$ and $\cup^m_{i=1} E_i$ meets with normal crossings
whenever they meet and no three distinct components of
$V^{(m)}(F)$ and $\cup E_i$ meet, where $V^{(m)}(F)$ and $\cup
E_i$ are defined just before Theorem $2.2$.
\endproclaim

\definition{Remark 2.3.1} If $F(y,z)$ of Corollary $2.3$ is
irreducible in $\BC\{y,z\}$ with isolated singularity at the origin,
then each exceptional curve $E_i$ of the first kind meets at most
three distinct intersection points with other exceptional curves and
the proper transform $V^{(m)}(F)$. It will be proved by Theorem
$3.6$ and Theorem $3.7$.
\enddefinition

\definition{Definition 2.4(The standard resolution, A homeomorphic resolution,
Having the same divisor under two standard resolutions)}

{\bf(I)} For a singularity of a plane curve, the smallest
resolution with normal crossings in the sense of Corollary $2.3$
is called the standard resolution of a given singularity. \ms

{\bf(II)} As in Definition $2.1$, let $V(F)$ and $V(G)$ be analytic
varieties at $0\in \BC^2$ where $F=F(y,z)$ and $G=G(y,z)$ are in
$\BC\{y,z\}$ with isolated singularity at the origin. $V(F)$ and
$V(G)$ are said to either have a homeomorphic resolution, or be
equisingular, if $(F\circ\tau_m)_{divisor}$ and
$(G\circ\tau_m)_{divisor}$ are equivalent in the sense of Definition
$2.1$ where $\tau_m=\pi_1\circ\pi_2\circ\cdots \circ\pi_m:M^{(m)}\to
\BC^2$ is the composition of the same number $m$ of successive
blow-ups at the origin, which is the standard resolution of the
singularity $(0,0)$ of both $V(F)$ and $V(G)$. Then, we denote this
relation by either $V(F) \buildrel \text{{\rm resol}} \over \sim
V(G)$ at $0\in \BC^2$ or $(F\circ\tau_m)_{divisor} =
(G\circ\tau_m)_{divisor}$ under the same standard resolution
$\tau_m$ of both $V(F)$ and $V(G)$. In other words, if we write
$(F\circ\tau_m)_{divisor}=V^{(m)}(F)+\sum^m_{i=1}e_iE_i$ and
$(G\circ\tau_m)_{divisor}=V^{(m)}(G)+\sum^m_{i=1}{\bar{e_i}}E_i$ in
the sense of Definition $2.1$, then it is said that $V(F) \buildrel
\text{{\rm resol}} \over \sim V(G)$ at $0\in \BC^2$ under the same
standard resolution $\tau_m$ of both $V(F)$ and $V(G)$ if and only
if $e_i={\bar e_i}$. \ms

{\bf(III)} In particular, let $f\in \text{\rm Family(0)}$ and $g\in
\text{\rm Family(0)}$ in the sense of Definition 1.2, that is, $f$
and $g$ be analytically irreducible in $\BC\{y,z\}$ with isolated
singularity at the origin. Let
$\tau_{\xi}=\pi_1\circ\pi_2\circ\cdots \circ\pi_{\xi}:M^{({\xi})}\to
\BC^2$ be the composition of a finite number ${\xi}$ of successive
blow-ups $\pi_i$ at $0\in \BC^2$, which is needed only to get the
standard resolution of the singularity of $V(f)$. Using the same
method as before, let
$\mu_{\eta}=\bar{\pi}_1\circ\bar{\pi}_2\circ\cdots
\circ\bar{\pi}_{\eta}:\bar{M}^{({\eta})}\to \BC^2$ be the
composition of a finite number ${\eta}$ of successive blow-ups at
$0\in \BC^2$, which is needed only to get the standard resolution of
the singularity of $V(g)$.  Following the same properties and
notations as in Definition 2.1, we write
$(f\circ\tau_m)_{divisor}=V^{(m)}(f)+\sum^m_{i=1}e_iE_i$ and
$(g\circ\mu_{\eta})_{divisor}=V^{({\eta})}(g)+\sum^{\eta}_{i=1}{\bar{e_i}}{\bar{E}}_i$
in the sense of Definition $2.1$.

It is said that either $V(f)$ and $V(g)$ have the same divisor under
the standard resolutions, or $(f\circ\tau_{\xi})_{divisor}$ and
$(g\circ\mu_{\eta})_{divisor}$ are equivalent,  denoted by either
$V(f) \buildrel \text{{\rm divisor}} \over \sim V(g)$ or ${f}
\buildrel \text{{\rm divisor}} \over \sim {g}$ under the standard
resolutions or
$(f\circ\tau_{\xi})_{divisor}=(g\circ\mu_{\eta})_{divisor}$ under
the standard resolutions, if the following condition is satisfied:
$$\align
(2.4.1) \quad  \text{\rm either} \quad
&\text{$\{(f\circ\tau_{\xi})_{divisor}\}_{seq.}\equiv
\{(g\circ\mu_{\eta})_{divisor}\}_{seq.}$} \quad \text{\rm as
sequence,} \\
\text{\rm or} \quad
&\text{$\{e_i:i=1,2,\dots,\xi\}\equiv\{{\bar{e}_i}\in
N:i=1,2,\dots,\eta\}$ \quad as an increasing sequence.}\qquad \qquad \\
\endalign$$

{\bf(IV)} $\underline{\text{\bf{Family(4)}}}$ is the fourth family,
consisting of all the divisors of $(f\circ\tau)$ defined by the
total transform of $V(f)$, denoted by $(f\circ\tau)_{divisor}$ where
$\tau:M\to \BC^2$ is the standard resolution of the singularity of
{$V(f)$} with $f\in \text{\rm Family(0)}$, denoted by

\noindent(2.4.2) \qquad
$\underline{\text{\rm{Family(4)}}=\text{\{$(f\circ\tau)_{divisor}$:
$f\in \text{\rm Family(0)}$ where $\tau:M\to \BC^2$ is}}$

\qquad \qquad \qquad \qquad \qquad $\underline{\text{the standard
resolution of the singularity of {$V(f)$}\}}}$. \enddefinition \ms

\definition{Remark 2.4.1} Let $f\in \text{\rm Family(0)}$ and $g\in
\text{\rm Family(0)}$. By Definition 2.4, $V(f) \buildrel \text{{\rm
resol}} \over \sim V(g)$ at $0\in \BC^2$ implies ${f} \buildrel
\text{{\rm divisor}} \over \sim {g}$ under two standard resolutions.
But, it remains to be proved that ${f} \buildrel \text{{\rm
divisor}} \over \sim {g}$ under two standard resolutions implies
$V(f) \buildrel \text{{\rm resol}} \over \sim V(g)$ at $0\in \BC^2$.
\enddefinition \ms

\definition{Definition 2.5} Suppose that $f\in\BC\{y,z\}$ has
an isolated singular point at $0\in \BC^2$. It is said that
\text{$f\in$the type $[j]$} or belongs to the type $[j]$ under the
standard resolution if $f$ satisfies the following properties (a)
and (b) after m iterations of blow-ups which is needed only to get
the standard resolution of the singularity of the curve defined by
$f$: Sometimes, we write \text{$V(f)\in$the type $[j]$}, instead of
\text{$f\in$the type $[j]$}.

(a) There are exactly j exceptional curves of the first kind with
$j\le m$, each of which has three distinct intersections with other
exceptional curves and the proper transform.

(b) Each of the remaining (m-j) exceptional curves rather than the
above j exceptional curves in (a) has at most two distinct
intersection points with other exceptional curves and the proper
transform.
\enddefinition
\ms

\definition{Remark 2.5.1}

(i) If $f$ has a nonsingular point at the origin, we say that $f\in$
the type $[0]$.

(ii) If $f\in$ the type $[j]$ for an integer $j\ge 1$, then $f$
may not be irreducible by the following example:

 Let $f=(z+y)(z+2y)(z+3y)+y^6$. Then $f\in$ the type $[1]$, but
$f$ is not irreducible in $\BC\{y,z\}$.

(iii) Using Definition $2.5$, then the following will proved by
Theorem $12.0$:

Let $V(f)$ be an analytic variety at $(y,z)=(0,0)$ where $f=f(y,z)$
is irreducible in $\BC\{y,z\}$ with isolated singularity at the
origin. As it have been in Definition $2.1$, and Definition $2.4$,
let $\tau_m=\pi_1\circ\pi_2\circ\cdots \circ\pi_m:M^{(m)}\to \BC^2$
be the standard resolution of the singular point $(0,0)$ of $V(f)$.
Then, $f\in$ the type $[j]$ under the standard resolution $\tau_m$
in the sense of Definition $2.5$ where $j$ is a positive integer
with $j<m$.
\enddefinition \ms

\noindent{\bf{Definition 2.6(Quasisingularity, homeomorphic
resolutions, Having the same divisor under two standard resolutions
and equivalence of multiplicity sequences for irreducible plane
curve singularities).}}

\noindent$\underline{\text{\bf (I) Definition for
quasisingularity:}}$ Let $V(F)$ be an analytic variety at
$(y,z)=(0,0)$ where $F=F(y,z)$ is in $\BC\{y,z\}$ with isolated
singularity at the origin. Assume that $V(F)$ satisfies the same
properties and notations as in Definition $2.1$. After $m$
iterations of blow-ups, let $\tau_m:M^{(m)}\to \BC^2$ be an
arbitrary composition of a finite number $m$ of successive blow-ups
at the origin in $\BC^2$, which is the singular point of $V(F)$. As
in Definition $2.1$, let $\tau^{-1}_m(0,0)=\cup^m_{i=1}E_i$ where
each $E_i$ is an called an exceptional curve of the first kind, and
$V^{(m)}(F)$ be the proper transform under $\tau_m$.

For a given $\tau_m$, let $P\in\tau^{-1}_m(0,0)$ be chosen
arbitrary. It is said that P is a quasisingular point of
$V^{(m)}(F)$ under $\tau_m$ if additional blow-ups at $P\in
V^{(m)}(F)$ are still necessary after $m$ iterations of blow-ups at
$(y,z)=(0,0)$, in order to get the standard resolution of the
singular point of $V(F)$. For notation, $qs(V^{(m)}(F))$ is called
the set of all quasisingular points of $V^{(m)}(F)$ under $\tau_m$.
Assuming that $qs(V^{(m)}(F))\not=\emptyset$, note that $P\in
qs(V^{(m)}(F))$ may not be a singular point of $V^{(m)}(F)$. \ms

\noindent$\underline{\text{\bf (II) Homeomorphic resolutions of some
plane curve singularities:}}$

Let $\phi(y,z)=a_0z^n+a_1y^{\alpha_1}z^{n-1}+ \dots+a_ny^{\alpha_n}$
be irreducible in $\BC\{y,z\}$, where $a_0$, $a_n$ are units in
$\BC\{y,z\}$, and also each $a_i$ is a unit in $\BC\{y,z\}$ if
exists, and the $\alpha_i$ are positive integers, and $1\le n<
\alpha_n$.

Let $\psi(y,z)=b_0z^{\ell}+b_1y^{\beta_1}z^{{\ell}-1}+
\dots+b_{\ell}y^{\beta_{\ell}}$ be irreducible in $\BC\{y,z\}$,
where $b_0$, $b_{\ell}$ are units in $\BC\{y,z\}$, and also each
$b_i$ is a unit in $\BC\{y,z\}$ if exists, and the $\beta_i$ are
positive integers, and $1\le {\ell}< \beta_{\ell}$.

For example, let $e$ be the multiplicity of $\phi$ at $0$. Then,
it was well-known by Hensel's Lemma or Lemma $3.1$ that $e=n$ and
$n<\alpha_i+n-i$ for all $i=1,2,\dots,n$.

Let $V(\Phi)$ and $V(\Psi)$ be analytic varieties at $(y,z)=(0,0)$
defined by
$$\align
(2.6.1) \qquad \qquad & \Phi=\Phi(y,z)=y^{\zeta}z^{\eta}\phi(y,z)
\quad \text{and}
\quad \Psi=\Psi(y,z)=y^{\zeta'}z^{\eta'}\psi(y,z), \qquad \qquad \qquad \\
& \text{where each of $\zeta$, $\eta$, $\zeta'$ and $\eta'$ is
either a positive integer or $0$.}
\endalign$$

In order for both $V(\Phi)$ and $V(\Psi)$ to have an isolated
singular point at the origin as reduced varieties, in addition we
may assume that the following property holds: Note by assumption
that $\phi=\phi(y,z)$ and $\psi=\psi(y,z)$ are irreducible in
$\BC\{y,z\}$.

\noindent {\rm{(2.6.2)} \qquad \quad {\rm(i)}  If $n=1$ from
$\Phi$, then either $\zeta$ or $\eta$ is positive.

\qquad \qquad \quad {\rm(ii)}  If ${\ell}=1$ from $\Psi$, then
either $\zeta'$ or $\eta'$ is positive. \ms

As in Definition $2.1$, let $V(\Phi)$ and $V(\Psi)$ be analytic
varieties at $(y,z)=(0,0)$ in $\BC^2$ where $\Phi=\Phi(y,z)$ and
$\Psi=\Psi(y,z)$ are in $\BC\{y,z\}$ with isolated singularity at
the origin. $V(\Phi)$ and $V(\Psi)$ are said to have a homeomorphic
resolution if $(\Phi\circ\tau_m)_{divisor}$ and
$(\Psi\circ\tau_m)_{divisor}$ are equivalent in the sense of
Definition $2.1$ where $\tau_m=\pi_1\circ\pi_2\circ\cdots
\circ\pi_m:M^{(m)}\to \BC^2$ is the composition of the same number
$m$ of successive blow-ups at the origin, which is the standard
resolution of the singularity $(0,0)$ of both $V(\Phi)$ and
$V(\Psi)$. Then, we denote this relation by either $V(\Phi)
\buildrel \text{{\rm resol}} \over \sim V(\Psi)$ under the same
standard resolution or $(\Phi\circ\tau_m)_{divisor} \buildrel
\text{{\rm resol}} \over \sim (\Psi\circ\tau_m)_{divisor}$ under the
same standard resolution $\tau_m$ of both $V(\Phi)$ and $V(\Psi)$.
In other words, if we write
$(\Phi\circ\tau_m)_{divisor}=V^{(m)}(\Phi)+\sum^m_{i=1}e_iE_i$ and
$(\Psi\circ\tau_m)_{divisor}=V^{(m)}(\Psi)+\sum^m_{i=1}{\bar{e_i}}E_i$
in the sense of Definition $2.1$, then it is said that $V(\Phi)
\buildrel \text{{\rm resol}} \over \sim V(\Psi)$ at the origin in
$\BC^2$ under the same standard resolution $\tau_m$ of both
$V(\Phi)$ and $V(\Psi)$ if and only if $e_i={\bar e_i}$ under the
same standard resolution $\tau_m$. \ms

As an application of Definition $2.4$, we have the following
definition for $V(\Phi)$ and $V(\Psi)$ to have a homeomorphic
resolution at the origin in $\BC^2$, denoted by $V(\Phi) \buildrel
\text{{\rm resol}} \over \sim V(\Psi)$ at the origin in $\BC^2$:
$$\align
  & \text{$V(y^{\zeta}z^{\eta}\phi) \buildrel \text{{\rm
resol}} \over \sim V(y^{\zeta'}z^{\eta'}\psi)$ under the same
standard resolution}  \tag 2.6.3\\
\text{$\Longleftrightarrow$} \quad \quad &\text{the following
conditions
are satisfied:} \\
&  \text{{\rm (i)}} \quad \text{$\zeta=\zeta'$ and $\eta=\eta'$.} \\
&  \text{{\rm (ii)}} \quad  \text{$e_{j}=\bar{e}_{j}$ \quad for each
$j=1,2,\dots,m$.}
\endalign$$

\noindent$\underline{\text{\bf (III) Having the same divisor under
two standard resolutions:}}$

By the same way as in (II), assume that $\phi=\phi(y,z)$ and
$\psi=\psi(y,z)$ are irreducible in $\BC\{y,z\}$,  satisfying the
same properties and notations as in (II). Also, let $V(\Phi)$ and
$V(\Psi)$ be analytic varieties at $(y,z)=(0,0)$ defined by
equalities in (2.6.1) and (2.6.2).

Let $\tau_{\lambda}=\pi_1\circ\pi_2\circ\cdots
\circ\pi_{\lambda}:M^{(\lambda)}\to \BC^2$ be the composition of a
finite number $\lambda$ of successive blow-ups $\pi_i$ at the origin
in $\BC^2$, which is needed only to get the standard resolution of
the singular point of an analytic variety $V(\Phi)$ where $\Phi$ is
defined as above. Let
$\tau^{-1}_{\lambda}(0,0)=\cup^{\lambda}_{i=1}E_{i}$ where each
$E_{i}$ is called an exceptional curve of the first kind.

By Definition $2.1$, let $(\Phi\circ\tau_{\lambda})_{divisor}$ be
the divisor of $\Phi\circ\tau_{\lambda}$ defined by
$$\align
(\Phi\circ\tau_{\lambda})_{divisor}
=V^{(\lambda)}(\Phi)+\sum^{\lambda}_{i=1}e_{i}E_{i},  \tag 2.6.4
\endalign$$ where each $e_{i}$ is the multiplicity of
$\Phi\circ\tau_{\lambda}$ along $E_{i}$ for $1\le i\le \lambda$ and
$V^{(\lambda)}(\Phi)$ is the proper transform of $V(\Phi)$ under
$\tau_{\lambda}$.

In more detail, let $\pi_{i}:M^{(i)}\to M^{(i-1)}$ be the {\rm
{i}}-th blow-up of $\tau_{\lambda}$ at a quasisingular point of
$V^{(i-1)}(\Phi)$ for $1\le i< \lambda$ in the sense of (I) of
Definition $2.6$ where $M^{(0)}=\BC^2$ and $V^{(0)}(\Phi)=V(\Phi)$.
Then, it is clear by Remark $2.6.0$, later that $V^{(i)}(\Phi)$ has
one and only one quasisingular point along $E_{i}$ for each
$i=1,2,\dots, {\lambda-1}$. So, it is well-defined that $E_{i}$ is
called the {\rm {i}}-th exceptional curve from
$\tau^{-1}_{\lambda}(0,0)=\cup^{\lambda}_{i=1}E_{i}$. Therefore, it
is clear from $(\Phi\circ\tau_{\lambda})_{divisor}$ of (2.6.4) that
$e_{i+1}>e_{i}$ for $1\le i\le \lambda-1$. \ms

Using the same method as before, let
$\mu_{\sigma}=\bar{\pi}_1\circ\bar{\pi}_2\circ\cdots
\circ\bar{\pi}_{\sigma}:\bar{M}^{(\sigma)}\to \BC^2$ be the
composition of a finite number ${\sigma}$ of successive blow-ups at
the origin in $\BC^2$, which is needed only to get the standard
resolution of the singular point of an analytic variety $V(\Psi)$
where $\Psi$ is defined as above. Let
$\mu^{-1}_{\sigma}(0,0)=\cup^{\sigma}_{j=1}\bar{E}_{j}$ for
$V(\Psi)$ where each $\bar{E}_{j}$ is called an exceptional curve of
the first kind.

By Definition $2.1$, let $(\Psi\circ\mu_{\sigma})_{divisor}$ be the
divisor of $\Psi\circ\mu_{\sigma}$ defined by
$$(\Psi\circ\mu_{\sigma})_{divisor}
=V^{({\sigma})}(\Psi)+\sum^{\sigma}_{j=1}\bar{e}_{j}\bar{E}_{j},
\tag 2.6.5$$ where each $\bar{e}_{j}$ is the multiplicity of
$\Psi\circ\mu_{\sigma}$ along $\bar{E}_{j}$ for $1\le j\le {\sigma}$
and $V^{({\sigma})}(\Psi)$ is the proper transform of $V(\Psi)$
under $\mu_{\sigma}$. \ms

As an application of Definition $2.4$, we have the following
definition for $V(\Phi)$ and $V(\Psi)$ to have the same divisor
under two standard resolutions, denoted by $V(\Phi) \buildrel
\text{{\rm divisor}} \over \sim V(\Psi)$ under two standard
resolutions:
$$\align
  & \text{$V(y^{\zeta}z^{\eta}\phi) \buildrel \text{{\rm
divisor}} \over \sim V(y^{\zeta'}z^{\eta'}\psi)$ under two
standard resolutions}  \tag 2.6.6\\
\text{$\Longleftrightarrow$} \quad \quad &\text{the following
conditions
are satisfied:} \\
&  \text{{\rm (i)}} \quad \text{$\zeta=\zeta'$, $\eta=\eta'$ and
$\lambda={\sigma}$.} \\
&  \text{{\rm (ii)}} \quad \text{We may assume that $E_{j}$ and
$\bar{E}_{j}$ are the same} \\
& \text{\quad \quad and so $e_{j}=\bar{e}_{j}$ for each
$j=1,2,\dots,\lambda$.}
\endalign$$

\noindent$\underline{\text{\bf (IV) Equivalence of multiplicity
sequences of irreducible plane curve singularities:}}$

As in (II), let $\Phi(y,z)=\phi(y,z)=a_0z^n+a_1y^{\alpha_1}z^{n-1}+
\dots+a_ny^{\alpha_n}$ be irreducible in $\BC\{y,z\}$ with isolated
singularity at the origin, where $a_0$, $a_n$ are units in
$\BC\{y,z\}$, and also each $a_i$ is a unit in $\BC\{y,z\}$ if
exists, and the $\alpha_i$ are positive integers, and $2\le n<
\alpha_n$. Let $\tau_{\lambda}=\pi_1\circ\pi_2\circ\cdots
\circ\pi_{\lambda}:M^{(\lambda)}\to \BC^2$ be the composition of a
finite number $\lambda$ of successive blow-ups $\pi_i$ at the origin
in $\BC^2$, which is needed only to get the standard resolution of
the singular point of an analytic variety $V(\Phi)$ where $\Phi$ is
defined as above. Let
$\tau^{-1}_{\lambda}(0,0)=\cup^{\lambda}_{i=1}E_{i}$ where each
$E_{i}$ is called an exceptional curve of the first kind.

Now, we may use the same notations and properties as in (II). Since
$\Phi(y,z)$ be irreducible in $\BC\{y,z\}$ with isolated singularity
at the origin, then $V^{(i)}(\Phi)$ with has one and only one
quasisingular point along $E_{i}$ for each $i=1,2,\dots,
{\lambda-1}$.

In more detail, let $\pi_{i}:M^{(i)}\to M^{(i-1)}$ be the {\rm
{i}}-th blow-up of $\tau_{\lambda}$ at a quasisingular point of
$V^{(i-1)}(\Phi)$ for $1\le i< \lambda$ in the sense of (I) of
Definition $2.6$ where $M^{(0)}=\BC^2$ and $V^{(0)}(\Phi)=V(\Phi)$.
Then, it may be assumed that $V^{(i)}(\Phi)$ has the multiplicity
$\nu_{i+1}$ at one and only one quasisingular point along $E_{i}$
for each \text{\rm i=1,2,\dots, ${\lambda}-1$}}. By definition,
$\{\nu_{i}:i=1,2,\dots,\lambda\}$ is called a multiplicity sequence
of $V(\Phi)$, denoted by $\text{\rm Multiseq(V($\Phi$))}$, where
$\nu_1=n$ is the multiplicity of $V(\Phi)$ at $(y,z)=(0,0)$.

Let $\psi(y,z)=b_0z^{\ell}+b_1y^{\beta_1}z^{{\ell}-1}+
\dots+b_{\ell}y^{\beta_{\ell}}$ be irreducible in $\BC\{y,z\}$ with
isolated singularity at the origin, where $b_0$, $b_{\ell}$ are
units in $\BC\{y,z\}$, and also each $b_i$ is a unit in $\BC\{y,z\}$
if exists, and the $\beta_i$ are positive integers, and $2\le
{\ell}< \beta_{\ell}$. Using the same method as before, let
$\mu_{\sigma}=\bar{\pi}_1\circ\bar{\pi}_2\circ\cdots
\circ\bar{\pi}_{\sigma}:\bar{M}^{(\sigma)}\to \BC^2$ be the
composition of a finite number ${\sigma}$ of successive blow-ups at
the origin in $\BC^2$, which is needed only to get the standard
resolution of the singular point of an analytic variety $V(\Phi)$
where $\Phi$ is defined as above. Let
$\mu^{-1}_{\sigma}(0,0)=\cup^{\sigma}_{j=1}\bar{E}_{j}$ for
$V(\Phi)$ where each $\bar{E}_{j}$ is called an exceptional curve of
the first kind.

In more detail, let $\bar{\pi}_{i}:\bar{M}^{(i)}\to \bar{M}^{(i-1)}$
be the {\rm {i}}-th blow-up of $\mu_{\sigma}$ at a quasisingular
point of $V^{(i-1)}(\Phi)$ for $1\le i< \lambda$ in the sense of (I)
of Definition $2.6$ where $\bar{M}^{(0)}=\BC^2$ and
$V^{(0)}(\Psi)=V(\Psi)$. Then, it may be assumed that
$V^{(i)}(\Psi)$ has the multiplicity $\eta_{i+1}$ at one and only
one quasisingular point along $\bar{E}_{i}$ for each $i=1,2,\dots,
{\sigma}-1$. By definition, $\{\eta_{i}:i=1,2,\dots,\sigma\}$ is
called a multiplicity sequence of $V(\Psi)$, denoted by $\text{\rm
Multiseq(V($\Psi$))}$, where $\eta_1={\ell}$ is the multiplicity of
$V(\Psi)$ at $(y,z)=(0,0)$.

If either $\text{\rm Multiseq(V($\Phi$))}=\text{\rm
Multiseq(V($\Psi$))}$ as sequence, or $\nu_i=\eta_i$ for each
$i=1,2,\dots,\lambda=\sigma$, then it is said that either
\text{$\Phi \buildrel \text{{\rm multiseq}} \over \sim \Psi$},  or
\text{$V(\Phi) \buildrel \text{{\rm multiseq}} \over \sim V(\Psi)$}
\text{at $0\in \BC^2$}. Otherwise, we write $\text{\rm
Multiseq(V($\Phi$))}\not =\text{\rm Multiseq(V($\Psi$))}$ as
sequence. \ms

\definition{Remark 2.6.0}
{\rm (3)} As in Definition $2.6$, let $V(\Phi)$ and $V(\Psi)$ have
a homeomorphic resolution at the origin in $\BC^2$.

{\rm (3a)} Then, $V^{(\lambda)}(\Phi)$ and $V^{(\lambda)}(\Psi)$
have the same number of components, that is, $1+\zeta+
\eta=1+\zeta'+ \eta'$. Also, $\Phi(y,z)=0$ and $\Psi(y,z)=0$ have
the same multiplicity
$e_1=\zeta+\eta+n=\zeta'+\eta'+{\ell}=\bar{e}_1$ at $(y,z)=(0,0)$ by
construction. So, $n={\ell}$. \ms

{\rm (3b)} Then, $V^{(\lambda)}(\Phi)\cap E_1$ and
$V^{(\lambda)}(\Psi)\cap \bar{E}_1$ have the same number of elements
as set, that is, ${\zeta}={\zeta'}$ because $1\le n< \alpha_n$ and
$1\le {\ell}< \beta_{\ell}$. So, $\eta=\eta'$ by (3a). \ms

{\rm (3c)} Whenever there are two integers $n$ and $\alpha_{n}$ with
$1\le n<\alpha_{n}$, there is a positive integer $s$ such that
$sn<\alpha_{n}\le (s+1)n$. Also, whenever there are two integers
${\ell}$ and $\beta_{\ell}$ with $1\le {\ell}<\beta_{\ell}$, there
is a positive integer $\bar{s}$ such that
$\bar{s}{\ell}<\beta_{\ell}\le (\bar{s}+1){\ell}$. Note that
$e_{s+1}=\zeta+ (s+1)\eta+(s+1){\alpha_{n}}$ and
$e_{\bar{s}+1}={\zeta}'+
(\bar{s}+1){\eta}'+(\bar{s}+1){\beta_{\ell}}$ with $s={\bar{s}}$ are
equal, and so ${\alpha_{n}}={\beta_{\ell}}$. Therefore,
$V^{(s)}(\Phi)=V^{(s)}(\phi)$ and $V^{({s})}(\Psi)=V^{({s})}(\psi)$.
Then, $V^{(s)}(\phi)=V^{({s})}(\psi)$ implies that
$V^{(s)}(\Phi)=V^{({s})}(\Psi)$, and so $V^{(s)}(\Phi) \buildrel
\text{{\rm resol}} \over \sim V^{({s})}(\Psi)$ at the origin in
$\BC^2$ implies that $V^{(s)}(\phi) \buildrel \text{{\rm resol}}
\over \sim V^{({s})}(\psi)$ at the origin in $\BC^2$.

Thus, Definition $2.6$ is well-defined by {\rm(3a)}, {\rm(3b)} and
{\rm(3c)}.
\enddefinition \ms

\definition{Remark 2.6.1} A quasisingular point of $V^{(m)}(F)$ may not be a
singular point of  $V^{(m)}(F)$ by the following example:

Let $V(F)=\{(y,z):F(y,z)=z^2+y^3=0\}$. Then $(y,z)=(0,0)$ is a
singular point of $V(F)$. Let $\pi:M\to \BC^2$ be a blow-up at
$(0,0)$ such that $\pi(v,u)=(y,z)=(v,vu)$ and
$\pi(v',u')=(y,z)=(v'u',v')$ where $u'=\frac{1}{u}$ and $v'=vu$.
\roster
\item "(a1)" Along $v=0$, $(F\circ\pi)_{total}=v^2F_1(v,u)=v^2(u^2+v)$
where $F_1(v,u)$ is in $\BC\{v,u\}$.

\item "(a2)" Along $v'=0$,
$(F\circ\pi)_{total}={v'}^2 F_2(v',u')={v'}^2(1+{u'}^3v')$ where
$F_2(v',u')$ is a unit in $\BC\{v',u'\}$.
\endroster

Then, $V^{(1)}(F)$ has no singularity everywhere, but
$(v,u)=(0,0)$ is only one quasisingular point of $V^{(1)}(F)$
because additional blow-ups at $(v,u)=(0,0)$ are still necessary
to get the standard resolutions of the singular point of $V(F)$,
noting that $V^{(1)}(F)$ and $\{v=0\}$ meet tangentially at
$(v,u)=(0,0)$.
\enddefinition
\ms

\proclaim{Lemma 2.7} Assume that $V(F)$ satisfies the same
properties as in Definition $2.6$. If $P$ is a quasisingular point
of $V^{(m)}(F)$ under $\tau_m$, then $P$ satisfies at least one of
the following three properties {\rm(i)}, {\rm(ii)} and {\rm(iii)}:
\roster

\item "(i)" $V^{(m)}(F)$ itself has a singular point at $P$ as a
reduced variety.

\item "(ii)" There exists an exceptional curve of the first kind
of
 $\tau^{-1}_m(0,0)$ and at least one irreducible component of
 $V^{(m)}(F)$ which meet tangentially at $P$.

\item "(iii)" There are two exceptional curves of the first kind
of
 $\tau^{-1}_m(0,0)$ and at least one irreducible component of
 $V^{(m)}(F)$ which meet at $P$.
\endroster
\endproclaim

\demo{\bf Proof of Lemma 2.7} It is trivial.
\enddemo \ms

\definition{Definition 2.8} Let $V(F)$ be an analytic variety at
the origin. Assume that $V(F)$ satisfies the same properties and
notations as in Definition $2.6$.  It is said that $F\in$ the type
$[j]$ or belongs to the type $[j]$ under $\tau_m$ if $F$ satisfies
the following four properties:
 Sometimes, we say that $V(F)\in$ the type $[j]$ or belongs to the
type $[j]$ under $\tau_m$.

\roster

\item "(1)" For any point of the the set $\cup^{m-1}_{i=1}E_i$, no
blow-ups are needed to get the standard resolution of the
singular point of $V(F)$. That is, there are no quasisingular
points on $\tau^{-1}_{m-1}(0,0)=\cup^{m-1}_{i=1}E_i$.

\item "(2)" At some point of $E_m-\cup^{m-1}_{i=1}E_i$, additional
blow-ups may be needed for the standard resolution of the singular
point of $V(F)$.

\item "(3)" On the set $\cup^m_{i=1}E_i$, there are $j$ exceptional
curves of the first kind with $j\le m$, each of which has three
distinct intersection points with other exceptional curves and
the proper transform $V^{(m)}(F)$.

\item "(4)" Each of the $(m-j)$ remaining exceptional curves of the
first kind from $\cup^m_{i=1}E_i$, rather than the above $j$
exceptional curves as just mentioned in $(3)$, has at most two
distinct intersection points with other exceptional curves of the
first kind and the proper transform $V^{(m)}(F)$.
\endroster
\enddefinition
\ms

\definition{Remark 2.8.1} Assume  by Definition $2.8$ that $F\in$ the type
$[j]$ or belongs to the type $[j]$ under $\tau_m$. Whenever
$\cup^{m-1}_{i=1}E_i$ and $V^{(m)}(F)$ meet, then they meet with
normal crossings, but no three distinct components of
$\cup^{m-1}_{i=1}E_i$ and $V^{(m)}(F)$ meet.
\enddefinition
\ms

\definition{Definition 2.9} Let
$V(F)$ be an analytic variety at $(y,z)=(0,0)$ where $F=F(y,z)$ is
in $\BC\{y,z\}$ with isolated singularity at the origin. Assume that
$V(F)$ satisfies the same properties and notations as in Definition
$2.1$ and Definition $2.6$.

(1) Then, each of $U_1=(v_1,u_1)$ and $U_2=(v'_1,u'_1)$ is called
one coordinate patch of the given local coordinates for a blow-up
$\pi_1$, respectively. Also, $E_1=E^+_1\cup\ E^{-}_1$ is called an
exceptional curve of the first kind where
$E^+_1=\{(v_1,u_1):v_1=0\}$ and $E^{-}_1=\{(v'_1,u'_1):v'_1=0\}$.

(2) For notation, let $qs(V^{(1)}(F))$ be the set of all
quasisingular points of $V^{(1)}(F)$ under $\tau_1=\pi_1$ in the
sense of Definition $2.6$. Observe that $qs(V^{(1)}(F))$ may be
either empty or nonempty. Let $P\in qs(V^{(1)}(F))$ be chosen
arbitrary. For the given two local coordinates for $M^{(1)}$ as
above, it is said that the coordinate of $P$ is written by one and
only one of three different types :

\roster \item "(2a)" $\{(v_1,u_1)=(0,0)\}=E^+_1-E^{-}_1$.

\item "(2b)" $\{(v'_1,u'_1)=(0,0)\}=E^{-}_1-E^+_1$.

\item "(2c)" either $(v_1,u_1)=(0,\xi)\in {E^+_1}\cap{E^{-}_1}$ or
$(v'_1,u'_1)=(0,\frac {1}{\xi})\in {E^+_1}\cap{E^{-}_1}$ for a
nonzero number $\xi$.
\endroster

\roster \item "(3)" Now, let $P$ be a quasisingular point of
$V^{(1)}(F)$ in the sense of Definition $2.6$.

\item "(3a)" It is said that $P$ satisfies the first type of
coordinate under $\tau_1=\pi_1$ if $P$ satisfies either (2a) or (2b)
in (2). Equivalently, $P$ satisfies the first type of coordinate
under $\pi_1$ if and only if the coordinate of $P$ coincides with
the origin in one and only one of the given two coordinate patches
for $M$.

\item "(3b)" It is said that $P$ satisfies the second type of
coordinate under $\pi_1$ if $P$ satisfies (2c).
\endroster
\enddefinition

\proclaim{Lemma 2.10} Assume that $V(F)$ satisfies the same
properties and notations as in Definition $2.9$. Suppose that
$qs(V^{(1)}(F))$ is nonempty.

{\rm(1)} Then, $qs(V^{(1)}(F))$ satisfies one and only one of the
following four properties {\rm(1a)}, {\rm(1b)}, {\rm(1c)} and
{\rm(1d)}:

{\rm(1a)} $qs(V^{(1)}(F))\subset E^+_1$ and
$qs(V^{(1)}(F))\not\subset E^{-}_1$.

{\rm(1b)} $qs(V^{(1)}(F))\subset E^{-}_1$ and
$qs(V^{(1)}(F))\not\subset E^+_1$.

{\rm(1c)} $qs(V^{(1)}(F))\subset E^+_1\cap E^{-}_1$.

{\rm(1d)} $qs(V^{(1)}(F))\subset E_1$, but
$qs(V^{(1)}(F))\not\subset E^+_1$ and $qs(V^{(1)}(F))\not\subset
E^{-}_1$.\ms

{\rm(2)} By {\rm (1)}, there are four subcases:

{\rm(2a)} Suppose that $qs(V^{(1)}(F))\subset E^+_1$ and
$qs(V^{(1)}(F))\not\subset E^{-}_1$. Then, there is one and only
one point in $qs(V^{(1)}(F))$ satisfying the first type of
coordinate.

{\rm(2b)} Suppose that $qs(V^{(1)}(F))\subset E^-_1$ and
$qs(V^{(1)}(F))\not\subset E^+_1$. Then, there is one and only one
point in $qs(V^{(1)}(F))$ satisfying the first type of coordinate.

{\rm(2c)} Suppose that $qs(V^{(1)}(F))\subset E^+_1\cap E^{-}_1$.
Then, any $P\in qs(V^{(1)}(F))$ satisfies the second type of
coordinate.

{\rm(2d)} Suppose that $qs(V^{(1)}(F))\subset E_1$, but
$qs(V^{(1)}(F))\not\subset E^+_1$ and $qs(V^{(1)}(F))\not\subset
E^{-}_1$. Then, there are two distinct points in $qs(V^{(1)}(F))$,
each of which satisfies the first type of coordinate, respectively.
\endproclaim
The proof of lemma is clear. \ms

\definition{Remark 2.10.1}
Let $V(G)$ be an analytic variety at $(y,z)=(0,0)$ where $G=G(y,z)$
is in $\BC\{y,z\}$ with isolated singularity at the origin. As we
have seen in Definition $2.1$ and Definition $2.4$, let
$\tau_m=\pi_1\circ\pi_2\circ\cdots \circ\pi_m:M^{(m)}\to \BC^2$ be
the composition of a finite number $m$ of successive blow-ups at the
origin in $\BC^2$, which is needed only to get the standard
resolution of the singular point $(0,0)$ of $V(G)$.

(1) First, let $G=G(y,z)$ be irreducible in $\BC\{y,z\}$ with
isolated singularity at the origin. Note that $qs(V^{(t)}(G))$ is
a one-point set for each $t=1,2,\dots,(m-1)$ and $qs(V^{(m)}(G))$
is empty, because $\tau_m$ is the standard resolution of the
singular point of $V(G)$ by assumption. In order to either study
$V^{(t)}(G)$ or find $qs(V^{(t)}(G))$ under $\tau_t$, it is
possible to choose just one coordinate patch of the local
coordinates for each blow-up $\pi_t:M^{(t)}\to M^{(t-1)}$, where
$1\le t\le m$ and $M^{(0)}=\BC^2$, by Definition $2.9$ and Lemma
$2.10$, as it has been well-known in the standard resolution of
irreducible plane curve singularities.

In order to study $\tau_{s}:M^{(s)}\to\BC^2$ with $1\le s\le m$,
it is clear that $\tau_{s}$ can be viewed as a composition of a
finite number $s$ of successive analytic mappings from one
coordinate patch of $M^{(t)}$ to another coordinate patch of
$M^{(t-1)}$ in the sense of Definition $2.9$, where $1\le t\le s$
and $M^{(0)}=\BC^2$. \ms

{\rm(2)} Next, without assuming that $G(y,z)$ is irreducible in
$\BC\{y,z\}$, for $1\le t\le m$ let
$\tau_t=\pi_1\circ\pi_2\circ\cdots\circ\pi_t:M^{(t)}\to\BC^2$ be
defined by the composition of a finite number $t$ of successive
blow-ups at the origin in $\BC^2$, which is in process of the
standard resolution of the singular point $(0,0)$ of $V(G)$.

Assuming that $qs(V^{(t)}(G))$ is a one-point set for each
$t=1,2,\dots,s-1$ if exists where $s$ is a positive integer with
$2\le s\le m$, then by the similar way as in (1) it is easy to
prove that $\tau_{t}:M^{(t)}\to\BC^2$ can be viewed as a
composition of a finite number $t$ of successive analytic mappings
from one coordinate patch of $M^{(j)}$ to another coordinate patch
of $M^{(j-1)}$ in the sense of Definition $2.9$, where $1\le j\le
t$ and $M^{(0)}=\BC^2$.

Now, whether or not $qs(V^{(s)}(G))$ is a one-point set, by using
the definition of the quasisingularity and also following the
definitions and lemmas in this section and the next section, we
are going to generalize the possibility of the terminology about
the choices of just one coordinate patch of the local coordinates
for each blow-up $\pi_t$, in order to study the local defining
equation for the total transform of $V(G)$ under $\tau_s$.
\enddefinition
\ms

\definition{Definition 2.11} Let $V(F)$ and $V(G)$ be analytic varieties
at $(y,z)=(0,0)$ where $F=F(y,z)$ and $G=G(y,z)$ are in
$\BC\{y,z\}$. Let $V(F)$ and $V(G)$ have an isolated singular point
at the origin as reduced varieties. Assuming that $V(F)$ and $V(G)$
satisfy the same kind of properties and notations as in Definition
$2.9$ and Lemma $2.10$, let $\pi_1:M^{(1)}\to \BC^2$ be the first
blow-up of $\BC^2$ at $(y,z)=(0,0)$ in process of the standard
resolution of the isolated singular point $(0,0)$ of both $V(F)$ and
$V(G)$. Let $U_1=(v_1,u_1)$ and $U_2=(v'_1,u'_1)$ be two coordinate
patches with $u'_1=1/u_1$ and $v'_1=v_1u_1$ for the blow-up $\pi_1$
as before. Let $qs(V^{(1)}(F))$ and $qs(V^{(1)}(G))$ be the set of
quasisingular points of $V^{(1)}(F)$ and $V^{(1)}(G)$, respectively.
Note by Definition $2.9$ that each of $U_1=(v_1,u_1)$ and
$U_2=(v'_1,u'_1)$ is called one coordinate patch of the given local
coordinates for a blow-up $\pi_1$, respectively. As we have seen in
Definition $2.9$, each of $U_1=(v_1,u_1)$ and $U_2=(v'_1,u'_1)$ is
called one coordinate patch of the given local coordinates for a
blow-up $\pi_1$, respectively. Also, $E_1=E^+_1\cup\ E^{-}_1$ is
called an exceptional curve of the first kind where
$E^+_1=\{(v_1,u_1):v_1=0\}$ and $E^{-}_1=\{(v'_1,u'_1):v'_1=0\}$.

Then, we have the following new terminologies [I] and [II]:

{\bf [I]} If $qs(V^{(1)}(F))\not= \emptyset$ satisfies either
$qs(V^{(1)}(F))\subset E^+_1$ or $qs(V^{(1)}(F))\subset E^{-}_1$,
then we say that just one coordinate patch covering is needed for
the blow-up $\pi_1$, in order to either study $V^{(1)}(F)$ or find
$qs(V^{(1)}(F))$, including the case that
$qs(V^{(1)}(F))=\emptyset$, whether or not $qs(V^{(1)}(F))$ is
nonempty.

In more detail, we say that if $qs(V^{(1)}(F))$ satisfies one and
only one of three properties $(1a)$, $(1b)$ and $(1c)$ except for
$(1d)$ in Lemma $2.10$ then it has one coordinate patch covering.
Equivalently, depending on $(1a)$, $(1b)$ and $(1c)$ in Lemma
$2.10$, then [I] can be divided into the following three subcases
{\rm[Ia]}, {\rm[Ib]} and {\rm[Ic]}:

\roster
\item "{\bf[Ia]}" If $qs(V^{(1)}(F))\not= \emptyset$
satisfies the property{\rm (1a)} in Lemma $2.10$, then it is said
by definition that just one coordinate patch covering $U_1$ of
$E^+_1$ is needed for the study of $V^{(1)}(F)$ under $\pi_1$,
because the other coordinate patch covering $U_2$ of $E^{-}_1$ is
not needed to find $qs(V^{(1)}(F))$, as a consequence.

\item "{\bf[Ib]}" If $qs(V^{(1)}(F))\not= \emptyset$ satisfies
{\rm (1b)} in Lemma $2.10$, then it is said by definition that
just one coordinate patch covering $U_2$ of $E^{-}_1$ is needed
for the study of $V^{(1)}(F)$ under $\pi_1$, because the other
coordinate patch covering of $U_1$ of $E^+_1$ is not needed to
find $qs(V^{(1)}(F))$, as a consequence.

\item "{\bf[Ic]}" If $qs(V^{(1)}(F))\not= \emptyset$ satisfies
{\rm (1c)} in Lemma $2.10$, then it is said by definition that
either one coordinate patch covering $U_1$ of $E^+_1$ or another
coordinate patch covering $U_2$ of $E^{-}_1$ may be chosen
arbitrary for the study of $V^{(1)}(F)$ under $\pi_1$, because any
one of two coordinate patches is needed enough to find
$qs(V^{(1)}(F))$, as a consequence.
\endroster
For example, if $qs(V^{(1)}(F))$ is a one-point set, then it is
clear that we can use just one coordinate patch covering of the
given local coordinates for the blow-up $\pi_1$, in order to
either study $V^{(1)}(F)$ or find $qs(V^{(1)}(F))$. \ms

{\bf[Id]} If $qs(V^{(1)}(F))$ satisfies (1d) in Lemma $2.10$, then
it is said that two coordinate patches $U_1$ and $U_2$ are all
needed without using any nonsingular change of coordinates, in order
to either study $V^{(1)}(F)$ or find $qs(V^{(1)}(F))$.

In other words, we say that if $qs(V^{(1)}(F))$ satisfies (1d) in
Lemma $2.10$ then it must have two coordinate patch coverings, but
it does not have one coordinate patch covering under the given
blow-up $\pi_1$, in order to study $V^{(1)}(F)$. \ms

{\bf [II]} It is said that we may use a common one coordinates patch
of the given local coordinates for the blow-up $\pi_1$, in order to
study both $V^{(1)}(F)$ and $V^{(1)}(G)$ simultaneously, if one of
the following properties holds :

(a) \quad $qs(V^{(1)}(F))\subset E^+_1$ and $qs(V^{(1)}(G))\subset
E^+_1$.

(b) \quad $qs(V^{(1)}(F))\subset E^{-}_1$ and $qs(V^{(1)}(G))\subset
E^{-}_1$. \ms

Equivalently, it is said that $qs(V^{(1)}(F))$ and $qs(V^{(1)}(G))$
may have the same one coordinate patch covering. Observe that either
$qs(V^{(1)}(F))$ or $qs(V^{(1)}(G))$ may be empty.

\enddefinition
\ms

\proclaim{Lemma 2.12} $\underline{\text{\bf Assumptions}}$ Let
$V(F)$ be an analytic variety at $(y,z)=(0,0)$ where $F=F(y,z)$ is
in $\BC\{y,z\}$. Assume that $V(F)$ has an isolated singular point
at the origin as a reduced variety. Let $\tau_{\ell}:M^{(\ell)}\to
\BC^2$ be the composition of the first finite number $\ell$ of
successive blow-ups at the origin which is in process of the
standard resolution of the singular point of $V(F)$. For each
$t=1,2,\dots,\ell$, write
$\tau_t=\pi_1\circ\pi_2\circ\cdots\circ\pi_t:M^{(t)}\to \BC^2$ where
$\{\pi_i:M^{(i)}\to M^{(i-1)}$ is a blow-up of $M^{(i-1)}$ at some
point of $M^{(i-1)}$ for $1\le i\le \ell$ with $M^{(0)}=\BC^2\}$.

Assume that $V(F)$ satisfies the following properties:

{\rm(i)} $F=F(y,z)$ may not be irreducible in $\BC\{y,z\}$.

{\rm(ii)} $qs(V^{(t)}(F))$ is a one-point set for each
$t=1,2,\dots,\ell-1$, which is denoted by $P_t$.

{\rm(iii)} $qs(V^{(\ell)}(F))$ may be empty. If not empty, then it
has a one coordinate patch covering in the sense of Definition
$2.11$.

$\underline{\text{\bf Conclusions}}$ Each $\tau_t$ with $1\le t\le
\ell$ can be also viewed as the composition of the same number of
analytic mappings from a two-dimensional complex analytic polydisc
centered at the origin in $\BC^2$ to itself, in order to either
study $V^{(t)}(F)$ or find $(F\circ\tau_t)_{proper}$.
\endproclaim

The Proof of Lemma 2.12 is clear. \ms

\definition{Definition 2.13}
Let $V(F)$ and $V(G)$ be analytic varieties at $(y,z)=(0,0)$ where
$F=F(y,z)$ and $G=G(y,z)$ are in $\BC\{y,z\}$. Let $V(F)$ and $V(G)$
have an isolated singular point at the origin as reduced varieties.
Assuming that $V(F)$ and $V(G)$ satisfy the same kind of properties
and notations as in Definition $2.9$ and Definition $2.11$, let
$\pi_1:M^{(1)}\to \BC^2$ be the first blow-up of $\BC^2$ at
$(y,z)=(0,0)$ in process of the standard resolution of the isolated
singular point $(0,0)$ of both $V(F)$ and $V(G)$. Let
$U_1=(v_1,u_1)$ and $U_2=(v'_1,u'_1)$ be two coordinate patches with
$u'_1=1/u_1$ and $v'_1=v_1u_1$ for the blow-up $\pi_1$ as before.
Let $qs(V^{(1)}(F))$ and $qs(V^{(1)}(G))$ be the set of
quasisingular points of $V^{(1)}(F)$ and $V^{(1)}(G)$, respectively.
Note by Definition $2.9$ that each of $U_1=(v_1,u_1)$ and
$U_2=(v'_1,u'_1)$ is called one coordinate patch of the given local
coordinates for a blow-up $\pi_1$, respectively.

It is said that we may use a common one coordinates patch of the
given local coordinates for the blow-up $\pi_1$, in order to study
both $V^{(1)}(F)$ and $V^{(1)}(G)$ simultaneously, if one of the
following properties holds : Recall that
$E^{+}_1=\{(v_1,u_1):v_1=0\}$ and
$E^{-}_1=\{(v'_1,u'_1):v'_1=0\}$.

(a) \quad $qs(V^{(1)}(F))\subset E^+_1$ and
$qs(V^{(1)}(G))\subset E^+_1$.

(b) \quad $qs(V^{(1)}(F))\subset E^{-}_1$ and
$qs(V^{(1)}(G))\subset E^{-}_1$. \ms

Equivalently, it is said that $qs(V^{(1)}(F))$ and
$qs(V^{(1)}(G))$ may have the same one coordinate patch covering.
Observe that either $qs(V^{(1)}(F))$ or $qs(V^{(1)}(G))$ may be
empty.
\enddefinition
\ms

\definition{Remark 2.13.1}

{\rm(i)} In particular, if $qs(V^{(1)}(F))=qs(V^{(1)}(G))$ is a
one-point set, then we may use a common one coordinate patch of
the local coordinates under the blow-up $\pi_1$, in order to study
both $V^{(1)}(F)$ and $V^{(1)}(G)$ at the same time. \ms

{\rm(ii)} Let $V(F)=\{(y,z):F(y,z)=0\}$, $V(G)=\{(y,z):G(y,z)=0\}$
and $V(H)=\{(y,z):H(y,z)=0\}$ be analytic varieties defined
respectively, as follows:
$$\align
&F(y,z)=z^2+y^3, \tag 2.13.1 \\
&G(y,z)=(z^2+y^2)^2+y^5, \\
&H(y,z)=z^3+y^2.
\endalign$$

Let $\pi_1:M^{(1)}\to \BC^2$ be the first blow-up of $\BC^2$ at
$(y,z)=(0,0)$ which is an isolated singular point of $V(F)$,
$V(G)$ and $V(H)$. Let $(v, u)$ and $(v', u')$ be the given local
coordinates for $M^{(1)}$ where $\pi_1(v,u)=(y,z)=(v,vu)$ and
$\pi_1(v',u')=(y,z)=(v'u',v')$ with $u'=1/u$ and $v'=vu$.

Note that $E_1=E^+_1\cup E^{-}_1$ where $E^+_1=\{(v, u):v=0\}$ and
$E^{-}_1=\{(v', u'):v'=0\}$. For notation of two coordinate
patches, we write $U_1=(v,u)$ and $U_2=(v',u')$.

Observe that $qs(V^{(1)}(F))\not=qs(V^{(1)}(G))$ and
$qs(V^{(1)}(H))\not=qs(V^{(1)}(G))$, because
$qs(V^{(1)}(F))\subset E^+_{1}-E^{-}_{1}$, $qs(V^{(1)}(G))\subset
E^+_{1}\cap E^{-}_{1}$ and $qs(V^{(1)}(H))\subset
E^{-}_{1}-E^+_{1}$.

For the study of $V^{(1)}(G)$, it is enough to choose either one
of two coordinate patches $U_1=(v,u)$ and $U_2=(v',u')$ under the
blow-up $\pi_1$. So, in order to study both $V^{(1)}(F)$ and
$V^{(1)}(G)$ simultaneously, it is possible to use one and only
one coordinate $U_1=(v,u)$ as a common one coordinate patch of the
local coordinates under the blow-up $\pi_1$. Also, in order to
study both $V^{(1)}(H)$ and $V^{(1)}(G)$ simultaneously, it is
possible to use one and only one coordinate $U_2=(v',u')$ as a
common one coordinate patch of the local coordinates under the
blow-up $\pi_1$. But, in order to find $V^{(1)}(F)$, $V^{(1)}(G)$
and $V^{(1)}(H)$ at the same time, we cannot choose a common one
coordinate patch of the given two local coordinates under the
blow-up $\pi_1$ without using a nonsingular change of coordinates.

Define a local nonsingular mapping $\phi$ from $(y,z)=(0,0)$ to
$(y',z')=(0,0)$ as follows:
$$\align
\phi(y,z)&=(y',z') \quad \text{with} \quad \phi(0,0)=(0,0), \tag 2.13.2 \\
y'&=y \quad \text{and} \quad z'=z+y. \\
\endalign$$

From (2.13.1) and (2.13.2), let $F^*=F\circ\phi$, $G^*=G\circ\phi$
and $H^*=H\circ\phi$, respectively. Then, $V(F^*)$, $V(G^*)$ and
$V(H^*)$ have still the same isolated singular point at
$(y,z)=(0,0)$, as follows:
$$\align
&F^*(y,z)=(z+y)^2+y^3, \tag 2.13.3 \\
&G^*(y,z)=\{(z+y)^2+y^2)\}^2+y^5, \\
&H^*(y,z)=(z+y)^3+y^2.
\endalign$$

Let $\pi_1:M^{(1)}\to \BC^2$ be the first blow-up of $\BC^2$ at
$(y,z)=(0,0)$ with the same local coordinates $(v,u)$ and
$(v',u')$ for $M^{(1)}$ just as above. Then, $F^*$, $G^*$ and
$H^*$ have a common one coordinate patch of the given two local
coordinates under the blow-up $\pi_1$, in order to find
$V^{(1)}(F^*)$, $V^{(1)}(G^*)$ and $V^{(1)}(H^*)$ at the same
time.
\enddefinition \ms

\proclaim{Lemma 2.14}  $\underline{\text{\bf Assumptions}}$ As in
Definition $2.1$, let $V(F)$ and $V(G)$ be analytic varieties at
$(y,z)=(0,0)$ where $F=F(y,z)$ and $G=G(y,z)$ are in $\BC\{y,z\}$.
Assume that $V(F)$ and $V(G)$ have an isolated singular point at the
origin in a two-dimensional manifold $\BC^2$, as reduced varieties.
Suppose that $V(F)$ satisfies the same assumptions and notations as
in Lemma $2.12$. Let $\omega_q:L^{(q)}\to \BC^2$ be an arbitrary
composition of a finite number $q$ of successive blow-ups at the
origin which is in process of the standard resolution of the
singular point of $V(G)$. For each $r=1,2,\dots,q$, write
$\omega_q=\sigma_1\circ\sigma_2\circ\cdots\circ\sigma_q:L^{(q)} \to
\BC^2$ with $\pi_1=\sigma_1$, where $\{\sigma_r:L^{(r)}\to
L^{(r-1)}$ is a blow-up of $L^{(r-1)}$ at some point of $L^{(r-1)}$
for $1\le r\le q \}$ with $L^{(0)}=\BC^2$ and $L^{(1)}=M^{(1)}$.

As we have seen in the assumption of Lemma $2.12$, assume that
$V(G)$ satisfies the following properties: \roster

\item "(i)" $G=G(y,z)$ may not be irreducible in $\BC\{y,z\}$.

\item "(ii)" $qs(V^{(r)}(G))$ is a one-point set for each
$r=1,2,\dots,q-1$, which is denoted by $Q_r$.

\item "(iii)" $qs(V^{(q)}(G))$ may be empty. If not empty, then it
has one coordinate patch covering in the sense of Definition
$2.11$.
\endroster

Define $s=\min\{\ell,q\}$, and then we assume the following
property:
$$\align
(2.14.1) \qquad \qquad &\text{There exists a positive integer~ $p$
~with~ $1\le p \le s$ ~~such that } \qquad \qquad \qquad \qquad  \\
&\text{\quad
$P_i=Q_i$\quad for each ~$i=1,2,\dots,p-1$.}
\endalign$$

$\underline{\text{\bf Conclusions}}$ We have the following: \roster

\item "(a)" For each $i=1,2,\dots,s$, it is possible to choose the
same local coordinates for both $\pi_i$ and $\sigma_i$
simultaneously, without any nonsingular change of coordinates.

\item "(b)" We may use a common one coordinate patch of the given
local coordinates for each blow-up $\pi_p=\sigma_p$, in order to
study both $V^{(p)}(F)$ and $V^{(p)}(G)$ simultaneously for each
$p=1,2,\dots,s$. But, a common one coordinate patch may not be
unique.

\item "(c)" In more detail, for each $p=1,2,\dots,s$,
$\tau_p=\omega_p$ can be considered as a local analytic mapping
from a two-dimensional complex analytic manifold to an analytic
polydisc centered at the origin in $\BC^2$, in order either to
study both $V^{(s)}(F)$ and $V^{(s)}(G)$ simultaneously, or to
find both $(F\circ\tau_s)_{proper}$ and
$(G\circ\omega_s)_{proper}$ at the same time.
\endroster
\endproclaim
\demo{\bf Proof of Lemma 2.14} If $qs(V^{(1)}(F))=qs(V^{(1)}(G))$
is a one-point set, then we may use a common one coordinate patch
of the local coordinates under the blow-up $\pi_1$ for the study
of both $V^{(1)}(F)$ and $V^{(1)}(G)$ at the same time. Then, the
proof can be just finished by induction on the positive integer
$s$ where $s=\min\{\ell,q\}$, by using Lemma $2.12$ and Definition
$2.13$ only.
\enddemo \ms

\vfill \pagebreak

{\bf \S3.  The representation for the local defining equations of
irreducible plane curve singularities which have either the same
multiplicity sequence or the homeomorphic resolution under the
standard resolution as the curve singularity
\text{(\{$z^n+y^k=0$\})} does and its generalizations } \bs

{\bf \S3.0. Introduction } \ms

In $\S 3.1$, by Theorem $3.2$ and Corollary $3.3$, we will find the
method how to construct the local defining equation for irreducible
plane curve singularities, which have the homeomorphic resolution as
an irreducible plane curve singularity defined by $z^{n}+y^{k}=0$
does in the sense of Definition 2.4 and its generalizations.

In $\S 3.2$, by Theorem $3.6$ and and Theorem $3.7$, we will find
the method how to compute the local defining equation for the total
transform of irreducible plane curve singularities, which have the
homeomorphic resolution as an irreducible plane curve singularity
defined by $z^{n}+y^{k}=0$ does, after a finite number of successive
blow-ups, which is just needed to get the standard resolution of the
above singularity and its generalizations. \ms

{\bf  {\S3.1.} How to find the necessary condition for any local
defining equation defining plane curve singularities to be
irreducible in $\BC\{y,z\}$ and its applications}

\proclaim{Lemma 3.1(Hensel's Lemma)} Let
$f(y,z)=a_0z^n+a_1y^{\alpha_1}z^{n-1}+ \dots+a_ny^{\alpha_n}$ be
irreducible in $\BC\{y,z\}$, where $n\ge 2$ and $\alpha_n\ge 2$, and
each $a_i$ is a unit in $\BC\{y,z\}$ if exists, and the $\alpha_i$
are positive integers. Let $m$ be the multiplicity of $f$ at $0$.
Then $m=n$ or $\alpha_n$. If $m=n=\alpha_i+n-i$ for some $i$, then
$n=\alpha_i+n-i$ for all $i=1,\dots,n$, and so $f$ can be written in
a power series as follows: \text{$f=f_n(y,z)+$ terms of degree
$>n$}, where $f_n=f_n(y,z)$ is a homogeneous polynomial of degree
$n$ with $f_n=(ay+bz)^n$ for some  ${a, b \in \BC}$.
\endproclaim \ms

\proclaim{Theorem 3.2(The generalized Hensel's lemma)}

$\underline{\text{\bf Assumptions}}$ Let $f(y,z)$ be in $\BC\{y,z\}$
with isolated singularity at $(0,0)$, defined by the following:
$$\align
(3.2.1) \qquad  \text{$f(y,z)=a_0z^n+a_1y^{\alpha_1}z^{n-1} +
\dots+a_ny^{\alpha_n}$ \quad with $n\ge 2$ and $k=\alpha_n\ge 2$,}
\qquad \qquad
\endalign$$
where each $a_i=a_i(y,z)$ is a unit in $\BC\{y,z\}$ if exists, and
the $\alpha_i$ are positive integers for all $i=1,2,\dots,n$, and
also $a_0=a_0(y,z)$ and $a_n=a_{n}(y,z)$ are units in $\BC\{y,z\}$.

If necessary, we may assume without need of the proof that $1\le
\alpha_1<\alpha_2<\cdots <\alpha_{n}$. \ms

$\underline{\text{\bf Conclusions}}$

{\rm(1)} The necessary condition for $f(y,z)$ to be irreducible in
$\BC\{y,z\}$ is as follows:
$$\align
\text{${\dfrac{\alpha_i}{i}} \ge {\dfrac{k}{n}}$ \quad for all
$i=1,2,\dots,n-1$.} \tag 3.2.2
\endalign$$

{\rm(2)} Equivalently, the inequality in {\rm(3.2.2)} holds if and
only if $f$ is represented as follows:
$$\align
(3.2.3) \qquad \qquad \qquad
\text{$f(y,z)=a_{0,0}z^n+a_{n,0}y^{k}+\sum_{\alpha,\beta\ge
0}c_{\alpha,\beta}y^{\alpha}z^{\beta}$} \quad \text{with
$n\alpha+k\beta\ge nk$,} \qquad \qquad
\endalign$$
where the $c_{\alpha,\beta}$ are nonzero complex numbers for some
nonnegative integers $\alpha$ and $\beta$ with $n\alpha+k\beta\ge
nk$, satisfying {\rm(i)} and {\rm(ii)}:

{\rm(i)} $a_0$ and $a_n$ are units in $\BC\{y,z\}$ with
$a_{0,0}=a_0(0,0)$ and $a_{n,0}=a_{n}(0,0)$.

{\rm(ii)} if $c_{0,\beta}\not=0$ for some integers $\beta>0$, then
$\beta>n$, and

\qquad if $c_{\alpha,0}\not=0$ for some integers $\alpha>0$, then
$\alpha>k$. \ms

{\rm(3)} In particular, as far as the inequality in $(3.2.2)$ is
concerned, if $k=np$ for some positive integer $p$,
$\dfrac{k}{n}=\dfrac{\alpha_i}{i}=p$ for $i=1,2,\dots,n-1$ where all
the $a_i$ are units in $\BC\{y,z\}$.
\endproclaim \ms

\demo{\bf Proof of Theorem 3.2} Assume that $f(y,z)$ is irreducible
in $\BC\{y,z\}$ with isolated singularity at $(0,0)$. The proof will
be by induction on the multiplicity of $f$. Then, it suffices to
prove that either the inequality in (3.2.2) or the representation in
(3.2.3) is true. For the induction proof, assuming that $2\le n\le
k$, it suffices to consider two facts, respectively:

$\underline{\text{\bf Fact(1)}}$ \quad Let $n=k\ge 2$. If $f(y,z)$
is irreducible in $\BC\{y,z\}$, by Lemma $3.1$ and by (3.2.1),
$\alpha_i+(n-i)=n$ for all $i=1,2,\dots,n-1$, and then $\alpha_i=i$.
So, there is nothing to prove for the inequality in (3.2.2). \ms

$\underline{\text{\bf Fact(2)}}$ To prove that the inequality in
(3.2.2) is true, by using Fact(1) it suffices to consider two cases
by induction on the multiplicity $n$ of the local defining equation
of $f(y,z)$ at the origin.

Case(I) $n=2<k$ and Case(II) $2<n<k$. \ms

$\underline{\text{\rm Case(I)}}$ Let $n=2<k$. The proof is clear.
\ms

$\underline{\text{\rm Case(II)}}$ Let $2<n<k$. Suppose we have shown
that the theorem is true whenever the multiplicity of the local
defining equation of $f(y,z)$ at the origin is less than $n$. Since
$2<n<k$, then there is a positive integer $s$ such that $0<k-sn\le
n$. For the induction proof, first of all, let $\tau_m$ be the
composition of a finite number $m$ of successive blow-ups which is
needed only to get the standard resolution of the singularity of
$V(f)$. Then, it will be shown by Sublemma $3.2.1$ that $s<m$, and
also if $V^{(s)}(f)$ is the proper transform under $\tau_s$ where
$\tau_s$ is the composition of a finite number $s$ of successive
blow-ups at $(y,z)=(0,0)$, then the local defining equation for
$V^{(s)}(f)$ is of the multiplicity $k-sn$ where $0<k-sn\le n$.
Then, there are two subcases {\rm(i)} and {\rm(ii)}.

{\rm(i)} If $0<k-sn<n$, apply the induction method to the
multiplicity of the local defining equation for $V^{(s)}(f)$ in
(3.2.7) of Sublemma $3.2.1$. After then, these will give the desired
proof of Case(II) by Sublemma $3.2.1$. \ms

{\rm(ii)} If $k-sn=n$, then the proof will be done by the same
method as we have used in the proof of Fact(1). \ms

\proclaim{Sublemma 3.2.1} $\underline{\text{\bf Assumptions}}$
Suppose that $f\in \BC\{y,z\}$ satisfies the same properties and
notations as in the assumption of Theorem $3.2$. In addition, let
$2\le n<k$. Then, there is a unique integer $s$ such that $sn<k\le
(s+1)n$.

$\underline{\text{\bf Conclusions}}$ \quad Let $\tau_m$ be the
composition of a finite number $m$ of successive blow-ups which is
needed only to get the standard resolution of the singular point of
$V(f)$. For each $t=1,2,\dots,m$, write $\tau_t=
\pi_1\circ\pi_2\circ\cdots \circ\pi_t:M^{(t)}\to \BC^2$ where
$\{\text{$ \pi_i:M^{(i)}\to M^{(i-1)}$ is a blow-up}$

\noindent $\text{of $M^{(i-1)}$ at some point of $M^{(i-1)}$ for
$1\le i\le t$} \}$ with $M^{(0)}=\BC^2$. Let $V^{(t)}(f)$ be the
proper transform under $\tau_t$ for $1\le t\le m$. Let
$E^{(m)}=\tau^{-1}_m(0,0)$, and let $E^{(m)}=\cup E_i$, $1 \le i \le
m$, be the decomposition of $E^{(m)}$ into irreducible components
where each $E_i$ is called an exceptional curve of the first kind.

As a conclusion, $1\le s<m$, and in order to study the proper
transform $V^{(t)}(f)$ for each $t=1,2,\dots,s$, we can use one and
only one coordinate patch of the local coordinates for each blow-up
$\pi_i:M^{(i)}\to M^{(i-1)}$ of $\tau_t$ with $1\le i\le t$,
satisfying the following properties:

\roster \item  "(a)" For each $t=1,2,\dots,s$, $qs(V^{(t)}(f))$ is a
one-point set in the sense of Definition $2.6$.

\item "(b)"  For each $i=1,2,\dots,t$, let $(v_i,u_i)$ and
$(v'_i,u'_i)$ be the local coordinates for $M^{(i)}$ where
$\pi_i:M^{(i)}\to M^{(i-1)}$ is a blow-up of $M^{(i-1)}$ at some
point of $M^{(i-1)}$ for $1\le i\le t$ where $u'_i=1/{u_i}$ and
$v'_i=v_iu_i$ and $M^{(0)}=\BC^2$ with $(v_0,u_0)=(y,z)$. We write
$E_i=\{(v_i,u_i):v_i=0\}\cup\{(v'_i,u'_i):v'_i=0\}$ for each
$i=1,2,\dots,t$.

\item "(b1)" Let $\pi_t(v_t,u_t)=(v_{t-1},u_{t-1})=(v_t,v_tu_t)$
and $\pi_t(v'_t,u'_t)=(v_{t-1},u_{t-1})=(v'_tu'_t,v'_t)$ for
$i=1,2,\dots,t$ where $(v_0,u_0)=(y,z)$. Then, $(v'_{t},u'_{t})$ is
not needed for the study of the proper transform $V^{(t)}(f)$ under
$\pi_{t}$ by $(3.2.6)$, and so $qs(V^{(t)}(f))\cap
E_t=\{(v_t,u_t)=(0,0)\}$, whose proof just follows from $(3.2.6)$
and $(3.2.7)$, below. \ms
\endroster

{\rm(c)} For $1\le t\le s$, let
$\tau_t=\pi_1\circ\pi_2\circ\cdots\circ\pi_t:M^{(t)}\to\BC^2$ be
defined by the local coordinates in $(b)$ where
$E_t=\{v_t=0\}\cup\{v'_t=0\}$.

In order to study the proper transform $V^{(t)}(f)$ for each
$t=1,2,\dots,s$, we may assume that $\tau_{t}:M^{(t)}\to\BC^2$
satisfies the same assumptions and notations as in $(b)$. For each
fixed $t=1,2,\dots,s$, $\tau_{t}$ can be rewritten in the form
$$\align
\noindent (3.2.4) \quad \quad \noindent
& \tau_t=\pi_1\circ\pi_2\circ\cdots\circ\pi_t,  \\
 & M^{(t)}@> \pi_{t} >> M^{(t-1)}@>
 \pi_{t-1} >> M^{(t-2)}\rightarrow\cdots \rightarrow
 M^{(1)}@> \pi_1 >> M^{(0)}=\BC^2, \\
& (v_t,u_t)@> \pi_{t} >>
 (v_{t-1},u_{t-1})@> \pi_{t-1} >>
 (v_{t-2},u_{t-2}) \rightarrow\cdots \rightarrow
 (v_1,u_1)@>\pi_1 >>
 (v_{0},u_{0}), \\
&(v'_t,u'_t)@> \pi_{t} >>
 (v_{t-1},u_{t-1})@> \pi_{t-1} >>
 (v_{t-2},u_{t-2}) \rightarrow\cdots \rightarrow
 (v_1,u_1)@>\pi_1 >> (v_{0},u_{0}), \qquad \qquad \\
 &  \text{where $(v_0,u_0)=(y,z)$}.
\endalign$$

By {\rm(3.2.4)}, along $E_t$ \ $\tau_{t}:M^{(t)}\to\BC^2$ can be
represented, as a composition of analytic mappings, as follows:
$$
\align
 \tau_t(v_t,u_t) &= (y,z)=(v_t,{v_t}^t u_t), \tag 3.2.5 \\
 \tau_t(v'_t,u'_t) &=(y,z)=(v'_t u'_t,{v'_t}^t {u'_t}^{t-1}),
\endalign
$$
where $u'_t={1/u_t}$ and $v'_t=v_tu_t$.

For each fixed $t=1,2,\dots,s$, along $v'_{t}=0$,
$(f\circ\tau_{t})_{total}$ can be  written in the form
$$\align
(3.2.6) \quad  (f\circ\tau_{t})_{total}
&=f(v'_{t}u'_{t},{v'_{t}}^{t}{u'_{t}}^{t-1})
={v'_{t}}^{tn}{u'_{t}}^{(t-1)n}f_{t}(v'_{t},u'_{t}), \\
f_{t}(v'_{t},u'_{t})&=a_{0,0}+ a_{\ell+1,0}
{v'_{t}}^{k-tn}{u'_{t}}^{k-(t-1)n} +\sum_{\alpha,\beta\ge
0}c_{\alpha,\beta}{v'_{t}}^{\alpha+t\beta-tn}
{u'_{t}}^{\alpha+(t-1)\beta-(t-1)n} \\
 \text{with}  & \quad  {\alpha+t\beta-tn>0.}
\endalign$$

For each $t=1,2,\dots,s$, $\{(v'_i,u'_i): i=1,2,\dots,t\}$ is not
needed for the study of the proper transform $V^{(t)}(f)$ under
$\tau_{t}$, because $\alpha+t\beta-tn>0$ by $(3.2.6)$ and $k-tn>0$
imply that $f_{t}(v'_t,u'_t)= \varepsilon'_t$ where $\varepsilon'_t$
is defined to be a unit in $\BC\{v'_t,u'_t\}$. \ms

For each fixed $t=1,2,\dots,s$, along $v_{t}=0$,
$(f\circ\tau_{t})_{total}$ can be written in the form
$$\align
(3.2.7) \qquad \qquad (f\circ\tau_{t})_{total}
&=f(v_{t},v_{t}^{t}u_{t})
={v_{t}}^{tn}f_{t}(v_{t},u_{t}),   \\
f_{t}(v_{t},u_{t}) &=a_{0,0}u_{t}^{n}+a_{\ell+1,0} {v_{t}}^{k-tn}
+\sum_{\alpha,\beta\ge
0}c_{\alpha,\beta}v_{t}^{\alpha+t\beta-tn}u_{t}^{\beta}, \qquad \qquad \qquad \\
\text{with}   \qquad & \alpha+t\beta-tn>0.
\endalign$$

So, $qs(V^{(t)}(f))\cap E_t=\{(v_t,u_t)=(0,0)\}$, and also the
multiplicity of the local defining equation for the proper
transform $V^{(s)}(f)$ under $\tau_s$ is $k-ns\le n$ at
$(v_s,u_s)=(0,0)$. \ms

{\rm(c-1)} If $1\le k-sn<n$ in $(3.2.7)$, then
$(k-sn)\beta+n(\alpha+s\beta-sn){\ge}(k-sn)n$ if and only if
$n\alpha+k\beta{\ge}nk$.

{\rm(c-2)} If $k-sn=n$ in $(3.2.7)$, then by Lemma $3.1$ we have the
following:
$$\align
(3.2.8) \qquad  & \text{There are nonzero complex numbers
$c_{{\alpha_i}{\beta_i}}$
for all $i=1,2,\dots,n-1$}  \qquad \qquad \\
& \text{such that $\alpha_i+s\beta_i-sn+\beta_i=n$, i.e.,
$\frac{\alpha_i}{n-\beta_i}=\frac{k}{n}=s+1$. {$\square$}}
\endalign$$
\endproclaim

Now, it is clear that Sublemma $3.2.1$ is true, and so we finished
the proof of the theorem.
\enddemo
\ms

\proclaim{Corollary 3.3} $\underline{\text{\bf Assumptions}}$
Suppose that $f\in \BC\{y,z\}$ with isolated singularity at $(0,0)$
satisfies the same properties and notations as in the assumption of
Theorem $3.2$.

In addition, we assume that the following equality holds:

\noindent {\rm (3.3.1)} \qquad \qquad $\gcd(n,k)=1.$ \ms

$\underline{\text{\bf Conclusions}}$

{\rm(1)} The necessary and sufficient condition for $f(y,z)$ to be
irreducible in $\BC\{y,z\}$ is as follows:
$$\align
\text{${\dfrac{\alpha_i}{i}} > {\dfrac{k}{n}}$ \quad for all
$i=1,2,\dots,n-1$.} \tag 3.3.2
\endalign$$

{\rm(2)} Equivalently, the inequality in {\rm(3.3.2)} holds if and
only if $f$ is represented as follows:
$$\align
(3.2.3) \qquad \qquad \qquad
\text{$f(y,z)=a_{0,0}z^n+a_{n,0}y^{k}+\sum_{\alpha,\beta\ge
0}c_{\alpha,\beta}y^{\alpha}z^{\beta}$} \quad \text{with
$n\alpha+k\beta> nk$,} \qquad \qquad
\endalign$$
where the $c_{\alpha,\beta}$ are nonzero complex numbers for some
nonnegative integers $\alpha$ and $\beta$ with $n\alpha+k\beta\ge
nk$, and $a_{0,0}$ and $a_{n,0}$ are nonzero constant.
\endproclaim \ms

\demo{\bf Proof of Corollary 3.3} It was already proved by Theorem
3.2 that the inequality in (3.3.2) is $\underline{\text{\bf the
necessary condition}}$ for $f$ to be irreducible in $\BC\{y,z\}$.
Now, in order to prove that  the inequality in (3.3.2) is
$\underline{\text{\bf a sufficient condition}}$ for $f(y,z)$ to be
irreducible in $\BC\{y,z\}$, assuming for convenience of proof that
$n<k$, the proof of this corollary can be easily finished by the
induction on the multiplicity of $f$ at the origin, using the same
method as we have seen in (3.2.6) and (3.2.7) of Sublemma $3.2.1$.
$\square$
\enddemo\ms

{\bf \S3.2. The algorithm for finding the total transform of an
irreducible plane curve singularity under the standard resolution
which has the homeomorphic resolution as the
curve{\text($z^n+y^k=0$)} does and its generalizations}

\proclaim{Lemma 3.4} $\underline{\text{\bf Assumptions}}$ \quad Let
$f(y,z)=a_nz^n+a_0y^k+\sum^{n-1}_{i=1}a_iy^{\alpha_i}z^i$ be
irreducible in $\BC\{y,z\}$ where for $0\le i\le n$, each
$a_i=a_i(y,z)$ is a unit in $\BC\{y,z\}$ if exists and the
$\alpha_i$ are positive integers. Let $d=\gcd(n,k)$ with $1\le n\le
k$, and write $n=n_1d$ and $k=k_1d$ with $\gcd(n_1,k_1)=1$.

$\underline{\text{\bf Conclusions}}$ \quad Then, $f$ can be
written in the form
$$\align
(3.4.1) \qquad \qquad f=A\prod^d_{i=1}(z^{n_1}+\xi_iy^{k_1})
+\sum_{\alpha,\beta\ge 0}c_{\alpha,\beta}y^{\alpha}z^{\beta} \quad
\text{with} \quad n_1\alpha+k_1\beta>n_1k_1d, \qquad \qquad
\endalign $$
 where $A=a_n(0,0)$ is a nonzero complex number, and for $1\le i
\le d$, all $\xi_i$ are nonzero complex numbers, and the
$c_{\alpha,\beta}$ are nonzero complex numbers for some nonnegative
integers $\alpha$ and $\beta$ with $n_1\alpha+k_1\beta>n_1k_1d$.
\endproclaim

\demo{\bf Proof of Lemma 3.4} Assume that $d=\gcd(n,k)>1$ with $2\le
n\le k$, otherwise there is nothing to prove. For any nonzero
monomial $y^{\alpha}z^{\beta}$ of $f(y,z)$ in the assumption, by
Theorem $3.2$ it suffices to consider two cases, denoted by Case(I)
$n\alpha+k\beta=nk$ and Case(II) $n\alpha+k\beta>nk$. Then, the
proof of the remainder is trivial.
\enddemo
\ms

\proclaim{Theorem 3.5(The representation theorem for the local
defining equations of irreducible plane curve singularities)}

$\underline{\text{\bf Assumptions}}$ \ Let
$f(y,z)=a_nz^n+a_0y^k+\sum^{n-1}_{i=1}a_iy^{\alpha_i}z^i$ be
irreducible in $\BC\{y,z\}$ where each $a_i=a_i(y,z)$ is a unit in
$\BC\{y,z\}$ for $0\le i\le n$, if exists and the $\alpha_i$ are
positive integers. Let $n$ be the multiplicity of $f$ at the origin
with $1\le n\le k$. Let $d=gcd(n,k)$, and write $n=n_1d$ and
$k=k_1d$ with $gcd(n_1,k_1)=1$. \ms

$\underline{\text{\bf Conclusions}}$ \ Then, f can be represented
as follows:
$$\align
f =A(z^{n_1}+ \xi y^{k_1})^d +\sum_{\alpha, \beta\ge 0}
   c_{\alpha,\beta}y^{\alpha}z^{\beta} \quad \text{with} \quad
   n_1\alpha+k_1\beta>n_1k_1d, \tag 3.5.1
 \endalign$$
where the $c_{\alpha,\beta}$ are nonzero complex numbers for some
nonnegative integers $\alpha$ and $\beta$, satisfying the following
properties :

\roster
\item "(i)" $A$ and $\xi$ are the unique nonzero complex
numbers such that $A=a_n(0,0)\ne 0$, $dA\xi=a_{n-n_1}(0,0)\ne 0$,
and \quad ${{d}\choose{i}}A\xi^i=a_{n-in_1}(0,0)$ for $1\le i\le
d$.

\item "(ii)" $\frac{\alpha_i}{n-i}\ge \frac{k}{n}=\frac{k_1}{n_1}$ for
$0\le i\le n-1$.
\endroster
\endproclaim
The proof of Theorem 3.5 will be done by Lemma $3.4$ and Theorem
3.6. \ms

\noindent{\bf Remark 3.5.1.} {\rm(a)} If $d=\gcd(n,k)=1$, then it is
clear by Corollary $3.3$ that the representation form in
{\rm(3.5.1)} is a sufficient condition for $f$ to be irreducible in
$\BC\{y,z\}$.

{\rm(b)} If $f$ is defined by $f=(z^2+y^3)^2+y^2z^4$, then note that
$f$ is not irreducible in $\BC\{y,z\}$, satisfying an equation in
{\rm(3.5.1)}.  \ms

\proclaim{Theorem 3.6}\ $\underline{\text{\bf Assumptions}}$ \quad
Let $V(g_i)=\{(y,z):g_i(y,z)=0\}$ for $1\le i\le d$, $V(f)=\{(y,z):
f(y,z)=0\}$ and $V(F)=\{(y,z): F(y,z)=0\}$ be analytic varieties at
$(0,0)$ in $\BC^2$, each of which is written respectively in the
form,
$$
\align
 g_i &=z^{n_1}+\xi_iy^{k_1}, \tag 3.6.1 \\
 f&=\prod^d_{i=1}g_i+\sum_{\alpha,\beta\ge 0}c_{\alpha,\beta}y^{\alpha}z^{\beta}
 \quad \text{with} \quad n_1\alpha+k_1\beta>n_1k_1d, \\
 F&=y^{\delta_1}z^{\delta_2}f,
\endalign
$$
satisfying the properties {\rm(i)}, {\rm(ii)}, {\rm(iii)},
{\rm(iv)}, {\rm(v)} and {\rm(vi)}: \roster \item "(i)"
$\gcd(n_1,k_1)=1$ with $1\le n_1<k_1$ and $d$ is a positive
integer.

\item "(ii)" The $\xi_i$ are nonzero complex numbers for $1\le
i\le d$, and $\mu$ is the counting number of distinct elements in
the set $\{\xi_1,\xi_2,\dots,\xi_d\}$, as a set, but not as a
sequence.

\item "(iii)" The $c_{\alpha,\beta}$ are nonzero complex numbers
for some nonnegative integers $\alpha$ and $\beta$ such that
$n_1\alpha+k_1\beta> n_1k_1d$, if exist.

\item "(iv)" Each $\delta_i$ is either a positive integer or $0$
for $i=1,2$.

\item "(v)" If $n_1=d=1$, assume additionally that $\delta_2>0$.

\item "(vi)" If $n_1=1$ and $\delta_2=0$, note by {\rm(v)} that
$d\ge 2$. If $d\ge 2$ and $n_1\ge 1$, then assume that $V(f)$ has an
isolated singular point at the origin as a reduced variety.
\endroster

$\underline{\text{\rm Note.}}$ By the property {\rm(v)}, if
$n_1=d=1$ then $V(f)$ does not have a singularity at $(0,0)$, but
$V(F)$ has an isolated singular point at $(0,0)$ as a reduced
variety.

In preparation for the construction of the statement in the
conclusion, let $V(G)=\{(y,z):G(y,z)=0\}$ be another analytic
variety with isolated singularity at the origin in $\BC^2$ defined
by the form
$$\align
G =z^{\gamma}g_1 \quad \text{and} \quad
g_1 =z^{n_1}+\xi_1y^{k_1}, \tag 3.6.2 \\
\endalign
$$
satisfying the properties {\rm(vii)} and {\rm(viii)}:

\roster \item "(vii)" If $n_1=1$, then $\gamma=1$.

\item "(viii)" If $n_1\ge 2$, then $\gamma=0$.
\endroster

Let $\tau_m$ be the composition of a finite number $m$ of
successive blow-ups which is needed only to get the standard
resolution of the singular point of $V(G)$. For each
$t=1,2,\dots,m$, write $\tau_t= \pi_1\circ\pi_2\circ\cdots
\circ\pi_t:M^{(t)}\to \BC^2$ where $\{\text{$ \pi_i:M^{(i)}\to
M^{(i-1)}$ is a blow-up}$

\noindent $\text{of $M^{(i-1)}$ at some point of $M^{(i-1)}$ for
$1\le i\le t$} \}$ with $M^{(0)}=\BC^2$. For brevity of notation,
let $V^{(t)}(G)$ be the proper transform under $\tau_t$ for $1\le
t\le m$.

Let $E^{(m)}=\tau^{-1}_m(0,0)$, and let $E^{(m)}=\cup E_i$, $1 \le
i \le m$, be the decomposition of $E^{(m)}$ into irreducible
components where each $E_i$ is called an exceptional curve of the
first kind.
\ms

$\underline{\text{\bf Conclusions}}$ \quad We have the following
facts.

{\bf Fact(1).} In order to study $V^{(t)}(G)$ under $\tau_t$, we can
find just one coordinate patch of the local coordinates for each
blow-up $\pi_t:M^{(t)}\to M^{(t-1)}$, where $1\le t\le m$ and
$M^{(0)}=\BC^2$, which will be proved in process of the standard
resolution of the singular point of $V(G)$ by the following lemma,
Lemma $4.1$.

{\bf Fact(2).} By {\rm Fact(1)}, we can use the same $\tau_m$ for
the composition of the first finite number $m$ of successive
blow-ups in preparation for finding the standard resolution of the
singular point $(0,0)$ of either $V(f)$ or $V(F)$ if exists, as a
reduced variety.

{\bf Fact(3).} In order to study each proper transform of either
$V(f)$ or $V(F)$ under $\tau_t$, without using a nonsingular
change of coordinates, we can use the common one coordinate patch
of the same local coordinates simultaneously, as it has been
already used for each blow-up $\pi_t:M^{(t)}\to M^{(t-1)}$ in {\rm
Fact(1)}, where $1\le t\le m$.

{\bf Fact(4).} If $f$ is irreducible in $\BC\{y,z\}$, then
$\mu=1$, that is, all the $\xi_i$ are the same. But, the converse
dose not hold. \ms

After $m$ iterations of blow-ups, let $(v_m,u_m)$ and
$(v'_m,u'_m)$ be the local coordinates for $M^{(m)}$ where by {\rm
Fact(1)} $\pi_m:M^{(m)}\to M^{(m-1)}$ was defined to be the $m$-th
blow-up at some point of $M^{(m-1)}$ with $u'_m=1/u_m$ and
$v'_m=v_mu_m$. Note that $E_m=\{v_m=0\}\cup \{v'_m=0\}$.

For brevity of the notation, we write
$$\align
\text{$(v,u)=(v_m,u_m)$ \quad and \quad $(v',u')=(v'_m,u'_m)$.} \tag
3.6.3
\endalign$$

Also, let $(F\circ\tau_m)_{divisor}$ be a divisor of
$F\circ\tau_m$ defined by
$$\align
(F\circ\tau_m)_{divisor}=V^{(m)}(F)+\sum^m_{i=1}e_iE_i, \tag 3.6.4
\endalign$$
where each $e_i$ is the multiplicity of $F\circ\tau_m$ along $E_i$
for $1\le i\le m$ and $V^{(m)}(F)$ is the proper transform of
$V(F)$ under $\tau_m$. \ms

Using the proofs for {\rm Fact(1)}, {\rm Fact(2)}, {\rm Fact(3)}
and {\rm Fact(4)} of the above, and also for the next statements
{\rm[I]}, {\rm[II]} and {\rm[III]} in the remainder of this
theorem, then we will find some elementary solutions for the
following problems:

{\bf Problem 1.} Along $E_m$, construct $\tau_m:M^{(m)}\to\BC^2$
as a local analytic mapping and the local defining equation
$(F\circ\tau_m)_{total}$ for the total transform of $V(F)$ under
$\tau_m$.

{\bf Problem 2.} Find the necessary condition for $f$ to be
irreducible in $\BC\{y,z\}$.

{\bf Problem 3.} Find the conditions under which $f$ and $F$
belong to the type $[1]$ under $\tau_m$ in the sense of Definition
$2.8$.

$$\text{\bf[I]}$$ \quad Along $v=0$, $\tau_m:M^{(m)}\to\BC^2$
as a composition of analytic mappings and $(F\circ\tau_m)_{total}$
can be written in the following form:
$$
\align
\tau_m(v,u)&=(y,z)=(v^{n_1}u^a,v^{k_1}u^b), \tag 3.6.5  \\
(F\circ\tau_m)_{total}&=v^{e_m}u^{\varepsilon}(f\circ\tau_m)_{proper}
\quad \text{with} \\
(f\circ\tau_m)_{proper}&=\prod^d_{i=1}(u+\xi_i)+\sum_{\alpha,\beta\ge
0}c_{\alpha,\beta}v^{n_1\alpha+k_1\beta-n_1k_1d}u^{\ve_{\alpha,\beta}},\\
(G\circ\tau_m)_{total}&=v^{k_1\gamma+n_1k_1}u^{b\gamma+ak_1}(u+\xi_1),
\endalign
$$
where {\rm(i)} $a$ and $b>0$ are nonnegative integers such that
$bn_1-ak_1=1$,

{\rm(ii)} $e_m=n_1\delta_1+k_1\delta_2+n_1k_1d$,
$\ve=a\delta_1+b\delta_2+ak_1d$ and
$\ve_{\alpha,\beta}=a\alpha+b\beta-ak_1d\ge 0$,

{\rm(iii)} by assumption of $V(G)$, if $n_1=1$ then $\gamma=1$, and
if $n_1\ge 2$ then $\gamma=0$. \ms

$$\text{\bf[II]}$$ \quad Along $v'=0$, $\tau_m:M^{(m)}\to\BC^2$ as a composition
of analytic mappings and $(F\circ\tau_m)_{total}$ can be written
in the following form:
$$
\align \tau_m(v',u') &=(y,z)=({v'}^{n_1}{u'}^{p},
{v'}^{k_1}{u'}^{q}), \tag 3.6.6  \\
(F\circ\tau_m)_{total} &={v'}^{e_m}{u'}^{\ve'}
(f\circ\tau_m)_{proper} \quad \text{with} \\
(f\circ\tau_m)_{proper} &=\prod^d_{i=1}(1+\xi_i
u')+\sum_{\alpha,\beta\ge 0}c_{\alpha,\beta} {v'}^{n_1\alpha + k_1
\beta - n_1 k_1 d}{u'}^{\ve'_{\alpha,\beta}}, \\
(G\circ\tau_m)_{total}
&={v'}^{k_1\gamma+n_1k_1}{u'}^{q\gamma+qn_1}(G\circ\tau_m)_{proper} \quad \text{with} \\
(G\circ\tau_m)_{proper} &=(1+\xi_1u'),
\endalign
$$
where {\rm(i)} $p$ and $q$ are positive integers such that
$pk_1-qn_1=1$,

{\rm(ii)} $e_m=n_1\delta_1+k_1\delta_2+n_1k_1d$,
$\ve'=p\delta_1+q\delta_2+qn_1d>1$ and
$\ve'_{\alpha,\beta}=p\alpha+q\beta-qn_1d>0$,

{\rm(iii)} $q\gamma+qn_1>1$ because $q>0$, noting that if $n_1=1$
then $\gamma=1$ by assumption of $V(G)$. \ms

$$\text{\bf[III]}$$

After $m$ iterations of blow-ups, denoted by $\tau_m$, we have the
following consequences: \roster \item "(i)" $V(G)\in the ~type[1]$
under $\tau_m$ in the sense of Definition $2.8$.

\item "(ii)" If $f$ is irreducible in $\BC\{y,z\}$, then all the $\xi_i$
are the same or $\mu=1$.

\item "(iii)" If $f$ is  irreducible in $\BC\{y,z\}$, then $F$ belongs
to the type $[1]$ under $\tau_m$ in the sense of Definition $2.8$
whether or not $f$ has a singularity at the origin.

\item "(iv)" Let $f$ be irreducible in $\BC\{y,z\}$ with isolated
singularity at the origin. If $n=1$, then $f\in the~ type[0]$ under
$\tau_m$ in the sense of Definition $2.8$, and if $n\ge 2$, then
$f\in the~ type[1]$ under $\tau_m$ in the sense of Definition $2.8$.
$\square$
\endroster
\endproclaim
The proof of Theorem 3.6 will be done in $\S4$. \ms

\definition{Remark 3.6.1} Let $\tau_m=\pi_1\circ\pi_2\circ\cdots\circ\pi_m$
be the composition of a finite number $m$ of successive blow-ups at
the origin in $\BC^2$ which is needed to get the standard resolution
of the singular point of $V(G)$, as we have seen in the assumption
of Theorem $3.6$. \ms

{\rm(a)} \ Let $f(y,z)=(z^2+y^3)^2+y^2z^4$ with $G(y,z)=z^2+y^3$.
Then, note that $G$ and $f$ satisfies the  same assumption as in
Theorem $3.6$. Then, $V(G)\in the \ type[1]$ under $\tau_m$, and
$V(f)\in the \ type[1]$ under $\tau_m$, too. But, it can be easily
shown that $f$ is not irreducible in $\BC\{y,z\}$. \ms

{\rm(b)} \  Let $\tau_t=\pi_1\circ\pi_2\circ\cdots\circ\pi_t$ for
$t=1,2,\dots,m$. Suppose that $G$, $f$ and $F$ satisfy the same
assumption as in Theorem $3.6$. For each $t=1,2,\dots,m$, consider
$qs(V^{(t)}(G))$, $qs(V^{(t)}(f))$ and $qs(V^{(t)}(F))$ under
$\tau_t$ in the sense of Definition $2.6$.

{\rm(b1)} \ For each $t=1,2,\dots,m-1$, $qs(V^{(t)}(G))$ is a
one-point set under $\tau_t$.

{\rm(b2)} \ If $t=m$, then $qs(V^{(m)}(G))$ is empty under $\tau_m$.

{\rm(b3)} \ For each $t=1,2,\dots,m-1$, $V^{(t)}(f)$ and
$V^{(t)}(F)$ has the same quasisingular point as $V^{(t)}(G)$ does
in $(b1)$, if exists.

{\rm(b4)} \ As in the conclusion of this theorem, after $m$
iterations of blow-ups, let $(v_m,u_m)$ and $(v'_m,u'_m)$ be the
local coordinates for $M^{(m)}$ where $\pi_m:M^{(m)}\to M^{(m-1)}$
was defined to be the $m-th$ blow-up at some point of $M^{(m-1)}$
with $u'_m=1/u_m$ and $v'_m=v_mu_m$.

Then, $qs(V^{(m)}(F))=qs(V^{(m)}(f))$ such that
$qs(V^{(m)}(F))=qs(V^{(m)}(F))\cap \{v_m=0\}=qs(V^{(m)}(F))\cap
\{v'_m=0\}$ and $qs(V^{(m)}(f))=qs(V^{(m)}(f))\cap
\{v_m=0\}=qs(V^{(m)}(f))\cap \{v'_m=0\}$. Note that
$E_m=\{v_m=0\}\cup \{v'_m=0\}$.
\enddefinition
 \ms

\proclaim{Theorem 3.7(The representation theorem for some total
transforms of the local defining equations of irreducible plane
curve singularities in process of the standard resolution of their
singularity)}

$\underline{\text{\bf Assumptions}}$ Let
$V(g_i)=\{(y,z):g_i(y,z)=0\}$ for $1\le i\le d$, $V(f)=\{(y,z):
f(y,z)=0\}$ and $V(F)=\{(y,z): F(y,z)=0\}$ be analytic varieties at
$(0,0)$ in $\BC^2$, satisfying the same properties and notations as
in the assumption of Theorem $3.6$.

In addition, for the necessity of the irreducibility of $f$ in
$\C\{y,z\}$, we may assume that the following holds:
$$\align
 \text{all the $\xi_i$ are the same numbers and so $g_i=g_1$
for all $i$.} \tag 3.7.0\\
\endalign$$

Let $V(\phi)=\{(y,z): \phi(y,z)=0\}$ and $V(\Psi)=\{(y,z):
\Psi(y,z)=0\}$ be analytic varieties at $(0,0)$ in $\BC^2$, each of
which is written respectively in the form,
 $$\align
 \phi&=g_1^{{\ell}d}+\sum_{p,q\ge
 0}a_{p,q}y^{p}z^{q}
 \qquad \qquad \text{with} \quad n_1p+k_1q>n_1k_1d{\ell}, \tag 3.7.1\\
 \Phi&=y^{\ell\delta_1}z^{\ell\delta_2}\phi,
\endalign$$
satisfying the properties {\rm(ix)} and {\rm(x)}:
\roster
\item "(ix)" $\ell$ is a positive integer.

\item "(x)" The $a_{p,q}$ are nonzero complex numbers for some nonnegative
integers $p$ and $q$ such that $n_1p+k_1q> n_1k_1d\ell$, if exist.
\endroster \ms

By {\rm Theorem 3.6}, we can use the same $\tau_m$ for the
composition of the first finite number $m$ of successive blow-ups in
the process of the standard resolution of the singular point $(0,0)$
of both $V(f)$ and $V(F)$ if exist, as reduced varieties. Let
$E^{(m)}=\tau^{-1}_m(0,0)$, and let $E^{(m)}=\cup E_i$, $1 \le i \le
m$, be the decomposition of $E^{(m)}$ into irreducible components
where each $E_i$ is called an exceptional curve of the first kind
under $\tau_m$.

As in $(3.6.4)$ of {\rm Theorem 3.6}, let
$(G\circ\tau_m)_{divisor}$, $(f\circ\tau_m)_{divisor}$,
$(F\circ\tau_m)_{divisor}$ and $(\Phi\circ\tau_m)_{divisor}$ be the
divisors of $G\circ\tau_{m}$, $f\circ\tau_{m}$, $F\circ\tau_{m}$ and
$\Phi\circ\tau_{m}$ under $\tau_m$, respectively, each of which can
be written as follows:
$$\align
(3.7.2)\quad  (G\circ\tau_m)_{divisor}
&=V^{(m)}(G)+\sum^m_{i=1}e_{0,i} E_i, \qquad
(f\circ\tau_m)_{divisor} =V^{(m)}(f)+\sum^m_{i=1}e_{1,i} E_i,  \qquad\\
(F\circ\tau_m)_{divisor} &=V^{(m)}(F)+\sum^m_{i=1}e_{2,i} E_i,
\qquad (\Phi\circ\tau_m)_{divisor}
=V^{(m)}(\Phi)+\sum^m_{i=1}e_{3,i} E_i,
\endalign$$
where $e_{0,i}$, $e_{1,i}$, $e_{2,i}$ and $e_{3,i}$ are the
multiplicities of $(G\circ\tau_m)_{divisor}$,
$(f\circ\tau_m)_{divisor}$, $(F\circ\tau_m)_{divisor}$ and
$(\Phi\circ\tau_m)_{divisor}$ along $E_{i}$, respectively, and
$V^{(m)}(G)$, $V^{(m)}(f)$, $V^{(m)}(F)$ and $V^{(m)}(\Phi)$ are the
proper transform of $V(G)$, $V(f)$, $V(F)$ and $V(\Phi)$ under
$\tau_{m}$, respectively. \ms

$\underline{\text{\bf Conclusions}}$ We have three statements, {\rm
[I], [II] and [III]}.

{\bf [I]} In the sense of {\rm Lemma 2.14}, using the same notations
and properties as in {\rm Theorem 3.6}, there is the composition of
a finite number $m$ of successive blow-ups, denoted by
$\tau_m=\pi_1\circ\pi_2\circ\cdots\circ\pi_m:M^{(m)}\to\BC^2$, which
can be viewed as a local analytic mapping from a polydisc
$\Delta(r_{11},r_{12})$ to a polydisc $\Delta(r_{21},r_{22})$, with
$$\align
&\Delta(r_{11},r_{12})=\{(v,u)\in\M^{(m)}:|v|<r_{11},
|u|<r_{12}\} \quad \text{and} \tag 3.7.3 \\
&\Delta(r_{21},r_{22})=\{(y,z)\in\BC^2:|y|<r_{21}, |z|<r_{22}\},
\endalign$$ satisfying the following properties :
$$\align
&\tau_m(v,u)=(y,z)=(v^{n_1}u^a,v^{k_1}u^b), \tag 3.7.4 \\
&(F\circ\tau_m)(v,u)=v^{e_m}u^{\varepsilon}(f\circ\tau_m)_{proper}, \\
&(f\circ\tau_m)_{proper}=(u+\xi)^d+\sum_{\alpha,\beta\ge
0}c_{\alpha,\beta}v^{n_1\alpha+k_1\beta-n_1k_1d}u^{\varepsilon_{\alpha,\beta}},
\endalign
$$
where \roster \item "(i)" $a$ and $b>0$ are nonnegative integers
such that $bn_1-ak_1=1$ and $E_m=\{(u,v):v=0\}$, \item "(ii)"
$e_m=n_1\delta_1+k_1\delta_2+n_1k_1d$, $\varepsilon
=a\delta_1+b\delta_2+ak_1d$ and $\varepsilon_{\alpha,\beta}=
a\alpha+b\beta-ak_1d\ge 0$.
\endroster

Also, if $n=1$, then $f\in the~ type[0]$ and $F\in the ~type[1]$
under $\tau_m$, but if $n\ge 2$, then $f\in the~ type[1]$ and $F\in
the ~type[1]$ under $\tau_m$ in the sense of Definition $2.8$.

Observe by the notation that $(f\circ\tau_m)_{proper}$ is the local
defining equation for the proper transform $V^{(m)}(f)$ and that
$(F\circ\tau_m)(v,u)$ is the local defining equation for the total
transform of $V(F)$ under $\tau_m$. \bs

{\bf [II]} Let ${\tau_m}^{-1}(0)=E^{(m)}=\cup E_i$, $1 \le i \le m$.
\ms

{\bf Fact(1) of [II]:} In case of $V(G)$, for $1\le i\le m$, each
$E_i$ belongs to exactly one of the following two types relative to
${{\tau_m}^{-1}(V(G))}$: For any $A\subset {\tau_m}^{-1}(V(G))$,
$\overline{A}$ is called the closure of $A$ in ${\tau_m}^{-1}(V(G))$
where ${\tau_m}^{-1}(V(G))=E^{(m)}\cup V^{(m)}(V(G))$. \ms

$\underline{\text{\rm Type(1)}}$ $E_{m}$ of ${{\tau_m}^{-1}(V(G))}$
and $\overline{{\tau_m}^{-1}(V(G))-E_{m}}$ of
${{\tau_m}^{-1}(V(G))}$ have three intersection points in
${{\tau_m}^{-1}(V(G))}$ in the sense of Definition $2.8$.

$\underline{\text{\rm Type(2)}}$ For any $j\not =m$, $E_{j}$ of
${{\tau_m}^{-1}(V(G))}$ and $\overline{{\tau_m}^{-1}(V(G))-E_{j}}$
of ${{\tau_m}^{-1}(V(G))}$ have at  most two intersection points in
${{\tau_m}^{-1}(V(G))}$ in the sense of {\rm Definition 2.8}.  \ms

By the same method as in {\rm Fact(1)}, each of $V(F)$ and $V(\Phi)$
have the same type of exceptional curves of the first kind relative
to the corresponding one of ${{\tau_m}^{-1}(V(F))}$ and
${{\tau_m}^{-1}(V(\Phi))}$ respectively, as $V(G)$ does in {\rm
Fact(1)}, as a reduced variety. Also, whenever $V(f)$ has an
isolated singularity as a reduced variety, then $V(f)$ has the same
type of exceptional curves of the first kind relative to
${{\tau_m}^{-1}(V(f))}$ as the above $V(G)$ does. \ms

{\bf Fact(2) of [II]:} In any case of $V(G)$, $V(f)$, $V(F)$ and
$V(\Phi)$, for $1\le i\le m$, each $E_i$ belongs to exactly one of
the following two types relative to $E^{(m)}=\cup E_i$, $1 \le i \le
m$: \ms

$\underline{\text{\rm Type(1)}}$ \quad For each $i=1,2,\dots,m$,
there are two distinct exceptional curves of the first kind in
$E^{(m)}$, denoted by $E_{1}$ and $E_{s}$ with $1<s\le m$, which
satisfy the same property in the following sense: Note that $s$ is
an integer such that $(s-1)n_1<k_1\le sn_1$.

{\rm (1a)} $E_{1}$ of $E^{(m)}$ and $\overline{E^{(m)}-E_{1}}$ of
$E^{(m)}$ has one and only one intersection point in $E^{(m)}$.

{\rm(1b)} For some integer $s\not= 1$, $E_{s}$ of $E^{(m)}$ and
$\overline{E^{(m)}-E_{s}}$ of $E^{(m)}$ has one and only one
intersection point in $E^{(m)}$.

$\underline{\text{\rm Type(2)}}$ \quad If $E_t$ of $E^{(m)}$ is
neither $E_1$ nor $E_s$ in {\rm Type(1)}, then $E_{t}$ of $E^{(m)}$
and $\overline{E^{(m)}-E_{t}}$ of $E^{(m)}$ have two distinct
intersection points in $E^{(m)}$. \bs

{\bf [III]} Then, we have the following:
$$
\text{$\frac{e_{3,i}}{e_{2,i}}=\ell$ \quad   for all
$i=1,2,\dots,m$.} \tag 3.7.5
$$

In particular, $e_{0,m}=k_1\gamma+n_1k_1$ and
$e_{2,m}=n_1\delta_1+k_1\delta_2+n_1k_1d$.

Moreover, whenever an analytic variety $V(F)$ satisfies the
conditions in $(3.6.1)$, then a finite sequence $\{e_i:
i=1,2,\dots,m\}$ is uniquely determined by
$\{\delta_1,\delta_2,n_1,k_1,d\}$, without depending on the nonzero
complex numbers $c_{\alpha\beta}$. $\square$
\endproclaim
The proof of Theorem $3.7$ will be done in $\S 4$. \bs

\definition{Corollary 3.8} $\underline{\text{\bf Assumptions}}$ As in
Corollary $3.3$ or Theorem $3.7$, assume that
$V(h_1)=\{(y,z):h_1(y,z)=0\}$ is represented as follows: Note that
$2\le n_1<k_1$.
$$\align
(3.8.1) \qquad \qquad  &
\text{$h_1(y,z)=z^{n_1}+y^{k_1}+\sum_{\alpha,\beta\ge
0}c_{\alpha,\beta}y^{\alpha}z^{\beta}$ } \quad \text{with
$n_1\alpha+k_1\beta>n_1k_1$,} \qquad \qquad \qquad
\endalign$$
where the $c_{\alpha,\beta}$ are nonzero complex numbers for some
integers $\alpha\ge 0$ and $\beta \ge 0$, if exist.  \ms

$\underline{\text{\bf Conclusions}}$

{\rm(i)} Let $g_1(y,z)=z^{n_1}+y^{k_1}$ with $\gcd(n_1,k_1)=1$.
Then, $g_1(y,z)=0$ and $h_1(y,z)=0$ have the homeomorphic resolution
at $(y,z)=(0,0)$. \ms

{\rm(ii)} If $f(y,z)$ of Theorem $3.7$ is irreducible in
$\BC\{y,z\}$, then we have the following:
$$\align
(3.8.2) \qquad \qquad & \text{\rm Multiseq(V(f))}=\text{\rm
Join}(\{[dn_1:dk_1]\},
\text{\rm Multiseq($V^{(m)}(f)$)}) \quad \text{by Definition $1.7$},  \\
& \text{\rm Multiseq(V(g))} =\{[n_1:k_1]\},   \\
\endalign$$
noting that if $\{[n_1:k_1]\}=\{c_i:i=1,2,\dots,w\}$, write
$\{[dn_1:dk_1]\}=\{dc_i:i=1,2,\dots,w\}$ for notation.
\enddefinition \ms

\vfill \pagebreak

{\bf \S{4}. The proofs of Theorem 3.6 and Theorem 3.7 in \S {3}} \bs

{\bf \S4.0. Introduction } \ms

In preparation for the proof of Theorem $3.6$, we prove three
lemmas, Lemma $4.1$, Lemma $4.2$ and Lemma $4.3$ in $\S4.1$. In
\S{4.2,} using three lemmas, it suffices to prove Theorem $3.6$
because Theorem $3.5$ and Theorem $3.7$ can be proved as an easy
corollary of Theorem $3.6$. \ms

\text{\bf \S{4.1}. Lemma 4.1, Lemma 4.2 and Lemma 4.3 with proofs}

\proclaim{Lemma 4.1} $\underline{\text{\bf Assumptions}}$ \quad
Let $V(G)=\{(y,z):G(y,z)=0\}$ be an analytic variety at $(0,0)$ in
$\BC^2$ defined by the form
$$\align
G =z^{\gamma}g \quad \text{and} \quad
g =z^{n_1}+y^{k_1} \quad \text{with} \quad \gcd(n_1,k_1)=1, \tag 4.1.1 \\
\endalign$$
satisfying the following properties:

{\rm(i)}  $1\le n_1 <k_1$. {\rm(ii)} If $n_1=1$, then $\gamma=1$.
{\rm(iii)} If $n_1 \ge 2$, then $\gamma=0$. \ms

$\underline{\text{\bf Conclusions}}$ \quad We can find the standard
resolution of the singular point $(0,0)$ of $V(G)$, that is, the
composition $\tau_m=\pi_1\circ\pi_2\circ\cdots\circ\pi_m$ of a
finite number $m$ of successive blow-ups at the origin in $\BC^2$
such that for each blow-up $\pi_j$ with $1\le j\le m$, just one
coordinate patch is needed in order to study the $j$-th proper
transform $V^{(j)}(G)$ under $\tau_j=\pi_1\circ\cdots \circ\pi_j$.
$\square$
\endproclaim

The proof of Lemma 4.1 is trivial. \ms

\proclaim{Lemma 4.2} $\underline{\text{\bf Assumptions}}$ \quad Let
$V(G)=\{(y,z):G(y,z)=0\}$ be an analytic variety at $(0,0)$ in
$\BC^2$, satisfying the same properties and notations as the
assumption of Lemma $4.1$. In addition, for given integers $n_1$ and
$k_1$ in $(4.1.1)$, let $s$ be a positive integer such that
$sn_1<k_1\le (s+1)n_1$. Since $gcd(n_1,k_1)=1$ with $1\le n_1<k_1$,
note that
$$\align
&\text{if \quad $n_1=1$, then \quad $k_1=s+1$,} \tag 4.2.1 \\
&\text{if \quad $n_1\ge 2$, then \quad $sn_1<k_1<(s+1)n_1$.}
\endalign$$

$\underline{\text{\bf Conclusions}}$ \quad Let $\tau_m$ be the
composition of a finite number $m$ of successive blow-ups which is
needed only to get the standard resolution of the singular point
of $V(G)$.

Then, $1\le s<m$, and in preparation for the proof of Theorem $3.6$,
we can construct the first $t$ iterations of blow-ups, in order to
study the $t-th$ proper transform $V^{(t)}(G)$ of $V(G)$ under
$\tau_t=\pi_1\circ\pi_2\circ\cdots \circ\pi_t$ for each
$t=1,2,\dots,s$, as follows.

{\rm(a)} For each $t=1,2,\dots,s$, $qs(V^{(t)}(G))$ is a one-point
set in the sense of Definition $2.6$.

{\rm(b)} In order to study the $t-th$ proper transform $V^{(t)}(G)$
under $\tau_t=\pi_1\circ\pi_2\circ\cdots \circ\pi_t$, we can use
just one coordinate patch of the local coordinates for each blow-up
$\pi_i:M^{(i)}\to M^{(i-1)}$ at some point of $M^{(i-1)}$,  where
$M^{(0)}=\BC^2$ with $(v_0,u_0)=(y,z)$, as follows: For each
$i=1,2,\dots,t$, let $(v_i,u_i)$ and $(v'_i,u'_i)$ be the local
coordinates for $M^{(i)}$ with $u'_i={1/u_i}$ and $v'_i=v_iu_i$.
Then, $\pi_t(v_i,u_i)=(v_{i-1},u_{i-1})=(v_i,v_iu_i)$ and
$\pi_t(v'_i,u'_i)=(v_{i-1},u_{i-1})=(v'_iu'_i,v'_i)$ where
$(v'_{i},u'_{i})$ is not needed for the study of $V^{(i)}(G)$ under
$\pi_{i}$.

{\rm(c)} For $1\le t\le s$, let
$\tau_t=\pi_1\circ\pi_2\circ\cdots\circ\pi_t:M^{(t)}\to\BC^2$ be
defined by the local coordinates in $(b)$ where
$E_t=\{v_t=0\}\cup\{v'_t=0\}$.

Suppose that $\tau_t:M^{(t)}\to\BC^2$ satisfies the same assumptions
and notations as in $(b)$. For each fixed $t$, along $E_t$, $\tau_t$
can be represented, as a composition of analytic mappings, as
follows: Note that $u'_t={1/u_t}$ and $v'_t=v_tu_t$.
$$
\align
 \tau_t(v_t,u_t) &= (y,z)=(v_t,{v_t}^t u_t) \quad \text{and} \tag \text{4.2.2} \\
 \tau_t(v'_t,u'_t) &=(y,z)=(v'_t u'_t,{v'_t}^t {u'_t}^{t-1}).
\endalign
$$

So, using the first $s$ iterations of blow-ups, denoted by
$\tau_s$, the sequence of coordinates
$\{(v'_1,u'_1),(v'_2,u'_2),\dots,(v'_{s},u'_{s})\}$ is not needed
for the study of $V^{(s)}(G)$ by {\rm(b)}.

In addition, $qs(V^{(t)}(G))=\{(v_t,u_t)=(0,0)\}$ for each
$t=1,2,\dots,s$.

Moreover, a divisor of $G\circ\tau_s$ under $\tau_s$, denoted by
$(G\circ\tau_s)_{divisor}$, can be written as follows:
$$\align
 (G\circ\tau_s)_{divisor}&=V^{(s)}(G)+\sum^s_{j=1}e_jE_j \tag 4.2.3 \\
 \quad \text{where} \quad &e_j=(\gamma+n_1)j \quad \text{for} \quad
 j=1,2,\dots,s.
\endalign$$
Note that $V^{(s)}(G)$ is the proper transform of $V(G)$ under
$\tau_s$. $\square$
\endproclaim

\demo{\bf Proof of Lemma 4.2} Following the same method as we have
used in the proof of Sublemma $3.2.1$ of Theorem $3.2$ for each
$t=1,2,\dots,s$, it is easy to find the proof of (a), (b) and (c) at
the same time. Also, there is nothing to prove for (4.2.3). Thus,
the proof can be finished. $\square$
\enddemo
\ms

\definition{Remark 4.2.1} Recall by Theorem $3.6$ that
$g_1(y,z)=z^{n_1}+\xi_1 y^{k_1}$. with $n_1<k_1$. Then,
$g(y,z)=g_1(y,z)$ with $\xi_1=1$ and $G(y,z)=z^{\gamma}g
=z^{\gamma}g_1$ where $\gamma$ is either $1$ or zero. In particular,
$G(y,z)$ can be defined by $F(y,z)$ of (3.6.1) with some
coefficients. Also, note that $\{z=0\}$ and $\{G(y,z)=0\}$ meet
tangentially at the origin.
\enddefinition \ms

\proclaim{Lemma 4.3} $\underline{\text{\bf Assumptions}}$ Suppose
that $V(g_i)$, $1\le i\le d$, $V(f)$, $V(F)$ and $V(G)$ satisfy the
same properties and notations as in the assumption of Theorem $3.6$
and Lemma $4.2$.

Let $s$ be a positive integer such that $sn_1<k_1\le (s+1)n_1$,
and $\tau_m$ be the composition of a finite number $m$ of
successive blow-ups which is needed only to get the standard
resolution of the singular point of $V(G)$, where $s<m$. Using the
same properties and notations as in the conclusion of Lemma $4.2$,
then we proved by Lemma $4.2$ that we can construct the first $t$
iterations of blow-ups, denoted by $\tau_t$, in order to study the
t-th proper transform $V^{(t)}(G)$ under
$\tau_t=\pi_1\circ\pi_2\circ\cdots \circ\pi_t$ for each
$t=1,\dots,s$, as follows:

{\bf Property(a).} For the study of the $t-th$ proper transform
$V^{(t)}(G)$ with $1\le t\le s$, we can use one and only one
coordinate patch of the local coordinates for each blow-up
$\pi_i:M^{(i)}\to M^{(i-1)}$ at some point of $M^{(i-1)}$, where
$1\le i\le t$ and $M^{(0)}=\BC^2$.

{\bf Property(b).} Then, $qs(V^{(t)}(G))=\{(v_t,u_t)=(0,0)\}$ for
each $t=1,2,\dots,s$. \ms

$\underline{\text{\bf Conclusions}}$ \quad We can get five facts:

{\bf Fact(1).} Whenever $V(F)$ has the singular point at the
origin as a reduced variety, then we can use the same $\tau_s$ for
the composition of the first finite number $s$ of successive
blow-ups in preparation for the standard resolution of the
singular point $(0,0)$ of $V(F)$, as we have seen in the
assumption of this lemma.

{\bf Fact(2).}  In order to study the $t-th$ proper transform
$V^{(t)}(F)$ of $V(F)$ under $\tau_t$, without using a nonsingular
change of coordinates, we can use one and only common one coordinate
patch of the same local coordinates simultaneously, as it was
already used for each blow-up $\pi_t:M^{(t)}\to M^{(t-1)}$ in the
above {\rm Property(a)}, where $1\le t\le s$.

{\bf Fact(3).} Each $t-th$ proper transform $V^{(t)}(F)$ under
$\tau_t$ has the same quasisingular point $(v_t,u_t)=(0,0)$ as
$V^{(t)}(G)$ does in the above {\rm Property(b)}, as a reduced
variety for $1\le t\le s$.

{\bf Fact(4).} Along $E_t$, we can construct
$\tau_t:M^{(t)}\to\BC^2$ as an analytic mapping in the sense of
$(4.2.2)$ and the local defining equation $(F\circ\tau_t)_{total}$
for the total transform of $V(F)$ under $\tau_t$, in order to
study the $t-th$ proper transform $V^{(t)}(F)$ where $1\le t\le
s$.

{\bf Fact(5).} Then, $F\in \text{\rm {the type}}[0]$ under $\tau_s$
in the sense of \text{\rm Definition 2.8}. \ms

In order to prove the above five facts, in the remainder of this
lemma, it suffices to construct the following three statements
{\rm[I]}, {\rm[II]} and {\rm[III]}, respectively: \ms

{\bf[I]}\quad For any $t=1,2,\dots,s$, along $v_t=0$,
$\tau_t:M^{(t)}\to\BC^2$ as a composition of analytic mappings in
the sense of $(4.2.2)$ and $(F\circ\tau_t)_{total}$ can be
written, respectively as follows:
$$ \align
(4.3.1) \qquad \qquad \tau_t(v_t,u_t) &=(y,z)=(v_t,{v_t}^t u_t),  \\
 (F\circ\tau_t)_{total}
 &={v_t}^{e_t}(F\circ\tau_t)_{proper} \quad \text{with} \quad
 (F\circ\tau_t)_{proper}={u_t}^{\delta_2}(f\circ\tau_t)_{proper}, \\
 (f\circ\tau_t)_{proper} &=\prod^d_{i=1}({u_t}^{n_1}+\xi_i{v_t}^{k_1-tn_1})
 + \sum_{\alpha,\beta\ge 0}
  c_{\alpha,\beta}{v_t}^{\lambda_t}{u_t}^{\beta},  \\
 (f\circ\tau_t)_{total} &={v_t}^{tn_1d}(f\circ\tau_t)_{proper}, \\
(G\circ\tau_t)_{total} &=v_t^{t\gamma+tn_1}(G\circ\tau_t)_{proper}
\quad \text{with} \quad (G\circ\tau_t)_{proper}
=u_t^{\gamma}(u_t^{n_1}+\xi_1v_t^{k_1-tn_1}),
 \endalign
$$
where

{\rm(i)} $e_t=\delta_1+t\delta_2+tn_1d$ and
$\lambda_t=\alpha+t\beta-tn_1d>0$ for $1\le t\le s$, respectively,

{\rm(ii)} $n_1\lambda_t+(k_1-tn_1)\beta>n_1(k_1-tn_1)d$, that is,
$n_1\alpha+k_1\beta>n_1k_1d$ for all $\alpha$ and $\beta$. \ms

Moreover, $V^{(t)}(G)$ has one and only one quasisingularity at
$(v_t,u_t)=(0,0)$ along $v_t=0$. Also, $V^{(t)}(F)$ has one and
only one quasisingularity at $(v_t,u_t)=(0,0)$ along $v_t=0$, as a
reduced variety. \bs

{\bf[II]} For any $t=1,2,\dots,s$, along $v'_t=0$,
$\tau_t:M^{(t)}\to\BC^2$ as a composition of analytic mappings in
the sense of $(4.2.2)$ and $(F\circ\tau_t)_{total}$ can be
written, respectively as follows:
$$ \align
 (4.3.2) \quad \quad \tau_t(v'_t,u'_t) &=(y,z)=(v'_tu'_t,{v'_t}^t{u'_t}^{t-1}),  \\
(F\circ\tau_t)_{total}&= \left\{\eqalign{& {v'_1}^{e_1}
(F\circ\tau_1)_{proper}, \qquad \qquad \text{if \ $t=1$ or $s=1$},
\cr &{v'_t}^{e_t}{u'_t}^{e_{t-1}}
(F\circ\tau_t)_{proper}, \qquad  \text{if \ $t\ge 2$}, \cr} \right.  \\
(F\circ\tau_t)_{proper}&= \left\{\eqalign{&
{u'_1}^{\delta_1}(f\circ\tau_1)_{proper}, \qquad \qquad \text{if \
$t=1$ or $s=1$}, \cr
&(f\circ\tau_t)_{proper},\qquad \qquad \qquad  \text{if \ $t\ge 2$}, \cr} \right.  \\
(f\circ\tau_t)_{proper}
&=\prod^d_{i=1}(1+\xi_i{v'_t}^{k_1-tn_1}{u'_t}^{k_1-(t-1)n_1})
+\sum_{\alpha,\beta\ge 0}
c_{\alpha,\beta}{v'_t}^{\lambda_t}{u'_t}^{\lambda_{t-1}} \quad \text{for \ $t\ge 1$,} \quad \\
(f\circ\tau_t)_{total}&={v'_t}^{tn_1d}(f\circ\tau_t)_{proper}
\quad \qquad \qquad \qquad \qquad  \text{for \ $t\ge 1$,} \\
(G\circ\tau_t)_{total}
&={v'_t}^{t\gamma+tn_1}{u'_t}^{(t-1)\gamma+(t-1)n_1}
(G\circ\tau_t)_{proper} \quad  \text{for \ $t\ge 1$,} \\
(G\circ\tau_t)_{proper}
&=1+\xi_1{v'_t}^{k_1-tn_1}{u'_t}^{k_1-(t-1)n_1} \quad \qquad
\qquad \text{for \ $t\ge 1$,}
\endalign$$
where $e_t=\delta_1+t\delta_2+tn_1d$ and $\lambda_0=\alpha$ are
well-defined, being compared with $\lambda_t=\alpha+t\beta-tn_1d>0$
for $1\le t\le s$, respectively. \ms

But, note that $(F\circ\tau_t)_{proper}$ is represented as follows
:

{\rm(i)} If $t=1$, then
$(F\circ\tau_1)_{proper}={u'_1}^{\delta_1}(f\circ\tau_1)_{proper}$.

{\rm(ii)} If $t\ge 2$, then
$(F\circ\tau_t)_{proper}=(f\circ\tau_t)_{proper}$.

Moreover, there is no quasisingular point of $V^{(t)}(G)$ along
$v'_t=0$. Also, there is no quasisingular point of $V^{(t)}(F)$
along $v'_t=0$, as a reduced variety. \ms

{\bf[III]} \quad After $t$ iterations of blow-ups for any
$t=1,2,\dots,s$, $V(G)\in the ~ type~ [0]$ and also $V(F)\in the ~
type~ [0]$ under $\tau_t$ whether or not $f$ is irreducible in
$\BC\{y,z\}$, satisfying the following properties : Let
$\cup^t_{i=1}E_i=\tau^{-1}_t(0,0)$ where each $E_i$ is called an
exceptional curve of the first kind under $\tau_t$.

{\rm(a)} $qs(V^{(t)}(G))=qs(V^{(t)}(F))=\{(v_t,u_t)=(0,0)\}$. \ms

{\rm(b)} For $1\le i\le t$, each $E_i$ has at most two distinct
intersection points with any other exceptional curves and
$V^{(t)}(F)$. Also, no three distinct components of
$\cup^t_{i=1}E_i$ meet.\ms

{\rm(c)} In more detail, for $2\le i\le t$, each $E_i$ has two
distinct intersection points with any other exceptional curves and
$V^{(t)}(F)$. But, $E_1$ has at most two distinct intersection
points with any other exceptional curves and $V^{(t)}(F)$. In
particular, if $\delta_1>0$, then $E_1$ has two distinct
intersection points with $\cup^t_{i=2}E_i$ and $V^{(t)}(F)$, and if
$\delta_1=0$, then $E_1$ has just one intersection point with
$\cup^t_{i=2}E_i$ and $V^{(t)}(F)$, whether or not $t=1$. $\square$
\endproclaim  \ms

The proofs of three statements [I], [II] and [III] in Lemma $4.3$
can be easily done by induction on $t$ with $1\le t\le s$, and so
the proof of the lemma is completely finished. \bs

\text{\bf \S{4.2}. The proofs of Theorem 3.6 and Theorem 3.7} \ms

\demo{\bf Proof of Theorem 3.6} For proof of this theorem, we are
going to follow continuously the same notations and indices as we
have used in the proof of Lemma $4.3$.

By a given integer $s$ with $sn_1<k_1\le (s+1)n_1$, recall by
$(4.3.1)$ of Lemma $4.3$ that along $v_s=0$,
$\tau_s:M^{(s)}\to\BC^2$ as a composition of analytic mappings, and
$(F\circ\tau_s)_{total}$ and $(G\circ\tau_s)_{total}$ can be
written, respectively in the form
$$\align
(3.6.7) \qquad \qquad \tau_s(v_s,u_s)&=(y,z)=(v_s,v_s^s u_s),  \\
(F\circ\tau_s)_{total}&=v_s^{e_s}(F\circ\tau_s)_{proper} \quad
\text{with} \quad
 (F\circ\tau_s)_{proper}=u_s^{\delta_2}(f\circ\tau_s)_{proper},
\\
(f\circ\tau_s)_{proper}&=\prod^d_{i=1}(u_s^{n_1}+\xi_i{v_s}^{k_1-sn_1})
+\sum_{\alpha,\beta\ge 0} c_{\alpha,\beta}v_s^{\lambda_s}u_s^{\beta}, \\
(f\circ\tau_s)_{total} &={v_s}^{sn_1d}(f\circ\tau_s)_{proper}, \\
(G\circ\tau_s)_{total} &=v_s^{s\gamma+sn_1}(G\circ\tau_s)_{proper}
\quad \text{with} \quad (G\circ\tau_s)_{proper}
=u_s^{\gamma}(u_s^{n_1}+\xi_1v_s^{k_1-sn_1}),
\endalign
$$
where (i) $e_s=\delta_1+s\delta_2+sn_1d$ and
$\lambda_s=\alpha+s\beta-sn_1d>0$ by $(4.3.1)$,

\quad (ii) $n_1\lambda_s+(k_1-sn_1)\beta>n_1(k_1-sn_1)d$, that is,
  $n_1\alpha+k_1\beta>n_1k_1d$.

Note that $qs(V^{(s)}(F))=qs(V^{(s)}(G)) =\{(v_{s},u_{s})=(0,0)\}$
along $E_s$, as a reduced variety. \ms

For the proof of the theorem in more detail, in order to apply the
induction method to the multiplicity $n_1$ of
$g_1=z^{n_1}+\xi_1y^{k_1}$ at $(y,z)=(0,0)$ where $1\le n_1<k_1$
and $G(y,z)=z^{\gamma}g_1=z^{\gamma}(z^{n_1}+\xi_1y^{k_1})$, it is
enough to consider two cases:

{\bf Case(A):}  Let $n_1=1$. Then $\gamma=1$ by assumption.

{\bf Case(B):}  Let $n_1\ge 2$. Then $\gamma=0$ by assumption.  \ms

Therefore, first we will show that Case(A) satisfies
$[\text{\rm{I}}]$, $[\text{\rm{II}}]$ and $[\text{\rm{III}}]$,
respectively, and next we will show that Case(B) satisfies
$[\text{\rm{I}}]$, $[\text{\rm{II}}]$ and $[\text{\rm{III}}]$,
respectively, where [I], [II] and [III] are three distinct
statements in the conclusion of this theorem. \ms

{\bf Case(A):} Let $n_1=1$. In order to prove that this case
satisfies $[\text{\rm{I}}]$, $[\text{\rm{II}}]$ and
$[\text{\rm{III}}]$, let $V(f)$, $V(F)$ and $V(G)$ satisfy the same
properties and notations as in $(3.6.1)$ and $(3.6.2)$ of the
assumption of Theorem $3.6$, each of which was already defined by
the following form:
$$\align
 f &=\prod^d_{i=1}(z+\xi_iy^{k_1})+ \sum_{\alpha,\beta\ge
 0}c_{\alpha,\beta}y^{\alpha}z^{\beta} \quad \text{with
  $\alpha+k_1\beta>k_1d$},  \tag 3.6.8  \\
 F &=y^{\delta_1}z^{\delta_2}f \quad
\text{and} \quad G =z(z+\xi_1y^{k_1}).
\endalign$$

Since $sn_1<k_1\le (s+1)n_1$ with $n_1=1$, note that
$k_1=(s+1)n_1=s+1$.

Along $v_s=0$, $(F\circ\tau_s)_{total}$ can be rewritten as follows,
using (3.6.7) and (3.6.8):
$$\align
(3.6.9) \qquad \qquad
(F\circ\tau_s)_{total}&=v_s^{e_s}(F\circ\tau_s)_{proper} \quad
\text{with} \quad
(F\circ\tau_s)_{proper}=u_s^{\delta_2}(f\circ\tau_s)_{proper}, \qquad \qquad\\
(f\circ\tau_s)_{proper} &=\prod^d_{i=1}(u_s+\xi_i v_s)
 +\sum_{\alpha,\beta\ge 0} c_{\alpha,\beta}v_s^{\lambda_s}u_s^{\beta},\\
(G\circ\tau_s)_{total} &=v_s^{2s}(G\circ\tau_s)_{proper}\quad
\text{with} \quad (G\circ\tau_s)_{proper} =u_s(u_s+\xi_1v_s),
\endalign $$
where $e_s=\delta_1+s\delta_2+sd$, and also
$\lambda_s=\alpha+s\beta-sd>0$ by (4.3.1), because $n_1=1$ implies
that $\gamma=1$ by assumption of the theorem.

It is interesting to observe from (3.6.9) that $n_1\lambda_s+(k_1-s
n_1)\beta=1(\alpha+s\beta-sd)+1\cdot \beta>1\cdot d$ if and only if
$\alpha+k_1\beta>k_1d$ for all $\alpha\ge 0$ and $\beta\ge 0$. \ms

For the proof,  let $\pi_{s+1}:M^{(s+1)}\to M^{(s)}$ be a blow-up
at $(v_s,u_s)=(0,0)$ defined by
$$
\align
\pi_{s+1}(v_{s+1},u_{s+1})&=(v_s,u_s)=(v_{s+1},v_{s+1}u_{s+1}), \tag 3.6.10-a \\
\pi_{s+1}(v'_{s+1},u'_{s+1})&=(v_s,u_s)=(v'_{s+1}u'_{s+1},v'_{s+1}),
\tag 3.6.10-b
\endalign$$
where $(v_{s+1},u_{s+1})$ and $(v'_{s+1},u'_{s+1})$ are the local
coordinates for $M^{(s+1)}$ with $u'_{s+1}=1/u_{s+1}$ and
$v'_{s+1}=v_{s+1}u_{s+1}$, and write $E_{s+1}=\{v_{s+1}=0\}\cup
\{v'_{s+1}=0\}$.

Let $\tau_{s+1}=\tau_s\circ\pi_{s+1}$, and so (3.6.7), (3.6.10-a)
and (3.6.10-b) imply that
$$\align
\text{(3.6.11-a)} \qquad \qquad \tau_{s+1}(v_{s+1},u_{s+1})
&=(\tau_s\circ\pi_{s+1})(v_{s+1},u_{s+1})=\tau_s(v_{s+1},v_{s+1}u_{s+1})
\qquad \qquad
\\
 &=(v_{s+1},v_{s+1}^{s+1}u_{s+1})=(y,z) \quad \text{and}  \\
\text{(3.6.11-b)} \qquad \qquad \tau_{s+1}(v'_{s+1},u'_{s+1})
&=(\tau_s\circ\pi_{s+1})(v'_{s+1},u'_{s+1})=\tau_s(v'_{s+1}u'_{s+1},v'_{s+1})
\qquad \qquad
\\ &=(v'_{s+1}u'_{s+1},{v'_{s+1}}^{s+1}{u'_{s+1}}^s)=(y,z).
\endalign$$

Since $k_1=s+1$, then (3.6.11-a) and (3.6.11-b) are rewritten as
follows: For brevity of notations, write $(v,u)=(v_{s+1},u_{s+1})$
and $(v',u')=(v'_{s+1},u'_{s+1})$.
$$\align
& \tau_{k_1}(v,u)=(y,z)=(v,v^{k_1}u), \tag 3.6.12-a \\
& \tau_{k_1}(v',u')=(y,z)=(v'u',{v'}^{k_1}{u'}^s). \tag 3.6.12-b \\
\endalign$$
Thus, it can be easily proved by (3.6.9) and (3.6.10-a) and
(3.6.10-b) that $\tau_{k_1}$ is the desired standard resolution of
the singular point of $V(G)$.

Now, we will prove [I], [II] and [III] for Case(A), respectively as
follows. \ms

{\bf The proof of [I] for Case(A):} By (3.6.5) and (3.6.12-a),
compare the following:
$$\align
\tau_m(v,u)&=(y,z)=(v^{n_1}u^a,v^{k_1}u^b), \tag 3.6.13 \\
\tau_{k_1}(v,u)&=(y,z)=(v,v^{k_1}u),
\endalign$$
where $m$ is, by definition, a finite number of successive
blow-ups which is needed only to get the standard resolution of
the singular point of $V(G)$. Then, $m=k_1$, $n_1=1$, $a=0$ and
$b=1$. So, $bn_1-ak_1=1-0=1$.

Assuming that $n_1=1$, then apply (3.6.12-a) to (3.6.1) and (3.6.2),
respectively. Or, we may apply (3.6.10-a) to (3.6.9), if convenient.

Then, along $v=0$,
$(F\circ\tau_{k_1})_{total}=(F\circ\tau_{k_1})(v,u)$ can be
written in the form
$$\align
(3.6.14) \qquad \qquad  (F\circ\tau_{k_1})_{total}
&=v^{\delta_1}(v^{k_1}u)^{\delta_2} \{\prod^d_{i=1}(v^{k_1}u+\xi_i
v^{k_1})+\sum_{\alpha,\beta\ge 0}
c_{\alpha,\beta}v^{\alpha}(v^{k_1}u)^{\beta}\}
\qquad \qquad  \\
&=v^{e_{k_1}}(F\circ\tau_{k_1})_{proper} \quad \text{with} \quad
(F\circ\tau_{k_1})_{proper}
=u^{\delta_2}(f\circ\tau_{k_1})_{proper},\\
(f\circ\tau_{k_1})_{proper} &=\prod^d_{i=1}(u+\xi_i)
 +\sum_{\alpha,\beta\ge 0} c_{\alpha,\beta}v^{\alpha+k_1\beta-k_1d}u^{\beta}, \\
(G\circ\tau_{k_1})_{total} &=v^{2k_1}(G\circ\tau_{k_1})_{proper}
\quad \text{with} \quad (G\circ\tau_{k_1})_{proper} =u(u+\xi_1),
 \endalign
 $$
 where
(i) $a=0$ and $b=1$ with $n_1=1$ imply that $bn_1-ak_1=1$,

(ii) $e_{k_1}=\delta_1+k_1{\delta_2}+k_1d=n_1\delta_1
+k_1\delta_2+n_1k_1d$, $\varepsilon=\delta_2$ and
$\varepsilon_{\alpha\beta}=\beta$, as compared with $(3.6.5)$.

Note by assumption in (3.6.1) that
$\alpha+k_1\beta-k_1d=n_1\alpha+k_1\beta-n_1k_1d>0$ because $n_1=1$.
So, the proof of [I] is finished for Case $(A)$. \ms

{\bf The proof of [II] for Case(A):} \ By $(3.6.6)$ and (3.6.12-b),
compare the following:
$$\align
\tau_{m}(v',u')&=(y,z)=({v'}^{n_1}{u'}^{p}, {v'}^{k_1}{u'}^{q}), \tag 3.6.15 \\
\tau_{k_1}(v',u')&=(y,z)=(v'u',{v'}^{k_1}{u'}^s),
\endalign$$
where $m$ is, by definition, a finite number of successive
blow-ups which is needed only to get the standard resolution of
the singular point of $V(G)$. Then, $m=k_1$, $n_1=1$, $p=1$ and
$q=s=k_1-1$ because $k_1=s+1$. So, $pk_1-qn_1=k_1-s=1$.

Since $n_1=1$, then apply (3.6.12-b) to (3.6.1) and (3.6.2),
respectively. Or, we may apply (3.6.10-b) to (3.6.9), if convenient.

Then, along $v'=0$,
$(F\circ\tau_{k_1})_{total}=(F\circ\tau_{k_1})(v',u')$ can be
written in the form
$$\align
(3.6.16) \quad (F\circ\tau_{k_1})_{total}
 &=(v'u')^{\delta_1}({v'}^{k_1}{u'}^s)^{\delta_2}
(\prod^d_{i=1}({v'}^{k_1}{u'}^s+\xi_i(v'u')^{k_1})
 +\sum_{\alpha,\beta\ge 0}c_{\alpha,\beta}(v'u')^{\alpha}({v'}^{k_1}{u'}^s)^{\beta})
 \\
 &={v'}^{e_{k_1}}{u'}^{\varepsilon_1'}(F\circ\tau_{k_1})_{proper}
 \quad \text{with} \quad (F\circ\tau_{k_1})_{proper} =(f\circ\tau_{k_1})_{proper},\\
(f\circ\tau_{k_1})_{proper} &=\prod^d_{i=1}(1+\xi_i u')
 +\sum_{\alpha,\beta\ge 0}
c_{\alpha,\beta}{v'}^{\alpha+k_1\beta-k_1d}{u'}^{\lambda_s},\\
 (G\circ\tau_{k_1})_{total}
&={v'}^{2 k_1}{u'}^{2s}(G\circ\tau_{k_1})_{proper} \quad \text{with}
\quad (G\circ\tau_{k_1})_{proper} =1+\xi_1u',
 \endalign $$
where (i) $p=1$ and $q=s=k_1-1$ imply that $pk_1-qn_1=1$,

(ii) $e_{k_1}=\delta_1+k_1\delta_2+k_1d$ and
 $\ve'_1=\delta_1+s\delta_2+sd=e_s=e_{(k_1-1)}\ge 2$,

(iii) $\lambda_s=\alpha+s\beta-sd>0$ by (4.3.1) or (3.6.9). \ms

From (ii) of $(3.6.16)$, if $d\ge 2$, then $\ve'_1\ge 2$. Also, if
$d=1$, then $\delta_2>0$ by assumption of the theorem, and so
$\ve'_1\ge 2$. Thus, note by (3.6.16) that
$(f\circ\tau_{k_1})_{proper}=(F\circ\tau_{k_1})_{proper}$. As
compared with $(3.6.6)$,
$\ve'=\delta_1+(k_1-1)\delta_2+(k_1-1)d=\delta_1+s\delta_2+sd
=\varepsilon'_1$  and ${\varepsilon'}_{\alpha\beta}
=\alpha+s\beta-sd=\lambda_s$, since $k_1=s+1$. Thus, the proof of
[II] is done for Case $(A)$. \ms

{\bf{The proof of [III] for Case(A):}} For brevity of notation in
the conclusion of Theorem $3.6$, recall that $\widehat{\delta}_2$ is
defined as follows:

If $\delta_2$ is positive, $\widehat{\delta}_2=1$, and if
$\delta_2$ is zero, $\widehat{\delta}_2=0$.

Now, to prove [III] for Case(A), note that $n_1=1$ and $m=k_1=s+1$
as in the proof of [I] and [II] for Case$(A)$ where $\tau_m$ is the
composition of a finite number $m$ of successive blow-ups which is
needed only to get the standard resolution of the singularity of
$V(G)$.

First, $G\in \text{\rm{the type}}[0]$ under $\tau_s$ by Lemma $4.3$,
and after $s+1$ iterations of blow-ups, $E_{s+1}=E_m$ has three
distinct intersection points with other exceptional curves and the
proper transform $V^{(m)}(G)$ by $(3.6.14)$ and $(3.6.16)$, and so
$G\in \text{\rm{the type}}[1]$ under $\tau_m$.

Next, as far as $V(F)$ is concerned, after $s$ iterations of
blow-ups, we proved by [III] of Lemma $4.3$ that any irreducible
component of $\cup^s_{i=1}E_i$ has at most two distinct intersection
points with any other exceptional curves and $V^{(s)}(F)$, and that
no three distinct components of $\cup^s_{i=1}E_i$ meet. Rigorously
speaking, note by Lemma $4.3$ or (3.6.9) that for any $s\ge 2$,
$E_s\cap E_{s-1}\ne \emptyset$, $E_s\cap V^{(s)}(F)$ is a one-point
subset and $E_s\cap E_{s-1}\cap V^{(s)}(F)=\emptyset$, and that for
any $s\ge 1$, $E_s$ and $V^{(s)}(F)$ meet transversely at $E_s\cap
V^{(s)}(F)$, noting that $qs(V^{(s)}(F))=qs(V^{(s)}(G))=E_s\cap
V^{(s)}(F)$ is a one-point subset under $\tau_s$ and that $F\in
\text{\rm{the type}}[0]$, as a reduced variety, and also $G\in
\text{\rm{the type}}[0]$ under $\tau_s$.

Now, after $k_1$ iterations of blow-ups, by $(3.6.14)$ and
$(3.6.16)$ $E_{k_1}=E_{s+1}$ has $(\mu+1+\widehat{\delta}_2)$
distinct intersection points with other exceptional curves and the
proper transform $V^{(k_1)}(F)$, because $E_{k_1}$ and
$V^{(k_1)}(f)$ have $\mu$ distinct intersection points by (ii) of
(3.6.1) in the assumption of this theorem, as a reduced variety. So,
we get that $F\in \text{\rm{the type}} [1]$ if and only if
$\mu+1+\widehat{\delta}_2=3$, as a reduced variety.  In particular,
if $\mu=1$, then $f\in \text{\rm{the type}}[0]$ under $\tau_m$,
assuming that $f$ has a singularity at the origin as a reduced
variety.

 (i) If $d=1$, then $\mu=1$ and also $\delta_2>0$ by assumption
because $n_1=1$. In this case, $F$ belongs to the
$\text{\rm{type}}[1]$ under $\tau_m$, as a reduced variety because
$\mu+1+\widehat{\delta}_2=3$. In particular, if $n_1=d=1$, then it
is clear that $f$ has no singularity at the origin, and so $f$ is
irreducible in $\BC\{y,z\}$.

(ii) Let $d\ge 2$. If $\mu=1$, then it is trivial that $F\in
\text{\rm{the type}}[1]$ under $\tau_m$ if and only if
$\widehat{\delta}_2=1$, as a reduced variety.  If $f$ is
irreducible in $\BC\{y,z\}$, then $\mu=1$, but the converse does
not hold, because $f=(z^2+y^3)^2+y^6z^2$ is reducible in
$\BC\{y,z\}$, for example.

Thus, the proof of [III] is done, and so we finished the proof for
Case$(A)$. \ms

{\bf Case(B):} Let $n_1\ge 2$. Using the same methods and results as
we have seen in the proof of Lemma $4.3$ and Case(A), it can be
shown computationally by induction on the integer $n_1$ that Case(B)
implies the truth of three statements [I], [II] and [III] where
$n_1$ is the multiplicity of $g_1(y,z)=z^{n_1}+\xi_1y^{k_1}$ at the
origin. Thus, the proof of Case(B) is done.

Therefore, the proof of theorem is completely finished. $\square$
\enddemo \ms
\bs

\vfill \pagebreak

{\bf Part[B2] How to construct the Puiseux convergent power
series of the recursive type in $\BC\{y,z\}$(irreducible W-polys of
two complex variables of recursive types in $\BC\{y,z\}$)} \bs

{\bf \S5. To find the necessary and sufficient condition for the
semi-quasi-Puiseux convergent power series in $\BC\{y,z\}$ of the
recursive type to be irreducible in $\C\{y,z\}$} \ms

{\bf \S 5.0. Introduction} \ms

In order to succeed in the computation of Explicit algorithm, in
this section it is very interesting and important to prove that we
can construct the new terminology, irreducible Weierstrass
polynomials of two complex variables of the recursive type, called
{\rm ``the standard Puiseux polynomial in $\BC[y,z]$ of the
recursive type"} throughout this paper, which will be shown to be
equivalent to the standard Puiseux expansion with Explicit Algorithm
of $\S11$, as far as the multiplicity sequences of irreducible plane
curve singularities are concerned.

In preparation, we need Definition $5.0.0$ and Theorem $5.0$. \ms

\definition{Definition 5.0.0} Let $N_0$ be the set of
nonnegative integers and $N^k_0$ be its $k$-dimensional copy. Let
$r$ be an arbitrary positive integer. \ms

\noindent{\bf [A]} $g_r\in\BC\{y,z\}$ is called
$\underline{\text{\rm a semi-quasi-Puiseux convergent power series
of the recursive r-type}}$ if there are sequences
$\{X_k:k=1,2,\dots,r\}$ with $X_k\subset N_0$,
$\{g_k:k=1,2,\dots,r\}$ with $g_k\in\BC\{y,z\}$ and
$\{\text{${{\Delta}_k}:N^k_0\to N_0$ is an integer-valued function
for}$ $\text{$k=1,2,\dots,r$}\}$ satisfying the following
$\underline{\text{\rm four conditions}}$: \ms

\noindent$\underline{\text{\bf Four conditions}}$ are denoted by
\text{\rm The 1st ${\text{\rm{Cond}}}^{\text{{\rm(0)}}}$}, $\dots$,
\text{\rm The 4-th ${\text{\rm{Cond}}}^{\text{{\rm(0)}}}$} for
notation. \ms

\noindent $\underline{\text{\rm The 1st
${\text{\rm{Cond}}}^{\text{{\rm(0)}}}$}}$ Let
$\{X_j:j=1,2,\dots,r\}$ with $X_j\subset N_0$ be defined as follows:

\roster
\item"(1) (1a)"
 $X_1=\{n_1,\beta_{1,1}\}$ with $n_1\ge 2$ and $\beta_{1,1}\ge 1$.
 \item"$\quad$ (1b)" $X_j=\{n_j,\beta_{j,1},\beta_{j,2},\dots,\beta_{j,j}\}$
 with $n_j\ge 2$ \quad where $j=2,\dots,r$.
 \endroster
If $j\ge 2$, then assume that at least one of
$\beta_{j,1},\beta_{j,2},\dots,\beta_{j,j}$ is nonzero. \ms

\noindent $\underline{\text{\rm The 2nd
${\text{\rm{Cond}}}^{\text{{\rm(0)}}}$}}$ For each $j=1,2,\dots,r$,
let $g_j=g_j(y,z)$ be in $\BC\{y,z\}$, each of which is defined by
the following way:
 \roster
 \item"(2) (2a)" $g_1=z^{n_1}+\ve_1y^{\beta_{1,1}}$.
  \item"$\quad$ (2b)"
 $g_j=g^{n_j}_{j-1}+\ve_jy^{\beta_{j,1}}z^{\beta_{j,2}}g^{\beta_{j,3}}_1\cdots
 g^{\beta_{j,j}}_{j-2}$ \quad with $g_{-1}=y$ and $g_0=z$, where $j=2,\dots,r$.
\endroster
Note that each $\ve_i=\ve_i(y,z)$ is a unit in $\BC\{y,z\}$ for
$1\le i\le r$. \ms

\noindent $\underline{\text{\rm The 3rd
${\text{\rm{Cond}}}^{\text{{\rm(0)}}}$}}$ Let $\{\Delta_k: N^k_0\to
N_0: k=1,2,\dots,r\}$ be a sequence such that each $\Delta_k$ is an
integer-valued function defined by the following:

\roster
\item"(3) (3a)" $\Delta_1(t)=t$ for each $t\in N_0$.
\item"$\quad$ (3b)"
$\Delta_j(t_j)^j_{k=1}=t_j\Delta_{j-1}(\beta_{j-1,k})^{j-1}_{k=1}
+n_{j-1}\Delta_{j-1}(t_k)^{j-1}_{k=1}$ for each $(t_k)^j_{k=1}\in
N^j_0$

\noindent where $j=2,\dots,r$.
\endroster \ms

\noindent $\underline{\text{\rm The 4-th
${\text{\rm{Cond}}}^{\text{{\rm(0)}}}$}}$ Then, the following
inequalities hold: Note that $2\le j\le r$. \roster
\item"(4) (4a)" $\Delta_1(\beta_{1,1})=\beta_{1,1}>0$ with $n_1\ge 2$.
\item"$\quad$ (4b)" $\Delta_j(\beta_{j,k})^j_{k=1}
>n_jn_{j-1}\Delta_{j-1}(\beta_{j-1,k})^{j-1}_{k=1}$
where $j=2,\dots,r$.
\endroster \ms

\noindent{\bf [B]} Let $g_r\in \BC\{y,z\}$ be $\underline{\text{\rm
a semi-quasi-Puiseux convergent power series of the recursive
r-type}}$ as in [A].

There are two additional conditions, denoted by \text{\rm The 5-th
${\text{\rm{Cond}}}^{\text{{\bf(0)}}}$} and \text{\rm The 6-th
${\text{\rm{Cond}}}^{\text{{\bf(0)}}}$}.

\noindent $\underline{\text{\rm The 5-th
${\text{\rm{Cond}}}^{\text{{\rm(0)}}}$}}$ \quad The following
inequalities hold:

\roster \item"(5)(5a)" $\gcd(n_j,\Delta_j(\beta_{j,k})^j_{k=1})=1$
for $1\le j\le r$.
\endroster

\noindent $\underline{\text{\rm The 6-th
${\text{\rm{Cond}}}^{\text{{\rm(0)}}}$}}$ \quad The following
inequalities hold: Note that $2\le j\le r$.

\noindent\rm{(6)(6a)}  $2\le n_1<\beta_{1,1}$.

\rm{(6b)} $n_{j}\ge 2$, $\beta_{j,1}>0$, and $0\le
\beta_{j,k}<n_{k-1}$ for $2\le j\le r$ and $2\le k\le j$. \bs

Now, we define the new terminology in [B1], [B2], [B3] and [B4] of
[B].\ms

\noindent{\bf [B1]} $g_r\in \BC\{y,z\}$ is called
$\underline{\text{\rm the quasi-Puiseux convergent power series of
the recursive r-type}}$ in $\BC\{y,z\}$ if $g_r$ in [A] satisfies an
additional condition, denoted by \text{\rm The 5-th
${\text{\rm{Cond}}}^{\text{{\bf(0)}}}$}.

Namely, it is said that $g_r\in \BC\{y,z\}$ is the quasi-Puiseux
convergent power series of the recursive r-type if $g_r$ satisfies
\text{\rm The 1st ${\text{\rm{Cond}}}^{\text{{\rm(0)}}}$}, $\dots$,
\text{\rm The 4-th ${\text{\rm{Cond}}}^{\text{{\rm(0)}}}$},\text{\rm
The 5-th ${\text{\rm{Cond}}}^{\text{{\rm(0)}}}$}. \ms

\noindent{\bf[B2]} $g_r\in \BC\{y,z\}$ is called
$\underline{\text{\rm the Puiseux convergent power series of the
recursive r-type}}$ in $\BC\{y,z\}$ if $g_r$ in [A] satisfies
\text{\rm The 5-th ${\text{\rm{Cond}}}^{\text{{\bf(0)}}}$} and an
inequality in (6a) of \text{\rm The 6-th
${\text{\rm{Cond}}}^{\text{{\bf(0)}}}$}. \ms

\noindent{\bf[B3]} $g_r\in \BC\{y,z\}$ is called
$\underline{\text{\rm the standard Puiseux convergent power series
of the recursive r-type}}$ in $\BC\{y,z\}$ if $g_r$ in [A] satisfies
\text{\rm The 5-th ${\text{\rm{Cond}}}^{\text{{\bf(0)}}}$} and
\text{\rm The 6-th ${\text{\rm{Cond}}}^{\text{{\bf(0)}}}$}. \ms

\noindent{\bf[B4]} Let $g_r\in \BC\{y,z\}$ be the standard Puiseux
convergent power series of the recursive r-type as in [B3]. Then,
$g_r$ is called $\underline{\text{\rm the standard Puiseux
polynomial of two complex variables of}}$ 

\noindent$\underline{\text{\rm the recursive r-type}}$ if each unit
$\ve_i=\ve_i(y,z)$ is equal to an integer one for $1\le i\le r$ in
\text{\rm The 2-th ${\text{\rm{Cond}}}^{\text{{\rm(0)}}}$} of
\text{\rm[A]}.
\enddefinition \ms

\proclaim{Theorem 5.0(To find the necessary and sufficient condition
for the semi-quasi-Puiseux convergent power series in $\BC\{y,z\}$
of the recursive type to be irreducible in $\C\{y,z\}$)}

$\underline{\text{\bf {Assumptions}}}$  Let $g_r$ be
$\underline{\text{\rm a semi-quasi-Puiseux convergent power series
of the recursive r-type}}$, satisfying the same properties and
notations as in {\rm [A]} of {\rm Definition 5.0.0}.\ms

$\underline{\text{\bf {Conclusions}}}$ For each $j=1,2,\dots,r$, let
$(0,0)$ be the singularity of an analytic variety
$V(g_j)=\{(y,z):g_j(y,z)=0\}$ except that
$g_1=z^{n_1}+\ve_1y^{\beta_{1,1}}$ with $\beta_{1,1}=1$. Then, we
get two statements {\rm[A]} and {\rm[B]} as follows:
$$ \split
[A] \qquad \qquad &\text{$g_r$ is irreducible in $\BC\{y,z\}$} \\
  \iff \quad &\text{$g_1,g_2,\dots,g_{r-1}$ are irreducible in $\BC\{y,z\}$ and
   $\gcd(n_r,\Delta_r(\beta_{r,k})^r_{k=1})=1$} \\
\iff \quad &\text{$\gcd(n_1,\beta_{1,1})=1$,
$\gcd(n_2,\Delta_2(\beta_{2,1},\beta_{2,2}))=1$,
$\dots,\gcd(n_r,\Delta_r(\beta_{r,k})^r_{k=1})=1$}. \qquad \qquad \\
\iff \quad &\text{$g_r$ is the quasi-Puiseux convergent power series
of the recursive r-type} \\
\endsplit
$$

\noindent{\rm[B]} Let $g_r$ be irreducible in $\BC\{y,z\}$.

\noindent{\rm[B1]} Let
$V(y^{\gamma}g_r)=\{(y,z):y^{\gamma}g_r(y,z)=0\}$ be an analytic
variety at $(0,0)$ in $\BC^2$ defined by
$$\align
  y^{\gamma}g_r(y,z) \quad \text{such
that} \quad \left\{\eqalign{& \text{$\gamma=1$, \quad if \quad
$\beta_{1,1}=1$}, \cr & \text{$\gamma=0$, \quad if
\quad $\beta_{1,1}>1$.} \cr} \right.  \tag 5.0.0 \\
\endalign$$
Then, $y^{\gamma}g_r\in$ the type $[r]$ under the standard
resolution, denoted by $\tau$, in the sense of Definition $2.5$.
Also, if $\beta_{1,1}=1$ then $g_r\in$ the type $[r-1]$ under the
standard resolution. \ms

\noindent{\rm[B2]} In particular, $z^{\delta}yg_r\in$ the type $[r]$
under the same standard resolution $\tau$, whether $\delta$ is
either one or zero. $\square$
\endproclaim \ms

{\bf \S5.1. In preparation for the proof of Theorem 5.0 } \ms

For the proof of Theorem $5.0$, in this section, we will prepare the
statements of five sublemmas, consisting of Sublemma 5.1, Sublemma
5.2, ..., Sublemma 5.5 without proofs. After then, we will finish
the proofs of these sublemmas and Theorem 5.0 in \S6.

Moreover, as a corollary of the above theorem and sublemmas, we can
easily get Corollary $5.6$ and Corollary $5.7$ with no need of
proofs. \ms

\proclaim{Sublemma 5.1} $\underline{\text{\bf Assumptions}}$ \quad
Suppose that the same properties and notations as in the assumptions
of Theorem $5.0$ hold.

For any integer $r\ge 2$, let
$\Delta^{\sharp}_2(\beta_{2,1},\beta_{2,2})$ and
$\Delta^{\sharp}_j(\beta_{j,k})^j_{k=1}$ with $3\le j\le r$ be the
notations defined as follows : Note that
$\Delta_2(t_1,t_2)=n_1t_1+\beta_{1,1}t_2$ for each $(t_1,t_2)\in
N^2_0$.
$$\align
  \Delta^{\sharp}_2(\beta_{2,1},\beta_{2,,2})
&=\Delta_2(\beta_{2,1},\beta_{2,2}),  \tag 5.1.1\\
  \Delta^{\sharp}_j(\beta_{j,k})^j_{k=1}
 &=\Delta_2(\beta_{j,1},\beta_{j,2})+n_1\beta_{1,1}\beta_{j,3}
 +n_1\beta_{1,1}n_2\beta_{j,4} \\
 &\quad +n_1\beta_{1,1} n_2 n_3 \beta_{j,5}
 + \cdots +n_1\beta_{1,1}n_2 \cdots n_{j-2}\beta_{j,j}.
\endalign$$

$\underline{\text{\bf Conclusions}}$ \quad  Then, we have the
following:
$$\align
\Delta^{\sharp}_2(\beta_{2,1},\beta_{2,2})
&> n_1\beta_{1,1}n_2 \quad \text{on $g_2$},  \tag 5.1.2\\
\Delta^{\sharp}_j(\beta_{j,k})^j_{k=1} &> n_1\beta_{1,1}n_2n_3\cdots
n_j \quad \text{on $g_j$}.\quad \text{$\square$}
\endalign
$$
\endproclaim
\ms

\proclaim{Sublemma 5.2} $\underline{\text{\bf Assumptions}}$ \quad
Suppose that the same properties and notations as in the assumptions
of Theorem $5.0$ hold. Let $r$ be any integer with either $r\ge 2$
or $r=1$.

$\underline{\text{\bf Conclusions}}$ \quad Then, we get the
following:

{\rm(a)} For any $r\ge 2$, $g_r=g_r(y,z)$ can be written in the
form
$$\align
(5.2.1) \qquad & g_r=(z^{n_1}+\ve_1 y^{\beta_{1,1}})^{n_2n_3\cdots
n_r}+\sum_{\alpha,\beta\ge 0}
c^{(r)}_{\alpha,\beta}y^{\alpha}z^{\beta} \quad \text{with $\ve_1=1$
and}
\qquad \qquad \qquad \\
& \quad \text{with
$n_1\alpha+\beta_{1,1}\beta>n_1\beta_{1,1}n_2n_3\cdots n_r$},
\endalign$$
where a unit $\ve_1=\ve_1(y,z)$ may be analytically assumed to be
one in $\BC\{y,z\}$, and the $c^{(r)}_{\alpha,\beta}$ are nonzero
complex numbers for some nonnegative integers $\alpha$ and $\beta$.
\ms

{\rm(b)} For each $r\ge 2$, we have the following:

{\rm(b1)} The multiplicity of $g_r(0,z)$ at $z=0$ is
$n_1\prod^r_{k=2}n_k$ when $g_r=g_r(y,z)$.

{\rm(b2)} The multiplicity of $g_r(y,0)$ at $y=0$ is
 $\beta_{11}\prod^r_{k=2}n_k$ when $g_r=g_r(y,z)$. \ms

{\rm(c)} For each $r\ge 2$, we have the following:

{\rm(c1)} If $n_1<\beta_{1,1}$ then $\alpha+\beta>n_1n_2\cdots n_r$,
and so the multiplicity of $g_r$ at $(y,z)=(0,0)$ is
$n_1\prod^r_{k=2}n_k$.

{\rm(c2)} If $n_1>\beta_{1,1}$ then
$\alpha+\beta>\beta_{1,1}n_2n_3\cdots n_r$, and so the multiplicity
of $g_r$ at $(y,z)=(0,0)$ is $\beta_{11}\prod^r_{k=2}n_k$. \ms

{\rm(d)} If $g_r$ is irreducible in $\BC\{y,z\}$ for any $r\ge 2$,
then either $\gcd(n_1,\beta_{1,1})=1$ or $g_1$ in \text{\rm The 2-th
${\text{\rm{Cond}}}^{\text{{\rm(0)}}}$} is irreducible in
$\BC\{y,z\}$.

{\rm(e)} In particular, from {\rm (5.2.1)} let $h_1=h_1(y,z)$ be
defined in the form
$$\align
(5.2.1.1) \qquad \qquad & h_1=(z^{n_1}+\ve_1
y^{\beta_{1,1}})+\sum_{\alpha,\beta\ge 0}
c^{(1)}_{\alpha,\beta}y^{\alpha}z^{\beta} \quad \text{with
$\ve_1=1$ and} \qquad \qquad \qquad  \quad \\
& \quad \text{with $n_1\alpha+\beta_{1,1}\beta>n_1\beta_{1,1}$},
\endalign$$
where a unit $\ve_1=\ve_1(y,z)$ may be analytically assumed to be
one in $\BC\{y,z\}$, and the $c^{(1)}_{\alpha,\beta}$ are nonzero
complex numbers for some nonnegative integers $\alpha$ and $\beta$,
if exist.

Then, $h_1$ satisfies the same kind of results as $g_r$ of
\text{\rm(a)} does in {\rm(b)}, {\rm(c)} and {\rm(d)}. $\square$
\endproclaim
\ms

\proclaim{Sublemma 5.3} $\underline{\text{\bf Assumptions}}$ \quad
Suppose that the same properties and notations as in the assumptions
of Theorem $5.0$ hold. Let $r$ be an arbitrary integer with $r\ge
2$.

In addition, we need the following assumptions:

\noindent {\rm(5.3.0)} \qquad \quad $\gcd(n_1,\beta_{1,1})=1$ or
$g_1$ is irreducible in $\BC\{y,z\}$.

Since $\gcd(n_1,\beta_{1,1})=1$ with $n_1\ge 2$ and $\beta_{1,1}\ge
1$, then there are two nonnegative integers $a>0$ and $b\ge 0$ such
that $a\beta_{1,1}-bn_1=1$.

For given two integers $a>0$ and $b\ge 0$, let $\Omega_2:N^2_0\to
N_0$ be a function defined by
$$\align
\Omega_2(t_1,t_2)=at_1+bt_2. \tag 5.3.1
\endalign$$

Let $\Omega^{\sharp}_2(\beta_{2,1},\beta_{2,2})$ and
$\Omega^{\sharp}_j(\beta_{j,k})^j_{k=1}$ with $3\le j\le r$ be the
notations defined as follows:
$$\align
(5.3.2) \qquad  \Omega^{\sharp}_2(\beta_{2,1},\beta_{2,2})
 &=\Omega_2(\beta_{2,1},\beta_{2,2}).   \\
  \Omega^{\sharp}_j(\beta_{j,k})^j_{k=1} &=\Omega_2(\beta_{j,1},\beta_{j,2})
 +bn_1\beta_{j,3}+bn_1n_2\beta_{j,4}
+\cdots +bn_1n_2\cdots n_{j-2}\beta_{j,j}. \qquad \qquad
\endalign
$$

$\underline{\text{\bf Conclusions}}$ \quad Then, we get the
following:
 $$
\align
 \Omega^{\sharp}_2(\beta_{2,1},\beta_{2,2}) &\ge bn_1n_2,  \tag 5.3.3 \\
 \Omega^{\sharp}_j(\beta_{j,k})^j_{k=1} &\ge bn_1n_2n_3\cdots n_j.
 \quad \text{$\square$}
 \endalign $$
\endproclaim
\ms

\proclaim{Sublemma 5.4} $\underline{\text{\bf Assumptions}}$ \quad
Suppose that the same properties and notations as in the assumptions
of Theorem $5.0$ hold. As in Sublemma $5.3$, additionally assume
that we have the following properties:
$$\align
(\ast1) \quad & \text{$g_1$ is irreducible in $\BC\{y,z\}$ or
$\gcd(n_1, \beta_{1,1})=1$}. \\
(\ast2)  \quad & \text{$g_{r}$ may not be irreducible in
$\BC\{y,z\}$ for some $r\ge 2$, but note by \text{\rm The 4-th
${\text{\rm{Cond}}}^{\text{{\rm(0)}}}$} in} \\
& \text{the assumption of this theorem that
$\Delta_{r}(\beta_{r,k})^{r}_{k=1}
>n_{r}n_{r-1}\Delta_{r-1}(\beta_{r-1,k})^{r-1}_{k=1}$.}  \\
\endalign$$

$\underline{\text{\bf {Conclusions}}}$ \quad For each
$j=1,2,\dots,r$, let $V(g_j)=\{(y,z):g_j(y,z)=0\}$  be an analytic
variety at the origin in $\BC^2$. For the construction of the
statement of the conclusion, let $V(G)=\{(y,z):G(y,z)=0\}$ be
another analytic variety with isolated singularity at the origin in
$\BC^2$ defined by the form
$$
\align g_1 &=z^{n_1}+\ve_{1}y^{\beta_{1,1}}
\quad \text{with a unit  $\ve_{1}\in \C\{y,z\}$,} \tag 5.4.0 \\
G &=y^{\gamma}g_1,
\endalign
$$
satisfying the properties {\rm(i)} and {\rm(ii)}:

\roster \item "(i)" If $\beta_{1,1}=1$, then $\gamma=1$.

\item "(ii)" If $\beta_{1,1}\ge 2$, then $\gamma=0$.
\endroster \ms

Let $\tau_m=\pi_1\circ\pi_2\circ\cdots\circ\pi_m:M^{(m)}\to\BC^2$ be
the compositions of a finite number $m$ of successive blow-ups
$\pi_i$ which is needed to get the standard resolution of the
singular point of $V(y^{\gamma}g_1)$. If $V(g_1)$ has the singular
point at the origin, then as compared with the above $\tau_m$,
exactly the same $\tau_m$ can be also used for the standard
resolution of the singular point of $V(yg_1)$ as far as the blow-ups
process is concerned.

{\rm(a)(a1)} We can use just one coordinate patch of the local
coordinates for each blow-up $\pi_i$ of $\tau_m$ with $1\le i\le
m$ in the sense of Lemma $2.12$.

{\rm \quad(a2)} Just as above, we can use the same $\tau_m$ for
the composition of the first finite number $m$ of successive
blow-ups in preparation for the standard resolution of the
singular point $(0,0)$ of $V(g_j)$ for all $j=2,3,\dots,r$.

{\rm \quad(a3)} Also, we can use just the common one coordinate
patch of the given local coordinates for each blow-up $\pi_i$ of
the above $\tau_m$ in {\rm (a1)}, in order to study any of
$V^{(i)}(g_j)$ for all $j=2,3,\dots,r$ and all $i=1,2,\dots,m$ in
the sense of Lemma $2.14$. \ms

{\rm(b)} For brevity of notations, let $(v,u)$ be the common one of
the local coordinates for the $m-th$ blow-up $\pi_m:M^{(m)}\to
M^{(m-1)}$ at $(0,0)$ which is the quasisingular point of
$V^{(m-1)}(y^{\gamma}g_1)$ in the sense of Definition $2.6$. Being
viewed as an analytic mapping, $\tau_m:M^{(m)}\to\BC^2$ can be
written in the form
 $$
 \tau_m(v,u)=(y,z)=(v^{n_1}u^a,v^{\beta_{1,1}}u^b), \tag 5.4.1
 $$
where

$(b_1)$ $a>0$ and $b\ge 0$ are some nonnegative integers such that
$a\beta_{1,1}-bn_1=1$,

$(b_2)$ $E_m=\{v=0\}$ is defined by the $m-th$ exceptional curve
of the first kind. \ms

{\rm(c)} Use the same notations for
$\Delta^{\sharp}_q(\beta_{q,k})^q_{k=1}$ and
$\Omega^{\sharp}_q(\beta_{q,k})^q_{k=1}$ as in Sublemma $5.1$, and
Sublemma $5.3$ where $\Omega_2:N^2_0\to N_0$ is a function defined
by $\Omega_2(t_1,t_2)=at_1+bt_2$  for given two nonnegative integers
a and b in $(b_1)$ of {\rm(5.4.1)}, and we may start with assuming
that $\varepsilon_1=1$ in
$V(y^{\gamma}g_1)=\{y^{\gamma}(z^{n_1}+\varepsilon_1y^{\beta_{1,1}})=0\}$,
in order to study $V^{(i)}(g_j)$ for all $i=1,2,\dots,m,$ and all
$j=1,2,\dots,r$. Whether $\beta_{1,1}\ge 2$ or $\beta_{1,1}=1$, we
may write that $(g_1\circ\tau_m)_{proper}=(1+\ve_1u)$ with
$\ve_1=1$, without complexity of the notation if necessary, noting
that if $\beta_{1,1}=1$ then $V(g_1)$ has no singularity at the
origin. \ms

Now, along $v=0$, $(g_j\circ\tau_m)_{total}$ with
$(g_j\circ\tau_m)_{proper}$ can be written as follows:
$$\align
(5.4.2) \quad ((y^{\gamma}g_1)\circ\tau_m)_{total}
&=v^{(\gamma+\beta_{1,1})n_1}u^{bn_1}
((y^{\gamma}g_1)\circ\tau_m)_{proper} \quad\text{with}  \\
((y^{\gamma}g_1)\circ\tau_m)_{proper} &=u^{a\gamma}(1+u),\\
 (g_1\circ\tau_m)_{total} &=v^{n_1\beta_{1,1}}u^{bn_1}
 (g_1\circ\tau_m)_{proper}\quad\text{with}  \\
 (g_1\circ\tau_m)_{proper} &=(1+u),\\
  (g_j\circ\tau_m)_{total} &=v^{n_1\beta_{1,1}n_2\cdots
 n_j}u^{bn_1n_2\cdots n_j}
 (g_j\circ\tau_m)_{proper}\quad\text{with}\\
 (g_j\circ\tau_m)_{proper} &=(g_{j-1}\circ\tau_m)^{n_j}_{proper}
 +\{\ve'_j v^{\Delta^{\sharp}_j(\beta_{j,k})^j_{k=1}-n_1\beta_{1,1}n_2\cdots
 n_j}{\times} \\
 & \quad u^{\Omega^{\sharp}_j(\beta_{j,k})^j_{k=1}-bn_1n_2\cdots n_j}
 (1+u)^{\beta_{j,3}} (g_2\circ\tau_m)^{\beta_{j,4}}_{proper}
 \cdots (g_{j-2}\circ\tau_m)^{\beta_{j,j}}_{proper}\},
\endalign
$$where each $\ve'_j=\ve_j\circ\tau_m(v,u)$ is a unit in
$\BC\{v,1+u\}$ for $2\le j\le r$.

Note by Sublemma $5.1$ and Sublemma $5.3$ that for
$j=2,3,\dots,r$,
$$\align
\Delta^{\sharp}_j(\beta_{j,k})^j_{k=1} &> n_1\beta_{1,1}n_2\cdots
n_j
\quad \text {and} \tag 5.4.3 \\
\Omega^{\sharp}_j(\beta_{j,k})^j_{k=1} &\ge bn_1n_2n_3\cdots n_j.
\endalign$$

Moreover,
$(y^{\delta_{r,1}}z^{\delta_{r,2}}g^{\delta_{r,3}}_1g^{\delta_{r,4}}_2\cdots
g^{\delta_{r,r}}_{r-2})\circ\tau_m(v,u)$ can be viewed as
$$
u^{\Omega^{\sharp}_r(\delta_{r,k})^r_{k=1}}v^{\Delta^{\sharp}_r(\delta_{r,k})^r_{k=1}}
(1+u)^{\delta_{r,3}}(g_2\circ\tau_m)^{\delta_{r,4}}_{proper} \cdots
(g_{r-2}\circ\tau_m)^{\delta_{r,r}}_{proper}, \tag 5.4.4
$$ where

$(c_1)$ the $\delta_{r,i}$ are nonnegative integers for $1\le i\le
r$,

$(c_2)$
$\Delta^{\sharp}_r(\delta_{r,k})^r_{k=1}=\Delta_2(\delta_{r,1},\delta_{r,2})
 +n_1\beta_{1,1}\delta_{r,3}+n_1\beta_{1,1}n_2\delta_{r,4}+\cdots
 +n_1\beta_{1,1}n_2\cdots n_{r-2}\delta_{r,r}$ as in the definition of
 $\Delta^{\sharp}_r(\beta_{r,k})^r_{k=1}$ of Sublemma $5.1$,

$(c_3)$ $\Omega^{\sharp}_r(\delta_{r,k})^r_{k=1}
 =\Omega_2(\delta_{r,1},\delta_{r,2})+bn_1\delta_{r,3}
 +bn_1n_2\delta_{r,4}+\cdots +bn_1n_2\cdots n_{r-2}\delta_{r,r}$
 as in the definition of $\Omega^{\sharp}_r(\beta_{r,k})^r_{k=1}$ of Sublemma $5.3$. \ms

{\rm(d)} Note again that $\tau_m$ is the composition of a finite
number $m$ of successive blow-ups, which is needed to get the
standard resolution of the singular point of $V(y^{\gamma}g_1)$ or
$V(g_1)$. Let ${\tau_m}^{-1}(0,0)=\cup^m_{i=1}E_i$ where $E_i$ is
an exceptional curve of the first kind under $\tau_m$. For
$j=1,2,\dots,r$, let
$$(g_j\circ\tau_m)_{divisor}=V^{(m)}(g_j)+\sum^m_{i=1}e_{j,i}E_i, \tag 5.4.5 $$
where $V^{(m)}(g_j)$ is the proper transform of $V(g_j)$ under
$\tau_m$.

Then we have the following:

{\rm(d1)} If $\beta_{1,1}\ge 2$, then $e_{j+1,i}=n_{j+1}e_{j,i}$ for
any $j\ge 1$ and for $i=1,2,\dots,m$.

\qquad If $\beta_{1,1}=1$, then $e_{j+1,i}=n_{j+1}e_{j,i}$ for any
$j\ge 2$ and  for $i=1,2,\dots,m$.

In particular, $e_{j,m}=n_1\beta_{1,1}n_2\cdots n_j$ for
$j=2,\dots,r$, and $e_{1,m}=n_1\beta_{1,1}$, if exists. \ms

{\rm(d2)} $V^{(m)}(g_j)\cap(\cup^m_{i=1}E_i)=V^{(m)}(g_j)\cap
E_m=\{(v,1+u)=(0,0)\}$ for any $j=2,\dots,r$. \ms

{\rm(d3)} If $\beta_{1,1}\ge 2$, then for any $j=1,2,\dots,r$,
$g_j\in$ the type $[1]$ under $\tau_m$.

\qquad If $\beta_{1,1}=1$, then for any $j=1,2,\dots,r$, $g_j\in$
the type $[0]$ under $\tau_m$.

In particular, if $\beta_{1,1}\ge 1$, note that for all
$j=1,2,\dots,r$, $z^{\delta}yg_j\in$ the type $[1]$ under $\tau_m$
whether $\delta=1$ or $\delta=0$, by Theorem $3.6$. \quad $\square$
\endproclaim
\ms

\proclaim{Sublemma 5.5} $\underline{\text{\bf {Assumptions}}}$
Suppose that the same properties and notations as in the assumptions
of Theorem $5.0$ hold. As in Sublemma $5.3$, additionally assume
that we have the following properties:

\noindent {\rm(5.5.0)} \qquad \quad $\gcd(n_1,\beta_{1,1})=1$ or
$g_1$ is irreducible in $\BC\{y,z\}$.

Let $r$ be an arbitrary positive integer with $r\ge 2$. Throughout
this sublemma, we will use the same notations and consequences as
in Sublemma $5.4$, in order to get the representation for the
conclusion of this sublemma. \ms

$\underline{\text{\bf {Conclusions}}}$ \quad As $\{g_k:
k=1,2,\dots,r\}$ with $g_k\in \BC\{y,z\}$ satisfies four conditions
in the assumptions of Theorem $5.0$, denoted by \text{\bf The 1st
${\text{\bf{Cond}}}^{\text{{\bf(0)}}}$}, \dots, \text{\bf The 4-th
${\text{\bf{Cond}}}^{\text{{\bf(0)}}}$}, then
$\{(g_k\circ\tau_m)_{proper}: k=2,3,\dots,r\}$ with
$(g_k\circ\tau_m)_{proper}$ in $\BC\{v,1+u\}$ satisfies the same
kind of four conditions, which will be denoted by \text{\bf The 1st
${\text{\bf{Cond}}}^{\text{{\bf(1)}}}$}, \dots, \text{\bf The 4-th
${\text{\bf{Cond}}}^{\text{{\bf(1)}}}$}. Note that
$\{(g_k\circ\tau_m)_{proper}: k=2,3,\dots,r\}$ has been already
well-defined by Sublemma $5.4$. \ms

In more detail, in order to construct four conditions recursively,
which will be denoted by \text{\bf The 1st
${\text{\bf{Cond}}}^{\text{{\bf(1)}}}$}, \dots, \text{\bf The 4-th
${\text{\bf{Cond}}}^{\text{{\bf(1)}}}$}, we prefer to add one more
condition to the above four conditions, denoted by \text{\bf The
5-th ${\text{\bf{Cond}}}^{\text{{\bf(1)}}}$}, for convenience of
representation.  By using the same kind of properties and notations
as in Theorem $5.0$, the desired construction is as follows:
$$\align
\text{There are}&\text{{\quad}sequences:} \\
& \{Y_k: k=1,2,\dots,r-1\} \quad
\text{with $Y_k\subset N_0$}, \\
& \{h_k: k=1,2,\dots,r-1\} \quad \text{with
$h_k=(g_{k+1}\circ\tau_m)_{proper}$ in $\BC\{v,u+1\}$ and} \qquad \qquad \\
&\{\text{$\Xi_k:N^k_0\to N_0$ is an integer-valued
function for $k=1,2,\dots,r-1$}\}, \\
& \text{satisfying the following four conditions for each k}:
\endalign$$

Such conditions are denoted by \text{\bf The 1st
${\text{\bf{Cond}}}^{\text{{\bf(1)}}}$}{\bf, \dots,} \text{\bf The
4-th ${\text{\bf{Cond}}}^{\text{{\bf(1)}}}$}. \ms

\noindent \text{\bf The 1st ${\text{\bf{Cond}}}^{\text{{\bf(1)}}}$:}
Let $\{Y_k:k=1,2,\dots,r-1\}$ with $Y_k\subset N_0$ be defined by
$$\align
(1a) \qquad \qquad \qquad  Y_1 &=\{s_1,\gamma_{1,1}\} \qquad \qquad
\qquad \text{with $s_1\ge 2$ and
$\gamma_{1,1}\ge 1$}, \\
\text{\rm (1b)} \qquad \qquad \qquad Y_{j}
&=\{s_{j},\gamma_{j,1},\gamma_{j,2},\dots,
 \gamma_{j,j}\} \quad \text{with $s_{j}\ge 2$}, \quad \text{where
 $2\le j\le r-1$}. \qquad \qquad
 \endalign$$
such that for each $j=1,2,\dots,r-1$,
$$\align
(5.5.1) \qquad \qquad  e_{1,m}&=n_1\Delta_1(\beta_{1,1})=n_1\beta_{1,1}, \\
  s_{j} &=n_{j+1}\ge 2,
 ~\gamma_{j,1}=\Delta^{\sharp}_{j+1}(\beta_{{j+1},k})^{j+1}_{k=1}
 -n_1\beta_{1,1}n_2n_3\cdots n_{j+1}>0, \qquad \qquad\\
  \gamma_{j,2}&=\beta_{{j+1},3},\quad \gamma_{j,3}=\beta_{{j+1},4},\quad \dots ,
\gamma_{j,j}=\beta_{{j+1},{j+1}},
\endalign
$$
noting that $\gamma_{1,1},\gamma_{2,1},\dots,\gamma_{r-1,1}$ are
 positive by Sublemma $5.1$. \ms

\noindent \text{\bf The 2nd ${\text{\bf{Cond}}}^{\text{{\bf(1)}}}$:}
Let $(g_2\circ\tau_m)_{proper},(g_3\circ\tau_m)_{proper},\ldots,
(g_r\circ\tau_m)_{proper}$ be denoted by $h_1,h_2,\dots,h_{r-1}$,
respectively in $\BC\{v,u+1\}$ as follows:
$$\align
(2a)\qquad \qquad
&(g_2\circ\tau_m)_{total}=v^{n_2e_{1,m}}(g_2\circ\tau_m)_{proper}
\quad \text{with}    \\
\qquad \qquad &(g_2\circ\tau_m)_{proper}=(u+1)^{s_1}+\eta_1
v^{\gamma_{1,1}}
\qquad \text{with { } $\eta_1=1$,} \qquad \qquad \\
\text{\rm(2b)} \qquad \qquad &(g_j\circ\tau_m)_{total}
=v^{n_jn_{j-1}\cdots n_2e_{1,m}}(g_j\circ\tau_m)_{proper}
\quad \text{with}\\
\quad  &(g_j\circ\tau_m)_{proper}
 =(g_{j-1}\circ\tau_m)^{s_{j-1}}_{proper}
+\eta_{j-1}v^{\gamma_{j-1,1}}(1+u)^{\gamma_{j-1,2}} \qquad \qquad \qquad \quad\\
 & \qquad \qquad \qquad \quad \times(g_2\circ\tau_m)^{\gamma_{j-1,3}}_{proper}\cdots
 (g_{j-2}\circ\tau_m)^{\gamma_{j-1,j-1}}_{proper},
 \endalign$$
where $\eta_i=\eta_i(v,u+1)$ is a unit in $\BC\{v,u+1\}$ for $1\le
i\le r-1$, noting that $\eta_i=\ve'_{i+1}
u^{\Omega^{\sharp}_{i+1}(\beta_{i+1,k})^{i+1}_{k=1}-bn_1n_2\cdots
n_{i+1}}$. Here, we may assume by a nonsingular change of
coordinates that $\eta_1$ can be equal to an integer one for the
standard resolution of the quasisingular point $(v,u+1)=(0,0)$ of
$V((g_k\circ\tau_m)_{proper})$ for $2\le k\le r$. \ms

\noindent \text{\bf The 3rd ${\text{\bf{Cond}}}^{\text{{\bf(1)}}}$:}
Let $\{\text{$\Xi_k:N^k_0\to N_0$ is an integer-valued function for
{\rm k=1,2,\dots,r-1}}\}$ be a sequence defined by the following:
$$\align
(3a)\quad \quad  &\text{$\Xi_1(t)=t$  for each $t\in N_0$.} \\
\text{\rm(3b)} \quad \quad
&\text{$\Xi_{j-1}(t_k)^{j-1}_{k=1}=t_{j-1}\Xi_{j-2}
(\gamma_{j-2,k})^{j-2}_{k=1}+s_{j-2}\Xi_{j-2}(t_k)^{j-2}_{k=1}$
for each $(t_k)^{j-1}_{k=1}\in  N^{j-1}_0.$}   \qquad \\
\endalign$$ \ms

\noindent \text{\bf The $(4\alpha)$-th
${\text{\bf{Cond}}}^{\text{{\bf(1)}}}$:} For each
$q=1,2,3,\dots,{r-1}$, we have the following:
$$\align
(5.5.4\alpha) \qquad & \Xi_1(\gamma_{1,1})=\gamma_{1,1}
=\Delta_2(\beta_{2,1},\beta_{2,2})-n_1\beta_{1,1}n_2>0, \\
& \Xi_q(\gamma_{q,k})^q_{k=1}
-s_qs_{q-1}\Xi_{q-1}(\gamma_{q-1,k})^{q-1}_{k=1} \\
& \quad  =\Delta_{q+1}(\beta_{q+1,k})^{q+1}_{k=1}
-n_{q+1}n_q\Delta_q(\beta_{q,k})^q_{k=1}>0 \quad \text{for $2\le
q\le {r-1}$.} \qquad \qquad \\
\endalign$$

\noindent \text{\bf The 4-th
${\text{\bf{Cond}}}^{\text{{\bf(1)}}}$:} For each
$q=1,2,3,\dots,{r-1}$, we have the following:
$$\align
(5.5.4) \qquad \qquad \Xi_1(\gamma_{1,1})&=\gamma_{1,1}>0, \\
\Xi_q(\gamma_{q,k})^q_{k=1}&>s_qs_{q-1}\Xi_{q-1}(\gamma_{q-1,k})^{q-1}_{k=1}
\quad \text{for each $q=2,3,\dots,{r-1}$.} \qquad \qquad \qquad
\endalign$$
\ms

\noindent \text{\bf The $(5\alpha)$-th
${\text{\bf{Cond}}}^{\text{{\bf(1)}}}$:} For each
$j=1,2,3,\dots,{r-1}$, we have the following:
$$\align
&\gcd(s_j,\Xi_j(\gamma_{j,k})^j_{k=1})
=\gcd(n_{j+1},\Delta_{j+1}(\beta_{j+1,k})^{j+1}_{k=1}). \tag
5.5.$5\alpha$
\endalign$$
\endproclaim \ms

\noindent{\bf Remark 5.5.1.} For brevity of the proof of {\bf
Sublemma 5}, suppose that \text{\bf The 1st
${\text{\bf{Cond}}}^{\text{{\bf(1)}}}$}, \text{\bf The 2nd
${\text{\bf{Cond}}}^{\text{{\bf(1)}}}$} and \text{\bf The 3rd
${\text{\bf{Cond}}}^{\text{{\bf(1)}}}$} have proved in the
conclusions of this sublemma. Then, for the remaining proof, it
suffices to show that \noindent \text{\bf The $(4\alpha)$-th
${\text{\bf{Cond}}}^{\text{{\bf(1)}}}$} is true, because of the
following two facts:

{\rm Fact(1):} If \text{\bf The $4\alpha$-th
${\text{\bf{Cond}}}^{\text{{\bf(1)}}}$} is true, then it is clear
that \text{\bf The 4-th ${\text{\bf{Cond}}}^{\text{{\bf(1)}}}$} is
true. \ms

{\rm Fact(2):} If \text{\bf The $4\alpha$-th
${\text{\bf{Cond}}}^{\text{{\bf(1)}}}$} is true, then we can easily
prove by \text{\bf The 1-th ${\text{\bf{Cond}}}^{\text{{\bf(1)}}}$}
that \text{\bf The $5\alpha$-th
${\text{\bf{Cond}}}^{\text{{\bf(1)}}}$} is true, by using the
following elementary computation:

{\rm(i)}
$\gcd(s_1,\gamma_{1,1})=\gcd(n_2,\Delta_2(\beta_{2,1},\beta_{2,2})-n_1\beta_{1,1}n_2)=
\gcd(n_2,\Delta_2(\beta_{2,1},\beta_{2,2}))$.

{\rm(ii)} For each $j=2,3,\dots,r-1$, $s_{j}=n_{j+1}$, and
$$\split
(5.5.5\alpha) \qquad
&\gcd(s_j,\Xi_j(\gamma_{j,k})^j_{k=1})=\gcd(s_j,\Xi_j(\gamma_{j,k})^j_{k=1}
-s_js_{j-1}\Xi_{j-1}(\gamma_{j-1,k})^{j-1}_{k=1}) \qquad \qquad \\
&=\gcd(n_{j+1},\Delta_{j+1}(\beta_{j+1,k})^{j+1}_{k=1}
-n_{j+1}n_j\Delta_j(\beta_{j,k})^j_{k=1}) \quad \text{by
$(5.5.4\alpha)$}\\
&=\gcd(n_{j+1},\Delta_{j+1}(\beta_{j+1,k})^{j+1}_{k=1}).\text{\quad
$\square$}
\endsplit$$ \ms

\proclaim{Corollary 5.6} $\underline{\text{\bf Assumptions}}$ Let
$g_r\in \BC\{y,z\}$ be $\underline{\text{\rm a semi-quasi-Puiseux
convergent power}}$ \noindent $\underline{\text{\rm series of the
recursive r-type}}$, as either in {\rm[A]} of {\rm Definition 5.0.0}
or in the assumption of {\rm Theorem 5.0}. \ms

$\underline{\text{\bf Conclusions}}$ \quad For any $r\ge 2$,
$g_r=g_r(y,z)$ can be written in the form
$$\align
 g_r=(z^{n_1}+\ve_1 y^{\beta_{1,1}})^{n_2n_3\cdots
n_r}+\sum_{\alpha,\beta\ge 0}
c^{(r)}_{\alpha,\beta}y^{\alpha}z^{\beta} \quad \text{with
$\ve_1=1$,} \tag 5.6.1
\endalign$$
where a unit $\ve_1=\ve_1(y,z)$ may be analytically assumed to be
one in $\BC\{y,z\}$, and the $c^{(r)}_{\alpha,\beta}$ are nonzero
complex numbers for some nonnegative integers $\alpha$ and $\beta$
such that $n_1\alpha+\beta_{1,1}\beta>n_1\beta_{1,1}n_2n_3\cdots
n_r$.
\endproclaim \ms

\proclaim{Corollary 5.7} $\underline{\text{\bf {Assumptions}}}$ As
we have seen in {\rm Definition 5.0.0}, let $f\in \BC\{y,z\}$ be
$\underline{\text{\rm a semi-quasi-Puiseux convergent power series
of the of the recursive r-type}}$. \ms

$\underline{\text{\bf {Conclusions}}}$ Then, $f$ is irreducible in
$\BC\{y,z\}$ if and only if $f$ is the quasi-Puiseux convergent
power series of the recursive r-type in $\BC\{y,z\}$.
\endproclaim \ms

\proclaim{Corollary 5.8} $\underline{\text{\bf {Assumptions}}}$ Let
$g_r\in \BC\{y,z\}$ be $\underline{\text{\rm a semi-quasi-Puiseux
convergent power}}$ \noindent $\underline{\text{\rm series of the
recursive r-type}}$, as either in {\rm[A]} of {\rm Definition 5.0.0}
or in the assumption of {\rm Theorem 5.0}. In addition, assume that
$g_r$ is the quasi-Puiseux convergent power series of the recursive
r-type in $\BC\{y,z\}$. Let
$\tau_m=\pi_1\circ\pi_2\circ\cdots\circ\pi_m:M^{(m)}\to\BC^2$ be the
compositions of a finite number $m$ of successive blow-ups $\pi_i$
which is needed to get the standard resolution of the singular point
of $V(y^{\gamma}g_1)$. \ms

$\underline{\text{\bf {Conclusions}}}$ Then,
$(g_{r+1}\circ\tau_m)_{proper}$ in $\BC\{v,1+u\}$ is

$\underline{\text{\rm a quasi-Puiseux convergent power series of the
recursive r-type}}$ satisfying $\underline{\text{\rm four
conditions}}$ as we have seen in Sublemma 5.5.
\endproclaim \bs

\vfill \pagebreak

{\bf \S{6}. The proofs of Theorem 5.0 with five sublemmas and
corollaries in \S {5}} \bs

{\bf \S{6.0}. Introduction}

In preparation for the proof of Theorem $5.0$, first  we prove five
sublemmas, Sublemma $5.1$, Sublemma $5.2$,\dots, Sublemma $5.5$,
respectively in $\S{6.1}$. After then, we will finish the proof of
Theorem $5.0$ with corollaries in $\S{6.2}$. \ms

{\bf \S{6.1}. The proofs of five sublemmas}

\demo{\bf Proof of Sublemma 5.1} If $r=2$, then it is trivial to
prove that $\Delta^{\sharp}_2(\beta_{2,1},\beta_{2,2})>
n_1\beta_{1,1}n_2 \quad \text{on $g_2$}$, because
$\Delta^{\sharp}_2(\beta_{2,1},\beta_{2,2})=\Delta_2(\beta_{2,1},\beta_{2,2})$
by (5.1.1) and $\Delta_2(\beta_{2,1},\beta_{2,2})
>n_1\beta_{1,1}n_2$ by \text{\rm The 4-th ${\text{\rm{Cond}}}^{\text{{\rm(0)}}}$}
in the assumptions of Theorem $5.0$.

Let $r\ge 3$. For any $\ell=3,4,\dots,r$, it is trivial to note by
definition of
$\Delta^{\sharp}_{\ell}(\beta_{{\ell},k})^{\ell}_{k=1}$ in (5.1.1)
that the following three equalities are the same, and so we can
write
$c=\Delta^{\sharp}_{\ell}(\beta_{{\ell},k})^{\ell}_{k=1}-n_1\beta_{1,1}n_2n_3\cdots
n_\ell$ for convenience of notation:
$$\align
(5.1.3)  \quad \qquad \qquad
&\text{$c=\Delta^{\sharp}_{\ell}(\beta_{{\ell},k})^{\ell}_{k=1}
-n_1\beta_{1,1}n_2n_3\cdots n_\ell$} \\
\iff \quad &  \\
(5.1.4)  \quad \qquad \qquad
&\text{$c=\Delta_2(\beta_{\ell,1},\beta_{\ell,2})+n_1\beta_{1,1}\beta_{\ell,3}
 +n_1\beta_{1,1}n_2\beta_{\ell,4}
 +n_1\beta_{1,1} n_2 n_3 \beta_{\ell,5}$}   \\
& \quad \text{\qquad $+ \cdots +n_1\beta_{1,1}n_2 \cdots
n_{\ell-2}\beta_{\ell,\ell}-n_1\beta_{1,1}n_2n_3\cdots n_\ell$} \\
\iff \quad &  \\
(5.1.5)  \quad \qquad \qquad
&\text{$c=\Delta_2(\beta_{\ell,1},\beta_{\ell,2})-(n_\ell
n_{\ell-1}\cdots n_2
-\beta_{\ell,\ell}n_{\ell-2}n_{\ell-3}\cdots n_2$}  \\
& \qquad \text{\qquad $
-\beta_{\ell,\ell-1}n_{\ell-3}n_{\ell-4}\cdots n_2 -\cdots
-\beta_{\ell,4}n_2-\beta_{\ell,3})n_1\beta_{1,1}$}. \qquad \qquad
\endalign$$ \ms

So, for any integer $\ell\ge 3,$ it suffices to show that the above
integer $c$ of (5.1.5) can be equal to an integer $\xi_{\ell-2}>0$
such that $\xi_{\ell-2}$ is the {\rm $(\ell-2)$-th} element of a
positive sequence $\{\xi_j:j=1,2,\dots,\ell-2\}$, where each $\xi_j$
satisfies the following properties:
$$\align
(5.1.6) \quad \xi_0 &=\Delta_{\ell}(\beta_{\ell, k})^{\ell}_{k=1}
 -\{n_{\ell}\}n_{\ell-1}\Delta_{\ell-1}(\beta_{\ell-1,k})^{\ell-1}_{k=1}>0,  \\
 \xi_1  &=\Delta_{\ell-1}(\beta_{\ell, k})^{\ell-1}_{k=1}
 -\{n_{\ell}n_{\ell-1}-\beta_{\ell,\ell}\}n_{\ell-2}
 \Delta_{\ell-2}(\beta_{\ell-2,k})^{\ell-2}_{k=1}>0,  \\
 \xi_2  &=\Delta_{\ell-2}(\beta_{\ell,k})^{\ell-2}_{k=1}-\{n_{\ell}n_{\ell-1}n_{\ell-2}
 -\beta_{\ell,\ell}n_{\ell-2}-\beta_{\ell,\ell-1}\}
 n_{\ell-3}\Delta_{\ell-3}(\beta_{\ell-3,k})^{\ell-3}_{k=1}>0, \quad\\
   \xi_j &=\Delta_{\ell-j}(\beta_{\ell,k})^{\ell-j}_{k=1}
  -\{n_{\ell}n_{\ell-1}\cdots n_{\ell-j}
  -\beta_{\ell,\ell}n_{\ell-2}n_{\ell-3}\cdots n_{\ell-j}  \\
 &\qquad -\beta_{\ell,\ell-1}n_{\ell-3}n_{\ell-4}
  \cdots n_{\ell-j}-\cdots -\beta_{\ell,\ell-j+2}n_{\ell-j} \\
  &\qquad -\beta_{\ell,\ell-j+1}\}
  \times n_{\ell-j-1}\Delta_{\ell-j-1}(\beta_{\ell-j-1,k})^{\ell-j-1}_{k=1}>0
  \quad \text{for $3\le j\le {\ell-2}$.} \\
\endalign$$

Let $\ell\ge 3$ be chosen arbitrary. Now, we will show by the
induction method on the nonnegative integer $j\le \ell-2$ that
$\xi_j$ is positive for all $j$.

It is trivial by (d) in the assumption of Theorem $5.0$ that
$\xi_0>0$.

To prove that $\xi_1$ is positive, first of all, it is easy to
observe the following by the definition of
$\Delta_{\ell}(\beta_{\ell,k})^{\ell}_{k=1}$:
$$\align
(5.1.7) \quad 0<\xi_0 &=\Delta_{\ell}(\beta_{\ell,k})^{\ell}_{k=1}
 -n_{\ell}n_{\ell-1}\Delta_{\ell-1}(\beta_{\ell-1,k})^{\ell-1}_{k=1}\\
&=\beta_{\ell,\ell}\Delta_{\ell-1}(\beta_{\ell-1,k})^{\ell-1}_{k=1}
+n_{\ell-1}\Delta_{\ell-1}(\beta_{\ell,k})^{\ell-1}_{k=1}
  -n_{\ell}n_{\ell-1}\Delta_{\ell-1}(\beta_{\ell-1,k})^{\ell-1}_{k=1} \\
&=n_{\ell-1}\Delta_{\ell-1}(\beta_{\ell,k})^{\ell-1}_{k=1}-(n_{\ell}n_{\ell-1}
-\beta_{\ell,\ell})\Delta_{\ell-1}(\beta_{\ell-1,k})^{\ell-1}_{k=1}.
\endalign$$

Since $\Delta_{\ell-1}(\beta_{\ell-1,k})^{\ell-1}_{k=1}
>n_{\ell-1}n_{\ell-2}\Delta_{\ell-2}(\beta_{\ell-2,k})^{\ell-2}_{k=1}$
by \text{\rm The 4-th ${\text{\rm{Cond}}}^{\text{{\rm(0)}}}$} in the
assumption of Theorem $5.0$, then the third inequality of (5.1.7)
and $\beta_{\ell,\ell}<n_{\ell-1}$ imply that
$$\align
(5.1.8) \qquad \qquad
n_{\ell-1}\Delta_{\ell-1}(\beta_{\ell,k})^{\ell-1}_{k=1}
&>(n_{\ell}n_{\ell-1}-\beta_{\ell,\ell})
\Delta_{\ell-1}(\beta_{\ell-1,k})^{\ell-1}_{k=1},
\quad \text{and so} \\
n_{\ell-1}\Delta_{\ell-1}(\beta_{\ell,k})^{\ell-1}_{k=1}
&>(n_{\ell}n_{\ell-1}-\beta_{\ell,\ell})n_{\ell-1}n_{\ell-2}
\Delta_{\ell-2}(\beta_{\ell-2,k})^{\ell-2}_{k=1}, \qquad \qquad
\endalign$$
whether or not $n_{\ell}n_{\ell-1}-\beta_{\ell,\ell}>0$. \ms

Dividing both sides on (5.1.8) by $n_{\ell-1}$, then
$$
\Delta_{\ell-1}(\beta_{\ell,k})^{\ell-1}_{k=1}
>(n_{\ell}n_{\ell-1}-\beta_{\ell,\ell})n_{\ell-2}
\Delta_{\ell-2}(\beta_{\ell-2,k})^{\ell-2}_{k=1}, \tag 5.1.9
$$
which is equivalent to the fact that $\xi_1>0$.

By the induction assumption on the positive integer $j\le \ell-2$,
suppose we have shown that $\xi_j$ is positive with $1\le j\le
\ell-3$. To prove that $\xi_{j+1}$ is positive, for convenience of
notations, let $\xi_j$ of $(5.1.6)$ be written again in the form
$$
\align (5.1.10) \qquad   \xi_j
&=\Delta_{\ell-j}(\beta_{\ell,k})^{\ell-j}_{k=1}
 -\omega_jn_{\ell-j-1}\Delta_{\ell-j-1}(\beta_{\ell-j-1,k})^{\ell-j-1}_{k=1}>0
 \quad \text{with} \qquad \qquad \quad \\
 \omega_j &=n_{\ell}n_{\ell-1}\cdots n_{\ell-j}-\beta_{\ell,\ell}n_{\ell-2}n_{\ell-3}\cdots
 n_{\ell-j}\\
& \quad -\beta_{\ell,\ell-1}n_{\ell-3}n_{\ell-4}\cdots
n_{\ell-j}-\cdots
-\beta_{\ell,\ell-j+2}n_{\ell-j}-\beta_{\ell,\ell-j+1}.
\endalign$$

Now, by (5.1.10) and the definition of
$\Delta_{\ell-j}(\beta_{\ell,k})^{\ell-j}_{k=1}$ only, it is easy
to prove that
$$\align
(5.1.11) \quad  0<\xi_j
 &=\Delta_{\ell-j}(\beta_{\ell,k})^{\ell-j}_{k=1}
 -\omega_jn_{\ell-j-1}\Delta_{\ell-j-1}(\beta_{\ell-j-1,k})^{\ell-j-1}_{k=1}
 \qquad \qquad \quad \\
 &=\beta_{\ell,\ell-j}\Delta_{\ell-j-1}(\beta_{\ell-j-1,k})^{\ell-j-1}_{k=1}
 +n_{\ell-j-1}\Delta_{\ell-j-1}(\beta_{\ell,k})^{\ell-j-1}_{k=1}\\
 & \quad -\omega_jn_{\ell-j-1}\Delta_{\ell-j-1}(\beta_{\ell-j-1,k})^{\ell-j-1}_{k=1}\\
 &=n_{\ell-j-1}\Delta_{\ell-j-1}(\beta_{\ell,k})^{\ell-j-1}_{k=1}
 -(\omega_jn_{\ell-j-1}-\beta_{\ell,\ell-j})
\Delta_{\ell-j-1}(\beta_{\ell-j-1,k})^{\ell-j-1}_{k=1}. \quad
\endalign$$

Since $\Delta_{r-j-1}(\beta_{r-j-1,k})^{r-j-1}_{k=1}
>n_{r-j-1}n_{r-j-2}\Delta_{r-j-2}(\beta_{r-j-2,k})^{r-j-2}_{k=1}$
by \text{\rm The 4-th ${\text{\rm{Cond}}}^{\text{{\rm(0)}}}$} in the
assumption of Theorem $5.0$, then we get the following from the last
equality in $(5.1.11)$:
$$\align
(5.1.12) {\quad} &
n_{\ell-j-1}\Delta_{\ell-j-1}(\beta_{\ell,k})^{\ell-j-1}_{k=1}
 >(\omega_jn_{\ell-j-1}-\beta_{\ell,\ell-j})
 n_{\ell-j-1}\Delta_{\ell-j-1}(\beta_{\ell,k})^{\ell-j-1}_{k=1}
\text{\rm and so} \\
\quad
&n_{\ell-j-1}\Delta_{\ell-j-1}(\beta_{\ell,k})^{\ell-j-1}_{k=1}>
(\omega_jn_{\ell-j-1}-\beta_{\ell,\ell-j})n_{\ell-j-1}n_{\ell-j-2}
\Delta_{\ell-j-2}(\beta_{\ell-j-2,k})^{\ell-j-2}_{k=1}, \\
\endalign$$
whether or not $\omega_jn_{\ell-j-1}-\beta_{\ell,\ell-j}>0$. \ms

Dividing both sides of $(5.1.12)$ by $n_{\ell-j-1}$, then we get
$$\align
(5.1.13) \qquad \qquad
\Delta_{\ell-j-1}(\beta_{\ell,k})^{\ell-j-1}_{k=1}
>(\omega_jn_{\ell-j-1}-\beta_{\ell,\ell-j})n_{\ell-j-2}
\Delta_{\ell-j-2}(\beta_{\ell-j-2,k})^{\ell-j-2}_{k=1}. \qquad
\endalign$$

Before we prove that $\xi_{j+1}>0$, then it is trivial to observe
by (5.1.10) that $\xi_{j+1}$ of (5.1.6) can be rewritten as
follows:
$$\align
(5.1.14) \quad \quad
\xi_{j+1}=\Delta_{\ell-j-1}(\beta_{\ell,k})^{\ell-j-1}_{k=1}
-(\omega_jn_{\ell-j-1}-\beta_{\ell,\ell-j})n_{\ell-j-2}
\Delta_{\ell-j-2}(\beta_{\ell-j-2,k})^{\ell-j-2}_{k=1}.
\endalign$$
Now, we can show by (5.1.13) that $\xi_{j+1}$ is positive. Thus,
the proof is done. $\square$
\enddemo
\ms

\demo{\bf Proof of Sublemma 5.2} We will prove (a), (b), (c), (d)
and (e), respectively. \ms

\text{\bf(a)} \quad In preparation for the proof of an equality in
(5.2.1), it suffices to show that for any integer $r\ge 2$,
$g_r=g_r(y,z)$ in the assumption of this theorem can be generally
represented in the following form:
$$\align
\text{\rm(5.2.2)} \qquad \qquad  g_r &=\Sigma_{r,0}+\Sigma_{r,1},   \\
\Sigma_{r,0} & =(z^{n_1}+\ve_{1}y^{\beta_{1,1}})^{n_2n_3\cdots
n_r}=(g_1)^{n_2n_3\cdots n_r}
\quad \text{with $\ve_{1}=1$ or a unit in $\BC\{y,z\}$},\\
\Sigma_{r,1} &=\sum_{\g,\de\ge 0} c^{(1)}_{\g,\de}y^{\g}z^{\de}
\quad \text{with
$n_1\g+\beta_{1,1}\de>n_1\beta_{1,1}n_2n_3\cdots n_r$}, \quad \text{and} \\
\text{\rm(5.2.3)} \qquad \qquad  g_r^{\ell}
&=\Sigma^{(\ell)}_{r,0}+\Sigma^{(\ell)}_{r,1} \quad
\text{for any integer $\ell\ge 2$},  \\
\Sigma^{(\ell)}_{r,0} &
=(z^{n_1}+\ve_{1}y^{\beta_{1,1}})^{n_2n_3\cdots n_r\ell}
=(g_1)^{n_2n_3\cdots n_r\ell}=(\Sigma_{r,0})^{\ell},\\
\Sigma^{(\ell)}_{r,1} &=\sum_{\g,\de\ge 0}
c^{(\ell)}_{\g,\de}y^{\g}z^{\de} \quad \text{with} \quad
n_1\g+\beta_{1,1}\de>n_1\beta_{1,1}n_2n_3\cdots n_r\ell, \qquad
\qquad
\endalign$$
where $\ve_{1}=\ve_{1}(y,z)$ is assumed to be one in $\BC\{y,z\}$ if
necessary, and the $c^{(1)}_{\g,\de}$ are nonzero complex numbers
for some nonnegative integers $\g$ and $\de$, if exists, and the
$c^{(\ell)}_{\g,\de}$ are nonzero complex numbers for some
nonnegative integers $\g$ and $\de$ if exists.

For the induction proof of the equalities in (5.2.2) and (5.2.3),
it suffices to consider the following two cases:

Case(I) $r=2$ and Case(II) $r>2$. \ms

$\underline{\text{\rm Case(I)}}$ \quad Assuming that $r=2$, recall
by \text{\rm The 2nd ${\text{\rm{Cond}}}^{\text{{\rm(0)}}}$} and
\text{\rm The 4-th ${\text{\rm{Cond}}}^{\text{{\rm(0)}}}$} in the
assumptions of this theorem and by (5.1.2) of Sublemma $5.1$ that
$g_2$ can be written in the form
$$\align
(5.2.4) \qquad \qquad
&\text{$g_2=g^{n_2}_1+\ve_{2}y^{\beta_{2,1}}z^{\beta_{2,2}}$
\quad and} \qquad \qquad \\
&\text{$\Delta^{\sharp}_2(\beta_{2,1},\beta_{2,2})=\Delta_2(\beta_{2,1},\beta_{2,2})
=n_1\beta_{2,1}+\beta_{1,1}\beta_{2,2}> n_1\beta_{1,1}n_2$ \quad {on
$g_2$}}, \qquad \qquad
\endalign$$
where $\ve_2$ is a unit in $\BC\{y,z\}$.

In order to prove the equalities in (5.2.2) and (5.2.3) for $r=2$,
it remains to show that the equality in (5.2.3) holds, because the
proof of the equality in (5.2.2) was already proved by (5.2.4).

So, in preparation for the proof of the equality in (5.2.3) with
$r=2$, then it is clear by (5.2.2) or (5.2.4) that $g_2(y,z)$ can
be rewritten in the form
$$\align
\text{\rm(5.2.5)} \qquad \qquad \qquad g_2 &=\Sigma_{2,0}+\Sigma_{2,1},   \\
\Sigma_{2,0} & =(z^{n_1}+\ve_{1}y^{\beta_{1,1}})^{n_2}
\qquad \text{with  $\ve_{1}=1$ or a unit in $\BC\{y,z\}$}, \qquad \qquad \qquad \\
\Sigma_{2,1} &=\sum_{\g,\de\ge 0} c^{(1)}_{\g,\de}y^{\g}z^{\de}
\quad \quad \text{with  $n_1\g+\beta_{1,1}\de>n_1\beta_{1,1}n_2$},  \\
\endalign$$
where the $c^{(1)}_{\g,\de}$ are nonzero complex numbers for some
nonnegative integers $\g$ and $\de$.

For the proof of an inequality in (5.2.3), note by (5.2.5) that
$g_2^{\ell}$ can be written in the form
$$\align
\text{\rm(5.2.6)} \qquad \quad g_2^{\ell}
&=(\Sigma_{2,0}+\Sigma_{2,1})^{\ell}
\quad \text{for any integer $\ell\ge 2$} \\
&=(\Sigma_{2,0})^{\ell}+\sum^{\ell-1}_{k=1}\binom{\ell}{k}
(\Sigma_{2,0})^{k}(\Sigma_{2,1})^{\ell-k}+(\Sigma_{2,1})^{\ell}\\
&=\Sigma^{(\ell)}_{2,0}+\Sigma^{(\ell)}_{2,1}, \\
\text{where} \quad \Sigma^{(\ell)}_{2,0} &=(\Sigma_{2,0})^{\ell}
=(z^{n_1}+\ve_{1}y^{\beta_{1,1}})^{n_2\ell}
\quad \text{with $\ve_{1}=1$ or a unit in $\BC\{y,z\}$}, \qquad \qquad \\
\Sigma^{(\ell)}_{2,1} &=\sum^{\ell-1}_{k=1}\binom{\ell}{k}
(\Sigma_{2,0})^{k}(\Sigma_{2,1})^{\ell-k}+(\Sigma_{2,1})^{\ell}
\quad \text{for notation.} \qquad
\endalign$$

To prove that an inequality in (5.2.3) is true for $r=2$, apply the
equalities of (5.2.5) and Sublemma $5.1$ to equalities of (5.2.6).
After then,  it suffices to show that for any nonzero monomial
$y^{\alpha}z^{\beta}\in \Sigma^{(\ell)}_{2,1}$,
$$
n_1\alpha+\beta_{1,1}\beta> n_1\beta_{1,1}n_2\ell. \tag 5.2.7
$$

To prove $(5.2.7)$, it is enough to show that the following two
claims hold by using the defining equation for
$\Sigma^{(\ell)}_{2,1}$ in $(5.2.6)$:
$$\align
(5.2.8) \qquad \qquad  & \text{Claim(i)} \qquad \text{For any
monomial $y^{\alpha}z^{\beta}\in (\Sigma_{2,1})^{\ell}$}, \qquad
\qquad \qquad \qquad \qquad \qquad\\
& \qquad \qquad \quad  n_1\alpha+\beta_{1,1}\beta> n_1\beta_{1,1}n_2\ell. \\
&\text{Claim(ii)} \qquad \text{For any monomial
$y^{\gamma}z^{\delta}\in (\Sigma_{2,0})^{k}(\Sigma_{2,1})^{\ell-k}$}, \\
& \qquad \qquad \quad n_1\gamma+\beta_{1,1}\delta>
n_1\beta_{1,1}n_2k+n_1\beta_{1,1}n_2(\ell-k)=n_1\beta_{1,1}n_2\ell.
\endalign$$
Note by (5.2.5) that the proof of two claims in (5.2.8) is trivial
because $\ell-k>0$, and so the proof of (5.2.3) for $r=2$ is done.
So, the proof of Case(I) is done. \ms

$\underline{\text{\rm Case(II)}}$ Let $r>2$. Now, suppose we have
proved by induction assumption on the positive integer $j<r$ that
the representation of $g^{\ell}_j$ in (5.2.2) and (5.2.3) is true
for $2\le j<r$ and for any integer $\ell\ge 1$. Then, recall by
\text{\rm The 2nd ${\text{\rm{Cond}}}^{\text{{\rm(0)}}}$} and
\text{\rm The 4-th ${\text{\rm{Cond}}}^{\text{{\rm(0)}}}$} in the
assumptions of this theorem and by (5.1.2) of Sublemma $5.1$ that
$g_{j+1}$ can be rewritten as follows:
$$\align
 & g_{j+1}=g^{n_{j+1}}_{j}+\ve_{j+1}y^{\beta_{j+1,1}}
z^{\beta_{j+1,2}}g^{\beta_{j+1,3}}_1\cdots g^{\beta_{j+1,j+1}}_{j-1}
\quad \text{with}  \tag 5.2.9 \\
& \Delta^{\sharp}_{j+1}(\beta_{j+1,k})^{j+1}_{k=1} >
n_1\beta_{1,1}n_2n_3\cdots
n_{j}n_{j+1} \quad \text{on $g_{j+1}$}, \\
\endalign$$
where $\ve_{j+1}$ is defined to be a unit in $\BC\{y,z\}$ by
(5.2.9).

Now, applying the induction assumption on the integer $j$, we may
assume that for each $k=1,2,\dots,j$, and any integer $\ell>0$, we
have
$$\align
 g_k^{\ell}&=\Sigma^{(\ell)}_{k,0}+\Sigma^{(\ell)}_{k,1} \quad
\text{for any integer $\ell\ge 2$}, \tag 5.2.10 \\
\Sigma^{(\ell)}_{k,0} &
=(z^{n_1}+\ve_{1}y^{\beta_{1,1}})^{n_2n_3\cdots n_k\ell}
=(g_1)^{n_2n_3\cdots n_k\ell}=(\Sigma^{(1)}_{k,0})^{\ell},\\
\Sigma^{(\ell)}_{k,1} &=\sum_{\g_,\de\ge 0}
c^{(\ell)}_{\g,\de}y^{\g}z^{\de} \quad \text{with} \quad
n_1\g+\beta_{1,1}\de>n_1\beta_{1,1}n_2n_3\cdots n_k\ell,
\endalign$$
where the $c^{(\ell)}_{\g,\de}$ are nonzero complex numbers for some
nonnegative integers $\g$ and $\de$.

To prove that an inequality in (5.2.2) is true for $k=j+1$, apply
the equalities of (5.2.10) and Sublemma $5.1$ to the equalities of
$g_{j+1}$ in (5.2.9). Since it is clear by (5.2.10) that
$g^{n_{j+1}}_{j}$ can be written in the form
$$\align
\text{\rm(5.2.11)} \qquad \qquad g^{n_{j+1}}_{j}
&=\Sigma^{(n_{j+1})}_{j,0}+\Sigma^{(n_{j+1})}_{j,1} \quad
\text{for any integer $n_{j+1}\ge 2$},  \\
\Sigma^{(n_{j+1})}_{j,0} &
=(z^{n_1}+\ve_{1}y^{\beta_{1,1}})^{n_2n_3\cdots n_jn_{j+1}}
=(g_1)^{n_2n_3\cdots n_jn_{j+1}},\\
\Sigma^{(n_{j+1})}_{j,1} &=\sum_{\g_,\de\ge 0}
c^{(n_{j+1})}_{\g,\de}y^{\g}z^{\de} \quad \text{with} \quad
n_1\g+\beta_{1,1}\de>n_1\beta_{1,1}n_2n_3\cdots n_jn_{j+1}, \qquad
\qquad
\endalign$$
where the $c^{(n_{j+1})}_{\g,\de}$ are nonzero complex numbers for
some nonnegative integers $\g$ and $\de$.

For the proof of an inequality in (5.2.2) for $k=j+1$, applying
(5.2.11) to (5.2.9), it remains to show the following: for any
nonzero monomial $y^{\alpha}z^{\beta}\in g_{j+1}-g^{n_{j+1}}_j$,
$$
n_1\alpha+\beta_{1,1}\beta> n_1\beta_{1,1}n_2n_3\cdots n_jn_{j+1}.
\tag 5.2.12
$$ \ms

First of all, it is clear by (5.2.10) that for any nonzero
monomial $y^{\alpha_k}z^{\beta_k}\in g^{\ell}_k$,
$$\align
 n_1\alpha_k+\beta_{1,1}\beta_k\ge
n_1\beta_{1,1}n_2n_3\cdots n_k\ell \quad \text{for $k=1,2,\dots,j$.}
\tag 5.2.13
\endalign$$

So, in order to prove that either (5.2.2) or (5.2.12) is true on
$g_{j+1}$, it is clear by (5.2.9), (5.2.10) and (5.2.13) that any
nonzero monomial $y^{\gamma}z^{\de}\in
\dfrac{1}{\ve_{j+1}}\{g_{j+1}-g^{n_{j+1}}_j\}$ can be represented
as follows:
$$\align
(5.2.14) \quad & y^{\gamma}z^{\de}=y^{\beta_{j+1,1}}
z^{\beta_{j+1,2}}\Pi^{j-1}_{k=1}y^{\alpha_k}z^{\beta_k} \quad
\text{such that}  \\
& \quad \text{$n_1\alpha_k+\beta_{11}\beta_k\ge
n_1\beta_{11}n_2n_3\cdots n_k\beta_{j+1,k+2}$ \quad for some
$y^{\alpha_k}z^{\beta_k}\in g^{\beta_{j+1,k+2}}_k$,} \qquad \qquad
\endalign$$
where $\ve_{j+1}$ is a unit in $\BC\{y,z\}$.

In other words, whenever $y^{\gamma}z^{\de}\in
\dfrac{1}{\ve_{j+1}}\{g_{j+1}-g^{n_{j+1}}_j\}$ is chosen
arbitrary, then it may be assumed by (5.2.14) that
$$\align
\text{\rm (5.2.15)} \qquad \qquad \gamma &=\beta_{j+1,1}+\alpha_1
+\cdots+\alpha_{j-1}\quad \text{and} \quad
\delta =\beta_{j+1,2}+\beta_1+\cdots+\beta_{j-1}, \qquad \qquad\\
\text{where} \quad & \text{for each $k=1,2,\dots,j-1$,
$y^{\alpha_k}z^{\beta_k}\in
g^{\beta_{j+1,k+2}}_k$} \\
& \text{with $n_1\alpha_k+\beta_{1,1}\beta_k\ge
n_1\beta_{1,1}n_2n_3\cdots n_k\beta_{j+1,k+2}$.} \qquad \qquad
\qquad \qquad
\endalign$$

Therefore, by (5.2.14) and (5.2.15) and by Sublemma $5.1$ again,
for any nonzero monomial $y^{\gamma}z^{\de}\in
\dfrac{1}{\ve_{j+1}}\{g_{j+1}-g^{n_{j+1}}_j\}$, we can prove the
following:
$$\align
\text{\rm (5.2.16)} \qquad \qquad n_1\gamma+\beta_{1,1}\delta & =
n_1\beta_{j+1,1}+\beta_{1,1}\beta_{j+1,2}+\sum^{j-1}_{k=1}
(n_1\alpha_k+\beta_{1,1}\beta_k)  \\
\ge &n_1\beta_{j+1,1}+\beta_{1,1}\beta_{j+1,2}+\sum^{j-1}_{k=1}
(n_1\beta_{1,1}n_2n_3\cdots n_k\beta_{j+1,k+2}) \qquad \qquad \\
=&\Delta^{\sharp}_{j+1}(\beta_{j+1,k})^{j+1}_{k=1}
>n_1\beta_{1,1}n_2n_3\cdots n_jn_{j+1}.
\endalign$$

Thus, if $r=j+1$, then the proof of (5.2.2) is done by (5.2.11)
and (5.2.16). Also, if r=j+1, then the proof of (5.2.3) is trivial
by the same method as we have used in the proof for (5.2.3) with
r=2. Then, the proof of Case(II) is done.

So, we finished the proof of (a) by Case(I) and Case(II). \ms

\text{\bf(b)} To prove {\rm(b1)}, it suffices to consider $g_r(0,z)$
from $g_r(y,z)$ of (5.2.1). Then
$$\align
g_r(0,z) &=z^{n_1n_2\cdots n_r}+\sum c^{(r)}_{0,\beta}z^{\beta}
\quad\text{with} \tag 5.2.17\\
 \beta_{1,1}\beta &>n_1\beta_{1,1}n_2n_3\cdots n_r.
\endalign$$
Thus, $\beta>n_1n_2\cdots n_r$, and so it is done. Also, the proof
of (b2) can be done similarly. \ms

\text{\bf(c)} \quad To prove (c1), suppose that $n_1<\beta_{11}$. If
then, by (a), $\beta_{1,1}(\alpha+\beta)\ge
n_1\alpha+\beta_{1,1}\beta>n_1\beta_{1,1}n_2n_3\cdots n_r$, which
implies that $\alpha+\beta>n_1\prod^r_{k=2}n_k$. Thus, the proof of
(c1) is done. Similarly, (c2) can be proved. \ms

\text{\bf(d)} \quad By (a) and Theorem $3.6$, it is clear. \ms

\text{\bf(e)} \quad There is nothing to prove. Thus, the proof of
this sublemma is finished. $\square$
\enddemo
\ms

\demo{\bf Proof of Sublemma 5.3} To prove (5.3.3) for any $r\ge
2$, it is enough to show by (5.3.2) that the following equation in
(5.3.4) is nonnegative:
$$\split
(5.3.4) \qquad
&\Omega^{\sharp}_r(\beta_{r,k})^r_{k=1}-bn_1n_2n_3\cdots n_r \\
 =~& a\beta_{r,1}+b\beta_{r,2}+bn_1\beta_{r,3}+bn_1n_2\beta_{r,4}
+\cdots +bn_1n_2\cdots n_{r-2}\beta_{r,r}-bn_1n_2n_3\cdots n_r \qquad \qquad \\
=~&a\beta_{r,1}+bD\ge 0 \quad \text{with} \\
& D=\beta_{r,2}+n_1\beta_{r,3}+n_1n_2\beta_{r,4}+\cdots
+n_1n_2\cdots
n_{r-2}\beta_{r,r}-n_1n_2n_3\cdots n_r,\\
\endsplit$$
by definition of $\Omega^{\sharp}_r(\beta_{r,k})^r_{k=1}$ in (5.3.2)
where $a>0$ and $b\ge 0$.

Then, it remains to prove by (5.3.4) that $a\beta_{r,1}+bD\ge 0$,
independently of $D$. For the proof, it suffices to consider two
cases, Case(i) $D\ge 0$ and Case(ii) $D<0$. \ms

\noindent$\underline{\text{\rm Case(i)}}$ \ Let $D\ge 0$. It is
clear that $a\beta_{r,1}+bD\ge 0$ because $a>0$, and also $b\ge 0$
with $\beta_{r,1}$ nonnegative. Thus, the proof of Case(i) is done.
\ms

\noindent$\underline{\text{\rm Case(ii)}}$ \ Let $D<0$. In
preparation for proof of the inequality, first of all, note that the
inequality
$\Delta^{\sharp}_r(\beta_{r,k})^r_{k=1}-n_1\beta_{1,1}n_2n_3\cdots
n_r>0$ of Sublemma $5.1$ can be equivalently rewritten as follows:
$$\align
(5.3.5) \qquad \qquad & \Delta^{\sharp}_r(\beta_{r,k})^r_{k=1}
-n_1\beta_{1,1}n_2n_3\cdots n_r \\
=& n_1\beta_{r,1}+\beta_{1,1}\beta_{r,2}+n_1\beta_{1,1}\beta_{r,3}
+n_1\beta_{1,1}n_2\beta_{r,4} \\
&\quad +\cdots +n_1\beta_{1,1}n_2\cdots n_{r-2}\beta_{r,r}
-n_1\beta_{1,1}n_2n_3\cdots n_r  \quad \text{by (5.1.1)} \\
=& n_1\beta_{r,1}+\beta_{1,1}D \quad \text{is \quad  positive}, \\
where \quad & D=\beta_{r,2}+n_1\beta_{r,3}+n_1n_2\beta_{r,4}+\cdots
+n_1n_2\cdots
n_{r-2}\beta_{r,r}-n_1n_2n_3\cdots n_r \ \text{by (5.3.4)}.\\
\endalign$$

Since $-D>0$ and $n_1\ge 2>0$, then the inequality
$n_1\beta_{r,1}+\beta_{1,1}D>0$ in $(5.3.5)$ can be equivalently
represented as follows:
$$
 \frac{\beta_{r,1}}{-D} > \frac{\beta_{1,1}}{n_1} \ . \tag
5.3.6
$$
Also, $a\beta_{1,1}-bn_1=1$ with $a>0$ implies that
$\dfrac{\beta_{1,1}}{n_1}
> \dfrac{b}{a}$. Therefore, we proved by $(5.3.6)$ that
$\dfrac{\beta_{r,1}}{-D}
>\dfrac{b}{a}$, that is, $a\beta_{r,1}+bD>0$. Thus, the proof of
Case(ii) is done.

Therefore, we showed by Case(i) and Case(ii) that the equation in
(5.3.4) is nonnegative, and so the proof of this sublemma is
finished. $\square$
\enddemo
\ms

\demo{\bf Proof of Sublemma 5.4} Following the same assumptions and
notations as in Sublemma 5.1, Sublemma 5.2 and Sublemma 5.3, then
(a) in the conclusion of Sublemma 5.2 is true, and so for each
$j=2,3,\dots,r$, \ $g_j=g_j(y,z)$ of (5.2.1) can be easily rewritten
as follows:
$$\split
\noindent(5.4.6) \qquad \qquad  g_j&=(z^{n_1}+\ve_1
y^{\beta_{1,1}})^{n_2n_3\cdots n_j}+\sum_{\alpha,\beta\ge 0}
c^{(j)}_{\alpha,\beta}y^{\alpha}z^{\beta} \quad \text{with $\ve_1=1$
and}
\qquad \qquad \qquad \\
&  \text{with \quad
$n_1\alpha+\beta_{1,1}\beta>n_1\beta_{1,1}n_2n_3\cdots n_j$,}
\endsplit$$
where $\ve_1=\ve_1(y,z)$ is assumed to be one in $\BC\{y,z\}$, and
the $c^{(j)}_{\alpha,\beta}$ are nonzero complex numbers for some
nonnegative integers $\alpha$ and $\beta$ such that
$n_1\alpha+\beta_{1,1}\beta>n_1\beta_{1,1}n_2n_3\cdots n_j$.

First, we will show how to apply Theorem $3.6$ to the proof of (a),
(b) and (d) in this sublemma, and after then, the remaining part (c)
of this sublemma will be proved computationally.

In preparation for the proof of (a), (b) and (d) in this sublemma,
it is clear that the equation of $g_j$ of (5.4.6) satisfies the same
kind of properties as $f$ does in the assumption of Theorem $3.6$,
which can be represented as follows: \ms

$\underline{\text{$g_j$ of (5.4.6) satisfies the same kind of
assumption as in Theorem 3.6}}$ \quad Let
$V(g_1)=\{(y,z):g_1(y,z)=0\}$, $V(f)=\{(y,z): f(y,z)=0\}$ and
$V(G)=\{(y,z): G(y,z)=0\}$ be analytic varieties at $(0,0)$ in
$\BC^2$, each of which is written respectively as follows: For
convenience of notation, substitute $g_j$ of (5.4.6) by $f$, for an
application of Theorem $3.6$.
$$\align
 g_1 &=z^{n_1}+\ve_1y^{\beta_{1,1}} \quad \text{with $\ve_1=1$},  \tag 5.4.7 \\
f&={g_1}^{d_j}+\sum_{\alpha,\beta\ge
0}c^{(j)}_{\alpha,\beta}y^{\alpha}z^{\beta} \quad \text{with} \quad
n_1\alpha+\beta_{1,1}\beta>n_1k_1d_j, \\
F&=y^{\delta_1}z^{\delta_2}f, \\
G&=y^{\gamma}g_1 \\
\endalign
$$
satisfying the properties {\rm(i)}, {\rm(ii)}, {\rm(iii)},
{\rm(iv)} and {\rm(v)}:

\roster
\item "(i)" $\gcd(n_1,\beta_{1,1})=1$ with $n_1\ge 2$ and
$\beta_{11}\ge 1$.

\item "(ii)" $d_j=n_2n_3\cdots n_j$ is a positive integer with
$d_j\ge 2$.

\item "(iii)" $\ve_1$ is assumed to be one in $\BC\{y,z\}$, and
the $c^{(j)}_{\alpha,\beta}$ are nonzero complex numbers for some
nonnegative integers $\alpha$ and $\beta$ such that
$n_1\alpha+k_1\beta> n_1k_1d_j$, if exist.

\item "(iv)" Assume that $V(f)$ has an isolated singular point at
the origin as a reduced variety.

\item "(v)" If $\beta_{1,1}=1$, then $\gamma=1$, and if
$\beta_{1,1}\ge 2$, then $\gamma=0$.

\item "(vi)" In addition, assume that each $\delta_i$ is either a
positive integer or $0$ for $i=1,2$, as far as $V(F)$ has an
isolated singular point at the origin as a reduced variety, even
if $d_j\ge 1$.

\endroster \ms

So, we have the same kind of conclusion as we have seen in Theorem
$3.6$, up to change of notations:

$\underline{\text{The same kind of conclusion as in Theorem 3.6}}$
\quad Let
$\tau_m=\pi_1\circ\pi_2\circ\cdots\circ\pi_m:M^{(m)}\to\BC^2$ be the
compositions of a finite number $m$ of successive blow-ups $\pi_i$
which is needed to get the standard resolution of the singular point
of $V(G)=V(y^{\gamma}g_1)$.

Therefore, by the conclusion of Theorem $3.6$, there is nothing to
prove for (a), (b), (d2) and (d3) in this sublemma. Also, the proof
of (d1) is trivial, applying Corollary $3.8$ to the defining
equation of $g_j$ and the defining equation of $g_{j+1}$ in (5.4.6).

So, in order to prove (c) throughout this sublemma we can use the
same kind of notations and properties as in (a), (b), (d2) and
(d3) as follows:

For (a), (b) and (d2),  along $v=0$ $\tau_m:M^{(m)}\to\BC^2$ as a
composition of analytic mappings and $(f\circ\tau_m)_{total}$ can
be rewritten in the following form: Note that $2\le j\le r$.
$$\align
(5.4.8) \quad \quad  \tau_m(v,u)&=(y,z)=(v^{n_1}u^a,v^{\beta_{1,1}}u^b),  \\
(f\circ\tau_m)_{total}&=(f\circ\tau_m)(v,u)
=v^{e_{j,m}}u^{\rho_{j,m}}(f\circ\tau_m)_{proper}\quad
\text{with  $g_j=f$}, \qquad \qquad \\
(f\circ\tau_m)_{proper}&=(1+\ve_1u)^{d_j}+\sum_{\alpha,\beta\ge
0}c^{(j)}_{\alpha,\beta}v^{n_1\alpha+\beta_{1,1}\beta-n_1\beta_{1,1}d_j}
u^{a\alpha+b\beta-bn_1d_j}, \qquad \qquad \\
\endalign
$$
where \roster \item "(i)" $a$ and $b$ are some nonnegative integers
such that $a\beta_{1,1}-bn_1=1$ and $\ve_1$ is assumed to be one in
$\BC\{y,z\}$,

\item "(ii)" $e_{j,m}=n_1\beta_{1,1}d_j$ and $\rho_{j,m}=bn_1d_j$
and $\rho_{\alpha,\beta}=a\alpha+b\beta-bn_1d_j\ge 0$,

\item "(iii)" $E_m=\{v=0\}$ is defined by the $m-th$ exceptional
curve of the first kind.

\item "(iv)" $V^{(m)}(g_j)\cap(\cup^m_{i=1}E_i)=V^{(m)}(g_j)\cap
E_m=\{(v,1+\ve_1 u)=(0,0)\}$ for any $j=2,\dots,r$.

\endroster \ms

For (d3), after $m$ iterations of blow-ups, denoted by $\tau_m$,
we have the following consequences:
$$\align
(5.4.9) \qquad \qquad &\text{(i) \quad If $\beta_{1,1}=1$, then
$f\in \text{\rm the type[0]}$ under
$\tau_m$, and }\\
&\text{\qquad \quad if $\beta_{1,1}\ge 2$, then $f\in \text{\rm the
type[1]}$ under
$\tau_m$.} \\
&\text{(ii) \quad Whether $\beta_{1,1}=1$ or $\beta_{1,1}\ge 2$,
then $F\in the~ type[1]$ under $\tau_m$.} \qquad \qquad
\endalign$$

$\underline{\text{Remark 5.4.1}}$

(i) In the assumption of Theorem $3.6$, the construction for
$G(y,z)=z^{\gamma}g_1$ with $g_1=z^{n_1}+y^{k_1}$ was defined as
follows: Note that $\gcd(n_1,k_1)=1$.

Let $1\le n_1<k_1$, and if $n_1=1$, then $\gamma=1$, and if
$n_1\ge 2$, then $\gamma=0$.

(ii) In the conclusion of Theorem $3.6$, whether $n_1=1$ or $2\le
n_1<k_1$, or $2\le k_1 <n_1$, there are some nonnegative integers
$a$ and $b$ such that $bn_1-ak_1=1$ because of $[\text{\rm{I}}]$ and
$[\text{\rm{II}}]$ in Theorem $3.6$. \ms

In preparation for the proof of (c), apply the above conclusion with
(5.4.8), to $g_{j+1}=(z^{n_1}+\ve_1 y^{\beta_{1,1}})^{n_2n_3\cdots
n_{j+1}}+\sum c^{(j+1)}_{\alpha,\beta}y^{\alpha}z^{\beta}$ in
$(5.4.6)$. Then for any $j=1,2,\dots,r-1$, we have the following:
$$\align
(5.4.10) \quad  &(g_{j+1}\circ \tau_m)_{total} \\
=&v^{e_{j+1,m}}u^{\rho_{j+1,m}}(g_1\circ
\tau_m)^{d_{j+1}}_{proper}+\sum_{\alpha,\beta\ge 0}
c^{(j+1)}_{\alpha,\beta}(v^{n_1}u^a)^{\alpha}(v^{\beta_{1,1}}u^b)^{\beta}
\qquad \qquad \\
=&v^{e_{j+1,m}}u^{\rho_{j+1,m}}\{(g_1\circ
\tau_m)^{d_{j+1}}_{proper} +\sum_{\alpha,\beta\ge 0}
c^{(j+1)}_{\alpha,\beta}v^{n_1\alpha+\beta_{1,1}\beta-n_1\beta_{1,1}d_{j+1}}
u^{a\alpha+b\beta-bn_1d_{j+1}}\} \\
=&v^{e_{j+1,m}}u^{\rho_{j+1,m}}(g_{j+1}\circ \tau_m)_{proper},
\endalign$$
where $d_{j+1}=n_2n_3\cdots n_{j+1}$,
$e_{j+1,m}=n_1\beta_{1,1}d_{j+1}$ and $\rho_{j+1,m}=bn_1d_{j+1}$,
noting by (5.4.7) and (ii) of (5.4.8) that
$n_1\alpha+\beta_{1,1}\beta-n_1\beta_{1,1}d_{j+1}>0$ and
$a\alpha+b\beta-bn_1d_{j+1}\ge 0$.

On the other hand, recall that for $j=2,3,\dots,r-1$,
$$\align (5.4.11) \qquad \qquad \qquad
&g_{j+1}=g_j^{n_{j+1}}+\ve_{j+1}y^{\beta_{j+1,1}}z^{\beta_{j+1,2}}
g_1^{\beta_{j+1,3}}\cdots g_{j-1}^{\beta_{j+1,j+1}}, \qquad \qquad
\qquad \\
& g_j =g_1^{n_2\cdots n_j} +\sum_{\gamma,\delta\ge 0}
c^{(j)}_{\gamma,\delta}y^{\gamma}z^{\delta}, \\
& g_1 =z^{n_1}+\ve_1y^{\beta_{1,1}} \quad \text{with $\ve_1=1$}, \\
\text{where} \quad \quad \text{\rm(i)}
\quad&\Delta_2(\gamma,\delta)=n_1\gamma+\beta_{1,1}\delta
>n_1\beta_{1,1}n_2\cdots n_j
\quad  \qquad  \text{by (5.4.6)}, \\
\text{\rm(ii)} \quad
&\Delta^{\sharp}_{j+1}(\beta_{j+1,k})^{j+1}_{k=1}>n_1\beta_{1,1}n_2n_3\cdots
n_{j+1} \qquad  \quad \text{by (5.1.2)}, \\
\text{\rm(iii)} \quad
&\Omega^{\sharp}_{j+1}(\beta_{j+1,k})^{j+1}_{k=1}\ge
bn_1n_2n_3\cdots n_{j+1}
\qquad \qquad   \text{by (5.3.3)}, \\
\text{\rm(iv)} \quad &\text{$\ve_{j+1}$ is a unit in
$\BC\{y,z\}$.}
\endalign$$

Now, apply (5.4.8) and (5.4.10) to
$g_{j+1}=g_j^{n_{j+1}}+\ve_{j+1}y^{\beta_{j+1,1}}z^{\beta_{j+1,2}}
g_1^{\beta_{j+1,3}}\cdots g_{j-1}^{\beta_{j+1,j+1}}$ in $(5.4.11)$.
Then, we have the following:
$$\align
(5.4.12)   & \quad (g_{j+1}\circ \tau_m)_{total}\\
 &=((g_j\circ
\tau_m)(v,u))^{n_{j+1}} +{(\ve_{j+1}\circ \tau_m)(v,u)}{((y\circ
\tau_m)(v,u))^{\beta_{j+1,1}}}  \qquad \qquad \\
& \quad \times {((z\circ
\tau_m)(v,u))^{\beta_{j+1,2}}}{(({g_1}\circ
\tau_m)(v,u))^{\beta_{j+1,3}}}\cdots
{((g_{j-1}\circ \tau_m)(v,u))^{\beta_{j+1,j+1}}} \qquad \qquad \\
 &=\{v^{e_{j,m}}u^{\rho_{j,m}}(g_{j}\circ
 \tau_m)_{proper}\}^{n_{j+1}}
 +{\ve'_{j+1}}{(v^{n_1}u^a)^{\beta_{j+1,1}}} {(v^{\beta_{1,1}}u^b)^{\beta_{j+1,2}}} \\
& \quad \times
\{v^{n_1\beta_{1,1}}u^{bn_1}(1+\ve_1u)\}^{\beta_{j+1,3}}
\{v^{e_{2,m}}u^{\rho_{2,m}}(g_{2}\circ
\tau_m)_{proper}\}^{\beta_{j+1,4}}\times\cdots \\
& \quad \times \{v^{e_{j-1,m}}u^{\rho_{j-1,m}}
(g_{j-1}\circ \tau_m)_{proper}\}^{\beta_{j+1,j+1}} \\
&=v^{e_{j+1,m}}u^{\rho_{j+1,m}}\{(g_{j}\circ
\tau_m)^{n_{j+1}}_{proper}+{\ve'_{j+1}}
v^{\Delta^{\sharp}_{j+1}(\beta_{j+1,k})^{j+1}_{k=1}-e_{j+1,m}} \\
& \quad \times u^{\Omega^{\sharp}_{j+1}(\beta_{j+1,k})^{j+1}_{k=1}-b
n_1d_{j+1}}(1+\ve_1u)^{\beta_{j+1,3}}\cdots (g_{j-1}\circ
\tau_m)^{\beta_{j+1,j+1}}_{proper}\}  \quad \text{by (5.4.13)} \\
&=v^{e_{j+1,m}}u^{\rho_{j+1,m}}(g_{j+1}\circ \tau_m)_{proper} \quad
\text{by (5.4.10)},
\endalign$$
where \ (i) $\ve'_{j+1}={(\ve_{j+1}\circ \tau_m)(v,u)}$ is a unit in
$\BC\{v,1+u\}$ because $n_1>0$ and $\beta_{1,1}>0$,

\quad \quad (ii) if we write $d_{j}=n_2n_3\cdots n_{j}$,
$e_{j,m}=n_1\beta_{1,1}d_{j}$ and $\rho_{j,m}=bn_1d_{j}$ for

\qquad \qquad $j=2,3,\dots,r$, then $d_{j}n_{j+1}=d_{j+1}$,
$e_{j,m}n_{j+1}=e_{j+1,m}$ and $\rho_{j,m}n_{j+1}=\rho_{j+1,m}$. \ms

The proof for the representation in (5.4.12) just follows from
$(*)$ and $(**)$ of (5.4.13):

\noindent$(5.4.13)-(*)$ It is clear by Sublemma 5.2 that
$(n_1\beta_{j+1,1}+\beta_{1,1}\beta_{j+1,2})+n_1\beta_{1,1}\beta_{j+1,3}
+e_{2,m}\beta_{j+1,4}+\cdots+e_{j-1,m}\beta_{j+1,j+1}$

$=\Delta_2(\beta_{j+1,1},\beta_{j+1,2})+n_1\beta_{1,1}\beta_{j+1,3}
+n_1\beta_{1,1}d_2\beta_{j+1,4}+\cdots+n_1\beta_{1,1}d_{j-1}\beta_{j+1,j+1}$

$=\Delta_2(\beta_{j+1,1},\beta_{j+1,2})+n_1\beta_{1,1}\beta_{j+1,3}
+n_1\beta_{1,1}n_2\beta_{j+1,4}+\cdots+n_1\beta_{1,1}n_2 \cdots
n_{j-1}\beta_{j+1,j+1}$

$=\Delta^{\sharp}_{j+1}(\beta_{j+1,k})^{j+1}_{k=1}>e_{j+1,m}=n_1\beta_{1,1}d_{j+1}$
by (5.1.2) and by (5.4.10). \ms

\noindent$(5.4.13)-(**)$ It is clear by Sublemma 5.3 that
$a\beta_{j+1,1}+b\beta_{j+1,2}+bn_1\beta_{j+1,3}+\rho_{2,m}\beta_{j+1,4}+
+\cdots+\rho_{j-1,m}\beta_{j+1,j+1}$

$=\Omega_2(\beta_{j+1,1},\beta_{j+1,2})+bn_1\beta_{j+1,3}+bn_1d_2\beta_{j+1,4}
+\cdots+bn_1d_{j-1}\beta_{j+1,j+1}$

$=\Omega_2(\beta_{j+1,1},\beta_{j+1,2})+bn_1\beta_{j+1,3}+bn_1n_2\beta_{j+1,4}
+\cdots+bn_1n_2n_3\cdots n_{j-1}\beta_{j+1,j+1}$

$=\Omega^{\sharp}_{j+1}(\beta_{j+1,k})^{j+1}_{k=1}\ge
\rho_{j+1,m}=bn_1d_{j+1}$ by (5.3.2) and (5.3.3) and by (5.4.10).

Thus, we can prove that (c) and so, the proof of the sublemma are
finished. $\square$
\enddemo \ms

\demo{\bf Proof of Sublemma 5.5} First of all, let $\{Y_k:
k=1,2,\dots,r-1\}$ with $Y_k\subset N_0$, $\{h_k: k=1,2,\dots,r-1\}$
with $h_k=(g_{k+1}\circ\tau_m)_{proper}$ in $\BC\{v,1+u\}$ and
$\{\Xi_k: N^k_0\to N_0: k=1,2,\dots,r-1\}$, where each $\Xi_k$ is an
integer-valued function, be three sequences satisfying the given
three conditions, denoted by \text{\bf The 1-th
${\text{\bf{Cond}}}^{\text{{\bf(1)}}}$}, \text{\bf The 2-th
${\text{\bf{Cond}}}^{\text{{\bf(1)}}}$}, \text{\bf The 3-th
${\text{\bf{Cond}}}^{\text{{\bf(1)}}}$}, in the conclusion of this
sublemma.

After the proof of Sublemma $5.4$ was done, it is easy to observe
without any more need of the proof that the above three sequences
with three conditions are well-constructed.

For the proof of this sublemma, it suffices to show that these three
sequences satisfy the remaining two conditions, {\rm (i)} \text{\bf
The $(4\alpha)$-th ${\text{\bf{Cond}}}^{\text{{\bf(1)}}}$} with
\text{\bf The 4-th ${\text{\bf{Cond}}}^{\text{{\bf(1)}}}$} and
{\rm(ii)} \text{\bf The $(5\alpha)$-th
${\text{\bf{Cond}}}^{\text{{\bf(1)}}}$} in $\underline {\text{\rm
Conclusions}}$ of this sublemma.  So, by Remark $5.5.1$, it is
enough to show that the  equality in (5.5.$4\alpha$) of \text{\bf
The $(4\alpha)$-th ${\text{\bf{Cond}}}^{\text{{\bf(1)}}}$} is true.

Now, we will prove that the equality in (5.5.$4\alpha$) is true,
using the following three steps: Let $\ell$ and $q$ be arbitrary
positive integers such that $r-1\ge \ell\ge q\ge 2$.
$$\align
\text{$\underline{\text{\rm Step(i)}}$} \qquad  & \Xi_q(\gamma_{\ell, k})^q_{k=1} \\
=& \Delta_{q+1}(\beta_{\ell+1,k})^{q+1}_{k=1}+n^2_qn^2_{q-1}\cdots
n^2_2n_1\beta_{1,1}\{\beta_{\ell+1,q+2}+n_{q+1}\beta_{\ell+1,q+3} \\
& +n_{q+1}n_{q+2}\beta_{\ell+1,q+4}+\cdots +n_{q+1}n_{q+2}\cdots
n_{\ell-1}\beta_{\ell+1,\ell+1}-n_{q+1}n_{q+2}\cdots n_{\ell+1}\}.
\qquad \qquad \\
\text{$\underline{\text{\rm Step(ii)}}$} \qquad  & \text{In
particular, if $\ell=q$
then} \\
& \Xi_q(\gamma_{q,k})^q_{k=1}
=\Delta_{q+1}(\beta_{q+1,k})^{q+1}_{k=1}
 -n_{q+1}n^2_qn^2_{q-1}\cdots n^2_2n_1\beta_{1,1} \quad \text{from {\rm Step(i)}}.
\qquad \qquad \\
\text{$\underline{\text{\rm Step(iii)}}$} \qquad &
\Xi_q(\gamma_{q,k})^q_{k=1}-s_qs_{q-1}\Xi_{q-1}(\gamma_{q-1,k})^{q-1}_{k=1}\\
 =& \Delta_{q+1}(\beta_{q+1,k})^{q+1}_{k=1}
 -n_{q+1}n_q\Delta_q(\beta_{q,k})^q_{k=1}>0 \quad \text{from {\rm
 Step(ii)}}.
\endalign$$
\ms

We prove {\rm Step(i)}, {\rm Step(ii)} and {\rm Step(iii)}
simultaneously by induction on the integer $q\ge 2$.

So, it is enough to consider two cases, respectively:

Case(I) $q=2$, and Case(II) $q\ge 2$. \ms

\noindent{\bf Case(I):} \ Let $q=2$. Note by \text{\bf The 3-th
${\text{\bf{Cond}}}^{\text{{\bf(1)}}}$} that
$\Xi_2(t_1,t_2)=t_2\Xi_1(\gamma_{1,1})+s_1\Xi_1(t_1)
=t_2\gamma_{1,1}+s_1t_1$ for each $(t_1,t_2)\in N^2_0$.
$$\align
\text{$\underline{\text{\rm Step(i)}}$} \quad  &
\Xi_2(\gamma_{\ell,1},\gamma_{\ell,2})
 = s_1\gamma_{\ell,1}+\gamma_{1,1}\gamma_{\ell,2}
 \quad \text{by definition of} \ {\Xi_{2}}\\
=& n_2
\{\Delta^{\sharp}_{\ell+1}(\beta_{\ell+1,k})^{\ell+1}_{k=1}-n_1\beta_{1,1}n_2\cdots
n_{\ell+1}\} \\
& +\{\Delta^{\sharp}_{2}(\beta_{2,1},\beta_{2,2})-n_1\beta_{1,1}n_2
\}\beta_{\ell+1,3}
\qquad \qquad \qquad \quad \text{by (5.5.1)}\qquad \\
=& n_2
\{\Delta_2(\beta_{\ell+1,1},\beta_{\ell+1,2})+n_1\beta_{1,1}\beta_{\ell+1,3}
+n_1\beta_{1,1}n_2\beta_{\ell+1,4}+\cdots  \\
&  +n_1\beta_{1,1}n_2\cdots n_{\ell-1}\beta_{\ell+1,\ell+1}
-n_1\beta_{1,1}n_2\cdots n_{\ell+1} \} \\
&+ \{\Delta_2(\beta_{2,1},\beta_{2,2})-n_1\beta_{1,1}n_2
\}\beta_{\ell+1,3}
\qquad \qquad \qquad \quad \text{by (5.1.1)}\\
=&\Delta_3(\beta_{\ell+1,1},\beta_{\ell+1,2},\beta_{\ell+1,3})+n^2_2n_1\beta_{1,1}
\{\beta_{\ell+1,4}+n_3\beta_{\ell+1,5}+n_3n_4\beta_{\ell+1,6} \qquad \qquad \\
& + \cdots +n_3n_4\cdots
n_{\ell-1}\beta_{\ell+1,\ell+1}-n_3n_4\cdots n_{\ell+1}\},
\endalign$$
by the definition of
$\Delta_3(\beta_{\ell+1,1},\beta_{\ell+1,2},\beta_{\ell+1,3})$
only, which implies the proof of {\rm Step(i)}. \ms

\noindent $\underline{\text{\rm Step(ii)}}$ \quad In particular,
if $\ell=2$ then by {\rm Step(i)}
$$
\Xi_2(\gamma_{2,1},\gamma_{2,2})=\Delta_3(\beta_{3,1},\beta_{3,2},\beta_{3,3})
-n_3n^2_2n_1\beta_{1,1}.
$$

Thus, the proof of {\rm Step(ii)} is done. \ms

\noindent $\underline{\text{\rm Step(iii)}}$ To prove that
$\Xi_2(\gamma_{2,1},\gamma_{2,2})-s_2s_1\Xi_1(\gamma_{1,1})>0$,
first note by (5.5.1) that
$$\align
&s_2s_1\Xi_1(\gamma_{1,1})=s_2s_1\gamma_{1,1}=n_3n_2\{\Delta_2(\beta_{2,1},\beta_{2,2})
-n_1\beta_{1,1}n_2\}. \\
\text{Then,} \qquad &\Xi_2(\gamma_{2,1},\gamma_{2,2})-s_2s_1\gamma_{1,1} \\
=&\Delta_3(\beta_{3,1},\beta_{3,2},\beta_{3,3})-n_3n^2_2n_1\beta_{1,1}
-n_3n_2\{\Delta_2(\beta_{2,1},\beta_{2,2})-n_1\beta_{1,1}n_2\}  \qquad \qquad \\
=&\Delta_3(\beta_{3,1},\beta_{3,2},\beta_{3,3})-n_3n_2\Delta_2(\beta_{2,1},\beta_{2,2})>0,
\endalign$$
by   in the assumption of Theorem $5.0$, which implies the proof
of {\rm Step(iii)}.

Thus, if $q=2$, then we proved that this sublemma is true. \ms

\noindent{\bf Case(II):} \ Let $q\ge 2$. By induction proof,
suppose that the sublemma is true on the integer $q\le r-2$ with
$r-1\ge \ell\ge q$. Then, it is enough to prove {\rm Step(i)},
{\rm Step(ii)} and {\rm Step(iii)} simultaneously on the integer
$(q+1)\le \ell$ as follows:
$$\align
\text{$\underline{\text{\rm Step(i)}}$} \qquad  &
\Xi_{q+1}(\gamma_{\ell, k})^{q+1}_{k=1}
=\gamma_{\ell,q+1}\Xi_q(\gamma_{q,k})^q_{k=1}+s_q\Xi_q(\gamma_{\ell,
k})^q_{k=1}
\quad \text{by definition of \ $\Xi_{q+1}$} \\
=&\beta_{\ell+1,q+2}
\{\Delta_{q+1}(\beta_{q+1,k})^{q+1}_{k=1}-n_{q+1}n^2_qn^2_{q-1}\cdots
n^2_2n_1\beta_{1,1} \} \\
& + n_{q+1} \{ \Delta_{q+1}(\beta_{\ell+1,k})^{q+1}_{k=1}
+n^2_qn^2_{q-1}\cdots n^2_2n_1\beta_{1,1}(\beta_{\ell+1,q+2}+n_{q+1}\beta_{\ell+1,q+3} \\
& +n_{q+1}n_{q+2}\beta_{\ell+1,q+4}+\cdots+n_{q+1}n_{q+2}\cdots
n_{\ell-1}\beta_{\ell+1,\ell+1}-n_{q+1}n_{q+2}\cdots n_{\ell+1})
\} \\
& \quad   \text{by the induction assumption on the integer $q$}  \qquad \qquad   \\
=&\Delta_{q+2}(\beta_{\ell+1,k})^{q+2}_{k=1}+n^2_{q+1}n^2_q
n^2_{q-1}\cdots
n^2_2n_1\beta_{1,1} \\
& \times\{\beta_{\ell+1,q+3}+n_{q+2}\beta_{\ell+1,q+4}
 +n_{q+2}n_{q+3}\beta_{\ell+1,q+4}+\cdots  \\
& \quad +n_{q+2}n_{q+3}\cdots
n_{\ell+1-2}\beta_{\ell+1,\ell+1}-n_{q+2}n_{q+3}\cdots n_{\ell+1}
\},
\endalign$$
by  the definition of $\Delta_{q+2}(\beta_{\ell+1,k})^{q+2}_{k=1}$
only, which implies the proof of Step(i). \ms

$\underline{\text{Step(ii)}}$ \quad In particular, if $\ell=q+1$,
then $\ell+1=q+2<q+3$ and so
$$
\Xi_{q+1}(\gamma_{q+1,k})^{q+1}_{k=1}=\Delta_{q+2}(\beta_{q+2,k})^{q+2}_{k=1}
-n_{q+2}n^2_{q+1}n^2_q\cdots n^2_2n_1\beta_{11},
$$
by Step(i) on the integer $q+1$, which implies the proof of Step(ii)
on $q+1$. \ms

$\underline{\text{Step(iii)}}$ \quad To prove that the equality in
$(5.5.4\alpha)$ is true, as an application of Step(ii) on two
integers $\ell=q+1$ and $\ell=q$,  then we have
$$\align \Xi_{q+1} &(\gamma_{q+1,k})^{q+1}_{k=1}
-s_{q+1}s_q\Xi_q(\gamma_{q,k})^q_{k=1} \\
&=\{\Delta_{q+2}(\beta_{q+2,k})^{q+2}_{k=1}-n_{q+2}n^2_{q+1}n^2_q\cdots
n^2_2n_1\beta_{1,1}\} \\
&\quad
-n_{q+2}n_{q+1}\{\Delta_{q+1}(\beta_{q+1,k})^{q+1}_{k=1}-n_{q+1}n^2_qn^2_{q-1}\cdots
n^2_2n_1\beta_{1,1}\} \\
&=\Delta_{q+2}(\beta_{q+2,k})^{q+2}_{k=1}-n_{q+2}n_{q+1}\Delta_{q+1}
(\beta_{q+1,k})^{q+1}_{k=1}>0,
\endalign$$
by  \text{\bf The $4\alpha$-th
${\text{\bf{Cond}}}^{\text{{\bf(0)}}}$} in the assumption of Theorem
$5.0$, which implies the proof of Step(iii) on the integer $q+1$.

Thus, we proved that the equality in $(5.5.4\alpha)$ is true, and so
the proof of this sublemma is finished. $\square$ \bs

{\bf \S{6.2}. The proofs of Theorem 5.0 with corollaries}

\demo{\bf Proof of Theorem 5.0} For the induction proof, it suffices
to consider two cases, respectively:

Case(1)  $r=1$, and Case(2) $r\ge 1$.

For each case, we shall find first, the proof of $[A]$, and after
then, the the proof of $[B]$. \ms

{\bf Case(1):} Let $r=1$. Then, there is nothing to prove for $[A]$.
For the proof of $[B]$, assuming that $g_1=z^{n_1}+y^{\beta_{1,1}}$
is irreducible in $\BC\{y,z\}$ with $n_1\ge 2$, then
$y^{\gamma}g_1\in \text{the type[1]}$ under the standard resolution
by Sublemma $5.4$ or Theorem $3.6$, because if $\beta_{1,1}=1$ then
$\gamma=1$ and if $\beta_{1,1}>1$ then $\gamma=0$. In particular,
$z^{\delta}yg_1\in \text{the type[1]}$ under the standard resolution
by Sublemma $5.4$ or Theorem $3.6$ whether $\delta=1$ or $\delta=0$.
Thus, $[B]$ can be easily proved. So, if $r=1$, the proofs of $[A]$
and $[B]$ are done. \ms

{\bf Case(2):} Let $r\ge 1$. For each integer $r\ge 2$, we shall
find first, the proof of $[A]$ and next, the the proof of $[B]$. If
$g_j$ is irreducible in $\BC\{y,z\}$ for any $j\ge 2$, then
$\gcd(n_1,\beta_{1,1})=1$ by (d) of Sublemma $5.2$ and Theorem
$3.6$. So, to find the proofs of $[A]$ and $[B]$ for each $r\ge 2$,
we may assume without proof that $\gcd(n_1,\beta_{1,1})=1$. Since
$\gcd(n_1,\beta_{1,1})=1$, for the proof of the theorem we may
follow the same notations and consequences as in Sublemma $5.4$.

Therefore, as we have done in Sublemma $5.4$, recall that
$\tau_m:M^{(m)}\to\BC^2$ and $(g_j\circ\tau_m)_{total}$ with
$(g_j\circ\tau_m)_{proper}$ satisfies all facts in {\rm(a)},
{\rm(b)}, {\rm(c)} and {\rm(d)} in the conclusion of Sublemma $5.4$:
Note that $\tau_m:M^{(m)}\to\BC^2$ is the composition of a finite
number $m$ of successive blow-ups $\pi_i$ which is needed to get the
standard resolution of the singular point of $V(y^{\gamma}g_1)$ in
the conclusion of Sublemma $5.4$. Also, for any positive integer
$r\ge 2$, observe by $(c)$ of Sublemma $5.4$ that $g_r$ is
irreducible in $\BC\{y,z\}$ if and only if $g_1$ is irreducible in
$\BC\{y,z\}$ and $(g_r\circ\tau_m)_{proper}$ is irreducible in
$\BC\{v,u+1\}$, noting that $g_1$ is irreducible in $\BC\{y,z\}$ if
and only if $\gcd(n_1,\beta_{1,1})=1$.

Now, we consider the proof of Case(2). For the proof of the theorem,
by the induction assumption, suppose we have shown that [A] and [B]
of the theorem are true on the integer $r\ge 1$. In order to prove
the theorem on the integer $(r+1)$, first we will find the proof of
[A] on the integer (r+1) and next, the proof of [B] on the integer
(r+1). \ms

$\underline{\text{\rm The proof of [A]}}$. To prove [A] on the
integer $r+1$, assume that $\gcd(n_1,\beta_{1,1})=1$ and
$g_{r+1}\in\BC\{y,z\}$ is $\underline{\text{\rm a semi-quasi-Puiseux
convergent power series of the recursive (r+1)-type}}$.

To prove the theorem, it suffices to consider the following defining
equation for $g_{r+1}$ with the additive assumptions:
$$
g_{r+1}=g^{n_{r+1}}_r+\ve_{r+1}y^{\beta_{r+1,1}}z^{\beta_{r+1,2}}
g^{\beta_{r+1,3}}_1\cdots g^{\beta_{r+1,r+1}}_{r-1}, \tag 5.0.1
$$
where

\noindent (i) \quad $g_1,g_2,\dots,g_r$ satisfies the same
assumptions and conclusions as in this theorem,

\noindent (ii) \quad
$X_{r+1}=\{n_{r+1},\beta_{r+1,1},\beta_{r+1,2},\dots,\beta_{r+1,r+1}\}\subset
N_0$ with $n_{r+1}\ge 2$ and $\ve_{r+1}=\ve_{r+1}(y,z)$ is a unit
in $\BC\{y,z\}$,

\noindent (iii) \quad
$\Delta_{r+1}(t_k)^{r+1}_{k=1}=t_{r+1}\Delta_r(\beta_{r,k})^r_{k=1}
+n_r\Delta_r(t_k)^r_{k=1}$ for each $(t_k)^{r+1}_{k=1}\in
N^{r+1}_0$,
$$\align
&\text{\rm(iv)} \qquad \quad
\Delta_{r+1}(\beta_{r+1,k})^{r+1}_{k=1}>n_{r+1}n_r\Delta_r(\beta_{r,k})^r_{k=1},
\qquad \qquad \\
&\text{\rm(v)} \qquad \quad
\Delta^{\sharp}_{r+1}(\beta_{r+1,k})^{r+1}_{k=1}
 =\Delta_2(\beta_{r+1,1},\beta_{r+1,2})+n_1\beta_{1,1}\beta_{r+1,3}
 +n_1\beta_{1,1}n_2\beta_{r+1,4} \qquad \qquad \\
 &\qquad \qquad \qquad \qquad \qquad \qquad \quad +n_1\beta_{1,1}
 n_2 n_3 \beta_{r+1,5}
 + \cdots +n_1\beta_{1,1}n_2 \cdots n_{r-1}\beta_{r+1,r+1}, \\
&\text{\rm(vi)} \qquad \quad
\Omega^{\sharp}_{r+1}(\beta_{r+1,k})^{r+1}_{k=1}
=\Omega_2(\beta_{r+1,1},\beta_{r+1,2})+bn_1\beta_{r+1,3}+bn_1n_2\beta_{r+1,4} \\
 &\qquad \qquad \qquad \qquad \qquad \qquad \quad
 +bn_1n_2 n_3 \beta_{r+1,5}
 +\cdots +bn_1n_2\cdots n_{r-1}\beta_{r+1,r+1}. \\
\endalign$$

\text{\bf Remark 5.0.1.1.} Note by (5.0.1) that the following
inequalities in \text{\rm(a)} and \text{\rm(b)} can be proved by
Sublemma $5.1$ and Sublemma $5.3$, respectively, using the same
methods as we have used in the proof of Sublemma $5.1$ and Sublemma
$5.3$.
$$\align
\Delta^{\sharp}_{r+1}(\beta_{r+1,k})^{j+1}_{k=1}
&> n_1\beta_{1,1}n_2n_3\cdots n_{r}n_{r+1}. \tag {a} \\
\Omega^{\sharp}_{r+1}(\beta_{r+1,k})^{r+1}_{k=1} &\ge
bn_1n_2n_3\cdots n_{r}n_{r+1}. \tag {b} \\
\endalign$$ \ms

Using the composition $\tau_m$ again of a finite number $m$ of
successive blow-ups which is needed to get the standard resolution
of the singular point of $V(y^{\gamma}g_1)$ as in the beginning of
the proof, then $h_r=(g_{r+1}\circ\tau_m)_{proper}$ can be written
in the form
$$
h_r=h^{s_r}_{r-1}+\eta_r
v^{\gamma_{r,1}}(1+u)^{\gamma_{r,2}}h^{\gamma_{r,3}}_1 \cdots
h^{\gamma_{r,r}}_{r-2}, \tag 5.0.2
$$
where

\noindent(i) each $h_j=(g_{j+1}\circ\tau_m)_{proper}$ is in
$\BC\{v,1+u\}$ for $j=1,2,\dots,r$, which has been already
represented by Sublemma $5.5$,

\noindent(ii)
$(g_{r+1}\circ\tau_m)_{total}=v^{n_1\beta_{1,1}n_2\cdots
 n_{r+1}}u^{bn_1n_2\cdots n_{r+1}}(g_{r+1}\circ\tau_m)_{proper}$,

\noindent(iii) $\eta_r=\ve'_{r+1}
u^{\Omega^{\sharp}_{r+1}(\beta_{r+1,k})^{r+1}_{k=1}-bn_1n_2\cdots
n_{r+1}}$ is a unit in $\BC\{v,u+1\}$, satisfying the same
notations and consequences as in Sublemma $5.5$,

\noindent(iv)  $s_{r} =n_{r+1}\ge 2$,
$\gamma_{r,1}=\Delta^{\sharp}_{r+1}(\beta_{r+1,k})^{r+1}_{k=1}
 -n_1\beta_{1,1}n_2n_3\cdots n_{r+1}>0$,
$\gamma_{r,2}=\beta_{r+1,3}$, $\gamma_{r,3}=\beta_{r+1,4}$,
\dots,$\gamma_{r,r}=\beta_{r+1,r+1}$. \ms

Because blow-ups process preserves irreducibility of plane curve
singularity, note by Sublemma 5.2 and Sublemma 5.4 that $g_{r+1}$ is
irreducible in $\BC\{y,z\}$ if and only if $g_1$ is irreducible in
$\BC\{y,z\}$ and $h_r=(g_{r+1}\circ\tau_m)_{proper}$ is irreducible
in $\BC\{v,u+1\}$.

Note that $g_1$ is irreducible in $\BC\{y,z\}$ if and only if
$\gcd(n_1,\beta_{1,1})=1$. First of all, whether or not $h_r$ is
irreducible in $\BC\{v,u+1\}$, we proved by Sublemma $5.5$ that
$h_r$ satisfies the same kind of assumptions as $g_r$ does the
assumptions in this theorem, up to change of notations.

So, by the induction assumption on the integer $r$ and by
following the same notations as in Sublemma $5.5$, then it is easy
to get the following:
$$\split
(5.0.3)\quad  &\text{$h_r$ is irreducible in $\BC\{v,u+1\}$} \iff
\gcd(s_j,\Xi_j(\gamma_{j,k})^j_{k=1})=1  \quad \text{for
$j=1,2,\dots, r$}.
\endsplit$$

Also, by \text{\bf The $(5\alpha)$-th
${\text{\bf{Cond}}}^{\text{{\bf(1)}}}$} of Sublemma $5.5$,  we have
the following:
$$\align
(5.0.4) \qquad &\gcd(s_j,\Xi_j(\gamma_{j,k})^j_{k=1})=1 \iff
\gcd(n_{j+1},\Delta_{j+1}(\beta_{j+1,k})^{j+1}_{k=1})=1. \quad
\text{for $j=1,2,\dots, r$}.
\endalign$$

Since  it is clear by Sublemma 5.4 that $g_{r+1}$ is irreducible in
$\BC\{y,z\}$ if and only if $\gcd(n_1,\beta_{11})=1$ and
$h_r=(g_{r+1}\circ\tau_m)_{proper}$ is irreducible in
$\BC\{v,u+1\}$, then we can prove by (5.0.3) and (5.0.4) and by the
induction assumption that
$$\align
(5.0.5) \qquad \qquad &\text{$g_{r+1}$ is irreducible in $\BC\{y,z\}$} \\
\iff \quad  & \text{$\gcd(n_1,\beta_{1,1})=1$ and \quad
$h_r$ is irreducible in $\BC\{v,u+1\}$} \\
\iff \quad  & \text{$\gcd(n_1,\beta_{1,1})=1$ and
$\gcd(n_{j+1},\Delta_{j+1}(\beta_{j+1,k})^{j+1}_{k=1})=1$ for $1\le j\le r$.}\\
\iff \quad &\text{$g_1, g_2,\dots,g_r$ are irreducible in
$\BC\{y,z\}$ and
$\gcd(n_{r+1},\Delta_{r+1}(\beta_{r+1,k})^{r+1}_{k=1})=1$.}
\endalign$$

Thus, we finished the proof of $[A]$ on the integer (r+1). \ms

$\underline{\text{\rm The proof of [B]}}$ \quad Next to prove $[B]$
on the integer $(r+1)$, let $g_{r+1}$ be irreducible in
$\BC\{y,z\}$. If $g_j$ is irreducible in $\BC\{y,z\}$ for any $j\ge
1$, note that $\gcd(n_1,\beta_{11})=1$ by (d) of Sublemma $5.2$ and
Theorem $3.6$. By Sublemma $5.4$, we have the following:

\noindent(5.0.6) \quad \text{(i) \  If $\beta_{1,1}>1$, then
$g_{r+1}\in$ the type $[1]$ under $\tau_m$.}

\qquad \quad \text{ (ii) \ If $\beta_{1,1}=1$, then $g_{r+1}\in$ the
type $[0]$ under $\tau_m$.}

\qquad \quad \text{ (iii) \ If $\beta_{1,1}\ge 1$, then
$z^{\delta}yg_{r+1}\in$ the type $[1]$ under $\tau_m$ whether
$\delta$ is $1$ or $0$.} \ms

Let $V_{r+1}=\{(y,z):g_{r+1}(y,z)=0\}$ and
$W_{r+1}=\{(y,z):yg_{r+1}(y,z)=0\}$ be analytic varieties at the
origin in $\BC\{y,z\}$, respectively.  Let $\tau_m$ be again the
composition of a finite number $m$ of successive blow-ups which is
needed to get the standard resolution of the singular point of
$V(yg_1)$.

At $(v,1+u)=(0,0)$, by (5.4.1) $\tau^{-1}_m(V_{r+1})$ and
$\tau^{-1}_m(W_{r+1})$ can be written as follows:
$$\align
(5.0.7) \quad \quad\qquad   \tau^{-1}_m(V_{r+1})
&=\{(g_{r+1}\circ\tau_m)_{total}=v^{n_1\beta_{1,1}d'_2}
u^{bn_1d'_2}(g_{r+1}\circ\tau_m)_{proper}=0\}, \qquad \qquad \\
  \tau^{-1}_m(W_{r+1}) &=
\{((yg_{r+1})\circ\tau_m)_{total}=v^{n_1\beta_{1,1}d'_2+n_1}
u^{bn_1d'_2+a}(g_{r+1}\circ\tau_m)_{proper}=0\}, \qquad \qquad \\
& \text{noting by (5.0.2) that} \\
(5.0.8) \quad
(g_{r+1}\circ\tau_m)_{proper}&=(g_r\circ\tau_m)^{n_{r+1}}_{proper}
+\eta_rv^{\Delta^{\sharp}_{r+1}(\beta_{r+1,k})^{r+1}_{k=1}-n_1\beta_{1,1}d'_2}\\
&\qquad
\times(g_1\circ\tau_m)^{\beta_{r+1,3}}_{proper}(g_2\circ\tau_m)^{\beta_{r+1,4}}_{proper}
\cdots (g_{r-1}\circ\tau_m)^{\beta_{r+1,r+1}}_{proper}. \qquad
\qquad
\endalign$$
where $d'_2=n_2n_3\cdots n_{r+1}$ and $\eta_r=\ve'_{r+1}
u^{\Omega^{\sharp}_{r+1}(\beta_{r+1,k})^{r+1}_{k=1}-bn_1d'_2}$ is a
unit in $\BC\{v,u+1\}$.

Let $V(\psi_r)=\{(v,u+1):\psi_r(v,u+1)=0\}$ be an analytic variety
at $(v,u+1)=(0,0)$ defined by
$$\align
(5.0.9) \qquad \qquad &\psi_r=\psi_r(v,u+1)=v^{\gamma_1}h_r \quad
\text{with \quad $h_r=(g_{r+1}\circ\tau_m)_{proper}$} \qquad \qquad \qquad \qquad \\
& \text{such that}  \quad \left\{\eqalign{& \text{$\gamma_1=1$,
\quad  if \quad $\Delta_2(\beta_{2,1},
\beta_{2,2})-n_1\beta_{1,1}n_2=1$}, \cr & \text{$\gamma_1=0$, \quad
if \quad $\Delta_2(\beta_{2,1},
\beta_{2,2})-n_1\beta_{1,1}n_2\ge 2$.} \cr} \right. \\
\endalign$$ \ms

Let $Z_r=\{(v,u+1):vh_r(v,u+1)=0\}$ be an analytic variety at
$(v,u+1)=(0,0)$. Since $(g_{r+1}\circ\tau_m)_{proper}$ is
irreducible in $\BC\{v,u+1\}$, then note that
$Z_{r}=\tau^{-1}_m(V_{r+1})=\tau^{-1}_m(W_{r+1})$ have the same
two irreducible components at $(v,u+1)=(0,0)$ as reduced
varieties, and so they have the same standard resolution of the
singular point $(v,u+1)=(0,0)$. Let $\omega_s$ be the composition
of $s$ iterations of blow-ups which is needed to get the standard
resolution of the singular point $(v,u+1)=(0,0)$ of $V(\psi_r)$.

Since we proved by Sublemma $5.5$ that $V(h_r)$ satisfies the same
kind of properties and notations as $V(g_r)$ does in the
assumption of this theorem, then by the induction assumption on
the integer $r$, both $V(\psi_r)$ and $Z_{r}$ belong to the type
$[r]$ under $\omega_s$. That is, both $\tau^{-1}_m(V_{r+1})$ and
$\tau^{-1}_m(W_{r+1})$ belong to the type $[r]$ under $\omega_s$
at $(v,u+1)=(0,0)$, as reduced varieties.

Now, the $m$-th exceptional curve of the first kind among
$\tau^{-1}_m(0,0)$, denoted by $E_m=\{v=0\}$, can be viewed as one
of two irreducible components of $Z_r$ at $(v,u+1)=(0,0)$. Then,
both $V_{r+1}$ and $W_{r+1}$ can have the same standard resolution
$\tau_m\circ\omega_s$ of the singular point $(y,z)=(0,0)$, using the
composition $\tau_m\circ\omega_s$ of a finite number $(m+s)$ of
successive blow-ups at $(y,z)=(0,0)$. Since $Z_r\in$ the type $[r]$
under $\omega_s$, then by (5.0.6) or Sublemma $5.4$ again, it is
trivial to get the followings:

\noindent(5.0.10) \quad \text{(i) \ If $\beta_{11}>1$, then
$g_{r+1}\in$ the type $[r+1]$ under $\tau_m\circ\omega_s$.}

\qquad \quad \text{ (ii) If $\beta_{11}=1$, then $g_{r+1}\in$ the
type $[r]$ under $\tau_m\circ\omega_s$.}

\qquad \quad \text{ (iii) If $\beta_{11}\ge 1$, then
$z^{\delta}yg_{r+1}\in$ the type $[r+1]$ under
$\tau_m\circ\omega_s$, whether $\delta$ is $1$ or $0$.} \ms

Thus, the proof for {\rm[B]} is done. Therefore, we completed the
proof of theorem. $\square$
\enddemo \ms

The proofs of Corollary 5.6, Corollary 5.7 and Corollary 5.8 just
follow from Theorem $5.0$ and Sublemmas. \ms

\vfill \pagebreak

\ms \vfill \pagebreak

{\bf \S7. For any two Puiseux convergent power series $g_r$ and
$\phi_\rho$ of the recursive types in $\BC\{y,z\}$, how to compute
the necessary and sufficient condition for \text{$\phi_\rho
\buildrel \text{{\rm divisor}} \over \sim g_r$} under the standard
resolutions and their classifications} \bs

{\bf \S7.0. Introduction}  \ms

As we have seen in Definition $1.2$ and Definition $2.4$, four
families with equivalence relations,
$\underline{\text{\rm{Family(1)}}}$, $\dots$,
$\underline{\text{\rm{Family(4)}}}$ with equivalence relations have
been already defined.

\definition{Definition 7.0} $\underline{\text{\rm{Quasi-Family(1)}}}$
is a family consisting of all the quasi-Puiseux convergent power
series \text{$f \in \BC\{y,z\}$} of the recursive type with isolated
singularity at the origin, denoted by $\underline{\text{\{\text{\rm
$f$ is arbitrary quasi-Puiseux convergent power series of the
recursive type}: $f\in \text{\rm {Family(0)}}$\}}}$ by Definition
$5.0.0$. Note that $\underline{\text{\rm{Family(1)}}}$ is a subset
of $\underline{\text{\rm{Quasi-Family(1)}}}$. \enddefinition

{\bf Problem[1].}\quad To solve the problem is to succeed in The 1st
algorithm and its application in $\S1.4$. The solution of the
problem can be divided into three small problems, denoted by
Problem[1-A], Problem[1-B] and Problem[1-C] in order.

{\bf Problem[1-A].} \quad In preparation, the first small problem is
to prove that there is a one-to-one function from
$\underline{\text{\rm{Family(1)}}}$ into
$\underline{\text{\rm{Family(j)}}}$ for any $j=2,3,4$. It was
already proved by Theorem A([K]) that there is a one-to-one function
$\phi$ from $\underline{\text{\rm{Family(2)}}}$ onto
$\underline{\text{\rm{Family(3)}}}$. Then, this problem can be
solvable by {\rm(i)} and {\rm(ii)}:

{\bf(i)} It will be proved by Theorem $7.3$ and Theorem $7.7$ in $\S
7$ that there exists a one-to-one function from
$\underline{\text{\rm{Family(1)}}}$ into
$\underline{\text{\rm{Family(4)}}}$.

{\bf(ii)} It will be proved by Theorem $10.2$ and Theorem $7.7$ that
there exists a one-to-one function from
$\underline{\text{\rm{Family(1)}}}$ into
$\underline{\text{\rm{Family(3)}}}$.

{\bf Problem[1-B].} \quad The second small problem is to prove that
we can compute a one-to-one map $\phi$ from
$\underline{\text{\rm{Family(1)}}}$ into
$\underline{\text{\rm{Family(2)}}}$. As an application, for given
any standard Puiseux polynomial $g_r\in \BC\{y\}[z]$ of the
recursive r-type, we compute an algorithm for finding the standard
Puiseux expansion $C(t)$ such that $g_r$ and $C(t)$ have the same
multiplicity sequence by Theorem $11.2$ and Corollary $11.3$.

{\bf Problem[1-C].} \quad Whenever the standard Puiseux expansion
$C(t)$ is chosen arbitrary, we can compute an algorithm for finding
the standard Puiseux polynomial $g_r\in \BC\{y\}[z]$ of the
recursive r-type such that $g_r$ and $C(t)$ have the same
multiplicity sequence by Theorem $11.4$, which implies that $\phi:
{\text{\rm{Family(1)}}}\rightarrow {\text{\rm{Family(2)}}}$ is onto.
\ms

Now, in this section the aim is to solve (i) of Problem[1-A] by
Theorem $7.3$ and Theorem $7.7$ completely. First of all, we are
preparing Theorem $7.1$. \ms

{\bf \S7.1. In preparation for solving (i) of Problem[1-A] by
Theorem 7.1}

In this section, the aim is to solve the following in terms of
Theorem 7.1:
$$\align
(*) \quad &\text{If $g_r\in \BC\{y,z\}$ is any quasi-Puiseux
convergent power series of the recursive r-type} \\
&\text{and $\psi_s\in \BC\{y,z\}$ is any quasi-Puiseux convergent
power series of the
recursive s-type, }\\
&\text{then we can compute a necessary condition for $V(g_r)
\buildrel \text{{\rm divisor}} \over \sim
V(\psi_s)$ under the standard} \\
&\text{resolutions in the sense of Definition 2.4.}
\endalign$$

For the study of Theorem $7.1$, it is clear by Sublemma $5.1$ and
Theorem $3.6$ that $f=g_r\in \BC\{y,z\}$ with
$F=y^{\zeta}z^{\eta}g_r\in \BC\{y,z\}$ and $\phi=\psi_s\in
\BC\{y,z\}$ with $\Phi=y^{\zeta'}z^{\eta'}\psi_s\in \BC\{y,z\}$ in
$(*)$ can be defined by the same properties and notations as in the
assumption of Theorem $7.1$. Then, we will solve the following in
terms of Theorem 7.1:

(a) The first aim is to find a necessary condition for
$V(y^{\zeta}z^{\eta}g_r) \buildrel \text{{\rm divisor}} \over \sim
V(y^{\zeta'}z^{\eta'}\psi_s)$ under the standard resolutions by
Theorem $7.1$ where $y^{\zeta}z^{\eta}g_r\in \BC\{y,z\}$ and
$y^{\zeta'}z^{\eta'}\psi_s\in \BC\{y,z\}$ are defined just as above,
by using Theorem $7.1$.

(b) As an application of Theorem $7.1$, the second aim is to find
the necessary and sufficient condition for $V(y^{\zeta}z^{\eta}g_r)
\buildrel \text{{\rm divisor}} \over \sim
V(y^{\zeta'}z^{\eta'}\psi_s)$ under the standard resolutions in an
elementary and concrete way, by using Theorem $7.3$. \ms

\proclaim{Theorem 7.1} $\underline{\text{\bf {Assumptions}}}$

{\rm (a)} Let $f=f(y,z)$ be irreducible in $\BC\{y,z\}$, and let
$V(F)=\{(y,z): F(y,z)=0\}$ be an analytic variety at $(0,0)$ in
$\BC^2$, each of which is written respectively in the form,
$$\align
\text{\rm (7.1.1)} \qquad \qquad f&=(z^{n_1}+\ve
y^{\beta_{1,1}})^{d} +\sum_{\alpha,\beta\ge
0}c_{\alpha,\beta}y^{\alpha}z^{\beta} \quad \text{with} \quad
n_1\alpha+\beta_{1,1}\beta>n_1\beta_{1,1}d,
\qquad \qquad \qquad \\
F&=y^{\zeta}z^{\eta}f
\endalign$$
satisfying the properties {\rm(i)}, {\rm(ii)},\dots, {\rm(vi)}:

\roster \item "(i)" $\gcd(n_1,\beta_{1,1})=1$ with $1\le
n_1<\beta_{1,1}$ and $d$ is a positive integer.

\item "(ii)" $\ve$ is a unit in $\BC\{y,z\}$.

\item "(iii)" The $c_{\alpha,\beta}$ are nonzero complex numbers
for some nonnegative integers $\alpha$ and $\beta$ such that
$n_1\alpha+\beta_{1,1}\beta>n_1\beta_{1,1}d$.

\item "(iv)" $\zeta$ and $\eta$ are nonnegative integers.

\item "(v)" If $n_1=1$, assume that $\eta$ is a positive integer.

\item "(vi)" If $d=n_1=1$, note that $V(f)$ has no singular point
at the origin.
\endroster \ms

{\rm (b)} Let $\phi=\phi(y,z)$ be irreducible in $\BC\{y,z\}$, and
let $V(\Phi)=\{(y,z):\Phi(y,z)=0\}$ be an analytic variety at
$(0,0)$ in $\BC^2$, each of which is written respectively in the
form,
$$\align
\text{\rm (7.1.2)} \qquad \qquad \phi&=(z^{\ell_1}+\bar{\ve}
y^{\delta_{1,1}})^{d'} +\sum_{p,q\ge 0}a_{p,q}y^{p}z^{q} \quad
\text{with} \quad \ell_1p+\delta_{1,1}q>\ell_1\delta_{1,1}d'
\qquad \qquad \qquad \\
\Phi&=y^{\zeta'}z^{\eta'}\phi \\
\endalign$$
satisfying the properties {\rm(i)}, {\rm(ii)},\dots, {\rm(vi)}:

\roster \item "(i)" $\gcd(\ell_1,\delta_{1,1})=1$ with $1\le
\ell_1<\delta_{11}$ and $d'$ is a positive integer.

\item "(ii)" $\bar{\ve}$ is a unit in $\BC\{y,z\}$.

\item "(iii)" The $a_{p,q}$ are nonzero complex numbers for some
nonnegative integers $p$ and $q$ such that
$\ell_1p+\delta_{11}q>\ell_1\delta_{1,1}d'$.

\item "(iv)" $\zeta'$ and $\eta'$ are nonnegative integers.

\item "(v)" If $\ell_1=1$, assume that $\eta'$ is a positive
integer.

\item "(vi)" If $d'=\ell_1=1$, note that $V(\phi)$ has no singular
point at the origin.

\endroster \ms

$\underline{\text{\bf {Conclusions}}}$ As a necessary condition for
\text{$V(F) \buildrel \text{{\rm divisor}} \over \sim V(\Phi)$}
under the standard resolutions as reduced varieties, we have the
following fact:
$$\align
  &\text{If \quad  \text{$V(F) \buildrel
\text{{\rm divisor}} \over \sim V(\Phi)$ under the standard
resolutions},} \tag 7.1.3 \\
& \text{then \quad  $\zeta=\zeta'$,  $\eta=\eta'$,  $n_1=\ell_1$,
$\beta_{1,1}=\delta_{1,1}$ \ {and} \ $d=d'$}.
\endalign$$

In particular, if $d=d'=1$, then the necessary condition for
\text{$V(F) \buildrel \text{{\rm divisor}} \over \sim V(\Phi)$}
under the standard resolutions, being given by {\rm(7.1.3)}, is also
sufficient.  $\square$
\endproclaim \ms

\definition{Remark 7.1.0} In preparation for the proof of
this theorem, first of all, we need the following lemma, that is,
Lemma $7.2$ with proof. After then, it will be found that there is
nothing to prove for this theorem, later. $\square$
\enddefinition \ms

\proclaim{Lemma 7.2} $\underline{\text{\bf {Assumptions}}}$ Let
$V(F)=\{(y,z): F(y,z)=0\}$ and $V(\Phi)=\{(y,z): \Phi(y,z)=0\}$ be
analytic varieties at $(0,0)\in \BC^2$, satisfying the same
assumption and notation as we have seen in the assumption of Theorem
$7.1$.

{\rm (a)(a1)} Let $V(G)=\{(y,z):G(y,z)=0\}$ be an analytic variety
at $(0,0)$ in $\BC^2$ defined by the form
$$\align
 G=& z^{\gamma}g \tag 7.2.1 \\
 g=& z^{n_1}+y^{\beta_{1,1}} \quad \text{with} \quad
\gcd(n_1,\beta_{1,1})=1,
\endalign$$
satisfying the following properties:

\roster

\item "(i)"  $1\le n_1 <\beta_{1,1}$.

\item "(ii)" If $n_1=1$, then $\gamma=1$.

\item "(iii)" If $n_1 \ge 2$, then $\gamma=0$.
\endroster \ms

{\rm (a2)} Let $\tau_m$ be the composition of a finite number $m$ of
successive blow-ups which is needed to get the standard resolution
of the singular point of $V(G)$ as we have used in {\rm(3.6.2)} of
the assumption of Theorem $3.6$. For each $t=1,2,\dots,m$, write
$\tau_t= \pi_1\circ\pi_2\circ\cdots \circ\pi_t:M^{(t)}\to \BC^2$
where \text{$ \pi_i:M^{(i)}\to M^{(i-1)}$ is a blow-up of
$M^{(i-1)}$ at a point of $M^{(i-1)}$} for $1\le i\le t$, and
$M^{(0)}=\BC^2$. For brevity of notation, let $V^{(t)}(G)$ be the
proper transform under $\tau_t$ for $1\le t\le m$. By the conclusion
of Theorem $3.6$, we can use the same $\tau_m$ for the composition
of the first finite number $m$ of successive blow-ups in process of
the standard resolution of the singular point $(0,0)$ of $V(F)$, as
a reduced variety.

Let $E^{(m)}=\tau^{-1}_m(0,0)$, and let $E^{(m)}=\cup E_i$, $1 \le
i \le m$, be the decomposition of $E^{(m)}$ into irreducible
components where each $E_i$ is called an exceptional curve of the
first kind under $\tau_m$. Then, $V^{(t)}(G)$ has one and only one
quasisingular point along $E_{t}$ for $1\le t\le m-1$. \ms

{\rm(a3)} Let $(F\circ\tau_m)_{divisor}$ be a divisor of
$F\circ\tau_m$ defined by
$$\align
(F\circ\tau_m)_{divisor}=V^{(m)}(F)+\sum^m_{i=1}e_iE_i, \tag 7.2.2
\endalign$$
where each $e_i$ is the multiplicity of $F\circ\tau_m$ along $E_i$
for $1\le i\le m$ and $V^{(m)}(F)$ is the proper transform of
$V(F)$ under $\tau_m$.

In more detail, let $\pi_{i}:M^{(i)}\to M^{(i-1)}$ is the {\rm
{i}}-th blow-up of $\tau_m$ at a quasisingular point of
$V^{(i-1)}(F)$ for $1\le i\le m-1$ in the sense of Definition
$2.6$ where $M^{(0)}=\BC^2$. Then, $V^{(i)}(F)$ has one and only
one quasisingular point along $E_{i}$ for $1\le i\le m$, if
exists.

Using {\rm(3.6.3)} of the conclusion of Theorem $3.6$, {\rm(a1)} and
{\rm(a2)}, then the equation in {\rm(7.2.2)} is well-defined,
satisfying the following properties:

{\rm(a3-1)} $e_i<e_{i+1}$ for $i=1,2,\dots,m$.

{\rm(a3-2)} $V(F)$ belongs to the \text{\rm type[1]} under
$\tau_m$ in the sense of Definition $2.8$ because $E_m$ is one and
only one exceptional curve of the first kind among
$E^{(m)}=\cup^{m}_{i=1} E_i$ which has three distinct intersection
points with other exceptional curves and the proper transform
$V^{(m)}(F)$ under $\tau_m$. \ms

{\rm (b)(b1)} Let $V(H)=\{(y,z):H(y,z)=0\}$ be an analytic variety
at $(0,0)$ in $\BC^2$ defined by the form
$$
\align
 H&=z^{\delta}h \tag 7.2.3  \\
 h&=z^{\ell_1}+y^{\delta_{1,1}} \quad \text{with} \quad
 \gcd(\ell_1,\delta_{1,1})=1,
\endalign
$$
satisfying the following properties:

\roster

\item "(i)"  $1\le \ell_1 <\delta_{1,1}$.

\item "(ii)" If $\ell_1=1$, then $\delta=1$.

\item "(iii)" If $\ell\ge 2$, then $\delta=0$.
\endroster \ms

{\rm (b2)} Let $\mu_\rho$ be the composition of a finite number
$\rho$ of successive blow-ups which is needed to get the standard
resolution of the singular point of $V(H)$ as we have used in the
assumption of Theorem $3.6$. For each $t=1,2,\dots,\rho$, write
$\mu_t= {\bar{\pi}_1}\circ{\bar{\pi}_2}\circ\cdots
\circ{\bar{\pi}_t}:{\bar{M}}^{(t)}\to \BC^2$ where \text{$
\bar\pi_i:\bar{M}^{(i)}\to \bar{M}^{(i-1)}$ is a blow-up of
$\bar{M}^{(i-1)}$ at some point of $\bar{M}^{(i-1)}$} for $1\le i\le
t$, and $\bar{M}^{(0)}=\BC^2$. For brevity of notation, let
$V^{(t)}(H)$ be the proper transform under $\mu_t$ for $1\le t\le
\rho$. By the conclusion of Theorem $3.6$, we can use the same
$\mu_\rho$ for the composition of the first finite number $\rho$ of
successive blow-ups in process of the standard resolution of the
singular point $(0,0)$ of $V(\Phi)$, as a reduced variety.

Let $\bar{E}^{(\rho)}={\mu_\rho}^{-1}(0,0)$, and let
$\bar{E}^{(\rho)}=\cup \bar{E}_i$, $1 \le i \le \rho$, be the
decomposition of $\bar{E}^{(\rho)}$ into irreducible components
where each $\bar{E}_i$ is called an exceptional curve of the first
kind under $\mu_\rho$. Then, $V^{(t)}(H)$ has one and only one
quasisingular point along $\bar{E}_{t}$ for $1\le t\le \rho-1$,
which is called the {\rm {t}}-th exceptional curve. \ms

{\rm(b3)} Let $(\Phi\circ\mu_\rho)_{divisor}$ be a divisor of
$\Phi\circ\mu_\rho$ defined by
$$\align
(\Phi\circ\mu_\rho)_{divisor}=V^{(\rho)}(\Phi)
+\sum^\rho_{i=1}\bar{e}_i\bar{E}_i, \tag 7.2.4
\endalign$$
where each $\bar{e}_i$ is the multiplicity of $\Phi\circ\mu_\rho$
along $\bar{E}_i$ for $1\le i\le \rho$ and $V^{\rho)}(\Phi)$ is
the proper transform of $V(\Phi)$ under $\mu_\rho$.

In more detail, let \text{$\bar\pi_i:\bar{M}^{(i)}\to
\bar{M}^{(i-1)}$} is the {\rm {i}}-th blow-up of $\mu_t$ at a
quasisingular point of $V^{(t-1)}(\Phi)$ for $1\le t\le \rho-1$ in
the sense of Definition $2.6$ where $M^{(0)}=\BC^2$. Then,
$V^{(t)}(\Phi)$ has one and only one quasisingular point along
$\bar{E}_{t}$ for $1\le t\le \rho$, if exists.

Using {\rm(3.6.3)} of the conclusion of Theorem $3.6$, {\rm(b1)} and
{\rm(b2)}, then the equation in {\rm(7.2.4)} is well-defined,
satisfying the following properties:

{\rm(b3-1)} $\bar{e}_i<\bar{e}_{i+1}$ for $1\le i\le \rho-1$.

{\rm(b3-2)} $V(\Phi)$ belongs to the \text{\rm type[1]} under
$\mu_\rho$ in the sense of Definition $2.8$ because $\bar{E}_\rho$
is one and only one exceptional curve of the first kind among
$\bar{E}^{(\rho)}=\cup^{\rho}_{i=1} \bar{E}_i$ which has three
distinct intersection points with other exceptional curves and the
proper transform $V^{(\rho)}(\Phi)$ under $\mu_\rho$. \ms

$\underline{\text{\bf {Conclusions}}}$ As a necessary condition or a
sufficient condition for \text{$V(F) \buildrel \text{{\rm divisor}}
\over \sim V(\Phi)$} under the standard resolutions as reduced
varieties, we have the following:

$\underline{\text{\rm Fact(1)}}$ If \text{$V(F) \buildrel \text{{\rm
divisor}} \over \sim V(\Phi)$} under the standard resolutions as
reduced varieties, then $\zeta=\zeta'$, $\eta=\eta'$ and
$n_1=\ell_1$, $\beta_{1,1}=\delta_{1,1}$ and $d=d'$.

$\underline{\text{\rm Fact(2)}}$ If \text{$V(F) \buildrel \text{{\rm
divisor}} \over \sim V(\Phi)$} under the standard resolutions as
reduced varieties, then $m=\rho$ with $\zeta=\zeta'$ and
$\eta=\eta'$, and $e_i=\bar{e}_i$ for $i=1,2,\dots,m=\rho$.

$\underline{\text{\rm Fact(3)}}$ Let $d\ge 1$ and $d'\ge 1$. Then
\text{$V(F) \buildrel \text{{\rm divisor}} \over \sim V(\Phi)$}
under the standard resolutions, as reduced varieties, if and only if
\text{$V(v(f\circ\tau_m)_{proper}) \buildrel \text{{\rm divisor}}
\over \sim V(v(\phi\circ\mu_\rho)_{proper})$} under the standard
resolutions of a quasisinguar point along $E_m=\{v=0\}$ or
$\bar{E}_{m}=\{{\bar{v}}=0\}$ with $m=\rho$, and either of {\rm
Fact(1)} and {\rm Fact(2)} is satisfied. Note that
$(F\circ\tau_m)_{proper}=(f\circ\tau_m)_{proper}$ and
$(\Phi\circ\mu_\rho)_{proper}=(\phi\circ\mu_\rho)_{proper}$.
$\square$
\endproclaim \ms

\definition{Remark 7.2.1} As far as the above assumptions of Lemma $7.2$
are concerned, whenever there exist the statements in both (a) and
(b) of the assumptions of Theorem $7.1$, it was already proved by
Theorem $3.6$ or Sublemma $5.4$ that the statements in (a1), (a2),
(a3), (b1), (b2) and (b3) of the assumptions of this lemma are true.
$\square$
\enddefinition \ms

{\bf \S7.2. The proofs of Lemma 7.2 and Theorem 7.1 }

In this section, we prove Lemma 7.2, and after then, it will be
found that there is nothing to prove for Theorem 7.1.

\demo{\bf{Proof of Lemma 7.2}} In preparation for the proof of
Fact(1), Fact(2) and Fact(3), first of all, we need to observe the
following:

Whenever there are two integers $n_1$ and $\beta_{1,1}$ with
$\gcd(n_1,\beta_{1,1})=1$ and $1\le n_1<\beta_{1,1}$, there is a
positive integer $s$ such that $sn_1<\beta_{1,1}\le (s+1)n_1$. Also,
whenever there are two integers $\ell_1$ and $\delta_{1,1}$ with
$\gcd(\ell_1,\delta_{1,1})=1$ and $1\le \ell_1<\delta_{1,1}$, there
is a positive integer $\bar{s}$ such that
$\bar{s}\ell_1<\delta_{1,1}\le (\bar{s}+1)\ell_1$.  Then, note that
$s<m$ with $\tau_m$ and $\bar{s}<\rho$ with $\mu_{\rho}$. \ms

$\underline{\text{\rm The proof of Fact(1)}}$ Since  \text{$V(F)
\buildrel \text{{\rm divisor}} \over \sim V(\Phi)$} under the
standard resolutions as reduced varieties, it is trivial to show
that the following are true:

{\rm(i)} $V(F)$ and $V(\Phi)$ have the same number of components,
that is, $1+\zeta+ \eta=1+\zeta'+ \eta'$.

{\rm(ii)} Then, $V^{(1)}(F)\cap E_1$ and $V^{(1)}(\Phi)\cap
\bar{E}_1$ have the same number of elements as set, that is,
${\zeta}={\zeta'}$ because $1\le n_1< \beta_{1,1}$ and $1\le \ell_1<
\delta_{1,1}$. So, $\eta=\eta'$ by (i).

{\rm(iii)}  Note that $e_1=\zeta+ \eta+n_1d$ and $\bar{e}_1=\zeta'+
\eta'+{\ell_1}d'$ are equal, and so $n_1d={\ell_1}d'$.

{\rm(iv)} Note that $e_{s+1}=\zeta+ (s+1)\eta+(s+1){\beta_{1,1}}d$
and $e_{\bar{s}+1}={\zeta}'+
(\bar{s}+1){\eta}'+(\bar{s}+1){\delta_{1,1}}d'$ with $s={\bar{s}}$
are equal, and so ${\beta_{1,1}}d={\delta_{1,1}}d'$.

{\rm(v)} Since $\gcd(n_1,\beta_{1,1})=1$ and
$\gcd(\ell_1,\delta_{1,1})=1$, then $n_1=\ell_1$,
${\beta_{1,1}}={\delta_{1,1}}$, and $d=d'$.

Thus, the proof of Fact(I) is done. \ms

$\underline{\text{\rm The proof of Fact(2)}}$ By (a3-2) and (b3-2)
in the assumption of this lemma, as reduced varieties, $E_m$ is one
and only one exceptional curve of the first kind among
$E^{(m)}=\cup^{m}_{i=1} E_i$ which has three distinct intersection
points with other exceptional curves and the proper transform
$V^{(m)}(F)$ under $\tau_m$, and also $\bar{E}_\rho$ is one and only
one exceptional curve of the first kind among
$\bar{E}^{(\rho)}=\cup^{\rho}_{i=1} \bar{E}_i$ which has three
distinct intersection points with other exceptional curves and the
proper transform $V^{(\rho)}(\Phi)$ under $\bar{\tau}_\rho$. Since
it is assumed that \text{$V(F) \buildrel \text{{\rm divisor}} \over
\sim V(\Phi)$} under the standard resolutions as reduced varieties,
then it is clear by {\rm Fact(1)} that
$e_m=n_1\zeta+\beta_{1,1}\eta+n_1\beta_{1,1}d$ and
$\bar{e}_\rho=\ell_1\zeta'+\delta_{1,1}\eta'+\ell_1\delta_{1,1}d'$
are equal, and that $G(y,z)$ of (7.2.1) and $H(y,z)$ of (7.2.3) are
equal, too. So, $m=\rho$ and $e_i=\bar{e}_i$ for
$i=1,2,\dots,m=\rho$. Thus, Fact(2) can be easily proved.

Also, it is trivial to prove by Theorem $3.6$ that Fact(3) is true.
$\square$
\enddemo \ms

\demo{\bf{Proof of Theorem 7.1}} If \text{$V(F) \buildrel \text{{\rm
divisor}} \over \sim V(\Phi)$} under the standard resolutions, then
it is clear by Definition $2.6$ that $\zeta=\zeta'$ and $\eta=\eta'$
and $V(f)$ and $V(\phi)$ have the same multiplicity
$n_1d={\ell_1}d'$ at $(y,z)=(0,0)$. For the proof of theorem, it
remains to show that $n_1=\ell_1$, $\beta_{1,1}=\delta_{1,1}$ and
$d=d'$, which can be easily proved by the induction method on the
multiplicity of $V(f)$ at $(y,z)=(0,0)$, using the same kind of
methods and notations as we have used in the proof of Theorem $3.6$
with Lemma $4.2$ and Lemma $4.3$. Thus, the proof can be finished.
$\square$
\enddemo \ms

{\bf \S7.3 The necessary and sufficient condition for any two
quasi-Puiseux convergent power series of recursive types in C\{y,z\}
to have the same divisor under the standard resolutions and The
solution of (i) of Problem[1-A]}

\proclaim{Theorem 7.3} Let $r$ and $\rho$ be arbitrary positive
integers.

$\underline{\text{\bf {Assumptions}}}$ By the same way as in {\rm
Definition 5.0.0}, let $g_r\in \BC\{y,z\}$ be the quasi-Puiseux
convergent power series of the recursive r-type defined by {\rm
{Sequences[I]}, and $\phi_{\rho}\in \BC\{y,z\}$ be the quasi-Puiseux
convergent power series of the recursive {$\rho$}-type defined by
{\rm Sequences[II]} such that {\rm {Sequences[I]} and {\rm
{Sequences[I]} are defined respectively, as follows: Note that {\rm
{Sequences[I]} and {\rm Sequences[II]} are the same up to the change
of notations. \ms

$\text{\bf Sequences[I]}$ Let $\{X_k:k=1,2,\dots,r\}$ with
$X_k\subset N_0$, $\{g_k:k=1,2,\dots,r\}$ with $g_k\in \BC\{y,z\}$
and $\{\text{$\Delta_k:N^k_0\to N_0$ is an integer-valued function
for $k=1,2,\dots,r$}\}$ be three sequences satisfying the following
$\underline{\text{\rm five conditions}}:$

{\noindent} Five conditions are denoted by \text{\bf The 1st
${\text{\bf{Cond}}}^{\text{{\bf(0)}}}$}, \dots, \text{\bf The 5-th
${\text{\bf{Cond}}}^{\text{{\bf(0)}}}$ of Sequences[I]}. \ms

{\bf[I]-(1)} \text{\bf The 1st
${\text{\bf{Cond}}}^{\text{{\bf(0)}}}$ of Sequences[I]:}

{\rm (1a)}  $X_1=\{n_1,\beta_{1,1}\}$ with $n_1\ge 2$ and
$\beta_{1,1}\ge 1$.

 {\rm (1b)} $X_j=\{n_j,\beta_{j,1},\beta_{j,2},\dots,\beta_{j,j}\}$
 with $n_j\ge 2$ \quad where $j=2,\dots,r$.

 If $j\ge 2$, then assume that at least one of
$\beta_{j,1},\beta_{j,2},\dots,\beta_{j,j}$ is nonzero. \ms

{\bf [I]-(2)} \text{\bf The 2nd
${\text{\bf{Cond}}}^{\text{{\bf(0)}}}$ of Sequences[I]:}

{\rm (2a)} $g_1=z^{n_1}+{\ve_1}y^{\beta_{1,1}}$.

{\rm (2b)} $g_j=g_{j-1}^{n_j}+{\ve_j}y^{\beta_{j.1}}z^{\beta_{j.2}}
g_1^{\beta_{j.3}}\cdots g_{j-2}^{\beta_{j,j}}$ \quad where
$j=2,\dots,r$.

Note that each $\ve_i=\ve_i(y,z)$ is a unit in $\BC\{y,z\}$ for
$1\le i\le r$.\ms

{\bf [I]-(3)} \text{\bf The 3rd
${\text{\bf{Cond}}}^{\text{{\bf(0)}}}$ of Sequences[I]:}

{\rm(3a)} $\Delta_1(t)=t$ for each $t\in N_0$.

{\rm(3b)}
$\Delta_j(t_j)^j_{k=1}=t_j\Delta_{j-1}(\beta_{j-1,k})^{j-1}_{k=1}
+n_{j-1}\Delta_{j-1}(t_k)^{j-1}_{k=1}$ for each $(t_k)^j_{k=1}\in
N^j_0$

\indent where $j=2,\dots,r$. \ms

{\bf [I]-(4)} \text{\bf The 4-th
${\text{\bf{Cond}}}^{\text{{\bf(0)}}}$ of Sequences[I]:}

{\rm(4)} $\Delta_j(\beta_{j,k})^j_{k=1}>n_jn_{j-1}
\Delta_{j-1}(\beta_{j-1,k})^{j-1}_{k=1}$ for $2\le j\le r$. \bs

{\bf [I]-(5)} \text{\bf The 5-th
${\text{\bf{Cond}}}^{\text{{\bf(0)}}}$ of Sequences[I]:}

{\rm(5)} $\gcd(n_j,\Delta_j(\beta_{j,k})^j_{k=1})=1$ for $1\le j\le
r$.  \bs

$\text{\bf Sequences[II]}$ Let $\{W_k:k=1,2,\dots,\rho\}$ with
$W_k\subset N_0$, $\{\phi_k:k=1,2,\dots,\rho\}$ with $\phi_k\in
\BC\{y,z\}$ and $\{\text{$\omega_k:N^k_0\to N_0$ is an
integer-valued function for $k=1,2,\dots,\rho$}\}$ be three
sequences satisfying the following $\underline{\text{\rm five
conditions}}:$

{\noindent} Five conditions are denoted by \text{\bf The 1st
${\text{\bf{Cond}}}^{\text{{\bf(0)}}}$}, \dots, \text{\bf The 5-th
${\text{\bf{Cond}}}^{\text{{\bf(0)}}}$ of Sequences[II]}. \ms

{\bf [II]-(1)} \text{\bf The 1st
${\text{\bf{Cond}}}^{\text{{\bf(0)}}}$ of Sequences[II]:}

{\rm (1a)} $W_1=\{\ell_1,\delta_{1,1}\}$ with $\ell_1\ge 2$ and
$\delta_{1,1}\ge 1$.

{\rm (1b)}
$W_{j}=\{\ell_{j},\delta_{j,1},\delta_{j,2},\dots,\delta_{j,j}\}$
with $\ell_{j}\ge 2$ \quad where $j=2,\dots,\rho$.

If $j\ge 2$, then assume that at least one of
$\delta_{j,1},\delta_{j,2},\dots,\delta_{j,j}$ is nonzero. \ms

{\bf [II]-(2)} \text{\bf The 2nd
${\text{\bf{Cond}}}^{\text{{\bf(0)}}}$ of Sequences[II]:}

{\rm (2a)} $\phi_1=z^{\ell_1}+{\bar{\ve}_1}y^{\delta_{1,1}}$.

{\rm $(2b)$}
$\phi_{j}=\phi_{j-1}^{\ell_{j}}+{\bar{\ve}_j}y^{\delta_{j,1}}z^{\delta_{j,2}}
      \phi_1^{\delta_{j,3}}\cdots \phi_{j-2}^{\delta_{j,j}}$ \quad
where $j=2,\dots,{\rho}$.

Note that each $\bar{\ve}_i=\bar{\ve}_i(y,z)$ is a unit in
$\BC\{y,z\}$ for $1\le i\le {\rho}$.\ms

{\bf [II]-(3)} \text{\bf The 3rd
${\text{\bf{Cond}}}^{\text{{\bf(0)}}}$ of Sequences[II]:}

{\rm (3a)} $\omega_1(t)=t$ for each $t\in N_0$.

{\rm (3b)}
$\omega_j(t_k)^{j}_{k=1}=t_j\omega_{j-1}(\delta_{j-1,k})^{j-1}_{k=1}
+\ell_{j-1}\omega_{j-1}(t_k)^{j-1}_{k=1}$ for each
$(t_k)^{j}_{k=1}\in N^{j}_0$

\indent where $j=2,\dots,{\rho}$. \ms

{\bf [II]-(4)} \text{\bf The 4-th
${\text{\bf{Cond}}}^{\text{{\bf(0)}}}$ of Sequences[II]:}

{\rm(4)} $\omega_j(\delta_{j,k})^j_{k=1}>\ell_j \ell_{j-1}
\omega_{j-1}(\delta_{j-1,k})^{j-1}_{k=1}$  for $2\le j\le \rho$. \ms

{\bf [II]-(5)} \text{\bf The 5-th
${\text{\bf{Cond}}}^{\text{{\bf(0)}}}$ of Sequences[II]:}

{\rm(5)} $\gcd(\ell_{j},\omega_{j}(\delta_{j,k})^{\rho}_{k=1})=1$
for $1\le j\le \rho$.   \bs

$\underline{\text{\bf {Conclusions}}}$

Let $V(y^{\zeta} g_r)$ and $V(y^{\eta}\phi_\rho)$ be analytic
varieties at $(y,z)=(0,0)$ with isolated singularity at the origin
defined as follows:
$$\align
 V(y^{\zeta}g_r) &=\{y,z):y^{\zeta}g_r(y,z)=0\},  \tag 7.3.0 \\
 V(y^{\eta}\phi_\rho) &=\{(y,z):y^{\eta}\phi_\rho(y,z)=0\}.
\endalign$$

Note by assumption that $g_{r}$ and $\phi_{\rho}$ are irreducible in
$\BC\{y,z\}$.

For convenience of notation, we need to assume in addition that the
following hold:

\roster \item"(i)"  if $\beta_{1,1}=1$ from $g_r$, then $\zeta$ is a
positive integer.

\item"(ii)"  if $\delta_{1,1}=1$ from $\phi_\rho$, then $\eta$ is a
positive integer.

\item"(iii)"  if $\zeta=0$,  then $\beta_{1,1}>1$ by {\rm (i)}, and
so we may assume without loss of generality that $2\le
n_1<\beta_{1,1}$, because if $n_1>\beta_{1,1}\ge 2$ then replace $z$
by $y$, and $y$ by $z$.

\item"(iv)" if $\eta=0$, then $\delta_{1,1}>1$ by {\rm (ii)}, and
so we may assume without loss of generality that $2\le \ell_1
<\delta_{1,1}$, because if $\ell_1
>\delta_{1,1}\ge 2$ then replace $z$ by $y$, and $y$ by $z$.
\endroster \ms

{\bf[A]} Then, we have the following:
 $$\align
(7.3.1)\qquad \quad  & \text{$V(y^{\zeta}g_r) \buildrel \text{{\rm
divisor}} \over \sim V(y^{\eta}\phi_\rho)$ under the standard
resolutions as reduced varieties} \qquad \qquad\\
 \iff \quad & \\
(7.3.2)\qquad \quad  & \text{$\zeta=\eta$, \quad $r=\rho$, \quad and
\quad
 $n_j=\ell_j$ \quad for $1\le j\le r=\rho$, \quad and} \\
 & \text{$\Delta_j(\beta_{j,k})^j_{k=1}=\omega_j(\delta_{j,k})^j_{k=1}$
 \quad for each $j=1,2,\dots,r$}.
\endalign$$ \ms

{\bf[B]} In particular, assume that the above $g_r\in \BC\{y,z\}$ in
{\rm {Sequences[I]} satisfies either the property of the Puiseux
convergent power series or $2\le n_1<\beta_{1,1}$. Also, assume that
the above $\phi_\rho\in \BC\{y,z\}$ in {\rm {Sequences[II]}
satisfies either the property of the Puiseux convergent power series
or $2\le \ell_1<\delta_{1,1}$.

Then, we have
$$\align
& \text{$g_r \buildrel \text{{\rm divisor}} \over \sim \phi_\rho$
\quad under the standard
resolutions} \tag 7.3.3\\
 \iff \quad & \\
& \text{$n_j=\ell_j$ \qquad \qquad for each $j=1,2,\dots,r=\rho$,
\quad  and}  \tag 7.3.4 \\
 & \text{$\Delta_j(\beta_{j,k})^j_{k=1}=\omega_j(\delta_{j,k})^j_{k=1}$
\qquad  for each $j=1,2,\dots,r$} \\
 \iff \quad & \\
& \text{$n_j=\ell_j$ \qquad \qquad  for each $j=1,2,\dots,r=\rho$,
\quad
  and}  \tag 7.3.5 \\
& \Delta_j(t_k)^j_{k=1}=\omega_j(t_k)^j_{k=1} \qquad \text{for
each $(t_k)^j_{k=1}\in N^j_0$, \quad  and} \\
 & \text{$\Delta_j(\beta_{j,k})^j_{k=1}=\Delta_j(\delta_{j,k})^j_{k=1}$
\qquad for each $j=1,2,\dots,r$}. \quad \quad \text{$\square$}
 \endalign$$
\endproclaim

\definition{Remark 7.3.1} By Definition $5.0.0$, $g_r\in \BC\{y,z\}$
of Sequences[I] is called a quasi-Puiseux convergent power series of
the recursive r-type. In addition, if $2\le n_1<\beta_{1,1}$ in
\text{\rm The 1-th ${\text{\rm{Cond}}}^{\text{{\rm(0)}}}$,} then
$g_r\in \BC\{y,z\}$ of Sequences[I] is a Puiseux convergent power
series of the recursive r-type. Also, by Definition $5.0.0$,
$\phi_{\rho}\in \BC\{y,z\}$ of Sequences[II] is called a
quasi-Puiseux series of the recursive $\rho$-type. In addition, if
$2\le \ell_1<\delta_{1,1}$ in \text{\rm The 1-th
${\text{\rm{Cond}}}^{\text{{\rm(0)}}}$,} then $\phi_{\rho}\in
\BC\{y,z\}$ of $\text{\rm Sequences[II]}$ is called a Puiseux
convergent power series of the recursive $\rho$-type. $\square$
\enddefinition \ms

{\bf \S7.4. In preparation for the proof of Theorem 7.3 } \ms

For the proof of Theorem $7.3$, first we will prepare two sublemmas,
Sublemma 7.4 and Sublemma 7.5 in $\S7.4$. After then, in $\S7.5$ we
will finish the proof of Theorem $7.3$, using these two sublemmas.
First, to construct the statement for Sublemma 7.4 with its proof,
it suffices to apply Sublemma 5.1, Sublemma 5.2, Sublemma 5.3 and
Sublemma 5.4 of Theorem 5.0 to the assumptions of Theorem $7.3$.
Next, to construct the statement for Sublemma 7.5 with its proof, it
suffices to apply Sublemma 5.5 of Theorem 5.0 to Sublemma 7.4. \ms

\proclaim{Sublemma 7.4} $\underline{\text{\bf {Assumptions}}}$
Suppose that the same properties and notations as in {\rm
Sequences[I]} and  {\rm Sequences[II]} of  the assumptions of
Theorem $7.3$ hold. Let $r$ be an arbitrary integer with $r\ge 1$ in
{\rm Sequences[I]}, and $\rho$ be an arbitrary integer with $\rho\ge
1$ in {\rm Sequences[II]}. \ms

$\underline{\text{\bf {Conclusions}}}$ Then, we have the following
three facts, {\rm Fact(1)}, {\rm Fact(2)} and {\rm Fact(3)}. \ms

{\bf Fact(1):} By {\rm(5.2.1)} in Sublemma $5.2$ of Theorem $5.0$,
we may assume without any need of proof that two properties, {\rm
Property(1) and Property(2)} of {\rm Fact(1)} are true: \ms

$\underline{\text{\rm {Property(1) of Fact(1)}}}$ For any $r\ge 1$,
$g_r$ of Theorem $7.3$ can be written in the form
$$\align
 & g_r=g_r(y,z)=(z^{n_1}+\ve_1 y^{\beta_{1,1}})^{n_2n_3\cdots
n_r}+\sum_{\alpha,\beta\ge 0}
c^{(r)}_{\alpha,\beta}y^{\alpha}z^{\beta} \tag 7.4.1 \\
 & \text{with \quad
$n_1\alpha+\beta_{1,1}\beta>n_1\beta_{1,1}n_2n_3\cdots n_r$,}
\endalign$$
where $\ve_1=\ve_1(y,z)$ is a unit in $\BC\{y,z\}$, and the
$c^{(r)}_{\alpha,\beta}$ are nonzero complex numbers for some
nonnegative integers $\alpha$ and $\beta$, if exist. \ms

$\underline{\text{\rm {Property(2) of Fact(1)}}}$ For any $\rho\ge
1$, $\phi_\rho$ of Theorem $7.3$ can be written in the form
$$\align
& \phi_\rho=\phi_\rho(y,z)=(z^{\ell_1}+\bar{\ve}_1
y^{\delta_{1,1}})^{\ell_2\ell_3\cdots
\ell_\rho}+\sum_{\gamma,\delta\ge 0}
a^{(\rho)}_{\gamma,\delta}y^{\gamma}z^{\delta} \tag 7.4.2 \\
&\text{with \quad
$\ell_1\gamma+\delta_{1,1}\delta>\ell_1\delta_{1,1}\ell_2\ell_3\cdots
\ell_\rho$}, \\
\endalign$$
where $\bar{\ve}_1=\bar{\ve}_1(y,z)$ is a unit in $\BC\{y,z\}$, and
the $a^{(\rho)}_{\gamma,\delta}$ are nonzero complex numbers for
some nonnegative integers $\gamma$ and $\delta$, if exist.  \ms

{\bf Fact(2):} Let $\tau_{m}=\pi_1\circ\pi_2\circ\cdots
\circ\pi_{m}:M^{(m)}\to \BC^2$ and
$\mu_{\lambda}=\bar{\pi}_1\circ\bar{\pi}_2\circ\cdots
\circ\bar{\pi}_{\lambda}:\bar{M}^{(\lambda)}\to \BC^2$ be the
composition of two finite numbers $m$ and $\lambda$ of successive
blow-ups which are needed only to get the standard resolutions of
the singular point $(y,z)=(0,0)$ of $V(y^{\zeta}g_1)$ and
$V(y^{\eta}\phi_1)$, respectively. For brevity of notation, write
$G_j=y^{\zeta}g_j$ for $1\le j\le r$, and let
$\Phi_i=y^{\eta}\phi_i$ for $1\le i\le \rho$, as we have seen in
{\rm (7.3.0)}.

Then, $y^{\zeta}g_j\in type [1]$ under the standard resolution
$\tau_{m}$ and $y^{\eta}\phi_i\in type [1]$ under the standard
resolution $\mu_{\lambda}$ by Theorem $5.0$. \ms

{\bf Fact(3):} As in Lemma $7.2$, let $(G_j\circ\tau_m)_{divisor}$
be a divisor of $G_j\circ\tau_m$ defined by
$$(G_j\circ\tau_m)_{divisor}=V^{(m)}(G_j)+\sum^m_{i=1}e_iE_i, \tag 7.4.3$$ where
each $e_i$ is the multiplicity of $G_j\circ\tau_m$ along $E_i$ for
$1\le i\le m$ and $V^{(m)}(G_j)$ is the proper transform of $V(G_j)$
under $\tau_m$.

Also, as in Lemma $7.2$, let $(\Phi_j\circ\mu_{\lambda})_{divisor}$
be a divisor of $\Phi_j\circ\mu_{\lambda}$ defined by
$$(\Phi_j\circ\mu_{\lambda})_{divisor}=V^{(\lambda)}(\Phi_j)
+\sum^{\lambda}_{i=1}{\bar{e_i}}{\bar{E_i}}, \tag 7.4.4 $$ where
each $\bar{e_i}$ is the multiplicity of $\Phi_j\circ\mu_{\lambda}$
along $\bar{E_i}$ for $1\le i\le {\lambda}$ and
$V^{({\lambda})}(\Phi_j)$ is the proper transform of $V(\Phi_j)$
under $\mu_{\lambda}$.

If \text{$V(y^{\zeta}g_r) \buildrel \text{{\rm divisor}} \over \sim
V(y^{\eta}\phi_\rho)$ under the standard resolutions as reduced
varieties}, we have the following by Lemma $7.2$:
$$\text{$\tau_m=\mu_{\lambda}$ and
$e_i=\bar{e_i}$ for $i=1,2,\dots,m=\lambda$}. \quad \text{$\square$}
\tag 7.4.5 $$
\endproclaim

\demo{\bf Proof of Sublemma 7.4} Applying Sublemma 5.1, Sublemma
5.2, Sublemma 5.3 and Sublemma 5.4 of Theorem 5.0, and Lemma $7.2$
to the assumptions of Theorem $7.3$, then there is nothing to prove
for this sublemma. $\square$
\enddemo \ms

\proclaim{Sublemma 7.5} $\underline{\text{\bf Assumptions}}$ \quad
For brevity of notation, let {\rm Sequences[I]} of the {\rm(r+1)-th}
type and {\rm Sequences[II]} of the {\rm($\rho+1$)-th} type be given
sequences, each of which satisfies the same kind of properties and
notations, as we have seen in {\rm Sequences[I]} of the r-th type
and {\rm Sequences[II]} of the $\rho$-th type under the assumptions
of Theorem $7.3$, respectively. Note by Theorem $7.3$ that $r$ is an
arbitrary integer with $r\ge 1$ in {\rm Sequences[I]}, and $\rho$ is
an arbitrary integer with $\rho\ge 1$ in {\rm Sequences[II]}. \ms

$\underline{\text{\bf Conclusions}}$ As in Sublemma $7.4$, let
$\tau_{m}:M^{(m)}\to \BC^2$ and
$\mu_{\lambda}:\bar{M}^{(\lambda)}\to \BC^2$ be the composition of
two finite numbers $m$ and $\lambda$ of successive blow-ups which
are needed only to get the standard resolutions of the singular
point $(y,z)=(0,0)$ of $V(y^{\zeta}g_1)$ and $V(y^{\eta}\phi_1)$,
respectively. Then, we have two properties, {\rm Property(1)} and
{\rm Property(2)}: \ms

\noindent$\underline{\text{\bf Property(1) of Conclusions}}$ As we
have seen in Sublemma $5.5$, using
$\{h_k=(g_{k+1}\circ\tau_m)_{proper}: k=1,2,\dots,r\}$ \text{with
$h_k\in \BC\{v,u+1\}$}, it has been already proved that we can
construct another new sequences generated by
\text{\rm{Sequences[I]}} of the {\rm(r+1)-th} type, denoted by {\rm
{Sequences[I]}$^{(1)}$} of the {\rm r-th} type, which is rewritten
as follows:
$$\align
&  \text{\bf Sequences[I]}^{\text{{\bf(1)}}} \text{\bf of the r-th
type}\text{\bf :} \quad \text{Let $\{Y_k: k=1,2,\dots,r\}$
with $Y_k\subset N_0$}, \\
& \text{$\{h_k: k=1,2,\dots,r\}$ with
$h_k=(g_{k+1}\circ\tau_m)_{proper}\in \BC\{v,u+1\}$, and} \\
&\{\text{$\Xi_k:N^k_0\to N_0$ is an integer-valued
function for $k=1,2,\dots,r$}\} \\
& \text{be three sequences, satisfying the following five conditions
for each k.} \qquad \qquad \\
&\text{Such five conditions are denoted by \text{\bf The 1st
${\text{\bf{Cond}}}^{\text{{\bf(1)}}}$}{\bf, \dots,} \text{\bf The
5-th ${\text{\bf{Cond}}}^{\text{{\bf(1)}}}$}}.
\endalign$$

\indent{\rm {\bf [I]-(1)}} \quad \text{\bf The 1st
${\text{\bf{Cond}}}^{\text{{\bf(1)}}}$} $\text{\bf of
Sequences[I]}^{\text{{\bf(1)}}}$ {\bf :}
$$\align
\quad (1a) \qquad  Y_1 &=\{s_1,\gamma_{1,1}\} \quad
\text{where} \\
s_1 &=n_2\ge 2 \quad \text{and} \quad
\gamma_{1,1}=\Delta^{\sharp}_2(\beta_{2,1},\beta_{2,2})-n_1\beta_{1,1}n_2>0. \\
\quad (1b) \qquad  Y_j
&=\{s_j,\gamma_{j,1},\gamma_{j,2},\dots,\gamma_{j,j}\}
\quad \text{\rm where for j=2,3,\dots,r,} \\
 s_j &=n_{j+1}\ge 2, \quad \text{and} \quad
\gamma_{j,1}=\Delta^{\sharp}_{j+1}(\beta_{j+1,k})^{j+1}_{k=1}
-n_1\beta_{1,1}n_2n_3\cdots n_{j+1}>0, \qquad \qquad\\
\gamma_{j,2} &=\beta_{j+1,3}, \gamma_{j,3}=\beta_{j+1,4},\dots,
\gamma_{j,j}=\beta_{j+1,j+1}.
\endalign$$

{\rm {\bf [I]-(2)}} \quad \text{\bf The 2nd
${\text{\bf{Cond}}}^{\text{{\bf(1)}}}$} $\text{\bf
 of Sequences[I]}^{\text{{\bf(1)}}}$ {\bf :} Let $2\le j\le r$. \ms

{\rm(7.5.1)} \quad {\rm(2a)}\quad $h_1
=(u+1)^{s_1}+\ve_{1,2}v^{\gamma_{1,1}}$.

\quad \quad \qquad {\rm(2b)}\quad $h_j
=h_{j-1}^{s_j}+\ve_{1,j+1}v^{\gamma_{j,1}}(u+1)^{\gamma_{j,2}}
  h_1^{\gamma_{j,3}}h_2^{\gamma_{j,4}}\cdots
h_{j-2}^{\gamma_{j,j}}$.

Note that $\ve_{1,i}=\ve_{1,i}(v,u+1)$ is a unit in $\BC\{v,u+1\}$
for $2\le i\le r+1$. \bs

{\rm {\bf [I]-(3)}} \quad \text{\bf The 3rd
${\text{\bf{Cond}}}^{\text{{\bf(1)}}}$} $\text{\bf
 of Sequences[I]}^{\text{{\bf(1)}}}$ {\bf :} Let $2\le j\le r$. \ms

{\rm(3a)}\quad $\Xi_1(t)=t$  for each $t\in N_0$.

{\rm(3b)}\quad
$\Xi_j(t_k)^j_{k=1}=t_j\Xi_{j-1}(\gamma_{j-1,k})^{j-1}_{k=1}
  +s_{j-1}\Xi_{j-1}(t_k)^{j-1}_{k=1}$ for $(t_k)^j_{k=1}\in N^j_0$.
\ms

\noindent{\rm {\bf [I]-(4$\alpha$)}} \quad \text{\bf The
(4$\alpha$)-th ${\text{\bf{Cond}}}^{\text{{\bf(1)}}}$} $\text{\bf
 of Sequences[I]}^{\text{{\bf(1)}}}$ {\bf :} \ms

\noindent{\rm(7.5.2)} \qquad  {\rm(4a)}\qquad
$\Xi_1(\gamma_{1,1})=\gamma_{1,1}
=\Delta_2(\beta_{2,1},\beta_{2,2})-n_1\beta_{1,1}n_2>0$.

\quad \quad \qquad {\rm(4b)}\qquad $\Xi_q(\gamma_{q,k})^q_{k=1}
-s_qs_{q-1}\Xi_{q-1}(\gamma_{q-1,k})^{q-1}_{k=1}$

\qquad \qquad \qquad \quad $=\Delta_{q+1}(\beta_{q+1,k})^{q+1}_{k=1}
-n_{q+1}n_q\Delta_q(\beta_{q,k})^q_{k=1}>0$ for  $q=2,3,\dots,r$.
\ms

\noindent{\rm {\bf [I]-(5$\alpha$)}} \quad \text{\bf The
(5$\alpha$)-th ${\text{\bf{Cond}}}^{\text{{\bf(1)}}}$} $\text{\bf
 of Sequences[I]}^{\text{{\bf(1)}}}$ {\bf :}
$$\align
(7.5.3) \qquad &\gcd(s_j,\Xi_j(\gamma_{j,k})^j_{k=1})
=\gcd(n_{j+1},\Delta_{j+1}(\beta_{j+1,k})^{j+1}_{k=1})=1 \quad
\text{for $j=1,2,\dots,r.$} \qquad
\endalign$$ \ms

\noindent$\underline{\text{\bf Property(2) of Conclusions}}$ As we
have seen in Sublemma $5.5$, using
$\{\psi_k=(\phi_{k+1}\circ\mu_\lambda)_{proper}: k=1,2,\dots,\rho\}$
\quad \text{with $\psi_k\in \BC\{\bar{v},\bar{u}+1\}$}, it was
already proved that we can construct another new sequences generated
by \text{\rm{Sequences[II]}} of the {\rm($\rho+1$)-th} type, denoted
by {\rm {Sequences[II]}$^{(1)}$} of the {\rm $\rho$-th} type, which
is rewritten as follows:
$$\align
&  \text{\bf Sequences[II]}^{\text{{\bf(1)}}} \text{\bf of the
{$\rho$}-th type}\text{\bf :}
 \quad  \text{Let $\{L_k: k=1,2,\dots,\rho\}$
with $L_k\subset N_0$}, \\
& \text{$\{\psi_k: k=1,2,\dots,\rho\}$ with
$\psi_k=(\phi_{k+1}\circ\mu_{\lambda})_{proper}\in
\BC\{\bar{v},\bar{u}+1\}$, and} \qquad \qquad \\
&\{\text{$\theta_k:N^k_0\to N_0$ is an integer-valued
function for $k=1,2,\dots,\rho$}\} \\
& \text{be three sequences, satisfying the following five
conditions for each k}: \\
&\text{Such conditions are denoted by \text{\bf The 1st
${\text{\bf{Cond}}}^{\text{{\bf(1)}}}$}{\bf, \dots,} \text{\bf The
5-th ${\text{\bf{Cond}}}^{\text{{\bf(1)}}}$}.}
\endalign$$

{\rm {\bf [II]-(1)}} \quad \text{\bf The 1st
${\text{\bf{Cond}}}^{\text{{\bf(1)}}}$} $\text{\bf of
Sequences[II]}^{\text{{\bf(1)}}}$ {\bf :}
$$\align
\quad (1a) \qquad  L_1 &=\{p_1,\nu_{1,1}\} \quad
\text{where} \\
p_1 &=\ell_2\ge 2 \quad \text{and} \quad
\nu_{1,1}=\omega^{\sharp}_2(\delta_{2,1},\delta_{2,2})-\ell_1\delta_{1,1}\ell_2>0. \\
\quad (1b) \qquad  L_j &=\{p_j,\nu_{j,1},\nu_{j,2},\dots,\nu_{j,j}\}
\quad \text{\rm where for j=2,3,\dots,$\rho$,} \\
 p_j &=\ell_{j+1}\ge 2, \quad \text{and} \quad
\nu_{j,1}=\omega^{\sharp}_{j+1}(\delta_{j+1,k})^{j+1}_{k=1}
-\ell_1\delta_{1,1}\ell_2\ell_3\cdots \ell_{j+1}>0, \qquad \qquad\\
\nu_{j,2} &=\delta_{j+1,3}, \nu_{j,3}=\delta_{j+1,4},\dots,
\nu_{j,j}=\delta_{j+1,j+1}.
\endalign$$

{\rm {\bf [II]-(2)}} \quad \text{\bf The 2nd
${\text{\bf{Cond}}}^{\text{{\bf(1)}}}$} $\text{\bf
 of Sequences[II]}^{\text{{\bf(1)}}}$ {\bf :} Let $2\le j\le \rho$.
\ms

{\rm(7.5.4)} \quad {\rm(2a)}\quad $\psi_1
=(\bar{u}+1)^{p_1}+\bar{\ve}_{1,2} \bar{v}^{\nu_{1,1}}$.

\quad \quad \qquad {\rm(2b)}\quad $\psi_j=\psi^{p_j}_{j-1}
+\bar{\ve}_{1,j+1}
v^{\nu_{j,1}}(\bar{u}+1)^{\nu_{j,2}}\psi^{\nu_{j,3}}_1\psi^{\nu_{j,4}}_2
 \cdots \psi^{\nu_{j,j}}_{j-2}$.

Note that $\bar{\ve}_{1,i}=\bar{\ve}_{1,i}(\bar{v},\bar{u}+1)$ is a
unit in $\BC\{\bar{v},\bar{u}+1\}$ for $2\le i\le \rho+1$. \bs

{\rm {\bf [II]-(3)}} \quad \text{\bf The 3rd
${\text{\bf{Cond}}}^{\text{{\bf(1)}}}$} $\text{\bf
 of Sequences[II]}^{\text{{\bf(1)}}}$ {\bf :} Let $2\le j\le \rho$.
\ms

{\rm(3a)} \quad $\theta_1(t) =t$ for each $t\in N_0$.

{\rm(3b)} \quad $\theta_j(t_k)^j_{k=1}
 =t_j\theta_{j-1}(\nu_{j-1,k})^{j-1}_{k=1}
 +p_{j-1}\theta_{j-1}(t_k)^{j-1}_{k=1}$ for  $(t_k)^j_{k=1}\in
N^j_0$. \ms

{\rm {\bf [II]-(4$\alpha$)}} \quad \text{\bf The (4$\alpha$)-th
${\text{\bf{Cond}}}^{\text{{\bf(1)}}}$} $\text{\bf
 of Sequences[II]}^{\text{{\bf(1)}}}$ {\bf :} \ms

{\rm(7.5.5)} \quad {\rm(4a)}\quad $\theta_1(\nu_{1,1})=\nu_{1,1}
=\omega_2(\delta_{2,1},\delta_{2,2})-\ell_1\delta_{1,1}\ell_2>0$.

\quad \quad \qquad {\rm(4b)}\quad $\theta_q(\nu_{q,k})^q_{k=1}-p_q
p_{q-1}\theta_{q-1}(\nu_{q-1,k})^{q-1}_{k=1}$

\qquad \qquad \qquad $=\omega_{q+1}(\delta_{q+1,k})^{q+1}_{k=1}
-\ell_{q+1}\ell_q\omega_q(\delta_{q,k})^q_{k=1}>0$ for
$q=2,3,\dots,\rho$. \ms

{\rm {\bf [II]-(5$\alpha$)}} \quad \text{\bf The (5$\alpha$)-th
${\text{\bf{Cond}}}^{\text{{\bf(1)}}}$} $\text{\bf
 of Sequences[II]}^{\text{{\bf(1)}}}$ {\bf :}
$$\align
\quad (7.5.6) \quad  &\gcd(p_q,\theta_q(\nu_{q,k})^q_{k=1})=
\gcd({\ell}_{q+1},\omega_{q+1}(\delta_{q+1,k})^{q+1}_{k=1})=1 \quad
\text{for $q=1,2,\dots, \rho$}. \quad   {\square} \qquad \qquad
\endalign$$
\endproclaim \ms

\demo{\bf Proof of Sublemma 7.5} The proof of this sublemma just
follows from Sublemma $5.5$ of Theorem $5.0$. $\square$
\enddemo \ms

{\bf \S7.5. For the proof of Theorem 7.3 }

In this section, we prove Theorem $7.3$ by Sublemma 7.4 and Sublemma
7.5 in {\rm \S7.4}.

\demo{\bf Proof of Theorem 7.3} The proof of the theorem is as
follows:

{\bf[A]} The first aim is to show that two statements in (7.3.1)
and (7.3.2) are equivalent.

{\bf[B]} As a consequence of [A], the second aim is to show easily
that three statements in (7.3.3), (7.3.4) and (7.3.5) are
equivalent. \ms

In preparation for the proof of theorem, for notation, we write
$G_j=y^{\zeta}g_j$ with $1\le j\le r$ and $\Phi_{s}=y^{\eta}\phi_s$
with $1\le s\le \rho$, as we have seen in {\rm(7.3.0)}. It is clear
by Sublemma $5.2$ or Corollary $5.6$ that $y^{\zeta}g_r$ and
$y^{\eta}\phi_\rho$ of (7.3.0) can be rewritten as follows:
$$\align
\text{\rm (7.3.6)} \qquad \qquad G_r &=y^{\zeta}g_r, \\
g_r&=(z^{n_1}+\ve_1 y^{\beta_{1,1}})^{n_2n_3\cdots
n_r}+\sum_{\alpha,\beta\ge 0}
c^{(r)}_{\alpha,\beta}y^{\alpha}z^{\beta} \quad \text{with
$\ve_1=1$} \\
& \text{and \quad
$n_1\alpha+\beta_{1,1}\beta>n_1\beta_{1,1}n_2n_3\cdots n_r$}, \\
\Phi_\rho &= y^{\eta}\phi_\rho, \\
\phi_\rho &=(z^{\ell_1}+\bar{\ve}_1
y^{\delta_{1,1}})^{\ell_2\ell_3\cdots
\ell_\rho}+\sum_{\gamma,\delta\ge 0}
a^{(\rho)}_{\gamma,\delta}y^{\gamma}z^{\delta} \quad \text{with
$\bar{\ve}_1=1$} \qquad \qquad \qquad \qquad \\
&  \text{and \quad
$\ell_1\gamma+\delta_{1,1}\delta>\ell_1\delta_{1,1}\ell_2\ell_3\cdots
\ell_\rho$}, \\
\endalign
$$
where {\rm(i)} $n_1\ge 2$ and $\gcd(n_1,\beta_{1,1})=1$,

{\rm(ii)}  $\ell_1\ge 2$ and $\gcd(\ell_1,\delta_{1,1})=1$,

{\rm(iii)} $\zeta$ is either a positive integer or $0$, and $\eta$
is either a positive integer or $0$,

{\rm(iv)} if $\beta_{1,1}=1$ then $\zeta$ is a positive integer, and
if $\delta_{1,1}=1$ then $\eta$ is a a positive integer,

{\rm(v)} if $\zeta=0$, then $\beta_{1,1}\ge 2$, and also we may
assume for brevity of notation that $n_1<\beta_{1,1}$ ; and if
$\eta=0$, then $\delta_{1,1}\ge 2$, and also we may assume for
brevity of notation that $\ell_1<\delta_{1,1}$,

{\rm(vi)} both $\ve_1=\ve_1(y,z)$ and $\bar{\ve}_1=\bar{\ve}_1(y,z)$
are units in $\BC\{y,z\}$, which may be analytically assumed to be
one, if necessary,

{\rm(vii)}  the $c^{(r)}_{\alpha,\beta}$ are nonzero complex numbers
for some nonnegative integers $\alpha$ and $\beta$ if exist, and the
$a^{(\rho)}_{\gamma,\delta}$ are nonzero complex numbers for some
nonnegative integers $\gamma$ and $\delta$ if exist. \ms

Now, we shall prove firstly the equivalence of the condition in
$[A]$, and secondly that of the condition in $[B]$, respectively.
\ms

{\bf [A]} For proof, we may assume without loss of generality that
$1\le r \le\rho$. Then, for the induction proof on $r\ge 1$, it
suffices to consider two cases, respectively:

Case[I])  $r=1$, and Case[II] $r\ge 1$. \ms

$\underline{\text{\bf Case[I]}}$ Let $r=1\le\rho$. Assume that
\text{$y^{\zeta}g_1 \buildrel \text{{\rm divisor}} \over \sim
y^{\eta}\phi_\rho$} under the standard resolutions as reduced
varieties. By Theorem $7.1$ and Lemma $7.2$, there is nothing to
prove for [A].

$\underline{\text{\bf Case[II]}}$ \quad Let $r\ge 1$. Now, suppose
we have shown for the induction assumption on $r$ that the
equivalence of the condition in {\rm (7.3.2)} for
\text{$V(y^{\zeta}g_r) \buildrel \text{{\rm divisor}} \over \sim
V(y^{\eta}\phi_\rho)$} under the standard resolutions as reduced
varieties in (7.3.1) is true.

Then, it suffices to prove that [A] of the theorem is true on the
positive integer $r+1$, which can be represented as follows: We may
assume that $r+1\le \rho+1$.
$$\align
\text{(7.3.7)} \quad &\text{\text{$G_{r+1}=y^{\zeta}g_{r+1}
\buildrel \text{{\rm divisor}} \over \sim
y^{\eta}\phi_{\rho+1}=\Phi_{\rho+1}$} under the standard resolutions
as reduced varieties} \quad\\
\text{$\iff$} \quad &\text{$\zeta=\eta$, $r+1=\rho+1$, $n_j=\ell_j$
and $\Delta_j(\beta_{j,k})^j_{k=1}=\omega_j(\delta_{j,k})^j_{k=1}$
for $1\le j\le r+1$.} \qquad \qquad
\endalign$$

In order to prove (7.3.7), let $(G_{r+1}\circ\tau_m)_{divisor}$ be a
divisor of $G_{r+1}\circ\tau_m$ defined by
$(G_{r+1}\circ\tau_m)_{divisor}=V^{(m)}(G_{r+1})+\sum^m_{i=1}e_iE_i$,
where each $e_i$ is the multiplicity of $G_{r+1}\circ\tau_m$ along
$E_i$ for $1\le i\le m$ and $V^{(m)}(G_{r+1})$ is the proper
transform of $V(G_{r+1})$ under $\tau_m$. Also, as in Lemma $7.2$,
let $(\Phi_{\rho+1}\circ\mu_{\lambda})_{divisor}$ be a divisor of
$\Phi_{\rho+1}\circ\mu_{\lambda}$ defined by
$(\Phi_{\rho+1}\circ\mu_{\lambda})_{divisor}=V^{(\lambda)}(\Phi_{\rho+1})
+\sum^{\lambda}_{i=1}{\bar{e_i}}{\bar{E_i}}$, where each $\bar{e_i}$
is the multiplicity of $\Phi_{\rho+1}\circ\mu_{\lambda}$ along
$\bar{E_i}$ for $1\le i\le {\lambda}$ and
$V^{({\lambda})}(\Phi_{\rho+1})$ is the proper transform of
$V(\Phi_{\rho+1})$ under $\mu_{\lambda}$.

Following the notations in Sublemma $7.5$, recall that as reduced
varieties,
$$\align
\text{(7.3.8)} \qquad &\text{
$\{vh_r=0\}=\{(G_{r+1}\circ\tau_m)_{total}=0\}$
at $(v,u+1)=(0,0)$}, \qquad \qquad \qquad \qquad\\
& \text{
$\{\bar{v}\psi_{\rho}=0\}=\{(\Phi_{\rho+1}\circ\mu_{\lambda})_{total}=0\}$
at $(\bar{v},\bar{u}+1)=(0,0)$,}
\endalign$$
where $h_r=(g_{r+1}\circ\tau_m)_{proper}$ in $\BC\{v,u+1\}$ and
$\psi_{\rho}=(\phi_{\rho+1}\circ\mu_{\lambda})_{proper}$ in
$\BC\{\bar{v},\bar{u}+1\}$.

It is clear by Lemma $7.2$ and by Sublemma $7.4$ that
$$\align
(7.3.9) \qquad \qquad  & \text{$y^{\zeta}g_{r+1} \buildrel
\text{{\rm divisor}} \over \sim y^{\eta}\phi_{\rho+1}$ under the
standard resolutions
as reduced varieties}  \\
 \quad \iff \quad & \\
(7.3.10) \qquad \qquad   &\text{$\zeta=\eta$, $n_1=\ell_1$,
$\beta_{1,1}=\delta_{1,1}$,
$n_2n_3\cdots n_{r+1}=\ell_2\ell_3\cdots \ell_{\rho+1}$ and}  \\
& \text{$\tau_m=\mu_{\lambda}$ and $e_i=\bar{e_i}$ for $1\le i\le
m=\lambda$,} \\
 & \text{$vh_r \buildrel \text{{\rm divisor}} \over \sim
\bar{v}\psi_{\rho}$ under the standard resolutions as reduced
varieties}
\endalign$$

It is clear that $e_m=n_1\beta_{1,1}n_2\cdots n_{r+1}+\zeta n_1=
\ell_1\delta_{1,1}\ell_2\cdots
\ell_{\rho+1}+\eta\ell_1=\bar{e}_{\lambda}$.

So, it remains to show by (7.3.7), (7.3.9) and (7.3.10) that
$r+1=\rho+1$, $n_j=\ell_j$ and
$\Delta_j(\beta_{j,k})^j_{k=1}=\omega_j(\delta_{j,k})^j_{k=1}$ for
$2\le j\le r+1$.

By \text{\rm The 1-th ${\text{\rm{Cond}}}^{\text{{\rm(1)}}}$}
$\text{\rm of Sequences[I]}^{\text{{\rm(1)}}}$ and by (7.5.1) of
Sublemma $7.5$, along $E_{m}=\{v=0\}$,
$(G_{r+1}\circ\tau_m)_{total}$ can be written in the form
$$\align
(7.3.11) \qquad (G_{r+1}\circ\tau_m)_{total}
&=v^{n_1\beta_{1,1}n_2\cdots n_{r+1}+\zeta
n_1}u^{bn_1n_2\cdots n_{r+1}}(g_{r+1}\circ\tau_m)_{proper}, \qquad \qquad\\
(g_{r+1}\circ\tau_m)_{proper}
&=h_{r-1}^{s_r}+\ve_{1,r+1}v^{\gamma_{r,1}}(u+1)^{\gamma_{r,2}}
h_1^{\gamma_{r,3}}h_2^{\gamma_{r,4}}\cdots
h_{r-2}^{\gamma_{r,r}}=h_r, \qquad \qquad
\endalign$$
where for notation, it may be assumed by Sublemma $5.4$ that
$(g_1\circ\tau_m)_{proper}=(u+1)$ and that $\ve_{1,r+1}$ is a unit
in $\BC\{v,u+1\}$. Note that
$(G_{r+1}\circ\tau_m)_{proper}=(g_{r+1}\circ\tau_m)_{proper}$. \ms

Also, by \text{\rm The 1-th ${\text{\rm{Cond}}}^{\text{{\rm(1)}}}$}
$\text{\rm of Sequences[II]}^{\text{{\rm(1)}}}$ and by (7.5.4) of
Sublemma $7.5$, along $\bar{E}_{\lambda}=\{\bar{v}=0\}$,
$(\Phi_\rho\circ\mu_{\lambda})_{total}$ can be written in the form
$$\align
(7.3.12) \quad \qquad (\Phi_{\rho+1}\circ\mu_{\lambda})_{total}
&=\bar{v}^{\ell_1\delta_{1,1}\ell_2\cdots
\ell_{\rho}+\eta\ell_1}\bar{u}^{\bar{b}\ell_1\ell_2\cdots \ell_\rho}
(\phi_{\rho+1}\circ\mu_{\lambda})_{proper},   \qquad \qquad \quad \\
(\phi_{\rho+1}\circ\mu_{\lambda})_{proper}
&=\psi^{p_{\rho}}_{\rho-1} +\bar{\ve}_{1,{\rho+1}}
v^{\nu_{{\rho},1}}(\bar{u}+1)^{\nu_{{\rho},2}}{\psi^{\nu_{{\rho},3}}_1}
{\psi^{\nu_{{\rho},4}}_2} \cdots
\psi^{\nu_{{\rho},{\rho}}}_{\rho-1}=\psi_{\rho}, \qquad \qquad
\endalign$$
where for notation, it may be assumed by Sublemma $5.4$ that
$(\phi_1\circ\mu_{\lambda})_{proper} =(\bar{u}+1)$ and that
$\bar{\ve}_{1,{\rho+1}}$ is a unit in $\BC\{\bar{v},\bar{u}+1\}$.
Note that
$(\Phi_{\rho+1}\circ\mu_{\lambda})_{proper}=(\phi_{\rho+1}\circ\mu_{\lambda})_{proper}$.

Apply the induction assumption on the integer $r$ to two sequences
in Sublemma $7.5$, which are denoted by {\rm {Sequences[I]}$^{(1)}$}
and {\rm {Sequences[II]}$^{(1)}$}. Because these two sequences, {\rm
{Sequences[I]}$^{(1)}$} and {\rm {Sequences[II]}$^{(1)}$} satisfy
the same kind of five conditions as we have seen in the assumptions
of Theorem $7.3$, then it is clear by (7.3.11) and (7.3.12), and by
induction assumption on the integer $r$ that the following are true:
$$\align
(7.3.13)\qquad \qquad \qquad  & \text{$vh_r \buildrel \text{{\rm
divisor}} \over \sim \bar{v}\psi_{\rho}$ under the standard
resolutions as reduced varieties}
\qquad \qquad \quad \quad \quad\\
 \iff \quad & \\
(7.3.14) \qquad \qquad  \qquad & (a) {\quad} s_i =p_i \quad
\text{for
$1\le i\le r=\rho$.} \qquad \qquad \quad\\
& (b) {\quad}  \gamma_{1,1} =\nu_{1,1} \quad \text{or} \quad
\Xi_1(\gamma_{1,1})=\theta_1(\nu_{1,1}). \\
& (c) {\quad} \Xi_q(\gamma_{q,k})^q_{k=1}
=\theta_q(\nu_{q,k})^q_{k=1} \quad \text{for each
$q=2,3,\dots,r=\rho$}.
\endalign$$
Note that $\gamma_{1,1}
=\Delta^{\sharp}_2(\beta_{2,1},\beta_{2,2})-n_1\beta_{1,1}n_2>0$ and
$\nu_{1,1}
=\omega^{\sharp}_2(\delta_{2,1},\delta_{2,2})-\ell_1\delta_{1,1}\ell_2>0$.

Noting by \text{\rm The 1-th ${\text{\rm{Cond}}}^{\text{{\rm(1)}}}$}
$\text{\rm of Sequences[I]}^{\text{{\rm(1)}}}$ and \text{\rm The
1-th ${\text{\rm{Cond}}}^{\text{{\rm(1)}}}$} $\text{\rm of
Sequences[II]}^{\text{{\rm(1)}}}$ that $s_i=n_{i+1}$ and
$p_i=\ell_{i+1}$ for $1\le i\le r=\rho$, it is clear by $(a)$ of
(7.3.14) that $n_{i+1}=\ell_{i+1}$ for $1\le i\le r=\rho$.

Now, since it is clear by $(c)$ of (7.3.14) that
$\Xi_q(\gamma_{q,k})^q_{k=1}=\theta_q(\nu_{q,k})^q_{k=1}$ for each
$q=1,2,\dots,r$, and also by (7.3.10) that $\zeta=\eta$,
$n_1=\ell_1$, $\beta_{1,1}=\delta_{1,1}$, then by Sublemma $5.5$ and
(7.3.14) we have the following: Note that $2\le q\le r$.
$$\align
(7.3.15)\qquad \qquad \qquad &
\Delta_{q+1}(\beta_{q+1,k})^{q+1}_{k=1}
-n_{q+1}n_q\Delta_{q}(\beta_{q,k})^{q}_{k=1}  \\
=\quad &\Xi_q(\gamma_{q,k})^q_{k=1}
-s_qs_{q-1}\Xi_{q-1}(\gamma_{q-1,k})^{q-1}_{k=1}
\quad \text{by (7.5.2) of Sublemma 7.5} \qquad\\
=\quad & \theta_q(\nu_{q,k})^q_{k=1}-
p_qp_{q-1}\theta_{q-1}(\nu_{q-1,k})^{q-1}_{k=1} \\
=\quad & \omega_{q+1}(\delta_{q+1,k})^{q+1}_{k=1}
-\ell_{q+1}\ell_q\omega_q(\delta_{q,k})^q_{k=1} \quad \text{by
(7.5.5) of Sublemma 7.5.}
\endalign$$

To finish the proof for {\rm Case[II]}, it remains to prove by
(7.3.7) that
$\Delta_{j}(\beta_{j,k})^{j}_{k=1}=\omega_{j}(\delta_{j,k})^{j}_{k=1}$
for each $j=2,3,\dots,r+1$ because
$\Delta_1(\beta_{1,1})=\beta_{1,1}=\omega_1(\delta_{1,1})=\delta_{1,1}$
and $n_j=\ell_j$ for $1\le j\le r+1$ and $\zeta=\eta$ by (7.3.10)
and (7.3.14). Since
$\Delta_{2}(\beta_{2,k})^{2}_{k=1}-n_2n_1\beta_{1,1}=\gamma_{1,1}
=\nu_{1,1}=\omega_2(\delta_{2,k})^2_{k=1}-{\ell}_2{\ell}_1\delta_{1,1}$
by $(b)$ of (7.3.14), then
$\Delta_{2}(\beta_{2,k})^{2}_{k=1}=\omega_2(\delta_{2,k})^2_{k=1}$,
which can be applicable to the equation of (7.3.15) by induction on
the integer $r=\rho$. Since $n_j=\ell_j$ for $1\le j\le r+1$ by
(7.3.14), then it is clear by (7.3.15) that
$\Delta_{j}(\beta_{j,k})^{j}_{k=1}=\omega_{j}(\delta_{j,k})^{j}_{k=1}$
for each $j=1,2,\dots,r+1$. Thus, we finished the proof for {\rm
Case[II]}. So, the proof of [A] is done by Case[I] and Case[II]. \ms

{\bf[B]} As an application of [A], it is clear that three statements
in (7.3.3), (7.3.4) and (7.3.5) are equivalent, and so the proof of
[B] is done.

Therefore, we have completed the proof of the theorem by [A] and
[B]. $\square$
\enddemo \ms

{\bf \S7.6. In preparation for construction of the standard Puiseux
series $\phi_\rho \in \BC\{y,z\}$ such that $\phi_\rho \buildrel
\text{{\rm divisor}} \over \sim g_r$ under the standard resolutions
for given any Puiseux series $g_r \in \BC\{y,z\}$} \ms

In preparation for finding the solution of a unique standard Puiseux
polynomial $\phi_\rho \in \BC[y,z]$ such that $\phi_\rho \buildrel
\text{{\rm divisor}} \over \sim g_r$ under the standard resolutions
for given any Puiseux series $g_r \in \BC\{y,z\}$, first of all, we
need to find Corollary $7.6$, which can be proved by Theorem
$7.6$(The Euclidean Algorithm). Note by Definition $5.0.0$ and
Corollary $5.7$ that $\phi_\rho$ is a Weierstrass polynomial.

After then, we will solve the above problem by Theorem 7.7, using
Corollary $7.6$ only.\ms

\proclaim{Theorem (A well-known theorem)} Let A and B be positive
integers.

{\rm(1)} It is well-known by {\rm Definition $1.14$ of the Euclidean
algorithm} that there are two integers $\gamma$ and $\delta$ such
that $\gcd(A,B)=\gamma A+\delta B$ where $\gcd(A,B)$ is a greatest
common divisor of $A$ and $B$. \ms

{\rm (2)} In addition, we have the following:

{\rm (i)} Let $\gamma$ and $\delta$ be two integers such that
$\gcd(A,B)=\gamma A+\delta B$. Whether $\gamma\delta\not =0$ or not,
we may assume without loss of generality that $\delta$ is chosen
negative.

{\rm (ii)} For example, if A and B are relatively prime, then
$1=\gamma A+\delta B$ for some integers $\gamma$ and $\delta$ where
$\gamma>0$ and $\delta<0$. $\square$
\endproclaim

\proclaim{Corollary 7.6} Let A and B be positive integers which are
relatively prime.

{\rm(1)} If $p$ is a positive integer such that $p>AB$, then there
are two positive integers $s$ and $t$ such that $p=s A+t B$. \ms

{\rm(2)} In particular, under the same assumption as in {\rm (1)},
there is a unique pair of two nonnegative integers $s_1$ and $t_1$
such that $p=s_1A+t_1 B$ with $0\le s_1 <B$ and $t_1>0$. \ms

{\rm(3)} For example, let $p$ be a positive integer such that
$p>nAB$ for some integer $n\ge 2$.

{\rm(3a)} Then it is clear  by {\rm(2)} that there is a unique pair
of two nonnegative integers $s_1$ and $t_1$ such that $p=s_1A+t_1 B$
with $0\le s_1<B$ and $t_1>A$.

{\rm(3b)} There is a finitely different pairs of two nonnegative
integers $(s_1+kB)\ge 0$ and $(t_1-kA)\ge 0$ with $1\le k<n$ such
that $p=(s_1+kB)A+(t_1-kA)B$ with $0\le (s_1+kB)$ and $0\le
(t_1-kA)$. $\square$
\endproclaim

The proof just follows from Theorem(A well-known theorem). \ms

{\bf \S7.7. \text{An algorithm for finding the standard Puiseux
polynomial $\phi_\rho$ such that}

\noindent \text{$\phi_\rho \buildrel \text{{\rm divisor}} \over \sim
g_r$ under the standard resolutions} for given any Puiseux
convergent power series $g_r$} \ms

\proclaim{Theorem 7.7} $\underline{\text{\bf {Assumptions}}}$ By the
same way as in \text{\rm Definition 1.1}, define arbitrary Puiseux
convergent power series $g_r\in \BC\{y,z\}$ of the recursive
$r$-type by {\rm {Sequences[I]}, which consists of three sequences
with the following five conditions(or, by the same way as in the
assumption of Theorem $7.3$, define arbitrary quasi-Puiseux
convergent power series $g_r\in \BC\{y,z\}$ of the recursive
$r$-type by \text{\rm Sequences[I]} with an additional inequality
$2\le n_1< \beta_{11}$ in \text{\rm The 1-th
${\text{\rm{Cond}}}^{\text{{\rm(0)}}}$}, which consists of three
sequences with the following five conditions): \ms

$\text{\bf Sequences[I]}$ Let $\{X_k:k=1,2,\dots,r\}$ with
$X_k\subset N_0$, $\{g_k:k=1,2,\dots,r\}$ with $g_k\in \BC\{y,z\}$
and $\{\text{$\Delta_k:N^k_0\to N_0$ is an integer-valued function
for $k=1,2,\dots,r$}\}$ be three different sequences satisfying the
following five conditions:

\noindent Five conditions are denoted by \text{\bf The 1st
${\text{\bf{Cond}}}^{\text{{\bf(0)}}}$}, \dots, \text{\bf The 5-th
${\text{\bf{Cond}}}^{\text{{\bf(0)}}}$ of Sequences[I]}. \ms

{\bf[I]-(1)} \text{\bf The 1st
${\text{\bf{Cond}}}^{\text{{\bf(0)}}}$ of Sequences[I]:}

{\rm (1a)}  $X_1=\{n_1,\beta_{1,1}\}$ with $2\le n_1< \beta_{1,1}$.

 {\rm (1b)} $X_j=\{n_j,\beta_{j,1},\beta_{j,2},\dots,\beta_{j,j}\}$
 with $n_j\ge 2$ \quad where $j=2,3,\dots,r$.

For each $j\ge 2$, assume that at least one of
$\beta_{j,1},\beta_{j,2},\dots,\beta_{j,j}$ is nonzero. \ms

{\bf [I]-(2)} \text{\bf The 2nd
${\text{\bf{Cond}}}^{\text{{\bf(0)}}}$ of Sequences[I]:}

{\rm (2a)} $g_1=z^{n_1}+y^{\beta_{1,1}}$.

{\rm (2b)} $g_j=g_{j-1}^{n_j}+y^{\beta_{j,1}}z^{\beta_{j,2}}
g_1^{\beta_{j,3}}\cdots g_{j-2}^{\beta_{j,j}}$ \quad where
$j=2,3,\dots,r$. \ms

{\bf [I]-(3)} \text{\bf The 3rd
${\text{\bf{Cond}}}^{\text{{\bf(0)}}}$ of Sequences[I]:}

{\rm(3a)} $\Delta_1(t)=t$ for each $t\in N_0$.

{\rm(3b)}
$\Delta_j(t_j)^j_{k=1}=t_j\Delta_{j-1}(\beta_{j-1,k})^{j-1}_{k=1}
+n_{j-1}\Delta_{j-1}(t_k)^{j-1}_{k=1}$ for each $(t_k)^j_{k=1}\in
N^j_0$

\indent where $j=2,3,\dots,r$. \ms

{\bf [I]-(4)} \text{\bf The 4-th
${\text{\bf{Cond}}}^{\text{{\bf(0)}}}$ of Sequences[I]:}

{\rm (4a)} $\Delta_j(\beta_{j,k})^j_{k=1}>n_jn_{j-1}
\Delta_{j-1}(\beta_{j-1,k})^{j-1}_{k=1}$ for $2\le j\le r$. \ms

{\bf [I]-(5)} \text{\bf The 5-th
${\text{\bf{Cond}}}^{\text{{\bf(0)}}}$ of Sequences[I]:}

{\rm (5a)} $\gcd(n_j,\Delta_j(\beta_{j,k})^j_{k=1})=1$ for $1\le
j\le r$.  \bs

$\underline{\text{\bf {Conclusions}}}$ \quad We have two statements,
denoted by {\rm Fact[I]} and {\rm Fact[II]}. It is clear that {\rm
Fact[I]} and {\rm Fact[II]} are equivalent, and so it remains to
show that {\rm Fact[II] is true.

$\underline{\text{\bf {Fact[I]:}}}$ We can find the computation
algorithm of a unique standard Puiseux polynomial $\phi_r$ in
$\BC[y,z]$ of the recursive $r$-type with $\phi_r \buildrel
\text{{\rm divisor}} \over \sim g_r$ under the standard resolutions,
satisfying the same properties and notations in {\rm Sequences[II]}
of the assumption of Theorem $7.3$, and also the following two
additional conditions. Such two additional conditions are denoted by
\text{\bf The 6-th ${\text{\bf{Cond}}}^{\text{{\bf(0)}}}$ of
Sequences[II]} and \text{\bf The 7-th
${\text{\bf{Cond}}}^{\text{{\bf(0)}}}$ of Sequences[II]} for brevity
of notation. \ms

{\bf [II]-(6)} \text{\bf The 6-th
${\text{\bf{Cond}}}^{\text{{\bf(0)}}}$ of Sequences[II]:}

{\rm(6a)} $2\le \ell_1<\delta_{1,1}$.

\rm{(6b)} $\ell_{j}\ge 2$, $\delta_{j,1}>0$, and $0\le
\delta_{j,k}<\ell_{k-1}$ for $2\le j\le r$ and $2\le k\le j$. \ms

{\bf [II]-(7)} \text{\bf The 7-th
${\text{\bf{Cond}}}^{\text{{\bf(0)}}}$ (for $\phi_r \buildrel
\text{{\rm divisor}} \over \sim g_r$ under the standard
resolutions):}

{\rm(7a)} \text{$n_j=\ell_j$ \quad for each $j=1,2,\dots,r$.}

{\rm(7b)}
\text{$\Delta_j(\beta_{j,k})^j_{k=1}=\omega_j(\delta_{j,k})^j_{k=1}$
\quad  for each $j=1,2,\dots,r$}. \bs

$\underline{\text{\bf {Fact[II]:}}}$ To find the computation
algorithm of $\phi_r(y,z)$ in Fact[I], it remains to apply (3) of
Corollary $7.6$ only to the following step by induction on the
positive integer, and then we can find the desired algorithm in
terms of Sublemma $7.9$. \ms

{\bf Subemma 7.8.} It is clear that $X_1=W_1$ with $2\le
n_1=\ell_1<\beta_{1,1}=\delta_{1,1}$. Whenever each given set $X_j$
satisfies the same conditions in {\rm Sequences[I]} of Theorem $7.3$
with $\Delta_j(\beta_{j,k})^j_{k=1}
>n_jn_{j-1}\Delta_{j-1}(\beta_{j-1,k})^{j-1}_{k=1}$ for
each $j=2,3,\dots,r$, then by (3) of Corollary $7.6$ we can find an
algorithm for computing one and only one solution set
$W_{j}=\{\ell_{j},\delta_{j,1},\delta_{j,2},\dots,\delta_{j,j}\}$
with $\omega_{j}(\delta_{{j},k})^j_{k=1}
>\ell_{j}\ell_{j-1}\omega_{j-1}(\delta_{j-1,k})^{j-1}_{k=1}$
for $W_j$ such that the following conditions are true: \ms

\noindent{\rm(7.8.1)} \qquad {\rm(a)} \quad $n_j=\ell_j\ge 2$ for
$1\le j\le r$.

\qquad \qquad {\rm(b)} \quad
$\Delta_j(\beta_{j,k})^j_{k=1}=\omega_j(\delta_{j,k})^j_{k=1}$ \quad
for each $j=1,2,\dots,r$.

\qquad \qquad {\rm(c) \quad $\delta_{j,1}>0$, and $0\le
\delta_{j,k}<\ell_{k-1}$ for $1\le j\le r$ and $2\le k\le j$. \ms

$\underline{\text{\rm Note}}$ In order to apply (3) of Corollary 7.6
to Subemma 7.8, we use the following:

For $j\ge 2$, $p=\Delta_j(\beta_{j,k})^j_{k=1}>n_jAB$ where
$A=n_{j-1}$ and $B=\Delta_{j-1}(\beta_{j-1,k})^{j-1}_{k=1}$ with
$\gcd(A,B)=1$, and $n_j\ge 2$. In general, $p$ may be chosen
arbitrary such that $p>n_jAB$. \ms

{\bf Subemma 7.9(The computation algorithm for finding $\phi_r(y,z)$
in Fact[I]).}

$\underline{\text{\bf Problem on The 1st step:}}$ Let
$\omega_1:N_0\to N_0$ be a function defined by
$$ \omega_1(t)=t. \tag 7.9.1 $$

For a given set $X_1=\{n_1, \beta_{1,1}\}$ with $2\le
n_1<\beta_{1,1}$, we can compute a unique solution set
$W_1=\{\ell_1, \delta_{1,1}\}$ such that
 $$\align
 \text{$\ell_1=n_1$ and
 $\Delta_1(\beta_{1,1})=\omega_1(\delta_{1,1})$ with $\gcd(\ell_1,
 \delta_{1,1})=1$.} \tag 7.9.2
 \endalign$$

 $\underline{\text{\rm The computation algorithm for the 1st
step:}}$ It is clear. \ms

$\underline{\text{\bf Problem on The 2nd step:}}$ For a solution set
$W_1=\{\ell_1,\delta_{1,1}\}$ in $(7.9.2)$ of the 1st step,  let
$\omega_2:N^2_0\to N_0$ be a function defined by
$$\align
\omega_2(t_1,t_2) &=t_2\delta_{1,1}+\ell_1t_1=\Delta_2(t_1,t_2).
\tag 7.9.3
\endalign$$

For a given set $X_2=\{n_2,\beta_{2,1},\beta_{2,2}\}$ with $n_2\ge
2$, by the computation algorithm we can find a unique solution set
$W_2=\{\ell_2,\delta_{2,1},\delta_{2,2}\}$ with $\ell_2\ge 2$ such
that
$$\align
\text{\rm (7.9.4)} \quad \quad &\text{$\ell_2=n_2$ and
$\Delta_2(\beta_{2,k})^2_{k=1}=\omega_2(\delta_{2,k})^2_{k=1}$ with
$\gcd(\ell_2,\omega_2(\delta_{2,1},\delta_{2,2}))=1$,}  \quad \quad\\
&\text{$\ell_2\ge 2$, $\delta_{2,1}>0$, and $0\le
\delta_{2,2}<\ell_1$.}
\endalign$$

$\underline{\text{\rm The computation algorithm for the 2nd step:}}$
\quad Let $n_2=\ell_2$. Then,
$\Delta_2(\beta_{2,k})^2_{k=1}>2\ell_1\omega_1(\delta_{1,1})=2\ell_1\delta_{1,1}$
with $\gcd(\ell_1,\delta_{1,1})=1$, because of the following:

$\Delta_2(\beta_{2,k})^2_{k=1}>n_2n_1\Delta_1(\beta_{1,1})=n_2n_1\beta_{1,1}$
and $n_2\ge 2$, $n_1=\ell_1$, $\beta_{1,1}=\delta_{1,1}$ and
$\gcd(n_1,\beta_{1,1})=1$.

By (3) of Corollary $7.6$, we can find a unique solution
$\{(\delta_{2,1},\delta_{2,2})\}\subset N^2_0$ such that
$$\align
(7.9.5) \qquad \quad
\omega_2(\delta_{2,k})^2_{k=1}=\Delta_2(\beta_{2,k})^2_{k=1} \quad
\text{with $\delta_{2,1}>\delta_{1,1}>0$ and $0\le
\delta_{2,2}<\ell_1$}, \qquad \qquad \qquad
\endalign$$
where
$\omega_2(\delta_{2,k})^2_{k=1}=\delta_{2,2}\delta_{1,1}+\ell_1\delta_{2,1}$.
\ms

The general case is proved by induction on the positive integer
$j<r$. By a finite number $\frac{j(j-1)}{2}$ iterations of (3) of
Corollary $7.6$, suppose we have shown by induction on the integer
$j<r$ that we can compute a unique set
$W_{j}=\{\ell_{j},\delta_{j,1},\delta_{j,2},\dots,\delta_{j,j}\}$
with $\ell_{j}\ge 2$ as a sequence such that the following
conditions are satisfied:

\noindent{\rm(7.9.6)} \qquad {\rm(a)} $n_j=\ell_j\ge 2$ .

\qquad \qquad {\rm(b)}
$\Delta_j(\beta_{j,k})^j_{k=1}=\omega_j(\delta_{j,k})^j_{k=1}$.

\qquad \qquad  {\rm(c) $\ell_{j}\ge 2$, $\delta_{j,1}>0$, and $0\le
\delta_{j,k}<\ell_{k-1}$ for $2\le k\le j$. \ms

Therefore, whenever $p$ is any integer such that
$p>\ell_{j}\ell_{j-1}\omega_{j-1}(\delta_{j-1,k})^{j-1}_{k=1}$, then
we may assume without loss of generality that by a finite number
$\frac{j(j-1)}{2}$ iteration of (3) of Corollary $7.6$, we can
compute $\{t_k:k=1,2,\dots,j\}\subset N$ such that $t_1>0$,
$t_k<\ell_{k-1}$ for $k=2,3,\dots,j$ and $p=\omega_j(t_k)^j_{k=1}$,
using the same method as we have done in the $j$-th step.

Then, it remain to prove that we can find the algorithm for the
\text{\rm problem on (j+1)-th} step, as follows:

$\underline{\text{\bf Problem on The {(j+1)}-th step:}}$ From a
unique solution set
$W_{j}=\{\ell_{j},\delta_{{j,1}},\delta_{{j},2},\dots,\delta_{{j},{j}}\}$
in the $j$-th step, let $\omega_{j+1}:N^{j+1}_0\to N_0$ be a
function defined by
$$
\omega_{j+1}(t_k)^{j+1}_{k=1}
=t_{j+1}\omega_{j}(\delta_{j,k})^{j}_{k=1}
+\ell_{j}\omega_{j}(t_k)^{j}_{k=1}. \tag 7.9.7
$$

For a given set
$X_{j+1}=\{n_{j+1},\beta_{j+1,1},\beta_{j+1,2},\dots,\beta_{j+1,j+1}\}$
with $n_{j+1}\ge 2$, we can compute a unique solution set
$W_{j+1}=\{\ell_{j+1},\delta_{j+1,1},\delta_{j+1,2},\dots,\delta_{j+1,j+1}\}$
with $\ell_{j+1}\ge 2$ such that
$$\align
&\text{$\ell_{j+1}=n_{j+1}$,
$\Delta_{j+1}(\beta_{j+1,k})^{j+1}_{k=1}=\omega_{j+1}(\delta_{{j+1,}k})^{j+1}_{k=1}$
with
$\gcd(\ell_{j+1},\omega_{j+1}(\delta_{{j+1,}k})^{j+1}_{k=1})=1$,} \tag 7.9.8\\
&\text{$\ell_{j+1}\ge 2$, $\delta_{{j+1,}1}>0$ and $0\le
\delta_{{j+1,}k}<\ell_{k-1}$ for $2\le k<j+1$.}
\endalign$$ \ms

$\underline{\text{\rm The computation algorithm for the (j+1)-th
step:}}$ \quad Let $n_{j+1}=\ell_{j+1}$.

\noindent Then, $\Delta_{j+1}(\beta_{{j+1},k})^{j+1}_{k=1}>2\ell_{j}
\omega_{j}(\delta_{{j,}k})^{j}_{k=1}$ with
$\gcd(\ell_{j},\omega_{j}(\delta_{{j,k}})^{j}_{k=1}))=1$, because of
the following:
$$\align
\text{\rm(7.9.9)} \quad
&\text{$\Delta_{j+1}(\beta_{{j+1},k})^{j+1}_{k=1}
>n_{j+1}n_{j}\Delta_{j}(\beta_{{j},k})^{j}_{k=1}
=\ell_{j+1}\ell_{j}\omega_{j}(\delta_{{j,}k})^{j}_{k=1}$,
$n_{j+1}=\ell_{j+1}\ge 2$,} \qquad \quad\\
&\text{$n_{j}=\ell_{j}$,
$\Delta_{j}(\beta_{{j},k})^{j}_{k=1}=\omega_{j}(\delta_{{j,}k})^{j}_{k=1}$
and $\gcd(n_{j},\Delta_{j}(\beta_{{j},k})^{j}_{k=1})=1$.}
\endalign$$

Applying (3) of Corollary $7.6$ to (7.9.9) once, we can find a
unique solution $\{p,\delta_{{j+1},{j+1}}\}\subset N^2_0$ such that
$$\align
& \Delta_{j+1}(\beta_{{j+1},k})^{j+1}_{k=1}
=\delta_{{j+1},{j+1}}\omega_j(\delta_{j,k})^j_{k=1}+p\ell_j \quad
\text{with} \tag 7.9.10 \\
&\text{$p>\omega_j(\delta_{j,k})^j_{k=1}$, $0\le
\delta_{{j+1},{j+1}}<\ell_j$}. \qquad \qquad \qquad
\endalign$$

Also,  note by(7.9.10) that
$$\align
 p &>\omega_j(\delta_{j,k})^j_{k=1}
=\Delta_{j}(\beta_{{j},k})^{j}_{k=1} \tag 7.9.11\\
&>n_jn_{j-1}\Delta_{j-1}(\beta_{j-1,k})^{j-1}_{k=1}=\ell_j\ell_{j-1}
\omega_{j-1}(\delta_{j-1,k})^{j-1}_{k=1}.
\endalign$$

Applying (3) of Corollary $7.6$ to (7.9.11), then  we assume by
induction method on the positive integer $j$-th step that by a
finite number iterations $\frac{j(j-1)}{2}$ of (3) of Corollary
$7.6$, we can find a unique solution $(t_1,t_2,\dots, t_j)\in
N^j_0$, denoted by
$(\delta_{j+1,1},\delta_{j+1,2},\dots,\delta_{j+1,j})$, such that
$$\align
(7.9.12) \qquad  p=\omega_j(\delta_{j+1,k})^j_{k=1} \quad \text{with
$\delta_{j+1,1}>0$, and $0\le \delta_{j+1,k}<\ell_{k-1}$ for $2\le
k\le j$}. \qquad \qquad \qquad \qquad \qquad
\endalign$$

By (7.9.10) and (7.9.12), we get the following:
$$\align
(7.9.13) \qquad   \Delta_{j+1}(\beta_{{j+1},k})^{j+1}_{k=1}
&= \delta_{{j+1},{j+1}}\omega_j(\delta_{j,k})^j_{k=1}+p\ell_j  \\
&= \delta_{{j+1},{j+1}}\omega_j(\delta_{j,k})^j_{k=1}
+\ell_j\omega_{j}(\delta_{j+1,k})^{j}_{k=1} =
\omega_{j+1}(\delta_{j+1,k})^{j+1}_{k=1}. \qquad \qquad
\endalign$$

Thus, with a finite number $\frac{(j+1)j}{2}$ iteration of (3) of
Corollary $7.6$, the computation algorithm for finding
$\phi_{j+1}(y,z)$ or $W_{j+1}$ on the (j+1)-th step can be
completely finished. $\square$
\endproclaim
\ms

\demo{\bf Proof of Theorem 7.7} It is clear that {\rm Fact[I]} and
{\rm Fact[II]} in the conclusion of Theorem $7.7$ are equivalent.
So, as far as the proof of Theorem $7.7$ are concerned, it suffices
to show that {\rm Fact[II]} is true. But, the proof of {\rm
Fact[II]} just follows from (3) of Corollary $7.6$. Therefore, there
is nothing to prove for the theorem. $\square$
\enddemo \bs

\vfill \pagebreak

{\bf Part[B3] Review on the definition of the Puiseux pairs and
the multiplicity sequences} \bs

{\bf \S8. The definition of the multiplicity and Puiseux
exponents(the Puiseux pairs) and the review of theorems about an
equivalence of the multiplicity and Puiseux exponents and the
multiplicity sequence for irreducible curves} \ms

{\bf \S8.0. Introduction } \ms

In this section, in preparation for a good success of the main
algorithm in this paper, first of all, we will review the complete
proof of the following theorem(Theorem A) only, without using the
well-known theorem(Theorem B):

{\bf Theorem A:} Whenever any two irreducible parametrization have
the same Puiseux pairs by a nonsingular change of the
parametrization(equivalently, the same multiplicity and Puiseux
exponents in the sense of Definition $8.1$, later), then they have
the same multiplicity sequences, and conversely.

{\bf Theorem B:} As far as arbitrary Puiseux expansion of
irreducible plane curve singularities is concerned, any two
irreducible plane curve singularities have the same topological
types if and only if they have the same Puiseux pairs. \ms

So, the first part is to find an equivalence of irreducible
parametrization for irreducible plane curve singularities by Theorem
$8.8$. After then, as an application, the second part is to review
the proof of Theorem A(Theorem 8.10) in an elementary way, without
using the well-known theorem(Theorem B or Theorem 8.5).

For example, assume that the standard Puiseux expansion of an
irreducible curve $C_1$ is given by $y(t)=t^n$ and
$z(t)=t^{\alpha_1}+t^{\alpha_2}+\cdots +t^{\alpha_r}$ where $2\le n
<\alpha_1<\cdots <\alpha_r$ and $n>d_1>\cdots >d_r=1$ with
$d_i=\gcd(n,\alpha_1,\dots,\alpha_i)$, $1\le i\le r$, and also that
the parametrization of another irreducible plane curve $C_2$ is
given by $y(t)=t^n(1+H(t))$ and $z(t)=t^{\alpha_1}$ where
$1+H(t)=1+t^{\alpha_2-\alpha_1}+\cdots+t^{\alpha_r-\alpha_1}$. Then,
it is easily shown by Theorem A that two irreducible curves $C_1$
and $C_2$ have the same multiplicity sequence, and also the same
Puiseux pairs by a nonsingular change of a parameter, without using
Theorem B. But, without using Theorem A, it has been not yet proved
by Theorem B only that these two irreducible curves $C_1$ and $C_2$
have the same multiplicity sequences. \ms

{\bf \S8.1. The definition and known preliminaries} \ms

Now, in order to avoid the complexity of the terminology in this
section, first of all, we can rewrite the statement of the
definition of the Puiseux pairs. \ms

\definition{Definition 8.1} Let the
parametrization for arbitrary irreducible curve $C$ be defined by
$$ \text{$y(t)=t^n$ and
$z(t)=c_1t^{k_1}+c_2t^{k_2}+\cdots=c_1t^{k_1}(1+H(t))$,} \tag
8.1.1$$ where $1<n$, $1<k_1<k_2<\cdots, $ and the $c_i$ are nonzero
complex numbers and $H(t)$ is just the substitution.

Moreover, note that the curve $C$ is irreducible in $\BC\{y,z\}$
$\iff$ $n\ge \gcd(n,k_1)\ge \gcd(n,k_1,k_2)\ge \cdots \ge
\gcd(n,k_1,k_2,\dots)=1$. \ms

Now, consider two cases, respectively.

Case[I] Let $n\le k_1$. Then, the parametrization for the curve
$C$ of $(8.1.1)$ is called the Puisuex expansion.

Case[II] Let $n>k_1$. Then, the parametrization for the curve $C$
of $(8.1.1)$ is not called the Puisuex expansion. \ms

{\bf Case[I]} Assume that $n\le k_1$. Now, we can define the
sequence $\{\gamma_1,\gamma_2,\dots,\gamma_p\}$ from the set
$\{k_i:i=1,2,\dots\}$, consisting of the exponents of the above
parameter $t$, as follows: Note that $n$ is the multiplicity of
the curve $C$ at the origin.

\noindent $(\ast)$ $\gamma_1$ is the smallest positive integer among
the exponents $k_i$ such that $n>\gcd(n,k_i)$; $\gamma_2$ is the
smallest positive integer among the exponents $k_i$ such that $n>
\gcd(n,\gamma_1)>\gcd(n,\gamma_1,k_i)$; \dots; $\gamma_p$ is the
smallest positive integer among the exponents $k_i$ such that $n>
\gcd(n,\gamma_1)> \gcd(n,\gamma_1,\gamma_2)>
\gcd(n,\gamma_1,\gamma_2,\gamma_3)
>\cdots>\gcd(n,\gamma_1,\gamma_2,\dots,\gamma_p)=1$.

{\rm (1)} By the uniqueness of construction of the set
$\{\gamma_i: 1\le i\le p\}$, $\gamma_i$ is called $i$-th Puiseux
exponent in this paper. \ms

{\rm (2)} By (1), let S be the set defined by
$\{n,\gamma_1,\gamma_2,\dots,\gamma_p \}$. Whenever the Puiseux
expansion for the curve $C$ is given, the set $S$ is uniquely
determined by the curve $C$.

{\rm (2a)} $S$ is called the multiplicity and Puiseux exponents for
the curve $C$, that is, a new terminology.

{\rm (2b)} As in (2a), the parametrization defined by $y=t^n$ and
$z=t^{\gamma_1}+t^{\gamma_2}+\dots+t^{\gamma_p}$ is called the
standard Puiseux expansion for the curve $C$. Note that $2\le
n<\gamma_1<\gamma_2<\cdots <\gamma_p$ and
$n>\gcd(n,\gamma_1)>\gcd(n,\gamma_1,\gamma_2)>\cdots
>\gcd(n,\gamma_1,\gamma_2,\dots,\gamma_p)=1$. \ms

{\rm (3)} By (2), let $d_i=\gcd(n,\gamma_1,\dots,\gamma_i)$ for
$1\le i\le p$, and write $d_0=n$ for brevity of notation.

Define $\lambda_i$ and $\mu_i$ by $\lambda_i=\dfrac{\gamma_i}{d_i}$
and $\mu_i=\dfrac{d_{i-1}}{d_i}$ for $1\le i\le p$, and let
$(\lambda_i,\mu_i)$ be defined by the Puiseux pair for each $i$.
Then, $\{(\lambda_i,\mu_i):i=1,2,\dots,p\}$ is called a finite
sequence of Puiseux pairs for the curve $C$.  \ms

{\bf Case[II]} Assume that $n>k_1$. The generalization of
definitions for the Puiseux pairs and the standard Puiseux expansion
for irreducible curves is as follows.

For the convenience of the notation, we may begin without loss of
generality that the parametrization of the pair $(y(t),z(t))$ for
the curve $C$ of $(8.1.1)$ is rewritten in the following:
$$\align
y(t)=t^m,\quad z(t)=b_1t^{\beta_1}+b_2t^{\beta_2}+\cdots, \quad
\text{with $m>\beta_1$} \tag 8.1.2
\endalign$$
where the $b_i$ are nonzero complex numbers, and $m>1$ and
$1<\beta_1<\beta_2<\cdots$, and $m \ge \gcd(m,\beta_1) \ge
\gcd(m,\beta_1,\beta_2) \ge \cdots \ge
\gcd(n,\beta_1,\beta_2,\dots)=1$.

By $(8.1.2)$, let $s$ be the new parameter defined by a conformal
mapping
$$\align
s(t)=t(b_1+\sum_{i\ge2}b_it^{\beta_i-\beta_1})^{\frac{1}{\beta_1}}
\tag 8.1.3
\endalign$$
of $t$ at the origin such that $z(t)=(s(t))^{\beta_1}$ and
$s(0)=0$, and let $t=\phi(s)$ be its inverse. \ms

Then, the Puiseux expansion defined by $y_1(s)=y(\phi(s))$ and
$z_1(s)=z(\phi(s))$, which is equivalent to the parametrization of
the pair $(y(t),z(t))$ in $(8.1.2)$, can be written as follows:
$$\align
z_1(s)=s^{\beta_1}, \quad
y_1(s)=c_1s^{\ell_1}+c_2s^{\ell_2}+\cdots, \quad \text{with
$\beta_1<\ell_1$} \tag 8.1.4
\endalign$$
where $1<m=\ell_1<\ell_2<\cdots $, and $\beta_1<\ell_1$, and the
$c_i$ are nonzero complex numbers. \ms

Therefore, if $m=\ell_1$ is greater than $\beta_1$, first we will
find the inverse $t=\phi(s)$ of a conformal mapping $s=s(t)$ in
(8.1.3) by using (8.8.3) of Theorem $8.8$ in this section, and next
compute an algorithm for the construction of the Puiseux expansion
in $(8.1.4)$ of Theorem $8.8$, that is, an equivalent
parametrization for the above curve $C$.

Next, by the same way as we have used in Case[I] for the definition
of the words in (8.1.5), and by Definition 8.9 and Theorem 8.10, we
can naturally generalize the definition of the words in (8.1.5) for
this curve $C$ of $(8.1.2)$ in Case[II], respectively: \ms

\noindent(8.1.5) \quad The multiplicity and Puiseux exponents; the
standard Puiseux expansion; a finite sequence of the Puiseux pairs.
$\square$
\enddefinition \ms

\definition{Remark 8.1.1 for Case[I]} {\rm(i)} It can be easily
shown that there is a one-to-one correspondence between the set of
the multiplicity and Puiseux exponents, and the set of Puiseux pairs
because (2) and (3) have the same type of definitions
arithmetically.

{\rm (ii)} If the parametrization defined by $(y(t),z(t))$ in
(8.1.1) is the Puiseux expansion, then it is said that this Puiseux
expansion have either the multiplicity and Puiseux exponents
$\{n,\gamma_1,\gamma_2,\dots,\gamma_p\}$ or the Puiseux pairs
$\{(\lambda_i,\mu_i):i=1,2,\dots,p\}$ where each $\lambda_i$ and
$\mu_i$ is defined as we have seen in {\rm(3)}.

{\rm (iii)} By (i) of the above remark, throughout this paper, we
prefer to choose the terminology in (2) rather than that in (3), if
necessary. $\square$ \enddefinition \ms

\definition{Definition 8.2} Definition 8.2 has
the same statement as Definition 1.2 of $\S1$ does in Family(3) with
respect to the definition of the multiplicity sequences of
irreducible plane curves with isolated singularity under the
standard resolution. $\square$
\enddefinition \ms

\proclaim{Theorem 8.3(Enriques-Chisini)}

{\rm(i)} For an irreducible curve with Puiseux expansion
$$\align
x &=t^m  \tag 8.3.1 \\
y &=a_1t^{k_1}+a_2t^{k_2}+\cdots +a_qt^{k_q},
\endalign$$ in which only essential (characteristic) term appear,
the multiplicity sequence is determined by the following chain of
$g$ Euclidean algorithms: Let $i=1,2,\dots$.
$$\align
\lambda_i &=\mu_{i,1}r_{i,1}+r_{i,2}, \\
r_{i,1} &=\mu_{i,2}r_{i,2}+r_{i,3}, \\
&\ldots\ldots \\
r_{i,w(i)-1} &=\mu_{i,w(i)}r_{i,w(i)} \qquad
\text{with \quad $0\le r_{i,j+1}<r_{i,j}$}, \\
\lambda_i &=k_i-k_{i-1} \quad \text{for ~
$1\le i\le g$, \quad and \quad $k_0=0$,} \\
r_{i,1} &=r_{i-1,w(i-1)} \quad \text{for ~ $i>1$, \quad and \quad
$r_{1,1}=m$.}
\endalign$$
In the multiplicity sequence, the multiplicity $r_{ij}$ then appears
$\mu_{ij}$ times,  where $i=1,\dots,g; j=1,\dots,w(i)$. {\rm(}If a
certain multiplicity arises from several successive algorithms, then
it is also counted multiply.{\rm)}

{\rm(ii)} For an arbitrary irreducible curve one obtains the
multiplicity  sequences by omitting all non-characteristic terms
from the Puiseux expansion and then applying the algorithm above.

{\rm(iii)} Conversely, one can reconstruct the exponents of the
characteristic terms of the Puiseux expansion of an irreducible
curve, i.e. the Puiseux pairs of the curve, from the multiplicity
sequence, by the chain of Euclidean algorithms. $\square$
\endproclaim
\demo{\bf Proof of Theorem 8.3}  See  [Bri-Kn].
\enddemo \ms

\definition{Definition 8.4} Let $V(f)=\{(y,z): f(y,z)=0\}$ and
$V(g)=\{(y,z): g(y,z)=0\}$ be analytic varieties at the origin in
$\BC^2$ where $f$ and $g$ are analytically irreducible in
$\BC\{y,z\}$  with isolated singularity at the origin $\BC^2$.

{\rm(i)} $f$ and $g$ are said to have the same topological type of
the singularity at the origin if there is a germ at the origin of
homeomorphisms $\phi:(U_1,0)\to (U_2,0)$ such that $\phi(V)=W$ and
$\phi(0)=0$ where $U_1$ and $U_2$ are open subsets in $\BC^{2}$. In
this case, denote this relation by $f\sim g$ or $V\sim W$.
Otherwise, we write $f\not \sim g$ or $V\not \sim W$.

{\rm(ii)} $f$ and $g$ are said to have the same analytic type of the
singularity at the origin if there is a germ at the origin of
biholomorphisms $\psi:(U_1,0)\to (U_2,0)$ such that $\psi(V)=W$ and
$\psi(0)=0$ where $U_1$ and $U_2$ are open subsets in $\BC^{2}$,
that is, $f\circ \psi =ug$ where $u$ is a unit in ${}_{2}\CO$, the
ring of germs of holomorphic functions at the origin in $\BC^{2}$.
In this case, denote this relation by $f\approx g$ or $V\approx W$.
Otherwise, we write $f\not \approx g$ or $V\not \approx W$.

{\rm(iii)} Following Definition $8.2$, for notation, \text{\rm
Multiseq(V(f))} or \text{\rm \{[Mult(V(f))]\}} is called the
multiplicity sequence of $V(f)$. If $V(f)$ and $V(g)$ have the same
multiplicity sequences, we write either \text{\text{\rm
Multiseq(V(f))}$\equiv$ \text{\rm Multiseq(V(g))}} as sequence or
$f\sim g$ (Multiseq). Otherwise, we write either \text{\text{\rm
Multiseq(V(f))}{$\not \equiv$}\text{\rm Multiseq(V(g))}} as sequence
or $f\not \sim g$ (Multiseq). $\square$
\enddefinition \ms

\proclaim{Theorem 8.5([Br],[Bu],[Z1])} Let $f(y,z)$ be irreducible
in $\BC\{y,z\}$ with isolated singularity at the origin in $\BC^2$.
Then the curve defined by $f$ at the origin can be described
topologically by $y=t^n$ and
$z=t^{\alpha_1}+t^{\alpha_2}+\cdots+t^{\alpha_p}$ where
$n<\alpha_1<\cdots<\alpha_p$ and
$n>\gcd(n,\alpha_1)>\cdots>\gcd(n_1,\alpha_1,\dots,\alpha_p)=1$. If
for a given $f$ there is another homeomorphic parametrization
defined by $y=t^m$ and $z=t^{\beta_1}+\cdots+t^{\beta_q}$ where
$m<\beta_1<\cdots<\beta_q$ and
$m>\gcd(m,\beta_1)>\cdots>\gcd(m,\beta_1,\dots,\beta_q)=1$, then
$n=m$, and $p=q$ and $\alpha_i=\beta_i$ for $1\le i\le p$.

Note by {\rm Definition $8.1$} and {\rm Remark $8.1.1$} that the
multiplicity and Puiseux exponents for the standard Puiseux
expansion determine the topological types of the irreducible plane
curve singularities, and conversely. $\square$
\endproclaim
 \ms

{\bf{\S8.2. How to get an equivalent parametrization from given any
irreducible parametrization by the inverse mapping theorem of one
complex variable}} \ms \ms

\definition{Definition 8.6} Let $\phi(t)$ be an analytic function in
a neighborhood of zero such that $\phi(0)=\phi'(0)=
\cdots=\phi^{(k)}(0)=0$, but $\phi^{(k+1)}(0) \not =0$. Then, it is
said that $\phi(t)$ has a multiplicity $k$ at $t=0$ and write
$\text{\rm Mult}(\phi(t),0)=k$ for notation. Let $f(y,z)$ be in
$\BC\{y,z\}$. It is said that $f(y,z)$ has a multiplicity $\nu$ at
$(y,z)=(0,0)$, denoted by $\text{\rm Mult}(f(y,z),(0,0))=\nu$, if
there is the least integer $\nu$ such that some partial derivative
of $f$ of order $\nu$ is nonzero at the origin. $\square$
\enddefinition

\proclaim{Lemma 8.7(The rearrangement of an irreducible
parametrization)}

$\underline{\text{\bf Assumptions}}$ Let the curve $V=\{f(y,z)=0\}$
with $f(y,z) \in \BC\{y,z\}$ have an irreducible parametrization at
the origin, which is defined by
$$\align
 \text{$y(t)=t^n$ and
$z(t)=c_1t^{k_1}+c_2t^{k_2}+\cdots=c_1t^{k_1}(1+H(t))$,} \tag 8.7.1
\endalign$$
where the $c_i$ are nonzero complex numbers and $1\le n$, $1\le
k_1<k_2<\cdots,$ and $n \ge \gcd(n,k_1) \ge \gcd(n,k_1,k_2) \ge
\cdots \ge \gcd(n,k_1,k_2,\dots)=1$.

To get a desired rearrangement of $y=t^n$ and
$z=\sum^{\infty}_{i=1}c_it^{k_i}$ in the conclusion of this lemma,
first we can define a finite sequence $\{\alpha_1,\alpha_2,
\dots,\alpha_{r+1} \}$ from the sequence $\{k_i:i=1,2,\dots\}$
consisting of the exponents $k_i$ in $(8.7.1)$ as follows: \ms

{\rm(1)} Let $\alpha_1=k_1$, and then note that $n \ge
\gcd(n,\alpha_1)$. That is, either $n=\gcd(n,\alpha_1)$ or $n>
\gcd(n,\alpha_1)$.

{\rm(2)} Let $\alpha_2$ be the smallest positive integer among the
exponents $k_i$ such that $n \ge
\gcd(n,\alpha_1)>\gcd(n,\alpha_1,k_i)$.

$\ldots\ldots$

{\rm(r+1)} Let $\alpha_{r+1}$ be the smallest positive integer among
the exponents $k_i$ such that $n \ge \gcd(n,\alpha_1)>
\gcd(n,\alpha_1,\alpha_2)> \cdots> \gcd(n,\alpha_1,\alpha_2,
\dots,\alpha_r)> \gcd(n,\alpha_1,\alpha_2, \dots,$
\text{$\alpha_r,k_i$}{\rm {)}}$=1$. \ms

Let $d$ and $k$ be arbitrary positive integers. For brevity of
notation, if $k$ is divisible by $d$, then we write $d \vert k$.
Otherwise, we write $d \not \vert k$.

Now, let $d_i=\gcd(n,\alpha_1,\dots,\alpha_i)$ for $1 \le i \le
r+1$,and then $n \ge d_1>d_2> \cdots>d_{r+1}$. Note that $d_i \vert
(\alpha_i-\alpha_1)$, $d_i \not \vert (\alpha_{i+1}-\alpha_1)$, and
$d_{i+1} \vert d_i$. \ms

$\underline{\text{\bf Conclusions}}$ The irreducible parametrization
of $V$ can be rearranged in $t$ as follows:
$$\align
(8.7.2) \qquad\qquad y=& t^n \quad \text{and}  \\
z=&
ct^{\alpha_1}\{(1+D_1(t))+t^{\alpha_2-\alpha_1}(c_{2,0}+D_2(t))+\cdots
\\
&+t^{\alpha_r-\alpha_1}(c_{r,0}+D_r(t))+t^{\alpha_{r+1}-\alpha_1}(c_{r+1,0}
+D_{r+1}(t))\}  \qquad \qquad \qquad \qquad\qquad\\
\endalign$$
satisfying the properties {\rm (i)},{\rm (ii)},\dots,{\rm (v)}.

{\rm(i)} $1\le n$ and $1\le\alpha_1 <\alpha_2<\cdots <\alpha_{r+1}$.

{\rm(ii)} $n\ge d_1> d_2>\cdots > d_{r+1}=1$ with
$\gcd(n,\alpha_1,\alpha_2,\dots,\alpha_i)=d_i$ for $1\le i\le r+1$.

{\rm(iii)} $p_1,p_2,\dots,p_r$ are nonnegative integers such that
$p_1d_1 <\alpha_2-\alpha_1<(p_1+1)d_1$,
$p_2d_2<\alpha_3-\alpha_2<(p_2+1)d_2$,\dots,
$p_rd_r<\alpha_{r+1}-\alpha_r<(p_r+1)d_r$.

{\rm(iv)} If $p_j\neq 0$ for some $j\le r$, write
$D_j(t)=\sum^{p_j}_{i=1}c_{j,i}t^{id_j}
 \in \BC[t]$, and if $p_j= 0$ for some $j\le r$, write
$D_j(t)=0$, and also $D_{r+1}(t)=\sum^\infty_{k=1}c_{r+1,k}t^k \in
\BC\{t\}$.

{\rm(v)} $c,c_{1,0}=1, c_{2,0},c_{3,0},\dots,c_{r+1,0}$
 are all nonzero complex numbers. $\square$
\endproclaim

For the proof of Lemma 8.7, see Lemma $3.3$([K2]). \ms

\proclaim{Theorem 8.8(An equivalence of irreducible
parametrization)}

$\underline{\text{\bf Assumptions}}$ \quad We may assume without
loss of generality that the curve $V=\{f(y,z)=0\}$ with $f(y,z) \in
\BC\{y,z\}$ at the origin has an irreducible parametrization as
follows:
$$\align
(8.8.1) \qquad y=& t^n \quad \text{and}  \\
z=&
ct^{\alpha_1}\{(1+D_1(t))+t^{\alpha_2-\alpha_1}(c_{2,0}+D_2(t))+\cdots
\\
&+t^{\alpha_r-\alpha_1}(c_{r,0}+D_r(t))+t^{\alpha_{r+1}-\alpha_1}(c_{r+1,0}
+D_{r+1}(t))\}=ct^{\alpha_1}(1+H(t)),   \\
\endalign$$
satisfying the properties {\rm (i)},{\rm (ii)},\dots,{\rm (v)}.

{\rm(i)} $1\le n$ and $1\le\alpha_1 <\alpha_2<\cdots <\alpha_{r+1}$.

{\rm(ii)} $n\ge d_1> d_2>\cdots > d_{r+1}=1$ with
$\gcd(n,\alpha_1,\alpha_2,\dots,\alpha_i)=d_i$ for $1\le i\le r+1$.

{\rm(iii)} $p_1,p_2,\dots,p_r$ are nonnegative integers such that
$\alpha_1+p_1d_1 <\alpha_2<\alpha_1+(p_1+1)d_1$,
$\alpha_2+p_2d_2<\alpha_3<\alpha_2+(p_2+1)d_2$,\dots,
$\alpha_r+p_rd_r<\alpha_{r+1}<\alpha_r+(p_r+1)d_r$.

{\rm(iv)} If $p_j\neq 0$ for some $j\le r$, write
$D_j(t)=\sum^{p_j}_{i=1}c_{j,i}t^{id_j}
 \in \BC[t]$, and if $p_j= 0$ for some $j\le r$, write $D_j(t)=0$,
and also $D_{r+1}(t)=\sum^\infty_{k=1}c_{r+1,k}t^k \in \BC\{t\}$.
Note that $1+H(t)=1+D_1(t)+t^{\alpha_2-\alpha_1}(c_{2,0}+D_2(t))+
\cdots
 \quad +t^{\alpha_r-\alpha_1}(c_{r,0}+D_r(t))+t^{\alpha_{r+1}
-\alpha_1}(c_{r+1,0}+D_{r+1}(t))$.

{\rm(v)} $c,c_{1,0}=1, c_{2,0},c_{3,0},\dots,c_{r+1,0}$
 are all nonzero complex numbers. \ms

$\underline{\text{\bf Conclusions}}$ \quad In preparation for the
construction of an equivalent irreducible parametrization of $V$,
let $s$ be the new parameter defined by
$$ s(t)=c^{\frac{1}{\alpha_1}}t(1+H(t))^{\frac{1}{\alpha_1}} \tag
8.8.2
$$
where {\rm(i)} $c^{\frac{1}{\alpha_1}}$ is a complex root such that
$\omega^{\alpha_1}=c$,

{\rm(ii)} $s=s(t)$ is a conformal mapping of $t$ at the origin,

{\rm (iii)} $z=s^{\alpha_1}$. \ms

Then, we have the following {\rm(I)} and {\rm(II)}:

{\bf (I)}
$t=c^{-\frac{1}{\alpha_1}}s(1+H(t))^{-\frac{1}{\alpha_1}}$, as
$t=\phi(s) \in \BC\{s\}$, can be written as follows: Note that
$y=(\phi(s))^n$.
$$\align
(8.8.3) \qquad t= &\phi(s) \\
=&c^{-\frac{1}{\alpha_1}}s\{1+Q_1(s)+s^{\alpha_2-\alpha_1}(B_{2,0}+Q_2(s))
\\ &+ \cdots
+s^{\alpha_r-\alpha_1}(B_{r,0}+Q_r(s))+s^{\alpha_{r+1}
-\alpha_1}(B_{r+1,0}+Q_{r+1}(s))\} \qquad \qquad \qquad \\
\endalign$$
satisfying the properties {\rm (i)} and {\rm (ii)}. \ms

\noindent{\rm (i)}
$B_{2,0}=\frac{c_{2,0}}{-\alpha_1}(c^{-\frac{1}{\alpha_1}})^{\alpha_2-\alpha_1}$,
$B_{3,0}=\frac{c_{3,0}}{-\alpha_1}(c^{-\frac{1}{\alpha_1}})^{\alpha_3-\alpha_1}$,
\dots,
$B_{r+1,0}=\frac{c_{r+1,0}}{-\alpha_1}(c^{-\frac{1}{\alpha_1}})
^{\alpha_{r+1}-\alpha_1}$.

\noindent{\rm (ii)} $Q_1(s)=B_{1,1}s^{d_1}+B_{1,2}s^{2d_1}+
\cdots+B_{1,p_1}s^{p_1d_1} \in \BC[s]$,
$Q_2(s)=B_{2,1}s^{d_2}+B_{2,2}s^{2d_2}+ \cdots+B_{2,p_2}s^{p_2d_2}
\in \BC[s]$, \dots, $Q_r(s)=B_{r,1}s^{d_r}+B_{r,2}s^{2d_r}+
\cdots+B_{r,p_r}s^{p_rd_r} \in \BC[s]$, $Q_{r+1}(s)=
\sum^{\infty}_{k=1}B_{r+1,k}s^{k} \in \BC \{s\}$ such that all the
$B_{i,j}$ are complex numbers and that in particular the $B_{i,0}$
are nonzero for $2 \le i \le r+1$. Note that $Q_i(0)=0$ for $1\le i
\le r+1$. \ms

{\bf(II)} Then, the equivalent parametrization with the new
parameter $s$ for $V$ can be analytically written in the following
form:
$$\align
 z=&s^{\alpha_1},  \tag 8.8.4 \\
y=&c^{-\frac{n}{\alpha_1}}s^n\{1+Q^*_1(s)
+s^{\alpha_2-\alpha_1}(b_{2,0}+Q^*_2(s)) \\
&+s^{\alpha_3-\alpha_1}(b_{3,0}+Q^*_3(s))+
\cdots+s^{\alpha_{r+1}-\alpha_1}(b_{r+1,0}+Q^*_{r+1}(s))\}  \\
\endalign$$
satisfying the properties {\rm (i)} and {\rm (ii)}. \ms

\noindent{\rm(i)}
$b_{2,0}=\frac{n}{-\alpha_1}c_{2,0}(c^{-\frac{1}{\alpha_1}})^{\alpha_2-\alpha_1}$,
$b_{3,0}=\frac{n}{-\alpha_1}c_{3,0}(c^{-\frac{1}{\alpha_1}})^{\alpha_3-\alpha_1}$,
\dots,
$b_{r+1,0}=\frac{n}{-\alpha_1}c_{r+1,0}(c^{-\frac{1}{\alpha_1}})^{\alpha_{r+1}-\alpha_1}$.

\noindent{\rm(ii)} $Q^*_1(s)=b_{1,1}s^{d_1}+b_{1,2}s^{2d_1}+
\cdots+b_{1,p_1}s^{p_1d_1} \in \BC[s]$,
$Q^*_2(s)=b_{2,1}s^{d_2}+b_{2,2}s^{2d_2}+ \cdots+b_{2,p_2}s^{p_2d_2}
\in \BC[s]$,\dots,$Q^*_{r+1}(s)=\sum^{\infty}_{k=1}b_{r+1,k}s^k \in
\BC \{s\}$, such that all the $b_{i,j}$ are complex numbers and that
in particular the $b_{i,0}$ are nonzero for $2 \le i \le r+1$. Note
that $Q^*_i(0)=0$ for all $i=1,2,\dots,r+1$. $\square$
\endproclaim

\definition{Remark 8.8.1} Observe by $(8.8.3)$ and $(8.8.4)$ that
$b_{20}=nB_{20},b_{30}=nB_{30},\dots,b_{r+1,0}=nB_{r+1,0}$. For the
proof of Theorem 8.8, see Theorem $3.4$([K2]). $\square$
\enddefinition \ms

{\bf \S 8.3. The generalization of definitions for the Puiseux pairs
and the standard Puiseux expansion for irreducible curves} \ms

Now, as an application of Theorem $8.8$, we are going to generalize
the words in (8.1.5) of Definition $8.1$ by the following.

\definition{Definition 8.9} Let $f(y,z)$ be analytically
irreducible in $\BC\{y,z\}$ with isolated singularity at the origin
in $\BC^2$. By Lemma $8.7$, we may assume without loss of generality
that the curve $V(f)$ defined by the above $f$ at the origin has an
irreducible parametrization as follows:
$$
(8.9.1) \quad \text{$V(f):=$} \left\{\eqalign{ y= &t^n \cr z=
&a_1t^{\alpha_1}(1+D_1(t))+a_2t^{\alpha_2}(1+D_2(t))+\cdots
+a_{r+1}t^{\alpha_{r+1}}(1+D_{r+1}(t)), \cr} \right.  $$

satisfying the properties, {\rm(1a)},{\rm(1b)},\dots, {\rm(1e)}.

\roster \item "(1a)" $2\le n$ and $2\le \alpha_1<\alpha_2<\cdots
<\alpha_{r+1}$.

\item "(1b)" $n\ge d_1>d_2>\cdots
>d_{r+1}=1$ with $d_i=\gcd(n,\alpha_1,\alpha_2,\dots,\alpha_i)$ for $1\le i\le
r+1$.
\item "(1c)" The $a_i$ are all nonzero numbers for
$i=1,2,\dots,r+1$.

\item "{\rm(1d)}" Define
$p_1,p_2,\dots,p_r$ to be nonnegative integers such that $p_1d_1
<\alpha_2-\alpha_1<(p_1+1)d_1$, $p_2d_2
<\alpha_3-\alpha_2<(p_2+1)d_2$, \dots, $p_rd_r
<\alpha_{r+1}-\alpha_r<(p_r+1)d_r$.

\item "{\rm(1e)}" ~ If $p_j\neq 0$ for some $j\le r$, write
$D_j(t)=\sum^{p_j}_{i=1}c_{j,i}t^{id_j}
 \in \BC[t]$, and if $p_j= 0$ for some $j\le r$, write
$D_j(t)=0$, and also $D_{r+1}(t)=\sum^\infty_{k=1}c_{r+1,k}t^k \in
\BC\{t\}$.

\item "{\rm(1f)}" ~ $c,c_{10}=1, c_{20},c_{30},\dots,c_{r+1,0}$
 are all nonzero complex numbers.
\endroster

Then, the multiplicity and Puiseux exponents for the curve $V(f)$
are defined as follows:

\roster \item "{\rm (A)}" If $n\le \alpha_1$ and $n$ is not a
divisor of $\alpha_1$, then note that the parametrization defined by
$y=t^n$ and $z=t^{\alpha_1}+t^{\alpha_2}+\cdots +t^{\alpha_{r+1}}$
is called the standard Puiseux expansion. Then, it is said that the
set $\{n,\alpha_1,\alpha_2,\dots,\alpha_{r+1}\}$ is a finite
sequence of the multiplicity and Puiseux exponents for the Puiseux
expansion of $V(f)$.

\item "{\rm (B)}" If $n\le \alpha_1$ and $n$ is a divisor of
$\alpha_1$, then note that the parametrization defined by $y=t^n$
and $z=t^{\alpha_2}+t^{\alpha_3}+\cdots +t^{\alpha_{r+1}}$ is called
the standard Puiseux expansion. Then, it is said that the set
$\{n,\alpha_2,\alpha_3,\dots,\alpha_{r+1}\}$ is a finite sequence of
the multiplicity and Puiseux exponents for the Puiseux expansion of
$V(f)$.
\endroster

In case $n>\alpha_1$, using the equation of (8.8.4) in the
conclusion of Theorem $8.8$, we can compute the Puiseux expansion
which is equivalent to the parametrization of $V(f)$, as follows:
$$
\text{$V(f):\approx$} \left\{\eqalign{ z=& s^{\alpha_1}\cr
y=&c_1^{-\frac{n}{\alpha_1}}s^n\{(1+Q^*_1(s))
+s^{\alpha_2-\alpha_1}(b_{20}+Q^*_2(s))\cr &+s^{\alpha_3-\alpha_1}
(b_{30}+Q^*_3(s)) +\cdots
+s^{\alpha_{r+1}-\alpha_1}(b_{r+1,0}+Q^*_{r+1}(s))\} \cr
=&c_1^{-\frac{n}{\alpha_1}}s^n\{1+L(s)\}, \cr} \right. \tag 8.9.2
$$

where \roster \item "{\rm(i)}" $\gcd
(n,\alpha_1,\alpha_2-\alpha_1,\dots,\alpha_i-\alpha_1)=
\gcd(n,\alpha_1,\alpha_2,\dots,\alpha_i)=d_i$ for $1\le i\le r+1$,

\item"{\rm (ii)}" $Q^*_j(s)=\sum^{p_j}_{i=1}b_{j,i}s^{id_j}\in
\BC[s]$ for $1\le j\le r$ and
$Q^*_{r+1}(s)=\sum^{\infty}_{i=1}b_{r+1,i}s^i\in \BC\{s\}$,

\item"{\rm(iii)}" all the $b_{j,i(j)}$ are complex numbers with
$1\le j\le r+1$ and $1\le i(j)\le p_j$, noting that $p_{r+1}$ may be
infinite,

\item "{\rm(iv*)}" the $b_{j,0}$ are all nonzero complex numbers
for $2\le j\le r+1$, noting by (8.8.3) that
$b_{j,0}=\frac{n}{-\alpha_1}c_{j0}c^{\frac{1}{-\alpha_1}
(\alpha_j-\alpha_1)}$ for $2\le j\le r+1$,

\item"{\rm (v)}" $L(s)$ is just the substitution.
\endroster \ms

By the same method as we have done in two cases (A) and (B), then it
is enough to consider the following cases:

\roster \item "{\rm(C)}" If $n>\alpha_1$ and
$\alpha_1>\gcd(n,\alpha_1)$, then it is said that $y=t^{\alpha_1}$
and $z=t^n+t^{n+\alpha_2-\alpha_1}+t^{n+\alpha_3-\alpha_1}+\cdots
+t^{n+\alpha_{r+1}-\alpha_1}$ is the standard Puiseux expansion.
Then, it is said that the set
$\{\alpha_1,n,n+\alpha_2-\alpha_1,n+\alpha_3-\alpha_1,\dots,
n+\alpha_{r+1}-\alpha_1\}$ is a finite sequence of the multiplicity
and Puiseux exponents for the curve $V(f)$.

\item "{\rm (D)}" If $n>\alpha_1$ and $\alpha_1$ is a divisor of
$n$, then it is said that $y=t^{\alpha_1}$ and
$z=t^{n+\alpha_2-\alpha_1}+t^{n+\alpha_3-\alpha_1}+\cdots
+t^{n+\alpha_{r+1}-\alpha_1}$ is the standard Puiseux expansion.
Then, it is said that the set
$\{\alpha_1,n+\alpha_2-\alpha_1,n+\alpha_3-\alpha_1,\dots,
n+\alpha_{r+1}-\alpha_1\}$ is a finite sequence of the multiplicity
and Puiseux exponents for the curve $V(f)$. $\square$
\endroster
\enddefinition \ms

{\bf \S 8.4. The review of theorems about an equivalence of the
Puiseux pairs and the multiplicity sequence for irreducible curves}
\ms

Finally, as mentioned in the beginning of this section, we will
review the statement of Theorem $5.1$[K] without its proof,
recalling that Theorem $5.1$[K] is equivalent to Theorem $A$. So, to
finish this section, we just rewrite Theorem $5.1$[K] by the
following theorem(Theorem $8.10$), which was already proved by
Theorem $8.8$(an equivalence of irreducible plane curve
singularities) and $\sigma-$ process only, without using the
well-known theorem(Theorem $B$). \ms

\proclaim{Theorem 8.10(Theorem A, An equivalence of the Puiseux
expansions with the multiplicity and Puiseux exponents in the sense
of Definition 8.1 and the multiplicity sequences for irreducible
parametrization by Theorem 5.1[K])}

$\underline{\text{\bf Assumptions}}$ \quad Let $f(y,z)$, $g(y,z)$
and $h(y,z)$ be analytically irreducible in $\BC\{y,z\}$ with
isolated singularity at the origin in $\BC^2$. Assume that three
curves $V(f)$, $V(g)$ and $V(h)$ defined by the above analytic
functions $f$, $g$ and $h$ at the origin have irreducible
parametrization, respectively as follows:

{\rm(1)} Let the parametrization of $V(f)$ be the Puiseux expansion
with the multiplicity and Puiseux exponents
$\{n,\alpha_1,\alpha_2,\dots,\alpha_{r+1}\}$, defined by the
following:
$$
(8.10.1) \quad \text{$V(f):=$} \left\{\eqalign{ y= &t^n \cr z=
&a_1t^{\alpha_1}(1+D_1(t))+a_2t^{\alpha_2}(1+D_2(t))+\cdots
+a_{r+1}t^{\alpha_{r+1}}(1+D_{r+1}(t)), \cr} \right.  $$

satisfying the properties, {\rm(1a)},{\rm(1b)},\dots, {\rm(1e)}.

\roster \item "{\rm(1a)}" $2\le n<\alpha_1<\alpha_2<\cdots
<\alpha_{r+1}$.

\item "{\rm (1b*)}" $n>d_1>d_2>\cdots
>d_{r+1}=1$ with $d_i=\gcd(n,\alpha_1,\alpha_2,\dots,\alpha_i)$ for $1\le i\le
r+1$.

\item "{\rm(1c)}" The $a_i$ are all nonzero numbers for
$i=1,2,\dots,r+1$.

\item "{\rm(1d)}" Define
$p_1,p_2,\dots,p_r$ to be nonnegative integers such that $p_1d_1
<\alpha_2-\alpha_1<(p_1+1)d_1$, $p_2d_2
<\alpha_3-\alpha_2<(p_2+1)d_2$, \dots, $p_rd_r
<\alpha_{r+1}-\alpha_r<(p_r+1)d_r$.

\item "{\rm(1e)}" If $p_j\neq 0$ for some $j\le r$, write
$D_j(t)=\sum^{p_j}_{i=1}a_{j,i}t^{id_j}
 \in \BC[t]$, and if $p_j= 0$ for some $j\le r$, write $D_j(t)=0$,
and also $D_{r+1}(t)=\sum^\infty_{k=1}a_{r+1,k}t^k \in \BC\{t\}$.
\endroster

{\rm{Remark}}: In the above condition {\rm(1b*)} of {\rm(8.10.1)},
if $n\ge \gcd(n,\alpha_1)$ and $n$ is a divisor of $\alpha_1$, then
it is clear that the singularity of $V(f)$ is analytically invariant
at the origin, whether or not $a_1$ is zero, and so from the
beginning we may assume without loss of generality that
$n>\gcd(n,\alpha_1)$ and $a_1\not =0$. \ms

{\rm(2)} Let the parametrization of $V(g)$ be the Puiseux expansion
with the multiplicity and Puiseux exponents
$\{m,\beta_1.\beta_2,\dots,\beta_{u+1}\}$, defined by the following:
$$
(8.10.2) \quad \text{$V(g):=$} \left\{\eqalign{ y= &t^m \cr z=
&b_1t^{\beta_1}(1+L_1(t))+b_2t^{\beta_2}(1+L_2(t))+\cdots
+b_{u+1}t^{\beta_{u+1}} (1+L_{u+1}(t)), \cr} \right.  $$

satisfying the properties, {\rm(2a)},{\rm(2b)},\dots, {\rm(2e)}.

\roster \item"{\rm(2a)}" $2\le m<\beta_1<\beta_2<\cdots
<\beta_{u+1}$,

\item"{\rm(2b)}" $m>e_1>e_2>\cdots >e_{u+1}=1$
with $e_i=\gcd(m,\beta_1,\beta_2,\dots,\beta_i)$ for $1\le i\le
u+1$.

\item"{\rm(2c)}" The $b_i$ are all nonzero numbers for
$i=1,2,\dots,u+1$.

\item"{\rm(2d)}" Define
$q_1$,$q_2$,\dots,$q_u$ to be nonnegative integers such that $q_1e_1
<\beta_2-\beta_1<(q_1+1)e_1$, $q_2e_2 <\beta_3-\beta_2<(q_2+1)e_2$,
\dots, $q_ue_u <\beta_{u+1}-\beta_u<(q_u+1)e_u$.

\item"{\rm(2e)}" If $q_j\neq 0$ for some $j\le u$, write
$L_j(t) =\sum^{q_j}_{i=1}b_{j,i}t^{ie_j}\in \BC[t]$, and if $q_j= 0$
for some $j\le u$, write $L_j(t) =0$,  and also $L_{u+1}(t)
=\sum^{\infty}_{i=1}b_{u+1,i}t^i \in \BC\{t\}$.
\endroster \ms

{\rm (3)} Let the parametrization of $V(h)$ be defined by the
following:
$$
(8.10.3) \quad \text{$V(h):=$} \left\{\eqalign{ y=
&c_1t^{\ell_1}(1+R_1(t)) +c_2t^{\ell_2}(1+R_2(t))+\cdots
+c_{v+1}t^{\ell_{v+1}}(1+R_{v+1}(t)) \cr z= &t^{\gamma}, \cr}
\right.  $$

satisfying properties, {\rm(3a)},{\rm(3b)},\dots, {\rm(3e)}.

\roster \item "{\rm(3a)}" $2\le \ell_1<\gamma$ and
$\ell_1<\ell_2<\dots <\ell_{v+1}$.

\item "{\rm (3b)}" $\ell_1\ge
\tau_1>\tau_2>\cdots
>\tau_{v+1}=1$  with $\tau_i=\gcd(\gamma,\ell_1,\ell_2,\dots,\ell_i)$
for $1\le i\le v+1$.

\item "{\rm(3c)}" the $c_i$ are all nonzero numbers for
$i=1,2,\dots,v+1$.

\item "{\rm(3d)}" Define $\ve_1$, $\ve_2$,\dots, $\ve_v$ to be
nonnegative integers such that $\ve_1\tau_1
<\ell_2-\ell_1<(\ve_1+1)\tau_1$, $\ve_2\tau_2
<\ell_3-\ell_2<(\ve_2+1)\tau_2$, \dots, $\ve_v\tau_v
<\ell_{v+1}-\ell_v<(\ve_v+1)\tau_v$.

\item "{\rm(3e)}" If $\ve_j\neq 0$ for some $j\le v$, write
write $R_j(t)=\sum^{\ve_j}_{i=1}c_{j,i}t^{i\tau_j}\in \BC[t]$, and
if $\ve_j= 0$ for some $j\le v$, write $R_j(t) =0$, and also
$R_{v+1}(t) =\sum^{\infty}_{i=1}c_{v+1,i}t^i\in \BC\{t\}$.
\endroster  \ms

$\underline{\text{\bf Conclusions}}$ \quad We get the following:

\noindent {{\bf (I)}} \quad Note that $n>\gcd(n,\alpha_1)$ and
$m>\gcd(m,\beta_1)$.
$$\align
(8.10.4) \quad \quad & \quad \text{\text{\rm Multiseq(V(f))} =
\text{\rm Multiseq(V(g))} as  sequence}    \\
\iff &\quad \text{ the multiplicity and Puiseux exponents are the
same,
by Definition $8.1$,} \qquad \qquad \qquad\\
& \quad \text{that is, $n=m$,~ $r+1=u+1$, ~ and ~ $\alpha_i=\beta_i$
\quad for all $i=1,2,\dots,r+1$} \\
\iff &\quad \text{the Puiseux pairs for both $V(f)$ and $V(g)$ are
the same.}
\endalign$$

\noindent {{\bf (II)}} \quad Let $\gamma>\ell_1\ge 2$. Then, it is
enough to consider two cases:

\quad{\rm{(IIa)}} $\ell_1>\gcd(\gamma,\ell_1)$. \quad {\rm{(IIb)}}
$\ell_1=\gcd(\gamma,\ell_1)$, that is, $\ell_1$ is a divisor of
$\gamma$. \ms

\noindent{{\bf(IIa)}}\quad Let $\ell_1>\gcd(\gamma,\ell_1)$.
$$\align
(8.10.5) \quad \quad  &  \quad \text{\text{\rm Multiseq(V(f))} =
\text{\rm Multiseq(V(h))} as sequence}  \qquad \qquad \\
\iff & \quad n=\ell_1, \alpha_1=\gamma, r+1=v+1 ~ and ~
\alpha_i=\gamma+\ell_i-\ell_1 \quad \text {for} ~ 1\le i \le
r+1 \qquad \qquad \qquad \\
\iff &\quad \text{the Puiseux pairs for both $V(f)$ and $V(h)$ are
the same.}
\endalign$$

\noindent{{\bf(IIb)}}\quad Let $\ell_1=\gcd(\gamma,\ell_1)$, that
is, $\ell_1$ is a divisor of $\gamma$.
$$\align
(8.10.6) \quad \quad & \quad \text{\text{\rm Multiseq(V(f))} =
\text{\rm Multiseq(V(h))} as sequence}  \qquad \qquad   \\
\iff & \quad n=\ell_1,\alpha_1=\gamma+\ell_2-\ell_1, r+1=v ~ and ~
\alpha_i=\gamma+\ell_{i+1}-\ell_1 ~ for~ 2\le i\le r+1 \qquad \qquad \qquad \\
\iff & \quad \text{the Puiseux pairs for both $V(f)$ and $V(h)$ are
the same. \qquad$\square$}
\endalign$$
\endproclaim

\demo{\bf Proof of Theorem 8.10} See Theorem $5.1$([K2]).
\enddemo \ms

\vfill \pagebreak

{\bf \S 9 New definition of the join of subsequences of a finite
sequence and its application to the representation of the
multiplicity sequences for irreducible plane curve singularities}
\ms

{\bf \S9.0. Introduction } \ms

In preparation for finding a computation algorithm for the
multiplicity sequences of all the irreducible plane curve
singularities, using the Euclidean algorithm of an integer ring $Z$,
first we are going to construct two new terminology in Definition
$9.0$ and Definition $9.1$. \ms

{\bf \S9.1. New Definitions} \ms

\definition{Definition 9.1}  Let $S=\{e_i \in {N}:i=1,2,\dots,q\}$ be
a finite sequence of positive integers. Then, it is said that {\rm
S} is the join of r subsequences of $S$ in order, denoted by
$S=\text{\rm Join}\{B_1,B_2\dots,B_r\}$ of r subsequences in order,
where each $B_i$ is a subsequence of $S$ for $i=1,2,\dots,r$, if the
following properties are satisfied:

{\rm(a)} Let $q=\lambda_r$. For each $j=1,2,\dots,r$, define the
number of elements of $B_j$ by $\lambda_j-\lambda_{j-1}$, which is
positive. Note that $\lambda_0=0$.

{\rm(b)} Each subsequence $B_i$ of $S$ can be written as follows:
Let $q=\lambda_r$.
$$\align
(9.1.1) \qquad \qquad
& B_1=\{b_{0,i}=e_i: i=1,2,\dots,\lambda_1 \},  \qquad \qquad  \\
&
B_j=\{b_{j-1,i}=e_{\lambda_{j-1}+i}:i=1,2,\dots,(\lambda_j-\lambda_{j-1})\}
\quad \text{for $j=2,3,\dots,r$.} \qquad \qquad
\endalign$$
\enddefinition \ms

\definition{Definition 9.2}
(1) Let $\beta$ and $m$ be two arbitrary positive integers such that
$\beta\ge m$. To find $\gcd(\beta,m)$, it suffices to consider the
following case except for that $\beta=m=1$:

Let $\beta\ge m\ge 1$ with ${\beta}m>1$. To find $\gcd(\beta,m)$
with ${\beta}m>1$, then the Euclidean algorithms for $\beta=\mu_0$
and $m=\mu_1$ can be written as follows:
$$\align
\beta&=q_1\mu_1+\mu_2 \quad \text{with $0\le
\mu_2<\mu_1$,} \tag 9.2.1 \\
\mu_1 &=q_2\mu_2+\mu_3 \quad
\text{with $0\le\mu_3<\mu_2$,} \\
\mu_2 &=q_3\mu_3+\mu_4 \quad \text{with
$0\le\mu_4<\mu_3$,}\\
\dots\\
 \mu_{w-2} &=q_{w-1}\mu_{w-1}+\mu_w \quad \text{with
$0\le\mu_w<\mu_{w-1}$,}\\
 \mu_{w-1} &=q_w\mu_w+0 \qquad \qquad
\text{~with
$\gcd(\beta,m)=\mu_w$,} \\
\endalign$$
where $S=\{\mu_i:i=1,2,\dots,w\}$ is a strictly decreasing finite
sequence.

Then, it is said that $\mu_w$ is the greatest common divisor of
$\beta$ and $m$, denoted by $\gcd(\beta,m)$. \ms

{\rm(2)} Let  $m$ and $\beta$ be any two positive integers such that
$\beta\ge m\ge 1$ and ${\beta}m>1$. By (9.2.1), {\text{\bf the
Euclidean multiplicity sequence for two positive integers $\beta$
and $m$}} with {$\gcd(\beta,m)$}, denoted by either one of four
notations \text{Ems$[\beta:m]$}, $\{[\beta:m]\}$,
\text{Ems$[m:\beta]$} and $\{[m:\beta]\}$, is defined by the
following:
$$\align
 \text{$\{[\beta:m]\}$}
&=\{\mu_1,\mu_1,\dots,\mu_1; \mu_2,\mu_2,\dots,\mu_2;\dots;
 \mu_{w},\mu_{w},\dots,\mu_{w}\} \tag 9.2.2\\
&=\text{Join}\{S_1,S_2, \dots,S_w\} \quad \text{by Definition 9.1},
\endalign$$
where the sequence $S_i=\{\mu_i:i=1,2,\dots,w\}$ is a subsequence of
$\{[\beta:m]\}$ for each $i=1,2,\dots,w$.

For example, we use the same kind of notations as in {\rm(9.2.3)},
as follows:

Let $\gcd(\beta,m)=1$ with $\beta\ge m\ge 1$ and ${\beta}m>1$ and
$\{[\beta,m]\}=\{c_1,c_2,\dots,c_t\}$ by (9.2.2). For brevity of
notation, the following may be rewritten by
$$
(9.2.3) \qquad \qquad
\text{$\{[d\beta,dm]\}=\{d[\beta,m]\}=\{dc_1,dc_2,\dots,dc_t\}$ for
any positive integer $d$.} \qquad \qquad
$$

{\rm(3)} Let $V(g)=\{(y,z):g(y,z)=z^m+y^{\beta}=0\}$ be irreducible
in $\BC\{y,z\}$ with isolated singularity at $0\in \BC^2$ where
$m\ge 2$ and $\beta\ge 2$ are positive integers and
$\gcd(m,\beta)=1$. In this case, it is said by (9.2.1) and (9.2.2)
that $\text{\rm Multiseq$(V(g))$}$, called the multiplicity sequence
of $V(g)$, can be rewritten by either $\{[\beta:m]\}$ or
$\{[m:\beta]\}$, as a sequence. \ms
\enddefinition
\ms

\definition{Remark 9.2.1}

(i) For any positive integers $\beta$ and $m$ such that $\beta\ge
m$, it is well-known by {\rm(1)} of Definition $9.2$ that there are
two integers $\gamma$ and $\delta$ such that $\gcd(\beta,m)={\gamma}
\beta+{\delta}m$.

(ii) Assuming that $\beta>m\ge 2$ for given any Euclidean
multiplicity sequence \text{Ems$[\beta:m]$} and that $m$ is not a
divisor of $\beta$, to count the number $\sigma$ of
$d=\min\{\text{Ems$[\beta:m]$}\}$ in \text{Ems$[\beta:m]$} as a
sequence, then it is enough to compute $\sigma=\dfrac{\min{\{a_i\in
E
:a_i>d\}}}{d}=\dfrac{\mu_{w-1,q_{w-1}}}{\mu_{w,q_w}}=\dfrac{\mu_{w-1}}{\mu_w}$
from (9.2.1) and (9.2.2), because if $m$ is a divisor of $\beta$
then \text{Ems$[\beta:m]$} is equal to a one-point set.

(iii)  Let \text{Eds$[\beta:m]$} and \text{Eds$[\beta':m']$} be two
arbitrary Euclidean multiplicity sequences where $m\le \beta$ and
$m'\le \beta'$. Then, \text{Eds$[\beta:m]$} and
\text{Eds$[\beta':m']$} are the same Euclidean multiplicity
sequences if and only if $m=m'$ and $\beta=\beta'$. Note by
definition that \text{Eds$[m:\beta]$} and \text{Eds$[\beta:m]$} are
the same Euclidean multiplicity sequences.
\enddefinition \ms

{{\bf \S9.2 Some examples for the representation of the multiplicity
sequence for any irreducible plane curve singularity in terms of a
collection of subsequences of $D$ with the subdivision}} \ms

\definition{Definition 9.3}
Let $V(f)=\{(y,z):f(y,z)=0\}$ be an analytic variety at
$(y,z)=(0,0)$ in $\BC^2$ where $f(y,z)$ is irreducible in
$\BC\{y,z\}$ with isolated singularity at the origin in $\BC^2$.

(i) Following Definition $8.2$, brevity for notation, either
\text{\rm Multiseq(V(f))} or \text{\rm \{[Mult(V(f))]\}} is said to
be the multiplicity sequence of $V(f)$.

(ii) As we have seen in the definitions just before Theorem $2.2$,
let $V^{(\sigma)}(f)$ be the $\sigma$-th proper transform of $V(f)$
at the origin, whenever any finite suitable number $\sigma$
iterations of blow-ups in process of the standard resolution of the
singular point $(0,0)$ of $V(f)$ are chosen arbitrary. Following the
notation in (i), \text{\rm Multiseq($V^{(\sigma)}(f)$)} or \text{\rm
\{[Mult($V^{(\sigma)}(f)$)]\}} is said to be the multiplicity
sequence of $V^{(\sigma)}(f)$.
\enddefinition \ms

In this section, in preparation for finding a solution of
Problem[2], by using Definition $9.2$ and Definition $9.3$, observe
the following proposition and corollaries in order.

\proclaim{Proposition 9.4} $\underline{\text{\bf Assumptions}}$ \
Let $V(f)=\{(y,z): f(y,z)=0\}$ be an analytic variety at $(0,0)$ in
$\BC^2$ defined by $f(y,z)=a_0z^n+a_1y^{\alpha_1}z^{n-1}+
\dots+a_ny^{\alpha_n}$ in $\BC\{y,z\}$, where each $a_i$ is a unit
in $\BC\{y,z\}$ if exists, and the $\alpha_i$ are positive integers.

Let $1\le n< k=\alpha_n$ with $d=\gcd(n,k)$, and write $n=n_1d$ and
$k=k_1d$. \ms

{\rm(1)}{\rm(1a)} By Theorem $3.6$ of $\S 3$, we may assume without
any need of proof that if $f$ is irreducible in $\BC\{y,z\}$, then f
can be represented as follows:
$$\align
(9.4.1) \qquad \qquad f =A(z^{n_1}+ \xi y^{k_1})^d +\sum_{\alpha,
\beta\ge 0}
   c_{\alpha\beta}y^{\alpha}z^{\beta} \quad \text{with} \quad
   n_1\alpha+k_1\beta>n_1k_1d, \qquad \qquad
 \endalign$$
where the $c_{\alpha\beta}$ are nonzero complex numbers for some
nonnegative integers $\alpha$ and $\beta$, and $A$ and $\xi$ are the
unique nonzero complex numbers.

{\rm(1b)} If $\gcd(n,\alpha_n)=1$, the necessary condition for $f$
to be irreducible in $\BC\{y,z\}$ in {\rm(9.4.1) is sufficient. Note
that $n=n_1$ and $k=k_1$ with $d=1$. \ms

{\rm (2)} As in the assumption of Theorem $3.7$ of $\S 3$, let
$V(G)=\{(y,z):G(y,z)=0\}$ be an analytic variety at $(0,0)$ in
$\BC^2$ defined by the form
$$
\align
 G=& z^{\gamma}g, \tag 9,4.2  \\
 g_1=& z^{n_1}+y^{k_1} \quad \text{with} \quad \gcd(n_1,k_1)=1,
\endalign
$$
satisfying the following properties:

\roster

\item "(i)"  $1\le n_1 <k_1$.

\item "(ii)" If $n_1=1$, then $\gamma=1$.

\item "(iii)" If $n_1 \ge 2$, then $\gamma=0$.
\endroster

Let $\tau_m$ be the composition of a finite number $m$ of successive
blow-ups which is needed to get the standard resolution of the
singular point of $V(G)$ in $(9.4.2)$ by the same way as we have
used in the assumption of Theorem $3.7$.

For brevity, let $V^{(t)}(G)$ be the proper transform under $\tau_t$
for $1\le t\le m$. For each $t=1,2,\dots,m$, suppose that
$\tau_{t}:M^{(t)}\to\BC^2$ satisfies the same properties and
notations as in the assumption of Theorem $3.7$. \ms

{\rm (3)} As we have seen in Definition $9.3$, let $V^{(\sigma)}(f)$
be the $\sigma-th$ proper transform of $V(f)$ at the origin, and
\text{\rm \{[Mult($V^{(\sigma)}(f)$)]\}} be the multiplicity
sequence of $V^{(\sigma)}(f)$, whenever any finite suitable number
$\sigma$ iterations of blow-ups in process of the standard
resolution of the singular point $(0,0)$ of $V(f)$ are chosen
arbitrary. \ms

$\underline{\text{\bf {Conclusions}}}$ \quad Assuming that $f$ is
analytically irreducible in $\BC\{y,z\}$  with isolated singularity
at the origin $\BC^2$, then we have the following:
$$\align
\text{\rm(9.4.3)} &\quad \text{If $1\le n_1<k_1$, then $\text{\rm
\{[Mult($V(f)$)]\}} =\text{\rm{Join}}(\{[dn_1:dk_1]\},\text{\rm
\{[Mult($V^{(m)}(f)$)]\}} )$.} \qquad \qquad \\
&\quad \text{If $2\le n_1<k_1$, then $\text{\rm
\{[Mult($V(G)$)]\}}=\text{\rm
\{[Mult($V(g_1)$)]\}}=\{[n_1:k_1]\}$.} \\
\endalign$$
For brevity of notation, we write $\{[dn_1:dk_1]\}=\{d[n_1:k_1]\}$,
if necessary. $\square$
\endproclaim
\ms

\demo{\bf Proof of Proposition 9.4} By Theorem $3.7$ or Corollary
$3.8$, the proof is done. $\square$
\enddemo \ms

\definition{Remark 9.4.1}
{\rm(a)} Let $f(y,z)=a_0z^n+a_1y^{\alpha_1}z^{n-1}+
\dots+a_ny^{\alpha_n}$ be in $\BC\{y,z\}$ where $2\le n<\alpha_n$,
each $a_i$ is a unit in $\BC\{y,z\}$ if exists, and the $\alpha_i$
are positive integers.

If $\gcd(n,\alpha_n)=1$, then one of the necessary and sufficient
condition for $f$ to be irreducible in $\BC\{y,z\}$ is that $f=0$
and $z^n+y^{\alpha_n}=0$ have the same multiplicity sequence,
denoted by $\{[n,\alpha_n]\}$, by Proposition $9.4$ or Corollary
$3.3$. $\square$
\enddefinition \ms

By Proposition $9.4$, it is easy to get the next corollary.

\proclaim{Corollary 9.5} $\underline{\text{\bf Assumptions}}$ Let
$g_r\in \BC\{y,z\}$ be $\underline{\text{\rm a semi-quasi-Puiseux
series of the recursive }}$

\noindent$\underline{\text{\rm r-type}}$, as either in {\rm[A]} of
{\rm Definition 5.0.0} or in the assumption of Theorem $5.0$. If
$g_r$ is irreducible in $\BC\{y,z\}$, for each $r\ge 2$ $g_r$ can be
written in the form
$$\align
 g_r=(z^{n_1}+\ve_1
y^{\beta_{1,1}})^{d}+\sum_{\alpha,\beta\ge 0}
c^{(r)}_{\alpha\beta}y^{\alpha}z^{\beta} \quad \text{with
$d=n_2n_3\cdots n_r$,} \tag 9.5.1
\endalign$$
by Sublemma $5.2$ of $\S 5$ where $\ve_1$ is a unit in $\BC\{y,z\}$,
and the $c^{(r)}_{\alpha\beta}$ are nonzero complex numbers for some
nonnegative integers $\alpha$ and $\beta$ such that
$n_1\alpha+\beta_{1,1}\beta>n_1\beta_{1,1}n_2n_3\cdots n_r$. \ms

$\underline{\text{\bf {Conclusions}}}$ \quad Using the same
properties and notations in Proposition $9.4$, we have the
following:
$$\align
\text{\rm \{[Mult($V(g_r)$)]\}}
=\text{\rm{Join}}\{[dn_1:dk_1]; \text{\rm [Mult($V^{(m)}(g_r)$)]}\} \tag 9.5.2 \\
\endalign$$
where $k_1=\beta_{1,1}$, $\{[dn_1:dk_1]\}=\{d[n_1:k_1]\}$ and
$d={n_2n_3\cdots n_r}$. $\square$
\endproclaim

\vfill \pagebreak

{\bf Part[B4] In preparation for the proof of The 1st Algorithm}
\bs

{\bf \S10. To find the necessary and sufficient condition for any
two Puiseux convergent power series of recursive types in
$\BC\{y,z\}$ to have the same multiplicity sequence and their
classifications} \ms

{\bf \S10.0. Introduction}  \ms

In this section, the problem is to prove by Theorem 10.2 that we can
compute a one-to-one function from
$\underline{\text{\rm{Family(1)}}}$ into
$\underline{\text{\rm{Family(2)}}}$, using Theorem $7.7$, which
gives a solution of Problem[1-B] in $\S7$. Therefore, it suffices to
prove by Theorem $7.7$ and Theorem $10.2$ that if any two Puiseux
convergent power series of recursive types in $\BC\{y,z\}$ have the
same multiplicity sequence then they have the same divisor under two
standard resolutions in the sense of Definition 2.4, and conversely.
Note by Definition $7.0$ in $\S7$ that
$\underline{\text{\rm{Family(1)}}}$ is the subset of Quasi-Family(1)
where $\underline{\text{Quasi-Family(1)=\{\text{\rm $f$ is arbitrary
quasi-Puiseux convergent power}}}$ $\underline{\text{\text{\rm
series of the recursive type}: $f\in \text{\rm {Family(0)}}$\}}}$ in
the sense of Definition $5.0.0$. \ms

{\bf \S10.1. In preparation for computing a one-to-one function from
Family(1) into Family(2)} \ms

\proclaim{Theorem 10.1} $\underline{\text{\bf {Assumptions}}}$ By
the same way as in {\rm Theorem 7.3} or {\rm Definition 5.0.0},
define a quasi-Puiseux convergent power series $g_r$ of recursive
$r$-type in $\BC\{y,z\}$ by {\rm {Sequences[I]} in either the
assumptions of Theorem $7.3$ or Definition $5.0.0$. Suppose that the
same properties and notations as either in the assumptions of
Theorem $7.3$ or in Definition $5.0.0$ hold. Let $r$ be an arbitrary
positive integer. \ms

\noindent $\underline{\text{\bf Conclusions}}$ \quad By the same
notations in Definition $9.1$ and Definition $9.2$, \text{\rm
Multiseq($V(g_r)$)}, called the multiplicity sequence of $V(g_r)$,
can be represented as follows:
$$\align
(10.1.1) \quad & \text{\rm Multiseq($V(g_r)$)} =
\text{Join}(\{[n:\alpha_1]\},
\text\{[\gcd(n,\alpha_1):\alpha_2-\alpha_1]\},\dots, \qquad\\
 &
 \text\{[\gcd(n,\alpha_1,\dots,\alpha_{r-2}):\alpha_{r-1}-\alpha_{r-2}]\},
\text\{[\gcd(n,\alpha_1,\dots,\alpha_{r-1}):\alpha_{r}-\alpha_{r-1}]\}), \qquad \qquad \\
&  \quad\text{\rm such  that}  \qquad \qquad n=n_1n_2\cdots n_r,
\qquad \alpha_1 =\beta_{1,1}n_2\cdots n_r \quad \text{and}\\
& \quad \quad \qquad \qquad \qquad  \alpha_j =\alpha_{j-1}+
\widehat{\Delta}_j(\beta_{j,k})^j_{k=1}n_{j+1}n_{j+2}\cdots  n_r
\quad \text{for $2\le j\le r$},\\
 \endalign$$
where $\widehat{\Delta}_j(\beta_{j,k})^j_{k=1}
=\Delta_j(\beta_{j,k})^j_{k=1}-n_jn_{j-1}\Delta_{j-1}(\beta_{j-1,k})^{j-1}_{k=1}$
for $2\le j\le r$. $\square$
\endproclaim

\proclaim{Remark 10.1.1} {\rm (a)} \text{\rm Multiseq($V(g_r)$)} of
{\rm(10.1.1)} can be rewritten as follows:
$$\align
(10.1.1^*) \qquad & \text{\rm Multiseq($V(g_r)$)} = \text{\rm
Join}(\{[n:\alpha_1]\}, \text\{[d_1:\alpha_2-\alpha_1]\},\dots,
\text\{[d_{r-1}:\alpha_{r}-\alpha_{r-1}]\}),  \qquad\\
&  \quad\text{\rm such  that}  \qquad \qquad n=n_1d_1,
\qquad \alpha_1 =\beta_{1,1}d_1 \quad \text{\rm {and}} \\
& \quad \quad \qquad \qquad \qquad  \alpha_j-\alpha_{j-1}=
\widehat{\Delta}_j(\beta_{j,k})^j_{k=1}d_j
\quad \text{for $2\le j\le r$},\\
 \endalign$$
where $d_j=n_{j+1}n_{j+2}\cdots n_r$ for $1\le j\le r-1$ and
$d_r=1$, and $d_1=\gcd(n,\alpha_1)$ and
$d_j=\gcd(n,\alpha_1,\alpha_2-\alpha_1,\dots,\alpha_j-\alpha_{j-1})
=\gcd(d_{j-1},\alpha_j-\alpha_{j-1})$ for $2\le j\le r$ because
$d_j=\gcd(n,\alpha_1,\dots,\alpha_j)$. \ms

{\rm (b)} Suppose that $g_r$ is irreducible in $\BC\{y,z\}$. Then,
it was already proved by Theorem $5.0$ that the following are true:
$$ \split
(10.1.1.0) \qquad &\text{$g_r$ is irreducible in $\BC\{y,z\}$} \\
  \iff &\text{$g_1,g_2,\dots,g_{r-1}$ are irreducible in $\BC\{y,z\}$ and
   $\gcd(n_r,\Delta_r(\beta_{r,k})^r_{k=1})=1$} \\
\iff  &\text{$\gcd(n_1,\beta_{1,1})=1$,
$\gcd(n_2,\widehat{\Delta}_2(\beta_{2,1},\beta_{2,2}))=1$,
$\dots,\gcd(n_r,\widehat{\Delta}_r(\beta_{r,k})^r_{k=1})=1$}. \qquad
\qquad
\endsplit$$

For each $j=1,2,\dots,r$, note that $(0,0)$ is an isolated singular
point of an analytic variety $V(g_j)=\{(y,z):g_j(y,z)=0\}$ except
the case that $V(g_1)$ with $n_1=1$.  $\square$
\endproclaim \ms

{\bf \S10.2. The proof of Theorem 10.1 } \ms

Note that the assumptions of Theorem $10.1$ satisfies the same
assumptions and notations as in Theorem $5.0$. For the proof of
Theorem $10.1$, we can use five sublemmas, that is, Sublemma $5.1$,
Sublemma $5.2$, $\dots$, Sublemma $5.5$ of Theorem $5.0$,
respectively. \ms

\demo{\bf Proof of Theorem 10.1} Now, the proof of the equality in
{\rm (10.1.1)} will be by induction on the positive integer $r\ge
1$. Note that $n_1\ge 2$ and $\beta_{1,1}\ge 1$ with
$\gcd(n_1,\beta_{1,1})=1$.

Then, it is enough to consider two cases, respectively.

Case(I) $r=1$, and Case(II) $r\ge 2$. \ms

$\underline{\text{\bf Case(I)}}$ \ Let $r=1$. To prove the equality
in {\rm (10.1.1)}, if $n_1\ge 2$ and $\beta_{1,1}=1$, then there is
nothing to prove. For the proof of theorem, if $r=1$ then it may be
assumed that $n_1\ge 2$ and $\beta_{1,1}>1$ with
$\gcd(n_1,\beta_{1,1})=1$. So, it suffices to prove the following:
$$\align
\text{\rm \{[Multiseq$(V(g_1))$]\}}
&=\{[n:\alpha_1]\}, \tag 10.1.2 \\
\endalign$$
where $n=n_1$ and $\alpha_1=\beta_{1,1}$. Then, there is nothing to
prove by Proposition $9.4$. \ms

$\underline{\text{\bf Case(II)}}$ \ Let $r\ge 2$. To prove the
equality in {\rm (10.1.1)}, suppose we have shown on the positive
integer $(r-1)$ by the induction method that for a given
$V(g_{r-1})=\{(y,z):g_{r-1}(y,z)=0\}$, \text{\rm
Multiseq($V(g_{r-1})$)=\{[Mult($V(g_{r-1})$)]\}}, called the
multiplicity sequence of $V(g_{r-1})$, can be represented as
follows:
$$\align
(10.1.3) \quad & \text{\rm Multiseq($V(g_{r-1})$)} =
\text{Join}(\{[n:\alpha_1]\},
\text\{[\gcd(n,\alpha_1):\alpha_2-\alpha_1]\},\dots, \qquad\\
&
\text\{[\gcd(n,\alpha_1,\dots,\alpha_{r-3}):\alpha_{r-2}-\alpha_{r-3}]\},
\text\{[\gcd(n,\alpha_1,\dots,\alpha_{r-2}):\alpha_{r-1}-\alpha_{r-2}]\}), \qquad \qquad \\
 &\text{such ~ that} \quad n =n_1n_2\cdots n_{r-1},
\quad \alpha_1 =\beta_{1,1}n_2\cdots n_{r-1} \quad \text{and}\\
 & \quad \quad \qquad   \alpha_j =\alpha_{j-1}+
\widehat{\Delta}_j(\beta_{j,k})^j_{k=1}n_{j+1}n_{j+2}\cdots  n_{r-1}
\quad \text{for $2\le j\le r-1$,} \\
  \endalign$$
where $\widehat{\Delta}_j(\beta_{j,k})^j_{k=1}
=\Delta_j(\beta_{j,k})^j_{k=1}-n_jn_{j-1}\Delta_{j-1}(\beta_{j-1,k})^{j-1}_{k=1}$
for \text{$2\le j\le r-1$} and \text{$\Delta_1(t)=t$}. \ms

In preparation for the proof of the equality in {\rm (10.1.1)}, it
was already proved by Proposition $9.4$(Corollary $9.5$) or Sublemma
$5.4$ that the following equality is true:
$$\align
&  \text{\{[Mult$(V(g_r))$]\}} =\text{Join}(\{[n:\alpha_1]\},
\text{\{[Mult$(V^{(\lambda_1)}(g_r))$]}\}), \tag 10.1.4 \\
& where \quad n =n_1n_2\cdots n_r \quad \text{and} \quad
 \alpha_1 =\beta_{1,1}n_2\cdots n_r. \\
\endalign$$

So, for the proof of the equality in (10.1.1), it suffices to show
by (10.1.3) and (10.1.4) that the following equality is true:
$$\align
(10.1.5)\quad & \{\text{[Mult$(V^{(\lambda_1)}(g_r))$]}\} =
\text{Join}(\{[\gcd(n,\alpha_1):\alpha_2-\alpha_1]\},
\{[\gcd(n,\alpha_1,\alpha_2):\alpha_3-\alpha_2]\},  \\
&
\dots,\{[\gcd(n,\alpha_1,\dots,\alpha_{r-2}):\alpha_{r-1}-\alpha_{r-2}]\},
\{[\gcd(n,\alpha_1,\dots,\alpha_{r-1}):\alpha_{r}-\alpha_{r-1}]\}). \\
\endalign$$

For the proof of the equality in (10.1.5), it remains to compute
$\{\text{[Mult$(V^{(\lambda_1)}(g_r))$]}\}$.

First of all, note by Sublemma $5.5$ that
$V^{(\lambda_1)}(g_r)=V(h_{r-1})$ is well-defined where
$h_{r-1}=(g_r\circ\tau_m)_{proper}
=h^{s_{r-1}}_{r-2}+\eta_{r-1}v^{\gamma_{r-1,1}}(u+1)^{\gamma_{r-1,2}}
h^{\gamma_{r-1,3}}_1\cdots h^{\gamma_{r-1,r-1}}_{r-3}$ with
$\lambda_1=m$. That is, it was already shown by Sublemma $5.5$ that
the local defining equation $h_{r-1}\in \BC\{v,u+1\}$ of
$V(h_{r-1})$ satisfies the same kind of assumptions and notations as
the local defining equation of $V(g_{r-1})$ in {\rm (10.1.3)} does,
which is a necessary and sufficient condition for the use of the
induction proof on the positive integer \text{\rm (r-1)}.\ms

So, by the induction assumption on the positive integer \text{\rm
(r-1)}, we can apply the same results of {\rm (10.1.3)} to
$V(h_{r-1})$, up to the same kind of notations and properties below.
\ms

Therefore, we can prove by (10.1.3) and by the same kind of
notations and properties as in Sublemma $5.5$ that for a given
$V(h_{r-1})=\{(u+1,v):h_{r-1}(u+1,v)=0\}$, \text{\rm
Multiseq($V(h_{r-1})$)}, called the multiplicity sequence of
$V(h_{r-1})$, is represented as follows:
$$\align
(10.1.6) \qquad\quad &  \{\text{[Mult$(V^{(\lambda_1)}(g_r))$]}\}
= \{\text{[Mult$(V(h_{r-1}))$]}\}   \\
\quad =&\text{Join}(\{[b:\delta_1];
[\gcd(b,\delta_1):\delta_2-\delta_1];
[\gcd(b,\delta_1,\delta_2):\delta_3-\delta_2];  \cdots; \\
 \qquad &
[\gcd(b,\delta_1,\delta_2,\dots,\delta_{r-3}):\delta_{r-2}-\delta_{r-3}];
[\gcd(b,\delta_1,\delta_2,\dots,\delta_{r-2}):\delta_{r-1}-\delta_{r-2}]\}),
\qquad \qquad \\
 &\text{such ~ that} \qquad b=s_1 s_2\cdots s_{r-1},
 \quad \delta_1 =\gamma_{1,1}s_2\cdots s_{r-1} \quad \text{and} \quad \\
& \qquad \qquad \qquad \delta_j =\delta_{j-1}+
\widehat{\Xi}_j(\gamma_{j,k})^j_{k=1}s_{j+1}s_{j+2}\cdots s_{r-1}
\quad \text{for $2\le j\le r-1$}, \\
 \endalign$$
where $\widehat{\Xi}_j(\gamma_{j,k})^j_{k=1}
=\Xi_j(\gamma_{j,k})^j_{k=1}-n_jn_{j-1}\Xi_{j-1}(\gamma_{j-1,k})^{j-1}_{k=1}$
for $2\le j\le r-1$ and $\Xi_1(t)=t$.

Here, note by Sublemma $5.5$ that
$$\align
(10.1.7) \qquad \quad s_1 &=n_2\ge 2,
 ~\gamma_{1,1}=\Delta^{\sharp}_2(\beta_{2,1},\beta_{2,2})-n_1\beta_{1,1}n_2>0,  \\
  s_{j-1} &=n_j\ge 2,
 ~\gamma_{j-1,1}=\Delta^{\sharp}_j(\beta_{j,k})^j_{k=1}
 -n_1\beta_{1,1}n_2n_3\cdots n_j>0, \\
 &\hskip 2 true cm  \gamma_{j-1,2}=\beta_{j,3},
\gamma_{j-1,3}=\beta_{j,4} , \dots ,
 \gamma_{j-1,j-1}=\beta_{j,j} \quad \text{for $2\le j\le r$,} \qquad \qquad \\
\endalign
$$
where $\gamma_{1,1},\gamma_{2,1},\dots,\gamma_{r-1,1}$ are positive
by Sublemma $5.1$, noting that
$\gamma_{1,1}=\Delta^{\sharp}_2(\beta_{2,1},\beta_{2,2})-n_1\beta_{1,1}n_2
=\Delta_2(\beta_{2,1},\beta_{2,2})-n_1\beta_{1,1}n_2
=\widehat{\Delta}_2(\beta_{2,1},\beta_{2,2})$. \ms

Moreover, as we have seen in a conclusion of Sublemma $5.5$, we have
the following representation, too: Let $q=2,3,\dots,r-1$.
$$\align
(10.1.8) \quad  &\widehat{\Xi}_q(\gamma_{q,k})^q_{k=1}
=\Xi_q(\gamma_{q,k})^q_{k=1}-s_qs_{q-1}\Xi_{q-1}(\gamma_{q-1,k})^{q-1}_{k=1}
\\
 =&\Delta_{q+1}(\beta_{q+1,k})^{q+1}_{k=1}
-n_{q+1}n_q\Delta_q(\beta_{q,k})^q_{k=1}
=\widehat{\Delta}_{q+1}(\beta_{q+1,k})^{q+1}_{k=1}>0. \quad
\text{(Sublemma $5.5$)} \qquad \qquad\\
\endalign$$

Now, for the complete proof of this theorem, comparing (10.1.5) with
(10.1.6), it suffices to show the following equations are true:
$$\align
(10.1.9)\qquad & \text{{\rm (i)}
$[\gcd(n,\alpha_1):\alpha_2-\alpha_1]=[b:\delta_1]$.} \\
& \text{{\rm (ii)}
$[\gcd(n,\alpha_1,\alpha_2,\dots,\alpha_{j-1}):\alpha_{j}-\alpha_{j-1}]=
[\gcd(b,\delta_1,\delta_2,\dots,\delta_{j-2}):\delta_{j-1}-\delta_{j-2}]$} \qquad \\
& \qquad   \text{for each $j=3,4,\dots,r$}.
\endalign$$

Recall by (10.1.1) that
$$\align
(10.1.10) \quad \qquad\qquad n &=n_1n_2\cdots n_r, \quad
 \alpha_1 =\beta_{1,1}n_2\cdots n_r \quad \text{and} \\
   \alpha_j &=\alpha_{j-1}+
\widehat{\Delta}_j(\beta_{j,k})^j_{k=1}n_{j+1}n_{j+2}\cdots
  n_r \quad \text{for $2\le j\le r$.} \qquad\qquad \qquad\qquad\\
 \endalign$$

In preparation for the  proof of the equality in (10.1.9), first of
all, we can get easily the following equations from (10.1.6),
(10.1.7), (10.1.8) and (10.1.10):
$$\align
(10.1.11) \quad \qquad \qquad b &=s_1 s_2\cdots s_{r-1}
=n_2n_3\cdots n_r=\gcd(n,\alpha_1),  \quad \qquad \qquad\\
 \delta_1 &=\gamma_{1,1}s_2\cdots s_{r-1}
=\widehat{\Delta}_2(\beta_{2,1},\beta_{2,2})n_3n_4\cdots
n_r=\alpha_2-\alpha_1>0,
\quad \qquad \qquad\\
\delta_j -\delta_{j-1} &=
\widehat{\Xi}_j(\gamma_{j,k})^j_{k=1}s_{j+1}s_{j+2}\cdots s_{r-1}
\quad \text{for $2\le j\le r-1$} \\
&=\widehat{\Delta}_{j+1}(\beta_{j+1,k})^{j+1}_{k=1}n_{j+2}n_{j+3}\cdots  n_r
=\alpha_{j+1}-\alpha_j>0, \\
\endalign$$
noting by assumption that  $\gcd(n,\alpha_1) =\gcd(n_1n_2\cdots
n_r,\beta_{1,1}n_2\cdots n_r)=n_2n_3\cdots n_r=b$ because
$\gcd(n_1,\beta_{1,1})=1$.

In order to finish the proof of the equality of (10.1.9), compare
(10.1.9) with (10.1.11), and then it remains to prove the following:
$$\align
(10.1.12) \quad  &\text{{\rm (i)} \quad $\gcd(n,\alpha_1)=b$.} \\
& \text{{\rm (ii)} \quad
$\gcd(n,\alpha_1,\alpha_2,\dots,\alpha_{j})=
\gcd(b,\delta_1,\delta_2,\dots,\delta_{j-1})$ for each
$j=2,3,\dots,r-1.$} \qquad \qquad
\qquad \\
\endalign$$

Using (10.1.11) again, the proof of (10.1.12) is as follows:

(a) Since it is clear by (10.1.11) that $\gcd(n,\alpha_1)=b$, then
there is nothing for the proof of (i) of (10.1.12).

(b) If $j=2$, then it is clear by (a) and (10.1.11) that
$\gcd(n,\alpha_1,\alpha_2)=\gcd(n,\alpha_1,\alpha_2-\alpha_1)=\gcd(b,\delta_1)$.
So, if $j=2$, then the proof of (ii) of (10.1.12) is done.

The general case in (ii) of (10.1.12) will be proved by induction.
Suppose we have shown on the positive integer $j<r-1$ that
$\gcd(n,\alpha_1,\alpha_2,\dots,\alpha_{j})=
\gcd(b,\delta_1,\delta_2,\dots,\delta_{j-1})$. Then,
$\gcd(n,\alpha_1,\alpha_2,\dots,\alpha_{j+1})
=\gcd(n,\alpha_1,\alpha_2-\alpha_1,\alpha_3-\alpha_2,\dots,\alpha_{j+1}-\alpha_j)
=\gcd(b,\delta_1,\delta_2-\delta_1,\dots,\delta_j-\delta_{j-1})
=\gcd(b,\delta_1,\delta_2,\dots,\delta_j)$ by (10.1.11) and (a).
Thus, the proof of (ii) of (10.1.12) is done, and so the proof of
(10.1.12) is finished.

Therefore, the proof of the theorem can be completely finished.
$\square$
\enddemo \ms

{\bf \S10.3. How to compute a 1-1 function from Family(1) into
Family(2)} \ms

\proclaim{Theorem 10.2} $\underline{\text{\bf {Assumptions}}}$ Let
$r$ and $\rho$ be arbitrary positive integers. By the same way as in
the assumption of Theorem $7.3$, define arbitrary quasi-Puiseux
series $g_r\in \BC\{y,z\}$ of the recursive $r$-type by {\rm
{Sequences[I]} and arbitrary quasi-Puiseux series $\phi_{\rho}\in
\BC\{y,z\}$ of the recursive $\rho$-type  by {\rm {Sequences[II]} in
Theorem $7.3$.

We may assume that the same properties and notations as in the
assumptions of Theorem $7.3$ hold. Note that {\rm {Sequences[I]} and
{\rm Sequences[II]} are the same up to the change of notations.
Assume in addition that $2\le n_1<\beta_{1,1}$ and $2\le
\ell_1<\delta_{1,1}$.  \ms

\noindent $\underline{\text{\bf Conclusions}}$ \quad Note by either
Definition 1.1 or Definition 5.0.0 that $g_r$ and $\phi_\rho$ are
called the Puiseux series in $\BC\{y,z\} $ because  $2\le
n_1<\beta_{1,1}$ and $2\le \ell_1<\delta_{1,1}$ by assumptions.
Then, we have the following:
$$\align
  & \text{$g_r$ and $\phi_{\rho}$ have the same multiplicity
sequence.} \tag 10.2.1 \\
  \iff \quad & \\
(10.2.2)\qquad \qquad \qquad \qquad & \text{$n_j=\ell_j$ and
$\Delta_j(\beta_{j,k})^j_{k=1}=\omega_j(\delta_{j,k})^j_{k=1}$ for
each $j=1,2,\dots,r=\rho$,}   \\
\endalign$$

Moreover, it can be easily proved by Theorem $7.3$ that the
following holds:
$$\align
& \text{$g_r$ and $\phi_{\rho}$ have the same multiplicity
sequence} \tag 10.2.3 \\
 \iff \quad & \\
& \text{$g_r \buildrel \text{{\rm divisor}} \over \sim \phi_\rho$
under the standard resolutions. \quad $\square$} \tag 10.2.4
\endalign$$
\endproclaim

\definition{Remark 10.2.1} Without assuming that both $2\le n_1<\beta_{1,1}$
and $2\le \ell_1<\delta_{1,1}$, it can be easily proved that the
conclusion of the theorem may not be true by the following example:
$$\align
& g_1=z^3+y^8 \qquad \text{and} \tag 10.2.1.1\\
& \phi_2=\phi^3_1+y^2z^3 \quad \text{with \quad $\phi_1=z+y^2$,}
\endalign$$
because $g_1$ and $\phi_{2}$ have the same multiplicity sequence,
and also they have the same divisor under two standard resolutions,
but the condition in (10.2.2) does not hold.
\enddefinition \ms

{\bf \S10.4. The proof of Theorem 10.2 } \ms

\demo{\bf Proof of Theorem 10.2} First of all, the first half for
the proof of Theorem 10.2 is to show that two statements in (10.2.1)
and (10.2.2) are equivalent, and then the second half for the proof
of Theorem 10.2 is to show that two statements in (10.2.3) and
(10.2.4) are equivalent. After the proof of the first half is done,
there is nothing to prove for the second half, because it was
already proved by (7.3.4) of Theorem $7.3$ that the condition in
(10.2.2) is necessary and sufficient for \text{$g_r \buildrel
\text{{\rm divisor}} \over \sim \phi_\rho$} under the standard
resolutions.

Now, for the proof of the first half, it suffices to consider two
cases:

{\rm Fact[I].} We prove the sufficiency of the condition in {\rm
(10.2.2)} for $V(g_r)$ and $V(\phi_\rho)$ to have the same
multiplicity sequence at $(0,0)\in \BC^2$.

{\rm Fact[II].} We prove the necessity of the condition in {\rm
(10.2.2)} for $V(g_r)$ and $V(\phi_\rho)$ to have the same
multiplicity sequence at $(0,0)\in \BC^2$. \ms

In preparation for solving the first half, since each of
Sequences[I] and Sequences[II] in Theorem $10.2$ satisfies the same
kind of assumption as in Sequences[I] of Theorem $10.1$, then we can
apply Theorem $10.1$ to each of Sequences[I] and Sequences[II] in
Theorem $10.2$, respectively. So, we can easily get the following
sublemma, which will be applicable for the proofs of Fact(I) and
Fact(II).

\proclaim{Sublemma 10.2.2} $\underline{\text{\bf Assumptions}}$
Suppose that the same properties and notations as in the assumption
of Theorem $10.2$ are true. Note that $g_r$ and $\phi_\rho$ are
irreducible in $\BC\{y,z\}$. \ms

\noindent $\underline{\text{\bf Conclusions}}$ Then, the
multiplicity sequences of $V(g_r)$ and $V(\phi_\rho)$ are as
follows:
$$\align
(10.2.2.1) \quad & \text{\rm{Multiseq}$(V(g_r))$} =
\text{\rm{Join}}(\{[n:\alpha_1]\},
\text\{[\gcd(n,\alpha_1):\alpha_2-\alpha_1]\},\dots, \qquad\\
 &
 \text\{[\gcd(n,\alpha_1,\dots,\alpha_{r-2}):\alpha_{r-1}-\alpha_{r-2}]\},
\text\{[\gcd(n,\alpha_1,\dots,\alpha_{r-1}):\alpha_{r}-\alpha_{r-1}]\}), \qquad \qquad \\
&  \quad\text{\rm such  that}  \qquad \quad n=n_1n_2\cdots n_r,
\qquad \alpha_1 =\beta_{1,1}n_2\cdots n_r \quad \text{and}\\
& \quad \quad \quad \qquad \qquad  \alpha_j =\alpha_{j-1}+
\widehat{\Delta}_j(\beta_{j,k})^j_{k=1}n_{j+1}n_{j+2}\cdots  n_r
\quad \text{for $2\le j\le r-1$},\\
&  \quad \quad \quad \qquad \qquad \alpha_r =\alpha_{r-1}+
\widehat{\Delta}_r(\beta_{r,k})^r_{k=1},
 \endalign$$
where $\widehat{\Delta}_j(\beta_{j,k})^j_{k=1}
=\Delta_j(\beta_{j,k})^j_{k=1}-n_jn_{j-1}\Delta_{j-1}(\beta_{j-1,k})^{j-1}_{k=1}$
for $2\le j\le r$ and $\Delta_1(t)=t$.
$$\align
\text{\rm(10.2.2.2)}\quad & \text{\rm{Multiseq}$(V(\phi_{\rho}))$} =
\text{\rm{Join}}(\{[b:\chi_1]\},
\text\{[\gcd(b,\chi_1):\chi_2-\chi_1]\},\dots, \qquad\\
&
\text\{[\gcd(b,\chi_1,\chi_2,\dots,\chi_{\rho-2}):\chi_{\rho-1}-\chi_{\rho-2}]\},
\text\{[\gcd(n,\chi_1,\chi_2,\dots,\chi_{\rho-1}):\chi_{\rho}-\chi_{\rho-1}]\}),  \\
&  \text{such ~ that} \quad \qquad  b ={\ell_1} {\ell_2}\cdots
{\ell_{\rho}},\qquad
\chi_1 =\delta_{1,1}{\ell_2} {\ell_3}\cdots {\ell_{\rho}} \quad \text{and}\\
& \qquad \qquad \qquad \chi_{j} =\chi_{j-1}
+\widehat{\omega}_{j}(\delta_{j,k})^{j}_{k=1}{\ell_{j+1}}{\ell_{j+2}}\cdots{\ell_{\rho}}
  \quad \text{for $2\le j\le {\rho}-1$}, \\
& \qquad \qquad \qquad  \chi_{\rho} =\chi_{\rho-1}+
\widehat{\omega}_{\rho}(\delta_{\rho,k})^{\rho}_{k=1},
 \endalign$$
where $\widehat{\omega}_j(\delta_{j,k})^j_{k=1}
=\omega_j(\delta_{j,k})^j_{k=1}-{\ell_j}{\ell_{j-1}}
\omega_{j-1}(\delta_{j-1,k})^{j-1}_{k=1}$ for $2\le j\le \rho$ and
$\omega_1(t)=t$. $\square$
\endproclaim

The proof of this sublemma just follows from Theorem $10.1$. \ms

${\underline{\text{\bf{Fact[I]}}}}$ For the proof of the sufficiency
of the condition in (10.2.2), suppose that the following are true:
$$\align
(10.2.5)\qquad \qquad   \text{$n_j=\ell_j$ and
$\Delta_j(\beta_{j,k})^j_{k=1}=\omega_j(\delta_{j,k})^j_{k=1}$ for
each $j=1,2,\dots,r=\rho.$} \qquad \qquad \qquad
\endalign$$

Since $n_j=\ell_j$ for each $j=1,2,\dots,r=\rho$ and
$\beta_{1,1}=\Delta_1(\beta_{1,1})=\omega_1(\delta_{1,1})=\delta_{1,1}$
by (10.2.5), then
\text{$\widehat{\Delta}_j(\beta_{j,k})^j_{k=1}=\widehat{\omega}_j(\delta_{j,k})^j_{k=1}$
for each $j=1,2,\dots,r$} by (10.2.5), again.

So, it is clear by (10.2.5) and Sublemma $10.2.2$ that the
following equalities are true:
$$\align
 & \text{$n=b$ and $\alpha_1=\chi_1$,\quad and}
\tag 10.2.6 \\
&\text{$\alpha_j-\alpha_{j-1}=\chi_j-\chi_{j-1}$ for each
$j=2,\dots,r=\rho$.} \\
\endalign$$

Noting by (10.2.6) that $\alpha_j=\chi_j$ for each
$j=1,2,\dots,r=\rho$, then it is clear by (10.2.2.1) and
(10.2.2.2) of Sublemma $10.2.2$ and by (10.2.6) that $g_r$ and
$\phi_\rho$ have the same multiplicity sequence. \ms

${\underline{\text{\bf{Fact[II]}}}}$ To prove the necessity of the
condition in (10.2.2), suppose that $g_r$ and $\phi_\rho$ have the
same multiplicity sequence and by assumption that $2\le
n_1<\beta_{1,1}$ and $2\le \ell_1<\delta_{1,1}$. For the proof, we
may assume that $1\le r\le \rho$.

Case(I) $r=1\le \rho$, and Case(II) $2\le r\le \rho$. \ms

${\underline{\text{\bf{Case(I) of Fact[II]}}}}$ Let $r=1\le \rho$.
By the definition of $\text{\rm{Multiseq}$(V(g_r))$}$ and
\text{\rm{Multiseq}$(V(\phi_{\rho}))$}, observe (i) and (ii):

(i) The largest element of $\{\text{\rm [Mult$(V(g_1))$]}\}$ is
$n_1$, and the largest element of
$\{\text{\rm[Mult$(V(\phi_{\rho}))$]}\}$ is $b ={\ell_1}
{\ell_2}\cdots {\ell_{\rho}}$, too. So, $n_1=b$.

(ii) Let $s$ be a positive integer such that
$sn_1<\beta_{1,1}<(s+1)n_1$. Then, we can define the second largest
element $K_2$ of $\{\text{[Mult$(V(g_1))$]}\}$ with $n_1>K_2$. So,
it is clear that $K_2=\beta_{1,1}-sn_1$, because
$K_2=\beta_{1,1}-sn_1$ is the (s+1)-th element of
$\{\text{[Mult$(V(g_1))$]}\}$, viewed as a finite sequence.

Applying the same method to (10.2.2.2) of Sublemma $10.2.2$ just as
above, let $s'$ be a positive integer such that
$s'b<\chi_1<(s'+1)b$. Then, we can define the second largest element
$L_2$ of $\{\text{[Mult$(V(\phi_{\rho}))$]}\}$ with $b>L_2$, noting
that $L_2=\chi_1-s'b$ is the $(s'+1)$-th element of
$\{\text{[Mult$(V(\phi_{\rho}))$]}\}$, as a finite sequence.

Since any multiplicity sequence is monotonically decreasing, then it
is trivial by (i) and (ii) that  $s=s'$ and
$\beta_{1,1}-sn_1=K_2=L_2=\chi_1-s'b$ by the definition of the
multiplicity sequence for two positive integers in Definition $9.1$.
So, $\beta_{1,1}=\chi_1=\delta_{1,1}\ell_2\ell_3\cdots \ell_\rho$
because $n_1=b$. So,
$\frac{\beta_{1,1}}{n_1}=\frac{\chi_1}{b}=\frac{\delta_{1,1}}{\ell_1}$
by (10.2.2.2) of Sublemma $10.2.2$. Since $\gcd(n_1,\beta_{1,1})=1$
and $\gcd(\ell_1,\delta_{1,1})=1$, then $n_1=\ell_1$ and
$\beta_{1,1}=\delta_{1,1}$. Since $n_1=b$, then
$\delta_{1,1}=\chi_1$ with $\rho=1$. Thus, the proof for Case(I) is
done. \ms

${\underline{\text{\bf{Case(II) of Fact[II]}}}}$ Let $2\le r\le
\rho$. For the proof for Case(II), we prove first that the assertion
in (10.2.7) is true, and after then, it suffices to show by (10.2.7)
that the assertion in (10.2.8) is true:
$$\align
(10.2.7) \qquad & \text{$\alpha_{j}=\chi_{j}$
\quad for each $j=1,2,\dots,r\le \rho$,} \\
& \text{and so $\gcd(n,\alpha_1,\alpha_2,\dots,\alpha_{j})
=\gcd(b,\chi_1,\chi_2,\dots,\chi_{j})$,}  \\
&\text{which is the smallest element of
$\{[\gcd(n,\alpha_1,\alpha_2,\dots,\alpha_{j-1}):\alpha_{j}-\alpha_{j-1}]\}$.}\\
(10.2.8) \qquad & \text{$n_j=\ell_j$ and
$\Delta_j(\beta_{j,k})^j_{k=1}=\omega_j(\delta_{j,k})^j_{k=1}$ \quad
for each $j=1,2,\dots,r=\rho$.} \qquad \qquad
\endalign$$

$\underline{\text{\rm The assertion in (10.2.7)}}$ The assertion in
(10.2.7) will be proved by induction on the positive integer $r$
where $1\le r\le \rho$, using the following steps:

Step(1). \quad $\alpha_1=\chi_1$.

Step(j). \quad $\alpha_j=\chi_j$ \quad for each $j=2,3,\dots,r$. \ms

$\underline{\text{\rm Step(1) for the proof of the assertion in
(10.2.7)}}$ As in Case(I), note from the definition of
$\{\text{[Mult$(V(g_r))$]}\}$ and
$\{\text{[Mult$(V(\phi_{\rho}))$]}\}$ that the largest element of
$\{\text{[Mult$(V(g_r))$]}\}$ is $n=n_1n_2\cdots n_r$ and the
largest element of $\{\text{[Mult$(V(\phi_{\rho}))$]}\}$ is $b
={\ell_1} {\ell_2}\cdots {\ell_{\rho}}$. So, $n=b$ by the necessity
of the condition in (10.2.2).

Let $s$ be a positive integer such that $sn<\alpha_1<(s+1)n$ because
$2\le n_1<\beta_{1,1}$ and $\gcd(n_1,\beta_{1,1})=1$ and
$\dfrac{\alpha_1}{n}=\dfrac{\beta_{1,1}}{n_1}$. Then, we can define
the second largest integer $K_2$ of $\{\text{[Mult$(V(g_r))$]}\}$
with $n>K_2$. So, it is clear that $K_2=\alpha_1-sn$, noting that
$K_2=\alpha_1-sn$ is the (s+1)-th element of
$\{\text{[Mult$(V(g_r))$]}\}$, as a sequence.

By the same way as above, let $s'$ be a positive integer such that
$s'b<\chi_1<(s'+1)b$. Then, we can define the second largest integer
$L_2$ of $\{\text{[Mult$(V(\phi_{\rho}))$]}\}$ with $b>L_2$, noting
that $L_2=\chi_1-s'b$ is the $(s'+1)$-th element of
$\{\text{[Mult$(V(\phi_{\rho}))$]}\}$, as a sequence. Since any
multiplicity sequence is monotonically decreasing, then it is
trivial that $s+1=s'+1$ and $\alpha_1-sn=K_2=L_2=\chi_1-s'b$ by the
definition of the multiplicity sequence for two positive integers in
Definition $9.1$, and therefore $\alpha_1=\chi_1$. Thus, the proof
of Step(1) is done. \ms

$\underline{\text{\rm Step(j+1) for the proof of the assertion in
(10.2.7)}}$ By induction assumption on the positive integer $r$,
suppose we have shown by Step(1), Step(2),\dots,Step(j) that
$\alpha_i=\chi_i$ for $i=2,3,\dots,j$ where $j$ is an arbitrary
integer such that $2\le j<r\le \rho$. Then, it remains to prove that
$\alpha_{j+1}=\chi_{j+1}$. Since
$\gcd(n,\alpha_1,\alpha_2,\dots,\alpha_{j})>\gcd(n,\alpha_1,\alpha_2,\dots,\alpha_{j+1})$
by assumption, and then
$\gcd(n,\alpha_1,\alpha_2,\dots,\alpha_{j})\not
=\alpha_{j+1}-\alpha_{j}$, for the proof, it suffices to consider
two subcases, respectively.

Subcase(A)
$\gcd(n,\alpha_1,\alpha_2,\dots,\alpha_{j})<\alpha_{j+1}-\alpha_{j}$,
and

Subcase(B)
$\gcd(n,\alpha_1,\alpha_2,\dots,\alpha_{j})>\alpha_{j+1}-\alpha_{j}$.
\ms

Subcase(A): Let
$\gcd(n,\alpha_1,\alpha_2,\dots,\alpha_{j})<\alpha_{j+1}-\alpha_{j}$.
Note that the largest element of
$\{[\gcd(n,\alpha_1,\alpha_2,\dots,\alpha_{j}):\alpha_{j+1}-\alpha_{j}]\}$
is $\gcd(n,\alpha_1,\alpha_2,\dots,\alpha_{j})$. Since
$\gcd(b,\chi_1,\chi_2,\dots,\chi_{j})>\gcd(b,\chi_1,\chi_2,\dots,\chi_{j+1})$
by assumption, then $\gcd(b,\chi_1,\chi_2,\dots,\chi_{j})\not
=\chi_{j+1}-\chi_{j}$. Then, for this case, it is enough to
consider two possibilities, respectively:

(A1) $\gcd(b,\chi_1,\chi_2,\dots,\chi_{j})<\chi_{j+1}-\chi_{j}$, and
(A2) $\gcd(b,\chi_1,\chi_2,\dots,\chi_{j})>\chi_{j+1}-\chi_{j}$. \ms

(A1) Let $\gcd(b,\chi_1,\chi_2,\dots,\chi_{j})<\chi_{j+1}-\chi_{j}$.
By the definition of the multiplicity sequence for the greatest
common divisor of two positive integers, the largest element of
$\{[\gcd(b,\chi_1,\dots,\chi_{j}):\chi_{j+1}-\chi_{j}]\}$ is
$\gcd(b,\chi_1,\dots,\chi_{j})$. Since
$\gcd(n,\alpha_1,\dots,\alpha_{j})<\alpha_{j+1}-\alpha_{j}$ and
$\gcd(b,\chi_1,\dots,\chi_{j})<\chi_{j+1}-\chi_{j}$, then
$\gcd(n,\alpha_1,\dots,\alpha_{j})=\gcd(b,\chi_1,\dots,\chi_{j})$
implies that $\alpha_{j+1}-\alpha_j=\chi_{j+1}-\chi_j$ by the same
techniques as we have used for the proof of Case(I). So,
$\alpha_{j+1}=\chi_{j+1}$.

(A2) Let $\gcd(b,\chi_1,\chi_2,\dots,\chi_{j})>\chi_{j+1}-\chi_{j}$.
By the definition of the multiplicity sequence for the greatest
common divisor of two positive integers, the largest element of
$\{[\gcd(b,\chi_1,\chi_2,\dots,\chi_{j}):\chi_{j+1}-\chi_{j}]\}$ is
$\chi_{j+1}-\chi_j$. So,
$\gcd(n,\alpha_1,\alpha_2,\dots,\alpha_{j})=\chi_{j+1}-\chi_j
<\gcd(b,\chi_1,\chi_2,\dots,\chi_{j})=\gcd(n,\alpha_1,\alpha_2,\dots,\alpha_{j})$,
which would be impossible. \ms

Subcase(B): Let
$\gcd(n,\alpha_1,\alpha_2,\dots,\alpha_{j})>\alpha_{j+1}-\alpha_{j}$.
By the definition of the multiplicity sequence for the greatest
common divisor of two positive integers, note that the largest
element of
$\{[\gcd(n,\alpha_1,\dots,\alpha_{j}):\alpha_{j+1}-\alpha_{j}]\}$ is
$\alpha_{j+1}-\alpha_j$. Since
$\gcd(b,\chi_1,\dots,\chi_{j})>\gcd(b,\chi_1,\chi_2,\dots,\chi_{j+1})$
by assumption, then $\gcd(b,\chi_1,\chi_2,\dots,\chi_{j})\not
=\chi_{j+1}-\chi_{j}$. Then, for this case, it is enough to consider
two possibilities, respectively:

(B1) $\gcd(b,\chi_1,\chi_2,\dots,\chi_{j})>\chi_{j+1}-\chi_{j}$, and
(B2) $\gcd(b,\chi_1,\chi_2,\dots,\chi_{j})<\chi_{j+1}-\chi_{j}$. \ms

\noindent(B1) Let
$\gcd(b,\chi_1,\chi_2,\dots,\chi_{j})>\chi_{j+1}-\chi_{j}$. By the
definition of the multiplicity sequence for the greatest common
divisor of two positive integers, the largest element of
$\{[\gcd(b,\chi_1,\chi_2,\dots,\chi_{j}):\chi_{j+1}-\chi_{j}]\}$ is
$\chi_{j+1}-\chi_j$. So, $\chi_{j+1}-\chi_j=\alpha_{j+1}-\alpha_j$.
So, $\alpha_{j+1}=\chi_{j+1}$ because $\alpha_{j}=\chi_{j}$ by the
induction assumption. \ms

\noindent(B2) Let
$\gcd(b,\chi_1,\chi_2,\dots,\chi_{j})<\chi_{j+1}-\chi_{j}$. By
definition of the multiplicity sequence for the greatest common
divisor of two positive integers, the largest element of
$\{[\gcd(b,\chi_1,\chi_2,\dots,\chi_{j}):\chi_{j+1}-\chi_{j}]\}$ is
$\gcd(b,\chi_1,\chi_2,\dots,\chi_{j})$. So,
$\alpha_{j+1}-\alpha_j=\gcd(b,\chi_1,\dots,\chi_{j})
=\gcd(n,\alpha_1,\dots,\alpha_{j})$, which would be impossible
because $\gcd(n,\alpha_1,\dots,\alpha_j)=n_{j+1}n_{j+2}\cdots
n_r>n_{j+2}n_{j+3}\cdots n_r=\gcd(n,\alpha_1,\dots,\alpha_{j+1})$.

By summarizing two subcases, that is, Subcase(A) and Subcase(B),
we proved that $\alpha_{j+1}=\chi_{j+1}$, and so, the proof of
Step(j+1) is done. Thus, the proof of (10.2.7) is finished. \ms

$\underline{\text{\rm The assertion in (10.2.8)}}$ In preparation
for the proof of the assertion in (10.2.8), first we prove that the
following are true:
$$\align
(10.2.9) \qquad \qquad &\text{ $n_j=\ell_j$, \quad
$\Delta_j(\beta_{j,k})^j_{k=1}=\omega_j(\delta_{j,k})^j_{k=1}$ \quad and} \\
&\text{$n_{j+1}n_{j+2}\cdots n_r={\ell}_{j+1}{\ell}_{j+2}\cdots
{\ell}_{\rho}$ \quad for $j=1,2,\dots,r\le {\rho}$.} \qquad \qquad
\qquad \qquad
\endalign$$

After the proof of (10.2.9) is done, then for the proof of (10.2.8)
it remains to show that $r=\rho$.

$\underline{\text{\rm For the proof of the assertion in (10.2.8)}}$
Now, in order to prove the equalities in (10.2.9), we will use
(10.2.7) together with (10.2.2.1) and (10.2.2.2) of Sublemma
$10.2.2$.

Now, the proof of (10.2.9) will be by induction on $j=1,2,\dots,r$,
as follows:

{\rm(1)} First, note by Sublemma $10.2.2$ and (10.2.7) that
$\dfrac{\beta_{1,1}}{n_1}=\dfrac{\alpha_1}{n}=\dfrac{\chi_1}{b}
=\dfrac{\delta_{1,1}}{\ell_1}$ and that
$\gcd(n_1,\beta_{1,1})=\gcd(\ell_1,\delta_{1,1})=1$. So,
$n_1=\ell_1$ and
$\Delta_1(\beta_{1,1})=\beta_{1,1}=\delta_{1,1}=\omega_1(\delta_{1,1})$.
Also, $n_1n_2n_3\cdots n_r=n=b=\ell_1\ell_2\ell_3\cdots \ell_{\rho}$
imply that $n_2n_3\cdots n_r=\ell_2\ell_3\cdots \ell_{\rho}$. \ms

{\rm(2)} Suppose we have shown by induction on the integer $r$ that
for any integer $j<r$, the following are true:
$$\align
(10.2.10) \qquad \qquad  & \text{$n_{j-1}=\ell_{j-1}$ \quad and
\quad $\Delta_{j-1}(\beta_{j-1,k})^{j-1}_{k=1}
=\omega_{j-1}(\delta_{j-1,k})^{j-1}_{k=1}$,} \qquad \qquad \qquad \\
& \text{$n_jn_{j+1}\cdots n_r=\ell_j\ell_{j+1}\cdots \ell_{\rho}$
\quad for $2\le j\le r$.}
\endalign$$

Then, by (10.2.7) and Sublemma 10.2.2,

$\dfrac{\widehat{\Delta}_j(\beta_{j,k})^j_{k=1}n_{j+1}n_{j+2}\cdots
n_r}{n_jn_{j+1}\cdots
n_r}=\dfrac{\alpha_j-\alpha_{j-1}}{n_{j}n_{j+1}\cdots
n_r}=\dfrac{\chi_j-\chi_{j-1}}{\ell_{j}\ell_{j+1}\cdots \ell_\rho}=
\dfrac{\widehat{\omega}_{j}(\delta_{j,k})^{j}_{k=1}\ell_{j+1}
\ell_{j+2}\cdots \ell_{\rho}}{\ell_j\ell_{j+1}\cdots \ell_\rho}$,

which implies that
$\dfrac{\widehat{\Delta}_j(\beta_{j,k})^j_{k=1}}{n_j}
=\dfrac{\widehat{\omega}_j(\delta_{j,k})^{j}_{k=1}}{\ell_j}$.

That is,
$\dfrac{{\Delta}_j(\beta_{j,k})^j_{k=1}-n_jn_{j-1}\Delta_{j-1}
(\beta_{j-1,k})^{j-1}_{k=1}}{n_j}
=\dfrac{{\omega}_j(\delta_{j,k})^{j}_{k=1}-
\ell_j\ell_{j-1}{\omega}_{j-1}(\delta_{j-1,k})^{j-1}_{k=1}}{\ell_j}$.

Since $\gcd(n_j,\Delta_j(\beta_{j,k})^j_{k=1})
=\gcd(\ell_j,\omega_j(\delta_{j,k})^j_{k=1})=1$ by Sublemma 10.2.2,
and $n_{j-1}=\ell_{j-1}$ and
$\Delta_{j-1}(\beta_{j-1,k})^{j-1}_{k=1}
=\omega_{j-1}(\delta_{j-1,k})^{j-1}_{k=1}$ by (10.2.10), then
$n_j=\ell_j$ and
${\Delta}_j(\beta_{j,k})^j_{k=1}={\omega}_j(\delta_{j,k})^{j}_{k=1}$.
Noting by (10.2.10) that $n_jn_{j+1}\cdots
n_r=\ell_j\ell_{j+1}\cdots \ell_{\rho}$, it is clear that
$n_{j+1}n_{j+2}\cdots n_r=\ell_{j+1}\ell_{j+2}\cdots \ell_{\rho}$.
Thus, the proof of (10.2.9) is done.

Now, to finish the proof of (10.2.8), it remains to prove that
$r=\rho$. Assume the contrary. Since it was assumed that $1\le j\le
r\le \rho$, then $r<\rho$. Then, $n_r=\ell_r$ and
${\Delta}_r(\beta_{r,k})^r_{k=1}={\omega}_r(\delta_{r,k})^{r}_{k=1}$.
Also, $n_r=\ell_r\ell_{r+1}\cdots \ell_{\rho}$ implies that
$1=\ell_{r+1}\ell_{r+2}\cdots \ell_{\rho}$, and so
$\ell_{r+1}=\ell_{r+2}=\cdots=\ell_\rho=1$, which would be
impossible. So, we proved that $r=\rho$.

Thus, the proof of (10.2.8) is done, and then the proof for Case(II)
is done. So, the proof for Fact[II] is done by Case(I) and Case(II),
and therefore the proof of this theorem is completely finished by
Fact[I] and Fact[II]. $\square$
\enddemo \ms

{\bf $\S$ 10.5. The difference between quasi-Puiseux convergent
power series of the recursive r-type and the Puiseux convergent
power series of the recursive r-type by Theorem 10.3}

\proclaim{Theorem 10.3} Let $\rho$ and $r$ be arbitrary positive
integers.

$\underline{\text{\bf Assumptions}}$

\noindent{\bf Assumptions[I]} By the same way as in either {\rm
Definition 5.0.0 or Theorem 7.3}, define a semi-quasi-Puiseux
convergent power series $g_r$ by {\rm Sequences[I]}, as follows:

{\bf Sequences[I]} Let $\{X_k:k=1,2,\dots,r\}$ with $X_k\subset
N_0$, $\{g_k:k=1,2,\dots,r\}$ with $g_k=g_k(y,z) \in \BC\{y,z\}$ and
\text\{$\Delta_k:N^k_0\to N_0$\} is an integer-valued function for
$1\le k\le r$ be three sequences satisfying the following four
conditions:

\noindent Four conditions are denoted by \text{\bf The 1st
${\text{\bf{Cond}}}^{\text{{\bf(0)}}}$}, \dots, \text{\bf The 4-th
${\text{\bf{Cond}}}^{\text{{\bf(0)}}}$} \text{\bf of Sequences[I]}.
\ms

{\bf[I]-(1)} \text{\bf The 1st
${\text{\bf{Cond}}}^{\text{{\bf(0)}}}$ of Sequences[I]:}

{\rm (1a)}  $X_1=\{n_1,\beta_{1,1}\}$ with $n_1=1< \beta_{1,1}$.

 {\rm (1b)} $X_j=\{n_j,\beta_{j,1},\beta_{j,2},\dots,\beta_{j,j}\}$
 with $n_j\ge 2$ \quad where $j=2,\dots,r$.

 If $j\ge 2$, then assume that at least one of
$\beta_{j,1},\beta_{j,2},\dots,\beta_{j,j}$ is nonzero. \ms

{\bf[I]-(2)} \text{\bf The 2nd
${\text{\bf{Cond}}}^{\text{{\bf(0)}}}$ of Sequences[I]:}

{\rm (2a)} $g_1=z^{n_1}+y^{\beta_{1,1}}$.

{\rm (2b)} $g_j=g_{j-1}^{n_j}+y^{\beta_{j,1}}z^{\beta_{j,2}}
g_1^{\beta_{j,3}}\cdots g_{j-2}^{\beta_{j,j}}$ \quad where
$j=2,\dots,r$. \ms

{\bf[I]-(3)} \text{\bf The 3rd
${\text{\bf{Cond}}}^{\text{{\bf(0)}}}$ of Sequences[I]:}

{\rm(3a)} $\Delta_1(t)=t$ for each $t\in N_0$.

{\rm(3b)}
$\Delta_j(t_j)^j_{k=1}=t_j\Delta_{j-1}(\beta_{j-1,k})^{j-1}_{k=1}
+n_{j-1}\Delta_{j-1}(t_k)^{j-1}_{k=1}$ for each $(t_k)^j_{k=1}\in
N^j_0$

\indent where $j=2,\dots,r$. \ms

{\bf[I]-(4)} \text{\bf The 4-th
${\text{\bf{Cond}}}^{\text{{\bf(0)}}}$ of Sequences[I]:}

{\rm(4)} $\Delta_j(\beta_{j,k})^j_{k=1}>n_jn_{j-1}
\Delta_{j-1}(\beta_{j-1,k})^{j-1}_{k=1}$ for $2\le j\le r$. \ms

\noindent{\bf Assumptions[II]} By using {\rm Sequences[I]} as above,
it is clear to define another convergent power series $h_{r-1}$ by
{\rm Sequences[II]}, as follows:

{\bf Sequences[II]} \quad Using the nonsingular analytic mapping
defined by $(g_1,y) \rightarrow (z,y)$ from {\bf[I]-(2)} of {\bf
Sequences[I]} where $g_1 =z+y^{\beta_{11}}$, construct {\bf
Sequences[II]} consisting of three sequences, that is,
$\{W_k:k=1,2,\dots,r-1\}$ with $W_k\subset N_0$,
$\{h_k:k=1,2,\dots,r-1\}$ with $h_k=h_k(y,z)\in \BC\{y,z\}$ and
\text\{$\omega_k:N^k_0\to N_0$ is an integer-valued function for
$k=1,2,\dots,r-1$\}, satisfying the following three conditions:

\noindent Three conditions are denoted by \text{\bf The 1-th
${\text{\bf{Cond}}}^{\text{{\bf(0)}}}$}, \dots, \text{\bf The 3-th
${\text{\bf{Cond}}}^{\text{{\bf(0)}}}$} \text{\bf of Sequences[II]}.
\ms

{\bf[II]-(1)} \text{\bf The 1st
${\text{\bf{Cond}}}^{\text{{\bf(0)}}}$ of Sequences[II]:}

{\rm(1a)} $W_1 =\{c_1,\sigma_{1,1}\}$ \quad with \quad $c_1\ge 2$,

{\rm(1b)} $W_j
=\{c_j,\sigma_{j,1},\sigma_{j,2},\dots,\sigma_{j,j}\}$ \quad with
\quad $c_j\ge 2$, \quad noting that $2\le j\le r-1$,

satisfying the following inequalities:

{\rm(1a-a)}  $c_1=n_2$, $\sigma_{1,1}=\Delta_2(\beta_{2,1},
\beta_{2,2})>n_2n_1\beta_{1,1}=n_2\beta_{1,1}>n_2$,

{\rm(1b-a)}  $c_j=n_{j+1}$, $\sigma_{j,1}=\Delta_2(\beta_{j+1,1},
\beta_{j+1,2})$, $\sigma_{j,2}=\beta_{j+1,3}$,
$\sigma_{j,3}=\beta_{j+1,4}$,\dots,$\sigma_{j,j}=\beta_{j+1,j+1}$.
\ms

{\bf[II]-(2)} \text{\bf The 2nd
${\text{\bf{Cond}}}^{\text{{\bf(0)}}}$ of Sequences[II]:} Note that
$2\le j\le r-1$.

{\rm(2a)} $h_1 =z^{c_1}+\eta_1y^{\sigma_{1,1}}$.

{\rm(2b)} $h_j
=h^{c_j}_{j-1}+\eta_jy^{\sigma_{j,1}}z^{\sigma_{j,2}}h^{\sigma_{j,3}}_1\cdots
 h^{\sigma_{j,j}}_{j-2}$ for $2\le j\le r-1$.

Note that $\eta_i=\eta_i(y,z)$ is a unit in $\BC\{y,z\}$ for $1\le
i\le r-1$. \ms

{\bf[II]-(3)} \text{\bf The 3rd
${\text{\bf{Cond}}}^{\text{{\bf(0)}}}$ of Sequences[II]:} Note that
$2\le j\le r-1$.

{\rm(3a)} $\omega_1(t)=t$ \quad for each $t\in N_0$.

{\rm(3b)}
$\omega_{j}(t_k)^{j}_{k=1}=t_{j}\omega_{j-1}(\sigma_{j-1,k})^{j-1}_{k=1}
+c_{j-1}\omega_{j-1}(t_k)^{j-1}_{k=1}$ \quad for each
$(t_k)^j_{k=1}\in N^j_0$. \ms

In the above assumptions, note that {\bf Sequences[I]} and {\bf
Sequences[II]} have the same kind of three conditions except for
\text{\bf The 4-th ${\text{\bf{Cond}}}^{\text{{\bf(0)}}}$ of
Sequences[I]}. \ms

$\underline{\text{\bf Conclusions}}$ Then, we have the following
three facts: \ms

$\underline{\text{\bf Fact(1)}}$ Then, $h_{r-1}$ of {\bf
Sequences[II]} either is a semi-quasi-Puiseux convergent power
series or satisfies the following equations.
$$\align
\text{\rm(10.3.1)} \quad \quad
&\text{$\omega_j(\sigma_{j,k})^j_{k=1}=\Delta_{j+1}(\beta_{j+1,k})^{j+1}_{k=1}$
\quad for each $j=1,2,\dots,r-1$,} \qquad \qquad \qquad \qquad \\
\text{\rm (10.3.2)} \quad \quad
&\text{$\omega_j(\sigma_{j,k})^j_{k=1}>c_jc_{j-1}
\omega_{j-1}(\sigma_{j-1,k})^{j-1}_{k=1}$ \quad for each
$j=2,3,\dots,r-1$}.
\endalign$$

In particular, the inequality in {\rm(10.3.2)} is equivalent to the
{\rm 4-th} condition of {\bf Sequences[II]}, denoted by \text{\bf
The 4-th ${\text{\bf{Cond}}}^{\text{{\bf(0)}}}$ of Sequences[II]},
noting that $\omega_1(t)=t$  for each $t\in N_0$. \ms

$\underline{\text{\bf Fact(2)}}$ $g_{r}$ is quasi-Puiseux convergent
power series by {\rm Sequence(I)} $\iff$ $h_{r-1}$ is quasi-Puiseux
convergent power series. \ms

$\underline{\text{\bf Fact(3)}}$ Let $g_r$ be irreducible in
$\BC\{y,z\}$. Then, we have the following:
$$\align
\text{$g_r$ and $h_{r-1}$ have the same
multiplicity sequence.} \quad  \text{$\square$}  \tag 10.3.3 \\
\endalign$$
\endproclaim \ms

\demo{\bf Proof of Theorem 10.3} We prove Fact(1), Fact(2) and
Fact(3), respectively.

$\underline{\text{\bf Fact(1)}}$ If the equation in (10.3.1) is
true, there is nothing to prove for the inequality in {\rm (10.3.2)}
because
$\omega_j(\sigma_{j,k})^j_{k=1}=\Delta_{j+1}(\beta_{j+1,k})^{j+1}_{k=1}
>n_{j+1}n_j\Delta_{j}(\beta_{j,k})^{j}_{k=1}
=c_jc_{j-1} \omega_{j-1}(\sigma_{j-1,k})^{j-1}_{k=1}$ for each
$j=2,3,\dots,r-1$, by using the equation in (10.3.1), \text{\rm The
4-th ${\text{\rm{Cond}}}^{\text{{\rm(0)}}}$ of Sequences[I]} and
\text{\rm The 1st ${\text{\rm{Cond}}}^{\text{{\rm(0)}}}$ of
Sequences[II]}. So, for the proof of Fact(1), it suffices to show
that the equality in {\rm(10.3.1)} is true.

In preparation for the proof of the equality in {\rm(10.3.1)}, using
$c_j=n_{j+1}$ and
$\sigma_{j,1}=\Delta_2(\beta_{j+1,1},\beta_{j+1,2})$ for $1\le j\le
r-1$ in \text{\rm The 1st ${\text{\rm{Cond}}}^{\text{{\rm(0)}}}$} of
{\rm Sequences[II]}, first it suffices to show by (10.3.5) that the
following is true: Let $(t_k)^j_{k=1}\in N^j_0$ be chosen arbitrary.
$$\align
(10.3.4) \quad & \text{Whenever $t_1$ satisfies the equation
$t_1=\Delta_2(\delta_1,\delta_2)$
for any $(\delta_1,\delta_2)\in   N^2_0$,} \\
& \text{ $\omega_j(t_1,t_2,t_3,\dots,t_j)
=\Delta_{j+1}(\delta_1,\delta_2,t_2,t_3,\dots,t_j)$ for all
$j=1,2,\dots,r-1$.} \qquad \qquad \qquad \\
\endalign$$

As an application of {\rm(10.3.4)}, if
$t_1=\sigma_{j,1}=\Delta_2(\beta_{j+1,1},\beta_{j+1,2})$,
$t_2=\sigma_{j,2}$, $t_3=\sigma_{j,3}$, $ \dots$,
$t_j=\sigma_{j,j}$, then there is nothing to prove for the equality
in (10.3.1), because
$$\align
(10.3.5) \qquad \qquad
\omega_j(\sigma_{j,1},\sigma_{j,2},\dots,\sigma_{j,j})
&=\Delta_{j+1}(\beta_{j+1,1},\beta_{j+1,2},\sigma_{j,2},\sigma_{j,3},\dots,\sigma_{j,j})
\qquad \qquad \qquad \qquad\\
&=\Delta_{j+1}(\beta_{j+1,k})^{j+1}_{k=1}, \qquad \qquad \qquad \qquad \\
\endalign$$
where $\sigma_{j,1}=\Delta_2(\beta_{j+1,1},\beta_{j+1,2})$,
$\sigma_{j,k}=\beta_{j+1,k}$ with $2\le k\le j$, by using the
equation in \text{\rm The 1st ${\text{\rm{Cond}}}^{\text{{\rm(0)}}}$
of Sequences[I]}. \ms

Now, the proof of $(10.3.4)$ will be by induction on the positive
integer $j\le r-1$.

If $j=1$, then there is nothing to prove because $\omega_1(t)=t$ for
each $t\in N_0$.

By the induction proof, suppose we have shown that $(10.3.4)$ is
true for any positive integer $j$ where $1\le j\le r-2$.

Then, it is enough to prove the following: Let $(t_k)^{j+1}_{k=1}\in
N^{j+1}_0$ be chosen arbitrary.
$$\align
(10.3.6) \qquad \qquad & \text{Whenever $t_1$ satisfies the equation
$t_1=\Delta_2(\delta_1,\delta_2)$
for any $(\delta_1,\delta_2)\in   N^2_0$,} \qquad \\
& \text{ $\omega_{j+1}(t_1,t_2,t_3,\dots,t_{j+1})
=\Delta_{j+2}(\delta_1,\delta_2,t_2,t_3,\dots,t_{j+1})$,}
\qquad \qquad \qquad \\
\endalign$$
for any $j+1$ where $2\le j\le r-2$.

For the proof of $(10.3.6)$, whenever $t_1$ satisfies the equation
$t_1=\Delta_2(\delta_1,\delta_2)$ for any $(\delta_1,\delta_2)\in
N^2_0$, then we may suppose by the induction assumption on \text{\rm
$j<r-1$} that
$\omega_j(t_1,t_2,\dots,t_j)=\Delta_{j+1}(\delta_1,\delta_2,t_2,\dots,t_j)$
with
$\omega_j(\sigma_{j,k})^j_{k=1}=\Delta_{j+1}(\beta_{j+1,k})^{j+1}_{k=1}$.
\ms

Then, it is trivial to show that the following equality is true:
$$\align
(10.3.7) \qquad \qquad \omega_{j+1}(t_k)^{j+1}_{k=1}
&=c_j\omega_j(t_1,t_2,\dots,t_j)+t_{j+1}\omega_j(\sigma_{j,k})^j_{k=1}
\\
 &=n_{j+1}\Delta_{j+1}(\delta_1,\delta_2,t_2,\dots,t_j)
+t_{j+1}\Delta_{j+1}(\beta_{j+1,k})^{j+1}_{k=1} \qquad \qquad \\
&=\Delta_{j+2}(\delta_1,\delta_2,t_2,\dots,t_j,t_{j+1}),
\endalign$$
by {\rm[II]-(3)} of {\rm Sequences[II]}, by the induction assumption
on \text{\rm $j<r$} and {\rm[I]-(3)} of {\rm Sequences[I]}, which
gives the proof for (10.3.4) on the integer $(j+1)$.

Thus, the proof of (10.3.4) is done, and so the proof of (10.3.1) is
done, too. Therefore, the proof of {\rm Fact(1)} is finished. \ms

$\underline{\text{\bf Fact(2)}}$ Since Sequences[II] has the same
kind of four conditions as we have seen in the assumption of Theorem
$5.0$ by Fact(1), then by Theorem $5.0$ we have the following:
$$\split
(10.3.8) \qquad \qquad &\text{$h_r$ is irreducible in $\BC\{y,z\}$} \\
\iff &
\text{$\gcd(c_1,\sigma_{1,1})=\gcd(c_2,\omega_2(\sigma_{2,1},\sigma_{2,2}))
=\cdots  =\gcd(c_r,\omega_r(\sigma_{r,k})^r_{k=1})=1$}. \qquad
\qquad \\
\endsplit$$

Since for each $j=1,2,\dots,r-1$,
$\omega_j(\sigma_{j,k})^j_{k=1}=\Delta_{j+1}(\beta_{j+1,k})^{j+1}_{k=1}$
by (10.3.1) and $c_j=n_{j+1}$ by (1) of {\rm Sequences(II)}, then
$\gcd(c_j,\omega_j(\sigma_{j,k})^j_{k=1})
=\gcd(n_{j+1},\Delta_{j+1}(\beta_{j+1,k})^{j+1}_{k=1})$. Therefore,
by (10.3.8) and Theorem $5.0$, we have the following:
$$
(10.3.9) \quad \qquad \text{$h_{r-1}$ is irreducible in
$\BC\{y,z\}$} \iff \text{$g_{r}$ is irreducible in $\BC\{y,z\}$}.
\qquad \qquad \qquad
$$
Thus, the proof of {\rm Fact(2)} is done. \ms

$\underline{\text{\bf Fact(3)}}$ Let $g_r$ be irreducible in
$\BC\{y,z\}$. By (10.1.1) of Theorem $10.1$,
$\{\text{[Mult$(V(g_r))$]}\}$, called the multiplicity sequence of
$V(g_r)$, can be represented as follows:
$$\align
(10.3.10) \qquad & \text{\rm Multiseq($V(g_r)$)} =
\text{Join}(\{[n:\alpha_1]\},
\text\{[\gcd(n,\alpha_1):\alpha_2-\alpha_1]\},\dots, \qquad\\
 &
 \text\{[\gcd(n,\alpha_1,\dots,\alpha_{r-2}):\alpha_{r-1}-\alpha_{r-2}]\},
\text\{[\gcd(n,\alpha_1,\dots,\alpha_{r-1}):\alpha_{r}-\alpha_{r-1}]\}), \qquad \qquad \\
&  \quad\text{\rm such  that}  \qquad \qquad n=n_1n_2\cdots n_r,
\qquad \alpha_1 =\beta_{1,1}n_2\cdots n_r \quad \text{and}\\
& \quad \quad \qquad \qquad \qquad  \alpha_j =\alpha_{j-1}+
\widehat{\Delta}_jn_{j+1}n_{j+2}\cdots  n_r
\quad \text{for $2\le j\le r$},\\
 \endalign$$
where $\widehat{\Delta}_j
=\Delta_j(\beta_{j,k})^j_{k=1}-n_jn_{j-1}\Delta_{j-1}(\beta_{j-1,k})^{j-1}_{k=1}>0$
for $2\le j\le r$ and $\Delta_1(t)=t$.
 \ms

Note by Fact(2) that $h_{r-1}$ is irreducible in $\BC\{y,z\}$. So,
as above, \text{\rm Multiseq($V(h_{r-1})$)} can be represented as
follows:
$$\align
(10.3.11) \quad & \text{\rm Multiseq($V(h_{r-1})$)} =
\text{Join}(\{[c:\xi_1]\},
\{[\gcd(c,\xi_1):\xi_2-\xi_1]\},\dots, \\
 & \{[\gcd(c,\xi_1,\xi_2,\dots,\xi_{r-3}):\xi_{r-2}-\xi_{r-3}]\},
\{[\gcd(c,\xi_1,\xi_2,\dots,\xi_{r-2}):\xi_{r-1}-\xi_{r-2}]\}), \qquad \qquad \\
&  \text{such ~ that}  \qquad \qquad c=c_1c_2\cdots c_{r-1}, \quad
\xi_1 =\sigma_{1,1}c_2\cdots c_{r-1} \quad \text{and} \\
& \quad \quad \qquad \qquad \qquad  \xi_j =\xi_{j-1}+
\widehat{\omega}_j(\sigma_{j,k})^j_{k=1}c_{j+1}c_{j+2}\cdots
c_{r-1} \quad \text{for} \quad 2\le j\le r-1,
 \endalign$$
where $\widehat{\omega}_j
=\omega_j(\sigma_{j,k})^j_{k=1}-c_jc_{j-1}\omega_{j-1}(\sigma_{j-1,k})^{j-1}_{k=1}$
for $2\le j\le r-1$ and $\omega_1(t)=t$. \ms

Then, by Fact(1), (10.3.10) and (10.3.11), we have the following
properties:
$$\align
(10.3.12) \qquad\quad
&\xi_j-\xi_{j-1}=\widehat{\omega}_{j}c_{j+1}c_{j+2}\cdots
c_{r-1} \qquad \qquad \text{for any $j=2,3,\dots,r-1$,} \qquad \qquad \\
&\qquad \qquad=\widehat{\Delta}_{j+1}n_{j+2}n_{j+3}\cdots
n_{r}=\alpha_{j+1}-\alpha_j>0. \qquad \qquad  \\
\endalign$$

Since $c=c_1c_2\cdots c_{r-1}$ and $\xi_1 =\sigma_{1,1}c_2\cdots
c_{r-1}$ by (10.3.11) and $\alpha_1=\beta_{1,1}n_2\cdots n_{n}$ by
(10.3.10), then
$$\align
(10.3.13) \qquad \quad
\xi_1-\beta_{1,1}c&=(\sigma_{1,1}-\beta_{1,1}c_1)c_2\cdots c_{r-1}
\\
&=(\Delta_{2}(\beta_{2,k})^{2}_{k=1} -n_1n_2\beta_{1,1})n_3\cdots
n_{r}=\alpha_2-\alpha_1>0 \qquad \qquad \qquad \qquad
\endalign$$ by (10.3.10)
where $n_{1}=1$ was in Sequences[I], and
$\sigma_{1,1}=\Delta_{2}(\beta_{2,k})^{2}_{k=1}$ and $c_j=n_{j+1}$
for $1\le j\le r-1$ were defined by Sequences[II].

Since $n=n_1n_2\cdots n_{r}$ and $\alpha_1 =\beta_{1,1}n_2n_3\cdots
n_{r}$ with $n_{1}=1$ is defined by (10.3.10), then
$\gcd(n,\alpha_1)=n_2\cdots n_{r}=c_1c_2\cdots c_{r-1}=c$, and so we
have the following:
$$\align
(10.3.14) \qquad \qquad & \gcd(c,\xi_1,\xi_2,\dots,\xi_{j})
\quad \text{for any $j=2,3,\dots,r-1$,} \\
=& \gcd(c,\xi_1-\beta_{1,1}c,\xi_2-\xi_1,\dots,\xi_{j}-\xi_{j-1})
 \quad \text{by (10.3.13)} \\
 =&
\gcd(\gcd(n,\alpha_1),\alpha_2-\alpha_1,\alpha_3-\alpha_2,\dots,\alpha_{j+1}-\alpha_j)
\quad \text{by (10.3.12)} \qquad \qquad \qquad \\
=& \gcd(n,\alpha_1,\alpha_2,\dots,\alpha_{j+1}). \\
\endalign$$

By (10.3.12) and (10.3.14), it can be proved that the following
equality is true in the sense of Definition 9.2: For any
$j=2,3,\dots,r-1$,
$$\align
(10.3.15) \qquad
\{[\gcd(n,\alpha_1,\alpha_2,\alpha_3,\dots,\alpha_{j+1}),\alpha_{j+1}-\alpha_j]\}=
\{[\gcd(c,\xi_1,\xi_2,\dots,\xi_{j}),\xi_{j}-\xi_{j-1}]\}. \qquad
\qquad
\endalign$$

To prove that (10.3.10) and (10.3.11) are the same, it remains to
show by (10.3.15) that the following equality is true:
$$\align
\text{Join}(\{[n:\alpha_1]\},\{[\gcd(n,\alpha_1):\alpha_2-\alpha_1]\})
=\{[c:\xi_1]\}. \tag 10.3.16
\endalign$$

Now, the proof of (10.3.16) is as follows:
$$\align
(10.3.17) \quad  \{[c:\xi_1]\} &=\{[c_1c_2\cdots
c_{r-1}:\sigma_{1,1}c_2c_3\cdots
c_{r-1}]\} \quad \text{by (10.3.11)} \\
&=\{[n_2n_3n_4\cdots n_{r}:
\sigma_{1,1}n_3n_4\cdots n_{r}]\} \quad \text{by Sequences[II]} \\
&=\{[n_1n_2n_3\cdots
n_{r}:\Delta_{2}(\beta_{2,k})^{2}_{k=1}n_3n_4\cdots
n_{r}]\} \quad \text{with $1=n_1<\beta_{1,1}$} \\
&= \text{Join}(\{[n_1n_2\cdots n_{r}\:\beta_{1,1}n_2\cdots n_{r}]\}
,\{[n_2n_3\cdots n_r: (\widehat{\Delta}_2n_3\cdots
n_r]\}) \qquad \qquad \qquad\\
&=\text{Join}(\{[n:\alpha_1]\};\{[\gcd(n,\alpha_1):\alpha_2-\alpha_1]\}),
\endalign$$
by Definition 9.2(the definition of the multiplicity sequence for
two positive integers) because $\beta_{1,1}>n_1=1$ and
$\widehat{\Delta}_2=\Delta_{2}(\beta_{2,k})^{2}_{k=1}-\beta_{1,1}n_1n_2>0$
by Sequences[II]. Thus, the proof of Fact(3) is finished.

Therefore, the proof of theorem is completely finished by the proofs
of Fact(1), Fact(2) and Fact(3). $\square$
\enddemo \bs

\vfill \pagebreak

{\bf Part[B5] The 1st Algorithm with proofs} \bs

{\bf {\S11.} The 1st Algorithm for computing
 a one-to-one function from Family(1) onto Family(2) with proofs}
\ms

{\bf \S11.0. Introduction}\ms

In this section, the aim is to find the 1st Algorithm for compute a
one-to-one function between Family(1) and Family(2). Then, the first
half of The 1st Algorithm is given by Theorem 11.2 of \S11.2, and
the second half of The 1st Algorithm is given by Theorem 11.4 of
\S11.4. \ms

{\bf \S11.1. Some definitions} \ms

\definition{Definition 11.1} Let $f(y,z)$ be irreducible in
$\BC\{y,z\}$ with isolated singularity at $0\in \BC^2$. By
Definition 1.2, let the standard Puiseux expansion $C_r(t)$ of the
$r$-type for any irreducible curve $C(t)$ be given as follows:
$$\align
(11.1.1) \quad \quad \text{$C_r(t)=:$} \left\{\eqalign{ y=&t^n, \cr
z=&t^{\alpha_1}+t^{\alpha_2}+\cdots +t^{\alpha_r}, \cr} \right. \\
\text{where} \quad
 2\le n <\alpha_1<\cdots <\alpha_r   \quad  \text{and} \quad & \\
\qquad  n >d_1>\cdots
 >d_r=1  \quad   \text{with}  \quad &
\text{$d_i=\gcd(n,\alpha_1,\dots,\alpha_i)$, $1\le i\le r$.} \qquad
\qquad \qquad
\endalign$$

If the above curve $C(t)$ and $V(f)$ have the same multiplicity
sequence at $(y,z)=(0,0)$, it is said by Definition 1.2 that
$$\align
(11.1.2) \quad \text{either ${C(t)\equiv V(f)}$ {\rm (Multiseq)}
\quad {or} \quad \text{\text{\rm
Multiseq($C(t)$)}{$\equiv$}\text{\rm Multiseq(V(f))}} as sequence.}
\endalign$$
\enddefinition \ms

\definition{Remark 11.1.1 for Definition 11.1} As in Definition 8.1,
let the parametrization $C(t)$ for the above irreducible curve $C$
of (11.1.2) be defined by
$$ \text{$y(t)=t^n$ and
$z(t)=c_1t^{k_1}+c_2t^{k_2}+\cdots=c_1t^{k_1}(1+H(t))$,} \tag
11.1.3$$ where $1<n$, $1<k_1<k_2<\cdots, $ and the $c_i$ are nonzero
complex numbers and $H(t)$ is just the substitution. Also, let the
parametrization of another irreducible plane curve $C'$ be given by
$y(t)=t^n(1+H(t))$ and $z(t)=t^{k_1}$. Note that either $n\le k_1$
or $k_1<n$.

Then, it was already directly shown by Theorem A(Theorem 8.10) that
two irreducible curves $C$ and $C'$ have the same multiplicity
sequence, and also the same Puiseux pairs by a nonsingular change of
a parameter, without using Theorem B. But, without using Theorem A,
it has been not yet proved by Theorem B only that these two
irreducible curves $C$ and $C'$ have the same multiplicity
sequences.

So, it was proved by Theorem 8.10 that an equivalence relation in
(11.1.2) is well-defined.
\enddefinition \ms

{\bf \S11.2. The first half of The 1st Algorithm(Theorem 11.2)}

\proclaim{Theorem 11.2(Theorem 1.4:Algorithm for finding a
one-to-one function from Family(1) into Family(2))}

$\underline{\text{\bf {Assumptions}}}$ Let $r$ be arbitrary positive
integer. By the same way as in {\rm Definition 5.0.0} or {\rm
Theorem 7.3}, define a quasi-Puiseux convergent power series $g_r$
of recursive $r$-type in $\BC\{y,z\}$ by {\rm {Sequences[I]} in the
assumptions of Theorem $7.3$. For notation, assume in addition that
$1\le n_1<\beta_{1,1}$ in \text{\rm The 1-th
${\text{\rm{Cond}}}^{\text{{\rm(0)}}}$} of {\rm {Sequences[I]} of
Theorem $7.3$. Now, we may use the same conditions and notations as
in the assumptions of Theorem $7.3$. \ms

$\underline{\text{\bf Conclusions}}$ It is very elementary to
compute the standard Puiseux expansion $C^*(t)$ of an irreducible
curve $C(t)$ such that \text{$V(g_r)\equiv C^*(t)$ \text{\rm
(Multiseq)}} with desired algorithms in {\rm Fact(A)} and {\rm
Fact(B)}. \ms

\noindent$\underline{\text{\rm {\bf Fact(A):} {\rm By explicit
algorithm(Algorithm 11.2.1),} we can compute the Puiseux }}$

\noindent$\underline{\text{\rm expansion $C(g_r:t)$ for the curve
$C(t)$ such that \text{$V(g_r)\equiv C(g_r:t)$ \text{\rm
(Multiseq)}}.}}$ \ms

\noindent$\underline{\text{\rm(Algorithm 11.2.1 for Theorem 11.2)}}$
$$\align
 \text{$C(g_r:t):=$} & \left\{\eqalign{y&=t^n \cr
z&=t^{\alpha_1}+t^{\alpha_2}+\cdots +t^{\alpha_r}, \cr} \right. \tag
11.2.1
 \\
 \text{such that}  \qquad
 n &=n_1n_2\cdots n_r \quad \text{and} \quad
 \alpha_1 =\beta_{1,1}n_2\cdots n_r, \\
 \alpha_j &=\alpha_{j-1}+
\widehat{\Delta}_jn_{j+1}n_{j+2}\cdots
 n_r, \\
 \endalign$$
where $\widehat{\Delta}_j
=\Delta_j(\beta_{j,k})^j_{k=1}-n_jn_{j-1}\Delta_{j-1}(\beta_{j-1,k})^{j-1}_{k=1}$
is a positive integer for $2\le j\le r$ and $\Delta_1(t)=t$. \ms

{\rm Remark for Fact(A).} Note that the parametrization in (11.2.1)
satisfies the following:

\noindent (11.2.2) \quad {\rm(i)} \quad
$n<\alpha_1<\alpha_2<\alpha_3<\cdots <\alpha_r$.

\qquad \quad \text{\rm(ii)} \quad $n\ge d_1>d_2>\cdots >d_r=1$ with
$d_i= \gcd(n,\alpha_1,\dots,\alpha_i)$,

\quad \quad\qquad where $d_i=n_{i+1}n_{i+2}\cdots n_r$ for
$i=1,2,\dots,r-1$. \ms

\noindent$\underline{\text{\rm {\bf Fact(B):} {\rm By explicit
algorithm(Algorithm 11.2.2),} we can compute the standard }}$

\noindent$\underline{\text{\rm Puiseux expansion $C^*(t)$ for the
curve $C(g_r:t)$ in Fact(A) such that \text{$V(g_r)\equiv C^*(t)$
\text{\rm (Multiseq)}}.}}$ \ms

\noindent$\underline{\text{\rm(Algorithm 11.2.2 for Theorem 11.2)}}$

{\rm(a)} If $n>\gcd(n,\alpha_1)$, $C^*(t)=C(g_r:t)$ is the standard
Puiseux expansion of the {\rm r-th} type where the curve $C(g_r:t)$
was already in Fact(A). \ms

{\rm(b)} If $n=\gcd(n,\alpha_1)$, $C^*(t)$ can be parametrized by
the standard Puiseux expansion $C_{r-1}(t)$ of the {\rm (r-1)-th}
type:
$$\text{$C_{r-1}(t):=$}\left\{\eqalign{ y &=t^n \cr
z &=t^{\alpha_2}+t^{\alpha_3}+\cdots +t^{\alpha_{r}}. \cr} \right.
\tag 11.2.3$$

Note that \text{$C(g_r:t)\equiv C_{r-1}(t)$} \text{\rm (Multiseq)}
where the curve $C(g_r:t)$ is in (11.2.1). \ms

\noindent{\bf Fact(C):} If for any $g_r$ in {\rm Family(1)} define a
function \text{\rm $\Psi$ :Family(1)$\rightarrow$Family(2)} by
$\Psi(g_r)=C^*(t)$, then $\Psi$ is a one-to-one function from {\rm
{Family(1)}} into {\rm Family(2)}. \quad $\square$
\endproclaim \ms

\definition{Remark 11.2.1} {\rm(a)} $n_1=1$ if and only if
$n=\gcd(n,\alpha_1)=d_1$. \ms

(b) For each $j=1,2,\dots,r$, note that $(0,0)$ is an isolated
singularity of an analytic variety $V(g_j)$ except for $V(g_1)$ with
$n_1=1$. \ms

(c) Whenever $V(g_r)$ is chosen arbitrary as in the assumption of
Theorem $11.2$, then it was already proved by {\rm(10.1.1)} of
Theorem $10.1$ that the multiplicity sequence of $V(g_r)$, that is,
\text{\rm Multiseq($V(g_r)$)}, can be represented as follows:
$$\align
\text{\rm(11.2.1-1)} \quad & \text{\rm Multiseq($V(g_r)$)} =
\text{Join}(\{[n:\alpha_1]\},
\text\{[\gcd(n,\alpha_1):\alpha_2-\alpha_1]\},\dots, \quad \\
&
\text\{[\gcd(n,\alpha_1,\dots,\alpha_{r-2}):\alpha_{r-1}-\alpha_{r-2}]\},
\text\{[\gcd(n,\alpha_1,\dots,\alpha_{r-1}):\alpha_{r}-\alpha_{r-1}]\}). \quad \text{$\square$} \\
 \endalign$$
\enddefinition

\noindent{\bf Example 11.2.2 for the first half of The 1st algorithm
in Theorem 11.2:} Example 11.2.2 and Example 1.4.1 of $\S 1$ are the
same. \quad \text{$\square$} \ms

Rigorously, the statement for Fact(C) of Theorem 11.2 can be
rewritten by Corollary $11.3$.

\proclaim{Corollary 11.3(The uniqueness of the solution for Fact(B)
of Theorem 11.2)}

$\underline{\text{\bf Assumptions}}$ First, suppose that the same
assumption as in {\rm Sequences[I]} of Theorem $7.3$ or Theorem
$11.2$ is satisfied with an additional condition that
$\underline{2\le n_1<\beta_{11}}$. Next, as we have seen in {\rm
Sequences[I]} of the assumption of  Theorem $7.3$, we define another
quasi-Puiseux convergent power series $\phi_{\rho}$ of recursive
$\rho$-type in $\BC\{y,z\}$ by {\rm Sequences[II]} in the
assumptions of Theorem $7.3$, with an additional condition that
$\underline{2\le \ell_1<\delta_{11}}$. \ms

$\underline{\text{\bf Conclusions}}$ Note that
$V(g_r)=\{(y,z):g_r(y,z)=0\}$ and
$V(\phi_\rho)=\{(y,z):\phi_\rho(y,z)=0\}$ be analytic varieties at
$(y,z)=(0,0)$ with isolated singularity at the origin. Then, we have
the following:
$$\align
& \text{$n_j=\ell_j$ \quad  and}  \tag 11.3.2 \\
& \text{$\Delta_j(\beta_{j,k})^j_{k=1}
=\omega_j(\delta_{j,k})^j_{k=1}$ \quad for each $j=1,2,\dots,r=\rho$} \\
\Longleftrightarrow \quad & \\
& \text{$V(g_r)\equiv V(\phi_\rho)$ \quad \text{\rm (Multiseq)}}  \tag 11.3.3\\
\Longrightarrow \quad & \\
& \text{$C(g_r:t)\equiv C(\phi_\rho:t)$ \quad \text{\rm (Multiseq)}}  \tag 11.3.4 \\
& \text{where $V(g_r)\equiv C(g_r:t)$ and $V(\phi_\rho)\equiv
C(\phi_\rho:t)$}
\quad \text{\rm (Multiseq)}\\
\Longrightarrow \quad & \\
& \text{$C(g_r:t)$ and $C(\phi_\rho:t)$ have the same
standard Puiseux expansion}  \tag 11.3.5 \\
\endalign$$

Moreover, if $g_{r}$ and $\phi_\rho$ is irreducible in $\BC\{y,z\}$,
it can be easily proved by Theorem $7.3$ and Theorem $10.2$ that
\text{$V(g_r)\equiv V(\phi_\rho)$ \text{\rm (Multiseq)}} if and only
if $g_r \buildrel \text{{\rm divisor}} \over \sim \phi_\rho$ under
the standard resolutions. \quad \text{$\square$}
\endproclaim \ms

{\bf \S11.3. The proof of Theorem 11.2 } \ms

\demo{\bf Proof of Theorem 11.2} First of all, note by Sublemma
$5.2$ of Theorem $5.0$ that the multiplicity of $g_r$ at
$(y,z)=(0,0)$ is $n=\prod^r_{k=1}n_k$ because $1 \le
n_1<\beta_{1,1}$.

For the proof of theorem, it suffices to prove that Fact(A) is true
because of the following:

{\rm(i)} The proof of Fact(A) with Theorem 8.10 and Definition 11.1
implies clearly that Fact(B) is true. \ms

{\rm(ii)}  It is trivial by Theorem $7.3$, Theorem 7.7, and Theorem
10.2 that Fact(C) is true. \ms

$\underline{\text{\bf Fact(A)}}$ For the proof, Fact(A) is divided
into two steps:

\text{\rm Fact(A-1)} Then, the curve $C(g_r:t)$ of (11.2.1) is the
Puiseux expansion at $(y,z)=(0,0)$.

\text{\rm Fact(A-2)} Then, \text{$V(g_r)\equiv C(g_r:t)$ \text{\rm
(Multiseq)}}. For notation, $D=D_{g_r}$ is defined by the curve
which is parametrized by $D(t)=D_{g_r}(t)=C(g_r:t)$ if necessary.
\ms

$\underline{\text{\bf Fact(A-1)}}$ There is nothing to prove for
Fact(A-1).

$\underline{\text{\bf Fact(A-2)}}$ Then, the proof will be by
induction on the multiplicity $n=\prod^r_{k=1}n_k$ of $g_r$ at
$(y,z)=(0,0)$ that \text{$V(g_r)\equiv C(g_r:t)$ \text{\rm
(Multiseq)}} where $1\le n_1<\beta_{1,1}$ and $n_i\ge 2$ for $2\le
i\le r$. For notation, we write $D(t)=D_{g_r}(t)=C(g_r:t)$ if
necessary.

For the induction proof of Fact(A-2), it suffices to consider the
following two cases:

{\rm Case[I]} $n=2$ and {\rm Case[II]} $n\ge 3$. \ms

Note by (11.2.1) that $\gcd(n,\alpha_1)=n_2n_3\cdots n_r\le
n=n_1n_2\cdots n_r$ with $\gcd(n_1,\beta_{1,1})=1$. So, for each of
Case[I] and Case[II], it suffices to consider two subcases,
respectively:
$$\align
 \text{(a) $\gcd(n,\alpha_1)<n$ \quad and \quad (b)
$\gcd(n,\alpha_1)=n$}. \tag 11.2.4
\endalign$$

$\underline{\text{\bf Case[I] of Fact(A-2)}}$ \quad Let $n=2$. By
(11.2.4), there are two subcases:

Case[Ia] $1=\gcd(n,\alpha_1)<n=2$ \quad and \quad Case[Ib]
$\gcd(n,\alpha_1)=n=2$. \ms

$\underline{\text{\bf Case[Ia] of Fact(A-2)}}$ \quad Let $n=2$ with
$d_1=\gcd(n,\alpha_1)=1$. Then, it is trivial that the local
defining equation $g_1$ for an analytic variety $V(g_1)$ at the
origin is analytically defined by  $g_1(y,z)=z^2+y^{\beta_{1,1}}$
where $2<\beta_{1,1}=\alpha_1$ and $\gcd(2,\alpha_1)=1$.

On the other hand, since $\gcd(2,\alpha_1)=1$ then the irreducible
parametrization of $C(g_1:t)$ for the above $g_1$ can be written by
$y=t^2$ and $z=t^{\alpha_1}$ with $\alpha_1=\beta_{1,1}$.

Then, it is clear that \text{$V(g_1)\equiv C(g_1:t)$ \text{\rm
(Multiseq)}}. Thus, the proof of Case[Ia] is done. \ms

$\underline{\text{\bf Case[Ib] of Fact(A-2)}}$ \quad Let $n=2$ with
$d_1=\gcd(n,\alpha_1)=n=2$. That, is, $d_1=\gcd(2,\alpha_1)=2$.

Then, it is trivial that $g_2(y,z)$ of Sequences[I] with five
conditions is written in the form
$$\align
(11.2.5) \qquad &g_2=(z^{n_1}+y^{\beta_{1,1}})^{n_2}+\ve
y^{\beta_{2,1}}z^{\beta{2,2}}
 =(z+y^{\beta_{1,1}})^2+\ve y^{\beta_{2,1}}z^{\beta{2,2}} \quad
 \text{with}  \qquad \qquad \qquad \\
 &\Delta_2(\beta_{2,1},\beta_{2,2})=n_1\beta_{2,1}
 +\beta_{1,1}\beta_{2,2}=\beta_{2,1}
 +\beta_{1,1}\beta_{2,2}>n_2n_1\beta_{1,1}=2\beta_{1,1}, \quad \text{and} \\
 &\gcd(n_2,\Delta_2(\beta_{2,1},\beta_{2,2}))
 =\gcd(2,\Delta_2(\beta_{2,1},\beta_{2,2}))=1,
\endalign$$
where $\Delta_2(t_1,t_2)$ is defined by
$n_1t_1+\beta_{1,1}t_2=t_1+\beta_{1,1}t_2$ and $\ve$ is a unit in
$\BC\{y,z\}$.

By Theorem $10.3$ or by a nonsingular change of coordinates at the
origin, it can be easily proved that \text{$V(g_2)\equiv V(h)$
\text{\rm (Multiseq)}} where
$h(y,z)=z^2+y^{\Delta_2(\beta_{2,1},\beta_{2,2})}$.

On the other hand, since $\alpha_1=2\beta_{1,1}$ by (11.2.5) and
$\alpha_2-\alpha_1=\widehat{\Delta}_2(\beta_{2,1},\beta_{2,2})
=\Delta_2(\beta_{2,1},\beta_{2,2})-n_2n_1\beta_{1,1}
=\Delta_2(\beta_{2,1},\beta_{2,2})-2\beta_{1,1}
=\Delta_2(\beta_{2,1},\beta_{2,2})-\alpha_1$ by (11.2.1), then
$C(g_2:t)$ of (11.2.1) can be represented by $y=t^2$ and
$z=t^{2\beta_{1,1}}+t^{2\beta_{1,1}+
\widehat{\Delta}_2(\beta_{2,1},\beta_{2,2})}$. Noting that
$2\beta_{1,1}+\widehat{\Delta}_2(\beta_{2,1},\beta_{2,2})=
\Delta_2(\beta_{2,1},\beta_{2,2})$, let $C^*$ be defined by $y=t^2$
and $z=t^{\Delta_2(\beta_{2,1},\beta_{2,2})}$.

Then, it is clear that \text{$C^*\equiv C(g_2:t)$} \text{\rm
(Multiseq)}, \text{$C^*\equiv V(h)$} \text{\rm (Multiseq)} and
\text{$V(g_2)\equiv V(h)$} \text{\rm (Multiseq)}. So,
\text{$V(g_2)\equiv C(g_2:t)$} \text{\rm (Multiseq)}. Thus, we
proved that Case[Ib] is true, and therefore the proof of Case[I] is
done. \ms

$\underline{\text{\bf Case[II] of Fact(A-2)}}$ \quad Let $n\ge 3$.
For the induction proof, suppose we have shown by {\rm Case[I]} that
Fact(A-2) is true when $g_r$ has a multiplicity at the origin which
is either less than $n$ or equal to two. Now, assuming that $g_r$
has a multiplicity $n$ at the origin, then it is enough to consider
two such subcases, respectively:

{\rm Case[IIa]} $\gcd(n,\alpha_1)<n(n_1>1)$. \quad {\rm Case[IIb]}
$\gcd(n,\alpha_1)=n(n_1=1)$. \ms

$\underline{\text{\bf Case[IIa] of Fact(A-2)}}$ \quad Let
$\gcd(n,\alpha_1)<n$. To avoid the complexity of notations for the
proof, we may assume that $V(g_r)$ and $D(t)=C(g_r:t)$ of Fact(A-1)
satisfy the corresponding properties, respectively in the assumption
of this theorem up to the same notations.

Since $\gcd(n,\alpha_1)<n$ with $2\le n<\alpha_1$, let $q$ be the
positive integer such that
$$\align
(11.2.6)\qquad \qquad qn_1<\beta_{1,1}<(q+1)n_1, \quad \text{that
is}, \quad qn<\alpha_1<(q+1)n, \qquad \qquad \qquad \qquad
\endalign$$
where $q\ge 1$, $d_1=\gcd(n,\alpha_1)=n_2n_3\cdots n_r$, $n=n_1d_1$
and $\alpha_1=\beta_{1,1}d_1$.

In preparation for the proof of {\rm Case[IIa]}, it suffices to
prove the following two sublemmas, Sublemma 11.2.1 for Case[IIa] and
Sublemma 11.2.2 for Case[IIa] respectively, because of the
following:

{\rm(i)} Firstly, we show by Sublemma $11.2.1$ for Case[IIa] that we
find the representation of $V^{(q)}(g_r)$ and $D^{(q)}(t)$ in
Sublemma $11.2.1$ for Case[IIa], using the same method as we have
done in Lemma $4.2$ where the {\rm q}-th proper transforms of
$V(g_r)$ and $D(t)=C(g_r:t)$ are denoted by $V^{(q)}(g_r)$ and
$D^{(q)}(t)$, respectively. In addition, we compute by Sublemma
$11.2.1$ for Case[IIa] that if we write $\text{\rm
Multiseq($V(g_r)$)}=\{a_1,a_2,\dots,a_m\}$ and $\text{\rm
Multiseq($D(t)$)}=\{b_1,b_2,\dots,b_s\}$ then $a_i=b_i=n$ for
$i=1,2,\dots,q$. \ms

{\rm(ii)} Secondly, using Sublemma $11.2.1$ for Case[IIa], it is
enough to show by Sublemma $11.2.2$ for Case[IIa] that
\text{$V^{(q)}(g_r)\equiv D^{(q)}(t)$} \text{\rm (Multiseq)}. \ms

{\rm(iii)} After then, Sublemma $11.2.1$ and Sublemma $11.2.2$ for
Case[IIa] imply by {\rm(i)} and Corollary $3.8$ that
$\{\text{[Mult$(V(g_r))$]}\}=\text{\rm {Join}}(\{[n:qn]\},
\{\text{[Mult$(V^{(q)}(g_r))$]}\})$ and
$\{\text{[Mult$(D(t))$]}\}=\text{\rm {Join}}(\{\{[n:qn]\},
\{\text{[Mult$(D^{(q)}(t))$]}\})$ are equal. \ms

$\underline{\text{\bf Sublemma 11.2.1 for Case[IIa] of Fact(A-2)}}$
Claim the following:

(a) Whenever we use the composition of $q$ iterations of blow-ups,
denoted by $\tau_q=\pi_1\circ\pi_2\circ\cdots \circ\pi_q:M^{(q)}\to
\BC^2$ where $q$ is a positive integer defined by {\rm (11.2.6)},
then it can be easily proved by Lemma $4.2$ that just one of the
local coordinates for each blow-up $\pi_i$ with $1\le i\le q$ is
needed for the study of the $i-th$ proper transforms of both
$V(g_r)$ and $D(t)$ simultaneously under
$\tau_i=\pi_1\circ\pi_2\circ\cdots \circ\pi_i$ because $V(g_r)$ and
$D(t)$ are irreducible in $\BC\{y,z\}$. Note that the $i-th$ proper
transforms of $V(g_r)$ and $D(t)=C(g_r:t)$ are denoted by
$V^{(i)}(g_r)$ and $D^{(i)}(t)$ respectively. Then by Lemma $4.2$,
$\tau_i$ can be represented, as a composition of analytic mappings,
as follows:
$$\align
 \tau_i(v_i,u_i)=(y,z)=(v_i,{v_i}^iu_i), \tag 11.2.7
\endalign$$
where for brevity of notation, $(v_i,u_i)$ is called one coordinate
patch of the local coordinates for the $i$-th blow-up
$\pi_i:M^{(i)}\to M^{(i-1)}$ at some point of $M^{(i-1)}$ with
$M^{0}=\BC^2$. If we write $\text{\rm
Multiseq($V(g_r)$)}=\{a_1,a_2,\dots,a_m\}$ and $\text{\rm
Multiseq($D(t)$)}=\{b_1,b_2,\dots,b_s\}$, then $a_i=b_i=n$ for
$i=1,2,\dots,q$.

{\rm(b)} For the proof, first construct the local defining equation
for the proper transform $D^{(q)}(t)$ of $D(t)$ under $\tau_q$, and
next construct the local defining equation for the proper transform
$V^{(q)}(g_r)$ of $V(g_r)$ under $\tau_q$, denoted by
$(g_r\circ\tau_q)_{proper}$, as follows:

{\rm(b1)} The local defining equation for the curve $D^{(q)}$ is as
follows:
$$\align
\text{$D^{(q)}(t):=$} \left\{\eqalign{ v=& t^n, \cr u=&
t^{\alpha_1-qn}+t^{\alpha_2-qn}+\cdots
             +t^{\alpha_r-qn}. \cr} \right. \tag 11.2.8 \\
\endalign$$

Since $0<\alpha_1-qn<n$, the parametrization in (11.2.8) is not the
Puiseux expansion.

Let $\Omega=\Omega(s)$ be another irreducible curve with isolated
singularity at $0\in\BC$, which is defined by the Puiseux expansion
with a parameter $s$, as follows:
$$\align
\text{$\Omega(s):=$} \left\{\eqalign{v=&
s^n+s^{n+\alpha_2-\alpha_1}+\cdots +s^{n+\alpha_r-\alpha_1},
\cr u=& s^{\alpha_1-qn}, \cr} \right. \tag 11.2.9 \\
\endalign$$
because $0<\alpha_1-qn<n$ and $0<\alpha_1<\alpha_2<\cdots<\alpha_r$.

Then, it is clear by Theorem $8.8$ and Theorem $8.10$ that
\text{$D^{(q)}(t)\equiv \Omega(t)$ \ \text{\rm (Multiseq)}}. \ms

{\rm(b2)} For each $t=1,2,\dots,q$, along $E_t=\{v_t=0\}$ the local
defining equation $(g_r\circ\tau_t)_{total}$ for the $t-th$ total
transform of $V(g_r)$ is as follows:
$$\align
(11.2.10) \quad &(g_j\circ\tau_t)_{total}=v^{tn_1n_2\cdots
n_j}_t(g_j\circ\tau_t)_{proper} \quad \text{for $1\le j\le r$,} \\
&(g_1\circ\tau_t)_{proper}=u^{n_1}_t+\ve_1v^{\beta_{1,1}-tn_1}_t, \dots,\\
  &(g_j\circ\tau_t)_{proper}=(g_{j-1}\circ\tau_t)^{n_j}_{proper}+\ve_ju^{\beta_{j,2}}_t
  v^{\omega_{t,j}}_t\prod^{j}_{k=3}
  (g_{k-2}\circ\tau_t)^{\beta_{j,k}}_{proper} \\
 &  \text{with
   $\omega_{t,j}=\beta_{j,1}+t\beta_{j,2}+tn_1\beta_{j,3}+\cdots
  +tn_1n_2\cdots n_{j-2}\beta_{j,j}-tn_1n_2\cdots n_j>0$}, \qquad
\endalign$$
such that $\omega_{t,j}>0$ for all $t=1,2,\dots,q$ and all
$j=2,3,\dots,r$ because $\beta_{1,1}-n_1q>0$. Note that each
$\ve_i=\ve_i(u_t,v_t)$ is a unit in $\BC\{u_t,v_t\}$ for $1\le i\le
r$.

{\rm(b3)} Moreover, the multiplicity of $(g_r\circ\tau_q)_{proper}$
at $(v_q,u_q)=(0,0)$ is $\alpha_1-qn$, which is less than the
multiplicity $n$ of $g_r$ at $(y,z)=(0,0)$.  $\square$ \ms

$\underline{\text{\bf Proof of Sublemma 11.2.1 for Case[IIa]}}$ As
some applications of Theorem $3.6$, Sublemma $5.1$, Sublemma $5.2$,
Sublemma $5.3$, Sublemma $5.4$, Theorem $8.8$ and Theorem $8.10$,
the proof of this sublemma can be easily done. $\square$ \ms

$\underline{\text{\bf Sublemma 11.2.2 for Case[IIa] of Fact(A-2) }}$
Claim the following:

For the proof of \text{$V^{(q)}(g_r)\equiv \Omega(s)$ \ \text{\rm
(Multiseq)}}, since the multiplicity of $(g_r\circ\tau_q)_{proper}$
at $(v_q,u_q)=(0,0)$ is $\alpha_1-qn<n$ by Sublemma 11.2.1 for
Case[IIa], it suffices to show by the induction method that the
following two facts are true:

$\underline{\text{\rm Fact(i) of Sublemma 11.2.2}}$ $V^{(q)}(g_r)$
satisfies the same kind of properties as $V(g_r)$ does in the
assumption of Fact(A-2).

$\underline{\text{\rm Fact(ii) of Sublemma 11.2.2}}$ $\Omega(s)$ of
(11.2.9) satisfies the same kind of properties relative to
$V^{(q)}(g_r)$ as $D(t)$ of (11.2.1) does relative to $V(g_r)$ in
the assumption of Fact(A-2).

In order to prove that Fact(i) and Fact(ii) are true, we will
represent Fact(i) and Fact(ii) in more detail.

For notation, rewrite $(g_1\circ\tau_q)_{proper}$,
$(g_2\circ\tau_q)_{proper}$,\dots,$(g_r\circ\tau_q)_{proper}$ of
(11.2.10) by $h_1,h_2,\dots,h_r$, respectively in $\C\{v,u\}$ as
follows: Note that $(v,u)=(v_q,u_q)$.
$$\align
(11.2.11) \quad  & h_1 =\ve_1v^{m_1}+u^{\gamma_{1,1}}\quad
{where} \quad 1\le m_1 =\beta_{1,1}-n_1q<\gamma_{1,1}=n_1, \\
& h_j =h_{j-1}^{m_j}+\ve_ju^{\gamma_{j,1}}v^{\gamma_{j,2}}
h_1^{\gamma_{j,3}}\cdots h_{j-2}^{\gamma_{j,j}} \quad \text{for $j=2,3,\dots,r,$} \\
{where} & \quad 2\le m_j =n_j, \gamma_{j,1}=\beta_{j,2}, \\
& \quad
\gamma_{j,2}=\omega_{q,j}=\beta_{j,1}+q\beta_{j,2}+qn_1\beta_{j,3}+\cdots
+qn_1n_2\cdots n_{j-2}\beta_{j,j}-qn_1n_2\cdots n_j>0, \\
& \quad \gamma_{j,3}=\beta_{j,3},\dots,~\gamma_{j,j}=\beta_{j,j}.
\endalign$$
Note that $0<m_1=\beta_{1,1}-n_1q<n_1$ and each $\ve_i=\ve_i(u,v)$
is a unit in $\BC\{u,v\}$ for $1\le i\le r$. \ms

Then, apply the notations in (11.2.11) to Fact(i) and Fact(ii), as
follows:

$\underline{\text{\bf Fact(i) of Sublemma 11.2.2 for Case[IIa] of
Fact(A-2)}}$ Applying the same kind of notations as we have used for
$\{X_k:k=1,2,\dots,r\}$ in the assumption of this theorem to the
proof of Fact(i) of this sublemma,  we claim that Fact(i) can be
rewritten as follows:
$$\align
&  \text{\bf Sequences[I]}^{\text{{\bf(1)}}} \text{\bf of the r-th
type}\text{\bf :} \quad \text{Let $\{Y_k: k=1,2,\dots,r\}$}
\quad \text{with $Y_k\subset N_0$}, \\
& \{h_k: k=1,2,\dots,r\} \quad \text{with
$h_k=(g_{k}\circ\tau_q)_{proper}$ in $\BC\{u,v\}$,} \qquad \qquad \\
&\{\text{$\Xi_k:N^k_0\to N_0$ is an integer-valued
function for $k=1,2,\dots,r$}\} \\
& \text{be three sequences, satisfying the following five conditions
for each k.} \qquad \qquad
\endalign$$

\noindent Five conditions are denoted by \text{\bf The 1st
${\text{\bf{Cond}}}^{\text{{\bf(1)}}}$}{\bf, \dots,} \text{\bf The
5-th ${\text{\bf{Cond}}}^{\text{{\bf(1)}}}$} of $\text{\bf
Sequences[I]}^{\text{{\bf(1)}}}$.

\text{\bf [I]-(1)} \quad \text{\bf The 1st
${\text{\bf{Cond}}}^{\text{{\bf(1)}}}$} $\text{\bf
 of Sequences[I]}^{\text{{\bf(1)}}}${\bf :}

{\rm (1a)} $Y_1=\{m_1,\gamma_{1,1}\}$ with $1\le
m_1=\beta_{1,1}-n_1q<\gamma_{1,1}=n_1$.

{\rm(1b)} $Y_j=\{m_j,\gamma_{j,1},\gamma_{j,2},\dots,\gamma_{j,j}\}$
with $m_j\ge 2$, where $j=2,\dots,r$. \ms

\text{\bf [I]-(2)} \quad \text{\bf The 2nd
${\text{\bf{Cond}}}^{\text{{\bf(1)}}}$} $\text{\bf
 of Sequences[I]}^{\text{{\bf(1)}}}$ {\bf :}

{\rm(2a)} $h_1 =v^{m_1}+\ve_1u^{\gamma_{1,1}}.$

{\rm(2b)} $h_j =h_{j-1}^{m_j}+\ve_ju^{\gamma_{j,1}}v^{\gamma_{j,2}}
h_1^{\gamma_{j,3}}\cdots h_{j-2}^{\gamma_{j,j}}$ for $j=2,\dots,r$.

Note that each $\ve_j=\ve_j(u,v)$ is a unit in $\BC\{u,v\}$ for
$1\le j\le r$. \ms

\text{\bf [I]-(3)} \quad \text{\bf The 3rd
${\text{\bf{Cond}}}^{\text{{\bf(1)}}}$} $\text{\bf of
Sequences[I]}^{\text{{\bf(1)}}}$ {\bf :}

{\rm (3a)} $\Xi_1(t)=t$ for each $t\in N_0$.

{\rm (3b)} $\Xi_j(t_k)^j_{k=1}=t_j\Xi_{j-1}(\gamma_{j-1,k})^j_{k=1}
 +m_{j-1}\Xi_{j-1}(t_k)^{j-1}_{k=1}$
for each $(t_k)^j_{k=1}\in N^j_0$. \ms

\text{\bf [I]-(4)} \quad \text{\bf The 4-th
${\text{\bf{Cond}}}^{\text{{\bf(1)}}}$} $\text{\bf
 of Sequences[I]}^{\text{{\bf(1)}}}$ {\bf :}

{\rm (4a)} $\Xi_1(\gamma_{1,1})=\gamma_{1,1}>m_1\ge 1$.

 {\rm (4b)} $\Xi_j(\gamma_{j,k})^j_{k=1}
 >m_jm_{j-1}\Xi_{j-1}(\gamma_{j-1,k})^{j-1}_{k=1}$
for $j=2,\dots,r$. \ms

\text{\bf [I]-(5)} \quad \text{\bf The 5-th
${\text{\bf{Cond}}}^{\text{{\bf(1)}}}$} $\text{\bf
 of Sequences[I]}^{\text{{\bf(1)}}}$ {\bf :}

{\rm (5)} $\gcd(m_j,\Xi_j(\gamma_{j,k})^j_{k=1})=1$ for
 $j=1,\dots,r$. \ms

$\underline{\text{\bf Fact(ii) of Sublemma 11.2.2 for Case[IIa] of
Fact(A-2)}}$ For the given $V(h_r)=\{(u,v):h_r(u,v)=0\}$ in Fact(i),
the Puiseux expansion $\Omega(s)$ at $s=0$ for the curve $\Omega$ is
given as follows:
$$\align
\text{$\Omega(s):=$} \left\{\eqalign{v=&
s^n+s^{n+\alpha_2-\alpha_1}+\cdots +s^{n+\alpha_r-\alpha_1},
\cr u=& s^{\alpha_1-qn}, \cr} \right. \tag 11.2.12 \\
\endalign$$
where $0<\alpha_1-qn<n$ and $0<\alpha_1<\alpha_2<\cdots<\alpha_r$.

As we have seen in $(11.2.1)$, $\Omega(s)$ of (11.2.12) satisfies
the following property:
$$\align
(11.2.13) \qquad  & m_1m_2\cdots m_r=\alpha_1-qn,  \quad
 \gamma_{1,1}m_2\cdots m_r=n, \quad \text{and} \quad
 \text{for $j=1,\dots,r,$} \qquad \qquad\\
&(n+\alpha_j-\alpha_1)-(n+\alpha_{j-1}-\alpha_1)
 =\widehat{\Xi}_j(\gamma_{j,k})^j_{k=1}m_{j+1}m_{j+2}\cdots m_r. \\
\endalign$$
Note that
$\alpha_j-\alpha_{j-1}=(n+\alpha_j-\alpha_1)-(n+\alpha_{j-1}-\alpha_1)$
for $j=1,\dots,r$.  $\square$ \ms

$\underline{\text{\bf Proof of Sublemma 11.2.2 for Case[IIa] of
Fact(A-2)}}$ For the proof of Sublemma 11.2.2 for Case[IIa] of
Fact(A-2), It suffices to prove Fact(i) of Sublemma 11.2.2 for
Case[IIa] of Fact(A-2) and Fact(ii) of Sublemma 11.2.2 for Case[IIa]
of Fact(A-2), respectively. \ms

$\underline{\text{\bf Proof of Fact(i) of Sublemma 11.2.2 for
Case[IIa] of Fact(A-2)}}$ For the proof of Fact(i), it suffices to
show that the properties in \text{\rm The 4-th
${\text{\rm{Cond}}}^{\text{{\rm(1)}}}$}  and \text{\rm The 5-th
${\text{\rm{Cond}}}^{\text{{\rm(1)}}}$} of $\text{\rm
Sequences[I]}^{\text{{\rm(1)}}}$ are true.

For the proofs of the properties in \text{\rm The 4-th
${\text{\rm{Cond}}}^{\text{{\rm(1)}}}$}  and \text{\rm The 5-th
${\text{\rm{Cond}}}^{\text{{\rm(1)}}}$}, it is trivial by (11.2.11)
that $\Xi_1(\gamma_{1,1})=\gamma_{1,1}>m_1\ge 1$ and
$1=\gcd(n_1,\beta_{1,1})=\gcd(n_1,\beta_{1,1}-qn_1)=
\gcd(m_1,\gamma_{1,1})$, and so it suffices to show that the
following assertion is true. \ms

$\underline{\text{Assertion}}$ \quad Let $w$ and $\ell$ be arbitrary
positive integers such that $2\le w\le \ell\le r$. Then, we have the
following inequalities:
$$ \align
(11.2.14) \quad \text{(i)} \quad &
 \Xi_w(\gamma_{\ell,k})^w_{k=1}=\Delta_w(\beta_{\ell,k})^w_{k=1}+qn^2_1\cdots
 n^2_{w-1}(\beta_{\ell,w+1}+n_w\beta_{\ell,w+2}+\cdots  \qquad \qquad \\
 & +n_wn_{w+1}\cdots
 n_{\ell-2}\beta_{\ell,\ell}-n_wn_{w+1}\cdots n_{\ell}). \\
 \text{(ii)}\quad &
 \Xi_w(\gamma_{w,k})^w_{k=1}=\Delta_w(\beta_{w,k})^w_{k=1}-qn^2_1\cdots
 n^2_{w-1}n_w. \\
 \text{(iii)}\quad & \widehat{\Xi}_w(\gamma_{w,k})^w_{k=1}
  =\Xi_w(\gamma_{w,k})^w_{k=1}-m_wm_{w-1}\Xi_{w-1}
 (\gamma_{w-1,k})^{w-1}_{k=1} \\
 & =\Delta_w(\beta_{w,k})^w_{k=1}
 -n_wn_{w-1}\Delta_{w-1}(\beta_{w-1,k})^{w-1}_{k=1}
  =\widehat{\Delta}_w(\beta_{w,k})^w_{k=1}>0.  \\
 \text{(iv)} \quad & \gcd(\Xi_w(\gamma_{w,k})^w_{k=1},m_w)=1.
\endalign$$

So, it is clear that the inequalities in (iii) and (iv) of (11.2.14)
imply the properties in \text{\rm The 4-th
${\text{\rm{Cond}}}^{\text{{\rm(1)}}}$}  and \text{\rm The 5-th
${\text{\rm{Cond}}}^{\text{{\rm(1)}}}$}. \ms

$\underline{\text{Proof of the Assertion}}$: The proof will be by
induction on the integer $w\ge 1$. So, for the proof of this
assertion, it suffices to consider two subcases, respectively:

{\rm  Subcase(1)}  $w=1$ and {\rm Subcase(2)}  $w\ge 2$. \ms

$\underline{\text{Subcase(1) for the Assertion}}$ If $w=1$, then,
there is nothing to prove. \ms

$\underline{\text{Subcase(2) for the Assertion}}$ Let \text{$w\ge
2$}. Next, by using the induction assumption on the integer $w\ge
1$, suppose we have shown that the above assertion is true on the
integer $w$ with $r-1\ge \ell\ge w\ge 1$. Assuming that $r\ge
\ell\ge w+1> 2$, then it suffices to show that the proof of four
equalities in (11.2.14) of the assertion just follows from the
proofs of (a), (b), (c) and (d):

(a) Compute $\Xi_{w+1}(\gamma_{\ell,{k}})^{w+1}_{k=1})$ by (11.2.11)
and \text{\rm The 3-th ${\text{\rm{Cond}}}^{\text{{\rm(1)}}}$}, as
follows:
$$\align
(11.2.15)\quad  &\Xi_{w+1}(\gamma_{\ell,{k}})^{w+1}_{k=1}
=\gamma_{\ell,w+1}\Xi_w(\gamma_{w,k})^w_{k=1}+m_w\Xi_w(\gamma_{\ell,{k}})^w_{k=1}  \\
 =&\beta_{\ell,w+1}(\Delta_w(\beta_{w,k})^r_{k=1}-qn^2_1\cdots
n^2_{w-1}n_w)+n_w \{\Delta_w(\beta_{\ell,{k}})^w_{k=1}+ qn^2_1\cdots n^2_{w-1}\cdot  \\
& (\beta_{\ell,w+1}+n_w\beta_{\ell,w+2}+n_w
n_{w+1}\beta_{\ell,w+3}+\cdots
  +n_w n_{w+1}\cdots n_{\ell-2}\beta_{\ell,\ell}-n_w n_{w+1}\cdots n_\ell)\} \quad\\
=&\Delta_{w+1}(\beta_{\ell,{k}})^{w+1}_{k=1}+qn^2_1\cdots
n^2_{w-1}n^2_w(\beta_{\ell,w+2}+n_{w+1}\beta_{\ell,w+3}+\cdots \\
& +n_{w+1}n_{w+2}\cdots
n_{\ell-2}\beta_{\ell,\ell}-n_{w+1}n_{w+2}\cdots n_\ell),
\endalign$$
which gives the proof of {\rm (i)} of (11.2.14). \ms

(b) From $(11.2.15)$ in $(a)$, if $\ell=w+1$, then
 $$\align
(11.2.16) \qquad \qquad \Xi_{w+1}(\gamma_{w+1,k})^{w+1}_{k=1}
=\Delta_{w+1}(\beta_{w+1,k})^{w+1}_{k=1} -qn^2_1\cdots
n^2_{w-1}n^2_{w} n_{w+1},\qquad \qquad \qquad \\
\endalign$$
and so the proof of {\rm (ii)} of (11.2.14) is done. \ms

(c) By $(b)$ and (11.2.11), we have
$$\align
(11.2.17) \qquad &\widehat{\Xi}_{w+1}(\gamma_{w+1,k})^{w+1}_{k=1}
=\Xi_{w+1}(\gamma_{w+1,k})^{w+1}_{k=1}-m_{w+1}m_w\Xi_w(\gamma_{w,k})^w_{k=1}
\qquad \qquad \\
 =&\Delta_{w+1}(\beta_{w+1,k})^{w+1}_{k=1}-qn^2_1\cdots n^2_wn_{w+1}
-n_{w+1}n_w\{\Delta_w(\beta_{w,k})^w_{k=1}
 -qn^2_1\cdots n^2_{w-1}n_w\} \\
=&\Delta_{w+1}(\beta_{w+1,k})^{w+1}_{k=1}-n_{w+1}n_w\Delta_w(\beta_{w,k})^w_{k=1}
=\widehat{\Delta}_{w+1}(\beta_{w+1,k})^{w+1}_{k=1}>0. \\
\endalign$$
Thus, the proof of {\rm (iii)} of (11.2.14) is done. \ms

(d) By (c) and (11.2.11) again, we have
$$\align
(11.2.18) \qquad
&\gcd(\Xi_{w+1}(\gamma_{w+1,k})^{w+1}_{k=1},m_{w+1})
=\gcd(\widehat{\Xi}_{w+1}(\gamma_{w+1,k})^{w+1}_{k=1},m_{w+1})  \qquad \qquad\qquad \\
=&\gcd(\widehat{\Delta}_{w+1}(\beta_{w+1,k})^{w+1}_{k=1},n_{w+1})
=\gcd(\Delta_{w+1}(\beta_{w+1,k})^{w+1}_{k=1},n_{w+1})=1. \\
\endalign$$
Thus, the proof of {\rm (iv)} of (11.2.14) is done, and so we proved
that {\rm Subcase(2)} is true. Therefore, the proof of the assertion
is done, and so we finished the proof of Fact(i) of Sublemma
$11.2.2$. \ms

$\underline{\text{\bf Proof of Fact(ii) of Sublemma 11.2.2 for
Case[IIa] of Fact(A-2)}}$ It suffices to show that $\Omega(s)$ of
(11.2.12) satisfies the equalities of (11.2.13). By use of (11.2.11)
and (11.2.14), and (11.2.1) of the assumption, the proof is as
follows:
$$\align
(11.2.19)\qquad \text{\rm(i)} \quad m_1m_2\cdots m_r
&=(\beta_{1,1}-qn_1)n_2n_3\cdots n_r   \\
&=\beta_{1,1}n_2n_3\cdots n_r-qn_1n_2\cdots n_r=\alpha_1-qn, \\
\text{\rm(ii)} \quad \gamma_{1,1}m_2\cdots m_r&=n_1n_2n_3\cdots n_r=n, \\
\text{\rm(iii)} \quad \quad  \alpha_{j+1}-\alpha_j
&=\widehat{\Delta}_{j+1}(\beta_{j+1,k})^{j+1}_{k=1}n_{j+2}n_{j+3}\cdots
n_r \quad \text{for $1\le j\le r-1$} \qquad \qquad\\
&=\widehat{\Xi}_{j+1}(\gamma_{j+1,k})^{j+1}_{k=1}m_{j+2}m_{j+3}\cdots
m_r.
\endalign$$

Thus, the proof of Fact(ii) of Sublemma $11.2.2$ for Case[IIa] of
Fact(A-2) is done, and so we finished the proof of Sublemma $11.2.2$
for Case[IIa] of Fact(A-2).  $\square$

Therefore, we can show by the proofs of Sublemma $11.2.1$ and
Sublemma $11.2.2$ that {\rm Case[IIa]} of Fact(A-2) is true.
$\square$ \ms

$\underline{\text{\bf Case[IIb] of Fact(A-2)}}$ \quad Let
$\gcd(n,\alpha_1)=n$. Then, $n=n_1d_1$ and $\alpha_1=\beta_{1,1}d_1$
with $d_1=\gcd(n,\alpha_1)=n>1$, and so $d_1=n_2n_3\cdots n_r$ by
(11.2.2), noting that $n_1=1<\beta_{1,1}$.

Since $1=n_1<\beta_{1,1}$ in Case[IIb], then the representation of
$g_r$ in the assumption of Case[IIb] and the representation of $g_r$
in the assumptions of Sequences[I] of Theorem $10.3$ are the same.
Since $g_r$ is irreducible in $\BC\{y,z\}$, then by the assumption
of Theorem $10.3$ we can use the same results and notations as in
the conclusions of Theorem $10.3$:
$$\align
\text{\rm(11.2.20)}  \qquad \qquad
&\text{$\omega_j(\sigma_{j,k})^j_{k=1}=\Delta_{j+1}(\beta_{j+1,k})^{j+1}_{k=1}$
for each $j=1,2,\dots,r-1$,} \qquad \qquad \qquad \qquad \\
(11.2.21)  \qquad \qquad &\widehat{\omega}_q(\sigma_{q,k})^q_{k=1}
=\omega_q(\sigma_{q,k})^q_{k=1}-c_qc_{q-1}\omega_{q-1}(\sigma_{q-1,k})^{q-1}_{k=1}
\qquad \qquad \qquad \qquad \qquad \\
\qquad &=\Delta_{q+1}(\beta_{q+1,k})^{q+1}_{k=1}
-n_{q+1}n_q\Delta_q(\beta_{q,k})^q_{k=1} \quad \text{for each} \quad {q=2,3,\dots,r-1},\\
\qquad &=\widehat{\Delta}_{q+1}(\beta_{q+1,k})^{q+1}_{k=1}>0, \\
\text{\rm(11.2.22)} \qquad \qquad &\text{$g_{r}$ is irreducible in
$\BC\{y,z\}$ if and only if  $h_{r-1}$ is irreducible in $\BC\{y,z\}$,} \\
\text{\rm(11.2.23)} \qquad \qquad &\text{$g_r$ and $h_{r-1}$ have
the same multiplicity sequence.}
\endalign$$

Noting by Theorem $10.3$ that $c_j=n_{j+1}$ for $1\le j\le r-1$, we
showed that $\gcd(n,\xi_1)<n$, $\gcd(n,\xi_{1})=1$ and $2\le
c_1<\sigma_{1,1}$ or $2\le n<\xi_1$ in Sequences[II] where
$n=c_1c_2\cdots c_{r-1}$, $\xi_1=\sigma_{1,1}c_2\cdots c_{r-1}$ and
$n$ is the multiplicity of $h_{r-1}$ at $(y,z)=(0,0)$, and so the
results of Case[IIa] can be applied to Case[IIb].\ms

Applying Case[IIa] to this case, then for a given
$\{(y,z):h_{r-1}(y,z)=0\}$ in {\rm [II]-(2)} of Sequences{\rm[II]}
of Theorem $10.3$, we can compute the Puiseux expansion
$C(h_{r-1}:t)$ for the curve such that  \text{$V(h_{r-1})\equiv
C(h_{r-1}:t)$ \text{\rm (Multiseq)}} where
$V(g_r)=\{(y,z):g_r(y,z)=0\}$ in the assumption, as follows:
$$\align
(11.2.24)\qquad \qquad \qquad  &\text{$C(h_{r-1}:t):=$}
\left\{\eqalign{ y&=t^c, \cr z&=t^{\xi_1}+t^{\xi_2}+\cdots
+t^{\xi_{r-1}}, \cr} \right.
 \\
 \qquad     \text{such that} & \qquad  c =c_1c_2\cdots
c_{r-1} \quad \text{and} \quad
\xi_1 =\sigma_{1,1}c_2c_3\cdots c_{r-1}, \qquad \qquad\\
& \qquad  \xi_j
=\xi_{j-1}+\widehat{\omega}_j(\sigma_{j,k})^j_{k=1}c_{j+1}c_{j+2}\cdots
c_{r-1} \quad \text{for $2\le i\le r-1$.} \qquad \qquad\\
\endalign$$

Now, in order to finish the proof of the subcase Case[IIb] of
Case[II], since $n$ is  not a divisor of $\alpha_2$ and $c=n$ is not
a divisor of $\xi_1$, then it remains to show by {\rm Case[IIa]}
that the following equalities are true:
$$\align
n =c \quad \text{and} \quad
\xi_j =\alpha_{j+1} \quad \text{for $j=1,2,\dots,r-1$}.\tag 11.2.25 \\
\endalign$$

For the proof, first compute $\{\xi_k:k=1,2,\dots,r-1\}$ as follows:

{\rm (i)} By (11.2.20) and (11.2.24), and by (11.2.1), we have
$$\align
(11.2.26)  \qquad \qquad \xi_1 &=\sigma_{1,1}c_2c_3\cdots c_{r-1}
=\Delta_2(\beta_{2,1},\beta_{2,2})n_3n_4\cdots n_r,  \qquad \qquad\\
\alpha_2
&=\alpha_1+\widehat{\Delta}_2(\beta_{2,1},\beta_{2,2})n_3n_4\cdots
n_r
\\ &=\beta_{1,1}n_2n_3\cdots n_r
+(\Delta_2(\beta_{2,1},\beta_{2,2})-n_2n_1\beta_{1,1})n_3n_4\cdots
n_r
\qquad \qquad \qquad\\
 &=\Delta_2(\beta_{2,1},\beta_{2,2})n_3n_4\cdots n_r,
\endalign$$
because $n_1=1$. So, $\xi_1=\alpha_2$. \ms

{\rm (ii)} By (11.2.20), (11.2.21) and (11.2.24), we have for each
$j=2,3,\dots,r$,
$$\align
\xi_{j-1}-\xi_{j-2}
&=\widehat{\omega}_{j-1}(\sigma_{j-1,k})^{j-1}_{k=1}c_jc_{j+1}\cdots
c_{r-1}  \tag 11.2.27 \\
&=\widehat{\Delta}_j(\beta_{j,k})^j_{k=1}n_{j+1}n_{j+2}\cdots
n_r =\alpha_j-\alpha_{j-1}. \\
\endalign$$

So, the proof of {\rm(11.2.25)} follows very easily from (11.2.26)
and {\rm (11.2.27)}. Thus, the proof of Case[IIb] is done, which
implies the completion for the proof of {\rm Case[II]}.

Therefore, we showed by Case[I] and Case[II] that Fact(A) is true.

$\underline{\text{\bf The proofs of Fact(B) and Fact(C)}}$ It is
trivial by Theorem $7.3$, Theorem 7.7, and Theorem 10.2 that Fact(B)
and Fact(C) are true.

Thus, we finished the proof of the theorem. $\square$
\enddemo \ms

The proof of Corollary 11.3 just follows from Theorem 11.2 with
Theorem 10.2 and Theorem 8.10. \bs

{\bf \S11.4. The second half of The 1st Algorithm(Theorem 11.4)} \ms

Now, in order to solve Problem[1-C] mentioned in \S7.0 of \S7, it
suffices to show how to use the Euclidean algorithm in {\rm (3)} of
Corollary $7.6$ only.

\proclaim{Theorem 11.4(Theorem 1.6:Algorithm for finding the unique
element of Family(1) corresponding to any given standard Puiseux
expansion of Family(2))} {\rm Theorem 11.4} has the same statement
as {\rm Theorem 1.6} does in $\S 1$. $\square$
\endproclaim

\demo{\bf Proof of Theorem 11.4} The proof of Theorem 11.4 follows
from Theorem 7.7 of $\S 7$. $\square$
\enddemo \ms

\noindent{\bf Example 11.4.1 for the second half of The 1st
algorithm in Theorem 11.4:} Example 11.4.1 and Example 1.6.3 of $\S
1$ are the same. \bs

{{\bf \S 11.5. An application of The 1st Algorithm }}

\proclaim{Theorem 11.5} $\underline{\text{\rm As an application of
The 1st Algorithm}}$, instead of using {\rm the 2nd, the 3rd and the
4-th algorithms}, we can show directly that there exists a
one-to-one correspondence between four families, {\rm Family(1),
Family(2), Family(3) and Family(4)} as we have seen in {\rm
Definition 1.2} and {\rm Definition 2.4} , which was already proved.
$\square$
\endproclaim \ms

\bs \vfill \pagebreak

\newpage

$$\align
\quad &\text{\bf \centerline{Part[C](Part[C1],\dots, Part[C4]}}\\
&\text{\bf A complete and explicit irreducibility algorithms for the
W-polys of two}\\
&\text{\bf complex variables with proofs and related topics in the
Puiseux expansions}
\endalign$$

{\bf Part[C1] The new terminology and notations in preparation
for finding irreducibility criterion of W-polys in
$\BC\{y,z\}$}

{\bf $\S$12. The generalized standard irreducible W-polys of the
recursive r-type {\indent}with Theorem 12.0} \ms

In order to succeed in the computation of Explicit algorithm, it is
very interesting and important to prove that we can define the new
terminology, irreducible Weierstrass polynomials of two complex
variables of the recursive type, called {\rm ``the generalized
standard Puiseux polynomial in $\BC[y,z]$ of the recursive type"}
throughout this paper, which will be well-defined. In preparation,
we need Definition $12.0.0$ and Theorem $12.0$. \ms

\definition{Definition 12.0.0} Let $N_0$ be the set of
nonnegative integers and $N^k_0$ be its $k$-dimensional copy, and
$N$ be the set of positive integers. Let $r$ be an arbitrary
positive integer. \ms

\noindent{\bf [A]} $g_r\in\BC\{y,z\}$ is called either
$\underline{\text{\rm a generalized semi-quasi-Puiseux germ of the
recursive r-type}}$  $\underline{\text{\rm or a generalized
semi-quasi Puiseux convergent power series of the recursive
r-type}}$ if there are sequences $\{X_k:k=1,2,\dots,r\}$ with
$X_k\subset N_0$, $\{g_k:k=1,2,\dots,r\}$ with $g_k\in\BC\{y,z\}$
and $\{\text{${{\Delta}_k}:N^k_0\to N_0$ is an integer-valued
function for}$ $\text{$k=1,2,\dots,r$}\}$ satisfying the following
$\underline{\text{\rm four conditions}}$:\ms

\noindent$\underline{\text{\bf Four conditions}}$ are denoted by
\text{\rm The 1-th ${\text{\rm{Cond}}}^{\text{{\rm(0)}}}$}, $\dots$,
\text{\rm The 4-th ${\text{\rm{Cond}}}^{\text{{\rm(0)}}}$} for
brevity. \ms

\noindent$\underline{\text{\rm The 1st
${\text{\rm{Cond}}}^{\text{{\rm(0)}}}$}}$ The family
$\{X_{\ell}:\ell=1,2,\dots,r\}$ with $X_{\ell}\subset N_0$ is as
follows: \ms

\noindent{\rm(1)(1a)} $X_1=\{n_1\}\cup \{\beta_{1,i,1}:0\le i<n_1
\}$ with $n_1\ge 2$ and $\beta_{1,0,1}\ge 1$ where $X_{1}\subset N$.

{\rm(1b)} $X_j=\{n_j\}\cup \{\beta_{j,i,1}:0\le i<n_j \} \cup
\{\beta_{j,i,2}:0\le i<n_j \}\cup \cdots \cup \{\beta_{j,i,j}:0\le
i<n_j \}$ with $n_j\ge 2$ where $j=2,3,\dots,r$.

For each $j=2,3,\dots,r$, assume that at least one of
$\beta_{j,0,1},\beta_{j,0,2},\dots,\beta_{j,0,j}$ is nonzero. \ms

\noindent $\underline{\text{\rm The 2nd
${\text{\rm{Cond}}}^{\text{{\rm(0)}}}$}}$ For each $j=1,2,\dots,r$,
let $g_j=g_j(y,z,c_j)$ be in $\BC\{y,z\}$ where all the $c_j$ are
complex numbers, each of which is defined by the following way: \ms

\noindent \text{\rm(2)(2a)} \quad
\text{$g_1=z^{n_1}+\ve_{1,0}y^{\beta_{1,0,1}}
+c_1\sum^{n_1-1}_{i=1}\ve_{1,i}y^{\beta_{1,i,1}}z^i$ with
$\ve_{1,0}=1$.}

\text{\rm(2b)} \quad \text{
$g_j=g^{n_j}_{j-1}+\ve_{j,0}y^{\beta_{j,0,1}}z^{\beta_{j,0,2}}
g^{\beta_{j,0,3}}_1g^{\beta_{j,0,4}}_2\cdots
g^{\beta_{j,0,j}}_{j-2}$}

\qquad \qquad \quad
\text{$+c_j\sum^{n_j-1}_{i=1}\ve_{j,i}y^{\beta_{j,i,1}}z^{\beta_{j,i,2}}
g^{\beta_{j,i,3}}_1\cdots g^{\beta_{j,i,j}}_{j-2}g^i_{j-1}$,} where
$j=2,3,\dots,r$. \ms

Note that each $\ve_{j,i}=\ve_{j,i}(y,z)$ is a unit in $\BC\{y,z\}$
for $1\le j\le r$ and $0\le i<n_j$, if exists. As far as analytic
equivalence of isolated plane curve singularities defined by all
$g_j$, $1\le j\le r$, is concerned, then we may assume that
$\ve_{1,0}$ is equal to one by a suitable nonsingular change of
coordinates at the origin in $\BC^2$. \ms

\noindent $\underline{\text{\rm The 3rd
${\text{\rm{Cond}}}^{\text{{\rm(0)}}}$}}$ \quad Let $\{\Delta_k:
N^k_0\to N_0: k=1,2,\dots,r\}$ be a sequence such that each
$\Delta_k$ is an integer-valued function defined by the following
way: \roster \item"(3) (3a)" $\Delta_1(t)=t$ for each $t\in N_0$.

\item"$\quad$ (3b)"
$\Delta_j(t_k)^j_{k=1}=t_j\Delta_{j-1}(\beta_{j-1,0,k})^{j-1}_{k=1}
+n_{j-1}\Delta_{j-1}(t_k)^{j-1}_{k=1}$ for each $(t_k)^j_{k=1}\in
N^j_0$.
\endroster \ms

\noindent $\underline{\text{\rm The 4-th
${\text{\rm{Cond}}}^{\text{{\rm(0)}}}$}}$ \quad Then, the following
inequalities hold: Note that $r\ge 2$.
$$\align
\text{\rm (4)(4a)}& \quad \Delta_1(\beta_{1,i,1})={\beta_{1,i,1}}>0
\quad \text {for $0\le i<n_1$.}   \qquad \qquad \\
\text{\rm(4b)}& \quad {\Delta_{j}(\beta_{j,i,k})^{j}_{k=1}}
> {(n_{j}-i)}n_{j-1}\Delta_{j-1}(\beta_{j-1,0,k})^{j-1}_{k=1} \quad
\text {for $0\le i<n_j$ where $j=2,\dots,r$}. \qquad \qquad \qquad
\qquad
\endalign$$ \ms

\noindent{\bf [B]} Let $g_r\in \BC\{y,z\}$ be $\underline{\text{\rm
a generalized semi-quasi-Puiseux germ of the recursive r-type}}$ as
in [A]. There are two additional conditions, denoted by \text{\rm
The 5-th ${\text{\rm{Cond}}}^{\text{{\bf(0)}}}$} and \text{\rm The
6-th} ${\text{\rm{Cond}}}^{\text{{\bf(0)}}}$. \ms

\noindent $\underline{\text{\rm The 5-th
${\text{\rm{Cond}}}^{\text{{\rm(0)}}}$}}$ \quad For each
$q=1,2,\dots,r$, the following inequalities hold:
$$\align
\noindent \text{\rm (5)}&\text{\rm(5a)} \qquad \qquad\qquad\qquad
\text{ $\gcd(n_q,\Delta_q(\beta_{q,0,k})^q_{k=1})=1$ \quad for $1\le
q\le
r$.}  \\
&\text{\rm (5b)} \qquad \qquad\qquad\qquad
\f{\Delta_{q}(\beta_{q,i,k})^{q}_{k=1}}{n_{q}-i}
> \f{\Delta_{q}(\beta_{q,0,k})^{q}_{k=1}}{n_{q}} \quad
\text {for $0< i<n_q$.} \qquad \qquad\qquad\qquad
\endalign$$

\noindent $\underline{\text{\rm The 6-th
${\text{\rm{Cond}}}^{\text{{\rm(0)}}}$}}$ \quad The following
inequalities hold: Note that $2\le j\le r$.

\noindent\rm{(6)(6a)}  $2\le n_1<\beta_{1,0,1}$.

\rm{(6b)} $n_{j}\ge 2$, $\beta_{j,i,1}>0$, and $0\le
\beta_{j,i,k}<n_{k-1}$ for $2\le j\le r$, $0\le i<n_j$ and $2\le
k\le j$. \ms

Now, we define the new terminology in [B1], [B2] and [B3] of [B].

\noindent{\bf [B1]} It is said that $g_r\in \BC\{y,z\}$ is
$\underline{\text{\rm a generalized quasi-Puiseux germ of the
recursive r-type}}$ if $g_r$ in [A] satisfies an additional
condition, denoted by \text{\rm The 5-th
${\text{\rm{Cond}}}^{\text{{\bf(0)}}}$}, equivalently, if $g_r\in
\BC\{y,z\}$ satifies the above five conditions, i.e., \text{\rm The
1-th ${\text{\rm{Cond}}}^{\text{{\rm(0)}}}$}, $\dots$, \text{\rm The
4-th ${\text{\rm{Cond}}}^{\text{{\rm(0)}}}$} and \text{\rm The 5-th
${\text{\rm{Cond}}}^{\text{{\bf(0)}}}$}. \ms

\noindent{\bf[B2]} $g_r\in \BC\{y,z\}$ is called
$\underline{\text{\rm a generalized Puiseux germ of the recursive
r-type}}$ in $\BC\{y,z\}$ if $g_r$ in [A] satisfies \text{\rm The
5-th ${\text{\rm{Cond}}}^{\text{{\bf(0)}}}$} and an inequality in
(6a) of \text{\rm The 6-th ${\text{\rm{Cond}}}^{\text{{\bf(0)}}}$},
without mentioning any other equalities in (6b). \ms

\noindent{\bf[B3]} $g_r\in \BC\{y,z\}$ is called
$\underline{\text{\rm a generalized standard Puiseux germ of the
recursive r-type}}$ in $\BC\{y,z\}$ if $g_r$ in [A] satisfies
\text{\rm The 5-th ${\text{\rm{Cond}}}^{\text{{\bf(0)}}}$} and
\text{\rm The 6-th ${\text{\rm{Cond}}}^{\text{{\bf(0)}}}$}. \ms

\noindent{\bf[C]} Let $g_r\in \BC\{y,z\}$ be a generalized standard
Puiseux germ of the recursive r-type in [B3]. If each unit
$\ve_{j,i}=\ve_{j,i}(y,z)$ is equal to $\ve_{j,i}(y,0)$ for $1\le
j\le r$ and $0\le i<n_j$ if exists in \text{\rm The 2-th
${\text{\rm{Cond}}}^{\text{{\rm(0)}}}$} of \text{\rm[A]}, $g_r$ is
called $\underline{\text{\rm a generalized standard Puiseux
Weierstrass polynomial }}$ $\underline{\text{\rm  of the recursive
r-type}}$. In particular, if each unit $\ve_{j,i}=\ve_{j,i}(y,z)$ is
equal to an integer one for $1\le j\le r$ and $0\le i<n_j$ if
exists, as in \text{\rm The 2-th
${\text{\rm{Cond}}}^{\text{{\rm(0)}}}$} of \text{\rm[A]}, then $g_r$
is called $\underline{\text{\rm the generalized standard Puiseux
polynomial of the recursive r-type.}}$
\enddefinition \ms

\proclaim{Theorem 12.0(To find the necessary and sufficient
condition for a generalized semi-quasi-Puiseux convergent power
series in $\BC\{y,z\}$ of the recursive type to be irreducible in
$\C\{y,z\}$)}

$\underline{\text{\bf {Assumptions}}}$  Let $g_r$ be
$\underline{\text{\rm a generalized semi-quasi-Puiseux germ of the
recursive r-type}}$, satisfying the same properties and notations as
in {\rm [A]} of {\rm Definition 12.0.0}.

In addition, assume that we have the following:
$\gcd(n_1,\beta_{1,0,1})=1$. \ms

$\underline{\text{\bf {Conclusions}}}$ For each $j=1,2,\dots,r$, let
$(0,0)$ be the singularity of an analytic variety
$V(g_j)=\{(y,z):g_j(y,z)=0\}$ except for $V(g_1)$ with
$\beta_{1,0,1}=1$.

$\underline{\text{For brevity of representation, we may assume that the above $g_r$
satisfies all the equalities}}$
$\underline{\text{in {\rm(5a)} of \text{\rm The 5-th
${\text{\rm{Cond}}}^{\text{{\bf(0)}}}$} of {\rm [B]} of {\rm
Definition 12.0.0}, without mentioning all the inequalities}}$
$\underline{\text{in {\rm(5b)}
of \text{\rm The 5-th ${\text{\rm{Cond}}}^{\text{{\bf(0)}}}$}}}$. \ms

Then, we get two statements {\rm[A]} and {\rm[B]} as follows:
$$ \align
\noindent & [A] \quad \text{$g_r$ is irreducible in $\BC\{y,z\}$} \\
&  \iff \text{$g_1,\dots,g_{r-1}$ are irreducible in $\BC\{y,z\}$
and
   $\f{\Delta_{r}(\beta_{r,i,k})^{r}_{k=1}}{n_{r}-i}
> \f{\Delta_{r}(\beta_{r,0,k})^{r}_{k=1}}{n_{r}}$ for $0< i<n_j$} \qquad \qquad\\
& \iff {\f{\Delta_{j}(\beta_{j,i,k})^{j}_{k=1}}{n_{j}-i}}
> {\f{\Delta_{j}(\beta_{j,0,k})^{j}_{k=1}}{n_{j}}} \quad
\text {for each $j=1,2,\dots,r$ and for $0< i<n_j$.}
\endalign
$$

\noindent{\rm[B]} Let $g_r$ be irreducible in $\BC\{y,z\}$.

\noindent{\rm[B1]} Let
$V(y^{\gamma}g_r)=\{(y,z):y^{\gamma}g_r(y,z)=0\}$ be an analytic
variety at $(0,0)$ in $\BC^2$ defined by
$$\align
  y^{\gamma}g_r(y,z) \quad \text{such
that} \quad \left\{\eqalign{& \text{$\gamma=1$, \quad if \quad
$\beta_{1,0,1}=1$}, \cr & \text{$\gamma=0$, \quad if
\quad $\beta_{1,0,1}>1$.} \cr} \right.  \tag 12.0.0 \\
\endalign$$
Then, $y^{\gamma}g_r\in$ the type $[r]$ under the standard
resolution, denoted by $\tau$, in the sense of {\rm Definition 2.5}.
Also, if $\beta_{1,0,1}=1$ then $g_r\in$ the type $[r-1]$ under the
standard resolution. \ms

\noindent{\rm[B2]} In particular, $z^{\delta}yg_r\in$ the type $[r]$
under the same standard resolution $\tau$, whether $\delta$ is
either one or zero. \ms

\noindent{\rm[C]} Let $g_r$ be irreducible in $\BC\{y,z\}$.

\noindent{\rm[C1]} For each $j=1,2,\dots,r$, let
$V(\psi_j)=\{(y,z):\psi_j(y,z)=0\}$ be an analytic variety at
$(0,0)$ in $\BC^2$ defined by
$$\align
\psi_j(y,z)=g_j(y,z,0). \tag {$**$}
\endalign$$

Then, the singularity of both $V(\psi_j)$ and $V(g_j(y,z,c_j))$ have
the same divisor under the standard resolutions. \ms

\noindent{\rm[C2]} In particular, for each $j=1,2,\dots,r$, let
$V(H_j)=\{(y,z):H_j(y,z)=0\}$ be an analytic variety at $(0,0)$ in
$\BC^2$, each of which is defined as follows:

{\rm(i)} \qquad $H_1=z^{n_1}+y^{\beta_{1,0,1}}$.

{\rm(ii)} \qquad
$H_j=H^{n_j}_{j-1}+y^{\beta_{j,0,1}}z^{\beta_{j,0,2}}H^{\beta_{j,0,3}}_1\cdots
H^{\beta_{j,0,j}}_{j-2}$ for $j=2,3,\dots,r.$ \ms

Then, the singularity of both $V(g_j(y,z,c_j))$ and $V(H_j)$ have
the same divisor under the standard resolutions. $\square$
\endproclaim \ms

{\bf \S12.1. In preparation for the proof of Theorem 12.0 } \ms

In this section, for the proof of Theorem $12.0$ we will prepare the
statements of five sublemmas without proofs, consisting of Sublemma
12.1, Sublemma 12.2, ..., Sublemma 12.5. In $\S 13$, we will finish
the proof of Theorem $12.0$ with Corollary $12.6$, using the proofs
of these five sublemmas. \ms

\proclaim{Sublemma 12.1} $\underline{\text{\bf {Assumptions}}}$
\quad Suppose that the same properties and notations as in {\bf
{Assumptions}} of Theorem $12.0$ hold.

For any integer $r\ge 2$, let
$\Delta^{\sharp}_2(\beta_{2,i,1},\beta_{2,i,2})$ and
$\Delta^{\sharp}_j(\beta_{j,i,k})^j_{k=1}$ with $3\le j\le r$ be the
notations defined as follows : Note that
$\Delta_2(t_1,t_2)=n_1t_1+\beta_{1,0,1}t_2$ for each $(t_1,t_2)\in
N^2_0$.
$$\align
(12.1.1) \quad   \Delta^{\sharp}_2(\beta_{2,i,1},\beta_{2,i,2})
&=\Delta_2(\beta_{2,i,1},\beta_{2,i,2}) \quad \text {for $0\le i<n_2$.}  \\
 \Delta^{\sharp}_j(\beta_{j,i,k})^j_{k=1}
 &=\Delta_2(\beta_{j,i,1},\beta_{j,i,2})+n_1\beta_{1,0,1}\beta_{j,i,3}
 +n_1\beta_{1,0,1}n_2\beta_{j,i,4} \\
 &\quad +n_1\beta_{1,0,1} n_2 n_3 \beta_{j,i,5}
 + \cdots +n_1\beta_{1,0,1}n_2 \cdots n_{j-2}\beta_{j,i,j}
 \quad \text {for $0\le i<n_j$.} \qquad
 \qquad \qquad
\endalign$$

$\underline{\text{\bf {Conclusions}}}$ \quad  Then, we have the
following:
$$\align
\Delta^{\sharp}_2(\beta_{2,i,1},\beta_{2,i,2})
&> n_1\beta_{1,0,1}(n_2-i) \quad \text{on $g_2$}.  \tag 12.1.2\\
\Delta^{\sharp}_j(\beta_{j,i,k})^j_{k=1} &>
n_1\beta_{1,0,1}n_2n_3\cdots n_{j-1}(n_j-i) \quad \text{on
$g_j$}.\quad \text{$\square$}
\endalign$$
\endproclaim \ms

\proclaim{Sublemma 12.2} $\underline{\text{\bf {Assumptions}}}$
\quad Suppose that the same properties and notations as in {\bf
{Assumptions}} of Theorem $12.0$ hold. Let $r$ be arbitrary integer
with $r\ge 1$.

In addition, we need the following assumption: Note that
$\gcd(n_1,\beta_{1,0,1})=1$.
$$\align
(12.2.0) \quad \text{$g_1$ is irreducible in $\BC\{y,z\}$, \quad
{equivalently,} \quad $\f{\beta_{1,i,1}}{n_1-i}>
\f{\beta_{1,0,1}}{n_1}>0$ \quad for $0< i<n_1$}.  \qquad \qquad
\qquad
\endalign$$

But, note that $g_{r}$ may not be irreducible in $\BC\{y,z\}$ for
some $r\ge 2$. \ms

$\underline{\text{\bf {Conclusions}}}$ \quad  Then, we get the
following:

{\rm(a)} $g_r=g_r(y,z)$ can be written in the form
$$\align
(12.2.1) \qquad \qquad & g_r=(z^{n_1}+\ve_{1,0}
y^{\beta_{1,0,1}})^{d_r}+\sum_{\alpha,\beta\ge 0}
c^{(r)}_{\alpha,\beta}y^{\alpha}z^{\beta} \quad \text{with
$\ve_{1,0}=1$
and} \qquad \qquad \qquad \\
& \quad \text{with
$n_1\alpha+\beta_{1,0,1}\beta>n_1\beta_{1,0,1}d_r$},
\endalign$$
where if $r\ge 2$ then $d_r=\prod^r_{k=2}n_k$ and $d_1=1$, and a
unit $\ve_{1,0}=\ve_{1,0}(y,z)$ may be analytically assumed to be
one in $\BC\{y,z\}$, and the $c^{(r)}_{\alpha,\beta}$ are nonzero
complex numbers for some nonnegative integers $\alpha$ and $\beta$
such that $n_1\alpha+\beta_{1,0,1}\beta>n_1\beta_{1,0,1}d_r$. \ms

{\rm(b)} For each $r\ge 1$, we have the following:

{\rm(b1)} The multiplicity of $g_r(0,z)$ at $z=0$ is
$n_1\prod^r_{k=2}n_k$ when $g_r=g_r(y,z)$.

{\rm(b2)} The multiplicity of $g_r(y,0)$ at $y=0$ is
 $\beta_{1,0,1}\prod^r_{k=2}n_k$ when $g_r=g_r(y,z)$. \ms

{\rm(c)} For each $r\ge 1$, we have the following:

{\rm(c1)} If $n_1<\beta_{1,0,1}$ then
$\alpha+\beta>n_1\prod^r_{k=2}n_k$, and so the multiplicity of $g_r$
at $(y,z)=(0,0)$ is $n_1\prod^r_{k=2}n_k$.

{\rm(c2)} If $n_1>\beta_{1,0,1}$ then
$\alpha+\beta>\beta_{1,0,1}\prod^r_{k=2}n_k$, and so the
multiplicity of $g_r$ at $(y,z)=(0,0)$ is
$\beta_{1,0,1}\prod^r_{k=2}n_k$.  $\square$
\endproclaim
\ms

\proclaim{Sublemma 12.3} $\underline{\text{\bf {Assumptions}}}$
\quad Suppose that the same properties and notations as in {\bf
{Assumptions}} of Theorem $12.0$ hold. Let $r\ge 2$ be an arbitrary
positive integer.

In addition, assume that we have the following: Note that
$\gcd(n_1,\beta_{1,0,1})=1$ by the assumption of Theorem $12.0$, but
that the condition in {\rm(12.3.0)} does not belong to the
assumption of Theorem $12.0$.
$$\align
(12.3.0) \quad \text{$g_1$ is irreducible in $\BC\{y,z\}$, \quad
{equivalently,} \quad $\f{\beta_{1,i,1}}{n_1-i}>
\f{\beta_{1,0,1}}{n_1}>0$ \quad for $0< i<n_1$}.  \qquad \qquad
\qquad
\endalign$$

Since $\gcd(n_1,\beta_{1,0,1})=1$ with $n_1\ge 2$ and
$\beta_{1,0,1}\ge 1$, then there are two nonnegative integers $a>0$
and $b\ge 0$ such that $a\beta_{1,0,1}-bn_1=1$.

For given two integers $a>0$ and $b\ge 0$, let $\Omega_2:N^2_0\to
N_0$ be a function defined by
$$\align
             \Omega_2(t_1,t_2)=at_1+bt_2. \tag 12.3.1
\endalign$$

Let $\Omega^{\sharp}_2(\beta_{2,i,1},\beta_{2,i,2})$ and
$\Omega^{\sharp}_j(\beta_{j,i,k})^j_{k=1}$ with $3\le j\le r$ be the
notations defined as follows:
$$\align
(12.3.2)\qquad \quad \Omega^{\sharp}_2(\beta_{2,i,1},\beta_{2,i,2})
 &=\Omega_2(\beta_{2,i,1},\beta_{2,i,2}) \quad \text {for $0\le i<n_2$.}
 \qquad \qquad \qquad \qquad \qquad \qquad \\
  \Omega^{\sharp}_j(\beta_{j,i,k})^j_{k=1} &=\Omega_2(\beta_{j,i,1},\beta_{j,i,2})
 +bn_1\beta_{j,i,3}+bn_1n_2\beta_{j,i,4}\\
 &\quad+\cdots +bn_1n_2\cdots n_{j-2}\beta_{j,i,j} \quad \text {for $0\le i<n_j$.}
\endalign
$$

$\underline{\text{\bf {Conclusions}}}$ \quad  Then, we get the
following:
$$\align
 \Omega^{\sharp}_2(\beta_{2i1},\beta_{2i2})
 & \ge bn_1n_2  \quad \text {for $0\le i<n_2$.} \tag 12.3.3 \\
 \Omega^{\sharp}_j(\beta_{jik})^j_{k=1} &\ge bn_1n_2n_3\cdots
n_{j-1}(n_j-i) \quad \text {for $0\le i<n_j$.} \text{$\square$}
\endalign$$
\endproclaim
\ms

\proclaim{Sublemma 12.4} $\underline{\text{\bf {Assumptions}}}$
\quad Suppose that the same properties and notations as in {\bf
{Assumptions}} of Theorem $12.0$ hold. As in Sublemma $12.3$,
additionally assume that we have the following:
$$\align
(*1) \quad \text{$g_1$ is irreducible in $\BC\{y,z\}$, \quad
{equivalently,} \quad $\f{\beta_{1,i,1}}{n_1-i}>
\f{\beta_{1,0,1}}{n_1}>0$ \quad for $0< i<n_1$}.  \qquad \qquad
\qquad
\endalign$$

But, note that $g_{r}$ may not be irreducible in $\BC\{y,z\}$ for
some $r\ge 2$. \ms

$\underline{\text{\bf {Conclusions}}}$ \quad For each
$j=1,2,\dots,r$, let $V(g_j)=\{(y,z):g_j(y,z)=0\}$  be an analytic
variety at the origin in $\BC^2$. For the construction of the
statement of the conclusion, let $V(G)=\{(y,z):G(y,z)=0\}$ be
another analytic variety with an isolated singularity at the origin
in $\BC^2$ defined by the form
$$
\align g_0 &=z^{n_1}+\ve_{1,0}y^{\beta_{1,0,1}}
\quad \text{with a unit  $\ve_{1,0}\in \C\{y,z\}$,} \tag 12.4.0 \\
G &=y^{\gamma}g_0,
\endalign
$$
satisfying the properties {\rm(i)} and {\rm(ii)}:

\roster \item "(i)" If $\beta_{1,0,1}=1$, then $\gamma=1$.

\item "(ii)" If $\beta_{1,0,1}\ge 2$, then $\gamma=0$.
\endroster \ms

Let $\tau_m=\pi_1\circ\pi_2\circ\cdots\circ\pi_m:M^{(m)}\to\BC^2$ be
the compositions of a finite number $m$ of successive blow-ups
$\pi_i$ which is needed to get the standard resolution of the
singular point of $V(G)$. If $V(g_1)$ in the assumption of Theorem
$12.0$ has the singular point at the origin, then as compared with
the above $\tau_m$, exactly the same $\tau_m$ can be also used for
the standard resolution of the singular point of $V(yg_1)$ as far as
the blow-ups process is concerned.

{\rm(a)(a1)} We can use just one coordinate patch of the local
coordinates for each blow-up $\pi_i$ of $\tau_m$ with $1\le i\le m$
in the sense of Lemma $2.12$.

{\rm \quad(a2)} Just as above, we can use the same $\tau_m$ for the
composition of the first finite number $m$ of successive blow-ups in
process of the standard resolution of the singular point $(0,0)$ of
$V(g_j)$ for all $j=2,3,\dots,r$.

{\rm \quad(a3)} Also, we can use just the common one coordinate
patch of the given local coordinates for each blow-up $\pi_i$ of the
above $\tau_m$ in {\rm (a1)}, in order to study any of
$V^{(i)}(g_j)$ for all $j=2,3,\dots,r$ and all $i=1,2,\dots,m$ in
the sense of Lemma $2.14$. \ms

{\rm(b)} For brevity of notation, by {\rm (a3)} let $(v,u)$ be the
common one of the local coordinates for the $m-th$ blow-up
$\pi_m:M^{(m)}\to M^{(m-1)}$ at $(0,0)$ which is the quasisingular
point of $V^{(m-1)}(G)$. Being viewed as an analytic mapping,
$\tau_m:M^{(m)}\to\BC^2$ can be written in the form
 $$
 \tau_m(v,u)=(y,z)=(v^{n_1}u^a,v^{\beta_{1,0,1}}u^b), \tag 12.4.1
 $$
where \text{\rm(i)} $a>0$ and $b\ge 0$ are some nonnegative integers
such that $a\beta_{1,0,1}-bn_1=1$,

\noindent{\text{\rm(ii)}}  $E_m=\{v=0\}$ is defined by the $m-th$
exceptional curve of the first kind. \ms

{\rm(c)} As in Sublemma $12.1$ and Sublemma $12.3$, use the same
notations for a sequence $\{\Delta^{\sharp}_i: N^i_0\to N_0,
\text{functions for i=2,\dots,r}\}$ and $\{\Omega^{\sharp}_i:
N^i_0\to N_0, \text{functions for i=2,\dots,r}\}$ where
$\Omega_2:N^2_0\to N_0$ is a function defined by
$\Omega_2(t_1,t_2)=at_1+bt_2$  for given two nonnegative integers a
and b in $(b_1)$ of {\rm(12.4.1)}, and we may start with assuming
that $\ve_{1,0}=1$ in
$V(y^{\gamma}g_0)=\{y^{\gamma}(z^{n_1}+\ve_{1,0}y^{\beta_{1,0,1}})=0\}$,
in order to study $V^{(i)}(g_j)$ for all $i=1,2,\dots,m,$ and all
$j=1,2,\dots,r$. Whether $\beta_{1,0,1}\ge 2$ or $\beta_{1,0,1}=1$,
we may write that
$((y^{\gamma}g_0)\circ\tau_m)_{proper}=u^{a,\gamma}(1+\ve_{1,0}u)$
with $\ve_{1,0}=1$, without complexity of the notation if necessary,
noting that if $\beta_{1,0,1}=1$, $V(g_0)$ and $V(g_1)$ have no
singularity at the origin.

Now, along $v=0$, $(g_j\circ\tau_m)_{total}$ with
$(g_j\circ\tau_m)_{proper}$ can be written as follows:
$$\align
(12.4.2)\qquad  ((y^{\gamma}g_0)\circ\tau_m)_{total}
&=v^{(\gamma+\beta_{1,0,1})n_1}u^{bn_1+a\gamma}
((y^{\gamma}g_0)\circ\tau_m)_{proper} \quad\text{with} \\
((y^{\gamma}g_0)\circ\tau_m)_{proper} &=(1+\ve_{1,0}u) \quad \text{with $\ve_{1,0}=1$,}\\
 (g_1\circ\tau_m)_{total} =&v^{n_1\beta_{1,0,1}}u^{bn_1}
 (g_1\circ\tau_m)_{proper}\quad \text{with}  \\
 (g_1\circ\tau_m)_{proper} =&(1+\ve_{1,0}u)
 +c_1\sum^{n_1-1}_{i=1}\ve'_{1,i}v^{\Delta^{\sharp}_2(\beta_{1,i,1},i)-n_1\beta_{1,0,1}}
 u^{\Omega^{\sharp}_2(\beta_{1,i,1},i)-bn_1}  \\
 =&1+\ve_{1,0}\bar{u}=1+\bar{u} \qquad \qquad \text{for brevity of notation}, \\
 (g_j\circ\tau_m)_{total} =&v^{n_1\beta_{1,0,1}d_j}u^{bn_1d_j}
 (g_j\circ\tau_m)_{proper}\quad\text{with}\\
 (g_j\circ\tau_m)_{proper} =&(g_{j-1}\circ\tau_m)^{n_j}_{proper}
 +\{\ve'_{j,0} v^{\Delta^{\sharp}_j(\beta_{j,0,k})^j_{k=1}-n_1\beta_{1,0,1}d_j}
 u^{\Omega^{\sharp}_j(\beta_{j,0,k})^j_{k=1}-bn_1d_j} \\
   \times ({g_1}&\circ\tau_m)_{proper}^{\beta_{j,0,3}}
 (g_2\circ\tau_m)^{\beta_{j,0,4}}_{proper}
 \cdots (g_{j-2}\circ\tau_m)^{\beta_{j,0,j}}_{proper} \} \\
  +{c_j}{\cdot}\{&\sum^{n_j-1}_{i=1}\ve'_{j,i} v^{\Delta^{\sharp}_j(\beta_{j,i,k})^j_{k=1}
 -n_1\beta_{1,0,1}d_{j-1}(n_j-i)}
 u^{\Omega^{\sharp}_j(\beta_{j,i,k})^j_{k=1}-bn_1d_{j-1}(n_j-i)} \\
 \times ({g_1}&\circ\tau_m)_{proper}^{\beta_{j,i,3}}
 (g_2\circ\tau_m)^{\beta_{j,i,4}}_{proper}
 \cdots (g_{j-2}\circ\tau_m)^{\beta_{j,i,j}}_{proper}(g_{j-1}\circ\tau_m)^{i}_{proper}\},
\endalign
$$where if $j\ge 2$ then $d_j=\prod^j_{k=2}n_k$ and $d_1=1$, and
each $\ve'_{j,i}=\ve_{j,i}\circ\tau_m(v,u)$ is a unit in
$\BC\{v,1+u\}$ for $2\le j\le r$ and $0\le i<n_j$.

Note by Sublemma $12.1$ and Sublemma $12.3$ that for $j=2,3,\dots,r$
and $i=0,1,\dots,n_j-1$,
$$\align
(12.4.3) \qquad \qquad \Delta^{\sharp}_j(\beta_{j,i,k})^j_{k=1} &>
n_1\beta_{1,0,1}n_2\cdots
 n_{j-1}(n_j-i)=n_1\beta_{1,0,1}d_{j-1}(n_j-i)
 \quad \text {and} \qquad \qquad  \\
 \Omega^{\sharp}_j(\beta_{j,i,k})^j_{k=1} &\ge bn_1n_2\cdots n_{j-1}(n_j-i)
 =bn_1d_{j-1}(n_j-i).
\endalign
$$
Moreover,
$(y^{\beta_{j,i,1}}z^{\beta_{j,i,2}}g^{\beta_{j,i,3}}_1g^{\beta_{j,i,4}}_2\cdots
g^{\beta_{j,i,j}}_{j-2})\circ\tau_m(v,u)$ can be viewed as
$$
u^{\Omega^{\sharp}_r(\beta_{j,i,k})^j_{k=1}}v^{\Delta^{\sharp}_j(\beta_{j,i,k})^j_{k=1}}
(g_1\circ\tau_m)^{\beta_{j,i,3}}_{proper}(g_2\circ\tau_m)^{\beta_{j,i,4}}_{proper}
\cdots (g_{j-2}\circ\tau_m)^{\beta_{j,i,j}}_{proper}, \tag 12.4.4
$$
where $(c_1)$ for $j=1,2,\dots,r$ and for $1\le i< n_j$,
$(\beta_{j,i,k})^j_{k=1}\in N^j_0$ as in {\rm [A]} of {\rm
Definition 12.0.0},

$(c_2)$
$\Delta^{\sharp}_j(\beta_{j,i,k})^j_{k=1}=\Delta_2(\beta_{j,i,1},\beta_{j,i,2})
 +n_1\beta_{1,0,1}\beta_{j,i,3}+n_1\beta_{1,0,1}n_2\beta_{j,i,4}+\cdots
 +n_1\beta_{1,0,1}n_2\cdots n_{j-2}\beta_{j,i,j}$ by the definition of
 $\Delta^{\sharp}_j(\beta_{j,i,k})^j_{k=1}$ in Sublemma $12.1$,

$(c_3)$ $\Omega^{\sharp}_j(\beta_{j,i,k})^j_{k=1}
 =\Omega_2(\beta_{j,i,1},\beta_{j,i,2})+bn_1\beta_{j,i,3}
 +bn_1n_2\beta_{j,i,4}+\cdots +bn_1n_2\cdots n_{j-2}\beta_{j,i,j}$
by the definition of $\Omega^{\sharp}_j(\beta_{j,i,k})^j_{k=1}$ in
Sublemma $12.3$. \ms

{\rm(d)} Let ${\tau_m}^{-1}(0,0)=\cup^m_{i=1}E_i$ where $E_i$ is an
exceptional curve of the first kind. For $j=1,2,\dots,r$, let
$$(g_j\circ\tau_m)_{divisor}=V^{(m)}(g_j)+\sum^m_{i=1}e_{j,i}E_i, \tag 12.4.5 $$
where $V^{(m)}(g_j)$ is the proper transform of $V(g_j)$ under
$\tau_m$. \ms

Then we have the following: Note again that $\tau_m$ is the
composition of a finite number $m$ of successive blow-ups, which is
needed to get the standard resolution of the singular point of
$V(y^{\gamma}g_1)$ or $V(g_1)$.

{\rm(d1)} If $\beta_{1,0,1}\ge 2$, then $e_{j+1,i}=n_{j+1}e_{j,i}$
for any $j\ge 1$ and for $i=1,2,\dots,m$. If $\beta_{1,0,1}=1$, then
$e_{j+1,i}=n_{j+1}e_{j,i}$ for any $j\ge 2$ and for $i=1,\dots,m$.

In particular, $e_{j,m}=n_1\beta_{1,0,1}n_2\cdots n_j$ for
$j=2,\dots,r$, and if $j=1$ with $\beta_{1,0,1}>1$, then
$e_{j,m}=n_1\beta_{1,0,1}$. \ms

{\rm(d2)} $V^{(m)}(g_j)\cap(\cup^m_{i=1}E_i)=V^{(m)}(g_j)\cap
E_m=\{(v,1+\ve_{1,0}u)=(0,0)\}$ for any $j=2,\dots,r$ where
$1+\ve_{1,0}u=(g_1\circ\tau_m)_{proper}$. \ms

{\rm(d3)} If $\beta_{1,0,1}\ge 2$, then for any $j=1,2,\dots,r$,
$g_j\in$ the type $[1]$ under $\tau_m$. If $\beta_{1,0,1}=1$, then
for any $j=1,2,\dots,r$, $g_j\in$ the type $[0]$ under $\tau_m$.

If $\beta_{1,0,1}\ge 1$, note that for all $j=1,2,\dots,r$,
$z^{\delta}yg_j\in$ the type $[1]$ under $\tau_m$ whether $\delta=1$
or $\delta=0$, by Theorem $3.6$.  $\square$
\endproclaim
\ms

\proclaim{Sublemma 12.5} $\underline{\text{\bf {Assumptions}}}$
Suppose that the same properties and notations as in {\bf
{Assumptions}} of Theorem $12.0$ hold. As in either Sublemma $12.3$
or Sublemma $12.4$, additionally assume that we have the following:
$$\align
(12.5.0) \quad \text{$g_1$ is irreducible in $\BC\{y,z\}$, \quad
{equivalently,} \quad $\f{\beta_{1,i,1}}{n_1-i}>
\f{\beta_{1,0,1}}{n_1}>0$ \quad for $0< i<n_1$}.  \qquad \qquad
\qquad
\endalign$$

Let $r$ be an arbitrary positive integer with $r\ge 2$. Throughout
this sublemma, we will use the same notations and consequences as in
Sublemma $12.4$, in order to get the representation for the
conclusion of this sublemma. \ms

$\underline{\text{\bf {Conclusions}}}$ \quad As $\{g_k:
k=1,2,\dots,r\}$ with $g_k\in \BC\{y,z\}$ satisfies four conditions
in the assumptions of Theorem $12.0$, denoted by \text{\bf The 1-th
${\text{\bf{Cond}}}^{\text{{\bf(0)}}}$}, \dots, \text{\bf The 4-th
${\text{\bf{Cond}}}^{\text{{\bf(0)}}}$} in Definition $12.0.0$, then
$(g_r\circ\tau_m)_{proper}\in \BC\{v,1+\bar{u}\}$ is a generalized
semi-quasi-Puiseux germ of the recursive $(r-1)$ type in the sense
of Definition $12.0.0$ where $\{(g_k\circ\tau_m)_{proper}:
k=2,3,\dots,r\}$ with $(g_k\circ\tau_m)_{proper}$ in
$\BC\{v,1+\bar{u}\}$ satisfies the same kind of four conditions,
which will be denoted by \text{\bf The 1-th
${\text{\bf{Cond}}}^{\text{{\bf(1)}}}$}, \dots, \text{\bf The 4-th
${\text{\bf{Cond}}}^{\text{{\bf(1)}}}$}. Note that
$\{(g_k\circ\tau_m)_{proper}: k=2,3,\dots,r\}$ has been already
well-defined by Sublemma $12.4$. \ms

In more detail, in order to construct four conditions recursively,
which will be denoted by \text{\bf The 1-th
${\text{\bf{Cond}}}^{\text{{\bf(1)}}}$}, \dots, \text{\bf The 4-th
${\text{\bf{Cond}}}^{\text{{\bf(1)}}}$}, we prefer to add one more
condition to the above four conditions, denoted by \text{\bf The
5-th ${\text{\bf{Cond}}}^{\text{{\bf(1)}}}$}, for convenience of
representation.  By using the same kind of properties and notations
as in Definition $12.0.0$, the desired construction is as follows:
$$\align
 \quad & Let \quad \{Y_k: k=1,2,\dots,r-1\} \quad
\text{with $Y_k\subset N_0$}, \\
& \{h_k: k=1,2,\dots,r-1\} \quad \text{with
$h_k=(g_{k+1}\circ\tau_m)_{proper}$ in $\BC\{v,1+\bar{u}\}$,} \qquad \qquad \\
&\{\text{$\Xi_k:N^k_0\to N_0$ is an integer-valued
function for $k=1,2,\dots,r-1$}\} \\
& \text{be three sequences, satisfying the following five conditions
for each k}:
\endalign$$

Such conditions are denoted by \text{\bf The 1-th
${\text{\bf{Cond}}}^{\text{{\bf(1)}}}$}{\bf, \dots,} \text{\bf The
5-th ${\text{\bf{Cond}}}^{\text{{\bf(1)}}}$}. \ms

\noindent \text{\bf The 1-th
${\text{\bf{Cond}}}^{\text{{\bf(1)}}}$:} Let
$\{Y_j:j=1,2,\dots,r-1\}$ with $Y_j\subset N_0$ be defined by
\roster

\item"(12.5.1)" $Y_1=\{s_1\}\cup \{\g_{1,i,1}:0\le i<s_1 \}$ with
$s_1\ge 2$ and $\g_{1,0,1}\ge 1$,

\item"{}" $Y_j=\{s_j\}\cup \{\g_{j,i,1}:0\le i<s_j \} \cup
\{\g_{j,i,2}:0\le i<s_j \}\cup \{\g_{j,i,j}:0\le i<s_j \}$

with $s_j\ge 2$,
\endroster
such that for each $j=1,2,\dots,r-1$,
$e_{1,m}=n_1\Delta_1(\beta_{1,0,1})=n_1\beta_{1,0,1}$, and
$$\align
(12.5.1.1) \quad s_1 &=n_2\ge 2,
 ~\gamma_{1,i,1}=\Delta^{\sharp}_2(\beta_{2,i,k})^2_{k=1}-n_1\beta_{1,0,1}(n_2-i)>0
 \text{ for $0\le i<s_1$},   \\
  s_j &=n_{j+1}\ge 2,
\gamma_{j,i,1}=\Delta^{\sharp}_{j+1}(\beta_{{j+1},i,k})^{j+1}_{k=1}
-n_1\beta_{1,0,1}n_2n_3\cdots n_{j}(n_{j+1}-i)>0 \text{ and} \\
  \gamma_{j,i,2}&=\beta_{{j+1},i,3},
\gamma_{j,i,3}=\beta_{{j+1},i,4}, \gamma_{j,i,4}=\beta_{{j+1},i,5} ,
\dots ,
 \gamma_{j,i,j}=\beta_{{j+1},i,{j+1}} \text{ for $0\le i<s_j$}, \qquad \quad\\
  \endalign
$$
noting that $\gamma_{1,i,1},\gamma_{2,i,1},\dots,\gamma_{r-1,i,1}$
are
 positive by Sublemma $12.1$. \ms

\noindent \text{\bf The 2-th
${\text{\bf{Cond}}}^{\text{{\bf(1)}}}$:} Let
$(g_2\circ\tau_m)_{proper},(g_3\circ\tau_m)_{proper},\ldots,
(g_r\circ\tau_m)_{proper}$ be denoted by $h_1,h_2,\dots,h_{r-1}$,
respectively in $\BC\{v,u+1\}$ as follows: For brevity of notation,
we can define a local biholomorphic mapping $\phi$ from
$(u,v)=(-1,0)$ to $(\bar{u},v)=(-1,0)$ such that
$1+\bar{u}=(g_1\circ\tau_m)_{proper}=(1+\ve_{1,0}u)
+c_1\sum^{n_1-1}_{i=1}\ve'_{1,i}v^{\Delta^{\sharp}_2(\beta_{1,i,1},i)-n_1\beta_{1,0,1}}
u^{\Omega^{\sharp}_2(\beta_{1,i,1},i)-bn_1}$ \text{ with
$\ve_{1,0}=1$} and $\phi(u,v)=(\bar{u},v)$.
$$\align
(12.5.2) \quad
(g_2\circ\tau_m)_{total}&=v^{n_2e_{1,m}}(g_2\circ\tau_m)_{proper}
=v^{n_2e_{1,m}}h_1 \quad \text{with} \\
 \quad \qquad h_1&=(1+\bar{u})^{s_1}+\eta_{1,0}
v^{\gamma_{1,0,1}}+c_2\sum^{s_1-1}_{i=1}\eta_{1,i}v^{\g_{1,i,1}}(1+\bar{u})^i
\text{ with $\eta_{1,0}=1$,} \qquad \qquad \\
 \quad \qquad (g_j\circ\tau_m)_{total}&=v^{n_jn_{j-1}\cdots
n_2e_{1,m}}(g_j\circ\tau_m)_{proper} =v^{n_jn_{j-1}\cdots
n_2e_{1,m}}h_{j-1} \quad \text{with} \\
  h_{j-1}&=
h^{s_{j-1}}_{j-2}+\eta_{j-1,0}v^{\gamma_{j-1,0,1}}(1+\bar{u})^{\gamma_{j-1,0,2}}
 h^{\gamma_{j-1,0,3}}_1\cdots h^{\gamma_{j-1,0,j-1}}_{j-3} \\
&\quad+c_{j}\sum^{s_{j-1}-1}_{i=1}\eta_{j-1,i}v^{\gamma_{j-1,i,1}}
(1+\bar{u})^{\gamma_{j-1,i,2}}\alpha
 h^{\gamma_{j-1,i,3}}_1\cdots
 h^{\gamma_{j-1,i,j-1}}_{j-3}h^{i}_{j-2}, \qquad \qquad
 \endalign$$
where $\eta_{j,i}=\eta_{j,i}(v,1+\bar{u})$ is a unit in
$\BC\{v,1+\bar{u}\}$ for $1\le j\le r-1$ and $1\le i\le s_j-1$,
noting that $\eta_{ji}=\ve'_{j+1,i}
u^{\Omega^{\sharp}_{j+1}(\beta_{j+1,i,k})^{j+1}_{k=1}-bn_1n_2\cdots
n_j(n_{j+1}-i)}$. Here, we may assume by a nonsingular change of
coordinates that $\eta_1$ can be equal to an integer one for the
standard resolution of the quasisingular point $(v,1+\bar{u})=(0,0)$
of $V(h_k)$ for $1\le k\le r-1$. \ms

\noindent \text{\bf The 3-th
${\text{\bf{Cond}}}^{\text{{\bf(1)}}}$:} Let $\{\text{\rm $\Xi_j:
N^j_0\to N_0$ is an integer-valued function for j=1,2,\dots,r-1}\}$
be a sequence defined by the following:
$$\align
(12.5.3)\quad \quad  &\text{$\Xi_1(t)=t$  for each $t\in N_0$.} \\
 \quad \quad
&\text{$\Xi_{j-1}(t_k)^{j-1}_{k=1}=t_{j-1}\Xi_{j-2}
(\gamma_{j-2,0,k})^{j-2}_{k=1}+s_{j-2}\Xi_{j-2}(t_k)^{j-2}_{k=1}$
for each $(t_k)^{j-1}_{k=1}\in  N^{j-1}_0.$}    \\
\endalign$$ \ms

\noindent \text{\bf The $(4{\alpha})$-th
${\text{\bf{Cond}}}^{\text{{\bf(1)}}}$:} For each $q=1,2,\dots,r-1$,
we have the following: Note that $r\ge 2$.
$$\align
(12.5.4{\alpha}) \qquad \quad
&\Xi_1(\gamma_{1,i,1})=\Delta^{\sharp}_2(\beta_{2,i,1},\beta_{2,i,2})
-n_1\beta_{1,0,1}(n_2-i)>0, \\
 \qquad \quad &
\Xi_q(\gamma_{q,i,k})^q_{k=1}-(s_q-i)s_{q-1}
\Xi_{q-1}(\gamma_{q-1,0,k})^{q-1}_{k=1} \quad \text{ for
$2\le q\le {r-1}$}\\
& \quad =\Delta_{q+1}(\beta_{q+1,i,k})^{q+1}_{k=1}
-(n_{q+1}-i)n_q\Delta_q(\beta_{q,0,k})^q_{k=1}>0 \quad \text{for
$0\le i<s_{q}$.} \qquad \qquad
\endalign$$ \ms

\noindent \text{\bf The 4-th
${\text{\bf{Cond}}}^{\text{{\bf(1)}}}$:} By \text{\rm The
$(4{\alpha})$-th ${\text{\rm{Cond}}}^{\text{{\rm(1)}}}$}, for
$q=1,2,\dots,{r-1}$, it is clear:
$$\align
(12.5.4) \qquad \qquad &  \Xi_1(\gamma_{1,i,1})=\gamma_{1,i,1}>0
\quad \text{for $0\le i<{s_1}$.}  \\
& \Xi_q(\gamma_{q,i,k})^q_{k=1}
>(s_q-i)s_{q-1}\Xi_{q-1}(\gamma_{q-1,0,k})^{q-1}_{k=1}
\quad \text{for $0\le i<s_{q}$.} \qquad \qquad \qquad \qquad
\endalign$$
\ms

\noindent \text{\bf The $(5{\alpha})$-th
${\text{\bf{Cond}}}^{\text{{\bf(1)}}}$:} By \text{\rm The
$(4{\alpha})$-th ${\text{\rm{Cond}}}^{\text{{\rm(1)}}}$}, we have
the following for $q=1,2,\dots,{r-1}$:
$$\align
(12.5.5{\alpha})(12.5.5{\alpha}.1) \qquad \qquad
\gcd(s_{q},\Xi_{q}(\gamma_{q,0,k})^{q}_{k=1})=
\gcd(n_{q+1},\Delta_{q+1}(\beta_{q+1,0,k})^{q+1}_{k=1}). \qquad
\qquad \qquad \qquad
\endalign$$

$$\align
(12.5.5{\alpha})(12.5.5{\alpha}.2) \qquad \qquad \quad &
\f{\Delta_{q+1}(\beta_{q+1,i,k})^{q+1}_{k=1}}{n_{q+1}-i}
> \f{\Delta_{q+1}(\beta_{q+1,0,k})^{q+1}_{k=1}}
{n_{q+1}} \quad \text {for $0<i<n_{q+1}$} \qquad \qquad\\
\qquad \iff \qquad &\\
&  \f{\Xi_{q}(\g_{q,i,k})^{q}_{k=1}}{s_{q}-i}
> \f{\Xi_{q}(\g_{q,0,k})^{q}_{k=1}}{s_{q}} \quad
\text {for $0< i<s_q$. \qquad $\square$} \qquad \qquad \\
\endalign$$
\endproclaim \ms

\newpage

\vfill \pagebreak

{\bf \S{13}. The proofs of Theorem 12.0 with five sublemmas and
corollaries in \S {12}} \bs

{\bf \S 13.1. For the proofs of five sublemmas} \ms

Using the same method as we have used in the proofs of Sublemma 5.1,
Sublemma 5.2,\dots, Sublemma 5.5 and following the same kind of
properties and notations as we have seen in the proof of Sublemma
5.1, Sublemma 5.2, \dots, Sublemma 5.5, the proofs of Sublemma 12.1,
Sublemma 12.2, \dots, Sublemma 12.5 can be generalized.  \ms

\demo{\bf Proof of Sublemma 12.1} If $r=2$, then it is trivial to
prove that $\Delta^{\sharp}_2(\beta_{2,i,1},\beta_{2,i,2})>
n_1\beta_{1,0,1}(n_2-i) \text{ on $g_2$}$, because
$\Delta^{\sharp}_2(\beta_{2,i,1},\beta_{2,i,2})=\Delta_2(\beta_{2,i,1},\beta_{2,i,2})$
by (12.1.1) and $\Delta_2(\beta_{2,i,1},\beta_{2,i,2})
>n_1\beta_{1,0,1}(n_2-i)$ by \text{\rm The 4-th
${\text{\rm{Cond}}}^{\text{{\rm(0)}}}$} in the assumptions of
Theorem $12.0$.

Let $r\ge 3$. For any $\ell=3,4,\dots,r$, it is trivial to note by
definition of
$\Delta^{\sharp}_{\ell}(\beta_{{\ell},i,k})^{\ell}_{l=1}$ in
(12.1.1) that the following three equalities are the same, and so we
can write
$c=\Delta^{\sharp}_{\ell}(\beta_{{\ell},i,k})^{\ell}_{k=1}-n_1\beta_{1,0,1}n_2n_3\cdots
n_{\ell-1}(n_{\ell}-i)$ for convenience of notation:
$$\align
(12.1.3) \qquad \qquad  &
c=\Delta^{\sharp}_{\ell}(\beta_{{\ell},i,k})^{\ell}_{k=1}
-n_1\beta_{1,0,1}n_2n_3\cdots n_{\ell-1}(n_{\ell}-i) \quad
\text{by $(c)$ of Sublemma $12.4$} \qquad \qquad\\
\iff &\\
(12.1.4) \qquad \qquad  &c=\Delta_2(\beta_{{\ell},
i,1},\beta_{{\ell}, i,2})+n_1\beta_{1,0,1}\beta_{{\ell},i,3}
 +n_1\beta_{1,0,1}n_2\beta_{{\ell},i,4}
 +n_1\beta_{1,0,1} n_2 n_3 \beta_{{\ell},i,5}   \\
& \quad \text{\qquad $+ \cdots +n_1\beta_{1,0,1}n_2 \cdots
n_{\ell-2}\beta_{{\ell},i,\ell}-n_1\beta_{1,0,1}n_2n_3\cdots n_{\ell-1}(n_{\ell}-i)$} \qquad \\
\iff &\\
(12.1.5) \qquad \qquad &\text{$c=\Delta_2(\beta_{{\ell},
i,1},\beta_{{\ell}, i,2})-\{(n_\ell-i) n_{\ell-1}\cdots n_2
-\beta_{{\ell},i,\ell}n_{\ell-2}n_{\ell-3}\cdots n_2$}   \\
& \quad \text{\qquad $
-\beta_{{\ell}i,\ell-1}n_{\ell-3}n_{\ell-4}\cdots n_2 -\cdots
-\beta_{{\ell},i,4}n_2-\beta_{{\ell},i,3}\}n_1\beta_{1,0,1}$}.
\endalign$$ \ms

So, for any integer $\ell\ge 3,$ it suffices to show that the above
integer $c$ of (12.1.5) can be equal to an integer $\xi_{\ell-2}>0$,
where $\xi_{\ell-2}$ is the {\rm $(\ell-2)$-th} element of a
positive sequence $\{\xi_j:j=1,2,\dots,\ell-2\}$ such that each
$\xi_j$ satisfies the following properties:
$$\align
(12.1.6) \quad \xi_0
&=\Delta_{\ell}(\beta_{{\ell},i,k})^{\ell}_{k=1}
 -(n_{\ell}-i)n_{\ell-1}\Delta_{\ell-1}(\beta_{\ell-1,0,k})^{\ell-1}_{k=1}>0,  \\
 \xi_1  &=\Delta_{\ell-1}(\beta_{{\ell},i, k})^{\ell-1}_{k=1}
 -\{(n_{\ell}-i)n_{\ell-1}-\beta_{{\ell},i,\ell}\}n_{\ell-2}
 \Delta_{\ell-2}(\beta_{\ell-2,0,k})^{\ell-2}_{k=1}>0,  \\
 \xi_2  =\Delta_{\ell-2}&(\beta_{{\ell},i,k})^{\ell-2}_{k=1}-
 \{(n_{\ell}-i)n_{\ell-1}n_{\ell-2}
 -\beta_{{\ell},i,\ell}n_{\ell-2}-\beta_{{\ell},i,\ell-1}\}
 n_{\ell-3}\Delta_{\ell-3}(\beta_{\ell-3,0,k})^{\ell-3}_{k=1}>0, \qquad\\
  & \ldots\ldots \\
  \xi_j &=\Delta_{\ell-j}(\beta_{{\ell},i,k})^{\ell-j}_{k=1}
  -\{(n_{\ell}-i)n_{\ell-1}\cdots n_{\ell-j}
  -\beta_{{\ell},i,\ell}n_{\ell-2}n_{\ell-3}\cdots n_{\ell-j}  \\
 &\qquad -\beta_{{\ell},i,\ell-1}n_{\ell-3}n_{\ell-4}
  \cdots n_{\ell-j}-\cdots -\beta_{{\ell},i,\ell-j+2}n_{\ell-j} \\
  &\qquad -\beta_{{\ell},i,\ell-j+1}\}
  \times n_{\ell-j-1}\Delta_{\ell-j-1}(\beta_{\ell-j-1,0,k})^{\ell-j-1}_{k=1}>0
\quad \text{for $3\le j\le \ell-2$}. \\
\endalign$$

Let $\ell\ge 3$ be chosen arbitrary. Now, we will show by the
induction method on the nonnegative integer $j\le \ell-2$ that
$\xi_j$ is positive for all $j$.

It is trivial by \text{\rm The 4-th
${\text{\rm{Cond}}}^{\text{{\rm(0)}}}$} in the assumption of Theorem
$12.0$ that $\xi_0>0$.

In order to prove that $\xi_1$ is positive, first of all, it is easy
to observe the following by the definition of
$\Delta_{\ell}(\beta_{{\ell},i,k})^{\ell}_{k=1}$:
$$\align
(12.1.7) \quad 0<\xi_0
&=\Delta_{\ell}(\beta_{{\ell},i,k})^{\ell}_{k=1}
 -(n_{\ell}-i)n_{\ell-1}\Delta_{\ell-1}(\beta_{\ell-1,0,k})^{\ell-1}_{k=1}\\
&=\beta_{{\ell},i,\ell}\Delta_{\ell-1}(\beta_{\ell-1,0,k})^{\ell-1}_{k=1}
+n_{\ell-1}\Delta_{\ell-1}(\beta_{{\ell},i,k})^{\ell-1}_{k=1}
  -(n_{\ell}-i)n_{\ell-1}\Delta_{\ell-1}(\beta_{\ell-1,0,k})^{\ell-1}_{k=1} \\
&=n_{\ell-1}\Delta_{\ell-1}(\beta_{{\ell},i,k})^{\ell-1}_{k=1}-((n_{\ell}-i)n_{\ell-1}
-\beta_{{\ell},i,\ell})\Delta_{\ell-1}(\beta_{\ell-1,0,k})^{\ell-1}_{k=1}.
\endalign$$

Since $\Delta_{\ell-1}(\beta_{\ell-1,0,k})^{\ell-1}_{k=1}
>n_{\ell-1}n_{\ell-2}\Delta_{\ell-2}(\beta_{\ell-2,0,k})^{\ell-2}_{k=1}$
by \text{\rm The 4-th ${\text{\rm{Cond}}}^{\text{{\rm(0)}}}$} in the
assumption of Theorem $12.0$, then the third inequality of (12.1.7)
implies that
$$\align
(12.1.8) \qquad
&n_{\ell-1}\Delta_{\ell-1}(\beta_{{\ell},i,k})^{\ell-1}_{k=1}
>((n_{\ell}-i)n_{\ell-1}-\beta_{{\ell},i,\ell})n_{\ell-1}n_{\ell-2}
\Delta_{\ell-2}(\beta_{\ell-2,0,k})^{\ell-2}_{k=1}, \qquad \qquad  \\
& \text{whether or not
$(n_{\ell}-i)n_{\ell-1}-\beta_{{\ell},i,\ell}>0$.}
\endalign$$

Dividing both sides on (12.1.8) by $n_{\ell-1}$, then
$$
\Delta_{\ell-1}(\beta_{{\ell},i,k})^{\ell-1}_{k=1}
>((n_{\ell}-i)n_{\ell-1}-\beta_{{\ell},i,\ell})n_{\ell-2}
\Delta_{\ell-2}(\beta_{\ell-2,0,k})^{\ell-2}_{k=1}, \tag 12.1.9
$$
which is equivalent to the fact that $\xi_1>0$.

By the induction assumption on the positive integer $j\le \ell-2$,
suppose we have shown that $\xi_j$ is positive with $1\le j\le
\ell-3$. To prove that $\xi_{j+1}$ is positive, for convenience of
notations, let $\xi_j$ of $(12.1.6)$ be written again in the form
$$
\align (12.1.10) \qquad \qquad  \xi_j
&=\Delta_{\ell-j}(\beta_{{\ell},i,k})^{\ell-j}_{k=1}
 -\omega_jn_{\ell-j-1}\Delta_{\ell-j-1}(\beta_{\ell-j-1,0,k})^{\ell-j-1}_{k=1}>0
 \quad \text{with} \qquad \qquad \quad \\
\omega_j &=(n_{\ell}-i)n_{\ell-1}\cdots
n_{\ell-j}-\beta_{{\ell},i,\ell}n_{\ell-2}n_{\ell-3}\cdots
n_{\ell-j}\\
& \quad -\beta_{{\ell},i,\ell-1}n_{\ell-3}n_{\ell-4}\cdots
n_{\ell-j}-\cdots
-\beta_{{\ell},i,\ell-j+2}n_{\ell-j}-\beta_{{\ell},i,\ell-j+1}.
\endalign$$

Now, by (12.1.10) and the definition of
$\Delta_{\ell-j}(\beta_{{\ell},i,k})^{\ell-j}_{k=1}$ only, it is
easy to prove that
$$\align
(12.1.11)  \quad 0<\xi_j
 &=\Delta_{\ell-j}(\beta_{{\ell},i,k})^{\ell-j}_{k=1}
 -\omega_jn_{\ell-j-1}\Delta_{\ell-j-1}(\beta_{\ell-j-1,0,k})^{\ell-j-1}_{k=1}
 \qquad \qquad \quad \\
 &=\beta_{{\ell},i,\ell-j}\Delta_{\ell-j-1}(\beta_{\ell-j-1,0,k})^{\ell-j-1}_{k=1}
 +n_{\ell-j-1}\Delta_{\ell-j-1}(\beta_{{\ell},i,k})^{\ell-j-1}_{k=1}\qquad \qquad \quad \\
 & \quad -\omega_jn_{\ell-j-1}\Delta_{\ell-j-1}(\beta_{\ell-j-1,0,k})^{\ell-j-1}_{k=1}\\
 &=n_{\ell-j-1}\Delta_{\ell-j-1}(\beta_{{\ell},i,k})^{\ell-j-1}_{k=1}
 -(\omega_jn_{\ell-j-1}-\beta_{{\ell},i,\ell-j})
\Delta_{\ell-j-1}(\beta_{\ell-j-1,0,k})^{\ell-j-1}_{k=1}.
\endalign$$

Since $\Delta_{\ell-j-1}(\beta_{\ell-j-1,0,k})^{\ell-j-1}_{k=1}
>n_{\ell-j-1}n_{\ell-j-2}\Delta_{\ell-j-2}(\beta_{\ell-j-2,0,k})^{\ell-j-2}_{k=1}$
by \text{\rm The 4-th ${\text{\rm{Cond}}}^{\text{{\rm(0)}}}$} in the
assumption of Theorem $12.0$, we get the following from the last
equality in $(12.1.11)$:
$$\align
& n_{\ell-j-1}\Delta_{\ell-j-1}(\beta_{{\ell},i,k})^{\ell-j-1}_{k=1}
 >(\omega_jn_{\ell-j-1}-\beta_{{\ell},i,\ell-j})
 \Delta_{\ell-j-1}(\beta_{\ell-j-1,0,k})^{\ell-j-1}_{k=1},
\quad \text{and so}  \tag 12.1.12\\
&n_{\ell-j-1}\Delta_{\ell-j-1}(\beta_{{\ell},i,k})^{\ell-j-1}_{k=1}>
(\omega_jn_{\ell-j-1}-\beta_{{\ell},i,\ell-j})n_{\ell-j-1}n_{\ell-j-2}
\Delta_{\ell-j-2}(\beta_{\ell-j-2,0,k})^{\ell-j-2}_{k=1}, \\
& \text{whether or not
$\omega_jn_{\ell-j-1}-\beta_{{\ell},i,\ell-j}>0$.}
\endalign$$

Dividing both sides of $(12.1.12)$ by $n_{\ell-j-1}$, then we get
$$\align
(12.1.13) \quad \quad
\Delta_{\ell-j-1}(\beta_{{\ell},i,k})^{\ell-j-1}_{k=1}
>(\omega_jn_{\ell-j-1}-\beta_{{\ell},i,{\ell-j}})n_{\ell-j-2}
\Delta_{\ell-j-2}(\beta_{\ell-j-2,0,k})^{\ell-j-2}_{k=1}. \qquad
\quad \quad
\endalign$$

Before we prove that $\xi_{j+1}>0$, then it is trivial to observe by
(12.1.10) that $\xi_{j+1}$ of (12.1.6) can be rewritten as follows:
$$\align
(12.1.14) \quad
\xi_{j+1}=\Delta_{\ell-j-1}(\beta_{{\ell},i,k})^{\ell-j-1}_{k=1}
-(\omega_jn_{\ell-j-1}-\beta_{{\ell},i,\ell-j})n_{\ell-j-2}
\Delta_{\ell-j-2}(\beta_{\ell-j-2,0,k})^{\ell-j-2}_{k=1}. \quad
\endalign$$

Then, it is clear by (12.1.13) that $\xi_{j+1}$ is positive for all
$j=1,2,\dots,\ell-3$. Since $c=\xi_{\ell-2}$ for some integer
$\ell\ge 3$ by (12.1.6), the proof is done. $\square$
\enddemo
\ms

\demo{\bf Proof of Sublemma 12.2} We prove (a), (b), (c) and (d),
respectively.

$\underline{\text{(a)}}$ \quad In preparation for the proof of an
equality in (12.2.1), it is clear by (2a) of \text{\rm The 2-th
${\text{\rm{Cond}}}^{\text{{\rm(0)}}}$} in the assumptions of this
theorem and by an inequality in (12.2.0) of this sublemma that the
following are true:

If $r=1$ then $g_1=g_1(y,z)$ can be written in the form
$$\align
g_1 &=\Sigma_{1,0}+\Sigma_{1,1}, \tag 12.2.2 \\
\text{where} \quad \Sigma_{1,0} &
=z^{n_1}+\ve_{1,0}y^{\beta_{1,0,1}} \quad \text{with $\ve_{1,0}=1$}, \\
\Sigma_{1,1} &=\sum_{\alpha,\beta\ge 0}
c^{(1)}_{\alpha,\beta}y^{\alpha}z^{\beta} \quad \text{with} \quad
n_1\alpha+\beta_{1,0,1}\beta>n_1\beta_{1,0,1}, \qquad
\endalign$$
where a unit $\ve_{1,0}=\ve_{1,0}(y,z)$ may be analytically assumed
to be one in $\BC\{y,z\}$ if necessary, and the
$c^{(1)}_{\alpha,\beta}$ are nonzero complex numbers for some
nonnegative integers $\alpha$ and $\beta$ such that
$n_1\alpha+\beta_{1,0,1}\beta>n_1\beta_{1,0,1}$, if exist.

Also, note by (12.2.2) that for any positive integer $\ell$,
$g_1^{\ell}$ can be written in the form
$$\align
 g_1^{\ell} &=(\Sigma_{1,0}+\Sigma_{1,1})^{\ell}
=(\Sigma_{1,0})^{\ell}+\sum^{\ell-1}_{k=1}\binom{\ell}{k}
(\Sigma_{1,0})^{k}(\Sigma_{1,1})^{\ell-k}+(\Sigma_{1,1})^{\ell} \tag 12.2.3\\
&=\Sigma^{(\ell)}_{1,0}+\Sigma^{(\ell)}_{1,1}
\quad \text{for any integer $\ell\ge 2$},  \\
\text{where write} \quad \Sigma^{(\ell)}_{1,0}
&=(\Sigma_{1,0})^{\ell}
=(z^{n_1}+\ve_{1,0}y^{\beta_{1,0,1}})^{\ell} \quad \text{with $\ve_{1,0}=1$}, \\
\text{\qquad and} \quad \Sigma^{(\ell)}_{1,1}
&=\sum^{\ell-1}_{k=1}\binom{\ell}{k}
(\Sigma_{1,0})^{k}(\Sigma_{1,1})^{\ell-k}+(\Sigma_{1,1})^{\ell}.
\qquad
\endalign$$

It is clear that for any monomial $y^{\alpha}z^{\beta}\in
(\Sigma_{1,0})^{\ell}$, $n_1\alpha+\beta_{1,0,1}\beta=
n_1\beta_{1,0,1}\ell$. So, to prove that
$n_1\alpha+\beta_{1,0,1}\beta> n_1\beta_{1,0,1}\ell$ for any
monomial $y^{\alpha}z^{\beta}\in g_1^{\ell}-(\Sigma_{1,0})^{\ell}$,
is equivalent to prove by (12.2.3) that  for any nonzero monomial
$y^{\alpha}z^{\beta}\in \Sigma^{(\ell)}_{1,1}$
$$
n_1\alpha+\beta_{1,0,1}\beta> n_1\beta_{1,0,1}\ell. \tag 12.2.4
$$

To prove $(12.2.4)$, it suffices to show that the following three
claims hold by using an equation of $\Sigma^{(\ell)}_{1,1}$ in
$(12.2.3)$: Let $\ell\ge 1$ be an arbitrary integer.
$$\align
(12.2.5) \quad  & \text{Claim(i)} \quad \text{For any monomial
$y^{\alpha}z^{\beta}\in (\Sigma_{1,0})^{\ell}$},  \quad
n_1\alpha+\beta_{1,0,1}\beta= n_1\beta_{1,0,1}\ell. \qquad \qquad \qquad \qquad\\
& \text{Claim(ii)} \quad \text{For any monomial
$y^{\alpha}z^{\beta}\in (\Sigma_{1,1})^{\ell}$},
\quad  n_1\alpha+\beta_{1,0,1}\beta> n_1\beta_{1,0,1}\ell. \\
&\text{Claim(iii)} \quad \text{For any monomial
$y^{\gamma}z^{\delta}\in (\Sigma_{1,0})^{k}(\Sigma_{1,1})^{\ell-k}$}, \\
& \qquad \qquad \quad n_1\gamma+\beta_{1,0,1}\delta>
n_1\beta_{1,0,1}k+n_1\beta_{1,0,1}(\ell-k)=n_1\beta_{1,0,1}\ell.
\qquad \qquad
\endalign$$
Note by (12.2.2) that the proof of three claims in (12.2.5) is
trivial, and so the proof of (12.2.4) is done. \ms

Now, for the proof of the sublemma, using equations in (12.2.2),
(12.2.3) and (12.2.4), it suffices to show that for any integer
$r\ge 2$, $g_r$ of $(12.2.1)$ can be generally represented in the
following form: Let $\ell\ge 2$ be an arbitrary integer.
$$\align
\text{\rm(12.2.6)(12.2.6.1)} \quad  g_r &=\Sigma_{r,0}+\Sigma_{r,1},   \\
\Sigma_{r,0} &=(\Sigma_{1,0})^{d_r} \quad \text{and
$\Sigma_{1,0}=z^{n_1}+\ve_{1,0}y^{\beta_{1,0,1}}$ with
$\ve_{1,0}=1$ and $d_r=n_2\cdots n_r$}, \qquad \qquad \\
\Sigma_{r,1} &=\sum_{\g,\de\ge 0} c^{(1)}_{\g,\de}y^{\g}z^{\de}
\quad \text{with
$n_1\g+\beta_{1,0,1}\de>n_1\beta_{1,0,1}d_r$}, \quad \text{and  then} \qquad\qquad \\
\text{\rm(12.2.6.2)} \quad g_r^{\ell}
&=\Sigma^{(\ell)}_{r,0}+\Sigma^{(\ell)}_{r,1} \quad
\text{for any integer $\ell\ge 2$},  \\
\Sigma^{(\ell)}_{r,0} &=(\Sigma_{1,0})^{d_r{\ell}}
=(\Sigma_{r,0})^{\ell},\\
\Sigma^{(\ell)}_{r,1} &=\sum_{\g,\de\ge 0}
c^{(\ell)}_{\g,\de}y^{\g}z^{\de} \quad \text{with} \quad
n_1\g+\beta_{1,1}\de>n_1\beta_{1,0,1}d_r\ell, \qquad \qquad
\endalign$$
where $\ve_{1,0}$ is assumed to be one in $\BC\{y,z\}$ if necessary,
and the $c^{(\ell)}_{\g,\de}$ are nonzero complex numbers for some
nonnegative integers $\g$ and $\de$ such that
$n_1\g+\beta_{1,0,1}\de>n_1\beta_{1,0,1}d_r\ell$.

For the induction proof of an equality in (12.2.6), it suffices to
consider the following two cases:

Case(I) $r=2$ and Case(II) $2<r$. \ms

$\underline{\text{\rm Case(I)}}$ Let $r=2$. Recall by (2b) of
\text{\rm The 2-th ${\text{\rm{Cond}}}^{\text{{\rm(0)}}}$} and (4b)
of \text{\rm The 4-th ${\text{\rm{Cond}}}^{\text{{\rm(0)}}}$} in the
assumptions of this theorem and by (12.1.2) of Sublemma $12.1$ that
$g_2$ can be written in the form
$$\align
(12.2.7)(12.2.7.1) \qquad \qquad
&\text{$g_2=g^{n_2}_1+\ve_{2,0}y^{\beta_{2,0,1}}z^{\beta_{2,0,2}}
+c_2\sum^{n_2-1}_{i=1}\ve_{2,i}y^{\beta_{2,i,1}}z^{\beta_{2,i,2}}g^{i}_1$
\quad and} \qquad\qquad \\
(12.2.7.2)\qquad \qquad
&\text{$\Delta^{\sharp}_2(\beta_{2,i,1},\beta_{2,i,2})
=\Delta_2(\beta_{2,i,1},\beta_{2,i,2})> n_1\beta_{1,0,1}(n_2-i)$
\quad {on} \quad $g_2$}.
\endalign$$

In order to prove an equality in (12.2.6) for $r=2$, since
$n_1\beta_{2,0,1}+\beta_{1,0,1}\beta_{2,0,2}> n_2n_1\beta_{1,0,1}$
by (4b) of \text{\rm The 4-th
${\text{\rm{Cond}}}^{\text{{\rm(0)}}}$} in the assumption of Theorem
$12.0$, then it suffices to show that the following two claims hold
by using equations in $(12.2.7)$ with (12.2.4):
$$\align
(12.2.8) \quad & \text{Claim(i)} \quad \text{For any monomial
$y^{\alpha}z^{\beta}\in (g_1)^{n_2}-(\Sigma_{1,0})^{n_2}$, \quad
$n_1\alpha+\beta_{1,0,1}\beta>n_1\beta_{1,0,1}n_2$.} \qquad\\
&\text{Claim(ii)} \quad \text{For any monomial
$y^{\gamma}z^{\delta}\in y^{\beta_{2,i,1}}z^{\beta_{2,i,2}}g^{i}_1$,
\quad $n_1\gamma+\beta_{1,0,1}\delta\ge $} \\
& \qquad \qquad \quad
\Delta^{\sharp}_2(\beta_{2,i,1},\beta_{2,i,2})+ n_1\beta_{1,0,1}i>
n_1\beta_{1,0,1}(n_2-i)+n_1\beta_{1,0,1}i=n_1\beta_{1,0,1}n_2.
\endalign$$

Thus, the proof of two claims in (12.2.8) is trivial by (12.2.4) and
(12.2.7.2), and so the proof of (12.2.6.1) is done for $r=2$. Also,
if $r=2$, the proof of (12.2.6.2) is trivial by the same method as
we have used in the proof for (12.2.4) with (12.2.3). Thus, the
proof of Case(I) is done. \ms

$\underline{\text{\rm Case(II)}}$ Let $2<r$. Now, suppose we have
proved by induction assumption on the positive integer $j<r$ that
the representation of $g_j$ in (12.2.6) is true for $2\le j<r$.
Then, recall by (2b) of \text{\rm The 2-th
${\text{\rm{Cond}}}^{\text{{\rm(0)}}}$} and (4b) of \text{\rm The
4-th ${\text{\rm{Cond}}}^{\text{{\rm(0)}}}$} in the assumptions of
this theorem and by (12.1.2) of Sublemma $12.1$ that $g_{j+1}$ can
be rewritten as follows:
$$\align
(12.2.9) \qquad \qquad   &
g_{j+1}=g^{n_{j+1}}_{j}+\ve_{j+1,0}y^{\beta_{j+1,0,1}}z^{\beta_{j+1,0,2}}
g^{\beta_{j+1,0,3}}_1\cdots g^{\beta_{j+1,0,j+1}}_{j-1} \qquad \qquad \qquad \qquad \qquad  \\
& \qquad \quad
+c_{j+1}\sum^{n_{j+1}-1}_{i=1}\ve_{j+1,i}y^{\beta_{j+1,i,1}}
z^{\beta_{j+1,i,2}}g^{\beta_{j+1,i,3}}_1\cdots
g^{\beta_{j+1,i,j+1}}_{j-1}g^i_{j}, \\
& \Delta^{\sharp}_{j+1}(\beta_{j+1,i,k})^{j+1}_{k=1} >
n_1\beta_{1,0,1}n_2n_3\cdots
n_{j}(n_{j+1}-i) \quad \text{on $g_{j+1}$}, \\
\endalign$$
where $\ve_{j+1,0}$ and $\ve_{j+1,i}$ are units in $\BC\{y,z\}$, and
the $c_{j+1}$ are nonzero complex numbers.

First of all, applying the induction assumption on $g_j$, then for
any $k=1,2,\dots,j$ and any integer $\ell>0$, suppose we have shown
that the following are true:
$$\align
\text{\rm(12.2.10)(12.2.10.1)} \quad  g_k &=\Sigma_{k,0}+\Sigma_{k,1},   \\
\Sigma_{k,0} &=(\Sigma_{1,0})^{d_k} \quad \text{and
$\Sigma_{1,0}=z^{n_1}+\ve_{1,0}y^{\beta_{1,0,1}}$ with
$\ve_{1,0}=1$ and $d_k=n_2\cdots n_k$}, \qquad \qquad \\
\Sigma_{k,1} &=\sum_{\g,\de\ge 0} c^{(1)}_{\g,\de}y^{\g}z^{\de}
\quad \text{with
$n_1\g+\beta_{1,0,1}\de>n_1\beta_{1,0,1}d_k$}, \quad \text{and  then} \qquad\qquad \\
\text{\rm(12.2.10.2)} \quad g_k^{\ell}
&=\Sigma^{(\ell)}_{k,0}+\Sigma^{(\ell)}_{k,1} \quad
\text{for any integer $\ell\ge 2$},  \\
\Sigma^{(\ell)}_{k,0} &=(\Sigma_{1,0})^{d_k{\ell}}
=(\Sigma_{k,0})^{\ell},\\
\Sigma^{(\ell)}_{k,1} &=\sum_{\g,\de\ge 0}
c^{(\ell)}_{\g,\de}y^{\g}z^{\de} \quad \text{with} \quad
n_1\g+\beta_{1,1}\de>n_1\beta_{1,0,1}d_k\ell, \qquad \qquad
\endalign$$
where $\ve_{1,0}$ is assumed to be one in $\BC\{y,z\}$ if necessary,
and the $c^{(\ell)}_{\g,\de}$ are nonzero complex numbers for some
nonnegative integers $\g$ and $\de$ such that
$n_1\g+\beta_{1,0,1}\de>n_1\beta_{1,0,1}d_k\ell$.

First, it is clear by (12.2.10) that for any nonzero monomial
$y^{\alpha_k}z^{\beta_k}\in g^{\ell}_k
=\Sigma^{(\ell)}_{k,0}+\Sigma^{(\ell)}_{k,1}$,
$$\align
n_1\alpha_k+\beta_{1,0,1}\beta_k\ge n_1\beta_{1,0,1}d_k\ell \quad
\text{on $g_k$.} \tag 12.2.11
\endalign$$

In order to prove that (12.2.6) is true on $g_{j+1}$, it is clear by
(12.2.9) and (12.2.11) that any nonzero monomial
$y^{\gamma}z^{\de}\in g_{j+1}-g^{n_{j+1}}_{j}$ can be represented as
follows:
$$\align
(12.2.12) \qquad \qquad & y^{\gamma}z^{\de}=y^{\beta_{j+1,i,1}}
z^{\beta_{j+1,i,2}}\Pi^{j}_{k=1}y^{\alpha_k}z^{\beta_k} \quad
\text{and}  \\
 &\text{$y^{\alpha_k}z^{\beta_k}\in
g^{\beta_{j+1,i,k+2}}_k$ with $n_1\alpha_k+\beta_{1,0,1}\beta_k\ge
n_1\beta_{1,0,1}d_k\beta_{j+1,i,k+2}$.} \qquad \qquad \qquad
\endalign$$

So, by (12.2.11), (12.2.12) and Sublemma $12.1$, for any nonzero
monomial $y^{\gamma}z^{\de}\in g_{j+1}-g^{n_{j+1}}_{j}$,
$$\align
(12.2.13) \qquad \qquad & n_1\gamma+\beta_{1,0,1}\delta   \\
\ge & n_1\beta_{j+1,i,1}+\beta_{101}\beta_{j+1,i,2}+\sum^{j-1}_{k=1}
(n_1\alpha_k+\beta_{1,0,1}\beta_k)+n_1\beta_{1,0,1}d_j{i}
\qquad \qquad \qquad \qquad\\
\ge
&n_1\beta_{j+1,i,1}+\beta_{1,0,1}\beta_{j+1,i,2}+\sum^{j-1}_{k=1}
(n_1\beta_{1,0,1}
d_k\beta_{j+1,i,k+2})+n_1\beta_{1,0,1}d_j{i} \qquad\\
=&\Delta^{\sharp}_{j+1}(\beta_{j+1,i,k})^{j+1}_{k=1}+
n_1\beta_{1,0,1}d_j{i} \\
>& n_1\beta_{1,0,1}d_j(n_{j+1}-i)+n_1\beta_{1,0,1}d_j{i}
=n_1\beta_{1,0,1}d_jn_{j+1}=n_1\beta_{1,0,1}d_{j+1}.
\endalign$$

Since any nonzero monomial $y^{\gamma}z^{\de}\in
g^{n_{j+1}}_{j}-\Sigma^{(n_{j+1})}_{j,0}$ implies that
$n_1\gamma+\beta_{1,0,1}\delta
> n_1\beta_{1,0,1}d_{j+1}$, then
if $r=j+1$, the proof of (12.2.6.1) is done. Also, if r=j+1, then
the proof of (12.2.6.2) is trivial by the same method as we have
used in the proof for (12.2.4) with (12.2.3). Then, the proof of
Case(II) is done. So, we finished the proof of (a). \ms

$\underline{\text{(b)}}$ \quad To prove {\rm(b1)}, it is enough to
consider $g_r(0,z)$ from $g_r(y,z)$ of (12.2.1). Then
$$\align
g_r(0,z) &=z^{n_1n_2\cdots n_r}+\sum c^{(r)}_{0,\beta}z^{\beta}
\quad\text{with} \tag 12.2.14\\
 \beta_{1,0,1}\beta &>n_1\beta_{1,0,1}n_2n_3\cdots n_r.
\endalign$$
Thus, $\beta>n_1n_2\cdots n_r$, and so it is done. Also, the proof
of (b2) can be done similarly. \ms

$\underline{\text{(c)}}$ \quad To prove (c1), suppose that
$n_1<\beta_{1,0,1}$. By (a),
$\beta_{1,0,1}(\alpha+\beta)>n_1\alpha+\beta_{1,0,1}\beta>n_1\beta_{1,0,1}n_2n_3\cdots
n_r$ implies that $\alpha+\beta>n_1\prod^r_{k=2}n_k$. Thus, the
proof of (c1) is done. Similarly, (c2) can be proved.

Thus, the proof of this sublemma is finished. $\square$
\enddemo
\ms

\demo{\bf Proof of Sublemma 12.3} To prove (12.3.3) for any $r\ge
2$, it is enough to show by (12.3.1) and (12.3.2) that the following
equation in (12.3.4) is nonnegative:
$$\align
(12.3.4) \quad
&\Omega^{\sharp}_r(\beta_{r,i,k})^r_{k=1}-bn_1n_2\cdots
n_{r-1}(n_r-i)=a\beta_{r,i,1}+bD\ge 0  \quad \text{with} \\
&  D=\beta_{r,i,2}+n_1\beta_{r,i,3}+n_1n_2\beta_{r,i,4}+\cdots
+n_1n_2\cdots
n_{r-2}\beta_{r,i,r}-n_1n_2\cdots n_{r-1}(n_r-i),\\
\endalign$$
by definition of $\Omega^{\sharp}_r(\beta_{r,i,k})^r_{k=1}$ in
(12.3.2) where $a>0$ and $b\ge 0$.

Then, it remains to prove by (12.3.4) that $a\beta_{r,i,1}+bD\ge 0$
depending on $D$, and so it suffices to consider the following two
cases:

\noindent$\underline{\text{\rm Case(i)}}$ \ Let $D\ge 0$. It is
clear that $a\beta_{r,i,1}+bD\ge 0$ because $a>0$, and also $b\ge 0$
with $\beta_{r,i,1}$ nonnegative. Thus, the proof of Case(i) is
done. \ms

\noindent$\underline{\text{\rm Case(ii)}}$ \ Let $D<0$. In order to
prove that $a\beta_{r,i,1}+bD\ge 0$, first note that the inequality
$\Delta^{\sharp}_r(\beta_{r,i,k})^r_{k=1}-n_1\beta_{1,0,1}n_2n_3\cdots
n_{r-1}(n_r-i)>0$ of Sublemma $12.1$ can be equivalently rewritten
as follows:
$$\align
(12.3.5) \quad 0&<\Delta^{\sharp}_r(\beta_{r,i,k})^r_{k=1}
-n_1\beta_{1,0,1}n_2n_3\cdots n_{r-1}(n_r-i)  \\
&=
n_1\beta_{r,i,1}+\beta_{1,0,1}\beta_{r,i,2}+n_1\beta_{1,0,1}\beta_{r,i,3}
+n_1\beta_{1,0,1}n_2\beta_{r,i,4} \\
&\quad +\cdots +n_1\beta_{1,0,1}n_2\cdots n_{r-2}\beta_{r,i,r}
-n_1\beta_{1,0,1}n_2n_3\cdots n_{r-1}(n_r-i)  \quad \text{by (12.1.1)} \\
&= n_1\beta_{r,i,1}+\beta_{1,0,1}D, \quad \text{where} \\
&  D=\beta_{r,i,2}+n_1\beta_{r,i,3}+\cdots +n_1n_2\cdots
n_{r-2}\beta_{r,i,r}-n_1n_2\cdots n_{r-1}(n_r-i) \ \text{by (12.3.4)}. \qquad\\
\endalign$$

Since $-D>0$ and $n_1\ge 2>0$, then the inequality
$n_1\beta_{r,i,1}+\beta_{1,0,1}D>0$ in $(12.3.5)$ can be
equivalently represented as follows:
$$
 \frac{\beta_{r,i,1}}{-D} > \frac{\beta_{1,0,1}}{n_1}.  \tag
12.3.6
$$

Also, $a\beta_{1,0,1}-bn_1=1$ with $a>0$ implies that
$\dfrac{\beta_{1,0,1}}{n_1}
> \dfrac{b}{a}$. Therefore, we proved by $(12.3.6)$ that
$\dfrac{\beta_{r,i,1}}{-D}
>\dfrac{b}{a}$, that is, $a\beta_{r,i,1}+bD>0$. Thus, the proof of
Case(ii) is done.

Therefore, we showed by Case(i) and Case(ii) that the equation in
(12.3.4) is nonnegative, and so the proof of this sublemma is
finished. $\square$
\enddemo
\ms

\demo{\bf Proof of Sublemma 12.4} Following the same assumptions and
notations as in Sublemma 12.1, Sublemma 12.2 and Sublemma 12.3, then
(a) of Sublemma 12.2 is true, and so for each $j=2,3,\dots,r$, \
$g_j=g_j(y,z)$ of (12.2.1) can be easily rewritten as follows:
$$\split
\noindent(12.4.6) \qquad \qquad  g_j&=(z^{n_1}+\ve_{1,0}
y^{\beta_{1,0,1}})^{d_j}+\sum_{\alpha,\beta\ge 0}
c^{(j)}_{\alpha,\beta}y^{\alpha}z^{\beta} \quad \text{with
$\ve_{1,0}=1$ and}
\qquad \qquad \qquad \\
&  \text{with \quad
$n_1\alpha+\beta_{1,0,1}\beta>n_1\beta_{1,0,1}d_j$, \quad where
$d_j=n_2n_3\cdots n_j$,}
\endsplit$$
where $\ve_{1,0}=\ve_{1,0}(y,z)$ is assumed to be one in
$\BC\{y,z\}$, and the $c^{(j)}_{\alpha,\beta}$ are nonzero complex
numbers for some nonnegative integers $\alpha$ and $\beta$ such that
$n_1\alpha+\beta_{1,0,1}\beta>n_1\beta_{1,0,1}d_j$. \ms

For the proof of this sublemma, it suffices to follow Step(I) and
Step(II) in order:

{\rm Step(I)} First, we will show how to apply Theorem $3.6$ to the
proof of (a), (b), and $(d_2)$ and $(d_3)$ of (d) in this sublemma.

{\rm Step(II)} Next, the remaining part of this sublemma will be
proved computationally. \ms

\noindent{\bf Step(I)} In preparation for the proof of (a), (b),
$(d_2)$ and $(d_3)$, it is clear that the equation of $g_j$ of
(12.4.6) satisfies the same kind of properties as $f$ does in the
assumption of Theorem $3.6$, which can be represented as follows:
\ms

$\underline{\text{$g_j$ of (12.4.6) satisfies the same kind of
assumption as in Theorem 3.6}}$ \quad Let $V(g_0)=
\{(y,z):g_0(y,z)=0\}$, $V(f)=\{(y,z): f(y,z)=0\}$ and $V(G)=\{(y,z):
G(y,z)=0\}$ be analytic varieties at $(0,0)$ in $\BC^2$, each of
which is written respectively as follows: For brevity of notation,
substitute $g_j$ of (12.4.6) by $f$, for an application of Theorem
$3.6$.
$$\align
(12.4.7) \qquad \qquad g_0 &=z^{n_1}+\ve_{1,0}y^{\beta_{1,0,1}}
\quad \text{with $\ve_{1,0}=1$}, \qquad \qquad \qquad \qquad \\
f&={g_0}^{d_j}+\sum_{\alpha,\beta\ge
0}c^{(j)}_{\alpha,\beta}y^{\alpha}z^{\beta} \quad \text{with} \quad
n_1\alpha+\beta_{1,0,1}\beta>n_1\beta_{1,0,1}d_j, \qquad \qquad \qquad \qquad \\
F&=y^{\delta_1}z^{\delta_2}f, \\
G&=y^{\gamma}g_0, \\
\endalign
$$
satisfying the properties {\rm(i)}, {\rm(ii)}, {\rm(iii)},
{\rm(iv)}, {\rm(v)} and {\rm(vi)}:

\roster \item "(i)" $\gcd(n_1,\beta_{1,0,1})=1$ with $n_1\ge 2$ and
$\beta_{1,0,1}\ge 1$.

\item "(ii)" $d_j=n_2n_3\cdots n_j$ is a positive integer with
$d_j\ge 2$, and $d_1=1$ if necessary.

\item "(iii)" $\ve_{1,0}$ is assumed to be one in $\BC\{y,z\}$, and
the $c_{\alpha,\beta}$ are nonzero complex numbers for some
nonnegative integers $\alpha$ and $\beta$ such that
$n_1\alpha+k_1\beta> n_1\beta_{1,0,1}d_j$, if exist.

\item "(iv)" Assume that $V(f)$ has an isolated singular point at
the origin as a reduced variety.

\item "(v)" If $\beta_{1,0,1}=1$, then $\gamma=1$, and if
$\beta_{1,0,1}\ge 2$, then $\gamma=0$.

\item "(vi)" In addition, assume that each $\delta_i$ is either a
positive integer or $0$ for $i=1,2$, as far as $V(F)$ has an
isolated singular point at the origin as a reduced variety, even if
$d_j\ge 1$.

\endroster \ms

So, $\underline{\text{$g_j$}}$ of (12.4.6) has the same kind of
conclusion as in Theorem $3.6$, up to change of notations: \ms

$\underline{\text{The same kind of conclusion as in Theorem 3.6}}$
\quad As we have seen in Theorem $3.6$, let
$\tau_m=\pi_1\circ\pi_2\circ\cdots\circ\pi_m:M^{(m)}\to\BC^2$ be the
compositions of a finite number $m$ of successive blow-ups $\pi_i$
which is needed to get the standard resolution of the singular point
of $V(G)=V(y^{\gamma}g_0)$. Therefore, by the conclusion of Theorem
$3.6$, there is nothing to prove for (a), (b), (d2) and (d3) in this
sublemma. \ms

\noindent{\bf Step(II)} To finish the proof of this sublemma, it
remains to prove (c) and the remaining part $(d_1)$ of (d). Now, we
will use the same kind of notations and properties as in (a), (b),
(d2) and (d3) as follows:

For (a), (b) and (d2),  along $v=0$ $\tau_m:M^{(m)}\to\BC^2$ as a
composition of analytic mappings and $(f\circ\tau_m)_{total}$ can be
rewritten in the following form: Note that $2\le j\le r$.
$$\align
(12.4.8) \qquad \qquad \tau_m(v,u)&=(y,z)=(v^{n_1}u^a,v^{\beta_{1,0,1}}u^b),  \\
(f\circ\tau_m)_{total}&=(f\circ\tau_m)(v,u)
=v^{e_{j,m}}u^{\rho_{j,m}}(f\circ\tau_m)_{proper}\quad \text{with \  $g_j=f$},\\
(f\circ\tau_m)_{proper}&=(1+\ve_{1,0}u)^{d_j}+\sum_{\alpha,\beta\ge
0}c^{(j)}_{\alpha,\beta}v^{n_1\alpha+\beta_{1,0,1}\beta-n_1\beta_{1,0,1}d_j}
u^{a\alpha+b\beta-bn_1d_j}, \qquad \qquad \\
\endalign
$$
where \roster \item "(i)" $a$ and $b$ are some nonnegative integers
such that $a\beta_{1,0,1}-bn_1=1$ and $\ve_{1,0}$ is assumed to be
one in $\BC\{y,z\}$,

\item "(ii)" $e_{j,m}=n_1\beta_{1,0,1}d_j$ and $\rho_{j,m}=bn_1d_j$
and $\rho_{\alpha,\beta}=a\alpha+b\beta-bn_1d_j\ge 0$,

\item "(iii)" $E_m=\{v=0\}$ is defined by the $m-th$ exceptional
curve of the first kind.

\item "(iv)" $V^{(m)}(g_j)\cap(\cup^m_{i=1}E_i)=V^{(m)}(g_j)\cap
E_m=\{(v,1+\ve_{1,0} u)=(0,0)\}$ for any $j=2,\dots,r$.
\endroster \ms

For (d3), after $m$ iterations of blow-ups, denoted by $\tau_m$, we
have the following consequences:
$$\align
(12.4.9) \qquad \qquad &\text{(i) \quad If $\beta_{1,0,1}=1$, then
$f\in \text{\rm the  type}[0]$ under
$\tau_m$, and }\\
&\text{\qquad \quad if $\beta_{1,0,1}\ge 2$, then $f\in \text{\rm
the type}[1]$ under
$\tau_m$.} \\
&\text{(ii) \quad Whether $\beta_{1,0,1}=1$ or $\beta_{1,0,1}\ge 2$,
then $F\in \text{\rm the  type}[1]$ under $\tau_m$.} \qquad \qquad
\endalign$$

$\underline{\text{Remark 12.4.1}}$

(a) In the assumption of Theorem $3.6$, the construction for
$G(y,z)=z^{\gamma}g_0$ with $g_0=z^{n_1}+y^{k_1}$ was defined as
follows: Note that $\gcd(n_1,k_1)=1$.

Let $1\le n_1<k_1$, and if $n_1=1$, then $\gamma=1$, and if $n_1\ge
2$, then $\gamma=0$.

(b) In the conclusion of Theorem $3.6$, whether $n_1=1$ or $2\le
n_1<k_1$, or $2\le k_1 <n_1$, there are some nonnegative integers
$a$ and $b$ such that $bn_1-ak_1=1$ because of $[\text{\rm{I}}]$ and
$[\text{\rm{II}}]$ in Theorem $3.6$. \ms

In preparation for the proof of (c), apply the above conclusion with
(12.4.8), to $g_{j+1}={g_0}^{d_{j+1}}+\sum
c^{(j+1)}_{\alpha,\beta}y^{\alpha}z^{\beta}$ in $(12.4.6)$ where
$g_0=z^{n_1}+\ve_{1,0} y^{\beta_{1,0,1}}$ with $\ve_{1,0}=1$. Then
for any $j=1,2,\dots,r-1$, we have the following:
$$\align
(12.4.10) \quad (g_{j+1}\circ \tau_m)_{total}
&=v^{e_{j+1,m}}u^{\rho_{j+1,m}}(g_0\circ
\tau_m)^{d_{j+1}}_{proper}+\sum_{\alpha,\beta\ge 0}
c^{(j+1)}_{\alpha,\beta}(v^{n_1}u^a)^{\alpha}(v^{\beta_{1,0,1}}u^b)^{\beta}
\qquad \\
&=v^{e_{j+1,m}}u^{\rho_{j+1,m}}\{(g_0\circ
\tau_m)^{d_{j+1}}_{proper}
\\
& \quad +\sum_{\alpha,\beta\ge 0}
c^{(j+1)}_{\alpha,\beta}v^{n_1\alpha+\beta_{1,0,1}\beta-n_1\beta_{1,0,1}d_{j+1}}
u^{a\alpha+b\beta-bn_1d_{j+1}}\} \\
&=v^{e_{j+1,m}}u^{\rho_{j+1,m}}(g_{j+1}\circ \tau_m)_{proper},
\endalign$$
where $d_{j+1}=n_2n_3\cdots n_{j+1}$,
$e_{j+1,m}=n_1\beta_{1,0,1}d_{j+1}$ and $\rho_{j+1,m}=bn_1d_{j+1}$,
noting by (12.4.7) and (ii) of (12.4.8) that
$n_1\alpha+\beta_{1,0,1}\beta-n_1\beta_{1,0,1}d_{j+1}>0$ and
$a\alpha+b\beta-bn_1d_{j+1}\ge 0$. \ms

On the other hand, recall by (2b) of \text{\rm The 2-th
${\text{\rm{Cond}}}^{\text{{\rm(0)}}}$} in the assumption of the
theorem and by (12.4.6) that
$$\align
&g_{j+1}=g_j^{n_{j+1}}+\ve_{j+1,0}y^{\beta_{j+1,0,1}}z^{\beta_{j+1,0,2}}
g_1^{\beta_{j+1,0,3}}\cdots g_{j-1}^{\beta_{j+1,0,j+1}} \tag
12.4.11 \\
& \qquad \quad +c_{j+1}\sum^{n_{j+1}-1}_{i=1}\ve_{j+1,i}
y^{\beta_{j+1,i,1}} z^{\beta_{j+1,i,2}} g^{\beta_{j+1,i,3}}_1\cdots
g^{\beta_{j+1,i,j+1}}_{j-1}g^i_{j}, \\
& g_j =g_0^{n_2\cdots n_j} +\sum_{\gamma,\delta\ge 0}
c^{(j)}_{\gamma,\delta}y^{\gamma}z^{\delta}, \\
& g_0 =z^{n_1}+\ve_{1,0}y^{\beta_{1,0,1}} \quad \text{with $\ve_{1,0}=1$}, \\
where \quad \quad \text{\rm(i)} \quad
&\Delta_2(\gamma,\delta)=n_1\gamma+\beta_{1,0,1}\delta
>n_1\beta_{1,0,1}n_2\cdots n_j
\qquad \qquad \quad  \text{by (12.4.6)}, \\
\text{\rm(ii)} \quad
&\Delta^{\sharp}_{j+1}(\beta_{j+1,i,k})^{j+1}_{k=1}>n_1\beta_{1,0,1}n_2n_3\cdots
n_j(n_{j+1}-i) \qquad \text{by (12.1.2)}, \\
\text{\rm(iii)} \quad
&\Omega^{\sharp}_{j+1}(\beta_{j+1,i,k})^{j+1}_{k=1}\ge
bn_1n_2n_3\cdots n_j(n_{j+1}-i)
\qquad \qquad   \text{by (12.3.3)}, \\
\text{\rm(iv)} \quad &\text{$\ve_{j+1,i}$ is a unit in $\BC\{y,z\}$
\qquad for $0\le i\le n_{j+1}-1$}.
\endalign$$

Now, apply (12.4.8) or (12.4.10), to $g_{j+1}$ in $(12.4.11)$.

Then, we have the following:
$$\align
(12.4.12)
&\quad (g_{j+1}\circ \tau_m)_{total} \\
 &=((g_j\circ
\tau_m)(v,u))^{n_{j+1}} +{(\ve_{j+1,0}\circ \tau_m)(v,u)}{((y\circ
\tau_m)(v,u))^{\beta_{j+1,0,1}}}  \qquad \qquad \\
&  \quad \times {((z\circ
\tau_m)(v,u))^{\beta_{j+1,0,2}}}{(({g_1}\circ
\tau_m)(v,u))^{\beta_{j+1,0,3}}}\cdots
{((g_{j-1}\circ \tau_m)(v,u))^{\beta_{j+1,0,j+1}}} \qquad \qquad \\
& \quad +c_{j+1}\sum^{n_{j+1}-1}_{i=1}\{{(\ve_{j+1,i}\circ
\tau_m)(v,u)} {((y\circ \tau_m)(v,u))^{\beta_{j+1,i,1}}} {((z\circ
\tau_m)(v,u))^{\beta_{j+1,i,2}}} \\
& \quad \times {(({g_1}\circ \tau_m)(v,u))^{\beta_{j+1,i,3}}}\cdots
{((g_{j-1}\circ \tau_m)(v,u))^{\beta_{j+1,i,j+1}}}
{((g_{j}\circ \tau_m)(v,u))^{i}}\}, \\
 &=\{v^{e_{j,m}}u^{\rho_{j,m}}(g_{j}\circ
 \tau_m)_{proper}\}^{n_{j+1}}
 +{\ve'_{j+1,0}}{(v^{n_1}u^a)^{\beta_{j+1,0,1}}} {(v^{\beta_{1,0,1}}u^b)^{\beta_{j+1,0,2}}} \\
& \quad \times \{v^{n_1\beta_{1,0,1}}u^{bn_1}(g_{1}\circ
\tau_m)_{proper}\}^{\beta_{j+1,0,3}}
\{v^{e_{2,m}}u^{\rho_{2,m}}(g_{2}\circ
\tau_m)_{proper}\}^{\beta_{j+1,0,4}}\cdots \\
& \quad \times \{v^{e_{j-1,m}}u^{\rho_{j-1,m}}
(g_{j-1}\circ \tau_m)_{proper}\}^{\beta_{j+1,0,j+1}} \\
& \quad
+c_{j+1}\sum^{n_{j+1}-1}_{i=1}\{{\ve'_{j+1,i}}{(v^{n_1}u^a)^{\beta_{j+1,i,1}}}
{(v^{\beta_{1,0,1}}u^b)^{\beta_{j+1,i,2}}} \\
& \quad \times [v^{n_1\beta_{1,0,1}}u^{bn_1}(g_{1}\circ
\tau_m)_{proper}]^{\beta_{j+1,i,3}}
[v^{e_{2,m}}u^{\rho_{2,m}}(g_{2}\circ
\tau_m)_{proper}]^{\beta_{j+1,i,4}}\cdots \\
& \quad \times [v^{e_{j-1,m}}u^{\rho_{j-1,m}} (g_{j-1}\circ
\tau_m)_{proper}]^{\beta_{j+1,i,j+1}} [v^{e_{j,m}}u^{\rho_{j,m}}
(g_{j}\circ \tau_m)_{proper}]^{i} \}\\
&=v^{e_{j+1,m}}u^{\rho_{j+1,m}}\{ (g_{j}\circ
\tau_m)^{n_{j+1}}_{proper}+\ve'_{j+1,0}
v^{\Delta^{\sharp}_2(\beta_{j+1,0,1},\beta_{j+1,0,2})-e_{j+1,m}} \\
& \quad \times
u^{\Omega^{\sharp}_2(\beta_{j+1,0,1},\beta_{j+1,0,2})-b
n_1n_2n_3\cdots n_{j+1}}(g_{1}\circ
\tau_m)_{proper}^{\beta_{j+1,0,3}}\cdots (g_{j-1}\circ
\tau_m)^{\beta_{j+1,0,j+1}}_{proper} \\
& \quad +c_{j+1}\sum^{n_{j+1}-1}_{i=1}{\ve'_{j+1,i}}
v^{\Delta^{\sharp}_{j+1}(\beta_{j+1,i,k})^{j+1}_{k=1})-e_{j+1,m}+e_{j,m}{i}}
u^{\Omega^{\sharp}_{j+1}(\beta_{j+1,i,k})^{j+1}_{k=1}-bn_1n_2\cdots n_j(n_{j+1}-i)} \\
& \quad \times (g_{1}\circ \tau_m)_{proper}^{\beta_{j+1,i,3}}\cdots
(g_{j-1}\circ
\tau_m)^{\beta_{j+1,i,j+1}}_{proper}(g_{j}\circ \tau_m)^{i}_{proper}\} \\
&=v^{e_{j+1,m}}u^{\rho_{j+1,m}}(g_{j+1}\circ \tau_m)_{proper} \quad
\text{because of the following two equations in (12.4.13)},
\endalign$$
where

(a) $\ve'_{j+1,i}={(\ve_{j+1,i}\circ \tau_m)(v,u)}$ is a unit in
$\BC\{v,1+u\}$ for $0\le i<n_{j+1}$ and

(b) noting by (12.4.7) and (12.4.8) that $d_{j}=n_2n_3\cdots n_{j}$,
$e_{j,m}=n_1\beta_{1,0,1}d_{j}$ and $\rho_{j,m}=bn_1d_{j}$ for
$j=2,3,\dots,r$, then $d_{j}n_{j+1}=d_{j+1}$,
$e_{j,m}n_{j+1}=e_{j+1,m}$ and $\rho_{j,m}n_{j+1}=\rho_{j+1,m}$. \ms

Thus, the proof for the representation in (12.4.12) just follows
from Equation(1) and Equation(2) of (12.4.13):
$$\align
&\text{$\underline{\text{\rm Equation(1) of (5.4.13).}}$} \\
&\text{$(n_1\beta_{j+1,i,1}+\beta_{1,0,1}\beta_{j+1,i,2})+n_1\beta_{1,0,1}\beta_{j+1,i,3}
+e_{2,m}\beta_{j+1,i,4}+\cdots+e_{j-1,m}\beta_{j+1,i,j+1}$} \qquad
\\
=&\text{$\Delta_2(\beta_{j+1,i,1},\beta_{j+1,i,2})+n_1\beta_{1,0,1}\beta_{j+1,i,3}
+n_1\beta_{1,0,1}d_2\beta_{j+1,i,4}+\cdots+n_1\beta_{1,0,1}d_{j-1}\beta_{j+1,i,j+1}
$}\\
=&\text{$\Delta_2(\beta_{j+1,i,1},\beta_{j+1,i,2})+n_1\beta_{1,0,1}\beta_{j+1,i,3}
+n_1\beta_{1,0,1}n_2\beta_{j+1,i,4}+\cdots+n_1\beta_{1,0,1}n_2
\cdots
n_{j-1}\beta_{j+1,i,j+1}$} \\
=&\text{$\Delta^{\sharp}_{j+1}(\beta_{j+1,i,k})^{j+1}_{k=1}>e_{j+1,m}-ie_{j,m}
=n_1\beta_{1,0,1}d_{j}(n_{j+1}-i)$
by (12.1.2) and (b) of (12.4.12).} \\
&\text{$\underline{\text{\rm Equation(2) of (12.4.13).}}$} \\
&\text{$a\beta_{j+1,i,1}+b\beta_{j+1,i,2}+bn_1\beta_{j+1,i,3}+\rho_{2,m}\beta_{j+1,i,4}+
+\cdots+\rho_{j-1,m}\beta_{j+1,i,j+1}$}  \\
=&\text{$\Omega_2(\beta_{j+1,i,1},\beta_{j+1,i,2})+bn_1\beta_{j+1,i,3}+bn_1d_2\beta_{j+1,i,4}
+\cdots+bn_1d_{j-1}\beta_{j+1,i,j+1}$} \\
=&\text{$\Omega_2(\beta_{j+1,i,1},\beta_{j+1,i,2})+bn_1\beta_{j+1,i,3}+bn_1n_2\beta_{j+1,i,4}
+\cdots+bn_1n_2n_3\cdots n_{j-1}\beta_{j+1,i,j+1}$} \\
=&\text{$\Omega^{\sharp}_{j+1}(\beta_{j+1,i,k})^{j+1}_{k=1}\ge
\rho_{j+1,m}-i\rho_{j,m}=bn_1d_j(n_{j+1}-i)$} \text{ by (12.3.2),
(12.3.3) and (b) of (12.4.12).}
\endalign$$

Thus, the proof of (c) is done.  Now, it suffices to show that the
remaining part (d1) of (d) in this sublemma is true, which just
follows from Corollary $3.8$. So, the proof of the sublemma is
finished. $\square$
\enddemo
\ms

\demo{\bf Proof of Sublemma 12.5} First of all, let $\{Y_k:
k=1,2,\dots,r-1\}$ with $Y_k\subset N_0$, $\{h_k: k=1,2,\dots,r-1\}$
with $h_k=(g_{k+1}\circ\tau_m)_{proper}$ in $\BC\{v,1+\bar{u}\}$ and
$\{\Xi_k: N^k_0\to N_0: k=1,2,\dots,r-1\}$ be three different
sequences where each $\Xi_k$ is an integer-valued function,
satisfying the given three conditions, denoted by \text{\bf The 1-th
${\text{\bf{Cond}}}^{\text{{\bf(1)}}}$}{\bf, \dots,} \text{\bf The
3-th ${\text{\bf{Cond}}}^{\text{{\bf(1)}}}$}, in the conclusion of
this sublemma. After the proof of Sublemma $12.4$ was done, it is
easy to observe without any need of the proof that this three
sequences with three conditions are well-constructed.

For the proof of this sublemma, it suffices to show that this three
sequences satisfy the remaining two conditions \text{\bf The
4{$\alpha$}-th ${\text{\bf{Cond}}}^{\text{{\bf(1)}}}$} and \text{\bf
The 5{$\alpha$}-th ${\text{\bf{Cond}}}^{\text{{\bf(1)}}}$} in
$\underline {\text{\rm Conclusions}}$ of this sublemma, because
there is nothing to prove that the truth of \text{\bf The
4{$\alpha$}-th ${\text{\bf{Cond}}}^{\text{{\bf(1)}}}$} implies that
of The \text{\bf 4-th ${\text{\bf{Cond}}}^{\text{{\bf(1)}}}$}. So,
for the proof, firstly we will prove by {\bf[I]} that \text{\bf The
4{$\alpha$}-th ${\text{\bf{Cond}}}^{\text{{\bf(1)}}}$} is true, and
secondly,  by {\bf[II]} that \text{\bf The 5-th
${\text{\bf{Cond}}}^{\text{{\bf(1)}}}$} is true. \ms

{\bf [I]} For the proof of the truth of \text{\bf The 4{$\alpha$}-th
${\text{\bf{Cond}}}^{\text{{\bf(1)}}}$}, it remains to prove the
second inequality in {\rm(12.5.4${\alpha}$)} by using the following
three steps: Let $\ell$ and $q$ be an arbitrary positive integer
such that $r-1\ge {\ell}\ge q\ge 2$.
$$\align
\text{$\underline{\text{\rm Step(i)}}$} \quad  & \Xi_q(\gamma_{\ell,
i,k})^q_{k=1}=
\Delta_{q+1}(\beta_{\ell+1,i,k})^{q+1}_{k=1}+n^2_qn^2_{q-1}\cdots
n^2_2n_1\beta_{1,0,1}\\
& \qquad \qquad \quad  \times
\{\beta_{\ell+1,i,q+2}+n_{q+1}\beta_{\ell+1,i,q+3}
+n_{q+1}n_{q+2}\beta_{\ell+1,i,q+4} \\
&\qquad \qquad \quad +\cdots +n_{q+1}n_{q+2}\cdots
n_{\ell-1}\beta_{\ell+1,i,\ell+1}-n_{q+1}n_{q+2}\cdots n_{\ell}
(n_{\ell+1}-i)\}.
\qquad \qquad \\
\text{$\underline{\text{\rm Step(ii)}}$} \quad  & \text{In
particular, if $\ell=q$
then} \\
& \Xi_q(\gamma_{q,i,k})^q_{k=1}
=\Delta_{q+1}(\beta_{q+1,i,k})^{q+1}_{k=1}
 -(n_{q+1}-i)n^2_qn^2_{q-1}\cdots n^2_2n_1\beta_{1,0,1} \quad \text{from {\rm Step(i)}}.
\qquad \qquad \\
\text{$\underline{\text{\rm Step(iii)}}$} \quad &
\Xi_q(\gamma_{q,i,k})^q_{k=1}-(s_q-i)s_{q-1}\Xi_{q-1}(\gamma_{q-1,0,k})^{q-1}_{k=1}\\
 =& \Delta_{q+1}(\beta_{q+1,i,k})^{q+1}_{k=1}
 -(n_{q+1}-i)n_q\Delta_q(\beta_{q,0,k})^q_{k=1}>0 \quad \text{from {\rm
 Step(ii)}}. \quad
\endalign$$

We will prove {\rm Step(i)}, {\rm Step(ii)} and {\rm Step(iii)} in
order, by induction on the integer $q\ge 2$.

So, it is enough to consider two cases, respectively:

Case(I) $q=2$, and Case(II) $q\ge 2$. \ms

\noindent{\bf Case(I):} \ Let $q=2$. Note by \text{\bf The 3-th
${\text{\bf{Cond}}}^{\text{{\bf(1)}}}$} that
$\Xi_2(t_1,t_2)=t_2\Xi_1(\gamma_{1,0,1})+s_1\Xi_1(t_1)
=t_2\gamma_{1,0,1}+s_1t_1$ for each $(t_1,t_2)\in N^2_0$.
$$\align
\text{$\underline{\text{\rm Step(i)}}$} \qquad \qquad
& \Xi_2(\gamma_{{\ell},i,1},\gamma_{{\ell},i,2}) \\
 =& s_1\gamma_{{\ell},i,1}+\gamma_{1,0,1}\gamma_{{\ell},i,2} \\
=& n_2
\{\Delta^{\sharp}_2(\beta_{\ell+1,i1},\beta_{\ell+1,i2})-n_1\beta_{1,0,1}n_2\cdots
n_{\ell}(n_{\ell+1}-i)\} \\
&
+\{\Delta^{\sharp}_2(\beta_{2,0,1},\beta_{2,0,2})-n_1\beta_{1,0,1}n_2\}
\beta_{\ell+1,i,3}
\qquad \qquad \qquad \quad \text{by (12.5.1)}\qquad \\
=& n_2
\{\Delta_2(\beta_{\ell+1,i,1},\beta_{\ell+1,i,2})+n_1\beta_{1,0,1}\beta_{\ell+1,i,3}
+n_1\beta_{1,0,1}n_2\beta_{\ell+1,i,4}+\cdots  \\
&  +n_1\beta_{1,0,1}n_2\cdots n_{\ell-1}\beta_{\ell+1,i,\ell+1}
-n_1\beta_{1,0,1}n_2\cdots n_{\ell}(n_{\ell+1}-i) \} \\
&+ \{\Delta_2(\beta_{2,0,1},\beta_{2,0,2})-n_1\beta_{1,0,1}n_2 \}
\beta_{\ell+1,i,3} \qquad \qquad \qquad \quad \text{by (12.1.1)}\\
=&\Delta_3(\beta_{\ell+1,i,1},\beta_{\ell+1,i,2},\beta_{\ell+1,i,3})
+n^2_2n_1\beta_{1,0,1}\{\beta_{\ell+1,i,4}+n_3\beta_{\ell+1,i,5}
+n_3n_4\beta_{\ell+1,i,6} \qquad \qquad \\
& +\cdots +n_3n_4\cdots
n_{\ell-1}\beta_{\ell+1,i,\ell+1}-n_3n_4\cdots
n_{\ell}(n_{\ell+1}-i)\},
\endalign$$
by the definition of
$\Delta_3(\beta_{\ell+1,i,1},\beta_{\ell+1,i,2},\beta_{\ell+1,i,3})$
only, which implies the proof of {\rm Step(i)}. \ms

\noindent $\underline{\text{\rm Step(ii)}}$ \quad In particular, if
$\ell=2$ then  an equation of {\rm Step(i)} gives
$$
\Xi_2(\gamma_{2,i,1},\gamma_{2,i,2})=\Delta_3(\beta_{3,i,1},\beta_{3,i,2},\beta_{3,i,3})
-n^2_2n_1\beta_{1,0,1}(n_3-i).
$$

Thus, the proof of {\rm Step(ii)} is done. \ms

\noindent $\underline{\text{\rm Step(iii)}}$ \quad  To prove that
$\Xi_2(\gamma_{2,i,1},\gamma_{2,i,2})-(s_2-i)s_1\Xi_1(\gamma_{1,0,1})
=\Delta_3(\beta_{3,i,1},\beta_{3,i,2},\beta_{3,i,3})
-(n_3-i)n_2\Delta_2(\beta_{2,0,1},\beta_{2,0,2})>0$, first note by
(12.5.1.1) that
$$\align
&(s_2-i)s_1\Xi_1(\gamma_{1,0,1})=(s_2-i)s_1\gamma_{1,0,1}
=(n_3-i)n_2\{\Delta_2(\beta_{2,0,1},\beta_{2,0,2})
-n_1\beta_{1,0,1}n_2\}. \\
\text{Then,} \qquad &\Xi_2(\gamma_{2,i,1},\gamma_{2,i,2})-(s_2-i)s_1\gamma_{1,0,1} \\
=&\Delta_3(\beta_{3,i,1},\beta_{3,i,2},\beta_{3,i,3})-(n_3-i)n^2_2n_1\beta_{1,0,1}
-(n_3-i)n_2\{\Delta_2(\beta_{2,0,1},\beta_{2,0,2})-n_1\beta_{1,0,1}n_2\}  \qquad \qquad \\
=&\Delta_3(\beta_{3,i,1},\beta_{3,i,2},\beta_{3,i,3})
-(n_3-i)n_2\Delta_2(\beta_{2,0,1},\beta_{2,0,2})>0,
\endalign$$
by {\rm Step(ii)} and by \text{\bf The 4-th
${\text{\bf{Cond}}}^{\text{{\bf(0)}}}$} in the assumption of Theorem
$12.0$, which implies the proof of {\rm Step(iii)}.

Thus, if $q=2$, then we proved that the second inequality in
{\rm(12.5.4$\alpha$)} holds. \ms

\noindent{\bf Case(II):} \ Let $q\ge 2$. By the induction proof,
suppose we have shown that all the equalities of {\rm Step(i)}, {\rm
Step(ii)} and {\rm Step(iii)} are true on the integer $q\le r-1$
with $r-1\ge \ell\ge q$.

Then, it is enough to prove {\rm Step(i)}, {\rm Step(ii)} and {\rm
Step(iii)} in order, on the integer $(q+1)\le \ell$ as follows:
$$\align
\text{$\underline{\text{\rm Step(i)}}$} \qquad \qquad  &
\Xi_{q+1}(\gamma_{\ell,i,k})^{q+1}_{k=1} \\
=&\gamma_{\ell,i,q+1}\Xi_q(\gamma_{q,0,k})^q_{k=1}+s_q\Xi_q(\gamma_{\ell,i,k})^q_{k=1}
\quad \text{by definition of \ $\Xi_{q+1}$} \qquad \qquad \\
=&\beta_{\ell+1,i,q+2}
\{\Delta_{q+1}(\beta_{q+1,0,k})^{q+1}_{k=1}-n_{q+1}n^2_qn^2_{q-1}\cdots
n^2_2n_1\beta_{1,0,1} \} \\
& + n_{q+1} \{ \Delta_{q+1}(\beta_{\ell+1,i,k})^{q+1}_{k=1}
+n^2_qn^2_{q-1}\cdots n^2_2n_1\beta_{1,0,1}[\beta_{\ell+1,i,q+2} \\
&+n_{q+1}\beta_{\ell+1,i,q+3}+n_{q+1}n_{q+2}\beta_{\ell+1,i,q+4}+\cdots
\\ &+n_{q+1}n_{q+2}\cdots
n_{\ell-1}\beta_{\ell+1,i,\ell+1}-n_{q+1}n_{q+2}\cdots
n_{\ell}(n_{\ell+1}-i)]\} \\
& \quad   \text{by the induction assumption on the integer $q$}  \qquad \qquad   \\
=&\Delta_{q+2}(\beta_{\ell+1,i,k})^{q+2}_{k=1}+n^2_{q+1}n^2_q
n^2_{q-1}\cdots
n^2_2n_1\beta_{1,0,1} \\
& \times\{\beta_{\ell+1,i,q+3}+n_{q+2}\beta_{\ell+1,i,q+4}
 +n_{q+2}n_{q+3}\beta_{\ell+1,i,q+5}+\cdots  \\
& \quad +n_{q+2}n_{q+3}\cdots
n_{\ell-1}\beta_{\ell+1,i,\ell+1}-n_{q+2}n_{q+3}\cdots n_{\ell}
(n_{\ell+1}-i) \},
\endalign$$
by  the definition of $\Delta_{q+2}(\beta_{\ell+1,i,k})^{q+2}_{k=1}$
only, which implies the proof of Step(i). \ms

$\underline{\text{Step(ii)}}$ \quad In particular, if $\ell=q+1$,
then $\ell+1=q+2<q+3$ and so
$$
\Xi_{q+1}(\gamma_{q+1,i,k})^{q+1}_{k=1}=\Delta_{q+2}(\beta_{q+2,i,k})^{q+2}_{k=1}
-(n_{q+2}-i)n^2_{q+1}n^2_q\cdots n^2_2n_1\beta_{1,0,1},
$$
by Step(i) on the integer $q+1$, which implies the proof of Step(ii)
on the integer $q+1$. \ms

$\underline{\text{Step(iii)}}$ \quad To prove that the equality in
(12.5.4$\alpha$) is true, we have
$$\align \Xi_{q+1} &(\gamma_{q+1,i,k})^{q+1}_{k=1}
-(s_{q+1}-i)s_q\Xi_q(\gamma_{q,0,k})^q_{k=1} \\
&=\{\Delta_{q+2}(\beta_{q+2,i,k})^{q+2}_{k=1}-(n_{q+2}-i)n^2_{q+1}n^2_q\cdots
n^2_2n_1\beta_{1,0,1}\} \\
&\quad
-(n_{q+2}-i)n_{q+1}\{\Delta_{q+1}(\beta_{q+1,0,k})^{q+1}_{k=1}-n_{q+1}n^2_qn^2_{q-1}\cdots
n^2_2n_1\beta_{1,0,1}\} \\
&=\Delta_{q+2}(\beta_{q+2,i,k})^{q+2}_{k=1}-(n_{q+2}-i)n_{q+1}\Delta_{q+1}
(\beta_{q+1,0,k})^{q+1}_{k=1}>0,
\endalign$$
by Step(ii) on the integer $q$ and $q+1$, and by \text{\bf The 4-th
${\text{\bf{Cond}}}^{\text{{\bf(0)}}}$} of Theorem $12.0$, which
implies the proof of Step(iii).

Thus, we proved that the second inequality in $(12.5.3)$ is true,
and so we can finish the proof of the truth of \text{\bf The
4{$\alpha$}-th ${\text{\bf{Cond}}}^{\text{{\bf(1)}}}$}. \ms

{\bf [II]} By \text{\rm The $(4{\alpha})$-th
${\text{\rm{Cond}}}^{\text{{\rm(1)}}}$} and \text{\rm The 5-th
${\text{\rm{Cond}}}^{\text{{\rm(0)}}}$}, the proof of the truth of
an equality in (12.5.5{$\alpha$}.1) or (12.5.6.1) of \text{\bf The
5-th ${\text{\bf{Cond}}}^{\text{{\bf(1)}}}$} can be easily done:
Note that $r\ge 2$.
$$\align
(12.5.6.1) \quad & \gcd(\gamma_{1,0,1},s_1)=
\gcd(\Delta_2(\beta_{2,0,1},\beta_{2,0,2})
-n_2n_1\beta_{1,0,1},n_2)=
\gcd(\Delta_2(\beta_{2,0,1},\beta_{2,0,2}),n_2)=1,\\
&\gcd(\Xi_q(\gamma_{q,0,k})^q_{k=1},s_q) \quad \text{ for $2\le q\le
{r-1}$}\\
 \qquad \quad
=&\gcd(\Xi_q(\gamma_{q,0,k})^q_{k=1}-s_qs_{q-1}
\Xi_{q-1}(\gamma_{q-1,0,k})^{q-1}_{k=1},s_q) \quad \text{ for
$2\le q\le {r-1}$}\\
=&\gcd(\Delta_{q+1}(\beta_{q+1,0,k})^{q+1}_{k=1}
-n_{q+1}n_q\Delta_q(\beta_{q,0,k})^q_{k=1},n_{q+1}) \quad \text{for
$0\le
i<s_{q}$.} \qquad \qquad \\
=&\gcd(\Delta_{q+1}(\beta_{q+1,0,k})^{q+1}_{k=1},n_{q+1})=1 \quad
\text{for $0\le
i<s_{q}$.} \qquad \qquad \\
\endalign$$ \ms

Moreover, by \text{\rm The $(4{\alpha})$-th
${\text{\rm{Cond}}}^{\text{{\rm(1)}}}$} and \text{\rm The 4-th
${\text{\rm{Cond}}}^{\text{{\rm(0)}}}$}, the proof of the truth of
an equality in (12.5.{$5\alpha$}.2) or (12.5.6.2) of \text{\bf The
5-th ${\text{\bf{Cond}}}^{\text{{\bf(1)}}}$} is easily done: For
$q=2,3,\dots,r-1$,
$$\align
(12.5.6.2)\qquad \qquad &\quad \f{\Xi_{q}(\g_{q,i,k})^{q}_{k=1}
-(s_{q}-i)s_{q-1}\Xi_{q-1}(\g_{q-1,0,k})^{q-1}_{k=1}}{s_{q}-i}
\\    &> \f{\Xi_{q}(\g_{q,0,k})^{q}_{k=1}
-s_{q}s_{q-1}\Xi_{q-1}(\g_{q-1,0,k})^{q-1}_{k=1}}{s_{q}}>0 \quad
\text {for $0< i<s_q$,} \qquad \qquad\qquad \qquad\\
\iff \qquad &\\
&\quad \f{\Delta_{q+1}(\beta_{q+1,i,k})^{q+1}_{k=1}
-(n_{q+1}-i)n_q\Delta_q(\beta_{q,0,k})^q_{k=1}}{n_{q+1}-i}
\\   & > \f{\Delta_{q+1}(\beta_{q+1,0,k})^{q+1}_{k=1}
-n_{q+1}n_q\Delta_q(\beta_{q,0,k})^q_{k=1}}{n_{q+1}}>0 \quad \text
{for $0< i<n_{q+1}$.}
\endalign$$

In particular, if $q=1$, then note that an inequality in
$(12.5.6.2)$ of \text{\bf The 5-th
${\text{\bf{Cond}}}^{\text{{\bf(1)}}}$} holds if and only if the
following inequality is true: Note that $s_1=n_2$.
$$\align
(12.5.6.2^*)\quad  &\f{\Xi_{1}(\g_{1,i,1})}{s_{1}-i}
 > \f{\Xi_{1}(\g_{1,0,1})}{s_{1}}>0 \quad
\text {for $0< i<s_1$,} \qquad \qquad\qquad \qquad\\
\iff \quad & \\
&\f{\Delta_{2}(\beta_{2,i,1},\beta_{2,i,2})-n_1\beta_{1,0,1}(n_2-i)}{n_2-i}>
\f{\Delta_{2}(\beta_{2,0,1},\beta_{2,0,2})-n_1\beta_{1,0,1}n_2}{n_2}
>0 \quad \text {for $0< i<n_2$,}  \\
\endalign$$
which is equivalent to an equality in {\rm (5b)} of \text{\bf The
5-th ${\text{\bf{Cond}}}^{\text{{\bf(0)}}}$} of the assumption of
this theorem. $\square$
\enddemo
\bs

{\bf \S 13.2. For the proof of Theorem 12.0} \ms

Now, we will prove Theorem $12.0$ by using five sublemmas. For
convenience of the proof of the theorem, we need another additional
sublemma as follows. \ms

\proclaim{Sublemma 12.7} $\underline{\text{\bf {Assumptions}}}$
\quad

{\rm(i)} Suppose that the same properties and notations as in
$\underline{\text{\rm {Assumptions of Theorem 12.0}}}$ hold.

{\rm(ii)} Let $r$ be arbitrary integer with $r\ge 1$. As we have
seen in an additional condition of $\underline{\text{\rm
{Conclusions of Theorem 12.0}}}$, we may assume that the above $g_r$
satisfies all the equalities in {\rm(5a)} of \text{\rm The 5-th
${\text{\rm{Cond}}}^{\text{{\bf(0)}}}$} of {\rm [B]} of {\rm
Definition 12.0.0}, without mentioning any inequality in {\rm(5b)}.

{\rm(iii)} In addition, assume that $g_s$ is irreducible in
$\BC\{y,z\}$ for some integer $s$ where $1\le s\le r$. \ms

$\underline{\text{\bf {Conclusions}}}$ \quad  Then, we have the
following:

{\rm(a)} $g_1$ is irreducible in $\BC\{y,z\}$.

{\rm(b)} For all $j=1,2,\dots,r$, $g_j=g_j(y,z)$ can be written in
the form
$$\align
(12.7.1) \qquad \qquad & g_j=(z^{n_1}+\ve_{1,0}
y^{\beta_{1,0,1}})^{d_j}+\sum_{\alpha,\beta\ge 0}
c^{(j)}_{\alpha,\beta}y^{\alpha}z^{\beta} \quad \text{with
$\ve_{1,0}=1$
and} \qquad \qquad \qquad \\
& \quad \text{with
$n_1\alpha+\beta_{1,0,1}\beta>n_1\beta_{1,0,1}d_j$},
\endalign$$
where if $j\ge 2$ then $d_j=\prod^j_{k=2}n_k$ with $d_1=1$, and a
unit $\ve_{1,0}=\ve_{1,0}(y,z)$ may be analytically assumed to be
one in $\BC\{y,z\}$, and the $c^{(j)}_{\alpha,\beta}$ are nonzero
complex numbers for some nonnegative integers $\alpha$ and $\beta$
such that $n_1\alpha+\beta_{1,0,1}\beta>n_1\beta_{1,0,1}d_j$.

Note that $g_{s+1}$ may not be irreducible in $\BC\{y,z\}$ for a
given $s< r$. $\square$
\endproclaim \ms

\demo{\bf Proof of Sublemma 12.7} We are going to prove (a) and (b), simultaneously.
Assuming that (a) is true, it was already proved by Sublemma 12.2 that (b) is true. 
For the proof of the remaining,  assuming that (a) is not true, it suffices to prove that 
(b) is not true.

Whether $g_1$ is irreducible in $\BC\{y,z\}$ or not, 
we may assume that the above $g_1$
satisfies all the equalities in {\rm(2a)} of \text{\rm The 2nd
${\text{\rm{Cond}}}^{\text{{\bf(0)}}}$} of {\rm [A]} of {\rm
Definition 12.0.0}, as follows:
$$\align
 g_1&=z^{n_1}+\ve_{1,0}y^{\beta_{1,0,1}}
+\sum^{n_1-1}_{i=1}\ve_{1,i}y^{\beta_{1,i,1}}z^i \quad \text{with} \quad
\ve_{1,0}=1, \tag 12.7.2 \\
&=g_0+\sum_1 \quad \text{with} \quad 
\sum_1=\sum^{n_1-1}_{i=1}\ve_{1,i}y^{\beta_{1,i,1}}z^i. \\ 
\endalign$$

Since $g_1$ is not irreducible in $\BC\{y,z\}$,  $g_1$ of (12.7.2) can be rewritten  
as follows:
$$\align
(12.7.3) \qquad &g_1=g_0+\sum_1={g_0}+\sum_{1,1}+\sum_{1,2} \quad \quad \text{with} \quad 
g_0=z^{n_1}+\ve_{1,0}y^{\beta_{1,0,1}}, \qquad \qquad \qquad \qquad \qquad \qquad \qquad\\
 \quad &\text{where} \quad \sum_{1,1}=\sum_{\alpha,\beta\ge 0}
c_{\alpha,\beta}y^{\alpha}z^{\beta}\not =0 \quad \text{and}
\quad \sum_{1,2}=\sum_{\gamma,\delta\ge 0} 
c_{\gamma,\delta}y^{\gamma}z^{\delta},  \\
&\text{such that} \quad  n_1\alpha+\beta_{1,0,1}\beta<n_1\beta_{1,0,1} \quad \text{for any nonzero monomial $c_{\alpha,\beta}y^{\alpha}z^{\beta}$ in $\sum_{1,1}$} \quad\\ 
& \text{and that} \quad n_1\gamma+\beta_{1,0,1}\delta>n_1\beta_{1,0,1} \quad \text{for any  nonzero monomial $c_{\gamma,\delta}y^{\gamma}z^{\delta}$ in $\sum_{1,2}$}. \\ 
\endalign$$

Let $m=\min\{n_1\alpha+\beta_{1,0,1}\beta: n_1\alpha+\beta_{1,0,1}\beta <n_1\beta_{1,0,1} \}$
for all nonzero monomials $c_{\alpha,\beta}y^{\alpha}z^{\beta}$ in $\sum_{1,1}$. 
Then it is clear  by (12.7.3)
that there is one and only one nonzero monomial $c_{\gamma,\delta}y^{\gamma}z^{\delta}$ 
in $\sum_{1,1}$ of $g_1$ because ${\gamma}<{\beta}_{1,0,1}$ and ${\delta}<n_1$, and 
$1=\gcd(n_1,{\beta}_{1,0,1})$, and $\sum_{1,1}$ must be nonempty.
$$\align
 & \gamma<{\beta}_{1,0,1} \quad \text{and} \quad {\delta}<n_1 
 \quad \text{with \quad $1=\gcd(n_1,{\beta}_{1,0,1})$}, \tag 12.7.4\\ 
& \text{$m=n_1\gamma+\beta_{1,0,1}\delta<n_1\beta_{1,0,1}$ \quad and \quad
$n_1\alpha+\beta_{1,0,1}\beta>m$},
\endalign$$
where if $r\ge 2$,  $d_r=\prod^r_{k=2}n_k$ and $d_1=1$, and a
unit $\ve_{1,0}=\ve_{1,0}(y,z)$ may be analytically assumed to be
one in $\BC\{y,z\}$, and the $c^{(r)}_{\alpha,\beta}$ are nonzero
complex numbers for some nonnegative integers $\alpha$ and $\beta$
such that $n_1\alpha+\beta_{1,0,1}\beta>n_1\beta_{1,0,1}d_r$. \ms

By all the equalities in {\rm(2a)} and {\rm(2b)} of \text{\rm The 2nd
${\text{\rm{Cond}}}^{\text{{\bf(0)}}}$} of {\rm [A]} of {\rm
Definition 12.0.0}, using all the equalities in {\rm(12.1.2)} of Sublemma 12.1,
$\{g_k:k=2,3,\dots,j\}$ in {\rm(2b)} of \text{\rm The 2nd
${\text{\rm{Cond}}}^{\text{{\bf(0)}}}$} of {\rm [A]} of {\rm
Definition 12.0.0}  can be reviewed as follows:

{\rm(a)} $g_2={g_1}^{n_2}+\sum_2
={g_1}^{n_2}+\sum_{2,1}+\sum_{2,2}$
\text{with
$n_1\alpha+\beta_{1,0,1}\beta>n_1\beta_{1,0,1}n_2$} \ms

{\rm(b)} 
$$\align
g_k&={g_{k-1}}^{n_k}+\sum_k
={g_{k-1}}^{n_k}+\sum_{k-1,1}+\sum_{k-1,2} \quad
\text{for $k=2,3,\dots,j$.} \\
&\sum_{k-1,1}=\ve_{j,0}y^{\beta_{j,0,1}}z^{\beta_{j,0,2}}
g^{\beta_{j,0,3}}_1g^{\beta_{j,0,4}}_2\cdots
g^{\beta_{j,0,j}}_{j-2}\\
&\sum_{k-1,2}=\text{$c_j\sum^{n_j-1}_{i=1}\ve_{j,i}y^{\beta_{j,i,1}}z^{\beta_{j,i,2}}
g^{\beta_{j,i,3}}_1\cdots g^{\beta_{j,i,j}}_{j-2}g^i_{j-1}$,}
\endalign$$

$$\align
(12.7.3) \qquad &g_1=g_0+\sum_1={g_0}+\sum_{1,1}+\sum_{1,2} \quad \quad \text{with} \quad 
g_0=z^{n_1}+\ve_{1,0}y^{\beta_{1,0,1}}, \qquad \qquad \qquad \qquad \qquad \qquad \qquad\\
 \quad &\text{where} \quad \sum_{1,1}=\sum_{\alpha,\beta\ge 0}
c_{\alpha,\beta}y^{\alpha}z^{\beta}\not =0 \quad \text{and}
\quad \sum_{1,2}=\sum_{\gamma,\delta\ge 0} 
c_{\gamma,\delta}y^{\gamma}z^{\delta},  \\
&\text{such that} \quad  n_1\alpha+\beta_{1,0,1}\beta<n_1\beta_{1,0,1} \quad \text{for any nonzero monomial $c_{\alpha,\beta}y^{\alpha}z^{\beta}$ in $\sum_{1,1}$} \quad\\ 
& \text{and that} \quad n_1\gamma+\beta_{1,0,1}\delta>n_1\beta_{1,0,1} \quad \text{for any  nonzero monomial $c_{\gamma,\delta}y^{\gamma}z^{\delta}$ in $\sum_{1,2}$}. \\ 
\endalign$$
\text{with
$n_1\alpha+\beta_{1,0,1}\beta>n_1\beta_{1,0,1}n_2\cdots n_k$} 
for $k=2,3,\dots,j$.

Therefore, For all $k=1,2,\dots,j$, $g_k=g_k(y,z)$ can be written in
the form
$$\align
(12.7.1) \qquad \qquad & g_k=(z^{n_1}+\ve_{1,0}
y^{\beta_{1,0,1}})^{d_k}+\sum_{\alpha,\beta\ge 0}
c^{(k)}_{\alpha,\beta}y^{\alpha}z^{\beta} \quad \text{with
$\ve_{1,0}=1$
and} \qquad \qquad \qquad \\
& \quad \text{with
$n_1\gamma+\delta_{1,0,1}\delta<n_1\beta_{1,0,1}d_k$},
\endalign$$
where if $j\ge 2$ then $d_j=\prod^j_{k=2}n_k$ with $d_1=1$, and a
unit $\ve_{1,0}=\ve_{1,0}(y,z)$ may be analytically assumed to be
one in $\BC\{y,z\}$, and the $c^{(k)}_{\gamma,\delta}$ a unique nonzero
complex numbers for some nonnegative integers $\gamma$ and $\delta$
such that $n_1\gamma+\beta_{1,0,1}\delta<n_1\beta_{1,0,1}d_k$.

Thus, it can be proved that $g_j$ is not irreducible in $\BC\{y,z\}$.
$\square$
\enddemo \ms

\noindent \text{\rm(2)(2a)} \quad
\text{$g_1=z^{n_1}+\ve_{1,0}y^{\beta_{1,0,1}}
+c_1\sum^{n_1-1}_{i=1}\ve_{1,i}y^{\beta_{1,i,1}}z^i$ with
$\ve_{1,0}=1$.}

\text{\rm(2b)} \quad \text{
$g_j=g^{n_j}_{j-1}+\ve_{j,0}y^{\beta_{j,0,1}}z^{\beta_{j,0,2}}
g^{\beta_{j,0,3}}_1g^{\beta_{j,0,4}}_2\cdots
g^{\beta_{j,0,j}}_{j-2}$}

\qquad \qquad \quad
\text{$+c_j\sum^{n_j-1}_{i=1}\ve_{j,i}y^{\beta_{j,i,1}}z^{\beta_{j,i,2}}
g^{\beta_{j,i,3}}_1\cdots g^{\beta_{j,i,j}}_{j-2}g^i_{j-1}$,} where
$j=2,3,\dots,r$. \ms

Note that each $\ve_{j,i}=\ve_{j,i}(y,z)$ is a unit in $\BC\{y,z\}$
for $1\le j\le r$, $0\le i<n_j$, if exists. $\square$

\bs \bs

\demo{\bf Proof of Theorem 12.0} For example, let $r=1$. By Theorem
$3.2$ and Corollary $3.3$, there is nothing to prove for $[A]$. For
the proof of $[B]$, if $g_1$ is irreducible in $\BC\{y,z\}$ with
$n_1\ge 2$, then $y^{\gamma}g_1\in \text{the type[1]}$ under the
standard resolution by either Sublemma $12.4$ or Theorem $3.6$,
where if $\beta_{101}=1$ then $\gamma=1$ and if $\beta_{101}>1$ then
$\gamma=0$. In particular, $z^{\delta}yg_1\in \text{the type[1]}$
under the standard resolution by Sublemma $12.4$ or Theorem $3.6$
whether $\delta=1$ or $\delta=0$. Thus, $[B]$ can be easily proved.
So, if $r=1$, the proofs of $[A]$ and $[B]$ are done.

For the main proof, let $r\ge 2$. If $g_j$ is irreducible in
$\BC\{y,z\}$ for any $j\ge 1$, then it was already proved by
Sublemma 12.7 that $g_1$ is irreducible in $\BC\{y,z\}$. So, to
prove the necessity and sufficiency of the condition in $[A]$ and
the condition in $[B]$ at the same time, we may start with assuming
that $g_1$ is irreducible in $\BC\{y,z\}$ because it was already
proved by Corollary 3.3 that $g_1$ is irreducible in $\BC\{y,z\}$ if
and only if $\f{\beta_{1,i,1}}{n_1-i}> \f{\beta_{1,0,1}}{n_1}>0$ for
$0< i<n_1$, that is, the statement for [A] with $r=1$. Since
$\f{\beta_{1,i,1}}{n_1-i}> \f{\beta_{1,0,1}}{n_1}>0$ for $0< i<n_1$,
we are going to follow the same notations and consequences as in
Sublemma $12.4$.

Therefore, as we have done in Sublemma $12.4$, recall that
$\tau_m:M^{(m)}\to\BC^2$ and $(g_j\circ\tau_m)_{total}$ with
$(g_j\circ\tau_m)_{proper}$ satisfy the following four facts, where
$\tau_m:M^{(m)}\to\BC^2$ is the composition of a finite number $m$
of successive blow-ups $\pi_i$ which is needed to get the standard
resolution of the singular point of $V(y^{\gamma}g_1)$ in the
conclusion of Sublemma $12.4$. \ms

$\underline{\text{\rm Fact[1]}}$ By {\rm(a)} and {\rm(b)} in the
conclusion of Sublemma $12.4$, let $(v,u)$ be the common one of the
local coordinates for the $m$-th blow-up $\pi_m:M^{(m)}\to
M^{(m-1)}$ at $(0,0)$ which is the quasisingular point of
$V^{(m-1)}(y^{\gamma}g_1)$, in order to study any of $V^{(m)}(g_j)$
for all $j=2,3,\dots,r$ in the sense of Lemma $2.14$. Being viewed
as an analytic mapping, $\tau_m:M^{(m)}\to\BC^2$ can be written in
the form
$$
\tau_m(v,u)=(y,z)=(v^{n_1}u^a,v^{\beta_{1,0,1}}u^b), \tag 12.0.1
$$
where \ (i) $a$ and $b$ are some nonnegative integers such that
 $a\beta_{1,0,1}-bn_1=1$,

\qquad (ii) $E_m=\{v=0\}$ is defined by the $m-th$ exceptional curve
of the first kind. \ms

$\underline{\text{\rm Fact[2]}}$ By {\rm(c)} of the conclusion in
Sublemma $12.4$, along $v=0$, $(g_j\circ\tau_m)_{total}$ with
$(g_j\circ\tau_m)_{proper}$ can be written as follows: Use the same
notations for $\Delta_2,\Delta^{\sharp}_2,\Omega^{\sharp}_2$ as in
Sublemma $12.1$, Sublemma $12.2$, and Sublemma $12.3$, and we may
start with assuming that $\varepsilon_{1,0}=1$ in
$V(y^{\gamma}g_0)=\{y^{\gamma}(z^{n_1}+\varepsilon_{1,0}y^{\beta_{1,0,1}})=0\}$,
in order to study $V^{(i)}(g_j)$ for all $i=1,2,\dots,m,$ and all
$j=1,2,\dots,r$. Whether $\beta_{1,0,1}\ge 2$ or $\beta_{1,0,1}=1$,
we may write that $(g_0\circ\tau_m)_{proper}=(1+\ve_{1,0}u)$ and
$(g_1\circ\tau_m)_{proper}=(1+\ve_{1,0}\bar{u})$ with $\ve_{1,0}=1$,
without complexity of the notation if necessary, noting that if
$\beta_{1,0,1}=1$ then $V(g_0)$ has no singularity at the origin.
$$\align
(12.0.2) \quad (g_1\circ\tau_m)_{total}
&=v^{n_1\beta_{1,0,1}}u^{bn_1}(g_1\circ\tau_m)_{proper} \quad \text{with} \\
 (g_1\circ\tau_m)_{proper} =&(1+\ve_{1,0}u)
 +c_1\sum^{n_1-1}_{i=1}\ve'_{1,i}v^{\Delta^{\sharp}_2(\beta_{1,i,1},i)-n_1\beta_{1,0,1}}
 u^{\Omega^{\sharp}_2(\beta_{1,i,1},i)-bn_1} \text{ with $\ve_{1,0}=1$} \\
 =&1+\ve_{1,0}\bar{u}=1+\bar{u} \qquad \qquad \text{for brevity of notation}, \\
 (g_j\circ\tau_m)_{total} =&v^{n_1\beta_{1,0,1}n_2\cdots
 n_j}u^{bn_1n_2\cdots n_j}
 (g_j\circ\tau_m)_{proper}\quad\text{with}\\
 (g_j\circ\tau_m)_{proper} =&(g_{j-1}\circ\tau_m)^{n_j}_{proper}
 +\ve'_{j,0} v^{\Delta^{\sharp}_j(\beta_{j,0,k})^j_{k=1}-n_1\beta_{1,0,1}n_2\cdots n_j}
 u^{\Omega^{\sharp}_j(\beta_{j,0,k})^j_{k=1}-bn_1n_2\cdots n_j} \\
 & \quad \times (g_1\circ\tau_m)_{proper}^{\beta_{j,0,3}}
 (g_2\circ\tau_m)^{\beta_{j,0,4}}_{proper}
 \cdots (g_{j-2}\circ\tau_m)^{\beta_{j,0,j}}_{proper} \\
 +& c_j\sum^{n_j-1}_{i=1}\ve'_{j,i} v^{\Delta^{\sharp}_j(\beta_{j,i,k})^j_{k=1}
 -n_1\beta_{1,0,1}n_2\cdots n_{j-1}(n_j-i)}
 u^{\Omega^{\sharp}_j(\beta_{j,i,k})^j_{k=1}-bn_1n_2\cdots n_{j-1}(n_j-i)} \\
 & \quad \times (g_1\circ\tau_m)_{proper}^{\beta_{j,i,3}}
 (g_2\circ\tau_m)^{\beta_{j,i,4}}_{proper}
 \cdots (g_{j-2}\circ\tau_m)^{\beta_{j,i,j}}_{proper}(g_{j-1}\circ\tau_m)^{i}_{proper},
\endalign
$$
for $2\le j\le r$ where each $\ve'_{j,i}=\ve_{j,i}\circ \tau_m(v,u)$
is a unit in $\BC\{v,1+\bar{u}\}$ for $2\le j\le r$ and $0\le
i<n_j$, noting by Sublemma $12.1$ and Sublemma $12.3$ that for $0\le
i<n_j$,
$$
(12.0.3) \quad \Delta^{\sharp}_j(\beta_{j,i,k})^j_{k=1}
>n_1\beta_{1,0,1}n_2\cdots n_{j-1}(n_j-i) \ \text{and} \
\Omega^{\sharp}_j(\beta_{j,i,k})^j_{k=1}\ge bn_1n_2\cdots
n_{j-1}(n_j-i). \qquad \qquad
$$

Here, we may write $\eta_{j-1,i}=\eta_{j-1,i}(v,1+\bar{u})$ is a
unit in $\BC\{v,1+\bar{u}\}$ for $2\le j\le r$ and $1\le i\le
n_j-1$, noting that $\eta_{j-1,i}=\ve'_{j,i}
u^{\Omega^{\sharp}_j(\beta_{j,i,k})^j_{k=1}-bn_1n_2\cdots
n_{j-1}(n_j-i)}$. Here, we may assume by a nonsingular change of
coordinates that $(g_1\circ\tau_m)_{proper}=(1+\bar{u})$ for brevity
of the notation in the standard resolution of the singular point
$(v,1+u)=(0,0)$ of $V((g_j\circ\tau_m)_{proper})$ for $2\le j\le r$.

$\underline{\text{\rm Fact[3]}}$ For any positive integer $r\ge 2$,
observe by $(c)$ of Sublemma $12.4$ that $g_r$ is irreducible in
$\BC\{y,z\}$ if and only if $g_1$ is irreducible in $\BC\{y,z\}$ and
$(g_r\circ\tau_m)_{proper}$ is irreducible in $\BC\{v,1+\bar{u}\}$,
noting that $g_1$ is irreducible in $\BC\{y,z\}$ if and only if
$\f{\beta_{1,i,1}}{n_1-i}> \f{\beta_{1,0,1}}{n_1}>0$ for $0< i<n_1$.
\ms

$\underline{\text{\rm Fact[4]}}$ Note that
$\gcd(n_1,\beta_{1,0,1})=1$ and let $j$ be an arbitrary positive
integer with $j\ge 1$. Using (d3) of Sublemma $12.4$ or Theorem
$3.6$, then we get the following:

(4a) If $\beta_{1,0,1}>1$, then $g_j\in$ the type $[1]$ under
 $\tau_m$.

(4b) If $\beta_{1,0,1}=1$, then $g_j\in$ the type $[0]$ under
 $\tau_m$.

(4c) If $\beta_{1,0,1}\ge 1$, then $z^{\delta}yg_j\in$ the type
 $[1]$ under $\tau_m$, whether $\delta$ is either $1$ or $0$. \ms

Now, for the induction proof, it is enough to consider two cases,
respectively:

Case(1)  $r=2$, and Case(2) $r\ge 2$. \ms

{\bf Case(1):} \quad Let $r=2$. As we have seen in (12.4.2), we may
write for brevity of the notation that
$(g_0\circ\tau_m)_{proper}=(1+\ve_{1,0}u)$ and
$(g_1\circ\tau_m)_{proper}=(1+\ve_{1,0}\bar{u})$ with $\ve_{1,0}=1$.
\ms

$\underline{\text{\rm The proof for [A]}}$ First to prove [A], by
(12.0.2) or Sublemma $12.4$, $(g_2\circ\tau_m)_{total}$ in
$\BC\{v,1+\bar{u}\}$ may be written in the form
$$\align
(12.0.4) \quad
(g_2\circ\tau_m)_{total}&=v^{n_1\beta_{1,0,1}n_2}u^{bn_1n_2}
 (g_2\circ\tau_m)_{proper}\quad \text{with} \\
 (g_2\circ\tau_m)_{proper}=&(g_1\circ\tau_m)^{n_2}_{proper}
 +\ve'_{2,0} v^{\Delta^{\sharp}_2(\beta_{2,0,1},\beta_{2,0,2})-n_1\beta_{1,0,1}n_2}
 u^{\Omega^{\sharp}_2(\beta_{2,0,1},\beta_{2,0,2})-bn_1n_2}\\
 {+c_2}\sum^{n_2-1}_{i=1}&\ve'_{2,i}
v^{\Delta^{\sharp}_2(\beta_{2,i,1},\beta_{2,i,2})-n_1\beta_{1,0,1}(n_2-i)}
 u^{\Omega^{\sharp}_2(\beta_{2,i,1},\beta_{2,i,2})-bn_1(n_2-i)}(g_1\circ\tau_m)^i_{proper},\\
\text{where} \quad  \Delta^{\sharp}_2(\beta_{2,i,1},\beta_{2,i,2})&
-n_1\beta_{1,0,1}(n_2-i)>0  \quad \text{and} \quad
\Omega^{\sharp}_2(\beta_{2,i,1},\beta_{2,i,2})-bn_1(n_2-i)>0 \\
\text{and} \quad &\text{$\ve'_{2,i}=\ve_{2,i}\circ \tau_m(v,u)$ is a
unit in  $\BC\{v,1+\bar{u}\}$
for $0\le i<n_i$.} \\
\endalign$$

By (d2) of \text{\rm The 4-th
${\text{\rm{Cond}}}^{\text{{\rm(0)}}}$} in the assumption of the
theorem,
$$\align
{\Delta_{2}(\beta_{2,i,k})^{2}_{k=1}}
> {(n_{2}-i)}n_1\beta_{1,0,1}\quad
\text {for $0\le i<n_1$. } \tag 12.0.4*\\
\endalign$$

Now, we are going to prove the necessity and sufficiency of the
condition in [A].

Since
$\Delta_2(\beta_{2,0,1},\beta_{2,0,2})=\Delta^{\sharp}_2(\beta_{2,0,1},\beta_{2,0,2})$
for notation, then it is clear by (12.0.4) and (12.0.4*) that
$(g_2\circ\tau_m)_{proper}$ is irreducible in $\BC\{v,\bar{u}+1\}$
under the standard resolution $\tau_m$ if and only if
$\f{\Delta_{2}(\beta_{2,i,k})^{2}_{k=1}}{n_{2}-i}
> \f{\Delta_{2}(\beta_{2,0,k})^{2}_{k=1}}{n_{2}}$ for $0< i<n_2$
because $\gcd(n_2,\Delta_{2}(\beta_{2,0,k})^{2}_{k=1})=1$. Note that
$v(g_2\circ\tau_m)_{proper}\in \text{the type} [1]$. \ms

Therefore, it can be easily proved by {\rm Fact[3]} or by $(c)$ of
Sublemma $12.4$ that
$$\align
(12.0.5) \quad &\text{$g_2$ is irreducible in $\BC\{y,z\}$ with
$yg_2\in \text{the type} [2]$ under the
standard resolution} \\
 \iff &\text{$g_1$ is irreducible in $\BC\{y,z\}$ and
$(g_2\circ\tau_m)_{proper}$ is irreducible in $\BC\{v,1+\bar{u}\}$}\\
\qquad &\text{with $v(g_2\circ\tau_m)_{proper}\in \text{the type}
[1]$ under the standard resolution} \\
 \iff &\text{$\f{\beta_{1,i,1}}{n_1-i}> \f{\beta_{1,0,1}}{n_1}>0$ for $0< i<n_1$ and
$\f{\Delta_{2}(\beta_{2,i,k})^{2}_{k=1}}{n_{2}-i}
> \f{\Delta_{2}(\beta_{2,0,k})^{2}_{k=1}}{n_{2}}$ for $0< i<n_2$}.
\endalign
$$
Thus, we proved $[A]$ for $r=2$. \ms

$\underline{\text{\rm The proof for [B]}}$ Next to prove the
condition in $[B]$, let $g_2$ be irreducible in $\BC\{y,z\}$ with
$yg_2\in \text{the type} [2]$ under the standard resolution. If
$r=2$, then by {\rm Fact[4]} or Sublemma $12.4$ we get the
following:

\noindent(12.0.6) \quad \text{(i) \ If $\beta_{1,0,1}>1$, then
$g_2\in$ the type $[1]$ under $\tau_m$.}

\qquad \quad  \text{(ii) \ If $\beta_{1,0,1}=1$, then $g_2\in$ the
type $[0]$ under $\tau_m$.}

\qquad \quad  \text{(iii) \ If $\beta_{1,0,1}\ge 1$, then
$z^{\delta}yg_2\in$ the type $[1]$ under $\tau_m$ whether $\delta$
is $1$ or zero.} \ms

Let $V_2=\{(y,z):g_2(y,z)=0\}$ and $W_2=\{(y,z):yg_2(y,z)=0\}$ be
analytic varieties at the origin in $\BC\{y,z\}$, respectively. At
$(v,1+\bar{u})=(0,0)$, $\tau^{-1}_m(V_2)$ and $\tau^{-1}_m(W_2)$ can
be written as follows:
$$\align
(12.0.7) \qquad \tau^{-1}_m(V_2)
&=\{(g_2\circ\tau_m)_{total}=v^{n_1\beta_{1,0,1}n_2}u^{bn_1n_2}
(g_2\circ\tau_m)_{proper}=0\},  \\
\tau^{-1}_m(W_2)
&=\{((yg_2)\circ\tau_m)_{total}=v^{n_1\beta_{1,0,1}n_2+n_1}
u^{bn_1n_2+a}(g_2\circ\tau_m)_{proper}=0\}, \qquad \qquad \\
\text{such that} &\quad \\
(g_2\circ\tau_m)_{proper}&=(1+\bar{u})^{n_2}+\ve'_{2,0}v^{\Delta^{\sharp}_2(\beta_{2,0,1},
\beta_{2,0,2})-n_1\beta_{1,0,1}n_2} \times
u^{\Omega^{\sharp}_2(\beta_{2,0,1},\beta_{2,0,2})-bn_1n_2} \\
&=(1+\bar{u})^{n_2}+\eta_1v^{\Delta^{\sharp}_2(\beta_{2,0,1},
\beta_{2,0,2})-n_1\beta_{1,0,1}n_2} \\
&=h_1 \quad \text{by following the notation in Sublemma $12.5$,}
\endalign$$
where $\ve'_{2,0}$ and
$\eta_1=\ve'_{2,0}u^{\Omega^{\sharp}_2(\beta_{2,0,1},\beta_{2,0,2})-bn_1n_2}$
are units in $\BC\{v,1+\bar{u}\}$.

Since $(g_2\circ\tau_m)_{proper}$ and
$(\bar{u}+1)^{n_2}+v^{\Delta^{\sharp}_2(\beta_{2,0,1},\beta_{2,0,2})-n_1\beta_{1,0,1}n_2}$
have the same analytic type of the singularity at
$(\bar{u}+1,v)=(0,0)$, we may start with assuming that they are the
same equations. Let $V(\phi)=\{(v,\bar{u}+1):\phi(v,\bar{u}+1)=0\}$
be an analytic variety at $(v,\bar{u}+1)=(0,0)$ defined by
$$\align
(12.0.8) \qquad \qquad
&\phi=\phi(v,\bar{u}+1)=v^{\gamma_1}(g_2\circ\tau_m)_{proper}=v^{\gamma_1}h_1
\qquad \qquad\qquad \qquad\qquad \qquad \qquad \qquad\\
& \text{such that}  \quad \left\{\eqalign{& \text{$\gamma_1=1$,
\quad  if \quad $\Delta^{\sharp}_2(\beta_{2,0,1},
\beta_{2,0,2})-n_1\beta_{1,0,1}n_2=1$}, \cr & \text{$\gamma_1=0$,
\quad if \quad $\Delta^{\sharp}_2(\beta_{2,0,1},
\beta_{2,0,2})-n_1\beta_{1,0,1}n_2\ge 2$.} \cr} \right. \\
\endalign$$

Let $Z_1=\{(v,\bar{u}+1):vh_1=v(g_2\circ\tau_m)_{proper}=0\}$ be
another analytic variety at $(v,\bar{u}+1)=(0,0)$. Since
$(g_2\circ\tau_m)_{proper}$ is irreducible in $\BC\{v,\bar{u}+1\}$,
then note that $Z_1=\tau^{-1}_m(V_2)=\tau^{-1}_m(W_2)$ have the same
two irreducible components at $(v,\bar{u}+1)=(0,0)$ as reduced
variety, and so they have the homeomorphic resolution, using the
composition of the same number of successive blow-ups at
$(v,\bar{u}+1)=(0,0)$.

Now, we apply consequences of Theorem $3.6$ to the proof of [B].
First of all, it is clear that the equation of
$V(\phi)=V(v^{\gamma_1}h_1)$ of (12.0.8) satisfies the same kind of
assumption as in Theorem $3.6$, which can be represented as follows:
\ms

$\underline{\text{$h_1$ of (12.0.7) satisfies the same assumption as
$g_1$ of (3.6.1) does in Theorem 3.6}}$ \quad Let

\noindent $V(h_1)=\{(v,u+1):h_1(v,\bar{u}+1)=0\}$,
$V(f)=\{(v,\bar{u}+1): f(y,z)=0\}$ and $V(\phi)=\{(v,\bar{u}+1):
v^{\gamma_1}h(v,\bar{u}+1)=0\}$ be analytic varieties at
$(v,\bar{u}+1)=(0,0)$ in $\BC^2$, each of which is written
respectively as follows: For convenience of notation, substitute
$g_0$ of (12.4.7) by $h_1$, for an application of Theorem $3.6$ or
Sublemma $12.4$.
$$\align
(12.0.9) \qquad \qquad h_1 &=(1+\bar{u})^{n_2}+\eta_1v^{k_2} \quad
\text{with
$k_2=\Delta^{\sharp}_2(\beta_{2,0,1},\beta_{2,0,2})-n_1\beta_{1,0,1}n_2$},
\qquad \qquad \qquad \\
f&={h_1}^d+\sum_{\alpha,\beta\ge
0}c^{(2)}_{\alpha,\beta}v^{\alpha}(1+\bar{u})^{\beta} \quad
\text{with} \quad
n_2\alpha+k_2\beta>n_2k_2d=n_2k_2, \\
F&=v^{\delta_1}(\bar{u}+1)^{\delta_2}f, \\
\phi&=v^{\gamma_1}h_1 \\
\endalign$$
satisfying the properties {\rm(i)}, {\rm(ii)}, {\rm(iii)}, {\rm(iv)}
and {\rm(v)}:

\roster

\item "(i)" $\gcd(n_2,k_2)=1$ and $d=1$.

\item "(ii)" If $k_2=1$, then $\gamma=1$, and if $k_2\ge 2$, then
$\gamma=0$.

\item "(iii)" $\eta_1$ is assumed to be one in
$\BC\{v,\bar{u}+1\}$, and the $c^{(2)}_{\alpha,\beta}$ are nonzero
complex numbers for some nonnegative integers $\alpha$ and $\beta$
such that $n_2\alpha+k_2\beta> n_2k_2$, if exist.

\item "(iv)"  Each $\delta_i$ is either a positive integer or $0$
for $i=1,2$.

\item "(v)" If $k_2=1$, assume additionally that $\delta_1>0$,
whether or not $\delta_2$ is zero.
\endroster \ms

So, we have the same kind of conclusion as we have seen in Theorem
3.6, up to change of notations. \ms

$\underline{\text{$h_1$ of (12.0.7) satisfies the same conclusion as
$g_1$ of (3.6.1) does in Theorem 3.6}}$ As we have seen in Theorem
$3.6$, let $\omega_s$ be the composition of $s$ iterations of
blow-ups which is needed to get the standard resolution for the
singular point at $(v,\bar{u}+1)=(0,0)$ of $V(\phi)$, by using the
same method as we have done with $\tau_m$, either in the beginning
of the proof of Theorem $3.6$ or in (12.4.8) of Sublemma $12.4$,
because $\phi=v^{\gamma_1}(g_2\circ\tau_m)_{proper}=v^{\gamma_1}h_1$
satisfies the same kind of properties as $y^{\gamma}g_1$ does.
Therefore, by the conclusion of Theorem $3.6$, there is nothing to
prove for the conclusion of Sublemma $12.4$.

As far as the above $\omega_s$ is concerned, we can always use one
of the local coordinates for each blow-up of $\omega_s$ in the sense
of Definition $2.11$, whenever it is necessary.

Whether
$k_2=\Delta_2(\beta_{2,0,1},\beta_{2,0,2})-n_1\beta_{1,0,1}n_2$ is
either greater than one or equal to one at $(v,\bar{u}+1)=(0,0)$,
then $\phi\in$ the type $[1]$ and $vh_1\in$ the type $[1]$ under
$\omega_s$ by Fact[4] or Theorem $3.6$. Also, the $m-th$ exceptional
curve of the first kind among $\tau^{-1}_m(0,0)$, denoted by
$E_m=\{v=0\}$, can be viewed as one of two irreducible components of
$V(vh_1)$ at $(v,\bar{u}+1)=(0,0)$. Thus, both $V_2$ and $W_2$ of
(12.0.7) can have the same standard resolution of the singular point
$(y,z)=(0,0)$, using the composition $\tau_m\circ\omega_s$ of a
finite number $(m+s)$ of successive blow-ups at $(y,z)=(0,0)$.

Since $\phi\in$ the type $[1]$ and $vh_1\in$ the type $[1]$ at
$(v,\bar{u}+1)=(0,0)$ under $\omega_s$ where
$Z_1=\{vh_1=0\}=\tau^{-1}_m(V_2)=\tau^{-1}_m(W_2)$, then by (12.0.6)
or Fact[4] again, we get the following:

\noindent(12.0.10) \quad \text{(i) \ If $\beta_{1,0,1}>1$, then
$g_2\in$ the type $[2]$ under $\tau_m\circ\omega_s$.}

\qquad \quad \text{ (ii) \ If $\beta_{1,0,1}=1$, then $g_2\in$ the
type $[1]$ under $\tau_m\circ\omega_s$.}

\qquad \quad \text{ (iii) \ If $\beta_{1,0,1}\ge 1$, then
$z^{\delta}yg_2\in$ the type $[2]$ under $\tau_m\circ\omega_s$,
whether $\delta$ is $1$ or $0$.}

Thus, we proved $[B]$ for $r=2$. \ms

{\bf Case(2):} Let $r\ge 2$. For the proof of the theorem, by the
induction assumption, suppose we have shown that $[A]$ and $[B]$ of
the theorem are true on the integer $r\ge 2$.

In order to prove the theorem on the integer $(r+1)$, it suffices to
consider the defining equation for $g_{r+1}$ with the additive
assumptions:
$$\align
(12.0.11) \qquad \qquad  &
g_{r+1}=g^{n_{r+1}}_{r}+\ve_{r+1,0}y^{\beta_{r+1,0,1}}z^{\beta_{r+1,0,2}}
g^{\beta_{r+1,0,3}}_1\cdots g^{\beta_{r+1,0,r+1}}_{r-1} \qquad \qquad \qquad \qquad  \\
& \qquad \quad
+c_{r+1}\sum^{n_{r+1}-1}_{i=1}\ve_{r+1,i}y^{\beta_{r+1,i,1}}
z^{\beta_{r+1,i,2}}g^{\beta_{r+1,i,3}}_1\cdots
g^{\beta_{r+1,i,r+1}}_{r-1}g^i_{r}, \\
& \Delta^{\sharp}_{j+1}(\beta_{j+1,i,k})^{j+1}_{k=1} >
n_1\beta_{1,0,1}n_2n_3\cdots
n_{r}(n_{r+1}-i) \quad \text{on $g_{r+1}$}, \\
\endalign$$
where

(i) $g_1,g_2,\dots,g_r$ satisfies the same properties and notations
in the assumptions and conclusions of the theorem by induction
assumption,

(ii) $X_{r+1}=\{n_{r+1}\}\cup \{\beta_{r+1,i,1}:0\le i<n_{r+1} \}
\cup \{\beta_{r+1,i,2}:0\le i<n_3 \}\cup \cdots \cup
\{\beta_{r+1,i,r+1,}:0\le i<n_{r+1} \}\subset N_0$ with $n_{r+1}\ge
2$, and $\ve_{r+1,i}=\ve_{r+1,i}(y,z)$ is a unit in $\BC\{y,z\}$ for
$0\le i<n_{r+1}$, and the $c_{j+1}$ are nonzero complex numbers,

(iii)
$\Delta_{r+1}(t_k)^{r+1}_{k=1}=t_{r+1}\Delta_r(\beta_{r,0,k})^r_{k=1}
+n_r\Delta_r(t_k)^r_{k=1}$ for each $(t_k)^{r+1}_{k=1}\in
N^{r+1}_0$,
$$
\quad \text{\rm (iv)} \qquad \qquad
\Delta_{r+1}(\beta_{r+1,i,k})^{r+1}_{k=1}
>(n_{r+1}-i)n_{r}\Delta_{r}(\beta_{r,0,k})^{r}_{k=1}
\quad \text {for $0\le i<n_{r+1}$.} \qquad \qquad\qquad \qquad $$

Now, it is needed to show that two statements $[A]$ and $[B]$ are
true for the integer $(r+1)$, respectively. In the beginning of the
proof, it was already known by Sublemma $12.7$ that if $g_j$ is
irreducible in $\BC\{y,z\}$ for any $j\ge 1$, then $g_1$ is
irreducible in $\BC\{y,z\}$. \ms

$\underline{\text{\rm The proof for [A]}}$ First, we are going to
prove the necessity and sufficiency of the condition in [A],
assuming that $\f{\beta_{1,i,1}}{n_1-i}> \f{\beta_{1,0,1}}{n_1}>0$
for $0< i<n_1$. Using the composition $\tau_m$ again of a finite
number $m$ of successive blow-ups which is needed to get the
standard resolution of the singular point of $V(y^{\gamma}g_1)$ in
the beginning of the proof, then $h_r=(g_{r+1}\circ\tau_m)_{proper}$
can be written in the form
$$\align
(12.0.12) \qquad \qquad h_r &= h^{s_r}_{r-1}+\eta_{r,0}
v^{\gamma_{r,0,1}}(1+\bar{u})^{\gamma_{r,0,2}}h^{\gamma_{r,0,3}}_1
\cdots
h^{\gamma_{r,0,r}}_{r-2},  \qquad \qquad\qquad \qquad\qquad \qquad\\
& \quad +c_{r+1}\sum^{s_r-1}_{i=1}{\eta_{r,i}} v^{\gamma_{r,i,1}}
(1+\bar{u})^{\gamma_{r,i,2}}h^{\gamma_{r,i,3}}_1 \cdots
h^{\gamma_{r,i,r}}_{r-2}h^{i}_{r-1}, \\
\endalign$$
 where

(i) each $h_j=(g_{j+1}\circ\tau_m)_{proper}$ is in
$\BC\{v,1+\bar{u}\}$ for $j=1,2,\dots,r$, which has been

\quad already represented by Sublemma $12.5$,

(ii) $(g_{r+1}\circ\tau_m)_{total}=v^{n_1\beta_{1,0,1}n_2\cdots
 n_r}u^{bn_1n_2\cdots n_r}(g_{r+1}\circ\tau_m)_{proper}$,

(iii) $\eta_{r,i}=\ve'_{r+1,i}
u^{\Omega^{\sharp}_{r+1}(\beta_{r+1,i,k})^{r+1}_{k=1}-bn_1n_2\cdots
n_{r+1}}$ is a unit in $\BC\{v,\bar{u}+1\}$, satisfying the

\qquad same notations and consequences as in Sublemma $12.5$,

(iv)  $s_{r} =n_{r+1}\ge 2,
 ~\gamma_{r,i,1}=\Delta^{\sharp}_{r+1}(\beta_{r+1,i,k})^{r+1}_{k=1}
 -n_1\beta_{1,0,1}n_2n_3\cdots n_r(n_{r+1}-i)>0,$

\qquad  $\gamma_{r,i,2}=\beta_{r+1,i,3},
\gamma_{r,i,3}=\beta_{r+1,i,4} ,
\dots,\gamma_{r,i,r}=\beta_{r+1,i,r+1}$. \ms

Because blow-ups process preserves irreducibility of plane curve
singularity, note by Fact[3] that $g_{r+1}$ is irreducible in
$\BC\{y,z\}$ if and only if $g_1$ is irreducible in $\BC\{y,z\}$ and
$h_r=(g_{r+1}\circ\tau_m)_{proper}$ is irreducible in
$\BC\{v,\bar{u}+1\}$. To prove the aim in [A], we may assume that
$\f{\beta_{1,i,1}}{n_1-i}> \f{\beta_{1,0,1}}{n_1}>0$ for $0< i<n_1$
because $g_1$ is irreducible in $\BC\{y,z\}$ if and only if
$\f{\beta_{1,i,1}}{n_1-i}> \f{\beta_{1,0,1}}{n_1}>0$ for $0< i<n_1$.
First of all, whether or not $h_r$ is irreducible in
$\BC\{v,\bar{u}+1\}$, we proved by Sublemma $12.5$ that $h_r$
satisfies the same kind of assumptions as $g_r$ does in this
theorem, up to change of notations.
 So, by induction assumption on the integer $r$ and by following
the same notations as in Sublemma $12.5$, then it is easy to get the
following: Note by \text{\rm The $(5{\alpha})$-th
${\text{\rm{Cond}}}^{\text{{\rm(1)}}}$} of Sublemma 12.5 that
$\gcd(s_{q},\Xi_{q}(\gamma_{q,0,k})^{q}_{k=1})=
\gcd(n_{q+1},\Delta_{q+1}(\beta_{q+1,0,k})^{q+1}_{k=1})=1$ for
$q=1,2,\dots,r$.
$$\split
(12.0.13) \quad & \text{$h_r$ is irreducible
in $\BC\{v,\bar{u}+1\}$} \\
& \text{with $v^{\gamma}h_r\in \text{the type} [r]$
under the standard resolution}\\
\iff \quad   &   \f{\g_{1,i,1}}{s_1-i}>
\f{\g_{1,0,1}}{s_1} \quad \text {for $0< i<s_1$,} \quad \text{and}   \\
&  \f{\Xi_{j-1}(\g_{j-1,i,k})^{j-1}_{k=1}}{s_{j-1}-i} >
\f{\Xi_{j-1}(\g_{j-1,0,k})^{j-1}_{k=1}}{s_{j-1}} \quad \text {for
each j=2,3,\dots,r+1 and for $0<i<s_{j-1}$.} \qquad \qquad
\endsplit$$

Therefore, we get by $(12.0.13)$ and by an equality in
(12.5.{$5\alpha$}.2) or (12.5.6.2) of \text{\rm The $(5{\alpha})$-th
${\text{\rm{Cond}}}^{\text{{\rm(1)}}}$} of Sublemma 12.5 that
$$\split
(12.0.14) \qquad \qquad &\text{$h_r=(g_{r+1}\circ\tau_m)_{proper}$
is irreducible in $\BC\{v,\bar{u}+1\}$} \\
& \text{with $v^{\gamma}h_r\in \text{the type} [r]$
under the standard resolution}\\
\iff \quad &{\f{\Delta_{j}(\beta_{j,i,k})^{j}_{k=1}}{n_{j}-i}}
> {\f{\Delta_{j}(\beta_{j,0,k})^{j}_{k=1}}{n_{j}}} \quad
\text {for each $j=2,3,\dots,r+1$ and for $0< i<n_j$.} \\
\endsplit$$

Note by Fact[3] or by $(c)$ of Sublemma $12.4$ that
$$\split
(12.0.15) \qquad &\text{$g_{r+1}$ is irreducible in $\BC\{y,z\}$} \\
& \text{with $y^{\gamma}g_{r+1}\in \text{the type} [r+1]$
under the standard resolution}\\
\iff \quad & \text{$\f{\beta_{1,i,1}}{n_1-i}>
\f{\beta_{1,0,1}}{n_1}>0$ for $0< i<n_1$ and
$h_r=(g_{r+1}\circ\tau_m)_{proper}$ is irreducible} \qquad \qquad
\\
& \text{in $\BC\{v,\bar{u}+1\}$ with $v^{\gamma}h_r\in \text{the
type} [r]$
under the standard resolution.} \\
\endsplit$$

Thus, we can prove by (12.0.14) and (12.0.15) and by the induction
assumption that
$$\align
(12.0.16) \quad &\text{$g_{r+1}$ is irreducible in $\BC\{y,z\}$} \\
& \text{with $y^{\gamma}g_{r+1}\in \text{the type} [r+1]$
under the standard resolution}\\
\iff &{\f{\Delta_{j}(\beta_{j,i,k})^{j}_{k=1}}{n_{j}-i}}
> {\f{\Delta_{j}(\beta_{j,0,k})^{j}_{k=1}}{n_{j}}} \quad
\text {for each $j=1,2,\dots,r+1$ and for $0< i<n_j$.} \\
\iff &\text{$g_1, g_2,\dots,g_r$ are irreducible in $\BC\{y,z\}$ and
${\f{\Delta_{r+1}(\beta_{r+1,i,k})^{r+1}_{k=1}}{n_{r+1}-i}}
> {\f{\Delta_{r+1}(\beta_{r+1,0,k})^{r+1}_{k=1}}{n_{r+1}}}$.}
\endalign$$
Therefore, we proved $[A]$. \ms

$\underline{\text{\rm The proof for [B]}}$ Next to prove $[B]$, let
$g_{r+1}$ be irreducible in $\BC\{y,z\}$. If $g_j$ is irreducible in
$\BC\{y,z\}$ for any $j\ge 1$, note that $\gcd(n_1,\beta_{1,0,1})=1$
by (d) of Sublemma $12.2$ and Theorem $3.7$. By Fact[4] or Sublemma
$12.4$ again, we have the following:

\noindent(12.0.17) \quad \rm{(i)} \quad If $\beta_{1,0,1}>1$, then
$g_{r+1}\in$ the type $[1]$ under $\tau_m$.

\qquad \quad \rm{(ii)} \quad If $\beta_{1,0,1}=1$, then $g_{r+1}\in$
the type $[0]$ under $\tau_m$.

\qquad \quad \rm{(iii)} \quad If $\beta_{1,0,1}\ge 1$, then
$z^{\delta}yg_{r+1}\in$ the type $[1]$ under $\tau_m$ whether
$\delta$ is $1$ or $0$. \ms

Let $V_{r+1}=\{(y,z):g_{r+1}(y,z)=0\}$ and
$W_{r+1}=\{(y,z):yg_{r+1}(y,z)=0\}$ be analytic varieties at the
origin in $\BC\{y,z\}$, respectively.  Let $\tau_m$ be again the
composition of a finite number $m$ of successive blow-ups which is
needed to get the standard resolution of the singular point of
$V(yg_1)$.

At $(v,1+\bar{u})=(0,0)$, $\tau^{-1}_m(V_{r+1})$ and
$\tau^{-1}_m(W_{r+1})$ can be written as follows:
$$\align
 \tau^{-1}_m(V_{r+1})
&=\{(g_{r+1}\circ\tau_m)_{total}=v^{n_1\beta_{1,0,1}n_2\cdots
n_{r+1}}u^{bn_1n_2\cdots n_r}(g_{r+1}\circ\tau_m)_{proper}=0\}
\quad \text{and} \tag 12.0.18 \\
\tau^{-1}_m(W_{r+1})
&=\{((yg_{r+1})\circ\tau_m)_{total}=v^{n_1\beta_{1,0,1}n_2\cdots
n_{r+1}+n_1}u^{bn_1n_2\cdots
n_r+a}(g_{r+1}\circ\tau_m)_{proper}=0\},
\endalign$$
noting by $(12.0.12)$ that $h_r=(g_{r+1}\circ\tau_m)_{proper}$ with
$$\align
(12.0.19) \qquad \qquad h_r &= h^{s_r}_{r-1}+\eta_{r,0}
v^{\gamma_{r,0,1}}(1+\bar{u})^{\gamma_{r,0,2}}h^{\gamma_{r,0,3}}_1
\cdots h^{\gamma_{r,0,r}}_{r-2} \qquad \qquad\qquad \qquad\qquad \qquad \\
& \quad +c_{r+1}\sum^{s_r-1}_{i=1}{\eta_{r,i}} v^{\gamma_{r,i,1}}
(1+\bar{u})^{\gamma_{r,i,2}}h^{\gamma_{r,i,3}}_1 \cdots
h^{\gamma_{r,i,r}}_{r-2}h^{i}_{r-1}, \\
\endalign$$
where $\eta_{r,i}=\ve'_{r+1,i}
u^{\Omega^{\sharp}_{r+1}(\beta_{r+1,i,k})^{r+1}_{k=1}-bn_1n_2\cdots
n_{r+1}}$ is a unit in $\BC\{v,\bar{u}+1\}$. \ms

Let $V(\phi_r)=\{(v,\bar{u}+1):\phi_r(v,\bar{u}+1)=0\}$ be an
analytic variety at $(v,\bar{u}+1)=(0,0)$ defined by
$$\align
(12.0.20) \qquad \qquad &\phi_r=\phi_r(v,\bar{u}+1)=v^{\gamma_1}h_r
\quad
\text{with \quad $h_r=(g_{r+1}\circ\tau_m)_{proper}$} \qquad \qquad \qquad \qquad \\
& \text{such that}  \quad \left\{\eqalign{& \text{$\gamma_1=1$ \quad
if \quad $\Delta^{\sharp}_2(\beta_{2,0,1},
\beta_{2,0,2})-n_1\beta_{1,0,1}n_2=1$}, \cr & \text{$\gamma_1=0$
\quad if \quad $\Delta^{\sharp}_2(\beta_{2,0,1},
\beta_{2,0,2})-n_1\beta_{1,0,1}n_2\ge 2$.} \cr} \right. \\
\endalign$$ \ms

Let $Z_r=\{(v,\bar{u}+1):vh_r(v,\bar{u}+1)=0\}$ be an analytic
variety at $(v,\bar{u}+1)=(0,0)$. Since
$(g_{r+1}\circ\tau_m)_{proper}$ is irreducible in
$\BC\{v,\bar{u}+1\}$, then note that
$Z_{r}=\tau^{-1}_m(V_{r+1})=\tau^{-1}_m(W_{r+1})$ have the same two
irreducible components at $(v,\bar{u}+1)=(0,0)$ as reduced variety,
and so they have the same standard resolution of the singular point
$(v,\bar{u}+1)=(0,0)$. Let $\omega_s$ be the composition of $s$
iterations of blow-ups which is needed to get the standard
resolution of the singular point $(v,\bar{u}+1)=(0,0)$ of
$V(\phi_r)$. Since we proved by Sublemma $12.5$ that $V(h_r)$
satisfies the same kind of properties and notations as $V(g_r)$ does
in the assumption of this theorem, then by the induction assumption
on the integer $r$, both $V(\phi_r)$ and $Z_{r}$ belong to the type
$[r]$ under $\omega_s$. That is, both $\tau^{-1}_m(V_{r+1})$ and
$\tau^{-1}_m(W_{r+1})$ belong to the type $[r]$ under $\omega_s$ at
$(v,\bar{u}+1)=(0,0)$, as reduced variety.

Now, the $m$-th exceptional curve of the first kind among
$\tau^{-1}_m(0,0)$, denoted by $E_m=\{v=0\}$, can be viewed as one
of two irreducible components of $Z_r$ at $(v,\bar{u}+1)=(0,0)$.
Then, both $V_{r+1}$ and $W_{r+1}$ can have the same standard
resolution $\tau_m\circ\omega_s$ of the singular point
$(y,z)=(0,0)$, using the composition $\tau_m\circ\omega_s$ of a
finite number $(m+s)$ of successive blow-ups at $(y,z)=(0,0)$.

Since $Z_r\in$ the type $[r]$ under $\omega_s$, then by (12.0.17) or
Sublemma $12.4$ again, it is trivial to get the following: Note that
$\delta$ is an integer.

\noindent(12.0.21) \quad {\rm(i)} If $\beta_{1,0,1}>1$,  $g_{r+1}\in$ the type $[r+1]$
under $\tau_m\circ\omega_s$.

\qquad \qquad {\rm(ii)} If $\beta_{1,0,1}=1$,  $g_{r+1}\in$ the type $[r]$
under $\tau_m\circ\omega_s$.

\qquad \qquad {\rm(iii)} If $\beta_{1,0,1}\ge 1$, $z^{\delta}yg_{r+1}\in$ the
type $[r]$ under $\tau_m\circ\omega_s$, whether $\delta$ is zero or
not.

Thus, the proof for $[B]$ is done. Therefore, we completed the proof
of theorem. $\square$
\enddemo \ms

\proclaim{Corollary 12.6}

$\underline{\text{\bf Assumptions}}$ \quad Suppose that the
assumptions of Theorem $12.0$ hold. Let $r$ be an arbitrary positive
integer with $r\ge 1$. \ms

$\underline{\text{\bf Conclusions}}$ \quad Then, we get the
following:

For any $r\ge 1$, $g_r=g_r(y,z)$ can be written in the form
$$\align
(12.6.1) \qquad & g_r=(z^{n_1}+\ve_1
y^{\beta_{1,0,1}})^{n_2n_3\cdots n_r}+\sum_{\alpha,\beta\ge 0}
c^{(r)}_{\alpha,\beta}y^{\alpha}z^{\beta} \quad \text{with $\ve_1=1$
and}
\qquad \qquad \qquad \\
& \quad \text{with
$n_1\alpha+\beta_{1,0,1}\beta>n_1\beta_{1,0,1}n_2n_3\cdots n_r$},
\endalign$$
where a unit $\ve_1=\ve_1(y,z)$ may be analytically assumed to be
one in $\BC\{y,z\}$, and the $c^{(r)}_{\alpha,\beta}$ are nonzero
complex numbers for some nonnegative integers $\alpha$ and $\beta$
such that $n_1\alpha+\beta_{1,0,1}\beta>n_1\beta_{1,0,1}n_2n_3\cdots
n_r$.  $\square$ \ms
\endproclaim

\demo{\bf Proof of Corollary} There is nothing to prove by Sublemma
$12.2$. $\square$
\enddemo \bs

\vfill \pagebreak

{\bf \S14. How to compute the divisor of the total transform of
irreducible plane curve singularities defined by the quasi-Puiseux
convergent power series of the recursive type in $\BC\{y,z\}$ under
the standard resolution} \ms

\proclaim{Theorem 14.0} $\underline{\text{\bf {Assumptions}}}$ Let
$r$ be an arbitrary positive integer.

\noindent{\bf [A]} Let $g_{r+1}\in\BC\{y,z\}$ be
$\underline{\text{\rm a generalized semi-quasi-Puiseux germ of the
recursive (r+1)-type}}$ which is defined by {\rm Definition 12.0.0}:
$$\align
\quad Let \quad &\{X_k:k=1,2,\dots,r+1\} \quad
\text{with $X_k\subset N_0$}, \\
&\{g_k:k=1,2,\dots,r+1\} \quad \text{with $g_k\in\BC\{y,z\}$,} \\
&\{\text{$\Delta_k:N^k_0\to N_0$ is an integer-valued
function for $k=1,2,\dots,r+1$}\} \\
&\text{be three sequences, satisfying four conditions for each k}:
\qquad
\endalign$$
Such conditions are denoted by \text{\bf The 1-th
${\text{\bf{Cond}}}^{\text{{\bf(0)}}}$}, \dots, \text{\bf The 4-th
${\text{\bf{Cond}}}^{\text{{\bf(0)}}}$} for brevity. \ms

\noindent$\underline{\text{\rm The 1st
${\text{\rm{Cond}}}^{\text{{\rm(0)}}}$}}$ The family
$\{X_{\ell}:\ell=1,2,\dots,r+1\}$ with $X_{\ell}\subset N_0$ is as
follows: \ms

\noindent{\rm(1)(1a)} $X_1=\{n_1\}\cup \{\beta_{1,i,1}:0\le i<n_1
\}$ with $n_1\ge 2$ and $\beta_{1,0,1}\ge 1$ where $X_{1}\subset N$.

{\rm(1b)} $X_j=\{n_j\}\cup \{\beta_{j,i,1}:0\le i<n_j \} \cup
\{\beta_{j,i,2}:0\le i<n_j \}\cup \cdots \cup \{\beta_{j,i,j}:0\le
i<n_j \}$ with $n_j\ge 2$ where $j=2,3,\dots,r+1$.

For each $j=2,3,\dots,r+1$, assume that at least one of
$\beta_{j,0,1},\beta_{j,0,2},\dots,\beta_{j,0,j}$ is nonzero. \ms

\noindent $\underline{\text{\rm The 2nd
${\text{\rm{Cond}}}^{\text{{\rm(0)}}}$}}$ For each
$j=1,2,\dots,r+1$, let $g_j=g_j(y,z,c_j)$ be in $\BC\{y,z\}$ where
all the $c_j$ are complex numbers, each of which is defined by the
following way: \ms

\noindent \text{\rm(2)(2a)} \quad
\text{$g_1=z^{n_1}+\ve_{1,0}y^{\beta_{1,0,1}}
+c_1\sum^{n_1-1}_{i=1}\ve_{1,i}y^{\beta_{1,i,1}}z^i$ with
$\ve_{1,0}=1$.}

\text{\rm(2b)} \quad \text{
$g_j=g^{n_j}_{j-1}+\ve_{j,0}y^{\beta_{j,0,1}}z^{\beta_{j,0,2}}g^{\beta_{j,0,3}}_1\cdots
g^{\beta_{j,0,j}}_{j-2}
+c_j\sum^{n_j-1}_{i=1}\ve_{j,i}y^{\beta_{j,i,1}}z^{\beta_{j,i,2}}g^{\beta_{j,i,3}}_1\cdots
g^{\beta_{j,i,j}}_{j-2}g^i_{j-1}$,} where $j=2,3,\dots,r+1$. \ms

Note that each $\ve_{j,i}=\ve_{j,i}(y,z)$ is a unit in $\BC\{y,z\}$
for $1\le j\le r+1$ and $0\le i<n_j$, if exists. As far as analytic
equivalence of isolated plane curve singularities defined by all
$g_j$, $1\le j\le r+1$, is concerned, then we may assume that
$\ve_{1,0}$ is equal to one by a suitable nonsingular change of
coordinates at the origin in $\BC^2$. \ms

\noindent $\underline{\text{\rm The 3rd
${\text{\rm{Cond}}}^{\text{{\rm(0)}}}$}}$ \quad Let $\{\Delta_k:
N^k_0\to N_0: k=1,2,\dots,r+1\}$ be a sequence such that each
$\Delta_k$ is an integer-valued function defined by the following
way: \roster \item"(3) (3a)" $\Delta_1(t)=t$ for each $t\in N_0$.

\item"$\quad$ (3b)"
$\Delta_j(t_k)^j_{k=1}=t_j\Delta_{j-1}(\beta_{j-1,0,k})^{j-1}_{k=1}
+n_{j-1}\Delta_{j-1}(t_k)^{j-1}_{k=1}$ for each $(t_k)^j_{k=1}\in
N^j_0$.
\endroster \ms

\noindent $\underline{\text{\rm The 4-th
${\text{\rm{Cond}}}^{\text{{\rm(0)}}}$}}$ \quad Then, the following
inequalities hold: Note that $r+1\ge 2$.
$$\align
\text{\rm (4)(4a)}& \quad \Delta_1(\beta_{1,i,1})={\beta_{1,i,1}}>0
\quad \text {for $0\le i<n_1$.}   \qquad \qquad \\
\text{\rm(4b)}& \quad {\Delta_{j}(\beta_{j,i,k})^{j}_{k=1}}
> {(n_{j}-i)}n_{j-1}\Delta_{j-1}(\beta_{j-1,0,k})^{j-1}_{k=1} \quad
\text {for $0\le i<n_j$ where $j=2,\dots,r+1$}. \qquad \qquad \qquad
\qquad
\endalign$$ \ms

\noindent{\bf [B]} Let $g_r\in \BC\{y,z\}$ be $\underline{\text{\rm
a generalized semi-quasi-Puiseux convergent power series of the
recursive}}$ $\underline{\text{\rm r-type}}$ as in {\rm [A]}. There
are an additional condition, denoted by \text{\rm The 5-th
${\text{\rm{Cond}}}^{\text{{\bf(0)}}}$}. \ms

\noindent $\underline{\text{\rm The 5-th
${\text{\rm{Cond}}}^{\text{{\rm(0)}}}$}}$ \quad For each
$q=1,2,\dots,r$, the following inequalities hold:
$$\align
\noindent \text{\rm (5)}&\text{\rm(5a)} \qquad \qquad\qquad\qquad
\text{ $\gcd(n_q,\Delta_q(\beta_{q,0,k})^q_{k=1})=1$ \quad for $1\le
q\le
r$.}  \\
&\text{\rm (5b)} \qquad \qquad\qquad\qquad
\f{\Delta_{q}(\beta_{q,i,k})^{q}_{k=1}}{n_{q}-i}
> \f{\Delta_{q}(\beta_{q,0,k})^{q}_{k=1}}{n_{q}} \quad
\text {for $0< i<n_q$.} \qquad \qquad\qquad\qquad
\endalign$$

In addition, for each $j=1,2,\dots,r+1$, let $(0,0)$ be the singular
point of analytic varieties $V(g_j)=\{(y,z):g_j(y,z)=0\}$ and
$V(yg_j)=\{(y,z):yg_j(y,z)=0\}$ except for $V(g_1)$ with
$\beta_{1,0,1}=1$, satisfying the following two properties $(\ast1)$
and $(\ast2)$:
$$\align
(\ast 1) \quad &\text{$g_r$ is irreducible in $\BC\{y,z\}$ for $r\ge
1$}.
 \\
 (\ast 2) \quad &\text{$g_{r+1}$ may not be irreducible in $\BC\{y,z\}$, \quad
 \text{but \quad $g_{r+1}$ satisfies} } \\
& \Delta_{r+1}(\beta_{r+1,i,k})^{r+1}_{k=1}
>(n_{r+1}-i)n_r\Delta_r(\beta_{r,0,k})^r_{k=1} \quad
\text {for $0\le i<n_j$ where $j=2,\dots,r+1$}.
\endalign$$

{\rm Remark 14.0.1.} It was already proved by Theorem $12.0$ and the
above assumptions that the following are true:
$$ \align
\noindent &  \quad \text{$g_r$ is irreducible in $\BC\{y,z\}$} \\
&  \iff \text{$g_1,\dots,g_{r-1}$ are irreducible in $\BC\{y,z\}$
and
   $\f{\Delta_{r}(\beta_{r,i,k})^{r}_{k=1}}{n_{r}-i}
> \f{\Delta_{r}(\beta_{r,0,k})^{r}_{k=1}}{n_{r}}$ for $0< i<n_r$} \qquad \qquad\\
& \iff {\f{\Delta_{j}(\beta_{j,i,k})^{j}_{k=1}}{n_{j}-i}}
> {\f{\Delta_{j}(\beta_{j,0,k})^{j}_{k=1}}{n_{j}}} \quad
\text {for each $j=1,2,\dots,r$ and for $0< i<n_j$.}
\endalign
$$

$\underline{\text{\bf {Conclusions}}}$

{\rm (1)} Let $V(g_j)=\{(y,z):g_j(y,z)=0\}$ and
$V(yg_j)=\{(y,z):yg_j(y,z)=0\}$ be analytic varieties at the origin
in $\BC^2$ for \text{\rm $1\le j\le r$}. For each $j=1,2,\dots,r$,
we may assume that
$\tau_{\lambda_j}=\pi_1\circ\pi_2\circ\cdots\circ\pi_{\lambda_j}
:M^{(\lambda_j)}\to\BC^2$ is the composition of a finite number
$\lambda_j$ of successive blow-ups which is needed to get the
standard resolution of the singular point $(0,0)$ of either $V(g_j)$
or $V(yg_j)$. If $r=1$ and $\beta_{1,0,1}=1$, then $V(g_1)$ has no
singularity at the origin, but $V(yg_1)$ has an isolated singular
point at the origin. \ms

{\rm (2)} Let $\tau^{-1}_{\lambda_j}(0,0)=\cup^{\lambda_j}_{i=1}E_i$
for each fixed $j=1,2,\dots,r$, where each $E_i$ is called an
exceptional curve of the first kind. We may assume by \text{\rm
Proposition 14.1} that $(g_j\circ\tau_{\lambda_j})_{divisor}$ and
$(g_{j+1}\circ\tau_{\lambda_j})_{divisor}$ are the divisors of
$g_j\circ\tau_{\lambda_j}$ and $g_{j+1}\circ\tau_{\lambda_j}$,
respectively, which can be defined by
$$\align
(g_j\circ\tau_{\lambda_j})_{divisor}
&=V^{(\lambda_j)}(g_j)+\sum^{\lambda_j}_{i=1}e_{j,i}E_i, \tag 14.0.2
 \\
(g_{j+1}\circ\tau_{\lambda_j})_{divisor}
&=V^{(\lambda_j)}(g_{j+1})+\sum^{\lambda_j}_{i=1}e_{j+1,i}E_i,
\endalign$$
where $V^{(\lambda_j)}(g_j)$ and $V^{(\lambda_j)}(g_{j+1})$ are the
proper transforms of $V(g_j)$ and $V(g_{j+1})$ under
$\tau_{\lambda_j}$, respectively.\ms

Then, we have the following:

{\rm(2a)} $e_{j+1,i}=n_{j+1}e_{j,i}$ for $i=1,2,\dots,\lambda_j$.
 \ms

{\rm(2b)} Let $d_j=n_{j+1}n_{j+2}\cdots n_r$ for $1\le j<r$, and
$d_r=1$. In particular, we have
 $$\align
 e_{r,\lambda_j}=n_jd_j\Delta_j(\beta_{j,0,k})^j_{k=1} \quad
 \text{with  \quad $e_{r,\lambda_1}=n_1d_1\beta_{1,0,1}$.} \tag 14.0.3 \\
\endalign$$

{\rm(3)} As far as the standard resolution of the singular point of
$V(g_r)$ is concerned, each $E_i$ of $\lambda_r$ exceptional curves
of the first kind has at most three distinct intersection points
with other exceptional curves and the proper transform under
$\tau_{\lambda_r}$:

{\rm(3a)} If $\beta_{1,0,1}=1$, then $g_r\in$ the type $[r-1]$ under
$\tau_{\lambda_r}$ in the sense of {\rm Definition 2.5} such that
$E_{\lambda_2},E_{\lambda_3},\dots,E_{\lambda_r}$ are the only
$(r-1)$ exceptional curves, each of which has three distinct
intersection points with other exceptional curves and the proper
transform.

{\rm(3b)} If $\beta_{1,0,1}\ge 2$, then $g_r\in$ the type $[r]$
under $\tau_{\lambda_r}$  in the sense of {\rm Definition 2.5} such
that $E_{\lambda_1},E_{\lambda_2},\dots,E_{\lambda_r}$ are the only
$r$ exceptional curves, each of which has three distinct
intersection points with other exceptional curves and the proper
transform. $\square$
\endproclaim \bs

{\bf \S14.1. In preparation for the proof of Theorem 14.0 } \ms

Let $\{g_j:j=1,2,\dots,r+1\}$ with $g_j\in\BC\{y,z\}$ be defined as
in Theorem $14.0$. For the proof of Theorem $14.0$, first of all, it
suffices to prepare the proposition(Proposition $14.1$) with its
proof, by the following observation:

\definition{Observation}

(1) For an easy proof of the theorem, Proposition $14.1$ is the
restatement of Theorem $14.0$.

(2) The proof of Proposition $14.1$ will be by induction on the
integer $j$, which says that $\tau_{\lambda_j}$ of Theorem $14.0$ is
the composition of a finite number $\lambda_j$ of successive
blow-ups, which is needed to get the standard resolution of the
singular point of $V(g_j)$ or $V(yg_j)$.

(3) Proposition $14.1$ is a generalization of Sublemma $12.4$ of
Theorem $12.0$, because it is clear that if the integer $j$ of
Proposition $14.1$ is equal to an integer one, then Sublemma $12.4$
and Proposition $14.1$ are the same statement.

(4) For the induction proof of Proposition $14.1$, we will prepare
three statements, denoted by Sublemma 14.2, Sublemma 14.3 and
Proposition $14.1$ in this section. After then, these three
statement will be proved in $\S14.2$, respectively. If then, we can
finish the proof of this theorem by Proposition $14.1$, without any
additional assumption.
\enddefinition \ms

\proclaim{Proposition 14.1(Theorem 14.0)}  $\underline{\text{\bf
Assumptions}}$ Suppose that the assumptions and notations of Theorem
$14.0$ hold. Note by Theorem $12.0$ and the above assumptions of
Theorem $14.0$ that the following are true:
$$\align
\text{\rm(14.1.1)} \qquad &\text{$g_r$ is irreducible in
$\BC\{y,z\}$ with $yg_r\in \text{the type}
[r]$ under}\\
&\text{the standard resolution, but $g_{r+1}$ may not be irreducible
in $\BC\{y,z\}$.} \qquad \qquad
\endalign$$

$\underline{\text{\bf Conclusions}}$ For each $j=1,2,\dots,r$, let
$V(g_j)=\{(y,z):g_j(y,z)=0\}$ and $V(yg_j)=\{(y,z):yg_j(y,z)=0\}$ be
analytic varieties at the origin in $\BC^2$.

Let $\tau_{\lambda_j}:M^{(\lambda_j)}\to\BC^2$ be the composition of
a finite number $\lambda_j$ of successive blow-ups which is needed
to get the standard resolution of the singular point of $V(g_j)$ or
$V(yg_j)$. If $j\ge 2$ or $\beta_{1,1}\ge 2$, then
$\tau_{\lambda_j}$ can be chosen with the same local coordinate
system, in order to study the resolution of the singularity of both
$V(g_j)$ and $V(yg_j)$, noting that both $V(g_j)$ and $V(yg_j)$ have
the singular point at the origin. But if $j=1$ and $\beta_{1,1}=1$,
then $V(g_1)$ has no singularity at the origin, but note that
$V(yg_1)$ still has the singular point at the origin. \ms

{\rm(a)} For any $s=1,2,\dots,\lambda_j-1$, $V^{(s)}(yg_{\ell})$
with $j\le {\ell}\le r$ has one and only one quasisingular point on
$\tau^{-1}_s(0,0)$ in the sense of \text{\rm Definition $2.6$}.

{\rm(a1)} Then, we can use just one coordinate patch of the local
coordinates for each blow-up $\pi_i$ of
$$
\tau_{\lambda_j}=\pi_1\circ\pi_2\circ\cdots \circ\pi_{\lambda_j}
\tag 14.1.2
$$
with $1\le i\le \lambda_j$ in the sense of \text{\rm Definition
$2.11$}, to study $V^{(s)}(yg_{\ell})$.

{\rm(a2)} Both $V^{(s)}(g_{r+1})$ and $V^{(s)}(yg_{r+1})$ have the
same quasisingular point on $\tau^{-1}_s(0,0)$ as $V^{(s)}(yg_j)$
does on $\tau^{-1}_s(0,0)$ for $1\le s\le \lambda_j-1$.

Then, we can use a common one coordinate patch of the local
coordinates for each blow-up $\pi_i$ of $\tau_{\lambda_j}$ with
$1\le \lambda_j$ in the sense of \text{\rm Definition $2.13$}, in
order to study both $V^{(s)}(g_{\ell})$ and $V^{(s)}(yg_{\ell})$
with $j\le {\ell}\le r+1$. \ms

{\rm(b)} Suppose that $\tau_{\lambda_j}:M^{(\lambda_j)}\to\BC^2$
satisfies the same assumptions and notations as in $(a)$. In more
detail, $\tau_{\lambda_j}$ can be written in the form
$$\align
\text{\rm(14.1.3)}\quad &
\tau_{\lambda_j}=\mu_{1,m_1}\circ\mu_{2,m_2}\circ\cdots
\circ\mu_{j,m_j}
\quad \text{with $m_1=\lambda_1$}, \\
& M^{(\lambda_j)} @> \mu_{j,m_j} >> M^{(\lambda_{j-1})}@>
\mu_{j-1,m_{j-1}} >>
 M^{(\lambda_{j-2})} @>  >> \cdots @>  >>
 M^{(\lambda_1)} @> \mu_{1,m_1} >> M^{(\lambda_0)}=\BC^2, \\
where \quad & \\
(14.1.4) \quad &\text{$\lambda_j=m_1+m_2+\cdots +m_j$
for $1\le j\le r$},   \\
 &\text{$\mu_{1,m_1}=\pi_1\circ\pi_2\circ\cdots
 \circ\pi_{\lambda_1}$
 with $\lambda_1=m_1$}, \\
 &\text{$\mu_{j,m_j}=\pi_{\lambda_{j-1}+1}\circ\pi_{\lambda_{j-1}+2}\circ\cdots
 \circ\pi_{\lambda_j}$
 with $\lambda_{j-1}+m_j=\lambda_j$, \quad $2\le j\le r$,} \qquad
 \qquad
\endalign$$
 satisfying the following properties:

{\rm (b1)} $\tau_{\lambda_1}=\mu_{1,\lambda_1}$ is the composition
of a finite number $\lambda_1$ of successive blow-ups which is
needed to get the standard resolution of the singular point of
either $V(yg_1)$ or $V(g_1)$.

{\rm (b2)} $\mu_{j,m_j}$ is the composition of a finite number $m_j$
of successive blow-ups which is needed to get the standard
resolution of one and only one quasisingular point of
$V^{(\lambda_{j-1})}(g_j)$ with $2\le j\le r$. That is,
$\tau_{\lambda_j}=\mu_{1,m_1}\circ\mu_{2,m_2}\circ\cdots
\circ\mu_{j,m_j}$. \ms

Now, we are going to compute the representation of the local
defining equation for $V^{(\lambda_j)}(g_{j+\ell})$, denoted by
$(g_{j+\ell}\circ\tau_{\lambda_j})_{proper}$, which is called the
proper transform of $V(g_{j+\ell})$ at $(y,z)=(0,0)$ under
$\tau_{\lambda_j}$, where $j+\ell\le r+1$ with $1\le j\le r$ and
$\ell\ge 0$.

Just for notation, let $(v_{\lambda_j},u_{\lambda_j})$ be a common
one coordinate patch of the same local coordinates for
$M^{(\lambda_j)}$, that is, also the $\lambda_j$-th common
coordinate patch such that $\pi_{\lambda_j}:M^{(\lambda_j)}\to
M^{(\lambda_j-1)}$ is the $\lambda_j$-th blow-up of
$\tau_{\lambda_j}$ at a quasisingular point of
$V^{(\lambda_j-1)}(g_j)$ in the sense of $(a2)$. Again, let
$(v_{\lambda_{j-1}},u_{\lambda_{j-1}})$ be a common one coordinate
patch of the same local coordinates for $M^{(\lambda_{j-1})}$ such
that $\pi_{\lambda_{j-1}}:M^{(\lambda_{j-1})}\to
M^{(\lambda_{j-1}-1)}$ is the $\lambda_{j-1}-th$ blow-up of
$\tau_{\lambda_{j-1}}$ at a quasisingular point of
$V^{(\lambda_{j-1}-1)}(g_j)$ in the sense of $(a2)$.

Being viewed as an analytic map by \text{\rm Sublemma $12.4$} and
\text{\rm Theorem $3.6$}, along $v_{\lambda_j}=0$,
$\mu_{j,m_j}:M^{(\lambda_j)}\to M^{(\lambda_{j-1})}$ can be written
in the form
$$
\align \mu_{j,m_j}(v_{\lambda_j},u_{\lambda_j})
&=(v_{\lambda_{j-1}},1+{\ve}_{j-1,1}u_{\lambda_{j-1}}) \tag 14.1.5 \\
&=(v^{p_j}_{\lambda_j}u^{a_j}_{\lambda_j},v^{q_j}_{\lambda_j}u^{b_j}_{\lambda_j}),
\endalign
$$

where

{\rm(i)} for $j=1$, $p_1=n_1$ and
$q_1=\beta_{1,1}=\Delta_1(\beta_{1,1})$ with
$(v_{\lambda_0},1+{\ve}_{0,1}u_{\lambda_0})=(y,z)$,

{\rm(ii)} for $j\ge 2$, $p_j=n_j$ and
$q_j=\Delta_j(\beta_{j,k})^j_{k=1}-n_jn_{j-1}\Delta_{j-1}
(\beta_{j-1,k})^{j-1}_{k=1}>0$,
$\ve_{j-1,1}=\ve_{j-1,1}(v_{\lambda_{j-1}},1+{\ve}_{j-1,1}u_{\lambda_{j-1}})$
is a unit in
$\BC\{v_{\lambda_{j-1}},1+{\ve}_{j-1,1}u_{\lambda_{j-1}}\}$, and for
brevity of notations $\ve_{j-1,1}$ may be defined by one, using a
nonsingular change of coordinates from
$\BC\{v_{\lambda_{j-1}},1+{\ve}_{j-1,1}u_{\lambda_{j-1}}\}$ to
itself, independently of analytic representation of
$V^{(\lambda_{j-1})}(g_i)$ for all $i\ge j-1$,

{\rm(iii)} $a_j>0$ and $b_j\ge 0$ are some nonnegative integers with
$a_jq_j-b_jp_j=1$,

{\rm(iv)} $E_{\lambda_j}=\{v_{\lambda_j}=0\}$ is defined by the
$\lambda_j-th$ exceptional curve of the first kind under
$\tau_{\lambda_j}$, for notation. \ms

For simplicity of the notations in {\rm(14.1.5)}, we may use the
following, if necessary:
$$\align
\text{\rm(14.1.6)} \qquad {\bar v}_j=v_{\lambda_j},~
\bar{u}_j=u_{\lambda_j},~ \bar{\ve}_{j,k}=\ve_{j,k} \quad \text{for
$1\le j\le r$ and $1\le k\le r-j+2$.} \qquad \qquad
\endalign$$
Note that $\bar{\ve}_{j,k}=\bar{\ve}_{j,k}({\bar
v}_j,1+\bar{\ve}_{j,1}{\bar u}_j)$ and
$\ve_{j,k}=\ve_{j,k}(v_{\lambda_j},1+\ve_{j,1} u_{\lambda_j})$ are
the same units in $\BC\{v_{\lambda_j},1+\ve_{j,1} u_{\lambda_j}\}$
for $1\le j\le r$. \ms

{\rm(c)} By induction on the integer $j=1,2,\dots,r$, at
$(v_{\lambda_{j}},1+{\ve}_{j,1}u_{\lambda_{j}})=(0,0)$ along
$v_{\lambda_j}=0$, $(g_j\circ\tau_{\lambda_j})_{total}$ with
$(g_j\circ\tau_{\lambda_j})_{proper}$ is written as follows:

\text{\rm For each {$\ell$}=1,2,\dots,r-j},
$$\align
(14.1.7) \qquad \qquad (g_j\circ\tau_{\lambda_j})_{total}
&=v^{e_{j,\lambda_j}}_{\lambda_j}(g_j\circ\tau_{\lambda_j})_{proper}
\quad \text{with}   \qquad \qquad \qquad \\
(g_j\circ\tau_{\lambda_j})_{proper} &=1+\ve_{j,1} u_{\lambda_j}
\quad
\text{with $\ve_{j,1}=1$},\\
(g_{j+1}\circ\tau_{\lambda_j})_{total}
&=v^{s^{(j)}_1e_{j,\lambda_j}}_{\lambda_j}(g_{j+1}\circ\tau_{\lambda_j})_{proper}
\quad \text{with}\\
(g_{j+1}\circ\tau_{\lambda_j})_{proper}
&=(1+\ve_{j,1}u_{\lambda_j})^{s^{(j)}_1}++\sum^{s^{(j)}_1-1}_{i=0}\eta^{(j)}_{1i}
v^{\gamma^{(j)}_{1i1}}_{\lambda_j}(1+\ve_{j,1}u_{\lambda_j})^{i},\\
(g_{j+{\ell}}\circ\tau_{\lambda_j})_{total}
&=v^{s^{(j)}_{{\ell}}s^{(j)}_{{\ell}-1}\cdots
s^{(j)}_1e_{j,\lambda_j}}
(g_{j+{\ell}}\circ\tau_{\lambda_{j}})_{proper} \quad \text{with}\\
(g_{j+{\ell}}\circ\tau_{\lambda_j})_{proper}
&=(g_{j+{\ell}-1}\circ\tau_{\lambda_j})^{s^{(j)}_{{\ell}}}_{proper} \\
& +\sum^{s^{(j)}_{{\ell}}-1}_{i=0}\eta^{(j)}_{{\ell},i}
v^{\gamma^{(j)}_{{\ell},i,1}}_{\lambda_j}
\{\prod^{\ell}_{k=2}(g_{j+k-2}\circ
\tau_{\lambda_j})^{\gamma^{(j)}_{{\ell},i,k}}_{proper}\}
(g_{j+{\ell}-1}\circ\tau_{\lambda_j})^{i}_{proper}, \\
\endalign$$
$$\align
\noindent  (14.1.8)&  \quad {where}  \quad  \text{\rm for each {$\ell$}=1,2,\dots,r-j}, \\
&\text{\rm (c0)}  \quad  e_{j,\lambda_j}= n_j\Delta_j(\beta_{j,0,k})^j_{k=1},  \\
&\text{\rm (c1)}  \quad    s^{(j)}_1 = n_{j+1},  \\
\qquad & \qquad \gamma^{(j)}_{1,i,1}=
\Delta_{j+1}(\beta_{j+1,i,k})^{j+1}_{k=1}
-(n_{j+1}-i)n_j\Delta_j(\beta_{j,i,k})^j_{k=1}>0, \qquad \qquad  \\
&\text{\rm (c,{$\ell$}+1)} \qquad s^{(j)}_{{\ell}+1} = n_{j+{\ell}+1},  \\
&\qquad \gamma^{(j)}_{{\ell}+1,i,1} =
\Delta_{j+1}(\beta_{j+{\ell}+1,i,k})^{j+1}_{k=1}
+\{\beta_{j+{\ell}+1,i,j+2}+n_{j+1}\beta_{j+{\ell}+1,i,j+3}+\cdots  \\
& \quad \qquad {+n_{j+1}n_{j+2}}\cdots
n_{j+{\ell}-1}\beta_{j+{\ell}+1,i,j+{\ell}+1}
-(\prod^{{\ell}+1}_{k=1} n_{j+k})\}n_j\Delta_j(\beta_{j,i,k})^j_{k=1}>0, \qquad \\
&\qquad \gamma^{(j)}_{{\ell}+1,i,2}= \beta_{j+{\ell}+1,i,j+2},
\gamma^{(j)}_{{\ell}+1,i,3}= \beta_{j+{\ell}+1,i,j+3}, ~ \dots, ~
\gamma^{(j)}_{{\ell}+1,i,{\ell}+1}=\beta_{j+{\ell}+1,i,j+{\ell}+1},
\qquad \qquad \qquad \qquad
\endalign$$
and
$\eta^{(j)}_{{{\ell}+1},i}=\eta^{(j)}_{{{\ell}+1},i}(v_{\lambda_j},1+\ve_{j,1}
u_{\lambda_j})$ is a unit in $\BC\{v_{\lambda_j},1+\ve_{j,1}
u_{\lambda_j}\}$ for \text{\rm $1\le {{\ell}}\le r-j$} and $0\le
i\le s_{\ell}$. \ms

For brevity of notation, $\ve_{j,1}$ may be defined by one, using a
nonsingular change of coordinates from
$\BC\{v_{\lambda_j},1+\ve_{j,1} u_{\lambda_j}\}$ to itself,
independently of analytic representation of $V^{(\lambda_{j})}(g_i)$
for all $i\ge j$.  \ms

{\rm(d)} Let $\tau^{-1}_{\lambda_j}(0,0)=\cup^{\lambda_j}_{i=1}E_i$
where each $E_i$ is an exceptional curve of the first kind. For
$0\le \ell\le r-j+1$, let
$(g_{j+\ell}\circ\tau_{\lambda_j})_{divisor}$ be the divisor of
$g_{j+\ell}\circ\tau_{\lambda_j}$ defined by
$$
(g_{j+\ell}\circ\tau_{\lambda_j})_{divisor}=V^{(\lambda_j)}(g_{j+\ell})
 +\sum^{\lambda_j}_{i=1}e_{j+\ell,i}E_i, \tag 14.1.9
$$where $V^{(\lambda_j)}(g_{j+\ell})$ is the proper transform of
$V(g_{j+\ell})$ under $\tau_{\lambda_j}$ and
$\lambda_j=m_1+m_2+\cdots +m_j$ as we have seen in $(14.1.4)$. \ms

Then, we have the following:

{\rm(d1)} For each $\ell=1,2,\dots,r-j+1$,
$$ e_{j+\ell,i}=n_{j+\ell}n_{j+\ell-1}\cdots n_{j+1}e_{j,i} \quad \text{for $1\le i\le
\lambda_{j}$.} \tag 14.1.10 $$

In particular,
$e_{{j},\lambda_{j}}=n_{j}\Delta_{j}(\beta_{j,k})^{j}_{k=1}$ by
{\rm(c0)} of $(14.1.7)$. \ms

{\rm(d2)} For any $\ell=0,1,2,\dots,r-j+1$, where $E_{\lambda_j}$ is
defined by $\{v_{\lambda_{j}}=0\}$,
$$\align
\text{\rm (14.1.11)} \quad  V^{(\lambda_j)}(g_{j+\ell})\cap
(\cup^{\lambda_j}_{i=1}E_i) &=V^{(\lambda_j)}(g_j)\cap E_{\lambda_j}
=\{(v_{\lambda_{j}},1+{\ve}_{j,1}u_{\lambda_{j}})=(0,0)\} \qquad
\qquad
\endalign$$

{\rm(d3)} We get the following in the sense of {\rm Definition
$2.6$}: Note that $\gcd(n_r,\Delta_r(\beta_{r,k})^r_{k=1})=1$
because $g_r$ is irreducible in $\BC\{y,z\}$ by assumption.

If $\beta_{1,1}\ge 2$, then $g_k\in$ the type $[j]$ under
$\tau_{\lambda_j}$ for any $k=j,j+1,\dots,r+1$.

If $\beta_{1,1}=1$, then $g_k\in$ the type $[j-1]$ under
$\tau_{\lambda_j}$ for any $k=j,j+1,\dots,r+1$.

Whether $\beta_{1,1}\ge 2$ or $\beta_{1,1}=1$, note that $yg_k\in$
the type $[j]$ under $\tau_{\lambda_j}$ for any $k\ge j$. $\square$
\endproclaim \ms

\definition{Remark 14.1.1} Moreover,
$(y^{\beta_{j+\ell,i,1}}z^{\beta_{j+\ell,i,2}}g^{\beta_{j+\ell,i,3}}_1
g^{\beta_{j+\ell,i,4}}_2\cdots
g^{\beta_{j+\ell,i,j+\ell}}_{j+\ell-2}g^{i}_{j+\ell-1})
\circ\tau_{\lambda_j}(v_j,u_j)$ can be viewed as
$$
 \ve_{j,\ell+1}v_{\lambda_j}^{\Delta^{\natural}_{j+1}(\beta_{j+\ell,k})^{j+1}_{k=1}}
 (\prod^{\ell}_{k=2}
 (g_{j+k-2}\circ\tau_{\lambda_j})^{\beta_{j+\ell,j+k}}_{proper})
 (g_{j+{\ell}}\circ\tau_{\lambda_j})^{i}_{proper}, \tag 14.1.12
$$
where {\rm(i)} \quad the $\beta_{j+\ell,j+k}$ are nonnegative
integers for $1\le k\le \ell$,
$$\align
\quad \text{\rm (ii)} \quad & \text{for brevity of notation, we
write $\Delta^{\natural}_{j+1}(\beta_{j+\ell,k})^{j+1}_{k=1}=
\Delta_{j+1}(\beta_{j+\ell,k})^{j+1}_{k=1}$} \\
& +(\beta_{j+\ell,j+2}+n_{j+1}\beta_{j+\ell,j+3}+\cdots
+n_{j+1}n_{j+2}\cdots
n_{j+\ell-2}\beta_{j+\ell,j+\ell})n_j\Delta_j(\beta_{j,k})^j_{k=1},
\endalign$$

{\rm(iii)} \quad
$\ve_{j,\ell+1}=\ve_{j,\ell+1}(v_{\lambda_j},1+\ve_{j,1}
u_{\lambda_j})$ is a unit in $\BC\{v_{\lambda_j},1+\ve_{j,1}
u_{\lambda_j}\}$. \ms

Using $e_{j,\lambda_j}=n_j\Delta_j(\beta_{j,0,k})^j_{k=1}$ in
{\rm(c0)} of $(14.1.8)$ and the notation in {\rm(ii)} of
$(14.1.12)$, for each $\ell=1,2,\dots,r-j+1$,
$\gamma^{(j)}_{\ell,i,1}$ in $(c,\ell)$ of $(14.1.8)$ can be
rewritten in the form
$$\align
\text{\rm (14.1.13)} \qquad
\gamma^{(j)}_{\ell,i,1}&=\Delta_{j+1}(\beta_{j+\ell,i,k})^{j+1}_{k=1}
+\{\beta_{j+\ell,i,j+2}
+n_{j+1}\beta_{j+\ell,i,j+3}+n_{j+1}n_{j+2}\beta_{j+\ell,i,j+4}+\cdots \qquad \\
&\quad +n_{j+1}n_{j+2}\cdots n_{j+\ell-2}\beta_{j+\ell,i,j+\ell}
-n_{j+1}n_{j+2}\cdots
n_{j+\ell-1}n_{j+\ell}\}e_{j,\lambda_j}\\
&=\Delta^{\natural}_{j+1}(\beta_{j+\ell,i,k})^{j+1}_{k=1}-n_{j+1}n_{j+2}\cdots
n_{j+\ell-1}n_{j+\ell}e_{j,\lambda_j} \\
&=\Delta^{\natural}_{j+1}(\beta_{j+\ell,i,k})^{j+1}_{k=1}-e_{j+\ell,\lambda_j}>0
\quad \text{\rm by (14.1.10).} \\
\endalign$$
\enddefinition

\proclaim{Sublemma 14.2} $\underline{\text{\bf {Assumptions}}}$
\quad Let $j$ be an arbitrary positive integer with $1\le j\le r$.
Suppose that the assumptions and notations of Proposition $14.1$
hold. Note by Theorem $12.0$ and the above assumptions of
Proposition $14.1$ that the following is true:
$$\align
\text{\rm(14.2.1)} \qquad &\text{$g_r$ is irreducible in
$\BC\{y,z\}$ with $yg_r\in \text{the type}
[r]$ under}\\
&\text{the standard resolution, but $g_{r+1}$ may not be irreducible
in $\BC\{y,z\}$.} \qquad \qquad \qquad
\endalign$$

$\underline{\text{\bf Conclusions}}$ \quad For each $j$, let
$\tau_{\lambda_j}:M^{(\lambda_j)}\to\BC^2$ be the composition of a
finite number $\lambda_j$ of successive blow-ups which is needed to
get the standard resolution of the singular point of an analytic
variety defined by $V(g_j)$ or $V(yg_j)$.

Then, we can construct three sequences, denoted by $\text{\bf
Sequences[I]}^{\text{{\bf(j)}}}$:
$$\align
(14.2.2)  \qquad &\{Y^{(j)}_\ell:\ell=1,2,\dots,r-j+1\} \quad
\text{with $Y^{(j)}_\ell\subset N_0$},  \\
&\{h^{(j)}_\ell:\ell=1,2,\dots,r-j+1\} \quad \text{with
$h^{(j)}_\ell=(g_{j+\ell}\circ\tau_{\lambda_{j}})_{proper} \in
\BC\{v_{\lambda_{j}},1+{\ve}_{j,1}u_{\lambda_{j}}\}$,} \qquad \qquad \\
&\{\text{$\Xi^{(j)}_\ell:N^\ell_0\to N_0$ is an integer-valued
function for $\ell=1,2,\dots,r-j+1$}\}, \\
&\text{satisfying the following five conditions for each k :} \\
\endalign$$

\noindent Five conditions are denoted by \text{\bf The 1st
${\text{\bf{Cond}}}^{\text{{\bf(j)}}}$}{\bf, \dots,} \text{\bf The
5-th ${\text{\bf{Cond}}}^{\text{{\bf(j)}}}$ of $\text{\bf
Sequences[I]}^{\text{{\bf(j)}}}$}. \ms

\noindent \text{\bf The 1st ${\text{\bf{Cond}}}^{\text{{\bf(j)}}}$:}
Let $\{Y^{(j)}_\ell:\ell=1,2,\dots,r-j+1\}$ with
$Y^{(j)}_\ell\subset N_0$ be defined by {\rm(14.2.3)}.

\roster
\item"(14.2.3)(1a)"$Y^{(j)}_1=\{s^{(j)}_1\}\cup
\{\gamma^{(j)}_{{1},i,1}:0\le i<s^{(j)}_1 \}$ with $s^{(j)}_1\ge 2$
and $\gamma^{(j)}_{1,0,1}\ge 1$ where $Y^{(j)}_1\subset N$,

\item"(1b)" $Y^{(j)}_{\ell}=\{s^{(j)}_{\ell}\}\cup
\{\gamma^{(j)}_{{\ell},i,1}:0\le i<s^{(j)}_{\ell} \} \cup
\{\gamma^{(j)}_{{\ell},i,2}:0\le i<s^{(j)}_{\ell} \}\cup \cdots \cup
\{\gamma^{(j)}_{{\ell},i,{\ell}}:0\le i<s^{(j)}_{\ell} \}$ with
$s^{(j)}_{\ell}\ge 2$ where $j=1,2,\dots,r$,
\endroster

such that for each $j=1,2,\dots,r$, assume that at least one of
$\gamma^{(j)}_{{\ell},0,1},\gamma^{(j)}_{{\ell},0,2},\dots,
\gamma^{(j)}_{{\ell},0,{\ell}}$ is nonzero by Sublemma $12.1$ and
\ms

such that each of {\rm(14.2.3)} satisfies either {\rm(a)} or
{\rm(b)} of {\rm(14.2.4)}:

\roster
\item"(14.2.4)(a)" if $j=1$, then the first family
$\{Y^{(1)}_\ell:\ell=1,2,\dots,r\}$ is defined as follows:

\item"(a0)" $e_{1,\lambda_1}=n_1\Delta_1(\beta_{1,0,1})=n_1\beta_{1,0,1}$,

\item"(a1)" $s^{(1)}_1=n_2\ge 2$,
$\gamma^{(1)}_{1,i,1}=\Delta^{\sharp}_2(\beta_{2,i,1},\beta_{2,i,2})
-(n_2-i)n_1\beta_{1,0,1}>0$ for $0\le i<n_2$,

\item"(a2)" $s^{(1)}_{\ell} =n_{{\ell}+1}\ge 2$,
$\gamma^{(1)}_{{\ell},i,1}=\Delta^{\sharp}_{{\ell}+1}(\beta_{{\ell}+1,i,k})^{{\ell}+1}_{k=1}
-(n_{{\ell}+1}-i)n_{\ell}\cdots n_2n_1\beta_{1,0,1}>0$, and

\noindent $\gamma^{(1)}_{{\ell},i,2}=\beta_{{\ell}+1,i,3}$,
$\gamma^{(1)}_{{\ell},i,3}=\beta_{{\ell}+1,i,4}$, $\dots $,
$\gamma^{(1)}_{{\ell},i,{\ell}}=\beta_{{\ell}+1,i,{\ell}+1}$ for
$\ell\ge 2$ and $0\le i<n_{\ell}$,
\endroster
noting that
$\gamma^{(1)}_{1,0,1},\gamma^{(1)}_{2,0,1},\dots,\gamma^{(1)}_{{\ell},0,1}$
are positive by {\rm Sublemma 12.1}. \ms

\roster
\item"(14.2.4)(b)" for each $j=2,3,\dots,r$,  the family
$\{Y^{(j)}_\ell:\ell=1,2,\dots,r-j+1\}$ is defined as follows:

\item"(b0)" $e_{j,\lambda_j}=n_j\Delta_j(\beta_{j,0,k})^j_{k=1}$,

\item"(b1)"$s^{(j)}_1=s^{(j-1)}_2=s^{(j-2)}_3=\cdots=s^{(1)}_j=n_{j+1}\ge 2$,
$\gamma^{(j)}_{1,i,1}=\Xi^{(j-1)^{\sharp}}_2(\gamma^{(j-1)}_{2,i,1},\gamma^{(j-1)}_{2,i,2})
-(s^{(j-1)}_2-i)s^{(j-1)}_1\gamma^{(j-1)}_{1,0,1}>0$ for $0\le
i<s^{(j-1)}_2$,

\item"(b2)"
$s^{(j)}_{\ell}=s^{(j-1)}_{\ell+1}=s^{(j-2)}_{\ell+2}=\cdots
=s^{({\ell})}_{j}=n_{j+{\ell}}$.

\noindent $\gamma^{(j)}_{{\text{\rm
$\ell$}},i,1}=\Xi^{(j-1)^{\sharp}}_{\ell+1}(\gamma^{(j-1)}_{{\text{\rm
$\ell$+1}},i,k})^{\text{\rm {$\ell$}+1}}_{k=1}
-(s^{(j-1)}_{{\text{\rm $\ell+1$}}}-i)s^{(j-1)}_{{\text{\rm
$\ell$}}} \cdots s^{(j-1)}_2 s^{(j-1)}_1\gamma^{(j-1)}_{1,0,1}>0$,
$\gamma^{(j)}_{{\ell},i,2}=\gamma^{(j-1)}_{{\ell+1},i,3}$,
$\gamma^{(j)}_{{\ell},i,3}=\gamma^{(j-1)}_{{\ell+1,}i,4}$,
$\gamma^{(j)}_{{\ell},i,4}=\gamma^{(j-1)}_{{\ell+1,}i,5}$,$\dots$,
$\gamma^{(j)}_{\ell,i,\ell}=\gamma^{(j-1)}_{\ell+1,i,\ell+1}$ for
$0\le i<s^{(j-1)}_{{\ell}+1}$,
\endroster \ms

\noindent \text{\bf The 2-th
${\text{\bf{Cond}}}^{\text{{\bf(j)}}}$:} For $\ell=1,2,\dots,
\text{\rm r-j+1}$, let
$h^{(j)}_\ell=(g_{j+\ell}\circ\tau_{\lambda_j})_{proper}$ be defined
by
$$\align
(14.2.5) \qquad \qquad (g_j\circ\tau_{\lambda_j})_{total}
&=v^{e_{j,\lambda_j}}_{\lambda_j}(g_j\circ\tau_{\lambda_j})_{proper}
\quad \text{with}   \qquad \qquad \qquad \\
(g_j\circ\tau_{\lambda_j})_{proper} &=1+\ve_{j,1} u_{\lambda_j}
\quad
\text{with $\ve_{j,1}=1$},\\
(g_{j+1}\circ\tau_{\lambda_j})_{total}
&=v^{s^{(j)}_1e_{j,\lambda_j}}_{\lambda_j}(g_{j+1}\circ\tau_{\lambda_j})_{proper}
\quad \text{with}\\
(g_{j+1}\circ\tau_{\lambda_j})_{proper}
&=(1+\ve_{j,1}u_{\lambda_j})^{s^{(j)}_1}++\sum^{s^{(j)}_1-1}_{i=0}\eta^{(j)}_{1i}
v^{\gamma^{(j)}_{1i1}}_{\lambda_j}(1+\ve_{j,1}u_{\lambda_j})^{i},\\
(g_{j+{\ell}}\circ\tau_{\lambda_j})_{total}
&=v^{s^{(j)}_{{\ell}}s^{(j)}_{{\ell}-1}\cdots
s^{(j)}_1e_{j,\lambda_j}}
(g_{j+{\ell}}\circ\tau_{\lambda_{j}})_{proper} \quad \text{with}\\
(g_{j+{\ell}}\circ\tau_{\lambda_j})_{proper}
&=(g_{j+{\ell}-1}\circ\tau_{\lambda_j})^{s^{(j)}_{{\ell}}}_{proper} \\
& +\sum^{s^{(j)}_{{\ell}}-1}_{i=0}\eta^{(j)}_{{\ell},i}
v^{\gamma^{(j)}_{{\ell},i,1}}_{\lambda_j}
\{\prod^{\ell}_{k=2}(g_{j+k-2}\circ
\tau_{\lambda_j})^{\gamma^{(j)}_{{\ell},i,k}}_{proper}\}
(g_{j+{\ell}-1}\circ\tau_{\lambda_j})^{i}_{proper}, \\
\endalign$$
where $\ve_{j,1}=\ve_{j,1}(v_{\lambda_j}, u_{\lambda_j})$ is a unit
in $\BC\{v_{\lambda_j}, u_{\lambda_j}\}$, and also
$\eta^{(j)}_{{\ell},i}=\eta^{(j)}_{{\ell},i({\ell})}(v_{\lambda_j},1+\ve_{j,1}u_{\lambda_j})$
is a unit in $\BC\{v_{\lambda_{j}},1+\ve_{j,1} u_{\lambda_{j}} \}$
for $\ell=1,2,\dots,r-j+1$ and
$i=i(\ell)=0,1,\dots,{s^{(j)}_{\ell}}-1$. For brevity of notation,
$\ve_{j,1}$ may be defined by one, using a nonsingular change of
coordinates from $\BC\{v_{\lambda_j},u_{\lambda_j}\}$ to itself,
independently of analytic representation of $V^{(\lambda_{j})}(g_i)$
for all $i\ge j$.  \ms

\noindent \text{\bf The 3-th
${\text{\bf{Cond}}}^{\text{{\bf(j)}}}$.} Let $\{\Xi^{(j)}_{\ell}:
N^{\ell}_0\to N_0, {\ell}=1,2,\dots,r-j+1\}$ be a sequence such that
each $\Xi^{(j)}_k$ is an integer-valued function defined by the
following: Note that $2\le{\ell}\le r-j+1$.
$$\align
(14.2.6) \qquad &(3a)\quad  \text{$\Xi^{(j)}_1(t)=t$  for each $t\in N_0$.}  \\
& \text{\rm(3b)} \quad
\text{$\Xi^{(j)}_{\ell}(t_k)^{\ell}_{k=1}=t_{\ell}\Xi^{(j)}_{{\ell-1}}
(\gamma^{(j)}_{{\ell-1},0,k})^{{\ell-1}}_{k=1}+s^{(j)}_{{\ell-1}}
\Xi^{(j)}_{{\ell-1}}(t_k)^{{\ell-1}}_{k=1}$}
\text{ for each $(t_k)^{\ell}_{k=1}\in  N^{\ell}_0.$}    \\
\endalign$$

Moreover, for brevity of notation, let $\{\Xi^{(j)^{\sharp}}_{\ell}:
N^{\ell}_0\to N_0, {\ell}=2,3,\dots,r-j+1\}$ be a sequence such that
each $\Xi^{(j)^{\sharp}}_{\ell}$ is an integer-valued function
defined by the following:
$$\align
\text{\rm (14.2.6-1)} \quad
&\Xi^{(j)^{\sharp}}_{2}(\gamma^{(j)}_{2,i,k})^{2}_{k=1}
=\Xi^{(j)}_{2}(\gamma^{(j)}_{2,i,k})^{2}_{k=1} \quad \text{for
$0\le i<s^{(j)}_2$,} \\
&\Xi^{(j)^{\sharp}}_{\ell}(\gamma^{(j)}_{{\ell},i,k})^{\ell}_{k=1}
=\Xi^{(j)}_{2}(\gamma^{(j)}_{{\ell},i,k})^{2}_{k=1}
+s^{(j)}_1\gamma^{(j)}_{1,0,1}\gamma^{(j)}_{{\ell},i,3}
+s^{(j)}_1\gamma^{(j)}_{1,0,1}s^{(j)}_2\gamma^{(j)}_{{\ell},i,4} \\
&+s^{(j)}_1\gamma^{(j)}_{1,0,1}s^{(j)}_2s^{(j)}_3\gamma^{(j)}_{{\ell},i,5}
+\cdots +s^{(j)}_1\gamma^{(j)}_{1,0,1}s^{(j)}_2\cdots
s^{(j)}_{{\ell}-2}\gamma^{(j)}_{{\ell},i,{\ell}} \quad \text{for
$0\le i<s^{(j)}_{\ell}$.} \qquad \qquad \\
\endalign$$ \ms

\noindent \text{\bf The $(4\alpha)$-th
${\text{\bf{Cond}}}^{\text{{\bf(j)}}}$.} Then, the following
equalities hold: Note that $r\ge 2$.

{\rm (a)} For $j=1$, the following equality holds with respect to
$\{Y^{(1)}_\ell:\ell=1,2,\dots,r\}$:
$$\align
\text{\rm(14.2.7)(14.2.7-1)} \qquad\qquad
&\Xi^{(1)}_1(\gamma^{(1)}_{1,i,1})=\gamma^{(1)}_{1,i,1}
=\Delta_2(\beta_{2,i,1},\beta_{2,i,2})-(n_2-i)n_1\beta_{1,0,1}>0, \qquad \qquad\\
&\Xi^{(1)}_{\ell}(\gamma^{(1)}_{{\ell},i,k})^{\ell}_{k=1}-(s^{(1)}_{\ell}-i)s^{(1)}_{{\ell}-1}
\Xi^{(1)}_{{\ell}-1}(\gamma^{(1)}_{{\ell}-1,0,k})^{{\ell}-1}_{k=1}
 \quad \text{for $2\le {\ell}\le {r}$} \qquad \qquad\\
& \qquad \qquad =
\Delta_{{\ell}+1}(\beta_{{\ell}+1,i,k})^{{\ell}+1}_{k=1}
-(n_{{\ell}+1}-i)n_{\ell}\Delta_{\ell}(\beta_{{\ell},0,k})^{\ell}_{k=1}>0. \\
\endalign$$

{\rm (b)} For each $j=2,3,\dots,r$, the following equality holds
with respect to the family $\{Y^{(j)}_\ell:\ell=1,2,\dots,r-j+1\}$:
$$\align
\text{\rm(14.2.7)(14.2.7-2)}  \quad \quad
&\Xi^{(j)}_1(\gamma^{(j)}_{1,i,1})=\gamma^{(j)}_{1,i,1}
=\Xi^{(j-1)}_2(\gamma^{(j-1)}_{2,i,1},\gamma^{(j-1)}_{2,i,2})
-(s^{(j-1)}_2-i)s^{(j-1)}_1\gamma^{(j-1)}_{1,0,1}>0, \\
&\Xi^{(j)}_{\ell}(\gamma^{(j)}_{{\ell},i,k})^{\ell}_{k=1}-(s^{(j)}_{\ell}-i)s^{(j)}_{{\ell}-1}
\Xi^{(j)}_{{\ell}-1}(\gamma^{(j)}_{{\ell}-1,0,k})^{{\ell}-1}_{k=1}
\qquad \text{ for $2\le {\ell}\le {r-j+1}$} \qquad \\
&\quad
=\Xi^{(j-1)}_{{\ell}+1}(\gamma^{(j-1)}_{{\ell}+1,i,k})^{{\ell}+1}_{k=1}
-(s^{(j-1)}_{{\ell}+1}-i)s^{(j-1)}_{\ell}
\Xi^{(j-1)}_{\ell}(\gamma^{(j-1)}_{{\ell},0,k})^{\ell}_{k=1}>0.
\\
\endalign$$

\noindent \text{\bf The 4-th
${\text{\bf{Cond}}}^{\text{{\bf(j)}}}$.} The following inequalities
hold: Note that $r\ge 2$.
$$\align
(14.2.8) \quad \quad  &\text{\rm(4a)} \qquad \quad \quad
\Xi^{(j)}_1(\g^{(j)}_{1,i,1})=\g^{(j)}_{1,i,1}>0. \\
&\text{\rm(4b)} \quad \quad
\Xi^{(j)}_{\ell}(\g^{(j)}_{{\ell},i,k})^{\ell}_{k=1}>(s^{(j)}_{\ell}-i)
s^{(j)}_{\ell-1}\Xi^{(j)}_{\ell-1}(\g^{(j)}_{\ell-1,i,k})^{\ell-1}_{k=1}
\quad \text{for $2\le {\ell} \le r-j+1$.} \qquad \qquad \\
\endalign$$

\noindent \text{\bf The 5-th
${\text{\bf{Cond}}}^{\text{{\bf(j)}}}$.} For \text{\rm
{$\ell$}=1,2,\dots,r-j+1}, the following equalities hold: Note that
$r\ge 2$.
$$\align
\text{\rm(14.2.9)(14.2.9-1)} \quad
&\text{$\gcd(s^{(1)}_{\ell},\Xi^{(1)}_{\ell}(\gamma^{(1)}_{{\ell},0,k})^{\ell}_{k=1})
=\gcd(n_{{\ell}+1},\Delta_{{\ell}+1}(\beta_{{\ell}+1,0,k})^{{\ell}+1}_{k=1})$,
\quad for $1\le{\ell}\le r$,} \\
&\text{$\gcd(s^{(j)}_{\ell},\Xi^{(j)}_{\ell}(\gamma^{(j)}_{{\ell},0,k})^{\ell}_{k=1})
=\gcd(s^{(j-1)}_{{\ell}+1},\Xi^{(j-1)}_{{\ell}+1}
(\gamma^{(j-1)}_{{\ell}+1,0,k})^{{\ell}+1}_{k=1})$ for $j\ge 2$ with
$j+{\ell}\le r+1$.}
\endalign$$

\noindent \text{\bf The 5-th
${\text{\bf{Cond}}}^{\text{{\bf(j)}}}$.} The following inequalities
hold: Note that $r\ge 2$.
$$\align
\text{\rm(14.2.9)(14.2.9-2)} \quad
&\f{\g^{(j)}_{1,i,1}}{s^{(j)}_1-i}>
\f{\g^{(j)}_{1,0,1}}{s^{(j)}_1}>0 \quad \text {for $0<
i<s^{(j)}_1$.}
\\
&\f{\Xi^{(j)}_{\ell}(\g^{(j)}_{{\ell},i,k})^{\ell}_{k=1}}{s^{(j)}_{\ell}-i}
> \f{\Xi^{(j)}_{\ell}(\g^{(j)}_{{\ell},0,k})^{\ell}_{k=1}}{s^{(j)}_{\ell}}
>s^{(j)}_{\ell-1}\Xi^{(j)}_{\ell-1}(\g^{(j)}_{{\ell-1},0,k})^{\ell-1}_{k=1}
\quad \text {for $0< i<s^{(j)}_{\ell}$.}  \qquad \qquad \qquad \qquad\\
\endalign$$

$\underline{\text{\rm Remark for \text{\rm The 5-th
${\text{\rm{Cond}}}^{\text{{\rm(j)}}}$.}}}$ For notation, the above
inequality in $\text{\rm(14.2.9-2)}$ may be equivalently rewritten
as follows: For each ${\ell}=1,2,\dots,r-j+1$,
$$\align
&\f{\Xi^{(j)}_{\ell}(\g^{(j)}_{{\ell},i,k})^{\ell}_{k=1}
-(s^{(j)}_{{\ell}}-i)s^{(j)}_{{\ell}-1}\Xi^{(j)}_{{\ell}-1}
(\g^{(j)}_{{\ell}-1,0,k})^{{\ell}-1}_{k=1}} {s^{(j)}_{{\ell}}-i}
\\   & > \f{\Xi^{(j)}_{{\ell}}(\g^{(j)}_{{\ell},0,k})^{{\ell}}_{k=1}
-s^{(j)}_{{\ell}}s^{(j)}_{{\ell}-1}\Xi^{(j)}_{{\ell}-1}
(\g^{(j)}_{{\ell}-1,0,k})^{{\ell}-1}_{k=1}} {s^{(j)}_{{\ell}}}>0
\quad \text {for $0< i<s^{(j)}_{\ell}$.}
\endalign$$
\endproclaim \ms

\definition{Remark 14.2.1.1}
\noindent $\underline{\text{\rm The 5-th
${\text{\rm{Cond}}}^{\text{{\rm(0)}}}$}\text{\rm in the assumptions
of Theorem 14.0}}$ For each ${\ell}=1,2,\dots,r$, the following
inequalities hold:
$$\align
\noindent \text{\rm (5)}&\text{\rm(5a)} \qquad \qquad\qquad\qquad
\text{
$\gcd(n_{\ell},\Delta_{\ell}(\beta_{{\ell},0,k})^{\ell}_{k=1})=1$
\quad for $1\le {\ell}\le r$.}  \\
&\text{\rm (5b)} \qquad \qquad\qquad\qquad
\f{\Delta_{{\ell}}(\beta_{{\ell},i,k})^{{\ell}}_{k=1}}{n_{{\ell}}-i}
> \f{\Delta_{{\ell}}(\beta_{{\ell},0,k})^{{\ell}}_{k=1}}{n_{{\ell}}} \quad
\text {for $0< i<n_{\ell}$.} \qquad \qquad\qquad\qquad
\endalign$$
\enddefinition \ms

\definition{Remark 14.2.1.2} For brevity of the proof of {\bf Sublemma
14.2}, suppose that \text{\bf The 1st
${\text{\bf{Cond}}}^{\text{{\bf(j)}}}$}, \text{\bf The 2nd
${\text{\bf{Cond}}}^{\text{{\bf(j)}}}$} and \text{\bf The 3rd
${\text{\bf{Cond}}}^{\text{{\bf(j)}}}$} have proved in the
conclusions of this sublemma. Then, for the remaining proof, it
suffices to show that \noindent \text{\bf The $(4\alpha)$-th
${\text{\bf{Cond}}}^{\text{{\bf(j)}}}$} is true, because of the
following two facts:

{\rm Fact(1):} If \text{\bf The $4\alpha$-th
${\text{\bf{Cond}}}^{\text{{\bf(j)}}}$} is true, then it is clear
that \text{\bf The 4-th ${\text{\bf{Cond}}}^{\text{{\bf(j)}}}$} is
true. \ms

{\rm Fact(2):} If \text{\bf The $4\alpha$-th
${\text{\bf{Cond}}}^{\text{{\bf(j)}}}$} is true, then we can easily
prove by \text{\bf The 1-th ${\text{\bf{Cond}}}^{\text{{\bf(j)}}}$}
that \text{\bf The $5$-th ${\text{\bf{Cond}}}^{\text{{\bf(j)}}}$} is
true, by using the following elementary computation:

Case(A): $j=1$,  and Case(B): $j\ge 2$. \ms

Case(A). Let $j=1$.

\noindent{\rm(i)} $\gcd(s^{(1)}_1,\gamma_{1,0,1})
=\gcd(n_2,\Delta_2(\beta_{2,0,1},\beta_{2,0,2})-n_2n_1\beta_{1,0,1})=
\gcd(n_2,\Delta_2(\beta_{2,0,1},\beta_{2,0,2}))=1$. \ms

\noindent{\rm(ii)} For each $q=2,3,\dots,r-1$,
$s^{(1)}_{q}=n_{q+1}$, and \text{by \rm The $(4\alpha)$-th
${\text{\rm{Cond}}}^{\text{{\rm(j)}}}$} we have
$$\split
\text{\rm(14.2.10-1)} \qquad
&\gcd(s^{(1)}_q,\Xi_q(\gamma_{q,0,k})^q_{k=1})
=\gcd(s^{(1)}_q,\Xi_q(\gamma_{q,0,k})^q_{k=1}
-s^{(1)}_qs^{(1)}_{q-1}\Xi_{q-1}(\gamma_{q-1,0,k})^{j-1}_{k=1}) \qquad \qquad \\
&=\gcd(n_{q+1},\Delta_{q+1}(\beta_{q+1,0,k})^{q+1}_{k=1}
-n_{q+1}n_q\Delta_q(\beta_{q,0,k})^q_{k=1}) \quad \text{by
$(12.5.4\alpha)$}\\
&=\gcd(n_{q+1},\Delta_{q+1}(\beta_{q+1,0,k})^{j+1}_{k=1})=1.\text{\quad
$\square$}
\endsplit$$

Case(B). Let $j\ge 2$.

{\rm(i)} $\gcd(s^{(j)}_1,\gamma^{(j)}_{1,0,1})=
\gcd(s^{(j-1)}_2,\Xi^{(j-1)}_2(\gamma^{(j-1)}_{2,0,k})^2_{k=1}
-s^{(j-1)}_2s^{(j-1)}_{1}\gamma^{(j-1)}_{1,0,1})$

\qquad
$=\gcd(s^{(j-1)}_2,\Xi^{(j-1)}_2(\gamma^{(j-1)}_{2,0,k})^2_{k=1})=1$.
\ms

{\rm(ii)} For each ${\ell}=2,3,\dots,r-1$,
$s^{(j)}_{\ell}=s^{(j-1)}_{{\ell}+1}$, and \text{by \rm The
$(4\alpha)$-th ${\text{\rm{Cond}}}^{\text{{\rm(j)}}}$} we have
$$\split
\text{\rm(14.2.10-2)}  \quad
&\gcd(s^{(j)}_{\ell},\Xi^{(j)}_{\ell}(\gamma^{(j)}_{{\ell},0,k})^{\ell}_{k=1})\\
&=\gcd(s^{(j)}_{\ell},\Xi^{(j)}_{\ell}(\gamma^{(j)}_{{\ell},0,k})^{\ell}_{k=1}
-s^{(j)}_{\ell}s^{(j)}_{{\ell}-1}\Xi^{(j)}_{{\ell}-1}
(\gamma^{(j)}_{{\ell}-1,0,k})^{{\ell}-1}_{k=1}) \\
&=\gcd(s^{(j-1)}_{{\ell}+1},\Xi^{(j-1)}_{{\ell}+1}
(\gamma^{(j-1)}_{{\ell}+1,0,k})^{{\ell}+1}_{k=1}
-s^{(j-1)}_{{\ell}+1}s^{(j-1)}_{{\ell}}\Xi^{(j-1)}_{{\ell}}
(\gamma^{(j-1)}_{{\ell},0,k})^{{\ell}}_{k=1})
\qquad \qquad \\
&=\cdots \\
&=\gcd(n_{j+{\ell}},\Delta_{j+{\ell}}(\beta_{j+{\ell},0,k})^{j+{\ell}}_{k=1}
-n_{j+{\ell}}n_{j+{\ell}-1}\Delta_{j+{\ell}-1}
(\beta_{j+{\ell}-1,0,k})^{j+{\ell}-1}_{k=1}) \\
&=\gcd(n_{j+{\ell}},\Delta_{j+{\ell}}(\beta_{j+{\ell},0,k})^{j+{\ell}}_{k=1})
=1 \quad \text{by (14.2.10-1)} \qquad \text{$\square$} \\
\endsplit$$
\enddefinition \ms

\proclaim{Sublemma 14.3} $\underline{\text{\bf {Assumptions}}}$ Let
$r$ be an arbitrary integer with $r\ge 2$. Let
$\tau_{\lambda_j}:M^{(\lambda_j)}\to\BC^2$ be the composition of a
finite number $\lambda_j$ of successive blow-ups which is needed to
get the standard resolution of the singular point of $V(g_j)$ or
$V(yg_j)$. Suppose that the assumptions of Sublemma $14.2$ hold. \ms

$\underline{\text{\bf {Conclusions}}}$ The aim is to prove the
following by induction on the integer $j>0$.

{\rm(1)} As we have seen in the conclusion of Proposition $14.1$,
$\tau_{\lambda_{j}}=\mu_{1,m_1}\circ\mu_{2,m_2}\circ\cdots\circ\mu_{j,m_j}$
can be represented as the composition of a finite number
$\lambda_{j}=m_1+m_2+\cdots+m_{j}$ of successive blow-ups which is
needed to get the standard resolution of the isolated singular point
of $V(yg_{j})$ or $V(g_{j})$, satisfying the following with the
desired properties and notations:

{\rm (i)} $\tau_{\lambda_1}=\mu_{1,\lambda_1}$ is the composition of
a finite number $\lambda_1$ of successive blow-ups which is needed
to get the standard resolution of the singular point of either
$V(yg_1)$ or $V(g_1)$.

{\rm (ii)} $\mu_{s,m_s}$ is the composition of a finite number $m_s$
of successive blow-ups which is needed to get the standard
resolution of one and only one quasisingular point of
$V^{(\lambda_{s-1})}(g_s)$ with $2\le s\le j\le r$. That is,
$\tau_{\lambda_s}=\mu_{1,m_1}\circ\mu_{2,m_2}\circ\cdots
\circ\mu_{s,m_s}$. \ms

{\rm(2)} Whenever the family
$\{(g_{j+{\ell}}\circ\tau_{\lambda_{j}})_{proper}:\ell=1,2,\dots,r-j+1\}$
with $(g_{j+{\ell}}\circ\tau_{\lambda_{j}})_{proper}\in
\BC\{1+\ve_{j,1}u_{\lambda_j},v_{\lambda_j}\}$ satisfies five
conditions in the assumptions of Sublemma $14.2$, denoted by
\text{\bf The 1-th ${\text{\bf{Cond}}}^{\text{{\bf(j)}}}$}, \dots,
\text{\bf The 5-th ${\text{\bf{Cond}}}^{\text{{\bf(j)}}}$}, then
without assuming irreducibility of
$(g_{j+{\ell}}\circ\tau_{\lambda_{j}})_{proper}\in
\BC\{1+\ve_{j,1}u_{\lambda_j},v_{\lambda_j}\}$ it was already proved
by either Sublemma $12.2$ of Theorem $12.0$ or Remark $14.2.1$ that
we get the following:

For any $\ell=1,2,\dots,r-j+1$,
$(g_{j+{\ell}}\circ\tau_{\lambda_{j}})_{total}$ with
$(g_{j+{\ell}}\circ\tau_{\lambda_{j}})_{proper}$ can be written in
the form
$$\align
(14.3.1) \quad (g_{j}\circ\tau_{\lambda_{j}})_{total}
&=v^{e_{j,\lambda_{j}}}_{\lambda_{j}}
(g_{j}\circ\tau_{\lambda_{j}})_{proper}
\quad \text{with}\\
(g_{j}\circ\tau_{\lambda_{j}})_{proper}
&=(1+\ve_{j,1}u_{\lambda_j})\\
 (g_{j+{\ell}}\circ\tau_{\lambda_{j}})_{total}
 &=v^{s^{(j)}_{\ell}s^{(j)}_{\ell-1}\cdots s^{(j)}_{2}
 s^{(j)}_{1}e_{j,\lambda_{j}}}_{\lambda_{j}}
 (g_{j+{\ell}}\circ\tau_{\lambda_{j}})_{proper}
\quad \text{with}\\
(g_{j+{\ell}}\circ\tau_{\lambda_{j}})_{proper}
&=\{(1+\ve_{j,1}u_{\lambda_j})^{s^{(j)}_1}
+v^{\gamma^{(j)}_{1,0,1}}_{\lambda_j}\}^{s^{(j)}_2s^{(j)}_3\cdots
s^{(j)}_{\ell}} +\sum_{\alpha,\beta\ge
0}B^{(j)}_{{\ell},\alpha,\beta}
v^{\alpha}_{\lambda_j}(1+\ve_{j,1}u_{\lambda_j})^{\beta},\\
\endalign$$
where a unit
$\ve_{j,1}=\ve_{j,1}(1+\ve_{j,1}u_{\lambda_j},v_{\lambda_j})$ may be
analytically assumed to be one in
$\BC\{1+\ve_{j,1}u_{\lambda_j},v_{\lambda_j}\}$, and the
$B^{(j)}_{{\ell},\alpha,\beta}$ are nonzero complex numbers for some
nonnegative integers $\alpha$ and $\beta$ such that
$s^{(j)}_1\alpha+\gamma^{(j)}_{1,0,1}\beta>s^{(j)}_{\ell}s^{(j)}_{\ell-1}\cdots
s^{(j)}_{1}\gamma^{(j)}_{1,0,1}$. \ms

{\rm(3)} Let
$\tau^{-1}_{\lambda_{j}}(0,0)=\cup^{\lambda_{j}}_{i=1}E_i$ where
each $E_i$ is an exceptional curve of the first kind. For $0\le
\ell\le r-j+1$, let $(g_{j}\circ\tau_{\lambda_{j}})_{divisor}$ and
$(g_{j+\ell}\circ\tau_{\lambda_{j}})_{divisor}$ be the divisors of
$g_{j}\circ\tau_{\lambda_{j}}$ and
$g_{j+\ell}\circ\tau_{\lambda_{j}}$ under $\tau_{\lambda_{j}}$,
respectively, which are defined by
$$\align
(g_{j}\circ\tau_{\lambda_{j}})_{divisor}&=V^{(\lambda_{j})}(g_{j})
 +\sum^{\lambda_{j}}_{i=1}e_{j,i}E_i, \tag 14.3.2 \\
 (g_{j+\ell}\circ\tau_{\lambda_{j}})_{divisor}&=V^{(\lambda_{j})}(g_{j+\ell})
 +\sum^{\lambda_{j}}_{i=1}e_{j+\ell,i}E_i,
\endalign$$
where $V^{(\lambda_{j})}(g_{j+\ell})$ is the proper transform of
$V(g_{j+\ell})$ under $\tau_{\lambda_{j}}$ and
$\lambda_{j}=m_1+m_2+\cdots +m_{j}$ as we have seen in $(14.1.4)$ of
Proposition $14.1$. \ms

Then, we have the following:

{\rm(3a)} For each $\ell=1,2,\dots,r-j+1$,

\noindent{\rm(14.3.3)} \qquad \qquad \qquad
$e_{j+\ell,i}=n_{j+\ell}n_{j+\ell}\cdots n_{j+1}e_{j,i} \quad
\text{for $1\le i\le \lambda_{j}$.}$

In particular,
$e_{{j},\lambda_{j}}=n_{j}\Delta_{j}(\beta_{j,0,k})^{j}_{k=1}$ by
{\rm(a0)} and {\rm(b0)} of $(14.2.4)$. \ms

{\rm(3b)} For any $\ell=0,1,2,\dots,r-j+1$, where $E_{\lambda_{j}}$
is defined by $\{v_{\lambda_{j}}=0\}$,

\noindent{\rm(14.3.4)} \quad \quad
$V^{(\lambda_{j})}(g_{j+\ell})\cap
(\cup^{\lambda_{j}}_{i=1}E_i)=V^{(\lambda_{j})}(g_{j})\cap
E_{\lambda_{j}}
=\{(v_{\lambda_{j}},1+{\ve}_{j,1}u_{\lambda_{j}})=(0,0)\}.$ \ms

{\rm(3c)} We get the following in the sense of {\rm Definition
$2.6$}.

If $\beta_{1,0,1}\ge 2$, then $g_k\in$ the type $[j]$ under
$\tau_{\lambda_{j}}$ for any $k=j,j+1,\dots,r+1$.

If $\beta_{1,0,1}=1$, then $g_k\in$ the type $[j-1]$ under
$\tau_{\lambda_{j}}$ for any $k=j,j+1,\dots,r+1$.

Whether $\beta_{1,0,1}\ge 2$ or $\beta_{1,0,1}=1$, note that
$yg_k\in$ the type $[j]$ under $\tau_{\lambda_{j}}$ for any $k\ge
j$. $\square$
\endproclaim \ms

{\bf \S14.2. The proof of Proposition 14.1 with Sublemma 14.2 and
Sublemma 14.3} \ms

In this section, Proposition 14.1 with Sublemma 14.2 and Sublemma
14.3 will be proved by induction on the integer $j$ where for $1\le
j\le r$, each $\tau_{\lambda_j}:M^{(\lambda_j)}\to\BC^2$ of
Proposition $14.1$ is defined to be the composition of a finite
number $\lambda_j$ of successive blow-ups which is needed to get the
standard resolution of the singular point of $V(g_j)$ or $V(yg_j)$.
So, for the proof of Proposition $14.1$, it is enough to consider
two cases, respectively:

Case(I) $j=1$, and Case(II) $2\le j\le r$. \ms

{\bf Case(I):} Let $j=1$. As far as $\tau_{\lambda_1}$ is concerned,
there is nothing to prove for Proposition $14.1$, because if $j=1$
for $\tau_{\lambda_j}$ then the proof of Proposition 14.1 with
Sublemma 14.2 and Sublemma 14.3 was already done by Sublemma $12.4$
and Sublemma $12.5$ of Theorem $12.0$. \ms

{\bf Case(II):} We are going to prove that Proposition 14.1 with
Sublemma 14.2 and Sublemma 14.3 is true for $2\le j\le r$. By the
induction proof on the positive integer $j$, suppose we have shown
by Case(I) that Sublemma $14.2$, Sublemma $14.3$ and Proposition
$14.1$ are true for $1\le j<r$. \ms

Then, it suffices to show that the following three statements are
true, respectively.

\noindent{\bf (1) Statement 14.4:} Firstly, we prove that Sublemma
$14.2$ on the integer {\rm {j+1}} is true.

\noindent{\bf (2) Statement 14.5:} Secondly, we prove that Sublemma
$14.3$ on the integer {\rm {j+1}} is true.

\noindent{\bf (3) Statement 14.6:} Thirdly, we prove that
Proposition 14.1 on the integer {\rm {j+1}} is true. \ms

\noindent$\underline{\text{\bf (1) Statement 14.4 with proof.}}$

In order to finish the proof of sublemma $14.2$ on the integer {\rm
{j+1}}, it suffices to prove the following, called Statement 14.4.
\ms

In preparation for the proof of Sublemma 14.2 on the integer (j+1),
first of all it is needed to show that the conclusion for Sublemma
$14.2$ on the integer \text{\rm (j)} implies the assumption of
Sublemma $12.5$ up to the change of notations, noting that instead
of an equality of (12.5.0) in the assumption of Sublemma $12.5$, we
need one and only one equality, being defined by
$\gcd(s^{(j)}_1,\gamma^{(j)}_{1,0,1})=1$ of (14.2.9), which was
already proved in \text{\bf The 5-th
${\text{\bf{Cond}}}^{\text{{\bf(j)}}}$} of Sublemma 14.2 on the
integer (j).

\proclaim{Statement 14.4} $\underline{\text{\bf {Assumptions}}}$
\quad Let $j$ be an arbitrary positive integer with $1\le j\le r$.
Suppose that the assumptions and notations of Proposition $14.1$
hold. Note by Theorem $12.0$ and the above assumptions of
Proposition $14.1$ that the following is true:
$$\align
\text{\rm(14.4.0)} \qquad &\text{$g_r$ is irreducible in
$\BC\{y,z\}$ with $yg_r\in \text{the type}
[r]$ under}\\
&\text{the standard resolution, but $g_{r+1}$ may not be irreducible
in $\BC\{y,z\}$.} \qquad \qquad \qquad
\endalign$$

By induction assumption on the integer $j$, suppose that the
statement on the integer $j$ in the conclusion of Sublemma $14.2$ is
true. In other words, we may assume by Sublemma $14.2$ that there
are three sequences of {\rm(14.2.2)} as we have seen, denoted by
$\text{\bf Sequences[I]}^{\text{{\bf(j)}}}$:
$$\align
(14.4.1)  \qquad &\{Y^{(j)}_\ell:\ell=1,2,\dots,r-j+1\} \quad
\text{with $Y^{(j)}_\ell\subset N_0$},  \\
&\{h^{(j)}_\ell:\ell=1,2,\dots,r-j+1\} \quad \text{with
$h^{(j)}_\ell=(g_{j+\ell}\circ\tau_{\lambda_{j}})_{proper} \in
\BC\{v_{\lambda_{j}},1+{\ve}_{j,1}u_{\lambda_{j}}\}$,} \qquad \qquad \\
&\{\text{$\Xi^{(j)}_\ell:N^\ell_0\to N_0$ is an integer-valued
function for $\ell=1,2,\dots,r-j+1$}\}, \\
&\text{satisfying five conditions in the conclusion of Sublemma $14.2.2$ :} \\
\endalign$$

\noindent Five conditions are denoted by \text{\bf The 1st
${\text{\bf{Cond}}}^{\text{{\bf(j)}}}$}{\bf, \dots,} \text{\bf The
5-th ${\text{\bf{Cond}}}^{\text{{\bf(j)}}}$ of $\text{\bf
Sequences[I]}^{\text{{\bf(j)}}}$}. \ms

$\underline{\text{\bf Conclusions}}$ \quad Then, we can prove that
the statement on the integer $(j+1)$ in the conclusion of Sublemma
$14.2$ is true.

Equivalently, for each $j$, let
$\tau_{\lambda_{j+1}}:M^{(\lambda_{j+1})}\to\BC^2$ be the
composition of a finite number $\lambda_{j+1}$ of successive
blow-ups which is needed to get the standard resolution of the
singular point of an analytic variety defined by $V(g_{j+1})$ or
$V(yg_{j+1})$.

Then, we can construct three sequences, denoted by $\text{\bf
Sequences[I]}^{\text{{\bf({j+1})}}}$:
$$\align
(14.4.2)  \qquad &\{Y^{({j+1})}_\ell:\ell=1,2,\dots,r-j\} \quad
\text{with $Y^{({j+1})}_\ell\subset N_0$},  \\
&\{h^{({j+1})}_\ell:\ell=1,2,\dots,r-j\} \quad \text{with
$h^{({j+1})}_\ell=(g_{{j+1}+\ell}\circ\tau_{\lambda_{{j+1}}})_{proper}
\in
\BC\{v_{\lambda_{{j+1}}},1+{\ve}_{{j+1},1}u_{\lambda_{{j+1}}}\}$,} \qquad \qquad \\
&\{\text{$\Xi^{({j+1})}_\ell:N^\ell_0\to N_0$ is an integer-valued
function for $\ell=1,2,\dots,r-j$}\}, \\
&\text{satisfying the following five conditions for each j:} \\
\endalign$$

\noindent Five conditions are denoted by \text{\bf The 1st
${\text{\bf{Cond}}}^{\text{{\bf({j+1})}}}$}{\bf, \dots,} \text{\bf
The 5-th ${\text{\bf{Cond}}}^{\text{{\bf({j+1})}}}$ of $\text{\bf
Sequences[I]}^{\text{{\bf({j+1})}}}$}. \ms

\noindent \text{\bf The 1st
${\text{\bf{Cond}}}^{\text{{\bf({j+1})}}}$:} Let
$\{Y^{({j+1})}_\ell:\ell=1,2,\dots,r-j\}$ with
$Y^{({j+1})}_\ell\subset N_0$ be defined by {\rm(14.2.3)}.

\roster
\item"(14.4.3)(1a)"$Y^{({j+1})}_1=\{s^{({j+1})}_1\}\cup
\{\gamma^{({j+1})}_{{1},i,1}:0\le i<s^{({j+1})}_1 \}$ with
$s^{({j+1})}_1\ge 2$ and $\gamma^{({j+1})}_{1,0,1}\ge 1$ where
$Y^{({j+1})}_1\subset N$,

\item"(1b)" $Y^{({j+1})}_{\ell}=\{s^{({j+1})}_{\ell}\}\cup
\{\gamma^{({j+1})}_{{\ell},i,1}:0\le i<s^{({j+1})}_{\ell} \} \cup
\{\gamma^{({j+1})}_{{\ell},i,2}:0\le i<s^{({j+1})}_{\ell} \}\cup
\cdots \cup \{\gamma^{({j+1})}_{{\ell},i,{\ell}}:0\le
i<s^{({j+1})}_{\ell} \}$ with $s^{({j+1})}_{\ell}\ge 2$ where
$j=1,2,\dots,r$,
\endroster

such that for each $j=1,2,\dots,r$, assume that at least one of
$\gamma^{({j+1})}_{{\ell},0,1},\gamma^{({j+1})}_{{\ell},0,2},\dots,
\gamma^{({j+1})}_{{\ell},0,{\ell}}$ is nonzero by Sublemma $12.1$
and \ms

such that each of {\rm(14.4.3)} satisfies {\rm(a)} of {\rm(14.4.4)}:

\roster
\item"(14.4.4)(a)" for each $j+1=2,3,\dots,r$,  the family
$\{Y^{({j+1})}_\ell:\ell=1,2,\dots,r-j\}$ is defined as follows:

\item"(a0)" $e_{{j+1},\lambda_{j+1}}=n_{j+1}\Delta_{j+1}(\beta_{{j+1},0,k})^{j+1}_{k=1}$,

\item"(a1)"$s^{({j+1})}_1=s^{(j)}_2=n_{j+2}\ge 2$,
$\gamma^{({j+1})}_{1,i,1}=\Xi^{(j)^{\sharp}}_2(\gamma^{(j)}_{2,i,1},\gamma^{(j)}_{2,i,2})
-(s^{(j)}_2-i)s^{(j)}_1\gamma^{(j)}_{1,0,1}>0$ for $0\le
i<s^{(j)}_2$,

\item"(a2)"
$s^{({j+1})}_{\ell}=s^{(j)}_{\ell+1}=n_{{j+1}+{\ell}}$.

\noindent $\gamma^{({j+1})}_{{\text{\rm
$\ell$}},i,1}=\Xi^{(j)^{\sharp}}_{\ell+1}(\gamma^{(j)}_{{\text{\rm
$\ell$+1}},i,k})^{\text{\rm {$\ell$}+1}}_{k=1} -(s^{(j)}_{{\text{\rm
$\ell+1$}}}-i)s^{(j)}_{{\text{\rm $\ell$}}} \cdots s^{(j)}_2
s^{(j)}_1\gamma^{(j)}_{1,0,1}>0$,
$\gamma^{({j+1})}_{{\ell},i,2}=\gamma^{(j)}_{{\ell+1},i,3}$,
$\gamma^{({j+1})}_{{\ell},i,3}=\gamma^{(j)}_{{\ell+1,}i,4}$,
$\gamma^{({j+1})}_{{\ell},i,4}=\gamma^{(j)}_{{\ell+1,}i,5}$,$\dots$,
$\gamma^{({j+1})}_{\ell,i,\ell}=\gamma^{(j)}_{\ell+1,i,\ell+1}$ for
$0\le i<s^{(j)}_{{\ell}+1}$,
\endroster \ms

\noindent \text{\bf The 2-th
${\text{\bf{Cond}}}^{\text{{\bf({j+1})}}}$:} For $\ell=1,2,\dots,
\text{\rm r-j}$, let
$h^{({j+1})}_\ell=(g_{{j+1}+\ell}\circ\tau_{\lambda_{j+1}})_{proper}$
be defined by
$$\align
(14.4.5) \qquad \qquad (g_{j+1}\circ\tau_{\lambda_{j+1}})_{total}
&=v^{e_{{j+1},\lambda_{j+1}}}_{\lambda_{j+1}}(g_{j+1}\circ\tau_{\lambda_{j+1}})_{proper}
\quad \text{with}   \qquad \qquad \qquad \\
(g_{j+1}\circ\tau_{\lambda_{j+1}})_{proper} &=1+\ve_{{j+1},1}
u_{\lambda_{j+1}} \quad
\text{with $\ve_{{j+1},1}=1$},\\
(g_{j+2}\circ\tau_{\lambda_{j+1}})_{total}
&=v^{s^{({j+1})}_1e_{{j+1},\lambda_{j+1}}}_{\lambda_{j+1}}
(g_{j+2}\circ\tau_{\lambda_{j+1}})_{proper}
\quad \text{with}\\
(g_{j+2}\circ\tau_{\lambda_{j+1}})_{proper}
&=(1+\ve_{{j+1},1}u_{\lambda_{j+1}})^{s^{({j+1})}_1}
+\sum^{s^{({j+1})}_1-1}_{i=0}\eta^{({j+1})}_{1,i}
v^{\gamma^{({j+1})}_{1,i,1}}_{\lambda_{j+1}}(1+\ve_{{j+1},1}u_{\lambda_{j+1}})^{i},\\
(g_{{j+1}+{\ell}}\circ\tau_{\lambda_{j+1}})_{total}
&=v^{s^{({j+1})}_{{\ell}}s^{({j+1})}_{{\ell}-1}\cdots
s^{({j+1})}_1e_{{j+1},\lambda_{j+1}}}
(g_{{j+1}+{\ell}}\circ\tau_{\lambda_{{j+1}}})_{proper} \quad \text{with}\\
(g_{{j+1}+{\ell}}\circ\tau_{\lambda_{j+1}})_{proper}
&=(g_{{j+1}+{\ell}-1}\circ\tau_{\lambda_{j+1}})^{s^{({j+1})}_{{\ell}}}_{proper} \\
 +\sum^{s^{({j+1})}_{{\ell}}-1}_{i=0}\eta^{({j+1})}_{{\ell},i}
v^{\gamma^{({j+1})}_{{\ell},i,1}}_{\lambda_{j+1}}&
\{\prod^{\ell}_{k=2}(g_{{j+1}+k-2}\circ
\tau_{\lambda_{j+1}})^{\gamma^{({j+1})}_{{\ell},i,k}}_{proper}\}
(g_{{j+1}+{\ell}-1}\circ\tau_{\lambda_{j+1}})^{i}_{proper}, \\
\endalign$$
where $\ve_{{j+1},1}=\ve_{{j+1},1}(v_{\lambda_{j+1}},
u_{\lambda_{j+1}})$ is a unit in $\BC\{v_{\lambda_{j+1}},
u_{\lambda_{j+1}}\}$, and also $\eta^{({j+1})}_{{\ell},i}
=\eta^{({j+1})}_{{\ell},i({\ell})}(v_{\lambda_{j+1}},1+\ve_{{j+1},1}u_{\lambda_{j+1}})$
is a unit in $\BC\{v_{\lambda_{{j+1}}},1+\ve_{{j+1},1}
u_{\lambda_{{j+1}}} \}$ for $\ell=1,2,\dots,r-{j+1}+1$ and
$i=i(\ell)=0,1,\dots,{s^{({j+1})}_{\ell}}-1$. For brevity of
notation, $\ve_{{j+1},1}$ may be defined by one, using a nonsingular
change of coordinates from
$\BC\{v_{\lambda_{j+1}},u_{\lambda_{j+1}}\}$ to itself,
independently of analytic representation of
$V^{(\lambda_{{j+1}})}(g_i)$ for all $i\ge {j+1}$.  \ms

\noindent \text{\bf The 3-th
${\text{\bf{Cond}}}^{\text{{\bf({j+1})}}}$.} Let
$\{\Xi^{({j+1})}_{\ell}: N^{\ell}_0\to N_0, {\ell}=1,2,\dots,r-j\}$
be a sequence such that each $\Xi^{({j+1})}_k$ is an integer-valued
function defined by the following: Note that $2\le{\ell}\le r-j$.
$$\align
(14.4.6) \qquad &(3a)\quad  \text{$\Xi^{({j+1})}_1(t)=t$  for each $t\in N_0$.}  \\
& \text{\rm(3b)} \quad
\text{$\Xi^{({j+1})}_{\ell}(t_k)^{\ell}_{k=1}=t_{\ell}\Xi^{({j+1})}_{{\ell-1}}
(\gamma^{({j+1})}_{{\ell-1},0,k})^{{\ell-1}}_{k=1}+s^{({j+1})}_{{\ell-1}}
\Xi^{({j+1})}_{{\ell-1}}(t_k)^{{\ell-1}}_{k=1}$}
\text{ for each $(t_k)^{\ell}_{k=1}\in  N^{\ell}_0.$}    \\
\endalign$$

Moreover, for brevity of notation, let
$\{\Xi^{({j+1})^{\sharp}}_{\ell}: N^{\ell}_0\to N_0,
{\ell}=2,3,\dots,r-j\}$ be a sequence such that each
$\Xi^{({j+1})^{\sharp}}_{\ell}$ is an integer-valued function
defined by the following:
$$\align
\text{\rm (14.4.6-1)} \quad
&\Xi^{({j+1})^{\sharp}}_{2}(\gamma^{({j+1})}_{2,i,k})^{2}_{k=1}
=\Xi^{({j+1})}_{2}(\gamma^{({j+1})}_{2,i,k})^{2}_{k=1} \quad
\text{for $0\le i<s^{({j+1})}_2$,} \\
&\Xi^{({j+1})^{\sharp}}_{\ell}(\gamma^{({j+1})}_{{\ell},i,k})^{\ell}_{k=1}
=\Xi^{({j+1})}_{2}(\gamma^{({j+1})}_{{\ell},i,k})^{2}_{k=1}
+s^{({j+1})}_1\gamma^{({j+1})}_{1,0,1}\gamma^{({j+1})}_{{\ell},i,3}
+s^{({j+1})}_1\gamma^{({j+1})}_{1,0,1}s^{({j+1})}_2\gamma^{({j+1})}_{{\ell},i,4} \\
&+s^{({j+1})}_1\gamma^{({j+1})}_{1,0,1}s^{({j+1})}_2s^{({j+1})}_3\gamma^{({j+1})}_{{\ell},i,5}
+\cdots +s^{({j+1})}_1\gamma^{({j+1})}_{1,0,1}s^{({j+1})}_2\cdots
s^{({j+1})}_{{\ell}-2}\gamma^{({j+1})}_{{\ell},i,{\ell}} \quad
\text{for
$0\le i<s^{({j+1})}_{\ell}$.} \qquad \qquad \\
\endalign$$ \ms

\noindent \text{\bf The $(4\alpha)$-th
${\text{\bf{Cond}}}^{\text{{\bf({j+1})}}}$.} Then, the following
equalities hold: Note that $r\ge 2$.

For each ${j+1}=2,3,\dots,r$, the following equality holds with
respect to the family $\{Y^{({j+1})}_\ell:\ell=1,2,\dots,r-j\}$:
$$\align
\text{\rm(14.4.7)} \quad \quad
&\Xi^{({j+1})}_1(\gamma^{({j+1})}_{1,i,1})=\gamma^{({j+1})}_{1,i,1}
=\Xi^{(j)}_2(\gamma^{(j)}_{2,i,1},\gamma^{(j)}_{2,i,2})
-(s^{(j)}_2-i)s^{(j)}_1\gamma^{(j)}_{1,0,1}>0, \\
&\Xi^{({j+1})}_{\ell}(\gamma^{({j+1})}_{{\ell},i,k})^{\ell}_{k=1}
-(s^{({j+1})}_{\ell}-i)s^{({j+1})}_{{\ell}-1}
\Xi^{({j+1})}_{{\ell}-1}(\gamma^{({j+1})}_{{\ell}-1,0,k})^{{\ell}-1}_{k=1}
\qquad \text{ for $2\le {\ell}\le {r-j}$} \qquad \\
&\quad
=\Xi^{(j)}_{{\ell}+1}(\gamma^{(j)}_{{\ell}+1,i,k})^{{\ell}+1}_{k=1}
-(s^{(j)}_{{\ell}+1}-i)s^{(j)}_{\ell}
\Xi^{(j)}_{\ell}(\gamma^{(j)}_{{\ell},0,k})^{\ell}_{k=1}>0.
\\
\endalign$$

\noindent \text{\bf The 4-th
${\text{\bf{Cond}}}^{\text{{\bf({j+1})}}}$.} The following
inequalities hold: Note that $r\ge 2$.
$$\align
(14.4.8) \quad \quad  &\text{\rm(4a)} \qquad \quad \quad
\Xi^{({j+1})}_1(\g^{({j+1})}_{1,i,1})=\g^{({j+1})}_{1,i,1}>0. \\
&\text{\rm(4b)} \quad \quad
\Xi^{({j+1})}_{\ell}(\g^{({j+1})}_{{\ell},i,k})^{\ell}_{k=1}>(s^{({j+1})}_{\ell}-i)
s^{({j+1})}_{\ell-1}\Xi^{({j+1})}_{\ell-1}(\g^{({j+1})}_{\ell-1,i,k})^{\ell-1}_{k=1}
\quad \text{for $2\le {\ell} \le r-j$.} \qquad \qquad \\
\endalign$$

\noindent \text{\bf The 5-th
${\text{\bf{Cond}}}^{\text{{\bf({j+1})}}}$.} For \text{\rm
{$\ell$}=1,2,\dots,r-j}, the following equalities hold: Note that
$r\ge 2$.
$$\align
\text{\rm(14.4.9-1)} \quad
&\text{$\gcd(s^{(1)}_{\ell},\Xi^{(1)}_{\ell}(\gamma^{(1)}_{{\ell},0,k})^{\ell}_{k=1})
=\gcd(n_{{\ell}+1},\Delta_{{\ell}+1}(\beta_{{\ell}+1,0,k})^{{\ell}+1}_{k=1})$,
\quad for $1\le{\ell}\le r$,} \\
&\text{$\gcd(s^{({j+1})}_{\ell},\Xi^{({j+1})}_{\ell}(\gamma^{({j+1})}_{{\ell},0,k})^{\ell}_{k=1})
=\gcd(s^{(j)}_{{\ell}+1},\Xi^{(j)}_{{\ell}+1}
(\gamma^{(j)}_{{\ell}+1,0,k})^{{\ell}+1}_{k=1})$ for $j\ge 1$ with
${j}+{\ell}\le r$.}
\endalign$$

\noindent \text{\bf The 5-th
${\text{\bf{Cond}}}^{\text{{\bf({j+1})}}}$.} The following
inequalities hold: Note that $r\ge 2$.
$$\align
\text{\rm(14.4.9-2)} \quad
&\f{\g^{({j+1})}_{1,i,1}}{s^{({j+1})}_1-i}>
\f{\g^{({j+1})}_{1,0,1}}{s^{({j+1})}_1}>0 \quad \text {for $0<
i<s^{({j+1})}_1$.}
\\
&\f{\Xi^{({j+1})}_{\ell}(\g^{({j+1})}_{{\ell},i,k})^{\ell}_{k=1}}{s^{({j+1})}_{\ell}-i}
> \f{\Xi^{({j+1})}_{\ell}(\g^{({j+1})}_{{\ell},0,k})^{\ell}_{k=1}}{s^{({j+1})}_{\ell}}
>s^{({j+1})}_{\ell-1}\Xi^{({j+1})}_{\ell-1}(\g^{({j+1})}_{{\ell-1},0,k})^{\ell-1}_{k=1}
\quad \text {for $0< i<s^{({j+1})}_{\ell}$.}  \qquad \qquad \qquad \qquad\\
\endalign$$

$\underline{\text{\rm Remark for \text{\rm The 5-th
${\text{\rm{Cond}}}^{\text{{\rm(j+1)}}}$}}}$ For notation, the above
inequality in $\text{\rm(14.4.9-2)}$ may be equivalently rewritten
as follows: For each ${\ell}=1,2,\dots,r-j$,
$$\align
&\f{\Xi^{({j+1})}_{\ell}(\g^{({j+1})}_{{\ell},i,k})^{\ell}_{k=1}
-(s^{({j+1})}_{{\ell}}-i)s^{({j+1})}_{{\ell}-1}\Xi^{({j+1})}_{{\ell}-1}
(\g^{({j+1})}_{{\ell}-1,0,k})^{{\ell}-1}_{k=1}}
{s^{({j+1})}_{{\ell}}-i}
\\   & > \f{\Xi^{({j+1})}_{{\ell}}(\g^{({j+1})}_{{\ell},0,k})^{{\ell}}_{k=1}
-s^{({j+1})}_{{\ell}}s^{({j+1})}_{{\ell}-1}\Xi^{({j+1})}_{{\ell}-1}
(\g^{({j+1})}_{{\ell}-1,0,k})^{{\ell}-1}_{k=1}}
{s^{({j+1})}_{{\ell}}}>0 \quad \text {for $0<
i<s^{({j+1})}_{\ell}$.}
\endalign$$
\endproclaim \ms

\noindent{\bf Proof of Statement 14.4.} For the proof, it suffices
to find the method how to apply Sublemma 12.4 and Sublemma 12.5 of
Theorem 12.0 to the proof of this statement. In order to prove that
the assumption of Sublemma of 12.5 and that of this statement are
same up to the change of notations, because the assumption of
Sublemma of 12.5 is the assumption of Theorem 12.0, which can be
rewritten by five conditions in Definition 12.0.0 and the assumption
of this statement is the same as the conclusion on Sublemma on the
integer $j$.

Therefore, as an application of Sublemma $12.5$, there is nothing to
prove for Sublemma $14.2$ on the integer {\rm {j+1}}, except for
proving the following equality:
$$
e_{j+1,\lambda_{j+1}}=n_{j+1}\Delta_{j+1}(\beta_{j+1,0,k})^{j+1}_{k=1}.
\tag 14.4.1
$$

To prove that an equality in (14.4.1) is true, we are going to apply
Sublemma $12.2$ and Corollary $3.8$ to this statement, by using the
following observations:

Let $\tau_{\lambda_{j}}:M^{(\lambda_{j})}\to\BC^2$ be the
composition of a finite number $\lambda_{j}$ of successive blow-ups
which is needed to get the standard resolution of the singular point
of $V(g_{j})$ or $V(yg_{j})$, satisfying the desired properties and
notations in Sublemma $14.2$ and Sublemma $14.3$ on the integer (j).

By applying Sublemma 12.2 of Theorem 12.0 to Sublemma $14.2$ on the
integer (j), for any $\ell=1,2,\dots,r-j+1$,
$(g_{j+{\ell}}\circ\tau_{\lambda_{j}})_{total}$ with
$(g_{j+{\ell}}\circ\tau_{\lambda_{j}})_{proper}$ can be written in
the form
$$\align
(g_{j}\circ\tau_{\lambda_{j}})_{total}
&=v^{e_{{j},\lambda_{{j}}}}_{\lambda_{{j}}}
(g_{j}\circ\tau_{\lambda_{{j}}})_{proper}
\quad \text{with}  \tag 14.4.2 \\
(g_{{j}}\circ\tau_{\lambda_{{j}}})_{proper}
&=(1+\ve_{{j},1}u_{\lambda_{j}}), \\
  (g_{j+{\ell}}\circ\tau_{\lambda_{j}})_{total}
  &=v^{s^{({j})}_{\ell}s^{({j})}_{\ell-1}\cdots
  s^{({j})}_{1}e_{{j},\lambda_{{j}}}}_{\lambda_{{j}}}
  (g_{j+{\ell}}\circ\tau_{\lambda_{{j}}})_{proper}
\quad \text{with}   \\
(g_{{j}+{\ell}}\circ\tau_{\lambda_{{j}}})_{proper}
&=\{(1+\ve_{{j},1}u_{\lambda_{j}})^{s^{({j})}_1}
+v^{\gamma^{({j})}_{1,1}}_{\lambda_{j}}\}^{s^{({j})}_2s^{({j})}_3\cdots
s^{({j})}_{\ell}} \\
&\quad +\sum_{\alpha,\beta\ge 0}B^{({j})}_{{\ell},\alpha,\beta}
v^{\alpha}_{\lambda_{j}}(1+\ve_{{j},1}u_{\lambda_{j}})^{\beta},\\
\endalign$$
where a unit
$\ve_{{j},1}=\ve_{{j},1}(u_{\lambda_{j}},v_{\lambda_{j}})$ may be
analytically assumed to be one in
$\BC\{1+\ve_{{j},1}u_{\lambda_{j}},v_{\lambda_{j}}\}$, and the
$B^{({j})}_{{\ell},\alpha,\beta}$ are nonzero complex numbers for
some nonnegative integers $\alpha$ and $\beta$ such that
$s^{({j})}_1\alpha+\gamma^{({j})}_{1,0,1}\beta>s^{({j})}_{\ell}s^{({j})}_{\ell-1}\cdots
s^{({j})}_{1}\gamma^{({j})}_{1,0,1}$. \ms

For the construction of the statement of the conclusion, let
$V(G^{(j+1)}_{0})=\{(v_{\lambda_j},1+\ve_{j,1}u_{\lambda_j})
:G^{(j+1)}_{0}=v^{\gamma}_{\lambda_j}g^{(j+1)}_{0}=0\}$ be another
analytic variety with isolated singularity at the origin in $\BC^2$
defined by the form
$$\align
g^{(j+1)}_{0}&=(1+\ve_{j,1}u_{\lambda_j})^{s^{(j)}_1} +
v^{\gamma^{(j)}_{1,0,1}}_{\lambda_j}, \tag 14.4.3
\\
G^{(j+1)}_{0} &=v^{\gamma}_{\lambda_j}g^{(j+1)}_{0},
\endalign
$$
satisfying the properties {\rm(i)} and {\rm(ii)}:

\roster \item "(i)" If $\gamma^{(j)}_{1,0,1}=1$, then $\gamma=1$.

\item "(ii)" If $\gamma^{(j)}_{1,0,1}\ge 2$, then $\gamma=0$.
\endroster \ms

Let $\mu_{j+1,m_{j+1}}$ be the composition of a finite number
$m_{j+1}$ of successive blow-ups which is needed to get the standard
resolution of the singular point of $V(G^{(j+1)}_{0})$. Then, as
compared with the above $\mu_{j+1,m_{j+1}}$, exactly the same
$\mu_{j+1,m_{j+1}}$ can be also used for the standard resolution of
one and only one quasisingular point of $V^{(\lambda_j)}(g_{j+1})$
being defined by $(g_{j+1}\circ\tau_{\lambda_j})_{proper}$ of
(14.4.2), as far as the blow-ups process is concerned. That is,
$\tau_{\lambda_{j+1}}=\mu_{1,m_1}\circ\mu_{2,m_2}\circ\cdots\circ\mu_{j+1,m_{j+1}}$
is the composition of a finite number
$\lambda_{j+1}=m_1+m_2+\cdots+m_{j+1}$ of successive blow-ups which
is needed to get the standard resolution of the isolated singular
point of $V(yg_{j+1})$ or $V(g_{j+1})$. Recall that for each
$j=1,2,\dots,r+1$, $V(g_j)=\{(y,z):g_j(y,z)=0\}$ is an analytic
variety at the origin in $\BC^2$ defined by \text{\bf The 2-th
${\text{\bf{Cond}}}^{\text{{\bf(0)}}}$,} as we have seen in the
assumption of the theorem. \ms

{\rm(b)} For simplicity of notations, let
$(v_{\lambda_{j+1}},u_{\lambda_{j+1}})$ be the common one of the
local coordinates for the $\lambda_{j+1}-th$ blow-up
$\pi_{\lambda_{j+1}}:M^{(\lambda_{j+1})}\to M^{(\lambda_{j+1}-1)}$
of $\tau_{\lambda_{j+1}}$ at a quasisingular point of
$V^{(m_{j+1}-1)}(G^{(j+1)}_{0})$. Being viewed as an analytic
mapping, $\mu_{j+1,m_{j+1}}:M^{(\lambda_{j+1})}\to
M^{(\lambda_{j})}$ can be written in the form
$$
(14.4.4)\qquad
\mu_{j+1,m_{j+1}}(v_{\lambda_{j+1}},u_{\lambda_{j+1}})
=(v_{\lambda_{j}},1+\ve_{j,1}u_{\lambda_j})=
(v^{s^{(j)}_1}_{\lambda_{j+1}}u^a_{\lambda_{j+1}},
v^{\gamma^{(j)}_{1,0,1}}_{\lambda_{j+1}}u^b_{\lambda_{j+1}}), \qquad
\qquad
$$
where {\rm(i)} $a>0$ and $b\ge 0$ are some integers such that
$a\gamma^{(j)}_{1,0,1}-bs^{(j)}_1=1$,

\noindent {\rm(ii)} $E_{\lambda_{j+1}}=\{v_{\lambda_{j+1}}=0\}$ is
defined by the $\lambda_{j+1}-th$ exceptional curve of the first
kind. \bs

Now, apply Corollary $3.8$ and
$\mu_{j+1,m_{j+1}}(v_{\lambda_{j+1}},v_{\lambda_{j+1}})$ in
(14.4.4), to the local defining equation of
$(g_{j+1}\circ\tau_{\lambda_{j}})_{total}$ of (14.4.2), and then by
(14.4.2) and (14.1.8), we have
$$\align
(14.4.5) \quad  e_{j+1,\lambda_{j+1}}
&=s^{(j)}_{1}s^{(j)}_{1}e_{j,\lambda_{j}}
+s^{(j)}_{1}\gamma^{(j)}_{1,0,1}    \\
&=n_{j+1}n_{j+1}n_{j}\Delta_{j}(\beta_{j,0,k})^{j}_{k=1}
+n_{j+1}(\Delta_{j+1}(\beta_{j+1,0,k})^{j+1}_{k=1}
-n_{j+1}n_{j}\Delta_{j}(\beta_{j,0,k})^{j}_{k=1})  \\
&=n_{j+1}\Delta_{j+1}(\beta_{j+1,0,k})^{j+1}_{k=1}. \\
\endalign$$

Thus, the proof of this statement is done. $\square$ \ms

For the proof of this sublemma, it suffices to show that this three
sequences satisfy the remaining two conditions \text{\bf The
4{$\alpha$}-th ${\text{\bf{Cond}}}^{\text{{\bf(1)}}}$} and \text{\bf
The 5{$\alpha$}-th ${\text{\bf{Cond}}}^{\text{{\bf(1)}}}$} in
$\underline {\text{\rm Conclusions}}$ of this sublemma, because
there is nothing to prove that the truth of \text{\bf The
4{$\alpha$}-th ${\text{\bf{Cond}}}^{\text{{\bf(1)}}}$} implies that
of The \text{\bf 4-th ${\text{\bf{Cond}}}^{\text{{\bf(1)}}}$}. So,
for the proof, firstly we will prove by {\bf[I]} that \text{\bf The
4{$\alpha$}-th ${\text{\bf{Cond}}}^{\text{{\bf(1)}}}$} is true, and
secondly,  by {\bf[II]} that \text{\bf The 5-th
${\text{\bf{Cond}}}^{\text{{\bf(1)}}}$} is true. \ms

{\bf [I]} For the proof of the truth of \text{\bf The 4{$\alpha$}-th
${\text{\bf{Cond}}}^{\text{{\bf(1)}}}$}, it remains to prove the
second inequality in {\rm(12.5.4${\alpha}$)} by using the following
three steps: Let $\ell$ and $q$ be an arbitrary positive integer
such that $r-1\ge {\ell}\ge q\ge 2$.
$$\align
\text{$\underline{\text{\rm Step(i)}}$} \quad  &
\Xi^{(j+1)}_{q}(\gamma^{(j+1)}_{\ell,i,k})^q_{k=1}=
\Xi^{(j)}_{q+1}(\gamma^{(j)}_{\ell+1,i,k})^{q+1}_{k=1}+(s^{(j)}_{q})^2(s^{(j)}_{q-1})^2\cdots
(s^{(j)}_{2})^2{s^{(j)}_{1}}\gamma^{(j)}_{1,0,1}\\
& \qquad \qquad \quad  \times
\{\gamma^{(j)}_{\ell+1,i,q+2}+s^{(j)}_{q+1}\gamma^{(j)}_{\ell+1,i,q+3}
+s^{(j)}_{q+1}s^{(j)}_{q+2}\gamma^{(j)}_{\ell+1,i,q+4} \\
&\qquad \qquad \quad +\cdots +s^{(j)}_{q+1}s^{(j)}_{q+2}\cdots
s^{(j)}_{\ell-1}\gamma^{(j)}_{\ell+1,i,\ell+1}-s^{(j)}_{q+1}s^{(j)}_{q+2}\cdots
s^{(j)}_{\ell} (s^{(j)}_{\ell+1}-i)\}.
\qquad \qquad \\
\text{$\underline{\text{\rm Step(ii)}}$} \quad  & \text{In
particular, if $\ell=q$
then} \\
& \Xi^{(j+1)}_{q}(\gamma^{(j+1)}_{q,i,k})^q_{k=1}
=\Xi^{(j)}_{q+1}(\gamma^{(j)}_{q+1,i,k})^{q+1}_{k=1}
-(s^{(j)}_{\ell+1}-i) (s^{(j)}_{q})^2(s^{(j)}_{q-1})^2\cdots
(s^{(j)}_{2})^2{s^{(j)}_{1}}\gamma^{(j)}_{1,0,1} \quad \text{from
{\rm Step(i)}}.
\qquad \qquad \\
\text{$\underline{\text{\rm Step(iii)}}$} \quad &
\Xi^{(j+1)}_{q}(\gamma^{(j+1)}_{q,i,k})^q_{k=1}
-(s^{(j+1)}_q-i)s^{(j+1)}_{q-1}\Xi^{(j+1)}_{q-1}(\gamma^{(j+1)}_{q-1,0,k})^{q-1}_{k=1}\\
 =& \Xi^{(j)}_{q+1}(\gamma^{(j)}_{q+1,i,k})^{q+1}_{k=1}
 -(s^{(j)}_{q+1}-i)s^{(j)}_q\Xi^{(j)}_{q}(\gamma^{(j)}_{q,0,k})^q_{k=1}>0
 \quad \text{from {\rm Step(ii)}}. \quad
\endalign$$

We will prove {\rm Step(i)}, {\rm Step(ii)} and {\rm Step(iii)} in
order, by induction on the integer $q\ge 2$.

So, it is enough to consider two cases, respectively:

Case(I) $q=2$, and Case(II) $q\ge 2$. \ms

\noindent{\bf Case(I):} \ Let $q=2$. Note by \text{\bf The 3-th
${\text{\bf{Cond}}}^{\text{{\bf(1)}}}$} that
$\Xi^{(j+1)}_2(t_1,t_2)=t_2\Xi^{(j+1)}_1(\gamma_{1,0,1})+s^{(j+1)}_1\Xi^{(j+1)}_1(t_1)
=t_2\gamma^{(j+1)}_{1,0,1}+s^{(j+1)}_1t_1$ for each $(t_1,t_2)\in
N^2_0$.
$$\align
\text{$\underline{\text{\rm Step(i)}}$} \qquad \qquad
& \Xi^{(j+1)}_2(\gamma^{(j+1)}_{{\ell},i,1},\gamma^{(j+1)}_{{\ell},i,2}) \\
 =& s^{(j+1)}_1\gamma^{(j+1)}_{{\ell},i,1}
 +\gamma^{(j+1)}_{1,0,1}\gamma^{(j+1)}_{{\ell},i,2}\\
=& s^{(j+1)}_1
\{\Xi^{(j)^{\sharp}}_{{\ell}+1}(\gamma^{(j)}_{{{\ell}+1,}k})^{{\ell}+1}_{k=1}
-(s^{(j)}_{{\ell}+2}-i)s^{(j)}_{{\ell}+1}\cdots
s^{(j)}_{2}s^{(j)}_{1}\gamma^{(j)}_{1,0,1} \} \\
& +\{\Xi^{(2)}_{2}(\gamma^{(j)}_{2,0,1},\gamma^{(j)}_{2,0,1})
-s^{(j)}_{2}s^{(j)}_{1}\gamma^{(j)}_{1,0,1}\}
\gamma^{(j)}_{\ell+1,i,3}
\qquad \qquad \qquad \quad \text{by (12.5.1)}\qquad \\
=& s^{(j)}_2
\{\Xi^{(j)}_2(\gamma^{(j)}_{{\ell+1},i,1},\gamma^{(j)}_{{\ell+1},i,2})
+s^{(j)}_{1}\gamma^{(j)}_{1,0,1}\gamma^{(j)}_{\ell+1,i,3}
+s^{(j)}_{1}\gamma^{(j)}_{1,0,1}s^{(j)}_{2}\gamma^{(j)}_{\ell+1,i,4}+\cdots  \\
&  +s^{(j)}_{1}\gamma^{(j)}_{1,0,1}s^{(j)}_{2}\cdots
s^{(j)}_{\ell-1}\gamma^{(j)}_{\ell+1,i,\ell+1}
-s^{(j)}_{1}\gamma^{(j)}_{1,0,1}s^{(j)}_2\cdots s^{(j)}_{\ell}(s^{(j)}_{\ell+1}-i) \} \\
&+
\{\Xi^{(j)}_2(\gamma^{(j)}_{2,0,1},\gamma^{(j)}_{2,0,2})-s^{(j)}_{2}s^{(j)}_{1}\gamma^{(j)}_{1,0,1}
\}
\gamma^{(j)}_{\ell+1,i,3} \qquad \qquad \qquad \quad \text{by (12.1.1)}\\
=&\{\Xi^{(j)}_3(\gamma^{(j)}_{\ell+1,i,1},\gamma^{(j)}_{\ell+1,i,2},\gamma^{(j)}_{\ell+1,i,3})
+(s^{(j)}_{2})^2(s^{(j)}_{1})\gamma^{(j)}_{1,0,1}\{\gamma^{(j)}_{\ell+1,i,4}
+s^{(j)}_3\gamma^{(j)}_{\ell+1,i,5}
+s^{(j)}_3s^{(j)}_4\gamma^{(j)}_{\ell+1,i,6} \qquad \qquad \\
& +\cdots +s^{(j)}_3s^{(j)}_4\cdots
s^{(j)}_{\ell-1}\gamma^{(j)}_{\ell+1,i,\ell+1}-s^{(j)}_3s^{(j)}_4\cdots
s^{(j)}_{\ell}(s^{(j)}_{\ell+1}-i)\},
\endalign$$

by the definition of
$\{\Xi^{(j)}_3(\gamma^{(j)}_{\ell+1,i,1},\gamma^{(j)}_{\ell+1,i,2},\gamma^{(j)}_{\ell+1,i,3})$
only, which implies the proof of {\rm Step(i)}. \ms

\noindent $\underline{\text{\rm Step(ii)}}$ \quad In particular, if
$\ell=2$ then  an equation of {\rm Step(i)} gives
$$
\Xi^{(j+1)}_2(\gamma^{(j+1)}_{2,i,1},\gamma^{(j+1)}_{2,i,2})
=\Xi^{(j)}_3(\gamma^{(j)}_{3,i,k})^3_{k=1}
-(s^{(j)}_{2})^2(s^{(j)}_{1})\gamma^{(j)}_{1,0,1}((s^{(j)}_3-i).
$$

Thus, the proof of {\rm Step(ii)} is done. \ms

\noindent $\underline{\text{\rm Step(iii)}}$ \quad  To prove that
$\Xi^{(j+1)}_2(\gamma^{(j+1)}_{2,i,1},\gamma^{(j+1)}_{2,i,2})
-(s^{(j+1)}_2-i)s^{(j+1)}_1\Xi^{(j+1)}_1(\gamma^{(j+1)}_{1,0,1})
=\Xi^{(j)}_3((\gamma^{(j)}_{3,i,k})^3_{k=1})
-(s^{(j)}_3-i)s^{(j)}_2\Xi^{(j)}_2(\gamma^{(j)}_{2,0,1},\gamma^{(j)}_{2,0,2})>0$,
first note by (12.5.1.1) that
$$\align
&(s^{(j+1)}_2-i)s^{(j+1)}_1\Xi^{(j+1)}_1(\gamma^{(j+1)}_{1,0,1}))
=(s^{(j+1)}_2-i)s^{(j+1)}_1\gamma^{(j+1)}_{1,0,1}\\
&=(s^{(j+1)}_2-i)s^{(j+1)}_1\{\Xi^{(j)}_2(\gamma^{(j)}_{2,0,1},\gamma^{(j)}_{2,0,2})
-s^{(j)}_1\gamma^{(j+1)}_{1,0,1}s^{(j)}_2\}. \\
\text{Then,} \qquad
&\Xi^{(j+1)}_2(\gamma^{(j+1)}_{2,i,1},\gamma^{(j+1)}_{2,i,2})
-(s^{(j+1)}_2-i)s^{(j+1)}_1\gamma^{(j+1)}_{1,0,1}\\
=&\Xi^{(j)}_3(\gamma^{(j)}_{3,i,k})^3_{k=1}
-(s^{(j)}_{2})^2(s^{(j)}_{1})\gamma^{(j)}_{1,0,1}((s^{(j)}_3-i)
-(s^{(j)}_3-i)s^{(j)}_2\{\Xi^{(j)}_2(\gamma^{(j)}_{2,0,1},\gamma^{(j)}_{2,0,2})
-s^{(j)}_1\gamma^{(j)}_{1,0,1}s^{(j)}_2\} \qquad \qquad \\
=&\Xi^{(j)}_3((\gamma^{(j)}_{3,i,k})^3_{k=1})
-(s^{(j)}_3-i)s^{(j)}_2\Xi^{(j)}_2(\gamma^{(j)}_{2,0,1},\gamma^{(j)}_{2,0,2})>0,
\endalign$$
by {\rm Step(ii)} and by \text{\bf The 4-th
${\text{\bf{Cond}}}^{\text{{\bf(0)}}}$} in the assumption of Theorem
$12.0$, which implies the proof of {\rm Step(iii)}.

Thus, if $q=2$, then we proved that the second inequality in
{\rm(12.5.4$\alpha$)} holds. \ms

\noindent{\bf Case(II):} \ Let $q\ge 2$. By the induction proof,
suppose we have shown that all the equalities of {\rm Step(i)}, {\rm
Step(ii)} and {\rm Step(iii)} are true on the integer $q\le r-1$
with $r-1\ge \ell\ge q$.

Then, it is enough to prove {\rm Step(i)}, {\rm Step(ii)} and {\rm
Step(iii)} in order, on the integer $(q+1)\le \ell$ as follows:

\text{$\underline{\text{\rm Step(i)}}$} Note that
$s^{(j+1)}_q=s^{(j)}_{q+1}$ and
$\gamma^{(j+1)}_{\ell,i,q+1}=\gamma^{(j)}_{\ell+1,i,q+2}$.
$$\align
\quad &\Xi^{(j+1)}_{q+1}(\gamma^{(j+1)}_{\ell,i,k})^{q+1}_{k=1} \\
=&s^{(j+1)}_q\Xi^{(j+1)}_q(\gamma^{(j+1)}_{\ell,i,k})^q_{k=1}
+\gamma^{(j+1)}_{\ell,i,q+1}\Xi^{(j+1)}_q(\gamma^{(j+1)}_{q,0,k})^q_{k=1}
\quad \text{by definition of \ $\Xi^{(j+1)}_{q+1}$} \qquad \qquad \\
= &s^{(j)}_{q+1} \{
\Xi^{(j)}_{q+1}(\gamma^{(j)}_{\ell+1,i,k})^{q+1}_{k=1}+(s^{(j)}_{q})^2(s^{(j)}_{q-1})^2\cdots
(s^{(j)}_{2})^2{s^{(j)}_{1}}\gamma^{(j)}_{1,0,1} \times
[\gamma^{(j)}_{\ell+1,i,q+2}+s^{(j)}_{q+1}\gamma^{(j)}_{\ell+1,i,q+3}\\
&+s^{(j)}_{q+1}s^{(j)}_{q+2}\gamma^{(j)}_{\ell+1,i,q+4} +\cdots
+s^{(j)}_{q+1}s^{(j)}_{q+2}\cdots
s^{(j)}_{\ell-1}\gamma^{(j)}_{\ell+1,i,\ell+1}-s^{(j)}_{q+1}s^{(j)}_{q+2}\cdots
s^{(j)}_{\ell} (s^{(j)}_{\ell+1}-i)]\}\\
&+\gamma^{(j)}_{\ell+1,i,q+2}
\{\Xi^{(j)}_{q+1}(\gamma^{(j)}_{q+1,0,k})^{q+1}_{k=1}
-(s^{(j)}_{q+1})(s^{(j)}_{q})^2(s^{(j)}_{q-1})^2\cdots
(s^{(j)}_{2})^2{s^{(j)}_{1}}\gamma^{(j)}_{1,0,1}\},\\
\endalign$$
by {\rm Step(i)}, {\rm Step(ii)} and {\rm Step(iii)} on the
induction integer $q$ because $s^{(j+1)}_q=s^{(j)}_{q+1}$ and
$\gamma^{(j+1)}_{\ell,i,q}=\gamma^{(j)}_{\ell,i,q+1}$.

Then, the definition of $\Xi^{(j)}_{q+2}
(\gamma^{(j)}_{\ell+1,i,k})^{q+2}_{k=1}$ implies the following:
$$\align
\Xi^{(j+1)}_{q+1}(\gamma^{(j+1)}_{\ell,i,k})^{q+1}_{k=1} =&
\Xi^{(j)}_{q+2} (\gamma^{(j)}_{\ell+1,i,k})^{q+2}_{k=1}
+\{(s^{(j)}_{q+1})^2(s^{(j)}_{q})^2(s^{(j)}_{q-1})^2\cdots
(s^{(j)}_{2})^2{s^{(j)}_{1}}\gamma^{(j)}_{1,0,1}\} \\
& \times
\{\gamma^{(j)}_{\ell+1,i,q+3}+s^{(j)}_{q+2}\gamma^{(j)}_{\ell+1,i,q+4}
+s^{(j)}_{q+2}s^{(j)}_{q+3}\gamma^{(j)}_{\ell+1,i,q+5} +\cdots \\
& \qquad +s^{(j)}_{q+2}s^{(j)}_{q+3}\cdots
s^{(j)}_{\ell-1}\gamma^{(j)}_{\ell+1,i,\ell+1}-s^{(j)}_{q+2}s^{(j)}_{q+3}\cdots
s^{(j)}_{\ell} (s^{(j)}_{\ell+1}-i)\},\\
\endalign$$
which implies the proof of Step(i). \ms

$\underline{\text{Step(ii)}}$ \quad In particular, if $\ell=q+1$,
then $\ell+1=q+2<q+3$ and so
$$\align
\Xi^{(j+1)}_{q+1}(\gamma^{(j+1)}_{q+1,i,k})^{q+1}_{k=1} =&
\Xi^{(j)}_{q+2} (\gamma^{(j)}_{q+2,i,k})^{q+2}_{k=1}
+\{(s^{(j)}_{q+2}-i)(s^{(j)}_{q+1})^2(s^{(j)}_{q})^2\cdots
(s^{(j)}_{2})^2{s^{(j)}_{1}}\gamma^{(j)}_{1,0,1}\} \\
\endalign$$
by Step(i) on the integer $q+1$, which implies the proof of Step(ii)
on the integer $q+1$. \ms

$\underline{\text{Step(iii)}}$ \quad To prove that the equality in
(12.5.4$\alpha$) is true, we have
$$\align \Xi^{(j+1)}_{q+1} &(\gamma^{(j+1)}_{q+1,i,k})^{q+1}_{k=1}
-(s^{(j+1)}_{q+1}-i)s^{(j+1)}_q\Xi^{(j+1)}_q(\gamma^{(j+1)}_{q,0,k})^q_{k=1} \\
&=  \Xi^{(j)}_{q+2} (\gamma^{(j)}_{q+2,i,k})^{q+2}_{k=1}
+\{(s^{(j)}_{q+2}-i)(s^{(j)}_{q+1})^2(s^{(j)}_{q})^2\cdots
(s^{(j)}_{2})^2{s^{(j)}_{1}}\gamma^{(j)}_{1,0,1}\} \\ &\quad
-((s^{(j)}_{q+2}-i)s^{(j)}_{q+1} \{\Xi^{(j)}_{q+1}
(\gamma^{(j)}_{q+1,0,k})^{q+1}_{k=1}-(s^{(j)}_{q+1})(s^{(j)}_{q})^2(s^{(j)}_{q-1})^2\cdots
(s^{(j)}_{2})^2{s^{(j)}_{1}}\gamma^{(j)}_{1,0,1}\} \\
&=\Xi^{(j)}_{q+2} (\gamma^{(j)}_{q+2,i,k})^{q+2}_{k=1}
-(s^{(j)}_{q+2}-i)(s^{(j)}_{q+1}) \Xi^{(j)}_{q+1}
(\gamma^{(j)}_{q+1,0,k})^{q+1}_{k=1},
\endalign$$
by Step(ii) on the integer $q$ and $q+1$, and by \text{\bf The 4-th
${\text{\bf{Cond}}}^{\text{{\bf(0)}}}$} of Theorem $12.0$, which
implies the proof of Step(iii).

Thus, we proved that the second inequality in $(12.5.3)$ is true,
and so we can finish the proof of the truth of \text{\bf The
4{$\alpha$}-th ${\text{\bf{Cond}}}^{\text{{\bf(1)}}}$}. \ms

\ms

\noindent$\underline{\text{\bf (2) Statement 14.5 with proof.}}$
Suppose that the assumptions of Sublemma $14.2$ are true and that
the proof of Sublemma $14.2$ is done.

In order to prove the proof of Statement 14.5, it suffices to prove
the following sublemma, called Sublemma 14.5.1. \ms

\noindent{\bf Sublemma 14.5.1.} $\underline{\text{\rm
{Assumptions}}}$ Under the same assumptions and notations of
Proposition $14.1$, suppose we have shown by Case(I) and by the
induction assumption on the positive integer $j$ that Sublemma
$14.2$, Sublemma $14.3$ and Proposition $14.1$ are true whenever $j$
is arbitrary positive integer with $1\le j<r$. Then, we showed by
proofs of Statement 14.4 and Statement 14.5 that Sublemma $14.2$ is
true. \ms

$\underline{\text{\rm {Conclusions}}}$ Since the assumptions of
Sublemma $14.2$ and the conclusions of Sublemma $14.2$ have the same
kind of conditions and notations, then by the same method as in
Sublemma 12.2 of Theorem 12.0, it is clear that the conclusions of
Sublemma $14.2$ have the same kind of representations as in the
proof of Statement 14.4 with (14.4.2), up to change of notations
only, as follows.

Whenever the family
$\{(g_{j+1+{\ell}}\circ\tau_{\lambda_{j+1}})_{proper}:\ell=1,2,\dots,r-j\}$
with $(g_{j+1+{\ell}}\circ\tau_{\lambda_{j+1}})_{proper}\in
\BC\{1+\ve_{j+1,1}u_{\lambda_{j+1}},v_{\lambda_{j+1}}\}$ satisfies
four conditions in the conclusions of Sublemma $14.2$, denoted by
\text{\bf The 1-th ${\text{\bf{Cond}}}^{\text{{\bf(j+1)}}}$}, \dots,
\text{\bf The 4-th ${\text{\bf{Cond}}}^{\text{{\bf(j+1)}}}$}, then
without assuming irreducibility of
$(g_{j+1+{\ell}}\circ\tau_{\lambda_{j+1}})_{proper}\in
\BC\{1+\ve_{j+1,1}u_{\lambda_{j+1}},v_{\lambda_{j+1}}\}$ the
conclusions of Sublemma $14.2$ have the following representations:

For any $\ell=1,2,\dots,r-j$,
$(g_{j+1+{\ell}}\circ\tau_{\lambda_{j+1}})_{total}$ with
$(g_{j+1+{\ell}}\circ\tau_{\lambda_{j+1}})_{proper}$ can be written
in the form
$$\align
(14.5.1) \quad (g_{j+1}\circ\tau_{\lambda_{j+1}})_{total}
&=v^{e_{j+1,\lambda_{j+1}}}_{\lambda_{j+1}}
(g_{j+1}\circ\tau_{\lambda_{j+1}})_{proper}
\quad \text{with}\\
(g_{j+1}\circ\tau_{\lambda_{j+1}})_{proper}
&=(1+\ve_{j+1,1}u_{\lambda_{j+1}})\\
 (g_{j+1+{\ell}}\circ\tau_{\lambda_{j+1}})_{total}
 &=v^{s^{(j+1)}_{\ell}s^{(j+1)}_{\ell-1}\cdots s^{(j+1)}_{2}
 s^{(j+1)}_{1}e_{j+1,\lambda_{j+1}}}_{\lambda_{j+1}}
 (g_{j+1+{\ell}}\circ\tau_{\lambda_{j+1}})_{proper}
\quad \text{with} \qquad \\
(g_{j+1+{\ell}}\circ\tau_{\lambda_{j+1}})_{proper}
&=\{(1+\ve_{j+1,1}u_{\lambda_{j+1}})^{s^{(j+1)}_1}
+v^{\gamma^{(j+1)}_{1,1}}_{\lambda_{j+1}}\}^{s^{(j+1)}_2s^{(j+1)}_3\cdots
s^{(j+1)}_{\ell}} \\
&\quad +\sum_{\alpha,\beta\ge 0}B^{(j+1)}_{{\ell},\alpha,\beta}
v^{\alpha}_{\lambda_{j+1}}(1+\ve_{j+1,1}u_{\lambda_{j+1}})^{\beta},\\
\endalign$$
where a unit
$\ve_{j+1,1}=\ve_{j+1,1}(1+\ve_{j+1,1}u_{\lambda_{j+1}},v_{\lambda_{j+1}})$
may be analytically assumed to be one in
$\BC\{1+\ve_{j+1,1}u_{\lambda_{j+1}},v_{\lambda_{j+1}}\}$, and the
$B^{(j+1)}_{{\ell},\alpha,\beta}$ are nonzero complex numbers for
some nonnegative integers $\alpha$ and $\beta$ such that
$s^{(j+1)}_1\alpha+\gamma^{(j+1)}_{1,1}\beta>s^{(j+1)}_{\ell}s^{(j+1)}_{\ell-1}\cdots
s^{(j+1)}_{1}\gamma^{(j+1)}_{1,1}$. $\square$

As an application of Theorem 3.8, there is nothing to prove for the
remaining for Sublemma 14.5.1.  \ms

\noindent$\underline{\text{\bf (3) Statement 14.6 with proof.}}$
\quad In order to prove Statement 14.6, it suffices to prove the
following sublemma, called Sublemma 14.6.1. \ms

\noindent{\bf Sublemma 14.6.1.}\quad $\underline{\text{\rm
{Assumptions}}}$ Under the same assumptions and notations of
Proposition $14.1$, suppose we have shown by Case(I) and by the
induction assumption on the positive integer $j$ that Sublemma
$14.2$, Sublemma $14.3$ and Proposition $14.1$ are true whenever $j$
is arbitrary positive integer with $1\le j<r$. Then, we showed by
proofs of Statement 14.4 and Statement 14.5 that Sublemma $14.2$ and
Sublemma $14.3$ on the integer $j+1$ are true.

$\underline{\text{\rm {Conclusions}}}$ To prove that Proposition
$14.1$ is true on the integer $(j+1)$, it remains to show that
either the equalities of (14.1.8) on the integer $(j+1)$ or the
following equalities are true for $1\le {\ell}\le r-j$:
$$\align
(14.6.1) \quad \text{\rm (0)}
\quad  & e_{{j+1},\lambda_{j+1}}=n_{j+1}\Delta_{j+1}(\beta_{{j+1},k})^{j+1}_{k=1}, \\
\text{\rm (1)} \quad  &s^{({j+1})}_1 =n_{j+2},  \\
&\gamma^{({j+1})}_{1,1}=\Delta_{j+2}(\beta_{j+2,k})^{j+2}_{k=1}
-n_{j+2}n_{j+1}\Delta_{j+1}(\beta_{{j+1},k})^{j+1}_{k=1}>0, \\
\text{\rm (2)} \quad  &s^{({j+1})}_2 =n_{j+3},  \\
&\gamma^{({j+1})}_{2,1}=\Delta_{j+2}(\beta_{j+3,k})^{j+2}_{k=1}
+\{\beta_{j+3,j+3}-n_{j+3}n_{j+2}\}n_{j+1}\Delta_{j+1}(\beta_{{j+1},k})^{j+1}_{k=1}>0,\\
&\gamma^{({j+1})}_{2,2}=\beta_{j+3,j+3}, \\
\text{\rm ($\ell$)} \quad &s^{({j+1})}_{\ell} =n_{j+1+{\ell}},  \\
&\gamma^{({j+1})}_{\ell,1}
=\Delta_{j+2}(\beta_{j+1+{\ell},k})^{j+2}_{k=1}
+\{\beta_{j+1+{\ell},j+3}+n_{j+2}\beta_{j+1+{\ell},j+4}+\cdots  \\
 &\cdots+ n_{j+2}n_{j+3}\cdots
 n_{j-1+{\ell}}\beta_{j+1+{\ell},j+1+{\ell}}
 -\prod^{1+{\ell}}_{k=2}n_{j+k}\}
n_{j+1}\Delta_{j+1}(\beta_{{j+1},k})^{j+1}_{k=1}>0, \qquad \qquad \qquad\\
&\gamma^{({j+1})}_{\ell,2}=\beta_{j+1+{\ell},j+3},
\gamma^{({j+1})}_{\ell,3}=\beta_{j+1+{\ell},j+4}, ~ \dots, ~
\gamma^{({j+1})}_{\ell,\ell}=\beta_{j+1+{\ell},j+1+{\ell}}. \qquad  \\
\endalign$$

\noindent{\bf Proof of Sublemma 14.6.1.} \quad Since it was shown by
Statement $14.4$ that
$e_{{j+1},\lambda_{j+1}}=n_{j+1}\Delta_{j+1}(\beta_{{j+1},k})^{j+1}_{k=1}$,
then for the proof of Sublemma 14.6.1 it suffices to show that the
following equalities are true:

(a) $s^{(j+1)}_q=n_{j+1+q}$ for each $q=1,2,\dots,r-j$. \ms

(b) $\gamma^{(j+1)}_{q,k}=\beta_{j+1+q,j+1+k}$ for any
$q=2,3,\dots,r-j$, and any $k=2,3,\dots,q$. \ms

(c) $\gamma^{(j+1)}_{q,1}=\Delta_{j+2}(\beta_{j+1+q,k})^{j+2}_{k=1}
+\{\beta_{j+1+q,j+3}+n_{j+2}\beta_{j+1+q,j+4}+\cdots $

$+ n_{j+2}n_{j+3}\cdots n_{j+q-1}\beta_{j+1+q,j+1+q}
-n_{j+2}n_{j+3}\cdots
n_{j+q}n_{j+1+q}\}n_{j+1}\Delta_{j+1}(\beta_{{j+1},k})^{j+1}_{k=1}$

for any $q=1,2,\dots,r-j$. \ms

Now, we prove (a), (b) and (c), respectively.

(a) Since it is clear by Statement 14.4 that
$s^{(j+1)}_q=s^{(j)}_{q+1}$ and by (14.1.8) that
$s^{(j)}_{q+1}=n_{j+1+q}$ for each $q=1,2,\dots,r-j$, then the proof
of (a) is done. \ms

(b) Since it is clear by Statement $14.4$ that
$\gamma^{(j+1)}_{q,k}=\gamma^{(j)}_{q+1,k+1}$ and by (14.1.8) that
$\gamma^{(j)}_{q+1,k+1}=\beta_{j+1+q,j+1+k}$ for any
$q=2,3,\dots,r-j$, and any $k=2,3,\dots,q$, then the proof of (b) is
done. \ms

(c) In preparation for the proof of (c), we use (c1), (c2), (c3) and
(c4).
$$\align
\text{\rm(c1)} \qquad \qquad \quad \gamma^{(j+1)}_{q,1}
&=\Xi^{(j)^{\sharp}}_{q+1}(\gamma^{(j)}_{{q+1},k})^{q+1}_{k=1}
-s^{(j)}_{q+1}s^{(j)}_{q}s^{(j)}_{q-1}\cdots
s^{(j)}_1\gamma^{(j)}_{1,1} \quad \text{\rm by Statement 14.4.} \\
\quad \text{\rm(c2)} \quad
\Xi^{(j)^{\sharp}}_{q+1}(\gamma^{(j)}_{{q+1,}k})^{q+1}_{k=1}
&=\Xi^{(j)}_{2}(\gamma^{(j)}_{{q+1,}k})^{2}_{k=1}
+s^{(j)}_1\gamma^{(j)}_{1,1}\{\gamma^{(j)}_{q+1,3}
+s^{(j)}_2\gamma^{(j)}_{q+1,4}\\
&  +s^{(j)}_2s^{(j)}_3\gamma^{(j)}_{q+1,5} +\cdots
+s^{(j)}_2s^{(j)}_3\cdots s^{(j)}_{{q+1}-2}\gamma^{(j)}_{q+1,q+1}\}
\quad \text{\rm by
(14.2.6-1).} \qquad \qquad \qquad \\
\quad \text{\rm(c3)}  \quad
\Xi^{(j)}_{2}(\gamma^{(j)}_{q+1,k})^{2}_{k=1}&=s^{(j)}_{1}\gamma^{(j)}_{q+1,1}
+\gamma^{(j)}_{q+1,2}\gamma^{(j)}_{1,1} \quad \text{by $(14.2.6)$}.
\endalign$$

(c4) For brevity of notation, put

$D=\gamma^{(j)}_{q+1,3} +s^{(j)}_2\gamma^{(j)}_{q+1,4}
+s^{(j)}_2s^{(j)}_3\gamma^{(j)}_{q+1,5} +\cdots
+s^{(j)}_2s^{(j)}_3\cdots
s^{(j)}_{{q+1}-2}\gamma^{(j)}_{q+1,q+1}-s^{(j)}_2s^{(j)}_3\cdots
s^{(j)}_{q+1}$.

Then, $D=\beta_{j+1+q,j+3}+n_{j+2}\beta_{j+1+q,j+4}+\cdots
+n_{j+2}n_{j+3}\cdots n_{j+q-1}\beta_{j+1+q,j+1+q}$

$-n_{j+2}n_{j+3}\cdots n_{j+q}n_{j+1+q}$ by (14.1.8). \ms

Also, applying $D$ to $\gamma^{(j)}_{q+1,1}$ of (14.1.8),
$\gamma^{(j)}_{q+1,1}$ can be rewritten as follows:
$$\align
(14.6.2) \qquad \gamma^{(j)}_{q+1,1}
&=\Delta_{j+1}(\beta_{j+q+1,k})^{j+1}_{k=1}
+\{\beta_{j+q+1,j+2}+n_{j+1}\beta_{j+q+1,j+3} \\
&\quad +n_{j+1}n_{j+2}\beta_{j+q+1,j+4}+\cdots
+ n_{j+1}n_{j+2}\cdots n_{j+q-1}\beta_{j+q+1,j+q+1} \qquad \quad  \\
& \quad -n_{j+1}n_{j+2}\cdots
n_{j+q}n_{j+q+1}\}n_j\Delta_{j}(\beta_{j,k})^{j}_{k=1} \\
&=\Delta_{j+1}(\beta_{j+q+1,k})^{j+1}_{k=1}
+\{\beta_{j+q+1,j+2}+n_{j+1}D\}n_j\Delta_{j}(\beta_{j,k})^{j}_{k=1}.
\endalign$$

Now, the truth of (c) can be just proved by the following
computation:
$$\align \gamma^{(j+1)}_{q,1} &=s^{(j)}_{1}\gamma^{(j)}_{q+1,1}
+\gamma^{(j)}_{q+1,2}\gamma^{(j)}_{1,1}+s^{(j)}_1\gamma^{(j)}_{1,1}D
\quad \text{by (c1), (c2),(c3), (c4),}\\
&=s^{(j)}_{1}\gamma^{(j)}_{q+1,1}
+(\gamma^{(j)}_{q+1,2}+s^{(j)}_1D)\gamma^{(j)}_{1,1} \\
&=n_{j+1}\{\Delta_{j+1}(\beta_{j+q+1,k})^{j+1}_{k=1}
+(\beta_{j+q+1,j+2}+n_{j+1}D)n_j\Delta_{j}(\beta_{j,k})^{j}_{k=1}\}
\\
&\quad+(\beta_{j+q+1,j+2}+n_{j+1}D)
(\Delta_{j+1}(\beta_{{j+1},k})^{j+1}_{k=1}-n_{j+1}n_j\Delta_j(\beta_{j,k})^{j}_{k=1})
\quad \text{by (14.6.2),}  \\
&=n_{j+1}\Delta_{j+1}(\beta_{j+q+1,k})^{j+1}_{k=1}
+\beta_{j+q+1,j+2}\Delta_{j+1}(\beta_{{j+1},k})^{j+1}_{k=1}
+Dn_{j+1}\Delta_{j+1}(\beta_{{j+1},k})^{j+1}_{k=1} \\
&= \Delta_{j+2}(\beta_{j+q+1,k})^{j+2}_{k=1}
+Dn_{j+1}\Delta_{j+1}(\beta_{{j+1},k})^{j+1}_{k=1}. \\
\endalign$$

Thus, we proved that Sublemma $14.6.1$ is true, and so the proof of
Statement 14.6 is finished. Therefore, the proof for Case(II) is
done, and then we finished the proof of Proposition $14.1$,
completely. $\square$ \ms

\demo{\bf Proof of Theorem 14.0 or Proposition 14.1} In preparation
for the proof of this theorem, it suffices to show that Proposition
$14.1$ is true , whenever $j$ is arbitrary with $1\le j\le r$. First
of all, as far as $\tau_{\lambda_1}$ is concerned, there is nothing
to prove for Proposition $14.1$, because if $j=1$ for
$\tau_{\lambda_j}$ then the proof of Proposition 14.1 was already
done by Sublemma $12.4$ and Sublemma $12.5$ of Theorem $12.0$. Also,
for $2\le j\le r$, the proof of Proposition $14.1$ is done by
Case(I) and Case(II). Thus, the proof of theorem is completely
finished. $\square$
\enddemo \ms

\newpage

{\bf Part[C2] The division algorithm for the W-polys} \ms

{\bf $\S$15. The division algorithm for W-polys} \bs

\noindent{\bf Notation 15.0.0.} Instead of Weierstrass polynomials,
the Weierstrass preparation theorem, and the Weierstrass division
theorem, we write $W$-polys, the WPT and the WDT respectively, for
brevity of notation. Recall the well-known theorems. \ms

\definition{Definition 15.0} Let $\C\{z_1,z_2,\dots,z_n\}$ or ${}_n\CO_0$
be the ring of convergent power series at the origin in $\C^n$ and
$\C\{z_1,z_2,\dots,z_m\}[z_{m+1},\dots,z_n]$ or
${}_m\CO[z_{m+1},\dots,z_n]$ be a polynomial ring in
$z_{m+1},z_{m+2},\dots,z_n$ with coefficients from the ring
${}_m\CO_0$.

(i) If $f\in {}_n\CO_0$ and $f(0,\dots,0,z_n)$ is not identically
zero as a function of $z_n$ in a neighborhood of $0 \in \C$, then
$f$ is regular at $0 \in \C^n$. Then, it is said that $f$ is regular
of order $\nu$ in $z_n$ at $0 \in \C^n$ if there is the least
integer $\nu$ such that ${\partial}^{\nu}f/{\partial z}^{\nu}_n$ is
nonzero at the origin.

(ii) If $f\in {}_n\CO_0$ is not identically zero, then the
multiplicity of $f$ at the origin in $\C^n$ is defined by the least
integer $\nu$ such that some partial derivative of $f$ of order
$\nu$ is nonzero at the origin.

(iii) If $h\in {}_{n-1}\CO[z_n]$ is a polynomial of degree $\nu$ in
$z_n$, then $h$ has the form $h=a_0z^\nu_n +a_1z^{\nu-1}_n+\cdots
+a_n$ where for each $i$, $a_i\in {}_{n-1}\CO_0$ and $a_0$ is a unit
in ${}_{n-1}\CO_0$. Also, if $h\in {}_{n-1}\CO[z_n]$ is a W-poly of
degree $\nu >0$ in $z_n$, then $h$ has the form
$h=z^\nu_n+a_1z^{\nu-1}+\cdots +a_n$ where for each $i$, $a_i$ is a
nonunit in ${}_{n-1}\CO_0$.

(iv) If $h\in {}_{n-1}\CO[z_n]$ is a Weierstrass polynomial of
degree $\nu >0$ in $z_n$, which has the multiplicity $\mu>0$ at the
origin in $\C^n$, then $h$ is said to be a Weierstrass polynomial of
degree $\nu >0$ in $z_n$ with the multiplicity $\mu>0$ at $0\in
\BC^n$.
\enddefinition \ms

\proclaim{Theorem 15.1 (The WPT and The WDT)}

{\bf (I)(The WPT):} If $f\in {}_n\CO_0$ is regular of order $\nu$ in
$z_n$, then there is a unique $W$-poly $h\in {}_{n-1}\CO [z_n]$ of
degree $\nu$ in $z_n$ such that $f=uh$ for some unit $u\in
{}_n\CO_0$. \ms

{\bf(II)(The WDT):} If $h\in {}_{n-1}\CO [z_n]$ is a $W$-poly of
degree $\nu$ in $z_n$, then any $f\in {}_n\CO_0$ can be written
uniquely in the form $f=gh+r$ where $g\in {}_n\CO_0$ and $r\in
{}_{n-1}\CO[z_n]$ is a polynomial of degree $<\nu$ in $z_n$. Also,
if $f\in {}_{n-1}\CO[z_n]$, then $g\in {}_{n-1}\CO[z_n]$.
\endproclaim

\ms \proclaim{Theorem 15.2 (The WDT for the W-polys)}

{\rm(1)} Let $h\in {}_{n-1}\CO[z_n]$ be a $W$-poly of degree $\nu
>0$ in $z_n$.

{\rm(1a)} Let $f\in {}_{n-1}\CO [z_n]$ be a $W$-poly of degree $\mu
\ge \nu$ in $z_n$. Then, $f$ can be written uniquely in the form
$$f=gh+r, \tag 15.2.1.1
$$
where if $\mu >\nu$ then $g\in {}_{n-1}\CO[z_n]$ is a $W$-poly of
degree $\mu -\nu>0$ in $z_n$ and if $\mu =\nu$ then $g$ is equal to
one, and if $\mu\ge\nu$ then $r\in {}_{n-1}\CO[z_n]$ is a polynomial
of degree $<\nu$ in $z_n$ with $r(0,\dots,0,z_n)$ identically zero.
\ms

{\rm(1b)} Let $f\in {}_{n-1}\CO[z_n]$ be a $W$-poly of degree
$\mu\ge \nu$ in $z_n$, and $\ell$ be a positive integer with
$\ell\nu \le \mu <(\ell+1)\nu$. Then, $f$ can be written uniquely in
the form
$$
f=\sum^{\ell}_{i=0} r_ih^i \quad \text{with \quad $h^{0}=1$}, \tag
15.2.1.2
$$
where if $\mu \ge \ell\nu$ then for each $i=0,1,\dots, \ell-1$,
$r_i\in {}_{n-1}\CO[z_n]$ is a polynomial of degree $<\nu$ in $z_n$
with $r_i(0,\dots, 0, z_n)$ identically zero, and if $\mu =\ell\nu$
then $r_{\ell}$ is equal to one and if $\mu
>\ell\nu$ then $r_{\ell}\in {}_{n-1}\CO[z_n]$ is a $W$-poly of degree
$\mu -\ell\nu <\nu$ in $z_n$. \ms

{\rm(1c)} Let $f\in {}_{n-1}\CO[z_n]$ be a polynomial of degree $\mu
\ge \nu$ in $z_n$ with $f(0,\dots, 0,z_n)$ identically zero. Then,
$f$ can be written uniquely in the form
$$f=gh+r \tag 15.2.2.1 $$
such that if $\mu >\nu$ then $g\in {}_{n-1}\CO[z_n]$ is a polynomial
of degree $\mu -\nu >0$ in $z_n$ and if $\mu =\nu$ then $g\in
{}_{n-1}\CO_0$ is a nonunit, and such that if $\mu \ge \nu$ then
$r\in {}_{n-1}\CO[z_n]$ is a polynomial of degree $<\nu$ in $z_n$
where $g(0,\dots, 0,z_n)$ and $r(0,\dots, 0,z_n)$ are identically
zero. \ms

{\rm(1d)} Let $f\in {}_{n-1}\CO[z_n]$ be a polynomial of degree $\mu
\ge \nu$ in $z_n$ with $f(0,\dots,0,z_n)$ identically zero, and
$\ell$ be a positive integer with $\ell\nu \le \mu <(\ell+1)\nu$.
Then, $f$ can be written uniquely in the form
$$
f=\sum^{\ell}_{i=0} r_ih^i \quad \text{with \quad $h^{0}=1$}, \tag
15.2.2.2
$$
where if $\mu \ge \ell\nu$ then for each $i=0,1,\dots,\ell-1$,
$r_i\in {}_{n-1}\CO[z_n]$ is a polynomial of degree $<\nu$ in $z_n$
with $r_i(0,\dots,0,z_n)$ identically zero, and if $\mu =\ell\nu$
then $r_{\ell}\in {}_{n-1}\CO_0$ is a nonunit and if $\mu
>\ell\nu$ then $r_{\ell}\in {}_{n-1}\CO[z_n]$ is a polynomial
of degree $\mu -\ell\nu <\nu$ in $z_n$ with $r_{\ell}(0,\dots,
0,z_n)$ identically zero. \ms

{\rm(2)} Let $h\in {}_{n-1}\CO[z_n]$ be a $W$-poly of degree $\nu
>0$ in $z_n$ with the multiplicity $\nu >0$ at $0\in \BC^n$.

{\rm(2a)} Let $f\in {}_{n-1}\CO[z_n]$ be a $W$-poly of degree $\mu
\ge \nu$ in $z_n$ with the multiplicity $\mu \ge \nu$ at $0\in
\BC^n$. Then, the above representation $f=\sum^{\ell}_{i=0}r_ih^i$
of {\rm (1b)} satisfies the property that for each
$i=0,1,\dots,{\ell}-1$, $r_i$ has a multiplicity $\ge \mu -i\nu$ at
$0\in \BC^n$ and that if $\mu = \ell\nu$ then $r_\ell$ is equal to
one and if $\mu
>\ell\nu$ then $r_\ell\in {}_{n-1}\CO[z_n]$ is a $W$-poly of degree
$\mu-\ell\nu$ in $z_n$ with the multiplicity $\mu -\ell\nu$ at $0\in
\BC^n$. \ms

{\rm(2b)} Let $f\in {}_{n-1}\CO[z_n]$ be a polynomial of degree $\mu
>\nu$ in $z_n$ with $f(0,\dots, 0,z_n)$ identically zero. If $f$ has
a multiplicity $m\ge \mu >\nu$ at $0 \in \C^n$, then the
representation $f=\sum^{\ell}_{i=0}r_ih^i$ of {\rm(1d)} satisfies
the property that for each $i=0,1,\dots,\ell$, $r_i$ has a
multiplicity $\ge m-i\nu$ at $0 \in \C^n$.
\endproclaim \ms

\demo{\bf Proof of Theorem 15.2} Recall by Theorem 15.1 that if
$h\in {}_{n-1}\CO[z_n]$ is a $W$-poly of degree $\nu >0$ in $z_n$
and $f\in {}_{n-1}\CO[z_n]$, then any $f$ can be written uniquely in
the form
$$f=gh +r \tag 15.2.3
$$
where $g\in {}_{n-1}\CO[z_n]$, and $r\in {}_{n-1}\CO[z_n]$ is a
polynomial of degree $<\nu$ in $z_n$. \ms

(1)(1a) Observe that $h(0,\dots, 0,z_n)=z^\nu_n$ and $f(0,\dots,
0,z_n)=z^\mu_n$ because $h$ and $f$ are $W$-polys in $z_n$. So,
(15.2.3) implies the following:
$$
z^\mu_n=z^\nu_ng(0,\dots,0,z_n)+r(0,\dots, 0,z_n). \tag 15.2.4
$$
Since $r\in {}_{n-1}\CO[z_n]$ is a polynomial of degree $<\nu$ in
$z_n$, then $r(0,\dots,0,z_n)$ is identically zero, and so
$g(0,\dots,0,z_n)=z_n^{\mu-\nu}$. If $\mu >\nu$ then $g\in
{}_{n-1}\CO[z_n]$ is a $W$-poly of degree $\mu -\nu$ in $z_n$
because $f$ and $h$ are $W$-polys in $z_n$, and if $\mu =\nu$ then
$g$ is equal to one. \ms

(1b) By (15.2.3) and (1a), rewrite $f=g_1h+r_0$ where $g_1=g$ and
$r_0=r$. Assume that $\mu -\nu\ge \nu$, otherwise it was already
done by (1a). Now, apply the WDT with a divisor $h$ to $g_1$, and
then $g_1$ can be written uniquely in the form $g_1=g_2h+r_1$ where
$r_1\in {}_{n-1}\CO[z_n]$ is a polynomial of degree $<\nu$ in $z_n$
with $r_1(0,\dots, 0,z_n)$ identically zero and $g_2\in
{}_{n-1}\CO[z_n]$ is either a $W$-poly of degree $\mu -2\nu
>0$ in $z_n$, or equal to one if $\mu =2\nu$. Again, assume that $\mu -2\nu
\ge \nu$, otherwise it can be finished by the same method as we have
done for $\mu -\nu< \nu$. Then apply the WDT with a divisor $h$ to
$g_2$ and so on, we can get the desired result. The uniqueness for
the above representation follows immediately from the WDT. \ms

(1c) Consider $f=gh+r$. By the assumption, $0=g(0,\dots,
0,z_n)z^\nu_n+r(0,\dots, 0,z_n)$. Since $r\in {}_{n-1}\CO[z_n]$ is a
polynomial of degree $<\nu$ in $z_n$, then $r(0,\dots,
0,z_n)=g(0,\dots, 0,z_n)=0$. So, if $\mu >\nu$ then $g\in
{}_{n-1}\CO[z_n]$ is a polynomial of degree $\mu -\nu >0$ in $z_n$
and if $\mu =\nu$ then $g\in {}_{n-1}\CO_0$ is a nonunit because
$f(0,\dots, 0,z_n)$ is zero. \ms

(1d) It just follows by the similar technique as we have used in the
proof of (1b) together with the result of (1c). \ms

(2)(2a) By the WDT and (1a), $f$ can be written uniquely with
$f=gh+r$ as follows:
$$\align
(15.2.5) {\qquad} f &=z^\mu_n+\sum^\mu_{i=1} a_iz^{\mu -i}_n,
{\quad} h=z^\nu_n +\sum^\nu_{j=1} b_jz^{\nu-j}_n {\quad} \text{and}
{\quad} g=z^{\mu-\nu}_n +\sum^{\mu-\nu}_{k=1} c_kz^{\mu-\nu-k}_n,
\qquad
\endalign$$
where $f$, $h$ and $g$ are $W$-polys in $z_n$.

Since $f$ and $h$ are $W$-polys of degree $\mu \ge \nu$ in $z_n$
with the multiplicity $\mu \ge \nu$ at $0\in \BC^n$, respectively,
then each $a_i$ is a nonunit in ${}_{n-1}\CO_0$ with multiplicity
$\ge i$ if exists, and also each $b_j$ is a nonunit in
${}_{n-1}\CO_0$ with multiplicity $\ge j$, if exists. If $\mu =\nu$
then $g$ is equal to one by (1a) and so it is clear by (15.2.5) that
$r$ has a multiplicity $\ge \mu$ at $0 \in \C^n$, because $f=gh+r$.

We claim that if $\mu>\nu$, each coefficient $c_k\in {}_{n-1}\CO_0$
of (15.2.5) has a multiplicity $\ge k$ if exists, and so $g$ is a
$W$-poly of degree $\mu -\nu$ in $z_n$ with the multiplicity $\mu
-\nu$ at $0 \in \C^n$.

For the proof of the claim, assume the contrary. Then there exists a
nonzero coefficient $c_k\in {}_{n-1}\CO_0$ with the multiplicity
$<k$, and so let $p$ be the smallest among the positive integers $k$
such that $c_k\in {}_{n-1}\CO_0$ has a multiplicity $<k$.

Now, rewrite $g$ in the following form
$$ g=\sum_1+\sum_2+\sum_3 , \tag 15.2.6 $$
where $\sum_1 =z^{\mu-\nu}_n+c_1z^{\mu-\nu-1}_n+\cdots
+c_{p-1}z^{\mu-\nu-(p-1)}_n$, $\sum_2=c_pz^{\mu-\nu-p}_n$ and
$\sum_3 = c_{p+1}z^{\mu-\nu-p-1}_n+c_{p+2}z^{\mu-\nu-p-2}_n+\cdots
+c_{\mu-\nu}$.

Consider $hg$ as $hg=h\Sigma_1+h\Sigma_2+h\Sigma_3$. Then $gh$
satisfies the following properties (2a-1), (2a-2) and (2a-3):

(2a-1) $h\Sigma_1$ has a multiplicity $\mu$ at at $0 \in \C^n$
because for $1\le k\le p-1$, $c_k$ has a multiplicity $\ge k$ at $0
\in \C^{n-1}$ by definition of $p$.

(2a-2) In $h\Sigma_2$, there exists a nonzero term $c_pz^{\nu}_n
z^{\mu-\nu-p}_n=c_pz^{\mu-p}_n$ with the multiplicity $<\mu$ because
$c_p$ has a multiplicity $<p$, and note that $h\Sigma_2\in
{}_{n-1}\CO[z_n]$ is a polynomial of degree $\mu-p\ge \nu$ in $z_n$.

(2a-3) For any nonzero term $dz^q_n\in h\Sigma_3$ with a nonunit
$d\in {}_{n-1}\CO_0$, $dz^q_n\in {}_{n-1}\CO[z_n]$ is a polynomial
of degree $q\le \nu +\mu-\nu-p-1<\mu-p$ in $z_n$. \ms

Therefore, $c_pz^{\mu-p}_n$ is a nonzero term in $f=gh+r$ with
multiplicity $<\mu$ by (2a-1), (2a-2) and (2a-3) because $r\in
{}_{n-1}\CO[z_n]$ is a polynomial of degree $<\nu$ and $\nu\le
\mu-p$. It would be impossible because $f$ has a multiplicity $\mu$
at $0 \in \C^n$. Thus we proved the claim, and so $r$ has a
multiplicity $\ge \mu$ at $0 \in \C^n$ because $hg$ has a
multiplicity $\mu$ at $0 \in \C^n$. \ms

Now, assume that $\mu-\nu \ge \nu$, otherwise it was done just
before. Since $g$ is a $W$-poly of degree $\mu-\nu$ in $z_n$ with
the multiplicity $\mu-\nu$ at $0 \in \C^n$, apply the WDT with a
divisor $h$ to $g$. Assuming by (15.2.3) that $f=gh+r$ as before,
then $g$ can be written uniquely in the form $g=g_2h+r_1$,
satisfying that either $g_2$ is a $W$-poly of degree $\mu -2\nu
>0$ in $z_n$ with the multiplicity $\mu -2\nu >0$ at $0 \in \C^n$ or
$g_2$ is equal to one by the similar technique as we have used just
before. Since $g_2h$ has a multiplicity $\mu-\nu$ at $0 \in \C^n$,
then $r_1$ has a multiplicity $\ge \mu-\nu$ at $0 \in \C^n$. Repeat
the above process by (1b), and then it can be done. \ms

(2b) It can be proved by the similar technique as in the proof of
(2a) together with the result of (1d).

Thus, we finished the proof of Theorem 15.2. $\square$
\enddemo \ms

\proclaim{Lemma 15.3} Let $f=z^n+\sum^{n-1}_{i=0} c_iy^{\beta_i}z^i$
be a $W$-poly of degree $n\ge 2$ in $z$ with the multiplicity $n\ge
2$ at $0 \in \C^n$ where for $0\le i\le n-1$, each $c_i=c_i(y)$ is a
unit in ${}_2\CO_0$ if exists and the $\beta_i$ are positive
integers. Assume that $f$ may not be irreducible in ${}_2\CO_0$. If
$c_{n-1}$ is not identically zero, then $f$ can be rewritten by the
\text{\rm WDT} only, without using a nonsingular change of
coordinates, as follows:
$$
f=(z+\f{c_{n-1}}ny^{\beta_{n-1}})^n+\sum^{n-2}_{i=0}
a_iy^{\a_i}(z+\f{c_{n-1}}ny^{\beta_{n-1}})^i, \tag 15.3.1
$$
where $a_{n-1}$ is identically zero and for $0\le i\le n-2$, each
$a_i=a_i(y)$ is a unit in ${}_2\CO_0$ if exists and the $\a_i$ are
positive integers.
\endproclaim

\demo {\bf Proof of Lemma 15.3} Apply the WDT to $f$ with a divisor
$z+\f{c_{n-1}}ny^{\beta_{n-1}}$. Then, it is trivial by (1b) and
(2a) of Theorem 15.2 that if $h=z+\f{c_{n-1}}ny^{\beta_{n-1}}$ and
$r_i=a_iy^{\a_i}$ then $\alpha_i\ge n-i$. $\square$
\enddemo \ms

Throughout this section, replace $z+\f{c_{n-1}}n y^{\beta_{n-1}}$ by
$z$ for brevity of notation, as we have done in Lemma 15.3. \ms

\vfill \pagebreak

{\bf \S15.1. The division algorithm for the W-polys in preparation
for the computations of The 2nd Algorithm and The 3rd Algorithm} \ms

As in Definition 15.0 with Notation 15.0.1, the Weierstrass
preparation theorem and the Weierstrass division theorem can be
written by The WPT and The WDT respectively, for brevity of
notation. In order to succeed in the computations of the 2nd and the
3rd algorithms in $\S1.6$, in this section it is very important to
say without any other proofs that The Division Algorithm for the
W-polys(Theorem 1.8(Theorem 15.4) with two sublemmas) can have an
important role of the 2nd and the 3rd algorithms in $\S1.6$, which
will be shown in $\S 16$, later.

\proclaim{Theorem 15.4 (The Division Algorithm for The W-polys)}

$\underline{\text{\bf Assumptions}}$ Let $f=z^n+\sum^{n-2}_{i=0}
a_iy^{\a_i}z^i$ be a $W$-poly of degree $n\ge 2$ in $z$ where for
$0\le i\le n-2$, each $a_i=a_i(y)$ is a unit in ${}_2\CO_0$ if
exists and the $\a_i$ are positive integers. Assume that $f$ may not
be irreducible in ${}_2\CO_0$, and note that $a_{n-1}$ is
identically zero for convenience. Write $n=\Pi^{\ell}_{k=1}n_k$ with
positive integers $n_k\ge 2$ for all $k$ where the $n_k$ may not be
the factorization of prime numbers. \ms

$\underline{\text{\bf Conclusions}}$ For each fixed $j$ with $1\le j\le
\ell-1$, there exists a unique sequence of $W$-polys in $z$,
$\{f_k: k=1,2,\dots,j\}$ with $f_k\in \C \{y\}[z]$, satisfying
the following two properties, called {\bf(1)} and {\bf(2)}, and notations:
Let $f_{-1}=y$, $f_0=z$, and write 
$$\align
f_k&=f^{n_k}_{k-1}+\sum^{n_k-2}_{i=0} R_{k,i}f^i_{k-1}
\quad \text{for each \  $k=1,2,\dots,j$,} \tag 15.4.1 \\
f&=f^{d_{j+1}}_j +\sum^{d_{j+1}-2}_{i=0} S_{j+1,i}f^i_j,
\endalign$$
such that
\roster
\item"(i)"
$n=d_{j+1}\Pi^j_{k=1}n_k$ if $\ell\ge 2$,
and $n=d_1$ if $\ell=1$,

\item"(ii)" for each $k\ge 1$,
$f_k=f_k(y,z)\in \C\{y\}[z]$ is a $W$-poly of degree
$\Pi^k_{t=1}n_t$ in $z$,

\item"(iii)" $f_k\in
\C\{f_{-1},f_0,\dots,f_{k-2}\}[f_{k-1}]$ is a $W$-poly of degree $n_k$ in
$f_{k-1}$ with coefficients from the ring $\C\{y,z,f_1,\dots,f_{k-2}\}$,

\item"(iv)"  $f\in \C\{f_{-1},f_0,\dots,f_{j-1}\}[f_j]$ is
a $W$-poly of degree $d_{j+1}$ in $f_j$ with coefficients from the ring
$\C\{y,z,f_1,\dots,f_{j-1}\}$, \endroster

{\noindent}considering $f_{-1}=y, f_0=z, f_1, \dots, f_j$ as independent complex
$(j+2)$-variables at $0\in \C^{j+2}$ if necessary, with two
properties {\rm(1)} and {\rm(2)}. \ms

\roster
\item"(1)(1a)" Let $k$ and $i$ be fixed with $1\le k\le
j$ and $0\le i\le n_k-2$. If exists, then $R_{1,i}=R_{1,i}(y)$ is a
nonunit in $\C\{y\}$ and for each $k\ge 2$, $R_{k,i}=R_{k,i}(y,z)\in
\C\{y\}[z]$ is a polynomial of degree $<\Pi^{k-1}_{t=1}n_t$ in $z$
with $R_{k,i}(0,z)=0$.

\item"(1)(1b)" Let $i$ be fixed with $0\le i\le d_{j+1}-2$. For each $j\ge
1$, $S_{j+1,i}=S_{j+1,i}(y,z)\in \C\{y\}[z]$ is a polynomial of
degree $<\Pi^j_{t=1}n_t$ in $z$ with $S_{j+1,i}(0,z)=0$. \ms

\item"(2)(2a)" Let $k$ and $i$ be fixed with $1\le
k\le j$ and $0\le i\le n_k-2$. For any nonzero monomial
$\Pi^k_{t=1}f^{\de_t}_{t-2}$ in $R_{k,i}=R_{k,i}(y,z,f_1,\dots,
f_{k-2})\in \C\{y\}[z,f_1,\dots, f_{k-2}]$, $\de_1>0$ and
$\de_t<n_{t-1}$ for $t=2,3,\dots, k$.

\item"(2b)" Let $i$ be fixed with $0\le i\le d_{j+1}-2$. For
any nonzero monomial $\Pi^{j+1}_{t=1}f^{\g_t}_{t-2}$ in
$S_{j+1,i}=S_{j+1,i}(y,z,f_1,\dots, f_{j-1})\in \C\{y\}[z,f_1,\dots,
f_{j-1}]$, $\g_1>0$ and $\g_t<n_{t-1}$ for $t=2,3,\dots,j+1$.
\endroster
\endproclaim \ms

\proclaim{Corollary 15.4.1}
$\underline{\text{\bf Assumptions}}$
Under the the same assumptions and notations as in Theorem 15.4,  in addition
we may assume that $f$ is a $W$-poly of degree $n\ge 2$ in $z$ with the multiplicity
$n$ at $0 \in \C^2$. \ms
$\underline{\text{\bf Conclusions}}$
Using the same properties and notations as in Theorem 15.4,
we can prove the following additional results:
\roster
\item"(i)" For each $k$, $f_k(y,z)$ is a $W$-poly of degree
$\Pi^k_{t=1}n_t$ in $z$ with the multiplicity $\Pi^k_{t=1}n_t$ at $0
\in \C^2$.

\item"(ii)" As in {\rm (2a)} of Theorem 15.4, $R_{k,i}\in \C\{y\}[z]$ of
{\rm(15.4.1)} has a multiplicity $\ge (n_k-i)\Pi^{k-1}_{t=1}n_t$ at $0
\in \C^2$.

\item"(iii)" As in $(2a)$ of Theorem 15.4, $S_{j+1,i}\in \C\{y\}[z]$
of $(15.4.1)$ has a multiplicity $\ge (d_{j+1}-i)\Pi^j_{t=1}n_t$ at $0
\in \C^2$.
\endroster
\endproclaim \bs

{\bf \S15.2. How to  prove Theorem 15.4 and Corollary 15.4.1 } \ms

In this section, for the proofs of Theorem $15.4$ and Corollary $15.4.1$
we are going to try how to construct the
statements of two sublemmas without proofs, consisting of Sublemma
15.5 and Sublemma 15.6. After then, in $\S 15.3$ $\S 15.4$, $\S 15.5$ and $\S 15.6$, 
we will give the proofs of two sublemmas, which finishes the proofs of
Theorem $15.4$ and Corollary $15.4.1$ completely. \ms

\noindent{\bf Proofs of Theorem 15.4 and Corollary 15.4.1.} First of all,
in this section it is trivial to show that Theorem $15.4$ and Corollary
$15.4.1$ with proofs are equivalent to finding two sublemmas with proofs, called Sublemma 15.5
(Existence of an algorithm) and Sublemma 15.6(Uniqueness of an algorithm),
which are defined by the following statements, under the same assumptions and notations as
it has been seen in Theorem 15.4 and Corollary 15.4.1: \ms

$\underline{\text{\rm How to find two sublemmas, Sublemma 15.5 and Sublemma 15.6 }}$ \
In $\S 15.3$, $\S 15.4$, $\S 15.5$ and $\S 15.6$, we are going to prove
Theorem $15.4$ and Corollary $15.4.1$ at the same time by induction on
the integer $j\ge 1$, as follows. 

{\bf (1)} In the process of the induction proof, first of all, if $j=1$ with 
$\ell=2$ in (15.4.1), it will be proved later that Theorem $15.4$ and Corollary $15.4.1$ are true.

That is, $f_{1}$ and $f$ can be uniquely constructed as follows:
$$
\cases f_{1} &=f^{n_{1}}_0 +\sum^{n_{1}-2}_{i=0}R_{1,i}f^i_0, \\
f &=f^{d_{2}}_{1}+\sum^{d_{2}-2}_{i=0} S_{2,i}f^i_{1},
\endcases  \tag 15.4.2
$$
where $n=d_{2}n_1$, and $f_{1}$ with
$R_{1,i}$ and $f$ with $S_{2,i}$ satisfy the same kind of
properties as $f_{1}$ with $R_{1,i}$ and $f$ with $S_{2,i}$ have
done in (1) and (2) of the conclusion of the theorem, respectively,
and also the same kind of properties in Corollary $15.4.1$ as before. \ms

{\bf(2)} Secondly, if $j\ge 1$ with $1\le j< \ell$, suppose we have shown by the induction assumption
on the integer $j\ge 1$ that
the theorem is true on the integer $1\le j< \ell$. Then on the positive
integer $(j+1)$, it suffices to show that given such a sequence
$\{f_1,f_2,\dots, f_j\}$ as we have seen in (15.4.1), then $f_{j+1}$
and $f$ can be uniquely constructed as follows:
$$
\cases f_{j+1} &=f^{n_{j+1}}_j +\sum^{n_{j+1}-2}_{i=0}R_{j+1,i}f^i_j, \\
f &=f^{d_{j+2}}_{j+1}+\sum^{d_{j+2}-2}_{i=0} S_{j+2,i}f^i_{j+1},
\endcases  \tag 15.4.2
$$
where $d_{j+2}=\Pi^{\ell}_{t=j+2}n_t$, and $f_{j+1}$ with
$R_{j+1,i}$ and $f$ with $S_{j+2,i}$ satisfy the same kind of
properties as $f_{j}$ with $R_{j,i}$ and $f$ with $S_{j+1,i}$ have
done in (1) and (2) of the conclusion of the theorem, respectively,
and also the same kind of properties in Corollary $15.4.1$ as before. \ms

For the proofs of Theorem 15.4 and Corollary $15.4.1$ simultaneously, using the Weierstrass
division theorem it suffices to show by induction
proof that Sublemma 15.5 and Sublemma 15.6 are true respectively, as follows:

\noindent{\bf Sublemma 15.5(Existence of an algorithm).} We prove that for any $j\ge 0$
such a sequence
$\{f_1,f_2,\dots,f_{j+1}\}$ for $f$ can be constructed in the sense
of (15.4.2), using {\bf The 1st half of Sublemma 15.5(Sublemma 15.5.{$\alpha$})} and {\bf
The 2nd half of Sublemma 15.5(Sublemma 15.5.{$\beta$})} for brevity of notation. \ms

\noindent{\bf Sublemma 15.6(Uniqueness of an algorithm).} We prove that for any $j\ge 0$
such a sequence $\{f_1,f_2,\dots,f_{j+1}\}$ for $f$ in the sense of (15.4.2)
constructed in Sublemma 15.5 must be unique. \ms

\noindent{\bf Remark 15.4.2.} Note that the assumptions in both Sublemma 15.5 and Sublemma 15.6 are the same as the assumptions in Theorem 15.4 and Corollary 15.4.1. In more detail, in preparation for finding such a sequence in (15.4.1) in Theorem 15.4 and Corollary 15.4.1, the proof is as follows:

1. To find such a sequence in (15.4.1) in Theorem 15.4 and Corollary 15.4.1, first of all,
it suffices to show by $\S15.3$ that we can find two small sublemmas of Sublemma 15.5,
called Sublemma 15.5.{$\alpha$}(the 1st half of Sublemma 15.5) and Sublemma 15.5.{$\beta$}
(the 2nd half of Sublemma 15.5), below.
Secondly, it is very interesting to show by $\S15.4$ and by $\S15.5$
that we can prove such two small sublemmas of $\S15.4$, which can give an algorithm
for finding such a sequence $\{f_1,f_2,\dots,f_{j+1}\}$ for $f$ in the sense of (15.4.2).

2. To find a uniqueness of such a sequence in (15.4.1) in Theorem 15.4 and Corollary 15.4.1,
it suffices to show that we can prove Sublemma 15.6. \bs

Therefore,
in order to prove Theorem $15.4$ and Corollary $15.4.1$,
it suffices to solve the following in order:
\roster
\item"(i)"  In $\S15.3$, the problem is how to write the 1st half of
Sublemma 15.5(Sublemma 15.5.{$\alpha$})
and the 2nd half of Sublemma 15.5(Sublemma 15.5.{$\beta$}) of Sublemma 15.5.

\item"(ii)"  In $\S15.4$, the problem is how to write two statements of Sublemma 15.5 for induction proof.

\item"(iii)"  In $\S15.5$, the problem is how to find the proofs of the 1st statement
and the 2nd statement of Sublemma 15.5 by the induction method.

\item"(iv)"  In $\S15.6$, the problem is how to find the proof of Sublemma 15.6
by the induction method.
\endroster \ms

{\bf \S15.3. How to write the 1st half of Sublemma 15.5(Sublemma 15.5.{$\alpha$})
and the 2nd half of Sublemma 15.5(Sublemma 15.5.{$\beta$}) of Sublemma 15.5} \ms

\noindent$\underline{\text{\bf Sublemma 15.5.{$\alpha$} of Sublemma 15.5}}$ \quad

$\underline{\text{\bf Assumptions}}$ Under the same induction
assumption on the integer $j\le \ell-1$ with $\ell \ge 2$, suppose
we have shown that Sublemma 15.5 is true on the integer $j\le \ell-1$,
just as above. \ms

$\underline{\text{\bf Conclusions}}$ We show that $h_{j,1}$ and
$f$ can be constructed as follows:
$$
\cases
h_{j,1} &=f^{n_{j}}_{j-1} +\sum^{n_{j}-2}_{i=0} R^{(1)}_{j,i}f^i_{j-1}, \\
f &=h^{d_{j+1}}_{j,1}+\sum^{d_{j+1}-1}_{i=0} T_{j+1,i}h^i_{j,1},
\endcases \tag 15.4.3
$$
where $h_{j,1}\in \C\{y\}[z]$ is a $W$-poly of degree
$\Pi^{j}_{t=1}n_t$ in $z$ with the multiplicity
$\Pi^{j}_{t=1}n_t$ at $0 \in \C^2$, satisfying the following
facts, Fact(A), Fact(B), Fact(C), Fact(D) and Fact(E).

$\underline{\text{\rm Fact(A)}}$ For each $i=0,1,\dots,n_{j}-2$,
$R^{(1)}_{j,i}=R_{j,i}^{(1)}(y,z)\in \C\{y\}[z]$ is a polynomial
of degree $<\Pi^{j-1}_{t=1}n_t$ in $z$  with $R^{(1)}_{j,i}(0,z)=0$
and has a multiplicity $\ge (n_{j}-i)\Pi^{j-1}_{t=1}n_t$ at $0 \in
\C^2$. \ms

$\underline{\text{\rm Fact(B)}}$ For each $i=0,1,\dots,n_{j}-2$,
and for any nonzero monomial $\Pi^{j}_{t=1}f^{\de_t}_{t-2}$ in
$R^{(1)}_{j,i}\in \C\{y\}[z,f_1,\dots, f_{j-2}]$, $\de_1>0$ and
$\de_t <n_{t-1}$ for $2\le t\le j$. \ms

$\underline{\text{\rm Fact(C)}}$ For each $i=0,1,\dots, d_{j+1}-1$,
$T_{j+1,i}=T_{j+1,i}(y,z)=\sum a_{p,q}y^pz^q$ with a nonzero
constant $a_{p,q}$ such that $p>0$ and $q<\Pi^{j}_{t=1}n_t$ and
that $T_{j+1,i}(0,z)=0$ and $T_{j+1,i}$ has a multiplicity $\ge
(d_{j+1}-i)\Pi^{j}_{t=1}n_t$ at $0 \in \C^2$. \ms

Moreover, consider $y,z,f_1,\dots,f_{j-1}$ as independent complex
$(j+1)$-variables at the origin in $\C^{j+1}$. Then, $T_{j+1,i}\in
\C\{y\}[z,f_1,\dots, f_{j-1}]$ satisfies two facts {\rm Fact(D) and
Fact(E)}.

$\underline{\text{\rm Fact(D)}}$  For each $i=0,1,\dots, d_{j+1}-1$,
and for any nonzero monomial $\Pi^{j+1}_{t=1}f^{\g_t}_{t-2}$ in
$T_{j+1,i}$, $\g_1>0$ and $\g_t<n_{t-1}$ for $2 \le t\le j+1$. \ms

$\underline{\text{\rm Fact(E)}}$ In particular, if $i=d_{j+1}-1$ for
$T_{j+1,i}$ of Fact(D), then $\g_{j+1}\le n_{j}-2$. $\square$ \ms

\noindent$\underline{\text{\bf Sublemma 15.5.{$\beta$ of Sublemma 15.5}}}$ \quad

$\underline{\text{\bf Assumptions}}$ By the
same way as we have done in {\rm (15.4.3)} of
the conclusion of Sublemma $15.5.\alpha$, we assume that $(h_{j,1},f)$
can be rewritten as follows:
$$
\cases
h_{j,1} &=f^{n_{j}}_{j-1} +\sum^{n_{j}-2}_{i=0} R^{(1)}_{j,i}f^i_{j-1}, \\
f &=h^{d_{j+1}}_{j,1}+\sum^{d_{j+1}-1}_{i=0} T_{j+1,i}h^i_{j,1},
\endcases \tag 15.4.3
$$
satisfying the facts, denoted by {\rm Fact(A), Fact(B), Fact(C),
Fact(D) and Fact(E)}. For brevity of notation, let $h_1=h_{j,1}$,
$R^{(1)}_i =R^{(1)}_{j,i}$ for $0\le i\le n_{j}-2$ and
$T^{(1)}_i=T^{(1)}_{j+1,i}=T_{j+1,i}$ for $0\le i\le d_{j+1}-1$,
respectively. \ms

$\underline{\text{\bf Conclusions}}$ Then, $(f_{j},f)$ for $f$ can
be constructed as follows:

{\bf Case [I]:} If $T^{(1)}_{d_{j+1}-1}$ in $(h_1,f)$ is zero, let
$f_{j}=h_1$, $R_{j,i}=R^{(1)}_i$ for $0\le i\le n_{j}-2$ and
$S_{j+1,i}=T^{(1)}_i$ for $0\le i\le d_{j+1}-2$, respectively. Then,
the construction of $(f_{j},f)$ has been already finished. \ms

{\bf Case [II]:} If $T^{(1)}_{d_{j+1}-1}$ is not zero, for finding
such a construction of $(f_{j},f)$, it suffices to follow two
steps, {\rm Step(1)} and {\rm Step(2)}.

$\underline{\text{\bf Step(1) for Case[II]}}$ Then, there is a
sequence of pairs, $H=\{(h_p,f):p=1,2,\dots\}$, each pair of which
can be constructed with five properties, called {\rm Property(1)},
{\rm Property(2)}, {\rm Property(3)}, {\rm Property(4)} and {\rm
Property(5)}, as follows:
$$
\cases h_{1} &=f^{n_{j}}_{j-1} +\sum^{n_{j}-2}_{i=0}R^{(1)}_if^i_{j-1}
\text{\quad with $R^{(1)}_i=R^{(1)}_{j,i}$ in $(15.4.3)$,}\\
f &=h^{d_{j+1}}_{1}+\sum^{d_{j+1}-1}_{i=0} T^{(1)}_ih^i_{1}
\text{\quad with $T^{(1)}_i=T^{(1)}_{j+1,i}$ in $(15.4.3)$,}
\endcases \tag 15.4.4)(15.4.4.1
$$
$$\cases h_{2} &=h_1 +\f 1{d_{j+1}}T^{(1)}_{d_{j+1}-1}=f^{n_{j}}_{j-1}
+\sum^{n_{j}-2}_{i=0}R^{(2)}_if^i_{j-1}, \qquad \qquad \quad\\
f &=h^{d_{j+1}}_{2}+\sum^{d_{j+1}-1}_{i=0} T^{(2)}_ih^i_{2},
\endcases \tag 15.4.4.2
$$
$$\cases h_{3} &=h_2 +\f 1{d_{j+1}}T^{(2)}_{d_{j+1}-1}=f^{n_{j}}_{j-1}
+\sum^{n_{j}-2}_{i=0}R^{(3)}_if^i_{j-1}, \qquad \qquad \quad\\
f &=h^{d_{j+1}}_{3}+\sum^{d_{j+1}-1}_{i=0} T^{(3)}_ih^i_{3},
\endcases \tag 15.4.4.3
$$
$$
 \dots, \qquad\qquad \qquad \qquad\qquad \qquad \qquad\qquad
 \qquad\qquad
$$
satisfying the following properties and notations: \ms

$\underline{\text{\rm Property(1)}}$ Let $p$ and $i$ be fixed with
$p\ge 1$ and $0\le i\le n_{j}-2$. Then
$R^{(p+1)}_i=R^{(p+1)}_i(y,z)\in \C\{y\}[z]$ is a polynomial of
degree $<\Pi^{j-1}_{t=1}n_t$ in $z$ and has a multiplicity $\ge
(n_{j}-i)\Pi^{j-1}_{t=1}n_t$ at $0 \in \C^2$.  \ms

$\underline{\text{\rm Property(2)}}$ Let $p$ and $i$ be fixed with
$p\ge 1$ and $0\le i\le d_{j+1}-1$. Then
$T^{(p+1)}_i=T^{(p+1)}_i(y,z)\in \C\{y\}[z]$ is a polynomial of
degree $<\Pi^{j}_{t=1}n_t$ in $z$ and has a multiplicity $\ge
(d_{j+1}-i)\Pi^{j}_{t=1}n_t$ at $0 \in \C^2$. \ms

Consider $y,z,f_1,\dots,f_{j-1}$ as independent complex
$(j+1)$-variables at the origin in $\C^{j+1}$.

$\underline{\text{\rm Property(3)}}$ Let $p$ and $i$ be fixed with
$p\ge 1$ and $0\le i\le n_{j}-2$. Then for any nonzero monomial
$\Pi^{j}_{t=1}f^{\de_t}_{t-2}$ in
$R^{(p+1)}_i=R^{(p+1)}_i(y,z,f_1,\dots,f_{j-2})\in
\C\{y\}[z,f_1,\dots,f_{j-2}]$, $\de_1>0$ and $\de_t <n_{t-1}$ for
$t=2,3,\dots, j$.  \ms

$\underline{\text{\rm Property(4)}}$ Let $p$ and $i$ be fixed with
$p\ge 1$ and $0\le i\le d_{j+1}-1$. Then for any nonzero monomial
$\Pi^{j+1}_{t=1}f^{\de_t}_{t-2}$ in
$T^{(p+1)}_i=T_i^{(p+1)}(y,z,f_1,\dots,f_{j-1})\in
\C\{y\}[z,f_1,\dots,f_{j-1}]$, $\de_1>0$ and $\de_t<n_{t-1}$ for
$t=2,3,\dots,j+1$. \ms

$\underline{\text{\rm Property(5)}}$ In particular, if $i=d_{j+1}-1$
for $T^{(p+1)}_i$ of {\rm Property(4)}, then $\de_{j+1}\le
n_{j}-2$. \bs

$\underline{\text{\bf Step(2) for Case[II]}}$ By {\rm Step(1)},
there is a pair $(h_{\nu+1},f)\in H$ which satisfies the following
property:

$\underline{\text{\rm Property(6)}}$ There is an integer $\nu\le
\f{n_{j}+1}2$ such that $T^{(p)}_{d_{j+1}-1}\ne 0$ for
$p=1,2,\dots,\nu$ and
$T^{(\nu+1)}_{d_{j+1}-1}=T^{(\nu+2)}_{d_{j+1}-1}=\cdots =0$. That
is, $(h_{\nu},f)\not= (f_{j},f)$ and $(h_{\nu+1},f)= (f_{j},f)$
for an integer $\nu \le \f{n_{j}+1}2$. \ms

{\bf Remark 15.4.3.} {\rm(a)} It is clear by Sublemma $15.4.\alpha$
that {\rm Property(1)}, {\rm Property(2)}, {\rm Property(3)}, {\rm
Property(4)} and {\rm Property(5)} for
\text{\rm($h_1$,f)=($h_{j+1,1}$,f)} in {\rm(15.5.2.1)} are equivalent
to {\rm Fact(A)}, {\rm Fact(C)}, {\rm Fact(B)}, {\rm Fact(D)} and
{\rm Fact(E)} for \text{\rm($h_{j+1,1}$,f)} in {\rm(15.4.3)},
respectively.

{\rm(b)} By {\rm(a)}, it is clear by Sublemma $15.4.\alpha$ that
$(h_1,f)$ was already constructed with five properties. $\square$
 \bs

{\bf \S15.4. How to write the statements of Sublemma 15.5 for induction proof} \ms

We are going to prove Sublemma 15.5 by induction on
the positive integer $j$. In case $j=0$ in (15.4.1), there
is nothing to prove.
In case $1\le j\le \ell-1$,  to find such an existence of a sequence
$\{f_1,f_2,\dots,f_j\}$ in the sense of Sublemma $15.5$, it is clear by induction
assumption on the positive integer $j$  that 
it suffices to consider the statement for Case(I)(called the 1st statement of Sublemma 15.5) and 
the statement for Case(II)(called the 2nd statement of Sublemma 15.5), as follows:

Case(I) $j=1\le {\ell-1}$ and  Case(II) $1\le j\le {\ell-1}$. \ms

\noindent$\underline{\text{\rm Case(I):}}$ Let $j=1\le {\ell-1}$. In $\S15.4.1$, the problem 
is how to write the 1st statement of Sublemma 15.5 for induction proof. 
\ms

\noindent$\underline{\text{\rm Case(II):}}$ Let $1\le j\le {\ell-1}$. In $\S15.4.2$, the problem  
is how to write the 2nd statement of Sublemma 15.5 for induction proof. 
\ms

\noindent{\bf Remark 15.4.0} To finish the proof of Sublemma 15.5,  it suffices to show by $\S15.5.1$ that
we can complete the proof of the 1st statement of Sublemma 15.5 in $\S15.4.1$
and to show by $\S15.5.2$ that applying directly the same properties and
notations as it can be used in the proof of the 1st statement of Sublemma 15.5
to the 2nd statement of Sublemma 15.5 in $\S15.4.2$ inductively, we can complete 
the proof of the 2nd statement of Sublemma 15.5 in $\S15.5.2$ 
very easily. \bs 

{\bf \S15.4.1. How to write the 1st statement of Sublemma 15.5 for induction proof} \ms
$\underline{\text{\bf The 1st statement for Sublemma 15.5}}$

$\underline{\text{\rm Assumptions for The 1st statement}}$
Suppose that this statement satisfies the same assumption as
we have seen in Theorem 15.4 on the integer $j=1\le {\ell-1}$.

$\underline{\text{\rm Conclusions for The 1st statement}}$ We can compute a unique
$W$-poly in $z$, $f_1\in \C \{y\}[z]$, satisfying the following
notations and properties: Let $f_{-1}=y$ and $f_0=z$.
$$\align
f_1&=f^{n_1}_{0}+\sum^{n_1-2}_{i=0} R_{1,i}f^i_{0}
\quad \text{and} \tag 15.4.5 \\
f&=f^{d_{2}}_1 +\sum^{d_{2}-2}_{i=0} S_{2,i}f^i_1,
\endalign$$
such that \quad{\rm(i)} \quad $n=d_{2}n_1$ with $n=d_1$,

\qquad \quad \quad {\rm(ii)} \quad $f_1=f_1(y,z)\in \C\{y\}[z]$ is a
$W$-poly of degree $n_1$ in $z$,

\qquad \quad \quad {\rm(iii)} \quad $f\in \C\{y\}[z,f_1]\subseteq
\C\{y,z\}[f_1]$ is a $W$-poly of degree $d_{2}$ in $f_1$, \ms

{\noindent}considering $f_{-1}=y,f_0=z,f_1$ as independent complex
$(3)$-variables at the origin in $\C^{3}$, with two properties {\rm
(1)} and {\rm (2)}: \ms

\roster
\item"(1)(1a)" Let $i$ be fixed with $0\le i\le n_1-2$. If exists, then
$R_{1,i}=R_{1,i}(y)$ is a nonunit in $\C\{y\}$.

\item"(1b)" Let $i$ be fixed with $0\le i\le d_{2}-2$. Then
$S_{2,i}=S_{2,i}(y,z)\in \C\{y\}[z]$ is a polynomial of degree
$<n_1$ in z and $S_{2,i}(0,z)=0$. \ms

\item"(2)(2a)" Let $i$ be fixed with $0\le i\le n_1-2$.
For any nonzero monomial $y^{\de_1}$ in $R_{1,i}=R_{1,i}(y)\in
\C\{y\}$, $\de_1>0$.

\item"(2b)" Let $i$ be fixed with $0\le i\le d_{2}-2$. For
any nonzero monomial $y^{\de_1}z^{\delta_2}$ in
$S_{2,i}=S_{2,i}(y,z)\in \C\{y\}[z]$, $\delta_1>0$ and
$\delta_2<n_{1}$. \quad $\square$
\endroster \bs

Under the the same assumptions and notations as in Corollary 15.4.1,  in addition
we may assume that $f$ is a $W$-poly of degree $n\ge 2$ in $z$ with the multiplicity
$n$ at $0 \in \C^2$. \ms
As a conclusion, using the same properties and notations as in Theorem 15.4,
we can prove the following additional results:
\roster
\item"(i)" $f_1(y,z)$ is a $W$-poly of degree
$n_1$ in $z$ with the multiplicity $n_1$ at $0
\in \C^2$.

\item"(ii)" As in {\rm (2a)} of the 1st statement of Sublemma 15.5, $R_{1,i}\in \C\{y\}[z]$ 
of {\rm(15.4.5)} has a multiplicity $\ge (n_1-i)$ at $0
\in \C^2$.

\item"(iii)" As in $(2b)$ the 1st statement of Sublemma 15.5, $S_{j+1,i}\in \C\{y\}[z]$
of $(15.4.5)$ has a multiplicity $\ge (d_{2}-i)n_1$ at $0
\in \C^2$.
\endroster \ms

For the proof the 1st statement of Sublemma 15.5,
it suffices to prove {\bf Sublemma 15.5.$\alpha$
of the 1st statement and Sublemma 15.5.$\beta$
of the 1st statement} respectively, as follows:

\proclaim{Sublemma 15.5.$\alpha$ of Sublemma 15.5 of the 1st statement}

$\underline{\text{\bf Assumptions}}$ Suppose that the same
properties and notations as in the assumption of Theorem 15.4 hold.
\ms

$\underline{\text{\bf Conclusions}}$ We show that $h_{1,1}$ and $f$
can be constructed as follows:
$$
\cases
h_{1,1} &=f^{n_{1}}_0 +\sum^{n_{1}-2}_{i=0} R^{(1)}_{1,i}f^i_0, \\
f &=h^{d_{2}}_{1,1}+\sum^{d_{2}-1}_{i=0} T_{2,i}h^i_{1,1},
\endcases \tag 15.5.1
$$
where $h_{1,1}\in \C\{y\}[z]$ is a $W$-poly of degree $n_1$ in z and
$n=n_1d_2$, satisfying the following facts, {\rm Fact(A), Fact(B),
Fact(C), Fact(D) and Fact(E)}.

\roster

\item"$\underline{\text{\rm Fact(A)}}$" For each
$i=0,1,\dots,n_{1}-2$, $R^{(1)}_{1,i}=R_{1,i}^{(1)}(y)\in \C\{y\}$
with $R^{(1)}_{1,i}(0)=0$, if exists, and
 has a multiplicity $\ge (n_1-i)$ at $0 \in \C^2$. \ms

\item"$\underline{\text{\rm Fact(B)}}$" For each
$i=0,1,\dots,n_{1}-2$, and for any nonzero monomial $y^{\de_1}$ in
$R^{(1)}_{1,i}\in \C\{y\}$, $\de_1>0$.

\item"$\underline{\text{\rm Fact(C)}}$" For each $i=0,1,\dots,
d_{2}-1$, $T_{2,i}=T_{2,i}(y,z)=\sum a_{p,q}y^pz^q$ with a nonzero
constant $a_{p,q}$ such that $p>0$ and $q<n_1$ and that
$T_{2,i}(0,z)=0$ and $T_{2,i}$ has a multiplicity $\ge
(d_2-i)n_1$ at $0 \in \C^2$.

Moreover, considering $y,z,f_1=h_{1,1}$ as independent complex
$(3)$-variables at the origin in $\C^{3}$, then, $T_{2,i}\in
\C\{y\}[z]\subseteq \C\{y,z\}$ satisfies two facts {\rm Fact(D) and
Fact(E)}.

\item"$\underline{\text{\rm Fact(D)}}$"  For each $i=0,1,\dots, d_{2}-1$,
and for any nonzero monomial $\Pi^{2}_{t=1}f^{\g_t}_{t-2}$ in
$T_{2,i}$, $\g_1>0$ and $\g_2<n_1$.

\item"$\underline{\text{\rm Fact(E)}}$" In particular, if $i=d_{2}-1$ for
$T_{2,i}$ of {\rm Fact(D)}, then $\g_{2}\le n_{1}-2$. $\square$
\endroster
\endproclaim \ms

\proclaim{Sublemma 15.5.$\beta$ of Sublemma 15.5 of the 1st statement}

$\underline{\text{\bf Assumptions}}$ Suppose that the same
properties and notations as in the assumption of Theorem 15.4 hold.
By the same way as we have seen in {\rm (15.5.1)} of the conclusion
of Sublemma 15.5.$\alpha$, we may assume that $(h_{1,1},f)$ can be
rewitten as follows:
$$
\cases
h_{1,1} &=f^{n_{1}}_0 +\sum^{n_{1}-2}_{i=0} R^{(1)}_{1,i}f^i_0, \\
f &=h^{d_{2}}_{1,1}+\sum^{d_{2}-1}_{i=0} T_{2,i}h^i_{1,1},
\endcases \tag 15.5.1
$$
satisfying the facts, denoted by {\rm Fact(A), Fact(B), Fact(C),
Fact(D) and Fact(E)}. For brevity of notation, let $h_1=h_{1,1}$,
$R^{(1)}_i =R^{(1)}_{1,i}$ for $0\le i\le n_{1}-2$ and
$T^{(1)}_i=T^{(1)}_{2,i}=T_{2,i}$ for $0\le i\le d_{2}-1$,
respectively. \ms

$\underline{\text{\bf Conclusions}}$ Then, $(f_{1},f)$ for $f$ in
the conclusion of Theorem 1.8 can be constructed as follows:

{\bf Case[I]:} If $T^{(1)}_{d_{2}-1}$ in $(h_1,f)$ is zero, let
$f_{1}=h_1$, $R_{1,i}=R^{(1)}_i$ for $0\le i\le n_{1}-2$ and
$S_{2,i}=T^{(1)}_i$ for $0\le i\le d_{2}-2$, respectively. Then, the
construction of $(f_{1},f)$ has been already finished. \ms

{\bf Case[II]:} If $T^{(1)}_{d_{2}-1}$ is not zero, for finding such
a construction of $(f_{1},f)$, it suffices to follow two steps, {\rm
Step(1)} and {\rm Step(2)}.

$\underline{\text{\bf Step(1) for Case[II]}}$ Then, there is a
sequence of pairs, $H=\{(h_p,f):p=1,2,\dots\}$, each pair of which
can be constructed with five properties, called {\rm Property(1)},
{\rm Property(2)}, {\rm Property(3)}, {\rm Property(4)} and {\rm
Property(5)}, as follows:
$$
(15.5.2)(15.5.2.1)\quad \quad  \cases h_{1} &=f^{n_{1}}_0
+\sum^{n_{1}-2}_{i=0}R^{(1)}_if^i_0
\text{\quad with $R^{(1)}_i=R^{(1)}_{1,i}$ in $(15.5.1)$,}\\
f &=h^{d_{2}}_{1}+\sum^{d_{2}-1}_{i=0} T^{(1)}_ih^i_{1} \text{\quad
with $T^{(1)}_i=T^{(1)}_{2,i}$ in $(15.5.1)$,} \qquad \qquad
\endcases $$
$$\cases h_{2} &=h_1 +\f 1{d_{2}}T^{(1)}_{d_{2}-1}=f^{n_{1}}_0
+\sum^{n_{1}-2}_{i=0}R^{(2)}_if^i_0, \qquad \qquad \quad\\
f &=h^{d_{2}}_{2}+\sum^{d_{2}-1}_{i=0} T^{(2)}_ih^i_{2},
\endcases \tag 15.15.2.2
$$
$$\cases h_{3} &=h_2 +\f 1{d_{2}}T^{(2)}_{d_{2}-1}=f^{n_{1}}_0
+\sum^{n_{1}-2}_{i=0}R^{(3)}_if^i_0, \qquad \qquad \quad\\
f &=h^{d_{2}}_{3}+\sum^{d_{2}-1}_{i=0} T^{(3)}_ih^i_{3},
\endcases \tag 15.15.2.3
$$
$$ {\dots,}\qquad \qquad \qquad \qquad \qquad \qquad
\qquad \qquad \quad  $$

satisfying the following properties and notations: \ms

$\underline{\text{\rm Property(1)}}$ Let $p$ be fixed with $p\ge 1$.
For each $i=0,1,\dots, d_{2}-2$, $R^{(p+1)}_i=R^{(p+1)}_i(y)\in
\C\{y\}$ with $R^{(p+1)}_{i}(0)=0$, if exists, and has a multiplicity $\ge
(n_1-i)$ at $0 \in \C^2$. \ms

$\underline{\text{\rm Property(2)}}$ Let $p$ and $i$ be fixed with
$p\ge 1$ and $0\le i\le {d_2}-1$. Then
$T^{(p+1)}_i=T^{(p+1)}_i(y,z)\in \C\{y\}[z]$ is a polynomial of
degree $< n_1$ in $z$ and has a multiplicity $\ge
({d_2}-i)n_1$ at $0 \in \C^2$. \ms

Consider $f_{-1}=y,f_0=z$ as independent complex $(2)$-variables at
the origin in $\C^{2}$.

$\underline{\text{\rm Property(3)}}$ Let $p$ and $i$ be fixed with
$p\ge 1$ and $0\le i\le n_{1}-2$. For any nonzero monomial
$f^{\de_1}_{-1}$ in $R^{(p+1)}_{i}=R^{(p+1)}_{i}(y)\in \C\{y\}$,
$\de_1>0$. \ms

$\underline{\text{\rm Property(4)}}$ Let $p$ and $i$ be fixed with
$p\ge 1$ and $0\le i\le d_{2}-1$. For any nonzero monomial
$f^{\de_1}_{-1}f^{\de_2}_{0}$ in
$T^{(p+1)}_i=T_i^{(p+1)}(f_{-1},f_{0})\in \C\{f_{-1},f_0\}$,
$\de_1>0$ and $\de_2<n_{1}$.

$\underline{\text{\rm Property(5)}}$ In particular, if $i=d_{2}-1$
for $T^{(p+1)}_i$ of {\rm Property(4)}, then $\de_{2}\le n_{1}-2$.
\ms

$\underline{\text{\bf Step(2) for Case[II]}}$ By {\rm Step(1)} for
{\rm Case[II]}, there is a pair $(h_{\nu+1},f)\in H$ which satisfies
the following property:

$\underline{\text{\rm Property(6)}}$ There is an integer $\nu\le
\f{n_{1}+1}2$ such that $T^{(p)}_{d_{2}-1}\ne 0$ for
$p=1,2,\dots,\nu$ and
$T^{(\nu+1)}_{d_{2}-1}=T^{(\nu+2)}_{d_{2}-1}=\cdots =0$. That is,
$(h_{\nu},f)\not= (f_{1},f)$ and $(h_{\nu+1},f)= (f_{1},f)$ for an
integer $\nu \le \f{n_{1}+1}2$. $\square$
\endproclaim

{\bf Remark 15.5.1.} {\rm(a)} It is clear by Sublemma $15.5$ that
{\rm Property(1)}, {\rm Property(2)}, {\rm Property(3)}, {\rm
Property(4)} and {\rm Property(5)} for \text{\rm($h_1$,f)} in
Sublemma 15.5.$\beta$ are equivalent to {\rm Fact(A)}, {\rm Fact(C)}, {\rm
Fact(B)}, {\rm Fact(D)} and {\rm Fact(E)} for
\text{\rm($h_{1,1}$,f)} in Sublemma 15.5.$\alpha$, respectively.

{\rm(b)} It is clear by Sublemma $15.5.\alpha$ that $(h_1,f)$ of
Sublemma 15.5 was already constructed with five properties. $\square$
\bs

{\bf \S15.4.2. How to write the 2nd statement of Sublemma 15.5 for induction proof} \ms

$\underline{\text{\bf The 2nd statement for Sublemma 15.5}}$

$\underline{\text{\rm Assumptions for the 2nd statement}}$
Suppose we have shown by induction method that Sublemma 15.5 is true on the integer
$j\le {\ell-1}$.

$\underline{\text{\rm Conclusions for the 2nd statement}}$
Then, on the integer $j+1$ it suffices to show that given such a sequence
$\{f_1,f_2,\dots, f_j\}$ as we have seen in (15.4.1), then $f_{j+1}$
and $f$ can be constructed as follows: Note that (15.4.6)=(15.4.2).
$$
\cases f_{j+1} &=f^{n_{j+1}}_j +\sum^{n_{j+1}-2}_{i=0}R_{j+1,i}f^i_j, \\
f &=f^{d_{j+2}}_{j+1}+\sum^{d_{j+2}-2}_{i=0} S_{j+2,i}f^i_{j+1},
\endcases  \tag 15.4.6
$$
where $d_{j+2}=\Pi^{\ell}_{t=j+2}n_t$, and $f_{j+1}$ with
$R_{j+1,i}$ and $f$ with $S_{j+2,i}$ satisfy the same kind of
properties as $f_{j}$ with $R_{j,i}$ and $f$ with $S_{j+1,i}$ have
done in (1) and (2) of the conclusion of the theorem, respectively,
and also the properties in Corollary $15.4.1$. $\square$ \ms

For the proof of the 2nd statement of Sublemma 15.5, it suffices to prove
Sublemma 15.5.$\alpha$ for the 2nd statement and Sublemma 15.5.$\beta$ for the 2nd statement
respectively,
as follows:
\ms

\proclaim{Sublemma 15.5.$\alpha$ for Sublemma 15.5 of the 2nd statement}

$\underline{\text{\bf Assumptions}}$ Under the same induction
assumption on the integer $j\le \ell-1$ with $\ell \ge 2$, suppose
we have shown that the theorem is true on the integer $j\le \ell-1$,
just as above. \ms

$\underline{\text{\bf Conclusions}}$ We show that $h_{j+1,1}$ and
$f$ can be constructed as follows:
$$
\cases
h_{j+1,1} &=f^{n_{j+1}}_j +\sum^{n_{j+1}-2}_{i=0} R^{(1)}_{j+1,i}f^i_j, \\
f &=h^{d_{j+2}}_{j+1,1}+\sum^{d_{j+2}-1}_{i=0} T_{j+2,i}h^i_{j+1,1},
\endcases \tag 15.5.1-2
$$
where $h_{j+1,1}\in \C\{y\}[z]$ is a $W$-poly of degree
$\Pi^{j+1}_{t=1}n_t$ in $z$ with the multiplicity
$\Pi^{j+1}_{t=1}n_t$ at $0 \in \C^2$, satisfying the following
facts, Fact(A), Fact(B), Fact(C), Fact(D) and Fact(E).

$\underline{\text{\rm Fact(A)}}$ For each $i=0,1,\dots,n_{j+1}-2$,
$R^{(1)}_{j+1,i}=R_{j+1,i}^{(1)}(y,z)\in \C\{y\}[z]$ is a polynomial
of degree $<\Pi^j_{t=1}n_t$ in $z$  with $R^{(1)}_{j+1,i}(0,z)=0$
and has a multiplicity $\ge (n_{j+1}-i)\Pi^j_{t=1}n_t$ at $0 \in
\C^2$. \ms

$\underline{\text{\rm Fact(B)}}$ For each $i=0,1,\dots,n_{j+1}-2$,
and for any nonzero monomial $\Pi^{j+1}_{t=1}f^{\de_t}_{t-2}$ in
$R^{(1)}_{j+1,i}\in \C\{y\}[z,f_1,\dots, f_{j-1}]$, $\de_1>0$ and
$\de_t <n_{t-1}$ for $2\le t\le j+1$. \ms

$\underline{\text{\rm Fact(C)}}$ For each $i=0,1,\dots, d_{j+2}-1$,
$T_{j+2,i}=T_{j+2,i}(y,z)=\sum a_{p,q}y^pz^q$ with a nonzero
constant $a_{p,q}$ such that $p>0$ and $q<\Pi^{j+1}_{t=1}n_t$ and
that $T_{j+2,i}(0,z)=0$ and $T_{j+2,i}$ has a multiplicity $\ge
(d_{j+2}-i)\Pi^{j+1}_{t=1}n_t$ at $0 \in \C^2$. \ms

Moreover, consider $y,z,f_1,\dots,f_{j}$ as independent complex
$(j+2)$-variables at the origin in $\C^{j+2}$. Then, $T_{j+2,i}\in
\C\{y\}[z,f_1,\dots, f_j]$ satisfies two facts {\rm Fact(D) and
Fact(E)}.

$\underline{\text{\rm Fact(D)}}$  For each $i=0,1,\dots, d_{j+2}-1$,
and for any nonzero monomial $\Pi^{j+2}_{t=1}f^{\g_t}_{t-2}$ in
$T_{j+2,i}$, $\g_1>0$ and $\g_t<n_{t-1}$ for $2 \le t\le j+2$. \ms

$\underline{\text{\rm Fact(E)}}$ In particular, if $i=d_{j+2}-1$ for
$T_{j+2,i}$ of {\rm Fact(D)}, then $\g_{j+2}\le n_{j+1}-2$. $\square$
\endproclaim \ms

\proclaim{Sublemma 15.5.$\beta$ for Sublemma 15.5 of the 2nd statement}

$\underline{\text{\bf Assumptions}}$ By the
same way as we have seen in {\rm (15.5.1-2)} of the conclusion of
Sublemma $15.5.\alpha$, we may assume that $(h_{j+1,1},f)$ can be
rewritten as follows:
$$
\cases
h_{j+1,1} &=f^{n_{j+1}}_j +\sum^{n_{j+1}-2}_{i=0} R^{(1)}_{j+1,i}f^i_j, \\
f &=h^{d_{j+2}}_{j+1,1}+\sum^{d_{j+2}-1}_{i=0} T_{j+2,i}h^i_{j+1,1},
\endcases \tag 15.5.1-2
$$
satisfying the facts, denoted by {\rm Fact(A), Fact(B), Fact(C),
Fact(D) and Fact(E)}. For brevity of notation, let $h_1=h_{j+1,1}$,
$R^{(1)}_i =R^{(1)}_{j+1,i}$ for $0\le i\le n_{j+1}-2$ and
$T^{(1)}_i=T^{(1)}_{j+2,i}=T_{j+2,i}$ for $0\le i\le d_{j+2}-1$,
respectively. \ms

$\underline{\text{\bf Conclusions}}$ Then, $(f_{j+1},f)$ for $f$ can
be constructed as follows:

{\bf Case [I]:} If $T^{(1)}_{d_{j+2}-1}$ in $(h_1,f)$ is zero, let
$f_{j+1}=h_1$, $R_{j+1,i}=R^{(1)}_i$ for $0\le i\le n_{j+1}-2$ and
$S_{j+1,i}=T^{(1)}_i$ for $0\le i\le d_{j+2}-2$, respectively. Then,
the construction of $(f_{j+1},f)$ has been already finished. \ms

{\bf Case [II]:} If $T^{(1)}_{d_{j+2}-1}$ is not zero, for finding
such a construction of $(f_{j+1},f)$, it suffices to follow two
steps, {\rm Step(1)} and {\rm Step(2)}.

$\underline{\text{\bf Step(1) for Case[II]}}$ Then, there is a
sequence of pairs, $H=\{(h_p,f):p=1,2,\dots\}$, each pair of which
can be constructed with five properties, called {\rm Property(1)},
{\rm Property(2)}, {\rm Property(3)}, {\rm Property(4)} and {\rm
Property(5)}, as follows:\ms

\noindent\text{\rm(15.5.2-2)}
$$
\cases h_{1} &=f^{n_{j+1}}_j +\sum^{n_{j+1}-2}_{i=0}R^{(1)}_if^i_j
\text{\quad with $R^{(1)}_i=R^{(1)}_{j+1,i}$ in \text{\rm(15.5.1-2)},}\\
f &=h^{d_{j+2}}_{1}+\sum^{d_{j+2}-1}_{i=0} T^{(1)}_ih^i_{1}
\text{\quad with $T^{(1)}_i=T^{(1)}_{j+2,i}$ in \text{\rm(15.5.1-2)},}
\endcases \tag 15.5.2-2.1
$$
$$\cases h_{2} &=h_1 +\f 1{d_{j+2}}T^{(1)}_{d_{j+2}-1}=f^{n_{j+1}}_j
+\sum^{n_{j+1}-2}_{i=0}R^{(2)}_if^i_j, \qquad \qquad \quad\\
f &=h^{d_{j+2}}_{2}+\sum^{d_{j+2}-1}_{i=0} T^{(2)}_ih^i_{2},
\endcases \tag 15.5.2-2.2
$$
$$\cases h_{3} &=h_2 +\f 1{d_{j+2}}T^{(2)}_{d_{j+2}-1}=f^{n_{j+1}}_j
+\sum^{n_{j+1}-2}_{i=0}R^{(3)}_if^i_j, \qquad \qquad \quad\\
f &=h^{d_{j+2}}_{3}+\sum^{d_{j+2}-1}_{i=0} T^{(3)}_ih^i_{3},
\endcases \tag 15.5.2-2.3
$$
$$
 \dots, \qquad\qquad \qquad \qquad\qquad \qquad \qquad\qquad
 \qquad\qquad
$$
satisfying the following properties and notations: \ms

$\underline{\text{\rm Property(1)}}$ Let $p$ and $i$ be fixed with
$p\ge 1$ and $0\le i\le n_{j+1}-2$. Then
$R^{(p+1)}_i=R^{(p+1)}_i(y,z)\in \C\{y\}[z]$ is a polynomial of
degree $<\Pi^j_{t=1}n_t$ in $z$ and has a multiplicity $\ge
(n_{j+1}-i)\Pi^j_{t=1}n_t$ at $0 \in \C^2$.  \ms

$\underline{\text{\rm Property(2)}}$ Let $p$ and $i$ be fixed with
$p\ge 1$ and $0\le i\le d_{j+2}-1$. Then
$T^{(p+1)}_i=T^{(p+1)}_i(y,z)\in \C\{y\}[z]$ is a polynomial of
degree $<\Pi^{j+1}_{t=1}n_t$ in $z$ and has a multiplicity $\ge
(d_{j+2}-i)\Pi^{j+1}_{t=1}n_t$ at $0 \in \C^2$. \ms

Consider $y,z,f_1,\dots,f_j$ as independent complex
$(j+2)$-variables at the origin in $\C^{j+2}$.

$\underline{\text{\rm Property(3)}}$ Let $p$ and $i$ be fixed with
$p\ge 1$ and $0\le i\le n_{j+1}-2$. Then for any nonzero monomial
$\Pi^{j+1}_{t=1}f^{\de_t}_{t-2}$ in
$R^{(p+1)}_i=R^{(p+1)}_i(y,z,f_1,\dots,f_{j-1})\in
\C\{y\}[z,f_1,\dots,f_{j-1}]$, $\de_1>0$ and $\de_t <n_{t-1}$ for
$t=2,3,\dots, j+1$.  \ms

$\underline{\text{\rm Property(4)}}$ Let $p$ and $i$ be fixed with
$p\ge 1$ and $0\le i\le d_{j+2}-1$. Then for any nonzero monomial
$\Pi^{j+2}_{t=1}f^{\de_t}_{t-2}$ in
$T^{(p+1)}_i=T_i^{(p+1)}(y,z,f_1,\dots,f_j)\in
\C\{y\}[z,f_1,\dots,f_j]$, $\de_1>0$ and $\de_t<n_{t-1}$ for
$t=2,3,\dots,j+2$. \ms

$\underline{\text{\rm Property(5)}}$ In particular, if $i=d_{j+2}-1$
for $T^{(p+1)}_i$ of {\rm Property(4)}, then $\de_{j+2}\le
n_{j+1}-2$. \bs

$\underline{\text{\bf Step(2) for Case[II]}}$ By {\rm Step(1)},
there is a pair $(h_{\nu+1},f)\in H$ which satisfies the following
property:

$\underline{\text{\rm Property(6)}}$ There is an integer $\nu\le
\f{n_{j+1}+1}2$ such that $T^{(p)}_{d_{j+2}-1}\ne 0$ for
$p=1,2,\dots,\nu$ and
$T^{(\nu+1)}_{d_{j+2}-1}=T^{(\nu+2)}_{d_{j+2}-1}=\cdots =0$. That
is, $(h_{\nu},f)\not= (f_{j+1},f)$ and $(h_{\nu+1},f)= (f_{j+1},f)$
for an integer $\nu \le \f{n_{j+1}+1}2$. \ms

{\bf Remark 15.5.2.} {\rm(a)} It is clear
that {\rm Property(1)}, {\rm Property(2)}, {\rm Property(3)}, {\rm
Property(4)} and {\rm Property(5)} for
\text{\rm($h_1$,f)=($h_{j+1,1}$,f)} in Sublemma 15.5.$\beta$ are equivalent
to {\rm Fact(A)}, {\rm Fact(C)}, {\rm Fact(B)}, {\rm Fact(D)} and
{\rm Fact(E)} for \text{\rm($h_{j+1,1}$,f)} in Sublemma 15.5.$\alpha$,
respectively.

{\rm(b)} By {\rm(a)}, it is clear by Sublemma $15.5.\alpha$ that
$(h_1,f)$ was already constructed with five properties. $\square$
\endproclaim \ms

{\bf \S15.5. How to find the proofs of the 1st statement
and the 2nd statement of Sublemma 15.5 by the induction method } \ms

\noindent{\bf Remark 15.5.0.} To finish the proof of Sublemma 15.5,  it suffices to show by $\S15.5.1$ that
we can complete the proof of the 1st statement of Sublemma 15.5 in $\S15.4.1$
and to show by $\S15.5.2$ that applying directly the same properties and
notations as it can be used in the proof of the 1st statement of Sublemma 15.5
to the 2nd statement of Sublemma 15.5 in $\S15.4.2$ inductively, we can complete 
the proof of the 2nd statement of Sublemma 15.5 in $\S15.5.2$ 
very easily. \ms 

{\bf \S15.5.1. how to find the proof of the 1st statement
of Sublemma 15.5 by the induction method} \ms

\demo{\bf Proof of Sublemma 15.5.$\alpha$ for Sublemma 15.5 of the 1st statement}

We prove the above facts, respectively.

$\underline{\text{\rm The Proofs of Fact(A) and Fact(B):}}$ It is clear. \ms

$\underline{\text{\rm The Proof of Fact(C):}}$ Just, apply the WDT
with a divisor $h_{1,1}$ to $f$ as an element of $\C\{y\}[z]$, and
so it is clear by (1b) and (2a) of Theorem 15.2. \ms

$\underline{\text{\rm The Proof of Fact(D):}}$ It is clear. \ms

$\underline{\text{\rm The Proof of Fact(E):}}$ Assume the contrary.
By Fact(D), we get
$$ \text{$\g'_{2}=n_{1}-1$ for some nonzero
monomial $\Pi^{2}_{t=1}f^{\g'_t}_{t-2}$ in $T_{2,d_2-1}$,} \tag 15.5.3
 $$ where $\g'_1>0$ and $0\le \g'_2<n_1$.

To find a contradiction, note that for any nonzero monomial $y^\a
z^\beta$ in $\Pi^{2}_{t=1}f^{\nu_t}_{t-2}\in \C\{y,z\}$ where
$\nu_1>0$ and $\nu_2<n_{1}$, then
$$\align
 \beta = \nu_2
<n_1=\ \text{the multiplicity of
$f_1(y,z)$ at $0 \in \C^2$} \tag 15.5.4
\endalign$$

So, for any positive integer $p$, we can get
$$(\Pi^{2}_{t=1}f^{\nu_t}_{t-2})^p=\sum
a_{e_1,e_2,e_{3}}\Pi^{3}_{t=1}f^{e_t}_{t-2}, \tag 15.5.5
$$
as a convergent power series in $\C\{y,z,f_1\}$ where each
$a_{e_1e_2,e_{3}}$ is a nonzero complex number, $e_1>0$,
$e_2<n_{1}$ and $e_3<n_{2}$, by using the same technique as in the proof of Fact (D)
and the inequality in (15.4.7).

To prove Fact (E), recall by Sublemma 15.5.$\alpha$ for Sublemma 15.5 of
the 1st statement that $f=h^{d_{2}}_{1,1}+\sum^{d_{2}-1}_{i=0} T_{2,i}h^i_{1,1}$
where for each fixed $i$, and any nonzero monomial
$\Pi^{2}_{t=1}f^{\g_t}_{t-2}$ in $T_{2,i}$ $\g_1>0$ and
$\g_2<n_{1}$.

Then for fixed $j$, and $0\le i\le d_{2}-1$,

$f-f^{d_{1}}_0 =h^{d_{2}}_{1,1}-f^{d_{1}}_0+\sum
^{d_{2}-1}_{i=0} T_{2,i}h^i_{1,1}=\sum_1 +\sum_2+\sum_3$ can
be written as follows:
$$
\cases \sum_1 &=h^{d_{2}}_{1}-f^{d_{1}}_0
=\{f^{n_{1}}_0+\sum^{n_{1}-2}_{i=0}
R^{(1)}_{1,i}f^i_0\}^{d_{2}}-f^{d_{1}}_0,\\
\sum_2 &=\sum^{d_{1}-2}_{i=0}T_{2,i}h^i_{1}, \\
\sum_3 &=T_{2,d_{2}-1}h^{d_{2}-1}_{1}.
\endcases \tag 15.5.6
$$
Note that $\sum_1$ can be rewritten as follows: With nonzero
constants $b_s$,

\noindent $\sum_1=\sum^{d_{2}}_{s=1}b_s(f^{n_{1}}_0)^{d_{2}-s}
(\sum^{n_{1}-2}_{i=0}R^{(1)}_{1,i}f^i_0)^s$.

First, by the assumption, observe by (15.5.3) that there is a nonzero
monomial
$$
(\Pi^{2}_{k=1}f^{\g'_k}_{k-2})(f^{n_{1}}_0)^{d_{2}-1}
=(\Pi^{1}_{k=1}f^{\g'_k}_{k-2})f^{d_{1}-1}_0, \tag 15.5.7
$$
in $\sum_3$ because
$\g'_{2}+n_{1}(d_{2}-1)=n_{1}-1+d_{1}-n_{1}=d_{1}-1$ by (15.5.3).
\ms

Next, we claim that
$(\Pi^{1}_{k=1}f^{\g'_k}_{k-2})f^{d_{1}-1}_j\in \sum_3$ does not
belong to both $\sum_1$ and $\sum_2$.

If the claim is proved, then the above monomial of (15.5.7) in
$\sum_3$ would belong to $f(y,z)-f^{d_{1}}_0$. But
it is impossible because such an monomial does not belong to
$f-f^{d_{1}}_0$ by induction assumption on
$f=f^{d_{1}}_0+\sum^{d_{1}-2}_{i=0} S_{2,i}f^i_0$.

Now, to prove the claim, it is enough to consider the following two
cases :

Case (1): Whenever $\Pi^{2}_{k=1}f^{\tau_k}_{k-2}$ is in $\sum_1$
with $\tau_1>0$ and $\tau_2<n_{1}$, then
$\tau_{2}\le n_{1}(d_{2}-1)+n_{1}-2=d_{1}-2$.

Case (2): Whenever $\Pi^{2}_{k=1}f^{\tau_k}_{k-2}$ is in $\sum_2$
with $\tau_1>0$ and $\tau_2<n_{1}$, then
$\tau_{2}<n_{1}(d_{2}-2)+n_{1}=d_{1}-n_{1}\le
d_{1}-2$. \ms

So, we proved the claim. Therefore, we can find a contradiction,
what we wanted. Thus, we proved Fact(E). $\square$
\enddemo \bs

\demo{\bf Proof of Sublemma 15.5.$\beta$ for Sublemma 15.5 of the 1st statement}

We prove two cases, respectively.

{\bf Case[I]}: By Sublemma $15.5.\alpha$, there is nothing to prove
for this sublemma.

{\bf Case[II]}: Assume that $T^{(1)}_{d_{j+2}-1}$ is not zero. For
the proof of this case, it suffices to follow two steps, Step(1) and
Step(2) in order.

$\underline{\text{\bf Step(1) for Case[II]}}$ We will prove this
case by induction on the integer $p\ge 1$ that each pair of $H$
satisfies {\rm Property(1)}, {\rm Property(2)}, \dots, {\rm
Property(5)} where $H=\{(h_p,f):p=1,2,\dots\}$ is called a sequence
of pairs, defined by (15.5.2.1),(15.5.2.2),\dots.

$\underline{\text{\bf Step(2) for Case[II]}}$ We will prove by
Step(1) that there is some integer $\nu \le \f{n_{1}+1}2$ which
satisfies {\rm Property(6)}, being equivalent to the fact that
$(h_{\nu},f)\not= (f_{1},f)$ and $(h_{\nu+1},f)= (f_{1},f)$. \ms

$\underline{\text{\bf The Proof of Step(1) for Case[II]}}$ As we
have seen in Remark 15.4.1, $(h_1,f)$ with $T^{(1)}_{d_{2}-1}\ne 0$
has already constructed by (15.5.1) and Sublemma $15.4.\alpha$. The
general case is proved by induction. Suppose we have shown that for
$p=1,2,\dots,w$, each pair $(h_p,f)$ satisfies the same kind of
properties corresponding to the properties {\rm Property(1)}, {\rm
Property(2)}, \dots, Property(5), which the pair $(h_1,f)$ does.

In order to construct such a pair $(h_{w+1},f)$, first define
$h_{w+1}$ by $h_{w+1}=h_w +\f 1{d_{2}}T^{(w)}_{d_{2}-1}$ where
$T^{(w)}_{d_{2}-1}$ was already defined from $(h_w,f)$. By
induction assumption on $(h_w,f)$, note that $T^{(w)}_{d_{2}-1}$
satisfies the corresponding properties, denoted by {\rm
Property(2)}, {\rm Property(4)} and Property(5) of Case [II], if
$p=w$. So, $h_{w+1}$ is rewritten in the form
$$\align
h_{w+1} &=f^{n_{1}}_0 +\sum^{n_{1}-2}_{i=0} R^{(w)}_if^i_0 +\f
1{d_{2}}T^{(w)}_{d_{2}-1} \tag 15.5.8 \\
&=f^{n_{1}}_{0}+\sum^{n_{1}-2}_{i=0} R^{(w+1)}_if^i_0,
\endalign
$$
and then $R^{(w+1)}_i$ satisfies the same kind of properties as
$R^{(w)}_i$ does in Property(1) and Property(3) in Case [II] of the
lemma, which can be proved by the following way:

Since for any nonzero monomial $\Pi^{2}_{t=1}f^{\g_t}_{t-2}$ in
$T^{(w)}_{d_{2}-1}\in \C\{y\}[z]$, $\g_1>0$ and $\g_{2}\le n_{1}-2$, then
we can put $T^{(w)}_{d_{2}-1}=\sum^{n_{1}-2}_{i=0}
U^{(w)}_if^i_0$ where for any nonzero monomial
$\Pi^{2}_{t=1}f^{\de_t}_{t-2}$ in $U^{(w)}_i\in
\C\{y\}[z]$, $\de_1>0$ and $\de_2=0$. So, note that for each $i$, $U^{(w)}_i=U^{(w)}_i(y)$
has a multiplicity $\ge n_1-i$ at $0 \in \C^2$ because $T^{(w)}_{d_{j+2}-1}$
has a multiplicity $\ge n_1$ at $0 \in \C^2$ by Property(2).

Now, apply the WDT to $f$ with a divisor $h_{w+1}$, and then it can
be proved by Theorem 15.2 that $f$ is written in the form
$$
f=h^{d_{2}}_{w+1}+\sum^{d_{2}-1}_{i=0} T^{(w+1)}_ih^i_{w+1},
\tag 15.5.9
$$
and then $T^{(w+1)}_i$ satisfies the same kind of properties as
$T^{(w)}_i$ does in Property(2), Property(4) and Property(5) in Case
[II] of the lemma. Thus, the proof of Step(1) is done. \ms

$\underline{\text{\bf The Proof of Step(2) for Case[II]}}$ Let
$T^{(1)}_{d_{2}-1}$ in $(h_1,f)$ be not zero. For the proof of
Step(2), it suffices to show by Step(1) that
$(f_{1},f)=(h_{\nu+1},f)$ with the following property:
$$\align
& \text{$T^{(\nu+1)}_{d_{2}-1}=0$ \quad for some integer \quad
$\nu\le \f{n_{1}+1}2$,} \tag 15.5.10 \\
& \text{$T^{(p)}_{d_{2}-1}\ne 0$ \quad  for \quad
$p=1,2,\dots,\nu$.}
\endalign$$

First of all, we have already proved by Step(1) that for given
$(h_p,f)$, $(h_{p+1},f)$ can be written in the form
$$ (15.5.11)\qquad \qquad \quad
\cases h_{p+1} &=h_p+\f 1{d_{2}}T^{(p)}_{d_{2}-1}=f^{n_{1}}_0
+\sum^{n_{1}-2}_{i=0}
R^{(p+1)}_if^i_0, \\
f &=h_{p+1}^{d_{2}}+\sum^{d_{2}-1}_{i=0} T^{(p+1)}_ih^i_{p+1},
\endcases \qquad \qquad \quad
$$
satisfying the properties {\rm Property(1)}, {\rm Property(2)},
\dots, Property(5) in this sublemma except possibly for the above
equation (15.5.5).

In preparation for the proof of the equation in (15.5.5), consider
$y=f_{-1}, z=f_0$ as independent complex $2$-variables at the
origin in $\C^{2}$ as we have seen in the proof of Step(1), if
necessary. For brevity of notation, if $T=T(f_{-1},f_0)$ is
not zero in $\C\{f_{-1}\}[f_0]$, define
$\text{deg}_0{T}$ by the degree of $T$ as a polynomial in $f_0$, when
$T$ is in $\C\{f_{-1}\}[f_0]$ as a polynomial in
$f_0$ with coefficients in $\C\{f_{-1}\}$. For such
notations, write $\tau_p=\text{deg}_j(T^{(p)}_{d_{2}-1})$ for all
$p\ge 1$, and so $0\le \tau_p\le n_{1}-2$ by Property(5).

First, in preparation for the proof of an equality in (15.5.10), it is
very interesting to show that if $T^{(w)}_{d_{2}-1}$ and
$T^{(w+1)}_{d_{2}-1}$ is not zero for any integer $w>0$, then the
following inequality is true:
$$\align
n_{1}-2\ge \tau_w
>\tau_{w+1}\ge 0 \quad \text{with} \quad \tau_w\ge \f{n_{1}-1}2, \tag 15.5.12
\endalign$$
because if $T^{(p)}_{d_{2}-1}$ is not zero for any $p\ge 1$ and an
inequality in (15.5.7) is true, then $\{\tau_p:p=1,2,\dots\}$ is an
infinite sequence which is strictly decreasing and bounded with
$n_{1}-2\ge \tau_p>\tau_{p+1}\ge 0$. This would be impossible,
which is equivalent to the fact that there is an integer
$\nu<n_{1}$ such that $T^{(\nu+1)}_{d_{2}-1}=0$.

Now, assuming that $T^{(w)}_{d_{2}-1}$ and $T^{(w+1)}_{d_{2}-1}$
is not zero for any integer $w>0$, then we are going to prove that
an inequality in (15.5.12) is true.

Since $h_{w+1}=h_w+\f 1{d_{2}}T^{(w)}_{d_{2}-1}$, then a
sequence in (15.5.11) implies that
$$\align
f&=h^{d_{2}}_w +\sum^{d_{2}-1}_{i=0} T^{(w)}_ih^i_w \tag 15.5.13 \\
&=h^{d_{2}}_{w+1} +\sum^{d_{2}-1}_{i=0} T^{(w+1)}_ih^i_{w+1}\\
&=\left(h_w+\f 1{d_{2}}T^{(w)}_{d_{2}-1}\right)^{d_{2}}
+\sum^{d_{2}-1}_{i=0} T^{(w+1)}_ih^i_{w+1}. \tag 15.5.14
\endalign
$$

\quad Therefore, \qquad
$f-h^{d_{2}}_w-T^{(w)}_{d_{2}-1}h^{d_{2}-1}_w =\sum^{(w)}_0
=\sum_{1,1}^{(w)}+\sum_{1,2}^{(w)}$ {\quad} such that
$$\align
 \sum^{(w)}_0 =\sum^{d_{1}-2}_{i=0}
T^{(w)}_ih^i_w \quad \text{by (15.5.13)}, \qquad \qquad \quad \tag 15.5.15 \\
\endalign$$
$$
(15.5.16) \qquad \qquad \qquad  \cases \sum^{(w)}_{1,1}
&=\sum^{d_{2}-2}_{q=0}
c_q(T^{(w)}_{d_{2}-1})^{d_{2}-q}h^q_w \quad \text{by (15.5.14)}, \\
\sum^{(w)}_{1,2} &=\sum^{d_{2}-1}_{i=0} T^{(w+1)}_ih^i_{w+1},
\endcases \qquad \qquad \qquad \qquad
$$
with some nonzero constant $c_q$ by (15.5.14).

Then, to prove that $\tau_w>\tau_{w+1}\ge 0$ and $\tau_w\ge
\f{n_{1}-1}2$, note that given any nonzero monomial
$\Pi^{2}_{t=1}f^{\de_t}_{t-2}$ in either $T^{(w)}_i$ or
$T^{(w+1)}_i$ for all $i$, $\de_1>0$ and $\de_2<n_{1}$, and
that $h_{w+1}\in\C\{y\}[z]$ is a $W$-poly of degree $n_{1}$ in $f_0$.

Also, note by (15.5.4) that for any
nonzero monomial $y^\a z^\beta$ in $\Pi^{2}_{t=1}f^{\g_t}_{t-2}$,
where $\g_1>0$ and $\g_2<n_{1}$,
$$
\beta <n_1=\ \text{the multiplicity of $f_1$ at $0 \in
\C^2$.} \tag 15.5.17
$$
Now, to finish the proof of (15.5.7), let $\mu_1,\mu_2$ and $\mu_3$
be defined by $\text{$\max$}\{\g_{2}\}$ among the degree of $f_0$
for all nonzero monomials $\Pi^{2}_{t=1}f^{\g_t}_{t-2}$ in
$\sum^{(w)}_0$, $\sum^{(w)}_{1,1}$ and $\sum^{(w)}_{1,2}$,
respectively, where $\g_1>0$ and $\g_2<n_{1}$.

Then, we need only to consider three subcases (i), (ii) and (iii),
respectively.

(i) In $\sum^{(w)}_0$,
$\mu_1<n_{1}+n_{1}(d_{2}-2)=d_{1}-n_{1}$ because
$d_{1}=n_{1}d_{2}$. \ms

(ii) In $\sum^{(w)}_{1,1}$, we have two possibilities by (15.5.12),
noting by Property(5) that $\tau_w
=\text{deg}_0(T^{(w)}_{d_{2}-1})\le n_{1}-2$.

(iia) If $\beta <\f {1}{2}$ for all $\beta$ in (15.5.12),
then
$\mu_2=2\de_{2}+n_{1}(d_{2}-2)=2\tau_w+d_{1}-2n_{1}\le
2(n_{1}-2)+d_{1}-2n_{1}=d_{1}-4$ for some nonzero monomials
$\Pi^{2}_{t=1}f^{\de_t}_{t-2}$ in $T^{(w)}_{d_{2}-1}$. \ms

(iib) If $\beta \ge \f {1}{2}$ for some $\beta$ in
(15.5.12), then
$\mu_2=2\de_{2}+1+n_{1}(d_{2}-2)=2\tau_w+1+d_{1}-2n_{1}\le
2(n_{1}-2)+1+d_{1}-2n_{1}=d_{1}-3$ for some nonzero
monomials $\Pi^{2}_{t=1}f^{\de_t}_{t-2}$ in $T^{(w)}_{d_{2}-1}$.
\ms

(iii) In $\sum^{(w)}_{1,2}$,
$\mu_3=\de'_{2}+n_{1}(d_{2}-1)=\tau_{w+1}+d_{1}-n_{1}\ge
d_{1}-n_{1}$ for some nonzero monomials
$\Pi^{2}_{t=1}f^{\de'_t}_{t-2}$ in $T^{(w+1)}_{d_{2}-1}$ where
$\tau_{w+1}=\text{deg}_j(T^{(w+1)}_{d_{2}-1})$.

By (i) and (iii), it is clear that $\mu_3>\mu_1$. \ms

Then, we claim that $\mu_2=\mu_3$. To prove the claim, assume the
contrary.

Then, it suffices to consider two possibilities (A) and (B),
respectively:

(A) If $\mu_3 >\mu_2$, then there is a nonzero monomial
$\left(\Pi^{2}_{t=1}f^{\de'_t}_{t-2}\right)f^{n_{1}(d_{2}-1)}_j$
which belongs to $\sum^{(w)}_{1,2}$, but does not belong to
$\sum^{(w)}_{1,1}$ and $\sum^{(w)}_0$ because $\mu_3>\mu_1$. It
would be a contradiction to $\sum^{(w)}_0
=\sum_{1,1}^{(w)}+\sum_{1,2}^{(w)}$.

(B) If $\mu_2>\mu_3$, then there is a nonzero monomial
$\left(\Pi^{1}_{t=1}f^{\de_t}_{t-2}\right)f^{\mu_2}_{j+2}$ which
belongs to $\sum^{(w)}_{1,1}$, but does not belong to $\sum^{(w)}_0$
and $\sum^{(w)}_{1,2}$ because $\mu_2>\mu_1$. It would be a
contradiction to $\sum^{(w)}_0 =\sum_{1,1}^{(w)}+\sum_{1,2}^{(w)}$.

Thus, we proved the claim. \ms

First, in preparation for the proof of an inequality in (15.5.7),
using (ii) and (iii) with $\mu_2=\mu_3$, then it suffices to
consider two possibilities (a) and (b), as follows:

(a) It is clear by (iia) and (iii) that $\mu_2=\mu_3$ implies
$2\tau_w=\tau_{w+1}+n_{j+1}$.

(b) It is clear by (iib) and (iii) that $\mu_2=\mu_3$ implies
$2\tau_w +1=\tau_{w+1}+n_{1}$. \ms

Recall by Property(5) of Case(II) in this sublemma that $\tau_p$ was
defined by $\tau_p=\text{deg}_j(T^{(p)}_{d_{2}-1})$ for all $p\ge
1$, and so $0\le \tau_p\le n_{1}-2$, if exists.

(a) First, if $2\tau_w=\tau_{w+1}+n_{1}$, then it is trivial that
$\tau_w-\tau_{w+1}=\f{\tau_{w+1}+n_{1}}{2}-\tau_{w+1}
=\f{n_{1}-\tau_{w+1}}{2}>0$ because $\tau_{w+1}\le n_{1}-2$.

(b) Next, if $2\tau_w+1=\tau_{w+1}+n_{1}$, then it is trivial that
$\tau_w-\tau_{w+1}=\f{\tau_{w+1}+n_{1}-1}{2}-\tau_{w+1}
=\f{n_{1}-\tau_{w+1}-1}{2}>0$ because $\tau_{w+1}\le n_{1}-2$.

Moreover, it is clear by (a) and (b) that $\tau_w\ge
\f{n_{1}-1}2$, which is an inequality in (15.5.7).

Therefore, we proved that $\tau_w>\tau_{w+1}\ge 0$ with $\tau_w\ge
\f{n_{1}-1}2$ for (15.5.7).

So, assuming that $T^{(p)}_{d_{2}-1}$ is not zero for any $p\ge
1$, since an inequality in (15.5.7) is true, then
$\{\tau_p:p=1,2,\dots\}$ is an infinite sequence which is strictly
decreasing and bounded with $n_{1}-2\ge \tau_p>\tau_{p+1}\ge 0$.
This would be impossible, and therefore there is an integer
$\nu<n_{1}$ such that $T^{(\nu+1)}_{d_{2}-1}=0$.

Finally, to prove that $\nu \le \f{n_{1}+1}2$, define a function
$\psi :\N\to \N\cup \{0\}$ by $\psi(p)=\tau_p$ where
$\tau_p=\text{deg}_j(T^{(p)}_{d_{2}-1})$ for all $p\ge 1$.

Then, $\nu-2\le (\psi(1)-\psi(2))+(\psi(2)-\psi(3))+\cdots
+(\psi(\nu-2)-\psi(\nu-1))=\psi(1)-\psi(\nu-1)
=\tau_1-\tau_{\nu-1}\le n_{1}-2-\f{n_{1}-1}2=\f{n_{1}-3}2$,
because $\psi(\nu-1)\ne 0$. So, $\nu\le \f{n_{1}+1}2$. Thus, we
proved Property(6), and so the proof of Step(2) is done.

Thus, we completed the proof of Sublemma 15.5. $\square$
\enddemo \bs

{\bf \S15.5.2. For the proof of the 2nd statement by Sublemma 15.5 } \ms

\demo{\bf Proof of Sublemma 15.5.$\alpha$ for Sublemma 15.5 of the 2nd statement}

We prove the above
facts, respectively.

$\underline{\text{\rm The Proofs of Fact(A) and Fact(B):}}$ We
construct $R^{(1)}_{j+1,i}\in \C\{y\}[z,f_1,\dots, f_{j-1}]$, which
satisfies Fact(A) and Fact(B), by the following method:

As an example, for each fixed $i=0,1,\dots,n_{j+1}-2$, and for any
nonzero monomial $\Pi^{j+1}_{t=1}f^{\de_t}_{t-2}$ in
$R^{(1)}_{j+1,i}$, $\de_1>0$ can be chosen sufficiently large such
that $\de_t<n_{t-1}$ for $2 \le t\le j+1$, which implies the
following:

(i) $R^{(1)}_{j+1,i}=R^{(1)}_{j+1,i}(y,z)$ is a polynomial of degree
$<\Pi^{j+1}_{t=1}n_t$ in $z$, because $\de_2+n_1\de_3 +n_1n_2\de_4
+\cdots +n_1n_2\cdots n_{j-1}\de_{j+1}\le
n_1-1+n_1(n_2-1)+n_1n_2(n_3-1)+\cdots +n_1n_2\cdots
n_{j-1}(n_j-1)<\Pi^j_{t=1}n_t,$ using by the induction assumption on
$j$ that each $f_k(y,z)$ is a polynomial of degree $\Pi^k_{t=1}n_t$
in $z$ for $0\le k\le j-1$. \ms

(ii) Also, it is clear that $R^{(1)}_{j+1,i}(0,z)=0$ and
$R^{(1)}_{j+1,i}(y,z)$ has the desired multiplicity by (2a) of
Theorem 15.2, because $\de_1>0$ can be chosen sufficiently large
such that $h_{j+1,1}\in \C\{y\}[z]$ is a $W$-poly of degree
$\Pi^{j+1}_{t=1}n_t$ in $z$ with the  multiplicity
$\Pi^{j+1}_{t=1}n_t$ at $(y,z)=(0,0)$.

Thus, Fact(A) and Fact(B) are easily proved. \ms

$\underline{\text{\rm The Proof of Fact(C):}}$ Just, apply the WDT
with a divisor $h_{j+1,1}$ to $f$ as an element of $\C\{y\}[z]$, and
so it is clear by (1b) and (2a) of Theorem 15.2. \ms

$\underline{\text{\rm The Proof of Fact(D):}}$ Recall that
$f_j(y,z)\in \C\{y\}[z]$ is a $W$-poly of degree $\Pi^j_{t=1}n_t$ in
$z$ with the multiplicity $\Pi^j_{t=1}n_t$ at $0 \in \C^2$ by the
induction assumption on $j$ and that $T_{j+2,i}(y,z)$ is a
polynomial of degree $<\Pi^{j+1}_{t=1}n_t$ in $z$ by Fact (C). For
each fixed $i=0,1,\dots,d_{j+2}-1$, apply the WDT with a divisor
$f_j$ to $T_{j+2,i}$ and then by (1d) and (2b) of Theorem 15.2,
$T_{j+2,i}$ can be written as follows:
$$\align
T_{j+2,i}=\sum^{n_{j+1}-1}_{k_1=0} Q_{k_1}f^{k_1}_j, \tag Eq.1
\endalign$$
where for $0\le k_1\le n_{j+1}-1$, each $Q_{k_1}\in \C\{y\}[z]$ is a
polynomial of degree $<\Pi^j_{k=1}n_t$ in $z$ with $Q_{k_1}(0,z)=0$,
if exists. Since $f_{j-1}\in \C\{y\}[z]$ is a $W$-poly of degree
$\Pi^{j-1}_{t=1}n_t$ in $z$ with the multiplicity
$\Pi^{j-1}_{t=1}n_t$ at $0 \in \C^2$, then for each fixed $k_1$
apply the WDT with a divisor $f_{j-1}$ to $Q_{k_1}$ again, and then
by (1d) and (2b) of Theorem 15.2 again, $Q_{k_1}$ may be written in
the form
$$
Q_{k_1}=\sum^{n_j-1}_{k_2=0} Q_{k_1,k_2}f^{k_2}_{j-1}, \tag Eq.2
$$
where for $0\le k_2\le n_{j}-1$, each $Q_{k_1,k_2}\in \C\{y\}[z]$ is
a polynomial of degree $<\Pi^{j-1}_{t=1}n_t$ in $z$ with
$Q_{k_1,k_2}(0,z)$ zero, if exists. Thus, continuing the same
process as above, then Fact (D) can be easily proved by Theorem
15.2. \ms

$\underline{\text{\rm The Proof of Fact(E):}}$ Assume the contrary.
By Fact(D), we get
$$ \text{\rm(15.5.3-2)} \qquad\text{$\g'_{j+2}=n_{j+1}-1$ for some nonzero
monomial $\Pi^{j+2}_{t=1}f^{\g'_t}_{t-2}$ in $T_{j+1,d_{j+2}-1}$,}
\qquad $$ where $\g'_1>0$ and $0\le \g'_t<n_{t-1}$ for $2\le t\le
j+2$.

To find a contradiction, note that for any nonzero monomial $y^\a
z^\beta$ in $\Pi^{j+1}_{t=1}f^{\nu_t}_{t-2}\in \C\{y,z\}$ where
$\nu_1>0$ and $\nu_t<n_{t-1}$ for $t=2,3,\dots, j+1$, then
$$\align
\text{\rm(15.5.4-2)} \qquad \beta &\le \nu_2+n_1\nu_3 +n_1n_2\nu_4 +\cdots
+n_1n_2\cdots
n_{j-1}\nu_{j+1} \\
&\le n_1-1+n_1(n_2-1)+n_1n_2(n_3-1)+\cdots +n_1n_2\cdots
n_{j-1}(n_j-1) \qquad \qquad \\
&= \Pi^j_{t=1}n_t-1<\Pi^j_{t=1}n_t=\ \text{the multiplicity of
$f_j(y,z)$ at $0 \in \C^2$}.
\endalign$$

So, for any positive integer $p$, we can get
$$(\Pi^{j+1}_{t=1}f^{\nu_t}_{t-2})^p=\sum
a_{e_1,e_2,\dots, e_{j+2}}\Pi^{j+2}_{t=1}f^{e_t}_{t-2}, \tag 15.5.5-2
$$
as a convergent power series in $\C\{y,z,f_1,\dots,f_j\}$ where each
$a_{e_1e_2\cdots e_{j+2}}$ is a nonzero complex number, and $e_1>0$,
$e_t<n_{t-1}$ for $t=2,3,\dots, j+1$ and $e_{j+2}<p$ if exists,
using the same technique as in the proof of Fact (D) and the
inequality in (15.5.4-2).

Now, to prove Fact (E), recall by Fact(4) that
$f=h^{d_{j+2}}_{j+1,1}+\sum^{d_{j+2}-1}_{i=0} T_{j+2,i}h^i_{j+1,1}$
where for each fixed $i$, and any nonzero monomial
$\Pi^{j+2}_{t=1}f^{\g_t}_{t-2}$ in $T_{j+2,i}$ $\g_1>0$ and
$\g_t<n_{t-1}$ for $t=2,3,\dots, j+2$.

Then for fixed $j$, and $0\le i\le d_{j+2}-1$,

$f-f^{d_{j+1}}_j =h^{d_{j+2}}_{j+1,1}-f^{d_{j+1}}_j+\sum
^{d_{j+2}-1}_{i=0} T_{j+2,i}h^i_{j+1,1}=\sum_1 +\sum_2+\sum_3$ can
be written as follows:
$$
\cases \sum_1 &=h^{d_{j+2}}_{j+1}-f^{d_{j+1}}_j
=\{f^{n_{j+1}}_j+\sum^{n_{j+1}-2}_{i=0}
R^{(1)}_{j+1,i}f^i_j\}^{d_{j+2}}-f^{d_{j+1}}_j,\\
\sum_2 &=\sum^{d_{j+1}-2}_{i=0}T_{j+2,i}h^i_{j+1}, \\
\sum_3 &=T_{j+2,d_{j+2}-1}h^{d_{j+2}-1}_{j+1}.
\endcases \tag 15.5.6-2
$$
Note that $\sum_1$ can be rewritten as follows: With nonzero
constants $b_s$,

\noindent $\sum_1=\sum^{d_{j+2}}_{s=1}b_s(f^{n_{j+1}}_j)^{d_{j+2}-s}
(\sum^{n_{j+1}-2}_{i=0}R^{(1)}_{j+1,i}f^i_j)^s$.

First, by the assumption, observe by (15.5.3-2) that there is a nonzero
monomial
$$
(\Pi^{j+2}_{k=1}f^{\g'_k}_{k-2})(f^{n_{j+1}}_j)^{d_{j+2}-1}
=(\Pi^{j+1}_{k=1}f^{\g'_k}_{k-2})f^{d_{j+1}-1}_j, \tag 15.5.7-2
$$
in $\sum_3$ because
$\g'_{j+2}+n_{j+1}(d_{j+2}-1)=n_{j+1}-1+d_{j+1}-n_{j+1}=d_{j+1}-1$.
\ms

Next, we claim that
$(\Pi^{j+1}_{k=1}f^{\g'_k}_{k-2})f^{d_{j+1}-1}_j\in \sum_3$ does not
belong to both $\sum_1$ and $\sum_2$.

If the claim is proved, then the above monomial of (15.5.7-2) in
$\sum_3$ would belong to $f(y,z,f_1,\dots,f_j)-f^{d_{j+1}}_j$. But
it is impossible because such an monomial does not belong to
$f-f^{d_{j+1}}_j$ by induction assumption on
$f=f^{d_{j+1}}_j+\sum^{d_{j+1}-2}_{i=0} S_{j+1,i}f^i_j$.

Now, to prove the claim, it is enough to consider the following two
cases :

Case (1): Whenever $\Pi^{j+2}_{k=1}f^{\tau_k}_{k-2}$ is in $\sum_1$
with $\tau_1>0$ and $\tau_k<n_{k-1}$ for $2\le k\le j+1$, then
$\tau_{j+2}\le n_{j+1}(d_{j+2}-1)+n_{j+1}-2=d_{j+1}-2$.

Case (2): Whenever $\Pi^{j+2}_{k=1}f^{\tau_k}_{k-2}$ is in $\sum_2$
with $\tau_1>0$ and $\tau_k<n_{k-1}$ for $2\le k\le j+1$, then
$\tau_{j+2}<n_{j+1}(d_{j+2}-2)+n_{j+1}=d_{j+1}-n_{j+1}\le
d_{j+1}-2$. \ms

So, we proved the claim. Therefore, we can find a contradiction,
what we wanted. Thus, we proved Fact(E). $\square$
\enddemo \ms

\demo{\bf Proof of Sublemma 15.5.$\beta$ for Sublemma 15.5 of the 2nd statement}

We prove two cases, respectively.

{\bf Case[I]}: By Sublemma $15.5.\alpha$, there is nothing to prove
for this sublemma.

{\bf Case[II]}: Assume that $T^{(1)}_{d_{j+2}-1}$ is not zero. For
the proof of this case, it suffices to follow two steps, Step(1) and
Step(2) in order.

$\underline{\text{\bf Step(1) for Case[II]}}$ We will prove this
case by induction on the integer $p\ge 1$ that each pair of $H$
satisfies {\rm Property(1)}, {\rm Property(2)}, \dots, {\rm
Property(5)} where $H=\{(h_p,f):p=1,2,\dots\}$ is called a sequence
of pairs, defined by (15.5.2-2.1),(15.5.2-2.2),\dots.

$\underline{\text{\bf Step(2) for Case[II]}}$ We will prove by
Step(1) that there is some integer $\nu \le \f{n_{j+1}+1}2$ which
satisfies {\rm Property(6)}, being equivalent to the fact that
$(h_{\nu},f)\not= (f_{j+1},f)$ and $(h_{\nu+1},f)= (f_{j+1},f)$. \ms

$\underline{\text{\bf The Proof of Step(1) for Case[II]}}$ As we
have seen in Remark 15.5.1, $(h_1,f)$ with $T^{(1)}_{d_{j+2}-1}\ne 0$
has already constructed by (15.5.1) and Sublemma $15.5.\alpha$. The
general case is proved by induction. Suppose we have shown that for
$p=1,2,\dots,w$, each pair $(h_p,f)$ satisfies the same kind of
properties corresponding to the properties {\rm Property(1)}, {\rm
Property(2)}, \dots, Property(5), which the pair $(h_1,f)$ does.

In order to construct such a pair $(h_{w+1},f)$, first define
$h_{w+1}$ by $h_{w+1}=h_w +\f 1{d_{j+2}}T^{(w)}_{d_{j+2}-1}$ where
$T^{(w)}_{d_{j+2}-1}$ was already defined from $(h_w,f)$. By
induction assumption on $(h_w,f)$, note that $T^{(w)}_{d_{j+2}-1}$
satisfies the corresponding properties, denoted by {\rm
Property(2)}, {\rm Property(4)} and Property(5) of Case [II], if
$p=w$. So, $h_{w+1}$ is rewritten in the form
$$\align
h_{w+1} &=f^{n_{j+1}}_j +\sum^{n_{j+1}-2}_{i=0} R^{(w)}_if^i_j +\f
1{d_{j+2}}T^{(w)}_{d_{j+2}-1} \tag 15.5.8-2 \\
&=f^{n_{j+1}}_{j}+\sum^{n_{j+1}-2}_{i=0} R^{(w+1)}_if^i_j,
\endalign
$$
and then $R^{(w+1)}_i$ satisfies the same kind of properties as
$R^{(w)}_i$ does in Property(1) and Property(3) in Case [II] of the
lemma, which can be proved by the following way:

Since for any nonzero monomial $\Pi^{j+2}_{t=1}f^{\g_t}_{t-2}$ in
$T^{(w)}_{d_{j+2}-1}\in \C\{y\}[z,f_1,\dots,f_j]$, $\g_1>0$,
$\g_t<n_{t-1}$ for $t=2,\dots,j+1$ and $\g_{j+2}\le n_{j+1}-2$, then
we can put $T^{(w)}_{d_{j+2}-1}=\sum^{n_{j+1}-2}_{i=0}
U^{(w)}_if^i_j$ where for any nonzero monomial
$\Pi^{j+1}_{t=1}f^{\de_t}_{t-2}$ in $U^{(w)}_i\in
\C\{y\}[z,f_1,\dots,f_{j-1}]$, $\de_1>0$ and $\de_t<n_{t-1}$ for
$t=2,3,\dots,j+1$, and for each $i$, $U^{(w)}_i$ has a multiplicity
$\ge \Pi^{j+1}_{t=1}n_t-i\Pi^j_{t=1}n_t =(n_{j+1}-i)\Pi^j_{t=1}n_t$
at $0 \in \C^2$ because $T^{(w)}_{d_{j+2}-1}$ has a multiplicity
$\ge \Pi^{j+1}_{t=1}n_t$ at $0 \in \C^2$ by Property(2).

Now, apply the WDT to $f$ with a divisor $h_{w+1}$, and then it can
be proved by Theorem 15.2 that $f$ is written in the form
$$
f=h^{d_{j+2}}_{w+1}+\sum^{d_{j+2}-1}_{i=0} T^{(w+1)}_ih^i_{w+1},
\tag 15.5.9-2
$$
and then $T^{(w+1)}_i$ satisfies the same kind of properties as
$T^{(w)}_i$ does in Property(2), Property(4) and Property(5) in Case
[II] of the lemma. Thus, the proof of Step(1) is done. \ms

$\underline{\text{\bf The Proof of Step(2) for Case[II]}}$ Let
$T^{(1)}_{d_{j+2}-1}$ in $(h_1,f)$ be not zero. For the proof of
Step(2), it is enough to show by Step(1) that
$(f_{j+1},f)=(h_{\nu+1},f)$ with the following property:
$$\align
& \text{$T^{(\nu+1)}_{d_{j+2}-1}=0$ \quad for some integer \quad
$\nu\le \f{n_{j+1}+1}2$,} \tag 15.5.10-2 \\
& \text{$T^{(p)}_{d_{j+2}-1}\ne 0$ \quad  for \quad
$p=1,2,\dots,\nu$.}
\endalign$$

First of all, we have already proved by Step(1) that for given
$(h_p,f)$, $(h_{p+1},f)$ can be written in the form
$$ \text{\rm(15.5.11-2)}\qquad \qquad \quad
\cases h_{p+1} &=h_p+\f 1{d_{j+2}}T^{(p)}_{d_{j+2}-1}=f^{n_{j+1}}_j
+\sum^{n_{j+1}-2}_{i=0}
R^{(p+1)}_if^i_j, \\
f &=h_{p+1}^{d_{j+2}}+\sum^{d_{j+2}-1}_{i=0} T^{(p+1)}_ih^i_{p+1},
\endcases \qquad \qquad \quad
$$
satisfying the properties {\rm Property(1)}, {\rm Property(2)},
\dots, Property(5) in this sublemma except possibly for the above
equation\text{\rm(15.5.10-2)}.

In preparation for the proof of the equation in \text{\rm(15.5.10-2)}, consider
$y,z,f_1,\dots,f_j$ as independent complex $(j+2)$-variables at the
origin in $\C^{j+2}$ as we have seen in the proof of Step(1), if
necessary. For brevity of notation, if $T=T(y,z,f_1,\dots,f_j)$ is
not zero in $\C\{y\}[z,f_1,\dots, f_j]$ for $j\ge 0$, then define
$\text{deg}_jT$ by the degree of $T$ as a polynomial in $f_j$, when
$T$ is in $\C\{y,z,f_1,\dots,f_{j-1}\}[f_j]$ as a polynomial in
$f_j$ with coefficients in $\C\{y,z,f_1,\dots,f_{j-1}\}$. For such
notations, write $\tau_p=\text{deg}_j(T^{(p)}_{d_{j+2}-1})$ for all
$p\ge 1$, and so $0\le \tau_p\le n_{j+1}-2$ by Property(5).

First, in preparation for the proof of an equality in \text{\rm(15.5.10-2)}, it is
very interesting to show that if $T^{(w)}_{d_{j+2}-1}$ and
$T^{(w+1)}_{d_{j+2}-1}$ is not zero for any integer $w>0$, then the
following inequality is true:
$$\align
n_{j+1}-2\ge \tau_w
>\tau_{w+1}\ge 0 \quad \text{with} \quad \tau_w\ge \f{n_{j+1}-1}2, \tag 15.5.12-2
\endalign$$
because if $T^{(p)}_{d_{j+2}-1}$ is not zero for any $p\ge 1$ and an
inequality in (15.5.12-2) is true, then $\{\tau_p:p=1,2,\dots\}$ is an
infinite sequence which is strictly decreasing and bounded with
$n_{j+1}-2\ge \tau_p>\tau_{p+1}\ge 0$. This would be impossible,
which is equivalent to the fact that there is an integer
$\nu<n_{j+1}$ such that $T^{(\nu+1)}_{d_{j+2}-1}=0$.

Now, assuming that $T^{(w)}_{d_{j+2}-1}$ and $T^{(w+1)}_{d_{j+2}-1}$
is not zero for any integer $w>0$, then we are going to prove that
an inequality in (15.5.12-2) is true.

Since $h_{w+1}=h_w+\f 1{d_{j+2}}T^{(w)}_{d_{j+2}-1}$, then a
sequence in (15.5.11-2) implies that
$$\align
 \text{\rm(15.5.13-2)}\qquad \qquad \qquad f&=h^{d_{j+2}}_w
+\sum^{d_{j+2}-1}_{i=0} T^{(w)}_ih^i_w  \\
&=h^{d_{j+2}}_{w+1} +\sum^{d_{j+2}-1}_{i=0} T^{(w+1)}_ih^i_{w+1},\\
 \text{\rm(15.5.14-2)}\qquad \qquad \qquad  f&=\left(h_w+\f 1{d_{j+2}}T^{(w)}_{d_{j+2}-1}\right)^{d_{j+2}}
+\sum^{d_{j+2}-1}_{i=0} T^{(w+1)}_ih^i_{w+1}. \qquad \qquad \qquad \qquad \qquad
\endalign
$$

\quad Therefore, \qquad
$f-h^{d_{j+2}}_w-T^{(w)}_{d_{j+2}-1}h^{d_{j+2}-1}_w =\sum^{(w)}_0
=\sum_{1,1}^{(w)}+\sum_{1,2}^{(w)}$ {\quad} such that
$$\align
 \sum^{(w)}_0 =\sum^{d_{j+1}-2}_{i=0}
T^{(w)}_ih^i_w \quad \text{by (15.5.13-2)}, \qquad \qquad \quad \tag 15.5.15-2 \\
\endalign$$
$$
\text{\rm(15.5.16-2)} \qquad \qquad \qquad  \cases \sum^{(w)}_{1,1}
&=\sum^{d_{j+2}-2}_{q=0}
c_q(T^{(w)}_{d_{j+2}-1})^{d_{j+2}-q}h^q_w \quad \text{by (15.5.14-2)}, \\
\sum^{(w)}_{1,2} &=\sum^{d_{j+2}-1}_{i=0} T^{(w+1)}_ih^i_{w+1},
\endcases \qquad \qquad \qquad \qquad
$$
with some nonzero constant $c_q$ by (15.5.14-2).

Then, to prove that $\tau_w>\tau_{w+1}\ge 0$ and $\tau_w\ge
\f{n_{j+1}-1}2$, note that given any nonzero monomial
$\Pi^{j+2}_{t=1}f^{\de_t}_{t-2}$ in either $T^{(w)}_i$ or
$T^{(w+1)}_i$ for all $i$, $\de_1>0$ and $\de_t<n_{t-1}$ for $2\le
t\le j+2$, and that $h_{w+1}\in\C\{y\}[z,f_1,\dots, f_j]$ is a
$W$-poly of degree $n_{j+1}$ in $f_j$.

Also, note by (15.4.12-2) that for each fixed $j\ge 1$ and for any
nonzero monomial $y^\a z^\beta$ in $\Pi^{j+1}_{t=1}f^{\g_t}_{t-2}$,
where $\g_1>0$ and $\g_t<n_{t-1}$ for $2\le t\le j+1$,
$$
\beta <\Pi^j_{t=1}n_t=\ \text{the multiplicity of $f_j$ at $0 \in
\C^2$.} \tag 15.5.17-2
$$
Now, to finish the proof of (15.5.12-2), let $\mu_1,\mu_2$ and $\mu_3$
be defined by $\text{$\max$}\{\g_{j+2}\}$ among the degree of $f_j$
for all nonzero monomials $\Pi^{j+2}_{t=1}f^{\g_t}_{t-2}$ in
$\sum^{(w)}_0$, $\sum^{(w)}_{1,1}$ and $\sum^{(w)}_{1,2}$,
respectively, where $\g_1>0$ and $\g_t<n_{t-1}$ for $2\le t\le j+1$.

Then, we need only to consider three subcases (i), (ii) and (iii),
respectively.

(i) In $\sum^{(w)}_0$,
$\mu_1<n_{j+1}+n_{j+1}(d_{j+2}-2)=d_{j+1}-n_{j+1}$ because
$d_{j+1}=n_{j+1}d_{j+2}$. \ms

(ii) In $\sum^{(w)}_{1,1}$, we have two possibilities by (15.5.17-2),
noting by Property(5) that $\tau_w
=\text{deg}_j(T^{(w)}_{d_{j+2}-1})\le n_{j+1}-2$.

(iia) If $\beta <\f 12\Pi^j_{t=1}n_t$ for all $\beta$ in (15.5.17-2),
then
$\mu_2=2\de_{j+2}+n_{j+1}(d_{j+2}-2)=2\tau_w+d_{j+1}-2n_{j+1}\le
2(n_{j+1}-2)+d_{j+1}-2n_{j+1}=d_{j+1}-4$ for some nonzero monomials
$\Pi^{j+2}_{t=1}f^{\de_t}_{t-2}$ in $T^{(w)}_{d_{j+2}-1}$. \ms

(iib) If $\beta \ge \f 12\Pi^j_{t=1}n_t$ for some $\beta$ in
(15.5.17-2), then
$\mu_2=2\de_{j+2}+1+n_{j+1}(d_{j+2}-2)=2\tau_w+1+d_{j+1}-2n_{j+1}\le
2(n_{j+1}-2)+1+d_{j+1}-2n_{j+1}=d_{j+1}-3$ for some nonzero
monomials $\Pi^{j+2}_{t=1}f^{\de_t}_{t-2}$ in $T^{(w)}_{d_{j+2}-1}$.
\ms

(iii) In $\sum^{(w)}_{1,2}$,
$\mu_3=\de'_{j+2}+n_{j+1}(d_{j+2}-1)=\tau_{w+1}+d_{j+1}-n_{j+1}\ge
d_{j+1}-n_{j+1}$ for some nonzero monomials
$\Pi^{j+2}_{t=1}f^{\de'_t}_{t-2}$ in $T^{(w+1)}_{d_{j+2}-1}$ where
$\tau_{w+1}=\text{deg}_j(T^{(w+1)}_{d_{j+2}-1})$.

By (i) and (iii), it is clear that $\mu_3>\mu_1$. \ms

Then, we claim that $\mu_2=\mu_3$. To prove the claim, assume the
contrary.

Then, it suffices to consider two possibilities (A) and (B),
respectively:

(A) If $\mu_3 >\mu_2$, then there is a nonzero monomial
$\left(\Pi^{j+2}_{t=1}f^{\de'_t}_{t-2}\right)f^{n_{j+1}(d_{j+2}-1)}_j$
which belongs to $\sum^{(w)}_{1,2}$, but does not belong to
$\sum^{(w)}_{1,1}$ and $\sum^{(w)}_0$ because $\mu_3>\mu_1$. It
would be a contradiction to $\sum^{(w)}_0
=\sum_{1,1}^{(w)}+\sum_{1,2}^{(w)}$.

(B) If $\mu_2>\mu_3$, then there is a nonzero monomial
$\left(\Pi^{j+1}_{t=1}f^{\de_t}_{t-2}\right)f^{\mu_2}_{j+2}$ which
belongs to $\sum^{(w)}_{1,1}$, but does not belong to $\sum^{(w)}_0$
and $\sum^{(w)}_{1,2}$ because $\mu_2>\mu_1$. It would be a
contradiction to $\sum^{(w)}_0 =\sum_{1,1}^{(w)}+\sum_{1,2}^{(w)}$.

Thus, we proved the claim. \ms

First, in preparation for the proof of an inequality in (15.5.12-2),
using (ii) and (iii) with $\mu_2=\mu_3$, then it suffices to
consider two possibilities (a) and (b), as follows:

(a) It is clear by (iia) and (iii) that $\mu_2=\mu_3$ implies
$2\tau_w=\tau_{w+1}+n_{j+1}$.

(b) It is clear by (iib) and (iii) that $\mu_2=\mu_3$ implies
$2\tau_w +1=\tau_{w+1}+n_{j+1}$. \ms

Recall by Property(5) of Case(II) in this sublemma that $\tau_p$ was
defined by $\tau_p=\text{deg}_j(T^{(p)}_{d_{j+2}-1})$ for all $p\ge
1$, and so $0\le \tau_p\le n_{j+1}-2$, if exists.

(a) First, if $2\tau_w=\tau_{w+1}+n_{j+1}$, then it is trivial that
$\tau_w-\tau_{w+1}=\f{\tau_{w+1}+n_{j+1}}{2}-\tau_{w+1}
=\f{n_{j+1}-\tau_{w+1}}{2}>0$ because $\tau_{w+1}\le n_{j+1}-2$.

(b) Next, if $2\tau_w+1=\tau_{w+1}+n_{j+1}$, then it is trivial that
$\tau_w-\tau_{w+1}=\f{\tau_{w+1}+n_{j+1}-1}{2}-\tau_{w+1}
=\f{n_{j+1}-\tau_{w+1}-1}{2}>0$ because $\tau_{w+1}\le n_{j+1}-2$.

Moreover, it is clear by (a) and (b) that $\tau_w\ge
\f{n_{j+1}-1}2$, which is an inequality in (15.5.12-2).

Therefore, we proved that $\tau_w>\tau_{w+1}\ge 0$ with $\tau_w\ge
\f{n_{j+1}-1}2$ for (15.5.12-2).

So, assuming that $T^{(p)}_{d_{j+2}-1}$ is not zero for any $p\ge
1$, since an inequality in (15.5.12-2) is true, then
$\{\tau_p:p=1,2,\dots\}$ is an infinite sequence which is strictly
decreasing and bounded with $n_{j+1}-2\ge \tau_p>\tau_{p+1}\ge 0$.
This would be impossible, and therefore there is an integer
$\nu<n_{j+1}$ such that $T^{(\nu+1)}_{d_{j+2}-1}=0$.

Finally, to prove that $\nu \le \f{n_{j+1}+1}2$, define a function
$\psi :\N\to \N\cup \{0\}$ by $\psi(p)=\tau_p$ where
$\tau_p=\text{deg}_j(T^{(p)}_{d_{j+2}-1})$ for all $p\ge 1$.

Then, $\nu-2\le (\psi(1)-\psi(2))+(\psi(2)-\psi(3))+\cdots
+(\psi(\nu-2)-\psi(\nu-1))=\psi(1)-\psi(\nu-1)
=\tau_1-\tau_{\nu-1}\le n_{j+1}-2-\f{n_{j+1}-1}2=\f{n_{j+1}-3}2$,
because $\psi(\nu-1)\ne 0$. So, $\nu\le \f{n_{j+1}+1}2$. Thus, we
proved Property(6), and so the proof of Step(2) is done.

Thus, we completed the proof of Sublemma 15.5. $\square$
\enddemo \ms

{\bf \S15.6. For the proof of Sublemma 15.6 } \ms

\proclaim{Sublemma 15.6} For a given $f$, the construction of a
sequence $\{f_k: k=1,2,\dots,j\}$ in the sense of Sublemma $15.5$ or
Theorem $15.4$ must be unique.
\endproclaim \ms

\noindent {\bf Proof of Sublemma 15.6.}
The proof will be by induction on
the integer $j$ with $0\le j\le \ell-1$.    In case $j=0$, there is nothing to prove.
In case $1\le j\le \ell-1$, if necessary, suppose we have shown by induction
assumption on the positive integer $j$ that such an existence of a sequence
$\{f_1,f_2,\dots,f_j\}$ in the sense of Sublemma $15.5$ is unique for a given $f$.
For the proof of the uniqueness on the positive integer $j$, it suffices to consider
the proofs for two cases, Case(I) and Case(II), respectively:

Case(I) $j=1\le {\ell-1}$ and  Case(II) $1\le j\le {\ell-1}$. \ms

\noindent$\underline{\text{\bf Case(I):}}$ Let $j=1\le {\ell-1}$.

$\underline{\text{\rm Assumptions for Case(I)}}$
Let $\phi_{-1}=\psi_{-1}=y$, $\phi_{0}=\psi_{0}=z$, and let
$$\cases
\phi_{1} &=\phi^{n_{1}}_{0}
+\sum^{n_{1}-2}_{i=0}R'_{1,i}\phi^i_{0},  \\
f
&=\phi^{d_{2}}_{1}+\sum^{d_{2}-2}_{p=0}S'_{2,p}\phi^p_{1}, \quad \text{and}
\endcases \tag 15.6.1 $$
$$\cases \psi_{1} &=\psi^{n_{1}}_{0}
+\sum^{n_{1}-2}_{i=0}R''_{1,i}\psi^i_{0}, \\
f
 &=\psi^{d_{2}}_{1}+\sum^{d_{2}-2}_{p=0}S''_{2,p}\psi^p_{1},
\endcases \tag 15.6.2
$$
where both \{$\phi_1$ and $f=f(\phi_{-1},\phi_0,\phi_1)$\}  in (15.6.1) and
\{$\psi_1$ and $f=f(\phi_{-1},\phi_0,\phi_1)$\} in (15.6.2) satisfy the same kind of
properties and notations as we have done in Assumptions and
Conclusions for Statement(1) for
Theorem 15.4  and Corollary 15.4.1 except possibly for the uniqueness. \ms

$\underline{\text{\rm Conclusions for Case(I)}}$
Then, $\phi_1(y,z)=\psi_1(y,z)$, and so  $f(\phi_{-1},\phi_0,\phi_1)=
f(\psi_{-1},\psi_0,\psi_1)$ because $S'_{2,p}=S''_{2,p}$
for each $p=0,1,\dots,d_2-2$, using the Weierstrass division theorem. \ms

\demo{$\underline{\text{\bf Proof of Case(I) for Proof of Sublemma 15.6}}$}
If $j=0$, the proof is trivial.
In preparation for the proof of the uniqueness on the integer
$(j=1)$, by (15.6.1) and (15.6.2), note that for any nonzero
monomial $y^s$ in $\phi^{\de'_1}_{-1}$ of
$R'_{1,i}=R'_{1,i}(\phi_{-1})$,
$$
s>0, \tag 15.6.3
$$
and for any nonzero monomial $y^\mu z^\nu$ in
$\Pi^{2}_{k=1}\phi^{\g'_k}_{k-2}$ of
$S'_{2,p}=S'_{2,p}(\phi_{-1},\phi_0)$,
$$
\nu <n_1=\ \text{the multiplicity of $\phi_{1}$ at
$0 \in \C^2$} \quad \text{and} \quad {\mu}>0. \tag 15.6.4
$$

{\noindent}Similarly, for any nonzero monomial $y^s$ in
$\psi^{\de''_1}_{-1}$ of $R''_{1,i}
=R''_{1,i}(\psi_{-1})$,
$$
 s>0, \tag 15.6.5
$$
and for any nonzero monomial $y^\mu z^\nu$ in $\Pi^{2}_{k=1}\psi^{\g''_k}_{k-2}$
of $S''_{2,p}=S''_{2,p}(\psi_{-1},\psi_0)$,
$$
\nu <n_1= \ \text{the multiplicity of  $\psi_{1}$
at $0 \in \C^2$} \quad \text{and} \quad {\mu}>0. \tag 15.6.6
$$

To prove the uniqueness on the integer $(j=1)$, it suffices to
follow two steps, Step(1) and Step(2). More rigorously, first we
will construct the detailed statement for Step(1) and the detailed
statement for Step(2), and next prove two statements, respectively.
\ms

$\underline{\text{\bf Step(1)}}$ If $j=0$, it is clear that the following is true:
$$
f(\phi_{-1},\phi_0)=f(\psi_{-1},\psi_0), \quad \text{because \
$\phi_{-1}=\psi_{-1}=y$ \ \text{and}\ $\phi_{0}=\psi_{0}=z$.}  \tag 15.6.7
$$

$\underline{\text{\bf Step(2)}}$ We prove by Step(1), (15.6.1) and
(15.6.2) that $\phi_{1}=\psi_{1}$, and by the uniqueness of The
WDT or Theorem $15.2$ that $f(\phi_{-1},\phi_0,\phi_1 )=f(\psi_{-1},\psi_0,\psi_1)$. \ms

In more detail, we write each step with proof, respectively.

$\underline{\text{\bf Step(1)}}$ \quad  To construct the statement
for Step(1), first it suffices to write two claims, i.e., Claim(I)
and Claim(II), and after then it remains to prove two claims for Step(1)
with (15.6.7). \ms

$\underline{\text{\rm Claim(I)}}$ \quad
$f(y,z)$ of (15.6.7) can be represented in the
form
$$\cases
\phi_0 &=z,\\
f &=\phi^{d_{1}}_{0}+\sum^{d_{1}-2}_{p=0}S'_{1,p}\phi^p_0,
\endcases \tag 15.6.8
$$
where $f$ of (15.6.8) is viewed as $f(\phi_{-1},\phi_0)$,
satisfying the following property:

{\rm(1)} Let $p$ be fixed with $0\le p\le d_{1}-2$. Then,
$S'_{1,p}=S'_{1,p}(\phi_{-1})\in \C\{y\}$ and $S'_{1,p}(0)=0$. \ms

{\rm(2)} Consider $y,z$ as independent
complex $2$-variables at the origin in $\C^{2}$. Let $p$ be
fixed with $0\le p\le d_{1}-2$. For any nonzero monomial
$y^{\de'_1}$ in $S'_{1,p}=S'_{1,p}(y)\in\C\{y\}$, $\de'_1>0$.  \ms

{\rm (3)} As in (iii) of Corollary 15.4.1, assuming that $f$ is a
$W$-poly of degree $n\ge 2$ in $z$ with the multiplicity $n$ at $0
\in \C^2$, then $S'_{1,p}\in \C\{y\}$ has a multiplicity $\ge
(d_{1}-p)$ at $0 \in \C^2$. \ms

$\underline{\text{\rm Claim(II)}}$ \quad
$f(y,z)$ of (15.6.2) can be represented in the
form
$$\cases
\psi_0 &=z,\\
f &=\psi^{d_{1}}_{0}+\sum^{d_{1}-2}_{p=0}S''_{1,p}\psi^p_0,
\endcases \tag 15.6.9
$$
where $f$ of (15.6.9) is viewed as $f(y,z)$,
satisfying the following property:

{\rm(1)} Let $p$ be fixed with $0\le p\le d_{1}-2$.
Then, $S''_{2,p}=S''_{2,p}(y)\in \C\{y\}$ and $S'_{2,p}(0)=0$. \ms

 {\rm(2)} Consider $y,z$ as independent
complex $2$-variables at the origin in $\C^{2}$. Let $i$ be
fixed with $0\le i\le d_{1}-2$. For any nonzero monomial
$y^{{\de}''_1}$ in $S''_{1,i}=S''_{1,i}(y)\in\C\{y\}$, ${\de_1}^{''}>0$. \ms

{\rm (3)} As in (iii) of Corollary 15.4.1, assuming that $f$ is a
$W$-poly of degree $n\ge 2$ in $z$ with the multiplicity $n$ at $0
\in \C^2$, then $S''_{1,p}\in \C\{y\}$ has a
multiplicity $\ge (d_{1}-p)$ at $0 \in \C^2$. \ms

$\underline{\text{\bf The Proof of Step(1)}}$ \quad
By assumptions of Theorem 15.4 and Corollary 15.4.1, $f=z^n+\sum^{n-2}_{i=0}
a_iy^{\a_i}z^i$ is a $W$-poly of degree $n\ge 2$ in $z$ with the multiplicity n
at $0\in \BC^{2}$ where for
$0\le i\le n-2$, each $a_i=a_i(y)$ is a unit in ${}_2\CO_0$. Then, $S'_{1,p}
=S''_{1,p}=a_py^{\a_p}$ for $p=0,1,\dots,d_1-2$, and so it is clear that
Claim(I) and Claim(II) are true. So, it can be easily
proved by (15.6.8) and (15.6.9) that the proof of Step(1) is true. \ms

$\underline{\text{\bf The Proof of Step(2)}}$ In preparation for the
proof of this step, using the equations in  (15.6.1) and (15.6.2) of Step(1),
then we get the following: Note that $n=n_1d_1$ with $n=d_1$.

Then, it is clear that $f$ can be rewritten in the form
$$\align
f &=\sum_1 +\sum_2\quad \text{for (15.6.1)}, \tag 15.6.10 \\
&=\sum_3+\sum_4 \quad \text{for (15.6.2)},
\endalign
$$
where
$$\cases
\sum_1
&=(\phi^{n_{1}}_{0}
+\sum^{n_{1}-2}_{i=0}R'_{1,i}\phi^i_{0})^{d_2}, \\
\sum_2 &=\sum^{d_{2}-2}_{p=0}S'_{2,p}
(\phi^{n_{1}}_{0}
+\sum^{n_{1}-2}_{i=0}R'_{1,i}\phi^i_{0})^p,
\endcases \tag 15.6.11
$$
$$\cases
\sum_3 &= (\psi^{n_{1}}_{0}
+\sum^{n_{1}-2}_{i=0}R''_{1,i}\psi^i_{0})^{d_{2}}, \\
\sum_4 &=\sum^{d_{2}-2}_{p=0}S''_{2,p}
(\psi^{n_{1}}_{0}
+\sum^{n_{1}-2}_{i=0}R''_{1,i}\psi^i_{0})^p.
\endcases \tag 15.6.12
$$

Now, consider $f-\phi^{n_{1}}_{0}=f-\psi^{n_{1}}_{0}$ as follows:
$$\align
f-\phi^{d_{1}}_{0} &=(\sum_1 - \phi^{d_{1}}_0)+\sum_2 \quad
\text{with (A) and (B)} \tag 15.6.13 \\
f-\psi^{d_{1}}_{0}&=(\sum_3 - \psi^{d_{1}}_{0})+\sum_4. \quad \text{with (C) and (D)}
\endalign$$

For the proof of this step, consider $y,z$ as
independent complex $2$-variables at the origin in $\C^{2}$.
Using the same method as we have seen in the proof of Fact(D) in
(15.4.3) of {\rm Sublemma 15.4.{$\alpha$}}, then $\sum_2$, $\sum_1
-\phi^{d_{1}}_0$, $\sum_4$ and $\sum_3-\psi^{d_{1}}_0$ which are
polynomials in $\C\{y\}[z]$, can be written by (A),
(B), (C) and (D) in order, as follows: \ms

$\underline{\text{\bf(A)}}$ Firstly, for any nonzero monomial
$y^rz^s\in \sum_2$ it is clear by (15.6.3) and (15.6.4) that $r>0$ and
$s<(d_{2}-2)n_{1}+n_1=d_{1}-n_1$.
$$\text{$r>0$ and
$0\le s< n_1(d_{2}-2)+n_{1}=d_{1}-n_1$.}$$

So, viewing $\sum_2$ as a polynomial in $\C\{y\}[z]$, by the same method
as we have seen in the proof of Fact(D)
in (15.4.3), $\sum_2$ can be rewritten as follows:
$$\align
(15.6.14) \qquad \quad  & \text{Whenever
$\Pi^{2}_{k=1}f^{\de_k}_{k-2}$ is in
$\sum_2$, \quad then $\delta_1>0$ and
$\delta_{2}<d_{1}-n_{1}$.} \qquad \qquad\\
\endalign$$

$\underline{\text{\bf(B)}}$ Secondly, assuming that $R'_{1,i}$ is
nonzero for some $i$, it is clear by (15.6.3) that for any nonzero
monomial $y^rz^s\in \sum_1-{\phi}^{d_{1}}_0$,
$$\text{$r>0$ and
$0\le s\le n_1(d_{2}-1)+n_{1}-2=d_{1}-2$.}$$

If $R'_{1,i}$ is nonzero for some $i$, there is a nonzero monomial
$y^{\alpha}z^{\beta}\in R'_{1,i}$, and so there is a nonzero
monomial $y^{r'}z^{s'}\in \sum_1 -{\phi}^{d_{1}}_0$ such that $r'>0$
and $s'\ge
(d_2-1)n_{1}=d_1-n_1$,
which implies that $y^{r'}z^{s'} \not \in\sum_2$.

So, viewing $\sum_1 -{\phi}^{d_{1}}_1$ as a polynomial in
$\C\{y\}[z]$, by the same method as we have seen in
the proof of Fact(D) in (15.4.3), if $R'_{1,i}$ is nonzero for some
$i$ then $\sum_1 -{\phi}^{d_{1}}_1$ can be rewritten as follows:
$$\align
(15.6.15) \qquad  & \text{Whenever $\Pi^{2}_{k=1}f^{\de_k}_{k-2}$
is in $\sum_1 -{\phi}^{d_{1}}_0$, then  $\delta_1>0$ \ and \ $\delta_{2}\le {d_{1}-2}$.} \\
& \text{Also, there is a monomial $\Pi^{2}_{k=1}f^{\de'_k}_{k-2}
\in \sum_1 -{\phi}^{d_{1}}_0$ \ with \ $\Pi^{2}_{k=1}f^{\de'_k}_{k-2}
\not \in \sum_2$ } \qquad \qquad \\
& \text{such that  $\delta'_1>0$
and $\delta'_{2}\ge d_{1}-n_{1}$.} \qquad \qquad \qquad
\endalign$$ \ms

$\underline{\text{\bf(C)}}$ Thirdly, for any nonzero monomial
$y^{r'}z^{s'}\in \sum_4$ it is clear by (15.6.5) and (15.6.6) that
$r>0$ and $s<(d_{2}-2)n_{1}+n_1=d_{1}-n_1$.
$$\text{$r>0$ and
$0\le s< n_1(d_{2}-2)+n_{1}=d_{1}-n_1$.}$$

So, viewing $\sum_4$ as a polynomial in $\C\{y\}[z]$, by the same method
as we have seen in the proof of Fact(D)
in (15.4.3), $\sum_4$ can be rewritten as follows:
$$\align
(15.6.16) \qquad \quad  & \text{Whenever
$\Pi^{2}_{k=1}f^{\de_k}_{k-2}$ is in
$\sum_4$, \quad then $\delta'_1>0$ and
$\delta'_{2}<d_{1}-n_{1}$.} \qquad \qquad\\
\endalign$$ \ms

$\underline{\text{\bf(D)}}$ Forthly, assuming that $R''_{1,i}$ is
nonzero for some $i$, it is clear by (15.6.5) that for any nonzero
monomial $y^rz^s\in \sum_3-{\psi}^{d_{1}}_0$,
$$\text{$r>0$ and
$0\le s\le n_1(d_{2}-1)+n_{1}-2=d_{1}-2$.}$$

If $R''_{1,i}$ is nonzero for some $i$, there is a nonzero monomial
$y^{\alpha}z^{\beta}\in R''_{1,i}$, and so there is a nonzero
monomial $y^{r''}z^{s''}\in \sum_3 -{\psi}^{d_{1}}_0$ such that $r''>0$
and $s''\ge
(d_2-1)n_{1}=d_1-n_1$,
which implies that $y^{r''}z^{s''} \not \in\sum_4$.

So, viewing $\sum_3 -{\psi}^{d_{1}}_0$ as a polynomial in
$\C\{y\}[z]$, by the same method as we have seen in
the proof of Fact(D) in (15.4.3), if $R''_{1,i}$ is nonzero for some
$i$ then $\sum_3 -{\psi}^{d_{1}}_0$ can be rewritten as follows:
$$\align
(15.6.17) \qquad  & \text{Whenever $\Pi^{2}_{k=1}f^{\de_k}_{k-2}$
is in $\sum_3 -{\psi}^{d_{1}}_0$, then  $\delta_1>0$ \ and \ $\delta_{2}\le {d_{1}-2}$.} \\
& \text{Also, there is a monomial $\Pi^{2}_{k=1}f^{\de''_k}_{k-2}
\in \sum_3 -{\psi}^{d_{1}}_0$ \ with \ $\Pi^{2}_{k=1}f^{\de''_k}_{k-2}
\not \in \sum_4$ } \qquad \qquad \\
& \text{such that  $\delta''_1>0$
and $\delta''_{2}\ge d_{1}-n_{1}$.} \qquad \qquad \qquad
\endalign$$ \ms

Recall by (15.6.1), (15.6.2) and Step(1) that
$$\align
\phi_{1} &=\phi^{n_{1}}_{0}
+\sum^{n_{1}-2}_{i=0}R'_{1,i}\phi^i_{0} \quad \text{and}  \tag 15.6.18 \\
\psi_{1}&=\psi^{n_{1}}_{0} +\sum^{n_{1}-2}_{i=0}R''_{1,i}\psi^i_0.
\endalign$$ \ms

In order to prove that $R'_{1,i}=R''_{1,i}$ for all $i=0,1,\dots,n_{1}-2$, first let
$$\align
(15.6.19) \qquad \qquad \qquad  m&=\text{$\max$}\{i:R'_{1,i}\
\text{is nonzero
in $\phi_{1}$ of (15.6.18)} \}, \quad \text{and}  \qquad \qquad \qquad \\
r &=\text{$\max$}\{i:R''_{1,i}\ \text{is nonzero in $\psi_{1}$
of (15.6.18)} \}.
\endalign$$

In order to prove that $m=r$ and $R'_{1,m}=R''_{1,r}$, it
suffices to consider two subcases: \ms

\noindent{\bf Subcases(1).} If there does not exist an integer $m$
satisfying (15.6.19), then $R'_{1,i}$ is identically zero for all
$i=0,1,\dots,n_{1}-2$, and so $R''_{1,i}$ is identically zero
for all $i=0,1,\dots,n_{1}-2$ by (15.6.13), (15.6.14),(15.6.15),
(15.6.16) and (15.6.17). Therefore, there is nothing to prove that
$\phi_{1}=\psi_{1}$. \ms

\noindent{\bf Subcase(2).} If there exist such integers $m\ge 0$ and
$r\ge 0$, $\sum_1 -{\phi}^{d_{1}}_0$ and $\sum_3-{\psi}^{d_{1}}_0$ can be
written as follows:

(a)(a1) $\sum_1-{\phi}^{d_{1}}_0=\sum a_{\g_1\g_2}
\Pi^{2}_{k=1}f^{\g_k}_{k-2}$ where each
$a_{\g_1\g_2}$ is a nonzero constant and $\g_1>0$ and
$\g_2\ge 0$.

(a2) Then, let $\tau$ be defined by $\text{$\max$}\{\g_{2}\}$ for
all nonzero monomials $\Pi^{2}_{k=1}f^{\g_k}_{k-2}$ in
$\sum_1-{\phi}^{d_{1}}_0$ where $\g_1>0$ and $\g_{2}\ge 0$. Then,
$\tau=n_{1}(d_{2}-1)+m$ by (15.6.15), (15.6.3) and (15.6.4). \ms

(b)(b1) $\sum_3-{\psi}^{d_{1}}_0=\sum b_{\de_1\de_2}
\Pi^{2}_{k=1}f^{\de_k}_{k-2}$ where each
$b_{\de_1\de_2}$ is a nonzero constant and $\de_1>0$ and
$\de_{2}\ge 0$.

(b2) Then, let $\omega$ be defined by $\text{Max}\{\de_{j+2}\}$ for
all nonzero monomials $\Pi^{j+2}_{k=1}f^{\de_k}_{k-2}$ in
$\sum_3-{\psi}^{d_{1}}_0$ where $\de_1>0$ and $\de_{j+2}\ge 0$.
Then, $\omega=n_{j+1}(d_{j+2}-1)+r$ by (15.6.17), (15.6.5) and (15.6.6). \ms

(c) Either $\tau=\omega$ or $m=r$ by (a2) and (b2). \ms

(d) To prove that $R'_{1,m}=R''_{1,r}$, first note that for any
nonzero monomial $\Pi^{2}_{k=1}f^{\g_k}_{k-2}\in
{\phi}^{{n_1}(d_{2}-1)+m}_0R'_{1,m}$, $\g_{2}=\tau={n_1}(d_{2}-1)+m$ by (15.6.3),
and that $\g'_{2}<\tau$ for any nonzero monomial
$\Pi^{2}_{k=1}f^{\g'_k}_{k-2}\in {\phi}^{n_{1}(d_{2}-1)+i}_0R'_{1,i}$
with $i<m$, if exists. By the similar result as in $\sum_3-{\psi}^{d_{1}}_0$, we
can get that
$$
{\phi}^{n_{1}(d_{2}-1)+m}_0R'_{1,m}={\psi}^{n_{1}(d_{2}-1)+m}_0
R''_{1,m}, \tag 15.6.20
$$
and so $R'_{1,m}=R''_{1,m}$. \ms
\enddemo \ms

If $m=r=0$, there is nothing to prove. To prove that
$R'_{1,i}=R''_{1,i}$ for $0\le i< m$, under the condition that
$m>0$, then it suffices to show that $R'_{1,i}=R''_{1,i}$ for
$0\le i<m$.

From (15.6.18) and (15.6.19) again, let $m_1$ and $r_1$ be defined by
$$\align
m_1 &=\text{Max}\{i:R'_{1,i}\ \text{is nonzero with}\ i<m\}, \tag 15.6.21 \\
r_1 &=\text{Max}\{i:R''_{1,i}\ \text{is nonzero with}\ i<r=m\}.
\endalign$$

For brevity of notation, let $R_{1,m}=R'_{1,m}=R''_{1,m}$ and $f_0={\phi}_0={\psi}_0$.
Noting that $\sum_1$,\dots,$\sum_4$ were already defined by (15.6.11)
and (15.6.12), then $\sum_{1,1}$ and $\sum_{1,3}$ can be defined by
$$\align
& \sum_{1,1}=\sum_1-({f}^{n_{1}}_0+R_{1,m}f^m_0)^{d_{2}} \quad
\text{and} \tag 15.6.22 \\
& \sum_{3,1} =\sum_3-({f}^{n_{1}}_0+R_{1,m}f^m_0)^{d_{2}}.
\endalign$$

By (15.6.10) and (15.6.22),
$f-(f^{n_{1}}_0+R_{1,m}f^m_0)^{d_{2}}$ can be rewritten as
follows:
$$\align
f-(f^{n_{1}}_0+R_{1,m}f^m_0)^{d_{2}}
=\sum_{1,1}+\sum_2=\sum_{3,1}+\sum_4. \tag 15.6.23
\endalign$$

In preparation for the proof of $m_1=r_1$ and
$R'_{1,m_1}=R''_{1,r_1}$, for brevity of notation,
$\sum_{1,1}$ and $\sum_{3,1}$ can be rewritten as follows:
$$\align
(15.6.24) \qquad \qquad
\sum_{1,1}&=\sum_1-(f^{n_{1}}_0+R_{1,m}f^m_0)^{d_{2}}
=\phi_{1}^{d_{2}}-(f^{n_{1}}_0+R_{1,m}f^m_0)^{d_{2}} \qquad \qquad\\
&=\{\sum^{m-1}_{i=0}R'_{1,i}f^i_0\}\cdot\{\sum^{d_{2}}_{i=1}
\phi_{1}^{d_{2}-i}(f^{n_{1}}_0+R_{1,m}f^m_0)^{i-1}\} \\
& \quad \text{where} \quad \sum_1 =\phi_{1}^{d_{2}} \quad
\text{with} \quad
\phi_{1}=(f_0^{n_{1}}+\sum^{n_{1}-2}_{i=0}R'_{1,i}f^i_0), \\
(15.6.25) \qquad \qquad
\sum_{3,1}&=\sum_3-(f^{n_{1}}_0+R_{1,m}f^m_0)^{d_{2}}
=\psi_{1}^{d_{2}}-(f^{n_{1}}_0+R_{1,m}f^m_0)^{d_{2}} \qquad \qquad\\
&=\{\sum^{m-1}_{i=0}R''_{1,i}f^i_0\}\cdot\{\sum^{d_{2}}_{i=1}
\psi_{1}^{d_{2}-i}(f^{n_{1}}_0+R_{1,m}f^m_0)^{i-1}\} \\
& \quad \text{where} \quad \sum_3 =\psi_{1}^{d_{2}} \quad
\text{with} \quad
\psi_{1}=(f_0^{n_{1}}+\sum^{n_{1}-2}_{i=0}R''_{1,i}f^i_0). \\
\endalign$$

For the proof, it suffices to consider two subcases:

\noindent$\underline{\text{\bf Subcases(1):}}$
If there does not exist an integer $m_1$
satisfying (15.6.21), then $R'_{1,i}$ is identically zero for all
$i=0,1,\dots,m-1$, and so $R''_{1,i}$ is identically zero for all
$i=0,1,\dots,m-1$, because otherwise $\sum_{3,1}$ satisfies the same
kind of properties as $\sum_3-f^{d_{1}}_j$ does in by (15.6.13),
(15.6.15) and (15.6.17). Therefore, there is nothing to prove that
$\phi_{1}=\psi_{1}$. \ms

\noindent$\underline{\text{\bf Subcases(2):}}$ If there exist such integers $m_1\ge 0$
and $r_1\ge 0$,
$\sum_{1,1}=\sum_1-(f^{n_{1}}_0+R_{1,m}f^m_0)^{d_{2}}$ and
$\sum_{3,1} =\sum_3-(f^{n_{1}}_0+R_{1,m}f^m_0)^{d_{2}}$ can be
written as follows:

(a) Let $\tau_1$ be defined by $\text{Max}\{\g'_{2}\}$ for all
nonzero monomials $\Pi^{2}_{k=1}f^{\g'_{k}}_{k-2}$ in
$\sum_{1,1}$, where $\g'_1>0$ and $\g'_{2}\ge 0$, Then,
$\tau_1=n_{1}(d_{2}-1)+m_1<\tau$ by (15.6.3) because
$\sum_{1,1}$ satisfies the same kind of properties
as $\sum_1 -f^{d_{1}}_0$ does in (15.6.15). \ms

(b) Let $\om_1$ be defined by $\text{Max}\{\g''_{2}\}$ for all
nonzero monomials $\Pi^{2}_{k=1}f^{\g''_{k}}_{k-2}$ in
$\sum_{3,1}$ where $\g''_1>0$ and $\g''_{2}\ge 0$. Then,
$\omega_1=n_{j+1}(d_{j+2}-1)+r_1<\omega$ by (15.6.5) because
$\sum_{3,1}$ satisfies the same kind of properties as $\sum_3
-f^{d_{1}}_0$ does in (15.6.17). \ms

(c) Either $\tau_1=\omega_1$ or $m_1=r_1$ by (a) and (b). \ms

(d) To prove that $R'_{1,m_1}=R''_{1,r_1}$, first note that for
any nonzero monomial $\Pi^{2}_{k=1}f^{\g_k}_{k-2}\in
f^{n_{1}(d_{2}-1)+m_1}_0R'_{1,m_1}$, $\g_{2}=\tau_1$ by
(15.6.3), and if $i<m_1$ then for any nonzero monomial
$\Pi^{2}_{k=1}f^{\g'_k}_{k-2}\in
f^{n_{1}(d_{2}-1)+i}_0R'_{1,i}$, $\g'_{2}<\tau_1$.

Note that whenever $\Pi^{2}_{k=1}f^{\g_k}_{k-2}\in
f^{n_{j+1}(d_{j+2}-1)+m_1}_0R'_{j+1,m_1}$ then
$\Pi^{2}_{k=1}f^{\g_k}_{k-2}\in \sum_{1,1}$.

By the similar result as in $\sum_{3,1}$, we can get that
$$
f^{n_{1}(d_{2}-1)+m_1}_0R'_{1,m_1}=f^{n_{1}(d_{2}-1)+m_1}_0
R''_{1,m_1}, \tag 15.6.26
$$
and so $R'_{1,m_1}=R''_{1,r_1}$ with $m_1=r_1$. \ms

If $m_1=r_1=0$, there is nothing to prove. To prove that
$R'_{1,i}=R''_{1,i}$ for $0\le i< m_1$, under the condition that
$m_1>0$, then it suffices to show that $R'_{1,i}=R''_{1,i}$ for
$0\le i<m_1$.

From (15.6.18) and (15.6.19) again, let $m_2$ and $r_2$ be defined by
$$\align
m_2 &=\text{Max}\{i:R'_{1,i}\ \text{is nonzero with}\ i<m_1\}, \tag 15.6.27 \\
r_2 &=\text{Max}\{i:R''_{1,i}\ \text{is nonzero with}\
i<r_1=m_1\}.
\endalign$$

For brevity of notation, let
$R_{1,m_1}=R'_{1,m_1}=R''_{1,r_1}$. Noting that
$\sum_1$,\dots,$\sum_4$ were already defined by (15.6.11) and
(15.6.12), then $\sum_{1,2}$ and $\sum_{3,2}$ can be defined by
$$\align
(15.6.28) \qquad \qquad \quad &
\sum_{1,2}=\sum_1-(f^{n_{1}}_0+R_{1,m}f^m_0+R_{1,m_1}f^{m_1}_0)^{d_{2}}
\quad
\text{and} \qquad \qquad \qquad \qquad \qquad\\
& \sum_{3,2}
=\sum_3-(f^{n_{1}}_0+R_{1,m}f^m_0+R_{1,m_1}f^{m_1}_0)^{d_{2}}.
\endalign$$

By (15.6.10) and (15.6.28),
$f-(f^{n_{1}}_0+R_{1,m}f^m_0+R_{1,m_1}f^{m_1}_0)^{d_{2}}$
can be rewritten as follows:
$$\align
(15.6.29) \qquad \qquad
f-(f^{n_{1}}_0+R_{1,m}f^m_0+R_{1,m_1}f^{m_1}_0)^{d_{2}}
=\sum_{1,2}+\sum_2=\sum_{3,2}+\sum_4. \qquad  \qquad \qquad
\endalign$$

By the same method as we have used in (15.6.24) and (15.6.25), it can
be easily shown that $R'_{1,m_2}=R''_{1,r_2}$ with $m_2=r_2$.
Thus, repeating the above process finitely many times, it can be
proved that $R'_{1,i}=R''_{1,i}$ for each
$i=0,1,\dots,n_{1}-2$, if exists. So, the proof of Step(2) is
done, and then we finished the proof for Case(I). $\square$
\enddemo \ms

\noindent$\underline{\text{\bf Case(II):}}$ Let $1\le j\le {\ell-1}$.

$\underline{\text{\rm Assumptions for Case(II)}}$
suppose we have shown by induction method that Sublemma 15.6 is true on the integer
$j\le {\ell-1}$.

Let $\phi_{-1}=\psi_{-1}=y$, $\phi_{0}=\psi_{0}=z$. In case $1\le j\le \ell-1$,
we may write by induction assumption on $\{j+1\}$ that for each
$k=1,2,\dots,j+1$,
$$\cases
\phi_{k} &=\phi^{n_{k}}_{k-1}
+\sum^{n_{k}-2}_{i=0}R'_{k,i}\phi^i_{k-1},  \\
f
&=\phi^{d_{j+2}}_{j+1}+\sum^{d_{j+2}-2}_{p=0}S'_{j+2,p}\phi^p_{j+1},
\endcases \tag 15.6.1-2  $$
and
$$\cases \psi_{k} &=\psi^{n_{k}}_{k-1}
+\sum^{n_{k}-2}_{i=0}R''_{k,i}\psi^i_{k-1}, \\
f
&=\psi^{d_{j+2}}_{j+1}+\sum^{d_{j+2}-2}_{p=0}S''_{j+2,p}\psi^p_{j+1},
\endcases \tag 15.6.2-2
$$
where each of $f=f(y,z,\phi_1,\dots,\phi_{j+1})$ and
$f=f(y,z,\psi_1,\dots,\psi_{j+1})$ satisfies the same kind of
properties and notations as we have done in the conclusion of
Theorem 15.4 and Corollary 15.4.1 except possibly for the uniqueness. \ms

$\underline{\text{\rm Conclusions for Case(II)}}$
Then, $\phi_i(y,z)=\psi_i(y,z)$ for each
$k=1,2,\dots,j+1$, and so  $f(\phi_{-1},\phi_0,\dots,\phi_{j+1})=f(\psi_{-1},\psi_0,\dots,\psi_{j+1})$ since $S'_{j+2,p}=S''_{j+2,p}$
for each $p=0,1,\dots,d_{j+2}-2$, using the Weierstrass division theorem. \bs

\noindent$\underline{\text{\bf Proof of Case(II) for Proof of Sublemma 15.6.}}$ \

The proof will be by induction on
the integer $j$ with $0\le j\le \ell-1$. If $j=0$, it is trivial.
So, for $1\le j\le \ell-1$, suppose we have shown by induction
assumption on $j$ that such an existence of a sequence
$\{f_1,f_2,\dots,f_j\}$ is unique for a given $f$.

In preparation for the proof of the uniqueness on the integer
$(j+1)$, by $f(\phi_{-1},\phi_0,\dots,\phi_{j+1})$ in (15.6.1-1) and by
$f(\psi_{-1},\psi_0,\dots,\psi_{j+1})$ in (15.6.2-2), note by the same method as in (15.4.7)
that for any nonzero
monomial $y^sz^t$ in $\Pi^{j+1}_{k=1}\phi^{\de'_k}_{k-2}$ of
$R'_{j+1,i}=R'_{j+1,i}(y,z,\phi_1,\dots,\phi_{j-1})$,
$$
t<\Pi^j_{k=1}n_k=\ \text{the multiplicity of $\phi_{j}$ at $0 \in
\C^2$} \quad \text{and} \quad s>0, \tag 15.6.3-2
$$
and for any nonzero monomial $y^\mu z^\nu$ in
$\Pi^{j+2}_{k=1}\phi^{\g'_k}_{k-2}$ of
$S'_{j+2,p}=S'_{j+2,p}(y,z,\phi_1,\dots,\phi_j)$,
$$
\nu <\Pi^{j+1}_{k=1}n_k=\ \text{the multiplicity of $\phi_{j+1}$ at
$0 \in \C^2$} \quad \text{and} \quad {\mu}>0. \tag 15.6.4-2
$$

{\noindent}Similarly, for any nonzero monomial $y^sz^t$ in
$\Pi^{j+1}_{k=1}\psi^{\de''_k}_{k-2}$ of $R''_{j+1,i}
=R''_{j+1,i}(y,z,\psi_1,\dots,\psi_{j-1})$,
$$
t<\Pi^j_{k=1}n_k=\ \text{the multiplicity of $\psi_j$ at $0 \in
\C^2$} \quad \text{and} \quad s>0, \tag 15.6.5-2
$$
and for any nonzero monomial $y^\mu z^\nu$ in
$\Pi^{j+2}_{k=1}\psi^{\g''_k}_{k-2}$ of
$S''_{j+2,p}=S''_{j+2,p}(y,z,\psi_1,\dots,\psi_j)$,
$$
\nu <\Pi^{j+1}_{k=1}n_k= \ \text{the multiplicity of  $\psi_{j+1}$
at $0 \in \C^2$} \quad \text{and} \quad {\mu}>0. \tag 15.6.6-2
$$

To prove the uniqueness on the integer $(j+1)$, it suffices to
follow two steps, Step(1) and Step(2). More rigorously, first we
will construct the detailed statement for Step(1) and the detailed
statement for Step(2), and next prove two statements, respectively.
\ms

$\underline{\text{\bf Step(1)}}$  We prove by induction assumption
on $j$ that
$$
f_k=\phi_k=\psi_k\quad \text{for}\quad 1\le k\le j, \tag 15.6.7-2
$$
and so $f(y,z,\phi_1,\dots,\phi_j)=f(y,z,\psi_1,\dots,\psi_j)$. \ms

$\underline{\text{\bf Step(2)}}$ We prove by Step(1), (15.6.1-2) and
(15.6.2-2) that $\phi_{j+1}=\psi_{j+1}$, and by the uniqueness of The
WDT or Theorem $15.2$ that there is nothing to prove for the
uniqueness. \ms

In more detail, we write each step with proof, respectively.

$\underline{\text{\bf Step(1)}}$ \quad  To construct the statement
for Step(1), first it suffices to write two claims, i.e., Claim(I)
and Claim(II), and after then it remains to prove two claims for Step(1)
with (15.6.7-2). \ms

$\underline{\text{\rm Claim(I)}}$ \quad
$f(y,z,\phi_1,\dots,\phi_{j})$ of (15.6.1-2) can be represented in the
form
$$\cases
\phi_j &=\phi^{n_j}_{j-1}+\sum^{n_j-2}_{i=0}R'_{j,i}\phi^i_{j-1},\\
f &=\phi^{d_{j+1}}_{j}+\sum^{d_{j+1}-2}_{i=0}S'_{j+1,i}\phi^i_j,
\endcases \tag 15.6.8-2
$$
where $f$ of (15.6.8-2) is viewed as $f(y,z,\phi_1,\dots,\phi_{j})$,
satisfying the following property:

{\rm(1)} Let $i$ be fixed with $0\le i\le d_{j+1}-2$. For each $j\ge
1$, $S'_{j+1,i}=S'_{j+1,i}(y,z)\in \C\{y\}[z]$ is a polynomial of
degree $<\Pi^j_{t=1}n_t$ in $z$ and $S'_{j+1,i}(0,z)=0$. \ms

{\rm(2)} Consider $y,z,\phi_1,\dots,\phi_{j-1}$ as independent
complex $(j+1)$-variables at the origin in $\C^{j+1}$. Let $i$ be
fixed with $0\le i\le d_{j+1}-2$. For any nonzero monomial
$\Pi^{j+1}_{t=1}\phi^{\de'_t}_{t-2}$ in
$S'_{j+1,i}=S'_{j+1,i}(y,z,\phi_1,\dots,\phi_{j-1})\in
\C\{y\}[z,\phi_1,\dots,\phi_{j-1}]$, $\de'_1>0$ and $\de'_t<n_{t-1}$
for $t=2,3,\dots, j+1$. \ms

{\rm (3)} As in (2) of Corollary 15.4.1, assuming that $f$ is a
$W$-poly of degree $n\ge 2$ in $z$ with the multiplicity $n$ at $0
\in \C^2$, then $S'_{j+1,i}\in \C\{y\}[z]$ has a multiplicity $\ge
(d_{j+1}-i)\Pi^j_{t=1}n_t$ at $0 \in \C^2$. \ms

$\underline{\text{\rm Claim(II)}}$ \quad
$f(y,z,\psi_1,\dots,\psi_{j})$ of (15.6.2-2) can be represented in the
form
$$\cases
\psi_j &=\psi^{n_j}_{j-1}+\sum^{n_j-2}_{i=0}R''_{j,i}\psi^i_{j-1},\\
f &=\psi^{d_{j+1}}_{j}+\sum^{d_{j+1}-2}_{i=0}S''_{j+1,i}\psi^i_j,
\endcases \tag 15.6.9-2
$$
where $f$ of (15.6.9-2) is viewed as $f(y,z,\psi_1,\dots,\psi_{j})$,
satisfying the following property:

{\rm(1)} Let $i$ be fixed with $0\le i\le d_{j+1}-2$. For each $j\ge
1$, $S''_{j+1,i}=S''_{j+1,i}(y,z)\in \C\{y\}[z]$ is a polynomial of
degree $<\Pi^j_{t=1}n_t$ in $z$ with $S''_{j+1,i}(0,z)=0$. \ms

{\rm(2)} Consider $y,z,\psi_1,\dots,\psi_{j-1}$ as independent
complex $(j+1)$-variables at the origin in $\C^{j+1}$. Let $i$ be
fixed with $0\le i\le d_{j+1}-2$. For any nonzero monomial
$\Pi^{j+1}_{t=1}\psi^{\de''_t}_{t-2}$ in
$S''_{j+1,i}=S''_{j+1,i}(y,z,\psi_1,\dots,\psi_{j-1})\in
\C\{y\}[z,\psi_1,\dots,\psi_{j-1}]$, $\de''_1>0$ and
$\de''_t<n_{t-1}$ for $t=2,3,\dots, j+1$. \ms

{\rm (3)} As in (iii) of Corollary 15.4.1, assuming that $f$ is a
$W$-poly of degree $n\ge 2$ in $z$ with the multiplicity $n$ at $0
\in \C^2$, then $S''_{j+1,i}\in \C\{y\}[z]$ of (15.4.1) has a
multiplicity $\ge (d_{j+1}-i)\Pi^j_{t=1}n_t$ at $0 \in \C^2$. \ms

If the proofs of Claim(I) and Claim(II) are done, then it can be
proved by (15.6.8-2) and (15.6.9-2) and by induction assumption on $j$
that $f_k=\phi_k=\psi_k$ for $1\le k\le j$ in (15.6.7-2), and so
$f(y,z,\phi_1,\dots,\phi_j)=f(y,z,\psi_1,\dots,\psi_j)$. \ms

$\underline{\text{\bf The Proof of Step(1)}}$ First, in order to
prove Claim(I), using the equations for $k=j$ and $k=j+1$ in
(15.6.1-2), then we get the following:
$$\align
f&=\phi^{d_{j+2}}_{j+1}+\sum^{d_{j+2}-2}_{p=0}S'_{j+2,p}\phi^p_{j+1}
\tag Eq.1\\
&=\sum_1 +\sum_2 ,  \tag Eq.2 \\
where \qquad \sum_1 &=\{\phi^{n_{j+1}}_j+\sum^{n_{j+1}-2}_{i=0}
R'_{j+1,i}\phi^i_{j}\}^{d_{j+2}} \quad \text{and}\\
\sum_2 &=\sum^{d_{j+2}-2}_{p=0}S'_{j+2,p}\{\phi^{n_{j+1}}_j+
\sum^{n_{j+1}-2}_{i=0}R'_{j+1,i}\phi^i_{j}\}^{p}.\\
\endalign$$

Now, consider $f-\phi^{d_{j+1}}_j$ as
$$\align
f-\phi^{d_{j+1}}_j &=(\sum_1 -\phi^{d_{j+1}}_j)+\sum_2. \tag Eq.3 \\
\endalign$$

Recall that $\phi_j(y,z)\in \C\{y\}[z]$ is a $W$-poly of degree
$\Pi^j_{t=1}n_t$ in $z$ with the multiplicity $\Pi^j_{t=1}n_t$ at
$0\in \BC^2$ by the induction assumption on $j$, and that
$R'_{j+1,i}(y,z)$ is a polynomial of degree $<\Pi^{j}_{t=1}n_t$ in
$z$ and $R'_{j+1,i}(0,z)=0$  by (15.6.3-2) and that $S'_{j+2,i}(y,z)$
is a polynomial of degree $<\Pi^{j+1}_{t=1}n_t$ in $z$ and
$S'_{j+2,i}(0,z)=0$ by (15.6.4-2).

In preparation for the proof of Claim(I), first of all, we will
prove the following:

(i) If we write $\sum_1 -\phi^{d_{j+1}}_j=\sum b_{p,q}y^pz^q$ with
some nonzero constant $b_{p,q}$, then $p>0$ and
$q<(d_{j+1}-1)\Pi^{j}_{t=1}n_t$.

(ii) If we write $\sum_2=\sum c_{r,s}y^rz^s$ with some nonzero
constant $c_{r,s}$, then $r>0$ and
$s<(d_{j+1}-n_{j+1})\Pi^{j}_{t=1}n_t$. \ms

For any nonzero monomial $y^pz^q\in \sum_1 -\phi^{d_{j+1}}_j$, it is
clear that $p>0$ and by (15.6.3-2) that
$q<(n_{j+1}(d_{j+2}-1)+n_{j+1}-2)\Pi^{j}_{t=1}n_t+\Pi^{j}_{t=1}n_t
=(d_{j+1}-1)\Pi^{j}_{t=1}n_t$. Thus, the proof of (i) is done.

For any nonzero monomial $y^rz^s\in \sum_2$, it is clear that $r>0$
and by (15.6.3-2) and (15.6.4-2) that
$s<(d_{j+2}-2)n_{j+1}\Pi^{j}_{t=1}n_t+\Pi^{j+1}_{t=1}n_t
=(d_{j+1}-2n_{j+1}+n_{j+1})\Pi^{j}_{t=1}n_t=(d_{j+1}-n_{j+1})\Pi^{j}_{t=1}n_t$.
Thus, the proof of (ii) is done.

Therefore, we proved that whenever any nonzero monomial
$y^{\alpha}z^{\beta}\in f-\phi^{d_{j+1}}_j$ then $\alpha>0$ and
$\beta<(d_{j+1}-1)\Pi^{j}_{t=1}n_t$, noting that $\phi_j(y,z)\in
\C\{y\}[z]$ is a $W$-poly of degree $\Pi^{j}_{t=1}n_t$ in $z$ with
the multiplicity $\Pi^{j}_{t=1}n_t$ at $0\in \BC^2$ by the induction
assumption on $j$.

Now, in order to prove Claim(I), apply the WDT with a divisor
$\phi_{j}$ to $f$ as an element of $\C\{y\}[z]$. Since for any
nonzero monomial $y^{\alpha}z^{\beta}\in f-\phi^{d_{j+1}}_j$
$\alpha>0$ and $\beta<(d_{j+1}-1)\Pi^{j}_{t=1}n_t$ by (i) and (ii),
then it is clear by (i) and (ii) of Theorem 15.2 that $f$ of (Eq.2)
can be rewritten as follows:
$$
f =\phi_{j}^{d_{j+1}}+\sum^{d_{j+1}-1}_{i=0} T^{(j+1)}_i\phi^i_{j}
\quad \text{with $T^{(j+1)}_{d_{j+1}-1}=0$}, \tag Eq.4
$$
where for each $i=0,1,\dots, d_{j+1}-2$,
$T^{(j+1)}_{i}=T^{(j+1)}_{i}(y,z)=\sum a_{p,q}y^pz^q$ with a nonzero
constant $a_{p,q}$ such that $p>0$ and $q<\Pi^{j}_{t=1}n_t$ and that
$T^{(j+1)}_{i}(0,z)=0$ and $T^{(j+1)}_{i}$ has a multiplicity $\ge
(d_{j+1}-i)\Pi^{j+1}_{t=1}n_t$ at $0 \in \C^2$. \ms

Now, to finish the proof of Claim(I), since $\phi_{j-1}\in
\C\{y\}[z]$ is a $W$-poly of degree $\Pi^{j-1}_{t=1}n_t$ in $z$ with
the multiplicity $\Pi^{j-1}_{t=1}n_t$ at $0\in \BC^2$, apply the WDT
with a divisor $\phi_{j-1}$ to $T^{(j+1)}_{i}$ for each
$i=0,1,\dots, d_{j+2}-2$. Then by Theorem $15.2$, each
$T^{(j+1)}_{i}$ can be written as
$$
T^{(j+1)}_{i}=\sum^{n_{j}-1}_{k_1=0} Q_{k_1}\phi^{k_1}_{j-1}, \tag
Eq.5
$$
where if exists, each $Q_{k_1}\in \C\{y\}[z]$ is a polynomial of
degree $<\Pi^{j-1}_{k=1}n_t$ in $z$ with $Q_{k_1}(0,z)=0$. Using the
similar technique as we have seen in the proof of Fact(D) in
(15.4.3), since $\phi_{j-2}\in \C\{y\}[z]$ is a $W$-poly of degree
$\Pi^{j-2}_{t=1}n_t$ in $z$ with the multiplicity
$\Pi^{j-2}_{t=1}n_t$ at $0\in \BC^2$, then for each fixed $k_1$
apply the WDT with a divisor $\phi_{j-2}$ to $Q_{k_1}$ for each
$k_1=0,1,\dots, n_{j}-1$, again. Then by Theorem 15.2, each
$Q_{k_1}$ may be written in the form
$$
Q_{k_1}=\sum^{n_{j-1}-1}_{k_2=0} Q_{k_1,k_2}\phi^{k_2}_{j-2}, \tag
Eq.6
$$
where if exists, each $Q_{k_1,k_2}\in \C\{y\}[z]$ is a polynomial of
degree $<\Pi^{j-2}_{t=1}n_t$ in $z$ with $Q_{k_1k_2}(0,z)$ zero. \ms

Thus, continue the same process as above, and consider
$y,z,\phi_1,\dots,\phi_{j-1}$ as independent complex
$(j+1)$-variables at the origin in $\C^{j+1}$. Let $i$ be fixed with
$0\le i\le d_{j+1}-2$. Then, we proved that for each fixed
$i=0,1,\dots,d_{j+1}-2$, and for any nonzero monomial
$\Pi^{j+1}_{t=1}\phi^{\de_t}_{t-2}$ in
$T_{j+1,i}=T_{j+1,i}(y,z,\phi_1,\dots,\phi_{j-1})\in
\C\{y\}[z,\phi_1,\dots,\phi_{j-1}]$, $\de_1>0$ and $\de_t<n_{t-1}$
for $t=2,3,\dots, j+1$. If we write $S'_{j+1,i}=T^{(j+1)}_i$ for
$0\le i\le d_{j+1}-2$, then the proof of Claim(I) with (15.6.8-2) is
done.

By the same method as we have used in the proof of Claim(I), the
proof of Claim(II) with (15.6.9-2) can be proved, too.

Therefore, by induction assumption on $j$ and (15.6.8-2) and (15.6.9-2),
we proved that
$$
f_k=\phi_k=\psi_k\quad \text{for}\quad 1\le k\le j, \tag Eq.7
$$
and then $f(y,z,\phi_1,\dots,\phi_j)=f(y,z,\psi_1,\dots,\psi_j)$.
Thus, the proof of Step(1) is done. \ms

$\underline{\text{\bf Step(2)}}$ For the proof, it suffices to show
that $R'_{j+1,i}=R''_{j+1,i}$ for $0\le i\le n_{j+1}-2$, because it
is clear by the uniqueness of The WDT or Theorem $15.2$ that
$\phi_{j+1}=\psi_{j+1}$ if and only if $R'_{j+1,i}=R''_{j+1,i}$ for
$0\le i\le n_{j+1}-2$ where $\phi_{j+1} =f^{n_{j}}_j
+\sum^{n_{j+1}-2}_{i=0} R'_{j+1,i}f^i_j$ by (15.6.1-2) and $\psi_{j+1}
=f^{n_{j}}_{j} +\sum^{n_{j+1}-2}_{i=0}R''_{j+1,i}f^i_j$ by (15.6.2-2),
noting by Step(1) that $f_j=\phi_j=\psi_j$. \ms

$\underline{\text{\bf The Proof of Step(2)}}$ In preparation for the
proof of this step, using the equation in either (15.6.7-2)
or (Eq.7) of Step(1), then apply (Eq.7) to (15.6.1-2) with $k=j+1$, and also
apply (Eq.7) to (15.6.2-2) with $k=j+1$. Then, it is clear that $f$
can be rewritten in the form
$$\align
f &=\sum_1 +\sum_2\quad \text{for (15.6.1-2)} \tag 15.6.10-2 \\
&=\sum_3+\sum_4 \quad \text{for (15.6.2-2)},
\endalign
$$
where
$$\cases
\sum_1
&=(f_j^{n_{j+1}}+\sum^{n_{j+1}-2}_{i=0}R'_{j+1,i}f^i_j)^{d_{j+2}}, \\
\sum_2 &=\sum^{d_{j+2}-2}_{p=0}S'_{j+2,p}
(f^{n_{j+1}}_j+\sum^{n_{j+1}-2}_{i=0}R'_{j+1,i}f^i_j)^p,
\endcases \tag 15.6.11-2
$$
$$\cases
\sum_3 &= (f^{n_{j+1}}_j
+\sum^{n_{j+1}-2}_{i=0}R''_{j+1,i}f^i_j)^{d_{j+2}}, \\
\sum_4 &=\sum^{d_{j+2}-2}_{p=0}S''_{j+2,p}
(f^{n_{j+1}}_j+\sum^{n_{j+1}-2}_{i=0}R''_{j+1,i}f^i_j)^p.
\endcases \tag 15.6.12-2
$$

Now, consider $f-f^{d_{j+1}}_j$ as
$$\align
f-f^{d_{j+1}}_j &=(\sum_1 - f^{d_{j+1}}_j)+\sum_2 \tag 15.6.13-2 \\
&=(\sum_3 - f^{d_{j+1}}_j)+\sum_4.
\endalign$$

For the proof of this step, consider $y,z,f_1,\dots,f_{j}$ as
independent complex $(j+2)$-variables at the origin in $\C^{j+2}$.
Using the same method as we have seen in the proof of Fact(D) in
(15.4.3) of {\rm Sublemma 15.4.{$\alpha$}}, then $\sum_2$, $\sum_1
-f^{d_{j+1}}_j$, $\sum_4$ and $\sum_3-f^{d_{j+1}}_j$ which are
polynomials in $\C\{y\}[z,f_1,\dots, f_{j}]$, can be written by (A),
(B), (C) and (D) in order, as follows: \ms

$\underline{\text{\bf(A)}}$ Firstly, for any nonzero monomial
$y^rz^s\in \sum_2$ it is clear by (15.6.3-2) and (15.6.4-2) that $r>0$ and
$s<(d_{j+2}-2)n_{j+1}\Pi^{j}_{t=1}n_t+\Pi^{j+1}_{t=1}n_t
=(d_{j+1}-2n_{j+1}+n_{j+1})\Pi^{j}_{t=1}n_t=(d_{j+1}-n_{j+1})\Pi^{j}_{t=1}n_t$.

So, considering $\sum_2$ as a polynomial in $\C\{y\}[z,f_1,\dots,
f_{j}]$, by the same method as we have seen in the proof of Fact(D)
in (15.4.3), $\sum_2$ can be rewritten as follows:
$$\align
\text{\rm(15.6.14-2)} \qquad \qquad \qquad  & \text{Whenever
$\Pi^{j+2}_{k=1}f^{\de_k}_{k-2}$ is in
$\sum_2$, then  $\delta_1>0$,} \qquad \qquad \qquad \qquad \qquad\\
& \text{$\delta_t<n_t$ for $2\le t\le j+1$ and
$\delta_{j+2}<d_{j+1}-n_{j+1}$.}
\endalign$$

$\underline{\text{\bf(B)}}$ Secondly, assuming that $R'_{j+1,i}$ is
nonzero for some $i$, for any nonzero monomial $y^rz^s\in \sum_1
-f^{d_{j+1}}_j$ it is clear by (15.6.3-2) that $r>0$ and
$s<\{(d_{j+2}-1)n_{j+1}+n_{j+1}-2\}\Pi^{j}_{t=1}n_t+\Pi^{j}_{t=1}n_t
=(d_{j+1}-1)\Pi^{j}_{t=1}n_t$.

If $R'_{j+1,i}$ is nonzero for some $i$, there is a nonzero monomial
$y^{\alpha}z^{\beta}\in R'_{j+1,i}$, and so there is a nonzero
monomial $y^{r'}z^{s'}\in \sum_1 -f^{d_{j+1}}_j$ such that $r'>0$
and $s'\ge
(d_{j+2}-1)n_{j+1}\Pi^{j}_{t=1}n_t=(d_{j+1}-n_{j+1})\Pi^{j}_{t=1}n_t$,
which implies that $y^{r'}z^{s'} \not \in\sum_2$.

So, considering $\sum_1 -f^{d_{j+1}}_j$ as a polynomial in
$\C\{y\}[z,f_1,\dots, f_{j}]$, by the same method as we have seen in
the proof of Fact(D) in (15.4.3), if $R'_{j+1,i}$ is nonzero for some
$i$ then $\sum_1 -f^{d_{j+1}}_j$ can be rewritten as follows:
$$\align
\text{\rm(15.6.15-2)} \qquad  & \text{Whenever $\Pi^{j+2}_{k=1}f^{\de_k}_{k-2}$
is in $\sum_1 -f^{d_{j+1}}_j$, then  $\delta_1>0$,} \\
& \text{$\delta_t<n_t$ for $2\le t\le j+1$ and
$\delta_{j+2}<d_{j+1}-1$.} \qquad \qquad \qquad \\
& \text{Also, there is a monomial $\Pi^{j+2}_{k=1}f^{\de'_k}_{k-2}
\in \sum_1 -f^{d_{j+1}}_j$ with $\Pi^{j+2}_{k=1}f^{\de'_k}_{k-2}
\not \in \sum_2$ } \\
& \text{such that  $\delta'_1>0$, $\delta'_t<n_t$ for $2\le t\le
j+1$ and $\delta'_{j+2}\ge d_{j+1}-n_{j+1}$.} \qquad \qquad \qquad
\endalign$$ \ms

$\underline{\text{\bf(C)}}$ Thirdly, for any nonzero monomial
$y^{r'}z^{s'}\in \sum_4$ it is clear by (15.6.5-2) and (15.6.6-2) that
$r>0$ and $s<(d_{j+2}-2)n_{j+1}\Pi^{j}_{t=1}n_t+\Pi^{j+1}_{t=1}n_t
=(d_{j+1}-2n_{j+1}+n_{j+1})\Pi^{j}_{t=1}n_t=(d_{j+1}-n_{j+1})\Pi^{j}_{t=1}n_t$.

So, considering $\sum_4$ as a polynomial in $\C\{y\}[z,f_1,\dots,
f_{j}]$, by the same method as we have seen in the proof of Fact(D)
in (15.4.3), $\sum_4$ can be rewritten as follows:
$$\align
\text{\rm(15.6.16-2)} \qquad \qquad \qquad
& \text{Whenever  $\Pi^{j+2}_{k=1}f^{\de'_k}_{k-2}$ is in
$\sum_4$, then  $\delta'_1>0$,}  \qquad \qquad \qquad \qquad \qquad\\
& \text{$\delta'_t<n_t$ for $2\le t\le j+1$ and
$\delta'_{j+2}<d_{j+1}-n_{j+1}$.}
\endalign$$ \ms

$\underline{\text{\bf(D)}}$ Forthly, assuming that $R''_{j+1,i}$ is
nonzero for some $i$, for any nonzero monomial $y^rz^s\in \sum_3
-f^{d_{j+1}}_j$ it is clear by (15.6.5-2) that $r>0$ and
$s<\{(d_{j+2}-1)n_{j+1}+n_{j+1}-2\}\Pi^{j}_{t=1}n_t+\Pi^{j}_{t=1}n_t
=(d_{j+1}-1)\Pi^{j}_{t=1}n_t$.

If $R''_{j+1,i}$ is nonzero for some $i$, there is a nonzero
monomial $y^{\alpha}z^{\beta}\in R''_{j+1,i}$, and so there is a
nonzero monomial $y^{r''}z^{s''}\in \sum_3 -f^{d_{j+1}}_j$ such that
$r''>0$ and $s''\ge
(d_{j+2}-1)n_{j+1}\Pi^{j}_{t=1}n_t=(d_{j+1}-n_{j+1})\Pi^{j}_{t=1}n_t$,
which implies that $y^{r''}z^{s''} \not \in\sum_2$.

So, considering $\sum_1 -f^{d_{j+1}}_j$ as a polynomial in
$\C\{y\}[z,f_1,\dots, f_{j}]$, by the same method as we have seen in
the proof of Fact(D) in (15.4.3), if $R''_{j+1,i}$ is nonzero for
some $i$ then $\sum_3 -f^{d_{j+1}}_j$ can be rewritten as follows:
$$\align
\text{\rm(15.6.17-2)} \qquad  & \text{Whenever $\Pi^{j+2}_{k=1}f^{\de_k}_{k-2}$
is in $\sum_3 -f^{d_{j+1}}_j$, then  $\delta_1>0$,} \\
& \text{$\delta_t<n_t$ for $2\le t\le j+1$ and
$\delta_{j+2}<d_{j+1}-1$.} \qquad \qquad \qquad \\
& \text{Also, there is a monomial $\Pi^{j+2}_{k=1}f^{\de''_k}_{k-2}
\in \sum_3 -f^{d_{j+1}}_j$ with $\Pi^{j+2}_{k=1}f^{\de''_k}_{k-2}
\not \in \sum_4$ } \\
& \text{such that  $\delta''_1>0$, $\delta''_t<n_t$ for $2\le t\le
j+1$ and $\delta''_{j+2}\ge d_{j+1}-n_{j+1}$.} \qquad \qquad \qquad
\endalign$$

Recall by (15.6.1-1), (15.6.2-2) and Step(1) that
$$\align
 \phi_{j+1}&=f^{n_{j}}_j +\sum^{n_{j+1}-2}_{i=0}
R'_{j+1,i}f^i_j  \quad \text{and}  \tag 15.6.18-2 \\
\psi_{j+1}&=f^{n_{j}}_{j} +\sum^{n_{j+1}-2}_{i=0}R''_{j+1,i}f^i_j.
\endalign$$ \ms

In order to prove that $R'_{j+1,i}=R''_{j+1,i}$ for $0\le i\le
n_{j+1}-2$, first let
$$\align
\text{\rm(15.6.19-2)} \qquad \qquad   m &=\text{$\max$}\{i:R'_{j+1,i}\
\text{is nonzero
in $\phi_{j+1}$ of (15.6.18-2)} \}, \quad \text{and}  \qquad \qquad \qquad \\
r &=\text{$\max$}\{i:R''_{j+1,i}\ \text{is nonzero in $\psi_{j+1}$
of (15.6.18-2)} \}.
\endalign$$

In order to prove that $m=r$ and $R'_{j+1,m}=R''_{j+1,r}$, it
suffices to consider two subcases: \ms

\noindent$\underline{\text{\bf Subcases(1):}}$ If there does not exist an integer $m$
satisfying (15.6.19-2), then $R'_{j+1,i}$ is identically zero for all
$i=0,1,\dots,{n_{j+1}-2}$, and so $R''_{j+1,i}$ is identically zero
for all $i=0,1,\dots,n_{j+1}-2$ by (15.6.13-2), (15.6.14-2),(15.6.15-2),
(15.6.16-2) and (15.6.17-2). Therefore, there is nothing to prove that
$\phi_{j+1}=\psi_{j+1}$. \ms

\noindent$\underline{\text{\bf Subcases(2):}}$ If there exist such integers $m\ge 0$ and
$r\ge 0$, $\sum_1-f^{d_{j+1}}_j$ and $\sum_3-f^{d_{j+1}}_j$ can be
written as follows:

(a)(a1) $\sum_1-f^{d_{j+1}}_j=\sum a_{\g_1\g_2\cdots
\g_{j+2}}\Pi^{j+2}_{k=1}f^{\g_k}_{k-2}$ where each
$a_{\g_1\g_2\cdots\g_{j+2}}$ is a nonzero constant and $\g_1>0$,
$\g_k<n_{k-1}$ for $2\le k\le j+1$, and $\g_{j+2}\ge 0$.

(a2) Then, let $\tau$ be defined by $\text{$\max$}\{\g_{j+2}\}$ for
all nonzero monomials $\Pi^{j+2}_{k=1}f^{\g_k}_{k-2}$ in
$\sum_1-f^{d_{j+1}}_j$ where $\g_1>0$, $\g_k <n_{k-1}$ for $2\le
k\le j+1$, and $\g_{j+2}\ge 0$. Then, $\tau=n_{j+1}(d_{j+2}-1)+m$ by
(15.6.15-2), (15.6.3-2) and (15.6.4-2). \ms

(b)(b1) $\sum_3-f^{d_{j+1}}_j=\sum b_{\de_1\de_2\cdots
\de_{j+2}}\Pi^{j+2}_{k=1}f^{\de_k}_{k-2}$ where each
$b_{\de_1\de_2\cdots\de_{j+2}}$ is a nonzero constant and $\de_1>0$,
$\de_k<n_{k-1}$ for $2\le k\le j+1$, and $\de_{j+2}\ge 0$.

(b2) Then, let $\omega$ be defined by $\text{Max}\{\de_{j+2}\}$ for
all nonzero monomials $\Pi^{j+2}_{k=1}f^{\de_k}_{k-2}$ in
$\sum_1-f^{d_{j+1}}_j$ where $\de_1>0$, $\de_k <n_{k-1}$ for $2\le
k\le j+1$, and $\de_{j+2}\ge 0$. Then, $\omega=n_{j+1}(d_{j+2}-1)+r$
by (15.6.17-2), (15.6.5-2) and (15.6.6-2). \ms

(c) Either $\tau=\omega$ or $m=r$ by (a2) and (b2). \ms

(d) To prove that $R'_{j+1,m}=R''_{j+1,r}$, first note that for any
nonzero monomial $\Pi^{j+2}_{k=1}f^{\g_k}_k\in
f^{n_{j+1}(d_{j+2}-1)+m}_jR'_{j+1,m}$, $\g_{j+2}=\tau$ by (15.6.3-2),
and that $\g'_{j+2}<\tau$ for any nonzero monomial
$\Pi^{j+2}_{k=1}f^{\g'_k}_k\in f^{n_{j+1}(d_{j+2}-1)+i}_jR'_{j+1,i}$
with $i<m$. By the similar result as in $\sum_3-f^{d_{j+1}}_j$, we
can get that
$$
f^{n_{j+1}(d_{j+2}-1)+m}_jR'_{j+1,m}=f^{n_{j+1}(d_{j+2}-1)+m}_j
R''_{j+1,m}, \tag 15.6.20-2
$$
and so $R'_{j+1,m}=R''_{j+1,m}$. \ms

If $m=r=0$, there is nothing to prove. To prove that
$R'_{j+1,i}=R''_{j+1,i}$ for $0\le i< m$, under the condition that
$m>0$, then it suffices to show that $R'_{j+1,i}=R''_{j+1,i}$ for
$0\le i<m$.

From (15.6.18-2) and (15.6.19-2) again, let $m_1$ and $r_1$ be defined by
$$\align
\text{\rm(15.6.21-2)} \qquad \qquad
m_1 &=\text{Max}\{i:R'_{j+1,i}\ \text{is nonzero with}\ i<m\},
\qquad \qquad \qquad \qquad \qquad \qquad \\
r_1 &=\text{Max}\{i:R''_{j+1,i}\ \text{is nonzero with}\ i<r=m\}.
\endalign$$

For brevity of notation, let $R_{j+1,m}=R'_{j+1,m}=R''_{j+1,m}$.
Noting that $\sum_1$,\dots,$\sum_4$ were already defined by (15.6.11-2)
and (15.6.12-2), then $\sum_{1,1}$ and $\sum_{1,3}$ can be defined by
$$\align
& \sum_{1,1}=\sum_1-(f^{n_{j+1}}_j+R_{j+1,m}f^m_j)^{d_{j+2}} \quad
\text{and} \tag 15.6.22-2 \\
& \sum_{3,1} =\sum_3-(f^{n_{j+1}}_j+R_{j+1,m}f^m_j)^{d_{j+2}}.
\endalign$$

By (15.6.10-2) and (15.6.22-2),
$f-(f^{n_{j+1}}_j+R_{j+1,m}f^m_j)^{d_{j+2}}$ can be rewritten as
follows:
$$\align
\text{\rm(15.6.23-2)} \qquad \qquad \quad
f-(f^{n_{j+1}}_j+R_{j+1,m}f^m_j)^{d_{j+2}}
=\sum_{1,1}+\sum_2=\sum_{3,1}+\sum_4. \qquad \qquad \qquad \qquad \qquad
\endalign$$

In preparation for the proof of $m_1=r_1$ and
$R'_{j+1,m_1}=R''_{j+1,r_1}$, for convenience of notation,
$\sum_{1,1}$ and $\sum_{3,1}$ can be rewritten as follows:
$$\align
\text{\rm(15.6.24-2)} \qquad
\sum_{1,1}&=\sum_1-(f^{n_{j+1}}_j+R_{j+1,m}f^m_j)^{d_{j+2}}
=\phi_{j+1}^{d_{j+2}}-(f^{n_{j+1}}_j+R_{j+1,m}f^m_j)^{d_{j+2}} \qquad \qquad \\
&=\{\sum^{m-1}_{i=0}R'_{j+1,i}f^i_j\}\cdot\{\sum^{d_{j+2}}_{i=1}
\phi_{j+1}^{d_{j+2}-i}(f^{n_{j+1}}_j+R_{j+1,m}f^m_j)^{i-1}\} \\
& \quad \text{where} \quad \sum_1 =\phi_{j+1}^{d_{j+2}} \quad
\text{with} \quad
\phi_{j+1}=(f_j^{n_{j+1}}+\sum^{n_{j+1}-2}_{i=0}R'_{j+1,i}f^i_j), \\
\text{\rm(15.6.25-2)} \qquad
\sum_{3,1}&=\sum_3-(f^{n_{j+1}}_j+R_{j+1,m}f^m_j)^{d_{j+2}}
=\psi_{j+1}^{d_{j+2}}-(f^{n_{j+1}}_j+R_{j+1,m}f^m_j)^{d_{j+2}} \qquad \qquad\\
&=\{\sum^{m-1}_{i=0}R''_{j+1,i}f^i_j\}\cdot\{\sum^{d_{j+2}}_{i=1}
\psi_{j+1}^{d_{j+2}-i}(f^{n_{j+1}}_j+R_{j+1,m}f^m_j)^{i-1}\} \\
& \quad \text{where} \quad \sum_3 =\psi_{j+1}^{d_{j+2}} \quad
\text{with} \quad
\psi_{j+1}=(f_j^{n_{j+1}}+\sum^{n_{j+1}-2}_{i=0}R''_{j+1,i}f^i_j). \\
\endalign$$

For the proof, it suffices to consider two subcases:

\noindent$\underline{\text{\bf Subcases(1):}}$ If there does not exist an integer $m_1$
satisfying (15.6.21-2), then $R'_{j+1,i}$ is identically zero for all
$i=0,1,\dots,m-1$, and so $R''_{j+1,i}$ is identically zero for all
$i=0,1,\dots,m-1$, because otherwise $\sum_{3,1}$ satisfies the same
kind of properties as $\sum_3-f^{d_{j+1}}_j$ does in by (15.6.13),
(15.6.15) and (15.6.17). Therefore, there is nothing to prove that
$\phi_{j+1}=\psi_{j+1}$. \ms

\noindent$\underline{\text{\bf Subcases(2):}}$ If there exist such integers $m_1\ge 0$
and $r_1\ge 0$,
$\sum_{1,1}=\sum_1-(f^{n_{j+1}}_j+R_{j+1,m}f^m_j)^{d_{j+2}}$ and
$\sum_{3,1} =\sum_3-(f^{n_{j+1}}_j+R_{j+1,m}f^m_j)^{d_{j+2}}$ can be
written as follows:

(a) Let $\tau_1$ be defined by $\text{Max}\{\g'_{j+2}\}$ for all
nonzero monomials $\Pi^{j+2}_{k=1}f^{\g'_{k-2}}_{k-2}$ in
$\sum_{1,1}$, where $\g'_1>0,\ \g'_k<n_{k-1}$ for $2\le k\le j+1$
and $\g'_{j+2}\ge 0$, Then, $\tau_1=n_{j+1}(d_{j+2}-1)+m_1<\tau$ by
(15.6.3-2) because $\sum_{1,1}$ satisfies the same kind of properties
as $\sum_1 -f^{d_{j+1}}_j$ does in (15.6.15-2). \ms

(b) Let $\om_1$ be defined by $\text{Max}\{\g''_{j+2}\}$ for all
nonzero monomials $\Pi^{j+2}_{k=1}f^{\g''_{k-2}}_{k-2}$ in
$\sum_{3,1}$ where $\g''_1>0$, $\g''_k<n_{k-1}$ for $2\le k\le j+1$
and $\g''_{j+2}\ge 0$. Then,
$\omega_1=n_{j+1}(d_{j+2}-1)+r_1<\omega$ by (15.6.5-2) since
$\sum_{3,1}$ satisfies the same kind of properties as $\sum_3
-f^{d_{j+1}}_j$ does in (15.6.17-2). \ms

(c) Either $\tau_1=\omega_1$ or $m_1=r_1$ by (a) and (b). \ms

(d) To prove that $R'_{j+1,m_1}=R''_{j+1,r_1}$, first note that for
any nonzero monomial $\Pi^{j+2}_{k=1}f^{\g_k}_k\in
f^{n_{j+1}(d_{j+2}-1)+m_1}_jR'_{j+1,m_1}$, $\g_{j+2}=\tau_1$ by
(15.6.3-2), and if $i<m_1$ then for any nonzero monomial
$\Pi^{j+2}_{k=1}f^{\g'_k}_k\in
f^{n_{j+1}(d_{j+2}-1)+i}_jR'_{j+1,i}$, $\g'_{j+2}<\tau_1$.

Note that whenever $\Pi^{j+2}_{k=1}f^{\g_k}_k\in
f^{n_{j+1}(d_{j+2}-1)+m_1}_jR'_{j+1,m_1}$ then
$\Pi^{j+2}_{k=1}f^{\g_k}_k\in \sum_{1,1}$.

By the similar result as in $\sum_{3,1}$, we can get that
$$
f^{n_{j+1}(d_{j+2}-1)+m_1}_jR'_{j+1,m_1}=f^{n_{j+1}(d_{j+2}-1)+m_1}_j
R''_{j+1,m_1}, \tag 15.6.26-2
$$
and so $R'_{j+1,m_1}=R''_{j+1,r_1}$ with $m_1=r_1$. \ms

If $m_1=r_1=0$, there is nothing to prove. To prove that
$R'_{j+1,i}=R''_{j+1,i}$ for $0\le i< m_1$, under the condition that
$m_1>0$, then it suffices to show that $R'_{j+1,i}=R''_{j+1,i}$ for
$0\le i<m_1$.

From (15.6.21-2) and (15.6.22-2) again, let $m_2$ and $r_2$ be defined by
$$\align
\text{\rm(15.6.27-2)} \qquad \qquad
 m_2 &=\text{Max}\{i:R'_{j+1,i}\ \text{is nonzero with}\ i<m_1\}, \qquad \qquad \qquad \qquad \qquad \qquad \\
r_2 &=\text{Max}\{i:R''_{j+1,i}\ \text{is nonzero with}\
i<r_1=m_1\}.
\endalign$$

For brevity of notation, let
$R_{j+1,m_1}=R'_{j+1,m_1}=R''_{j+1,r_1}$. Noting that
$\sum_1$,\dots,$\sum_4$ were already defined by (15.6.11-2) and
(15.6.12-2), then $\sum_{1,2}$ and $\sum_{3,2}$ can be defined by
$$\align
\text{\rm(15.6.28-2)} \qquad \qquad &
\sum_{1,2}=\sum_1-(f^{n_{j+1}}_j+R_{j+1,m}f^m_j+R_{j+1,m_1}f^{m_1}_j)^{d_{j+2}}
\quad
\text{and} \qquad \qquad \qquad \qquad \\
& \sum_{3,2}
=\sum_3-(f^{n_{j+1}}_j+R_{j+1,m}f^m_j+R_{j+1,m_1}f^{m_1}_j)^{d_{j+2}}.
\endalign$$

By (15.6.10-2) and (15.6.28-2),
$f-(f^{n_{j+1}}_j+R_{j+1,m}f^m_j+R_{j+1,m_1}f^{m_1}_j)^{d_{j+2}}$
can be rewritten as follows:
$$\align
\text{\rm(15.6.29-2)} \qquad
f-(f^{n_{j+1}}_j+R_{j+1,m}f^m_j+R_{j+1,m_1}f^{m_1}_j)^{d_{j+2}}
=\sum_{1,2}+\sum_2=\sum_{3,2}+\sum_4. \qquad  \qquad \qquad
\endalign$$

By the same method as we have used in (15.6.24-2) and (15.6.25-2), it can
be easily shown that $R'_{j+1,m_2}=R''_{j+1,r_2}$ with $m_2=r_2$.
Thus, repeating the above process finitely many times, it can be
proved that $R'_{j+1,i}=R''_{j+1,i}$ for each
$i=0,1,\dots,n_{j+1}-2$, if exists. So, the proof of Step(2) is
done, and then we finished the proof of Sublemma 15.6. $\square$
\ms

By Sublemma 15.5 and Sublemma 15.6, we completed the proof of Theorem 15.4 and
Corollary 15.4.1.
$\square$ \bs

\qquad \qquad

\newpage

\ms

{\bf Part[C3] Irreducibility criterion of W-polys of two complex
variables} \bs

{\bf \S {16}. In preparation for irreducibility criterion of W-polys
of two complex variables} \bs

{\bf $\S$16.1. Known preliminaries on irreducibility criterion of
germs of analytic functions of two complex variables} \bs

 \proclaim{Lemma 16.0} $\underline{\text{\bf
{Assumptions}}}$ Let $f(y,z)=z^n+a_{n-2}y^{\a_{n-2}}z^{n-2}+\cdots
+a_1y^{\a_1}z+a_0y^{\a_0}$ be irreducible in ${}_2\CO_0$ with
multiplicity $n\ge 2$ at $(0,0)\in \C^2$ where for $0\le i\le n-2$,
each $a_i=a_i(y,z)$ is a unit in ${}_2\CO_0$ if exists and the
$\a_i$ are positive integers. Note that $a_{n-1}$ is identically
zero. Assume that $d=\gcd(n,\a_{0})>1$. Then, we can write $n=dn_1$
and $\a_0=d\a_{1,0,1}$ with $\gcd(n_1,\a_{1,0,1})=1$ and $2\le
n_1<\a_{1,0,1}$.

In particular, if $a_i(y,z)=a_i(y)$ for all $i$, then $f(y,z)$ is
called a $W$-poly in $z$. \ms

$\underline{\text{\bf {Conclusions}}}$ Then, $f$ can be written
uniquely in the form
$$\align
(16.0.1) \qquad \qquad f=(z^{n_1}+\xi y^{\a_{1,0,1}})^d +\sum_{p,q\ge
0} c_{p,q}y^pz^q \quad \text{with} \quad n_1p+\a_{1,0,1}q>n_1\a_{1,0,1}d,
\qquad \qquad
\qquad \\
\endalign$$
where the $c_{p,q}$ are nonzero complex numbers for some nonnegative
integers $p$ and $q$ such that $n_1p+\a_{1,0,1}q>n_1\a_{1,0,1}d$,
satisfying the following properties:

\roster \item"{\rm(i)}" $\xi$ is a unique nonzero number such that
$$
{}_dC_i\xi^i = a_{n-in_1}(0,0)\quad \text{for}\quad 1\le i\le d.
\tag 16.0.2
$$
\item"{\rm(ii)}" $\f{\a_0}n =
\f{\a_{1,0,1}}{n_1}=\f{\a_{n-in_1}}{in_1}$, i.e.,
$\a_{n-in_1}=i\a_{1,0,1}$ for $1\le i\le d$, and
$$
(16.0.3) \qquad \qquad \f{\a_{n-j}}{j} >\f{\a_{1,0,1}}{n_1}\quad
\text{for any $j\ne n_1, 2n_1,\dots, (d-1)n_1$} \quad \text{where
$1\le j\le n$}. \qquad \qquad
$$
\item"{\rm(iii)}" If $q<n$, then $n_1p+\a_{1,0,1}q
>n_1\a_{1,0,1}d$ if and only if $\f p{n-q}>\f{\a_0}n
=\f{\a_{1,0,1}}{n_1}$.
\endroster

Moreover, if $f(y,z)$ is a $W$-poly in $z$, note that $p>0$ and
$q\le n-2$.
\endproclaim \ms

\demo{\bf Proof of Lemma 16.0} It just follows from Hensel's lemma or
Theorem 3.2. \enddemo \ms

\definition{Remark 16.0.1} Note by {\rm(i)} that $\xi =\f 1da_{n-n_1}(0,0)$ and
the $a_{n-in_1}(0,0)$ are nonzero for all $i=1,\dots,d$.
\enddefinition \ms

\proclaim{Theorem 16.1(By Theorem 3.6 and Theorem 3.7)}

$\underline{\text{\bf {Assumptions}}}$

{\rm (a)} Let $V(f)=\{(y,z): f(y,z)=0\}$, and $V(F)=\{(y,z):
F(y,z)=0\}$ be analytic varieties at $(0,0)$ in $\BC^2$, each of
which is written respectively in the form,
$$\align
\text{\rm (16.1.1)} \qquad \qquad f&=(z^{n_1}+\ve y^{k_1})^{d}
+\sum_{\alpha,\beta\ge 0}c_{\alpha,\beta}y^{\alpha}z^{\beta} \quad
\text{with} \quad n_1\alpha+k_1\beta>n_1k_1d,
\qquad \qquad \qquad \\
F&=y^{\de_1}z^{\de_2}f,
\endalign$$
satisfying the properties {\rm(i)}, {\rm(ii)},\dots, {\rm(vi)}:

\roster \item "(i)" $\gcd(n_1,k_1)=1$ with $1\le n_1<k_1$ and $d$ is
a positive integer.

\item "(ii)" $\ve$ is a unit in $\BC\{y,z\}$.

\item "(iii)" The $c_{\alpha,\beta}$ are nonzero complex numbers
for some nonnegative integers $\alpha$ and $\beta$ such that
$n_1\alpha+k_1\beta>n_1k_1d$.

\item "(iv)" $\de_1$ and $\de_2$ are nonnegative integers.

\item "(v)" If $n_1=1$, assume that $\de_2$ is a positive integer.

\item "(vi)" If $d\ge 2$ and $n_1=1$, then assume in addition that
$V(f)$ has an isolated singular point at the origin as a reduced
variety.
\endroster

{\rm(b)} Let $V(G)=\{(y,z):G(y,z)=0\}$ be another analytic variety
with isolated singularity at the origin in $\BC^2$ defined by the
form
$$\align
g &=z^{n_1}+y^{k_1} \quad
\text{with $\gcd(n_1,k_1)=1$ and $1\le n_1<k_1$}, \tag 16.1.2 \\
G &=z^{\gamma}g,
\endalign$$
satisfying the properties {\rm(i)} and {\rm(ii)}:

\roster \item "(i)" If $n_1=1$, then $\gamma=1$.

\item "(ii)" If $n_1\ge 2$, then $\gamma=0$, and so
$G(y,z)=g(y,z)$.
\endroster

$\underline{\text{\bf {Conclusions}}}$ \quad Let
$\tau_m=\pi_1\circ\pi_2\circ\cdots\circ\pi_m:M^{(m)}\to\BC^2$ be the
compositions of a finite number $m$ of successive blow-ups $\pi_i$
which is needed to get the standard resolution of the singular point
of $V(G)$. If $V(g)$ has the singular point at the origin, then as
compared with the above $\tau_m$, exactly the same $\tau_m$ can be
also used for the standard resolution of the singular point of
$V(g)$ as far as the blow-ups process is concerned.

{\rm(a)(a1)} We can use just one coordinate patch of the local
coordinates for each blow-up $\pi_i$ of $\tau_m$ with $1\le i\le m$
in the sense of {\rm Lemma 2.12}.

{\rm \quad(a2)} Just as above, we can use the same $\tau_m$ for the
composition of the first finite number $m$ of successive blow-ups in
preparation for the standard resolution of the singular point
$(0,0)$ of $V(F)$.

{\rm \quad(a3)} Also, we can use just the common one coordinate
patch of the given local coordinates for each blow-up $\pi_i$ of the
above $\tau_m$ in {\rm (a1)}, in order to study any of $V^{(i)}(F)$
for all $i=1,2,\dots,m$ in the sense of {\rm Lemma 2.14}. \ms

{\rm(b)} For simplicity of notations, let $(v,u)$ be the common one
of the local coordinates for the $m-th$ blow-up $\pi_m:M^{(m)}\to
M^{(m-1)}$ at $(0,0)$ which is the quasisingular point of
$V^{(m-1)}(G)$ in the sense of Definition 2.6. Being viewed as an
analytic mapping, $\tau_m:M^{(m)}\to\BC^2$ can be written in the
form
 $$
 \tau_m(v,u)=(y,z)=(v^{n_1}u^a,v^{k_1}u^b), \tag 16.1.3
 $$
where

$(b_1)$ $a$ and $b>0$ are nonnegative integers such that
 $bn_1-ak_1=1$,

$(b_2)$ $E_m=\{v=0\}$ is defined by the $m-th$ exceptional curve of
the first kind. \ms

{\rm(c)} By {\rm(b)}, along $v=0$, $(F\circ\tau_m)_{total}$ can be
written in the following form:
$$\align
(F\circ\tau_m)_{total}&=v^{e_m}u^{\varepsilon}(f\circ\tau_m)_{proper}
\quad \text{with} \tag 16.1.4 \\
(f\circ\tau_m)_{proper}&=(u+\ve')^d+\sum_{\alpha,\beta\ge
0}c_{\alpha,\beta}v^{n_1\alpha+k_1\beta-n_1k_1d}u^{\ve_{\alpha,\beta}},\\
(G\circ\tau_m)_{total}&=v^{k_1\gamma+n_1k_1}u^{b\gamma+ak_1}(u+1),
\endalign$$

where \roster \item "(i)" $e_m=n_1\delta_1+k_1\delta_2+n_1k_1d$,
$\ve=a\delta_1+b\delta_2+ak_1d$ and
$\ve_{\alpha,\beta}=a\alpha+b\beta-ak_1d\ge 0$,

\item "(ii)" by assumption of $V(G)$, if $n_1=1$ then $\gamma=1$,
and if $n_1\ge 2$ then $\gamma=0$,

\item "(iii)" $\ve'$ is a unit in $\C\{u+1,v\}$.
\endroster \ms

{\rm(d)} By {\rm(c)}, we have the following:

{\rm (d1)} Whether $n\ge 2$ or not, $G\in the ~type[1]$ under
$\tau_m$.

{\rm (d2)} If $n=1$ then $f\in the~ type[0]$ and $F\in the ~type[1]$
under $\tau_m$, but if $n\ge 2$, then $f\in the~ type[1]$ and $F\in
the ~type[1]$ under $\tau_m$ in the sense of Definition $2.8$. \ms

Observe by the notation that $(f\circ\tau_m)_{proper}$ is the local
defining equation for the proper transform $V^{(m)}(f)$ and that
$(F\circ\tau_m)_{total}$ is the local defining equation for the
total transform of $V(F)$ under $\tau_m$. \bs

\endproclaim \ms

{\bf $\S$16.2. An algorithm for computing a complete irreducibility
criterion of any W-poly in $\BC\{y,z\}$ which has the same
multiplicity sequence as the standard Puiseux expansion($y=t^{n}$
and $z=t^{\alpha}+t^{\beta}$) does} \ms

In $\S16.2$, the first aim is to find an algorithm for computing a
complete irreducibility criterion of any W-poly in $\BC\{y,z\}$
which has the same multiplicity sequence as the standard Puiseux
expansion($y=t^{n}$ and $z=t^{\alpha}+t^{\beta}$) does by Theorem
16.5 and Theorem 16.4 with Proposition 16.2 and Proposition 16.3,
whose proofs can be represented by $\S16.3$.

In $\S17$, the main aim is to finish a complete irreducibility
criterion of any W-poly in $\BC\{y,z\}$ having the same multiplicity
sequence as the standard Puiseux expansion does, which can be
represented by Theorem 16.6 with Proposition 16.7 and Proposition
16.8. \ms

\proclaim{Proposition 16.2} $\underline{\text{\bf {Assumptions}}}$
Let $f(y,z)=z^n+a_{n-2}y^{\a_{n-2}}z^{n-2}+\cdots
+a_1y^{\a_1}z+a_0y^{\a_0}$ be an irreducible $W$-poly in $z$ with
multiplicity $n\ge 2$ at $(0,0)\in \C^2$ where for $0\le i\le n-2$,
each $a_i=a_i(y)$ is a unit in ${}_2\CO_0$ if exists and the $\a_i$
are positive integers. Note that $a_{n-1}$ is identically zero. Let
$d_2=\text{\rm gcd}(n,\a_0)>1$ with $n=d_2n_1$ and $\a_0
=d_2\a_{1,0,1}$. Note that $2\le n_1<\a_{1,0,1}$. \ms

$\underline{\text{\bf {Conclusions}}}$ Then, $(g_1,f)$ can be
written in the form
$$\cases
g_1 &=z^{n_1}+\xi_1 y^{\sigma_{1}} \qquad \text{with
$\si_1=\alpha_{1,0,1}$,  $f_{-1}=y$  and  $f_{0}=z$},\\
f &= g^{d_2}_1 +\sum^{d_2-1}_{i=0}T_{2,i}g^i_1,
\endcases \tag 16.2.1
$$

where, considering $f_{-1},f_0,g_{1}$ as independent complex
3-variables at $0\in\C^{3}$,

{\rm(i)} $n=d_{2}n_1$ with $d_{2}\ge 2$ and $n_{1}\ge 2$, and
$n=d_1$ if necessary;

{\rm(ii)} $\si_1=\alpha_{n-n_{1}}=\alpha_{1,0,1}$ and $\xi_{1}=\f
{1}{d_2}a_{n-n_{1}}(0)$;

{\rm(iii)} $T_{2,i}=T_{2,i}(f_{-1},f_0)\in \C \{f_{-1},f_0\}$ of $f$
in {\rm(16.2.1)} for $i=0,1,\dots,d_{2}-1$;

{\rm(iv)} $g_{1}=g_{1}(f_{-1},f_0)\in \C\{f_{-1}\}[f_0]$;

{\rm(v)} $f=f(f_{-1},f_0,g_1)\in \C \{f_{-1},f_0\}[g_1]$ of $f$ in
{\rm(16.2.1)},

satisfying two conditions, denoted by The Necessary and Sufficient
Condition$[A]$ for $g_{1}(y,z)\in$the type$[1]$ and The Necessary
Condition$[B]$ for $f(y,z)\in$the type$[\ell]$ with $\ell\ge 2$,
each of which is represented respectively, as follows: \ms

$\underline{ \text{\bf [1]  The Necessary and Sufficient
Condition[A] for $g_{1}(y,z)\in$the type$[1]$:}}$

\noindent$\underline{ \text{\bf{\bf $g_{1}\in \BC\{f_{-1}\}[f_{0}]$}
is an irreducible W-poly  of degree $n_1$ in $f_0$ with a
coefficient of }}$

\noindent$\underline{ \text{\bf ${f_{0}}^{n_1-1}$ zero in
$\C\{f_{-1}\}$, and $g_{1}(y,z)\in$ the type[1] in the sense of
Definition 2.5}}$ \ms

To find \text{\rm The Necessary and Sufficient Condition[A] for
$g_{1}(y,z)\in$the type$[1]$}, it is clear that $g_{1}(y,z)\in
\C\{y\}[z]$ of {\rm(16.2.1)} itself is an irreducible $W$-poly of
degree $n_1$ in $z$ with coefficients in $\BC\{y\}$ and with
multiplicity $n_1$ at $0\in\C^{2}$, and $g_{1}(y,z) \in$ the type
$[1]$. \ms

$\underline{ \text{\bf [2] The Necessary Condition[B] for
$f(y,z)\in$the type$[\ell]$ with $\ell\ge 2$:}}$

\noindent$\underline{ \text{\bf $f\in \C\{f_{-1},f_0\}[g_{1}]$ is an
irreducible W-poly of degree $d_2$ in $g_{1}$ with a coefficient
of}}$

\noindent$\underline{ \text{\bf ${g_{1}}^{d_2-1}$ either zero or
nonzero in $\C\{f_{-1},f_0\}$, and  $f(y,z)\in $the type[${\ell}$]
with $\ell\ge 2$ }}$

\noindent$\underline{ \text{\bf in the sense of Definition 2.5}}$
\ms

To find \text{\rm the Necessary Condition[B] for $f(y,z)\in$the
type$[\ell]$ with $\ell\ge 2$}, it suffices to show that $f
=g^{d_{2}}_1 +\sum^{d_{2}-1}_{i=0} T_{2,i}g^i_1  \in
\BC\{f_{-1},f_0\}[g_{1}]$ of {\rm(16.2.1)} satisfies two properties
{\rm(1)} and {\rm(2)}: Note that either $\ell=2$ or $\ell\ge3$ and
that $T_{2,d_{2}-1}$ may not be zero. \ms

{\rm (1)} Each $T_{2,i}\in \C \{f_{-1},f_0\}$ of $f$ in
{\rm(16.2.1)} satisfies the properties {\rm(1a)}, {\rm(1b)},
{\rm(1c)} and {\rm(1d)} for $i=0,1,\dots,d_{2}-1$.

{\rm(1a)}{\rm (1a-1)} For any nonzero monomial
$\Pi^{2}_{k=1}f^{\delta_k}_{k-2}$ in $T_{2,i}$, $\de_1>0$ and
$\de_2<n_1$.

{\rm(1a-2)} In particular, for any nonzero monomial
$\Pi^{2}_{k=1}f^{\delta_k}_{k-2}$ in $T_{2,d_{2}-1}$, $\de_1>0$ and
$\de_2 \le n_1-2$. \ms

{\rm(1b)} Let $\N_0$ be the set of nonnegative integers and $\N^2_0$
be its two dimensional copy.

{\rm(1b-1)} Let $\th_1: \N_0\to \N_0$ and $\th_2: \N^2_0\to \N_0$ be
integer-valued functions, each of which is defined respectively, as
follows:

\qquad $\th_1(t)=t$ for each $t\in \N_0$, and

\qquad $\ol\th_2(t_1,t_2)=t_2\th_1(\a_{1,0,1}) + n_1\th_1(t_1) =
t_2\a_{1,0,1}+n_1t_1$ for each $(t_1,t_2)\in \N^2_0$. \ms

{\rm(1b-2)} For any two nonzero monomials
$\Pi^{2}_{k=1}f^{\beta_k}_{k-2}$ and $\Pi^{2}_{k=1}f^{\g_k}_{k-2}$
in $T_{2,i}$ with $i$ fixed,
$$\align
(16.2.2) \qquad \qquad
&\text{$\ol\th_{2}(\beta_k)^{2}_{k=1}=\ol\th_{2}(\g_k)^{2}_{k=1}$ if
and only if $\beta_k=\g_k$ for $k=1,2$.} \qquad \qquad \qquad \qquad
\\
&\text{So, there is a unique nonzero monomial
$C_{2,i}\Pi^{2}_{k=1}f^{{\beta}_{2,i,k}}_{k-2}$ in
$T_{2,i}$} \\
&\text{with a constant $C_{2,i}$ such that
$\ol\th_{2}(\beta_{2,i,k})^{2}_{k=1}
=\text{$\min$}\{\ol\th_{2}(\g_k)^{2}_{k=1}\}$}\\
&\text{for any nonzero monomial $\Pi^{2}_{k=1}f^{\g_k}_{k-2}$ in
$T_{2,i}$.}
\endalign$$

{\rm(1c)} For each $i=0,1,\dots,d_2-1$,
$$\align
{(16.2.3)}\qquad \qquad
\ol\th_{2}(\beta_{2,i,1},\beta_{2,i,2})>({d_2}-i)
n_{1}\th_1(\a_{1,0,1}).\qquad \qquad\qquad \qquad\qquad \qquad\qquad
\qquad
\endalign$$

$\underline{\text{\rm(1d)}}$ For all $i= 0,1,\dots,d_2-1$, the
following hold:
$$\align
&\gcd(d_2,\ol\th_{2}(\beta_{2,0,1},\beta_{2,0,2}))\ge 1 \quad \text{and}  \tag 16.2.4\\
& \dfrac{\ol\th_{2}(\beta_{2,i,1},\beta_{2,i,2})}{d_2-i}\ge
\dfrac{\ol\th_{2}(\beta_{2,0,1},\beta_{2,0,2})}{d_2}.
\endalign$$

Then, either $\gcd(d_{2},\ol\th_{2}(\beta_{2,0,k})^{2}_{k=1})=1$ or
$1<\gcd(d_{2},\ol\th_{2}(\beta_{2,0,k})^{2}_{k=1})\le d_{2}$. \ms

$\underline{\text{\rm(1d-1)}}$ Suppose
$\gcd(d_2,\ol\th_{2}(\beta_{2,0,k})^{2}_{k=1})=1$. Then $f$ is
irreducible in ${}_2\CO_0$ with $f\in $ the type $[2]$ in the sense
of Definition 2.5 if and only if the inequality in {\rm(16.2.4)}
holds and $g_1$ is irreducible in ${}_2\CO_0$ with $g_{1} \in$ the
type $[1]$ in the sense of Definition 2.5.

$\underline{\text{\rm(1d-2)}}$ Suppose
$1<\gcd(d_2,\ol\th_{2}(\beta_{2,0,k})^{2}_{k=1})\le d_{2}$ in
{\rm(16.2.4)}. To find an irreducible criterion, it remains to study
two subcases respectively:

$\underline{\text{\rm Subcase(i) of (1d-2)}}$ Let
$\gcd(d_2,\ol\th_2(\beta_{2,0,1},\beta_{2,0,2}))= d_{2}$ in
{\rm(16.2.4)}. Then, $f$ is either irreducible or not in
${}_2\CO_0$. If $f$ is irreducible in ${}_2\CO_0$ then $f \in$ the
type $[\ell]$ with $\ell\ge 2$ in the sense of Definition 2.5. \ms

$\underline{\text{\rm Subcase(ii) of (1d-2)}}$ Let
$1<\gcd(d_2,\ol\th_2(\beta_{2,0,1},\beta_{2,0,2}))< d_{2}$ in
{\rm(16.2.4)}. Then, $f$ is either irreducible or not in
${}_2\CO_0$. If $f$ is irreducible in ${}_2\CO_0$ then $f \in$ the
type $[\ell]$ with $\ell\ge 3$ in the sense of Definition 2.5. \ms

{\rm(2)}{\rm(2a)} $f=f(y,z)\in \C\{y\}[z]$ is an irreducible
$W$-poly of degree $n$ in $z$ with coefficients in $\BC\{y\}$ and
with multiplicity $n=d_{2}n_1$ at $0\in\C^2$. Also, either $f \in$
the type $[2]$ or $f \in$ the type $[\ell]$ with $\ell\ge 3$ in the
sense of Definition 2.5.

{\rm(2b)} $f=f(f_{-1},f_0,g_1)\in \C\{f_{-1},f_0\}[g_{1}]$ of
{\rm(16.2.1)} is an irreducible $W$-poly of degree $d_2$ in $g_1$
with coefficients in $\C\{f_{-1},f_0\}$ and with multiplicity
$d_{2}$ at $0\in\C^{3}$  where $f_{-1},f_0,g_{1}$ are viewed as
independent complex three variables at the origin in $\C^{3}$.
\endproclaim
\ms

\noindent$\underline{\text{\bf Remark 16.2.1.1}}$ {\rm(1)} It is
clear that $f=f(y,z)\in \C\{y\}[z]$ is a $W$-poly of degree $n$ in
$z$ with coefficients in $\BC\{y\}$ and with multiplicity
$n=d_{2}n_1$ at $0\in\C^2$.

{\rm (2)} It is clear that $g_{1}$ of (16.2.1) satisfies The
Necessary and Sufficient Condition$[A]$ for $g_{1}(y,z)\in$the
type$[1]$.

{\rm(3)} Whether $T_{2,d_{2}-1}=0$ or not, then it is said by
Proposition 16.2 that $f$ of (16.2.1) satisfies The Necessary
Condition[B] for $f(y,z)\in$the type[${\ell}$] with $\ell\ge 2$.

{\rm(4)} If $\gcd(d_{2},\ol\th_{2}(\beta_{2,0,k})^{2}_{k=1})=1$, it
is said by (1d) of Proposition 16.2 that $f$ satisfies The Necessary
and Sufficient Condition$[A]$ for $f(y,z)\in$the type[${2}$].

{\rm(5)} Let $1<\gcd(d_2,\ol\th_2(\beta_{2,0,1},\beta_{2,0,2}))<
d_{2}$ in {\rm(16.2.4)}. If $f$ is irreducible in ${}_2\CO_0$ with
$T_{2,d_{2}-1}$ either zero or not, it is proved by this proposition
that $f \in$ the type $[\ell]$ with $\ell\ge 3$ in the sense of
Definition 2.5. But, if $f$ is irreducible in ${}_2\CO_0$ and
$T_{2,d_{2}-1}\not=0$,
$\gcd(d_2,\ol\th_2(\beta_{2,0,1},\beta_{2,0,2}))$ may be equal to
$d_2$.

{\rm(6)} Let $\gcd(d_2,\ol\th_2(\beta_{2,0,1},\beta_{2,0,2}))=
d_{2}$ in Subcase(i) of (1d-2). If $f$ is irreducible in ${}_2\CO_0$
then $T_{2,d_{2}-1}\not=0$, otherwise if $T_{2,d_{2}-1}=0$,
$\gcd(d_2,\ol\th_2(\beta_{2,0,1},\beta_{2,0,2}))< d_{2}$ by Hensel's
lemma(Theorem 3.2).

{\rm(7)} Let
$f=g^{5}_1+\binom{5}{1}y^{6}zg^{4}_1+\binom{5}{2}y^{12}z^{2}g^{3}_1
+\binom{5}{3}y^{18}z^{3}g^{2}_1+\{\binom{5}{4}y^{24}z^{4}+y^{30}\}g_1
+y^{35}z^{2}-y^{37}$ where $g_1=z^{5}+y^{7}$. Note that
$g_1+y^{6}z=0$ and $g_1=0$ have the same multiplicity sequence.
Then, it is easy to compute that the above pair $(g_1,f)$ satisfies
a pair of (16.2.1), which belongs to $\underline{\text{\rm
Subcase(i) of (1d-2)}}$ of this proposition with $f(y,z)\in$the
type[${2}$]. \ms

\definition{Remark 16.2.1.2} Assuming that $f$ is irreducible in in
${}_2\CO_0$, let $h_2=h_1+\f{1}{d_2}T_{2,d_2-1}$ where $h_1=g_1$ and
$T_{2,d_2-1}$ were defined by (16.2.1). Then, $h_2$ is an
irreducible $W$-poly of degree $n_1$ in $z$, and also two curves
defined by $h_2=0$ and $g_1=0$ have the same multiplicity sequence
by Theorem 12.1 because for any nonzero monomial
$y^{\de_1}z^{\de_2}\in T_{2,d_2-1}$,
$\theta_2(\de_1,\de_2)>n_1\a_{1,0,1}$ by (1c) of Proposition $16.2$.
So, it is said by Definition of 2.4 that $V(h_2)$ and $V(h_1)$ have
the same multiplicity sequence under two standard resolutions,
denoted by ${h_2} \buildrel \text{{\rm multiseq}} \over \sim {h_1}$
under two standard resolutions.
\enddefinition \ms

\definition{Definition 16.2.2} If $h(y,z)\in \C\{y\}[z]$ and
$f(y,z)\in \C\{y\}[z]$ are $W$-polys in $z$ with coefficients in
$\C\{y\}$, then $(h,f)$ is called a pair of W-polys in $z$. Also,
assuming that $\ell(y,z)\in \C\{y\}[z]$ and $g(y,z)\in \C\{y\}[z]$
are $W$-polys in $z$ with coefficients in $\C\{y\}$, it is said for
brevity of notation that $(h,f)=(\ell,g)$ if and only if
$h(y,z)=\ell(y,z)$ and $f(y,z)=g(y,z)$, for the proof of Proposition
16.3 and its applications.
\enddefinition \ms

\proclaim{Proposition 16.3} $\underline{\text{\bf {Assumptions}}}$
\quad Let $f(y,z)=z^n+\sum^{n-2}_{i=0}a_iy^{\a_i}z^i$ be an
irreducible $W$-poly in $z$ with multiplicity $n\ge 2$ at $(0,0)\in
\C^2$ where for $0\le i\le n-2$, each $a_i=a_i(y)$ is a unit in
${}_2\CO_0$ if exists and the $\a_i$ are positive integers. Note
that $a_{n-1}$ is identically zero. Assume that
$d_{2}=\gcd(n,\a_0)>1$. Write $n=d_{2}n_1$ and
$\a_0=d_{2}\a_{1,0,1}$ with $\gcd(n_1,\a_{1,0,1})=1$. Note that
$2\le n_1<\a_{1,0,1}$. Without any need of proof, we may assume by
{\rm Proposition $16.2$} and {\rm Definition 16.2.2} that there
exists a pair of W-polys in $z$ with coefficients in $\C\{y\}$,
denoted by $(g_1,f)$, which can be written in the following form:
$$\cases
g_1 &=z^{n_1}+\xi_1 y^{\sigma_1}\ \text{with}\ \xi_1 ={\f 1{d_2}} a_{n-n_1}(0)\neq {0},\\
f &=g^{d_2}_1+\sum^{d_{2}-1}_{i=0}T_{2,i}g^i_1,
\endcases \tag 16.3.0
$$
$\si_1=\alpha_{n-n_1}=a_{1,0,1}$, satisfying two conditions, denoted
by The Necessary and Sufficient Condition$[A]$ for
$g_{1}(y,z)\in$the type$[1]$ and The Necessary Condition$[B]$ for
$f(y,z)\in$the type$[\ell]$ with $\ell\ge 2$, as we have seen in
Proposition $16.2$. \ms

$\underline{\text{\bf {Conclusions}}}$ The main aim is to construct
a unique pair $(f_{1},f)$ such that $(f_{1},f)$ can be written in
the form
$$\cases
f_{1} &=z^{n_{1}}+\sum^{n_{1}-2}_{i=0} R_{1,i}z^i \quad \text{with}
\ f_{-1}=y\ \text{and}\ f_0=z,
\\
f &=f^{d_{2}}_{1} +\sum^{d_{2}-2}_{i=0} S_{2,i}f^i_{1},
\endcases \tag 16.3.1
$$
where $f_{-1},f_0,f_1$ are considered as independent complex three
variables at the origin in $\C^{3}$ if necessary, satisfying the
following properties:

{\rm(i)} The first problem is how to construct $f_{1}=f_{1}(y,z)$
satisfying the condition in \text{\rm
$\widehat{\widehat{\text{\rm[1]}}}$} such that
\text{$f_{1}(y,z)\buildrel \text{{\rm multiseq}} \over \sim
g_{1}(y,z)$ \quad under the standard resolutions}.

$\underline{ \text{\bf $\widehat{\widehat{\text{\bf[1]}}}$  The
Necessary and Sufficient Condition[A] for $f_{1}(y,z)\in$the
type$[1]$:}}$

\noindent$\underline{ \text{\bf{\bf $f_{1}\in \BC\{f_{-1}\}[f_{0}]$}
is an irreducible W-poly  of degree $n_1$ in $f_0$ with a
coefficient of }}$

\noindent$\underline{ \text{\bf ${f_{0}}^{n_1-1}$ zero in
$\C\{f_{-1}\}$, and $f_{1}(y,z)\in$ the type[1] in the sense of
Definition 2.5}}$ \ms

{\rm(ii)} The second problem is to prove that $f=f(f_{-1},f_0,f_1)$
satisfies the condition in \text{\rm
$\widehat{\widehat{\text{\rm[2]}}}$} which is defined by the same
kind of property as $f(f_{-1},f_0,g_1)$ have done in The Necessary
Condition$[B]$ for $f(y,z)\in$the type$[\ell]$ with $\ell\ge 2$.

$\underline{ \text{\bf $\widehat{\widehat{\text{\bf[2]}}}$ The
Necessary Condition[B] for $f(y,z)\in$the type$[\ell]$ with $\ell\ge
2$:}}$

\noindent$\underline{ \text{\bf {\bf $f\in \C\{f_{-1},f_0\}[f_{1}]$}
is an irreducible W-poly of degree $d_2$ in $f_{1}$ with a
coefficient}}$

\noindent$\underline{ \text{\bf of ${f_{1}}^{d_2-1}$ zero in
$\C\{f_{-1},f_0\}$, and $f(y,z)\in$the type[${\ell}$] with $\ell\ge
2$}}$

\noindent$\underline{ \text{\bf in the sense of Definition 2.5}}$
\ms

In order to find the construction of a pair $(f_{1},f)$ in
{\rm(16.3.1)}, first of all, it suffices to consider the following
two cases, depending on the fact that either $T_{2,d_{2}-1}$ of
{\rm(16.3.0)} is zero or not. For brevity of notations, let
$h_1=g_1$ and $T^{(1)}_{2,i}=T_{2,i}$. \ms

$\underline{\text{\bf Case(1)}}$ Let $T^{(1)}_{2,d_{2}-1}=0$. Put
$f_{1}=g_{1}$, $R_{1,0}=\xi_1 y^{\sigma_1}$, and
$S_{2,i}=T^{(1)}_{2,i}$ for $0\le i\le d_{2}-2$. Then, it was
already proved by Proposition $16.2$ that $(f_{1},f)$ of the main
aim and $(g_{1},f)$ are the same pairs in the sense of Definition
$16.2.2$. \ms

$\underline{\text{\bf Case(2)}}$ Let $T^{(1)}_{2,d_{2}-1}\not=0$.
Then, there is a sequence of pairs of $W$-polys in $z$,
$\{(h_p,f):p=1,2,\dots\}$ with $h_1=g_1$ such that $(h_{\nu},f)\not
=(h_{\nu+1},f)=(h_{\nu+2},f)=\cdots$ for some integer $\nu\le
\f{n_1+1}2$, each pair of which can be written in the form
$$\cases
h_1 &=g_1=z^{n_1}+\xi_1 y^{\a_{1,0,1}}, \\
f &=h^{d_2}_1+\sum^{{d_2}-1}_{i=0}T^{(1)}_{2,i}h^i_1,
\endcases \tag 16.3.2
$$
and for each $p=2,3,\dots$
$$\cases
h_{p} &=h_{p-1}+\f 1{d_2}T^{({p-1})}_{2,d_2-1}
=z^{n_1}+\sum^{n_1-2}_{i=0}R^{(p)}_{1,i}z^i, \\
f &=h^{d_2}_{p} +\sum^{d_2-1}_{i=0}T^{(p+1)}_{2,i}h^i_{p},
\endcases \tag 16.3.3
$$
with $T^{(p)}_{2,d_{2}-1}\ne 0$ for $1\le p\le \nu$ and $T^{(\nu
+1)}_{2,d_{2}-1}=T^{(\nu +2)}_{2,d_{2}-1}=\cdots =0$,

where, considering $f_{-1},f_{0}, h_p$ as independent complex
$3$-variables at $0\in\C^{2}$,

{\rm(i)} $n=d_{2}n_1$ with $d_{2}\ge 2$ and $n_1\ge 2$;

{\rm(ii)} $R^{(p)}_{1,i}=R^{(p)}_{1,i}(f_{-1})\in\C\{f_{-1}\}$ for
$p\ge 1$ and $0\le i\le n_1-2$;

{\rm(iii)}
$T^{(p)}_{2,i}=T^{(p)}_{2,i}(f_{-1},f_0)\in\C\{f_{-1},f_0\}$ for
$p\ge 1$ and $0\le i\le d_{2}-1$;

{\rm(iv)} $h_p=h_p(f_{-1},f_0)\in \C\{f_{-1}\}[f_0]$ for $p\ge 1$;

{\rm(v)} $f=f(f_{-1},f_0, h_p)\in \C\{f_{-1},f_0\}[h_p]$ for each
fixed $(h_p,f)$ of either {\rm(16.3.2)} or {\rm(16.3.3)},

{\noindent}satisfying two conditions, denoted by The Necessary and
Sufficient Condition$[A]$ for $h_{p}(y,z)\in$the type$[1]$ and The
Necessary Condition$[B]$ for $f(y,z)\in$the type$[\ell]$ with
$\ell\ge 2$, each of which is represented respectively, as follows:

In more detail, for any fixed $(h_{p},f)$ of {\rm(16.3.3)}, $h_{p}$
of $(h_{p},f)$ satisfies The Necessary and Sufficient Condition$[A]$
for $h_{p}(y,z)\in$the type$[1]$ and $f(y,z)$ of $(h_{p},f)$
satisfies The Necessary Condition$[B]$ for $f(y,z)\in$the
type$[\ell]$ with $\ell\ge 2$. In particular, if $p_1$ is a given
integer with $p_1\ge \f{n_{1}+1}{2}$, define $(f_1,f)=(h_{p_1},f)$
by Definition $16.2.2$. Then, $f_1$ of $(f_1,f)$ of $(16.3.1)$
satisfies The Necessary and Sufficient Condition$[A]$ for
$f_1(y,z)\in$the type$[1]$ in \text{\rm
$\widehat{\widehat{\text{\rm[1]}}}$}, and $f$ of $(f_1,f)$ of
$(16.3.1)$ satisfies The Necessary Condition$[B]$ for $f(y,z)\in$the
type$[\ell]$ with $\ell\ge 2$ in \text{\rm
$\widehat{\widehat{\text{\rm[2]}}}$} by the same way  as $(h_p,f)$
of $(16.3.3)$ does up to the change of notations: Recall that
$f_{-1}=y$ and $f_0=z$. \ms

$\underline{ \text{\bf [1] The Necessary and Sufficient Condition[A]
for $h_{p}(y,z)\in$the type$[1]$:}}$

\noindent$\underline{ \text{\bf{\bf $h_{p}\in \BC\{f_{-1}\}[f_{0}]$}
is an irreducible W-poly  of degree $n_1$ in $f_0$ with a
coefficient of }}$

\noindent$\underline{ \text{\bf ${f_{0}}^{n_1-1}$ zero in
$\C\{f_{-1}\}$, and $h_{p}(y,z)\in$ the type[1] in the sense of
Definition 2.5}}$ \ms

To find \text{\rm The Necessary and Sufficient Condition[A] for
$h_{p}(y,z)\in$the type$[1]$}, it suffices to show that
$h_{p}=z^{n_1}+\sum^{n_1-2}_{i=0}R^{(p)}_{1,i}z^i$ of {\rm(16.3.3)}
satisfies two properties {\rm(1)} and {\rm(2)}: \ms

{\rm (1)} Let $p$ be fixed with $p\ge 1$. Each $R^{(p)}_{1,i}\ne 0$
satisfies the properties {\rm(1a)}, {\rm(1b)}, {\rm(1c)} and
{\rm(1d)} for $i=0,1,\dots,n_{1}-2$. Also, for each $p\ge 1$,
\text{$h_p \buildrel \text{{\rm multiseq}} \over \sim h_1$} and
$h_p\in \C\{y\}[z]$ is an irreducible $W$-poly in $z$ with
coefficients in $\BC\{y\}$ and with multiplicity $n_1$ at $0\in
\C^2$. \ms

{\rm(1a)} We write $R^{(p)}_i=b^{(p)}_iy^{\alpha^{(p)}_{1,i,1}}$
with $a$ unit $b^{(p)}_i$ in $\C\{y\}$ and a positive integer
$\alpha^{(p)}_{1,i,1}$, if exists. For all $p\ge 1$,
$\alpha^{(p)}_{1,0,1}=\a_{1,0,1}$ and $\xi_1= b^{(p)}_0(0)$ where
$\xi_1$ was found to be ${\f 1{d_2}} a_{n-n_1}(0)$ as in
$(h_1,f)=(g_1,f)$ of {\rm(16.3.2)}. \ms

{\rm(1b)} Define a function $\th_1:\N_0\to\N_0$ by $\th_1(t)=t$
where $\N_0$ is the set of nonnegative integers, by the same way as
in Proposition $16.2$. \ms

{\rm(1c)} For all $i=0,1,\dots,n_{j+1}-2$,
$$\align
\th_{1}(\a^{(p)}_{1,i,1})> {n_1-i}. \tag 16.3.4
\endalign$$

{\rm(1d)} For all $i=1,\dots,n_{j+1}-2$,
$$\align
 &\gcd(n_{1},\th_{1}(\a^{(p)}_{1,0,1}))= 1 \quad
 \text{with} \quad \si_1=\a^{(p)}_{1,0,1},
\tag 16.3.5\\
&\text{$\f{\th_{1}(\a^{(p)}_{1,i,1})}{n_{1}-i}>
\f{\th_{1}(\a^{(p)}_{1,0,1})}{n_{1}}$ \quad {or} \quad
$n_1\a^{(p)}_{1,i,1}+\a^{(p)}_{1,0,1}i>n_1\a^{(p)}_{1,0,1}$.}
\endalign$$

{\rm(2)}{\rm(2a)} For each $p\ge 1$, $h_p=h_{p}(y,z)\in \C\{y\}[z]$
is an irreducible $W$-poly in $z$ with coefficients in $\BC\{y\}$
and with $h_p\buildrel \text{{\rm multiseq}} \over \sim h_1=g_{1}$,
and $h_{p} \in$ the type $[1]$ in the sense of Definition 2.5. \ms

{\rm(2b)} $h_p=h_p(f_{-1},f_0)\in \C\{f_{-1}\}[f_{0}]$ of
{\rm(16.3.3)} is an irreducible $W$-poly in $f_{0}$ with
coefficients in $\C\{f_{-1}\}$ and with multiplicity $n_{1}$ at
$0\in\C^{2}$. \ms

$\underline{ \text{\bf [2] The Necessary Condition[B] for
$f(y,z)\in$the type$[\ell]$ with $\ell\ge 2$: }}$

\noindent$\underline{ \text{\bf $f\in \C\{f_{-1},f_0\}[h_p]$ is an
irreducible W-poly of degree $d_2$ in $h_p$ with a coefficient }}$

\noindent$\underline{ \text{\bf of ${h_p}^{d_2-1}$ either zero or
nonzero in $\C\{f_{-1},f_0\}$, and $f(y,z)\in$the type[${\ell}$]
with $\ell\ge 2$ }}$

\noindent$\underline{ \text{\bf in the sense of Definition 2.5}}$
\ms

To find \text{\rm the Necessary Condition[B] for $f(y,z)\in$the
type$[\ell]$ with $\ell\ge 2$}, it suffices to show that for each
$p\ge 1$, $f =h^{d_2}_{p}+\sum^{d_2-1}_{i=0} T^{(p)}_{2,i}h^i_{p},$
of {\rm(16.3.3)} satisfies two properties {\rm(1)} and {\rm(2)}:
Note that either $\ell=2$ or $\ell> 2$. \ms

{\rm(1)} Each $T^{(p)}_{2,i}(f_{-1},f_0)\in \C \{f_{-1},f_0\}$ of
$f$ in {\rm(16.3.3)} satisfies {\rm(1a)}, {\rm(1b)}, {\rm(1c)} and
{\rm(1d)} for $i=0,1,\dots,d_{2}-1$.

{\rm(1a)} For any nonzero monomial $\Pi^{2}_{k=1}f^{\g_k}_{k-2}$
 in $T^{(p)}_{2,i}$,
$$\align
\text{$\g_1>0$ \quad and \quad $\g_2<n_1$.}  \tag 16.3.6
\endalign$$

In particular, if $i=d_{2}-1$ for $T^{(p)}_{2,i}$ then $\g_1>0$ and
$\g_{2}\le n_{1}-2$. \ms

{\rm(1b)} Define $\ol\th_{2}(t_k)^{2}_{k=1}= t_{2}\si_1+n_{1}t_1$
for any $(t_k)^{2}_{k=1}\in N^{2}_0$ by the same way as we have seen
in \text{\rm {Proposition $16.2$}, (1b) of The Necessary
Condition[B] for $f(y,z)\in$the type$[\ell]$ with $\ell\ge 2$}.

For any two nonzero monomials $\Pi^{2}_{k=1}f^{\beta_k}_{k-2}$ and
$\Pi^{2}_{k=1}f^{\g_k}_{k-2}$ in $T^{(p)}_{2,i}$,
$$\align
(16.3.7) \qquad \qquad
&\text{$\ol\th_{2}(\beta_k)^{2}_{k=1}=\ol\th_{2}(\g_k)^{2}_{k=1}$ if
and only if $\beta_k=\g_k$ for $k=1,2$.} \qquad \qquad \qquad \qquad
\\
&\text{So, there is a unique nonzero monomial
$C^{(p)}_{2,i}\Pi^{2}_{k=1}f^{\beta^{(p)}_{2,i,k}}_{k-2}$ in
$T^{(p)}_{2,i}$} \\
&\text{with a constant $C^{(p)}_{2,i}$ such that
$\ol\th_{2}(\beta^{(p)}_{2,i,k})^{2}_{k=1}
=\text{$\min$}\{\ol\th_{2}(\g_k)^{2}_{k=1}\}$}\\
&\text{for any nonzero monomial $\Pi^{2}_{k=1}f^{\g_k}_{k-2}$ in
$T^{(p)}_{2,i}$.}
\endalign$$

{\rm(1c)} For all $i=0,1,\dots,d_{2}-1$,
$$\align
\ol\th_{2}(\beta^{(p)}_{2,i,k})^{2}_{k=1}>(d_{2}-i)
n_{1}\th_{1}(\si_1). \tag 16.3.8
\endalign$$

{\rm(1d)} For all $i=0,1,\dots,d_{2}-1$,
$$\align
&\gcd(d_2,\ol\th_{2}(\beta^{(p)}_{2,0,k})^{2}_{k=1})\ge 1,  \tag 16.3.9\\
& \f{\ol\th_{2}(\beta^{(p)}_{2,i,k})^{2}_{k=1}}{d_{2}-i}\ge
\f{\ol\th_{2}(\beta^{(p)}_{2,0,k})^{2}_{k=1}}{d_{2}}.
\endalign$$

Then, either
$\gcd(d_{2},\ol\th_{2}(\beta^{(p)}_{2,0,k})^{2}_{k=1})=1$ or
$1<\gcd(d_{2},\ol\th_{2}(\beta^{(p)}_{2,0,k})^{2}_{k=1})\le d_{2}$.
\ms

{\rm(1d-1)} Let
$\gcd(d_{2},\ol\th_{2}(\beta^{(p)}_{2,0,k})^{2}_{k=1})=1$. Then, $f$
is irreducible in ${}_2\CO_0$ if and only if the inequality in
$(16.3.9)$ holds. In this case, $f \in$ the type $[2]$ in the sense
of Definition 2.5, but note that $T^{(p)}_{2,d_{2}-1}$ may not be
zero where
$$
h_{p}=h_{p-1}+\f 1{d_{2}}T^{(p-1)}_{2,d_{2}-1} \quad \text{and}
\quad f=h^{d_{2}}_p+\sum^{d_{2}-1}_{i=0}T^{(p)}_{2,i}h^i_{p}. \tag
16.3.10
$$

{\rm(1d-2)}  Let
$1<\gcd(d_{2},\ol\th_{2}(\beta^{(p)}_{2,0,k})^{2}_{k=1})\le d_{2}$.
There is a positive integer $\nu$ with $\nu\le \f{n_{1}+1}2$ such
that $T^{(\nu +1)}_{2,d_{2}-1}=T^{(\nu +2)}_{2,d_{2}-1}=\cdots=0$
and $T^{(p)}_{2,d_{2}-1}\ne 0$ for $1\le p\le \nu$. In this case, $f
\in$ the type $[\ell]$ with $\ell\ge 2$ in the sense of Definition
2.5 and note that
$$
1\le \gcd(d_{2},\ol\th_{2}(\beta^{(\nu+1)}_{2,0,k})^{2}_{k=1}<d_{2}.
\tag 16.3.11
$$

{\rm(2)}{\rm(2a)} $f=f(y,z)\in \C\{y\}[z]$ is an irreducible
$W$-poly of degree $n$ in $z$ with coefficients in $\BC\{y\}$ and
with multiplicity $n=d_2n_1$ at $0\in\C^2$. Also, either $f \in$ the
type $[2]$ or $f \in$ the type $[\ell]$ with $\ell\ge 3$ in the
sense of Definition 2.5.

{\rm(2b)} $f=f(f_{-1},f_0,h_{\nu+1})\in \C\{f_{-1},f_0\}[h_{\nu+1}]$
of {\rm(16.3.3)} is an irreducible $W$-poly of degree $d_2$ in
$h_{\nu+1}$ with coefficients in $\C\{f_{-1},f_0\}$ and with
multiplicity $d_{2}$ at $0\in\C^{3}$ where $y,z,h_{\nu+1}$ are
viewed as independent complex three variables at the origin in
$\C^{3}$.
\endproclaim \ms

\noindent$\underline{\text{\bf Remark 16.3.0.}}$ {\bf(a)} As a
corollary of Proposition 16.3, following the same properties and
notations in $\underline{\text{\bf Case(2)}}$ of a conclusion of
Proposition 16.3, let $f_1=h_r$ for any $r\ge \nu+1$ where $\nu+1$
is the smallest positive integer such that $T^{(\nu+1)}_{2,d_2-1}$
is equal to zero.

For notation, if $p=\nu+1$, it is clear that $f=f(f_{-1},f_0,h_p)\in
\C \{f_{-1},f_0\}[h_p]$ satisfies the following up to the change of
notations:

{\rm(i)} $h_p$ of $(h_p,f)$ in (16.3.3) satisfies the same kind of
The Necessary and Sufficient Condition$[A]$ for $h_p(y,z)\in$the
type$[1]$ as $f_1$ of $(f_1,f)$ in $(16.3.1)$ does in \text{\rm
$\widehat{\widehat{\text{\rm[1]}}}$} up to the change of notations.

{\rm(ii)}  $f$ of $(h_p,f)$ in $(16.3.3)$ satisfies the same kind of
The Necessary Condition[B] for $f(y,z)\in$the type$[\ell]$ with
$\ell\ge 2$ as $f$ of $(f_1,f)$ in $(16.3.1)$ does in \text{\rm
$\widehat{\widehat{\text{\rm[2]}}}$} up to the change of notations:
Recall that $f_{-1}=y$ and $f_0=z$ up to the change of notations.
\ms

{\bf(b)} Let $(h_p,f)$ of (16.3.3) be defined as above, as we have
seen in $\underline{\text{\bf Case(2)}}$. In addition, assuming that
$p=\nu+1$ and
$\gcd(d_2,\ol\th_{2}(\beta^{(p)}_{2,0,k})^{2}_{k=1})=1$, then the
statement of [2] in (16.3.3) can be replaced by the statement of
\text{\rm $\widehat{\widehat{\text{\rm[2]}}}$} in (16.3.1), as
follows:

$\underline{ \text{\rm [2] The Necessary Condition[B] for
$f(y,z)\in$the type$[\ell]$ with $\ell\ge 2$: }}$

\noindent$\underline{ \text{\rm $f\in \C\{f_{-1},f_0\}[h_p]$ is an
irreducible W-poly of degree $d_2$ in $h_p$ with a coefficient }}$

\noindent$\underline{ \text{\rm of ${h_p}^{d_2-1}$ either zero or
nonzero in $\C\{f_{-1},f_0\}$, and $f(y,z)\in$the type[${\ell}$]
with $\ell\ge 2$ }}$

\noindent$\underline{ \text{\rm in the sense of Definition 2.5}}$
\ms

$\underline{ \text{\rm $\widehat{\widehat{\text{\rm[2]}}}$  The
Necessary and Sufficient Condition for $f(y,z)\in$the type$[2]$:}}$

\noindent$\underline{ \text{\rm $f\in \C\{f_{-1},f_0\}[h_p]$ is an
irreducible W-poly of degree $d_2$ in $h_p$ with a coefficient }}$

\noindent$\underline{ \text{\rm of ${h_p}^{d_2-1}$ zero in
$\C\{f_{-1},f_0\}$, and $f(y,z)\in$the type[${2}$] in the sense of
Definition 2.5}}$ \ms

\proclaim{Theorem 16.4(A complete algorithm for finding an
irreducibility criterion of arbitrary W-polys in $\BC\{Y,Z\}$ which
have the same multiplicity sequences as the standard Puiseux
expansion($y=t^{n}$ and $z=t^{\alpha}+t^{\beta}$) does)}

$\underline{\text{\bf Assumptions}}$ \quad Let
$f(y,z)=z^n+\sum^{n-2}_{i=0}a_iy^{\a_i}z^i$ be an irreducible
$W$-poly in $z$ with multiplicity $n\ge 2$ at $0\in\C^2$ where for
$0\le i\le n-2$, each $a_i=a_i(y)$ is a unit in ${}_2\CO_0$ if
exists, and the $\a_i$ are positive integers. Note that $a_{n-1}$ is
identically zero. Assume that $d_{2}=gcd(n,\a_0)>1$. Write
$n=d_{2}n_1$ and $\a_0=d_{2} \a_{1,0,1}$ with
$gcd(n_1,\a_{1,0,1})=1$. Note that $2\le n_1<\a_{1,0,1}$. Without
any need of proof, we may assume by Proposition $16.2$ and
Definition $16.2.2$ that $(g_1,f)$ can be written in the form
$$\cases
g_1 &=z^{n_1}+\xi_1 y^{\a_{1,0,1}}\ \text{with}\ \xi_1 =\f {1}{d_{2}} a_{n-n_1}(0), \\
f &=g^{d_{2}}_1+\sum^{d_{2}-1}_{i=1}T_{2,i}g^i_1,
\endcases \tag 16.4.1
$$
satisfying the same properties and notations as in the conclusion of
Proposition $16.2$. \ms

$\underline{\text{\bf Conclusions}}$ \quad Then, $(f_1,f)$ can be
uniquely written in the form
$$\cases
f_{1} &=z^{n_{1}}+\sum^{n_{1}-2}_{i=0} R_{1,i}z^i \quad \text{with}
\ f_{-1}=y\ \text{and}\ f_0=z,
\\
f &=f^{d_{2}}_{1} +\sum^{d_{2}-2}_{i=0} S_{2,i}f^i_{1},
\endcases \tag 16.4.2
$$
where, considering $f_{-1},f_{0},f_{1}$ as independent complex
3-variables at the origin in $\C^{3}$,

{\rm(i)} $n=d_{2}n_1$ with $d_{2}\ge 2$ and $n_{1}\ge 2$, and
$n=d_1$ if necessary;

{\rm(ii)} $R_{1,i}=R_{1,i}(y)\in \C\{y\}$ for each
\text{$i=0,1,\dots,n_1-2$}, and

{\rm(iii)} $S_{2,i}=S_{2,i}(y,z) \in \C\{y\}[z]$ for each
\text{$i=0,1,\dots,d_{2}-2$}, and

{\rm(iv)} $f_{1}=f_{1}(y,z)\in \C\{y\}[z]$;

{\rm(v)} \text{$f=f(y,z,f_1)\in
\C\{y\}[z,f_1]\subseteq\C\{y,z\}[f_1]$} in {\rm(16.4.2)},\ms
{\noindent} satisfying two algorithms, \text{\rm an algorithm for
finding an irreducibility criterion for} \text{\rm $f_{1}(y,z)\in$}
\text{the type$[1]$} in the sense of Definition 2.5 and \text{\rm an
algorithm for finding an irreducibility criterion} \noindent
\text{for $f(y,z)\in$the type$[2]$} in the sense of Definition 2.5,
each of which is represented as follows: Write $f_{-1}=y$ and
$f_0=z$, if necessary. \ms

$\underline{ \text{\bf [1] An algorithm for finding an
irreducibility criterion for $f_{1}(y,z)\in$the type$[1]$:}}$

\noindent$\underline{ \text{\bf $f_{1}\in \BC\{f_{-1}\}[f_{0}]$ is
an irreducible W-poly of degree $n_1$ in $f_0$ with a coefficient
}}$

\noindent$\underline{ \text{\bf of  ${f_{0}}^{n_1-1}$ zero, and
$f_{1}(y,z)\in$ the type[1] in the sense of Definition 2.5}}$

Note that $f \in$ the type$[\ell]$ with $\ell\ge 2$.\ms

Let $d_2=\gcd(n,\a_0)$ with $n=d_2n_1$ and $\a_0=d_2\a_{1,0,1}$ for
some integers $n_1$ and $\a_{1,0,1}$.

Then, $f_1 =f^{n_1}_{0}+\sum^{n_2}_{i=0} R_{1,i}f^i_{0}\in
\C\{y\}[f_0]$ of {\rm (16.4.2)} can be viewed as an element in
$\BC\{y,z\}$ if necessary, satisfying the following {\rm(1)} and
{\rm(2)}: \ms

{\rm(1)} Each $R_{1,i}\in\C\{y\}$ in $f_1$ of {\rm(16.4.2)}
satisfies {\rm(1a)}, {\rm(1b)}, {\rm(1c)} and {\rm(1d)} for
\text{\rm i=0,1,\dots,$n_1$-2}.

{\rm(1a)} For $0\le i\le n_1-2$, each $R_{1,i}=b_iy^{\a_{1,i,1}}$
with a unit $b_i\in\C \{y\}$ and a positive integer $\a_{1,i,1}$ if
exists. Denote $A_{1,i}$ by $b_i(0)$ for convenience of notations.

{\rm(1b)} Define a function $\th_1:\N_0\to\N_0$ by $\th_1(t)=t$
where $\N_0$ is the set of nonnegative integers.

{\rm(1c)} ${\th_1(\a_{1,i,1})}>({n_1}-i)$ for all $i=0,1,\dots,
n_1-2$. \ms

{\rm(1d)} For all $i=0,1,\dots, n_1-2$,
$$\align
&\gcd(n_1,\a_{1,0,1})=1 \quad \text{and} \tag 16.4.3\\
&\f{\th_1(\a_{1,i,1})}{n_1-i}=\f{\a_{1,i,1}}{n_1-i}>
\f{\a_{1,0,1}}{n_1}=\f{\th_1(\a_{1,0,1})}{n_1}.
\endalign$$

{\rm(2)}{\rm(2a)} $f_1=f_1(y,z)\in \C\{y\}[z]$ is an irreducible
$W$-poly in $z$ with coefficients in $\BC\{y\}$ and with
$f_1\buildrel \text{{\rm multiseq}} \over \sim h_1=g_{1}$, and $f_1
\in$ the type $[1]$ in the sense of Definition 2.5.

{\rm(2b)}$f_1 =f^{n_1}_{0}+\sum^{n_2}_{i=0} R_{1,i}f^i_{0}\in
\C\{f_{-1}\}[f_{0}]$ of $(16.4.2)$ is an irreducible $W$-poly in
$f_0$ with coefficients in $\C\{f_{-1}\}$ and with multiplicity
$n_{1}$ at $0\in\C^{2}$.\ms

$\underline{ \text{\bf [2] An algorithm for finding an
irreducibility criterion for $f(y,z)\in$the type$[2]$:}}$

\noindent$\underline{ \text{\bf  $f\in \C\{f_{-1},f_0\}[f_{1}]$ is
an irreducible W-poly of degree $d_2$ in $f_{1}$ with a coefficient
}}$

\noindent$\underline{ \text{\bf of ${f_1}^{d_2-1}$ zero in
$\C\{f_{-1},f_0\}$, and $f=f(y,z)\in$ the type[${\ell}$] with
$\ell\ge 2$}}$

\noindent$\underline{ \text{\bf in the sense of Definition 2.5}}$
\ms

Note that ${2}\le{\ell}$ where $j$ was already given by
{\rm(16.4.1)} and that $n=\Pi^{r}_{i=1}q_i$.

Then, $f =f^{d_{2}}_1+\sum^{{d_{2}}-2}_{i=0}S_{2,i}f^i_1$ of
{\rm(16.4.2)} can be viewed as an element in $\BC\{y,z\}$ if
necessary, satisfying the following {\rm(1)} and {\rm(2)}: \ms

{\rm(1)} Each $S_{2,i}\in\C\{y,z\}$ of $f$ in {\rm(16.4.2)}
satisfies {\rm(1a)}, {\rm(1b)}, {\rm(1c)} and {\rm(1d)} for
$i=0,1,\dots,d_{2}-2$.

{\rm(1a)} For any nonzero monomial $y^{\de_1}z^{\de_2}$ in
$S_{2,i}$,
$$
\de_1>0\quad \text{and}\quad \de_2<n_1. \tag 16.4.4
$$

{\rm(1b)} Define a function $\th_2 :\N^2_0\to \N_0$ by
$\th_2(t_k)^2_{k=1}=t_2\th_1(\a_{1,0,1})+n_1\th_1(t_1)
=t_2\a_{1,0,1}+n_1t_1$ for each $(t_1,t_2)\in \N^2_0$, by the same
way as in Proposition $16.4.2$.

Then, for any two nonzero monomials
$\Pi^{j+1}_{k=1}f^{\beta_k}_{k-2}$ and
$\Pi^{j+1}_{k=1}f^{\de_k}_{k-2}$ in $S_{2,i}$ with $i$ fixed,
$$\align
\text{\rm(16.4.5)} \qquad  &
\th_{2}(\beta_k)^{2}_{k=1}=\th_{2}(\de_k)^{2}_{k=1}\ \text{if and
only if} \ \beta_k=\de_k\ \text{for}\ k=1,2.
\qquad \qquad \qquad \qquad \qquad \\
&\text{So, there exists a unique nonzero monomial
$B_{2,i}\Pi^{2}_{k=1}f^{\beta_{2,i,k}}_{k-2}$ in
$S_{2,i}$}\\
&\text{with a nonzero constant $B_{2,i}$ such that
$\th_{2}(\beta_{2,i,k})^{2}_{k=1}
=\text{$\min$}\{\th_{2}(\de_k)^{2}_{k=1}\}$}\\
&\text{for any nonzero monomial $\Pi^{2}_{k=1}f^{\de_k}_{k-2}$ in
$S_{2,i}$ with $i$ fixed.}
\endalign$$

{\rm(1c)} For all $i=0,1,\dots, d_{2}-2$,
$$\align
(16.4.\text{\rm 6-1}) \qquad
{\th_{2}(\beta_{2,i,k})^{2}_{k=1}}>(d_{2}-i)
n_2\th_2(\beta_{2,0,k})^2_{k=1} \quad \text{for all} \quad
i=0,1,\dots, d_{2}-2. \qquad \qquad
\endalign$$

{\rm(1d)} For all $i=0,1,\dots,d_{2}-2$,
$$\align
& \gcd(d_{2},\th_{2}(\beta_{2,0,k})^{2}_{k=1})\ge 1 \quad \text{and} \tag 16.4.6-2\\
& \f{\th_{2}(\beta_{2,i,k})^{2}_{k=1}}{d_{2}-i}\ge
\f{\th_{2}(\beta_{2,0,k})^{2}_{k=1}}{d_{2}}.
\endalign$$

Then, either $\gcd(d_{2},\th_{2}(\beta_{2,0,k})^{2}_{k=1})=1$ or
$1<\gcd(d_{2},\th_{2}(\beta_{2,0,k})^{2}_{k=1})<d_{2}$. \ms

{\rm(1d-1)} Let $\gcd(d_{2},\th_{2}(\beta_{2,0,k})^{2}_{k=1})=1$.
Then, $f$ is irreducible in ${}_2\CO_0$ if and only if the
inequality in $(16.4.6.2)$ holds and $f \in$ the type $[2]$ in the
sense of Definition 2.5. \ms

{\rm(1d-2)}  Let $1<\gcd(d_{2},\th_{2}(\beta_{2,0,k})^{2}_{k=1})<
d_{2}$. If $f$ is irreducible in ${}_2\CO_0$, then $f \in$ the type
$[\ell]$ with $\ell\ge 3$ in the sense of Definition 2.5. \ms

{\rm(2)}{\rm(2a)} By {\rm(1d)} we can find a complete algorithm for
finding an irreducibility criterion of Weierstrass polynomials in
$\BC\{y,z\}$ which have the same multiplicity sequences as the
standard Puiseux expansion{\rm($y=t^{n}$ and
$z=t^{\alpha}+t^{\beta}$)} does. \ms

{\rm(2b)} $f=f(y,z)\in \C\{y\}[z]$ is an irreducible $W$-poly of
degree $n$ in $z$ with a coefficient of ${z}^{n-1}$ zero in
$\BC\{y\}$ and with multiplicity $n=d_{2}n_1$ at $0\in\C^2$. Also,
either $f \in$ the type $[2]$ or $f \in$ the type $[\ell]$ with
$\ell\ge 3$ in the sense of Definition 2.5. \ms

{\rm(2c)} $f=f(y,z,f_1)\in \C\{y,z\}[f_{1}]$ of {\rm(16.4.2)} is an
irreducible $W$-poly of degree $d_{2}$ in $f_1$ with a coefficient
of ${f_1}^{d_2-1}$ zero in $\C\{y,z\}$ and with multiplicity $d_{2}$
at $0\in\C^{3}$.
\endproclaim \ms

{\bf \S{16.3}. The proofs of Theorem 16.4 with Proposition 16.2 and
Proposition 16.3} \ms

\noindent{\bf Proof of Proposition 16.2.} We prove by (16.2.1) that
$g_{1}$ of $(g_{1},f)$ satisfies The Necessary and Sufficient
Condition$[A]$ for $g_{1}(y,z)\in$the type$[1]$ and that $f$ of
$(g_{1},f)$ satisfies The Necessary Condition$[B]$ for
$f(y,z)\in$the type$[\ell]$ with $\ell\ge 2$, respectively.\ms

$\underline{ \text{\rm The proof of The Necessary and Sufficient
Condition[A] for $g_{1}(y,z)\in$the type$[1]$:}}$ It is clear. \ms

$\underline{ \text{\rm The proof of The Necessary Condition[B] for
$f(y,z)\in$the type$[\ell]$ with $\ell\ge 2$:}}$ First we show that
(1a), (1b), (1c) and (1d) of (1) are true, and next that (2a) and
(2b) of (2) are true.

(1) We show that (1a), (1b), (1c) and (1d) are true, respectively.

{\rm(1a)} Apply the WDT with a divisor $g_1=z^{n_1}+\xi_1
y^{\a_{1,0,1}}$ to $f$ where $\xi_1 = \f {1}{d_2} a_{n-n_1}(0)$.

{\rm(1a-1)} By (1) of Theorem 15.2, $f$ may be written uniquely in
the form
$$
f=g^{d_2}_1 +\sum^{{d_2}-1}_{i=0} T_{2,i}g^i_1, \tag 16.2.5
$$
where each $T_{2,i}=T_{2,i}(y,z)\in \C\{y\}[z]$ is a polynomial in
$z$ of degree $<n_1$ such that $T_{2,i}(0,z)=0$. Thus, the proof of
(1a-1) is done.

{\rm(1a-2)} Assume the contrary. Then, we get that $\de_2 =n_1-1$
for a nonzero monomial $y^{\de_1}z^{\de_2}\in T_{2,d_2-1}$, and so
it would be a nonzero monomial
$y^{\de_1}z^{\de_2}(z^{n_1})^{d_2-1}\in T_{2,d_2-1}g^{d_2-1}_1$ as
in $\C\{y\}[z]$. But, note that
$y^{\de_1}z^{\de_2}(z^{n_1})^{d_2-1}=y^{\de_1}z^{n-1}$ because
$\de_2+n_1(d_2-1)=n_1-1+n_1(d_2-1)=n-1$. Now, we claim that there is
a nonzero monomial $cy^{\de_1}z^{n-1}\in f(y,z)$ for some nonzero
number $c$. If the claim is proved, then it would be a contradiction
to the assumption that $a_{n-1}$ is zero, and so there is nothing to
prove for (1a-2).

For the proof of the claim, first it is needed to compute any
nonzero monomial $y^{\g_1}z^{\g_2}\in T_ig^i_1$  for $1\le i\le
d_2-2$, and also any nonzero monomial $y^{\g'_1}z^{\g'_2}\in
g^{d_2}_1$, as in $\C\{y\}[z]$.

For any nonzero monomial $y^{\g_1}z^{\g_2}\in T_ig^i_1$ with $1\le
i\le d_2-2$,
$$
\g_2\le n_1-1+n_1i=n_1(i+1)-1\le n_1(d_2-1)-1<n-n_1\le n-2. \tag
16.2.6
$$

If a nonzero monomial $y^{\g'_1}z^{\g'_2}$ with $\g'_1>0$ is in
$g^{d_2}_1=(z^{n_1}+\xi y^{\a_{101}})^{d_2}$, it is clear that
$\g'_2\le n-2$ because $2\le n_1$. Thus, the claim is proved, and so
the proof of (1a-2) is done.\ms

(1b) It suffices to show that if
$\th_2(\g_1,\g_2)=\th_2(\de_1,\de_2)$ where $0<\g_1$, $0<\de_1$,
$0\le \g_2<n_1$ and $0\le \de_2<n_1$, then $\g_i=\de_i$ for $i=1,2$.
By definition, $n_1\g_1 +\a_{1,0,1}\g_2=n_1\de_1+\a_{1,0,1}\de_2$,
i.e., $n_1(\de_1-\g_1)=\a_{1,0,1}(\g_2-\de_2)$. Since
$\gcd(n_1,\a_{1,0,1})=1$ and $|\g_2-\de_2|<n_1$, it is clear that
$\g_1=\de_1$ and $\g_2=\de_2$. Thus, the proof of (1b) is done. \ms

(1c) By (1b), let $C_{2,i}y^{\g_{2,i,1}}z^{\g_{2,i,2}}$ be a unique
monomial in $T_{2,i}$ with a constant $C_{2,i}$ such that
$\th_2(\g_{2,i,1}, \g_{2,i,2})=\text{$\min$}\{\th_2(\de_1,\de_2)\}$
for any nonzero monomial $y^{\de_1}z^{\de_2}$ in $T_{2,i}$ with $i$
fixed.

To prove an inequality in (16.2.3), it suffices to prove that the
following are true:

{\rm(1c-1)} For each $i=0,1,\dots,d_2-1$,
$y^{\g_{2,i,1}}z^{\g_{2,i,2}}z^{in_1}$ belongs to $T_{2,i}g^i_1$.

{\rm(1c-2)} $\th_2(\g_{2,i,1},
\g_{2,i,2})+in_1\a_{1,0,1}>d_2n_1\a_{1,0,1}$ for all
$i=0,1,\dots,d_2-1$. \ms

To prove {\rm(1c-1)}, note that for any nonzero monomial
$y^{\g}z^{\de} \in T_ig^i_1$, $n_1\g+\a_{1,0,1}\de \ge
n_1\g_{2,i,1}+\a_{1,0,1}\g_{2,i,2}+in_1\a_{1,0,1}=\th_2(\g_{2,i,1},
\g_{2,i,2})+in_1\a_{1,0,1}$, because
$y^{\g}z^{\de}=y^{\de_1}z^{\de_2}y^{q\a_{1,0,1}}z^{pn_1}$ for a
nonzero monomial $y^{\de_1}z^{\de_2}\in T_{2,i}$ where
$g_1=z^{n_1}+\xi y^{\a_{1,0,1}}$ with $p+q=i$. Since $T_ig^i_1\in
\C\{y\}[z]$ is a polynomial in $z$ of degree $(in_1+\alpha)$ where
$0\le \alpha<n_1$, it is clear that
$y^{\g_{2,i,1}}z^{\g_{2,i,2}}z^{in_1}$ belongs to $T_ig^i_1$. Thus,
the proof of {\rm(1c-1)} can be done. Also, it is trivial by
{\rm(1c-1)} that {\rm(1c-2)} is true because
$y^{\g_{2,i,1}}z^{\g_{2,i,2}}z^{in_1}\in \sum^{d_2-1}_{i=0}
T_{2,i}g^i_1$ by Lemma $16.0$ and Theorem $16.1$. Thus, we proved an
inequality in (16.2.3). Then, the proof of (1c) is done. \ms

{\rm(1d)} In preparation for the proof of an inequality in (16.2.4),
whenever a nonzero monomial $y^{\g}z^{\de} \in T_{2,i}g^i_1$ is
chosen arbitrary for each $i=0,1,\dots,d_2-1$, then it was already
proved by (16.2.3) that $n_1\g+\a_{1,0,1}\de>n_1\a_{1,0,1}d_2$.

In order to apply Theorem $16.1$ to the proof of an inequality in
(16.2.4), recall by (16.2.5) that $V(g_1)=\{(y,z):g_1(y,z)=0\}$ is
an analytic variety with isolated singularity at the origin in
$\BC^2$ defined by the form
$$\align
(16.2.7) \quad \qquad g_1 &=z^{n_1}+\xi_1 y^{\a_{1,0,1}} \quad
\text{with $\gcd(n_1,k_1)=1$ and $2\le n_1<k_1=\a_{1,0,1}$}. \qquad \qquad \\
\endalign$$

Let $\tau_{\la_1}$ be the composition of a finite number $\la_1$ of
successive blow-ups which is needed only to get the standard
resolution of the singular point of $V(g_1)$. For each
$t=1,2,\dots,\la_1$, write $\tau_t= \pi_1\circ\pi_2\circ\cdots
\circ\pi_t:M^{(t)}\to \BC^2$ where $\{\text{$ \pi_i:M^{(i)}\to
M^{(i-1)}$ is a blow-up}$ of $M^{(i-1)}$ $\text{at some point for
$1\le i\le t$} \}$ with $M^{(0)}=\BC^2$. For brevity of notation,
let $V^{(t)}(g_1)$ be the proper transform under $\tau_t$ for $1\le
t\le \la_1$.

Let $E^{(\la_1)}=\tau^{-1}_{\la_1}(0,0)$, and let $E^{(\la_1)}=\cup
E_i$, $1 \le i \le \la_1$, be the decomposition of $E^{(\la_1)}$
into irreducible components where each $E_i$ is called an
exceptional curve of the first kind. By Theorem $16.1$, we can use
the following consequences. \ms

{\bf Consequence(1).} In order to study $V^{(t)}(g_1)$ under
$\tau_t$, we can find just one coordinate patch of the local
coordinates for each blow-up $\pi_t:M^{(t)}\to M^{(t-1)}$, where
$1\le t\le \la_1$ and $M^{(0)}=\BC^2$.

{\bf Consequence(2).} By {\rm Consequence(1)}, we can use the same
$\tau_{\la_1}$ for the composition of the first finite number
$\la_1$ of successive blow-ups in preparation for the standard
resolution of the singular point $(0,0)$ of $V(f)$.

{\bf Consequence(3).} In order to study each proper transform of
both $V(g_1)$ and $V(f)$ under $\tau_t$, without using a nonsingular
change of coordinates, we can use the common one coordinate patch of
the same local coordinates simultaneously, as it has been already
used for each blow-up $\pi_t:M^{(t)}\to M^{(t-1)}$ in {\rm
Consequence(1)}, where $1\le t\le \la_1$. \ms

After $\la_1$ iterations of blow-ups, let $(v_{\la_1},u_{\la_1})$
and $(v'_{\la_1},u'_{\la_1})$ be the local coordinates for
$M^{(\la_1)}$ where $\pi_{\la_1}:M^{(\la_1)}\to M^{(\la_1-1)}$ was
defined to be the $\la_1$-th blow-up at some point of
$M^{(\la_1-1)}$ with $u'_{\la_1}=1/u_{\la_1}$ and
$v'_{\la_1}=v_{\la_1}u_{\la_1}$. Note that
$E_{\la_1}=\{v_{\la_1}=0\}\cup \{v'_{\la_1}=0\}$. For brevity of
notation, write $(v,u)=(v_{\la_1},u_{\la_1})$ and
$(v',u')=(v'_{\la_1},u'_{\la_1})$. \ms

Being viewed as an analytic mapping,
$\tau_{\la_1}:M^{(\la_1)}\to\BC^2$ can be written in the form
 $$
 \tau_{\la_1}(v,u)=(y,z)=(v^{n_1}u^a,v^{k_1}u^b)
 \quad \text{with $k_1=\a_{1,0,1}$.} \tag 16.2.8
 $$
where

{\rm(i)} $a$ and $b>0$ are nonnegative integers such that
 $bn_1-ak_1=1$,

{\rm(ii)} $E_{\la_1}=\{v=0\}$ is defined by the ${\la_1}-th$
exceptional curve of the first kind. \ms

By Theorem $16.1$, along $v=0$, $(f\circ\tau_{\la_1})_{total}$ can
be written in the following form:
$$\align
(g_1\circ\tau_{\la_1})_{total}&=v^{n_1k_1}u^{ak_1}(u+\xi_1), \tag 16.2.9\\
 (f\circ\tau_{\la_1})_{total}&=v^{e_{\la_1}}u^{\varepsilon}
(f\circ\tau_{\la_1})_{proper},\\
(f\circ\tau_{\la_1})_{proper}&=(u+\xi_1)^{d_2}+\sum^{d_2-1}_{i=0}T'_i(u+\xi_1)^i,
\quad \text{with}\\
T'_i &=T'_i(u,v)=b_i v^{M'_i}, \\
M'_i &=\th_2(\g_{2,i,1}, \g_{2,i,2})+in_1k_1-d_2n_1k_1>0,
\endalign$$

where \roster \item "(i)" $e_{\la_1}=n_1k_1d_2$ and $\ve=ak_1d_2$,

\item "(ii)" $\xi_1$ is a nonzero constant and each $b_i$ is a unit
in $\C \{u+\xi_1,v\}$.
\endroster \ms

Since $(f\circ\tau_{\la_1})_{proper}$ is irreducible in
$\C\{u+\xi,v\}$, by Theorem $3.2$ and (16.2.9) we have an inequality
for each $i=1,2,\dots,d_2-1$:
$$
\f{M'_i}{d_2-i}\ge \f{M'_0}{d_2},\quad \text{that is}, \quad
\f{\th_2(\g_{2,i,1}, \g_{2,i,2})}{d_2-i}\ge
\f{\th_2(\g_{2,0,1},\g_{2,0,2})}{d_2}. \tag 16.2.10
$$

Thus, we proved the inequality in (16.2.4).

(1d-1) Suppose $\gcd(d_2,\th_2(\g_{2,0,1},\g_{2,0,2}))=1$. It is
clear by Corollary 3.3 and Theorem 3.7.

(1d-2) Suppose $1<\gcd(d_2,\th_2(\g_{2,0,1},\g_{2,0,2}))\le d_2$. In
order to prove that if $T_{2,d_2-1}$ is zero then
$\gcd(d_2,\th_2(\g_{2,0,1},\g_{2,0,2}))<d_2$, it suffices to
consider the following three subcases (i), (ii) and (iii) of (1d-2):
Note that $\gcd(d_2,\th_2(\g_{2,0,1}, \g_{2,0,2}))=\gcd(d_2,M'_0)$
where $M'_0 =\th_2(\g_{2,0,1}, \g_{2,0,2})-d_2n_1k_1>0$ by (16.2.9).

(i) of (1d-2). Let $d_2<M'_0$. Apply Theorem 3.2 to
$(f\circ\tau_{\la_1})_{proper}$ in (16.2.10). Since $T_{2,d_2-1}$ is
zero, it is clear that $\gcd(d_2,M'_0)<d_2$.

(ii) of (1d-2). Let $d_2=M'_0$. By Theorem 3.2, $T_{d_2-1}$ cannot
be zero.

(iii) of (1d-2). Let $d_2>M'_0$. It is clear.

Also, note by Theorem 3.2 that $\gcd(d,M'_0)$ may be equal to $M'_0$
whether $T_{d-1}$ is zero or not. Thus, the proof of (1d) is done.
\ms

(2)(2a) The proof is trivial by (1d) of (1).

(2b) Note that $T_{2,i}g^i_1\in
\C\{y\}[z,g_1]\subseteq\C\{y,z\}[g_1]$ for $0\le i\le d_2-1$ and
$f\in \C\{y,z\}[g_1]$, considering $y,z,g_1$ as independent three
complex variables. Let $y^{\de_1}z^{\de_2}g^i_1\in T_{2,i}g^i_1$ be
arbitrary for $0\le i\le d_2-1$ where $\de_1>0$ and $\de_2<n_1$.
Then, $n_1\a_{1,0,1}(\de_1+\de_2+i)>n_1\de_1 +\a_{1,0,1}\de_2 +
in_1\a_{1,0,1} \ge \th_2(\g_{2,i,1},
\g_{2,i,2})+in_1\a_{1,0,1}>(d_2-i)n_1\a_{1,0,1}+in_1\a_{1,0,1}=d_2n_1\a_{1,0,1}$
by (b) and (1c). Thus, $\de_1+\de_2+i>d_2$. So, $f(y,z,g_1)$ is an
irreducible $W$-poly of degree $d_2$ in $g_1$ with coefficients in
$\BC\{y,z\}$ and with multiplicity $d$ at the origin in $\C^3$
because ${}_3\CO_0$ is a unique factorization domain and $f\in
\C\{y\}[z]$ is irreducible in $\C\{y,z\}$. Thus, the proof of (2b)
is done. Therefore, the proof of this proposition is done. $\square$
\ms

\demo{\bf Proof of Proposition 16.3} For the construction of a pair
$(f_{1},f)$ in (16.3.1),  it suffices to consider the following two
cases, depending on the fact that $T_{2,d_{2}-1}$ of (16.3.0) is
either zero or not. For brevity of notations, let $h_1=g_{1}$ and
$T^{(1)}_{2,i}=T_{2,i}$ for $0\le i\le d_2-1$. \ms

$\underline{\text{\bf Case(1):}}$ Let $T^{(1)}_{2,d_{2}-1}$ be zero.
It is clear. \ms

$\underline{\text{\bf Case(2):}}$ Let $T^{(1)}_{2,d_{2}-1}$ be
nonzero. It has been already shown by Sublemma 15.5 and Sublemma
15.6 in the proof of Theorem 15.4 that the following assertion is
true:

There is a sequence of W-polys in z of pairs,
$\{(h_p,f):p=1,2,\dots\}$ such that
$$\align
(16.3.12) \qquad \qquad (h_{\nu+1},f)=(h_{\nu+2},f)=\cdots \quad
\text{for some integer} \quad \nu\le \f{n_1+1}2, \qquad \qquad
\qquad \qquad
\endalign$$
each pair of which can be written in the form
$$\cases
h_1 &=g_1=z^{n_1}+\xi_1 y^{\a_{1,0,1}}, \\
f &=h^{d_2}_1+\sum^{d_2-1}_{i=0}T^{(1)}_{2,i}h^i_1,
\endcases \tag 16.3.13
$$
and for each $p=2,3,\dots$
$$\cases
h_{p} &=h_{p-1}+\f1{d_2}T^{(p-1)}_{2,{d_2}-1}
=z^{n_1}+\sum^{n_1-2}_{i=0}R^{(p)}_{1,i}z^i, \\
f &=h^{d_2}_{p} +\sum^{{d_2}-1}_{i=0}T^{(p)}_{2,i}h^i_{p},
\endcases \tag 16.3.14
$$
with $T^{(p)}_{2,{d_2}-1}\ne 0$ for $1\le p\le \nu$ and
$T^{(\nu+1)}_{2,{d_2}-1}=T^{(\nu+2)}_{2,{d_2}-1}=\cdots =0$ where
$T^{(p)}_{2,i}=T^{(p)}_{2,i}(y,z)\in \C\{y\}[z]$ for $p\ge 1$ and
$0\le i\le {d_2}-1$, and $R^{(p)}_{2,i} =R^{(p)}_{2,i}(y)\in
\C\{y\}$ for $p\ge 1$ and $0\le i\le n_1-2$, if exist, satisfying
the same kind of the properties and notations as we have seen in
Sublemma $15.5$ of Theorem $15.4$, as follows: \ms

\noindent$\underline{\text{\rm (16.3.15-1) Property(1) and
Property(3)}}$ Let $p$ and $i$ be fixed with $p\ge 1$ and $0\le i\le
n_{1}-2$. Then, $R^{(p)}_{1,i}=R^{(p)}_{1,i}(y)\in \C\{y\}$ with
$R^{(p)}_{2,i}(0)=0$ and has a multiplicity $\ge (n_{1}-i)$ at $0
\in \C^2$. \ms

\noindent$\underline{\text{\rm (16.3.15-2)  Property(2) and
Property(4)}}$ Let $p$ and $i$ be fixed with $p\ge 1$ and $0\le i\le
d_2-1$. Then, $T^{(p)}_{2,i}=T^{(p)}_{2,i}(y,z)\in \C\{y\}[z]$ is a
polynomial of degree $\de<n_1$ in $z$ and has a multiplicity $\ge
(d_2-i)n_1$ at $0 \in \C^2$.

Also, for any nonzero monomial $y^{\de_1}z^{\de_2}$ in
$T^{(p)}_{2,i}=T_{2,i}^{(p)}(y,z)\in \C\{y\}[z]$, $\de_1>0$ and
$\de_2<n_1$. \ms

\noindent$\underline{\text{\rm (16.3.15-3) Property(5)}}$ In
particular, if $i=d_2-1$ for $T^{(p)}_{2,i}$ of {\rm Property(4)},
then $\de_{2}\le n_1-2$. \ms

\noindent$\underline{\text{\rm (16.3.15-4) Property(6)}}$ There is
an integer $\nu <n_{1}$ such that $T^{(p)}_{2,d_2-1}\ne 0$ for
$p=1,2,\dots,\nu$ and
$T^{(\nu+1)}_{2,d_2-1}=T^{(\nu+2)}_{2,d_2-1}=\cdots =0$. Also, $\nu
\le \f{n_{1}+1}2$. \ms

{\bf Remark.} Without any need of proof, Property(1),
Property(2),\dots, Property(6), which are mentioned just above,
follow clearly from Sublemma $15.5$ of Theorem $15.4$, which belongs
to Case[II] with $j=0$ in the conclusion of Sublemma $15.5$ of
Theorem $15.4$. In Sublemma $15.5$, note that $f_{-1}=y$ and
$f_{0}=z$. \ms

For the proof of this proposition in Case(2), it suffices to show
that two properties, denoted by, $\underline{ \text{\rm The
Necessary and Sufficient Condition[A] for $h_{p}(y,z)\in$the
type$[1]$}}$ and

\noindent{$\underline{ \text{\rm The Necessary Condition[B] for
$f(y,z)\in$the type$[\ell]$ with $\ell\ge 2$}}$} can be satisfied,
respectively. Then, the proof will be by induction on the integer
$p\ge 1$.

Now, it is enough to consider the following two subcases for
Case(2), respectively:

Subcase(A) $p=1$, and Subcase(B) $p\ge 1$. \ms

$\underline{\text{\bf Subcase(A) of Case(2):}}$ \quad Let $p=1$.
Then it suffices to show that $(h_1,f)$ given by (16.3.2), satisfies
$\underline{ \text{\rm The Necessary and Sufficient Condition[A] for
$g_{1}(y,z)\in$the type$[1]$}}$ and

\noindent{$\underline{ \text{\rm The Necessary Condition[B] for
$f(y,z)\in$the type$[\ell]$ with $\ell\ge 2$}}$}, which was already
proved by Proposition $16.2$. \ms

$\underline{\text{\bf Subcase(B) of Case(2):}}$ \quad Let $p\ge 1$.
For the proof of this subcase, it suffices to show by Subcase(A) of
Case(2) that the following sublemma is true: \ms

\noindent\text{\bf Sublemma 16.3.1 for Subcase(B) of Case(2).} \quad

$\underline{\text{\bf {Assumptions}}}$ \quad For the induction
proof, suppose we have shown on the integer $p\ge 1$ that
$\underline{ \text{\rm  The Necessary and Sufficient Condition[A]
for $h_{p}(y,z)\in$the type$[1]$}}$ and

$\underline{ \text{\rm The Necessary Condition[B] for $f(y,z)\in$the
type$[\ell]$ with $\ell\ge 2$ }}$ are true for $(h_p,f)$, following
the same notations and properties as we have seen in (16.3.3),
(16.3.4), \dots, (16.3.11).
 \ms

$\underline{\text{\bf {Conclusions}}}$ \quad Then, $(h_{p+1},f)$ can
be written by
$$\cases
h_{p+1} &=h_p+\f1{d_2}
T^{(p)}_{2,{d_2}-1}=z^{n_1}+\sum^{n_1-2}_{i=0}R^{(p+1)}_{1,i}z^i, \\
f &=h^{d_2}_{p+1} +\sum^{{d_2}-1}_{i=0}T^{(p+1)}_{2,i}h^i_{p+1},
\endcases \tag 16.3.16
$$
with $T^{(p)}_{2,{d_2}-1}\ne 0$ for $1\le p\le \nu$ and
$T^{(\nu+1)}_{2,{d_2}-1}=T^{(\nu+2)}_{2,{d_2}-1}=\cdots =0$ where
$T^{(p)}_{2,i}=T^{(p)}_{2,i}(y,z)\in \C\{y\}[z]$ for $p\ge 1$ and
$0\le i\le {d_2}-1$, and $R^{(p)}_{1,i} =R^{(p)}_{1,i}(y)\in
\C\{y\}$ for $p\ge 1$ and $0\le i\le n_1-2$, if exist, satisfying
the following properties, denoted by The Necessary and Sufficient
Condition[A] for $h_{p+1}(y,z)\in$the type$[1]$ and The Necessary
Condition[B] for $f(y,z)\in$the type$[\ell]$ with $\ell\ge 2$
inductively as we have seen in the conclusion of this proposition:
\ms

$\underline{ \text{\bf [1] The Necessary and Sufficient Condition[A]
for $h_{p+1}(y,z)\in$the type$[1]$:}}$

\noindent$\underline{ \text{\bf{\bf $h_{p+1}\in
\BC\{f_{-1}\}[f_{0}]$} is an irreducible W-poly  of degree $n_1$ in
$f_0$ with a coefficient of }}$

\noindent$\underline{ \text{\bf ${f_{0}}^{n_1-1}$ zero in
$\C\{f_{-1}\}$, and $h_{p+1}(y,z)\in$ the type[1] in the sense of
Definition 2.5}}$ \ms

To find \text{\bf The Necessary and Sufficient Condition[A] for
$h_{p+1}(y,z)\in$the type$[1]$}, as an element in $\BC\{y,z\}$ if
necessary, it suffices to show that
$h_{p+1}=z^{n_1}+\sum^{n_1-2}_{i=0}R^{(p+1)}_{1,i}z^i$ of
{\rm(16.3.16)} satisfies the following {\rm(1)} and {\rm(2)}: \ms

(1) Let $p$ be fixed with $p\ge 1$. Each $R^{({p+1})}_{1,i}\ne 0$
satisfies the properties {\rm(1a)}, {\rm(1b)}, {\rm(1c)} and
{\rm(1d)} for $i=0,1,\dots,n_{1}-2$. Also, for each $p\ge 1$,
\text{$h_{p+1} \buildrel \text{{\rm multiseq}} \over \sim h_1$} and
$h_{p+1}\in \C\{y\}[z]$ is an irreducible $W$-poly in $z$ with
coefficients in $\BC\{y\}$ and with multiplicity $n_1$ at $0\in
\C^2$. \ms

(1a) We write
$R^{({p+1})}_{1,i}=b^{({p+1})}_iy^{\alpha^{({p+1})}_{1,i,1}}$ with
$a$ unit $b^{({p+1})}_i$ in $\C\{y\}$ and a positive integer
$\alpha^{({p+1})}_{1,i,1}$, if exists. For all $p\ge 1$,
$\alpha^{({p+1})}_{1,0,1}=\a_{1,0,1}$ and $\xi_1= b^{(p)}_0(0)$
where $\xi_1$ was found to be ${\f 1{d_2}} a_{n-n_1}(0)$ as in
$(h_1,f)=(g_1,f)$ of {\rm(16.3.2)}. \ms

{\rm(1b)} Define a function $\th_1:\N_0\to\N_0$ by $\th_1(t)=t$
where $\N_0$ is the set of nonnegative integers, by the same way as
in Proposition $16.2$. \ms

{\rm(1c)} For all $i=0,1,\dots,n_{j+1}-2$,
$$\align
\th_{1}(\a^{({p+1})}_{1,i,1})> {n_1-i}. \tag 16.3.17
\endalign$$

(1d) For all $i=0,1,\dots,n_{1}-2$,
$$\align
 &\gcd(n_{1},\th_{1}(\a^{({p+1})}_{1,0,1}))= 1 \quad
 \text{with} \quad \si_1=\a^{({p+1})}_{1,0,1},
\tag 16.3.18\\
&\f{\th_{1}(\a^{({p+1})}_{1,i,1})}{n_{1}-i}\ge
\f{\th_{1}(\a^{({p+1})}_{1,0,1})}{n_{1}}.
\endalign$$
Note that ${\a^{({p+1})}_{1,i,1}}n_1+i\a^{({p+1})}_{1,0,1}>
n_1\a^{({p+1})}_{1,0,1}$ for all $i=1,\dots,n_{1}-2$.

{\rm(2)}{\rm(2a)} For each $p\ge 1$, $h_{p+1}=h_{{p+1}}(y,z)\in
\C\{y\}[z]$ is an irreducible $W$-poly in $z$ with coefficients in
$\BC\{y\}$ and with $h_{p+1}\buildrel \text{{\rm multiseq}} \over
\sim h_1=g_{1}$, and $h_{{p+1}} \in$ the type $[1]$ in the sense of
Definition 2.5. \ms

{\rm(2b)} $h_{p+1}=h_{p+1}(f_{-1},f_0)\in \C\{f_{-1}\}[f_{0}]$ of
{\rm(16.3.16)} is an irreducible $W$-poly in $f_{0}$ with
coefficients in $\C\{f_{-1}\}$ and with multiplicity $n_{1}$ at
$0\in\C^{2}$. \ms

$\underline{ \text{\bf [2] The Necessary Condition[B] for
$f(y,z)\in$the type$[\ell]$ with $\ell\ge 2$: }}$

\noindent$\underline{ \text{\bf $f\in \C\{f_{-1},f_0\}[h_{p+1}]$ is
an irreducible W-poly of degree $d_2$ in $h_{p+1}$ with a
coefficient }}$

\noindent$\underline{ \text{\bf of ${h_{p+1}}^{d_2-1}$ either zero
or nonzero in $\C\{f_{-1},f_0\}$, and $f(y,z)\in$the type[${\ell}$]
with $\ell\ge 2$ }}$

\noindent$\underline{ \text{\bf in the sense of Definition 2.5}}$
\ms

To find \text{\bf The Necessary Condition[B] for $f(y,z)\in$the
type$[\ell]$ with $\ell\ge 2$}, as an element in $\BC\{y,z\}$ if
necessary, it is enough to show that for each $p\ge 1$, $f
=h^{d_2}_{p+1}+\sum^{d_2-1}_{i=0} T^{(p+1)}_{2,i}h^i_{p+1},$ of
{\rm(16.3.16)} satisfies the following {\rm(1)} and {\rm(2)}: Note
that either $\ell=2$ or $\ell> 2$. \ms

{\rm(1)} Each $T^{(p+1)}_{2,i}\in \C \{f_{-1},f_0\}$ of $f$ in
{\rm(16.3.3)} satisfies {\rm(1a)}, {\rm(1b)}, {\rm(1c)} and
{\rm(1d)} for $i=0,1,\dots,d_{2}-1$.

(1a) For any nonzero monomial $\Pi^{2}_{k=1}f^{\g_k}_{k-2}$
 in $T^{(p+1)}_{2,i}$,
$$\align
\text{$\g_1>0$ \quad and \quad $\g_2<n_1$.}  \tag 16.3.19
\endalign$$

In particular, if $i=d_{2}-1$ for $T^{(p+1)}_{2,i}$ then $\g_1>0$
and $\g_{2}\le n_{1}-2$. \ms

(1b) Define $\ol\th_{2}(t_k)^{2}_{k=1}= t_{2}\si_1+n_{1}t_1$ for any
$(t_k)^{2}_{k=1}\in N^{2}_0$ by the same way as we have seen in the
definition of the integer valued-function $\ol\th_{2}$ in the
conclusion of Proposition $16.2$.

For any two nonzero monomials $\Pi^{2}_{k=1}f^{\beta_k}_{k-2}$ and
$\Pi^{2}_{k=1}f^{\g_k}_{k-2}$ in $T^{(p+1)}_{2,i}$,
$$\align
(16.3.20) \qquad \qquad
&\text{$\ol\th_{2}(\beta_k)^{2}_{k=1}=\ol\th_{2}(\g_k)^{2}_{k=1}$ if
and only if $\beta_k=\g_k$ for $k=1,2$.} \qquad \qquad \qquad \qquad
\\
&\text{So, there is a unique nonzero monomial
$C^{(p+1)}_{2,i}\Pi^{2}_{k=1}f^{\beta^{(p+1)}_{2,i,k}}_{k-2}$ in
$T^{(p+1)}_{2,i}$} \\
&\text{with a constant $C^{(p+1)}_{2,i}$ such that
$\ol\th_{2}(\beta^{(p+1)}_{2,i,k})^{2}_{k=1}
=\text{$\min$}\{\ol\th_{2}(\g_k)^{2}_{k=1}\}$}\\
&\text{for any nonzero monomial $\Pi^{2}_{k=1}f^{\g_k}_{k-2}$ in
$T^{(p+1)}_{2,i}$.}
\endalign$$

{\rm(1c)} For all $i=0,1,\dots,d_{2}-1$,
$$\align
\ol\th_{2}(\beta^{(p+1)}_{2,i,k})^{2}_{k=1}>(d_{2}-i)
n_{1}\th_{1}(\si_1). \tag 16.3.21
\endalign$$

(1d) For all $i=0,1,\dots,d_{2}-1$,
$$\align
&\gcd(d_2,\ol\th_{2}(\beta^{(p+1)}_{2,0,k})^{2}_{k=1})\ge 1,  \tag 16.3.22\\
& \f{\ol\th_{2}(\beta^{(p+1)}_{2,i,k})^{2}_{k=1}}{d_{2}-i}\ge
\f{\ol\th_{2}(\beta^{(p+1)}_{2,0,k})^{2}_{k=1}}{d_{2}}.
\endalign$$

Then, either
$\gcd(d_{2},\ol\th_{2}(\beta^{(p+1)}_{2,0,k})^{2}_{k=1})=1$ or
$1<\gcd(d_{2},\ol\th_{2}(\beta^{(p+1)}_{2,0,k})^{2}_{k=1})\le
d_{2}$. \ms

(1d-1) Let
$\gcd(d_{2},\ol\th_{2}(\beta^{(p+1)}_{2,0,k})^{2}_{k=1})=1$. Then,
$f$ is irreducible in ${}_2\CO_0$ if and only if the inequality in
(16.3.22) holds. In this case, $f \in$ the type $[2]$ in the sense
of Definition 2.5, but note that $T^{(p+1)}_{2,d_{2}-1}$ may not be
zero where
$$
h_{p+1}=h_{p}+\f 1{d_{2}}T^{(p)}_{2,d_{2}-1} \quad \text{and} \quad
f=h^{d_{2}}_{p+1}+\sum^{d_{2}-1}_{i=0}T^{(p+1)}_{2,i}h^i_{p+1}. \tag
16.3.23
$$

(1d-2)  Let
$1<\gcd(d_{2},\ol\th_{2}(\beta^{(p+1)}_{2,0,k})^{2}_{k=1})\le
d_{2}$. There is a positive integer $\nu$ with $\nu\le \f{n_{1}+1}2$
such that $T^{(\nu +1)}_{2,d_{2}-1}=0$ and $T^{(p+1)}_{2,d_{2}-1}\ne
0$ for $p+1=1,2,\dots,\nu$. In this case, $f \in$ the type $[\ell]$
with $\ell\ge 3$ in the sense of Definition 2.5 and note that
$$
1<\gcd(d_{2},\ol\th_{2}(\beta^{(\nu+1)}_{2,0,k})^{2}_{k=1}<d_{2}.
\tag 16.3.24
$$

{\rm(2)}{\rm(2a)} $f=f(y,z)\in \C\{y\}[z]$ is an irreducible
$W$-poly of degree $n$ in $z$ with coefficients in $\BC\{y\}$ and
with multiplicity $n=d_2n_1$ at $0\in\C^2$. Also, either $f \in$ the
type $[2]$ or $f \in$ the type $[\ell]$ with $\ell\ge 3$ in the
sense of Definition 2.5.

{\rm(2b)} $f=f(f_{-1},f_0,h_{p+1})\in \C\{f_{-1},f_0\}[h_{p+1}]$ of
{\rm(16.3.16)} is an irreducible $W$-poly of degree $d_2$ in
$h_{p+1}$ with coefficients in $\C\{f_{-1},f_0\}$ and with
multiplicity $d_{2}$ at $0\in\C^{3}$. \ms

By $\underline{\text{\rm The Necessary Condition[A] for
$f(y,z)\in$the type$[\ell]$ with $\ell\ge 2$ }}$, let $f_1=h_r$ for
any $r\ge \nu+1$ where $\nu+1$ is the smallest positive integer such
that $T^{(\nu+1)}_{2,d_2-1}$ is equal to zero. \ms

$\underline{ \text{\bf $\widehat{\widehat{\text{\bf[3]}}}$ The
Necessary Condition[A] for $f(y,z)\in$the type$[\ell]$ with $\ell\ge
2$: }}$

\noindent$\underline{ \text{\bf $f\in \C\{f_{-1},f_0\}[f_1]$ is an
irreducible W-poly of degree $d_2$ in $f_1$ with a coefficient }}$

\noindent$\underline{ \text{\bf of ${f_1}^{d_2-1}$ zero in
$\C\{f_{-1},f_0\}$, and $f=f(y,z)\in$the type[${\ell}$] with
$\ell\ge 2$ }}$

\noindent$\underline{ \text{\bf in the sense of Definition 2.5}}$
$\square$ \ms

\noindent$\underline{\text{\bf The proof of Sublemma 16.3.1 for
Subcase(B) of Case(2)}}$ Let $p\ge 1$ with $T^{(p)}_{2,d_2-1}\not
=0$. Suppose we have shown that $\underline{ \text{\rm The Necessary
and Sufficient Condition[A] for $h_{p}(y,z)\in$the type$[1]$}}$ and
$\underline{ \text{\rm The Necessary Condition[B] for $f(y,z)\in$the
type$[\ell]$ with $\ell\ge 2$}}$ are true for $(h_p,f)$.

We will prove that $h_{p+1}$ of $(h_{p+1},f)$ satisfies The
Necessary and Sufficient Condition[A] for \text{$h_{p+1}(y,z)\in$the
type$[1]$} and that $f$ of $(h_{p+1},f)$ satisfies The Necessary
Condition[B] for \text{$f(y,z)\in$the type$[\ell]$} with $\ell\ge
2$, respectively. \ms

$\underline{\text{\rm [1] The proof of  The Necessary and Sufficient
Condition[A] for \text{$h_{p+1}(y,z)\in$the type$[1]$}}}$

(1) (1a) It is clear by either (16.3.15-1) or Sublemma 15.5.5 of
Theorem 15.4.

(1b) There is nothing to prove.

(1c) Since $h_{p+1}=h_p+\f1{d_2}T^{(p)}_{2,d_2-1}$, it suffices to
show by induction on $p$ that for any nonzero monomial
$y^{\alpha}z^{\beta} \in T^{(p)}_{2,{d_2}-1}$, $\alpha+\beta>n_1$.
Since  $T^{(p)}_{2,d_2-1}$ of $(h_p,f)$ satisfies (16.3.7) and an
inequality in (16.3.8) by induction on $p$, $n_1<\a_{1,0,1}$ implies
that
$$\align
(16.3.25) \quad \a_{1,0,1}(\alpha+\beta)>
n_1\alpha+\a_{1,0,1}\beta=\th_2(\alpha,\beta)\ge
\th_2(\beta^{(p)}_{2,d_2-1,1},\beta^{(p)}_{2,d_2-1,2})>
n_1\a_{1,0,1}, \qquad \qquad
\endalign$$
for any nonzero monomials $y^{\alpha}z^{\beta}$ in
$T^{(p+1)}_{2,d_2-1}$. Thus, it is clear that $\alpha+\beta>n_1$.

(1d) By induction on $p$, it is clear by (16.3.25) that for any
nonzero monomial $y^{\alpha}z^{\beta} \in T^{(p+1)}_{2,{d_2}-1}$
with $\beta=0$, $\alpha>\a_{1,0,1}$. So,
$\a^{({p+1})}_{1,0,1}=\a_{1,0,1}$.

For any nonzero monomial $y^{\alpha}z^{\beta} \in
T^{(p+1)}_{2,{d_2}-1}$, write $i=\beta$. Then the proof of (1d) can
be done because it is clear by (16.3.25) that for all
$i=1,\dots,n_{1}-2$,
$$\align
 &\gcd(n_{1},\th_{1}(\a_{1,0,1}))= 1 \quad
 \text{with} \quad \a^{({p+1})}_{1,0,1}=\a_{1,0,1}=\sigma_1,
\tag 16.3.26\\
&\text{$\f{\alpha}{n_{1}-\beta}> \f{{\a}_{1,0,1}}{n_{1}}$ \quad {or}
\quad  $n_1\alpha+\a_{1,0,1}\beta>n_1\a_{1,0,1}$}. \\
\endalign$$

{\rm(2)} The proofs of (2a) and (2b) are clear by (1). Thus, the
proof of  The Necessary and Sufficient Condition[A] for
\text{$h_{p+1}(y,z)\in$the type$[1]$} is done. \ms

$\underline{ \text{\rm [2] The proof of The Necessary Condition[B]
for \text{$f(y,z)\in$the type$[\ell]$} with $\ell\ge 2$}}$

(1) (1a) It is clear by either (16.3.15-2) and (16.3.15-3), or
Sublemma 15.5.5 of Theorem 15.4.

(1b) Using the same method as we have used in the proof of (16.2.2)
in (1b) of The Necessary Condition[B] for \text{$f(y,z)\in$the
type$[\ell]$} with $\ell\ge 2$, it is clear that (16.3.7) in (1b) of
The Necessary Condition[B] for \text{$f(y,z)\in$the type$[\ell]$}
with $\ell\ge 2$ is true. So, it is clear that there exists a unique
nonzero monomial
$C^{(p+1)}_{2,i}y^{\beta^{(p+1)}_{2,i,1}}z^{\beta^{(p+1)}_{2,i,2}}\in
T^{(p+1)}_{2,i}$ with a constant $C^{(p+1)}_{2,i}$ such that
$\th_2(\beta^{(p+1)}_{2,i,1},\beta^{(p+1)}_{2,i,2})=\min\{\th_2(\g_1,\g_2)\}$
for any nonzero monomial $y^{\g_1}z^{\g_2}$ in $T^{(p+1)}_{2,i}$.

(1c) In preparation for the proof of (16.3.21), note by (16.0.1) and
(16.3.26) that for any nonzero monomial $y^{\a}z^{\beta} \in
h_{p+1}-h_1$ and for any nonzero monomial $y^{\a'}z^{\beta'} \in
h_{p+1}^i$ with $p\ge 1$,
$$\align
\th_2(\a,\beta)>n_1\a_{1,0,1} \quad \text{and} \quad
\th_2(\a',\beta')\ge in_1\a_{1,0,1}. \tag 16.3.27
\endalign$$

First note that $f(y,z)$ can be rewritten as follows:
$$\align
(16.3.28) \qquad \qquad f&=h^{d_2}_{p+1}
+\sum^{{d_2}-1}_{i=0}T^{(p+1)}_{2,i}h^i_{p+1}=(h_{p+1}-h_1+h_1)^{d_2}
+\sum^{{d_2}-1}_{i=0}T^{(p+1)}_{2,i}h^i_{p+1}, \qquad \qquad \\
&=(h_{p+1}-h_1)^{d_2}+h_1^{d_2}
+\sum^{{d_2}-1}_{i=1}{}_{d_2}C_i(h_{p+1}-h_1)^{i}h^{{d_2}-i}_{1}
+\sum^{{d_2}-1}_{i=0}T^{(p+1)}_{2,i}h^i_{p+1}
\endalign$$

Since $\ol\th_{2}(\beta^{(p+1)}_{2,i,k})^{2}_{k=1}
=\text{$\min$}\{\ol\th_{2}(\g_k)^{2}_{k=1}\}$ for any nonzero
monomial $\Pi^{2}_{k=1}f^{\g_k}_{k-2}$ in $T^{(p+1)}_{2,i}$ by
(16.3.20), note by (16.3.27) that for any nonzero monomial
$y^{\g}z^{\de} \in T^{(p+1)}_{2,i}h^i_{p+1}$,
$$\align
(16.3.29) \quad n_1\g+\a_{1,0,1}\de \ge
n_1\beta_{2,i,1}+\a_{1,0,1}\beta_{2,i,2}+in_1\a_{1,0,1}=\th_2(\beta_{2,i,1},
\beta_{2,i,2})+in_1\a_{1,0,1}, \qquad \qquad
\endalign$$
because
$y^{\g}z^{\de}=y^{\delta_1}z^{\delta_2}y^{t\a_{1,0,1}}z^{sn_1}$ for
a nonzero monomial $y^{\de_1}z^{\de_2}\in T^{(p+1)}_{2,i}$ where
$s+t=i$.

By (16.0.1) of Lemma 16.0, (16.3.27), (16.3.28) and (16.3.29), to
prove an inequality in (16.3.21), we show by the induction on the
integer $q$ that the following is true:
$$\align
\th_2(\beta_{2,d_2-q,k})^2_{k=1}+qn_1\a_{1,0,1}>d_2n_1\alpha_{1,0,1}
\tag 16.3.30 \endalign$$

Let $q=1$ for the proof of (16.3.30). It is clear by (16.3.19) that
$y^{\beta^{(p+1)}_{2,d_2-1,1}}z^{\beta^{(p+1)}_{2,d_2-1,2}}z^{n_1}\in
T^{(p+1)}_{2,d_2-1}h^{d_2-1}_{p+1}$ belongs to $f-h^{d_2}_{p+1}$.
So, it is clear that an inequality in (16.3.21) is true. Let $q=2$
for the proof of (16.3.30). To prove an inequality in (16.3.30), if
$y^{\beta^{(p+1)}_{2,d_2-2,1}}z^{\beta^{(p+1)}_{2,d_2-2,2}}z^{2n_1}\in
T^{(p+1)}_{2,d_2-2}h^2_{p+1}$ belongs to $f-h^{d_2}_{p+1}$, there is
nothing to prove by (16.0.1) of Lemma 16.0. Otherwise, it is clear
by (16.3.19) that
$y^{\beta^{(p+1)}_{2,d_2-2,1}}z^{\beta^{(p+1)}_{2,d_2-2,2}}z^{2n_1}$
would belong to $T^{(p+1)}_{2,d_2-1}h_{p+1}$, which implies by
(16.3.29) that $\th_2(\beta_{2,d_2-2,k})^2_{k=1}+2n_1\a_{1,0,1}\ge
\th_2(\beta_{2,d_2-1,k})^2_{k=1}+n_1\a_{1,0,1}
>d_2n_1\alpha_{1,0,1}$. Thus, the proof of inequality for $q=2$ is
done. If $q\ge 3$, by the same method as we have used for $q=1$ and
$q=2$ and by induction on the integer $q$, the proof of (16.3.21) is
easily done.

{\rm(1d)} For the proof of an inequality in (16.3.22), by the same
method as we have used in the proof of Proposition $16.2$, there is
nothing to prove.

{\rm(2)} The proof is trivial by (1d) of (1). Thus, the proof of The
Necessary Condition[B] for \text{$f(y,z)\in$the type$[\ell]$} with
$\ell\ge 2$ is done, and then the proof of Sublemma 16.3.1 for
Subcase(B) of Case(2) is done. So, the proof for Case(2) is finished
and therefore, the proof of the proposition can be completed.
$\square$
\enddemo \ms

\demo{\bf Proof of Theorem 16.4} The remaining proof of the theorem
just follows from Proposition 16.2 and Proposition 16.3 with Remark
16.3.0. $\square$
\enddemo \bs

\vfill \pagebreak

{\bf $\S$ {17.} Irreducibility criterion of W-polys of two complex
variables(A generalized representation of irreducible {W}-polys of
two complex variables)} \ms

Observe by Theorem 1.13(Irreducibility criterion of W-polys of two
complex variables) of $\S$ {1.7} in Part[A] that the necessary and
sufficient condition for $f(y,z)$ to be irreducible in $\BC\{y,z\}$
with $f\in {\text{\rm the type}}[{\ell}]$ in the sense of Definition
2.5 can be represented without proof. In this section, the first aim
is to prove that the sufficiency of the condition in Theorem 1.13
can be proved by Theorem 16.5. The second aim is to prove that the
necessity of the condition in Theorem 1.13(the converse of Theorem
16.5) can be represented by Theorem 16.6 with proof. Then, we can
find The 2nd Algorithm for computing irreducibility criterion of all
the W-polys of two complex variables in the process of the proof of
Theorem 16 with Proposition 16.7 and Proposition 16.8 completely and
rigorously, using Theorem 15.4. \ms

\proclaim{Theorem 16.5(A sufficient condition for any W-poly
$f(y,z)$ to be irreducible in $\BC\{y,z\}$ in terms of The Division
Algorithm for W-polys(Irreducibility criterion of W-polys of two
complex variables))}

$\underline{\text{\bf {Assumptions}}}$ Let
$f(y,z)=z^n+\sum^{n-2}_{i=0}a_iy^{\a_i}z^i$ be a $W$-poly in $z$
with multiplicity $n\ge 2$ at $0\in\C^2$ where for $0\le i\le n-2$,
each $a_i=a_i(y)$ is a unit in ${}_2\CO_0$ if exists, and the $\a_i$
are positive integers. Note that $a_{n-1}$ is identically zero. Note
that $a_{n-1}$ is identically zero for convenience. Write
$n=\Pi^{j+1}_{k=1}n_k$ with positive integers $n_k\ge 2$ for all $k$
where the $n_k$ may not be the factorization of prime numbers.

In addition, the following inequality holds:
$$  2\le n\le \alpha_0. \tag 16.5.0 $$

Suppose by {\rm {Theorem 15.4(The Division Algorithm for
$W$-polys)}} that for each fixed positive integer $j$ and for each
$k=1,2,\dots,j$, $f_k$ and $f$ can be written in the form
$$ \cases
f_k &=f^{n_k}_{k-1}+\sum^{n_k-2}_{i=0} R_{k,i}f^i_{k-1}  \\
f& =f^{n_{j+1}}_j +\sum^{n_{j+1}-2}_{i=0} R_{j+1,i}f^i_j
\endcases \tag 16.5.1
$$
where, considering $y,z,f_1,\dots, f_j$ as independent complex
$(j+2)$-variables at the origin in $\C^{j+2}$ with $f_{-1}=y$ and
$f_{0}=z$,

{\rm(i)} $n=\Pi^{j+1}_{k=1}n_k$ with $n_k\ge 2$ for $1\le k\le j+1$;

{\rm(ii)} for each fixed $k$ and for each $i$ with $0\le i\le
n_k-2$, $R_{k,i}\in \C\{y,z,f_1,\dots,f_{k-2}\}$;

{\rm(iii)} for each $k=1,2,\dots,j$, $f_k=f_k(y,z,f_1,\dots,
f_{k-1})\in \C\{y,z,f_1,\dots, f_{k-2}\}[f_{k-1}]$ and

{\rm(iv)} $f=f(y,z,f_1,\dots, f_j)\in \C\{y,z,f_1,\dots,
f_{j-1}\}[f_j]$ with $f=f_{j+1}$,

satisfying a finite number of conditions, each of which is
represented respectively, as follows: \ms

$\underline{ \text{\bf (1) Condition[A] for $f_1(y,z)\in$the type[1]
in the sense of Definition 2.5:}}$

$R_{1,i}\in \C\{y\}$ satisfies the properties {\rm(1a)}, {\rm(1b)}
and {\rm(1c)} for each $i=0,1,\dots,n_1-2,$ and then
$f_1=f_1(y,z)\in \C\{y,z\}$ satisfies the properties {\rm(1d)}. \ms

{\rm(1a)} For $0\le i\le n_1-2$, each $R_{1,i}=b_iy^{\a_{1,i,1}}$
with a unit $b_i\in\C \{y\}$ and a positive integer $\a_{1,i,1}$ if
exists. Denote $A_{1,i}$ by $b_i(0)$ for convenience of notations.

{\rm(1b)} Define a function $\th_1:\N_0\to\N_0$ by $\th_1(t)=t$
where $\N_0$ is the set of nonnegative integers.

{\rm(1c)} ${\th_1(\a_{1,i,1})}>({n_1}-i)$ for all $i=0,1,\dots,
n_1-2$.

{\rm(1d)} For all $i=0,1,\dots, n_1-2$,
$$\align
&\gcd(n_1,\a_{1,0,1})=1 \quad \text{and} \tag 16.5.2\\
&\f{\th_1(\a_{1,i,1})}{n_1-i}=\f{\a_{1,i,1}}{n_1-i}\ge
\f{\a_{1,0,1}}{n_1}=\f{\th_1(\a_{1,0,1})}{n_1}.
\endalign$$

$\underline{ \text{\bf (2) Condition[A] for $f_2(y,z)\in$the type[2]
in the sense of Definition 2.5:}}$

$R_{2,i}\in \C\{y\}[z]$ satisfies the properties {\rm(2a)},
{\rm(2b)} and {\rm(2c)} for each $i=0,1,\dots,n_2-2,$ and then
$f_2=f_2(y,z,f_1)\in \C\{y,z\}[f_1]$ satisfies the properties
{\rm(2d)}. \ms

{\rm(2a)} For any nonzero monomial $y^{\de_1}z^{\de_2}$ in
$R_{2,i}$, $\de_1>0$ and $\de_2<n_1$.

{\rm(2b)} Let $\N^2_0$ be two-dimensional cartesian product of
$\N_0$. For given integers $n_1,\a_{1,0,1}$ and a function $\th_1$
in {\rm Cond[A] of the 1st type}, define a function $\th_2:\N^2_0\to
\N_0$ by
$$
\th_2(t_1,t_2)=t_2\th_1(\a_{1,0,1})+n_1\th_1(t_1)=t_2\a_{1,0,1}+n_1t_1
\quad \text{for each $(t_1,t_2)\in \N^2_0$}. \tag 16.5.3 $$

Then, for any two nonzero monomials $y^{\a_1}z^{\a_2}$ and
$y^{\de_1}z^{\de_2}$ in $R_{2i}$ with $i$ fixed,
$$\align
\text{\rm(16.5.3-1)} \qquad \qquad &
\th_2(\a_1,\a_2)=\th_2(\de_1,\de_2)\ \text{if and only if}\
\a_1=\de_1\ \text{and}\ \a_2=\de_2. \qquad \qquad\qquad \qquad \\
 &\text{So, there exists a unique nonzero monomial
$A_{2,i}y^{\a_{2,i,1}}z^{\a_{2,i,2}}$ in $R_{2,i}$} \\
&\text{with a nonzero constant $A_{2,i}$ such that
$\th_2(\a_{2,i,1},
\a_{2,i,2})=\text{$\min$}\{\th_2(\de_1,\de_2)\}$} \\
&\text{for any nonzero monomial $y^{\de_1}z^{\de_2}$ in $R_{2,i}$
with $i$ fixed.}
\endalign$$

{\rm(2c)} For all $i=0,1,\dots, n_2-2$,
$$\align
{\th_2(\a_{2,i,k})^2_{k=1}}>({n_2}-i) n_{1}{\a_{1,0,1}}. \tag 16.5.4
\endalign$$

{\rm(2d)} For all $i=0,1,\dots, n_2-2$,
$$\align
&\gcd(n_2,\th_2(\a_{2,0,k})^2_{k=1})=1 \quad \text{and} \tag 16.5.5\\
&\f{\th_2(\a_{2,i,1}, \a_{2,i,2})}{n_2-i}\ge
\f{\th_2(\a_{2,0,1},\a_{2,0,2})}{n_2}.
\endalign$$

$\underline{ \text{\bf (3) Condition[A] for $f_m(y,z)\in$the type[m]
in the sense of Definition 2.5:}}$

For each fixed $m$ with $3\le m\le j+1$, $R_{m,i}\in
\C\{y,z,f_1,\dots,f_{m-2}\}$ satisfies the properties {\rm(3a)},
{\rm(3b)} and {\rm(3c)} for each $i=0,1,\dots,n_m-2$, and then
$f_m=f_m(y,z,f_1,\dots,f_{m-1})\in
\C\{y,z,f_1,\dots,f_{m-2}\}[f_{m-1}]$ satisfies the properties
{\rm(3d)}. \ms

{\rm(3a)} For any nonzero monomial $\Pi^m_{k=1}f^{\de_k}_{k-2}$ in
$R_{m,i}$ with $f_{-1}=y$ and $f_0=z$, $\de_1>0$ and $\de_k<n_{k-1}$
for $k=2,3,\dots,m$. \ms

{\rm(3b)} By induction assumption on the integer $(m-1)\le j$, there
exists a sequence $\{f_3,f_4,\dots,f_{m-1}\}$, each of which
satisfies the same kind of properties and notations as we have seen
in $\underline{ \text{\rm Condition[A] for $f_2(y,z)\in$the type[2]
in the sense of Definition 2.5}}$. Then inductively, define
$\th_m:\N^m_0\to \N_0$ where $\N^m_0$ is its $m$-dimensional
cartesian product by
$$
(16.5.6) \quad
\th_m(t_k)^m_{k=1}=t_m\th_{m-1}(\a_{m-1,0,k})^{m-1}_{k=1}
+n_{m-1}\th_{m-1}(t_k)^{m-1}_{k=1} \quad \text{for each
$(t_k)^m_{k=1}\in \N^k_0$},
$$
where recall by induction assumption that for a fixed $i$,
$A_{m-1,i}\Pi^{m-1}_{k=1}f^{\a_{m-1,i,k}}_{k-2}$ is a unique nonzero
monomial in $R_{m-1,i}$ with a constant $A_{m-1,i}$ such that
$$
\th_{m-1}(\a_{m-1,i,k})^{m-1}_{k=1}=\text{$\min$}\{\th_{m-1}(\de_k)^{m-1}_{k=1}\},
\tag 16.5.7
$$
for any nonzero monomial $\Pi^{m-1}_{k=1}f^{\de_k}_{k-2}$ in
$R_{m-1,i}$.

Then, for any two nonzero monomials $\Pi^m_{k=1}f^{\a_k}_{k-2}$ and
$\Pi^m_{k=1}f^{\de_k}_{k-2}$ in $R_{m,i}$ with $i$ fixed,
$$\align
\text{\rm(16.5.7-1)} \qquad \qquad
&\text{$\th_m(\a_k)^m_{k=1}=\th_m(\de_k)^m_{k=1}$
if and only if $\a_k=\de_k$ for $k=1,2,\dots, m.$} \qquad \qquad\qquad \qquad\\
&\text{So, there exists a unique nonzero-monomial
$A_{m,i}\Pi^m_{k=1}f^{\a_{m,i,k}}_{k-2}$ in $R_{m,i}$} \\
&\text{with a constant $A_{m,i}$ such that
$\th_m(\a_{m,i,k})^m_{k=1}=\text{$\min$}\{\th_m(\de_k)^m_{k=1}\}$}\\
&\text{for any nonzero monomial $\Pi^m_{k=1}f^{\de_k}_{k-2}$ in
$R_{m,i}$.}\\
\endalign$$

{\rm(3c)} For all $i=0,1,\dots, n_m-2$,
$$\align
{\th_m(\a_{m,i,k})^m_{k=1}}>({n_m}-i)
n_{m-1}\th_{m-1}(\a_{m-1,0,k})^{m-1}_{k=1}.\tag 16.5.8
\endalign$$

{\rm(3d)} For all $i=0,1,\dots,n_m-2$,
$$\align
&\gcd(n_m,\th_m(\a_{m,0,k})^m_{k=1})=1 \quad \text{and} \tag 16.5.9\\
&\f{\th_m(\a_{m,i,k})^m_{k=1}}{{n_m}-i}\ge
\f{\th_m(\a_{m,0,k})^m_{k=1}}{n_m}.
\endalign$$

$\underline{\text{\bf {Conclusions}}}$ Consider the sequence
$S=\{f_k: 1\le k\le j+1\}$ with $f_{j+1}=f$ where
$f_k=f_k(y,z,\dots,f_{k-1})\in \C\{y,z,f_1,\dots, f_{k-1}\}$, which
have the same properties and notations as we have seen in
{\rm(16.5.1)} of the assumptions of Theorem $16.5$. \ms

{\rm(1)} $f=f_{j+1}(y,z,f_1,\dots,f_{j})\in
\C\{y,z,f_1,\dots,f_{j-1}\}[f_{j}]$ is an irreducible $W$-poly of
degree $n_{j+1}$ in $f_{j}$ with a coefficient of
${f_j}^{n_{j+1}-1}$ zero in $\BC\{y,z,f_1,\dots,f_{j-1}\}$ and with
multiplicity $n_{j+1}$ at the origin in $\C^{j+1}$. \ms

{\rm(2)} $f\in \BC\{y,z\}$ is an irreducible $W$-poly of degree $n$
in $z$ with coefficients in $\BC\{y\}$ and with multiplicity $n$ at
the origin in $\C^{2}$, which can be written as follows:

$f=z^{n}+\sum^{n-2}_{i=0} a_iy^{\a_i}z^i$ is an irreducible $W$-poly
in $z$ with $f_{j+1}\in \text{\rm the type[j+1]}$ in the sense of
Definition of $2.5$ and with multiplicity $n$ at the origin in
$\C^{2}$ where for $0\le i\le n-2$, each $a_i=a_i(y)$ is a unit in
${}_2\CO_0$ if exists, and the $\a_i$ are positive integers, and
$n=\Pi^{j+1}_{i=1}n_i$. Note that $a_{n-1}$ is zero.
\endproclaim \bs

\noindent{\bf Proof of Theorem 16.5.} The proof of the theorem just
follows from Theorem 12.0. \bs

\proclaim{Corollary 16.5.1} $\underline{\text{\bf {Assumptions}}}$
\quad Under the same assumptions and notations as in Theorem $16.5$,
note that $f_k$ is irreducible in $\BC\{y,z\}$ with isolated
singularity at $(0,0)$ in $\BC^2$ for $r\ge 1$. In particular, for
each $k=1,2,\dots,j+1$, let $V(H_{j+1})=\{(y,z):H_{j+1}(y,z)=0\}$ be
an analytic variety at $(0,0)$ in $\BC^2$, each of which is defined
as follows:
$$\align
(16.5.1.1) \qquad \qquad \text{\rm(i)} \qquad
H_1&=z^{n_1}+y^{\alpha_{1,0,1}} \quad \text{with $n_1\ge 2$ and
$\alpha_{1,0,1}\ge 2$}. \qquad \qquad \qquad \qquad \\
\text{\rm(ii)} \qquad H_2&=H^{n_2}_1+y^{\alpha_{2,0,1}}z^{\alpha_{2,0,2}}.\\
\text{\rm(iii)} \qquad
H_3&=H^{n_3}_2+y^{\alpha_{3,0,1}}z^{\alpha_{3,0,2}}H^{\alpha_{3,0,3}}_1. \\
\qquad  &\ldots\ldots \\
\text{\rm(j+1)} \quad
H_{j+1}&=H^{n_{j+1}}_{j}+y^{\alpha_{j+1,0,1}}z^{\alpha_{j+1,0,2}}H^{\alpha_{j+1,0,3}}_1\cdots
H^{\alpha_{j+1,0,j+1}}_{j-1}. \\
\endalign$$

$\underline{\text{\bf {Conclusions}}}$ Then, $f_{j+1}\in \BC\{y,z\}$
and $H_{j+1}\in \BC\{y,z\}$ have the same multiplicity sequence at
$(0,0)$ in $\BC^2$. $\square$
\endproclaim \ms

\definition{Remark 16.5.2}
{\bf(I)} Note by Theorem $5.0$ that $H_{j+1}$ is irreducible in
$\BC\{y,z\}$ with $H_{j+1}\in \text{\rm the type[j+1]}$ in the sense
of Definition of $2.5$ $\iff$ $\gcd(n_1,\alpha_{1,0,1})=1$,
$\gcd(n_2,\theta_2(\alpha_{2,0,k})^{2}_{k=1})=1$,\dots,
$\gcd(n_{j+1},\theta_{j+1}(\alpha_{{j+1},0,k})^{j+1}_{k=1})=1$. \ms

\noindent{\bf(II)} Note that $f_{j+1}$ is irreducible in
$\BC\{y,z\}$ with $f_{j+1}\in \text{\rm the type[j+1]}$ in the sense
of Definition of $2.5$ $\iff$ the following hold:

(1) \qquad $\gcd(n_1,\alpha_{1,0,1})=1$ and
$\dfrac{\th_1(\a_{1,i,1})}{n_1-i}=\dfrac{\a_{1,i,1}}{n_1-i}\ge
\dfrac{\a_{1,0,1}}{n_1}=\dfrac{\th_1(\a_{1,0,1})}{n_1}$.

(2) \qquad $\gcd(n_2,\theta_2(\alpha_{2,0,1},\alpha_{2,0,2}))=1$ and
$\dfrac{\th_2(\a_{2,i,1}, \a_{2,i,2})}{n_2-i}\ge
\dfrac{\th_2(\a_{2,0,1},\a_{2,0,2})}{n_2}$. \ms

\qquad \qquad \qquad $\ldots\ldots$

(j+1) \quad
$\gcd(n_{j+1},\theta_{j+1}(\alpha_{{j+1},0,k})^{j+1}_{k=1})=1$ and
$\dfrac{\th_m(\a_{m,i,k})^m_{k=1}}{{n_m}-i}\ge
\dfrac{\th_m(\a_{m,0,k})^m_{k=1}}{n_m}$.
\enddefinition \bs

\proclaim{Theorem 16.6(How to find an algorithm for computing
irreducibility criterion of all the W-polys defining plane curve
singularities at $0\in \C^2$(The converse of Theorem 16.5))}

Let $N_0$ be the set of nonnegative integers and $N^k_0$ be its
$k$-dimensional copy. Let $r$ be an arbitrary positive integer.

$\underline{\text{\bf {Assumptions}}}$ Let $f=z^n+\sum^{n-2}_{i=0}
a_iy^{\a_i}z^i$ be a $W$-poly in $z$ with multiplicity $n\ge 2$ at
$0\in\C^2$ where for $0\le i\le n-2$, each $a_i=a_i(y)$ is a unit in
${}_2\CO_0$ if exists and the $\a_i$ are positive integers. Note
that $a_{n-1}$ is zero and that $n$ can be viewed as
$n=\Pi^{r}_{i=1}q_i$ where each integer $q_i$ is a prime number. \ms

$\underline{\text{\bf Conclusions}}$ To find a complete algorithm
for computing the irreducibility criterion of the Weierstrass
polynomials $f(y,z)$ in $\BC\{y,z\}$, it suffices to prove the
following:\ms

Whenever $f$ is irreducible in $\BC\{y,z\}$, we can find a positive
integer $\ell$ with $\ell\le r$ such that $f\in $ the type $[\ell]$
in the sense of {\rm Definition 2.5} with the following property,
and conversely: \ms

We may start to assume for brevity of representation that $\ell\ge
2$\text{\rm(that is, $\gcd(n,\alpha_0)>1$)}, because it was already
proved by {\rm Corollary 3.3 and Proposition 16.2} that $f$ is
irreducible in $\BC\{y,z\}$ with $f\in $ the type $[1]$ if and only
if $(16.6.0)$ holds.
$$
\text{$\f{\alpha_i}{n-i}\ge \f{\alpha_0}{n}$ \quad {for} \quad $0\le
i\le n-2$ \quad {and} \quad $\gcd(n,\alpha_0)=1$.} \tag 16.6.0
$$

For each fixed $j$ with $0\le j\le \ell-1$, there exists a sequence
of irreducible $W$-polys in $z$, $\{f_k: k=1,2,\dots,j\}$ where each
$f_k\in \C\{y\}[z]$ is a $W$-polys in $z$ with coefficients in
$\BC\{y\}$ and $f_k\in $ the type $[k]$ and $f_j\ne f$, satisfying
the following properties and notations: Note that $f_{-1}=y$,
$f_0=z$.

\text{For each fixed $j$ with $1\le k\le j$}, $f$ can be uniquely
written in the form
$$
\cases
f_k &=f^{n_k}_{k-1}+\sum^{n_k-2}_{i=0} R_{k,i}f^i_{k-1} \in
\BC\{y,z,f_1,\dots,f_{k-2}\}[f_{k-1}] \\
f &=f^{d_{j+1}}_j +\sum^{d_{j+1}-2}_{i=0} S_{j+1,i}f^i_j  \in
\BC\{y,z,f_1,\dots,f_{j-1}\}[f_{j}],
\endcases \tag 16.6.1
$$
where, considering $y,z,f_1,\dots, f_j$ as independent complex
$(j+2)$-variables at $0\in\C^{j+2}$,

{\rm(i)} $n=d_{j+1}\Pi^j_{k=1}n_k$ with $d_{j+1}\ge 2$ and $n_k\ge
2$ for $1\le k\le j$, and $n=d_1$ if $j=0$;

{\rm(ii)} for each fixed $k$ and for each $i$ with $0\le i\le
n_k-2$, $R_{k,i}\in \C\{y\}[z,f_1,\dots, f_{k-2}]$;

{\rm(iii)} for each $i$ with $0\le i\le d_{j+1}-2$, $S_{j+1,i} \in
\C\{y\}[z,f_1,\dots, f_{j-1}]$;

{\rm(iv)} for each $k=1,2,\dots,j$, $f_k=f_k(y,z,f_1,\dots,
f_{k-1})\in \C\{y\}[z,f_1,\dots, f_{k-1}]$ and

{\rm(v)} $f=f(y,z,f_1,\dots, f_j)\in \C\{y\}[z,f_1,\dots, f_j]$,

\noindent satisfying a finite number of algorithms, each of which is
represented respectively, as follows: \ms

$\underline{ \text{\bf [1] An algorithm for computing an
irreducibility criterion for $f_{1}(y,z)\in$the type$[1]$:}}$

\noindent$\underline{ \text{\bf $f_{1}\in \BC\{f_{-1}\}[f_{0}]$ is
an irreducible W-poly of degree $n_1$ in $f_0$ with a coefficient
}}$

\noindent$\underline{ \text{\bf of  ${f_{0}}^{n_1-1}$ zero, and
$f_{1}(y,z)\in$ the type[1] in the sense of Definition 2.5}}$

Note that $f=f(f_{-1},f_0,f_1) \in$ the type $[\ell]$ with $\ell\ge
2$, as an element in $\BC\{y,z\}$.\ms

Let $d_2=\gcd(n,\a_0)$ with $n=d_2n_1$ and $\a_0=d_2\a_{1,0,1}$ for
some integers $n_1$ and $\a_{1,0,1}$.

Then, $f_1 =f^{n_1}_{0}+\sum^{n_2}_{i=0} R_{1,i}f^i_{0}\in
\C\{y\}[f_0]$ of {\rm (16.6.1)} can be viewed as an element in
$\BC\{y,z\}$ if necessary, satisfying two properties {\rm(1)} and
{\rm(2)}: \ms

{\rm(1)} Each $R_{1,i}\in\C\{y\}$ in $f_1$ of {\rm(16.6.1)}
satisfies {\rm(1a)}, {\rm(1b)}, {\rm(1c)}, {\rm(1d)} for \text{\rm
i=0,1,\dots,$n_1-2$}.

{\rm(1a)} For $0\le i\le n_1-2$, each $R_{1,i}=b_iy^{\a_{1,i,1}}$
with a unit $b_i\in\C \{y\}$ and a positive integer $\a_{1,i,1}$ if
exists. Denote $A_{1,i}$ by $b_i(0)$ for convenience of notations.

{\rm(1b)} Define a function $\th_1:\N_0\to\N_0$ by $\th_1(t)=t$
where $\N_0$ is the set of nonnegative integers.

{\rm(1c)} ${\th_1(\a_{1,i,1})}>({n_1}-i)$ for all $i=0,1,\dots,
n_1-2$. \ms

{\rm(1d)} For all $i=0,1,\dots, n_1-2$,
$$\align
&\gcd(n_1,\a_{1,0,1})=1 \quad \text{and} \tag 16.6.2\\
&\f{\th_1(\a_{1,i,1})}{n_1-i}=\f{\a_{1,i,1}}{n_1-i}\ge
\f{\a_{1,0,1}}{n_1}=\f{\th_1(\a_{1,0,1})}{n_1}.
\endalign$$

{\rm(2)}{\rm(2a)} $f_1=f_1(y,z)\in \C\{y\}[z]$ is an irreducible
$W$-poly of degree $n_1$ in $z$ with coefficients in $\BC\{y\}$ and
with multiplicity $n_1$ at $0\in\C^2$, and $f_1 \in$ the type $[1]$.

{\rm(2b)}$f_1 =f^{n_1}_{0}+\sum^{n_2}_{i=0} R_{1,i}f^i_{0}$ of
$(16.6.1)$ is an irreducible $W$-poly of degree $n_1$ in $f_0$ with
coefficients in $\BC\{f_{-1}\}$ and with multiplicity $n_1$ at
$0\in\C^2$. \ms

$\underline{\text{\rm Remark.}}$ Note that an inequality in
{\rm(16.6.2)} is the necessary and sufficient condition for $f_1$ to
be irreducible in $\C\{y,z\}$ with $f_1 \in$ the type $[1]$. So,
$\text{$f_1 \buildrel \text{{\rm multiseq}} \over \sim H_1$}$ where
$H_1=z^{n_1}+y^{\a_{1,0,1}}$. \ms

$\underline{ \text{\bf [2] An algorithm for computing an
irreducibility criterion for $f_2(y,z)\in$the type$[2]$:}}$

\noindent$\underline{ \text{\bf  $f_2\in \C\{f_{-1},f_0\}[f_{1}]$ is
an irreducible W-poly of degree $n_2$ in $f_{1}$ with a coefficient
of}}$

\noindent$\underline{ \text{\bf ${f_1}^{n_2-1}$ zero in
$\C\{f_{-1},f_0\}$, and $f_2=f_2(y,z)\in$ the type[${2}$] in the
sense of Definition 2.5}}$

Note that $f=f(f_{-1},f_0,f_1,f_2) \in$ the type $[\ell]$ with
$\ell\ge 3$, as an element in $\BC\{y,z\}$. \ms

Then, $f_2 =f^{n_2}_{1}+\sum^{n_2-2}_{i=0} R_{2,i}f^i_{1}\in
\C\{y,z\}[f_1]$ of {\rm (16.6.1)} can be viewed as an element in
$\BC\{y,z\}$ if necessary, satisfying two properties {\rm(1)} and
{\rm(2)}: \ms

{\rm(1)} Each $R_{2,i}\in \C\{y\}[z]$ in $f_2$ of {\rm(16.6.1)}
satisfies {\rm(1a)}, {\rm(1b)}, {\rm(1c)}, and {\rm(1d)} for
\text{\rm i=0,1,\dots,$n_2$-2}.

{\rm(1a)} For any nonzero monomial $y^{\de_1}z^{\de_2}$ in
$R_{2,i}$, $\de_1>0$ and $\de_2<n_1$.

{\rm(1b)} Let $\N^2_0$ be two-dimensional cartesian product of
$\N_0$. For given integers $n_1,\a_{1,0,1}$ and a function $\th_1$
in {\rm Cond[A] of the 1st type}, define a function $\th_2:\N^2_0\to
\N_0$ by
$$
\th_2(t_1,t_2)=t_2\th_1(\a_{1,0,1})+n_1\th_1(t_1)=t_2\a_{1,0,1}+n_1t_1
\quad \text{for each $(t_1,t_2)\in \N^2_0$}. \tag 16.6.3 $$

Then, for any two nonzero monomials $y^{\a_1}z^{\a_2}$ and
$y^{\de_1}z^{\de_2}$ in $R_{2,i}$ with $i$ fixed,
$$\align
\text{\rm(16.6.3-1)} \qquad \qquad &
\th_2(\a_1,\a_2)=\th_2(\de_1,\de_2)\ \text{if and only if}\
\a_1=\de_1\ \text{and}\ \a_2=\de_2. \qquad \qquad\qquad \qquad \\
 &\text{So, there exists a unique nonzero monomial
$A_{2,i}y^{\a_{2,i,1}}z^{\a_{2,i,2}}$ in $R_{2,i}$} \\
&\text{with a nonzero constant $A_{2,i}$ such that
$\th_2(\a_{2,i,1},
\a_{2,i,2})=\text{$\min$}\{\th_2(\de_1,\de_2)\}$} \\
&\text{for any nonzero monomial $y^{\de_1}z^{\de_2}$ in $R_{2,i}$
with $i$ fixed.}
\endalign$$

{\rm(1c)} For all $i=0,1,\dots, n_2-2$,
$$\align
{\th_2(\a_{2,i,k})^2_{k=1}}>({n_2}-i) n_{1}{\a_{1,0,1}}. \tag 16.6.4
\endalign$$

{\rm(1d)} For all $i=0,1,\dots, n_2-2$,
$$\align
&\gcd(n_2,\th_2(\a_{2,0,k})^2_{k=1})=1 \quad \text{and} \tag 16.6.5\\
&\f{\th_2(\a_{2,i,1}, \a_{2,i,2})}{n_2-i}\ge
\f{\th_2(\a_{2,0,1},\a_{2,0,2})}{n_2}.
\endalign$$

Assuming that $f_1$ is an irreducible W-poly with {\bf $f_1$} $\in$
the $type[1]$, then an inequality of {\rm(16.6.5)} is the necessary
and sufficient condition for $f_2$ to be irreducible in $\C\{y,z\}$
with {\rm $f_2$} $\in$ the type $[2]$. \ms

{\rm(2)}{\rm(2a)} $f_2=f_2(y,z)\in \C\{y\}[z]$ is an irreducible
$W$-poly of degree $n_1n_2$ in $z$ with coefficients in $\BC\{y\}$
and with multiplicity $n_1n_2$ at $0\in\C^2$, and $f_2 \in$ the type
$[2]$.

{\rm(2b)} $f_2=f_2(f_{-1},f_0,f_1)\in \C\{f_{-1},f_0\}[f_1]$ of
$(16.6.1)$ is an irreducible $W$-poly of degree $n_2$ in $f_1$ with
coefficients in $\BC\{f_{-1},f_0\}$ and with multiplicity $n_2$ at
$0\in\C^3$. \ms

$\underline{\text{\rm Remark.}}$ By {\rm(1d)}, $\text{$f_2 \buildrel
\text{{\rm multiseq}} \over \sim
H_2=H^{n_2}_1+y^{\a_{2,0,1}}z^{\a_{2,0,2}}$}$ where
$H_1=z^{n_1}+y^{\a_{1,0,1}}$, if irreducible. \ms

$\underline{ \text{\bf [3] An algorithm for computing an
irreducibility criterion for $f_m(y,z)\in$the type$[\text{\bf
{m}}]$:}}$

\noindent$\underline{ \text{\bf  $f_m\in
\C\{f_{-1},f_0,\dots,f_{m-2}\}[f_{m-1}]$ is an irreducible W-poly of
degree $n_m$ in $f_{m-1}$ }}$

\noindent$\underline{ \text{\bf with a coefficient  of
${f_{m-1}}^{n_m-1}$ zero in $\BC\{f_{-1},f_0,\dots,f_{m-2}\}$ }}$

\noindent$\underline{ \text{\bf and $f_m=f_m(y,z)\in$ the
type[${m}$] in the sense of Definition 2.5}}$\ms

Let $m$ be arbitrary with $3\le m\le j$ and with
$j\le{\ell-1}\le{r-1}$. Note that $f=f(f_{-1},f_0,\dots,f_m) \in$
the type $[\ell]$ with $\ell\ge m+1$, as an element in $\BC\{y,z\}$.
\ms

Then, $f_m =f^{n_m}_{m-1}+\sum^{n_m-2}_{i=0} R_{m,i}f^i_{m-1} \in
\BC\{y,z,f_1,\dots,f_{m-2}\}[f_{m-1}]$ of {\rm(16.6.1)} can be
viewed as an element in $\BC\{y,z\}$ if necessary, satisfying two
properties {\rm(1)} and {\rm(2)}: \ms

{\rm(1)} Each $R_{m,i}\in \C \{y\}[z,f_1,\dots,f_{m-2}]$ of $f_m$ in
{\rm(16.6.1)} satisfies {\rm(1a)}, {\rm(1b)}, {\rm(1c)} and
{\rm(1d)} for $i=0,1,\dots,n_m-2$.

{\rm(1a)} For any nonzero monomial $\Pi^m_{k=1}f^{\de_k}_{k-2}$ in
$R_{m,i}$ with $f_{-1}=y$ and $f_0=z$, $\de_1>0$ and $\de_k<n_{k-1}$
for $k=2,3,\dots,m$. \ms

{\rm(1b)} By induction assumption on the integer $(m-1)\le j$,
suppose there exists a sequence $\{f_3,f_4,\dots,f_{m-1}\}$, each of
which satisfies the same kind of properties and notations as $f_2$
does in $\underline{\text{\rm An algorithm for finding an
irreducibility criterion for $f_2(y,z)\in$the type$[2]$}}$.

Inductively, define $\th_m:\N^m_0\to \N_0$ where $\N^m_0$ is its
$m$-dimensional cartesian product by
$$
(16.6.6) \quad
\th_m(t_k)^m_{k=1}=t_m\th_{m-1}(\a_{m-1,0,k})^{m-1}_{k=1}
+n_{m-1}\th_{m-1}(t_k)^{m-1}_{k=1} \quad \text{for each
$(t_k)^m_{k=1}\in \N^k_0$},
$$
where recall by induction assumption that for a fixed $i$,
$A_{m-1,i}\Pi^{m-1}_{k=1}f^{\a_{m-1,i,k}}_{k-2}$ is a unique nonzero
monomial in $R_{m-1,i}$ with a nonzero constant $A_{m-1,i}$ such
that
$$
\th_{m-1}(\a_{m-1,i,k})^{m-1}_{k=1}=\text{$\min$}\{\th_{m-1}(\de_k)^{m-1}_{k=1}\},
\tag 16.6.7
$$
for any nonzero monomial $\Pi^{m-1}_{k=1}f^{\de_k}_{k-2}$ in
$R_{m-1,i}$.

Then, for any two nonzero monomials $\Pi^m_{k=1}f^{\a_k}_{k-2}$ and
$\Pi^m_{k=1}f^{\de_k}_{k-2}$ in $R_{mi}$ with $i$ fixed,
$$\align
\text{\rm(16.6.7-1)} \quad \quad &
\th_m(\a_k)^m_{k=1}=\th_m(\de_k)^m_{k=1}\ \text{if and only if}\
\a_k=\de_k\ \text{for}\ k=1,2,\dots, m. \qquad \qquad \qquad \qquad\\
&\text{So, there exists a unique nonzero-monomial
$A_{m,i}\Pi^m_{k=1}f^{\a_{m,i,k}}_{k-2}$ in $R_{m,i}$}\\
&\text{with a nonzero constant $A_{m,i}$ such that
$\th_m(\a_{m,i,k})^m_{k=1}=\text{$\min$}\{\th_m(\de_k)^m_{k=1}\}$}\\
&\text{for any nonzero monomial $\Pi^m_{k=1}f^{\de_k}_{k-2}$ in
$R_{m,i}$.}
\endalign$$

{\rm(1c)} For all $i=0,1,\dots, n_m-2$,
$$\align
{\th_m(\a_{m,i,k})^m_{k=1}}>({n_m}-i)
n_{m-1}\th_{m-1}(\a_{m-1,0,k})^{m-1}_{k=1}.\tag 16.6.8
\endalign$$

{\rm(1d)} For all $i=0,1,\dots,n_m-2$,
$$\align
&\gcd(n_m,\th_m(\a_{m,0,k})^m_{k=1})=1 \quad \text{and} \tag 16.6.9\\
&\f{\th_m(\a_{m,i,k})^m_{k=1}}{{n_m}-i}\ge
\f{\th_m(\a_{m,0,k})^m_{k=1}}{n_m}.
\endalign$$

Assuming that {\bf $f_{k}$} is an irreducible W-poly with {\bf
$f_{k}$} $\in$ the type$[k]$ for $1\le k\le m-1$, then an inequality
of {\rm(16.6.9)} is the necessary and sufficient condition for $f_m$
to be irreducible in $\C\{y,z\}$ with {\bf $f_{m}$} $\in$ the
type$[m]$. \ms

{\rm(2)}{\rm(2a)} For each $m$ with $3\le m\le j$, $f_m=f_m(y,z)\in
\C\{y\}[z]$ is an irreducible $W$-poly of degree $\Pi^m_{k=1}n_k$ in
$z$ with coefficients in $\BC\{y\}$ and with multiplicity
$\Pi^m_{k=1}n_k$ at $0\in\C^{2}$, and $f_m \in$ the type $[m]$.

{\rm(2b)} $f_m=f_m(f_{-1},f_0,\dots,f_{m-1})\in
\C\{f_{-1},f_0,\dots,f_{m-2}\}[f_{m-1}]$ of {\rm(16.6.1)} is an
irreducible $W$-poly of degree $n_m$ in $f_{m-1}$ with coefficients
in $\C\{f_{-1},f_0,\dots,f_{m-2}\}$ and with multiplicity $n_m$ at
$0\in\C^{m+1}$. \ms

$\underline{\text{\rm Remark.}}$ By {\rm(1d)}, $\text{$f_m \buildrel
\text{{\rm multiseq}} \over \sim H_m$}$ where
$H_m=H^{n_m}_{m-1}+y^{\alpha_{m,0,1}}z^{\alpha_{m,0,2}}H^{\alpha_{m,0,3}}_1\cdots
H^{\alpha_{m,0,m}}_{m-2}$. \ms

$\underline{ \text{\bf [4] An algorithm for computing an
irreducibility criterion for $f(y,z)\in$the type$[j+1]$:}}$

\noindent$\underline{ \text{\bf  $f\in
\C\{f_{-1},f_0,\dots,f_{j-1}\}[f_{j}]$ is an irreducible W-poly of
degree $d_{j+1}$ in $f_{j}$ }}$

\noindent$\underline{ \text{\bf with a coefficient  of
${f_{j}}^{d_{j+1}-1}$ zero in $\BC\{f_{-1},f_0,\dots,f_{j-1}\}$ }}$

\noindent$\underline{ \text{\bf and $f=f(y,z)\in$ the type[${\ell}$]
with $\ell\ge{j+1}$ in the sense of Definition 2.5}}$\ms

Note that ${j+1}\le{\ell}$ where $j$ was already given by
{\rm(16.6.1)} and that $n=\Pi^{r}_{i=1}q_i$.

Then, $f =f^{d_{j+1}}_j +\sum^{d_{j+1}-2}_{i=0} S_{j+1,i}f^i_j  \in
\BC\{y,z,f_1,\dots,f_{j-1}\}[f_{j}]$ of {\rm(16.6.1)} can be viewed
as an element in $\BC\{y,z\}$ if necessary, satisfying two
properties {\rm(1)} and {\rm(2)}: \ms

{\rm(1)} Each $S_{j+1,i}\in\C\{y,z,f_1,\dots,f_{j-1}\}$ of $f$ in
{\rm(16.6.1)} satisfies {\rm(1a)}, {\rm(1b)}, {\rm(1c)} and
{\rm(1d)} for $i=0,1,\dots,d_{j+1}-2$.

{\rm(1a)} For any nonzero monomial $\Pi^{j+1}_{k=1}f^{\de_k}_{k-2}$
in $S_{j+1,i}$ with $f_{-1}=y$ and $f_0=z$,
$$
\de_1>0\ \text{and}\ \de_k<n_{k-1}\quad \text{for}\ k=2,3,\dots,
j+1. \tag 16.6.10
$$

{\rm(1b)} Define a function $\th_{j+1}:\N_0^{j+1}\to \N_0$ where
$\N_0^{j+1}$ is its $(j+1)$-dimensional cartesian product by
$$
(16.6.11) \qquad
\th_{j+1}(t_k)^{j+1}_{k=1}=t_{j+1}\th_j(\a_{j,0,k})^j_{k=1}+n_j\th(t_k)^j_{k=1}
\quad \text{for each $(t_k)^{j+1}_{k=1}\in \N^{j+1}_0$}. \qquad
\qquad
$$

Then, for any two nonzero monomials
$\Pi^{j+1}_{k=1}f^{\beta_k}_{k-2}$ and
$\Pi^{j+1}_{k=1}f^{\de_k}_{k-2}$ in $S_{j+1,i}$ with $i$ fixed,
$$\align
\text{\rm(16.6.11-1)} \quad &
\th_{j+1}(\beta_k)^{j+1}_{k=1}=\th_{j+1}(\de_k)^{j+1}_{k=1}\
\text{if and only if} \ \beta_k=\de_k\ \text{for}\ k=1,2,\dots, j+1.
\qquad \qquad\\
&\text{So, there exists a unique nonzero monomial
$B_{j+1,i}\Pi^{j+1}_{k=1}f^{\beta_{j+1,i,k}}_{k-2}$ in
$S_{j+1,i}$}\\
&\text{with a nonzero constant $B_{j+1,i}$ such that
$\th_{j+1}(\beta_{j+1,i,k})^{j+1}_{k=1}
=\text{$\min$}\{\th_{j+1}(\de_k)^{j+1}_{k=1}\}$}\\
&\text{for any nonzero monomial $\Pi^{j+1}_{k=1}f^{\de_k}_{k-2}$ in
$S_{j+1,i}$ with $i$ fixed.}
\endalign$$

{\rm(1c)} For all $i=0,1,\dots, d_{j+1}-2$,
$$\align
(16.6.12) \qquad
{\th_{j+1}(\beta_{j+1,i,k})^{j+1}_{k=1}}>(d_{j+1}-i)
n_j\th_j(\a_{j,0,k})^j_{k=1} \quad \text{for all} \quad i=0,1,\dots,
d_{j+1}-2. \qquad \qquad
\endalign$$

{\rm(1d)} For all $i=0,1,\dots,d_{j+1}-2$,
$$\align
& \gcd(d_{j+1},\th_{j+1}(\beta_{j+1,0,k})^{j+1}_{k=1}))<d_{j+1} \quad \text{and} \tag 16.6.13\\
& \f{\th_{j+1}(\beta_{j+1,i,k})^{j+1}_{k=1}}{d_{j+1}-i}\ge
\f{\th_{j+1}(\beta_{j+1,0,k})^{j+1}_{k=1}}{d_{j+1}}.
\endalign$$

Then, either
$\gcd(d_{j+1},\th_{j+1}(\beta_{j+1,0,k})^{j+1}_{k=1})=1$ or
$1<\gcd(d_{j+1},\th_{j+1}(\beta_{j+1,0,k})^{j+1}_{k=1})$. \ms

{\rm(1d-1)} Let
$\gcd(d_{j+1},\th_{j+1}(\beta_{j+1,0,k})^{j+1}_{k=1})=1$. Assuming
that {\bf $f_{j}\in \C\{y\}[z]$} is an irreducible W-poly with {\bf
$f_{j}$} $\in$ the type$[j]$, then an inequality of {\rm(16.6.13)}
is the necessary and sufficient condition for $f$ to be irreducible
in $\C\{y,z\}$ with {\bf $f$} $\in$ the type$[j+1]$. In this case,
$f=f_{j+1}$ with {\rm $f_{j+1}\in$ the type$[j+1]$}, after replacing
$d_{j+1}$, $B_{j+1,i}$, $\beta_{j+1,ik}$ and $S_{j+1,i}$ by
$n_{j+1}$, $A_{j+1,i}$, $\a_{j+1,i,k}$ and $R_{j+1,i}$ for
$i=0,1,\dots,n_{j+1}-2$ and $k=1,2,\dots,j+1$, respectively. \ms

{\rm(1d-2)} Let
$1<\gcd(d_{j+1},\th_{j+1}(\beta_{j+1,0k})^{j+1}_{k=1})<d_{j+1}$.
Assuming that {\bf $f_{j}$} is an irreducible W-poly with {\bf
$f_{j}$} $\in$ the type$[j]$, then $f$ is irreducible in $\C\{y,z\}$
with {\bf $f$} $\in$ the type$[\ell]$ with $\ell\ge j+2$. Note that
$\ell\le {r}$.

Therefore, if {\bf $f$} $\in$ the type$[\ell]$ with $\ell\ge j+2$,
by using the induction method on the positive integer $\ell$ as we
have used in {\rm(16.6.1)}, there exists a sequence of irreducible
$W$-polys in $z$, $\{\text{\rm $f_k\in \C\{y\}[z]$: k=1,2,\dots,
j+p}\}$ for some integer $p>1$ such that $f_k\in $ the type $[k]$
for $k=1,2,\dots, j+p$ and $f_{j+p}=f$, satisfying the same kind of
the properties and notations in {\rm the necessary and sufficient
condition for {\rm $f_{m}$} to be an irreducible W-poly with {\rm
$f_{m}$} $\in$ the type[m]} for $m=1,2,\dots,j+p$, by using the same
kind of the statements as we have used in {\rm(1d)} of {\rm[3]} with
finitely many times because $\ell\le {r}$ and $\ell=j+p$. \ms

{\rm(2)}{\rm(2a)} $f=f(y,z)\in \C\{y\}[z]$ is an irreducible
$W$-poly of degree $n$ in $z$ with coefficients in $\BC\{y\}$ and
with multiplicity $n=d_{j+1}\Pi^j_{k=1}n_k$ at $0\in\C^2$. Also,
either $f \in$ the type $[j+1]$ or $f \in$ the type $[\ell]$ with
$\ell\ge j+2$. \ms

{\rm(2b)} $f=f(f_{-1},f_0,\dots,f_j)\in
\C\{f_{-1},f_0,\dots,f_{j-1}\}[f_{j}]$ of {\rm(16.6.1)} is an
irreducible $W$-poly of degree $d_{j+1}$ in $f_j$ with a coefficient
of ${f_{j}}^{d_{j+1}-1}$ zero in $\C\{f_{-1},f_0,\dots,f_{j-1}\}$
and with multiplicity $d_{j+1}$ at $0\in\C^{j+2}$.
\endproclaim \ms

\newpage

\vfill \pagebreak

{\bf \S{18.} The proofs of Theorem 16.6 with Proposition 16.7 and
Proposition 16.8} \bs

\noindent{\bf Proof of Theorem 16.6.} For the proof of the theorem,
let $j$ be chosen arbitrary with $0\le j\le \ell-1$. Then, it
suffices to show that the following are true:

(i) First, for each $m=1,2,\dots,j$, {\bf An algorithm for finding
an irreducibility criterion for $f_m(y,z)\in$the type$[m]$} which
satisfies two properties (1) and (2), is true, as we have seen in
the conclusion of this theorem. \ms

(ii) Next, for $j+1\le \ell$, {\bf An algorithm for finding an
irreducibility criterion for $f(y,z)\in$the type$[j+1]$} which
satisfies two properties (1) and (2), is true, as we have seen in
the conclusion of this theorem. \ms

Now, the proof of the theorem will be by induction on the integer
$j$ with $0\le j\le \ell-1$.

If $\ell=1$, there is nothing to prove by Proposition $3.2$ and
Theorem $3.7$ because $f$ is irreducible in $\C\{y,z\}$ if and only
if $\gcd(n,\a_{0})=1$ and $\f{\a_{i}}{n-i}\ge \f{\a_{0}}{n}$ for all
$i=1,2,\dots,n-2$.

If $\ell\ge 2$ and $j=0,1$, then it was already proved by an
equation in (16.3.3) of Proposition $16.3$ and Theorem $16.4$.

For the induction proof with $\ell\ge 2$, suppose we have shown that
the theorem holds on the integer $j$, with $1\le j\le \ell-1$. Then,
we may assume without any need of proof that $f$ can be uniquely
written in the form
$$\cases
f_k &=f^{n_k}_{k-1}+\sum^{n_k-2}_{i=0} R_{k,i}f^i_{k-1}
\quad \text{for $k=1,2,\dots,j,$} \\
f &=f^{d_{j+1}}_j +\sum^{d_{j+1}-2}_{i=0} S_{j+1,i}f^i_j,
\endcases \tag 16.6.14
$$
satisfying the same properties and notations as in conclusion of the
theorem. \ms

For the remaining proof of the theorem, on the integer $j+1$, it
suffices to show that there exists such a $W$-poly $f_{j+1}\in \C
\{y\}[z]$ of (16.6.15) with $f_{j+1}\in $ the type $[j+1]$
satisfying the same kind of two properties as we have seen in {\bf
An algorithm for finding an irreducibility criterion for
$f_j(y,z)\in$the type$[j]$} which satisfies two properties (1) and
(2) in the theorem, and that $f=f(y,z,f_1,\dots, f_{j+1})$ of
(16.6.15) can be constructed with $f(y,z)\in $ the type $[j+2]$,
satisfying the same kind of two properties as we have seen in {\bf
An algorithm for finding an irreducibility criterion for
$f(y,z)\in$the type$[j+1]$} which satisfies two properties (1) and
(2) in the conclusion of this theorem for $f=f(y,z,f_1,\dots,f_j)$,
according to the following notations, if possible,
$$
\cases
f_{j+1} &=f^{n_{j+1}}_j +\sum^{n_{j+1}-2}_{i=0} R_{j+1,i}f^i_j, \\
f &=f^{d_{j+2}}_{j+1}+\sum^{d_{j+2}-2}_{i=0} S_{j+2,i}f^j_{j+1},
\endcases \tag 16.6.15
$$
where $n=d_{j+2}\Pi^{j+1}_{k=1}n_k$ with $d_{j+2}\ge 2$ and
$n_{k}\ge 2$ for $1\le k\le j+1$. \ms

Now, by the above induction assumption on the integer $j$ and by
(16.6.14), in preparation for applying Theorem $14.0$ and
Proposition $14.1$ to the proof of this theorem, we need to
substitute a new notation, and so let $\{g_k: k=1,2,\dots,j+1\}$ be
in $\BC\{y,z\}$ such that $g_k=f_k$ for $1\le k\le j$ and
$g_{j+1}=f$ where $f=f(y,z,f_1,\dots,f_j)\in
\C\{y\}[z,f_1,\dots,f_j]$, and note that $g_{j+1}$ of Theorem $14.0$
may not be irreducible in $\BC\{y,z\}$ by construction. Then,
$\{g_k: k=1,2,\dots,j+1\}$ of the above substitution satisfies the
same kind of properties and notations as in the assumption of
Theorem $14.0$ and Proposition $14.1$. For the proof of this
theorem, we can use the same kind of conclusion as in Theorem $14.0$
and Proposition $14.1$ without any need of proof, which can be
represented in the following sublemma. \ms

\noindent $\underline{\text{\bf Sublemma 16.6.1}}$ \quad
$\underline{\text{\bf Assumptions}}$ \quad By the above induction
assumption on the integer $j$ and by (16.6.1), \dots,(16.6.13) in
this theorem, for brevity of notation let $\{g_k: k=1,2,\dots,j+1\}$
be in $\BC\{y,z\}$ such that $g_k=f_k$ for $1\le k\le j$ and
$g_{j+1}=f$ where $f_k=f(y,z,f_1,\dots,f_{k-1})\in
\C\{y\}[z,f_1,\dots,f_{k-1}]$ for $k=1,2,\dots,j$ and
$f=f(y,z,f_1,\dots,f_j)\in \C\{y\}[z,f_1,\dots,f_j]$ by
construction. Then, it is clear without any need of proof that for
all $k=1,2,\dots,j+1$, $g_1,g_2,\dots,g_j$ are irreducible in
$\BC\{y,z\}$ but $g_{j+1}$ may not be irreducible in $\BC\{y,z\}$,
which satisfy the same kind of properties and notations as in the
assumption of Theorem $14.0$ and Proposition $14.1$. \ms

$\underline{\text{\bf Conclusions}}$ \quad As a consequence, we can
use the same kind of results and notations as in Theorem $14.0$ and
Proposition $14.1$, as follows:

Let $\tau_{\la_{j}}$ be the composition of a finite number $\la_{j}$
of successive blow-ups which is needed only to get the standard
resolution of the singular point of $V(f_j)$. For each
$t=1,2,\dots,\la_{j}$, write $\tau_t= \pi_1\circ\pi_2\circ\cdots
\circ\pi_t:M^{(t)}\to \BC^2$ where $\pi_i:M^{(i)}\to M^{(i-1)}$ is a
blow-up of $M^{(i-1)}$ at some point of $M^{(i-1)}$ for $1\le i\le
t$ with $M^{(0)}=\BC^2$. For brevity of notation, let $V^{(t)}(f_j)$
be the proper transform under $\tau_t$ for $1\le t\le \la_{j}$.

Let $E^{(\la_{j})}=\tau^{-1}_{\la_{j}}(0,0)$, and let
$E^{(\la_{j})}=\cup E_i$, $1 \le i \le \la_{j}$, be the
decomposition of $E^{(\la_{j})}$ into irreducible components where
each $E_i$ is called an exceptional curve of the first kind.  \ms

$\underline{\text{\rm Consequence(1)}}$ In order to study
$V^{(t)}(f_j)$ under $\tau_t$, we can find just one coordinate patch
of the local coordinates for each blow-up $\pi_t:M^{(t)}\to
M^{(t-1)}$, where $1\le t\le \la_{j}$ and $M^{(0)}=\BC^2$.

$\underline{\text{\rm Consequence(2)}}$ By {\rm Consequence(1)}, we
can use the same $\tau_{\la_{j}}$ for the composition of the first
finite number $\la_{j}$ of successive blow-ups in preparation for
the standard resolution of the singular point $(0,0)$ of both
$V(f_j)$ and $V(f)$.

$\underline{\text{\rm Consequence(3)}}$ In order to study each
proper transform of both $V(f_j)$ and $V(f)$ under $\tau_t$, without
using a nonsingular change of coordinates, we can use the common one
coordinate patch of the same local coordinates simultaneously, as it
has been already used for each blow-up $\pi_t:M^{(t)}\to M^{(t-1)}$
in {\rm Consequence(1)}, where $1\le t\le \la_{j}$. \ms

After $\la_{j}$ iterations of blow-ups, let
$(v_{\la_{j}},u_{\la_{j}})$ and $(v'_{\la_{j}},u'_{\la_{j}})$ be the
local coordinates for $M^{(\la_{j})}$ where by {\rm Consequence(3)}
$\pi_{\la_{j}}:M^{(\la_{j})}\to M^{(\la_{j}-1)}$ was defined to be
the $\la_{j}$-th blow-up at some point of $M^{(\la_{j}-1)}$ with
$u'_{\la_{j}}=1/u_{\la_{j}}$ and
$v'_{\la_{j}}=v_{\la_{j}}u_{\la_{j}}$. Note that
$E_{\la_{j}}=\{v_{\la_{j}}=0\}\cup \{v'_{\la_{j}}=0\}$. For brevity
of notation, write $(v,u)=(v_{\la_{j}},u_{\la_{j}})$ and
$(v',u')=(v'_{\la_{j}},u'_{\la_{j}})$. \ms

Note that $E_{\la_j}=\{(v,u):v=0\}\cup \{(v',u'):v'=0\}$ is the
$j$-th exceptional curve of the first kind. For notation, along
$E_{\la_j}$, $(f_j\circ\tau_{\la_j})_{total}=0$ and
$(f_j\circ\tau_{\la_j})_{proper}=0$ are called  the local defining
equations for the $\la_j$-th total transform of $f_j=0$ and the
$\la_j$-th proper transform of $f_j=0$ under $\tau_{\la_j}$,
respectively.

At $(v,u+a)=(0,0)$ along $v=0$, $(f_j\circ\tau_{\la_j})_{total}=0$
and $(f_j\circ\tau_{\la_j})_{proper}=0$ can be written in the form,
satisfying the following property:
$$\align
(f_j\circ\tau_{\la_j})_{total}&=
v^{n_j\th_j(\a_{j,0,k})^j_{k=1}}(f_j\circ\tau_{\la_j})_{proper}, \tag 16.6.16\\
(f_j\circ\tau_{\la_j})_{proper}&=(u+a+\ve), \\
(\Pi^{j+1}_{k=1}f^{\de_k}_{k-2}\circ\tau_{\la_j})_{total}&=
v^{\th_{j+1}({\de_k})^{j+1}_{k=1}}b(\de_1,\dots,\de_{j+1}),\\
\endalign$$
where $a$ is a nonzero constant, $\ve$ is a nonunit along $v=0$ and
$b(\de_1,\dots,\de_{j+1})$ is a unit in $\C\{v,u+a\}$. \ms

Moreover, as an application of the above consequences we have the
following:

$\underline{\text{\rm Consequence(4)}}$ \quad By (16.6.12) or
(16.6.14) and by the definition of a unique nonzero monomial
$B_{j+1,i}\Pi^{j+1}_{k=1}f^{\beta_{j+1,i,k}}_{k-2}$ in $S_{j+1,i}$
with a constant $B_{j+1,i}$, it can be easily shown that
$(S_{j+1,i}\circ\tau_{\la_j})_{total}=
v^{\th_{j+1}(\beta_{j+1,i,k})^{j+1}_{k=1}}b_{j+1,i}$ where
$b_{j+1,i}$ is a unit in $\C\{v,u+a\}$.

Using $f=f^{d_{j+1}}_j +\sum^{d_{j+1}-2}_{i=0} S_{j+1,i}f^i_j$ in
(16.6.14), then at $(v,u+a)=(0,0)$ along $v=0$,
$(f\circ\tau_{\la_j})_{total}=0$ and
$(f\circ\tau_{\la_j})_{proper}=0$ can be written as follows:
$$\align
\text{\rm(16.6.17)} \qquad \qquad (f\circ\tau_{\la_j})_{total}&=
v^{n_jd_{j+1}\th_j(\a_{j,0,k})^j_{k=1}}(f\circ\tau_{\la_j})_{proper},\\
(f\circ\tau_{\la_j})_{proper}&
=(u+a+\ve)^{d_{j+1}}+\sum^{d_{j+1}-2}_{i=0}
W_{j+1,i}(u+a+\ve)^i \quad \text{with}  \qquad \qquad \\
 W_{j+1,i}&=W_{j+1,i}(u,v)=b_{j+1,i}v^{M_{j+1,i}}\quad \text{and}\\
 M_{j+1,i}&=\th_{j+1}(\beta_{j+1,i,k})^{j+1}_{k=1}-n_j(d_{j+1}-i)
\th_j(\a_{j,0,k})^j_{k=1}>0,
\endalign$$
where each $b_{j+1,i}$ is a unit in $\C\{u+a,v\}$ for $0\le i\le
d_{j+1}-2$, if exists, noting that
$(S_{j+1,i}\circ\tau_{\la_j})_{total}=v^{n_j(d_{j+1}-i)
\th_j(\a_{j,0,k})^j_{k=1}}W_{j+1,i}$. $\square$

There is nothing to prove Sublemma 16.6.1 by Theorem 14.0 and
Proposition 14.1.

In preparation for the induction proof of the theorem on the integer
$(j+1)$, we must find the method for the construction of
$(f_{j+1},f)$ of (16.6.15). For brevity of the proof, we may assume
without loss of generality that
$\gcd(d_{j+1},\th_{j+1}(\beta_{j+1,0,k})^{j+1}_{k=1})>1$ and $j\le
\ell-2$ from (16.6.14), which will be proved by three cases in the
conclusion of the following sublemma. \ms

\noindent $\underline{\text{\bf Sublemma 16.6.2}}$ \quad
$\underline{\text{\bf Assumptions}}$ \quad Suppose that the same
assumption and notations of Sublemma 16.6.1 hold.

$\underline{\text{\bf Conclusions}}$ \quad For the proof of Theorem
16.6, it suffices to consider three cases, respectively. Moreover,
for the induction proof of the theorem on the integer $j+1$, it
remains to prove this theorem in Case(3) except for Case(1) and
Case(2):

$\underline{\text{\rm Case(1)}}$  If
$\gcd(d_{j+1},\th_{j+1}(\beta_{j+1,0,k})^{j+1}_{k=1})=1$, we may
assume that $f(y,z,f_1,\dots,f_j)=f_{j+1}(y,z,f_1,\dots,f_j)$ up to
the notations, and so there is nothing to prove for the theorem
because {\rm An algorithm for finding an irreducibility criterion
for $f(y,z)\in$the type$[j+1]$} with $j+1=\ell$ which satisfies two
properties (1) and (2), is true, as we have seen in the conclusion
of this theorem. \ms

$\underline{\text{\rm Case(2)}}$ If $j=\ell-1$, we may assume that
$f(y,z,f_1,\dots,f_j)=f_{j+1}(y,z,f_1,\dots,f_j)$ up to the
notations, and so there is nothing to prove for the theorem because
{\rm An algorithm for finding an irreducibility criterion for
$f(y,z)\in$the type$[j+1]$} with $j+1=\ell$ which satisfies two
properties (1) and (2), is true, as in the conclusion of this
theorem. \ms

$\underline{\text{\rm Case(3)}}$ After the proofs of Case(1) and
Case(2) are done, in preparation for the induction proof of the
theorem on the integer $j+1$, we may assume by Case(1) and Case(2)
that
$$
\gcd(d_{j+1},\th_{j+1}(\beta_{j+1,0,k})^{j+1}_{k=1})=\gcd(d_{j+1},M_{j+1,0})>1
\quad \text{and} \quad j\le \ell-2, \tag 16.6.18
$$
where
$M_{j+1,i}=\th_{j+1}(\beta_{j+1,i,k})^{j+1}_{k=1}-n_j(d_{j+1}-i)
\th_j(\a_{j,0,k})^j_{k=1}>0$.  $\square$  \ms

{\bf Proof of Sublemma 16.6.2.} \quad For the proof of the sublemma,
it suffices to prove three cases, respectively.

$\underline{\text{\rm The proof of Case(1).}}$ Let
$\gcd(d_{j+1},\th_{j+1}(\beta_{j+1,0,k})^{j+1}_{k=1})=1$. Then,
$\gcd(d_{j+1},M_{j+1,0})=1$ where
$M_{j+1,0}=\th_{j+1}(\beta_{j+1,0,k})^{j+1}_{k=1}-n_j{d_{j+1}}
\th_j(\a_{j,0,k})^j_{k=1}>0$ by (16.6.17).

So, $(f\circ\tau_{\la_j})_{proper}$ in (16.6.17) is irreducible in
$\C\{u+a,v\}$ if and only if
$$\align
\text{$\f{M_{j+1,i}}{d_{j+1}-i}\ge \f{M_{j+1,0}}{d_{j+1}}$ \quad for
$0\le i\le d_{j+1}-2$,} \tag 16.6.19
\endalign$$
by Corollary $3.3$. Also, an equality in (16.6.19) is equivalently
rewritten as follows:
$$
\f{\th_{j+1}(\beta_{j+1,i,k})^{j+1}_{k=1}}{d_{j+1}-i}\ge
\f{\th_j(\beta_{j+1,0,k})^{j+1}_{k=1}}{d_{j+1}}. \tag 16.6.20
$$

First, note that $\text{$(f\circ\tau_{\la_j})_{proper} \buildrel
\text{{\rm multiseq}} \over \sim
(u+a+\ve)^{d_{j+1}}+v^{M_{j+1,0}}$}$ where
$W_{j+1,0}=b_{j+1,0}v^{M_{j+1,0}}$ by (16.6.17), because
$\gcd(d_{j+1},M_{j+1,0})=1$ and $b_{j+1,0}$ is a unit in
$\C\{u+a,v\}$.

Observe that if $M_{j+1,0}=1$, that is,
$(f\circ\tau_{\la_j})_{proper}$ has no singular point at
$(u+a,v)=(0,0)$, then in order to have a resolution for the curve
$f=0$ at $(0,0)$ in the sense of Corollary 2.3, we need exactly one
more exceptional curve which has three distinct intersection points
with three components among additional exceptional curves including
the exceptional curve $v=0$ and the new proper transform, because
the curve $v=0$ and the last proper transform
$(f\circ\tau_{\la_j})_{proper}$ in (16.6.17) meet tangentially at
$(u+a,v)=(0,0)$. Also, if $M_{j+1,0}>1$, it is clear by Corollary
2.3 that $(f\circ\tau_{\la_j})_{proper}$ satisfies the same kind of
properties as we have done in case that $M_{j+1,0}=1$, as far as a
resolution is concerned in the sense of Corollary 2.3. Now, replace
$d_{j+1}$, $B_{j+1,i},\beta_{j+1,i,k}$, $S_{j+1,i}$ by $n_{j+1}$,
$A_{j+1,i}$, $\a_{j+1,i,k}$, $R_{j+1,i}$, respectively. Then, it is
trivial that $f=f_{j+1}$ satisfies the desired properties up to the
notations. Thus, the proof of Case(1) is done.

$\underline{\text{\rm The proof of Case(2).}}$ If $j=\ell-1$, then
we will prove that $\gcd(d_{j+1},M_{j+1,0})=1$. Assume the contrary,
i.e., $\gcd(d_{j+1},M_{j+1,0})>1$. Then, we would have three
subcases (i), (ii) and (iii).

(i) $d_{j+1}<M_{j+1,0}$ : Note that
$1<\gcd(d_{j+1},M_{j+1,0})<d_{j+1}$ and then
$(f\circ\tau_{\la_j})_{proper}\in $ the type $[k]$ for some positive
integer $k\ge 2$ by Hensel's lemma or Lemma $3.1$ because
$S_{j+1,i}=0$ for $i=d_{j+1}-1$, and so $f\in $ the type $[\ell']$
with $\ell'>\ell$. It would be impossible.

(ii) $d_{j+1}=M_{j+1,0}$ : Then, $(f\circ\tau_{\la_j})_{proper}$
would not be irreducible in $\C\{u+a,v\}$ by Lemma $3.1$, because
$S_{j+1,i}=0$ for $i=d_{j+1}-1$.

(iii) $d_{j+1}>M_{j+1,0}$ : It is enough to consider the following
two cases.

(iiia) If $\gcd(d_{j+1},M_{j+1,0})<M_{j+1,0}$, then
$(f\circ\tau_{\la_j})_{proper}\in $ the type $[k]$ for some positive
integer $k\ge 2$ because $\gcd(d_{j+1},M_{j+1,0})>1$, and so $f\in $
the type $[\ell']$ with $\ell'>\ell$. It would be impossible.

(iiib) If $d_{j+1}$ is a positive multiple of $M_{j+1,0}$, then note
that the last exceptional curve $v=0$ and the proper transform
$(f\circ\tau_{\la_j})_{proper}=0$ meet tangentially and that
$(f\circ\tau_{\la_j})_{proper}\in $ the type $[k]$ for $k\ge 1$,
since $M_{j+1,0}>1$. But, to resolve the total transform defined by
$v(f\circ\tau_{\la_j})_{proper}=0$ in the sense of Corollary 2.3, we
need at least two more exceptional curves, each of which has three
distinct intersection points with three components among additional
exceptional curves including the exceptional curve $v=0$ and the new
proper transform because the curve $v=0$ and the curve
$(f\circ\tau_{\la_j})_{proper}=0$ meet tangentially at
$(u+a,v)=(0,0)$. Therefore, $f\in $ the type $[\ell']$ with
$\ell'>\ell$, which would imply a contradiction.

Thus, the proof of Case(2) is done, and so by the same kind of
replacement as we have used in the proof of Case(1), $f$ can be
defined by $f_{j+1}$ with the desired properties. Therefore, the
proof of the sublemma is done.

$\underline{\text{\rm The proof of Case(3).}}$ Using Case(1) and
Case(2), there is nothing to prove. $\square$ \ms

\noindent{\bf Remark 16.6.2.1 for Sublemma 16.6.2.} Hereafter, using
Case(1), Case(2) and Case(3), in preparation for the induction proof
of the theorem on the integer $j+1$, we may assume without loss of
generality that
$$
\gcd(d_{j+1},\th_{j+1}(\beta_{j+1,0,k})^{j+1}_{k=1})=\gcd(d_{j+1},M_{j+1,0})>1
\quad \text{and} \quad j\le \ell-2, \tag 16.6.21.1
$$
where
$M_{j+1,i}=\th_{j+1}(\beta_{j+1,i,k})^{j+1}_{k=1}-n_j(d_{j+1}-i)
\th_j(\a_{j,0,k})^j_{k=1}>0.$ \ms

{\rm(iii)} Since $f$ is irreducible in $\C\{y,z\}$, then by Hensel's
lemma or Lemma 16.0, the equation $(f\circ\tau_{\la_j})_{proper}$ in
(16.6.17) can be rewritten in the form
$$
(f\circ\tau_{\la_j})_{proper}=[(u+a+\ve)^{n_{j+1}}+\zeta
v^{M'}]^{d_{j+2}}+\sum_{s,t\ge 0}a_{s,t}(u+a+\ve)^sv^t, \tag
16.6.21.2
$$
where $d_{j+2}=\gcd(d_{j+1}, M_{j+1,0})=\gcd(d_{j+1},
\th_{j+1}(\beta_{j+1,0,k})^{j+1}_{k=1})>1$ with
$d_{j+1}=n_{j+1}d_{j+2}$ and $M_{j+1,0}=M'd_{j+2}$ for some integers
$n_{j+1}\ge 2$ and $M'\ge 1$, and $\zeta$ is a nonzero constant, and
the $a_{s,t}$ are some nonzero constant such that
$sM'+tn_{j+1}>n_{j+1}d_{j+2}M'=d_{j+1}M'=n_{j+1}M_{j+1,0}$ for $s\ge
0$ and $t\ge 0$.

Observe that $W_{j+1,d_{j}-n_{j+1}}$ is not zero, applying Hensel's
lemma or Lemma 3.1 to the equation in (16.6.17) because
$n_{j+1}(d_{j+2}-1)=d_{j+1}-n_{j+1}$ in (16.6.21.2).  $\square$ \ms

For the completeness of the proofs of this theorem, by the induction
method on the integer $j$, $0\le j\le {\ell}-1$, the next aim is to
construct three steps with proofs, denoted by Proposition 16.7(Step
I), Proposition 16.8(Step II), and Proposition 16.9(Step III) called
Theorem 16.6, by using the same kind of methods and notations as we
have seen in Proposition 16.2(Step I), Proposition 16.3(Step II),
and Theorem 16.4 called Proposition 16.4(Step III).

\proclaim{Proposition 16.7(Step I)}

$\underline{\text{\bf Assumptions}}$ By induction assumption on the
integer $j$, suppose that $(f_j,f)$ of (16.6.1) can be uniquely
written in the form
$$
\cases f_j &=f^{n_j}_{j-1}+\sum^{n_j-2}_{i=0} R_{j,i}f^i_{j-1} \in
\BC\{y,z,f_1,\dots,f_{j-2}\}[f_{j-1}] \\
f &=f^{d_{j+1}}_j +\sum^{d_{j+1}-2}_{i=0} S_{j+1,i}f^i_j  \in
\BC\{y,z,f_1,\dots,f_{j-1}\}[f_{j}],
\endcases \tag 16.7.1
$$
satisfying the same properties and notations as in the conclusions
of Theorem 16.6.

By induction assumption on the integer $j$, recall by (16.6.11-1) in
the conclusion of this theorem that there exists a unique nonzero
monomial
$$\align
B_{j+1,d_{j+1}-n_{j+1}}\Pi^{j+1}_{k=1}f^{\beta_{j+1,d_{j+1}-n_{j+1},k}}_{k-2}
\tag 16.7.2
\endalign$$
in $S_{j+1,d_{j+1}-n_{j+1}}$ with a constant
$B_{j+1,d_{j+1}-n_{j+1}}$ such that
$\th_{j+1}(\beta_{j+1,d_{j+1}-n_{j+1},k})^{j+1}_{k=1}
=\text{$\min$}\{\th_{j+1}(\de_k)^{j+1}_{k=1}\}$ for any nonzero
monomial $\Pi^{j+1}_{k=1}f^{\de_k}_{k-2}\in
S_{j+1,d_{j+1}-n_{j+1}}$, if exists. \ms

For the induction proof on the integer $j+1$, we may assume without
proof that

$d_{j+2}=\gcd(d_{j+1},\th_{j+1}(\beta_{j+1,0,k})^{j+1}_{k=1})>1$
with $d_{j+1}=n_{j+1}d_{j+2}$ and $j\le \ell-2$ as we have seen in
Remark 16.6.2.1 for Sublemma 16.6.2. \ms

$\underline{\text{\bf Conclusions}}$ Then, $(g_{j+1},f)$ can be
uniquely written in the following form:
$$
\cases g_{j+1} &=f^{n_{j+1}}_j
+\xi_{j+1}\Pi^{j+1}_{k=1}f^{\si_k}_{k-2}\quad
\text{with} \ f_{-1}=y\ \text{and}\ f_0=z, \\
f &=g^{d_{j+2}}_{j+1}+\sum^{d_{j+2}-1}_{i=0} T_{j+2,i}g^i_{j+1},
\endcases \tag 16.7.3
$$
where, considering $y,z,f_1,\dots, f_j$ as independent complex
$(j+2)$-variables at $0\in\C^{j+2}$,

{\rm(i)} $n=d_{j+2}\Pi^{j+1}_{k=1}n_k$ with $d_{j+2}\ge 2$ and
$n_{k}\ge 2$ for $1\le k\le j+1$, and $n=d_1$ if $j=0$;

{\rm(ii)} $\si_k=\beta_{j+1,d_{j+1}-n_{j+1},k}$ for $1\le k\le j+1$
and $\xi_{j+1}=\f 1{d_{j+2}}B_{j+1,d_{j+1}-n_{j+1}}$ by (16.6.11-1);

{\rm(iii)} for each $i$ with $0\le i\le d_{j+2}-1$, $T_{j+2,i} \in
\C\{y,z,f_1,\dots, f_{j}\}$;

{\rm(iv)} $g_{j+1}=g_{j+1}(y,z,f_1,\dots, f_{j})\in
\C\{y,z,f_1,\dots, f_{j-1}\}[f_j]$;

{\rm(v)} $f=f(y,z,f_1,\dots, f_j,g_{j+1})\in \C\{y,z,f_1,\dots,
f_j\}[g_{j+1}]$,

satisfying two conditions, denoted by The Necessary and Sufficient
Condition$[A]$ for $g_{j+1}(y,z)\in$the type$[j+1]$ and The
Necessary Condition$[B]$ for $f(y,z)\in$the type$[\ell]$ with
$\ell\ge j+2$, each of which is represented respectively, as
follows: \ms

$\underline{ \text{\bf [1]  The Necessary and Sufficient
Condition[A] for $g_{j+1}(y,z)\in$the type$[j+1]$:}}$

\noindent$\underline{ \text{\bf{\bf $g_{j+1}\in
\BC\{f_{-1},f_0,\dots,f_{j-1}\}[f_{j}]$} is an irreducible W-poly of
degree $n_{j+1}$ in $f_j$  }}$

\noindent$\underline{ \text{\bf with a coefficient of
${f_{j}}^{n_{j+1}-1}$ zero in $\BC\{f_{-1},f_0,\dots,f_{j-1}\}$, and
$g_{j+1}(y,z)\in$ the type[j+1] }}$

\noindent$\underline{ \text{\bf in the sense of Definition 2.5}}$
\ms

Let $j$ be arbitrary with $3\le {j}$, noting that
$n=\Pi^{r}_{i=1}q_i$.

To find \text{\rm The Necessary and Sufficient Condition[A] for
$g_{j+1}(y,z)\in$the type$[j+1]$}, as an element in $\BC\{y,z\}$ if
necessary, it suffices to show that $g_{j+1} =f^{n_{j+1}}_j
+\xi_{j+1}\Pi^{j+1}_{k=1}f^{\si_k}_{k-2}$ of {\rm(16.7.3)} satisfies
two properties {\rm(1)} and {\rm(2)}: \ms

{\rm(1)} Each $R_{j+1,0}\in \C \{y\}[z,f_1,\dots,f_{j-1}]$ of
$g_{j+1}$ in {\rm(16.7.3)} satisfies {\rm(1a)}, {\rm(1b)}, {\rm(1c)}
and {\rm(1d)} for $i=0,1,\dots,n_{j+1}-2$.

{\rm(1a)} For a nonzero monomial $\Pi^{j+1}_{k=1}f^{\sigma_k}_{k-2}$
in $g_{j+1}$ with $f_{-1}=y$ and $f_0=z$,
$$
\sigma_1>0 \quad \text{and} \quad \sigma_k<n_{k-1}\quad \text{for}
\quad k=2,3,\dots, j+1. \tag 16.7.4
$$

{\rm(1b)} Define a function $\th_{j+1}:\N_0^{j+1}\to \N_0$ where
$\N_0^{j+1}$ is its $(j+1)$-dimensional cartesian product by
$$
(16.7.5) \qquad \qquad
\th_{j+1}(t_k)^{j+1}_{k=1}=t_{j+1}\th_j(\a_{j,0,k})^j_{k=1}+n_j\th(t_k)^j_{k=1}
\quad \text{for each $(t_k)^{j+1}_{k=1}\in \N^{j+1}_0$}. \qquad
\qquad
$$

For any two nonzero monomials $\Pi^{j+1}_{k=1}f^{\beta_k}_{k-2}$ and
$\Pi^{j+1}_{k=1}f^{\de_k}_{k-2}$ where $\beta_1>0$,
$\beta_k<n_{k-1}$ for $2\le k\le j+1$, $\delta_1>0$ and
$\delta_k<n_{k-1}$ for $2\le k\le j+1$, we have
$$
\th_{j+1}(\beta_k)^{j+1}_{k=1}=\th_{j+1}(\de_k)^{j+1}_{k=1}\
\text{if and only if} \ \beta_k=\de_k\ \text{for}\ k=1,2,\dots, j+1.
\tag 16.7.6
$$

{\rm(1c)} Then, $g_{j+1}$ of {\rm(16.7.3)} satisfies the following
properties:

Let $\tau_{\la_{j}}$ be the composition of a finite number $\la_{j}$
of successive blow-ups which is needed only to get the standard
resolution of the singular point of $V(f_j)$, as in Sublemma
$16.6.1$. After repeating the same number $\la_j$ of blow-ups with
the same local coordinates as we have used in the standard
resolution of the singular point of $V(f_j)$, the local defining
equation for the $\la_j$-th proper transform of the curve defined by
$g_{j+1}=0$ can be written in the form
$$\align
(16.7.7) \qquad &\text{$(u+a+\ve)^{n_{j+1}}+c_{j+1}v^{M'}$
\quad{with} $M'=\th_{j+1}(\sigma_k)^{j+1}_{k=1}
-n_jn_{j+1}\th_j(\a_{j,0,k})^j_{k=1}>0$} \qquad \qquad \\
& \quad\text{{where} \quad $\sigma_k=\beta_{j+1,d_{j+1}-n_{j+1},k}$ \quad {for} $1\le k\le j+1$.} \qquad \qquad \\
\endalign$$
from a uniquely defined nonzero monomial
$B_{j+1,d_{j+1}-n_{j+1}}\Pi^{j+1}_{k=1}f^{\beta_{j+1,d_{j+1}-n_{j+1},k}}_{k-2}$
of $S_{j+1,d_{j+1}-n_{j+1}}$ as in
$f(y,z,\dots,f_j)=f^{d_{j+1}}_j+\sum^{d_{j+1}-2}_{i=0}S_{j+1,i}f^i_j$
by induction assumption on the integer $j$ and
$c_{j+1}=c_{j+1}(u,v)$ is a unit in $\C\{u+a,v\}$ with
$c_{j+1}(-a,0)=\zeta$. \ms

{\rm(1d)} Also, an equation of $(16.7.7)$ satisfies the following:
$$\align
(16.7.8)\quad \quad \gcd(n_{j+1},\th_{j+1}(\si_k)^{j+1}_{k=1})=1
\quad \text{{with} \quad $\sigma_k=\beta_{j+1,d_{j+1}-n_{j+1},k}$ \quad {for} $1\le k\le j+1$.} \qquad \qquad\\
\endalign$$

{\rm(2)}{\rm(2a)} $g_{j+1}=g_{j+1}(y,z)\in \C\{y\}[z]$ is an
irreducible $W$-poly in $z$ of multiplicity $\Pi^{j+1}_{k=1}n_k$ at
$0\in\C^{2}$ with coefficients in $\BC\{y\}$, and $g_{j+1} \in$ the
type $[j+1]$. \ms

{\rm(2b)} $g_{j+1}=g_{j+1}(y,z,f_1,\dots,f_j)\in
\C\{y,z,f_1,\dots,f_{j-1}\}[f_{j}]$ of {\rm(16.7.3)} is an
irreducible $W$-poly in $f_{j}$ with coefficients in
$\C\{y,z,f_1,\dots,f_{j-1}\}$ and with multiplicity $n_{j+1}$ at
$0\in\C^{j+2}$. \ms

$\underline{ \text{\bf [2] The Necessary Condition[B] for
$f(y,z)\in$the type$[\ell]$ with $\ell\ge j+2$:}}$

\noindent$\underline{ \text{\bf $f\in
\C\{f_{-1},f_0,\dots,f_{j}\}[g_{j+1}]$ is an irreducible W-poly of
degree $d_{j+2}$ in $g_{j+1}$ with a coefficient of}}$

\noindent$\underline{ \text{\bf ${g_{j+1}}^{d_{j+2}-1}$ either zero
or nonzero in $\C\{f_{-1},f_0,\dots,f_{j}\}$, and $f(y,z)\in$the
type[${\ell}$] with $\ell\ge j+2$}}$ \ms

Note that ${j+1}\le{\ell}$ where $j$ was already given by
{\rm(16.6.1)} and that $n=\Pi^{r}_{i=1}q_i$.

To find \text{\rm the Necessary Condition[B] for $f(y,z)\in$the
type$[\ell]$ with $\ell\ge j+2$}, as an element in $\BC\{y,z\}$ if
necessary, it is enough to show that
$f=g^{d_{j+2}}_{j+1}+\sum^{d_{j+2}-1}_{i=0} T_{j+2,i}g^i_{j+1}$ of
{\rm(16.7.3)} satisfies two properties {\rm(1)} and {\rm(2)}: Note
that either $\ell=j+2$ or $\ell> j+2$ and that $T_{j+2,d_{j+2}-1}$
may be nonzero.\ms

{\rm(1)} Each $T_{j+2,i}\in \C \{y\}[z,f_1,\dots,f_j]$ of $f$ in
{\rm(16.7.8)} satisfies {\rm(1a)}, {\rm(1b)}, {\rm(1c)} and
{\rm(1d)} for $i=0,1,\dots,d_{j+2}-1$.

{\rm(1a)} For any nonzero monomial $\Pi^{j+2}_{k=1}f^{\de_k}_{k-2}$
in $T_{j+2,i}$,
$$
\de_1>0\quad \text{and}\quad \de_k<n_{k-1}\quad \text{for}\
k=2,3,\dots,j+2. \tag 16.7.9
$$

In particular, if $i=d_{j+2}-1$ for $T_{j+2,i}$, then $\de_{j+2}\le
n_{j+1}-2$. \ms

{\rm(1b)} Define a function $\ol \th_{j+2}:\N^{j+2}_0\to \N_0$ by
$$
\ol\th_{j+2}(t_k)^{j+2}_{k=1}=t_{j+2}\th_{j+1}(\si_k)^{j+1}_{k=1}
+n_{j+1}\th_{j+1}(t_k)^{j+1}_{k=1} \quad \text{for each
$(t_k)^{j+2}_{k=1}\in \N^2_0$}. \tag 16.7.10
$$

For any two nonzero monomials $\Pi^{j+2}_{k=1}f^{\g_k}_{k-2}$ and
$\Pi^{j+2}_{k=1}f^{\de_k}_{k-2}$ in $T_{j+2,i}$,
$$\align
(16.7.10^*) \qquad \qquad
&\ol\th_{j+2}(\g_k)^{j+2}_{k=1}=\ol\th_{j+2}(\de_k)^{j+2}_{k=1} \ \
\text{if and only if} \ \ \g_k=\de_k ~\text{for $k=1,2,\dots,j+2$}.
\qquad \qquad\qquad \qquad \\
&\text{So, there exists a unique nonzero monomial
$C_{j+2,i}\Pi^{j+2}_{k=1}f^{\beta_{j+2,i,k}}_{k-2}$ in $T_{j+2,i}$} \\
&\text{with a nonzero constant $C_{j+2,i}$ such that
$\ol\th_{j+2}(\beta_{j+2,i,k})^{j+2}_{k=1}
=\text{Min}\{\ol\th_{j+2}(\de_k)^{j+2}_{k=1}\}$}\\
&\text{for any nonzero monomial $\Pi^{j+2}_{k=1}f^{\de_k}_{k-2}$ in
$T_{j+2,i}$ with $i$ fixed.}\\
\endalign$$

{\rm(1c)} For all $i=0,1,\dots, d_{j+2}-1$,
$$\align
\ol\th_{j+2}(\beta_{j+2,i,k})^{j+2}_{k=1}>(d_{j+2}-i)
n_{j+1}\th_{j+1}(\sigma_k)^{j+1}_{k=1}. \tag 16.7.11
\endalign$$

{\rm(1d)} For all $ i=0,1,\dots,d_{j+2}-1$,
$$\align
 &\gcd(d_{j+2},\ol\th_{j+2}(\beta_{j+2,0,k})^{j+2}_{k=1})\le
d_{j+2} \quad \text{and}
\tag 16.7.12 \\
&\f{\ol\th_{j+2}(\beta_{j+2,i,k})^{j+2}_{k=1}}{d_{j+2}-i}\ge
\f{\ol\th_{j+2}(\beta_{j+2,0,k})^{j+2}_{k=1}}{d_{j+2}}.
\endalign$$

Then, either
$\gcd(d_{j+2},\ol\th_{j+2}(\beta_{j+2,0,k})^{j+2}_{k=1})=1$ or
$1<\gcd(d_{j+2},\ol\th_{j+2}(\beta_{j+2,0,k})^{j+2}_{k=1})$. \ms

{\rm(1d-1)} Let
$\gcd(d_{j+2},\ol\th_{j+2}(\beta_{j+2,0,k})^{j+2}_{k=1})=1$. Then
$f$ is irreducible in ${}_2\CO_0$ with $f \in$ the type $[j+2]$ if
and only if the inequality in (16.7.12) holds.  and
$f_1,f_2,\dots,f_j, g_{j+1}$ are irreducible in $\C\{y,z\}$ as
above.

{\rm(1d-2)} Let
$1<\gcd(d_{j+2},\ol\th_{j+2}(\beta_{j+2,0,k})^{j+2}_{k=1})\le
d_{j+2}$. If $f$ is irreducible in ${}_2\CO_0$ and
$T_{j+2,d_{j+2}-1}=0$, then
$1<\gcd(d_{j+2},\ol\th_{j+2}(\beta_{j+2,0,k})^{j+2}_{k=1})<d_{j+2}$,
and so $f \in$ the type $[\ell]$ with $\ell\ge j+3$. But, if $f$ is
irreducible in ${}_2\CO_0$ and $T_{j+2,d_{j+2}-1}\not=0$, then
$\gcd(d_{j+2},\ol\th_{j+2}(\beta_{j+2,0,k})^{j+2}_{k=1})$ may be
equal to $d_{j+2}$. \ms

{\rm(2)}{\rm(2a)} $f=f(y,z)\in \C\{y\}[z]$ is an irreducible
$W$-poly in $z$ with coefficients in $\BC\{y\}$ and with
multiplicity $n=d_{j+2}\Pi^{j+1}_{k=1}n_k$ at $0\in\C^2$. Also,
either $f \in$ the type $[j+2]$ or $f \in$ the type $[\ell]$ with
$\ell\ge j+3$.

{\rm(2b)} $f=f(y,z,f_1,\dots,f_j,g_{j+1})\in
\C\{y,z,f_1,\dots,f_{j}\}[g_{j+1}]$ of {\rm(16.7.8)} is an
irreducible $W$-poly of degree $d_{j+2}$ in $g_{j+1}$ with
coefficients in $\C\{y,z,f_1,\dots,f_{j}\}$ and  with multiplicity
$d_{j+2}$ at $0\in\C^{j+3}$. $\square$ \endproclaim \ms

\noindent$\underline{\text{\bf Remark 16.7.1.}}$ {\rm (1)} Whenever
$g_{j+1}$ of (16.7.1) satisfies The Necessary and Sufficient
Condition[A] for $g_{j+1}(y,z)\in$the type$[j+1]$, then it is said
that $g_{j+1}(y,z)\in$the type$[j+1]$ in the sense of Definition
$2.5$.

{\rm(2)} Assume that $f$ of (16.7.1) satisfies The Necessary
Condition[B] for $f(y,z)\in$the type$[\ell]$ with $\ell\ge j+2$.

{\rm(i)} If $f$ satisfies $T_{j+2,d_{j+2}-1}=0$ in addition, then it
is said that $f$ satisfies The Necessary Condition[A] for
$f(y,z)\in$the type$[\ell]$ with $\ell\ge j+2$.

{\rm(ii)} If $T_{j+2,d_{j+2}-1}=0$, and also
$\gcd(d_{j+2},\ol\th_{j+2}(\beta_{j+2,0,k})^{j+2}_{k=1})=1$, then it
is said that $f(y,z)\in$the type$[j+2]$ in the sense of Definition
$2.5$, noting that $f$ of (16.7.1) satisfies The Necessary and
Sufficient Condition[A] for $f(y,z)\in$the type$[j+2]$ in the sense
of Definition $2.5$. \ms

\proclaim{Proposition 16.8(Step II)} $\underline{\text{\bf
Assumptions}}$ For the induction proof on the integer $j+1$, we may
assume without loss of generality that
$\gcd(d_{j+1},\th_{j+1}(\beta_{j+1,0,k})^{j+1}_{k=1})>1$ and $j\le
\ell-2$ as we have seen in either Sublemma $16.6.2$ or Sublemma
$16.7$.

Then, we may assume by Proposition $16.7$ that $(g_{j+1},f)$ can be
written as follows:
$$
\cases g_{j+1} &=f^{n_{j+1}}_j
+\xi_{j+1}\Pi^{j+1}_{k=1}f^{\si_k}_{k-2} \quad
\text{with} \ f_{-1}=y\ \text{and}\ f_0=z, \\
f &=g^{d_{j+2}}_{j+1}+\sum^{d_{j+2}-1}_{i=0} T_{j+2,i}g^i_{j+1},
\endcases \tag 16.8.1
$$
$\si_k=\beta_{j+1,d_{j+1}-n_{j+1},k}$ for $1\le k\le j+1$,
satisfying two conditions as we have seen in Proposition $16.7$,
denoted by The Necessary and Sufficient Condition$[A]$ for
$g_{j+1}(y,z)\in$the type$[j+1]$ and The Necessary Condition$[B]$
for $f(y,z)\in$the type$[\ell]$ with $\ell\ge j+2$. \ms

$\underline{\text{\bf Conclusions}}$ The main aim is to construct a
unique pair $(f_{j+1},f)$ such that $(f_{j+1},f)$ can be written in
the form
$$\cases
f_{j+1} &=f^{n_{j+1}}_{j}+\sum^{n_{j+1}-2}_{i=0} R_{j+1,i}f^i_{j}
\quad
\text{with} \ f_{-1}=y\ \text{and}\ f_0=z, \\
f &=f^{d_{j+2}}_{j+1} +\sum^{d_{j+2}-2}_{i=0} S_{j+2,i}f^i_{j+1},
\endcases \tag 16.8.2
$$
where $y,z,f_1,\dots, f_{j+1}$ are considered as independent complex
$(j+2)$-variables at the origin in $\C^{j+3}$ if necessary,
satisfying the following properties:

{\rm(i)} The first problem is how to construct
$f_{j+1}=f_{j+1}(y,z)$ satisfying the condition in \text{\rm
$\widehat{\widehat{\text{\rm[1]}}}$} such that
\text{$f_{j+1}(y,z)\buildrel \text{{\rm multiseq}} \over \sim
g_{j+1}(y,z)$ \quad under the standard resolutions}.

$\underline{ \text{\bf $\widehat{\widehat{\text{\bf[1]}}}$  The
Necessary and Sufficient Condition[A] for $f_{j+1}(y,z)\in$the
type$[j+1]$:}}$

\noindent$\underline{ \text{\bf{\bf $f_{j+1}\in
\BC\{f_{-1},\dots,f_{j-1}\}[f_{j}]$} is an irreducible W-poly  of
degree $n_{j+1}$ in $f_j$ with a coefficient of }}$

\noindent$\underline{ \text{\bf ${f_{j}}^{n_{j+1}-1}$ zero in
$\C\{f_{-1},\dots,f_{j-1}\}$, and $f_{j+1}(y,z)\in$ the type[j+1] in
the sense of Definition 2.5}}$ \ms

{\rm(ii)} The second problem is to prove that $f=
f(y,z,f_1,\dots,f_{j},f_{j+1})$ satisfies the condition in \text{\rm
$\widehat{\widehat{\text{\rm[2]}}}$} which is defined by the same
kind of property as $f(y,z,f_1,\dots,f_{j},g_{j+1})$ have done in
The Necessary Condition$[B]$ for $f(y,z)\in$the type$[\ell]$ with
$\ell\ge j+2$.

$\underline{ \text{\bf $\widehat{\widehat{\text{\bf[2]}}}$ The
Necessary Condition[A] for $f(y,z)\in$the type$[\ell]$ with $\ell\ge
j+2$:}}$

\noindent$\underline{ \text{\bf {\bf $f\in
\C\{f_{-1},f_0,\dots,f_{j}\}[f_{j+1}]$} is an irreducible W-poly of
degree $d_{j+2}$ in $f_{j+1}$ with a coefficient }}$

\noindent$\underline{ \text{\bf of ${f_{j+1}}^{d_{j+2}-1}$ zero in
$\C\{f_{-1},f_0,\dots,f_{j}\}$, and $f(y,z)\in$the type[${\ell}$]
with $\ell\ge j+2$ }}$ \ms

For the construction of a pair $(f_{j+1},f)$ in $(16.8.2)$,  it
suffices to consider the following two cases, depending on the fact
that $T_{j+2,d_{j+2}-1}$ in a pair $(g_{j+1},f)$ of $(16.8.1)$ is
zero or not. For brevity of notations, let $h_1=g_{j+1}$ and
$T^{(1)}_{j+2,i}=T_{j+2,i}$. \ms

$\underline{\text{\bf Case(1):}}$ Let $T^{(1)}_{j+2,d_{j+2}-1}=0$.
By $(16.8.1)$, put $f_{j+1}=g_{j+1}$,
$R_{j+1,0}=\xi_{j+1}\Pi^{j+1}_{k=1}f^{\si_k}_{k-2}$, and
$S_{j+2,i}=T^{(1)}_{j+2,i}$ for $0\le i\le d_{j+2}-2$. Then, it is
clear by Proposition $16.7$ that $(f_{j+1},f)$ of the main aim and
$(g_{j+1},f)$ of $(16.8.1)$ are the same pairs in the sense of
Definition $16.2.2$. \ms

$\underline{\text{\bf Case(2):}}$ Let $T^{(1)}_{j+2,d_{j+2}-1}\neq
0$. Then, there is a sequence of pairs of $W$-polys in $z$,
$\{(h_p,f):p=1,2,\dots\}$, such that $h_1=g_{j+1}$ and
$(f_{j+1},f)=(h_{\nu +1},f)=(h_{\nu +2}, f)=\cdots $ for some
integer $\nu \le\f{n_{j+1}+1}2$, each pair of which can be written
in the form
$$\cases
h_1 &=f^{n_{j+1}}_j+\xi_{j+1}\Pi^{j+1}_{k=1}f^{\si_k}_{k-2}
=f^{n_{j+1}}_j+R^{(1)}_{j+1,0},\\
f &=h^{d_{j+2}}_1+\sum^{d_{j+2}-1}_{i=0} T^{(1)}_{j+2,i}h^i_1,
\endcases \tag 16.8.3
$$
and for $p=2,3,\dots $
$$\cases
h_{p}& =h_{p-1}+\f
1{d_{j+2}}T^{(p-1)}_{j+2,d_{j+2}-1}=f^{n_{j+1}}_j+\sum^{n_{j+1}-2}_{i=0}
R^{(p)}_{j+1,i}f^i_j,\\
f& =h^{d_{j+2}}_{p}+\sum^{d_{j+2}-1}_{i=0} T^{(p)}_{j+2,i}h^i_{p},
\endcases \tag 16.8.4
$$
with $T^{(p)}_{j+2,d_{j+2}-1}\ne 0$ for $1\le p\le \nu$ and $T^{(\nu
+1)}_{j+2,d_{j+2}-1}=T^{(\nu +2)}_{j+2,d_{j+2}-1}=\cdots =0$,

where, considering $f_{-1},f_{0},\dots,f_j,h_{p}$ as independent
complex $(j+3)$-variables at $0\in\C^{j+3}$,

{\rm(i)} $n=d_{j+2}\Pi^{j+1}_{k=1}n_k$ with $d_{j+2}\ge 2$ and
$n_k\ge 2$ for $1\le k\le j+1$, and $n=d_1$ if $j=0$;

{\rm(ii)}
$R^{(p)}_{j+1,i}=R^{(p)}_{j+1,i}(f_{-1},f_0,\dots,f_{j-1})\in\C\{f_{-1},f_0,\dots,f_{j-1}\}$
for $p\ge 1$ and $0\le i\le n_{j+1}-2$;

{\rm(iii)}
$T^{(p)}_{j+2,i}=T^{(p)}_{j+2,i}(f_{-1},f_0,\dots,f_j)\in\C\{f_{-1},f_0,\dots,f_j\}$
for $p\ge 1$ and $0\le i\le d_{j+2}-1$;

{\rm(iv)} $h_p=h_p(f_{-1},f_0,\dots, f_{j})\in \C\{f_{-1},f_0,\dots,
f_{j-1}\}[f_j]$ for $p\ge 1$;

{\rm(v)} $f=f(f_{-1},f_0,\dots, f_j,h_p)\in \C\{f_{-1},f_0,\dots,
f_j\}[h_p]$,

{\noindent}satisfying two conditions, denoted by The Necessary and
Sufficient Condition$[A]$ for $h_{p}(y,z)\in$the type$[j+1]$ and The
Necessary Condition$[B]$ for $f(y,z)\in$the type$[\ell]$ with
$\ell\ge j+2$, each of which is represented respectively, as
follows:

In more detail, for any fixed $(h_{p},f)$ of {\rm(16.8.4)}, $h_{p}$
of $(h_{p},f)$ satisfies The Necessary and Sufficient Condition$[A]$
for $h_{p}(y,z)\in$the type$[j+1]$ and $f$ of $(h_{p},f)$ satisfies
The Necessary Condition$[B]$ for $f(y,z)\in$the type$[\ell]$ with
$\ell\ge j+2$. In particular, if $p_1$ is a given integer with
$p_1\ge \f{n_{j+1}+1}{2}$, define $(f_{j+1},f)=(h_{p_1},f)$ by
Definition $16.2.2$. Then, $f_{j+1}$ of $(f_{j+1},f)$ of $(16.8.2)$
satisfies The Necessary and Sufficient Condition$[A]$ for
$f_{j+1}(y,z)\in$the type$[j+1]$ in \text{\rm
$\widehat{\widehat{\text{\rm[1]}}}$}, and $f$ of $(f_{j+1},f)$ of
$(16.8.2)$ satisfies The Necessary Condition[B] for $f(y,z)\in$the
type$[\ell]$ with $\ell\ge j+2$ in \text{\rm
$\widehat{\widehat{\text{\rm[2]}}}$} by the same way  as $(h_p,f)$
of $(16.8.4)$ does up to the change of notations: Recall that
$f_{-1}=y$ and $f_0=z$. \ms

$\underline{ \text{\bf [1] The Necessary and Sufficient Condition[A]
for $h_{p}(y,z)\in$the type$[j+1]$:}}$

\noindent$\underline{ \text{\bf $h_{p}\in
\BC\{f_{-1},\dots,f_{j-1}\}[f_{j}]$ is an irreducible W-poly of
degree $n_{j+1}$ in $f_j$ with a coefficient }}$

\noindent$\underline{ \text{\bf  of ${f_{j}}^{n_{j+1}-1}$ zero in
$\BC\{f_{-1},\dots,f_{j-1}\}$, and $h_{p}\in$ the type[j+1] in the
sense of Definition 2.5} }$ \ms

Let $j$ be arbitrary with $j\ge 0$, noting that
$n=d_{j+2}\Pi^{j+1}_{k=1}n_k$ with $d_{j+2}\ge 2$ and $n_{k}\ge 2$
for $1\le k\le j+1$, and $n=d_1$ if $j=0$.

To find \text{\rm The Necessary and Sufficient Condition[A] for
$h_{p}(y,z)\in$the type$[j+1]$}, it suffices to show that
$h_{p}=f^{n_{j+1}}_j+\sum^{n_{j+1}-2}_{i=0} R^{(p)}_{j+1,i}f^i_j$ of
{\rm(16.8.4)} satisfies two properties {\rm(1)} and {\rm(2)}: \ms

{\rm(1)} Let $p$ be fixed with $p\ge 1$. Each $R^{(p)}_{j+1,i}\ne 0$
satisfies the properties {\rm(1a)}, {\rm(1b)}, {\rm(1c)} and
{\rm(1d)} for $i=0,1,\dots,n_{j+1}-2$. Also, for each $p\ge 1$,
\text{$h_p \buildrel \text{{\rm multiseq}} \over \sim h_1$} and
$h_p\in \C\{y\}[z]$ is an irreducible $W$-poly in $z$ of
multiplicity $\Pi^{j+1}_{k=1}n_k$.

For notation, it may be said that for each $p\ge 1$, $h_p$ satisfies
The Necessary and Sufficient Condition[A] for $h_{p}(y,z)\in$the
type$[j+1]$.\ms

{\rm(1a)} For any nonzero monomial $\Pi^{j+1}_{k=1}f^{\de_k}_{k-2}$
in $R^{(p)}_{j+1,i}$,
$$
\de_1>0\ \text{and}\ \de_k<n_{k-1}\ \text{for}\ k=2,3,\dots,j+1.
\tag 16.8.5
$$

{\rm(1b)} For any two nonzero monomials
$\Pi^{j+1}_{k=1}f^{\a_k}_{k-2}$ and $\Pi^{j+1}_{k=1}f^{\de_k}_{k-2}$
in $R^{(p)}_{j+1,i}$,
$$\align
(16.8.6) \qquad &
\text{$\th_{j+1}(\a_k)^{j+1}_{k=1}=\th_{j+1}(\de_k)^{j+1}_{k=1}$ if
and only if $\a_k=\de_k$ \quad for $k=1,2,\dots,j+1$.} \qquad \qquad
\\
&\text{So, there exists a unique nonzero monomial
$A^{(p)}_{j+1,i}\Pi^{j+1}_{k=1}f^{\a^{(p)}_{j+1,i,k}}_{k-2}$ in
$R^{(p)}_{j+1,i}$}\\
&\text{with a constant $A^{(p)}_{j+1,i}$ such that
$\th_{j+1}(\a^{(p)}_{j+1,i,k})^{j+1}_{k=1}
=\text{$\min$}\{\th_{j+1}(\de_k)^{j+1}_{k=1}\}$}\\
&\text{for any nonzero monomial $\Pi^{j+1}_{k=1}f^{\de_k}_{k-2}$ in
$R^{(p)}_{j+1,i}$.} \\
\endalign$$

{\rm(1c)} For all $i=0,1,\dots,n_{j+1}-2$,
$$\align
\th_{j+1}(\a^{(p)}_{j+1,i,k})^{j+1}_{k=1}>({n_{j+1}-i})
n_j\th_j(\a_{j,0,k})^{j}_{k=1}. \tag 16.8.7
\endalign$$

{\rm(1d)} For all $i=0,1,\dots,n_{j+1}-2$ and for $k=1,2,\dots,j+1$,
$$\align
 &\gcd(n_{j+1},\th_{j+1}(\a^{(p)}_{j+1,0,k})^{j+1}_{k=1})= 1 \quad
 \text{with} \quad \si_k=\a^{(p)}_{j+1,0,k},
\tag 16.8.8\\
&\f{\th_{j+1}(\a^{(p)}_{j+1,i,k})^{j+1}_{k=1}}{n_{j+1}-i}\ge
\f{\th_{j+1}(\a^{(p)}_{j+1,0,k})^{j+1}_{k=1}}{n_{j+1}}.
\endalign$$

{\rm(2)}{\rm(2a)} For each $p\ge 1$, $h_p=h_{p}(y,z)\in \C\{y\}[z]$
is an irreducible $W$-poly in $z$ with coefficients in $\BC\{y\}$
and with $h_p\buildrel \text{{\rm multiseq}} \over \sim
h_1=g_{j+1}$, and $h_{p} \in$ the type $[j+1]$. \ms

{\rm(2b)} $h_p=h_p(f_{-1},f_0,\dots,f_j)\in
\C\{f_{-1},f_0,\dots,f_{j-1}\}[f_{j}]$ of {\rm(16.8.4)} is an
irreducible $W$-poly in $f_{j}$ with coefficients in
$\C\{f_{-1},f_0,\dots,f_{j-1}\}$ and with multiplicity $n_{j+1}$ at
$0\in\C^{j+2}$. \ms

$\underline{ \text{\bf [2] The Necessary Condition[B] for
$f(y,z)\in$the type$[\ell]$ with $\ell\ge j+2$: }}$

\noindent$\underline{ \text{\bf $f\in
\C\{f_{-1},f_0,\dots,f_{j}\}[h_p]$ is an irreducible W-poly of
degree $d_{j+2}$ in $h_p$ with a coefficient }}$

\noindent$\underline{ \text{\bf of ${h_p}^{d_{j+2}-1}$ either zero
or nonzero in $\C\{f_{-1},f_0,\dots,f_{j}\}$, and $f(y,z)\in$the
type[${\ell}$] with $\ell\ge j+2$ }}$ \ms

Note that ${j+1}\le{\ell}$ where $j$ was already given by
{\rm(16.6.1)} and that $n=\Pi^{r}_{i=1}q_i$.

To find \text{\rm The Necessary Condition[B] for $f(y,z)\in$the
type$[\ell]$ with $\ell\ge j+2$}, it suffices to show that $f
=h^{d_{j+2}}_{p}+\sum^{d_{j+2}-1}_{i=0} T^{(p)}_{j+2,i}h^i_{p},$ of
{\rm(16.8.4)} satisfies two properties {\rm(1)} and {\rm(2)}: Note
that either $\ell=j+2$ or $\ell> j+2$ and that
$T^{(p+1)}_{j+2,d_{j+2}-1}$ may be nonzero.\ms

{\rm(1)} Each $T^{(p)}_{j+2,i}\in \C \{f_{-1},f_0,\dots,f_j\}$ of
$f$ in {\rm(16.8.4)} satisfies {\rm(1a)}, {\rm(1b)}, {\rm(1c)} and
{\rm(1d)} for $i=0,1,\dots,d_{j+2}-1$.

{\rm(1a)} For any nonzero monomial $\Pi^{j+2}_{k=1}f^{\g_k}_{k-2}$
 in $T^{(p)}_{j+2,i}$,
$$\align
\text{$\g_1>0$ \quad and \quad $\g_k<n_{k-1}$ \quad for
$k=2,3,\dots,j+2$.} \tag 16.8.9
\endalign$$

In particular, if $i=d_{j+2}-1$ for $T^{(p)}_{j+2,i}$ then
$\g_{j+2}\le n_{j+1}-2$. \ms

{\rm(1b)} Define $\ol\th_{j+2}(t_k)^{j+2}_{k=1}=
t_{j+2}\th_{j+1}(\si_k)^{j+1}_{k=1}+n_{j+1}\th_{j+1}(t_k)^{j+1}_{k=1}$
for any $(t_k)^{j+2}_{k=1}\in N^{j+2}_0$ by the same way as we have
seen in Proposition $16.7$, {\rm(1b)} of {\rm The Necessary
Condition[B] for $f(y,z)\in$the type$[\ell]$ with $\ell\ge j+2$}.

For any two nonzero monomials $\Pi^{j+2}_{k=1}f^{\beta_k}_{k-2}$ and
$\Pi^{j+2}_{k=1}f^{\g_k}_{k-2}$ in $T^{(p)}_{j+2,i}$,
$$\align
(16.8.10) \qquad
&\text{$\ol\th_{j+2}(\beta_k)^{j+2}_{k=1}=\ol\th_{j+2}(\g_k)^{j+2}_{k=1}$
if and only if $\beta_k=\g_k$ for $k=1,2,\dots,j+2$.} \qquad \qquad
\\
&\text{So, there is a unique nonzero monomial
$C^{(p)}_{j+2,i}\Pi^{j+2}_{k=1}f^{\beta^{(p)}_{j+2,i,k}}_{k-2}$ in
$T^{(p)}_{j+2,i}$} \\
&\text{with a constant $C^{(p)}_{j+2,i}$ such that
$\ol\th_{j+2}(\beta^{(p)}_{j+2,i,k})^{j+2}_{k=1}
=\text{$\min$}\{\ol\th_{j+2}(\g_k)^{j+2}_{k=1}\}$}\\
&\text{for any nonzero monomial $\Pi^{j+2}_{k=1}f^{\g_k}_{k-2}$ in
$T^{(p)}_{j+2,i}$.}
\endalign$$

{\rm(1c)} For all $i=0,1,\dots,d_{j+2}-1$,
$$\align
\ol\th_{j+2}(\beta^{(p)}_{j+2,i,k})^{j+2}_{k=1}>(d_{j+2}-i)
n_{j+1}\th_{j+1}(\si_k)^{j+1}_{k=1}. \tag 16.8.11
\endalign$$

{\rm(1d)} For all $i=0,1,\dots,d_{j+2}-1$,
$$\align
&\gcd(d_{j+2},\ol\th_{j+2}(\beta^{(p)}_{j+2,0,k})^{j+2}_{k=1})\le
d_{j+2}, \tag 16.8.12\\
& \f{\ol\th_{j+2}(\beta^{(p)}_{j+2,i,k})^{j+2}_{k=1}}{d_{j+2}-i}\ge
\f{\ol\th_{j+2}(\beta^{(p)}_{j+2,0,k})^{j+2}_{k=1}}{d_{j+2}}.
\endalign$$

Then, either
$\gcd(d_{j+2},\ol\th_{j+2}(\beta^{(p)}_{j+2,0,k})^{j+2}_{k=1})=1$ or
$1<\gcd(d_{j+2},\ol\th_{j+2}(\beta^{(p)}_{j+2,0,k})^{j+2}_{k=1})$.
\ms

{\rm(1d-1)} Let
$\gcd(d_{j+2},\ol\th_{j+2}(\beta^{(p)}_{j+2,0,k})^{j+2}_{k=1})=1$.
Then, $f$ is irreducible in ${}_2\CO_0$ if and only if the
inequality in $(16.8.12)$ holds. In this case, $f \in$ the type
$[j+2]$, but note that $T^{(p)}_{j+2,d_{j+2}-1}$ may not be zero
where
$$
h_p=h_{p-1}+\f 1{d_{j+2}}T^{(p-1)}_{j+2,d_{j+2}-1} \quad \text{and}
\quad f=h^{d_{j+2}}_p+\sum^{d_{j+2}-1}_{i=0}T^{(p)}_{j+2,i}h^i_p.
\tag 16.8.13
$$

{\rm(1d-2)} Let
$1<\gcd(d_{j+2},\ol\th_{j+2}(\beta^{(p)}_{j+2,0,k})^{j+2}_{k=1})\le
d_{j+2}$. There is a positive integer $\nu$ with $\nu\le
\f{n_{j+1}+1}2$ such that $T^{(\nu +1)}_{j+2,d_{j+2}-1}=0$ and
$T^{(p)}_{j+2,d_{j+2}-1}\ne 0$ for $p=1,2,\dots,\nu$. In this case,
$f \in$ the type $[\ell]$ with $\ell\ge j+3$ and note that
$$
1<\gcd(d_{j+2},\ol\th_{j+2}(\beta^{(\nu+1)}_{j+2,0,k})^{j+2}_{k=1}<d_{j+2}.
\tag 16.8.14
$$

{\rm(2)}{\rm(2a)} $f=f(y,z)\in \C\{y\}[z]$ is an irreducible
$W$-poly in $z$ with coefficients in $\BC\{y\}$ and with
multiplicity $n=d_{j+2}\Pi^{j+1}_{k=1}n_k$ at $0\in\C^2$. Also,
either $f \in$ the type $[j+2]$ or $f \in$ the type $[\ell]$ with
$\ell\ge j+3$.

{\rm(2b)} $f=f(f_{-1},f_0,\dots,f_j,h_{p})\in
\C\{f_{-1},f_0,\dots,f_{j}\}[h_p]$ of {\rm(16.8.4)} is an
irreducible $W$-poly in $h_p$ with coefficients in
$\C\{f_{-1},f_0,\dots,f_{j}\}$ and with multiplicity $d_{j+2}$ at
$0\in\C^{j+3}$. \ms

$\underline{ \text{\bf [3] The Necessary Condition[A] for
$f(y,z)\in$the type$[\ell]$ with $\ell\ge j+2$: }}$

\noindent$\underline{ \text{\bf $f\in
\C\{f_{-1},\dots,f_{j}\}[f_{j+1}]$ is an irreducible W-poly of
degree $d_{j+2}$ in $f_{j+1}$ with a}}$

\noindent$\underline{ \text{\bf coefficient of
${f_{j+1}}^{d_{j+2}-1}$ zero in $\C\{f_{-1},\dots,f_{j}\}$, and
$f(y,z)\in$the type[${\ell}$] with $\text{\rm{$\ell\ge$ j+2}}$}}$
\ms

Note that ${j+1}\le{\ell}$ where $j$ was already given by
{\rm(16.6.1)} and that $n=\Pi^{r}_{i=1}q_i$. $\square$ \ms

$\underline{\text{\bf Remark 16.8.0.}}$ \quad {\rm(i)} For notation,
if $p=\nu+1$, it is said that $f=f(f_{-1},f_0,\dots,f_j,h_p)\in \C
\{f_{-1},f_0,\dots,\dots,f_j\}[h_p]$ satisfies {\rm  The Necessary
Condition[A] for $f(y,z)\in$the type$[\ell]$ with $\ell\ge j+2$}
because $h_p$ of $(h_p,f)$ in {\rm(16.8.4)} satisfies The Necessary
Condition and sufficient condition for $h_p\in$the type$[j+1]$.
Assuming that $p=\nu+1$, if
$\gcd(d_{j+2},\th_{j+2}(\beta_{j+2,0,k})^{j+2}_{k=1})=1$, then {\rm
The Necessary Condition[A] for $f(y,z)\in$the type$[\ell]$ with
$\ell\ge j+2$} may be replaced by {\rm  The Necessary and sufficient
Condition[A] for $f(y,z)\in$the type$[\ell]$ with $\ell= j+2$.} \ms

{\rm(ii)}  Assuming that
$\gcd(d_{j+1},\th_{j+1}(\beta_{j+1,0,k})^{j+1}_{k=1})=1$,
$\gcd(d_{j+2},\th_{j+2}(\beta_{j+2,0,k})^{j+2}_{k=1})=1$ with
$p=\nu+1$, and $\gcd(n_j,\th_j(\a_{j,0,k})^j_{k=1})=1$ for each
$j=1,2,\dots,m$, then note that an equality in {\rm(16.6.9)}, an
equality in {\rm(16.8.8)} and an equality in {\rm(16.8.12)} are the
necessary and sufficient condition for $f$ to be irreducible in
$\C\{y,z\}$. \ms

{\rm(iii)} As soon as the construction in {\rm {Case(1)}} and {\rm
{Case(2)}} is done, then $f$ of $(f_{j+1},f)$ in {\rm(16.8.2)}
satisfies the same kind of condition for {\rm  The Necessary
Condition[A] for $f(y,z)\in$the type$[\ell]$ with $\ell\ge j+2$}, as
$f$ of $(f_{j},f)$ in {\rm(16.7.1)} does for {\rm  The Necessary
Condition[A] for $f(y,z)\in$the type$[\ell]$ with $\ell\ge j+1$}},
as we have seen in the proof of the theorem. \ms

\noindent{\bf Corollary 16.8.1.} For a given integer $j+1$ in the
theorem, we can prove by Theorem 15.4, Proposition $16.7$ and
Proposition $16.8$ that there exists a unique sequence of
irreducible $W$-polys in $z$, $\{f_0=z,\ f_1,\dots, f_{j+1}\}$ with
$f_k\in \C\{y\}[z]$ and $f_k\in $ the type $[k]$ for $1\le k\le j+1$
and $f_{j+1}\ne f$, satisfying the desired properties and notations
as in the conclusion of the theorem.  $\square$ \endproclaim \ms

\proclaim{Proposition 16.9(Step III)} $\underline{\text{\bf
Assumptions}}$ For the induction proof on the integer $(j+1)$ of
this proposition, we may assume by the induction on the integer
$(j+1)$ that the assumptions and the conclusions of {\rm Proposition
16.7(Step(I)) and Proposition 16.8(Step(II))} have been already
shown. \ms

$\underline{\text{\bf Conclusions}}$ We are going to prove the
remaining part of {\rm Theorem 16.6}, equivalently, the following:

$\underline{\text{\rm Step III-1}}$ We prove by {\rm Step I and Step
II} that for a given integer $j+1\le \ell-1$ in this theorem we can
construct a sequence of irreducible $W$-polys in $z$, $\{f_0=z,\
f_1,\dots, f_{j+1}\}$ with $f_k\in \C\{y\}[z]$ and $f_k\in $ the
type $[k]$ for $1\le k\le j+1$ and $f_{j+1}\ne f$, such that
{\rm(i)} each $f_k$ satisfies An algorithm for finding an
irreducibility criterion for $f_k(y,z)\in$the type$[\text{\rm {k}}]$
and {\rm(ii)} $f$ satisfies An algorithm for finding an
irreducibility criterion for $f(y,z)\in$the type$[\ell]$ with
$\ell\ge j+2$ as we have seen in the conclusion of the theorem,
respectively.

$\underline{\text{\rm Step III-2}}$ A sequence of irreducible
$W$-polys satisfying the desired properties and notations as in {\rm
Step III-1} must be unique. $\square$ \endproclaim \ms

\noindent{\bf The Proof of Proposition 16.7.} We prove that
$g_{j+1}$ of $(g_{j+1},f)$ in (16.7.3) satisfies The Necessary and
Sufficient Condition[A] for $g_{j+1}(y,z)\in$the type$[j+1]$ and
that $f$ of $(g_{j+1},f)$ in (16.7.3) satisfies The Necessary
Condition[B] for $f(y,z)\in$the type$[\ell]$ with $\ell\ge j+2$,
respectively.

$\underline{\text{\rm The proof of The Necessary and Sufficient
Condition[A] for $g_{j+1}(y,z)\in$the type$[j+1]$:}}$

(1)(1a) The proof is clear by (16.6.10).  \ms

(1b) Let $\Pi^{j+1}_{k=1}f^{\beta_k}_{k-2}$ and
$\Pi^{j+1}_{k=1}f^{\de_k}_{k-2}$ be arbitrary nonzero monomials in
$\C\{y,z,f_1,\dots,f_{j-1}\}$ where $\beta_1>0$, $\beta_k<n_{k-1}$
for $2\le k\le j+1$, $\delta_1>0$ and $\delta_k<n_{k-1}$ for $2\le
k\le j+1$. Then, note by the definition of $\th_{j+1}$ that
$\th_{j+1}(\beta_k)^{j+1}_{k=1}=\th_{j+1}(\delta_k)^{j+1}_{k=1}$ if
and only if
$$
\beta_{j+1}\th_j(\a_{j,0,k})^{j}_{k=1} +n_j\th_j(\beta_k)^j_{k=1}
=\delta_{j+1}\th_j(\a_{j,0,k})^j_{k=1}+n_j\th_j(\delta_k)^j_{k=1}.
\tag 16.7.13
$$

Since $\gcd(n_j,\th_j(\a_{j,0,k})^j_{k=1})=1$, $0\le
\beta_{j+1}<n_j$ and $0\le \delta_{j+1}<n_j$, then \noindent
$(\beta_{j+1}-\delta_{j+1})\th_j(\a_{j,0,k})^j_{k=1} =
n_j(\th_j(\delta_k)^j_{k=1}-\th_j(\beta_k)^j_{k=1})$ if and only if
$\beta_{j+1}=\delta_{j+1}$ and
$\th_j(\beta_k)^j_{k=1}=\th_j(\delta_k)^j_{k=1}$. Next,
$\th_j(\beta_k)^j_{k=1}=\th_j(\delta_k)^j_{k=1}$ if and only if
$\beta_j\th_{j-1}(\a_{j-1,0,k})^{j-1}_{k=1}+n_{j-1}\th_{j-1}(\beta_k)^{j-1}_{k=1}
=\delta_j\th_{j-1}(\a_{j-1,0,k})^{j-1}_{k=1}+n_{j-1}\th_{j-1}(\delta_k)^{j-1}_{k=1}$
by the definition of $\th_{j-1}$.

Since $\gcd(n_{j-1},\th_{j-1}(\a_{j-1,0,k})^j_{k=1})=1$, $0\le
\beta_j<n_{j-1}$ and $0\le \delta_j<n_{j-1}$, then
$(\beta_{j}-\delta_{j})\th_{j-1}(\a_{j-1,0,k})^{j-1}_{k=1}=
n_{j-1}(\th_{j-1}(\delta_k)^{j-1}_{k=1}-\th_{j-1}(\beta_k)^{j-1}_{k=1})$
if and only if $\beta_j=\delta_j$ and
$\th_{j-1}(\beta_k)^{j-1}_{k=1}=\th_{j-1}(\delta_k)^{j-1}_{k=1}$.

Continuing the above process inductively, it can be easily shown
that
$\th_{j+1}(\beta_k)^{j+1}_{k=1}=\th_{j+1}(\delta_k)^{j+1}_{k=1}$ if
and only if for each $k=1,2,\dots, j+1$, $\beta_k=\delta_k$, by
using the fact that $\gcd(n_p,\th_p(\a_{j,0,k})^p_{k=1})=1$ for
$p=1,2,\dots,j$. Thus the proof of (1b) is done. \ms

(1c) For brevity of notation let $\{g_k: k=1,2,\dots,j+1\}$ be in
$\BC\{y,z\}$ such that $g_k=f_k$ for $1\le k\le j$ where
$f_k=f(y,z,f_1,\dots,f_{k-1})\in
\C\{y,z,f_1,\dots,f_{k-2}\}[f_{k-1}]$ for $k=1,2,\dots,j$ was
already defined in this theorem and
$g_{j+1}=g_{j+1}(y,z,f_1,\dots,f_j)\in
\C\{y,z,f_1,\dots,f_{j-1}\}[f_j]$ was already constructed by an
equation in (16.7.1). Then, it is clear without any need of proof
that $g_1,g_2,\dots,g_j$ are irreducible elements in $\BC\{y,z\}$,
which satisfy the same kind of properties and notations as in the
assumption of Theorem $14.0$ and Proposition $14.1$. In order to
apply Sublemma $16.6.1$ or Theorem $14.0$ to a sequence $\{g_k:
k=1,2,\dots,j+1\}$, it suffices to show that $g_{j+1}$ satisfies the
following inequality:
$$
\f{\th_{j+1}(\sigma_k)^{j+1}_{k=1}}{n_{j+1}}>n_j\th_j(\a_{j,0,k})^{j}_{k=1}.
\tag 16.7.14
$$
which is equivalently rewritten as follows:
$$
\f{\th_{j+1}(\beta_{j+1,d_{j+1}-n_{j+1},k})^{j+1}_{k=1}}{d_{j+1}-(d_{j+1}-n_{j+1})}
>n_j\th_j(\a_{j,0,k})^j_{k=1}.
\tag 16.7.15
$$

Since the inequality in (16.7.15) was already proved by (16.6.12) in
the induction assumption of this theorem, there is nothing to prove
for the inequality in (16.7.14). Following the same method and
notations as we have used in the conclusion of Sublemma $16.6.1$,
let $\tau_{\la_{j}}$ be the composition of a finite number $\la_{j}$
of successive blow-ups which is needed only to get the standard
resolution of the singular point of $V(f_j)$.

After the same number $\la_j$ of blow-ups with the same local
coordinates as we have done in (16.6.17) and (16.7.7), the local
defining equation $(g_{j+1}\circ\tau_{\la_j})_{proper}$ for the
$\la_j$-th proper transform of the curve defined by $g_{j+1}=0$ is
given analytically by
$$\align
 (g_{j+1}\circ\tau_{\la_j})_{total}
&=v^{n_jn_{j+1}\th_j(\a_{j,0,k})^j_{k=1}}(g_{j+1}\circ\tau_{\la_j})_{proper}
\tag 16.7.16 \\
(g_{j+1}\circ\tau_{\la_j})_{proper} &=(u+a+\ve)^{n_{j+1}}
+c_{j+1}(u,v)v^{M'},
\endalign$$
where $a$ is a nonzero constant, $\ve$ is a nonunit along $v=0$ and
$c_{j+1}(u,v)$ is a unit in $\C\{u+a,v\}$ with $c_{j+1}(-a,0)=\zeta$
and
$M'=\th_{j+1}(\sigma_k)^{j+1}_{k=1}-n_{j+1}n_j\th_j(\a_{j,0,k})^{j}_{k=1}$.
\ms

(1d) It is clear by (16.6.17) and (16.6.21.2) that
$\gcd(n_{j+1},\th_{j+1}(\si_k)^{j+1}_{k=1})=\gcd(n_{j+1},M')=1$ and
$$\align
(16.7.17) \qquad \qquad  M' &=\th_{j+1}(\si_k)^{j+1}_{k=1}
-n_jn_{j+1}\th_j(\a_{j,0,k})^j_{k=1}\\
&=\th_{j+1}(\beta_{j+1,d_{j+1}-n_{j+1},k})^{j+1}_{k=1}-
n_j(d_{j+1}-(d_{j+1}-n_{j+1}))\th_j(\a_{j,0,k})^j_{k=1} \qquad \qquad\\
&=M_{j+1,d_{j+1}-n_{j+1}}.\\
\endalign
$$

(2)(2a) Since $\gcd(n_{j+1},M')=1$ with
$M'=M_{j+1,d_{j+1}-n_{j+1}}$, then it is trivial by (16.7.16) that
$(g_{j+1}\circ\tau_{\la_j})_{proper}$ is irreducible in
$\C\{u+a,v\}$, and so $g_{j+1}\in \C\{y\}[z]$ is an irreducible
$W$-poly in $z$ by construction as we have done in (16.7.1) and
Theorem $15.4$. Now, for any given nonzero monomial
$\Pi^{j+1}_{k=1}f^{\beta_{j+1,i,k}}_{k-2}$ in $S_{j+1,i}\in
\C\{y\}[z,f_1,\dots,f_{j-1}]$, apply (ii) of Theorem 15.2 to
$f=f^{d_{j+1}}_j+\Sigma^{d_{j+1}-2}_{i=0} S_{j+1,i}f_j^i$. Then,
$S_{j+1,i}\in \C\{y\}[z]$ has a multiplicity $m_i\ge
(d_{j+1}-i)\Pi^j_{k=1}n_k$ at $0\in \C^2$ for each $i$. If
$i=d_{j+1}-n_{j+1}$, then $m_i\ge \Pi^{j+1}_{k=1}n_k$, and so the
proof is done because $f_j\in \C\{y,z\}$ has a multiplicity
$\Pi^j_{k=1}n_j$ at $0\in \C^2$. \ms

(2b) Observe that $0<n_1<\a_{1,0,1}$ and $2\le n_k$ for $1\le k\le
j+2$, and $\si_1>0$ and $0\le \si_k<n_{k-1}$ for $2\le k\le j+1$.
Since $g_{j+1}$ is irreducible in $\C\{y,z\}$, then by Lemma $3.1$
$g_{j+1}(y,z)$ of (16.7.1) can be represented in the form
$$
g_{j+1}=(z^{n_1}+\xi y^{\a_{1,0,1}})^{d} +\sum_{p,q\ge 0}
c_{p,q}y^pz^q \quad \text{with $n_1p+\a_{1,0,1}q>n_1\a_{1,0,1}d$},
\tag 16.7.18
$$
where the $c_{p,q}$ are some nonzero complex numbers, and $p\ge 0$
and $q\ge 0$ are integers such that $n_1p+\a_{1,0,1}q
>n_1\a_{1,0,1}d$ with $d=n_2n_3\cdots n_{j+1}$, and $\xi$ is a nonzero number.

Whenever $y,z,f_1,\dots,f_j$ is viewed as independent complex
$(j+2)$-variables, in order to prove for $g_{j+1}\in
\C\{y\}[z,f_1,\dots,f_j]$ to have multiplicity $n_{j+1}$ at the
origin in $\C^{j+2}$, it suffices to observe the following
computation:
$$\align
(16.7.19) \qquad \qquad  & \a_{1,0,1}n_1n_2\cdots n_j(\si_1+\si_2+\cdots +\si_{j+1}) \\
>\ & n_1\si_1 +\a_{1,0,1}(\si_2+n_1\si_3+n_1n_2\si_4+\cdots +n_1n_2\cdots
n_{j-1}\si_{j+1}) \qquad  \qquad \qquad  \qquad\\
=\ &n_1\si_1+\a_{1,0,1}\tau \ge \a_{1,0,1}n_1n_2\cdots n_{j+1},
\endalign$$
by (16.7.18) because $y^{\si_1}z^{\tau}$ belongs to
$\Pi^{j+1}_{k=1}f^{\si_k}_{k-2}\in\C\{y\}[z]$  where
$\tau=\si_2+n_1\si_3+n_1n_2\si_4+\cdots +n_1n_2\cdots
n_{j-1}\si_{j+1}$.

Thus, we proved by (16.7.19) that $\si_1+\si_2+\cdots
+\si_{j+1}>n_{j+1}$, and so $g_{j+1}\in \C\{y\}[z,f_1,\dots,f_j]$
has a multiplicity $n_{j+1}$ at the origin in $\C^{j+2}$. Also, it
is trivial by (2a) that $g_{j+1}\in \C \{y,z,f_1,\dots,f_{j-1}\}
[f_j]$ is an irreducible $W$-poly in ${f_{j}}$ with coefficients in
$\C \{y,z,f_1,\dots,f_{j-1}\}$ because ${}_n\CO_0$ is a unique
factorization domain.

Therefore, we finished the proof of The Necessary and Sufficient
Condition[A] for $g_{j+1}(y,z)\in$the type$[j+1]$. \ms

$\underline{\text{\rm The proof of The Necessary Condition[B] for
$f(y,z)\in$the type$[\ell]$ with $\ell\ge j+2$}}$

(1) (1a) To prove it, apply the WDT with a divisor $g_{j+1}$ to $f$.
By Theorem 15.2, $f$ may be written uniquely in the form
$$
f=g^{d_{j+2}}_{j+1}+\sum^{d_{j+2}-1}_{i=0} T_{j+2,i}g^i_{j+1} \quad
\text{with $T_{j+2,i}\in \C\{y,z\}$}, \tag 16.7.20
$$
where  $g_{j+1}\in \C\{y,z\}$ and $T_{j+2,i}\in \C\{y,z\}$. \ms

Whenever $y,z,f_1,\dots, f_j$ are viewed as independent complex
$(j+2)$-variables at the origin in $\C^{j+2}$, in order to prove
that $T_{j+2,i}\in \C\{y\}[z,f_1,\dots,f_j]$ and $T_{j+2,i}$
satisfies the property (1a), first of all note that $g_{j+1}$
satisfies the facts corresponding to the facts Fact(A) and Fact(B)
in the conclusions of Sublemma 15.4.{$\alpha$} of Theorem 15.4 which
have been already satisfied by $h_{j+1,1}$ of (15.4.3) in  the
conclusions of Sublemma 15.4.{$\alpha$}.

Also, as we have seen in the proof of Sublemma 15.4.{$\alpha$},
whenever a pair $(h_{j+1,1},f)$ in (15.4.3) satisfies the facts
Fact(A) and Fact(B), recall that the pair $(h_{j+1,1},f)$ satisfies
the facts Fact(C), Fact(D) and Fact(E) as we have seen in the proof
of Sublemma 15.4.{$\alpha$}. So, replace $(h_{j+1,1},f)$ by
$(g_{j+1},f)$ in the sense of Definition 16.2.2, and then there is
nothing to prove for (1a). \ms

(1b) The proof of the uniqueness of a function
$\bar{\th}_{j+2}:\N_0^{j+2}\to \N_0$ of (16.7.10) can be easily
done, by using the same method as we have used in the proof of the
uniqueness of a function $\th_{j+1}:\N_0^{j+1}\to \N_0$ of (16.7.5)
for (1b) of The Necessary and Sufficient Condition[A] for
$g_{j+1}(y,z)\in$the type$[j+1]$. \ms

(1c) In preparation for the proof of an inequality in (16.7.11), it
is needed to show by (1a) and (1b) that for each fixed $i$,
$T_{j+2,i}$ of (16.7.3) can be rewritten in the form
$$
T_{j+2,i}=\sum^{n_{j+1}-1}_{q=0} T^{(1)}_{j+2,i,q}f^{q}_j, \tag
16.7.21
$$
with two properties (i) and (ii):

(i) Let $i$ be fixed with $0\le i\le d_{j+2}-1$. For each
$q=0,1,\dots,n_{j+1}-1$ and for any nonzero monomial
$\Pi^{j+1}_{k=1}f^{\de_k}_{k-2}$ in $T^{(1)}_{j+2,i,q}\subseteq
\C\{y\}[z,f_1,\dots,f_{j-1}]$, $\de_1>0$ and $0\le\de_k<n_{k-1}$ for
$k=2,3,\dots,{j+1}$.

In particular, if $i=d_{j+2}-1$, then ${q}\le n_{j+1}-2$.

(ii) Let $i$ and $q$ be fixed. For any two nonzero monomials
$\Pi^{j+1}_{k=1}f^{\g_k}_{k-2}$ and $\Pi^{j+1}_{k=1}f^{\de_k}_{k-2}$
in $T^{(1)}_{j+2,i,q}\subseteq \C\{y\}[z,f_1,\dots,f_{j-1}]$, then
$$\align
(16.7.22) \qquad \qquad &
\text{$\th_{j+1}(\g_k)^{j+1}_{k=1}=\th_{j+1}(\de_k)^{j+1}_{k=1}$
if and only if $\g_k=\de_k$ for $k=1,2,\dots,j+1$.} \qquad \qquad\qquad \qquad \\
& \text{So, there exists a unique nonzero monomial
$C^{(1)}_{j+2,i,q}\Pi^{j+1}_{k=1}f^{\beta^{(1)}_{i,q,k}}_{k-2}$ in
$T^{(1)}_{j+2,i,q}$}\\
& \text{with a constant $C^{(1)}_{j+2,i,q}$ such that
$\th_{j+1}(\beta^{(1)}_{j+1,q,k})^{j+1}_{k=1}
=\text{$\min$}\{\th_{j+1}(\de_k)^{j+1}_{k=1}\}$}\\
& \text{for any nonzero monomial $\Pi^{j+1}_{k=1}f^{\de_k}_{k-2}$ in
$T^{(1)}_{j+2,i,q}$ with $q$ fixed.}
\endalign$$

It is trivial that the above summation for $T^{(1)}_{j+2,i,q}$ in
(16.7.21) and (16.7.22) can be clearly constructed by the similar
method as we have seen in the proof of (1a) and (1b). Thus, $f$ can
be rewritten in the form
$$\align
(16.7.23) \quad f &=g^{d_{j+2}}_{j+1}+\sum^{d_{j+2}-1}_{i=0} T_{j+2,i}g^i_{j+1}\\
&=(f^{n_{j+1}}_j+\xi_{j+1}\Pi^{j+1}_{k=1}f^{\si_k}_{k-2})^{d_{j+2 }}
+\sum^{d_{j+2}-1}_{i=0}\sum^{n_{j+1}-1}_{q=0}T^{(1)}_{j+2,i,q}f^{q}_j(f^{n_{j+1}}_j
+\xi_{j+1}\Pi^{j+1}_{k=1}f^{\si_k}_{k-2})^i \qquad
\endalign$$
for $0\le i\le d_{j+2}-1$ and $0\le q\le n_{j+1}-1$, as an element
of $\C\{y\}[z,f_1,\dots,f_j]$, noting that if $i=d_{j+2}-1$ then
$q=n_{j+1}-2$.

By (16.6.17) and (16.7.16), recall that the equation
$(f\circ\tau_{\la_j})_{total}$ for $f(y,z)$ of (16.7.23) can be
written in the form
$$\align
\text{\rm(16.7.24)} \qquad \qquad \qquad
(f\circ\tau_{\la_j})_{total}&=
v^{n_jd_{j+1}\th_j(\a_{j,0,k})^j_{k=1}}(f\circ\tau_{\la_j})_{proper},
\qquad \qquad\qquad \qquad \\
(f\circ\tau_{\la_j})_{proper}&
=(u+a+\ve)^{d_{j+1}}+\sum^{d_{j+1}-2}_{i=0}
W_{j+1,i}(u+a+\ve)^i \quad \text{with}  \qquad \qquad \\
 W_{j+1,i}&=W_{j+1,i}(u,v)=b_{j+1,i}v^{M_{j+1,i}}\quad \text{and}\\
 M_{j+1,i}&=\th_{j+1}(\beta_{j+1,i,k})^{j+1}_{k=1}-n_j(d_{j+1}-i)
\th_j(\a_{j,0,k})^j_{k=1}>0, \\
(g_{j+1}\circ\tau_{\la_j})_{total}
&=v^{n_jn_{j+1}\th_j(\a_{j,0,k})^j_{k=1}}(g_{j+1}\circ\tau_{\la_j})_{proper},
\\
(g_{j+1}\circ\tau_{\la_j})_{proper} &=(u+a+\ve)^{n_{j+1}}
+c_{j+1}(u,v)v^{M'},
\endalign$$
where $d_{j+2}=\gcd(d_{j+1}, M_{j+1,0})=\gcd(d_{j+1},
\th_{j+1}(\beta_{j+1,0,k})^{j+1}_{k=1})>1$ with
$d_{j+1}=n_{j+1}d_{j+2}$ and $M_{j+1,0}=M'd_{j+2}$ for some integers
$n_{j+1}\ge 2$ and $M'\ge 1$, and $a$ is a nonzero constant, $\ve$
is a unit along $v=0$ and $c_{j+1}(u,v)$ is a unit in $\C\{u+a,v\}$
with $c_{j+1}(-a,0)=\zeta$.

By (16.6.16), at $(v,u+a)=(0,0)$ along $v=0$,
$(f_j\circ\tau_{\la_j})_{total}=0$ with
$(f_j\circ\tau_{\la_j})_{proper}=0$,
$(\Pi^{j+1}_{k=1}f^{\de_k}_{k-2}\circ\tau_{\la_j})_{total}=0$ for
any nonzero monomial $\Pi^{j+1}_{k=1}f^{\de_k}_{k-2}\in
T^{(1)}_{j+2,i,q}$, and so
$(T^{(1)}_{j+2,i,q}\circ\tau_{\la_j})_{total}=0$ can be written as
follows:
$$\align
(f_j\circ\tau_{\la_j})_{total}&=
v^{n_j\th_j(\a_{j,0,k})^j_{k=1}}(f_j\circ\tau_{\la_j})_{proper}, \tag 16.7.25\\
(f_j\circ\tau_{\la_j})_{proper}&=(u+a+\ve), \\
(\Pi^{j+1}_{k=1}f^{\de_k}_{k-2}\circ\tau_{\la_j})_{total}&=
v^{\th_{j+1}({\de_k})^{j+1}_{k=1}}b(\de_1,\dots,\de_{j+1}),  \\
(T^{(1)}_{j+2,i,q}\circ\tau_{\la_j})_{total}&=
v^{\th_{j+1}(\beta^{(1)}_{i,q,k})^{j+1}_{k=1}}b_{j+2,i,q},
\endalign$$
where $a$ is a nonzero constant, $\ve$ is a nonunit along $v=0$ and
$b(\de_1,\dots,\de_{j+1})$ is a unit in $\C\{v,u+a\}$ and
$b_{j+2,i,q}$ is a unit in $\C\{v,u+a\}$.

In preparation for the proof of (16.7.11), whenever any nonzero
monomial $\Pi^{j+2}_{k=1}f^{\de_k}_{k-2}\in T_{j+2,i}$ is chosen
arbitrary, for brevity of notation put
$\Pi^{j+1}_{k=1}f^{\de_k}_{k-2}f^{q}_{j}\in
T^{(1)}_{j+2,i,q}f^{q}_j$ with $q=\de_{j+2}$. For any
$\Pi^{j+1}_{k=1}f^{\de_k}_{k-2}f^{q}_{j}\in
T^{(1)}_{j+2,i,q}f^{q}_j$, note by (16.6.16) and (16.6.17) that for
a unit $b=b(\de_1,\dots,\de_{j+1})$ in $\C\{v,u+a\}$,
$$ ((\Pi^{j+1}_{k=1}f^{\de_k}_{k-2}f^{q}_{j})\circ\tau_{\la_j})_{total}=
bv^{\th_{j+1}({\de_k})^{j+1}_{k=1}}v^{qn_j\th_j(\a_{j,0,k})^j_{k=1}}(u+a+\ve)^q.
\tag 16.7.26
$$

Then, it can be easily shown by (16.7.22), (16.7.23), \dots,
(16.7.26) that for each $i=0,1,\dots,d_{j+2}-1$ and for any nonzero
monomial $\Pi^{j+2}_{k=1}f^{\de_k}_{k-2}\in T_{j+2,i}$ with
$q=\de_{j+2}$,
$$
\th_{j+1}(\beta^{(1)}_{i,q,k})^{j+1}_{k=1}
+n_j(q+in_{j+1}-d_{j+1})\th_j(\a_{j,0,k})^j_{k=1}>0. \tag 16.7.27
$$

So, by the induction assumption on the integer $j$ and by Sublemma
$16.6.1$, after the same number $\la_j$ of blow-ups with the same
coordinates as we have used in (16.6.16) and (16.6.17), then the
local defining equation $(f\circ\tau_{\la_j})_{total}=0$ for the
$\la_j$-th total transform of $f(y,z)$ in (16.7.23) can be rewritten
in the following form:
$$\align
(16.7.28) \qquad (f\circ\tau_{\la_j})_{total}&=
v^{n_jd_{j+1}\th_j(\a_{j,0,k})^j_{k=1}}(f\circ\tau_{\la_j})_{proper}\\
(f\circ\tau_{\la_j})_{proper}
&=[(u+a+\ve)^{n_{j+1}}+b_{j+1}v^{M'}]^{d_{j+2}}, \\
&\quad +\sum_i\sum_qb_{i,q}v^{M'_{i,q}}(u+a+\ve)^q
[(u+a+\ve)^{n_{j+1}}+b_{j+1}v^{M'}]^i, \qquad \qquad \\
M'
&=\th_{j+1}(\si_k)^{j+1}_{k=1}-n_jn_{j+1}\th_j(\a_{j,0,k})^j_{k=1}
\quad
\text{with $\gcd(n_{j+1},M')=1$}, \\
M'_{iq} &=\th_{j+1}(\beta^{(1)}_{i,q,k})^{j+1}_{k=1}
+n_j(q+in_{j+1}-d_{j+1})\th_j(\a_{j,0,k})^j_{k=1}>0,
\endalign$$
noting by (16.7.27) that $M'_{i,q}>0$ where $0\le i\le d_{j+2}-1$
and $0\le q=\beta^{(1)}_{i,q,j+2}\le n_{j+1}-1$, and $a$ is a
nonzero constant, and $\ve$ is a nonunit along $v=0$, and $b_{j+1}$
and $b_{i,q}$ are units in $\C\{u+a,v\}$ if exist. Note that
$c_{j+1}$ is a nonunit in $\C\{u+a,v\}$ with $b_{j+1}=c_{j+1}$.

So, when $(f\circ\tau_{\la_j})_{proper}$ is viewed as an element in
$\C\{u+a,v\}$, then
$n_{j+1}M'_{i,q}+qM'+in_{j+1}M'>n_{j+1}M'd_{j+2}=M'd_{j+1}$ because
$(f\circ\tau_{\la_j})_{proper}$ is irreducible in $\C\{u+a,v\}$,
which is equivalently rewritten as
$$
\f{M'_{i,q}}{d_{j+1}-(in_{j+1}+q)}>\f{M'd_{j+2}}{n_{j+1}d_{j+2}}=\f{M'}{n_{j+1}},
\tag 16.7.29
$$
by Lemma 16.0 because $0\le q<n_{j+1}$ and $\gcd(n_{j+1},M')=1$.

By $M'_{i,q}$ and $M'$ of (16.7.28), an inequality in (16.7.29) can
be easily simplified as
$$
\f{\th_{j+1}(\beta^{(1)}_{i,q,k})^{j+1}_{k=1}}{d_{j+1}-(in_{j+1}+q)}
>\f{\th_{j+1}(\si_k)^{j+1}_{k=1}}{n_{j+1}}. \tag 16.7.30
$$
Thus, by the definition of $\ol\th_{j+2}$ in (16.7.10) with
$q=\beta^{(1)}_{i,q,j+2}$ and by (16.7.30), we have
$$\align
(16.7.31) \qquad \qquad
\ol\th_{j+2}(\beta^{(1)}_{i,q,k})^{j+2}_{k=1}
&=q\th_{j+1}(\si_k)^{j+1}_{k=1}+n_{j+1}\th_{j+1}(\beta^{(1)}_{i,q,k})^{j+1}_{k=1}
\qquad \qquad \qquad \qquad \qquad\\
&>q\th_{j+1}(\si_k)^{j+1}_{k=1}
+(d_{j+1}-in_{j+1}-q)\th_{j+1}(\si_k)^{j+1}_{k=1} \qquad \\
&=n_{j+1} (d_{j+2}-i)\th_{j+1}(\si_k)^{j+1}_{k=1},
\endalign$$
because $d_{j+1}=n_{j+1}d_{j+2}$. Since $q$ was chosen arbitrary for
any nonzero monomial $\Pi^{j+2}_{k=1}f^{\de_k}_{k-2} \in
T^{(1)}_{j+2,i,q}f^{\de_{j}}_j$ with $q=\de_{j+2}<n_{j+1}$, we
proved an inequality in (16.7.11). \ms

(1d) For the proof, consider the curve defined by
$v(f\circ\tau_{\la_j})_{proper}=0$ in $\C\{u+a,v\}$ as we have seen
in (16.7.28). As far as the standard resolution of the singular
point $(u+a,v)=(0,0)$ of $v(f\circ\tau_{\la_j})_{proper}=0$ is
concerned, for brevity of notation, we can easily construct a local
biholomorphic mapping $\phi$ from $(u,v)=(-a,0)$ to
$(\bar{u},\bar{v})=(0,0)$ as follows:
$$\align
&\phi(u,v)=(\bar{u},\bar{v})  \tag 16.7.32\\
\text{such that} \quad &\bar{u}=u+a+\ve \quad
\text{and} \quad \bar{v}=(b_{j+1})^{\f{1}{M'}}v. \\
\endalign$$

Let $V(h)$ and $V(H)$ be defined by $h(\bar{u},\bar{v})
=(f\circ\tau_{\la_j})_{proper}\circ{\phi}^{-1}$ and
$H(\bar{u},\bar{v})=(f\circ\tau_{\la_j})_{total}\circ{\phi}^{-1}$,
respectively in terms of the new coordinates.

Then, $(f\circ\tau_{\la_j})_{total}$ of (16.7.28) can be rewritten
in terms of the new coordinates, as follows:
$$\align
\text{(16.7.33)} \qquad \qquad H&=
b_0{\bar v}^{n_jd_{j+1}\th_j(\a_{j,0,k})^j_{k=1}}h,\\
h&=[{\bar u}^{n_{j+1}}+{\bar v}^{M'}]^{d_{j+2}}, \\
&\quad +\sum_i\sum_qb'_{i,q}{\bar v}^{M'_{i,q}}{\bar u}^q
[{\bar u}^{n_{j+1}}+{\bar v}^{M'}]^i, \qquad \qquad \\
M'
&=\th_{j+1}(\si_k)^{j+1}_{k=1}-n_jn_{j+1}\th_j(\a_{j,0,k})^j_{k=1}
\quad
\text{with $\gcd(n_{j+1},M')=1$}, \qquad \qquad\\
M'_{i,q} &=\th_{j+1}(\beta^{(1)}_{i,q,k})^{j+1}_{k=1}
-n_j(d_{j+1}-in_{j+1}-q)\th_j(\a_{j,0,k})^j_{k=1}>0, \qquad \qquad
\endalign$$
where $0\le i\le d_{j+2}-1$ and $0\le q=\beta^{(1)}_{i,q,j+2}\le
n_{j+1}-1$, and $b_0$ and $b'_{i,q}$ are units in
$\C\{\bar{u},\bar{v}\}$ if exist.

As an application of Theorem $16.1$(Theorem $3.7$), we can easily
get the following consequences:

Let $\mu_{\omega}=\pi_{\lambda_j+1}\circ\pi_{\lambda_j+2}
\circ\cdots\circ\pi_{\lambda_j+\omega}:M^{(\lambda_j+\omega)}\to
M^{(\lambda_j)}$ be the composition of a finite number $\omega$ of
successive blow-ups which is needed only to get the standard
resolution of the singular point $(\bar{u},\bar{v})=(0,0)$ of
$V(G)=\{(\bar{u},\bar{v})\in M^{(\lambda_j)}:G=0\}$ where $G={\bar
v}({\bar u}^{n_{j+1}}+{\bar v}^{M'})$, because $M'$ may be equal to
one.

{\rm(i)(ia)} We can use just one coordinate patch of the local
coordinates for each blow-up $\pi_i$ of $\mu_{\omega}$ with $1\le
i\le \omega$ in the sense of {\rm Lemma 2.12}.

{\rm \quad(ib)} Just as above, we can use the same $\mu_{\omega}$
for the composition of the first finite number $\omega$ of
successive blow-ups in preparation for the standard resolution of
the singular point $(0,0)$ of $V(h)$ where
$h(\bar{u},\bar{v})=(f\circ\tau_{\la_j})_{proper}\circ{\phi}^{-1}$.

{\rm \quad(ic)} Also, we can use just the common one coordinate
patch of the given local coordinates for each blow-up $\pi_i$ of the
above $\mu_{\omega}$ in {\rm (ia)}, in order to study any of
$V^{(i)}(h)$ for all $i=1,2,\dots,\omega$ in the sense of Lemma
2.14. \ms

{\rm(ii)} For simplicity of notations, let $(r,s)$ be the common one
of the local coordinates for the $\omega-th$ blow-up
$\pi_m:M^{(\omega)}\to M^{(\omega-1)}$ at $(0,0)$ which is the
quasisingular point of $V^{(\omega-1)}(G)$. Being viewed as an
analytic mapping, $\mu_\omega:M^{(\omega)}\to M^{0}$ can be written
in the form
$$
\mu_\omega(s,r)=(\bar{u},\bar{v})=(s^{M'}r^{\zeta},s^{n_{j+1}}r^{\eta}),
\tag 16.7.34
$$
where

{\rm (iia)} $\zeta$ and $\eta>0$ are nonnegative integers such that
 ${\eta}M'-{\zeta}n_{j+1}=1$,

{\rm (iib)} $E_\omega=\{s=0\}$ is defined by the $\omega-th$
exceptional curve of the first kind. \ms

{\rm(iii)} By {\rm(ii)}, along $s=0$, $(h\circ\mu_\omega)_{total}$
for the $\omega$-th total transform of the curve defined by $h =0$
is given analytically as follows:
$$\align
\text{(16.7.35)} \qquad \qquad \qquad (H\circ\mu_{\omega})_{total}&=
s^{d_{j+1}\th_{j+1}(\si_k)^{j+1}_{k=1}}(h\circ\mu_{\omega})_{proper}
\qquad \qquad \qquad \qquad\\
(h\circ\mu_{\omega})_{total}&=
s^{n_{j+1}d_{j+2}M'}(h\circ\mu_{\omega})_{proper}  \\
(h\circ\mu_{\omega})_{proper} &=(1+r)^{d_{j+2}}
+\sum_i\sum_q b''_{i,q}s^{M''_{i,q}}(1+r)^i,  \qquad \qquad \qquad \\
M''_{i,q} &=n_{j+1}M'_{i,q}-M'(d_{j+1}-in_{j+1}-q)>0,
\endalign$$
noting that $M''_{i,q}>0$, where for $0\le i\le d_{j+2}-1$ and $0\le
q\le n_{j+1}-1$, the $b''_{i,q}$ are units in $\C\{r+1,s\}$, if
exist.

It is clear by (16.7.28), (16.7.33) and the definition of
$\ol\th_{j+2}$ with $q=\beta_{i,q,j+2}$ that
$$\align
(16.7.36) \qquad \qquad M''_{i,q} &=n_{j+1}M'_{i,q}-M'(d_{j+1}-in_{j+1}-q) \\
&=n_{j+1}\{\th_{j+1}(\beta^{(1)}_{i,q,k})^{j+1}_{k=1}
-n_j(d_{j+1}-in_{j+1}-q)\th_{j}(\a_{j,0,k})^{j}_{k=1}\} \qquad \qquad \qquad \\
&\quad -\{\th_{j+1}(\si_k)^{j+1}_{k=1}
-n_jn_{j+1}\th_{j}(\a_{j,0,k})^{j}_{k=1}\}(d_{j+1}-in_{j+1}-q)\\
&=n_{j+1}\th_{j+1}(\beta^{(1)}_{i,q,k})^{j+1}_{k=1}
-(d_{j+1}-in_{j+1}-q)\th_{j+1}(\si_k)^{j+1}_{k=1} \\
&=\ol\th_{j+2}(\beta^{(1)}_{i,q,k})^{j+2}_{k=1}
-n_{j+1}(d_{j+2}-i)\th_{j+1}(\si_k)^{j+1}_{k=1},
\endalign$$
because $\ol\th_{j+2}(\beta^{(1)}_{i,q,k})^{j+2}_{k=1}
=n_{j+1}\th_{j+1}(\beta^{(1)}_{i,q,k})^{j+1}_{k=1}
+q\th_{j+1}(\si_k)^{j+1}_{k=1}$ with $q=\beta_{i,q,j+2}$ and
$d_{j+1}=n_{j+1}d_{j+2}$.

So, it is clear by (16.7.22) that $M''_{i,q}=M''_{i,q'}$ if and only
if $\beta^{(1)}_{i,q,k}=\beta^{(1)}_{i,q',k}$ for $1\le k\le j+2$
with $q=\beta^{(1)}_{i,q,j+2}$. For each fixed $i$ with $0\le i\le
d_{j+2}-1$, let $M''_i$ be defined by $\text{$\min$}\{M''_{i,q}:0\le
q\le n_{j+1}-1\}$. Then, by the definition of $M''_i$, for each
$i=0,1,\dots,d_{j+2}-1$,
$$
M''_i=\ol\th_{j+2}(\beta^{(1)}_{j+2,i,k})^{j+2}_{k=1}
-n_{j+1}(d_{j+2}-i)\th_{j+1}(\si_k)^{j+1}_{k=1}, \tag 16.7.37
$$
because $\ol\th_{j+2}(\beta^{(1)}_{i,k})^{j+2}_{k=1}=
\text{$\min$}\{\ol\th_{j+2}(\beta^{(1)}_{i,q,k})^{j+2}_{k=1}:0\le
\beta_{i,q,j+2}=q\le n_{j+1}-1\}$.

Now, since $(h\circ\mu_{\omega})_{proper}$ is irreducible in
$\C\{r+1,s\}$, then by Theorem 3.2
$$
\f{M''_i}{d_{j+2}-i}\ge \f{M''_0}{d_{j+2}}\quad \text{for $1\le i\le
d_{j+2}-1$}. \tag 16.7.38
$$

Therefore, by (16.7.37) and (16.7.38), we get
$$
\f{\ol\th_{j+2}(\beta^{(1)}_{j+2,i,k})^{j+2}_{k=1}}{d_{j+2}-i}\ge
\f{\ol\th_{j+2}(\beta^{(1)}_{j+2,0,k})^{j+2}_{k=1}}{d_{j+2}}\quad
\text{for $1\le i\le d_{j+2}-1$}. \tag 16.7.39
$$

Thus, the proof of (16.7.12) is done.

(1d-1) If
$\gcd(d_{j+2},\ol\th_{j+2}(\beta^{(1)}_{j+2,0,k})^{j+2}_{k=1})=1$,
then it is trivial by Corollary 3.3 that $f$ is irreducible in
$\C\{y,z\}$ with $f \in$ the type $[j+2]$ if and only if (16.7.39)
holds and $f_1,f_2,\dots,f_j, g_{j+1}$ are irreducible in
$\C\{y,z\}$ just as above, because $g_{j+1} \in$ the type $[j+1]$.

(1d-2) If
$1<\gcd(d_{j+2},\ol\th_{j+2}(\beta^{(1)}_{j+2,0,k})^{j+2}_{k=1})\le
d_{j+2}$, the proof is trivial by Theorem 3.2, Theorem 3.5 and
Theorem 3.6.

Thus we proved (1a), (1b), (1c) and (1d). \ms

(2)(2a) The proof is trivial by (1d) of (1).

(2b) To show that $f=g^{d_{j+2}}_{j+1}+\sum^{d_{j+2}-1}_{i=0}
T_{j+2,i}g^i_{j+1}\in \C\{f_{-1},f_0,\dots,f_{j}\}[g_{j+1}]$ has a
multiplicity $d_{j+2}$ at the origin in $\C^{j+3}$, consider any
nonzero monomial $\Pi^{j+2}_{k=1}f^{\de_k}_{k-2}$ in $T_{j+2,i}$
with $\de_1>0$ and $\de_k<n_{k-1}$ for $k=2,3,\dots,j+2$. Since
$\a_{1,0,1}>n_1>0$ and $\de_1>0$, then
$$\align
(16.7.40) \qquad \qquad & \a_{1,0,1}n_1n_2\cdots n_{j+1}
(\de_1+\de_2+\cdots +\de_{j+2}+i) \\
>\ &n_1\de_1+\a_{1,0,1}(\de_2+n_1\de_3+\cdots +n_1n_2\cdots
n_j\de_{j+2}+n_1n_2\cdots n_{j+1}i) \qquad \qquad \\
\ge \ &\a_{1,0,1}n_1n_2\cdots n_{j+1}d_{j+2}=n_1\a_{1,0,1}d_2,
\endalign$$
by Lemma 3.1 or (16.7.18) because $y^{\de_1}z^\tau$ belongs to
either $\sum^{d_{j+2}-1}_{i=0} T_{j+2,i}g^i_{j+1}$ or
$g^{d_{j+2}}_{j+1}$ as a convergent power series expansion in
$\C\{y,z\}$ where $\tau =\de_2+n_1\de_3+\cdots +n_1n_2\cdots
n_j\de_{j+2}+n_1n_2\cdots n_{j+1}i$, and so
$n_1\de_1+\a_{1,0,1}\tau\a_{1,0,1}\ge n_1n_2\cdots n_{j+1}d_{j+2}$.
Thus, we proved that $\de_1+\de_2+\cdots +\de_{j+1}+i>d_{j+2}$, what
we wanted. So, $f\in \C\{y\}[z,f_1,\dots,f_j,g_{j+1}]$ has a
multiplicity $d_{j+2}$ at the origin in $\C^{j+3}$. Also, it is
trivial that $f\in \C \{y\}[z,f_1,\dots, f_j,g_{j+1}]$ is an
irreducible $W$-poly of degree $d_{j+2}$ in $g_{j+1}$ because
${}_n\CO_0$ is a unique factorization domain and $f$ is irreducible
in $\C\{y,z\}$.

Therefore, we finished the proof of The Necessary Condition[B] for
$f(y,z)\in$the type$[\ell]$ with $\ell\ge j+2$, and so the proof of
this proposition is completed. $\square$ \ms

\noindent$\underline{\text{\bf Remark 16.7.2.}}$

{\rm(1)} Note by (16.7.16) and (16.7.26) that the local defining
equation $(g_{j+1}\circ\tau_{\la_j})_{proper}$ for the $\la_j$-th
proper transform of the curve defined by $g_{j+1}=0$ is analytically
written in the form
$$\align
\text{(16.7.43)} \qquad \qquad (g_{j+1}\circ\tau_{\la_j})_{total}
&=v^{n_jn_{j+1}\th_j(\a_{j,0,k})^j_{k=1}}(g_{j+1}\circ\tau_{\la_j})_{proper}
\qquad \qquad \qquad \qquad\\
(g_{j+1}\circ\tau_{\la_j})_{proper} &=(u+a+\ve)^{n_{j+1}}
+c_{j+1}(u,v)v^{M'}.
\endalign$$

By (16.7.34) and (16.7.36), it is easy to compute that
$$\align
\text{(16.7.44)} \qquad \qquad (g_{j+1}\circ\tau_{\la_j}\circ
\mu_\omega)_{total}
&=s^{n_{j+1}\th_{j+1}(\si_k)^{j+1}_{k=1}}(g_{j+1}\circ\tau_{\la_j}\circ
\mu_\omega)
\qquad \qquad \qquad \qquad\\
(g_{j+1}\circ\tau_{\la_j}\circ \mu_\omega) &=(r+1).
\endalign$$

Let $\tau_{\la_{j+1}}=\tau_{\la_j}\circ \mu_\omega$, and then
$\tau_{\la_{j+1}}$ is the standard resolution of the singular point
$(0,0)$ of $V(g_{j+1})=\{(y,z)\in \C^2:g_{j+1}(y,z)=0\}$.

{\rm(2)} By (1), at $(s,r)=(0,0)$ along $s=0$,
$(\Pi^{j+2}_{k=1}f^{\de_k}_{k-2}\circ\tau_{\la_{j+1}})_{total}=0$
and $(f_j\circ\tau_{\la_j})_{proper}=0$ can be written in the form,
satisfying the following property:
$$\align
(16.7.45) \qquad
(\Pi^{j+2}_{k=1}f^{\de_k}_{k-2}\circ\tau_{\la_{j+1}})_{total}&=
[(\Pi^{j+1}_{k=1}f^{\de_k}_{k-2}\circ\tau_{\la_{j+1}})]
[(f_j)^{\de_{j+2}}\circ\tau_{\la_{j+1}}]\\
&=[(\bar{v}^{\th_{j+1}(\de_k)^{j+1}_{k=1}}b)\circ\mu_\omega]
[(\bar{v}^{\de_{j+2}n_j\th_j(\a_{j,0,k})^j_{k=1}}\bar{u}^{\de_{j+2}})\circ\mu_\omega], \\
&=s^{n_{j+1}\th_{j+1}(\de_k)^{j+1}_{k=1}}s^{\de_{j+2}{n_j}n_{j+1}\th_j(\a_{j,0,k})^j_{k=1}}
s^{\de_{j+2}M'} b'\\
\endalign$$
where $b'$ is a unit in $\{\bar{u},\bar{v}\}$.

Since $M'
=\th_{j+1}(\si_k)^{j+1}_{k=1}-n_jn_{j+1}\th_j(\a_{j,0,k})^j_{k=1}$,
then it is easy to prove that
$$\align
&\text{$n_{j+1}\th_{j+1}(\de_k)^{j+1}_{k=1}+\de_{j+2}{n_j}n_{j+1}\th_j(\a_{j,0,k})^j_{k=1}
+\de_{j+2}M'$} \tag 16.7.46 \\
&\text{$=n_{j+1}\th_{j+1}(\de_k)^{j+1}_{k=1}+\de_{j+2}\th_{j+1}(\si_k)^{j+1}_{k=1}$}\\
&\text{$=\th_{j+2}(\de_k)^{j+2}_{k=1}$.} \\
\endalign$$

{\rm(3)} From the proof of the next proposition, denoted by
Proposition $16.8$, it will be shown that $\tau_{\la_{j+1}}$, which
is the standard resolution of the singular point $(0,0)$ of
$V(g_{j+1})$, is also the standard resolution of the singular point
$(0,0)$ of $V(f_{j+1})$ with \text{$f_{j+1}(y,z)\buildrel \text{{\rm
multiseq}} \over \sim g_{j+1}(y,z)$}. If proved, then by using (2),
it can be said that $\tau_{\la_{j+1}}$ is well-defined, satisfying
an additional property (d) for Fact[4] on the integer $j+1$. \ms

{\rm(4)} By (16.7.37),
$$\align
(f\circ\tau_{\la_{j+1}})_{total}
&=s^{d_{j+1}\th_{j+1}(\si_k)^{j+1}_{k=1}}(f\circ\tau_{\la_{j+1}})_{proper},
\tag 16.7.47 \\
(f\circ\tau_{\la_{j+1}})_{proper} &=(h\circ\mu_{\omega})_{proper}.
\endalign$$ \bs

{\bf Proof of Proposition 16.8.} For the construction of a pair
$(f_{j+1},f)$ in (16.8.2),  it suffices to consider the following
two cases, depending on the fact that $T_{j+2,d_{j+2}-1}$ of
(16.8.1) is either zero or not. For brevity of notations, let
$h_1=g_{j+1}$ and $T^{(1)}_{j+2,i}=T_{j+2,i}$ for $0\le i\le
d_{j+2}-1$.

$\underline{\text{\bf Case(1)}}$ Let $T^{(1)}_{j+2,d_{j+2}-1}$ be
zero. It is clear. \ms

$\underline{\text{\bf Case(2)}}$ Let $T^{(1)}_{j+2,d_{j+2}-1}$ be
nonzero. It has been already shown by Sublemma 15.5 and Sublemma
15.6 in the proof of Theorem 15.4 that the following assertion is
true:

There is a sequence of W-polys in z of pairs,
$\{(h_p,f):p=1,2,\dots\}$ such that
$$\align
(16.8.15) \qquad\qquad (h_{\nu+1},f)=(h_{\nu+2},f)=\cdots \quad
\text{for some integer} \quad \nu\le \f{n_{j+1}+1}2,
\qquad\qquad\qquad\qquad
\endalign$$
each pair of which can be written in the form
$$\cases
h_1 &=f^{n_{j+1}}_j+\xi_{j+1}\Pi^{j+1}_{k=1}f^{\si_k}_{k-2}
=f^{n_{j+1}}_j+R^{(1)}_{j+1,0},\\
f &=h^{d_{j+2}}_1+\sum^{d_{j+2}-1}_{i=0} T^{(1)}_{j+2,i}h^i_1,
\endcases \tag 16.8.16
$$
and for $p=2,3,\dots $
$$\cases
h_{p}& =h_{p-1}+\f
1{d_{j+2}}T^{(p-1)}_{j+2,d_{j+2}-1}=f^{n_{j+1}}_j+\sum^{n_{j+1}-2}_{i=0}
R^{(p)}_{j+1,i}f^i_j,\\
f& =h^{d_{j+2}}_{p}+\sum^{d_{j+2}-1}_{i=0} T^{(p)}_{j+2,i}h^i_{p},
\endcases \tag 16.8.17
$$
with $T^{(p)}_{{j+2},d_{j+2}-1}\ne 0$ for $1\le p\le \nu$ and
$T^{(\nu+1)}_{{j+2},d_{j+2}-1}=T^{(\nu+2)}_{{j+2},d_{j+2}-1}=\cdots
=0$ where $T^{(p)}_{{j+2},i}=T^{(p)}_{{j+2},i}(y,z)\in \C\{y\}[z]$
for $p\ge 1$ and $0\le i\le {d_{j+2}}-1$, and $R^{(p)}_{{j+1},i}
=R^{(p)}_{{j+1},i}(y)\in \C\{y\}$ for $p\ge 1$ and $0\le i\le
n_{j+1}-2$, if exist, satisfying the same kind of the properties and
notations as we have seen in Sublemma $15.5$ of Theorem $15.4$, as
follows: \ms

$\underline{\text{\rm (16.8.18)(16.8.18-1) Property(1)}}$ Let $p$
and $i$ be fixed with $p\ge 1$ and $0\le i\le n_{j+1}-2$. Then
$R^{(p)}_{j+1,i}=R^{(p)}_{j+1,i}(y,z)\in \C\{y\}[z]$ is a polynomial
in $z$ of degree $<\Pi^j_{t=1}n_t$ and has a multiplicity $\ge
(n_{j+1}-i)\Pi^j_{t=1}n_t$ at $0 \in \C^2$.  \ms

$\underline{\text{\rm (16.8.18-2) Property(2)}}$ Let $p$ and $i$ be
fixed with $p\ge 1$ and $0\le i\le d_{j+2}-1$. Then
$T^{(p)}_{j+2,i}=T^{(p)}_{j+2,i}(y,z)\in \C\{y\}[z]$ is a polynomial
in $z$ of degree $<\Pi^{j+1}_{t=1}n_t$ and has a multiplicity $\ge
(d_{j+2}-i)\Pi^{j+1}_{t=1}n_t$ at $0 \in \C^2$. \ms

Consider $y,z,f_1,\dots,f_j$ as independent complex
$(j+2)$-variables at the origin in $\C^{j+2}$.

$\underline{\text{\rm (16.8.18-3) Property(3)}}$ Let $p$ and $i$ be
fixed with $p\ge 1$ and $0\le i\le n_{j+1}-2$. Then for any nonzero
monomial $\Pi^{j+1}_{t=1}f^{\de_t}_{t-2}$ in
$R^{(p+1)}_{j+1,i}=R^{(p+1)}_{j+1,i}(y,z,f_1,\dots,f_{j-1})\in
\C\{y\}[z,f_1,\dots,f_{j-1}]$, $\de_1>0$ and $\de_t <n_{t-1}$ for
$t=2,3,\dots, j+1$.  \ms

$\underline{\text{\rm (16.8.18-4) Property(4)}}$ Let $p$ and $i$ be
fixed with $p\ge 1$ and $0\le i\le d_{j+2}-1$. Then for any nonzero
monomial $\Pi^{j+2}_{t=1}f^{\gamma_t}_{t-2}$ in
$T^{(p)}_{j+2,i}=T^{(p)}_{j+2,i}(y,z,f_1,\dots,f_j)\in
\C\{y\}[z,f_1,\dots,f_j]$, $\gamma_1>0$ and $\gamma_t<n_{t-1}$ for
$t=2,3,\dots,j+2$. \ms

$\underline{\text{\rm (16.8.18-5) Property(5)}}$ In particular, if
$i=d_{j+2}-1$ for $T^{(p)}_{j+2,i}$ of {\rm Property(4)}, then
$\gamma_{j+2}\le n_{j+1}-2$. \bs

$\underline{\text{\rm (16.8.18-6) Property(6)}}$ There is a pair
$(h_{\nu+1},f)\in H$ which satisfies the following:

There is an integer $\nu\le \f{n_{j+1}+1}2$ such that
$T^{(p)}_{j+2,d_{j+2}-1}\ne 0$ for $p=1,2,\dots,\nu$ and
$T^{(\nu+1)}_{j+2,d_{j+2}-1}=T^{(\nu+2)}_{j+2,d_{j+2}-1}=\cdots =0$.
That is, $(h_{\nu},f)\not= (f_{j+1},f)$ and $(h_{\nu+1},f)=
(f_{j+1},f)$ for an integer $\nu \le \f{n_{j+1}+1}2$. \ms

{\bf Remark.} Without any need of proof, Property(1),
Property(2),\dots, Property(6), which are mentioned just above,
follow clearly from Sublemma $15.5$ of Theorem $15.4$, which belongs
to Case (II) in the conclusion of Sublemma $15.5$ of Theorem $15.4$.
In Sublemma $15.5$, note that $f_{-1}=y$ and $f_{0}=z$. \ms

For the proof of this proposition in Case(2), it suffices to show
that two properties, denoted by, $\underline{ \text{\rm The
Necessary and Sufficient Condition[A] for $h_{p}(y,z)\in$the
type$[j+1]$}}$ and

\noindent{$\underline{ \text{\rm The Necessary Condition[B] for
$f(y,z)\in$the type$[\ell]$ with $\ell\ge j+2$}}$,} can be satisfied
respectively. Then, the proof will be by induction on the integer
$p\ge 1$.

Now, it is enough to consider the following two subcases for
Case(2), respectively:

Subcase(A) $p=1$, and Subcase(B) $p\ge 1$. \ms

$\underline{\text{\bf Subcase(A) of Case(2)}}$ \quad Let $p=1$. Then
it suffices to show that $(h_1,f)$, given by (16.8.3), satisfies
$\underline{ \text{\rm The Necessary and Sufficient Condition[A] for
$h_{1}(y,z)\in$the type$[j+1]$}}$, and

\noindent{$\underline{ \text{\rm The Necessary Condition[B] for
$f(y,z)\in$the type$[\ell]$ with $\ell\ge j+2$}}$,}

{\noindent}which was already proved by Proposition $16.7$. \ms

$\underline{\text{\bf Subcase(B) of Case(2)}}$ \quad Let $p\ge 1$.
For the proof of this subcase, it suffices to show by Subcase(A) for
Case(2) that the following sublemma is true: \bs

\noindent\text{\bf Sublemma 16.8.1 for Subcase(B) of Case(2).} \quad

$\underline{\text{\bf {Assumptions}}}$ \quad For the induction
proof, suppose we have shown on the integer $p\ge 1$ that
$\underline{ \text{\rm  The Necessary and Sufficient Condition[A]
for $h_{p}(y,z)\in$the type$[j+1]$}}$ and

$\underline{ \text{\rm The Necessary Condition[B] for $f(y,z)\in$the
type$[\ell]$ with $\ell\ge j+2$ }}$ are true for $(h_p,f)$,
following the same notations and properties as we have seen in
(16.8.4), (16.8.5), \dots, (16.8.14) of the statement of this
propositon.
 \ms

$\underline{\text{\bf {Conclusions}}}$ \quad Then, $(h_{p+1},f)$ can
be written by
$$\cases
h_{p+1} &=h_p+\f1{d_{j+2}}
T^{(p)}_{j+2,{d_{j+2}-1}}=f^{n_{j+1}}_j+\sum^{n_{j+1}-2}_{i=0}R^{(p+1)}_{j+1,i}f^i_j, \\
f &=h^{d_{j+2}}_{p+1}
+\sum^{{d_{j+2}}-1}_{i=0}T^{(p+1)}_{j+2,i}h^i_{p+1},
\endcases \tag 16.8.19
$$
with $T^{(p)}_{j+2,{d_{j+2}}-1}\ne 0$ for $1\le p\le \nu$ and
$T^{(\nu+1)}_{j+2,{d_{j+2}}-1}=T^{(\nu+2)}_{j+2,{d_{j+2}}-1}=\cdots
=0$ where $T^{(p)}_{j+2,i}\in \C\{f_{-1},f_0,\dots,f_{j}\}$ for
$p\ge 1$ and $0\le i\le {d_{j+2}}-1$, and $R^{(p+1)}_{j+1,i} \in
\C\{f_{-1},f_0,\dots,f_{j-1}\}$ for $p\ge 1$ and $0\le i\le
n_{j+1}-2$, if exist, satisfying the same kind of properties and
notations, denoted by The Necessary and Sufficient Condition[A] for
$h_{p+1}(y,z)\in$the type$[j+1]$ and The Necessary Condition[B] for
$f(y,z)\in$the type$[\ell]$ with $\ell\ge j+2$, inductively as we
have seen in the conclusion of Proposition 16.8. $\square$ \ms

\noindent$\underline{\text{\bf Proof of Sublemma 16.8.1 for
Subcase(B) of Case(2)}}$ Let $p\ge 1$ with $T^{(p)}_{j+2,d_2-1}\not
=0$.

$\underline{\text{\rm [1] The proof of  The Necessary and Sufficient
Condition[A] for \text{$h_{p+1}(y,z)\in$the type$[j+1]$}}}$

(1) The proof can be easily proved by Sublemma 15.5 of Theorem 15.4
and the induction assumption on the integer $p$.

{\rm(2)} The proofs of (2a) and (2b) are clear by (1), noting that
${}_n\CO$ is a unique factorization domain.  Thus, the proof of  The
Necessary and Sufficient Condition[A] for \text{$h_{p+1}(y,z)\in$the
type$[j+1]$} is done. \ms

$\underline{ \text{\rm [2] The proof of The Necessary Condition[B]
for \text{$f(y,z)\in$the type$[\ell]$} with $\ell\ge j+2$}}$

Using the same method as we have used in the proof of Proposition
16.7 and following the same kind of propertis and notations as we
have seen in (16.7.11), (16.7.12), (16.7.13),\dots, (16.7.37),
(16.7.38), (16.7.39), the proof of The Necessary Condition[B] for
\text{$f(y,z)\in$the type$[\ell]$} with $\ell\ge j+2$ can be easily
done. Thus, the proof of Sublemma 16.8.1 can be done. \ms

$\underline{ \text{\rm [3] The proof of The Necessary Condition[A]
for \text{$f(y,z)\in$the type$[\ell]$} with $\ell\ge j+2$}}$

It is trivial to prove. Therefore, the proof of Proposition 16.8 is
completely finished. $\square$ \ms

{\bf Proof of Proposition 16.9(Step III).}

$\underline{\text{The Proof of Step III-1}}$ \quad For a given
integer $j+1\le \ell-1$ and $(g_{j+1},f)$ in Proposition $16.7$,
the main aim is to construct a unique pair $(f_{j+1},f)$ which can
be uniquely written in the form
$$\cases
f_{j+1} &=f^{n_{j+1}}_{j}+\sum^{n_{j+1}-2}_{i=0} R_{j+1,i}f^i_{j}
\\
f &=f^{d_{j+2}}_{j+1} +\sum^{d_{j+2}-2}_{i=0} S_{j+2,i}f^i_{j+1},
\endcases \tag 16.9.1
$$
where $y,z,f_1,\dots, f_{j+1}$ are considered as independent complex
$(j+3)$-variables at the origin in $\C^{j+3}$ if necessary,
satisfying the following properties: \ms

\noindent (16.9.1-a)  The first problem is to prove that we can
construct $f_{j+1}=f_{j+1}(y,z,f_1,\dots,f_j)$ satisfying The
Necessary and Sufficient Condition[A] for $f_{j+1}(y,z)\in$the
type$[j+1]$ such that \text{$f_{j+1}(y,z)\buildrel \text{{\rm
multiseq}} \over \sim g_{j+1}(y,z)$} in the sense of Proposition
16.8. \ms

\noindent (16.9.1-b)  The second problem is to prove that
$f=f(y,z,f_1,\dots,f_{j+1})$ satisfies The Necessary Condition[A]
for $f(y,z)\in$the type$[\ell]$ with $\ell\ge j+2$, using the same
way as $f=f(y,z,f_1,\dots,f_{j})$ does The Necessary Condition[A]
for $f(y,z)\in$the type$[\ell]$ with $\ell\ge j+1$ in the sense of
Proposition 16.8. \ms

In preparation for finding the proofs in (16.9.1-a) and (16.9.1-b),
it was already shown by (16.7.3) of Proposition $16.7$ that
$(g_{j+1},f)$ can be uniquely written in the following form:
$$
\cases g_{j+1} &=f^{n_{j+1}}_j
+\xi_{j+1}\Pi^{j+1}_{k=1}f^{\si_k}_{k-2}\quad
\text{with} \ f_{-1}=y\ \text{and}\ f_0=z, \\
f &=g^{d_{j+2}}_{j+1}+\sum^{d_{j+2}-1}_{i=0} T_{j+2,i}g^j_{j+1},
\endcases \tag 16.9.2
$$
where $n=d_{j+2}\Pi^{j+1}_{k=1}n_k$ with $d_{j+2}\ge 2$ and
$n_{k}\ge 2$ for $1\le k\le j+1$, and
$\si_k=\beta_{j+1,d_{j+1}-n_{j+1},k}$ for $1\le k\le j+1$ and
$\xi_{j+1}=\f 1{d_{j+2}}B_{j+1,d_{j+1}-n_{j+1}}$, such that
$g_{j+1}$ of $(g_{j+1},f)$ satisfies The Necessary and Sufficient
Condition[A] for $g_{j+1}(y,z)\in$the type$[j+1]$ and $f$ of
$(g_{j+1},f)$ does The Necessary Condition[B] for $f(y,z)\in$the
type$[\ell]$ with $\ell\ge j+2$ in the sense of Proposition $16.7$
or Proposition $16.8$.

To solve the problems in (16.9.1-a) and (16.9.1-a), it suffices to
consider two cases:

\noindent{\rm Case(1)} Let $T_{j+2,d_{j+2}-1}=0$. Then there is
nothing to solve, because $g_{j+1}$ of $(g_{j+1},f)$ satisfies The
Necessary and Sufficient Condition[A] for $g_{j+1}(y,z)\in$the
type$[j+1]$ and $f$ of $(g_{j+1},f)$ does The Necessary Condition[A]
for $f(y,z)\in$the type$[{\ell}]$ with $\ell\ge j+2$ in the sense of
Proposition $16.8$. \ms

\noindent {\rm Case(2)} Let $T_{j+2,d_{j+2}-1}\not=0$. note that
$g_{j+1}$ of $(g_{j+1},f)$ satisfies The Necessary and Sufficient
Condition[A] for $g_{j+1}(y,z)\in$the type$[j+1]$, but note that $f$
of $(g_{j+1},f)$ does not satisfies The Necessary Condition[A] for
$f(y,z)\in$the type$[{\ell}]$ with $\ell\ge j+2$ in the sense of
Proposition $16.8$.

So, for finding the proofs in (16.9.1-a) and (16.9.1-b), it remains
to cosider Case(2). Since $T_{d_{j+2}-1}\not=0$ from (16.9.2), then
$f$ of $(g_{j+1},f)$ of (16.9.2) does not satisfies The Necessary
Condition[A] for $f(y,z)\in$the type$[j+2]$ in the sense of
Proposition $16.7$.

In preparation for the construction of the pair $(f_{j+1},f)$ in
(16.9.1), which satisfies the corresponding properties as we have
seen in Theorem 16.6, it was already known by Proposition $16.8$ and
Theorem $15.4$ that the pair $(h_{\nu+1},f)$ can be rewritten in the
form
$$\cases
h_{\nu+1} &=f^{n_{j+1}}_j+\sum^{n_{j+1}-2}_{i=0}
R^{(\nu+1)}_if^i_j, \\
f &=h^{d_{j+2}}_{\nu+1}+\sum^{d_{j+2}-2}_{i=0}
T^{(\nu+1)}_ih^i_{\nu+1},
\endcases \tag 16.9.3
$$
such that $h_{\nu+1}$ of $(h_{\nu+1},f)$ in (16.9.3) satisfies The
Necessary and Sufficient Condition[A] for $h_{\nu+1}(y,z)\in$the
type$[j+1]$ and $f$ of $(h_{\nu+1},f)$ in (16.9.3) satisfies The
Necessary Condition[A] for $f(y,z)\in$the type$[{\ell}]$ with
$\ell\ge j+2$. Note that $T^{(\nu+1)}_{d_{j+2}-1}$ was shown to be
zero.

Now, to construct such a pair $(f_{j+1},f)$, define $f_{j+1}$ by
$h_{\nu+1}$ of (16.9.3). Then, to find a representation of
$(f_{j+1},f)$ in terms of Theorem 16.6, let $R_{j+1,i}$ be
$R^{(\nu+1)}_i$ for $0\le i\le n_{j+1}-2$, and $S_{j+2,i}$ be
$T^{(\nu+1)}_i$ for $0\le i\le d_{j+1}-2$, following the notations
in Proposition $16.8$.

Then, by all the definitions of
$A_{j+1,i}\Pi^{j+1}_{k=1}f^{\a_{j+1,i,k}}_{k-2}$ and
$A^{(\nu+1)}_i\Pi^{j+1}_{k=1}f^{\a^{(\nu+1)}_{i,k}}_{k-2}$ for $0\le
i\le n_{j+1}-2$, we get that $\a_{j+1,i,k}=\a^{(\nu+1)}_{i,k}$ with
$\a^{(\nu+1)}_{0,k}=\si_k$ for $1\le k\le j+1$ and
$A_{j+1,i}=A^{(\nu+1)}_i$ for $0\le i\le n_{j+1}-2$. So,
$\ol\th_{j+2}$ is well defined, that is, $\ol\th_{j+2}=\th_{j+2}$.
Similarly, by all the definitions of
$B_{j+2,i}\Pi^{j+2}_{k=1}f^{\beta_{j+2,i,k}}_{k-2}$ and
$C^{(\nu+1)}_i\Pi^{j+2}_{k=1}f^{\beta^{(\nu+1)}_{i,k}}_{k-2}$ for
$0\le i\le n_{j+1}-2$, it is clear that
$\beta_{j+2,i,k}=\beta^{(\nu+1)}_{i,k}$ and
$B_{j+2,i}=C^{(\nu+1)}_i$ for $0\le i\le d_{j+2}-2$ and $1\le k\le
j+2$. Also, we already proved that
$\gcd(d_{j+2},\th_{j+2}(\beta_{j+2,0,k})^{j+2}_{k=1})=\gcd(d_{j+2},
\ol\th_{j+2}(\beta^{(\nu+1)}_{0,k})^{j+2}_{k=1})<d_{j+2}$.

Thus, the main problem in this theorem, equivalently, to compute
such a pair $(f_{j+1},f)$ is completely solved, whether
$T_{j+2,d_{j+2}-1}$ of (16.9.2) is zero or not. \ms

$\underline{\text{Step III-2}}$ We are going to prove that for a
given integer $j+1\le \ell-1$ in this theorem we can construct a
unique sequence of irreducible $W$-polys in $z$, $\{f_0=z,
f_1,\dots, f_{j+1}\}$ with $f_k\in \C\{y\}[z]$ for $1\le k\le j+1$
and $f_{j+1}\ne f$, such that each $f_k$ satisfies The Necessary and
Sufficient Condition[A] for $f_{k}(y,z)\in$the type$[k]$ and
$f=f(y,z,f_1,\dots,f_{j+1})$ satisfies The Necessary Condition[A]
for $f(y,z)\in$the type$[\ell]$ with $\ell\ge j+2$ as we have seen
in the conclusion of the theorem.  The proof just follows from
Theorem $15.4$, Step I, Step II and Step III-1. Thus, the proof of
Step III(Proposition 16.9) is completely finished.

Therefore, we can complete the proof of Theorem 16.6. $\square$
 \ms

\proclaim{Corollary 16.10} $\underline{\text{\bf Assumptions}}$
\quad Under the same assumptions and results in Theorem $16.6$, for
coincidence of notations $f$ may be rewritten as follows: Note that
$f_0=z$.
$$
f=z^n+\sum^{n-2}_{i=0} a_iy^{\a_i}z^i=f^n_0+\sum^{n-2}_{i=0}
S_{1i}z^i \tag 16.10.1
$$
where $S_{1,i}=a_iy^{\a_i}$ with $a_i$ units in $\C\{y\}$,
$B_{1i}=a_i(0)$ and $\beta_{1,i,1}=\a_i$. \bs

$\underline{\text{\bf Conclusions}}$ \quad  Then, we find the
following computational algorithm with $d_1=n$: \ms

{\rm(1)} Let $d_2=\gcd(d_1,\beta_{1,0,1})$ with $d_1=n$ and
$\beta_{1,0,1}=\a_0$.

Then we can find $n_1,\a_{1,0,1}$ and $A_{1,0}$ such that
$n=d_1=n_1d_2$, $\beta_{1,0,1}=\a_{1,0,1}d_2$,
$\a_{1,0,1}=\beta_{1,d_1-n_1,1}$ and $d_2A_{1,0}=B_{1,d_1-n_1}$.
Note that
$d_2\th_1(\a_{1,0,1})=d_2\a_{1,0,1}=\beta_{1,0,1}=\th_1(\beta_{1,0,1})=\a_0$.
\ms

{\rm(2)} Recall that $B_{2,i}y^{\beta_{2,i,1}}z^{\beta_{2,i,2}}$ is
a unique nonzero monomial in $S_{2,i}$ with a nonzero constant
$B_{2,i}$ such that
$\th_2(\beta_{2,i,k})^2_{k=1}=\text{$\min$}\{\th_2(\de_k)^2_{k=1}\}$
for any nonzero monomial $y^{\de_1}z^{\de_2}$ in $S_{2,i}$, and then
let $d_3=\gcd(d_2,\th_2(\beta_{2,0,k})^2_{k=1})$.

Then we can find $n_2$, $(\a_{2,0,k})^2_{k=1}$, and $A_{2,0}$ such
that $d_2=n_2d_3$, $\a_{2,0,k}=\beta_{2,d_2-n_2,k}$ for $k=1,2$ and
$d_3A_{2,0}=B_{2,d_2-n_2}$. Note that
$d_3\th_2(\a_{2,0,k})^2_{k=1}=\th_2(\beta_{2,0,k})^2_{k=1}$. \ms

{\rm(3)(i)} Let $m$ be such that $3\le m\le j$. Recall that
$B_{m,i}\Pi^m_{k=1}f^{\beta_{m,i,k}}_{k-2}$ is a unique nonzero
monomial in $S_{m,i}$ with a nonzero constant $B_{m,i}$ such that
$\th_m(\beta_{m,i,k})^m_{k=1}=\text{$\min$}\{\th_m(\de_k)^m_{k=1}\}$
for any nonzero monomial $\Pi^m_{k=1}f^{\de_k}_{k-2}$ in $S_{mi}$,
and let $d_{m+1}=\gcd(d_m,\th_m(\beta_{m,0,k})^m_{k=1})$.

Then we can find $n_m,\ (\a_{m,0,k})^m_{k=1}$ and $A_{m,0}$ such
that $d_m=n_md_{m+1}$, $\a_{m,0,k}=\beta_{m,d_m-n_m,k}$ for $1\le
k\le m$, and $d_{m+1}A_{m,0}=B_{m,d_m-n_m}$.

Moreover, $\{\a_{m,0,k}:k=1,2,\dots\}$ is a unique sequence such
that $d_m\th_m(\a_{m,0,k})^m_{k=1}=\th_m(\beta_{m,0,k})^m_{k=1}$,
and $\a_{m,0,1}>0$, $\a_{m,0,k}<n_{k-1}$ for $2\le k\le m$, and
$\beta_{m,0,1}>0$, $\beta_{m,0,k}<n_{k-1}$ for $2\le k\le m$. \ms

{\rm(ii)} Note that $n=d_1>d_2=\gcd(d_1,\beta_{1,0,1})>d_3 =
\gcd(d_2,\th_2(\beta_{2,0,k})^2_{k=1})>d_4 =
\gcd(d_3,\th_3(\beta_{3,0,k})^3_{k=1})>\cdots >d_{j+1} =
\gcd(d_j,\th_j(\beta_{j,0,k})^j_{k=1})$. Thus,
$\{d_1,d_2,\dots,d_{j+1}\}$ is a strictly decreasing positive
integer sequence. \ms

{\rm(4)} Let $H_1=z^{n_1}+y^{\a_{1,0,1}}$,
$H_2=H_1^{n_2}+y^{\a_{2,0,1}}z^{\a_{2,0,2}}$, $H_3=H^{n_3}_2
+\Pi^3_{k=1}H^{\a_{3,0,k}}_{k-2},\cdots,$ and
$H_m=H^{n_m}_{m-1}+\Pi^m_{k=1}H^{\a_{m,0,k}}_{k-2}$ for $1\le m\le
j$ where $H_{-1}=y$ and $H_0=z$.

Then, \text{$f_p\buildrel \text{{\rm multiseq}} \over \sim H_p$} for
$1\le p\le j$. $\square$ \endproclaim \ms

{\bf Remark 16.11.} Without assuming that $f$ is irreducible in
${}_2\CO_0$, suppose that $f$ satisfies the same assumption as in
Theorem 16.6. For a fixed $j$ with $0\le j\le \ell-1$ suppose also
that there exists such a unique sequence of irreducible $W$-polys in
$z$, $\{f_0,f_1,\dots,f_j\}$, satisfying the same conclusion as in
Theorem 16.6. Again, except that $F_{\la_j}$ is irreducible in
$\C\{u+a,v\}$ as in additional condition (d), assume that
$F_{\la_j}$ with the property $(16.6)$ satisfies the construction of
$g_{j+1}(y,z,f_1,\dots,f_j)$ and $f(y,z,f_1,\dots,f_j,g_{j+1})$ as
in Step I and II. Then we can prove by Step III that each $h_p$ is
irreducible in ${}_2\CO_0$ with $h_1=g_{j+1}$ and $p\ge 1$. Note
that $f$ may not be irreducible in ${}_2\CO_0$.

\vfill \pagebreak

{\bf Part[C4] The 2nd Algorithm for computing irreducible W-polys
from all {\indent}the W-polys of two complex variables and The 3rd
algorithm for computing {\indent}the corresponding standard Puiseux
expansion from any irreducible W-poly {\indent}of two complex
variables} \bs

{\bf\S {19}. The 2nd Algorithm for computing completely irreducible
W-polys from all the W-polys of two complex variables with proofs}
\ms

\proclaim{Theorem 19.1(The 2nd Algorithm)}

$\underline{\text{\bf Assumptions}}$ Suppose that Theoorem 19.1 and
Theorem 1.15 have the same assumption.

$\underline{\text{\bf Conclusions}}$ Then, they have have the same
conclusion. $\square$
\endproclaim \ms
\demo{Proof} There is nothing to prove for Theorem 19.1 because it
was already proved by Theorem 16.5 and by Theorem 16.6 with
Proposition 16.7 and Proposition 16.8. $\square$
\enddemo \ms

\noindent{\bf Example 19.1.1 for The 2nd algorithm in Theorem 19.1:}
Example 19.1.1 and Example 1.10.1 of $\S 1.10$ are the same. \quad
\text{$\square$} \bs

{\bf \S {20}. The 3rd Algorithm for computing the corresponding
standard Puiseux expansion from any irreducible W-poly of two
complex variables with respect to the multiplicity sequences with
proofs} \ms

\proclaim{Theorem 20.1(The 3rd Algorithm)}

$\underline{\text{\bf Assumptions}}$ Suppose that Theoorem 20.1 and
Theorem 1.16 have the same assumption.

$\underline{\text{\bf Conclusions}}$ Then, they have have the same
conclusion. $\square$
\endproclaim \ms

\demo{Proof} There is nothing to prove for Theorem 20.1 because it
was already proved by Theorem 16.5 and by Theorem 16.6 with
Proposition 16.7 and Proposition 16.8. $\square$
\enddemo \ms

\noindent{\bf Example 20.1.1 for The 3rd algorithm in Theorem 20.1:}
Example 20.1.1 and Example 1.10.2 of $\S 1.10$ are the same. \quad
\text{$\square$} \bs

\vfill \pagebreak

$$\text{\bf Appendix}$$ \bs

It is very interesting to study what problems can be computed in
irreducible plane curve singularities in algebraic geometry.
Appendix consists of Appendix A, Appendix B and Appendix C. As an
application of The 1st Algorithm, The 2nd Algorithm and The 3rd
Algorithm of Part[A], the aim is to find three algorithms, called
The $\alpha$-algorithm, The $\beta$-algorithm, The
$\gamma$-algorithm in Appendix A, Appendix B and Appendix C,
respectively, completely and rigorously in an elementary way, as
follows: \ms

$\underline{\text{\bf [1] In Appendix A,}}$ the aim is to find an
explicit algorithm for finding a one-to-one function from
{{Family(2)}(the family of the standard Puiseux expansion)} onto
{{Family(3)}(the family of all the multiplicity sequences of
irreducible plane curves with isolated singularity under the
standard resolution)}, called The $\alpha$-algorithm.
$\underline{\text{\rm By Theorem A.4 and Theorem A.5,}}$ we can find
The $\alpha$-algorithm completely and rigorously in an elementary
way. \ms

$\underline{\text{\bf [2] In Appendix B,}}$ the aim is to find an
explicit algorithm for finding a one-to-one function from Family(2)
onto {{Family(4)}}(the family of the divisors defined by the total
transforms of irreducible plane curves with isolated singularity
under the standard resolution), called The $\beta$-algorithm.
$\underline{\text{\rm By Theorem B.2 and Theorem B.3,}}$ we can find
The $\beta$-algorithm completely and rigorously in an elementary
way. \ms

$\underline{\text{\bf [3] In Appendix C,}}$ the aim is to find an
explicit algorithm for finding a one-to-one function from Family(2)
onto {{Family(5)}}(the family of the singular parts of the divisors
defined by the total transforms of irreducible plane curves with
isolated singularity under the standard resolution), called The
$\gamma$-algorithm. $\underline{\text{\rm By Theorem C.2 and Theorem
C.3,}}$ we can find The $\gamma$-algorithm completely and rigorously
in an elementary way. \ms

\noindent{\bf Remark.} Whenever any element of Family(4) is viewed
as a finite sequence of positive integers which is strictly
increasing by either Definition 1.2 of Part[A] or Definition A.2, it
will be proved by Remark A.2.2 that each element of Family(5) can be
viewed a proper subsequence of some element of Family(4), satisfying
an additional property.

\bs
\newpage

$$\align
 &\qquad \qquad\qquad  \text{\bf {Appendix A}}  \\
\quad &\text{\bf The $\alpha$-algorithm for finding a one-to-one}\\
 & \text{\bf function from Family(2) onto Family(3) }\\
\endalign$$ \ms

{\bf \S{A.0.} Introduction}

It was already proved by The 1st Algorithm of PART[A] that we can
find a one-to-one function $\phi$ from
$\underline{\text{\bf{Family(1)}}}$ onto
$\underline{\text{\bf{Family(2)}}}$. In Appendix A, the aim is to
find an algorithm for computing a one-to-one function from
$\underline{\text{\bf{Family(2)}}}$ onto
$\underline{\text{\bf{Family(3)}}}$ in the beginning of Appendix,
called the $\alpha$-algorithm, for notation. \ms

{\bf \S{A.1.} The terminology and notations in preparation for
finding the $\alpha$-algorithm}

In preparation for the explicit and rigorous representation of the
$\alpha$-algorithm, the $\beta$-algorithm and the $\gamma$-algorithm
as we have seen in the beginning of Appendix, recall the definitions
of four Families, that is, Family(k)($1\le k\le 4$) with equivalence
relations, as we have seen in Definition 1.2 and Definition 2.4 of
Part[A]. Moreover, Family(4) with some additional properties will be
said to be a new family in Definition A.2, called Family(5). It will
be proved by Appendix C later that Family(5) is well-defined and
that we can find the $\gamma$-algorithm for finding a one-to-one
function from Family(2) onto Family(5). \ms

\definition{Definition A.1} For brevity, let
$\underline{\text{\bf{Family(0)}}}$ be the $0$-th family, consisting
of all the convergent power series $f \in \BC\{y,z\}$ such that $f$
is irreducible in $\BC\{y,z\}$ with isolated singularity at $0\in
\BC^2$, as in Definition 1.2 and Definition 2.4 of Part[A]. \ms

\noindent{\bf(1)} $\underline{\text{\bf{Family(1)}}}$ is the
1st family, consisting of all the standard Puiseux W-polys $f \in
\BC\{y,z\}$ of the recursive type with isolated singularity at $0\in
\BC^2$, denoted by \ms

\noindent$\text{\rm(A.1.1)}  \quad
\underline{\text{\rm{Family(1)}}=\{\text{\rm $f$ is arbitrary
standard Puiseux W-poly of the recursive r-type:}}$

\quad\quad $\underline{\text{$f\in \text{\rm Family(0) and $r$ are
arbitrary positive integers}$}}\}$, \ms

{\noindent} with an equivalence relation for any two standard
Puiseux W-polys of the recursive type in Family(1), as in Definition
1.2 of Part[A]. \ms

\noindent{\bf(2)} $\underline{\text{\bf{Family(2)}}}$ is the 2nd
family, consisting of all the irreducible plane curve singularities
with the standard Puisuex expansions, denoted by \ms

\noindent(A.1.2)\quad \text{Family(2)=\{$C_r(t)$:$C_r(t)$ is the
standard Puiseux expansion of the $r$-type for any $r\in \N$\}}, \bs

{\noindent} with an equivalence relation for any two standard
Puiseux expansions in Family(2), as in Definition 1.2 of Part[A].
\ms

\noindent{\bf(3)} $\underline{\text{\bf{Family(3)}}}$ is the 3rd
family, consisting of all the multiplicity sequences of irreducible
plane curves with isolated singularity under the standard
resolution, denoted by \ms

\noindent(A.1.3) \qquad
\text{\rm{Family(3)}}=\text{\{Multiseq(V(f)):$f\in \text{\rm
{Family(0)}}$ and f is irreducible in in ${}_2\CO$\}} \ms

{\noindent}with an equivalence relation for any two multiplicity
sequences in Family(3), as in Definition 1.2 of Part[A].

\noindent{\bf(4)} $\underline{\text{\bf{Family(4)}}}$ is the 4-th
family, consisting of all the divisors of $(f\circ\tau)$ defined by
the total transform of $V(f)$, denoted by $(f\circ\tau)_{divisor}$
where $\tau:M\to \BC^2$ is the standard resolution of the
singularity of {$V(f)$} with $f\in \text{\rm Family(0)}$. Then,
$\text{\rm {Family(4)}}$ is denoted by \ms

\noindent(A.1.4) \quad
$\underline{\text{{Family(4)}={\{$(f\circ\tau)_{divisor}$ under
${\tau}$: $f\in \text{\rm Family(0)}$ and}}}$

\qquad \qquad \qquad \qquad  $\underline{\text{$\tau:M\to \BC^2$ is
the standard resolution of the singularity of {$V(f)$}\}}}$.
\qquad\ms

As in Definition 1.2 of Part[A], in more detail, we define
$(f\circ\tau)_{divisor}$, and also an equivalence relation for any
two divisors $(f\circ\tau)_{divisor}$ and $(g\circ\mu)_{divisor}$ in
Family(4), by the following: \ms

\noindent{\bf(4-a)} For any $f\in \text{\rm Family(0)}$, let
$\tau_{\xi}=\pi_1\circ\pi_2\circ\cdots \circ\pi_{\xi}:M^{(\xi)}\to
\BC^2$ with $\tau=\tau_{\xi}$ be the composition of a finite number
$\xi$ of successive blow-ups $\pi_i$ at $0\in \BC^2$, which is
needed only to get the standard resolution of the singularity of
$V(f)$, the zero set of $f$. Let $V(\widetilde {f})$ be the zero set
of $\widetilde{f}$ where $\widetilde{f}=f\circ\tau_{\xi}$. As a
subset of $M^{(\xi)}$, $V(\widetilde {f})$ consists of the
exceptional curve $E=\tau^{-1}_{\xi}(0,0)$ and the proper transform
$V^{(\xi)}(f)$ of $V({f})$. Let
$\tau^{-1}_{\xi}(0,0)=\cup^{\xi}_{i=1}E_{i}$ be the decomposition
into irreducible components, where each $E_{i}$ is called an
exceptional curve of the first kind. Then, $V(\widetilde {f})$
vanishes along each irreducible component $E_i$ of $E$ with a
certain multiplicity $e_i$. Define the total transform of $V(f)$ by
the divisor of $f\circ\tau_{\xi}$, denoted by
$(f\circ\tau_{\xi})_{divisor}$, which is written in the form

\noindent(A.1.5) \qquad \qquad
$(f\circ\tau_{\xi})_{divisor}=V^{(\xi)}(f)+\sum^{\xi}_{i=1}e_iE_i$,
\qquad \qquad

\noindent where each $e_i$ is the multiplicity of $f\circ\tau_{\xi}$
along $E_i$ for $1\le i\le {\xi}$ and $e_{i+1}>e_i$. \ms

\noindent{\bf(4-b)} For any $f$ of (A.1.5) of {\rm (4-a)},
\text{$\{(f\circ\tau_{\xi})_{divisor}\}_{seq.}=\{e_i:i=1,2,\dots,\xi\}$}
is called a sequence whose element consists of a coefficient of each
irreducible component of the exceptional curve
$E=\tau^{-1}_{\xi}(0,0)$. Then it is clear by (A.1.5) that
\text{$(f\circ\tau_{\xi})_{divisor}$} in Family(4) can be viewed as
a finite sequence \text{$\{(f\circ\tau_{\xi})_{divisor}\}_{seq.}$}
of positive integers strictly increasing.

So, $\text{\rm{Family(4)}}$ of (A.1.4) can be equivalently rewritten
by \ms

\noindent(A.1.6) \quad $\underline{\text{\rm{$\text{\rm
Family(4)}_{seq.}$}}=\text{\{$\{(f\circ\tau)_{divisor}\}_{seq.}:f\in
\text{\rm {Family(0)}}$ is arbitrary where $\tau:M\to \BC^2$ is}}$

\qquad \quad $\underline{\text{the standard resolution of the
singularity of {$V(f)$}\}}}$, viewed as a family of sequences. \ms

\noindent{\bf(4-c)} For any $g\in \text{\rm Family(0)}$, using the
same method as in both (4-a) and (4-b), let
$\mu_{\eta}=\bar{\pi}_1\circ\bar{\pi}_2\circ\cdots
\circ\bar{\pi}_{\eta}:\bar{M}^{(\eta)}\to \BC^2$ be the composition
of a finite number $\eta$ of successive blow-ups at the origin in
$\BC^2$, which is needed only to get the standard resolution of the
singularity of $V(g)$. Define the total transform of $V(g)$ by the
divisor of $g\circ\mu_{\eta}$, denoted by
$(g\circ\mu_{\eta})_{divisor}$, which is written in the form
$(g\circ\mu_{\eta})_{divisor}=V^{(\eta)}(g)+\sum^{\eta}_{i=1}\bar{e}_i\bar{E}_i$,
where each $\bar{e}_i$ is the multiplicity of $g\circ\mu_{\eta}$
along $\bar{E}_i$ for $1\le i\le {\eta}$ and
$\bar{e}_{i+1}>\bar{e}_i$. Then, we can define
$\{(g\circ\mu_{\eta})_{divisor}\}_{seq.}$  by a sequence of
$(g\circ\mu_{\eta})_{divisor}$ such that
$\{(g\circ\mu_{\eta})_{divisor}\}_{seq.}=\{{\bar{e}_i}\in
N:i=1,2,\dots,\eta\}$. \ms

\noindent{\bf(4-d)} For any two $f$ of (4-a) and $g$ of (4-c), it is
said that $(f\circ\tau_{\xi})_{divisor}$ and
$(g\circ\mu_{\eta})_{divisor}$ are equivalent(denoted by either ${f}
\buildrel \text{{\rm divisor}} \over \sim {g}$ under the standard
resolutions or
$(f\circ\tau_{\xi})_{divisor}=(g\circ\mu_{\eta})_{divisor}$ under
the standard resolutions) if either of the following condition is
satisfied:
$$\align
(A.1.7) \qquad \qquad \text{\rm either} \quad
&\text{$\{(f\circ\tau_{\xi})_{divisor}\}_{seq.}\equiv
\{(g\circ\mu_{\eta})_{divisor}\}_{seq.}$} \quad \text{\rm as
sequence,} \\
\text{\rm or} \qquad
&\text{$\{e_i:i=1,2,\dots,\xi\}\equiv\{{\bar{e}_i}\in
N:i=1,2,\dots,\eta\}$ \quad as sequence,}\qquad \qquad \\
\text{\rm or} \qquad &\text{\{$e_i={\bar{e}}_i$:
$i=1,2,\dots,\xi=\eta$\}. \quad {$\square$} }
\endalign$$
\enddefinition  \ms

\definition{Definition A.2}
In preparation for construction of the new family, called {\rm
Family(5)}, we can use the same properties and notations as in the
definition of Family(4). Firstly, for convenience of notation, we
define a subset of Family(4), which may be rewritten by
$$\align
(A.2.0) \qquad \qquad  &
\text{Subfamily(4)=\{$(g_r\circ\tau)_{divisor}$: either $g_r\in
\text{\rm Family(1)}$ or $g_r$ is a Puiseux} \qquad \qquad\\
& \text{convergent power series of the recursive r-type in
$\BC\{y,z\}\}
\subseteq \text{ \rm Family(4)}$\},}\\
\endalign$$
where $\tau=\tau_{\xi}=\pi_1\circ\pi_2\circ\cdots
\circ\pi_{\xi}:M^{(\xi)}\to \BC^2$ is the composition of a finite
number $\xi$ of successive blow-ups $\pi_i$ at the origin in
$\BC^2$, which is needed only to get the standard resolution of the
singularity of $V(g_r)$, the zero set of $g_r$, noting by Theorem
1.4 and Theorem 1.6 of Part[A] and by Theorem 7.3 and Theorem 10.2
of Part[B] that \text{\rm{$\text{\rm Family(4)}_{seq.}$}} and
\text{\rm{$\text{\rm Subfamily(4)}_{seq.}$}} are the same in the
sense of (A.1.7) of Definition A.1.

Secondly, we will construct the following three definitions, denoted
by Subdefinition A.2.2, Subdefinition A.2.4 and Subdefinition A.2.5,
which will be proved to be well-defined, later. After then, it
suffices to apply the above three definitions to the construction of
Family(5) with an equivalence relation, as follows: \ms

\noindent{\bf(a)} $\underline{\text{\bf Sublemma A.2.1}}$
$\underline{\text{\rm {(Assumptions)}}}$ By the definition of
\text{\rm{Subfamily(4)}}, let $g_r\in \text{\rm Family(1)}$ be such
that
$(g_r\circ\tau_{\xi})_{divisor}=V^{(\xi)}(g_r)+\sum^{\xi}_{i=1}e_iE_i$
where $\tau_{\xi}:M^{(\xi)}\to \BC^2$ is the composition of a finite
number $\xi$ of successive blow-ups at the origin in $\BC^2$, which
is the standard resolution of the singularity of $V(g_r)$ as we have
seen in (A.1.5) and (A.2.0) and
$\tau^{-1}_{\xi}(0,0)=E=\cup^{\xi}_{i=1}E_{i}$ is the decomposition
into irreducible components.

$\underline{\text{\rm {(Conclusions)}}}$ It was already proved by
Theorem 14.0 of $\S {14}$ of Part[C] that the standard
representation of the divisor of $g_r\circ\tau_{\xi}$, denoted by
$(g_r\circ\tau_{\xi})_{divisor}$, is defined by the following: For
any $A\subset M^{(\xi)}$, $\overline{A}$ is called the closure of
$A$ in $M^{(\xi)}$.
$$\align
(A.2.1) \qquad \qquad
(g_r\circ\tau_{\xi})_{divisor}&=V^{(\xi)}(g_r)+\sum^{\lambda_r}_{i=1}e_iE_i
\quad \text{with} \quad \lambda_r=\xi  \\
 &=V^{(\xi)}(g_r)+\sum^{\lambda_1}_{i=1}e_iE_i
+\sum^{\lambda_2}_{i=\lambda_1+1}e_iE_i+\dots+
\sum^{\lambda_{r}}_{i=\lambda_{r-1}+1}e_iE_i, \qquad \qquad \qquad
\qquad
\endalign$$
where each $e_i$ is the multiplicity of $g_r\circ\tau_{\xi}$ along
$E_i$ for $1\le i\le \xi=\lambda_r$ and $V^{(\xi)}(g_r)$ is the
proper transform of $V(g_r)$ under $\tau_{\xi}$ and write
$\lambda_0=0$ for notation, satisfying two properties:

Write $L=V^{(\xi)}(g_r)\bigcup ({\bigcup}^{\xi}_{i=1} E_i)$, and for
any $A\subset M^{(\xi)}$ $\overline{A}$ is the closure of $A$ in
$M^{(\xi)}$.

\noindent{\rm(i)} For each $i=1,\dots,\xi$, $E_i\bigcap
\overline{(L-E_i)}$ has at most three distinct points under
$\tau_{\xi}$ in $L$.

\noindent{\rm(ii)}{\rm(iia)} There is a strictly increasing finite
sequence $\{\lambda_i: 1\le i\le r\}$ such that
$E_{\lambda_i}\bigcap \overline{(L-E_{\lambda_i})}$ has exactly
three distinct points under $\tau_{\lambda_r}$ in $L$ for each
$i=1,2,\dots,r$.

{\rm(iib)} For any $j\not\in \{\lambda_i: 1\le i\le r\}$ where $1\le
j\le {\xi}$, $E_j\bigcap \overline{(L-E_j)}$ has at most two
distinct points under $\tau_{\lambda_r}$ in $L$.
 \ms

\noindent{\bf(a-1)} $\underline{\text{\bf Subdefinition A.2.2}}$
Suppose that the same properties and notations as in the assumptions
and conclusions of Sublemma $A.2.1$ hold. It is said that
$E_{\lambda_j}$ is called the $j$-th Puiseux exceptional curve of
the first kind. Note that
$1<\lambda_1<\lambda_2<\cdots<\lambda_r=\xi$. It can be proved by
Theorem 7.3 and Theorem 10.2 of Part[B] that a finite number $r$ of
Puiseux exceptional curves of the first kind for any $g_r\in
\text{\rm Family(1)}$ is invariant under an equivalence relation of
Family(4).
\enddefinition

\noindent{\bf(b)} $\underline{\text{\bf Sublemma A.2.3}}$
$\underline{\text{\rm {(Assumptions)}}}$ By Subdefinition A.2.1 with
Subfamily(4), study all the exceptional curves of the first kind of
$E=\bigcup^{\xi}_{i=1}{E_i}$ with $\xi=\lambda_r$ in more detail,
and let $E_{\lambda_j}$ be the $j$-th Puiseux exceptional curve of
the first kind. For convenience of notation, let
${\Omega}^{(1)}=\bigcup^{\lambda_1}_{i=1}{E_{i}}$,
${\Omega}^{(2)}=\bigcup^{\lambda_2}_{i={{\lambda_1}+1}}{E_{i}}$,\dots,
and ${\Omega}^{(r)}=\bigcup^{\lambda_r}_{i={{\lambda_{r-1}}+1}}
{E_{i}}$. Applying the conclusion of Theorem $3.7$ of Part[B] to the
equations in both $(14.2.5)$ and $(14.3.1)$ of $\S {14}$ of Part[C],
it is trivial to prove that the following in the conclusions are
true: For any $A\subset M^{(\xi)}$, $\overline{A}$ is the closure of
$A$ in $M^{(\xi)}$. \ms

$\underline{\text{\rm {(Conclusions)}}}$ Let $w$ be with $0\le w\le
r-1$. Then
${\Omega}^{(w+1)}=\bigcup^{\lambda_{w+1}}_{i=\lambda_{w}+1}E_i$
satisfies two properties : If necessary, we write $\lambda_0=0$.

Let $E_t$ be an arbitrary exceptional curve of the first kind for
$1\le t\le\xi=\lambda_r$.

$\underline{\text{\rm Property(1)}}$ \quad For any $E_t\subset
{\Omega}^{(w+1)}$ with $w+1\le r$, $E_t\bigcap
\overline{{\Omega}^{(w+1)}-E_{t}}$ \quad have at most two distinct
points in ${\Omega}^{(w+1)}$. \ms

$\underline{\text{\rm Property(2)}}$ Let $w$ be with $w+1\le r$, and
${\Omega}^{(w+1)}=\bigcup^{\lambda_{w+1}}_{i=\lambda_{w}+1}E_i$.

{\rm(i)} There are two distinct exceptional curves of the first kind
in ${\Omega}^{(w+1)}$, denoted by $E_{{\lambda_w}+1}$ and
$E_{{\lambda_w}+s_w}$ with $1<s_w\le\lambda_{w+1}-{\lambda_w}$, each
of which satisfies the following property: Note that
${\Omega}^{(w+1)}=\bigcup^{\lambda_{w+1}}_{i=\lambda_{w}+1}E_i$.

{\rm(ia)} $E_{\lambda_w+1}\bigcap
\overline{{\Omega}^{(w+1)}-E_{\lambda_w+1}}$ {\quad} has one and
only one intersection point in ${\Omega}^{(w+1)}$.

{\rm(ib)} $E_{{\lambda_w}+s_w}\bigcap
\overline{{\Omega}^{(w+1)}-E_{{\lambda_w}+s_w}}$ {\quad} has one and
only one intersection point in ${\Omega}^{(w+1)}$.

{\rm(ii)} $E_{{\lambda_w}+j}\bigcap
\overline{{\Omega}^{(w+1)}-E_{{\lambda_w}+j}}$ {\quad} has two
distinct intersection points in ${\Omega}^{(w+1)}$ for any $j$ where
$1<j\le\lambda_{w+1}-{\lambda_w}$ and $j\neq {s_w}$. \ms

\noindent{\bf(b-1)} $\underline{\text{\bf Subdefinition A.2.4.}}$
Suppose that the same properties and notations as in the assumptions
and conclusions of Sublemma $A.2.3$ hold.

(a) If $E_j$ satisfies {\rm(i)} of Property(2) in Sublemma A.2.3,
then  $E_j$ is called a singular exceptional curve of the first kind
for $E=\bigcup^{\xi}_{i=1}{E_i}$.

(b) If $E_j$ satisfies {\rm(ii)} of Property(2), then $E_j$ is
called a regular exceptional curve of the first kind for
$E=\bigcup^{\xi}_{i=1}{E_i}$.

Note that if $r\ge 2$, then $1<s_w\le \lambda_{w+1}-{\lambda_w}$ for
each $w=0,1,\dots,r-1$, and $1<s_0<\lambda_1$. Also, if $r=1$ then
$1<s_0<\lambda_1$. \ms

\noindent$\underline{\text {\bf Remark}}$ Note by Subdefinition
A.2.3 that the Puiseux exceptional curve of the first kind in {\rm
Subdefinition A.2.2} may be a singular exceptional curve of the
first kind in Subdefinition A.2.4, by computing the exceptional
curves for the standard resolution of the singularity of
$\{g(y,z)=(z^2+y^3)^2+y^5z=0\}$ at $0\in \BC^2$. \ms

\noindent{\bf(c)} $\underline{\text{\bf Subdefinition A.2.5
(Definition of singular part of the divisor).}}$ Suppose that the
same properties and notations as in the assumptions and conclusions
of Sublemma $A.2.3$ hold. Using the same properties and notations as
in Sublemma $A.2.3$ and Subdefinition A.2.4,
$(g_r\circ\tau_{\xi})_{\text{\rm singular part of the divisor}}$ is
defined by the following:
$$\align
 & (g_r\circ\tau_{\xi})_{\text{\rm singular part of the divisor}}
= V^{(\xi)}(g_r)+\sum_{1,1} \tag A.2.2\\
& \text{with} \qquad  \text{$\sum_{1,1}=\sum^{r-1}_{i=0}
\{e_{\lambda_i+1}E_{\lambda_i+1}+e_{\lambda_i+s_i}E_{\lambda_i+s_i}\}$},
\qquad\qquad\qquad
\endalign$$
where $\lambda_0=0$ and each $e_j$ is the multiplicity of
$g_r\circ\tau_{\xi}$ along $E_j$ for $1\le j\le \xi$ and
$V^{(\xi)}(g_r)$ is the proper transform of $V(g_r)$ under the
standard resolution $\tau_{\xi}$ of the singularity of $V(g_r)$.

Note that $E_{\lambda_i+1}$ and $E_{\lambda_i+s_i}$ are singular
exceptional curves of the first kind for each $i=1,2,\dots,r$. \ms

\noindent\noindent{\bf(d)} $\underline{\text{\bf{Family(5)}}}$ is
the 5-th family, consisting of the singular part of the divisor
defined by the total transform of $V(g_r)$ as we have seen in
(A.2.2), denoted by $(g_r\circ\tau)_{\text{\rm singular part of the
divisor}}$,
$$\align
(A.2.3) \quad \quad  &
\text{Family(5)=\{$(g_r\circ\tau_\xi)_{\text{\rm singular part of
the divisor}}$: $g_r\in
\text{\rm Family(1)}$ where {$\tau_\xi:M\to \BC^2$} is} \qquad \qquad\\
& \text{the standard resolution of the singularity of
$V(g_r)$, the zero set of $g_r$\}}.\\
\endalign$$

As we have seen in $\{\text{\rm Family(4)}\}_{seq.}$ of (A.1.6),
$(g_r\circ\tau_{\xi})_{\text{\rm singular part of the divisor}}$ can
be equivalently defined by $\underline{\text{\bf a sequence of 2r
elements}}$, denoted by \text{$\{(g_r\circ\tau_{\xi})_{\text{\rm
singular part of the divisor}}\}_{seq.}$},
$$\align
(A.2.4)  \qquad\qquad&
\text{$\{(g_r\circ\tau_{\xi})_{\text{\rm singular part of the divisor}}\}_{seq.}$} \\
=&\{e_{\lambda_0+1},e_{\lambda_0+s_0};e_{\lambda_1+1},e_{\lambda_1+s_1};
\text{$e_{\lambda_2+1},e_{\lambda_2+s_2};\dots;
e_{\lambda_{r-1}+1},e_{\lambda_{r-1}+s_{r-1}}\}$}, \qquad\qquad
\endalign$$
which is strictly increasing where  $\lambda_0=0$ and
$1<s_{j-1}\le{\lambda_j-\lambda_{j-1}}$ \text{for $1\le j\le r$}.
\ms

Also, \text{$\{(g_r\circ\tau_{\xi})_{\text{\rm singular part of the
divisor}}\}_{seq.}$} is called a sequence of
$(g_r\circ\tau)_{\text{\rm singular part of the divisor}}$. So,
viewed as a family of sequences, $\text{\bf {Family(5)}}$ of (A.2.3)
can be identified with $\text{\bf Family(5)}_{seq.},$
$$\align
(A.2.5) \quad \quad &\text{$\underline{\text{\rm{$\text{\bf
Family(5)}_{seq.}$}}=\text{\{$\{(g_r\circ\tau_\xi)_{{\text{\rm
singular part of the divisor}}}\}_{seq.}:g_r\in \text{\rm
{Family(1)}}$ }}$} \\
& \text{$\underline{\text{where $\tau_\xi:M\to \BC^2$ is the
standard resolution of the singularity of {$V(g_r)$}\}}}$.} \qquad
\qquad
\endalign$$

Observe that \text{$\{(g_r\circ\tau_{\xi})_{\text{\rm singular part
of the divisor}}\}_{seq.}$} is a proper subsequence of
\text{$\{(g_r\circ\tau)_{divisor}\}_{seq.}$} where
\text{$\{(g_r\circ\tau_{\xi})_{\text{\rm singular part of the
divisor}}\}_{seq.}$} consists of $2r$ elements. As far as the
computation for algorithms is concerned, it may be proved that the
computation method of \text{$\{(g_r\circ\tau_{\xi})_{\text{\rm
singular part of the divisor}}\}_{seq.}$} is much simpler than that
of \text{$\{(g_r\circ\tau)_{divisor}\}_{seq.}$}. \ms

\noindent{\bf(e)} {\bf In order to define an equivalence relation on
Family(5),} let $\phi_\rho\in \BC\{y,z\}$ be the standard Puiseux
polynomial of the recursive $\rho$-type as we have used in the
Family(1) of Definition 1.2 of Part[A]. By the same method as we
have used in the construction of
\text{$\{(g_r\circ\tau_{\xi})_{\text{\rm singular part of the
divisor}}\}_{seq.}$}, as we have seen in (A.2.2), let
$\mu_{\eta}=\bar{\pi}_1\circ\bar{\pi}_2\circ\cdots
\circ\bar{\pi}_{\eta}:\bar{M}^{(\eta)}\to \BC^2$ be the composition
of a finite number $\eta$ of successive blow-ups at $0\in \BC^2$,
which is needed only to get the standard resolution of the singular
point of $V(\phi_\rho)$.

First, by the same method as we have used in (A.2.1), first we can
define $\{(\phi_\rho\circ\mu_{\eta})_{divisor}\}_{seq.}$  by a
sequence of $(\phi_\rho\circ\mu_{\eta})_{divisor}$ where

$(\phi_\rho\circ\mu_{\eta})_{divisor}=V^{(\eta)}(\phi_\rho)
+\sum^{\eta}_{i=1}\bar{e_i}\bar{E}_i$ and

$\{(\phi_\rho\circ\mu_{\eta})_{divisor}\}_{seq.}=\{{\bar{e}_i}\in
N:i=1,2,\dots,\eta=\psi_{\rho}\}$, satisfying two properties:

{\rm(i)} Each $\bar{E_i}$ of $\psi_{\rho}$ exceptional curves of the
first kind has at most three distinct intersection points with other
exceptional curves of the first kind and the proper transform under
$\mu_{\eta}$.

{\rm(ii)} There exists a sequence $\{\psi_{i}: 1\le i\le \rho\}$
under $\mu_{\eta}$ such that
$\bar{E}_{\psi_1},\bar{E}_{\psi_2},\dots,\bar{E}_{\psi_\rho}$ are
the only $\rho$ exceptional curves of the first kind, each of which
has three distinct intersection points with other exceptional curves
of the first kind and the proper transform.

Next, by the same method as we have used in (A.2.2) and by (A.2.4),
define the following subsequence:
$$\align
(A.2.6)  \qquad & (\phi_\rho\circ\mu_{\eta})_{\text{\rm singular
part of the divisor}} = V^{(\xi)}(\phi_\rho)+\sum^{\rho}_{i=1}
\{e_{\psi_i+1}\bar{E}_{\psi_i+1}+e_{\psi_i+t_i}\bar{E}_{\psi_i+t_i}\}, \\
&\text{$\{(\phi_\rho\circ\mu_{\eta})_{\text{\rm singular part of the divisor}}\}_{seq.}$} \\
=&\{\bar{e}_{\psi_0+1},\bar{e}_{\psi_0+t_0};\bar{e}_{\psi_1+1},\bar{e}_{\psi_1+t_1};
\bar{e}_{\psi_2+1},e_{\psi_2+t_2};\dots;
\bar{e}_{\psi_{\rho-1}+1},\bar{e}_{\psi_{\rho-1}+t_{\rho-1}}\}, \qquad\qquad \\
\endalign$$
where  $\psi_0=0$ and $1<t_{i-1}\le{\psi_i-\psi_{i-1}}$ \text{for
$1\le i\le \rho$}.

As in (A.2.2) of Subdefinition A.2.5, note that $\bar{E}_{\psi_i+1}$
and $\bar{E}_{\psi_i+t_i}$ are singular exceptional curves of the
first kind for each $i=1,2,\dots,{\rho}$. \ms

Then it is said that $(g_r\circ\tau_{\xi})_{\text{\rm singular part
of the divisor}}$ and $(\phi_\rho\circ\mu_{\eta})_{\text{\rm
singular part of the divisor}}$ are equivalent if either (a) or (b)
is satisfied:
$$\align
\noindent (A.2.7) \quad  & \text{(a)
$\{(g_r\circ\tau_{\xi})_{\text{\rm singular part of the
divisor}}\}_{seq.}\}{\equiv} \{(\phi_\rho\circ\mu_{\eta})_{\text{\rm
singular part of the divisor}}\}_{seq.}\}$  as
sequence.} \qquad \qquad\\
& \text{(b)\quad(b1)\quad $\lambda_j=\psi_j$ and $s_j=t_j$ for
$j=0,\dots,r-1={\rho-1}$
with $\lambda_r=\psi_{\rho}$, and } \\
 &\text{\qquad(b2)\quad $e_{\lambda_j+1}=\bar{e}_{\lambda_j+1}$ and
$e_{\lambda_j+s_j}=\bar{e}_{\lambda_j+s_j}$ for $j=0,\dots,r-1$.
\quad {$\square$}}
\endalign$$

\definition {Remark A.2.6}
We may assume that \text{\rm{$\text{\rm Family(4)}_{seq.}$}} and
\text{\rm{$\text{\rm Subfamily(4)}_{seq.}$}} are the same in the
sense of (A.1.7) of Definition A.1. Since any element of Family(4)
is viewed as a finite sequence of positive integers which is
strictly increasing by Definition $A.1$, note that each element of
Family(5) can be viewed a proper subsequence of some element of
Family(4), satisfying an additional property. \quad {$\square$}
\enddefinition \ms

The solution of The $\alpha$-algorithm can be given by Theorem $A.4$
and Theorem $A.5$. In preparation for the representation of The
$\alpha$-algorithm and other algorithms in the beginning of this
appendix, we need the following definition with the new notation.
\ms

\definition{Definition A.3}
The representation of an equivalence of any two elements of both
{Family(1)} and {Family(2)}, and also the representation of an
equivalence of any two elements of both {Family(1)} and {Family(4)}
in terms of multiplicity sequences of irreducible plane curve
singularities can be defined as follows. \ms

\noindent$\underline{\text{\bf Definition A.3.0(The Euclidean
multiplicity sequence for two positive integers).}}$ Definition
A.3.0 has the same statement as Definition 9.2 of Part[B] does in
$\S 9$ of Part[B]. \ms

\noindent$\underline{\text{\bf Definition A.3.1}}$ {\bf[I]} Let
$S=\{t_i \in {\BC}:i=1,2,\dots,\lambda\}$ be a finite sequence of
complex numbers with $\lambda>0$. It is said that {\rm S} is the
join of a finite number $r$ of subsequences $B_i$ of $S$ in order,
denoted by $S=\text{\rm Join}(B_1,B_2,\dots,B_r)$, if the following
are satisfied:
$$\align
(A.3.1) \quad  &\text{$B_1=\{t_i: i=1,2,\dots,q_1 \}$,
\quad $B_2=\{t_{q_{1}+i}:i=1,2,\dots,q_2\}$,\dots,} \\
&\text{$B_r=\{t_{q_{r-1}+i}:i=1,2,\dots,q_r\}$} \\
&\text{where each $B_i$ is a subsequence of $S$ such that $q_i>0$
and $q_1+q_2+\cdots +q_r=\lambda$.} \qquad
\endalign$$
Namely, whenever $t_{\alpha}\in B_i$ and $t_{\beta}\in B_j$ for any
two $B_i$ and $B_j$ with $1\le i<j\le r$, then $1\le \alpha<\beta\le
\lambda$. \ms

\noindent{\bf [II]} Let $g_r\in \text{\rm Family(1)}$ be given
arbitrary. By the same properties and notations as in (A.2.1) of
Sublemma A.2.1 of Definition A.2 and by (A.2.0) of Definition A.2,
$\{(g_r\circ\tau_{\lambda_r})_{divisor}\}_{seq.}\in \text{\rm
Family(2)}$ is well-defined by
$$\align
(A.3.2) \qquad
\text{$\{(g_r\circ\tau_{\lambda_r})_{divisor}\}_{seq.}
\equiv\{e_i:i=1,2,\dots,\lambda_r\}\equiv\text{\rm
Join}(B_1,B_2,\dots,B_r)$, as sequence} \qquad
\endalign$$
for $j=1,2,\dots,r$, the j-th subsequence $B_j$ of which is written
respectively as follows:
$$\align
(A.3.3) \qquad \qquad & B_1=\{e_i: i=1,2,\dots,\lambda_1 \}  \quad
\text{with $1<\lambda_1<\lambda_2<\cdots<\lambda_r$,}   \\
& B_j=\{e_{\lambda_{j-1}+i}:i=1,2,\dots,(\lambda_j-\lambda_{j-1})\}
\quad \text{for $j=2,3,\dots,r$.} \qquad \qquad \qquad \qquad
\endalign$$

\noindent$\underline{\text{\bf Definition A.3.2}}$ {\bf[I]} Without
any need of computation, we may assume by {\rm Theorem 1.4} and {\rm
Theorem 1.6} of Part[A] that we can find the standard Puiseux
polynomial $g_r\in \BC\{y,z\}$ of the recursive $r$-type such that
\text{$V(g_r)\equiv C_r(t)$ \text{\rm (Multiseq)}} for given the
standard Puiseux expansion $C_r(t)$ of the $r$-type of the curve,
and conversely. \ms

\noindent{\bf [II]} For notation,  we may assume without loss of
generality that $g_r\in \text{\rm Family(1)}$ satisfies the same
notations and properties as in Definition 1.1 of Part[A], and also
$C_r(t)\in \text{\rm Family(2)}$ satisfies the same notations and
properties as in Definition 1.2 of Part[A]. Then, it was already
proved by Theorem 1.4 and Theorem 1.6 of Part[A] that
\text{$V(g_r)\equiv C_r(t)$ \text{\rm (Multiseq )}} if and only if
(i) and (ii) of (A.3.4) are true, and in addition (iii) is just the
computational formula for finding $e_{\lambda_{j}}$ in this case.
Thus, the following may be viewed as a definition for
\text{$V(g_r)\equiv C_r(t)$ \text{\rm (Multiseq)}}:
$$\align
(A.3.4) \qquad\qquad\qquad &\text{$V(g_r)\equiv C_r(t)$ \text{\rm (Multiseq)}} \\
\iff \quad \text{\rm(i)} \quad &\text{$n=n_1d_1$ and
$\alpha_1=\beta_{1,1}d_1$ with $d_1=\gcd(n,\alpha_1)$ and,} \\
\text{\rm{(ii)}} \quad &\text{$d_{j-1}=n_jd_{j}$ and
$\alpha_j-\alpha_{j-1}= \widehat{\Delta}_jd_j$ with
$d_j=\gcd(d_{j-1},\alpha_j-\alpha_{j-1})$}, \qquad\qquad\\
&\text{where $\widehat{\Delta}_j=\Delta_j(\beta_{j,k})^j_{k=1}
-n_jn_{j-1}\Delta_{j-1}(\beta_{j-1,k})^{j-1}_{k=1}$ for $2\le j\le
r$.}\\
\endalign$$

\noindent{\bf Remark A.3.2.1.} \text{\rm{\quad(iii)}}
$e_{\lambda_1}=d_1n_1\beta_{1,1}$ \quad and \quad
\text{$e_{\lambda_j}
=e_{\lambda_{j-1}}n_j+d_jn_j\text{$\widehat{\Delta}$}_j$ \quad for
$2\le j\le r$.} \ms

\noindent$\underline{\text{\bf Definition A.3.3}}$ For convenience
of notation, we may assume that $g_r\in \text{\rm Family(1)}$
satisfies the same properties and notations as in Definition A.3.1.
Note that $\{\lambda_j: 1<\lambda_1<\lambda_2<\cdots<\lambda_r, 1\le
j\le r \}$ is a finite integer sequence in Subdefinition A.2.2. By
the same way as we have used in (A.3.2) of Definition A.3.1,
$\text{\rm {Multiseq(V($g_r$))}}$ in \text{\rm Family(3)} is
well-defined by
$$\align
(A.3.5)\qquad  \text{\rm
Multiseq(V($g_r$))}\equiv\{c_i:i=1,2,\dots,\lambda_r\}
\equiv\text{\rm Join}\{P_1,P_2,\dots,P_r\}, \quad\text{as sequence}
\qquad
\endalign$$
for $j=1,2,\dots,r$, each subsequence $P_j$ of which can be uniquely
written as follows:
$$\align
\text{\rm (A.3.6)}  \qquad \qquad  & P_1=\{c_i:
i=1,2,\dots,\lambda_1 \}
 \quad \text{and}  \qquad \qquad\qquad \qquad \\
& P_j=\{c_{\lambda_{j-1}+i}:i=1,2,\dots,(\lambda_j-\lambda_{j-1})\}
\quad \text{for $j=2,3,\dots,r$. \quad {$\square$}} \qquad\qquad
\qquad
\endalign$$
\enddefinition \ms

{\bf {\S} A.2 and {\S} A.3. The $\alpha$-algorithm(The algorithm for
finding a one-to-one function between Family(2) and Family(3)) and
its examples} \bs

{\bf {\S} A.2. The 1st half of the $\alpha$-algorithm(Theorem A.4)}
\ms

\proclaim{Theorem A.4(The 1st half of the $\alpha$-algorithm: an
algorithm for finding a one-to-one function from Family(2) into
Family(3))}

$\underline{\text{\bf Assumptions}}$ Let the standard Puiseux
expansion $C_r(t)$ of the $r$-type for the curve $C$ be given as
follows:
$$\align
(A.4.1) \qquad \qquad \text{$C_r(t):=$} \left\{\eqalign{ y=&t^n, \cr
z=&t^{\alpha_1}+t^{\alpha_2}+\cdots +t^{\alpha_r}, \cr} \right. \\
\text{where} \quad
 2\le n <\alpha_1<\alpha_2<\cdots <\alpha_r  & \quad  \text{and} \\
 n >d_1>d_2>\cdots
 >d_r=1  \quad   \text{with} & \quad
\text{$d_i=\gcd(n,\alpha_1,\dots,\alpha_i)$, $1\le i\le r$.} \qquad
\endalign$$

$\underline{\text{\bf Conclusions}}$ \ms

{\rm(1)}{\rm(1a)} For a good representation, use {\rm Definition
A.3.3} with {\rm(A.3.5)} and {\rm Definition A.3.0}. Then, \text{\rm
Multiseq$(C_r(t))$}, called the multiplicity sequence of $C_r(t)$,
can be clearly represented by the following algorithm: \ms

\noindent$\underline{\text{\bf(The 1st half of the
$\alpha$-algorithm(Theorem A.4)}}$
$$\align
(A.4.2)^* \qquad \qquad \qquad   \text{\rm
Multiseq$(C_r(t))$}&=\text{\rm Join}(P_1,
P_2,\dots,P_r), \qquad \qquad \\
\text{where $P_1=\{[n:\alpha_1]\}$},
&\text{$P_2=\{[d_1:\alpha_2-\alpha_1]\}$,
\dots,$P_r=\{[d_{r-1}:\alpha_{r}-\alpha_{r-1}]\}$}. \qquad \qquad \qquad\\
\endalign$$

\noindent{\bf{Note}.} As we have seen in either {\rm (9.2.3)} or
{\rm (9.4.3)} of {\rm Part[B]}, we use the same kind of notations as
in $(A.4.2)^*$:

Assuming that $\gcd(n_1,k_1)=1$ and
$\{[n_1,k_1]\}=\{c_1,c_2,\dots,c_t\}$ by {\rm Definition A.3.0}, for
brevity of notation we may define
$\{[dn_1,dk_1]\}=\{d[n_1,k_1]\}=\{dc_1,dc_2,\dots,dc_t\}$ for any
positive integer $d$. \ms

{\rm(1b)} The representation of \text{\rm Multiseq({${C_r}$}(t))} in
{\rm(A.4.2)} is unique. That is, if \text{\rm Multiseq($C_r$(t))}
has another representation with $\text{\rm
Multiseq($C_r$(t))}=\text{\rm Join}(Q_1, Q_2,\dots,Q_s)$ for an
integer $s>0$ in the sense of {\rm Definition A.3.3}, then
$P_i{\equiv} Q_i$ as sequence for $1\le i\le r=s$ by {\rm Definition
A.2}. \ms

{\rm(2)} Assuming for notation that any standard Puiseux expansion
has the same multiplicity sequence as we have seen in {\rm(A.4.2)},
it is trivial to compute the corresponding standard Puiseux
expansion with the same multiplicity sequence. \quad {$\square$}
\endproclaim \ms

\definition{Remark A.4.1}
{\rm(i)} For any $C_r(t)\in \text{\rm Family}(2)$, we may assume
that \text{$V(g_r)\equiv C_r(t)$ \text{\rm (multi. seq.)}} for some
$g_r\in \text{\rm Family}(1)$. Using the properties and notations of
Definition $A.3$ for $(g_r\circ\tau_{\xi})_{divisor}\in \text{\rm
Family}(4)$, then  the number of elements in
$\{[d_{j-1}:\alpha_j-\alpha_{j-1}]\}$, as a sequence, is equal to
$\lambda_j-\lambda_{j-1}$ for $2\le j\le r$, which must be uniquely
determined.

{\rm(ii)} The proof of Theorem A.4 is clear.
\enddefinition \ms

\noindent{\bf Example A.4.2 for Theorem A.4:} As in Example 1.6.3
for Theorem 1.6 of Part[A], let the standard Puiseux expansion
$C_4(t)$ of the 4-th type for the curve be given by
$$
\text{$C_4(t):=$}  \left\{\eqalign{ y= &t^{100} \cr z=
&t^{250}+t^{375}+t^{410}+t^{417}. \cr } \right. \tag A.4.3
$$

By {\bf(The 1st half of the $\alpha$-algorithm(Theorem A.4)},
$\text{\rm Multiseq($C_4(t)$)}$ can be written as follows:
$$\align
(A.4.4) \qquad &\text{\rm Multiseq($C_4$(t))} =\text{\rm
Join}(P_1,P_2,
P_3,P_4),  \\
\text{where} \quad  &\text{$P_1=\{[100:250]\}$, $P_2=\{[50:125]\}$,
$P_3=\{[25:35]\}$ and $P_4=\{[5:7]\}$.} \qquad \qquad
\endalign$$
because of the following elementary computations:

{\rm(i)} $n=100$ and $\alpha_1=250\Longrightarrow
P_1=\{[n:\alpha_1]\}=\{[100:250]\}$. \ms

{\rm(ii)} $P_2=\{[d_1:\alpha_2-\alpha_1]\}$ with
$d_1=\gcd(n,\alpha_1)=\gcd(100,250)=50$ and $\alpha_2-\alpha_1=125$.
\ms

{\rm(iii)} $P_3=\{[d_2:\alpha_3-\alpha_2]\}$ with
$d_2=\gcd(d_1,\alpha_2-\alpha_1)=\gcd(50,125)=25$ and
$\alpha_3-\alpha_2=35$. \ms

{\rm(iv)} $P_4=\{[d_3:\alpha_4-\alpha_3]\}$ with
$d_3=\gcd(d_2,\alpha_3-\alpha_2)=\gcd(25,35)=5$ and
$\alpha_4-\alpha_3=7$. \bs

{\bf \S{A.3}  The 2nd half of The $\alpha$-algorithm(Theorem
A.5)}\ms

The 2nd half of the $\alpha$-algorithm is an algorithm for finding a
function from Family(2) onto Family(3). We will compute The 2nd half
of the $\alpha$-algorithm in more detail.

\proclaim{Theorem A.5(The 2nd half of The $\alpha$-algorithm: an
algorithm for finding the standard Puiseux expansion from either any
multiplicity sequence of irreducible plane curve singularities or
any finite sequence of positive integers)}

$\underline{\text{\bf Assumptions}}$ Assume that $E$ is a finite
sequence of positive integers as follows:
$$E=\{a_1,a_2,\dots,a_t\}, \tag A.5.1$$
where $a_1\ge a_2\ge \cdots\ge a_{t-1}\ge a_t=1$ are all positive
integers. \ms

{\rm(i)} $\underline{\text{\rm The 1st problem}}$ is to determine
whether or not $E$ is equal to the multiplicity sequence of an
irreducible plane curve with isolated singularity.

{\rm(ii)} $\underline{\text{\rm The 2nd problem}}$ is to find the
standard Puiseux expansion $C_r(t)$ of the r-type for the
irreducible plane curve such that $\text{\rm
Multiseq$(C_r(t))$}\equiv E$ as sequence if the solution for The 1st
problem is positive. \ms

To solve $\underline{\text{\rm The 1st problem}}$  and
$\underline{\text{\rm The 2nd problem}}$ at the same time, we may
assume by the same kind of notations as in $(A.4.2)^*$ that
$\text{\rm Multiseq$(C(t))$}$ can be represented as follows:
$$\align
(A.5.2) \qquad \qquad  \text{\rm Multiseq(C(t))} &=\text{\rm
Join}(P_1,P_2,\dots,P_{r}),  \\
with \quad  P_1=&\{[n:\alpha_1]\}  \quad \text{and \quad
$P_i=\{[d_{i-1}:\alpha_i-\alpha_{i-1}]\}$
for $2\le i\le r$,} \qquad \qquad\\
\endalign$$
where the standard Puiseux expansion of $C(t)$ is given by $y=t^n$
and $z=t^{\alpha_1}+t^{\alpha_2}+\cdots +t^{\alpha_{r}}$ for some
integer $r\ge 1$, with the following properties:

{\rm(a)} $2\le n<\alpha_1<\alpha_2<\cdots <\alpha_{r}$.

{\rm (b)} $n>d_1>\cdots
>d_{r}=1$,
and $d_{i}=\gcd(n,\alpha_1,\dots,\alpha_{i})$ for $1\le i\le r$. \ms

$\underline{\text{\bf Conclusions}}$ \quad For each $k=1,2,\dots,r$,
the aim is to find an elementary computational algorithm for the
construction of the $k-th$ Euclidean multiplicity subsequence
$P_k=\{[d_{k-1}:\alpha_k-\alpha_{k-1}]\}$ of $E$ for $1\le k\le
r${\bf (The 2nd half of The $\alpha$-algorithm)}. For notation write
$P_1=\{[n:\alpha_1]\}$ with $d_0=n$ and $\alpha_0=0$, if necessary.
\ms

\noindent $\underline{\text{\bf The first step for The 2nd half of
The $\alpha$-algorithm}}$ An elementary computational algorithm
formula for $P_1=\{[n:\alpha_1]\}$ with $d_1$ can be represented as
follows: \ms \noindent{\rm(A.5.3)} \qquad
$$\vbox
 {\offinterlineskip\tabskip=0pt\halign{\strut \vrule~ # \vrule \cr
\noalign{\hrule}
 {}\hfill  \cr
 \noalign{\vskip -8pt}
 \hskip 30pt~~$n=a_1=$\text{~the largest element in
$E$}\quad \text{and} \hfill \cr \hskip 30pt ~~ $\tau_{0}=\text{the
finite number of the element $n$ in $E$}$, \hfill \cr \hskip
30pt~~$\alpha_1=\tau_0n+ \max{\{a_i\in E :a_i<n\}} \quad
\text{with}$\hfill \cr
 \hskip 30pt~~$d_1=\gcd(n,\alpha_1)\in E$,\hfill \cr
  \noalign{\vskip -8pt}
 {}\hfill \cr
\noalign{\hrule}}}$$

where note by definition that $\max{\{a_i\in E :a_i<n\}}$ is equal
to the second largest positive integer in $E$. Then, observe that
$n$ and $\alpha_1$ were already computed by {\rm(A.5.3)}. \ms

Now, it remains only to determine whether or not the problem is
completely solved.

{\rm(i)} If $d_1=\gcd(n,\alpha_1)=1$, then the standard Puiseux
expansion of $C(t)$ is completely solved by $y=t^n$ and
$z=t^{\alpha_1}$.

{\rm (ii)} If $d_1=\gcd(n,\alpha_1)>1$, then take the next step. \ms

Assuming that $d_1>1$,  for $k=2,3,\dots,r$, an elementary
computational algorithm formula with {\rm($\alpha_k-\alpha_{k-1}$)}
and $d_k$ can be represented as follows:

\noindent $\underline{\text{\bf The k-th step for The 2nd half of
The $\alpha$-algorithm}}$ Let
$d_{k-1}=\gcd(n,\alpha_1,\dots,\alpha_{k-1})>1$ with $k\ge 2$, and
let $\alpha_{k-1}-\alpha_{k-2}$ be given with $\alpha_0=0$. Then, an
elementary computational algorithm formula for
$P_k=\{[d_{k-1}:\alpha_k-\alpha_{k-1}]\}$ with $d_k$ can be
represented as follows: \ms

\noindent {\rm(A.5.4)} \qquad
$$\vbox{\offinterlineskip\tabskip=0pt\halign{\strut \vrule~ #
\vrule\cr \noalign{\hrule}
 {}\hfill \cr
 \noalign{\vskip -8pt}
\hskip 20pt \quad  $\sigma_{k-1}=\dfrac{\min{\{a_i\in E
:a_i>d_{k-1}\}}}{d_{k-1}}$ \quad \text{and}, \hfill \cr \hskip 20pt
\quad $\tau_{k-1}=\text{  the finite number of the element $d_{k-1}$
in $E$}$,\hfill \cr \hskip 30pt~~ $\alpha_k-\alpha_{k-1}=
(\tau_{k-1}-\sigma_{k-1})d_{k-1}+ \max{\{a_i\in E :a_i<d_{k-1}\}}
\quad \text{with}$\hfill \cr \hskip 30pt \quad
$d_k=\gcd(d_{k-1},\alpha_k-\alpha_{k-1})\in E$. \hfill \cr
\noalign{\vskip -8pt}
 {}\hfill \cr
\noalign{\hrule}}}$$

Then, we can compute the following:
$\alpha_k=\alpha_{k-1}+(\alpha_k-\alpha_{k-1})$. \ms

Now, it remains only to determine whether or not the problem is
completely solved.

{\rm(i)} If $d_k=\gcd(d_{k-1},\alpha_k-\alpha_{k-1})=1$, then the
standard Puiseux expansion of $C(t)$ is completely solved by $y=t^n$
and $z=t^{\alpha_1}+t^{\alpha_2}+\cdots +t^{\alpha_{k}}$.

{\rm(ii)} If $d_k=\gcd(d_{k-1},\alpha_k-\alpha_{k-1})>1$, then take
the next step. \quad $\square$
\endproclaim \ms

\definition{Remark A.5.1} Note that $\sigma_{k-1}=\dfrac{\min{\{a_i\in E
:a_i>d_{k-1}\}}}{d_{k-1}}=\dfrac{\text{ \rm the 2nd smallest element
in $P_{k-1}$}}{d_{k-1}}$ is the finite number of the element
$d_{k-1}$ in $P_{k-1}$ because $d_{k-1}=\min{P_{k-1}}=\min\{a\in
P_{k-1}: \text{\rm a is arbitrary}\}$ for each fixed $k\ge 2$, and
so $\tau_{k-1}-\sigma_{k-1}\ge 0$.
\enddefinition \ms

\noindent{\bf Example A.6 and Example A.7 for The 2nd half of The
$\alpha$- algorithm(Theorem A.5):} Here are two finite sequences
$G_i$ of positive integers for $i=1,2$, as it has been seen in
(A.5.1). By The 2nd half of The $\alpha$- algorithm, we show that
the following are true:

{\rm Example A.6:} We prove that (i) $G_1$ is a multiplicity
sequence, which is equal to that of some irreducible plane curve
singularity defined by the standard Puiseux expansion $C_r(t)$ of
the r-type and (ii) compute $C_r(t)$. \ms

{\rm Example A.7:} We can prove that $G_2$ cannot be a multiplicity
sequence defining any irreducible plane curve singularity. \ms

{\bf Example A.6:}  If the finite sequence $G_1$ is given below,
find the standard Puiseux expansion for an irreducible curve $C(t)$
which has the same multiplicity sequence as in $G_1$, if exists:
$$\align
G_1&=\{a_i: i=1,2,\dots,31\} \tag A.6.1 \\
&=\{12600,12600;7200;5400;1800,1800,1800,1800;   \\
&\qquad500,500,500;300;200;100,100,100,100; \\
& \qquad 55;45;10,10,10,10;5,5,5,5;3;2;1,1\}. \\
\endalign$$ \ms

To get the solution, suppose that $G_1$ is equal to the multiplicity
sequence which is defined by the standard Puiseux expansion for the
curve $C(t)$, and then take the following steps.

$\underline{\text{The first step}}$ \quad  Find the construction for
$P_1=\{[n:\alpha_1]\}$, the Euclidean multiplicity sequence for two
given positive integers $n$ and $\alpha_1$. Note that
$\min\{P_1\}=d_1=\gcd(n,\alpha_1)$.

(I) By (A.5.3), an easy computation says that \noindent
$$\vbox{\offinterlineskip\tabskip=0pt\halign{\strut
\vrule~ # \vrule  \cr \noalign{\hrule}
 {}\hfill  \cr
 \noalign{\vskip -8pt}
\hskip 30pt ~~ $n=12600=$\text{~the largest element in $G_1$}\quad
\text{and} \hfill \cr \hskip 30pt ~~ $\tau_{0}=\text{ the finite
number of the element $n$ in $G_1$=2}$,\hfill \cr \hskip 30pt ~~
$\alpha_1=\tau_0n+ \max{\{a_i\in G_1 :a_i<n\}}=2\cdot
12600+7200=32400 \quad \text{with}$\hfill \cr \hskip 30pt ~~
$d_1=\gcd(n,\alpha_1)=\gcd(12600,32400)=\gcd(12600,7200)=1800 \in
G_1$. \hfill \cr \noalign{\vskip -8pt}
 {}\hfill \cr
\noalign{\hrule}}}$$

By Definition $A.3.0$, an easy computation says that $d_1>1$ and
$$\align
(A.6.2)\quad {P_1}&=\{[n:\alpha_1]\}= \{[12600:32400]\}
=\{12600,12600;7200;5400;1800,1800,1800\}, \\
G_1&=\text{{\rm Join}$(P_1,Q)$ for some subsequence $Q$ of $G_1$.}
\endalign$$
So, we can take the next step. \ms

$\underline{\text{The second step}}$ \quad  Find the construction
for $P_2=\{[d_1:\alpha_2-\alpha_1]\}$, as the second Euclidean
multiplicity subsequence of $G_1$. Note that
$\min{P_2}=d_2=\gcd(d_1,\alpha_2-\alpha_1)$.

(I) By (A.5.4), an easy computation says that

$$\vbox{\offinterlineskip\tabskip=0pt\halign{\strut
\vrule~ # \vrule\cr \noalign{\hrule}
 {}\hfill \cr
 \noalign{\vskip -8pt}
\hskip 20pt \quad $\sigma_{1}=\dfrac{\min{\{a_i\in G_1
:a_i>d_{1}\}}}{d_{1}}=\dfrac{5400}{1800}=3$ \quad \text{and} \hfill
\cr \hskip 20pt \quad  $\tau_{1}=\text{  the counting number of the
element $d_{1}$ in $G_1$}=4$,\hfill \cr \hskip 30pt ~~
$\alpha_2-\alpha_{1}= (\tau_{1}-\sigma_{1})d_{1}+ \max{\{a_i\in G_1
:a_i<d_{1}\}}=(4-3)1800+500=2300 \quad \text{with}$\hfill \cr \hskip
20pt \quad
$d_2=\gcd(d_1,\alpha_2-\alpha_1)=\gcd(1800,2300)=\gcd(1800,500)=100\in
G_1$.\hfill \cr \noalign{\vskip -8pt}
 {}\hfill \cr
\noalign{\hrule}}}$$ \ms

Thus, $\alpha_2=\alpha_1+(\alpha_2-\alpha_1)=32400+2300=34700$ with
$d_2=100>1$. \ms

(II) By (I), note that $\alpha_2-\alpha_1=2300$ and $d_1=1800$ with
$d_2=100>1$.

For the necessity of the next step with $d_2>1$, by Definition
$A.3.0$, it is clear that
$$\align
{P_2}&=\{[d_1:\alpha_2-\alpha_1]\}= \{[1800:2300]\}  \tag A.6.3 \\
&=\{1800;500,500,500;300;200;100,100\}.
\endalign$$

Then, $G_1=\text{\rm{Join}}(P_1,P_2,Q)$ for a subsequence $Q$ of
$G_1$, and so we can take the next step. \ms

$\underline{\text{The third step}}$ \quad  Find the construction for
${P_3}=\{[d_2:\alpha_3-\alpha_2]\}$, as the third Euclidean
multiplicity subsequence of $G_1$. Note that
$\min{P_3}=d_3=\gcd(d_2,\alpha_3-\alpha_2)$.

(I) By the same method as in (A.5.4), we can compute the following:

$\alpha_3=\alpha_2+(\alpha_3-\alpha_2)=34700+255=34955$ with
$d_3=5>1$. \ms

(II) By (I), note that $\alpha_3-\alpha_2=255$ and $d_2=100$ with
$d_3=5>1$.

For the necessity of the next step with $d_3>1$, by Definition
$A.3.0$, it is clear that
$$\align
(A.6.4) \qquad{P_3}&=\{[d_2:\alpha_3-\alpha_2]\}= \{[100:255]\}
=\{100,100; 55;45;10,10,10,10;5,5\}.
\endalign$$

Then, $G_1=\text{\rm{Join}}(P_1,P_2,P_3,Q)$ for some subsequence $Q$
of $G_1$, and so we can take the next step.  \ms

$\underline{\text{The fourth step}}$ \quad  Find the construction
for ${P_4}=\{[d_3:\alpha_4-\alpha_3]\}$, as the fourth Euclidean
multiplicity subsequence of $G_1$. Note that
$\min{P_4}=d_4=\gcd(d_3,\alpha_4-\alpha_3)$.

(I) \quad By the same method as in (A.5.4), we can compute the
following:

$\alpha_4=\alpha_3+(\alpha_4-\alpha_3)=34955+13=34968$ with $d_4=1$.
\ms

(II) By (I), note that $\alpha_4-\alpha_3=13$ and $d_3=5$ with
$d_4=1$.

Since $d_4=1$, then by Definition $A.3.0$, an easy computation says
that
$$\align
(A.6.5) \qquad  {P_4}&=\{[d_3:\alpha_4-\alpha_3]\}= \{[5:13]\}
=\{5,5;3;2,2;1,1\}~ \text{with $P_4=Q$, and so}  \qquad \qquad \\
\quad  G_1&=\text{\rm{Join}}(P_1,P_2,P_3,P_4) \quad \text{by
Theorem A.5.} \\
\endalign$$

Summarizing four steps, we proved that $G_1$ is a multiplicity
sequence being equivalent to some irreducible plane curve $C(t)$
with $G_1=\text{\rm{Join}}(P_1,P_2,P_3,P_4)$ where $n=12600$,
$\alpha_1=32400$, $\alpha_2=34700$, $\alpha_3=34955$ and
$\alpha_4=34968$.

So, the standard Puiseux expansion $C_4(t)$ of the 4-type for the
curve $C(t)$ is given as follows:
$$
\text{$C_4(t):=$} \left\{\eqalign{ y= &t^{12600} \cr z=
&t^{32400}+t^{34700}+t^{34955}+t^{34968}. \cr} \right. \tag A.6.6
$$ \ms

{\bf Example A.7:}  If the finite sequence $G_2$ is given below,
find the standard Puiseux expansion for an irreducible plane curve
$C(t)$ whose multiplicity sequence is equal to $G_2$, if exists:
$$\align
G_2&=\{b_i: i=1,2,\dots,30\}  \tag A.7.1 \\
&=\{12600,12600;7200;5400;1800,1800,1800,1800; \\
&\qquad500,500;300;200;100,100,100,100; \\
& \qquad 55;45;10,10,10,10;5,5,5,5;3;2;1,1\}. \\
\endalign$$

In order to find whether $G_2$ is equivalent to a multiplicity
sequence defining an irreducible plane curve singularity or not,
apply the same method and notations, as we have used in {\rm Example
A.6}, to {\rm Example A.7}. Then, we can compute that $G_2$ cannot
be a multiplicity sequence for any irreducible plane curve
singularity. \ms

{\bf \S{A.4.}  The proof for Theorem A.5(The 2nd half of the
$\alpha$-algorithm)} \ms

\demo{\bf Proof of Theorem A.5} It is clear that Theorem A.4 is
true. For the proof, we use the same properties and notations as in
{\rm Theorem A.4} of ${\S} A.2$ and {\rm Theorem A.5} of ${\S} A.3$.
Moreover, we may start to assume by Definition $A.3.0$, and Theorem
$11.2$ and Corollary $11.3$ of Part[B] that $E$ in the assumption of
Theorem A.5 satisfies two properties {\rm (i)} and {\rm (ii)}:

{\rm(i)} As a multiplicity sequence, $E$ can be written as follows:
$$\align
(A.5.5) \quad \quad E &=\{\mu_{1,1},\mu_{1,2},\dots,\mu_{1,q_1};
\mu_{2,1},\mu_{2,2},\dots,\mu_{2,q_2};\dots;  \\
& \qquad \mu_{w-1,1},\mu_{w-1,2},\dots,\mu_{w-1,q_{w-1}};
\mu_{w,1},\mu_{w,2},\dots,\mu_{w,q_w}\}  \\
& =\{Q_1;Q_2;\cdots;Q_{w-1};Q_w\}
\quad \text{with \quad $Q_i=P_i$ \quad for $1\le i\le w=r+1$,} \qquad \qquad\\
\endalign$$
where each subsequence
$\{Q_i\}=\{\mu_{i,1},\mu_{i,2},\dots,\mu_{i,q_i} \}$ is called the
$i-th$ Euclidean divisor subsequence of $E$, consisting of $q_i$
elements for $i=1,2,\dots,w$, and $q_1+q_2+\dots+q_w=t$. \ms

{\rm(ii)} We have the following:
$$\align
&\text{$\mu_{1,j}=a_j$ \quad  for $1\le j\le q_1$.}  \tag A.5.6 \\
&\text{$\mu_{2,j}=a_{{q_1}+j}$ \quad  for $1\le j\le q_2$.} \\
& \dots\dots \\
&\text{$\mu_{w,j}=a_{{q_1}+{q_2}+\dots +{q_{w-1}}+j}$ \quad  for
$1\le j\le q_w$.}
\endalign$$
\ms

Now, by induction on the positive integer, we are going to find an
elementary computational algorithm for the solution of the main
problem. First of all, note that $n=a_1=\mu_{1,1}$.

$\underline{\text{\bf [I] The 1st step:}}$ \ Let $n=a_1$. In this
step, the aim is to find an elementary computational algorithm for
$\alpha_1$, in construction of the Euclidean divisor sequence
$\{P_1\}=\{[n:\alpha_1]\}$.

Then, it is enough to consider the following in order:

\roster

\item "(a)" Let $\tau_0$ be the counting number of $n$ in $E$.
First, find $\tau_0$.

\item "(b)" Next, compute $\alpha_1$ by (a) and
by Definition $A.3.0$.

\item "(c)" If $\gcd(n,\alpha_1)=1$, then the algorithm is
completely finished.

\noindent If $\gcd(n,\alpha_1)>1$, then take the 2nd step.
\endroster

After we solve (a) and (b), then we have the following algorithm for
the 1st step.

\noindent $\underline{\text{\bf The algorithm for the 1st step:}}$
An elementary computational algorithm formula for $P_1=[n:\alpha_1]$
with $\alpha_1$ can be represented as follows:

\noindent {\rm(A.5.7)}
$$ \vbox
{\offinterlineskip\tabskip=0pt\halign{\strut \vrule~ # \vrule  \cr
\noalign{\hrule} \quad {}\hfill  \cr
 \noalign{\vskip -8pt}
 \hskip 30pt~~$n=a_1=$\text{~the largest element in
$E$}\quad \text{and} \hfill \cr \hskip 30pt ~~ $\tau_{0}=\text{ the
counting number of the element $n$ in $E$}$, \hfill \cr \hskip
30pt~~$\alpha_1=\tau_0n+ \max{\{a_i\in E :a_i<n\}} \quad
\text{with}$\hfill \cr
 \hskip 30pt~~$d_1=\gcd(n,\alpha_1)\in E$,\hfill \cr
  \noalign{\vskip -8pt}
 {}\hfill \cr
\noalign{\hrule}}}$$ where $\max{\{a_i\in E :a_i<n\}}$ is equal to
the second largest positive integer in $E$.

Thus, observe that $n$ and $\alpha_1$ were already computed by
{\rm(A.5.7)}. \ms

\roster

\item "(c)" Moreover, to finish this step, it is enough to
consider the following: Note that
$\gcd(n,\nu_2)=\gcd(n,\alpha_1)=d_1$ where $\alpha_1=\tau_0
n+\nu_2$.

\item "(c1)" If $\gcd(n,\alpha_1)=1$, then the problem is
completely solved.

\item "(c2)" If $\gcd(n,\alpha_1)>1$, take the next step.
\endroster \ms

$\underline{\text{\bf [II] The 2nd step:}}$ Let
$d_1=\gcd(n,\alpha_1)>1$. Note by  the first step that
$a_{q_1}=\mu_{1,q_1}=\gcd(n,\alpha_1)=\min\{P_1\}$ from (A.5.5) and
(A.5.6).

In this step, the aim is to find an elementary computational
algorithm for $\alpha_2-\alpha_1$, in construction of the Euclidean
divisor sequence $\{P_2\}=
\{[\gcd(n,\alpha_1):\alpha_2-\alpha_1]\}$. \ms

Then, it is enough to consider the following in order: \roster

\item "(a)" Let $\sigma_1$ be the counting number of
$\min\{P_1\}=\gcd(n,\alpha_1)$ in $\{[n,\alpha_1]\}=\{P_1\}$, and
$\tau_1$ be the counting number of the element $\gcd(n,\alpha_1)$ in
$E$. First, find $\sigma_1$ and $\tau_1$, respectively. It is clear
that $\tau_1\geqq\sigma_1$.

\item "(b)" Next, compute $\alpha_2-\alpha_1$, using the result of
(a) and Definition $A.3.0$.

\item "(c)" If $\gcd(n,\alpha_1,\alpha_2)=1$, then the algorithm
is completely finished.

\noindent If $\gcd(n,\alpha_1,\alpha_2)>1$, then take the third
step.
\endroster \ms

Now, let us solve (a), (b) and (c) in order.

(a) We are going to compute $\sigma_1$ and $\tau_1$, respectively.

(a1) Then, $\gcd(n,\alpha_1)=\mu_{1,q_1}$ in
$\{[n,\alpha_1]\}=\{P_1\}$ for coincidence of notation in $(A.5.5)$.
Then, by Definition $A.3$, we can compute $\sigma_1$ as follows:

Let $\mu_{1,s_1}$ be the second smallest element in $\{P_1\}$.
First, it is trivial to compute $\sigma_1$ such that
$\mu_{1,s_1}=\sigma_1\mu_{1,q_1}$ for $1\le s_1< q_1$.

(a2) It is clear that $\tau_1$ is directly countable in $E$. \ms

(b) In order to compute $\alpha_2-\alpha_1$, there are two cases:
Note that $\gcd(n,\alpha_1)=a_{q_1}$.

(b1) $\tau_1=\sigma_1$ and (b2) $\tau_1> \sigma_1$. \ms

(b1) Let $\tau_1=\sigma_1$. Then, observe by Definition $A.3$ that
$\alpha_2-\alpha_1<\gcd(n,\alpha_1)$.

We compute $\alpha_2-\alpha_1$ as follows:
$$\align
&\text{$\alpha_2-\alpha_1$= the
largest among the elements in $E$} \tag{$*$}\\
& \qquad \quad \quad \text{any of which is less than
$\gcd(n,\alpha_1)=d_1=\mu_{2,1}=a_{q_1+1}$}\\
\endalign$$

Therefore, $a_{q_1+1}\in E$ can be easily found from (A.5.5) and
(A.5.6). \ms

(b2) Let $\tau_1>\sigma_1$. Observe by Definition $A.3$ that
$\alpha_2-\alpha_1>\gcd(n,\alpha_1)$. So, by Definition $A.3$,
$$\align
&\text{$\alpha_2-\alpha_1$=$\gamma_1\gcd(n,\alpha_1)+$the second
largest element in $P_2$,}  \\
\endalign$$
where $\max\{P_2\}=\gcd(n,\alpha_1)$ and $\gamma_1$ is the number of
the element $\gcd(n,\alpha_1)$ in $\{P_2\}$. \ms

Observe that $\gamma_1=\tau_1-\sigma_1$, and also that the second
largest among the elements in $P_2$ is the same as the largest among
the elements in $E$, any of which is less than $\gcd(n,\alpha_1)$.
\ms

Find the largest among the elements in $E$, each of which is less
than $\gcd(n,\alpha_1)=\min\{[n,\alpha_1]\}$.

Therefore, $\alpha_2-\alpha_1$ can be computed as follows:
$$\align
(**) \qquad
&\text{$\alpha_2-\alpha_1$=$(\tau_1-\sigma_1)\gcd(n,\alpha_1)+$ the
largest among the elements in $E$} \qquad \qquad \\
& \qquad \text{any of which is less than $\gcd(n,\alpha_1)$.}
\endalign$$

Note that \text{$\gcd(n,\alpha_1)$= the largest among the elements}
in $\{[\gcd(n,\alpha_1):\alpha_2-\alpha_1]\}$ with $
\gcd(n,\alpha_1)<\alpha_2-\alpha_1$. \ms

$\underline{\text{\bf The algorithm for the second step:}}$ For
either of two cases(b1) and (b2), an elementary computational
algorithm formula for $P_2=[\gcd(n,\alpha_1):\alpha_2-\alpha_1]$
with $\alpha_2-\alpha_1$ can be represented as follows, at the same
time:

\noindent {\rm(A.5.8)}
$$\vbox{\offinterlineskip\tabskip=0pt\halign{\strut
\vrule~ # \vrule\cr \noalign{\hrule}
 {}\hfill \cr
 \noalign{\vskip -8pt}
\hskip 20pt \quad $\sigma_{1}=\dfrac{\min{\{a_i\in E
:a_i>d_{1}\}}}{d_{1}}$ \quad \text{and} \hfill   \cr \hskip 20pt
\quad  $\tau_{1}=\text{  the counting number of the element $d_{1}$
in $E$}$,\hfill \cr \hskip 30pt ~~ $\alpha_2-\alpha_{1}=
(\tau_{1}-\sigma_{1})d_{1}+ \max{\{a_i\in E :a_i<d_{1}\}} \quad
\text{with}$\hfill \cr \hskip 20pt \quad
$d_2=\gcd(d_1,\alpha_2-\alpha_1)\in E$.\hfill \cr \noalign{\vskip
-8pt}
 {}\hfill \cr
\noalign{\hrule}}}$$

Note that $\sigma_1=\dfrac{\min{\{a_i\in E
:a_i>d_1\}}}{d_1}=\dfrac{\text{the 2nd smallest element in
$\{P_1\}$}}{\text{$\min\{P_1\}$}}$ is the counting number of the
smallest element $d_1=\min\{P_1\}$ in $\{[n,\alpha_1]\}=\{P_1\}$,
and so $\tau_1-\sigma_1\ge 0$. \ms

Next, it is easy to compute $\alpha_2$ and
$\gcd(n,\alpha_1,\alpha_2)$ as follows:

{\rm(b3)} $\alpha_2=\alpha_1+(\alpha_2-\alpha_1)$.

{\rm(b4)}
$\gcd(n,\alpha_1,\alpha_2)=\gcd(n,\alpha_1,\alpha_2-\alpha_1)$. \ms

(c) Moreover, to finish this step, it is enough to consider the
following:

(c1) If $d_2=\gcd(n,\alpha_1,\alpha_2)=1$, then the problem is
completely solved.

(c2) If $d_2=\gcd(n,\alpha_1,\alpha_2)>1$, then take the next step.
\ms

The general case is proved by induction. Suppose we have shown that
the algorithm for the {\rm k-th} step is given by (A.5.4). In order
to find such an algorithm for {\rm (k+1)-th} step, it suffices to
prove the algorithm given by (A.5.9) later in the next case [III],
which can be computable as follows. For the proof of the {\rm
(k+1)-th} step, we may assume that $d_k>1$, otherwise there is
nothing to prove. \ms

$\underline{\text{\bf [III] The {(k+1)}-th step:}}$ Let
$d_{k}=\gcd(n,\alpha_1,\alpha_2,\dots,\alpha_{k})>1$. Note by the
{\rm k-th} step that $\mu_{{k},q_{k}}=d_{k}$ from (A.5.5) and
(A.5.6).

In this step, the aim is to find an elementary computational
algorithm for $\alpha_{k+1}-\alpha_{k}$, in construction of the
Euclidean divisor sequence
$\{P_{k+1}\}=\{[d_{k}:\alpha_{k+1}-\alpha_{k}]\}$. \ms

Then, it is enough to consider the following in order: \roster

\item "(a)" Let $\sigma_{k}$ be the counting number of
$d_{k}=\min\{P_{k}\}$ in
$\{[d_{k-1}:\alpha_{k}-\alpha_{k-1}]\}=\{P_{k}\}$, and $\tau_{k}$ be
the counting number of the element $d_{k}$ in $E$. First, find
$\sigma_{k}$ and $\tau_{k}$, respectively. It is clear that
$\tau_k\geqq\sigma_k$.

\item "(b)" Next, compute $\alpha_{k+1}-\alpha_{k}$ by the result
of (a) and the fundamental algorithm.

\item "(c)" If $d_{k+1}=1$, then the algorithm is completely
finished.

\noindent If $d_{k+1}>1$, then take the next step.
\endroster \ms

Now, let us solve (a), (b) and (c) in order.

(a) We are going to compute $\sigma_{k}$ and $\tau_{k}$,
respectively.

(a1)  Then,
$d_{k}=\gcd(n,\alpha_1,\alpha_2,\dots,\alpha_{k})=\mu_{{k},q_{k}}$
in $\{[d_{k-1}:\alpha_{k}-\alpha_{k-1}]\}=\{P_{k}\}$, for
coincidence of notation in (A.5.5). Then, by the fundamental
algorithm, we can compute $\sigma_{k}$ as follows:

Let $\mu_{k,s_k}$ be defined by the second smallest element in
$\{P_{k}\}$, which is easy to be found, for coincidence of notation
in (A.5.5). First, it is trivial to compute $\sigma_{k}$ such that
$\mu_{k,s_k}=\sigma_{k}\mu_{k,q_{k}}$ for $1\le s_k<q_{k}$.

(a2) It is clear that $\tau_{k}$ is directly countable in $E$. \ms

(b) In order to compute $\alpha_{k+1}-\alpha_{k}$, there are two
cases:

(b1) $\tau_{k}=\sigma_{k}$ and (b2) $\tau_{k}> \sigma_{k}$. \ms

(b1) Let  $\tau_{k}=\sigma_{k}$. Then, observe by Definition $A.3$
that $\alpha_{k+1}-\alpha_{k}<d_k$.

We compute $\alpha_{k+1}-\alpha_{k}$ as follows:
$$\align
&\text{$\alpha_{k+1}-\alpha_{k}$= the
largest among the elements in $E$} \tag $*$ \\
& \qquad \quad \quad \text{any of which is less than
$d_k=\mu_{k+1,1}=a_{q_1+q_2+\dots+q_{k}+1}$.}\\
\endalign$$

(b2) Let $\tau_{k}> \sigma_{k}$. Then, observe by Definition $A.3$
that $\alpha_{k+1}-\alpha_{k}>d_k$.

So, by the fundamental algorithm,
$$\align
&\text{$\alpha_{k+1}-\alpha_{k}$= $\gamma_{k}d_{k}+$the
second largest element in $P_{k+1}$,} \\
\endalign$$
where $\max\{P_{k+1}\}=d_{k}$ and $\gamma_{k}$ is the number of the
element $d_{k}=\min\{P_{k}\}$ in $\{P_{k}\}$. \ms

Observe that $\gamma_{k}=\tau_{k}-\sigma_{k}$, and also that the
second largest among the elements in $P_{k+1}$ is the same as the
largest among the elements in $E$, any of which is less than $d_{k}
=\min\{P_{k}\}$. \ms

Find the largest among the elements in $E$, each of which is less
than $d_{k} =\min\{P_{k}\}$.

Therefore, $\alpha_{k+1}-\alpha_{k}$ can be computed as follows:
$$\align
(**) \qquad \qquad &\text{$\alpha_{k+1}-\alpha_{k}$
=$(\tau_{k}-\sigma_{k})d_{k}+$
the largest among the elements in $E$,} \qquad \qquad \\
& \qquad \text{any of which is less than $d_{k}$}.
\endalign$$ \ms

Note that \text{$d_{k}$= the largest among the elements} in
$\{P_{k+1}\}=\{[d_{k}:\alpha_{k+1}-\alpha_{k}]\}$ with
$d_{k}<\alpha_{k+1}-\alpha_{k}$. \ms

$\underline{\text{\bf The algorithm for the (k+1)-th step:}}$ For
either of two cases (b1) and (b2), an elementary computational
algorithm formula for
$\{P_{k+1}\}=\{[d_{k}:\alpha_{k+1}-\alpha_{k}]\}$ with
$\alpha_{k+1}-\alpha_{k}$ can be represented as follows, at the same
time:

\noindent {\rm(A.5.9)}
$$\vbox{\offinterlineskip\tabskip=0pt\halign{\strut
\vrule~ # \vrule\cr \noalign{\hrule}
 {}\hfill \cr
 \noalign{\vskip -8pt}
\hskip 20pt \quad  $\sigma_{k}=\dfrac{\min{\{a_i\in E
:a_i>d_{k}\}}}{d_{k}}$ \quad \text{and} \hfill   \cr \hskip 20pt
\quad $\tau_{k}=\text{  the counting number of the element $d_{k}$
in $E$}$,\hfill \cr \hskip 30pt~~ $\alpha_{k+1}-\alpha_{k}=
(\tau_{k}-\sigma_{k})d_{k}+ \max{\{a_i\in E :a_i<d_{k}\}} \quad
\text{with}$\hfill \cr \hskip 30pt \quad
$d_{k+1}=\gcd(d_{k},\alpha_{k+1}-\alpha_{k})\in E$.\hfill \cr
\noalign{\vskip -8pt}
 {}\hfill \cr
\noalign{\hrule}}}$$

Note that $\sigma_{k}=\dfrac{\min{\{a_i\in E
:a_i>d_{k}\}}}{d_{k}}=\dfrac{\text{the 2nd smallest element in
$\{P_{k}\}$}}{\text{$\min\{P_{k}\}$}}$ is the counting number of the
smallest element $d_{k}=\min\{P_{k}\}$ in $\{P_{k}\}$ for $k\ge 1$,
and so $\tau_{k}-\sigma_{k}\ge 0$. \ms

Next, it is easy to compute $\alpha_{k+1}$ and
$d_{k+1}=\gcd(n,\alpha_1,\alpha_2,\dots,\alpha_{k+1})$ as follows:

{\rm(b3)} $\alpha_{k+1}=\alpha_{k}+(\alpha_{k+1}-\alpha_{k})$.

{\rm(b4)} $d_{k+1}=\gcd(n,\alpha_1,\alpha_2,\dots,\alpha_{k+1})
=\gcd(d_{k},\alpha_{k+1}-\alpha_{k})$ for $k\ge 1$. \ms

(c) Moreover, to finish this step, it is enough to consider the
following:

(c1) If $d_{k+1}=1$, then the problem is completely solved.

(c2) If $d_{k+1}>1$, then take the next step.

Thus, the proof of the theorem is completely finished. $\square$
\enddemo \ms

\definition{Remark A.5.2}
{\rm(a)} Assuming that $k=1$ in the above computational algorithm
formula (A.5.4) for the generalized k-th step, by the convenience of
notation, write $d_0=n$, and then $\sigma_0=0$. As we have seen in
(A.5.3), let $\tau_{0}$ be the counting number of the element
$d_0=n$ in $E$ when $E$ is viewed as a finite sequence. Thus, if
$k=1$, then the elementary computational algorithm formula for
$\alpha_1-\alpha_0$ with $d_1$ can be viewed as follows: Note that
$d_0=n$ and $\sigma_0=0$. Also, we write $\alpha_0=0$, if necessary.

\noindent (A.5.4)*
$$\vbox{\offinterlineskip\tabskip=0pt\halign{\strut
\vrule~ # \vrule\cr \noalign{\hrule}
 {}\hfill \cr
 \noalign{\vskip -8pt}
\hskip 20pt \quad  $\sigma_{0}=0 \quad \text{and}$ \hfill   \cr
\hskip 20pt \quad $\tau_{0}=\text{  the counting number of the
element $d_{0}$ in $E$}$,\hfill \cr \hskip 30pt~~
$\alpha_1-\alpha_{0}= (\tau_{0}-\sigma_{0})d_{0}+ \max{\{a_i\in E
:a_i<d_{0}\}} \quad \text{with}$\hfill \cr \hskip 30pt
$d_1=\gcd(d_{0},\alpha_1-\alpha_{0})\in E$. \hfill \cr
\noalign{\vskip -8pt}
 {}\hfill \cr
\noalign{\hrule}}}$$

{\rm(b)} A computational method for solving an existence of the
standard Puiseux expansion being equivalent to any given finite
sequence of all positive integers can be found by the same algorithm
as we have used in Theorem $A.5$.
\enddefinition \ms

\newpage

$$\align
 &\qquad \qquad\qquad  \text{\bf {Appendix B}}  \\
\quad &\text{\bf The $\beta$-algorithm for finding a one-to-one}\\
 & \text{\bf function from Family(2) onto Family(4) }\\
\endalign$$ \ms

{\bf \S{ B.0.} Introduction}

In Appendix B, the aim is to find an algorithm for computing a
one-to-one function from Family(2) onto Family(4) in Definition A.1,
called the $\beta$-algorithm for notation. \ms

{\bf \S{B.1.} The terminology and notations in preparation for
finding the $\beta$-algorithm}

We use the same terminology and notations as it has been in
Definition A.1 and Definition A.2 of Appendix A. \ms

{\bf {\S} B.2 and {\S} B.3. The $\beta$-algorithm(The algorithm for
finding a one-to-one function between Family(2) and Family(4)) and
its examples}\ms

{\bf {\S} B.2. The 1st half of the ${\beta}$ algorithm(Theorem B.2)}

\proclaim{Theorem B.2(The 1st half of the $\beta$ algorithm: an
algorithm for finding a one-to-one function from Family(4) into
Family(2))}

$\underline{\text{\bf {Assumptions}}}$ As in {\rm Family(4)} of {\rm
Definition A.1}, for any $f \in \text{\rm Family(0)}$ let
$(f\circ\tau_{\xi})_{divisor}$ be the divisor of $f\circ\tau_{\xi}$
given by
$$
(f\circ\tau_{\xi})_{divisor}=V^{(\xi)}(f)+\sum^{\xi}_{i=1}e_iE_i,
\tag B.2.1 $$ where each $e_i$ is the multiplicity of
$f\circ\tau_{\xi}$ along $E_i$ for $1\le i\le {\xi}$ and
$e_{i+1}>e_i$. Note that $\tau_{\xi}=\pi_1\circ\pi_2\circ\cdots
\circ\pi_{\xi}:M^{(\xi)}\to \BC^2$ is the composition of a finite
number $\xi$ of successive blow-ups $\pi_i$ at the origin in
$\BC^2$, which is needed only to get the standard resolution of the
singularity of $V(f)$. \ms

$\underline{\text{\bf Conclusions}}$ To compute the standard Puiseux
expansion $C(t)$ with \text{$f \buildrel \text{{\rm multiseq}} \over
\sim C(t)$} at the origin in $\BC^2$ is to find {\bf The 1st half of
the $\beta$ algorithm}, assuming that the standard Puiseux expansion
$C_r(t)$ of the $r$-type satisfies an equation for some r in
{\rm(B.2.2)}:
$$\text{$\rm C(t)=C_r(t):=$} \left\{\eqalign{ y& =t^n  \cr
 z &=t^{\alpha_1}+t^{\alpha_2}+\cdots +t^{\alpha_r}, \cr}
 \right. \tag B.2.2 $$
where  {\rm(i)} \quad $2\le n<\alpha_1<\alpha_2<\cdots
 <\alpha_r$ for a positive integer $r$ and

\noindent {\rm(ii)}\quad $n> d_1>d_2>\cdots
>d_{r}=1$, and
$d_{i}=\gcd(n,\alpha_1,\dots,\alpha_{i})$ for $1\le i\le r$. \ms

$\underline{\text{\bf (The 1st half of the $\beta$ algorithm)}}$
\quad For notation, write $\lambda_r=\xi$. Then, by {\rm Sublemma
A.2.1 of Definition A.2(Theorem 14.0 of Part[C])} and by
{\rm(B.2.2)}, there are exactly $r$ exceptional curves of the first
kind, written by $\{ E_{\lambda_w}:
1=\lambda_0<\lambda_1<\cdots<\lambda_w<\cdots<\lambda_r=\xi \}$,
each of which has three distinct intersection points with other
exceptional curves and the proper transform $V^{(\lambda_r)}(f)$
under $\tau_{\lambda_r}$. It is clear that $n=e_1$. So, to find the
desired solution, it suffices to compute $\alpha_w$ and
$e_{\lambda_w}$ for each {\rm w-th} step on the positive integer
$w=1,2,\dots,r$. \ms

$\underline{\text{\bf Step 1 for {\bf(The 1st half of the $\beta$
algorithm)}}}$ For this step, the problem is how to compute the
Puiseux exponent $\alpha_1$ and the coefficient $e_{\lambda_1}$.

Note that $n=e_1$, and let $S_1$ be the set defined by
$$\align
S_1 &=\{\text{$e_i\in S_0$ : $e_i$ cannot be divisible by $e_1$ for
$i>1$}\}, \tag B.2.3 \\
\endalign$$
where $S_0=\{\text{$e_i$: $i=1,2,\dots,\xi$}\}$ by
\text{\rm(B.2.1)}. \ms

Then, $\underline{\text{\rm the computational algorithm formula for
$\alpha_1$ and the coefficient $e_{\lambda_1}$ of $E_{\lambda_1}$}}$
can be represented as follows: \text{\rm{lcm(n,{$\alpha_1$})}=the
least common multiple of $n$ and $\alpha_1$}.
$$\align
\text{\rm (B.2.4)} \quad \quad &  n=e_1, \quad \alpha_1=\min S_1
\quad \text{and} \quad e_{\lambda_1}= \text{\rm $\text{{\rm
lcm}}(n,\alpha_1)$}=n_1\beta_{1,1}d_1,  \qquad \qquad \qquad \qquad\\
 \quad &\text{where \quad $d_1=\gcd(n,\alpha_1)$\quad with \quad $n=n_1d_1$ \quad
 {and} \quad $\alpha_1=\beta_{1,1}d_1$.} \\
\endalign$$

Then, there are two subcases: \ms

\noindent$\underline{\text{\rm Subcase(1)}}$ If $d_1=1$, then the
desired standard Puiseux expansion $C(t)=C_{1}(t)$ can be defined by
$y=t^n$ and $z=t^{\alpha_1}$.

\noindent$\underline{\text{\rm Subcase(2)}}$ If $d_1>1$, then take
the next step. \bs

Let $d_{w-1}=\gcd(n,\alpha_1,\alpha_2,\dots,\alpha_{w-1})>1$ with
$w\ge 2$. Then, an elementary computational algorithm formula for
$\alpha_w$ can be represented as follows: \ms

$\underline{\text{\bf Step w for (The 1st half of the $\beta$
algorithm)}}$ \quad Let $d_{w-1}>1$ with $2\le w\le r$.

Consider $(f\circ\tau_{\lambda_r})_{divisor}-V^{(\lambda_r)}(f)
-\sum^{\lambda_{w-1}}_{i=1}e_iE_i=\sum^{\lambda_r}_{i=\lambda_{w-1}+1}e_iE_i$.

For this step, the problem is how to compute the $\alpha_w$ and the
coefficient $e_{\lambda_w}$, assuming by the induction method that
each $\alpha_k$ and the coefficients $e_{\lambda_k}$ of
$E_{\lambda_k}$ have been already computed for any
$k=1,2,\dots,w-1$. \ms

Let $S_{w}=\{e_i-e_{\lambda_{w-1}}: \text{$ \lambda_{w-1}+1\le
i\le\lambda_r$}\}$.

For {\rm Step w}, we have exactly two cases:

{\rm Case(i)} $e_{(\lambda_{w-1}+1)}-e_{(\lambda_{w-1})}=d_{w-1}$.

{\rm Case(ii)} $e_{(\lambda_{w-1}+1)}-e_{(\lambda_{w-1})}<d_{w-1}$.
\ms

$\underline{\text{\bf Case(i) of Step w for (The 1st half of the
$\beta$ algorithm)}}$ \quad Let
$e_{(\lambda_{w-1}+1)}-e_{(\lambda_{w-1})}=d_{w-1}$.

Then, $\underline{\text{\rm the computational algorithm formula for
$\alpha_w$, $d_w$, and the coefficient $e_{\lambda_w}$ of
$E_{\lambda_w}$}}$ can be represented as follows:  It can be proved
by \text{\rm(14.0.3)} of \text{\rm Theorem 14.0} of \text{\rm
Part[C]} that $e_{(\lambda_{w-1})}$ is divisible by $d_{w-1}$ for
each $w\ge 2$.
$$\align
\text{\rm(B.2.5)} \quad  &\text{First, compute a unique positive
integer $s$
with $e_{(\lambda_{w-1}+s)}$, as follows:} \qquad \qquad \qquad\\
&\text{{\rm(i)}\quad$e_{(\lambda_{w-1}+1)}-e_{(\lambda_{w-1})}=d_{w-1}$.}\\
&\text{{\rm(ii)}\quad
$e_{(\lambda_{w-1}+t)}-e_{(\lambda_{w-1})}=td_{w-1}$
for each $t=1,2,\dots,s-1$.} \\
&\text{{\rm(iii)} \quad
$e_{(\lambda_{w-1}+s)}-e_{(\lambda_{w-1})}\not =sd_{w-1}$.}
\endalign$$
$$\align
\text{\rm(B.2.6)} \qquad \text{By {\rm (B.2.5)},}\quad &
\text{compute $\alpha_w$ and $d_w$, and $e_{\lambda_w}$(if
necessary), as follows:} \qquad \qquad \qquad \\
 \alpha_w
-\alpha_{w-1}&=e_{(\lambda_{w-1}+s)}-e_{(\lambda_{w-1})} \quad
\text{with} \quad
e_{(\lambda_{w-1}+1)}-e_{(\lambda_{w-1})}=d_{w-1}, \qquad \qquad \\
 {noting~ that}~ d_w &=\gcd(d_{w-1},\alpha_w -\alpha_{w-1})
\quad \text{with $d_{w-1}=d_wn_w$  and $\alpha_w
-\alpha_{w-1}=d_w\text{$\widehat{\Delta}$}_w$}, \\
\text{and} \quad e_{\lambda_w}&=e_{\lambda_{w-1}}n_w
 +d_wn_w\text{$\widehat{\Delta}$}_w \quad \text{\rm by (A.3.4)}.\\
\endalign$$

After the computation is done for {\rm Case(i) of Step w}, there are
two subcases for {\rm Case(i)}:\ms

\noindent$\underline{\text{\rm Subcase(1)}}$ \quad If $d_w=1$, then
the standard Puiseux expansion can be defined by $y=t^n$ and
$z=t^{\alpha_1}+t^{\alpha_2}+\cdots +t^{\alpha_w}$.

\noindent$\underline{\text{\rm Subcase(2)}}$ \quad If $d_w>1$, take
the next step, $\underline{\text{\rm Step (w+1) for (The 1st half of
the $\beta$ algorithm)}}$. \bs

$\underline{\text{\bf Case(ii) of Step w for (The 1st half of the
$\beta$ algorithm)}}$ \quad Let
$e_{(\lambda_{w-1}+1)}-e_{(\lambda_{w-1})}<d_{w-1}$.

Then, $\underline{\text{\rm the computational algorithm formula for
$\alpha_w$, $d_w$ and the coefficient $e_{\lambda_w}$ of
$E_{\lambda_w}$}}$ can be represented as follows:
$$\align
\text{\rm (B.2.7)} \quad \alpha_w-\alpha_{w-1}
&=e_{(\lambda_{w-1}+1)}-e_{(\lambda_{w-1})} \quad \text{with} \quad
e_{(\lambda_{w-1}+1)}-e_{(\lambda_{w-1})}<d_{w-1}, \qquad \qquad\\
{noting~ that}~ d_w &=\gcd(d_{w-1},\alpha_w -\alpha_{w-1}) \quad
\text{with $d_{w-1}=d_wn_w$  and $\alpha_w
-\alpha_{w-1}=d_w\text{$\widehat{\Delta}$}_w$}, \\
\text{and} \quad e_{\lambda_w}&=e_{\lambda_{w-1}}n_w
 +d_wn_w\text{$\widehat{\Delta}$}_w \quad \text{\rm by (A.3.4)}.\\
\endalign$$

After the computation is done for {\rm Case(ii) of Step w}, there
are two subcases for {\rm Case(ii)}: \ms

\noindent$\underline{\text{\rm Subcase(1)}}$ \quad If $d_w=1$, then
the standard Puiseux expansion can be defined by $y=t^n$ and
$z=t^{\alpha_1}+t^{\alpha_2}+\cdots+t^{\alpha_w}$.

\noindent$\underline{\text{\rm Subcase(2)}}$ If $d_w>1$, take the
next step, $\underline{\text{\rm Step (w+1) for (The 1st half of the
$\beta$ algorithm)}}{\square}$.
\endproclaim \ms

{\bf Remark.} Note by (B.2.4) of Step 1 and by Theorem 3.6 of
Part[B] that $g_r=(z^{n_1}+\ve_1y^{\beta_{1,1}})^{d_1}
+\sum_{\alpha,\beta\ge 0} c^{(0)}_{\alpha,\beta}y^{\alpha}z^{\beta}$
where $\ve_1$ is a unit $\BC\{y,z\}$ and the
$c^{(0)}_{\alpha,\beta}$ are some nonzero complex numbers such that
$n_1\alpha+\beta_{1,1}\beta>n_1\beta_{1,1}d_1$.

\proclaim{Remark B.2.1 for Theorem B.2}

\noindent{\bf(1)} After $\lambda_{w-1}$ iterations of blow-ups, let
$(v_{\lambda_{w-1}},u_{\lambda_{w-1}})$ and
$(v'_{\lambda_{w-1}},u'_{\lambda_{w-1}})$ be the local coordinates
for $M^{({\lambda_{w-1}})}$ where
$\pi_{\lambda_{w-1}}:M^{({\lambda_{w-1}})}\to
M^{({\lambda_{w-1}}-1)}$ was defined to be the  $\lambda_{w-1}$-th
blow-up at some point of $M^{(\lambda_{w-1}-1)}$ with
$u'_{\lambda_{w-1}}=1/u_{\lambda_{w-1}}$ and
$v'_{\lambda_{w-1}}=v_{\lambda_{w-1}}u_{\lambda_{w-1}}$. Note that
$E_{\lambda_{w-1}}=\{v_{\lambda_{w-1}}=0\}\cup
\{v'_{\lambda_{w-1}}=0\}$ is the $\lambda_{w-1}$-th exceptional
curve of the first kind.
$$\align
\text{\rm(B.2.7.1)} \qquad (g_{r}\circ\tau_{\lambda_{w-1}})_{total}
&=v^{e_{(\lambda_{w-1})}}_{\lambda_{w-1}}
(g_{r}\circ\tau_{\lambda_{w-1}})_{proper}
\quad \text{with $f=g_r \in \text{\rm Family(1)}$,} \qquad \qquad     \qquad \qquad  \\
(g_{r}\circ\tau_{\lambda_{w-1}})_{proper}
&=\{(1+\ve_{w-1}u_{\lambda_{w-1}})^{n_{w}}
+v^{\widehat{\Delta}_{w}}_{\lambda_{w-1}}\}^{d_{w}}\\
&\quad +\sum_{\alpha,\beta\ge 0}B^{(w-1)}_{\alpha,\beta}
v^{\alpha}_{\lambda_{w-1}}(1+\ve_{w-1}u_{\lambda_{w-1}})^{\beta},\\
\endalign$$
where {\rm(i)}  \text{$e_{\lambda_w}
={e_{\lambda_{w-1}}}n_w+d_wn_w\text{$\widehat{\Delta}$}_w$} for
$2\le j\le r$, and  $e_{\lambda_1}=d_1n_1\beta_{1,1}$,

{\rm(ii)} $\ve_{w-1}=\ve_{w-1}
(1+\ve_{w-1}u_{\lambda_{w-1}},v_{\lambda_{w-1}})$ is a unit in
$\BC\{1+\ve_{w-1}u_{\lambda_{w-1}},v_{\lambda_{w-1}}\}$ and the
$B^{(w-1)}_{\alpha,\beta}$ are some nonzero complex numbers,

{\rm(iii)}
$d_{w}=\gcd(n,\alpha_1,\dots,\alpha_{w})=n_{w+1}n_{w+2}\cdots n_r\ge
2$,

{\rm(iv)} $\widehat{\Delta}_{w}=\Delta_{w}(\beta_{w,k})^{w}_{k=1}
-n_{w}n_{w-1}\Delta_{w-1}(\beta_{w-1,k})^{w-1}_{k=1}>0$,

{\rm(v)} $\gcd(n_{w},\widehat{\Delta}_{w})=1$ and
$n_{w}\alpha+\widehat{\Delta}_{w}\beta>n_{w}\widehat{\Delta}_{w}d_{w}$.
\ms

\noindent{\bf(2)} Observe by {\rm(B.2.7.1)} that the following
remark may be necessary for the construction of the statement of
{\rm Theorem B.2}.

$\underline{\text{\rm Case(i) of Step w for Solution}:}$ Note that
there exists an integer $s>0$ such that
$(s-1)n_w<\widehat{\Delta}_w< sn_w$ because
$2\le{n_w}<\widehat{\Delta}_w$, and then
$e_{(\lambda_{w-1}+s)}=e_{(\lambda_{w-1})}+d_{w}\widehat{\Delta}_w$.
\ms

$\underline{\text{\rm Case(ii) of Step w for Solution}:}$ Note that
there exists an integer $s>0$ such that
$(s-1)\widehat{\Delta}_w<n_w\le s\widehat{\Delta}_w$ because $1
\le\widehat{\Delta}_w<n_w$, and then
$e_{(\lambda_{w-1}+s)}=se_{(\lambda_{w-1})}+d_{w-1}$.
\endproclaim \ms

\noindent{\bf Example B.2.2 for Theorem B.2.} As in Definition A.1,
let $(f\circ\tau_{\xi})_{divisor}$ be the divisor of
$f\circ\tau_{\xi}$ given by the following:
$$(f\circ\tau_{\xi})_{divisor}
=V^{(\xi)}(f)+\sum^{\xi}_{i=1}e_iE_i \quad \text{with} \quad \xi=12,
\tag B.2.8 $$ where $e_1=45$, $e_2=60$, $e_3=120$, $e_4=180$,
$e_5=186$, $e_6=372$, $e_7=555$, $e_8=930$, $e_9=933$, $e_{10}=935$,
$e_{11}=1869$, $e_{12}=2805$. \ms

To find the standard Puiseux expansion $C_r(t)$ for the curve $C$
such that \text{$C_r(t)\equiv V(f)$} \text{\rm (multi. seq.)}, by
{\bf (The 1st half of the $\beta$ algorithm for Theorem B.2)}, just
compute the following:

$\underline{\text{\rm Step 1:}}$ It is clear that $n=e_{1}=45$, and
$\alpha_1=e_{2}=60$.

In preparation for the computation in Step 2,
$d_1=\gcd(n,\alpha_1)=\gcd(45,60)=15>1$ because $n=n_1d_1$ and
$\alpha_1=\beta_{1,1}d_1$ imply that $n_1=3$, $\beta_{1,1}=4$, and
$e_{\lambda_1}=n_1\beta_{1,1}d_1=180$ by (A.3.4). \ms

$\underline{\text{\rm Step 2:}}$ Since
$e_{\lambda_1+1}-e_{\lambda_1}=186-180=6<15=d_1$, by Case(ii) of
Step 2, ${\alpha_2}-{\alpha_1}=e_{{\lambda_1}+1}-e_{\lambda_1}=6$.
So, it is clear by Step 1 that ${\alpha_2}={\alpha_1}+6=60+6=66$.

In preparation for the computation in Step 3, $d_2>1$ and
$e_{\lambda_2}=1155$ because of the following computations:

$d_2=\gcd(d_1,\alpha_2-\alpha_1)=\gcd(15,6)=3>1$ because
$d_1=n_2d_2$ and $\alpha_2-\alpha_1=\text{$\widehat{\Delta}$}_2d_2$
imply that $n_2=5$ and $\text{$\widehat{\Delta}$}_2=2$, and so
$e_{\lambda_2}=n_2\cdot{e_{\lambda_{1}}}+d_2n_2\text{$\widehat{\Delta}$}_2
=5\cdot{180}+{3}\cdot{5}\cdot{2}=930$ by (A.3.4). \ms

$\underline{\text{\rm Step 3:}}$ Since
$e_{\lambda_2+1}-e_{\lambda_2}=933-930=3=d_2$, by Case(i) of Step 3,
$\alpha_3-\alpha_2=e_{(\lambda_{2}+s_2)}-e_{\lambda_{2}}=5$ because
$e_{(\lambda_{2}+s_2)}$ with $s_2>0$ is defined to be the first
appearing integer which cannot be divisible by $d_2$. Note by
(14.0.3) of Theorem 14.0 of Part[C] that $e_{(\lambda_{j})}$ is
divisible by $d_j$. So, it is clear by Step 2 that
${\alpha_3}={\alpha_2}+5=66+5=71$. \ms

So, the standard Puiseux expansion for $C(t)$ such that
\text{$C(t)\equiv V(f)$} \text{\rm (multi. seq.)} can be given by
$y=t^{45}$ and $z=t^{60}+t^{66}+t^{71}$, noting that $d_3=1$ because
$d_2=3$. \bs

{\bf {\S} B.3. The 2nd half of the $\beta$ algorithm(Theorem B.3)}

\proclaim{Theorem B.3(The 2nd half of the $\beta$ algorithm: an
algorithm for finding a function from Family(4) onto Family(2))}

\noindent $\underline{\text{\bf {Assumptions}}}$ {\rm(1)} Let the
standard Puiseux expansion $C_r(t)$ of the $r$-th type for the curve
$C$ in {\rm Family(2)} be given by
$$\align
\text{\rm(B.3.1)} \quad \quad \text{$C(t)=C_r(t):=$}
\left\{\eqalign{ y=&t^n, \cr z=&t^{\alpha_1}+t^{\alpha_2}+\cdots
+t^{\alpha_r}, \cr} \right.
\\
\text{where} \quad
 2\le n <\alpha_1<\alpha_2<\cdots <\alpha_r  & \quad  \text{and} \\
 n >d_1>d_2>\cdots
 >d_r=1  \quad   \text{with} & \quad
\text{$d_i=\gcd(n,\alpha_1,\dots,\alpha_i)$, $1\le i\le r$.} \qquad
\qquad
\endalign$$

{\rm(2)} Without any need of proof, we may assume by {\rm Theorem
1.6} of {\rm Part[A]} that there exists the standard Puiseux
polynomial $g_r\in \BC\{y,z\}$ of the recursive $r$-type in {\rm
Family(1)} such that \text{$V(g_r)\equiv C_r(t)$ \text{\rm (multi.
seq.)}}. \ms

{\rm(3)} To compute $(g_r\circ\tau_{\lambda_r})_{divisor}$ of {\rm
Family(2)} with \text{$V(g_r)\equiv C_r(t)$} \text{\rm (multi.
seq.)} directly, for brevity of notation we may assume by {\rm
Theorem 1.6} of {\rm Part[A]} that
$(g_r\circ\tau_{\lambda_r})_{divisor}$ satisfies the same kind of
properties and notations in {\rm Definition A.3.2}.\ms

\noindent $\underline{\text{\bf Conclusions}}$ By \text{\rm(A.3.2)}
and \text{\rm(A.3.3)} of {\rm Definition A.3.1}, the problem is how
to write explicit algorithm for finding the sequence
$\{(g_r\circ\tau_{\lambda_r})_{divisor}\}_{seq.}\equiv\{e_i:i=1,2,\dots,\lambda_r\}$
as sequence, which consists of the coefficients of $E_i$ for
$(g_r\circ\tau_{\lambda_r})_{divisor}$.

Let $S=\{e_i:i=1,2,\dots,\lambda_r\}$. By \text{\rm Subdefinition
A.2.1} of \text{\rm Definition A.2} and by \text{\rm Definition
A.3}, $S$ is well-defined by $S=\text{\rm Join}(B_1,B_2\dots,B_r)$
with the following properties:
$$\align
\text{\rm (B.3.2)} \qquad
& B_1=\{b_{0,i}=e_i: i=1,2,\dots,\lambda_1 \},  \qquad \qquad  \\
&
B_w=\{b_{w-1,i}=e_{\lambda_{w-1}+i}:i=1,2,\dots,(\lambda_w-\lambda_{w-1})\}
\quad \text{for $w=2,3,\dots,r$.} \qquad \qquad \qquad \qquad
\endalign$$

{\bf A desired algorithm(The 2nd half of the $\beta$ algorithm for
Theorem B.3)} for finding the unique sequence $S=\text{\rm
Join}\{B_1,B_2\dots,B_r\}$ is as follows:

In preparation for finding such a desired algorithm by {\rm(A.3.4)},
first the computation formulas for \{$d_1$, $n_1$, $\beta_{1,1}$\}
by \text{\rm (i)}, for \{$d_w$, $n_w$, $\widehat{\Delta}_w$: $2\le
w\le r$\} by \text{\rm (ii)}, and for \{$e_{\lambda_w}$: $1\le w\le
r$\} by \text{\rm (iii)} respectively, will be most applicable:
$$\align
\text{\rm (B.3.3)} \quad \text{\rm(i)} \quad
&\text{$d_1=\gcd(n,\alpha_1)$ with
$n=n_1d_1$ and $\alpha_1=\beta_{1,1}d_1$.} \\
\text{\rm{(ii)}} \quad
&\text{$d_w=\gcd(d_{w-1},\alpha_w-\alpha_{w-1})$ with
$d_{w-1}=n_wd_w$ and $\alpha_w-\alpha_{w-1}= \widehat{\Delta}_wd_w$}. \qquad\qquad\\
\text{\rm{\quad(iii)}} \quad
&\text{\text{$e_{\lambda_1}=d_1n_1\Delta_1(\beta_{1,1})$ and
$e_{\lambda_w}=e_{\lambda_{w-1}}n_w+d_wn_w\text{$\widehat{\Delta}$}_w$}
for $2\le w\le r$}.
\endalign$$

{\bf(The 2nd half of the $\beta$ algorithm)}

$\underline{\text{\bf Step 1 for {\bf(The 2nd half of the $\beta$
algorithm)}}}$ To find $B_1$, let $h_0=(z^{n_1}+y^{\beta_{1,1}})$
and $H_0=\prod^{d_1}_{i=1}(z^{n_1}+ c_iy^{\beta_{1,1}})$ where all
the $c_i$ are the nonzero distinct complex numbers for $1\le i\le
d_1$.

Let $\tau_{\lambda_1}$ be the composition of a finite number
$\lambda_1$ of successive blow-ups at the origin in $\BC^2$, which
is the standard resolution of the isolated singularity of $V(h_0)$.

To find $B_1$, it suffices to compute $(H_0\circ
\tau_{\lambda_1})_{divisor}$ under $\tau_{\lambda_1}$ because of the
following:
$$\align
\text{\rm (B.3.4)} \qquad
\text{$(H_0\circ\tau_{\lambda_1})_{divisor}=V^{(\lambda_1)}(H_0)
+\sum^{\lambda_1}_{i=1}{\ell}_iE_i$} \text{\quad where
$B_1=\{e_i={\ell}_i: i=1,2,\dots,\lambda_1 \}$.} \qquad
\endalign$$

$\underline{\text{\bf Step w for {\bf(The 2nd half of the $\beta$
algorithm)}}}$ To find $B_w(2\le w\le r)$, let $h_{w-1}=v(u^{n_w}+
v^{\widehat{\Delta}_w})$, and let
$H_{w-1}=\prod^{e_{\lambda_{w-1}}}_{j=1}(v+a_ju^{2n_w})\prod^{d_w}_{i=1}(u^{n_w}+
{c_i}v^{\widehat{\Delta}_w})$ where all the $a_j$ are the nonzero
distinct complex numbers for $1 \le j\le e_{\lambda_{w-1}}$ and all
the $c_i$ are the nonzero distinct complex numbers for $1 \le i\le
d_w$.

Let $\mu_{w,m_w}$ be the composition of a finite number $m_w$ of
successive blow-ups at $(u,v)=(0,0)$ which are the standard
resolution of the isolated singularity of both $V(h_{w-1})$ and
$V(H_{w-1})$.

For a better representation, $H_{w-1}$ can be also replaced by
$H_{w-1}=v^{e_{\lambda_{w-1}}}\{(u^{n_w}+
v^{\widehat{\Delta}_w})^{d_{w}}+\sum_{\alpha,\beta\ge 0}
c^{(w-1)}_{\alpha,\beta}v^{\alpha}u^{\beta}\}$ where the
$c^{(w-1)}_{\alpha,\beta}$ are some nonzero complex numbers such
that $n_w\alpha+\widehat{\Delta}_w\beta>n_w\widehat{\Delta}_wd_w$,
if necessary.

To find $B_w$, it suffices to compute $(H_{w-1}\circ
\mu_{m_w})_{divisor}$ by \text{\rm Theorem 3.6} of \text{\rm
Part[B]} and \text{\rm(B.2.7.1)} because of the following:
$$\align
\text{\rm (B.3.5)} \quad  &\text{$(H_{w-1}\circ
\mu_{m_w})_{divisor}=V^{(m_w)}(H_{w-1})
+\sum^{m_w}_{i=1}{\ell}_iE_i$, } \\
&\text{{where} \quad $B_w=\{e_{\lambda_{w-1}+i}={\ell}_i:
i=1,2,\dots,\lambda_w-\lambda_{w-1}\}$ \quad {with} \quad
$m_w=\lambda_w-\lambda_{w-1}$.} \qquad \qquad\\
\endalign$$

Note by \text{\rm Proposition 14.1} of \text{\rm Part[C]} that
$\tau_{\lambda_r}=\mu_{1,m_1}\circ\mu_{2,m_2}\circ\cdots
\circ\mu_{r,m_r}$ with $\lambda_r =m_1+m_2+\cdots +m_r$ and that
$e_{\lambda_{w-1}}$, ${n_w}$, ${\widehat{\Delta}_w}$ and ${d_{w}}$
were already known.
\endproclaim \ms

\noindent{\bf Example B.3.1 for Theorem B.3.} Let the
parametrization for the Puiseux expansion be given by the following:
$$
\text{$C(t):=$}  \left\{\eqalign{ y= &t^{45} \cr z=
&t^{60}+t^{65}+t^{71}. \cr } \right. \tag B.3.6
$$

This is the standard Puiseux expansion $C_3(t)$ of the 3rd type
because of the  computations {\rm(a)} and {\rm(b)}:

{\rm(a)} $n<\alpha_1<\alpha_2<\alpha_3$ where $n=45$, $\alpha_1=60$,
$\alpha_1=65$ and $\alpha_3=71$.

{\rm(b)} $n=45>d_1=15>d_2=5>d_3=1$ where
$d_1=\gcd(n,\alpha_1)=\gcd(45,60)=15$,
$d_2=\gcd(d_1,\alpha_2-\alpha_1)=\gcd(15,5)=5$ and
$d_3=\gcd(d_2,\alpha_3-\alpha_2)=\gcd(5,6)=1$.\ms

Then, {\rm(i)} Since $d_1=\gcd(n,\alpha_1)=15$, then $n=n_1d_1$ and
$\alpha_1=\beta_{1,1}d_1$ imply that $n_1=3$ and $\beta_{11}=4$.

\qquad {\rm(ii)} Since
$d_2=\gcd(d_1,\alpha_2-\alpha_1)=\gcd(15,5)=5$, then $d_1=d_2n_2$
and $\alpha_2-\alpha_1=d_2\widehat{\Delta}_2$ imply that $n_2=3$ and
$\widehat{\Delta}_2=1$.

\qquad {\rm(iii)} Since
$d_3=\gcd(d_2,\alpha_3-\alpha_2)=\gcd(5,6)=1$, then $d_2=d_3n_3$ and
$\alpha_3-\alpha_2=d_3\widehat{\Delta}_3$ imply that $n_3=5$ and
$\widehat{\Delta}_3=6$.

Also, $e_{\lambda_1}=n_1\beta_{1,1}d_1=180$,
$e_{\lambda_2}=n_2e_{\lambda_1}+n_2\widehat{\Delta}_2d_2=540+15=555$,
$e_{\lambda_3}=n_3e_{\lambda_2}+n_3\widehat{\Delta}_3d_3=2775+30=2805$.
\ms

Therefore, we have the following:

{\rm(i)} $H_0=(z^{n_1}+ y^{\beta_{1,1}})^{d_1}+\sum_{\alpha,\beta\ge
0} c^{(0)}_{\alpha,\beta}y^{\alpha}z^{\beta}=(z^{3}+
y^{4})^{15}+\sum_{\alpha,\beta\ge 0}
c^{(0)}_{\alpha,\beta}y^{\alpha}z^{\beta}$ where the
$c^{(0)}_{\alpha,\beta}$ are some nonzero complex numbers such that
$3\alpha+4\beta=n_1\alpha+\beta_{1,1}\beta>n_1\beta_{1,1}d_1=180$.

{\rm(ii)} $H_{1}=v_1^{e_{\lambda_{1}}}\{(u_1^{n_2}+
v_1^{\widehat{\Delta}_2})^{d_{2}}+\sum_{\alpha,\beta\ge 0}
c^{(1)}_{\alpha,\beta}v_1^{\alpha}u_1^{\beta}\}=v_1^{180}\{(u_1^{3}+
v_1)^{5}+\sum_{\alpha,\beta\ge 0}
c^{(1)}_{\alpha,\beta}v^{\alpha}u^{\beta}\}$ where the
$c^{(1)}_{\alpha,\beta}$ are nonzero complex numbers such that
$3\alpha+\beta=n_2\alpha+{\widehat{\Delta}_2}\beta>n_2{\widehat{\Delta}_2}d_2=15$.

{\rm(iii)} $H_{2}=v_2^{e_{\lambda_{2}}}\{(u_2^{n_3}+
v_2^{\widehat{\Delta}_3})^{d_{3}}+\sum_{\alpha,\beta\ge 0}
c^{(2)}_{\alpha,\beta}v^{\alpha}u^{\beta}\}=v_2^{555}\{(u_2^{5}+
v_2^{6})+\sum_{\alpha,\beta\ge 0}
c^{(2)}_{\alpha,\beta}v^{\alpha}u^{\beta}\}$ where the
$c^{(2)}_{\alpha,\beta}$ are nonzero complex numbers such that
$5\alpha+6\beta=n_3\alpha+{\widehat{\Delta}_3}\beta>n_3{\widehat{\Delta}_3}d_3=30$.
\bs

Following the same notations and properties as in {\bf (The 2nd half
of the $\beta$ algorithm for Theorem B.3)}, compute

\qquad \text{$(H_0\circ\tau_{\lambda_1})_{divisor}$ \quad {and}
\quad $(H_{j-1}\circ \mu_{m_j})_{divisor}$ \quad  {for $j=2,3$},
\quad respectively.} \ms

Using the blow-up process directly, then we can compute the
following:

$B_1=\{e_1,e_2,e_3,e_4\}=\{45,60,120,180\}$,
$B_2=\{e_5,e_5,e_7\}=\{185,370,555\}$, and
$B_3=\{e_8,e_9,e_{10},e_{11},e_{12},e_{13}\}=\{560,561,1122,1683,2244,2805\}$.
\bs

\vfill \pagebreak

$$\align
 &\qquad \qquad\qquad  \text{\bf {Appendix [C]}}  \\
\quad &\text{\bf The $\gamma$-algorithm for finding a one-to-one}\\
 & \text{\bf function from Family(2) onto Family(5) }\\
\endalign$$ \ms

{\bf \S{ C.0.} Introduction}

In Appendix [C], the aim is to find an algorithm for computing a
one-to-one function from Family(2) of Definition A.1 onto Family(5)
of Definition A.2, called the $\gamma$-algorithm for notation, which
must be proved to be well-defined by Theorem C.2 and Theorem C.3.
\ms

{\bf \S{C.1.} The terminology and notations in preparation for
finding the $\beta$-algorithm}

We use the same terminology and notations as it has been in
Definition A.1, Definition A.2 and Definition A.3 of Appendix A. \ms

\proclaim{Theorem C.1} As we have seen in the definition of a
singular exceptional curve of the first kind in {\rm Subdefinition
A.2.4 and Subdefinition A.2.5 of Definition A.2}, a singular part of
the divisor in {\rm Subdefinition A.2.5} is well-defined. Rigorously
speaking, it will be proved by two theorems, {\rm Theorem C.2} and
{\rm Theorem C.3} that we can compute The $\gamma$-algorithm for
finding a one-to-one function between {\rm Family(2)} and {\rm
Family(5)}.
\endproclaim \ms

{\bf \S {C.2} and \S {C.3.} The $\gamma$-algorithm(The algorithm for
finding a one-to-one function between Family(2) and Family(5)) and
its examples)} \ms

{\bf {\S} C.2.  The 1st half of The 3rd algorithm(Theorem C.2)}\ms

\proclaim{Theorem C.2(The 1st half of the $\gamma$-algorithm: an
algorithm for finding a one-to-one function from Family(2) into
Family(5))}

\noindent $\underline{\text{\bf {Assumptions}}}$ {\rm(1)} Let the
standard Puiseux expansion $C_r(t)$ of the $r$-th type for the curve
$C$ in {\rm Family(2)} be given by
$$\align
\text{\rm (C.2.1)} \qquad \text{$C_r(t):=$} \left\{\eqalign{ y=&t^n,
\cr z=&t^{\alpha_1}+t^{\alpha_2}+\cdots +t^{\alpha_r}, \cr} \right.
\\
\text{where} \quad
 2\le n <\alpha_1<\alpha_2<\cdots <\alpha_r  & \quad  \text{and} \\
 n >d_1>d_2>\cdots
 >d_r=1  \quad   \text{with} & \quad
\text{$d_i=\gcd(n,\alpha_1,\dots,\alpha_i)$, $1\le i\le r$.} \qquad
\qquad
\endalign$$

{\rm(2)} Without any need of proof, we may assume by {\rm Theorem
1.6} of \text{\rm Part[A]} that we can compute the standard Puiseux
polynomial $g_r\in \BC\{y,z\}$ of the recursive $r$-type in {\rm
Family(1)} such that \text{$V(g_r)\equiv C_r(t)$ \text{\rm (multi.
seq.)}}. \ms

\noindent $\underline{\text{\bf Conclusions}}$ By {\rm Definition
A.2}, it is clear that ${Family(4)}$ and ${Family(4)}_{seq}$ can be
identified and  that ${Family(5)}$ and ${Family(5)}_{seq}$ can be
identified.

Also, $S_1=\{(g_r\circ\tau_{\lambda_r})_{\text{\rm singular part of
the divisor}}\}_{seq.}\in {Family(4)}_{seq}$ is a proper subsequence
of a finite sequence
$S=\{(g_r\circ\tau_{\lambda_r})_{divisor}\}_{seq.}\in
{Family(5)}_{seq}$.

{\bf To solve the problem is how to find an algorithm for computing
{$S_1$} directly, without computing {$S$} completely.}

Let $S_1=\text{\rm Join}(T_1,T_2,\dots,T_r)$, satisfying the
following: Write $\lambda_0=1$, if necessary.
$$\align
\text{\rm(C.2.2)} \qquad \qquad   & T_1=\{e_i: i=1,s_0 \}
\quad \text{with $1<s_0\le{\lambda_1}$},  \qquad \qquad\qquad \qquad \qquad  \\
& T_j=\{e_{\lambda_{j-1}+i}:i=1,s_{j-1}\} \quad \text{for
$2\le j\le r$ and $1<s_{j-1}\le{\lambda_j-\lambda_{j-1}}$.} \qquad \\
\endalign$$

{\bf A desired algorithm(The 1st half of the $\gamma$-algorithm for
Theorem C.2)} for finding the sequence $S_1=\text{\rm
Join}(T_1,T_2\dots,T_r)$ is as follows:

The arithmetic computation formula for \{$d_1$, $n_1$,
$\beta_{1,1}$\} by \text{\rm (i)}, and \{$d_w$, $n_w$,
$\widehat{\Delta}_w$: $2\le w\le r$\} by \text{\rm (ii)}, and
\{$e_{\lambda_w}$: $1\le w\le r$\} by \text{\rm (iii)} in
{\rm(C.2.2)} will be most helpful for finding such an algorithm:
$$\align
(C.2.3) \quad \text{\rm(i)} \quad &\text{$d_1=\gcd(n,\alpha_1)$ with
$n=n_1d_1$ and $\alpha_1=\beta_{1,1}d_1$.} \\
\text{\rm{(ii)}} \quad
&\text{$d_w=\gcd(d_{w-1},\alpha_w-\alpha_{w-1})$ with
$d_{w-1}=n_wd_w$ and $\alpha_w-\alpha_{w-1}= \widehat{\Delta}_wd_w$}. \qquad\qquad\\
\text{\rm{\quad(iii)}} \quad
&\text{\text{$e_{\lambda_1}=d_1n_1\Delta_1(\beta_{1,1})$ and
$e_{\lambda_w}=e_{\lambda_{w-1}}n_w+d_wn_w\text{$\widehat{\Delta}$}_w$}
for $2\le w\le r$}.
\endalign$$

\noindent$\underline{\text{\bf(The 1st half of the
$\gamma$-algorithm for Theorem C.2)}}$

\noindent$\underline{\text{\bf Step 1 for the 1st half of the
$\gamma$-algorithm}}$ Compute $e_1$ and $e_{s_0}$ such that $e_1=n$
and $e_{s_0}=\alpha_1$.

Compute $d_1=\gcd(n,\alpha_1)$ with $n=n_1d_1$ and
$\alpha_1=\beta_{1,1}d_1$, and then
$e_{\lambda_1}=n_1\beta_{1,1}d_1$.\ms

Then, there are two subcases:

$\underline{\text{\rm Subcase(1)}}$ If $d_1=1$, it is clear that
$S_1=T_1=\{e_1,e_{s_0}\}=\{n,\alpha_1\}$. \ms

$\underline{\text{\rm Subcase(2)}}$ If $d_1>1$, take the next step.
\ms

$\underline{\text{\rm Remark}}$ We can compute an integer $s_0$ such
that $(s_0-1)n<\alpha_1<s_0n$, if necessary. \ms

Assuming that $d_1>1$, for $w=2,3,\cdots,r$, an elementary
computational algorithm formula with $\alpha_w$ can be represented
as follows:

\noindent$\underline{\text{\bf Step w for the 1st half of the
$\gamma$-algorithm}}$ \quad Let $d_{w-1}>1$ with $2\le w\le r$.

Then, there are two cases:

{\rm Case(i) of Step w} \quad
$d_{w-1}<\alpha_w-\alpha_{w-1}$\text{\rm($\iff$
$e_{(\lambda_{w-1}+1)}-e_{(\lambda_{w-1})}=d_{w-1}$)}.

{\rm Case(ii) of Step w} \quad $d_{w-1}>\alpha_w-\alpha_{w-1}$
\text{\rm($\iff$
$e_{(\lambda_{w-1}+1)}-e_{(\lambda_{w-1})}<d_{w-1}$)}. \ms

\noindent$\underline{\text{\bf Case(i) of Step w for the 1st half of
the $\gamma$-algorithm}}$ \quad Let $d_{w-1}<\alpha_w-\alpha_{w-1}$.

Then, the algorithm for finding a subsequence
$T_w=\{e_{(\lambda_{w-1}+1)}, e_{(\lambda_{w-1}+s_{w-1})}\}$ of
$S_1$ can be uniquely represented as follows:
$$\align
(C.2.4) \qquad\qquad e_{(\lambda_{w-1}+1)}&=e_{(\lambda_{w-1})}+d_{w-1},\\
e_{(\lambda_{w-1}+s_{w-1})}&=e_{(\lambda_{w-1})}+(\alpha_w-\alpha_{w-1})
\quad \text{without computing an integer $s_{w-1}$}.\qquad\qquad\\
\endalign$$ \ms

After the computation is done for {\rm Case(i) of Step w}, there are
two subcases for {\rm Case(i)}:\ms

$\underline{\text{\rm Subcase(1)}}$ \quad If $d_w=1$, then the
computation of $S_1=\text{\rm Join}(T_1,T_2\dots,T_w)$ is finished
by $T_w=\{e_{(\lambda_{w-1}+1)},e_{(\lambda_{w-1}+s_{w-1})}\}
=\{e_{(\lambda_{w-1})}+d_{w-1},e_{(\lambda_{w-1})}+(\alpha_w-\alpha_{w-1})\}$.
\ms

$\underline{\text{\rm Subcase(2)}}$ \quad If $d_w>1$, take the next
step, $\underline{\text{\rm Step (w+1) for the 1st half of the
$\gamma$-algorithm}}$. \ms

$\underline{\text{\rm Remark}}$ Note by {\rm(C.2.3)} that $d_w
=\gcd(d_{w-1},\alpha_w-\alpha_{w-1})$ with $d_{w-1}=n_wd_w$ and
$\alpha_w-\alpha_{w-1}=\widehat{\Delta}_wd_w$. We can compute a
unique integer $s_{w-1}$ such that
$(s_{w-1}-1)n_w<\widehat{\Delta}_w<s_{w-1}n_w$ because
$d_{w-1}<\alpha_w-\alpha_{w-1}$. \bs

\noindent$\underline{\text{\bf Case(ii) of of Step w for the 1st
half of the $\gamma$-algorithm}}$ Let
$d_{w-1}>\alpha_w-\alpha_{w-1}$.

Then, the algorithm for finding a subsequence
$T_w=\{e_{(\lambda_{w-1}+1)}, e_{(\lambda_{w-1}+s_{w-1})}\}$ of
$S_1$ can be uniquely represented as follows:
$$\align
(C.2.5) \qquad e_{(\lambda_{w-1}+1)}&=e_{(\lambda_{w-1})}+(\alpha_w-\alpha_{w-1}),\\
e_{(\lambda_{w-1}+s_{w-1})}&=s_{w-1}e_{(\lambda_{w-1})}+d_{w-1}
\quad \text{with an integer $s_{w-1}$ in $(C.2.5.1)$ }. \\
\text{Noting that} \quad {d_w}&=\gcd(d_{w-1},\alpha_w-\alpha_{w-1})
\quad \text{with $d_{w-1}=n_wd_w$ and
$\alpha_w-\alpha_{w-1}=\widehat{\Delta}_wd_w$,} \\
(C.2.5.1) \quad (s_{w-1}-1)&\widehat{\Delta}_w<n_w\le
s_{w-1}\widehat{\Delta}_w \quad \text{for some}\quad s_{w-1}\in \N
\quad  \text{because} \quad
1 \le\widehat{\Delta}_w<n_w. \\
\endalign$$

After the computation is done for {\rm Case(i) of Step w}, there are
two subcases for {\rm Case(ii)}:\ms

$\underline{\text{\rm Subcase(1)}}$ \quad If $d_w=1$, then the
computation of $S_1=\text{\rm Join}(T_1,T_2\dots,T_w)$ is finished
by $T_w=\{e_{(\lambda_{w-1}+1)}, e_{(\lambda_{w-1}+s_{w-1})}\}=
\{e_{(\lambda_{w-1})}+(\alpha_w-\alpha_{w-1}),s_{w-1}e_{(\lambda_{w-1})}+d_{w-1}\}$.
\ms

$\underline{\text{\rm Subcase(2)}}$ \quad If $d_w>1$, take the next
step, $\underline{\text{\rm Step (w+1) for the 1st half of the
$\gamma$-algorithm}}$.
\endproclaim \ms

\noindent{\bf Example C.2.1 for the 1st half of the
$\gamma$-algorithm in Theorem C.2:}

Let the Puiseux expansion for $C_3(t)$ be given by
$$
\text{$C_3(t):=$}  \left\{\eqalign{ y= &t^{45} \cr z=
&t^{60}+t^{65}+t^{71}. \cr } \right. \tag C.2.6
$$

Prove by [I] that the above $C_3(t)$ is the standard Puiseux
expansion, and by [II] that compute a finite sequence
$S_1=\{(g_3\circ\tau_{\lambda_3})_{\text{\rm singular part of the
divisor}}\}_{seq.}$ directly. Note that $S_1$ is a subsequence of
$S=\{(g_3\circ\tau_{\lambda_3})_{divisor}\}_{seq.}$ where $g_3$ is
the standard Puiseux polynomial $g_3\in \BC[y,z]$ of the recursive
$3$-type in {\rm Family(1)} such that \text{$V(g_3)\equiv C_3(t)$
\text{\rm (multi. seq.)}} \ms

First, we compute The Arithmetic Computation Formula by [I], and
next apply [I] to solve the problem by [II].

{\bf [I]} By The Arithmetic Computation Formula, we can compute a
finite sequence $\{(n_j,\widehat{\Delta}_j)\in N^2: j=1,2,3,\}$ of
pairs, each of which satisfies the following: Recall that
$d_j=\gcd(d_{j-1},\alpha_j-\alpha_{j-1})$ for $1\le j\le 4$ where
$d_{0}=n$ and $\alpha_{0}=0$.

{\rm(1)} Let $45=n=d_1n_1$ and $60=\alpha_1=d_1\widehat{\Delta}_1$
with $d_1=\gcd(n,\alpha_1)$. Then, $d_1=15$, $n_1=3$ and
$\beta_{1,1}=\widehat{\Delta}_1=4$. Also,
$e_{\lambda_{1}}=n_1\beta_{1,1}d_1=180$. \ms

 {\rm(2)} Let $15=d_1=d_2n_2$
and $5=\alpha_2-\alpha_1=d_2\widehat{\Delta}_2$ with
$d_2=\gcd(d_1,\alpha_2-\alpha_1)$. Then, $d_2=5$, $n_2=3$ and
$\widehat{\Delta}_2=1$. Also,
$e_{(\lambda_{2})}=e_{(\lambda_{1})}n_2+d_2n_2\text{$\widehat{\Delta}$}_2
=180\cdot{3}+{5}\cdot{3}\cdot{1}=555$. \ms

{\rm(3)} Let $5=d_2=d_3n_3$ and
$6=\alpha_3-\alpha_2=d_3\widehat{\Delta}_3$ with
$d_3=\gcd(d_2,\alpha_3-\alpha_2)$. Then, $d_3=1$, $n_3=5$ and
$\widehat{\Delta}_3=6$. Also,
$e_{(\lambda_{3})}=e_{(\lambda_{2})}n_3+d_3n_3\text{$\widehat{\Delta}$}_3
=555\cdot{5}+{1}\cdot{5}\cdot{6}=2805$. \ms

By (1), (2) and (3), it is clear that $n<\alpha_1<\alpha_2<\alpha_3$
and $n>d_1>d_2>d_3=1$, and so $C_3(t)$ is the standard Puiseux
expansion. \ms

Using an algorithm in Theorem C.2, to solve the problem is just to
find $S_1=\text{\rm Join}(T_1,T_2,T_3)$ in the sense of Definition
$A.3$ where $T_j=\{e_{\lambda_{j-1}+i}:i=1,s_{j-1}\}$ \quad
\text{for $1\le j\le 3$} and
$1<s_{j-1}\le{\lambda_j-\lambda_{j-1}}$. For notation, we write
$\lambda_{0}=0$. \ms

{\bf[II]} For a complete solution of (C.2.2), using {\bf the 1st
half of the $\gamma$-algorithm} then it suffices to compute
$S_1=\text{\rm Join}(T_1,T_2,T_3)$ easily where
$T_j=\{e_{\lambda_{j-1}+i}:i=1,s_{j-1}\}$ \text{for $1\le j\le 3$}
and $1<s_{j-1}\le{\lambda_j-\lambda_{j-1}}$ by the following three
steps: Note that $\lambda_{0}=0$ for notation.

$\underline{\text{\rm Step 1:}}$ It is clear that
$T_1=\{e_1,e_{s_0}\}=\{45,60\}$ where $e_1=n=45$ and
$e_{s_0}=\alpha_1=60$. \ms

$\underline{\text{\rm Step 2:}}$ Note that
$d_1=15>6=\alpha_2-\alpha_1$. So, by Case(ii) of Step 2,
$e_{(\lambda_{1}+1)}=e_{\lambda_{1}}+\alpha_2-\alpha_1=180+5=185$
and $e_{(\lambda_{1}+s_{1})}=s_1\cdot{180}+d_1=3\cdot{180}+15=555$
because $(s_1-1)\widehat{\Delta}_2<n_2\le s_1\widehat{\Delta}_2$
with $n_2=3$ and $\widehat{\Delta}_2=1$ implies that $s_1=3$. So,
$T_2=\{e_{(\lambda_{1}+1)},e_{(\lambda_{1}+s_{1})}\}=\{185,555\}$.\ms

$\underline{\text{\rm Step 3:}}$ Note that
$d_2=5<6=\alpha_3-\alpha_2$. So, by Case(i) of Step 3,
$e_{(\lambda_{2}+1)}=e_{\lambda_{2}}+d_2=555+5=560$ and
$e_{(\lambda_{2}+s_{2})}=e_{\lambda_{2}}+\alpha_3-\alpha_2=555+6=561$.
So,
$T_3=\{e_{(\lambda_{2}+1)},e_{(\lambda_{2}+s_{2})}\}=\{560,561\}$.\ms

Therefore, it can be found by the above three steps that
$S_1=\text{\rm Join}(T_1,T_2,T_3)$. \bs

{\bf {\S} C.3.  The 2nd half of the $\gamma$-algorithm(Theorem C.3)}

\proclaim{Theorem C.3(The 2nd half of the $\gamma$-algorithm: an
algorithm for finding a function from Family(2) onto Family(5))}

\noindent $\underline{\text{\bf Assumptions}}$ By using the same
properties and notations as in \text{\rm Definition A.1}, \text{\rm
Definition A.2} and \text{\rm Definition A.3}, let
$S_1=\{(g_r\circ\tau_{\xi})_{\text{\rm singular part of
divisor}}\}_{seq.}$ be a given finite sequence of
$(g_r\circ\tau_{\xi})_{\text{\rm singular part of divisor}}$, which
is a subsequence of
$S=\{(g_r\circ\tau_{\lambda_r})_{divisor}\}_{seq.}=\{e_i:1\le i\le
\lambda_r\}$.

Let $S_1=\text{\rm Join}(T_1,T_2\dots,T_r)$ where each $T_i$ is a
subsequence of two elements in $S_1$, satisfying the following
properties: Write $\lambda_0=1$, if necessary.
$$\align
\text{\rm(C.3.1)} \quad \quad   & T_1=\{e_i: i=1,s_0 \}
\quad \text{with $1<s_0\le{\lambda_1}$},  \qquad \qquad\qquad \qquad  \\
& T_j=\{e_{\lambda_{j-1}+i}:i=1,s_{j-1}\} \quad \text{for
$1\le j\le r$ and $2<s_{j-1}\le{\lambda_j-\lambda_{j-1}}$.} \qquad\qquad \\
\endalign$$

\noindent $\underline{\text{\bf Conclusions}}$ Given the sequence of
{\text{\rm singular part of divisor}} in $(C.3.1)$,  the problem is
how to compute {\bf a desired algorithm(The 2nd half of the
$\gamma$-algorithm for Theorem C.3)} for finding the standard
Puiseux expansion $C_r(t)$ in $(C.2.1)$ such that $C_r(t)$ and
$(g_r\circ\tau_{\lambda_r})_{\text{\rm singular part of divisor}}$
have the same multiplicity sequence.

The arithmetic computation formula for \{$d_1$, $n_1$,
$\beta_{1,1}$\} by \text{\rm (i)}, and \{$d_w$, $n_w$,
$\widehat{\Delta}_w$: $2\le w\le r$\} by \text{\rm (ii)}, and
\{$e_{\lambda_w}$: $1\le w\le r$\} by \text{\rm (iii)} in
{\rm(C.3.2)} will be most helpful for finding such an algorithm:
$$\align
\text{\rm(C.3.2)} \quad \text{\rm(i)} \quad
&\text{$d_1=\gcd(n,\alpha_1)$ with
$n=n_1d_1$ and $\alpha_1=\beta_{1,1}d_1$.} \\
\text{\rm{(ii)}} \quad
&\text{$d_w=\gcd(d_{w-1},\alpha_w-\alpha_{w-1})$ with
$d_{w-1}=n_wd_w$ and $\alpha_w-\alpha_{w-1}= \widehat{\Delta}_wd_w$}. \qquad\qquad\\
\text{\rm{\quad(iii)}} \quad
&\text{\text{$e_{\lambda_1}=d_1n_1\Delta_1(\beta_{1,1})$ and
$e_{\lambda_w}=e_{\lambda_{w-1}}n_w+d_wn_w\text{$\widehat{\Delta}$}_w$}
for $2\le w\le r$}.
\endalign$$

\noindent$\underline{\text{\bf(The 2nd half of the
$\gamma$-algorithm for Theorem C.3)}}$

\noindent$\underline{\text{\bf Step 1 for  The 2nd half of the
$\gamma$-algorithm}}$ Let $T_1=\{e_1,e_{s_0}\}$ by {\rm(C.3.1)}. It
is trivial that $n=e_1$ and $\alpha_1=e_{s_0}$.

Compute $d_1=\gcd(n,\alpha_1)$ with $n=n_1d_1$ and
$\alpha_1=\beta_{11}d_1$, and then $e_{\lambda_1}=n_1\beta_{11}d_1$.
\ms

Then, there are two subcases:

$\underline{\text{\rm Subcase(1)}}$ If $d_1=1$, then the standard
Puiseux expansion $C_1(t)$ can be computed by $y=t^n$ and
$z=t^{\alpha_1}$.

$\underline{\text{\rm Subcase(2)}}$ If $d_1>1$, take the next step,
$\underline{\text{\rm Step (2) for The 2nd half of the
$\gamma$-algorithm}}$. \ms

$\underline{\text{\rm Remark}}$ Note that there exists an integer
$s_0$ such that $(s_0-1)n<\alpha_1<{s_0}n$ or
$(s_0-1)n_1<\beta_{11}<s_0n_1$. \ms

Assuming that $d_1>1$, for $w=2,3,\dots,r$, an elementary
computational algorithm formula with $\alpha_w$ can be represented
as follows:

\noindent$\underline{\text{\bf Step w for The 2nd half of the
$\gamma$-algorithm}}$ \quad Let $d_{w-1}>1$ with $2\le w\le r$. By
$(C.3.1)$, consider
$$\align
\text{$T_w=\{e_{\lambda_{w-1}+1},e_{\lambda_{w-1}+s_{w-1}}\}$, the
$w$-th subsequence of $S_1$.} \tag C.3.3
\endalign$$

Then, there are two cases:

{\rm Case(i) for Step w} \quad
$e_{(\lambda_{w-1}+1)}-e_{(\lambda_{w-1})}=d_{w-1}$.

{\rm Case(ii) for Step w} \quad
$e_{(\lambda_{w-1}+1)}-e_{(\lambda_{w-1})}<d_{w-1}$. \ms

\noindent$\underline{\text{\bf Case(i) of Step w for The 2nd half of
the $\gamma$-algorithm}}$ Let
$e_{(\lambda_{w-1}+1)}-e_{(\lambda_{w-1})}=d_{w-1}$.

Then, the algorithm for finding the exponent $\alpha_w$ of the
standard Puiseux expansion and the coefficient $e_{\lambda_w}$ of
$E_{\lambda_w}$ can be represented as follows:
$$\align
(C.3.4) \quad  \quad d_{w-1}&=e_{(\lambda_{w-1}+1)}-e_{(\lambda_{w-1})},\\
\alpha_w-\alpha_{w-1}&=e_{(\lambda_{w-1}+s_{w-1})}-e_{(\lambda_{w-1})},\\
\text{noting that} \quad d_w&=\gcd(d_{w-1},\alpha_w-\alpha_{w-1})
\quad \text{with $d_{w-1}=n_wd_w$ and
$\alpha_w-\alpha_{w-1}=\widehat{\Delta}_wd_w$}, \qquad\\
e_{\lambda_w}&=n_we_{\lambda_{w-1}}+n_wd_w\widehat{\Delta}_w
\quad \text{\rm {by (C.3.2)}}. \\
\endalign$$

$\underline{\text{\rm Remark}}$ Note that there exists an integer
$s_{w-1}$ such that $(s_{w-1}-1)n_w<\widehat{\Delta}_w<s_{w-1}n_w$.
\ms

After the computation is done for {\rm Case(i) of Step w}, there are
two subcases for  {\rm Case(i)}:\ms

$\underline{\text{\rm Subcase(1) }}$ \quad If $d_w=1$, then the
standard Puiseux expansion can be defined by $y=t^n$ and
$z=t^{\alpha_1}+t^{\alpha_2}+\cdots +t^{\alpha_w}$ with
$\alpha_w-\alpha_{w-1}=e_{(\lambda_{w-1}+s_{w-1})}-e_{(\lambda_{w-1})}$.
\ms

$\underline{\text{\rm Subcase(2)}}$ \quad If $d_w>1$, take the next
step, $\underline{\text{\rm Step (w+1) for The 2nd half of the
$\gamma$-algorithm}}$. \bs

\noindent$\underline{\text{\bf Case(ii) of Step w for The 2nd half
of the $\gamma$-algorithm}}$ Let
$e_{(\lambda_{w-1}+1)}-e_{(\lambda_{w-1})}<d_{w-1}$.

Then, the algorithm for finding the exponent $\alpha_w$ of the
standard Puiseux expansion and the coefficient $e_{\lambda_w}$ of
$E_{\lambda_w}$ can be represented as follows:
$$\align
(C.3.5) \quad \alpha_w-\alpha_{w-1}&=e_{(\lambda_{w-1}+1)}-e_{(\lambda_{w-1})},\\
d_{w-1}&=e_{(\lambda_{w-1}+s_{w-1})}-s_{w-1}e_{(\lambda_{w-1})}
\quad \text{with an integer $s_{w-1}$ in $(C.3.5.1)$}.\\
\text{Note that} \quad d_w&=\gcd(d_{w-1},\alpha_w-\alpha_{w-1})
\quad \text{with $d_{w-1}=n_wd_w$ and
$\alpha_w-\alpha_{w-1}=\widehat{\Delta}_wd_w$,} \\
(C.3.5.1) \quad (s_{w-1}-1)&\widehat{\Delta}_w <n_w\le
s_{w-1}\widehat{\Delta}_w, \quad \text{and} \quad
e_{\lambda_w}=n_we_{\lambda_{w-1}}+n_wd_w\widehat{\Delta}_w
\quad \text{\rm {by (1.13.2)}}. \\
\endalign$$

$\underline{\text{\rm Remark}}$ Note that there exists an integer
$s$ such that $(s_{w-1}-1)\widehat{\Delta}_w<n_w\le
s_{w-1}\widehat{\Delta}_w$ because $1 \le\widehat{\Delta}_w<n_w$.
\ms

After the computation is done for {\rm Case(ii) of Step w}, there
are two subcases for {\rm Case(ii)}:\ms

$\underline{\text{\rm Subcase(1)}}$ \quad If $d_w=1$, then the
standard Puiseux expansion can be defined by $y=t^n$ and
$z=t^{\alpha_1}+t^{\alpha_2}+\cdots +t^{\alpha_w}$ with
$\alpha_w-\alpha_{w-1}=e_{(\lambda_{w-1}+1)}-e_{(\lambda_{w-1})}$.
\ms

$\underline{\text{\rm Subcase(2)}}$ \quad If $d_w>1$, take the next
step, $\underline{\text{\rm Step (w+1) for (The 2nd half of the
$\gamma$-algorithm)}}$.
\endproclaim \ms

{\bf Example C.3.1 for The 2nd half of the $\gamma$-algorithm in
Theorem C.3:}\ms

Using the same properties and notation as in (C.3.1), let
$S_1=\text{\rm Join}(T_1,T_2,T_3)$

where $T_1=\{e_{1},e_{s_0}\}=\{45,75\}$,
$T_2=\{e_{\lambda_1+1},e_{\lambda_1+s_1}\}=\{231,690\}$,
$T_3=\{e_{\lambda_2+1},e_{\lambda_2+s_2}\}=\{1158,1159\}$(with
$\lambda_i<\lambda_i+s_i\le \lambda_{i+1}$ for $i=0,1,2$). \ms

(1) First, Prove that there is $g_3\in \text{\rm Family(1)}$ such
that $S_1$ is given by
$$\align
(C.3.6) \qquad \qquad S_1=\{(g_r\circ\tau_{\lambda_r})_{\text{\rm
singular part of the divisor}}\}_{seq.}= \text{\rm
Join}(T_1,T_2,T_3). \qquad \qquad \qquad \qquad
\endalign$$

(2) Next, using an algorithm in Theorem C.3, the problem is to find
the standard Puiseux expansion $C_3(t)$ of the $3$-th type for the
above irreducible plane curve $C$ such that \text{$C_3(t)\equiv
V(g_3)$} \text{\rm (multi. seq.)} where $C_3(t)$ is given by $y=
t^{n}$ and $z=t^{\alpha_1}+t^{\alpha_2}+t^{\alpha_3}$. \bs

For a complete solution for (C.3.1), using {\bf The 2nd half of the
$\gamma$-algorithm for Theorem C.3} then it suffices to compute
$S_1=\text{\rm Join}(T_1,T_2,T_3)$ easily where
$T_j=\{e_{\lambda_{j-1}+i}:i=1,s_{j-1}\}$ \text{for $1\le j\le 3$}
and $1<s_{j-1}\le{\lambda_j-\lambda_{j-1}}$ by the following three
steps: Note that $\lambda_{0}=0$ for notation.

$\underline{\text{\rm Step 1:}}$ It is clear that $n=e_{1}=45$, and
$\alpha_1=e_{s_{0}}=75$.

In preparation for the computation in Step 2, $e_{\lambda_1}=225$
because $d_1=\gcd(n,\alpha_1)=15>1$ with $n=n_1d_1$ and
$\alpha_1=\beta_{1,1}d_1$ implies that $n_1=3$, $\beta_{1,1}=5$ and
$e_{\lambda_1}=n_1\beta_{1,1}d_1=225$ by (C.3.2). \ms

$\underline{\text{\rm Step 2:}}$ Since
$e_{\lambda_1+1}-e_{\lambda_1}=231-225=6<15=d_1$, by Case(ii) of
Step 2, ${\alpha_2}-{\alpha_1}=e_{{\lambda_1}+1}-e_{\lambda_1}=6$.
So, it is clear by Step 1 that ${\alpha_2}={\alpha_1}+6=75+6=81$.

In preparation for the computation in Step 3,
$e_{\lambda_2}=n_2\cdot{e_{\lambda_{1}}}+d_2n_2\text{$\widehat{\Delta}$}_2
=5\cdot{225}+{3}\cdot{5}\cdot{2}=1125+30$ by (C.3.2) because of the
following computations:

$d_2=\gcd(d_1,\alpha_2-\alpha_1)=\gcd(15,6)=3>1$ with $d_1=n_2d_2$
and $\alpha_2-\alpha_1=\text{$\widehat{\Delta}$}_2d_2$ implies that
$n_2=5$ and $\text{$\widehat{\Delta}$}_2=2$, and so
$e_{\lambda_2}=1155$. \ms

$\underline{\text{\rm Step 3:}}$ Since
$e_{\lambda_2+1}-e_{\lambda_2}=1158-1155=3=d_2$, by Case(i) of Step
3,
$\alpha_3-\alpha_2=e_{(\lambda_{2}+s_2)}-e_{\lambda_{2}}=1159-1155=4$.
Then it is clear by Step 2 that ${\alpha_3}={\alpha_2}+4=81+4=85$.
\ms

So, the standard Puiseux expansion $C_3(t)$ such that
\text{$C_3(t)\equiv V(f)$} \text{\rm (multi. seq.)} can be given by
$y=t^{45}$ and $z=t^{75}+t^{81}+t^{85}$, because $d_2=1$ implies
that $d_3=\gcd(d_2,\alpha_3-\alpha_2)=1$. \bs

{\bf Acknowledgement.} To finish my paper, I must be very
appreciated for a good help of my three brothers, Chunghoo
Kang(Medical doctor of Physician), Chungin Kang(Professor of
Politics) and Chungmin Kang(Doctor of Philosopher) because of the
following. Whenever they had a good health then they had given a
better service to our parents for a long time than I had given when
our parents with no good health were alive, because my brothers have
already known that I have had no good health up to now. This has been very
helpful to study my paper by myself completely as far as the time is
concerned, because I really need a lot of time to finish writing a
long paper. Moreover, I am very appreciated for Dr. Young Rock Kim's
PDF Tex typewriting help, sometimes.

\newpage

\vfill \pagebreak

\ms\ms\bs

\newpage
  
Index 
\bs\bs

Puiseux expansion, 3, 23, 28
    
standard Puisuex expansions, 3,23, 28, 128 

the divisors under the standard
resolutions, 

multiplicity sequences, 23, 28

standard Puisuex
polynomial of two complex variables of the recursive type,  25,

standard Puiseux(irreducible)
W-poly of two complex variables of the recursive type,  3, 13, 23,24,25 

Euclidean multiplicity sequence for two positive integers,

quasisingularity,

reducible in the ring of convergent
power series, 24,27 

irreducible in the ring of convergent
power series, 24,27

Weierstrass polynomial(W-poly), 23,24

Weierstrass Preparation Theorem(WPT), 204 

Weierstrass Division Theorem(WDT), 204

Weierstrass division theorem for the W-polys, 13,34

Weierstrass division algorithm for the W-polys, 13, 35    
      
exceptional curve of the first kind, 84

Family(1), 27

Family(2), 28

Family(3), 28

 \bs\bs
 
\newpage                   
     
Chunghyuk Kang

\ms\ms\bs

\newpage                                                                                          
                                                                                                                                                                                    

Department of Mathematics, Seoul National
University, Seoul 151--742, KOREA

E-mail address: chkang\@snu.ac.kr

\end
\enddocument